\documentclass[11pt,oneside,english,refpage,intoc,openany]{book}
\usepackage{lmodern}

\usepackage[T1]{fontenc}
\usepackage[latin9]{inputenc}
\usepackage{geometry}
\usepackage{color}
\usepackage{babel}
\usepackage{varioref}
\usepackage{textcomp}
\usepackage{mathrsfs}
\usepackage{amsmath}
\usepackage{amsthm}
\usepackage{amssymb}
\usepackage{cancel}
\usepackage{makeidx}
\makeindex
\usepackage{graphicx}
\PassOptionsToPackage{normalem}{ulem}
\usepackage{ulem}
\usepackage{nomencl}
\providecommand{\printnomenclature}{\printglossary}
\providecommand{\makenomenclature}{\makeglossary}
\makenomenclature
\usepackage[unicode=true,
 bookmarks=true,bookmarksnumbered=true,bookmarksopen=false,
 breaklinks=false,pdfborder={0 0 1},backref=false,colorlinks=true]
 {hyperref}
\hypersetup{pdftitle={$\left(p,q\right)$-Adic Analysis and the Collatz Conjecture},
 pdfauthor={LyX Team},
 pdfsubject={LyX},
 pdfkeywords={LyX},
 linkcolor=black,citecolor=black,urlcolor=blue,filecolor=blue,pdfpagelayout=OneColumn,pdfnewwindow=true,pdfstartview=XYZ,plainpages=false}

\makeatletter

\providecommand{\tabularnewline}{\\}

\theoremstyle{definition}
    \ifx\thechapter\undefined
      \newtheorem{defn}{\protect\definitionname}
    \else
      \newtheorem{defn}{\protect\definitionname}[chapter]
    \fi
\theoremstyle{plain}
    \ifx\thechapter\undefined
      \newtheorem{prop}{\protect\propositionname}
    \else
      \newtheorem{prop}{\protect\propositionname}[chapter]
    \fi
\theoremstyle{plain}
    \ifx\thechapter\undefined
      \newtheorem{lem}{\protect\lemmaname}
    \else
      \newtheorem{lem}{\protect\lemmaname}[chapter]
    \fi
\theoremstyle{plain}
    \ifx\thechapter\undefined
	    \newtheorem{thm}{\protect\theoremname}
	  \else
      \newtheorem{thm}{\protect\theoremname}[chapter]
    \fi
\theoremstyle{definition}
    \ifx\thechapter\undefined
      \newtheorem{example}{\protect\examplename}
    \else
      \newtheorem{example}{\protect\examplename}[chapter]
    \fi
\theoremstyle{remark}
    \ifx\thechapter\undefined
      \newtheorem{rem}{\protect\remarkname}
    \else
      \newtheorem{rem}{\protect\remarkname}[chapter]
    \fi
\theoremstyle{plain}
    \ifx\thechapter\undefined
      \newtheorem{fact}{\protect\factname}
    \else
      \newtheorem{fact}{\protect\factname}[chapter]
    \fi
\theoremstyle{plain}
    \ifx\thechapter\undefined
  \newtheorem{cor}{\protect\corollaryname}
\else
      \newtheorem{cor}{\protect\corollaryname}[chapter]
    \fi
\theoremstyle{plain}
\newtheorem*{assumption*}{\protect\assumptionname}
\theoremstyle{plain}
    \ifx\thechapter\undefined
      \newtheorem{assumption}{\protect\assumptionname}
    \else
      \newtheorem{assumption}{\protect\assumptionname}[chapter]
    \fi
\theoremstyle{plain}
    \ifx\thechapter\undefined
      \newtheorem{conjecture}{\protect\conjecturename}
    \else
      \newtheorem{conjecture}{\protect\conjecturename}[chapter]
    \fi
\theoremstyle{remark}
    \ifx\thechapter\undefined
      \newtheorem{claim}{\protect\claimname}
    \else
      \newtheorem{claim}{\protect\claimname}[chapter]
    \fi
\theoremstyle{remark}
    \ifx\thechapter\undefined
      \newtheorem{notation}{\protect\notationname}
    \else
      \newtheorem{notation}{\protect\notationname}[chapter]
    \fi
\theoremstyle{definition}
    \ifx\thechapter\undefined
      \newtheorem{problem}{\protect\problemname}
    \else
      \newtheorem{problem}{\protect\problemname}[chapter]
    \fi
\theoremstyle{plain}
    \ifx\thechapter\undefined
      \newtheorem{question}{\protect\questionname}
    \else
      \newtheorem{question}{\protect\questionname}[chapter]
    \fi

\@ifundefined{date}{}{\date{}}
\usepackage{babel}
\addto\captionsenglish{}
\usepackage[figure]{hypcap}

\usepackage{fancyhdr}
\renewcommand{\headrulewidth}{2pt}

\@ifundefined{extrarowheight}{\usepackage{array}}{}
\providecommand{\assumptionname}{Assumption}
\providecommand{\claimname}{Claim}
\providecommand{\conjecturename}{Conjecture}
\providecommand{\corollaryname}{Corollary}
\providecommand{\definitionname}{Definition}
\providecommand{\examplename}{Example}
\providecommand{\factname}{Fact}
\providecommand{\lemmaname}{Lemma}
\providecommand{\notationname}{Notation}
\providecommand{\problemname}{Problem}
\providecommand{\propositionname}{Proposition}
\providecommand{\questionname}{Question}
\providecommand{\remarkname}{Remark}
\providecommand{\theoremname}{Theorem}

\makeatother

\providecommand{\assumptionname}{Assumption}
\providecommand{\claimname}{Claim}
\providecommand{\conjecturename}{Conjecture}
\providecommand{\corollaryname}{Corollary}
\providecommand{\definitionname}{Definition}
\providecommand{\examplename}{Example}
\providecommand{\factname}{Fact}
\providecommand{\lemmaname}{Lemma}
\providecommand{\notationname}{Notation}
\providecommand{\problemname}{Problem}
\providecommand{\propositionname}{Proposition}
\providecommand{\questionname}{Question}
\providecommand{\remarkname}{Remark}
\providecommand{\theoremname}{Theorem}

\begin{document}
\fancypagestyle{plain}{
\fancyhf{}
\fancyhead[R]{\thepage}
\renewcommand{\headrulewidth}{0pt}
\renewcommand{\footrulewidth}{0pt}
}\frontmatter

\global\long\def\headrulewidth{0pt}%

\thispagestyle{fancy}\fancyfoot[L]{May 2022}\fancyfoot[R]{Maxwell Charles Siegel}\pagenumbering{gobble} 
\begin{center}
$\left(p,q\right)$-ADIC ANALYSIS AND THE COLLATZ CONJECTURE
\par\end{center}

\vphantom{}
\begin{center}
by 
\par\end{center}

\vphantom{} 
\begin{center}
Maxwell Charles Siegel 
\par\end{center}

\vphantom{}

\vphantom{}

\vphantom{}
\begin{center}
A Dissertation Presented to the 
\par\end{center}

\begin{center}
FACULTY OF THE USC DORNSIFE COLLEGE OF LETTERS AND SCIENCES
\par\end{center}

\begin{center}
UNIVERSITY OF SOUTHERN CALIFORNIA 
\par\end{center}

\begin{center}
In Partial Fulfillment of the
\par\end{center}

\begin{center}
Requirements for the Degree 
\par\end{center}

\begin{center}
DOCTOR OF PHILOSOPHY 
\par\end{center}

\begin{center}
(MATHEMATICS) 
\par\end{center}

\vphantom{}

\vphantom{}

\vphantom{} 
\begin{center}
May 2022 
\par\end{center}

\newpage\pagenumbering{roman}\setcounter{page}{2}

\addcontentsline{toc}{chapter}{Dedication}

\ 

\ 

\ 

\ 

\ 

\ 

\ 

\ 

\ 

\ 

\ 

\ 

\ 

\ 
\begin{center}
{\large{}{}I dedicate this dissertation to those who listened, and
to those who cared. Also, to my characters; I'm sorry it took so long.}{\large\par}
\par\end{center}

\begin{center}
{\large{}{}I'll try to do better next time.}{\large\par}
\par\end{center}

\tableofcontents\pagebreak 

\listoftables

\newpage{}

\section*{Abstract}

\pagestyle{fancy}\fancyfoot{}\fancyhead[L]{\sl ABSTRACT}\fancyhead[R]{\thepage}\addcontentsline{toc}{chapter}{Abstract}

What use can there be for a function from the $p$-adic numbers to
the $q$-adic numbers, where $p$ and $q$ are distinct primes? The
traditional answer\textemdash courtesy of the half-century old theory
of non-archimedean functional analysis: \emph{not much}. It turns
out this judgment was premature. ``$\left(p,q\right)$-adic analysis''
of this sort appears to be naturally suited for studying the infamous
Collatz map and similar arithmetical dynamical systems. Given such
a map $H:\mathbb{Z}\rightarrow\mathbb{Z}$, one can construct a function
$\chi_{H}:\mathbb{Z}_{p}\rightarrow\mathbb{Z}_{q}$ for an appropriate
choice of distinct primes $p,q$ with the property that $x\in\mathbb{Z}\backslash\left\{ 0\right\} $
is a periodic point of $H$ if and only if there is a $p$-adic integer
$\mathfrak{z}\in\left(\mathbb{Q}\cap\mathbb{Z}_{p}\right)\backslash\left\{ 0,1,2,\ldots\right\} $
so that $\chi_{H}\left(\mathfrak{z}\right)=x$. By generalizing Monna-Springer
integration theory and establishing a $\left(p,q\right)$-adic analogue
of the Wiener Tauberian Theorem, one can show that the question ``is
$x\in\mathbb{Z}\backslash\left\{ 0\right\} $ a periodic point of
$H$'' is essentially equivalent to ``is the span of the translates
of the Fourier transform of $\chi_{H}\left(\mathfrak{z}\right)-x$
dense in an appropriate non-archimedean function space?'' This presents
an exciting new frontier in Collatz research, and these methods can
be used to study Collatz-type dynamical systems on the lattice $\mathbb{Z}^{d}$
for any $d\geq1$.

\newpage{}

\section*{Preface}

\fancyhead[L]{\sl PREFACE}\addcontentsline{toc}{chapter}{Preface}

Like with so many other of the world's wonders, my first encounter
with the Collatz Conjecture was on Wikipedia. I might have previously
stumbled across it as an undergraduate, but it wasn't until around
the time I started graduate school (Autumn 2015) that I gave it any
serious consideration. I remember trying to understand what $p$-adic
numbers were solely through what Wikipedia said they were, hoping
to apply it to the Conjecture, but then gave up after a day or two.
I spent a while fooling around with trying to prove the holomorphy
of tetration, and then fiddled with fractional differentiation and
the iterated Laplace transform, mostly mindlessly. But then, in March
of 2017, the Collatz ``bug'' bit\textemdash and bit \emph{hard}.
I realized I could encode the action of the Collatz map as a transformation
of holomorphic functions on the disk, and\textemdash with complex
analysis being near and dear to my heart\textemdash I was hooked.
Obsession was \emph{instantaneous}.

As of the writing of this preface some five years later, I suppose
the ardor of my obsession has dimmed, somewhat, primarily under the
weight of external pressures\textemdash angst, doubt, apprehension
and a healthy dose of circumspection, among them. Although all PhD
dissertations are uphill battles\textemdash frequently Sisyphean\textemdash the
lack of a real guiding figure for my studies made my scholarly journey
particularly grueling. The mathematics department of my graduate school\textemdash the
University of Southern California\textemdash has many strengths. Unfortunately,
these did not include anything in or adjacent to the wide variety
of sub-disciplines I drew from in my exploration of Collatz: harmonic
analysis, boundary behavior of power series, analytic number theory,
$p$-adic analysis, and\textemdash most recently\textemdash non-archimedean
functional analysis. I was up the proverbial creek without a paddle:
I did not get to benefit from the guidance of a faculty advisor who
could give substantive feedback on my mathematical work. The resulting
intellectual isolation was vaster than anything I'd previously experienced\textemdash and
that was \emph{before }the COVID-19 pandemic hit.

Of course, I alone am to blame for this. I \emph{chose }to submit
to my stubborn obsession, rather than fight it. The consequences of
this decision will most likely haunt me for many years to come. as
much trouble as my doltish tenacity has caused me, I still have to
thank it. I wouldn't have managed to find a light at the end of the
tunnel without it. As I write these words, I remain apprehensive\textemdash even
fearful\textemdash of what the future holds for me, mathematical or
otherwise. No course of action comes with a guarantee of success,
not even when you follow your heart. But, come what may, I can say
with confidence that I am proud of what I have accomplished in my
dissertation, and prouder to still be able to share it with others.
Connectedness is success all its own.

But enough about \emph{me}.

In the past five years, my research into the Collatz Conjecture can
be grouped into two-and-a-half different categories. The first of
these was the spark that started it all: my independent discovery
that the Collatz Conjecture could be reformulated in terms of the
solutions of a \emph{functional equation} defined over the open unit
disk $\mathbb{D}$ in $\mathbb{C}$. Specifically, the equation: 
\begin{equation}
f\left(z\right)=f\left(z^{2}\right)+\frac{z^{-1/3}}{3}\sum_{k=0}^{2}e^{-\frac{4k\pi i}{3}}f\left(e^{\frac{2k\pi i}{3}}z^{2/3}\right)\label{eq:Collatz functional equation}
\end{equation}
where the unknown $f\left(z\right)$ is represented by a power series
of the form: 
\begin{equation}
f\left(z\right)=\sum_{n=0}^{\infty}c_{n}z^{n}\label{eq:Power series}
\end{equation}
In this context, the Collatz Conjecture becomes equivalent to the
assertion that the set of holomorphic functions $\mathbb{D}\rightarrow\mathbb{C}$
solving equation (\ref{eq:Collatz functional equation}) is precisely:
\begin{equation}
\left\{ \alpha+\frac{\beta z}{1-z}:\alpha,\beta\in\mathbb{C}\right\} 
\end{equation}
This discovery was not new. The observation was first made by Meinardus
and Berg \cite{Berg =00003D000026 Meinardus,Meinardus and Berg} (see
also \cite{Meinardus,Opfer}). What made my take different was its
thrust and the broadness of its scope. I arrived at equation (\ref{eq:Collatz functional equation})
in the more general context of considering a linear operator (what
I call a \textbf{permutation operator}) induced by a map $H:\mathbb{N}_{0}\rightarrow\mathbb{N}_{0}$;
here $\mathbb{N}_{0}$ is the set of integers $\geq0$). The permutation
operator $\mathcal{Q}_{H}$ induced by $H$ acts upon the space of
holomorphic functions $\mathbb{D}\rightarrow\mathbb{C}$ by way of
the formula: 
\begin{equation}
\mathcal{Q}_{H}\left\{ \sum_{n=0}^{\infty}c_{n}z^{n}\right\} \left(z\right)\overset{\textrm{def}}{=}\sum_{n=0}^{\infty}c_{H\left(n\right)}z^{n}\label{eq:Definition of a permutation operator}
\end{equation}
It is not difficult to see that the space of functions fixed by $\mathcal{Q}_{H}$
is precisely: 
\begin{equation}
\left\{ \sum_{v\in V}z^{v}:V\textrm{ is an orbit class of }H\textrm{ in }\mathbb{N}_{0}\right\} \label{eq:orbit class set-series}
\end{equation}
A simple computation shows that the functional equation (\ref{eq:Collatz functional equation})
is precisely the equation $\mathcal{Q}_{T_{3}}\left\{ f\right\} \left(z\right)=f\left(z\right)$,
where $T_{3}$ is the \textbf{Shortened Collatz map}, with the general
\textbf{Shortened}\index{$qx+1$ map} \textbf{$qx+1$ map} being defined
by \index{$3x+1$ map}
\begin{equation}
T_{q}\left(n\right)\overset{\textrm{def}}{=}\begin{cases}
\frac{n}{2} & \textrm{if }n=0\mod2\\
\frac{qn+1}{2} & \textrm{if }n=1\mod2
\end{cases}\label{eq:qx plus 1 shortened}
\end{equation}
where $q$ is any odd integer $\geq3$. Maps like these\textemdash I
call them \textbf{Hydra maps}\textemdash have been at the heart of
my research in the past five years. Indeed, they are the central topic
of this dissertation.

My investigations of functional equations of the form $\mathcal{Q}_{H}\left\{ f\right\} \left(z\right)=f\left(z\right)$
eventually led me to a kind of partial version of a famous lemma\footnote{\textbf{Theorem 3.1 }in \cite{On Wiener's Lemma}.}
due to Norbert Wiener. I used this to devise a novel kind of Tauberian
theory, where the boundary behavior of a holomorphic function on $\mathbb{D}$
was represented by a singular measure on the unit circle\footnote{This viewpoint is intimately connected with the larger topic of representing
holomorphic functions on $\mathbb{D}$ in terms of their boundary
values, such as by the Poisson Integral Formula \cite{Bounded analytic functions,function classes on the unit disc}
or the Cauchy transform and its fractional analogues \cite{Ross et al,Fractional Cauchy transforms}.}. Of particular import was that these measures could also be interpreted
as \emph{functions} in the space $L^{2}\left(\mathbb{Q}/\mathbb{Z},\mathbb{C}\right)$\textemdash all
complex-valued functions on the additive group $\mathbb{Q}/\mathbb{Z}=\left[0,1\right)\cap\mathbb{Q}$
with finite $L^{2}$-norm with respect to that group's Haar measure
(the counting measure):

\begin{equation}
L^{2}\left(\mathbb{Q}/\mathbb{Z},\mathbb{C}\right)=\left\{ f:\mathbb{Q}/\mathbb{Z}\rightarrow\mathbb{C}\mid\sum_{t\in\mathbb{Q}/\mathbb{Z}}\left|f\left(t\right)\right|^{2}<\infty\right\} \label{eq:L2 of Q/Z}
\end{equation}
Because the Pontryagin dual of $\mathbb{Q}/\mathbb{Z}$ is the additive
group of profinite integers\footnote{The direct product $\prod_{p\in\mathbb{P}}\mathbb{Z}_{p}$ of the
rings of $p$-adic integers, over all primes $p$.}, I decided to interpret these measures/functions as the Fourier transforms
of complex-valued functions on the profinite integers. Detailed explanations
of these can be found in my paper \cite{Dreancatchers for Hydra Maps}
on arXiv. The notation in that paper is a mess, and I will likely
have to re-write it at some point in the future.

During these investigations, my \emph{idée fixe }was to use Fourier
analysis and exploit the Hilbert space structure of $L^{2}\left(\mathbb{Q}/\mathbb{Z},\mathbb{C}\right)$
to prove what I called \emph{rationality theorems}. These were generalizations
of the notion of a \emph{sufficient set} \index{sufficient set}(see
\cite{Andaloro - first paper on sufficient sets,Monks' sufficiency paper for 3x+1});
a set $S\subseteq\mathbb{N}_{0}$ is said to be ``sufficient'' for
the Collatz Conjecture whenever the statement ``the Collatz map iterates
every element of $S$ to $1$'' is sufficient to prove the Collatz
Conjecture in full. \cite{Remmert} showed that, for example, the
set $\left\{ 16n+1:n\geq0\right\} $ was sufficient for the Collatz
map; Monks et. al. \cite{The many Monks paper} extended this, showing
that any infinite arithmetic\index{arithmetic progression} progression
(``IAP'') (a set of the form $\left\{ an+b:n\geq0\right\} $ for
integers $a,b$, with $a\neq0$) was sufficient for the Collatz map.
Unfortunately, all I could manage to prove was a weaker form of this
conclusion: namely, that for any Hydra map satisfying certain simple
conditions, if an orbit class $V\subseteq\mathbb{N}_{0}$ of the Hydra
map could be written as the union of a finite set and finitely many
IAPs, then $V$ contained all but at most finitely many elements of
$\mathbb{N}_{0}$. The one advantage of my approach was that it \emph{appeared}
to generalize to ``multi-dimensional'' Hydra maps\textemdash a Collatz-type
map defined on the ring of integers a finite-degree field extension
of $\mathbb{Q}$, or, equivalently, on the lattice $\mathbb{Z}^{d}$.
I say ``appeared'' because, while the profinite integer Hilbert
space part of the generalization worked out without any difficulty,
it is not yet clear if (or how) the crux of the original argument
could be extended to power series of several complex variables. For
the curious, this crux was a beautiful theorem\footnote{See Section 4 in Chapter 11 of \cite{Remmert}.}
due to Gábor Szeg\H{o}\index{SzegH{o}, Gábor@Szeg\H{o}, Gábor} regarding
the analytic continuability of (\ref{eq:Power series}) in the case
where the set: 
\begin{equation}
\left\{ c_{n}:n\in\mathbb{N}_{0}\right\} 
\end{equation}
is finite. Taking the time to investigate this issue would allow for
the results of \cite{Dreancatchers for Hydra Maps} to be extended
to the multi-dimensional case, the details of which\textemdash modulo
the Szeg\H{o} gap\textemdash I have written up but not yet published
or uploaded in any form. I will most likely re-write \cite{Dreancatchers for Hydra Maps}
at a future date\textemdash assuming I have not already done so by
the time you are reading this.

The \emph{first-and-a-halfth} category of my work branched off from
my quest for rationality theorems. I turned my focus to using $L^{2}\left(\mathbb{Q}/\mathbb{Z},\mathbb{C}\right)$
to establish the possible value(s) of the growth constant $r>0$ \index{growth exponent}
for orbit classes of Hydra maps. The essentials are as follows. Let
$H$ be a Hydra map. Then, letting $\omega$ be a positive real number
and letting $f\left(z\right)=\left(1-z\right)^{-\omega}$, an elementary
computation with limits establishes that the limit
\begin{equation}
\lim_{x\uparrow1}\left(1-x\right)^{\omega}\mathcal{Q}_{H}\left\{ f\right\} \left(x\right)=1\label{eq:Eigenrate Equation}
\end{equation}
occurs if and only if $\omega$ satisfies a certain transcendental
equation. For the case of the Shortened $qx+1$ map, the values of
$\omega$ satisfying (\ref{eq:Eigenrate Equation}) are precisely
the solutions of the equation: 
\begin{equation}
2^{\omega}-q^{\omega-1}=1
\end{equation}
My hope was to prove that for: 
\begin{equation}
\omega_{0}\overset{\textrm{def}}{=}\min\left\{ \omega\in\left(0,1\right]:2^{\omega}-q^{\omega-1}=1\right\} 
\end{equation}
given any orbit class $V\subseteq\mathbb{N}_{0}$ of $T_{q}$ such
that the iterates $T_{q}\left(v\right),T_{q}\left(T_{q}\left(v\right)\right),\ldots$
were bounded for each $v\in V$, we would have: 
\begin{equation}
\left|\left\{ v\in V:v\leq N\right\} \right|=o\left(N^{\omega_{0}+\epsilon}\right)\textrm{ as }N\rightarrow\infty
\end{equation}
for any $\epsilon>0$. While initially quite promising, obtaining
this asymptotic result hinged on justifying an interchange of limits,
a justification which appears to be at least as difficult as establishing
this asymptotic directly, \emph{without appealing} to $\mathcal{Q}_{H}$
and complex analysis as a middleman. Although still conjecture as
of the writing of this dissertation, the $2^{\omega}-q^{\omega-1}=1$
formula would provide an elegant resolution to a conjecture of Kontorovich
and Lagarias regarding the growth exponent for the orbit class of
all integers that $T_{5}$ iterates to $1$ \cite{Lagarias-Kontorovich Paper}.
I have also studied this issue from profinite integer. This turns
out to be equivalent to studying the iterations of $H$ on the profinite
integer and subgroups thereof, such as $\mathbb{Z}_{2}\times\mathbb{Z}_{3}$,
the ring of ``$\left(2,3\right)$-adic integers''. I have obtained
minor results in that regard by exploiting the known ergodicity of
the $qx+1$ maps on $\mathbb{Z}_{2}$ (see \textbf{Theorem \ref{thm:shift map}}\vpageref{thm:shift map}).
I will likely write this work up for publication at a later date.

Whereas the above approach was grounded on the unit disk and harmonic
analysis on subgroups of the circle, the second-and-a-halfth branch
of my research\textemdash the one considered in this dissertation\textemdash concerns
a completely different line of reasoning, from a completely different
subspecialty of analysis, no less: $p$-adic and non-archimedean analysis.
By considering the different branches of a given, suitably well-behaved
Hydra map, I discovered a function $\chi_{H}:\mathbb{Z}_{p}\rightarrow\mathbb{Z}_{q}$
(where $p$ and $q$ are distinct primes which depend on $H$) with
the remarkable property\footnote{I call this the \textbf{Correspondence Principle}; see \textbf{Theorem
\ref{thm:CP v1} }and its other three variants, \textbf{Corollaries
\ref{cor:CP v2}}, \textbf{\ref{cor:CP v3}}, and \textbf{\ref{cor:CP v4}},
starting from page \pageref{thm:CP v1}.} that the set of periodic points of $H$ in $\mathbb{Z}$ was \emph{completely
determined} by the elements of the image of $\chi_{H}$ which lay
in $\mathbb{Z}_{q}\cap\mathbb{Z}$. While I developed this approach
in late 2019 and early 2020, and it was only in mid-2020 that I realized
that the $\chi_{H}$ corresponding to the Shortened Collatz map played
a pivotal role in a paper on the Collatz Conjecture published by Terence
Tao at the end of 2019 \cite{Tao Probability paper}. However, due
to Tao's \emph{deliberate choice} not\emph{ }to use a fully $p$-adic
formalism\textemdash one which would have invoked $\chi_{H}$ in full\textemdash it
appears he failed to notice this remarkable property. The relation
between Tao's work an my own is discussed in detail starting on page
\pageref{subsec:2.2.4 Other-Avenues}. It would be interesting to
see my approach investigated from Tao's perspective. As of the writing
of this dissertation, I have yet to grok Tao's probabilistic arguments,
and have therefore refrained from attempting to employ or emulate
them here.

To the extent there even is such a thing as ``Collatz studies'',
the subject is bazaar filled with odd and ends\textemdash \emph{ad
hoc}s\emph{ }of the moment formed to deal with the latest line of
thought about the Collatz map \emph{itself}. On the other hand, there
have been few, if any, noteworthy attempts to center the study of
Collatz-type arithmetical dynamical systems on a sound \emph{theoretical}
underpinning\textemdash and by ``sound \emph{theoretical }underpinning'',
I mean ``non-probabilistic investigations of multiple Collatz-type
maps which \emph{do not }''. While there \emph{are }some works that
address larger families of Collatz-type maps (\cite{dgh paper,Matthews' slides,Matthews and Watts,Moller's paper (german)},
for instance), they tend to be probabilistic, and are usually replete
with conjectures, heuristic results, or specific findings that do
not inspire much confidence in the possibility of a broader underlying
non-probabilistic theory. My greatest aspiration for this dissertation\textemdash my
labor of love, madness, obsession, and wonder\textemdash is that it
will help pull Collatz studies out of the shadows to to emerge as
a commendable, \emph{recommendable }field of study, complete with
a theoretical foundation worth admiring. If I can do even half of
that, I will consider my efforts successful, and all the angst I had
to wade through along the way will, at last, be vindicated.

\vphantom{}Los Angeles, California, USA

March 14, 2022

\vphantom{}

\vphantom{}

\vphantom{}

P.S. If you have questions, or have found typos or grammatical errors,
please e-mail me at siegelmaxwellc@ucla.edu (yes, I went to UCLA as
an undergraduate).

\newpage{}

\pagebreak\pagestyle{headings}

\mainmatter

\chapter{Introduction \label{chap:Introduction}}

\thispagestyle{headings}

\includegraphics[scale=0.45]{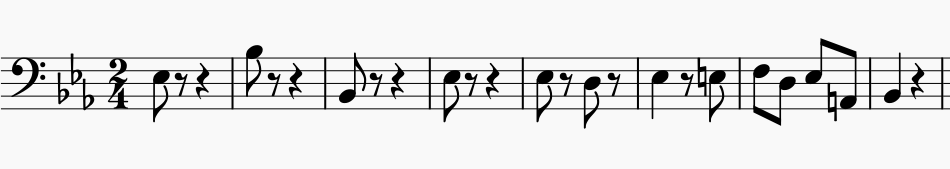}

\vphantom{}

Just like the title says, this is the introductory chapter of my PhD
dissertation. Section \ref{sec:1.1 A-Guide-for-the-Perplexed} begins
with an introduction the general notational conventions and idiosyncrasies
used throughout this tome (Subsection \ref{subsec:1.1.1 Notational-Idiosyncrasies}).
Subsection \ref{subsec:1.1.2 Some-Much-Needed-Explanations} provides
the eponymous much needed explanations regarding my main results,
and the agenda informing my (and our) pursuit of them\textemdash not
to mention a \emph{mea culpa }regarding my dissertation's prodigious
length. Speaking of length, Subsection \ref{subsec:1.1.3 An-Outline-of}
contains a chapter-by-chapter outline of the dissertation as a whole,
one I hope the reader will find useful.

Section \ref{sec:1.2 Dynamical-Systems-Terminology} gives a brief
exposition of basic notions from the theory of discrete dynamical
systems\textemdash orbit classes and the like. The slightly longer
Section \ref{sec:1.3 Crash course in ultrametric analysis} is a crash
course in $p$-adic numbers and ultrametric analysis, written under
the assumption that the reader has never heard of either topic. That
being said, even an expert should take a gander at Subsections \ref{subsec:1.3.3 Field-Extensions-=00003D000026},
\ref{subsec:1.3.4 Pontryagin-Duality-and}, and \ref{subsec:1.3.5 Hensel's-Infamous-Blunder},
because my approach will frequently bring us into dangerous territory
where we mix-and-match $p$-adic topologies, either with the real
or complex topology, or for different primes $p$.

\section{\emph{\label{sec:1.1 A-Guide-for-the-Perplexed}A Guide for the Perplexed}}

\subsection{\label{subsec:1.1.1 Notational-Idiosyncrasies}Notational Idiosyncrasies}

In order to resolve the endless debate over whether or not $0$ is
a natural number, given any real number $x$, I write $\mathbb{N}_{x}$\nomenclature{$\mathbb{N}_{k}$}{the set of integers greater than or equal to $k$ \nopageref}
to denote the set of all elements of $\mathbb{Z}$ which are $\geq x$.
An expression like $\mathbb{N}_{x}^{r}$, where $r$ is an integer
$\geq1$, then denotes the cartesian product of $r$ copies of $\mathbb{N}_{x}$.
Likewise, $\mathbb{Z}^{d}$ denotes the cartesian product of $d$
copies of $\mathbb{Z}$.

As a matter of personal preference, I use $\mathfrak{fraktur}$ (fraktur,
\textbackslash mathfrak) font ($\mathfrak{z},\mathfrak{y}$, etc.)
to denote $p$-adic (or $q$-adic) variables, and use the likes of
$\mathfrak{a},\mathfrak{b},\mathfrak{c}$ to denote $p$-adic (or
$q$-adic) constants. I prefer my notation to be as independent of
context as possible, and therefore find it useful and psychologically
comforting to distinguish between archimedean and non-archimedean
quantities. Also, to be frank, I think the fraktur letters look \emph{cool}.

I write $\left[\cdot\right]_{p^{n}}$ to denote the ``output the
residue of the enclosed object modulo $p^{n}$'' operator. The output
is \emph{always }an element of the set of integers $\left\{ 0,\ldots,p^{n}-1\right\} $.
$\left[\cdot\right]_{1}$ denotes ``output the residue of the enclosed
object modulo $1$''; hence, any integer (rational or $p$-adic)
is always sent by $\left[\cdot\right]_{1}$ to $0$. $\left|\cdot\right|_{p}$
and $\left|\cdot\right|_{q}$ denote the standard $p$-adic and $q$-adic
absolute values, respectively, while $v_{p}$ and $v_{q}$ denote
the standard $p$-adic and $q$-adic valuations, respectively, with
$\left|\cdot\right|_{p}=p^{-v_{p}\left(\cdot\right)}$, and likewise
for $\left|\cdot\right|_{q}$ and $v_{q}$. Similarly, I write $\mathbb{Z}/p\mathbb{Z}$
as a short-hand for the set $\left\{ 0,\ldots,p-1\right\} $, in addition
to denoting the group thereof formed by addition modulo $1$.

I write $\hat{\mathbb{Z}}_{p}$ to denote the Pontryagin dual of the
additive group of $p$-adic integers. I identify $\hat{\mathbb{Z}}_{p}$
with the set $\mathbb{Z}\left[1/p\right]/\mathbb{Z}$ of rational
numbers in $\left[0,1\right)$ of the form $k/p^{n}$ where $n$ is
an integer $\geq0$ and where $k$ is a non-negative integer which
is either $0$ or co-prime to $p$. $\hat{\mathbb{Z}}_{p}=\mathbb{Z}\left[1/p\right]/\mathbb{Z}$
is made into an additive group with the operation of addition modulo
$1$. All functions defined on $\hat{\mathbb{Z}}_{p}$ are assumed
to be $1$-periodic; that is, for such a function $\hat{\phi}$, $\hat{\phi}\left(t+1\right)=\hat{\phi}\left(t\right)$
for all $t\in\hat{\mathbb{Z}}_{p}$. $\left\{ \cdot\right\} _{p}$
denotes the $p$-adic fractional part, viewed here as a homomorphism
from the additive group $\mathbb{Q}_{p}$ of $p$-adic rational numbers
to the additive group $\mathbb{Q}/\mathbb{Z}$ of rational numbers
in $\left[0,1\right)$, equipped with the operation of addition modulo
$1$. Because it will be frequently needed, I write $\mathbb{Z}_{p}^{\prime}$
to denote $\mathbb{Z}_{p}\backslash\mathbb{N}_{0}$ (the set of all
$p$-adic integers which are not in $\left\{ 0,1,2,3,\ldots\right\} $).
I use the standard convention of writing $\mathcal{O}_{\mathbb{F}}$
to denote the ring of integers of a number field $\mathbb{F}$\nomenclature{$\mathcal{O}_{\mathbb{F}}$}{the ring of $\mathbb{F}$-integers \nopageref}.

With regard to embeddings (and this discussion will be repeated in
Subsections \ref{subsec:1.3.3 Field-Extensions-=00003D000026} and
\ref{subsec:1.3.4 Pontryagin-Duality-and}), while an algebraist might
be comfortable with the idea that for a primitive third root of unity
$\xi$ in $\mathbb{Q}_{7}$, the expression $\left|2-3\xi\right|_{7}$
is not technically defined, for my analytic purposes, this simply
\emph{will not do}. As such, throughout this dissertation, given an
odd prime $q$, in addition to writing $e^{2\pi i/\left(q-1\right)}$
to denote the complex number $\cos\left(2\pi/\left(q-1\right)\right)+i\sin\left(2\pi/\left(q-1\right)\right)$,
I also write $e^{2\pi i/\left(q-1\right)}$ to denote the \emph{unique}
primitive $\left(q-1\right)$th root of unity $\xi$ in $\mathbb{Z}_{q}^{\times}$
so that the value of $\xi$ modulo $q$ (that is, the first digit
in the $q$-adic representation of $\xi$) is the smallest integer
in $\left\{ 2,\ldots,q-2\right\} $ which is a primitive root modulo
$q-1$.

By far, the most common convention for expressing numerical congruences
is of the form: 
\begin{equation}
x=k\mod p
\end{equation}
However, I have neither any intention nor any inclination of using
this cumbersome notation. Instead, I write $\overset{p}{\equiv}$
to denote congruence modulo $p$. This is a compact and extremely
malleable notation. For example, given two elements $s,t\in\hat{\mathbb{Z}}_{p}$,
I write $s\overset{1}{\equiv}t$ \nomenclature{$\overset{1}{\equiv}$}{congruence modulo $1$ \nopageref}to
indicate that $s$ is congruent to $t$ modulo $1$; i.e., $s-t\in\mathbb{Z}$,
which\textemdash of course\textemdash makes $s$ and $t$ different
representatives of the same element of $\hat{\mathbb{Z}}_{p}$ (example:
$s=1/p$ and $t=\left(p+1\right)/p$). This congruence notation is
also compatible with $p$-adic numbers (integer or rational): $\mathfrak{z}\overset{p^{n}}{\equiv}k$\nomenclature{$\overset{p^{n}}{\equiv}$}{congruence modulo $p^{n}$ \nopageref}
holds true if and only if $\mathfrak{z}$ is an element of $k+p^{n}\mathbb{Z}_{p}$.

My congruence notation is particularly useful in conjunction with
the indispensable Iverson bracket notation. Given a statement or expression
$S$, the Iverson bracket is the convention of writing $\left[S\right]$
to denote a quantity which is equal to $1$ when $S$ is true, and
which is $0$ otherwise. For example, letting $k$ be an integer constant
and letting $\mathfrak{z}$ be a $p$-adic integer variable, $\left[\mathfrak{z}\overset{p^{n}}{\equiv}k\right]$
is then the indicator function for the set $k+p^{n}\mathbb{Z}_{p}$,
being equal to $1$ if $\mathfrak{z}$ is congruent to $k$ mod $p^{n}$
(that is, $\left[\mathfrak{z}\right]_{p^{n}}=k$), and being $0$
otherwise. Given a point $s\in\hat{\mathbb{Z}}_{p}$, I write $\mathbf{1}_{s}\left(t\right)$
to denote the function $\left[t\overset{1}{\equiv}s\right]$ of the
variable $t\in\hat{\mathbb{Z}}_{p}$. That is, $\mathbf{1}_{s}\left(t\right)$
is $1$ if and only if $t-s$ is in $\mathbb{Z}$, and is $0$ for
all other values of $t$. I write $\mathbf{1}_{\mathfrak{a}}\left(\mathfrak{z}\right)$
to denote the function $\left[\mathfrak{z}=\mathfrak{a}\right]$ of
the variable $\mathfrak{z}\in\mathbb{Z}_{p}$, which is $1$ if and
only if $\mathfrak{z}=\mathfrak{a}$ and is $0$ otherwise.

Because the absolute-value-induced-metric-topology with which we will
make sense of limits of sequences or sums of infinite series will
often frequently change, I have cultivated a habit of writing the
space in which the convergence occurs over the associated equal sign.
Thus, $\lim_{n\rightarrow\infty}f\left(x_{n}\right)\overset{\mathbb{Z}_{p}}{=}c$
means that the sequence $f\left(x_{n}\right)$ converges to $c$ in
$p$-adic absolute value. Note that this necessarily forces $c$ and
the $f\left(x_{n}\right)$s to be quantities whose $p$-adic absolute
values are meaningfully defined. Writing $\overset{\overline{\mathbb{Q}}}{=}$
or $\overset{\mathbb{Q}}{=}$ indicates that no limits are actually
involved; the sums in question are not infinite, but consist of finitely
many terms, all of which are contained in the field indicated above
the equals sign.

I write $\overset{\textrm{def}}{=}$\nomenclature{$\overset{\textrm{def}}{=}$}{"by definition" \nopageref}
to mean ``by definition''.

The generalizations of these notations to the multi-dimensional case
are explained on pages \pageref{nota:First MD notation batch}, \pageref{nota:Second batch},
and \pageref{nota:third batch}. These also include extensions of
other notations introduced elsewhere to the multi-dimensional case,
as well as the basic algebraic conventions we will use in the multi-dimensional
case (see page \pageref{nota:First MD notation batch}). Also, though
it will be used infrequently, I write \nomenclature{$\ll$}{Vinogradov notation (Big O)}$\ll$
and $\gg$ to denote the standard Vinogradov notation (an alternative
to Big $O$), indicating a bound on a real quantity which depends
on an implicit constant ($f\left(x\right)\ll g\left(x\right)$ iff
$\exists$ $K>0$ so that $f\left(x\right)\leq Kg\left(x\right)$).

Given a map $T$, I write $T^{\circ n}$ to denote the composition
of $n$ copies of $T$.

Finally, there is the matter of self-reference. Numbers enclosed in
(parenthesis) refer to the numeric designation assigned to one of
the many, many equations and expressions that occur in this document.
Whenever a Proposition, Lemmata, Theorem, Corollary or the like is
invoked, I do so by writing it in bold along with the number assigned
to it; ex: \textbf{Theorem 1.1}. Numbers in parenthesis, like (4.265),
refer to the equation with that particular number as its labels. References
are cited by their assigned number, enclosed in {[}brackets{]}. I
also sometimes refer to vector spaces as linear spaces; ergo, a $\overline{\mathbb{Q}}$-linear
space is a vector space over $\overline{\mathbb{Q}}$.

\subsection{Some Much-Needed Explanations \label{subsec:1.1.2 Some-Much-Needed-Explanations}}

Given that this document\textemdash my PhD Dissertation\textemdash is
over four-hundred fifty pages long, I certainly owe my audience an
explanation. I have explanations in spades, and those given in this
incipit are but the first of many. I have spent the past five years
studying the Collatz Conjecture from a variety of different analytical
approaches. My journey eventually led me into the \emph{ultima Thule
}of non-Archimedean analysis, which\textemdash to my surprise\textemdash turned
out to be far less barren and far more fruitful than I could have
ever suspected. With this monograph, I chronicle my findings.

To set the mood\textemdash this \emph{is }a study of the Collatz Conjecture,
after all\textemdash let me give my explanations in threes. Besides
the obvious goal of securing for myself a hard-earned doctoral degree,
the purpose of this dissertation is three-fold: 
\begin{itemize}
\item To present a analytical programme for studying a wide range of Collatz-type
arithmetical dynamical systems in a unified way; 
\item To detail new, unforeseen aspects of the subspecialty of non-archimedean
analysis involving the study of functions from the $p$-adic integers
to the $q$-adic (complex) numbers (I call this ``\emph{$\left(p,q\right)$-adic
analysis}''); 
\item To use $\left(p,q\right)$-adic analysis to study the dynamics of
Collatz-type maps\textemdash primarily the question of whether or
not a given integer is a periodic point of such a map. 
\end{itemize}
In this respect\textemdash as indicated by its title\textemdash this
dissertation is as much about the novel methods of non-archimedean
analysis as it is about using those methods to describe and better
understand Collatz-type dynamical systems. The reader has a winding
road ahead of them, monotonically increasing in complexity over the
course of their six-chapter journey.

To me, among the many ideas we will come across, three innovations
stand out among the rest: 
\begin{itemize}
\item Given a kind of a Collatz-type map $H:\mathbb{Z}\rightarrow\mathbb{Z}$
which I call a \textbf{Hydra map} (see page \pageref{def:p-Hydra map}),
provided $H$ satisfies certain simple qualitative properties, there
are distinct primes $p,q$ and a function $\chi_{H}:\mathbb{Z}_{p}\rightarrow\mathbb{Z}_{q}$\textemdash the
\textbf{Numen} of $H$ (see page \pageref{def:Chi_H on N_0 in strings})\textemdash with
the property that an integer $x\in\mathbb{Z}\backslash\left\{ 0\right\} $
is a periodic point of $H$ \emph{if and only if }there is a $p$-adic
integer $\mathfrak{z}_{0}\in\left(\mathbb{Q}\cap\mathbb{Z}_{p}\right)\backslash\left\{ 0,1,2,3,\ldots\right\} $
so that $\chi_{H}\left(\mathfrak{z}_{0}\right)=x$. I call this phenomenon
\textbf{the} \textbf{Correspondence Principle}\footnote{Actually, we can say more: modulo certain simple qualitative conditions
on $H$, if there is an ``irrational'' $p$-adic integer $\mathfrak{z}_{0}\in\mathbb{Z}_{p}\backslash\mathbb{Q}$
so that $\chi_{H}\left(\mathfrak{z}_{0}\right)\in\mathbb{Z}$, then
the sequence of iterates $\chi_{H}\left(\mathfrak{z}_{0}\right),H\left(\chi_{H}\left(\mathfrak{z}_{0}\right)\right),H\left(H\left(\chi_{H}\left(\mathfrak{z}_{0}\right)\right)\right),\ldots$
is unbounded with respect to the standard archimedean absolute value
on $\mathbb{R}$. I have not been able to establish the converse of
this result, however.} (see page \pageref{cor:CP v4}). 
\item My notion of \textbf{frames} (presented in Subsection \ref{subsec:3.3.3 Frames}),
which provide a formalism for dealing with functions $\chi:\mathbb{Z}_{p}\rightarrow\mathbb{Z}_{q}$
(where $p$ and $q$ are distinct primes) defined by infinite series
such that the \emph{topology of convergence} used to sum said series
\emph{varies from point to point}. The archetypical example is that
of a series in a $p$-adic integer variable $\mathfrak{z}$ which
converges in the $q$-adic topology whenever $\mathfrak{z}\in\mathbb{Z}_{p}\backslash\left\{ 0,1,2,3,\ldots\right\} $
and converges in the topology of the reals whenever $\mathfrak{z}\in\left\{ 0,1,2,3,\ldots\right\} $.
Using frames, we can significantly enlarge the class of $\left(p,q\right)$-adic
functions which can be meaningfully ``integrated''\textemdash I
call such functions \textbf{quasi-integrable }functions. With these
tools, I then establish corresponding $\left(p,q\right)$-adic generalizations
of the venerable \textbf{Wiener Tauberian Theorem }(\textbf{WTT})
from harmonic analysis (Subsection \ref{subsec:3.3.7 -adic-Wiener-Tauberian}). 
\item A detailed $\left(p,q\right)$-adic Fourier-analytical study of $\chi_{H}$.
I show that, for many Hydra maps $H$, the numen $\chi_{H}$ is in
fact quasi-integrable. Not only does this come with explicit non-trivial
$\left(p,q\right)$-adic series representations of $\chi_{H}$, and
formulae for its Fourier transforms ($\hat{\chi}_{H}$), but\textemdash in
conjunction with the Correspondence Principle and the $\left(p,q\right)$-adic
Wiener Tauberian Theorem\textemdash I show that the question ``is
$x\in\mathbb{Z}\backslash\left\{ 0\right\} $ a periodic point of
$H$?'' is essentially equivalent to ``is the span of the translates
of $\hat{\chi}_{H}\left(t\right)-x\mathbf{1}_{0}\left(t\right)$ dense
in the appropriate function space?'' In particular, the answer to
the former is ``no'' whenever the density of the span of the translates
occurs (\textbf{Theorem \ref{thm:Periodic Points using WTT}}\vpageref{thm:Periodic Points using WTT}).
I call this the \textbf{Tauberian Spectral Theorem }for $H$.\footnote{If the spans of the translates are \emph{not }dense, then, at present,
my methods then guarantee that $x$ is either a periodic point of
$H$ or the sequence of iterates $x,H\left(x\right),H\left(H\left(x\right)\right),\ldots$
is unbounded with respect to the standard archimedean absolute value
on $\mathbb{R}$.} 
\end{itemize}
I consider the above my ``main results''. It is worth mentioning
that these results hold not just for Collatz-type maps $H:\mathbb{Z}\rightarrow\mathbb{Z}$,
but for Collatz-type maps $H:\mathbb{Z}^{d}\rightarrow\mathbb{Z}^{d}$
for any integer $d\geq1$\textemdash the $d\geq2$ case being the
``\textbf{multi-dimensional case}'' (Chapters 5 \& 6). These arise
from studying generalizations of Collatz-type maps from $\mathbb{Z}$
to $\mathcal{O}_{\mathbb{F}}$, the ring of integers of $\mathbb{F}$,
a finite-dimensional field extension of $\mathbb{Q}$; such maps were
first considered by Leigh (1981) and subsequently explored by his
colleague K. R. Matthews \cite{Leigh's article,Matthews' Leigh Article,Matthews' slides}.
I briefly corresponded with Matthews in 2017, and again in 2019\textemdash the
man is a veteran Collatz scholar\textemdash and he informed me that
other than the works just cited, he knew of no others that have pursued
Collatz-type maps on rings of algebraic integers.

My final explanatory triplet addresses the forty-foot dragon in the
corner of the room with page after page of \emph{pages }clutched between
his claws: \emph{why is this dissertation over four-hundred-fifty
($450$) pages long?} 
\begin{itemize}
\item Unusually for ``hard'' analysis\textemdash generally a flexible,
well-mixed, free-ranged discipline\textemdash the \emph{sub}-subspecialty
of $\left(p,q\right)$-adic analysis is incredibly arcane, even within
the already esoteric setting of non-archimedean analysis as a whole.
While it is not anywhere near the level of, say, algebraic geometry,
non-archimedean analysis behooves the reader a good deal of legwork.
Many readers will find in these pages at least one (and perhaps many)
ideas they never would have thought of considering; common courtesy
and basic decency demands I give a thorough accounting of it. That
my work is left field even by the standards of non-archimedean analysis
only strengths my explanatory resolve. 
\item It is a general conviction of mine\textemdash and one strongly held,
at that\textemdash that thought or knowledge without warmth, patience,
openness, charity, and accessibility is indistinguishable from Brownian
motion or Gaussian noise. Mathematics is difficult and frustrating
enough as it is. I refuse to allow my writing\textemdash or, well,
any potential \emph{lack} thereof\textemdash to cause or perpetuate
pain, suffering, or discomfort which could otherwise have been avoided.
Also\textemdash as befits a budding new theory\textemdash many of
the most important methods and insights in this document arise directly
from demanding computational considerations. The formulas provide
the patterns, so it would make little sense (nor kindness, for that
matter) not to explore and discuss them at length. 
\item I am a \emph{wordy} person, and\textemdash in addition to my passion
for mathematics\textemdash I have a passion for writing. The Collatz
obsession infected me just as I was finishing up my second novel (woefully
stillborn), and persisted throughout the writing of my \emph{third}
novel, coeval with the research and writing that led to my dissertation. 
\end{itemize}
All I can say is, I hope you enjoy the ride.

\subsection{\label{subsec:1.1.3 An-Outline-of}An Outline of this Dissertation}

\subsubsection*{Chapter 1}

If this dissertation were a video game, Chapter 1 would be the Tutorial
level. Aside from my verbal excess, Chapter 1 provides the first significant
portion of necessary background material. Section \ref{sec:1.2 Dynamical-Systems-Terminology}
reviews elementary facts and terminology from the theory of dynamical
systems (periodic points, divergent trajectories, orbit classes, etc.).
The largest section of Chapter 1 is \ref{sec:1.3 Crash course in ultrametric analysis},
which is basically a crash course in ultrametric analysis, with a
special focus on the $p$-adic numbers in Subsection \ref{subsec:1.3.1. The-p-adics-in-a-nutshell}.
Most of the material from Subsection \ref{subsec:1.3.2. Ultrametrics-and-Absolute}
is taken from Schikhof's \emph{Ultrametric Calculus }\cite{Ultrametric Calculus},
which I highly recommend for anyone interested in getting a first
taste of the subject as a whole. Subsection \ref{subsec:1.3.3 Field-Extensions-=00003D000026}
briefly mentions $p$-adic field extensions, Galois actions, and spherical
completeness. Subsection \ref{subsec:1.3.4 Pontryagin-Duality-and}
introduces Pontryagin duality and the Fourier transform of a \emph{complex-valued
}function on a locally compact abelian group, again laser-focused
on the $p$-adic case. While Subsections \ref{subsec:1.3.1. The-p-adics-in-a-nutshell}
through \ref{subsec:1.3.4 Pontryagin-Duality-and} can be skipped
by a reader already familiar with their contents, I must recommend
that \emph{everyone }read \ref{subsec:1.3.5 Hensel's-Infamous-Blunder}.
Throughout this dissertation, we will skirt grave danger by intermixing
and interchanging archimedean and non-archimedean topologies on $\overline{\mathbb{Q}}$.
Subsection \ref{subsec:1.3.5 Hensel's-Infamous-Blunder} explains
this issue in detail, using Hensel's infamous $p$-adic ``proof''
of the irrationality of $e$ as a springboard. The reason we will
not suffer a similar, ignominious fate is due to the so-called ``universality''
of the geometric series.

\subsubsection*{Chapter 2}

Chapter 2 presents the titular Hydra maps, my name for the class of
Collatz-type maps studied in this dissertation, as well as the other
Collatz-related work discussed in the Preface. Similar generalizations
have been considered by \cite{Matthews and Watts,Moller's paper (german),dgh paper};
the exposition given here is an attempt to consolidate these efforts
into a single package. In order to avoid excessive generality, the
maps to be considered will be grounded over the ring of integers ($\mathbb{Z}$)
and finite-degree extensions thereof. In his slides, Matthews considers
Collatz-type maps defined over rings of polynomials with coefficients
in finite fields \cite{Matthews' slides}. While, in theory, it may
be possible to replicate the $\chi_{H}$ construction by passing from
rings of polynomials with coefficients in a finite field $\mathbb{F}$
to the local field of formal Laurent series with coefficients in $\mathbb{F}$,
in order to avoid the hassle of dealing with fields of positive characteristic,
this dissertation will defer exploration of Hydra maps on fields of
positive characteristic to a later date.

My presentation of Hydra maps deals with two types: a ``one-dimensional''
variant\textemdash Hydra maps defined on $\mathbb{Z}$\textemdash and
a ``multi-dimensional'' variant; given a number field $\mathbb{F}$,
a Hydra map defined on the ring of integers $\mathcal{O}_{\mathbb{F}}$
is said to be of dimension $\left[\mathbb{F}:\mathbb{Q}\right]=\left[\mathcal{O}_{\mathbb{F}}:\mathbb{Z}\right]$.
Equivalently, the multi-dimensional variant can be viewed as maps
on $\mathbb{Z}^{d}$. The study of the multi-dimensional variant is
more or less the same as the study of the one-dimensional case, albeit
with a sizable dollop of linear algebra.

After presenting Hydra maps, Chapter 2 proceeds to demonstrate how,
to construct the function $\chi_{H}$\textemdash the ``numen''\footnote{My original instinct was to call $\chi_{H}$ the ``characteristic
function'' of $H$ (hence the use of the letter $\chi$). However,
in light of Tao's work\textemdash where, in essence, one works with
the characteristic function (in the Probability-theoretic sense) of
$\chi_{H}$, this name proved awkward. I instead chose to call it
the ``numen'' of $H$ (plural: ``numina''), from the Latin, meaning
``the spirit or divine power presiding over a thing or place''.} of the Hydra map $H$\textemdash and establishes the conditions required
of $H$ in order for $\chi_{H}$ to exist. Provided that $H\left(0\right)=0$,
we will be able to construct $\chi_{H}$ as a function from $\mathbb{N}_{0}$
to $\mathbb{Q}$. Additional conditions (see page \pageref{def:Qualitative conditions on a p-Hydra map})
are then described which will guarantee the existence of a unique
extension/continuation of $\chi_{H}$ to a $q$-adic valued function
of a $p$-adic integer variable for appropriate $p$ and $q$ determined
by $H$ and the corresponding generalization for multi-dimensional
Hydra maps.

This then leads us to the cornerstone of this dissertation: the \textbf{Correspondence
Principle}, which I give in four variants, over pages \pageref{thm:CP v1}-\pageref{cor:CP v4}.
This principle is the fact that\textemdash provided $H$ satisfies
an extra condition beyond the one required for $\chi_{H}$ to have
a non-archimedean extension\textemdash there will be a one-to-one
correspondence between the periodic points of $H$ and the rational
integer values attained by the non-archimedean extension of $\chi_{H}$.
In certain cases, this principle even extends to cover the points
in divergent orbits of $H$!

Understandably, the Correspondence Principle motivates all the work
that follows. Chapter 2 concludes with a brief discussion of the relationship
between $\chi_{H}$ and other work, principally Tao's paper \cite{Tao Probability paper}
(discussed in this dissertation starting \vpageref{subsec:Connections-to-Tao}),
transcendental number theory, the resolution of \textbf{Catalan's
Conjecture }by P. Mih\u{a}ilescu in 2002 \cite{Cohen Number Theory},
and the implications thereof for Collatz studies\textemdash in particular,
the fundamental inadequacy of the current state of transcendental
number theory that precludes it from being able to make meaningful
conclusions on the Collatz Conjecture.

\subsubsection*{Chapter 3}

This is the core theory-building portion of this dissertation, and
is dedicated to presenting the facets of non-archimedean analysis
which I believe will finally provide Hydra maps and the contents of
Chapter 2 with a proper theoretical context. After a historical essay
on the diverse realizations of ``non-archimedean analysis'' (Subsection
\ref{subsec:3.1.1 Some-Historical-and}), the first section of Chapter
3 is devoted to explaining the essentials of what I call ``$\left(p,q\right)$-adic
analysis'', the specific sub-discipline of non-archimedean analysis
which concerns functions of a $p$-adic variable taking values in
$\mathbb{Q}_{q}$ (or a field extension thereof), where $p$ and $q$
are distinct primes. Since much of the literature in non-archimedean
analysis tends to be forbiddingly generalized\textemdash much more
so than the ``hard'' analyst is usually comfortable with\textemdash I
have taken pains to minimize unnecessary generality as much as possible.
As always, the goal is to make the exposition of the material as accessible
as possible without compromising rigor.

The other subsections of \ref{sec:3.1 A-Survey-of} alternate between
the general and the specific. \ref{subsec:3.1.2 Banach-Spaces-over}
addresses the essential aspects of the theory of non-archimedean Banach
spaces, with an emphasis on which of the major results of the archimedean
theory succeed or fail. Subsection \ref{subsec:3.1.3 The-van-der}
presents the extremely useful \textbf{van der Put basis }and \textbf{van
der Put series representation }for the space of continuous $\mathbb{F}$-valued
functions defined on the $p$-adic integers, where $\mathbb{F}$ is
a metrically complete valued field (either archimedean or non-archimedean).
The existence of this basis makes computations in $\left(p,q\right)$-adic
analysis seem much more like its classical, archimedean counterpart
than is usually the case in more general non-archimedean analysis.

Subsections \ref{subsec:3.1.4. The--adic-Fourier} and \ref{subsec:3.1.5-adic-Integration-=00003D000026}
are devoted to the $\left(p,q\right)$-adic Fourier(-Stieltjes) transform.
The contents of these subsections are a mix of W. M. Schikhof's 1967
PhD dissertation,\emph{ Non-Archimedean Harmonic Analysis} \cite{Schikhof's Thesis},
and my own independent re-discovery of that information prior to learning
that \emph{Non-Archimedean Harmonic Analysis} already existed. For
the classical analyst, the most astonishing feature of this subject
is the near-total absence of \emph{nuance}: $\left(p,q\right)$-adic
functions are integrable (with respect to the $\left(p,q\right)$-adic
Haar probability measure, whose construction is detailed in Subsection
\ref{subsec:3.1.5-adic-Integration-=00003D000026}) if and only if
they are continuous, and are continuous if and only if they have a
$\left(p,q\right)$-adic Fourier series representation which converges
uniformly \emph{everywhere}.

Whereas Subsections \ref{subsec:3.1.4. The--adic-Fourier} and \ref{subsec:3.1.5-adic-Integration-=00003D000026}
are primarily oriented toward practical issues of computation, \ref{subsec:3.1.6 Monna-Springer-Integration}
gives a thorough exposition of the chief theory of integration used
in non-archimedean analysis as a whole, that of\textbf{ }the\textbf{
Monna-Springer Integral}. My primary motivation for introducing Monna-Springer
theory is to make it better known, with an eye toward equipping other
newcomers to the subject with a guide for decoding the literature.

The novel content of Chapter 3 occurs in Sections \ref{sec:3.2 Rising-Continuous-Functions}
and \ref{sec:3.3 quasi-integrability}, particularly the latter. Section
\ref{sec:3.2 Rising-Continuous-Functions} presents and expands upon
material which is given but a couple exercise's\footnote{Exercise 62B, page 192, and its surroundings\cite{Ultrametric Calculus}.}
worth of attention in Schikhof's \emph{Ultrametric Calculus}. Instead,
in Subsection \ref{subsec:3.2.1 -adic-Interpolation-of}, we will
develop into a theory of what I call \textbf{rising-continuous functions}.
In brief, letting $\mathbb{F}$ be a metrically complete valued field,
these are functions $\chi:\mathbb{Z}_{p}\rightarrow\mathbb{F}$ satisfying
the point-wise limit condition:
\begin{equation}
\chi\left(\mathfrak{z}\right)\overset{\mathbb{F}}{=}\lim_{n\rightarrow\infty}\chi\left(\left[\mathfrak{z}\right]_{p^{n}}\right),\textrm{ }\forall\mathfrak{z}\in\mathbb{Z}_{p}
\end{equation}
where, as indicated by the $\overset{\mathbb{F}}{=}$, the convergence
occurs in the topology of $\mathbb{F}$. These functions arise naturally
when one considers interpolating functions on $\mathbb{N}_{0}$ to
ones on $\mathbb{Z}_{p}$. After this exposé, Subsection \ref{subsec:3.2.2 Truncations-=00003D000026-The}
demonstrates that the set of rising-continuous functions forms a non-archimedean
Banach algebra which extends the Banach algebra of continuous $\left(p,q\right)$-adic
functions. The conditions required for a rising-continuous function
to be a unit\footnote{I.E., the reciprocal of the function is also rising-continuous.}
of this algebra are investigated, with the \textbf{Square Root Lemma}
(see page \pageref{lem:square root lemma}) of Subsection \ref{subsec:3.2.2 Truncations-=00003D000026-The}
providing a necessary and sufficient condition. In the \emph{sub}-subsection
starting \vpageref{subsec:3.2.2 Truncations-=00003D000026-The}, the
theory of Berkovitch spaces from $p$-adic (algebraic) geometry is
briefly invoked to prove some topological properties of rising-continuous
functions. Subsection \ref{subsec:3.2.2 Truncations-=00003D000026-The}
also introduces the crucial construction I call \textbf{truncation}.
The $N$th truncation of a $\left(p,q\right)$-adic function $\chi$,
denoted $\chi_{N}$, is the function defined by: 
\begin{equation}
\chi_{N}\left(\mathfrak{z}\right)\overset{\textrm{def}}{=}\sum_{n=0}^{p^{N}-1}\chi\left(n\right)\left[\mathfrak{z}\overset{p^{N}}{\equiv}n\right]
\end{equation}
Even if $\chi$ is merely rising-continuous, its $N$th truncation
is continuous\textemdash in fact, locally constant\textemdash for
all $N$, and converges point-wise to $\chi$ everywhere as $N\rightarrow\infty$.
We also explore the interplay between truncation, van der Put series,
and the Fourier transform as it regards continuous $\left(p,q\right)$-adic
functions.

The heart of section \ref{sec:3.3 quasi-integrability}, however,
is the exposition of my chief innovations non-archimedean analysis:
\textbf{frames }and \textbf{quasi-integrability}. Because an arbitrary
rising-continuous function is not going to be continuous, it will
fail to be integrable in the sense of Monna-Springer theory. This
non-integrability precludes a general rising-continuous function from
possessing a well-defined Fourier transform. Quasi-integrability is
the observation that this difficulty can, in certain cases, be overcome,
allowing us to significantly enlarge the class of $\left(p,q\right)$-adic
functions we can meaningfully integrate. Aside from facilitating a
more detailed analysis of $\chi_{H}$, these innovations demonstrate
that $\left(p,q\right)$-adic analysis is far more flexible\textemdash and,
one would hope, \emph{useful}\textemdash than anyone had previously
thought.

Instead of diving right into abstract axiomatic definitions, in Subsection
\ref{subsec:3.3.1 Heuristics-and-Motivations}, I take the Arnoldian
approach (otherwise known as the \emph{sensible} approach) of beginning
with an in-depth examination of the specific computational examples
that led to and, eventually, necessitated the development of frames
and quasi-integrability. To give a minor spoiler, the idea of quasi-integrability,
in short, is that even if $\chi:\mathbb{Z}_{p}\rightarrow\mathbb{C}_{q}$
is not continuous, \emph{if we can find }a function $\hat{\chi}$
(defined on $\hat{\mathbb{Z}}_{p}$) so that the partial sums: 
\begin{equation}
\sum_{\left|t\right|_{p}\leq p^{N}}\hat{\chi}\left(t\right)e^{2\pi i\left\{ t\mathfrak{z}\right\} _{p}}
\end{equation}
converge to $\chi$ as $N\rightarrow\infty$ in an appropriate sense
for certain $\mathfrak{z}$, we can then use that $\hat{\chi}$ to
define \emph{a}\footnote{Alas, it is not \emph{the} Fourier transform; it is only unique modulo
the Fourier-Stieltjes transform of what we will call a \textbf{degenerate
measure}.}\textbf{ Fourier transform }of $\chi$. In this case, we say that
$\chi$ is \textbf{quasi-integrable}. The name is due to the fact
that the existence of $\hat{\chi}$ allows us to define the integral
$\int_{\mathbb{Z}_{p}}f\left(\mathfrak{z}\right)\chi\left(\mathfrak{z}\right)d\mathfrak{z}$
for any continuous $\left(p,q\right)$-adic function $f$. We view
$\chi$ as the ``derivative'' of the $\left(p,q\right)$-adic measure
$\chi\left(\mathfrak{z}\right)d\mathfrak{z}$; this makes $\hat{\chi}$
the Fourier-Stieltjes transform of this measure. The success of this
construction is of particular interest, given that, in general \cite{van Rooij - Non-Archmedean Functional Analysis},
the Radon-Nikodym theorem fails to hold in non-archimedean analysis.
As a brief\textemdash but fascinating\textemdash tangent, starting
from page \pageref{eq:Definition of (p,q)-adic Mellin transform},
I showcase how quasi-integrability can be used to meaningfully define
the \textbf{Mellin Transform }of a continuous $\left(p,q\right)$-adic
function.

Subsection \ref{subsec:3.3.2 The--adic-Dirichlet} shows how the limiting
process described above can be viewed as a $\left(p,q\right)$-adic
analogue of summability kernels, as per classical Fourier analysis.
All the most important classical questions in Fourier Analysis of
functions on the unit interval (or circle, or torus) come about as
a result of the failure of the \textbf{Dirichlet Kernel }to be a well-behaved
approximate identity. On the other hand, $\left(p,q\right)$-adic
analysis has no such troubles: the simple, natural $\left(p,q\right)$-adic
incarnation of the Dirichlet kernel is, in that setting, an approximate
identity\textemdash no questions asked!

In subsection \pageref{subsec:3.3.3 Frames}, the reader will be introduced
to the notion of \textbf{frames}.\textbf{ }A frame is, at heart, an
organizing tool that expedites and facilitates clear, concise discussion
of the ``spooky'' convergence behaviors of the partial sums of Fourier
series generated by $q$-adically bounded functions $\hat{\mathbb{Z}}_{p}\rightarrow\overline{\mathbb{Q}}$.
This subsection is primarily conceptual, as is its counterpart \ref{subsec:3.3.5 Quasi-Integrability},
where quasi-integrability with respect to a given frame is formally
defined. Subsection \ref{subsec:3.3.4 Toward-a-Taxonomy} makes the
first steps toward what will hopefully one day be an able-bodied theory
of $\left(p,q\right)$-adic measures, and introduces several (re-)summation
formulae for $\left(p,q\right)$-adic Fourier series that will take
central stage in our analysis of $\chi_{H}$. Subsection \ref{subsec:3.3.6 L^1 Convergence}
also shows that the hitherto neglected topic of absolute convergence
of $\left(p,q\right)$-adic integrals\textemdash that is, $L^{1}$
convergence\textemdash in non-archimedean analysis seems to provide
the logical setting for studying quasi-integrability at the theoretical
level. Lastly, Subsection \ref{subsec:3.3.7 -adic-Wiener-Tauberian}
introduces \textbf{Wiener's Tauberian Theorem }and states and proves
two new $\left(p,q\right)$-adic analogues of this chameleonic result,
one for continuous functions, and one for measures. Sub-subsection
\ref{subsec:A-Matter-of} then explores how the question of the existence
and continuity of the reciprocal of a $\left(p,q\right)$-adic function
can be stated and explored in terms of the spectral theory of \textbf{circulant
matrices}.

\subsubsection*{Chapter 4}

In Chapter 4, the techniques developed in Chapter 3 are mustered in
a comprehensive $\left(p,q\right)$-adic Fourier analysis of $\chi_{H}$.
After a bit of preparatory work (Subsection \ref{sec:4.1 Preparatory-Work--}),
Subsection \ref{sec:4.2 Fourier-Transforms-=00003D000026} launches
into the detailed, edifying computation needed to establish the quasi-integrability
of $\chi_{H}$ for a wide range of Hydra maps. The method of proof
is a kind of asymptotic analysis: by computing the the Fourier transforms
of the \emph{truncations} of $\chi_{H}$, we can identify fine structure
which exists independently of the choice of truncation. From this,
we derive an explicit formula for Fourier transform of $\chi_{H}$
and demonstrate its convergence with respect to the appropriate frame.
A recurring motif of this analysis is the use of functional equation
characterizations to confirm that our formulae do, in fact, represent
the desired functions. As a capstone, my $\left(p,q\right)$-adic
generalizations of Wiener's Tauberian Theorem are then invoked to
establish a near-equivalence of the problem of characterizing the
periodic points and divergent points of $H$ to determining those
values of $x$ for which the span of the translates of $\hat{\chi}_{H}\left(t\right)-x\mathbf{1}_{0}\left(t\right)$
are dense in $c_{0}\left(\hat{\mathbb{Z}}_{p},\mathbb{C}_{q}\right)$;
this is the content of the \textbf{Tauberian Spectral Theorem }for
$\chi_{H}$ (\textbf{Theorem \ref{thm:Periodic Points using WTT}},
on page \pageref{thm:Periodic Points using WTT}).

Section \ref{sec:4.3 Salmagundi} is a potpourri of miscellaneous
results about $\chi_{H}$. Some of these build upon the work of Section
\ref{sec:4.2 Fourier-Transforms-=00003D000026}; others which do not.
I suspect the lines of inquiry covered in Section \ref{sec:4.3 Salmagundi}
will be useful for future explorations of this subject. Of particular
interest is what I call the ``$L^{1}$ method'' (\textbf{Theorem
\ref{thm:The L^1 method}} on \pageref{thm:The L^1 method}), a means
for obtaining upper bounds on the \emph{archimedean }absolute value
of $\chi_{H}\left(\mathfrak{z}\right)$ over certain appropriately
chosen subsets of $\mathbb{Z}_{p}$.

\subsubsection*{Chapter 5}

Although my methods are applicable to Collatz-type maps on rings of
algebraic integers\textemdash equivalently $\mathbb{Z}^{d}$\textemdash I
prefer to cover the multi-dimensional case only after the one-dimensional
case, so as to let the underlying ideas shine through, unobscured
by technical considerations of multi-dimensional linear algebra. Chapter
5 contains the multi-dimensional case of the content of the one-dimensional
case as given in Chapters 2 and 3. Section \ref{sec:5.1 Hydra-Maps-on}
deals with the issue of converting a Collatz-type on some ring of
algebraic integers $\mathcal{O}_{\mathbb{F}}$ to an analogous map
acting on a lattice $\mathbb{Z}^{d}$ of appropriate dimension, with
Subsection \ref{subsec:5.1.2 Co=00003D0000F6rdinates,-Half-Lattices,-and}
providing the working definition for a Hydra map on $\mathbb{Z}^{d}$.
Section \ref{sec:5.2 The-Numen-of} is the multi-dimensional counterpart
of Chapter 2, and as such constructs $\chi_{H}$ and establishes the
multi-dimensional version of the \textbf{Correspondence Principle}.

Sections \ref{sec:5.3 Rising-Continuity-in-Multiple} and \ref{sec:5.4. Quasi-Integrability-in-Multiple}
treat rising-continuity and $\left(p,q\right)$-adic Fourier theory,
respectively, as they occur in the multi-dimensional context, with
the tensor product being introduced in Section \ref{sec:5.3 Rising-Continuity-in-Multiple}.
\ref{sec:5.4. Quasi-Integrability-in-Multiple} presents multi-dimensional
frames and quasi-integrability

\subsubsection*{Chapter 6}

Chapter 6 is a near-verbatim copy of Chapter 5, extended to the multi-dimensional
case, though without a corresponding extension of the salmagundi of
Section \ref{sec:4.3 Salmagundi}\textemdash I simply did not have
enough time.

\newpage{}

\section{\label{sec:1.2 Dynamical-Systems-Terminology}Elementary Theory of
Discrete Dynamical Systems}

In this brief section, we give an overview of the essential terminology
from the theory of discrete dynamical systems\index{discrete dynamical systems}.
A good reference for a beginner is \cite{Dynamical Systems}; it also
covers continuous systems.

For our purposes, it suffices to consider dynamical systems on $\mathbb{Z}^{d}$,
where $d$ is an integer $\geq1$, which is to say, there is some
map $T:\mathbb{Z}^{d}\rightarrow\mathbb{Z}^{d}$ whose iterates we
wish to study. Recall that we write $T^{\circ n}$ to denote the $n$-fold
iteration of $T$: 
\begin{equation}
T^{\circ n}=\underbrace{T\circ\cdots\circ T}_{n\textrm{ times}}
\end{equation}

\begin{defn}
Given an $\mathbf{x}=\left(x_{1},\ldots,x_{d}\right)\in\mathbb{Z}^{d}$,
the \textbf{forward orbit }of $\mathbf{x}$ under $T$ is the sequence
of iterates $\mathbf{x},T\left(\mathbf{x}\right),T\left(T\left(\mathbf{x}\right)\right),\ldots,T^{\circ n}\left(\mathbf{x}\right),\ldots$.
This is also called the \textbf{trajectory }of $\mathbf{x}$ under
$T$. More generally, a \textbf{trajectory }of $T$ refers to the
forward orbit of some $\mathbf{x}$ under $T$. 
\end{defn}
\vphantom{}

For us, the most important classical construction of the theory of
dynamical systems is that of the \textbf{orbit class}\index{orbit class}: 
\begin{defn}
\label{def:orbit class}Define a relation $\sim$ on $\mathbb{Z}^{d}$
by $\mathbf{x}\sim\mathbf{y}$ if and only if there are integers $m,n\geq0$
so that $T^{\circ m}\left(\mathbf{x}\right)=T^{\circ n}\left(\mathbf{y}\right)$,
where, recall, $T^{\circ0}$ is the identity map. It can be verified
that $\sim$ is an equivalence relation on $\mathbb{Z}^{d}$. The
distinct equivalence classes of $\mathbb{Z}^{d}$ under $\sim$ are
called the \textbf{irreducible orbit classes }of $T$ in $\mathbb{Z}^{d}$.
More generally, an \textbf{orbit class }is the union of irreducible
orbit classes of $T$; we say an orbit class is \textbf{reducible
}precisely when it can be written as a union of two or more non-empty
orbit classes, and say it is \textbf{irreducible }when no such decomposition
exists. 
\end{defn}
\vphantom{}

For those who have not seen this concept before, it is equivalent
to the idea of a watershed in earth sciences and hydrology: all elements
of a single irreducible orbit class of $T$ have the same long-term
behavior under iteration by $T$. 
\begin{prop}
A set $V\subseteq\mathbb{Z}^{d}$ is an orbit class of $T$ (possibly
reducible) if and only if $T^{-1}\left(V\right)=V$. 
\end{prop}
\vphantom{}

Using orbit classes, we can distinguish between sets according to
how $T$ behaves on them. 
\begin{defn}
We say $\mathbf{x}\in\mathbb{Z}^{d}$ is a \textbf{periodic point}
of $T$ whenever there is an integer $n\geq1$ so that $T^{\circ n}\left(\mathbf{x}\right)=\mathbf{x}$.
The \textbf{period }of $\mathbf{x}$ is the smallest integer $n\geq1$
for which $T^{\circ n}\left(\mathbf{x}\right)=\mathbf{x}$ holds true.
The set $\left\{ T^{\circ m}\left(\mathbf{x}\right):m\geq0\right\} $
(the forward orbit of $\mathbf{x}$) is then called the \textbf{cycle
generated by $\mathbf{x}$}. More generally, we say a set $\Omega\subseteq\mathbb{Z}^{d}$
is a \textbf{cycle }of $T$ whenever there is an $\mathbf{x}\in\Omega$
which is a periodic point of $T$ such that $\Omega=\left\{ T^{\circ m}\left(\mathbf{x}\right):m\geq0\right\} $. 
\end{defn}
\vphantom{}

The following facts are easily verified directly from the above definitions: 
\begin{prop}
Let $\Omega\subseteq\mathbb{Z}^{d}$ be a cycle of $\Omega$. Then:

\vphantom{}

I. $\Omega$ is a finite set.

\vphantom{}

II. Every element of $\Omega$ is a periodic point of $T$ of period
$\left|\Omega\right|$.

\vphantom{}

III. $T\mid_{\Omega}$ (the restriction of $T$ to $\Omega$) is a
bijection of $\Omega$.

\vphantom{}

IV. $\Omega\subseteq T^{-1}\left(\Omega\right)$

\vphantom{}

V. Either there exists an $\mathbf{x}\in\mathbb{Z}^{d}\backslash\Omega$
so that $T\left(\mathbf{x}\right)\in\Omega$ or $T^{-1}\left(\Omega\right)=\Omega$. 
\end{prop}
(IV) justifies the next bit of terminology: 
\begin{defn}
Given a cycle $\Omega$ of $T$, we say $\Omega$ is \textbf{attracting
}if $T^{-1}\left(\Omega\right)\backslash\Omega$ is non-empty\footnote{That is, if there are points not in $\Omega$ which $T$ sends into
$\Omega$.}; we say $\Omega$ is \textbf{isolated }(or, in analogy with continuous
dynamics, \textbf{repelling}) whenever $T^{-1}\left(\Omega\right)=\Omega$.

Additionally, we say $\mathbf{x}\in\mathbb{Z}^{d}$ is a \textbf{pre-periodic
point }of $T$ whenever there is an $n\geq0$ so that $T^{\circ n}\left(\mathbf{x}\right)$
is a periodic point of $T$ (where, recall, $T^{\circ0}$ is defined
to be the identity map). Given a cycle $\Omega$, we say $\mathbf{x}$
is \textbf{pre-periodic into $\Omega$ }if $T^{\circ n}\left(\mathbf{x}\right)\in\Omega$
occurs for some $n\geq0$. 
\end{defn}
\begin{prop}
Let $\mathbf{x}\in\mathbb{Z}^{d}$ be arbitrary. Then, either $\mathbf{x}$
is a pre-periodic point of $T$, or the forward orbit of $\mathbf{x}$
is unbounded. In particular, if $\mathbf{x}$ is not a pre-periodic
point, then: 
\begin{equation}
\lim_{n\rightarrow\infty}\left\Vert T^{\circ n}\left(\mathbf{x}\right)\right\Vert _{\infty}=+\infty\label{eq:Divergent Trajectory Definition}
\end{equation}
where for any $\mathbf{y}=\left(y_{1},\ldots,y_{d}\right)\in\mathbb{Z}^{d}$:
\begin{equation}
\left\Vert \mathbf{y}\right\Vert _{\infty}\overset{\textrm{def}}{=}\max\left\{ \left|y_{1}\right|,\ldots,\left|y_{d}\right|\right\} \label{eq:Definition of ell infinity norm on Z^d}
\end{equation}
\end{prop}
Proof: If $\mathbf{x}$ is not pre-periodic, then each element in
the forward orbit of $\mathbf{x}$ is unique; that is: for any integers
$m,n\geq0$, we have that $T^{\circ n}\left(\mathbf{x}\right)=T^{\circ m}\left(\mathbf{x}\right)$
occurs if and only if $m=n$. So, letting $N>0$ be arbitrary, note
that there are only finitely many $\mathbf{y}\in\mathbb{Z}^{d}$ with
$\left\Vert \mathbf{y}\right\Vert _{\infty}\leq N$. As such, if $\mathbf{x}$
is not pre-periodic, there can exist at most finitely many $n\geq0$
for which $\left\Vert T^{\circ n}\left(\mathbf{x}\right)\right\Vert _{\infty}\leq N$.
So, for any $N$, $\left\Vert T^{\circ n}\left(\mathbf{x}\right)\right\Vert _{\infty}$
will be strictly greater than $N$ for all sufficiently large $n$.
This establishes the limit (\ref{eq:Divergent Trajectory Definition}).

Q.E.D. 
\begin{defn}
We say $\mathbf{x}$ is a \textbf{divergent point }of $T$, or is
\textbf{divergent under $T$ }whenever the limit (\ref{eq:Divergent Trajectory Definition})\textbf{
}occurs. We use the terms \textbf{divergent trajectory}\index{divergent!trajectory}\textbf{
}and\textbf{ divergent orbit }to refer to a set of the form $\left\{ T^{\circ n}\left(\mathbf{x}\right):n\geq0\right\} $
where $\mathbf{x}$ is divergent under $T$. 
\end{defn}
\begin{defn}[\textbf{Types of orbit classes}]
Let $V\subseteq\mathbb{Z}^{d}$ be an orbit class of $T$. We say
$V$ is:

\vphantom{}

I. \textbf{attracting},\textbf{ }whenever it contains an attracting
cycle of $T$.

\vphantom{}

II. \textbf{isolated}, whenever it contains an isolated cycle of $T$.

\vphantom{}

III. \index{divergent!orbit class}\textbf{divergent}, whenever it
contains a divergent trajectory of $T$.
\end{defn}
\vphantom{}

It is a straight-forward and useful exercise to prove the following: 
\begin{lem}
If $V\subseteq\mathbb{Z}^{d}$ is an irreducible orbit class of $T$,
then one and only one of the following occurs:

\vphantom{}

I. $V$ is attracting. Additionally, there exists a unique attracting
cycle $\Omega\subset V$ so that every $\mathbf{x}\in V$ is pre-periodic
into $\Omega$.

\vphantom{}

II. $V$ is isolated. Additionally, $V$ itself is then an isolated
cycle of $T$.

\vphantom{}

III. $V$ is divergent. Additionally, $\lim_{n\rightarrow\infty}\left\Vert T^{\circ n}\left(\mathbf{x}\right)\right\Vert _{\infty}=+\infty$
for all $\mathbf{x}\in V$.
\end{lem}
\vphantom{}

Since the irreducible orbit classes of $T$ partition $\mathbb{Z}^{d}$,
they furnish a complete characterization of $T$'s dynamics on $\mathbb{Z}^{d}$. 
\begin{thm}
\label{thm:orbit classes partition domain}Let $T:\mathbb{Z}^{d}\rightarrow\mathbb{Z}^{d}$
be a map. Then, the irreducible orbit classes constitute a partition
of $\mathbb{Z}^{d}$ into at most countably infinitely many pair-wise
disjoint sets. In particular, for any $\mathbf{x}\in\mathbb{Z}^{d}$,
exactly one of the following occurs:

\vphantom{}

I. $\mathbf{x}$ is a periodic point of $T$ in an isolated cycle.

\vphantom{}

II. $\mathbf{x}$ is a pre-periodic point of $T$ which is either
in or eventually iterated \emph{into }an attracting cycle of $T$.

\vphantom{}

III. $\mathbf{x}$ has a divergent trajectory under $T$.
\end{thm}
\newpage{}

\section{\label{sec:1.3 Crash course in ultrametric analysis}$p$-adic Numbers
and Ultrametric Analysis}

THROUGHOUT THIS SECTION, $p$ DENOTES A PRIME NUMBER.

\subsection{\label{subsec:1.3.1. The-p-adics-in-a-nutshell}The $p$-adics in
a Nutshell}

We owe the conception of the $p$-adic numbers to the mathematical
work of\emph{ }Kurt Hensel\index{Hensel, Kurt} in 1897. Although
there are many ways of defining these numbers, I find Hensel's own
approach to be the most enlightening: the $p$-adics as \emph{power
series in the ``variable'' $p$} \cite{Gouvea's introudction to p-adic numbers book}.
In fact, many of the fundamental concepts associated with power series
turn out to be the inspiration for much of the paradigm-shifting developments
in number theory and algebraic geometry during the twentieth century.
Likely the most important of these is the simple observation that,
in general, a power series for an analytic function $f:\mathbb{C}\rightarrow\mathbb{C}$
about a given point $z_{0}\in\mathbb{C}$ gives us a formula for the
function which is valid only in a \emph{neighborhood} of $z_{0}$.
Quite often, we need to compute power series expansion about different
points (a process known as \emph{analytic continuation}) in order
to get formulae for the function on all the points in its domain.
That is to say, the function itself is a \emph{global }object, which
we study \emph{locally }(near the point $z_{0}$) by expanding it
in a power series about a point. 
\begin{defn}
The set of $p$-\textbf{adic integers}\index{$p$-adic!integers},
denoted \nomenclature{$\mathbb{Z}_{p}$}{the set of $p$-adic integers \nopageref}$\mathbb{Z}_{p}$,
is the set of all (formal) sums $\mathfrak{z}$ of the form: 
\begin{equation}
\mathfrak{z}=c_{0}+c_{1}p+c_{2}p^{2}+\ldots\label{eq:Definition of a p-adic integer}
\end{equation}
where the $c_{n}$s are elements of the set $\left\{ 0,1,\ldots,p-1\right\} $. 
\end{defn}
\vphantom{}

We can think of $p$-adic integers as ``power series in $p$''.
Note that every non-negative integer is automatically a $p$-adic
integer, seeing as every non-negative integer $x$ can be uniquely
written as a \emph{finite }sum $\mathfrak{z}=\sum_{n=0}^{N}c_{n}p^{n}$,
where the $c_{n}$s are the\index{$p$-adic!digits} $p$\textbf{-ary/adic
digits }of $\mathfrak{z}$. The $p$-adic representation of $\mathfrak{z}\in\mathbb{Z}_{p}$
is the expression: 
\begin{equation}
\mathfrak{z}=\centerdot_{p}c_{0}c_{1}c_{2}\ldots\label{eq:Definition of the p-adic digit representation}
\end{equation}
where the subscript $p$ is there to remind us that we are in base
$p$; this subscript can (and will) be dropped when there is no confusion
as to the value of $p$.

$\mathbb{Z}_{p}$ becomes a commutative, unital ring when equipped
with the usual addition and multiplication operations, albeit with
the caveat that 
\begin{equation}
cp^{n}=\left[c\right]_{p}p^{n}+\left(\frac{c-\left[c\right]_{p}}{p}\right)p^{n+1}
\end{equation}
for all $n,c\in\mathbb{N}_{0}$. 
\begin{example}
As an example, for $a\in\left\{ 0,p-1\right\} $: 
\begin{equation}
\centerdot_{p}00\left(p+a\right)=\left(p+a\right)p^{2}=ap^{2}+1p^{3}=\centerdot_{p}00a1
\end{equation}
where, in $\centerdot_{p}00\left(p+a\right)$, $p+a$ is the value
of the third $p$-adic digit, and where $\left[c\right]_{p}\in\left\{ 0,\ldots,p-1\right\} $
denotes the residue class of $c$ mod $p$. That is to say, like an
odometer, \emph{we carry over to the next $p$-adic digit's place
whenever a digit reaches $p$}.

Thus, in the $2$-adics: 
\begin{eqnarray*}
3 & = & 1\times2^{0}+1\times2^{1}=\centerdot_{2}11\\
5 & = & 1\times2^{0}+0\times2^{1}+1\times2^{2}=\centerdot_{2}101
\end{eqnarray*}
When we add these numbers, we ``carry over the $2$'': 
\[
3+5=\centerdot_{2}11+\centerdot_{2}101=\centerdot_{2}211=\centerdot_{2}021=\centerdot_{2}002=\centerdot_{2}0001=1\times2^{3}=8
\]
Multiplication is done similarly. 
\end{example}
\vphantom{}

Using this arithmetic operation, we can write negative integers $p$-adically.
\begin{example}
The $p$-adic number $\mathfrak{y}$ whose every $p$-adic digit is
equal to $p-1$: 
\[
\mathfrak{y}=\centerdot_{p}\left(p-1\right)\left(p-1\right)\left(p-1\right)\left(p-1\right)\ldots=\sum_{n=0}^{\infty}\left(p-1\right)p^{n}
\]
satisfies: 
\begin{eqnarray*}
1+\mathfrak{y} & = & \centerdot_{p}\left(p\right)\left(p-1\right)\left(p-1\right)\left(p-1\right)\ldots\\
 & = & \centerdot_{p}0\left(p-1+1\right)\left(p-1\right)\left(p-1\right)\ldots\\
 & = & \centerdot_{p}0\left(p\right)\left(p-1\right)\left(p-1\right)\ldots\\
 & = & \centerdot_{p}00\left(p-1+1\right)\left(p-1\right)\ldots\\
 & = & \centerdot_{p}00\left(p-1+1\right)\left(p-1\right)\ldots\\
 & = & \centerdot_{p}000\left(p-1+1\right)\ldots\\
 & \vdots\\
 & = & \centerdot_{p}000000\ldots\\
 & = & 0
\end{eqnarray*}
and thus, $\mathfrak{y}=-1$ in $\mathbb{Z}_{p}$.
\end{example}
\vphantom{}

The beauty of the power series conception of these numbers is that
it makes such formulae explicit: 
\begin{eqnarray*}
1+y & = & 1+\sum_{n=0}^{\infty}\left(p-1\right)p^{n}\\
 & = & 1+\sum_{n=0}^{\infty}p^{n+1}-\sum_{n=0}^{\infty}p^{n}\\
 & = & p^{0}+\sum_{n=1}^{\infty}p^{n}-\sum_{n=0}^{\infty}p^{n}\\
 & = & \sum_{n=0}^{\infty}p^{n}-\sum_{n=0}^{\infty}p^{n}\\
 & = & 0
\end{eqnarray*}
Consequently, just as when given a power series $\sum_{n=0}^{\infty}a_{n}z^{n}$,
we can compute its multiplicative inverse $\sum_{n=0}^{\infty}b_{n}z^{n}$
recursively by way of the equations: 
\begin{eqnarray*}
1 & = & \left(\sum_{m=0}^{\infty}a_{m}z^{m}\right)\left(\sum_{n=0}^{\infty}b_{n}z^{n}\right)=\sum_{k=0}^{\infty}\left(\sum_{n=0}^{k}a_{n}b_{k-n}\right)z^{n}\\
 & \Updownarrow\\
1 & = & a_{0}b_{0}\\
0 & = & a_{0}b_{1}+a_{1}b_{0}\\
0 & = & a_{0}b_{2}+a_{1}b_{1}+a_{2}b_{0}\\
 & \vdots
\end{eqnarray*}
we can use the same formula to compute the multiplicative inverse
of a given $p$-adic integer\textemdash \emph{assuming} it exists.
Working through the above equations, it can be seen that the $p$-adic
integer $\mathfrak{z}=\sum_{n=0}^{\infty}c_{n}p^{n}$ has a multiplicative
inverse if and only if $c_{0}\in\left\{ 0,\ldots,p-1\right\} $ is
multiplicatively invertible mod $p$ (i.e., $c_{0}$ is co-prime to
$p$). This is one of several reasons why we prefer to study $p$-adic
integers for prime $p$: every non-zero residue class mod $p$ is
multiplicatively invertible modulo $p$.
\begin{defn}
Just as we can go from the ring of power series to the field of Laurent
series, we can pass from the ring of $p$-adic integers $\mathbb{Z}_{p}$
to the field of \textbf{$p$-adic (rational) numbers}\index{$p$-adic!rational numbers}
\nomenclature{$\mathbb{Q}_{p}$}{the set of $p$-adic rational numbers \nopageref}$\mathbb{Q}_{p}$\textbf{
}by considering ``Laurent series'' in $p$. Every $\mathfrak{z}\in\mathbb{Q}_{p}$
has a unique representation as: 
\begin{equation}
\mathfrak{z}=\sum_{n=n_{0}}^{\infty}c_{n}p^{n}\label{eq:Laurent series representation of a p-adic rational number}
\end{equation}
for some $n_{0}\in\mathbb{Z}$. We write the ``fractional part''
of $\mathfrak{z}$ (the terms with negative power of $p$) to the
right of the $p$-adic point $\centerdot_{p}$. 
\end{defn}
\begin{example}
Thus, the $3$-adic number: 
\begin{equation}
\frac{1}{3^{2}}+\frac{2}{3^{1}}+0\times3^{0}+1\times3^{1}+1\times3^{2}+1\times3^{3}+\ldots
\end{equation}
would be written as: 
\begin{equation}
12\centerdot_{3}0\overline{1}=12\centerdot_{3}0111\ldots
\end{equation}
where, as usual, the over-bar indicates that we keep writing $1$
over and over again forever. 
\end{example}
\begin{defn}
We can equip $\mathbb{Q}_{p}$ with the structure of a metric space
by defining the \index{$p$-adic!valuation}\textbf{$p$-adic valuation}
\textbf{$v_{p}$ }and $p$\textbf{-absolute value}\footnote{Although this text and most texts on non-archimedean analysis by native-English-speaking
authors make a clear distinction between the $p$-adic absolute value
and the $p$-adic valuation, this distinction is not universal adhered
to. A non-negligible proportion of the literature (particularly when
authors of Russian extraction are involved) use the word ``valuation''
interchangeably, to refer to either $v_{p}$ or $\left|\cdot\right|_{p}$,
depending on the context.}\textbf{ }$\left|\cdot\right|_{p}$: 
\begin{equation}
v_{p}\left(\sum_{n=n_{0}}^{\infty}c_{n}p^{n}\right)\overset{\textrm{def}}{=}n_{0}\label{eq:Definition of the p-adic valuation}
\end{equation}
\begin{equation}
\left|\sum_{n=n_{0}}^{\infty}c_{n}p^{n}\right|_{p}\overset{\textrm{def}}{=}p^{-n_{0}}\label{eq:Definition of the p-adic absolute value}
\end{equation}
that is, $\left|\mathfrak{z}\right|_{p}=p^{-v_{p}\left(\mathfrak{z}\right)}$.
We also adopt the convention that $v_{p}\left(0\right)\overset{\textrm{def}}{=}\infty$.
The \textbf{$p$-adic metric}\index{p-adic@\textbf{$p$}-adic!metric}\textbf{
}is then the distance formula defined by the map: 
\begin{equation}
\left(\mathfrak{z},\mathfrak{y}\right)\in\mathbb{Z}_{p}\times\mathbb{Z}_{p}\mapsto\left|\mathfrak{z}-\mathfrak{y}\right|_{p}\label{eq:Definition of the p-adic metric}
\end{equation}
This metric (also known as an \textbf{ultrametric})\textbf{ }is said
to \textbf{non-archimedean} because, in addition to the triangle inequality
we all know and love, it also satisfies the \textbf{strong triangle
inequality} (also called the \textbf{ultrametric inequality}): 
\begin{equation}
\left|\mathfrak{z}-\mathfrak{y}\right|_{p}\leq\max\left\{ \left|\mathfrak{z}\right|_{p},\left|\mathfrak{y}\right|_{p}\right\} \label{eq:the Strong Triangle Inequality}
\end{equation}
\emph{Crucially}\textemdash and I \emph{cannot} emphasize this enough\textemdash (\ref{eq:the Strong Triangle Inequality})
holds \emph{with equality} whenever $\left|\mathfrak{z}\right|_{p}\neq\left|\mathfrak{y}\right|_{p}$.
The equality of (\ref{eq:the Strong Triangle Inequality}) when $\left|\mathfrak{z}\right|_{p}\neq\left|\mathfrak{y}\right|_{p}$
is one of the most subtly powerful tricks in ultrametric analysis,
especially when we are trying to contradict an assumed upper bound.
Indeed, this very method is at the heart of my proof of th $\left(p,q\right)$-adic
Wiener Tauberian Theorem.
\end{defn}
\begin{rem}
Although it might seem unintuitive that a \emph{large} power of $p$
should have a very small $p$-adic absolute value, this viewpoint
becomes extremely natural when viewed through Hensel's original conception
of $p$-adic integers as number theoretic analogues of power series
\cite{Hensel-s original article,Gouvea's introudction to p-adic numbers book,Journey throughout the history of p-adic numbers}.
It is a well-known and fundamental fact of complex analysis that a
function $f:U\rightarrow\mathbb{C}$ holomorphic on an open, connected,
non-empty set $U\subseteq\mathbb{C}$ which possesses a zero of infinite
degree in $U$ must be identically zero on $U$. In algebraic terms,
if $f$ has a zero at $z_{0}$ and $f$ is not identically zero, then
$f\left(z\right)/\left(z-z_{0}\right)^{n}$ can only be divided by
the binomial $z-z_{0}$ finitely many times before we obtain a function
which is non-zero at $z_{0}$. As such, a function holomorphic on
an open neighborhood of $z_{0}$ which remains holomorphic after being
divided by $z-z_{0}$ arbitrarily many times is necessarily the constant
function $0$. Hensel's insight was to apply this same reasoning to
numbers. Indeed, the uniform convergence to $0$ of the function sequence
$\left\{ z^{n}\right\} _{n\geq0}$ on any compact subset of the open
unit disk $\mathbb{D}\subset\mathbb{C}$ is spiritually equivalent
to the convergence of the sequence $\left\{ p^{n}\right\} _{n\geq0}$
to $0$ in the $p$-adics: just as the only holomorphic function on
$\mathbb{D}$ divisible by arbitrarily large powers of $z$ is the
zero function, the only $p$-adic integer divisible by arbitrarily
high powers of $p$ is the integer $0$.

Indeed, in the context of power series and Laurent series\textemdash say,
a series $\sum_{n=n_{0}}^{\infty}c_{n}\left(z-z_{0}\right)^{n}$ about
some $z_{0}\in\mathbb{C}$\textemdash the $p$-adic valuation corresponds
to the \textbf{zero degree }of the series at $z_{0}$. If $n_{0}$
is negative, then the function represented by that power series has
a pole of order $-n_{0}$ at $z_{0}$; if $n_{0}$ is positive, then
the function represented by that power series has a \emph{zero }of
order $n_{0}$ at $z_{0}$. Thus, for all $n\in\mathbb{Z}$: 
\begin{equation}
\left|p^{n}\right|_{p}=p^{-n}
\end{equation}
This is especially important, seeing as the $p$-adic absolute value
is a multiplicative group homomorphism from $\mathbb{Z}_{p}$ to $\mathbb{R}^{+}$:
\begin{equation}
\left|xy\right|_{p}=\left|x\right|_{p}\left|y\right|_{p}\label{eq:Multiplicativity of p-adic absolute value}
\end{equation}
as such, the level sets of $\mathbb{Z}_{p}$ are: 
\begin{equation}
\left\{ x\in\mathbb{Z}_{p}:\left|x\right|_{p}=p^{-n}\right\} =p^{n}\mathbb{Z}_{p}\overset{\textrm{def}}{=}\left\{ p^{n}y:y\in\mathbb{Z}_{p}\right\} \label{eq:Definition of p^n times Z_p}
\end{equation}
for all $n\in\mathbb{N}_{0}$.
\end{rem}
\vphantom{}

Returning to the matter at hand, by using absolute values, we see
that the $p$-adic integers are precisely those elements of $\mathbb{Q}_{p}$
with $p$-adic absolute value $\leq1$. Consequentially, the metric
space obtained by equipping $\mathbb{Z}_{p}$ with the $p$-adic metric
is \emph{compact}. $\mathbb{Q}_{p}$, meanwhile, is locally compact,
just like $\mathbb{R}$ and $\mathbb{C}$. 
\begin{defn}
We write \nomenclature{$\mathbb{Z}_{p}^{\times}$}{the group of multiplicatively invertible $p$-adic integers \nopageref}$\mathbb{Z}_{p}^{\times}$
to denote the set of all units of $\mathbb{Z}_{p}$\textemdash that
is, elements of $\mathbb{Z}_{p}$ whose reciprocals are contained
in $\mathbb{Z}_{p}$. This is an abelian group under multiplication,
with $1$ as its identity element. Note, also, that:
\begin{equation}
\mathbb{Z}_{p}^{\times}=\left\{ \mathfrak{z}\in\mathbb{Z}_{p}:\left|\mathfrak{z}\right|_{p}=1\right\} 
\end{equation}
\end{defn}
\begin{defn}
\textbf{Congruences }are extremely important when working with $p$-adic
integers; this is an intrinsic feature of the ``projective limit''
more algebraic texts frequently use to define the $p$-adic integers.
As we saw, every $p$-adic integer $\mathfrak{z}$ can be written
uniquely as: 
\begin{equation}
\mathfrak{z}=\sum_{n=0}^{\infty}c_{n}p^{n}\label{eq:p-adic series representation of a p-adic integer}
\end{equation}
for constants $\left\{ c_{n}\right\} _{n\geq0}\subseteq\left\{ 0,\ldots,p-1\right\} $.
The series representation (\ref{eq:p-adic series representation of a p-adic integer})
of a $p$-adic integer $\mathfrak{z}$is sometimes called the \textbf{Hensel
series}\index{Hensel series}\index{series!Hensel}series\textbf{
}or \textbf{Henselian series }of $\mathfrak{z}$; Amice uses this
terminology, for example \cite{Amice}. With this representation,
given any integer $m\geq0$, we can then define \nomenclature{$\left[\mathfrak{z}\right]_{p^{n}}$}{Projection of a $p$-adic integer modulo $p^n$ nomnorefpage}
$\left[\mathfrak{z}\right]_{p^{m}}$, the \textbf{projection of $\mathfrak{z}$
mod $p^{m}$} like so: 
\begin{equation}
\left[\mathfrak{z}\right]_{p^{m}}\overset{\textrm{def}}{=}\sum_{n=0}^{m-1}c_{n}p^{n}\label{eq:Definition of the projection of z mod p to the m}
\end{equation}
where the right-hand side is defined to be $0$ whenever $m=0$. Since
$\left[\mathfrak{z}\right]_{p^{m}}$ is a finite sum of integers,
it itself is an integer. Given $\mathfrak{z},\mathfrak{y}\in\mathbb{Z}_{p}$,
we write $\mathfrak{z}\overset{p^{m}}{\equiv}\mathfrak{y}$ if and
only if $\left[\mathfrak{z}\right]_{p^{m}}=\left[\mathfrak{y}\right]_{p^{m}}$.
Moreover, there is an equivalence between congruences and absolute
values: 
\begin{equation}
\mathfrak{z}\overset{p^{m}}{\equiv}\mathfrak{y}\Leftrightarrow\left|\mathfrak{z}-\mathfrak{y}\right|_{p}\leq p^{-m}
\end{equation}

There is a very useful notation for denoting subsets (really, ``open
neighborhoods'') of $p$-adic integers. Given $\mathfrak{z}\in\mathbb{Z}_{p}$
and $n\in\mathbb{N}_{0}$, we write: 
\begin{equation}
\mathfrak{z}+p^{n}\mathbb{Z}_{p}\overset{\textrm{def}}{=}\left\{ \mathfrak{y}\in\mathbb{Z}_{p}:\mathfrak{y}\overset{p^{n}}{\equiv}\mathfrak{z}\right\} +\left\{ \mathfrak{y}\in\mathbb{Z}_{p}:\left|\mathfrak{z}-\mathfrak{y}\right|_{p}\leq p^{-n}\right\} \label{eq:Definition of co-set notation for p-adic neighborhoods}
\end{equation}
In particular, note that: 
\begin{align*}
\mathfrak{z}+p^{n}\mathbb{Z}_{p} & =\mathfrak{y}+p^{n}\mathbb{Z}_{p}\\
 & \Updownarrow\\
\mathfrak{z} & \overset{p^{n}}{\equiv}\mathfrak{y}
\end{align*}
Additionally\textemdash as is crucial for performing integration over
the $p$-adics\textemdash for any $n\geq0$, we can partition $\mathbb{Z}_{p}$
like so: 
\begin{equation}
\mathbb{Z}_{p}=\bigcup_{k=0}^{p^{n}-1}\left(k+p^{n}\mathbb{Z}_{p}\right)
\end{equation}
where the $k+p^{n}\mathbb{Z}_{p}$s are pair-wise disjoint with respect
to $k$. 
\end{defn}
\vphantom{}

Finally, it is worth mentioning that $\mathbb{N}_{0}$ is dense in
$\mathbb{Z}_{p}$, and that $\mathbb{Q}$ is dense in $\mathbb{Q}_{p}$.
This makes both $\mathbb{Z}_{p}$ and $\mathbb{Q}_{p}$ into separable
topological spaces. Moreover, as topological spaces, they\textemdash and
any field extension thereof\textemdash are \emph{totally disconnected}.
This has profound implications for $p$-adic analysis, and for ultrametric
analysis in general.

\subsection{\label{subsec:1.3.2. Ultrametrics-and-Absolute}An Introduction to
Ultrametric Analysis}

While the $p$-adic numbers are surely the most well-known non-archimedean
spaces, they are far from the only ones. The study of function theory,
calculus, and the like on generic non-archimedean spaces is sometimes
called \textbf{Ultrametric analysis}\index{ultrametric!analysis}.
As will be addressed at length in the historical essay of Subsection
\ref{subsec:3.1.1 Some-Historical-and}, the sub-disciplines that
go by the names of non-archimedean analysis, $p$-adic analysis, ultrametric
analysis, and the like are sufficiently divers that it is worth being
aware of the undercurrent of common terminology.
\begin{defn}
Let $X$ be a set, and let $d:X\times X\rightarrow\left[0,\infty\right)$
be a metric; that is, a function satisfying:

\vphantom{}

I. $d\left(x,y\right)\geq0$ $\forall x,y\in X$, with equality if
and only if $x=y$;

\vphantom{}

II. $d\left(x,y\right)=d\left(y,x\right)$ $\forall x,y\in X$;

\vphantom{}

III. $d\left(x,y\right)\leq d\left(x,z\right)+d\left(y,z\right)$
$\forall x,y,z\in X$.

\vphantom{}

We say $d$ is a \textbf{non-archimedean} \textbf{metric }or \textbf{ultrametric}
whenever it satisfies the \textbf{Strong Triangle Inequality }(a.k.a.
\textbf{Ultrametric Inequality}): 
\begin{equation}
d\left(x,y\right)\leq\max\left\{ d\left(x,z\right),d\left(y,z\right)\right\} ,\forall x,y,z\in X\label{eq:Generic Strong Triangle Inequality}
\end{equation}
Like with the specific case of the $p$-adic ultrametric inequality,
the general ultrametric inequality holds \emph{with equality} whenever
$d\left(x,z\right)\neq d\left(y,z\right)$. An \index{ultrametric!space}\textbf{ultrametric
space }is a pair $\left(X,d\right)$, where $X$ is a set and $d$
is an ultrametric on $X$. 
\end{defn}
\vphantom{}

The most important ultrametrics are those that arise from \textbf{absolute
values} on abelian groups, particularly fields.
\begin{defn}
Let $K$ be an abelian group, written additively, and with $0$ as
its identity element. An \textbf{absolute value}\index{absolute value (on a field)}\textbf{
}on $K$ is a function $\left|\cdot\right|_{K}:K\rightarrow\left[0,\infty\right)$
satisfying the following properties for all $x,y\in K$:

\vphantom{}

I. $\left|0\right|_{K}=0$;

\vphantom{}

II. $\left|x+y\right|_{K}\leq\left|x\right|_{K}+\left|y\right|_{K}$;

\vphantom{}

III. $\left|x\right|_{K}=0$ if and only if $x=0$;

\vphantom{}

IV. If $K$ is a ring, we also require $\left|x\cdot y\right|_{K}=\left|x\right|_{K}\cdot\left|y\right|_{K}$.

\vphantom{}

Finally, we say that $\left|\cdot\right|_{K}$ is a \textbf{non-archimedean
absolute value} if, in addition to the above, $\left|\cdot\right|_{K}$
satisfies the \textbf{Ultrametric Inequality}\index{ultrametric!inequality}\index{triangle inequality!strong}:

\vphantom{}

V. $\left|x+y\right|_{K}\leq\max\left\{ \left|x\right|_{K},\left|y\right|_{K}\right\} $
(with equality whenever $\left|x\right|_{K}\neq\left|y\right|_{K}$),.

\vphantom{}

If (V) is not satisfied, we call $\left|\cdot\right|_{K}$ an \textbf{archimedean
absolute value}.

Finally, note that if $K$ is a field, any absolute value $\left|\cdot\right|_{K}$
on a field induces a metric $d$ on $K$ by way of the formula $d\left(x,y\right)=\left|x-y\right|_{K}$.
We call the pair $\left(K,\left|\cdot\right|_{K}\right)$ a \textbf{valued
group }(resp. \textbf{valued ring}; resp. \index{valued field}\textbf{valued
field}) whenever $K$ is an abelian group (resp. ring\footnote{The ring need not be commutative.};
resp., field), $\left|\cdot\right|_{K}$ is an absolute value on and
say it is \textbf{archimedean }whenever $\left|\cdot\right|_{K}$
is archimedean, and say that it is \textbf{non-archimedean }whenever
$\left|\cdot\right|_{K}$. 

Let $K$ be a non-archimedean valued ring.

\vphantom{}

I. Following Schikhof \cite{Schikhof Banach Space Paper}, let $B_{K}\overset{\textrm{def}}{=}\left\{ x\in K:\left|x\right|_{K}\leq1\right\} $
and let $B_{K}^{-}\overset{\textrm{def}}{=}\left\{ x\in K:\left|x\right|_{K}<1\right\} $.
Both $B_{K}$ and $B_{K}^{-}$ are rings under the addition and multiplication
operations of $K$, with $B_{K}^{-}$ being an ideal\footnote{In fact, $B_{K}^{-}$ is a maximal ideal in $B_{K}$, and is the unique
non-zero prime ideal of $B_{K}$. If $K$ is a ring, $K$ is then
called a \textbf{local ring}; if $K$ is a field, it is then called
a \textbf{local field}, in which case $B_{K}$ is the ring of $K$-integers,
and $K$ is the field of fractions of $B_{K}$.} in $B_{K}$. The ring $B_{K}/B_{K}^{-}$ obtained by quotienting
$B_{K}$ out by $B_{K}^{-}$ is called the \textbf{residue field }/
\textbf{residue class field }of $K$.

We say $K$ is\textbf{ $p$-adic} (where $p$ is a prime) when the
residue field of $K$ has characteristic $p$. In the case $K$ is
$p$-adic, in an abuse of notation, we will write $\left|\cdot\right|_{p}$
to denote the absolute value on $K$. Also, note that if $K$ is a
field, $B_{K}$ is then equal to $\mathcal{O}_{K}$, the \textbf{ring
of integers }of $K$.

\vphantom{}

II. The set: 
\begin{equation}
\left|K\backslash\left\{ 0\right\} \right|_{K}\overset{\textrm{def}}{=}\left\{ \left|x\right|_{K}:x\in K\backslash\left\{ 0\right\} \right\} \label{eq:Definition of the value group of a field}
\end{equation}
is called the \index{value group}\textbf{value group}\footnote{Not to be confused with a \emph{valued }group.}\textbf{
}of $K$. This group is said to be \textbf{dense }if it is dense in
the interval $\left(0,\infty\right)\subset\mathbb{R}$ in the standard
topology of $\mathbb{R}$, and is said to be \textbf{discrete }if
it is not dense. 
\end{defn}
\vphantom{}

Both the residue field and the value group of $K$ play a crucial
role in some of the fundamental properties of analysis on $K$ (or
on normed vector spaces over $K$). For example, they completely determine
whether or not $K$ is locally compact. 
\begin{thm}
A non-archimedean valued field $K$ is locally compact if and only
if its residue field is finite and its value group is discrete\footnote{\textbf{Theorem 12.2} from \cite{Ultrametric Calculus}.}. 
\end{thm}
\vphantom{}

The term ``non-archimedean'' comes from the failure of the \textbf{archimedean
property} of classical analysis, that being the intuitive notion that
\emph{lengths add up}. Non-archimedean spaces are precisely those
metric spaces where lengths \emph{need not} add up. In particular,
we have the following delightful result: 
\begin{thm}
Let $\left(K,\left|\cdot\right|_{K}\right)$ be a valued field, and
let $1$ denote the multiplicative identity element of $K$. Then,
$\left(K,\left|\cdot\right|_{K}\right)$ is non-archimedean if and
only if: 
\begin{equation}
\left|1+1\right|_{K}\leq\left|1\right|_{K}
\end{equation}
\end{thm}
Proof: Exercise.

Q.E.D.

\vphantom{}

The following list, adapted from Robert's book \cite{Robert's Book},
gives an excellent summary of the most important features of ultrametric
analysis. 
\begin{fact}[\textbf{The Basic Principles of Ultrametric Analysis}]
\label{fact:Principles of Ultrametric Analysis}Let $\left(K,\left|\cdot\right|_{K}\right)$
be a non-archimedean valued group with additive identity element $0$.
Then:

\vphantom{}

I. \textbf{\emph{The Strongest Wins}}\emph{:}

\vphantom{} 
\begin{equation}
\left|x\right|_{K}>\left|y\right|_{K}\Rightarrow\left|x+y\right|_{K}=\left|x\right|_{K}\label{eq:The Strongest Wins}
\end{equation}

II. \textbf{\emph{Equilibrium}}\emph{:} All triangles are isosceles
(or equilateral): 
\begin{equation}
a+b+c=0\textrm{ \& }\left|c\right|_{K}<\left|b\right|_{K}\Rightarrow\left|a\right|_{K}=\left|b\right|_{K}\label{eq:Equilibrium}
\end{equation}

\vphantom{}

III. \textbf{\emph{Competition}}\emph{:} If: 
\[
a_{1}+\cdots+a_{n}=0
\]
then there are distinct $i,j$ so that $\left|a_{i}\right|_{K}=\left|a_{j}\right|_{K}=\max_{k}\left|a_{k}\right|_{K}$.

\vphantom{}

IV. \textbf{\emph{The Freshman's Dream}}\emph{:} If the metric space
$\left(K,\left|\cdot\right|_{K}\right)$ is complete, a series $\sum_{n=0}^{\infty}a_{n}$
converges in $K$ if and only if $\left|a_{n}\right|_{K}\rightarrow0$
as $n\rightarrow\infty$. Consequently:

\vphantom{}

i. (The infinite version of (III) holds) $\sum_{n=0}^{\infty}a_{n}=0$
implies there are distinct $i,j$ so that $\left|a_{i}\right|_{K}=\left|a_{j}\right|_{K}=\max_{k}\left|a_{k}\right|_{K}$.

\vphantom{}

ii. The convergence of $\sum_{n=0}^{\infty}\left|a_{n}\right|_{K}$
in $\mathbb{R}$ implies the convergence of $\sum_{n=0}^{\infty}a_{n}$
in $K$, but the converse is not true: $\sum_{n=0}^{\infty}a_{n}$
converging in $K$ need not imply $\sum_{n=0}^{\infty}\left|a_{n}\right|_{K}$
converges in $\mathbb{R}$.

\vphantom{}

iii. For any $K$-convergent series $\sum_{n=0}^{\infty}a_{n}$: 
\begin{equation}
\left|\sum_{n=0}^{\infty}a_{n}\right|_{K}\leq\sup_{n\geq0}\left|a_{n}\right|_{K}=\max_{n\geq0}\left|a_{n}\right|_{K}\label{eq:Series Estimate}
\end{equation}
The right-most equality indicates that there is going to be an $n$
for which the absolute value is maximized.

\vphantom{}

V. \textbf{\emph{The Sophomore's Dream}}\emph{:} A sequence $\left\{ a_{n}\right\} _{n\geq0}$
is Cauchy if and only if $\left|a_{n+1}-a_{n}\right|_{K}\rightarrow0$
as $n\rightarrow\infty$.

\vphantom{}

\vphantom{}

VI. \textbf{\emph{Stationarity of the absolute value}}\emph{:} If
$a_{n}$ converges to $a$ in $K$, and if $a\neq0$, then there is
an $N$ so that $\left|a_{n}\right|_{K}=\left|a\right|_{K}$ for all
$n\geq N$. 
\end{fact}
\vphantom{}

We also have the following results regarding infinite series: 
\begin{prop}[\index{series!re-arrangement}\textbf{Series re-arrangement}\footnote{Given on page 74 of \cite{Robert's Book}.}]
\label{prop:series re-arrangement}Let $\left(K,\left|\cdot\right|_{K}\right)$
be a complete non-archimedean valued group, and let $\left\{ a_{n}\right\} _{n\geq0}$
be a sequence in $K$ which tends to $0$ in $K$, so that $\sum_{n=0}^{\infty}a_{n}$
converges in $K$ to $s$. Then, no matter how the terms of the sum
are grouped or rearranged, the resultant series will still converge
in $K$ to $s$. Specifically:

\vphantom{}

I. For any bijection $\sigma:\mathbb{N}_{0}\rightarrow\mathbb{N}_{0}$,
$\sum_{n=0}^{\infty}a_{\sigma\left(n\right)}$ converges in $K$ to
$s$.

\vphantom{}

II. For any partition of $\mathbb{N}_{0}$ into sets $I_{1},I_{2},\ldots$,
the series: 
\begin{equation}
\sum_{k}\left(\sum_{n\in I_{k}}a_{n}\right)
\end{equation}
converges in $K$ to $s$. 
\end{prop}
\begin{prop}[\index{series!interchange}\textbf{Series interchange}\footnote{Given on page 76 of \cite{Robert's Book}.}]
\label{prop:series interchange} Let $\left(K,\left|\cdot\right|_{K}\right)$
be a complete non-archimedean valued group, and let $\left\{ a_{m,n}\right\} _{m,n\geq0}$
be a double-indexed sequence in $K$. If, for any $\epsilon>0$, there
are only finitely many pairs $\left(m,n\right)$ so that $\left|a_{m,n}\right|_{K}>\epsilon$,
then the double sum: 
\begin{equation}
\sum_{\left(m,n\right)\in\mathbb{N}_{0}^{2}}a_{m,n}
\end{equation}
converges in $K$, and, moreover: 
\begin{equation}
\sum_{\left(m,n\right)\in\mathbb{N}_{0}^{2}}a_{m,n}=\sum_{m=0}^{\infty}\left(\sum_{n=0}^{\infty}a_{m,n}\right)=\sum_{n=0}^{\infty}\left(\sum_{m=0}^{\infty}a_{m,n}\right)
\end{equation}
where all equalities are in $K$. 
\end{prop}
\vphantom{}

The topological properties of ultrametric spaces are drastically different
from Euclidean spaces (ex: $\mathbb{R}^{n}$).
\begin{defn}
Let $\left(X,d\right)$ be an ultrametric space.

\vphantom{}

I. A closed ball\index{ball} in $X$ of radius $r$ (where $r$ is
a positive real number) centered at $x\in X$, written $B\left(x,r\right)$,
is the set: 
\begin{equation}
B\left(x,r\right)\overset{\textrm{def}}{=}\left\{ y\in X:d\left(x,y\right)\leq r\right\} \label{eq:Definition of a closed ball}
\end{equation}
Open balls are obtained by making the inequality $\leq r$ strict
($<r$). However, as we will see momentarily, this doesn't actually
amount to much of a distinction.

\vphantom{}

II. Given any non-empty subset $Y\subseteq X$, the \textbf{diameter
}of $Y$, denoted $d\left(Y\right)$, is defined by: 
\begin{equation}
d\left(Y\right)\overset{\textrm{def}}{=}\sup\left\{ d\left(a,b\right):a,b\in Y\right\} \label{eq:Definition of the diameter of an ultrametric set}
\end{equation}
\end{defn}
\begin{rem}
For any ball $B\left(x,r\right)\subseteq X$: 
\begin{equation}
d\left(B\right)=\inf\left\{ r^{\prime}>0:B\left(x,r^{\prime}\right)=B\left(x,r\right)\right\} \label{eq:Diameter of an ultrametric ball in terms of its radius}
\end{equation}
\end{rem}
\begin{prop}
In\footnote{\textbf{Propositions 18.4} and \textbf{18.5} from \cite{Ultrametric Calculus}.}
an ultrametric space\emph{ \index{ultrametric!balls}}$\left(X,d\right)$:

\vphantom{}

I. All open balls are closed, and all closed balls are open\footnote{However, it is not necessarily the case that $\left\{ y\in X:d\left(x,y\right)\leq r\right\} $
and $\left\{ y\in X:d\left(x,y\right)<r\right\} $ will be the same
set. For example, if $X=\mathbb{Z}_{3}$, then $\left\{ \mathfrak{y}\in\mathbb{Z}_{3}:\left|\mathfrak{z}-\mathfrak{y}\right|_{3}\leq\frac{1}{2}\right\} =\left\{ \mathfrak{y}\in\mathbb{Z}_{3}:\left|\mathfrak{z}-\mathfrak{y}\right|_{3}\leq\frac{1}{3}\right\} $,
because the value group of $\mathbb{Z}_{3}$ is $\left\{ 3^{-n}:n\in\mathbb{N}_{0}\right\} $.
On the other hand, for the same reason, $\left\{ \mathfrak{y}\in\mathbb{Z}_{3}:\left|\mathfrak{z}-\mathfrak{y}\right|_{3}<\frac{1}{3}\right\} =\left\{ \mathfrak{y}\in\mathbb{Z}_{3}:\left|\mathfrak{z}-\mathfrak{y}\right|_{3}\leq\frac{1}{9}\right\} $.}.

\vphantom{}

II. All points inside a ball are at the center of the ball; that is,
given $x,y\in K$ and $r>0$ such that $d\left(x,y\right)\leq r$,
we have that $B\left(x,r\right)=B\left(y,r\right)$.

\vphantom{}

III. Given two balls $B_{1},B_{2}\subseteq X$, either $B_{1}\cap B_{2}=\varnothing$
or either $B_{1}\subseteq B_{2}$ or $B_{2}\subseteq B_{1}$.

\vphantom{}

IV. Given any ball $B\left(x,r\right)$, there are infinitely many
real numbers $r^{\prime}$ for which $B\left(x,r^{\prime}\right)=B\left(x,r\right)$. 
\end{prop}
\begin{fact}
The \textbf{Heine-Borel Property}\footnote{A set is compact if and only if it is both closed and bounded.}\index{Heine-Borel property}
also holds in $\mathbb{Q}_{p}$, as well as in any finite-dimensional,
locally compact field extension thereof. 
\end{fact}
\vphantom{}

We also have the following result regarding open covers in an arbitrary
ultrametric space: 
\begin{thm}[\textbf{Open Set Decomposition Theorem}\footnote{\textbf{Theorem 18.6} on page 48 of \cite{Ultrametric Calculus}.}]
Let $X$ be an ultrametric space. Then, every non-empty open set
$U\subseteq X$ can be written as the union of countably many pair-wise
disjoint balls. 
\end{thm}
\begin{cor}
Let $X$ be a compact ultrametric space (such as $\mathbb{Z}_{p}$)
can be written as the union of finitely many pair-wise disjoint balls. 
\end{cor}
Proof: Let $U\subseteq\mathbb{Z}_{p}$ be non-empty and clopen. Since
$U$ is clopen, it is closed, and since $U$ is in $\mathbb{Z}_{p}$,
it is bounded in $\mathbb{Q}_{p}$. Since $\mathbb{Q}_{p}$ possess
the Heine-Borel property, the closedness and boundedness of $U$ then
force $U$ to be compact.

Now, by the theorem, since $U$ is clopen and non-empty, it is open
and non-empty, and as such, $U$ can be written as a union of countably
many disjoint balls. Since the balls are clopen sets, this collection
of balls forms an open cover for $U$. By the compactness of $U$,
this open cover must contain a finite sub-cover, which shows that
$U$ can, in fact, be written as the union of finitely many clopen
balls.

Q.E.D.

\vphantom{}

By this point, the reader might have noticed that this exposition
of basic ultrametric analysis has yet to say a word about functions.
This is intentional. Non-archimedean function theory is one of the
central overarching concerns of this dissertation. Nevertheless, there
are two important pieces of non-archimedean function theory I can
introduce here and now.
\begin{defn}[\textbf{Locally constant functions}\footnote{Taken from \cite{Ultrametric Calculus}.}]
Let $X$ be an ultrametric space, and let $K$ be any valued field.
We say a function $f:X\rightarrow K$ is \textbf{locally constant
}if, for each $x\in X$ there is an open neighborhood $U\subseteq X$
containing $x$ so that $f$ is constant on $U\cap X$. 
\end{defn}
\begin{example}
The most important example we will be working with are functions $\mathbb{Z}_{p}\rightarrow K$
of the form: 
\begin{equation}
\sum_{n=0}^{p^{N}-1}a_{n}\left[\mathfrak{z}\overset{p^{N}}{\equiv}n\right]
\end{equation}
where the $a_{n}$s are scalars in $K$. Such a function's value at
any given $\mathfrak{z}$ is entirely determined by the value of $\mathfrak{z}$
modulo $p^{N}$. 
\end{example}
\vphantom{}

Our second appetizer in non-archimedean function theory is, in all
honesty, a matter of non-archimedean \emph{functional }theory. It
is also one of several reasons why integration of $p$-adic-valued
functions of one or more $p$-adic variables is ``cursed'' as the
kids these days like to say.
\begin{defn}
Letting $C\left(\mathbb{Z}_{p},K\right)$ denote the space of continuous
functions from $\mathbb{Z}_{p}$ to $K$, where $K$ is a metrically
complete field extension of $\mathbb{Q}_{p}$. We say a linear functional
$\varphi:C\left(\mathbb{Z}_{p},K\right)\rightarrow K$ is\index{translation invariance}
\textbf{translation invariant }whenever $\varphi\left(\tau_{1}\left\{ f\right\} \right)=\varphi\left(f\right)$
for all $f\in C\left(\mathbb{Z}_{p},K\right)$ where: 
\begin{equation}
\tau_{1}\left\{ f\right\} \left(\mathfrak{z}\right)\overset{\textrm{def}}{=}f\left(\mathfrak{z}+1\right)
\end{equation}
\end{defn}
\begin{thm}
Let $K$ be a metrically complete field extension of $\mathbb{Q}_{p}$.
Then, the only translation invariant linear functional $\varphi:C\left(\mathbb{Z}_{p},K\right)\rightarrow K$
is the zero functional ($\varphi\left(f\right)=0$ for all $f$)\footnote{\textbf{Corollary 3 }on page 177 of Robert's book \cite{Robert's Book}.}. 
\end{thm}

\subsection{\label{subsec:1.3.3 Field-Extensions-=00003D000026}Field Extensions
and Spherical Completeness}

THROUGHOUT THIS SUBSECTION, WE WORK WITH A FIXED PRIME $p$.

\vphantom{}

It is often said that students of mathematics will not truly appreciate
the beauty of analytic functions until after they have endured the
pathological menagerie of functions encountered in a typical introductory
course in real analysis. In my experience, $p$-adic analysis does
much the same, only for fields instead than functions. The fact that
adjoining a square root of $-1$ to the real numbers results in a
metrically complete algebraically closed field is one of the greatest
miracles in all mathematics. Working with the $p$-adic numbers, on
the other hand, gives one an appreciation of just how \emph{awful}
field extensions can be.

The problem begins before we even take our first step up the tower
of $p$-adic field extensions. Anyone who has ever struggled with
the concept of negative numbers should be thankful that we do not
exist in a macroscopically $p$-adic universe. In real life, the existence
of negative integers is equivalent to the fact that $\mathbb{Z}^{\times}$\textemdash the
group of multiplicative;y invertible elements of $\mathbb{Z}$, a.k.a.
$\left\{ 1,-1\right\} $\textemdash is isomorphic to the cyclic group
$\mathbb{Z}/2\mathbb{Z}$. The isomorphism $\mathbb{Z}^{\times}\cong\mathbb{Z}/2\mathbb{Z}$
is also responsible for the existence of a meaningful ordering on
$\mathbb{Z}$ given by $<$ and friends.

Unfortunately, for any grade-school students who happen to live in
a $p$-adic universe, things are much less simple, and this has grave
implications for the study of field extensions of $\mathbb{Q}_{p}$.
To see why, for a moment, let us proceed naïvely. Let $K$ be a finite-degree
Galois extension\footnote{Recall, this means that the set:
\[
\left\{ \mathfrak{y}\in K:\sigma\left(\mathfrak{y}\right)=\mathfrak{y},\textrm{ }\forall\sigma\in\textrm{Gal}\left(K/\mathbb{Q}_{p}\right)\right\} 
\]
is equal to $\mathbb{Q}_{p}$.} of $\mathbb{Q}_{p}$. Because we are interested in doing analysis,
our first order of business is to understand how we might extend the
$p$-adic absolute value on $\mathbb{Q}_{p}$ to include elements
of $K$.

To do this, observe that since any $\sigma\in\textrm{Gal}\left(K/\mathbb{Q}_{p}\right)$
is a field automorphism of $K$, any candidate for a \emph{useful
}absolute value $\left|\cdot\right|_{K}$ on $K$ ought to satisfy:
\begin{equation}
\left|\sigma\left(\mathfrak{y}\right)\right|_{K}=\left|\mathfrak{y}\right|_{K},\textrm{ }\forall\mathfrak{y}\in K
\end{equation}
Moreover, because $\mathbb{Q}_{p}$ is a subset of $K$, it is not
unreasonable to require $\left|\mathfrak{y}\right|_{K}$ to equal
$\left|\mathfrak{y}\right|_{p}$ for all $\mathfrak{y}\in K$ that
happen to be elements of $\mathbb{Q}_{p}$. To that end, let us get
out the ol' determinant:
\begin{defn}
Define $N_{K}:K\rightarrow\mathbb{Q}_{p}$ by:
\begin{equation}
N_{K}\left(\mathfrak{y}\right)\overset{\textrm{def}}{=}\prod_{\sigma\in\textrm{Gal}\left(K/\mathbb{Q}_{p}\right)}\sigma\left(\mathfrak{y}\right)\label{eq:Definition of N_K}
\end{equation}
\end{defn}
\begin{rem}
$N_{K}$ is $\mathbb{Q}_{p}$-valued here precisely because $K$ is
a Galois extension.
\end{rem}
\vphantom{}

With this, we have that for all $\mathfrak{y}\in K$:
\[
\left|N_{K}\left(\mathfrak{y}\right)\right|_{p}=\left|N_{K}\left(\mathfrak{y}\right)\right|_{K}=\left|\prod_{\sigma\in\textrm{Gal}\left(K/\mathbb{Q}_{p}\right)}\sigma\left(\mathfrak{y}\right)\right|_{K}=\prod_{\sigma\in\textrm{Gal}\left(K/\mathbb{Q}_{p}\right)}\left|\sigma\left(\mathfrak{y}\right)\right|_{K}=\left|\mathfrak{y}\right|_{K}^{d}
\]
where $d$ is the order of $\textrm{Gal}\left(K/\mathbb{Q}_{p}\right)$\textemdash or,
equivalently (since $K$ is a Galois extension), the degree of $K$
over $\mathbb{Q}_{p}$. In this way, we can compute $\left|\mathfrak{y}\right|_{K}$
by taking the $d$th root of $\left|N_{K}\left(\mathfrak{y}\right)\right|_{p}$,
which is already defined, seeing as $N_{K}\left(\mathfrak{y}\right)$
is an element of $\mathbb{Q}_{p}$. With a little work, making this
argument rigorous leads to an extremely useful theorem:
\begin{thm}
\label{thm:Galois conjugates}Let $K$ be a Galois extension of $\mathbb{Q}_{p}$
of degree $d$, and let $\left|\cdot\right|_{p}$ be the $p$-adic
absolute value. Then, the map:
\begin{equation}
\mathfrak{y}\mapsto\left|N_{K}\left(\mathfrak{y}\right)\right|_{p}^{1/d}
\end{equation}
defines the unique absolute value $\left|\cdot\right|_{K}$ on $K$
which extends the $p$-adic absolute value of $\mathbb{Q}_{p}$.
\end{thm}
\begin{rem}
In fact, with a slight modification of $N_{K}$, one can show that
this method gives the absolute value on $K$ for \emph{any }finite-degree
extension $K$ of $\mathbb{Q}_{p}$, not only Galois extensions. Robert
does this on page 95 of \cite{Robert's Book}.
\end{rem}
\begin{example}
\label{exa:incorrect galois}As a sample application, consider $\mathbb{Q}_{5}$.
Let $\xi$ denote a primitive fourth root of unity, which we adjoin
to $\mathbb{Q}_{5}$. Then, since:
\begin{align*}
\left(1+2\xi\right)\left(1+2\xi^{3}\right) & =1+2\xi^{3}+2\xi+4\xi^{4}\\
\left(\xi+\xi^{3}=-1-\xi^{2}\right); & =5+2\left(-1-\xi^{2}\right)\\
\left(\xi^{2}=-1\right); & =5
\end{align*}
and so, by \textbf{Theorem \ref{thm:Galois conjugates}}: 
\begin{equation}
\left|1+2\xi\right|_{5}=\left|1+2\xi^{3}\right|_{5}=\sqrt{\left|5\right|_{5}}=\frac{1}{\sqrt{5}}
\end{equation}
Alas, this\emph{ }is \emph{wrong!} While this argument would be perfectly
valid if our base field was $\mathbb{Q}$ or $\mathbb{R}$, it fails
for $\mathbb{Q}_{5}$, due to the fact that $\xi$ was \emph{already
an element of $\mathbb{Q}_{5}$.}
\end{example}
\vphantom{}

With the sole exception of $\mathbb{Z}_{2}$, whose only roots of
unity are $1$ and $-1$, given any odd prime $p$, $\mathbb{Z}_{p}$
contains more roots of unity than just $\left\{ 1,-1\right\} $, and
this ends up complicating the study of field extensions. Fortunately,
the situation isn't \emph{too }bad; it is not difficult to determine
the roots of unity naturally contained in $\mathbb{Z}_{p}^{\times}$
(and hence, in $\mathbb{Z}_{p}$ and $\mathbb{Q}_{p}$):
\begin{thm}
When $p$ is an odd prime, $\mathbb{Z}_{p}^{\times}$ contains all
$\left(p-1\right)$th roots of unity. $\mathbb{Z}_{2}^{\times}$,
meanwhile, contains only second roots of unity ($-1$ and $1$).
\end{thm}
Proof: (Sketch) If $\mathfrak{z}\in\mathbb{Z}_{p}^{\times}$ satisfies
$\mathfrak{z}^{n}=1$ for some $n\geq1$, it must also satisfy $\mathfrak{z}^{n}\overset{p}{\equiv}1$,
so, the first $p$-adic digit of $\mathfrak{z}$ must be a primitive
root of unity modulo $p$. Conversely, for each primitive root $u\in\left\{ 1,2,\ldots,p-1\right\} $
of unity modulo $p$, with a little work (specifically, \textbf{Hensel's
Lemma }(see Section 6 of Chapter 1 of \cite{Robert's Book} for details)),
we can construct a unique $\mathfrak{z}\in\mathbb{Z}_{p}^{\times}$
with $u$ as its first digit. Because the group of multiplicative
units of $\mathbb{Z}/p\mathbb{Z}$ is isomorphic to $\mathbb{Z}/\left(p-1\right)\mathbb{Z}$,
this then yields the theorem.

Q.E.D.

\vphantom{}

Note that in order to compute the $p$-adic absolute value of an expression
involving a native root of unity of $\mathbb{Z}_{p}$, we need to
use said root's $p$-adic digit expansion.
\begin{example}
The roots of unity in $\mathbb{Z}_{7}$ are the $6$th roots of unity.
Of these, note that only two are primitive. By the argument used to
prove the previous theorem, note that these two roots of unity will
have $3$ and $5$, respectively, as their first $7$-adic digits,
seeing as those are the only two primitive roots of unity in $\mathbb{Z}/7\mathbb{Z}$.
If we let $\xi$ and $\zeta$ denote the the $3$-as-first-digit and
$5$-as-first-digit roots of unity, we then have that:
\begin{equation}
2\xi+1\overset{7}{\equiv}2\cdot3+1\overset{7}{\equiv}0
\end{equation}
and so, $\left|2\xi+1\right|_{7}\leq1/7$. By using \textbf{Hensel's
Lemma} to compute more and more digits of $\xi$, one can determine
the number of $0$s that occur at the start of $2\xi+1$'s $7$-adic
expansion and thereby determine its $7$-adic absolute value. On the
other hand:
\[
2\zeta+1\overset{7}{\equiv}2\cdot5+1\overset{7}{\equiv}4
\]
so, $\left|2\zeta+1\right|_{7}=1$, because $2\zeta+1$ is a $7$-adic
integer which is not congruent to $0$ modulo $7$.

In addition to Hensel's Lemma, there is a more algebraic approach
to computing the $7$-adic absolute values of these quantities. Since
$\xi$ is a primitive $6$th root of unity, the other primitive root
of unity (the one we denote by $\zeta$) must be $\xi^{5}$; $\xi^{2}$
and $\xi^{4}$ will be primitive $3$rd roots of unity, while $\xi^{3}$
will be the unique primitive square root of unity otherwise known
as $-1$. As such, writing $\zeta$ as $\xi^{5}$, observe that:
\begin{align*}
\left(2\xi+1\right)\left(2\xi^{5}+1\right) & =4\xi^{6}+2\xi+2\xi^{5}+1\\
 & =4+2\left(\xi+\xi^{5}\right)+1\\
 & =5+2\xi\left(1+\xi^{4}\right)\\
\left(\textrm{let }\omega=\xi^{2}\right); & =5+2\xi\left(1+\omega^{2}\right)
\end{align*}
Because $\omega$ is a primitive $3$rd root of unity, $1+\omega^{2}=-\omega=-\xi^{2}$,
and so:

\begin{align*}
\left(2\xi+1\right)\left(2\xi^{5}+1\right) & =5+2\xi\left(1+\omega^{2}\right)\\
 & =5+2\xi\left(-\xi^{2}\right)\\
 & =5-2\xi^{3}\\
\left(\xi^{3}=-1\right); & =5+2\\
 & =7
\end{align*}
Hence:
\[
\left|2\xi+1\right|_{7}\underbrace{\left|2\xi^{5}+1\right|_{7}}_{=\left|2\zeta+1\right|_{7}=1}=\left|7\right|_{7}=\frac{1}{7}
\]
So, $\left|2\xi+1\right|_{7}$ is precisely $1/7$.

What this shows in order to of said root in order to use complex exponential
notation to write roots of unity for any given odd prime $p$, we
need to specify the $p$-adic digits of $\xi_{p-1}\overset{\textrm{def}}{=}e^{2\pi i/\left(p-1\right)}$
so that we can then have a uniquely defined $p$-adic absolute value
for any polynomial in $\xi_{p-1}$ with coefficients in $\mathbb{Q}_{p}$.
\end{example}
\vphantom{}

In general, field extensions of $\mathbb{Q}_{p}$ are obtained in
the usual way, by adjoining to $\mathbb{Q}_{p}$ the root of a polynomial
with coefficients in $\mathbb{Q}_{p}$. Note that since every algebraic
number $\alpha$ is, by definition, the root of a polynomial with
coefficients in $\mathbb{Q}$, given any algebraic number $\alpha$,
we can create a $p$-adic field containing $\alpha$ by adjoining
$\alpha$ to $\mathbb{Q}_{p}$, seeing as $\mathbb{Q}_{p}$ contains
$\mathbb{Q}$ as a subfield. As described above, the absolute value
on $p$ extends with these extensions in the natural way, maintaining
its multiplicativity and ultrametric properties in the extension.
Thus, for example, in any extension $K$ of $\mathbb{Q}_{p}$, all
roots of unity will have a $p$-adic absolute value of $1$. More
generally, in any finite-degree extension $K$ of $\mathbb{Q}_{p}$,
the $p$-adic absolute value of an arbitrary element of $K$ will
be a number of the form $p^{r}$, where $r\in\mathbb{Z}\cup V$, where
$V$ is a finite set whose elements are non-integer rational numbers.

Far more problematic, however, is the issue of algebraic closure. 
\begin{defn}
We write \nomenclature{$\overline{\mathbb{Q}}_{p}$}{$p$-adic algebraic numbers}
$\overline{\mathbb{Q}}_{p}$, the \textbf{algebraic closure of $\mathbb{Q}_{p}$},
the field obtained by adjoining to $\mathbb{Q}_{p}$ the roots of
every polynomial with coefficients in $\mathbb{Q}_{p}$.\index{algebraic closure} 
\end{defn}
\vphantom{}

It is not an overgeneralization to say that algebraic number theory
is the study of the Galois group $\textrm{Gal}\left(\overline{\mathbb{Q}}/\mathbb{Q}\right)$,
and that the continued interest of the subject hinges upon the fact
that\textemdash unlike with $\mathbb{R}$\textemdash the algebraic
closure of $\mathbb{Q}$ is not a finite dimensional extension of
$\mathbb{Q}$. The same is true in the $p$-adic context: $\overline{\mathbb{Q}}_{p}$\textemdash the
algebraic closure of $\mathbb{Q}_{p}$\textemdash is an infinite-dimensional
extension of $\mathbb{Q}_{p}$. However, the parallels between the
two algebraic closures run deeper still. Just as $\overline{\mathbb{Q}}$
is incomplete as a metric space with respect to the complex absolute
value, so too is $\overline{\mathbb{Q}}_{p}$ incomplete as a metric
spaces with respect to the extension of $\left|\cdot\right|_{p}$.
Moreover, as infinite-dimensional extensions of their base fields,
neither $\overline{\mathbb{Q}}$ nor $\overline{\mathbb{Q}}_{p}$
is locally compact.

Just as taking the metric closure of $\overline{\mathbb{Q}}$ with
respect to the usual absolute value yields $\mathbb{C}$, the field
of complex numbers, taking the metric closure of $\overline{\mathbb{Q}}_{p}$
with respect to the $p$-adic absolute value yields the $p$-adic
complex numbers:
\begin{defn}
We write \nomenclature{$\mathbb{C}_{p}$}{the set of $p$-adic complex numbers \nomnorefpage}
$\mathbb{C}_{p}$ to denote the field of \index{$p$-adic!complex numbers}\textbf{$p$-adic
complex numbers}, defined here as the metric closure of $\overline{\mathbb{Q}}_{p}$
\end{defn}
\vphantom{}.

One of the many reasons complex analysis is the best form of analysis
is because $\mathbb{R}$ is locally compact and has $\mathbb{C}$
as a one-dimensional extension of $\mathbb{R}$; this then guarantees
that $\mathbb{C}$ is also locally compact. Unfortunately, the same
does not hold for $\mathbb{C}_{p}$; $\mathbb{C}_{p}$ is \emph{not
}locally compact. Worse still, $\mathbb{C}_{p}$ lacks a property
most analysts take for granted: \textbf{spherical completeness}.

Ordinarily, an (ultra)metric space is complete if and only if any
nested sequence of balls $B_{1}\supseteq B_{2}\supseteq\cdots$ whose
diameters/radii tend to $0$, the intersection $\bigcap_{n=1}^{\infty}B_{n}$
is non-empty. In the real/complex world, removing the requirement
that the diameters of the $B_{n}$s tend to $0$ would not affect
the non-emptiness of the intersection. Bizarrely, this is \emph{not
}always the case in an ultrametric space.
\begin{defn}
An ultrametric space $X$ is said to be \index{spherically complete}\textbf{spherically
complete }whenever every nested sequence of balls in $X$ has a non-empty
intersection. We say $X$ is \textbf{spherically incomplete }whenever
it is not spherically complete. 
\end{defn}
\begin{prop}
Any locally compact field is spherically complete. Consequently, for
any prime number $p$, $\mathbb{Q}_{p}$ is spherically complete,
as is any \emph{finite dimensional }field extension of $\mathbb{Q}_{p}$. 
\end{prop}
Proof: Compactness.

Q.E.D.

\vphantom{}

Despite this, there exist non-locally compact fields which are spherically
complete (see \cite{Robert's Book} for more details; Robert denotes
the spherical completion of $\mathbb{C}_{p}$ by $\Omega_{p}$).
\begin{thm}
Let $K$ be a non-archimedean valued field. Then, $K$ is spherically
complete whenever it is metrically complete, and has a discrete value
group. In particular, all locally compact metrically complete valued
fields (archimedean or not) are spherically complete\footnote{\textbf{Corollary 20.3} from \cite{Ultrametric Calculus}.}. 
\end{thm}
As it regards $\mathbb{C}_{p}$, the key properties are as follows: 
\begin{thm}
$\mathbb{C}_{p}$ is \textbf{not}\emph{ }spherically complete\emph{
\cite{Robert's Book}}. 
\end{thm}
\begin{thm}
$\mathbb{C}_{p}$ is not\emph{ }separable\footnote{That is, it contains no dense \emph{countable} subset.}\emph{
\cite{Robert's Book}}. 
\end{thm}
\vphantom{}

As will be mentioned in Subsection \ref{subsec:3.1.2 Banach-Spaces-over},
whether or not a non-archimedean field is spherically complete has
significant impact on the nature of the Banach spaces which can be
built over it. At present, it suffices to say that in most ``classical''
applications of $p$-adic analysis (especially to number theory and
algebraic geometry), spherical incompleteness is considered an undesirable
property, particularly because of its destructive interaction with
non-archimedean analogues of the \textbf{Hahn-Banach Theorem}\footnote{The thrust of this theorem, recall, is that we can take a continuous
linear functional $\phi$ defined on a subspace $\mathcal{Y}$ of
a given Banach space $\mathcal{X}$ and extend $\phi$ to a continuous
linear functional on $\mathcal{X}$ itself.}\textbf{ }of functional analysis \cite{Robert's Book}. Despite this,
the spherical \emph{incompleteness} of the fields in which our $\left(p,q\right)$-adic
functions will take their values plays a pivotal role in our work,
because of the effect this spherical incompleteness has upon the Fourier
analytic properties of our functions. See, for instance, \textbf{Theorem
\ref{thm:FST is an iso from measures to ell infinity}} on page \pageref{thm:FST is an iso from measures to ell infinity}
about the $\left(p,q\right)$-adic Fourier-Stieltjes transform.

Finally, we must give a brief mention of ``complex exponentials''
in the $p$-adic context. Although one can proceed to define a \textbf{$p$-adic
exponential function}\index{$p$-adic!exponential function} by way
of the usual power series $\exp_{p}:\mathbb{Q}_{p}\rightarrow\mathbb{Q}_{p}$:
\begin{equation}
\exp_{p}\left(\mathfrak{z}\right)\overset{\textrm{def}}{=}\sum_{n=0}^{\infty}\frac{\mathfrak{z}^{n}}{n!}\in\mathbb{Q}_{p}\label{eq:Definition of the p-adic exponential function}
\end{equation}
unlike its classical counterpart, this construction is quite unruly,
converging only for those $\mathfrak{z}$ with $\left|\mathfrak{z}\right|_{p}<p^{-\frac{1}{p-1}}$.
Even worse, thanks to the properties of $p$-adic power series (specifically,
due to \textbf{Strassman's Theorem} \cite{Robert's Book}), $\exp_{p}\left(\mathfrak{z}\right)$
is not a periodic function of $\mathfrak{z}$. As such, instead of
working natively, when $p$-adic analysis wants to make use of complex
exponentials in a manner comparable to their role in complex analysis,
we need to take a different approach, as we discuss in the next subsection.

\subsection{\label{subsec:1.3.4 Pontryagin-Duality-and}Pontryagin Duality and
Embeddings of Fields}

For those who haven't had the pleasure of meeting it, \textbf{Pontryagin
duality}\index{Pontryagin duality}\textbf{ }is that lovely corner
of mathematics where the concept of representing functions as Fourier
series is generalized to study complex-valued functions on spaces\textemdash specifically,
\textbf{locally compact abelian groups} (\textbf{LCAG}s)\footnote{Such a space is an abelian group that has a topology in which the
group's operations are continuous; moreover, as a topological space,
the group is locally compact.}\textemdash other than the circle ($\mathbb{R}/\mathbb{Z}$), the
line ($\mathbb{R}$), and their cartesian products. Essentially, given
a LCAG $G$, you can decompose functions $G\rightarrow\mathbb{C}$
as series or integrals of ``elementary functions'' (\textbf{characters}\index{character})
which transform in a simple, predictable way when they interact with
the group's operation. In classical Fourier analysis, the characters
are exponential functions $e^{2\pi ixt}$, $e^{-i\theta}$, and the
like, and their ``well-behavedness'' is their periodicity\textemdash their
invariance with respect to translations/shifts of the circle / line.

In terms of the literature, Folland's book \cite{Folland - harmonic analysis}
gives an excellent introduction and is very much worth reading in
general. For the reader on the go, Tao's course notes on his blog
\cite{Tao Fourier Transform Blog Post} are also quite nice. However,
because this dissertation will approach Haar measures and integration
theory primarily from the perspective of functional analysis, we will
not need to pull out the full workhorse of Pontryagin duality and
abstract harmonic analysis. Still, we need to spend at least \emph{some
}amount of time discussing characters and our conventions for using
them.

A character\textbf{ }is, in essence, a generalization of the complex
exponential. Classically, given an abelian group $G$, a \textbf{(complex-valued)
character }on $G$ is a group homomorphism $\chi:G\rightarrow\mathbb{T}$,
where $\left(\mathbb{T},\times\right)$ is the unit circle in $\mathbb{C}$,
realized as an abelian group with the usual complex-number multiplication
operation. Every character can be written as $\chi\left(g\right)=e^{2\pi i\xi\left(g\right)}$,
where $\xi:G\rightarrow\mathbb{R}/\mathbb{Z}$ is a \textbf{frequency},
a group homomorphism from $G$ to $\mathbb{R}/\mathbb{Z}$\textemdash the
group of real numbers in $\left[0,1\right)$\textemdash equipped with
the operation of addition modulo $1$. The set of all characters of
$G$ forms an abelian group under point-wise multiplication, and the
set of all frequencies of $G$ forms an abelian group under point-wise
addition. These two groups isomorphic to one another, and are both
realizations of the \textbf{(Pontryagin) dual }of $G$, denoted $\hat{G}$.

Depending on $G$, it can take some work to classify the characters
and frequencies and figure out explicit formulae to use for them.
Thankfully, since all of our work will be done over the additive group
$\mathbb{Z}_{p}$ (where $p$ is a prime), we only need to familiarize
ourselves with one set of formulae.
\begin{defn}[$\hat{\mathbb{Z}}_{p}$]
\ 

\vphantom{}

I. We write \nomenclature{$\hat{\mathbb{Z}}_{p}$}{$\overset{\textrm{def}}{=}\mathbb{Z}\left[\frac{1}{p}\right]/\mathbb{Z}$}$\hat{\mathbb{Z}}_{p}$
to denote the group $\mathbb{Z}\left[\frac{1}{p}\right]/\mathbb{Z}$.
This is the set of all $p$-ary rational numbers in $\left[0,1\right)$:
\begin{equation}
\hat{\mathbb{Z}}_{p}\overset{\textrm{def}}{=}\mathbb{Z}\left[\frac{1}{p}\right]/\mathbb{Z}=\left\{ \frac{k}{p^{n}}:n\in\mathbb{N}_{0},k\in\left\{ 0,\ldots,p^{n}-1\right\} \right\} \label{eq:Definition of Z_p hat}
\end{equation}
and is a group under the operation of addition modulo $1$. This is
the \textbf{additive}\footnote{The \emph{multiplicative} $p$-Prüfer group, generally denoted $\mathbb{Z}\left(p^{\infty}\right)$,
is the set of all $p$-power roots of unity in $\mathbb{T}$\textemdash that
is, all $\xi\in\mathbb{T}$ so that $\xi^{p^{n}}=1$ for some $n\in\mathbb{N}_{0}$\textemdash made
into a group by multiplication.}\textbf{ $p$-Prüfer group}\index{$p$-Prüfer group}. As the notation
$\hat{\mathbb{Z}}_{p}$ suggests, this group is one realization of
the \textbf{Pontryagin dual} of the additive group of $p$-adic integers:
$\left(\mathbb{Z}_{p},+\right)$.

\vphantom{}

II. Recall that every $p$-adic rational number $\mathfrak{z}\in\mathbb{Q}_{p}$
can be written as: 
\begin{equation}
\mathfrak{z}=\sum_{n=n_{0}}^{\infty}c_{n}p^{n}
\end{equation}
for constants $n_{0}\in\mathbb{Z}$ and $\left\{ c_{n}\right\} _{n\geq-n_{0}}\subseteq\left\{ 0,\ldots,p-1\right\} $.
The\textbf{ $p$-adic fractional part}\index{$p$-adic!fractional part},
denoted \nomenclature{$\left\{ \cdot\right\} _{p}$}{$p$-adic fractional part}
$\left\{ \cdot\right\} _{p}$, is the map from $\mathbb{Q}_{p}$ to
$\hat{\mathbb{Z}}_{p}$ defined by: 
\begin{equation}
\left\{ \mathfrak{z}\right\} _{p}=\left\{ \sum_{n=-n_{0}}^{\infty}c_{n}p^{n}\right\} _{p}\overset{\textrm{def}}{=}\sum_{n=n_{0}}^{-1}c_{n}p^{n}\label{eq:Definition of the p-adic fractional part}
\end{equation}
where the sum on the right is defined to be $0$ whenever $n_{0}\geq0$.
The \textbf{$p$-adic integer part}\index{$p$-adic!integer part},
denoted $\left[\cdot\right]_{1}$ \nomenclature{$ \left[\cdot\right]_{1}$}{$p$-adic integer part},
is the map from $\mathbb{Q}_{p}$ to $\mathbb{Z}_{p}$ defined by:
\begin{equation}
\left[\mathfrak{z}\right]_{1}=\left[\sum_{n=-n_{0}}^{\infty}c_{n}p^{n}\right]_{1}\overset{\textrm{def}}{=}\sum_{n=0}^{\infty}c_{n}p^{n}\label{eq:Definition of the p-adic integer part}
\end{equation}
\end{defn}
\vphantom{}

Note that $\left[\mathfrak{z}\right]_{1}=\mathfrak{z}-\left\{ \mathfrak{z}\right\} _{p}$.
Indeed, the $p$-adic fractional and integer parts are projections
from $\mathbb{Q}_{p}$ onto $\mathbb{Q}_{p}/\mathbb{Z}_{p}$ and $\mathbb{Z}_{p}$,
respectively. In fact, $\hat{\mathbb{Z}}_{p}$ and $\mathbb{Q}_{p}/\mathbb{Z}_{p}$
are isomorphic as additive groups. Moreover, we can express the group
$\left(\mathbb{Q}_{p},+\right)$ as the sum:
\begin{equation}
\mathbb{Q}_{p}=\mathbb{Z}\left[\frac{1}{p}\right]+\mathbb{Z}_{p}\label{eq:Sum decomposition of Q_p}
\end{equation}
where the right-hand side is group whose elements are of the form
$t+\mathfrak{z}$, where $t\in\mathbb{Z}\left[1/p\right]$ and $\mathfrak{z}\in\mathbb{Z}_{p}$.
Note, this sum is \emph{not }direct\footnote{For people who enjoy staring at diagrams of (short) exact sequences,
Subsection 5.4 (starting at page 40) of Robert's book \cite{Robert's Book}
gives a comprehensive discussion of (\ref{eq:Sum decomposition of Q_p})'s
failure to be a direct sum. Even though the \emph{set-theoretic }bijection
$\iota:\mathbb{Q}_{p}\rightarrow\hat{\mathbb{Z}}_{p}\times\mathbb{Z}_{p}$
given by $\iota\left(\mathfrak{y}\right)\overset{\textrm{def}}{=}\left(\left\{ \mathfrak{y}\right\} _{p},\left[\mathfrak{y}\right]_{1}\right)$
does more than your ordinary set-theoretic bijection (for example,
its restriction to $\mathbb{Z}_{p}$ \emph{is }an isomorphism of groups),
$\iota$ itself is not an isomorphism of groups. For example: 
\[
\iota\left(\frac{1}{p}+\frac{p-1}{p}\right)=\iota\left(1\right)=\left(0,1\right)
\]
\[
\iota\left(\frac{1}{p}\right)+\iota\left(\frac{p-1}{p}\right)=\left(\frac{1}{p},0\right)+\left(\frac{p-1}{p},0\right)=\left(1,0\right)=\left(0,0\right)
\]
}. For example:
\begin{equation}
\frac{1}{p}+2=\left(\frac{1}{p}+1\right)+1
\end{equation}
Despite this, every $p$-adic number is uniquely determined by its
fractional and integer parts, and we can view $\hat{\mathbb{Z}}_{p}$
as a subset of $\mathbb{Q}_{p}$ by identifying it with the subset
of $\mathbb{Q}_{p}$ whose elements all have $0$ as the value of
their $p$-adic integer part. In this way we can add and multiply
$t$s in $\hat{\mathbb{Z}}_{p}$ with $\mathfrak{z}$s in $\mathbb{Z}_{p}$
or $\mathbb{Q}_{p}$.
\begin{example}
For $t=\frac{1}{p}$ and $\mathfrak{z}=-3$, $t\mathfrak{z}=-\frac{3}{p}\in\mathbb{Q}_{p}$. 
\end{example}
\vphantom{}

Since every $p$-ary rational number $t\in\hat{\mathbb{Z}}_{p}$ is
an element of $\mathbb{Q}_{p}$, we can consider $t$'s $p$-adic
valuation and its $p$-adic absolute value. Specifically: 
\begin{fact}[\textbf{Working with $v_{p}$ and $\left|\cdot\right|_{p}$ on $\hat{\mathbb{Z}}_{p}$}]
\ 

\vphantom{}

I. $v_{p}\left(t\right)$ will be $\infty$ when $t=0$, and, otherwise,
will be $-n$, where $n$ is the unique integer $\geq1$ so that $t$
can be written as $t=k/p^{n}$, where $k\in\left\{ 1,\ldots,p^{n}-1\right\} $
is co-prime to $p$.

\vphantom{}

II. For non-zero $t$, $\left|t\right|_{p}=p^{-v_{p}\left(t\right)}$
will be the value in the denominator of $t$. Specifically, $\left|t\right|_{p}=p^{n}$
if and only if the irreducible fraction representing $t$ is $t=k/p^{n}$,
where $k\in\left\{ 1,\ldots,p^{n}-1\right\} $ is co-prime to $p$.

\vphantom{}

III. For non-zero $t$, $t\left|t\right|_{p}$ is the numerator $k$
of the irreducible fraction $k/p^{n}$ representing $t$.

\vphantom{}

IV. The $p$-adic fractional part also encodes information about congruences.
In particular, we have the following identity: 
\begin{equation}
\left\{ t\mathfrak{z}\right\} _{p}\overset{1}{\equiv}t\left[\mathfrak{z}\right]_{\left|t\right|_{p}},\textrm{ }\forall t\in\hat{\mathbb{Z}}_{p},\textrm{ }\forall\mathfrak{z}\in\mathbb{Z}_{p}\label{eq:p-adic fractional part identity modulo 1}
\end{equation}
where, $\left[\mathfrak{z}\right]_{\left|t\right|_{p}}$ is the projection
of $\mathfrak{z}$ modulo the power of $p$ in the denominator of
$t$, and where, recall, $\overset{1}{\equiv}$ means that the left
and the right hand side are equal to one another as real numbers modulo
$1$\textemdash that is, their difference is an integer.
\end{fact}
\begin{example}
Let $\mathfrak{z}\in\mathbb{Z}_{2}$. Then, we have the following
equality in $\mathbb{C}$
\begin{equation}
e^{2\pi i\left\{ \frac{\mathfrak{z}}{4}\right\} _{2}}\overset{\mathbb{C}}{=}e^{2\pi i\frac{\left[\mathfrak{z}\right]_{4}}{4}}
\end{equation}
 More generally, for any $\mathfrak{y}\in\mathbb{Q}_{p}$, we can
define: 
\begin{equation}
e^{2\pi i\left\{ \mathfrak{y}\right\} _{p}}\overset{\mathbb{C}}{=}\sum_{n=0}^{\infty}\frac{\left(2\pi i\left\{ \mathfrak{y}\right\} _{p}\right)^{n}}{n!}
\end{equation}
because $\left\{ \mathfrak{y}\right\} _{p}\in\hat{\mathbb{Z}}_{p}=\mathbb{Z}\left[\frac{1}{p}\right]/\mathbb{Z}\subset\left[0,1\right)\subset\mathbb{R}\subset\mathbb{C}$. 
\end{example}
\vphantom{}

There is also the issue of what happens when we ``mix our $p$s and
$q$s''.
\begin{example}
Let $\mathfrak{z}\in\mathbb{Z}_{2}$. Then, since $\left|\frac{1}{3}\right|_{2}=1$,
we have that $\frac{1}{3}\in\mathbb{Z}_{2}$, and thus, that $\frac{\mathfrak{z}}{3}\in\mathbb{Z}_{2}$.
However: 
\begin{equation}
\left\{ \frac{\mathfrak{z}}{24}\right\} _{2}=\left\{ \frac{1}{8}\frac{\mathfrak{z}}{3}\right\} _{2}\neq\frac{1}{24}\left[\mathfrak{z}\right]_{8}
\end{equation}
Rather, we have that: 
\begin{equation}
\left\{ \frac{\mathfrak{z}}{24}\right\} _{2}=\left\{ \frac{1}{8}\frac{\mathfrak{z}}{3}\right\} _{2}\overset{1}{\equiv}\frac{1}{8}\left[\frac{\mathfrak{z}}{3}\right]_{8}\overset{1}{\equiv}\frac{1}{8}\left[\mathfrak{z}\right]_{8}\left[\frac{1}{3}\right]_{8}\overset{1}{\equiv}\frac{\left[3\right]_{8}^{-1}}{8}\left[\mathfrak{z}\right]_{8}\overset{1}{\equiv}\frac{3}{8}\left[\mathfrak{z}\right]_{8}
\end{equation}
where $\left[3\right]_{8}^{-1}$ denotes the unique integer in $\left\{ 1,\ldots,7\right\} $
which is the multiplicative inverse of $3$ modulo $8$\textemdash in
this case, $3$, since $3$ is a square root of unity modulo $8$.
Lastly, to return to our discussion of characters, it turns out that
every complex-valued character on $\mathbb{Z}_{p}$ is a function
of the form $\mathfrak{z}\in\mathbb{Z}_{p}\mapsto e^{2\pi i\left\{ t\mathfrak{z}\right\} _{p}}\in\mathbb{C}$
for some unique $t\in\hat{\mathbb{Z}}_{p}$; this is an instance of
the isomorphism between the character group of $\mathbb{Z}_{p}$ and
$\hat{\mathbb{Z}}_{p}$. 
\end{example}
\vphantom{}

Note that all of the above discussion was predicated on the assumption
that our characters were \emph{complex-valued}. Importantly, because
$\left\{ \cdot\right\} _{p}$ outputs elements of $\hat{\mathbb{Z}}_{p}$,
each of which is a rational number, note also that the output of any
complex-valued character $\mathfrak{z}\in\mathbb{Z}_{p}\mapsto e^{2\pi i\left\{ t\mathfrak{z}\right\} _{p}}\in\mathbb{C}$
of $\mathbb{Z}_{p}$ is a $p$-power root of unity. Since every $p$-power
root of unity is the root of a polynomial of the form $X^{p^{n}}-1$
for some sufficiently large $n\geq1$, every $p$-power root of unity
is an \emph{algebraic }number. As such, our $p$-adic characters are
not \emph{true} $\mathbb{C}$-valued functions, but $\overline{\mathbb{Q}}$-valued
functions.

When working with functions taking values in a non-archimedean field
like $\mathbb{Q}_{q}$, where $q$ is a prime, it is generally more
useful to study them not by using Fourier series involving complex-valued
($\mathbb{C}$-valued) characters, but \textbf{$q$-adic (complex)
valued characters}\textemdash those taking values in $\mathbb{C}_{q}$.
However, regardless of whether we speak of the complex-valued or $q$-adic-valued
characters of $\mathbb{Z}_{p}$, our characters will still end up
taking values in $\overline{\mathbb{Q}}$. The reason there even exists
a verbal distinction between complex-valued and $q$-adic-valued characters
on $\mathbb{Z}_{p}$ is due to the needs of (algebraic) number theorists.

Ordinarily, in algebraic number theory, the question of how a given
field $\mathbb{F}$ is embedded in a larger field $K$ is very important,
and it is not uncommon to allow the embeddings to vary, or to even
do work by considering all possible embeddings \emph{simultaneously}.
However, in this dissertation, we will do exactly the opposite.

Our goal is to use the easy, comforting complex-exponential $e$-to-the-power-of-$2\pi i$-something
notation for roots of unity. Because the $\mathbb{C}$-valued characters
of $\mathbb{Z}_{p}$ are $\overline{\mathbb{Q}}$-valued, all of the
algebraic rules for manipulating expressions like $\left\{ t\mathfrak{z}\right\} _{p}$
and $e^{2\pi i\left\{ t\mathfrak{z}\right\} _{p}}$ will remain valid
when we view $e^{2\pi i\left\{ t\mathfrak{z}\right\} _{p}}$ as an
element of $\mathbb{C}_{q}$, because we will be identifying $e^{2\pi i\left\{ t\mathfrak{z}\right\} _{p}}$
with its counterpart in the copy of $\overline{\mathbb{Q}}$ that
lies in $\mathbb{C}_{q}$\textemdash a copy which is isomorphic to
the version of $\overline{\mathbb{Q}}$ that lies in $\mathbb{C}$.
If the reader is nervous about this, they can simply think of $e^{2\pi i\left\{ t\mathfrak{z}\right\} _{p}}$
as our notation of writing an abstract unitary character on $\mathbb{Z}_{p}$,
without reference to an embedding.
\begin{assumption*}[Embedding Convention]
For any primes $p$ and $q$, and any integer $n\geq1$, we will
write $e^{2\pi i/p^{n}}$ to denote a specific choice of a primitive
$p^{n}$th root of unity in $\mathbb{C}_{q}$, and shall then identify
$e^{2\pi ik/p^{n}}$ as the $k$th power of that chosen root of unity.
We also make these choices so that the are compatible in the natural
way, with 
\begin{equation}
\left(e^{2\pi i/p^{n+1}}\right)^{p}=e^{2\pi i/p^{n}},\textrm{ }\forall n\geq1
\end{equation}
More generally, when $q\geq5$, since the set of roots of unity in
$\mathbb{Q}_{q}$ is generated by a primitive $\left(q-1\right)$th
root of unity, we define $e^{2\pi i/\left(q-1\right)}$ to be the
primitive $\left(q-1\right)$th root of unity in $\mathbb{Z}_{q}^{\times}$
whose value mod $q$ is $r_{q}$, the smallest element of $\left(\mathbb{Z}/q\mathbb{Z}\right)^{\times}$
which is a primitive $\left(q-1\right)$th root of unity mod $q$.
Then, for every other root of unity $\xi$ in $\mathbb{Q}_{q}$, there
is a unique integer $k\in\left\{ 0,\ldots,q-2\right\} $ so that $\xi=\left(e^{2\pi i/\left(q-1\right)}\right)^{k}$,
and as such, we define $\xi$'s value mod $q$ to be the value of
$r_{q}^{k}$ mod $q$.
\end{assumption*}
In this convention, any \emph{finite }$\overline{\mathbb{Q}}$-linear
combination of roots of unity will then be an element of $\overline{\mathbb{Q}}$,
and, as such, we can and will view the linear combination as existing
in $\mathbb{C}$ and $\mathbb{C}_{q}$ \emph{simultaneously}. When
working with infinite sums, the field to which the sum belongs will
be indicated, when necessary, by an explicit mention of the field
in which the convergence happens to be occurring. \emph{The practical
implications of all this is that the reader will not need to worry
about considerations of local compactness, spherical completeness,
or field embeddings when doing or reading through computations with
complex exponentials in this dissertation. }Except for when limits
and infinite series come into play (and even then, only rarely), it
will be as if everything is happening in $\mathbb{C}$.

The reason everything works out is thanks to the fundamental identity
upon which all of our Fourier analytic investigations will be based:
\begin{equation}
\left[\mathfrak{z}\overset{p^{n}}{\equiv}\mathfrak{a}\right]\overset{\overline{\mathbb{Q}}}{=}\frac{1}{p^{n}}\sum_{\left|t\right|_{p}\leq p^{n}}e^{-2\pi i\left\{ t\mathfrak{a}\right\} _{p}}e^{2\pi i\left\{ t\mathfrak{z}\right\} _{p}},\textrm{ }\forall\mathfrak{z},\mathfrak{a}\in\mathbb{Z}_{p}
\end{equation}
where the sum is taken over all $t\in\hat{\mathbb{Z}}_{p}$ whose
denominators are at most $p^{n}$. This is the Fourier series for
the indicator function of the co-set $\mathfrak{a}+p^{n}\mathbb{Z}_{p}$.
As indicated, the equality holds in $\overline{\mathbb{Q}}$. Moreover,
observe that this identity is invariant under the action of $\textrm{Gal}\left(\overline{\mathbb{Q}}/\mathbb{Q}\right)$,
seeing as the left-hand side is always rational, evaluating to $1$
when $\mathfrak{z}$ is in $\mathfrak{a}+p^{n}\mathbb{Z}_{p}$ and
evaluating to $0$ for all other $\mathfrak{z}\in\mathbb{Z}_{p}$.
Because of this, nothing is lost in choosing a particular embedding,
and the theories we develop will be equally valid regardless of our
choice of an embedding. So, nervous number theorists can set aside
their worries and breathe easily. For staying power, we will repeat
this discussion in Subsection \ref{subsec:1.3.3 Field-Extensions-=00003D000026}.

One of the quirks of our approach is that when we actually go about
computing Fourier transforms and integrals for $q$-adic valued functions
of one or more $p$-adic variables, the computational formalisms will
be identical to their counterparts for real- or complex-valued functions
of one or more $p$-adic variables. As such, the overview of non-archimedean
Fourier analysis given in Subsections \ref{subsec:3.1.4. The--adic-Fourier}
and \ref{subsec:3.1.5-adic-Integration-=00003D000026} will include
all the necessary computational tools which will be drawn upon in
this dissertation.

That being said, I highly recommend the reader consult a source like
\cite{Bell - Harmonic Analysis on the p-adics} in case they are not
familiar with the nitty-gritty details of doing computations with
(real)-valued Haar measures on the $p$-adics and with the Fourier
transform of real- or complex-valued functions of a $p$-adic variable.
A more comprehensive\textemdash and rigorous\textemdash expositions
of this material can be found in the likes of \cite{Taibleson - Fourier Analysis on Local Fields,Automorphic Representations};
\cite{Automorphic Representations} is a book on representation theory
which deals with the methods of $p$-adic integration early on because
of its use in that subject, whereas \cite{Taibleson - Fourier Analysis on Local Fields}
is dedicated entirely to the matter of the Fourier analysis of complex-valued
functions on $\mathbb{Q}_{p}$ and other local fields.

\subsection{\label{subsec:1.3.5 Hensel's-Infamous-Blunder}Hensel's Infamous
Blunder and Geometric Series Universality}

In discussions with strangers over the internet, it became apparent
that some of the number-theoretic liberties I take in my work are
a source of controversy\textemdash and for good reason.

Despite the extraordinary importance that $p$-adic numbers would
eventually take in number theory, the $p$-adics' immediate reception
by the mathematical community of the \emph{fin de siècle }was one
of reservation and concerned apprehension. There are a variety of
reasons for this: the cultural attitudes of his contemporaries; Hensel\index{Hensel, Kurt}'s
unabashed zeal for his discovery; the nascent state of the theory
of topology\footnote{Poincaré published the installments of \emph{Analysis Situs }from
1895 through to 1904, contemporaneous with Hensel's introduction of
the $p$-adics in 1897; \emph{Analysis Situs }is to topology what
Newton's \emph{Principia }was to physics.} at that particular day and age; and, especially, Hensel's own fatal
blunder: a foundational flawed 1905 $p$-adic ``proof'' that $e$
was a transcendental number \cite{Gouvea's introudction to p-adic numbers book,Gouvea's p-adic number history slides,Koblitz's book,Conrad on p-adic series}.
Gouvea presents this blunder like so:
\begin{quotation}
Hensel's incorrect proof goes like this. Start from the equation:

\[
e^{p}=\sum_{n=0}^{\infty}\frac{p^{n}}{n!}
\]
Hensel checks that this series converges in $\mathbb{Q}_{p}$ and
concludes that it satisfies an equation of the form $y^{p}=1+p\varepsilon$
with $\varepsilon$ a $p$-adic unit. If $e$ is algebraic, it follows
that $\left[\mathbb{Q}\left(e\right):\mathbb{Q}\right]\geq p$. But
$p$ was arbitrary, so we have a contradiction, and $e$ is transcendental
\cite{Gouvea's p-adic number history slides}. 
\end{quotation}
Hensel's error is as much \emph{topological} it is number theoretic.
Specifically, he \textbf{\emph{incorrectly}} assumes the following
statement is true: 
\begin{assumption}
\emph{Let $\mathbb{F}$ and $\mathbb{K}$ be any field extensions
of $\mathbb{Q}$, equipped with absolute values $\left|\cdot\right|_{\mathbb{F}}$
and $\left|\cdot\right|$$_{\mathbb{K}}$, respectively, which make
$\left(\mathbb{F},\left|\cdot\right|_{\mathbb{F}}\right)$ and $\left(\mathbb{K},\left|\cdot\right|_{\mathbb{K}}\right)$
into complete metric spaces. Let $\left\{ x_{n}\right\} _{n\geq1}$
be a sequence in $\mathbb{Q}$ such that there are elements $\mathfrak{a}\in\mathbb{F}$
and $\mathfrak{b}\in\mathbb{K}$ so that $\lim_{n\rightarrow\infty}\left|\mathfrak{a}-x_{n}\right|_{\mathbb{F}}=0$
and $\lim_{n\rightarrow\infty}\left|\mathfrak{b}-x_{n}\right|_{\mathbb{F}}=0$.
Then $\mathfrak{a}=\mathfrak{b}$.} 
\end{assumption}
\vphantom{}

So, Hensel presumed that if $\mathbb{F}$ is a field extension of
$\mathbb{Q}$ in which the series $\sum_{n=0}^{\infty}\frac{p^{n}}{n!}$
converges, the sum of the series was necessarily $e^{p}$. This, of
course, is wrong. One might object and say ``but this is obvious:
not every $p$-adic integer can be realized as an ordinary real or
complex number!'', however, the situation is even more subtle than
that. For example, the statement in the previous paragraph remains
false even if we required $\mathfrak{a}\in\mathbb{F}\cap\mathbb{Q}$
and $\mathfrak{b}\in\mathbb{K}\cap\mathbb{Q}$. Consider the following: 
\begin{example}[\textbf{A Counterexample for Kurt Hensel}\footnote{Taken from \cite{Conrad on p-adic series}.}]
Let $a_{n}=\frac{1}{1+p^{n}}$, and let $r_{n}=a_{n}-a_{n-1}$, so
that $\sum_{n}r_{n}=\lim_{n\rightarrow\infty}a_{n}$. Then: 
\begin{align}
\sum_{n=0}^{\infty}r_{n} & \overset{\mathbb{R}}{=}0\\
\sum_{n=0}^{\infty}r_{n} & \overset{\mathbb{Q}_{p}}{=}1
\end{align}
because $a_{n}$ tends to $0$ in $\mathbb{R}$ and tends to $1$
in $\mathbb{Q}_{p}$. 
\end{example}
\vphantom{}

This is, quite literally, the oldest mistake in the subject, and still
worth pointing out to neophytes over a century later; I, myself, fell
for it in an earlier iteration of this dissertation. Nevertheless,
in this work, we are going to do what Hensel \emph{thought }he could
do, but in our case, \emph{it will be justified}. Such shenanigans
will occur multiple times in this dissertation, so it is worth taking
a moment to explain in detail what we are going to do, and why it
actually \emph{works}. And the exemplary exception that saves us is
none other than the tried and true geometric series formula. 
\begin{fact}[\textbf{\emph{Geometric Series Universality}}]
\textbf{\label{fact:Geometric series universality}} Let\index{geometric series universality}
$r\in\mathbb{Q}$, and suppose there is a valued field $\left(\mathbb{F},\left|\cdot\right|_{\mathbb{F}}\right)$
extending $\mathbb{Q}$ so that the series $\sum_{n=0}^{\infty}r^{n}$
converges in $\mathbb{F}$ (i.e., for which $\left|r\right|_{\mathbb{F}}<1$),
with the sum converging in $\mathbb{F}$ to $R\in\mathbb{F}$. Then,
$R\in\mathbb{Q}$, and, for any valued field $\left(\mathbb{K},\left|\cdot\right|_{\mathbb{K}}\right)$
extending $\mathbb{Q}$ for which the series $\sum_{n=0}^{\infty}r^{n}$
converges in $\mathbb{K}$ (i.e., for which $\left|r\right|_{\mathbb{K}}<1$),
the sum converges in $\mathbb{K}$ to $R$. \emph{\cite{Conrad on p-adic series}} 
\end{fact}
\vphantom{}

In practice, in our use of this universality, we will have a formal
sum $\sum_{n=0}^{\infty}a_{n}$ of rational numbers $\left\{ a_{n}\right\} _{n\geq0}\subseteq\mathbb{Q}$
such that $\sum_{n=0}^{\infty}a_{n}$ can be written as the sum of
finitely many geometric series: 
\begin{equation}
\sum_{n=0}^{\infty}b_{1}r_{1}^{n}+\sum_{n=0}^{\infty}b_{2}r_{2}^{n}+\sum_{n=0}^{\infty}b_{M}r_{M}^{n}
\end{equation}
where the $b_{m}$s and $r_{m}$s are rational numbers. In particular,
there will exist an integer $c\geq1$ so that: 
\begin{equation}
\sum_{n=0}^{cN-1}a_{n}\overset{\mathbb{Q}}{=}\sum_{n=0}^{N-1}b_{1}r_{1}^{n}+\sum_{n=0}^{N-1}b_{2}r_{2}^{n}+\sum_{n=0}^{N-1}b_{M}r_{M}^{n},\textrm{ }\forall N\geq1\label{eq:Partial sum decomposed as sum of partial sums of geometric series}
\end{equation}
Then, if there is a prime $p$ so that $\left|r_{m}\right|<1$ and
$\left|r_{m}\right|_{p}<1$ for all $m\in\left\{ 1,\ldots,M\right\} $,
the series: 
\begin{equation}
\sum_{n=0}^{N-1}b_{m}r_{m}^{n}
\end{equation}
will converge to $\frac{b_{m}}{1-r_{m}}$ in both $\mathbb{R}$ \emph{and
}$\mathbb{Q}_{p}$ as $N\rightarrow\infty$.
\begin{example}
The geometric series identity:
\[
\sum_{n=0}^{\infty}\left(\frac{3}{4}\right)^{n}=\frac{1}{1-\frac{3}{4}}=4
\]
holds in both $\mathbb{R}$ and $\mathbb{Q}_{3}$, in that the series
converges in both fields' topologies, and the limit of its partial
sum is $4$ in both of those topologies. Indeed:
\begin{equation}
\sum_{n=0}^{N-1}\left(\frac{3}{4}\right)^{n}=\frac{1-\left(\frac{3}{4}\right)^{N}}{1-\frac{3}{4}}=4-4\left(\frac{3}{4}\right)^{N}
\end{equation}
and so:
\[
\lim_{N\rightarrow\infty}\left|4-\sum_{n=0}^{N-1}\left(\frac{3}{4}\right)^{n}\right|=\lim_{N\rightarrow\infty}4\left(\frac{3}{4}\right)^{N}\overset{\mathbb{R}}{=}0
\]
and:

\begin{equation}
\lim_{N\rightarrow\infty}\left|4-\sum_{n=0}^{N-1}\left(\frac{3}{4}\right)^{n}\right|_{3}=\lim_{N\rightarrow\infty}\left|4\left(\frac{3}{4}\right)^{N}\right|_{3}=\lim_{N\rightarrow\infty}3^{-N}\overset{\mathbb{R}}{=}0
\end{equation}
\end{example}
\vphantom{}

Our second technique for applying universality is the observation
that, \emph{for any prime $p$, if a geometric series $\sum_{n=0}^{\infty}r^{n}$
(where $r\in\mathbb{Q}\cap\mathbb{Q}_{p}$) converges in $\mathbb{Q}_{p}$,
then its sum, $\frac{1}{1-r}$, is }\textbf{\emph{also}}\emph{ a rational
number. }In light of these two tricks, we will consider situations
like the following. Here, let $p$ and $q$ be primes (possibly with
$p=q$), and consider a function $f:\mathbb{Z}_{p}\rightarrow\mathbb{Q}_{q}$.
Then, if there is a subset $U\subseteq\mathbb{Z}_{p}$ such that,
for each $\mathfrak{z}\in U$, $f\left(\mathfrak{z}\right)$ is expressible
a finite linear combination of geometric series: 
\begin{equation}
f\left(\mathfrak{z}\right)\overset{\mathbb{Q}_{q}}{=}\sum_{m=1}^{M}\sum_{n=0}^{\infty}b_{m}r_{m}^{n}
\end{equation}
where the $b_{m}$s and the $r_{m}$s are elements of $\mathbb{Q}$,
we can view the restriction $f\mid_{U}$ as a rational-valued function
$f\mid_{U}:U\rightarrow\mathbb{Q}$. More generally, we can view $f\mid_{U}$
as a real-, or even complex-valued function on $U$. Depending on
the values of $\left|r\right|$ and $\left|r\right|_{q}$, we may
be able to compute the sum $\sum_{m=1}^{M}\sum_{n=0}^{\infty}b_{m}r_{m}^{n}$
solely in the topology of $\mathbb{R}$ or $\mathbb{C}$, solely in
the topology of $\mathbb{Q}_{q}$, or in both, simultaneously.

\chapter{\label{chap:2 Hydra-Maps-and}Hydra Maps and their Numina}

\includegraphics[scale=0.45]{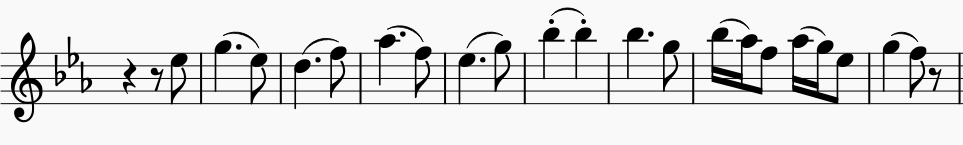}

\vphantom{}

Chapter 2 splits neatly (though lop-sidedly) into two halves. Section
\ref{sec:2.1 Hydra-Maps-=00003D000026 history} introduces Hydra maps\textemdash my
attempt to provide a taxonomy for some of the most important Collatz-type
maps, accompanied with numerous examples. Subsection \ref{subsec:2.1.2 It's-Probably-True}
is a historical essay on the Collatz Conjecture and most of the most
notable efforts to lay bare its mysteries.

Section \ref{sec:2.2-the-Numen}, meanwhile, is where we begin with
new results. The focus of \ref{sec:2.2-the-Numen} is the construction
of the titular \textbf{numen}, a $\left(p,q\right)$-adic function
we can associate to any sufficiently well-behaved Hydra map, and the
elucidation of its most important properties, these being its characterization
as the unique solution of a system of $\left(p,q\right)$-adic functional
equations\textemdash subject to a continuity-like condition\textemdash as
well as the all-important \textbf{Correspondence Principle}, the first
major result of my dissertation. Presented in Subsection \ref{subsec:2.2.3 The-Correspondence-Principle}\textemdash and
in four different variants, no less\textemdash the Correspondence
Principle establishes a direct correspondence between the values taken
by a numen and the periodic points and divergent points of the Hydra
map to which it is associated. \ref{sec:2.2-the-Numen} concludes
with Subsection \ref{subsec:2.2.4 Other-Avenues}, in which further
avenues of exploration are discussed, along with connections to other
recent work, principally Mih\u{a}ilescu's resolution of \textbf{Catalan's
Conjecture }and the shiny new state-of-the-art result on Collatz due
to Tao in 2019-2020 \cite{Tao Probability paper}.

\section{\label{sec:2.1 Hydra-Maps-=00003D000026 history}History and Hydra
Maps}

Much has been written about the Collatz Conjecture, be it efforts
to solve it, or efforts to dissuade others from trying. Bibliographically
speaking, I would place Lagarias' many writings on the subject\footnote{Such as his papers \cite{Lagarias-Kontorovich Paper,Lagarias' Survey},
and particularly his book \cite{Ultimate Challenge}.}, Wirsching's Book \cite{Wirsching's book on 3n+1}, and\textemdash for
this dissertation's purposes, Tao's 2019 paper \cite{Tao Probability paper}\textemdash near
the top of the list, along with K.R. Matthews maintains a dedicated
$3x+1$ page on his website, \href{http://www.numbertheory.org/3x\%2B1/}{page of Matthews' website},
filled with interesting links, pages, and, of course, his fascinating
powerpoint slides \cite{Matthews' slides,Matthews' Leigh Article,Matthews and Watts}.
The interested reader can survey these and other sources referenced
by this dissertation to get the ``lay of the land''\textemdash to
the extent that ``the land'' exists in this particular context.

The purpose of this section is two-fold: to introduce this dissertation's
fundamental objects of study\textemdash Hydra maps\textemdash and
to provide a short, likely inadequate history of the work that has
been done the Collatz Conjecture.

\subsection{\emph{\label{subsec:2.1.1 Release-the-Hydras!}Release the Hydras!
}- An Introduction to Hydra Maps}

We begin with the definition: 
\begin{defn}
\label{def:p-Hydra map}\nomenclature{$H$}{a Hydra map}\index{$p$-Hydra map}
\index{Hydra map!one-dimensional}\index{Hydra map}Fix an integer
$p\geq2$. A \textbf{$p$-Hydra map} is a map $H:\mathbb{Z}\rightarrow\mathbb{Z}$
of the form: 
\begin{equation}
H\left(n\right)=\begin{cases}
\frac{a_{0}n+b_{0}}{d_{0}} & \textrm{if }n\overset{p}{\equiv}0\\
\frac{a_{1}n+b_{1}}{d_{1}} & \textrm{if }n\overset{p}{\equiv}1\\
\vdots & \vdots\\
\frac{a_{p-1}n+b_{p-1}}{d_{p-1}} & \textrm{if }n\overset{p}{\equiv}p-1
\end{cases},\textrm{ }\forall n\in\mathbb{Z}\label{eq:Def of a Hydra Map on Z}
\end{equation}
and where $a_{j}$, $b_{j}$, and $d_{j}$ are integer constants (with
$a_{j},d_{j}\geq0$ for all $j$) so that the following two properties
hold:

\vphantom{}

I. (``co-primality'') $a_{j},d_{j}>0$ and $\gcd\left(a_{j},d_{j}\right)=1$
for all $j\in\left\{ 0,\ldots,p-1\right\} $.

\vphantom{}

II. (``one-sidedness'') $H\left(\mathbb{N}_{0}\right)\subseteq\mathbb{N}_{0}$
and $H\left(-\mathbb{N}_{1}\right)\subseteq-\mathbb{N}_{1}$, where
$-\mathbb{N}_{1}=\left\{ -1,-2,-3,\ldots\right\} $.

\vphantom{}

I call $H$ an \textbf{integral }Hydra map\index{Hydra map!integral},
if, in addition, it satisfies:

\vphantom{}

III. (``integrality''\footnote{For future or broader study, this condition might be weakened, or
even abandoned. It is not needed for the construction of the numen,
for example.}) For each $j\in\left\{ 0,\ldots,p-1\right\} $, $\frac{a_{j}n+b_{j}}{d_{j}}$
is an integer if and only if $n\overset{p}{\equiv}j$.

\vphantom{}

Hydra maps not satisfying (III) are said to be \textbf{non-integral}.\index{Hydra map!non-integral}
Finally, $H$ is said to be \textbf{prime }whenever $p$ is prime.
\end{defn}
\vphantom{}

When speaking of or working with these maps, the following notions
are quite useful:
\begin{defn}[\textbf{Branches of $H$}]
For each $j\in\mathbb{Z}/p\mathbb{Z}$, we write $H_{j}:\mathbb{R}\rightarrow\mathbb{R}$
to denote the \textbf{$j$th branch}\index{Hydra map!$j$th branch}
of $H$, defined as the function: 
\begin{equation}
H_{j}\left(x\right)\overset{\textrm{def}}{=}\frac{a_{j}x+b_{j}}{d_{j}}\label{eq:Definition of H_j}
\end{equation}
\end{defn}
Next, as a consequence of (\ref{eq:Def of a Hydra Map on Z}), observe
that for all $m\in\mathbb{N}_{0}$ and all $j\in\left\{ 0,\ldots,p-1\right\} $,
we can write: 
\[
\underbrace{H\left(pm+j\right)}_{\in\mathbb{N}_{0}}=\frac{a_{j}\left(pm+j\right)+b_{j}}{d_{j}}=\frac{pa_{j}m}{d_{j}}+\frac{ja_{j}+b_{j}}{d_{j}}=\frac{pa_{j}}{d_{j}}m+\underbrace{H\left(j\right)}_{\in\mathbb{N}_{0}}
\]
Setting $m=1$, we see that the quantity $\frac{pa_{j}}{d_{j}}=H\left(pm+j\right)-H\left(j\right)$
must be an integer for each value of $j$. This quantity occurs often
enough in computations to deserve getting a symbol of its own.
\begin{defn}[$\mu_{j}$]
We write \nomenclature{$\mu_{j}$}{$\overset{\textrm{def}}{=}\frac{pa_{j}}{d_{j}}$}$\mu_{j}$
to denote: 
\begin{equation}
\mu_{j}\overset{\textrm{def}}{=}\frac{pa_{j}}{d_{j}}\in\mathbb{N}_{1},\textrm{ }\forall j\in\mathbb{Z}/p\mathbb{Z}\label{eq:Def of mu_j}
\end{equation}
In particular, since $a_{j}$ and $d_{j}$ are always co-prime, \textbf{note
that $d_{j}$ must be a divisor of $p$}. In particular, this then
forces $\mu_{0}/p$ (the derivative of $H_{0}$) to be a positive
rational number which is not equal to $1$. 
\end{defn}
\begin{rem}
With this notation, we then have that: 
\begin{equation}
H\left(pm+j\right)=\mu_{j}m+H\left(j\right),\textrm{ }\forall m\in\mathbb{Z},\textrm{ }\forall j\in\mathbb{Z}/p\mathbb{Z}\label{eq:Value of rho m + j under H}
\end{equation}
\end{rem}
\vphantom{}

For lack of a better or more descriptive name, I have decided to call
maps of this form $p$-Hydra maps, in reference to the many-headed
monster of Greek myth. Much like their legendary namesake, Hydra maps
are not easily conquered. Questions and conjectures regarding these
maps' dynamical properties\textemdash such as the number of map's
periodic points, the existence of divergent trajectories for a given
map, or the total number of the map's orbit classes\textemdash rank
among the most difficult problems in all mathematics \cite{Lagarias' Survey}.

Case in point, the quintessential example of a Hydra map\textemdash for
the handful of people who aren't already aware of it\textemdash is
the infamous Collatz map: 
\begin{defn}
\label{def:The-Collatz-map}The\index{Collatz!map}\index{$3x+1$ map}
\textbf{Collatz map} $C:\mathbb{Z}\rightarrow\mathbb{Z}$ is the function
defined by: 
\begin{equation}
C\left(n\right)\overset{\textrm{def}}{=}\begin{cases}
\frac{n}{2} & \textrm{if }n\overset{2}{\equiv}0\\
3n+1 & \textrm{if }n\overset{2}{\equiv}1
\end{cases}\label{eq:Collatz Map}
\end{equation}
\end{defn}
First proposed by Lothar Collatz\footnote{Actually, Collatz came up with the map in 1929-1930; Matthews' website
has a translation of a letter by Collatz himself in which the man
states as such. It arose, apparently, as a diversion during time Collatz
spent attending classes by Perron and Schur, among others \cite{Collatz Letter}.} in 1937 \cite{Collatz Biography}, the eponymous map's well-deserved
mathematical infamy stems from the egregiously simple, yet notoriously
difficult conjecture of the same name \cite{Lagarias' Survey,Collatz Biography}: 
\begin{conjecture}[\textbf{Collatz Conjecture}\footnote{Also known as the $3n+1$ problem, the $3x+1$ problem, the Syracuse
Problem, the Hailstone problem, to give but a few of its many monickers.} \index{Collatz!Conjecture}]
\label{conj:Collatz}For every $n\in\mathbb{N}_{1}$, there exists
a $k\geq0$ so that $C^{\circ k}\left(n\right)=1$. 
\end{conjecture}
\begin{rem}
An equivalent reformulation of this result is the statement: \emph{$\mathbb{N}_{1}$
is an irreducible orbit class of $C$.} 
\end{rem}
\vphantom{}

A moment's thought reveals that the Collatz Conjecture actually consists
of \emph{two} separate statements: 
\begin{conjecture}[\textbf{Weak Collatz Conjecture}\footnote{Also known as the ``(No) Finite Cycles'' Conjecture.}]
\label{conj:Weak Collatz}The only periodic points of $C$ in $\mathbb{N}_{1}$
are $1$, $2$, and $4$. 
\end{conjecture}
\begin{conjecture}[\textbf{Divergent Trajectories Conjecture}]
$C$ has no divergent trajectories in $\mathbb{N}_{1}$. That is
to say, for each $n\in\mathbb{N}_{1}$: 
\begin{equation}
\sup_{k\geq0}C^{\circ k}\left(n\right)<\infty\label{eq:No divergent trajectories}
\end{equation}
\end{conjecture}
\vphantom{}

Since $3n+1$ will be even whenever $n$ is odd, it is generally more
expedient to alter $C$ by inserting an extra division by $2$ to
obtain the so-called Shortened Collatz map: 
\begin{defn}
\label{def:The-Shortened-Collatz Map}\nomenclature{$T_{3}$}{shortened Collatz map}The\index{$3x+1$ map}
\index{Collatz!map}\textbf{Shortened Collatz map }$T_{3}:\mathbb{Z}\rightarrow\mathbb{Z}$
is defined by: 
\begin{equation}
T_{3}\left(n\right)\overset{\textrm{def}}{=}\begin{cases}
\frac{n}{2} & \textrm{if }n\overset{2}{\equiv}0\\
\frac{3n+1}{2} & \textrm{if }n\overset{2}{\equiv}1
\end{cases}\label{eq:Definition of T_3}
\end{equation}
\end{defn}
\vphantom{}

The general dynamics of which are identical to that of $C$, with
$T_{3}\left(n\right)=C\left(n\right)$ for all even $n$ and $T_{3}\left(n\right)=C\left(C\left(n\right)\right)$
for all odd $n$. The use of $T$ or $T_{3}$ to denote this map is
standard in the literature (ex. \cite{Lagarias-Kontorovich Paper},
\cite{Wirsching's book on 3n+1}).

After the Collatz map $C$ and its shortened counterpart, $T_{3}$,
the next most common level of generalization comes from considering
the larger family of ``shortened $qx+1$ maps'': 
\begin{defn}
\label{def:Shortened qx plus 1 map}\nomenclature{$T_{q}$}{shortened $qx+1$ map}Let
$q$ be an odd integer $\geq1$. Then, the \textbf{Shortened $qx+1$
map}\index{$qx+1$ map} $T_{q}:\mathbb{Z}\rightarrow\mathbb{Z}$ is
the function defined by: 
\begin{equation}
T_{q}\left(n\right)\overset{\textrm{def}}{=}\begin{cases}
\frac{n}{2} & \textrm{if }n\overset{2}{\equiv}0\\
\frac{qn+1}{2} & \textrm{if }n\overset{2}{\equiv}1
\end{cases}\label{eq:Definition of T_q}
\end{equation}
\end{defn}
\vphantom{}

Of the (shortened) $qx+1$ maps, the most well-known are $T_{3}$
and $T_{5}$, with $T_{5}$ being of interest because its behavior
seems to be diametrically opposed to that of $T_{3}$: one the one
hand, the set of positive integers iterated by $T_{3}$ to $\infty$
has zero density \cite{Lagarias' Survey,Terras 76,Terras 79}; on
the other hand, $T_{5}$ iterates almost every positive integer to
$\infty$ \cite{Lagarias-Kontorovich Paper}. Nevertheless, despite
the probabilistic guarantee that almost every positive integer should
belong to a divergent trajectory, it has yet to be proven that any
\emph{specific} positive integer belongs to a divergent trajectory
\cite{Lagarias-Kontorovich Paper}!

A notable feature of Collatz Studies is its overriding single-mindedness.
Discussions of $C$ and its various shortened versions\footnote{Tao, for example, considers a variant which sends an $n$ odd integer
to the next odd integer $n^{\prime}$ in the forward orbit of $n$
under $C$; that is $n\mapsto\left(3n+1\right)\left|3n+1\right|_{2}$,
where $\left|\cdot\right|_{2}$ is the $2$-adic absolute value \cite{Tao Blog}.} predominate, with most in-depth number theoretic studies of Hydra
maps being focused on $C$. While other Hydra maps do receive a certain
degree of attention, it is generally from a broader, more holistic
viewpoint.

Some of the best examples of these broader studies can be found works
of \index{Matthews, K. R.}K. R. Matthews and his colleagues and collaborators,
going all the way back to the 1980s. The \emph{leitmotif} of this
body of work is the use of Markov chains to model the maps' dynamics.
These models provide a particularly fruitful heuristic basis for conjectures
regarding the maps' ``typical'' trajectories, and the expected preponderance
of cycles relative to divergent trajectories and vice-versa \cite{dgh paper,Matthews' Conjecture,Matthews' Leigh Article,Matthews' slides,Matthews and Watts}.
This line of investigation dates back to the work of H. Möller in
the 1970s \cite{Moller's paper (german)}, who formulated a specific
type of generalized multiple-branched Collatz map and made conjectures
about the existence of cycles and divergent trajectories. Although
probabilistic techniques of this sort will not on our agenda here,
Matthews and his colleagues' investigations provide ample ``non-classical''
examples of what I call Hydra maps.
\begin{example}
Our first examples are due to Leigh (\cite{Leigh's article,Matthews' Leigh Article});
we denote these as $L_{1}$ and $L_{2}$, respectively: 
\begin{equation}
L_{1}\left(n\right)\overset{\textrm{def}}{=}\begin{cases}
\frac{n}{4} & \textrm{if }n\overset{8}{\equiv}0\\
\frac{n+1}{2} & \textrm{if }n\overset{8}{\equiv}1\\
20n-40 & \textrm{if }n\overset{8}{\equiv}2\\
\frac{n-3}{8} & \textrm{if }n\overset{8}{\equiv}3\\
20n+48 & \textrm{if }n\overset{8}{\equiv}4\\
\frac{3n-13}{2} & \textrm{if }n\overset{8}{\equiv}5\\
\frac{11n-2}{4} & \textrm{if }n\overset{8}{\equiv}6\\
\frac{n+1}{8} & \textrm{if }n\overset{8}{\equiv}7
\end{cases}
\end{equation}
\begin{equation}
L_{2}\left(n\right)\overset{\textrm{def}}{=}\begin{cases}
12n-1 & \textrm{if }n\overset{4}{\equiv}0\\
20n & \textrm{if }n\overset{4}{\equiv}1\\
\frac{3n-6}{4} & \textrm{if }n\overset{4}{\equiv}2\\
\frac{n-3}{4} & \textrm{if }n\overset{4}{\equiv}3
\end{cases}
\end{equation}
Note that neither $L_{1}$ nor $L_{2}$ satisfy the integrality condition,
seeing as they have branches which are integer-valued for all integer
inputs. Also, while $L_{1}$ fixes $0$, $L_{2}$ sends $0$ to $-1$,
and as such, $L_{2}$ is not a $4$-Hydra map \emph{as written}; this
can be rectified by conjugating $L_{2}$ by an appropriately chosen
affine linear map\textemdash that is, one of the form $an+b$, where
$a$ and $b$ are integers; here, $a$ needs to be co-prime to $p=4$.
Setting $a=1$ gives: 
\begin{equation}
L_{2}\left(n-b\right)+b=\begin{cases}
12n-11b-1 & \textrm{if }n\overset{4}{\equiv}b\\
20n-19b & \textrm{if }n\overset{4}{\equiv}b+1\\
\frac{3n+b-6}{4} & \textrm{if }n\overset{4}{\equiv}b+2\\
\frac{n+3b-3}{4} & \textrm{if }n\overset{4}{\equiv}b+3
\end{cases}
\end{equation}
Setting $b=1$ then yields the map $\tilde{L}_{2}:\mathbb{N}_{0}\rightarrow\mathbb{N}_{0}$:
\begin{equation}
\tilde{L}_{2}\left(n\right)\overset{\textrm{def}}{=}L_{2}\left(n-1\right)+1=\begin{cases}
\frac{n}{4} & \textrm{if }n\overset{4}{\equiv}0\\
12n-12 & \textrm{if }n\overset{4}{\equiv}1\\
20n-19 & \textrm{if }n\overset{4}{\equiv}2\\
\frac{3n-5}{4} & \textrm{if }n\overset{4}{\equiv}3
\end{cases}\label{eq:Conjugated Hydra map}
\end{equation}
which, unlike $L_{2}$, satisfies conditions (I) and (II) for in (\ref{eq:Def of a Hydra Map on Z}).
Conjugating maps in this way is a useful tool to put them in a ``standard
form'' of sorts, one in which the conjugated map sends the non-negative
integers $\mathbb{N}_{0}$ to $\mathbb{N}_{0}$ and, preferably, also
fixes $0$. 
\end{example}
\vphantom{}

Our next example of a map is also non-integral, yet, it is significant
because it possess \emph{provably} divergent trajectories\index{divergent!trajectory}.
Moreover, it may be useful to broaden the notion of Hydra maps to
include maps like these, seeing as, despite its integrality failure,
we can still construct its numen using the methods of Section \ref{sec:2.2-the-Numen}.
\begin{example}[\textbf{Matthews' Map \& Divergent Trajectories}]
\label{exa:Matthews' map}The following map is due to Matthews \cite{Matthews' Conjecture}:
\begin{equation}
M\left(n\right)\overset{\textrm{def}}{=}\begin{cases}
7n+3 & \textrm{if }n\overset{3}{\equiv}0\\
\frac{7n+2}{3} & \textrm{if }n\overset{3}{\equiv}1\\
\frac{n-2}{3} & \textrm{if }n\overset{3}{\equiv}2
\end{cases}\label{eq:Matthews' Conjecture Map}
\end{equation}
That this map fails to satisfy the integrality requirement can be
seen in the fact that its $0$th branch $7n+3$ outputs integers even
when $n$ is \emph{not }congruent to $0$ mod $3$.

That being said, as defined, note that $-1$ is the only fixed point
of $M$ in $\mathbb{Z}$. For our purposes, we will need to work with
maps that fix $0$. To do so, we can conjugate $M$ by $f\left(n\right)=n-1$
($f^{-1}\left(n\right)=n+1$) to obtain: 
\begin{equation}
\tilde{M}\left(n\right)\overset{\textrm{def}}{=}\left(f^{-1}\circ M\circ f\right)\left(n\right)=M\left(n-1\right)+1=\begin{cases}
\frac{n}{3} & \textrm{if }n\overset{3}{\equiv}0\\
7n-3 & \textrm{if }n\overset{3}{\equiv}1\\
\frac{7n-2}{3} & \textrm{if }n\overset{3}{\equiv}2
\end{cases}\label{eq:Matthews' Conjecture Map, conjugated}
\end{equation}
Conjugating this map by $g\left(n\right)=-n$ ($g^{-1}=g$) then allows
us to study $\tilde{T}$'s behavior on the negative integers: 
\begin{equation}
\tilde{T}_{-}\left(n\right)\overset{\textrm{def}}{=}-\tilde{T}\left(-n\right)=\begin{cases}
\frac{n}{3} & \textrm{if }n\overset{3}{\equiv}0\\
\frac{7n+2}{3} & \textrm{if }n\overset{3}{\equiv}1\\
7n+3 & \textrm{if }n\overset{3}{\equiv}2
\end{cases}\label{eq:Matthews' conjecture map, conjugated to the negatives}
\end{equation}
In general, after conjugating a map (Hydra or not) so that it fixes
$0$, conjugating by $g\left(n\right)=-n$ is how one studies the
behavior of the maps on the non-positive integers.
\end{example}
\begin{rem}
It may be of interest to note that the $1$st and $2$nd branches
$\tilde{M}_{1}\left(x\right)=7x-3$ and $\tilde{M}_{2}\left(x\right)=\frac{7x-2}{3}$
commute with one another:
\begin{equation}
\tilde{M}_{1}\left(\tilde{M}_{2}\left(x\right)\right)=7\left(\frac{7x-2}{3}\right)-3=\frac{49x-14}{3}-3=\frac{49x-23}{3}
\end{equation}

\begin{equation}
\tilde{M}_{2}\left(\tilde{M}_{1}\left(x\right)\right)=\frac{7\left(7x-3\right)-2}{3}=\frac{49x-21-2}{3}=\frac{49x-23}{3}
\end{equation}
\end{rem}
\vphantom{}

As was mentioned above, the import of $M$ lies in the fact that it
has provably divergent trajectories: 
\begin{prop}
\label{prop:Matthews' map}$M$ has at least two divergent trajectories
over $\mathbb{Z}$:

\vphantom{}

I. The set of all non-negative integer multiples of $3$ (which are
iterated to $+\infty$);

\vphantom{}

II. The set of all negative integer multiples of three (which are
iterated to $-\infty$). 
\end{prop}
Proof: Since $7n+3$ is congruent to $0$ mod $3$ whenever $n$ is
congruent to $0$ mod $3$, any non-negative (resp. negative) integer
multiple of $3$ is iterated by $M$ to positive (resp. negative)
infinity.

Q.E.D.

\vphantom{}

For extra incentive, it is worth mentioning that there is money on
the table. \index{Matthews, K. R.}Matthews has put a bounty of $100$
Australian dollars, to be given to anyone who can bring him the severed,
fully-proven head of the following conjecture: 
\begin{conjecture}[\textbf{Matthews' Conjecture}\index{Matthews' Conjecture} \cite{Matthews' Conjecture}]
\label{conj:Matthews conjecture}The irreducible orbit classes of
$M$ in $\mathbb{Z}$ are the divergent orbit classes associated to
the trajectories (I) and (II) described in \textbf{\emph{Proposition
\ref{prop:Matthews' map}}}, and the orbit classes corresponding to
the attracting cycles $\left\{ -1\right\} $ and $\left\{ -2,-4\right\} $.
\end{conjecture}
\vphantom{}

Matthews' powerpoint slides \cite{Matthews' slides} (freely accessible
from his website) contain many more examples, both along these lines,
and in further, more generalized forms. The slides are quite comprehensive,
detailing generalizations and conjectures due to Hasse and Herbert
Möller\index{Möller, Herbert}, as well as delving into the details
of the aforementioned Markov-chain-based analysis of the maps and
their dynamics. The slides also detail how the Markov methods can
be extended along with the maps from $\mathbb{Z}$ to spaces of $p$-adic
integers, as well as to the ring of profinite integers\index{profinite integers}
(a.k.a,\emph{ polyadic integers}) given by the direct product $\prod_{p\in\mathbb{P}}\mathbb{Z}_{p}$.
However, for our purposes, the most significant offering in Matthews'
slides are the generalizations given in which simple Collatz-type
maps are defined over discrete spaces other than $\mathbb{Z}$. 
\begin{example}[\textbf{A ``Multi-Dimensional'' Hydra Map}]
The first example of such an extension appears to from Leigh in 1983\cite{Matthews' slides}.
Leigh considered the map\footnote{The notation $T_{3,\sqrt{2}}$ is my own.}
$T_{3,\sqrt{2}}:\mathbb{Z}\left[\sqrt{2}\right]\rightarrow\mathbb{Z}\left[\sqrt{2}\right]$
defined by: 
\begin{equation}
T_{3,\sqrt{2}}\left(z\right)\overset{\textrm{def}}{=}\begin{cases}
\frac{z}{\sqrt{2}} & \textrm{if }z\overset{\sqrt{2}}{\equiv}0\\
\frac{3z+1}{\sqrt{2}} & \textrm{if }z\overset{\sqrt{2}}{\equiv}1
\end{cases},\textrm{ }\forall z\in\mathbb{Z}\left[\sqrt{2}\right]\label{eq:Definition of T_3,root 2}
\end{equation}
which acts on the ring of quadratic integers of the form $a+b\sqrt{2}$,
where $a,b\in\mathbb{Z}$. Here, as indicated, a congruence $z\overset{\sqrt{2}}{\equiv}w$
for $z,w\in\mathbb{Z}\left[\sqrt{2}\right]$ means that $z-w$ is
of the form $v\sqrt{2}$, where $v\in\mathbb{Z}\left[\sqrt{2}\right]$.
Thus, for example, $3+2\sqrt{2}\overset{\sqrt{2}}{\equiv}1$, since:
\begin{equation}
3+2\sqrt{2}-1=2+2\sqrt{2}=\left(\sqrt{2}+2\right)\sqrt{2}
\end{equation}

Leigh's \index{multi-dimensional!Hydra map}map is the prototypical
example of what I call a \index{Hydra map!multi-dimensional}\textbf{Multi-Dimensional
Hydra Map}. The ``multi-dimensionality'' becomes apparent after
applying a bit of linear algebra. We can establish a ring isomorphism
of $\mathbb{Z}\left[\sqrt{2}\right]$ and $\mathbb{Z}^{2}$ by associating
$z=a+b\sqrt{2}\in\mathbb{Z}\left[\sqrt{2}\right]$ with the column
vector $\left(a,b\right)$. A bit of algebra produces: 
\begin{equation}
T_{3,\sqrt{2}}\left(a+b\sqrt{2}\right)=\begin{cases}
b+\frac{a}{2}\sqrt{2} & \textrm{if }a=0\mod2\\
3b+\frac{3a+1}{2}\sqrt{2} & \textrm{if }a=1\mod2
\end{cases},\textrm{ }\forall\left(a,b\right)\in\mathbb{Z}^{2}\label{eq:T_3,root 2 of a plus b root 2}
\end{equation}
Expressing the effect of $T_{3,\sqrt{2}}$ on $\left(a,b\right)$
in terms of these coordinates then gives us a map $\tilde{T}_{3,\sqrt{2}}:\mathbb{Z}^{2}\rightarrow\mathbb{Z}^{2}$
defined by: 
\begin{eqnarray}
\tilde{T}_{3,\sqrt{2}}\left(\left[\begin{array}{c}
a\\
b
\end{array}\right]\right) & \overset{\textrm{def}}{=} & \begin{cases}
\left[\begin{array}{cc}
1 & 0\\
0 & 2
\end{array}\right]^{-1}\left(\left[\begin{array}{cc}
0 & 1\\
1 & 0
\end{array}\right]\left[\begin{array}{c}
a\\
b
\end{array}\right]\right) & \textrm{if }a=0\mod2\\
\left[\begin{array}{cc}
1 & 0\\
0 & 2
\end{array}\right]^{-1}\left(\left[\begin{array}{cc}
0 & 3\\
3 & 0
\end{array}\right]\left[\begin{array}{c}
a\\
b
\end{array}\right]+\left[\begin{array}{c}
0\\
1
\end{array}\right]\right) & \textrm{if }a=1\mod2
\end{cases}\label{eq:Definition of the Lattice analogue of T_3,root 2}
\end{eqnarray}
which is a ``multi-dimensional'' analogue of a Hydra map; the branch:
\begin{equation}
n\mapsto\frac{an+b}{d}
\end{equation}
\textemdash an affine linear map on $\mathbb{Q}$\textemdash has been
replaced by an affine linear map: 
\begin{equation}
\mathbf{n}\mapsto\mathbf{D}^{-1}\left(\mathbf{A}\mathbf{n}+\mathbf{b}\right)
\end{equation}
on $\mathbb{Q}^{2}$, for $2\times2$ invertible matrices $\mathbf{A}$
and $\mathbf{D}$ and a $2$-tuple $\mathbf{b}$, all of which have
integer entries. In this way, notion of Hydra map defined at the beginning
of this subsection is really but a ``one-dimensional'' incarnation
of the more general notion of a map on a finitely generated module
over a principal ideal domain.

Note in (\ref{eq:Definition of the Lattice analogue of T_3,root 2})
the presence of the permutation matrix: 
\begin{equation}
\left[\begin{array}{cc}
0 & 1\\
1 & 0
\end{array}\right]
\end{equation}
which swaps the coordinate entries of a given $2$-tuple. As will
be seen in Section \ref{sec:5.1 Hydra-Maps-on}, this intertwining
of an affine linear map involving diagonal matrices with permutation
matrices is characteristic of the matrix representations of multi-dimensional
Hydra maps, illustrating the much greater degrees of freedom and flexibility
engendered by stepping out of $\mathbb{Z}$ and into spaces of higher
dimensionality.
\end{example}
\vphantom{}

Other than one additional example in Leigh's vein\textemdash presented
at the very end of Subsection \ref{subsec:5.1.2 Co=00003D0000F6rdinates,-Half-Lattices,-and}\textemdash Matthews
told me in our 2019 correspondence that he knew of no other examples
of Hydra maps of this type in the literature. For all intents and
purposes, the subject is simply unstudied. All the more reason, then,
for someone to try to whip it into shape.

\cite{Matthews' slides} goes even further, presenting near the end
examples of Hydra-like maps defined on a ring of polynomials with
coefficients in a finite field. While it may be possible to extend
my techniques to study such maps, I have yet to give consideration
to how my methods will apply when working with Hydra-like maps on
spaces of positive characteristic.

Finally, it should be mentioned that, for simplicity's sake, this
dissertation only concerns prime Hydra maps\textemdash those for which
$p$ is prime. While it remains to be seen whether or note my construction
of $\chi_{H}$ holds water in the case where $p$ is an arbitrary
composite integer, it seems pretty straightforward to modify the theory
presented here to deal with the case where $p$ is the power of a
prime\textemdash so, a $p^{n}$-Hydra map, where $p$ is prime and
$n\geq1$.

\newpage{}

\subsection{\emph{\label{subsec:2.1.2 It's-Probably-True}}A Short History of
the Collatz Conjecture}
\begin{quote}
\begin{flushright}
\emph{As you know, the $3x+1$ problem is a notorious problem, with
perhaps my favorite quote calling it a Soviet conspiracy to slow down
American mathematics.} 
\par\end{flushright}
\begin{flushright}
\textemdash Stephen J. Miller\footnote{Professor Miller made this remark to me as part of an e-mail correspondence
with me back in June of 2021.} 
\par\end{flushright}

\end{quote}
\index{Collatz!Conjecture} The Collatz Conjecture is infamous for
playing hard to get. With other difficult problems, lack of progress
toward a solution often comes with new ideas and insights worth studying
in their own right. On the other hand, for Collatz and its ilk, the
exact opposite appears to hold: a tool or insight of use for Collatz
seems to have little if any applicability to the outside mathematical
world. For people\textemdash like myself\textemdash with a vested
interest in mainstreaming Collatz studies this is especially problematic.
It isn't much help with conquering Collatz, either, given how much
mathematical progress owes to the discovery of intriguing new structures,
relations, and techniques, even if those novelties might not be immediately
applicable to the ``big'' questions we would like to solve.

The purpose of this essay is to give an overview of the primary ways
in which the mathematical community at large has attempted to frame
and understand the Collatz Conjecture and related arithmetical dynamical
systems. That my approach is essentially orthogonal to all of the
ones we are about to encounter is actually of the utmost relevance
to this discussion; I hope it will inspire readers to find other ways
of looking beyond the ``standard'' methods, so as to improve our
chances of finding exactly the sort of unexpected connections which
might be needed to better understand problems like Collatz.

To the extent that Collatz studies exists as a coherent subject, however,
the most preponderant trends are statistical in nature, where studies
center around defining a handful of interesting quantities\textemdash hoping
to bottle a sliver of the lightning animating Collatz's behavior\textemdash and
then demonstrating that interesting conclusions can be drawn. Arguably
the most impactful techniques of the statistical school have their
roots in the groundbreaking work done by Riho Terras\index{Terras, Riho}
in the 1970s (\cite{Terras 76,Terras 79}; see also Lagarias' survey
\cite{Lagarias' Survey}). Terras' paradigm rests upon two statistics:
the \textbf{stopping time} and the \textbf{parity sequence} (or \textbf{parity
vector}). 
\begin{defn}[\textbf{Stopping times}\footnote{Adapted from \cite{Lagarias-Kontorovich Paper}.}]
Let $\lambda$ be any positive real number, and let $T:\mathbb{Z}\rightarrow\mathbb{Z}$
be a map.

\vphantom{}

I. For any $n\in\mathbb{Z}$, the \textbf{$\lambda$-decay stopping
time}\index{stopping time}\textbf{ }of $n$ under $T$, denoted $\sigma_{T,\lambda}^{-}\left(n\right)$,
is the smallest integer $k\geq0$ so that $T^{\circ k}\left(n\right)<\lambda n$:
\begin{equation}
\sigma_{T,\lambda}^{-}\left(n\right)\overset{\textrm{def}}{=}\inf\left\{ k\geq0:\frac{T^{\circ k}\left(n\right)}{n}<\lambda\right\} \label{eq:Definition of the lambda decay stopping time}
\end{equation}

\vphantom{}

II. For any $n\in\mathbb{Z}$, the \textbf{$\lambda$-growth stopping
time }of $n$ under $T$, denoted $\sigma_{T,\lambda}^{+}\left(n\right)$,
is the smallest integer $k\geq0$ so that $T^{\circ k}\left(n\right)>\lambda n$:
\begin{equation}
\sigma_{T,\lambda}^{+}\left(n\right)\overset{\textrm{def}}{=}\inf\left\{ k\geq0:\frac{T^{\circ k}\left(n\right)}{n}>\lambda\right\} \label{eq:Definition of the lambda growth stopping time}
\end{equation}
\end{defn}
\vphantom{}

Terras' claim to fame rests on his ingenious implementation of these
ideas to prove: 
\begin{thm}[\textbf{Terras' Theorem} \cite{Terras 76,Terras 79}]
Let $S$ be set of positive integers with finite \index{Terras' Theorem}$1$-decay\footnote{i.e., $\lambda$-decay, with $\lambda=1$.}
stopping time under the Collatz map $C$. Then, $S$ has a natural
density of $1$; that is to say: 
\begin{equation}
\lim_{N\rightarrow\infty}\frac{\left|S\cap\left\{ 1,2,3,\ldots,N\right\} \right|}{N}=1\label{eq:Terras' Theorem}
\end{equation}
\end{thm}
In words, Terras' Theorem is effectively the statement \emph{almost
every positive integer does not belong to a divergent trajectory of
$C$}; or, equivalently, \emph{the set of divergent trajectories of
$C$ in $\mathbb{N}_{1}$ has density $0$}.\emph{ }This result is
of double importance: not only does it prove that $S$ \emph{possesses
}a density, it also shows that density to be $1$. One of the frustrations
of analytic number theory is that the limit (\ref{eq:Terras' Theorem})
defining the natural density of a set $S\subseteq\mathbb{N}_{1}$
does not exist for every\footnote{As an example, the set $S=\bigcup_{n=0}^{\infty}\left\{ 2^{2n},\ldots,2^{2n+1}-1\right\} $
does not have a well-defined natural density. For it, the $\limsup$
of (\ref{eq:Terras' Theorem}) is $2/3$, while the $\liminf$ is
$1/3$.} $S\subseteq\mathbb{N}_{1}$.

Tao's breakthrough was to improve ``almost every positive integer
does not go to infinity'' to ``almost all orbits of the Collatz
map attain almost bounded values''\textemdash indeed, that is the
very title of his paper \cite{Tao Probability paper}.

It is worth noting that by doing little more than replacing every
instance of a $3$ in Terras' argument with a $5$\textemdash and
adjusting all inequalities accordingly\textemdash an analogous result
can be obtained for the Shortened $5x+1$ map\index{$5x+1$ map} (or
$T_{q}$, for any odd $q\geq5$) \cite{Lagarias-Kontorovich Paper}: 
\begin{thm}[\textbf{Terras' Theorem for $qx+1$}]
Let $S$ be set of positive integers with finite $1$-growth\footnote{i.e., $\lambda$-growth, with $\lambda=1$.}
stopping time under the Shortened $qx+1$ map $T_{q}$. Then, $S$
has a natural density of $1$. 
\end{thm}
\vphantom{}

That is to say, while the set of positive integers which belong to
divergent trajectories of $T_{3}$ has density $0$, the set of divergent
trajectories of $T_{q}$ for $q\geq5$ has density $1$. Surreally,
for $q\geq5$, even though we know that almost every positive integer
should be iterated to $\infty$ by $T_{q}$, mathematics has yet to
prove that any specific positive integer \emph{suspected}\footnote{Such as $7$, for the case of the $5x+1$ map.}
of membership in a divergent trajectory actually belongs to one!

Next, we have Terras' notion of the \index{parity sequence}\textbf{parity
sequence} of a given input $n$. The idea is straightforward: since
the behavior of an integer $n$ under iterations of one of the $T_{q}$
maps is determined by the parity of those iterates\textemdash that
is, the value mod $2$\textemdash it makes sense to keep track of
those values. As such, the parity sequence of $n$ under $T_{q}$
is the sequence $\left[n\right]_{2},\left[T_{q}\left(n\right)\right]_{2},\left[T_{q}^{\circ2}\left(n\right)\right]_{2},\ldots$
consisting of the values (mod $2$) of the iterates of $n$ under
the map $T_{q}$ (\ref{eq:Definition of T_q}) \cite{Lagarias' Survey}.
More generally, since the iterates of $n$ under an arbitrary $p$-Hydra
map $H$ are determined not by their iterates' parities, but by their
iterates' values mod $p$, it will be natural to generalize Terras'
parity sequence into the \textbf{$p$-ity sequence}\index{$p$-ity sequence}\textbf{
}$\left[n\right]_{p},\left[H\left(n\right)\right]_{p},\left[H^{\circ2}\left(n\right)\right]_{p},\ldots$
consisting of the values mod $p$ of the iterates of $n$ under $H$;
doing so will play a crucial role early on in our construction of
the numen of $H$.

It should not come as much of a surprise that the parity sequence
of a given $n$ under $T_{q}$ ends up being intimately connected
to the behavior of $T_{q}$'s extension to the $2$-adic integers
(see for instance \cite{Parity Sequences}), where we apply the even
branch of $T_{q}$ to those $\mathfrak{z}\in\mathbb{Z}_{2}$ with
$\mathfrak{z}\overset{2}{\equiv}0$ and apply the odd branch of $T_{q}$
to those $\mathfrak{z}\in\mathbb{Z}_{2}$ with $\mathfrak{z}\overset{2}{\equiv}1$
Extending the environment of study from $\mathbb{Z}$ to $\mathbb{Z}_{2}$
in this manner is another one of Terras' innovations, one which is
relevant to our approach for all the obvious reasons.

Despite these promising connections, trying to conquer Collatz simply
by extending it to the $2$-adics is a quixotic undertaking, thanks
to the following astonishing result\footnote{Adapted from \cite{Lagarias-Kontorovich Paper}.}: 
\begin{thm}
\label{thm:shift map}For any odd integer $q\geq1$, the $2$-adic
extension of $T_{q}$ is a homeomorphism of $\mathbb{Z}_{2}$, and
is topologically and measurably conjugate to the $2$-adic shift map
$\theta:\mathbb{Z}_{2}\rightarrow\mathbb{Z}_{2}$ defined by: 
\begin{equation}
\theta\left(c_{0}+c_{1}2^{1}+c_{2}2^{2}+\cdots\right)\overset{\textrm{def}}{=}c_{1}+c_{2}2^{1}+c_{3}2^{2}+\cdots\label{eq:Definition of the 2-adic shift map}
\end{equation}
meaning that for every odd integer $q\geq1$, there exists a homeomorphism
$\Phi_{q}:\mathbb{Z}_{2}\rightarrow\mathbb{Z}_{2}$ which preserves
the (real-valued) $2$-adic Haar probability measure\footnote{$\mu\left(\Phi_{q}^{-1}\left(U\right)\right)=\mu\left(U\right)$,
for all measurable $U\subseteq\mathbb{Z}_{2}$.} so that: 
\begin{equation}
\left(\Phi_{q}\circ T_{q}\circ\Phi_{q}^{-1}\right)\left(\mathfrak{z}\right)=\theta\left(\mathfrak{z}\right),\textrm{ }\forall\mathfrak{z}\in\mathbb{Z}_{2}\label{eq:Conjugacy of T_q and the shift map}
\end{equation}
\end{thm}
\vphantom{}

As this theorem shows, in both a topological and measure theoretic
perspective, for any odd integers $p,q\geq1$, the dynamics $T_{p}$
and $T_{q}$ on $\mathbb{Z}_{2}$ are utterly indistinguishable from
one another, and from the shift map, which is ergodic on $\mathbb{Z}_{2}$
in a very strong sense \cite{Lagarias-Kontorovich Paper}. Because
$\mathbb{Z}$ has zero real-valued Haar measure in $\mathbb{Z}_{2}$,
it was already a given that employing tools of \index{ergodic theory}ergodic
theory and the measure-theoretic branches of dynamical systems would
be something of a wash; those methods cannot detect phenomena occurring
in sets of measure zero. However, the conjugate equivalence of the
different $T_{q}$ maps would appear to suggest that, at least over
the $2$-adics, probabilistic methods will useless for understanding
the dynamics of the $qx+1$ maps, let alone the Collatz map. Nevertheless,
as shown in \cite{Parity Sequences} and related works, there may
yet be fruitful work to be done in the $2$-adic setting, not by directly
studying the dynamics of $T_{3}$, but by how certain functions associated
to $T_{3}$ affect certain subsets of $\mathbb{Z}_{2}$.

Prior to Tao's innovations \cite{Tao Probability paper}, among the
motley collection of techniques used to study the Collatz Conjecture,
the least \emph{in}effective arguments in support of the Conjecture's
truth were based on the so-called ``difference inequalities'' of
I. Krasikov\index{Krasikov, Ilia} \cite{Krasikov}. In the 1990s,
Applegate and Lagarias\index{Lagarias, Jeffery} used Krasikov's inequalities
(along with a computer to assist them with a linear program in $\left(3^{9}-1\right)/2$
variables! \cite{Applegate and Lagarias - Difference inequalities})
to establish the asymptotic: 
\begin{thm}[\textbf{Applegate-Krasikov-Lagarias} \cite{Applegate and Lagarias - Difference inequalities}]
For any integer $a$, define $\pi_{a}:\mathbb{R}\rightarrow\mathbb{N}_{0}$
by: 
\begin{equation}
\pi_{a}\left(x\right)\overset{\textrm{def}}{=}\left|\left\{ n\in\mathbb{Z}:\left|n\right|\leq x\textrm{ and }T^{\circ k}\left(n\right)=a\textrm{ for some }k\in\mathbb{N}_{0}\right\} \right|\label{eq:Definition of Pi_a}
\end{equation}
Then, for any $a\in\mathbb{Z}$ which is not a multiple of $3$, there
is a positive real constant $c_{a}$ so that: 
\begin{equation}
\pi_{a}\left(x\right)\geq c_{a}x^{0.81},\textrm{ }\forall x\geq a\label{eq:Apple-Lag-Kra asymptotic theorem}
\end{equation}
\end{thm}
\vphantom{}

This school of work \cite{Applegate and Lagarias - Trees,Applegate and Lagarias - Difference inequalities,Krasikov,Wirsching's book on 3n+1}
intersects with graph-theoretic considerations (such as \cite{Applegate and Lagarias - Trees}),
a well-established player in probabilistic number theory and additive
combinatorics in the twentieth century and beyond.

With the exception of Rozier's work in \cite{Parity Sequences}, all
of the avenues of inquiry discussed above\textemdash even Tao's!\textemdash suffer
the same weakness, inherent to any \emph{probabilistic }approach:
they can only establish truths modulo ``small'' sets (zero density,
zero measure, etc.). Although results in this tradition provide often
deeply insightful, psychologically satisfying substitutes for resolutions
of the Collatz Conjecture, in part or in whole, they nevertheless
fail to capture or address whatever number-theoretic mechanisms are
ultimately responsible for the Collatz's dynamics, and the dynamics
of Hydra maps in general. In this respect, the Markov-chain\index{Markov chains}-based
work of Möller, Matthews, and their school offers an intriguing contrast,
if only because their methods show\textemdash or, at least, \emph{suggest}\textemdash how
the dynamical properties of Hydra maps might be related to the constants
used to define them.

Möller's generalization of Collatz began by fixing co-prime integers
$m,d>1$, with $m>d$ and considering the map $T_{\textrm{Möller}}:\mathbb{Z}\rightarrow\mathbb{Z}$
defined by \cite{Moller's paper (german),Matthews' slides}: 
\begin{equation}
T_{\textrm{Möller}}\left(x\right)\overset{\textrm{def}}{=}\begin{cases}
\frac{x}{d} & \textrm{if }x\overset{d}{\equiv}0\\
\frac{mx-1}{d} & \textrm{if }x\overset{d}{\equiv}\left[m\right]_{d}^{-1}\\
\frac{mx-2}{d} & \textrm{if }x\overset{d}{\equiv}2\left[m\right]_{d}^{-1}\\
\vdots & \vdots\\
\frac{mx-\left(d-1\right)}{d} & \textrm{if }x\overset{d}{\equiv}\left(d-1\right)\left[m\right]_{d}^{-1}
\end{cases},\textrm{ }\forall x\in\mathbb{Z}\label{eq:Moller's map}
\end{equation}
Möller's analysis led him to conjecture: 
\begin{conjecture}
$T_{\textrm{Möller}}$ eventually iterates every $x\in\mathbb{Z}$
if and only if $m<d^{d/\left(d-1\right)}$. Also, regardless of whether
or not $m<d^{d/\left(d-1\right)}$, $T_{\textrm{Möller}}$ has finitely
many cycles. 
\end{conjecture}
\vphantom{}

Matthews \cite{Matthews' slides} generalizes Möller's work like so.
Let $d\geq2$, and let $m_{0},\ldots,m_{d-1}$ be non-zero integers
(not necessarily positive). Also, for $j\in\left\{ 0,\ldots,d-1\right\} $,
let $r_{j}\in\mathbb{Z}$ satisfy $r_{j}\overset{d}{\equiv}jm_{j}$.
Then, define $T_{\textrm{Matthews}}:\mathbb{Z}\rightarrow\mathbb{Z}$
by: 
\begin{equation}
T_{\textrm{Matthews}}\left(x\right)\overset{\textrm{def}}{=}\begin{cases}
\frac{m_{0}x-r_{0}}{d} & \textrm{if }x\overset{d}{\equiv}0\\
\vdots & \vdots\\
\frac{m_{d-1}x-r_{d-1}}{d} & \textrm{if }x\overset{d}{\equiv}d-1
\end{cases},\textrm{ }\forall x\in\mathbb{Z}\label{eq:Matthews' Moller-map}
\end{equation}
He calls this map \textbf{relatively prime }whenever $\gcd\left(m_{j},d\right)=1$
for all $j\in\mathbb{Z}/d\mathbb{Z}$, and then makes the following
conjectures (given in \cite{Matthews' slides}): 
\begin{conjecture}
For the relatively prime case:

\vphantom{}

I. Every trajectory of $T_{\textrm{Matthews}}$ on $\mathbb{Z}$ is
eventually periodic whenever: 
\begin{equation}
\prod_{j=0}^{d-1}\left|m_{j}\right|<d^{d}
\end{equation}

\vphantom{}

II. If: 
\begin{equation}
\prod_{j=0}^{d-1}\left|m_{j}\right|>d^{d}
\end{equation}
then the union of all divergent orbit classes of $T_{\textrm{Matthews}}$
has density $1$ in $\mathbb{Z}$.

\vphantom{}

III. Regardless of the inequalities in (I) and (II), $T_{\textrm{Matthews}}$
always has a finite, non-zero number of cycles in $\mathbb{Z}$.

\vphantom{}

IV. Regardless of the inequalities in (I) and (II), for any $x\in\mathbb{Z}$
which belongs to a divergent orbit class of $T_{\textrm{Matthews}}$,
the iterates of $x$ under $T$ are uniformly distributed modulo $d^{n}$
for every $n\geq1$; that is: 
\begin{equation}
\lim_{N\rightarrow\infty}\frac{1}{N}\left|\left\{ k\in\left\{ 0,\ldots,N-1\right\} :T^{\circ k}\left(x\right)\overset{d^{n}}{\equiv}j\right\} \right|=\frac{1}{d^{n}}\label{eq:Matthews' Conjecture on uniform distribution modulo powers of rho of iterates of elements of divergent trajectories}
\end{equation}
\end{conjecture}
\vphantom{}

Even though these are conjectures, rather than results, they have
the appeal of asserting a connection between the maps' parameters
and the maps' dynamics. Even though probabilistic methods are not
my strong suit, I am fond of approaches like these because of their
open-mindedness, as well as the similarities between them and some
of my currently unpublished investigations (principally ergodic theoretic)
into the actions of Hydra maps on $\check{\mathbb{Z}}=\prod_{p\in\mathbb{P}}\mathbb{Z}_{p}$,
the ring of profinite integers (the product is taken over all primes
$p$) and its subgroups (ex: for Collatz, $\mathbb{Z}_{2}$ and $\mathbb{Z}_{2}\times\mathbb{Z}_{3}$).
As Lagarias has frequently told me (or intimated) in my on-and-off
correspondences with him over the past five years, any significant
progress on the Collatz Conjecture should either shed light on the
dynamics of close relatives like the $5x+1$ map, or give insight
into why the arguments used depend on the $3$ in $3x+1$. As such,
if we're going to try to surmount an insurmountable obstacle like
the Collatz Conjecture, we might as well cast wide nets, in search
of broader structures, avoiding any unnecessary dotage over $3x+1$
until we think there might be a way to tackle it. That, of course,
is exactly what I intend to do here.

\newpage{}

\section{\label{sec:2.2-the-Numen}$\chi_{H}$, the Numen of a Hydra Map}

THROUGHOUT THIS SECTION, WE ASSUME $p$ IS A PRIME AND THAT $H$ IS
A $p$-HYDRA MAP WITH $H\left(0\right)=0$.

\vphantom{}

In 1978, Böhm and Sontacchi \cite{Bohm and Sontacchi} (see also \cite{Wirsching's book on 3n+1})
established what has since become a minor (and oft-rediscovered) cornerstone
of Collatz literature: 
\begin{thm}[\textbf{Böhm-Sontacchi Criterion}\footnote{A monicker of my own invention.}]
\label{thm:Bohm-Sontacchi}Let\index{Böhm-Sontacchi criterion}\index{diophantine equation}
$x\in\mathbb{Z}$ be a periodic point of the Collatz map $C$ (\ref{eq:Collatz Map}).
Then, there are integers $m,n\geq1$ and a strictly increasing sequence
of positive integers $b_{1}<b_{2}<\cdots<b_{n}$ so that: 
\begin{equation}
x=\frac{\sum_{k=1}^{n}2^{m-b_{k}-1}3^{k-1}}{2^{m}-3^{n}}\label{eq:The Bohm-Sontacchi Criterion}
\end{equation}
\end{thm}
\vphantom{}

(\ref{eq:The Bohm-Sontacchi Criterion}) is arguably the clearest
expression of how the Collatz Conjecture reaches out of the gyre of
recreational mathematics and up toward more refined spheres of study.
Not only is (\ref{eq:The Bohm-Sontacchi Criterion}) an exponential
diophantine equation\textemdash an area of number theory as classical
as it is difficult\textemdash but is, worse yet, a \emph{family }of
exponential diophantine equations. Especially noteworthy is the denominator
term, $2^{m}-3^{n}$, which evidences (\ref{eq:The Bohm-Sontacchi Criterion})
as falling under the jurisdiction of of transcendental number theory\footnote{A slightly more in-depth explanation of this connection is given beginning
from page \pageref{subsec:Baker,-Catalan,-and}.}, in which the study of lower bounds for quantities such as $\left|2^{m}-3^{n}\right|$
is of the utmost importance.

The most significant result of the present section\textemdash and,
one of this dissertation's three principal achievements\textemdash is
the \textbf{Correspondence Principle }(proved in Subsection \ref{subsec:2.2.3 The-Correspondence-Principle}),
which significantly generalizes the \textbf{Böhm-Sontacchi Criterion},
and in two respects. First, it shows that the integers $m$, $n$,
and the $b_{k}$s can all be parameterized as functions of a single
variable in $\mathbb{Z}$. Secondly, it generalizes the criterion
to apply to a much larger family of Hydra maps. Not only that, but,
we will be able to do the same for Hydra maps on $\mathbb{Z}^{d}$
in Chapter 5! It should be noted that a generalized form of the Criterion
is given in Matthews' slides \cite{Matthews' slides}, but the version
there is not as broadly minded as mine, and lacks the single-variable
parameterization of the integer parameters in (\ref{eq:The Bohm-Sontacchi Criterion})
produced my $\chi_{H}$ formalism.

The ideas needed to make all this work emerged from considering the
$p$-ity vector associated to a given integer $n$ under the particular
$p$-Hydra map $H$ we are studying. To illustrate this, let us consider
the Shortened Collatz map. We will denote the two branches of $T_{3}$
($x/2$ and $\left(3x+1\right)/2$) by $H_{0}\left(x\right)$ and
$H_{1}\left(x\right)$, respectively.

Letting $n\geq1$, note that as we successively apply $T_{3}$ to
$n$, the iterates $T_{3}\left(n\right),T_{3}\left(T_{3}\left(n\right)\right),\ldots$
can be expressed in terms of compositions of $H_{0}$ and $H_{1}$:
\begin{example}
\ 
\begin{align*}
T_{3}\left(1\right) & =2=\frac{3\left(1\right)+1}{2}=H_{1}\left(1\right)\\
T_{3}\left(T_{3}\left(1\right)\right) & =1=\frac{T_{3}\left(1\right)}{2}=H_{0}\left(T_{3}\left(1\right)\right)=H_{0}\left(H_{1}\left(1\right)\right)\\
 & \vdots
\end{align*}
\end{example}
\vphantom{}

Given integers $m,n\geq1$, when we express the $m$th iterate of
$n$ under $T_{3}$: 
\begin{equation}
T_{3}^{\circ m}\left(n\right)=\underbrace{\left(T_{3}\circ T_{3}\circ\cdots\circ T_{3}\right)}_{m\textrm{ times}}\left(n\right)
\end{equation}
as a composition of $H_{0}$ and $H_{1}$, the sequence of the subscripts
$0$ and $1$ that occur in this sequence of maps is clearly equivalent
to the parity vector generated by $T_{3}$ over $n$. In this ``classical''
view of a parity vector, the parity vector is \emph{subordinated }to
$n$, being obtained \emph{from }$n$. My approach is to turn this
relationship on its head: what happens if we reverse it? That is,
what happens if we let the \emph{parity vector} (or, more generally,
$p$-ity vector) determine $n$?

For $T_{3}$, we can see this reversal in action by considering the
kinds of composition sequences of $H_{0}$ and $H_{1}$ that occur
when we are given an integer $x$ which is a periodic point $x$ of
$T_{3}$. For such an $x$, there is an integer $m\geq1$ so that
$T_{3}^{\circ m}\left(x\right)=x$, which is to say, there is some
length-$m$ composition sequence of $H_{0}$ and $H_{1}$ that has
$x$ as a fixed point. Instead of starting with $x$ as the given,
however, we treat it as an unknown to be solved for.
\begin{example}
Suppose $x$ is fixed by the sequence: 
\begin{equation}
x=H_{1}\left(H_{1}\left(H_{0}\left(H_{1}\left(H_{0}\left(x\right)\right)\right)\right)\right)
\end{equation}
Expanding the right-hand side, we obtain: 
\begin{equation}
x=\frac{1}{2}\left(3\left[\frac{1}{2}\left(3\left[\frac{1}{2}\left(\frac{3\left(\frac{x}{2}\right)+1}{2}\right)\right]+1\right)\right]+1\right)=\frac{3^{3}}{2^{5}}x+\frac{3^{2}}{2^{4}}+\frac{3}{2^{2}}+\frac{1}{2}\label{eq:Bohm-Sontacci Example}
\end{equation}
and hence: 
\[
x=\frac{3^{2}\cdot2+3\cdot2^{3}+2^{4}}{2^{5}-3^{3}}=\frac{58}{5}
\]
Since $58/5$ is not an integer, we can conclude that there is no
integer fixed point of $T_{3}$ whose iterates have the parity sequence
generated by the string $0,1,0,1,1$.
\end{example}
\vphantom{}

In general, to any \emph{finite }sequence $j_{1},j_{2},j_{3},\ldots,j_{N}$
consisting of $0$s and $1$s, we can solve the equation: 
\begin{equation}
x=\left(H_{j_{1}}\circ H_{j_{2}}\circ\cdots\circ H_{j_{N}}\right)\left(x\right)\label{eq:Periodic point set up for T3}
\end{equation}
for $x$, and, in doing so, obtain a map\footnote{When we generalize to consider an arbitrary $p$-Hydra map $H$ instead
of $T_{3}$, we will obtain a particular $x$ for each $H$.} from the space of all such sequences to the set of rational numbers.
That this works is because $H_{0}$ and $H_{1}$ are bijections of
$\mathbb{Q}$. Thus, any finite composition sequence $S$ of these
two branches is also a bijection, and as such, there is a unique $x\in\mathbb{Q}$
so that $x=S\left(x\right)$, which can be determined in the manner
shown above.

To put it another way, we can \textbf{\emph{parameterize}} the quantities
$m$, $n$, $x$, and $b_{k}$ which appear in the Böhm-Sontacchi
Criterion (\ref{eq:Bohm-Sontacci Example}) in terms of the sequence
$j_{1},\ldots,j_{N}$ of branches of $T_{3}$ which we applied to
map $x$ to itself. Viewing (\ref{eq:Bohm-Sontacci Example}) as a
particular case of this formula, note that in the above example, the
values of $m$, $n$, the $b_{k}$s in (\ref{eq:Bohm-Sontacci Example}),
and\textemdash crucially\textemdash $x$ itself were determined by
our choice of the composition sequence: 
\begin{equation}
H_{1}\circ H_{1}\circ H_{0}\circ H_{1}\circ H_{0}
\end{equation}
and our assumption that this sequence mapped $x$ to $x$. In this
way, we will be able to deal with all the possible variants of (\ref{eq:Bohm-Sontacci Example})
in the context of a single object\textemdash to be called $\chi_{3}$.

Just as we used sequences $j_{1},\ldots,j_{N}$ in $\left\{ 0,1\right\} $
to deal with the $2$-Hydra map $T_{3}$, for the case of a general
$p$-Hydra map, our sequences $j_{1},\ldots,j_{N}$ will consist of
integers in the set $\left\{ 0,\ldots,p-1\right\} $, with each sequence
corresponding to a particular composition sequences of the $p$ distinct
branches of the $p$-Hydra map $H$ under consideration. In doing
so, we will realize these expressions as specific cases of a more
general functional equation satisfied the function I denote by $\chi_{H}$.
We will initially define $\chi_{H}$ as a rational-valued function
of finite sequences $j_{1},j_{2},\ldots$. However, by a standard
bit of identifications, we will re-contextualize $\chi_{H}$ as a
rational function on the non-negative integers $\mathbb{N}_{0}$.
By placing certain qualitative conditions on $H$, we can ensure that
$\chi_{H}$ will remain well-defined when we extend its inputs from
$\mathbb{N}_{0}$ to $\mathbb{Z}_{p}$, thereby allowing us to interpolate
$\chi_{H}$ to a $q$-adic valued function over $\mathbb{Z}_{p}$
for an appropriate choice of primes $p$ and $q$.

\subsection{\label{subsec:2.2.1 Notation-and-Preliminary}Notation and Preliminary
Definitions}

Fix a $p$-Hydra map $H$, and \textbf{suppose that $b_{0}=0$} so
as to guarantee that $H\left(0\right)=0$; note also that this then
makes $H_{0}\left(0\right)=0$. In our study, \emph{it is imperative
that }\textbf{\emph{$b_{0}=0$}}.

We begin by introducing formalism\textemdash \textbf{strings}\textemdash to
used to denote sequences of the numbers $\left\{ 0,\ldots,p-1\right\} $.
\begin{defn}[\textbf{Strings}]
We write \nomenclature{$\textrm{String}\left(p\right)$}{ }$\textrm{String}\left(p\right)$
to denote the set of all finite sequences whose entries belong to
the set $\left\{ 0,\ldots,p-1\right\} $; we call such sequences \textbf{strings}.
We also include the empty set ($\varnothing$) as an element of $\textrm{String}\left(p\right)$,
and refer to it as the \textbf{empty string}. Arbitrary strings are
usually denoted by $\mathbf{j}=\left(j_{1},\ldots,j_{n}\right)$,
where $n$ is the \textbf{length }of the string, denoted by \nomenclature{$\left|\mathbf{j}\right|$}{ }$\left|\mathbf{j}\right|$.
In this manner, any $\mathbf{j}$ can be written as $\mathbf{j}=\left(j_{1},\ldots,j_{\left|\mathbf{j}\right|}\right)$.
We define the empty string as having length $0$, making it the unique
string of length $0$. We say the string $\mathbf{j}$ is \textbf{non-zero
}if $\mathbf{j}$ contains at least one non-zero entry.

\vphantom{}

We write \nomenclature{$\textrm{String}_{\infty}\left(p\right)$}{ }$\textrm{String}_{\infty}\left(p\right)$
to denote the set of all all sequences (finite \emph{or }infinite)
whose entries belong to the set $\left\{ 0,\ldots,p-1\right\} $.
A \textbf{finite }string is one with finite length; on the other hand,
any string in $\textrm{String}_{\infty}\left(p\right)$ which is not
finite is said to be \textbf{infinite}, and its length is defined
to be $+\infty$. We also include $\varnothing$ in $\textrm{String}_{\infty}\left(p\right)$,
again referred to as the empty string.
\end{defn}
\begin{defn}
Given any $\mathbf{j}\in\textrm{String}\left(p\right)$, we define
the \index{composition sequence}\textbf{composition sequence }$H_{\mathbf{j}}:\mathbb{Q}\rightarrow\mathbb{Q}$
as the affine linear map: \nomenclature{$H_{\mathbf{j}}\left(x\right)$}{$\overset{\textrm{def}}{=}\left(H_{j_{1}}\circ\cdots\circ H_{j_{\left|\mathbf{j}\right|}}\right)\left(x\right)$}
\begin{equation}
H_{\mathbf{j}}\left(x\right)\overset{\textrm{def}}{=}\left(H_{j_{1}}\circ\cdots\circ H_{j_{\left|\mathbf{j}\right|}}\right)\left(x\right),\textrm{ }\forall\mathbf{j}\in\textrm{String}\left(p\right)\label{eq:Def of composition sequence}
\end{equation}
When writing by hand, one can use vector notation: $\vec{j}$ instead
of $\mathbf{j}$. 
\end{defn}
\begin{rem}
Note the near-equivalence of $\mathbf{j}$ with $p$-ity vectors.
Indeed, given any $m,n\in\mathbb{N}_{1}$, there exists a unique $\mathbf{j}\in\textrm{String}\left(p\right)$
of length $m$ so that $H^{\circ m}\left(n\right)=H_{\mathbf{j}}\left(n\right)$.
Said $\mathbf{j}$ is the $p$-ity vector for the first $m$ iterates
of $n$ under $H$, albeit written in reverse order. While, admittedly,
this reverse-order convention is somewhat unnatural, it does have
one noteworthy advantage. In a moment, we will identify $\mathbf{j}$
with a $p$-adic integer $\mathfrak{z}$ by viewing the entries of
$\mathbf{j}$ as the $p$-adic digits of $\mathfrak{z}$, written
left-to-right in order of increasing powers of $p$. Our ``reverse''
indexing convention guarantees that the left-to-right order in which
we write the $H_{j_{k}}$s, the left-to-right order in which we write
the entries of $\mathbf{j}$, and the left-to-write order in which
we write the $p$-adic digits of the aforementioned $\mathfrak{z}$
are all the same. 
\end{rem}
\vphantom{}

Strings will play a vital role in this chapter. Not only do they simplify
many arguments, they also serve as a go-between for the dynamical-system
motivations (considering arbitrary composition sequences of $H_{j}$s)
and the equivalent reformulations in terms of non-negative and/or
$p$-adic integers. As was mentioned above and is elaborated upon
below, we will identify elements of $\textrm{String}\left(p\right)$
with the sequences of base $p$ digits of non-negative integers. Likewise,
we will identify elements of $\textrm{String}_{\infty}\left(p\right)$
with the sequences of $p$-adic digits of $p$-adic integers. We do
this by defining the obvious map for converting strings into ($p$-adic)
integers: 
\begin{defn}
We write \nomenclature{$\textrm{Dig}\textrm{Sum}_{p}$}{ }$\textrm{Dig}\textrm{Sum}_{p}:\textrm{String}_{\infty}\left(p\right)\rightarrow\mathbb{Z}_{p}$
by: 
\begin{equation}
\textrm{DigSum}_{p}\left(\mathbf{j}\right)\overset{\textrm{def}}{=}\sum_{k=1}^{\left|\mathbf{j}\right|}j_{k}p^{k-1}\label{eq:Definition of DigSum_p of bold j}
\end{equation}
\end{defn}
\begin{defn}[\textbf{Identifying Strings with Numbers}]
Given $\mathbf{j}\in\textrm{String}_{\infty}\left(p\right)$ and
$\mathfrak{z}\in\mathbb{Z}_{p}$, we say $\mathbf{j}$ \textbf{represents
}(or \textbf{is} \textbf{associated to})\textbf{ }$\mathfrak{z}$
(and vice-versa), written $\mathbf{j}\sim\mathfrak{z}$ or $\mathfrak{z}\sim\mathbf{j}$
whenever $\mathbf{j}$ is the sequence of the $p$-adic digits of
$n$; that is: 
\begin{equation}
\mathfrak{z}\sim\mathbf{j}\Leftrightarrow\mathbf{j}\sim\mathfrak{z}\Leftrightarrow\mathfrak{z}=j_{1}+j_{2}p+j_{3}p^{2}+\cdots\label{eq:Definition of n-bold-j correspondence.}
\end{equation}
equivalently: 
\begin{equation}
\mathbf{j}\sim\mathfrak{z}\Leftrightarrow\textrm{DigSum}_{p}\left(\mathbf{j}\right)=\mathfrak{z}
\end{equation}
As defined, $\sim$ is then an equivalence relation on $\textrm{String}\left(p\right)$
and $\textrm{String}_{\infty}\left(p\right)$. We write: 
\begin{equation}
\mathbf{i}\sim\mathbf{j}\Leftrightarrow D_{p}\left(\mathbf{i}\right)=D_{p}\left(\mathbf{j}\right)\label{eq:Definition of string rho equivalence relation, rational integer version}
\end{equation}
and we have that $\mathbf{i}\sim\mathbf{j}$ if and only if both $\mathbf{i}$
and $\mathbf{j}$ represent the same $p$-adic integer.

Note that in both $\textrm{String}\left(p\right)$ and $\textrm{String}_{\infty}\left(p\right)$,
the shortest string representing the number $0$ is then the empty
string.

Finally, we write \nomenclature{$\textrm{String}\left(p\right)/\sim$}{ }$\textrm{String}\left(p\right)/\sim$
and \nomenclature{$\textrm{String}_{\infty}\left(p\right)/\sim$}{ }$\textrm{String}_{\infty}\left(p\right)/\sim$
to denote the set of equivalence classes of $\textrm{String}\left(p\right)$
and $\textrm{String}_{\infty}\left(p\right)$ under this equivalence
relation.
\end{defn}
\begin{prop}
\label{prop:string number equivalence}\ 

\vphantom{}

I. $\textrm{String}\left(p\right)/\sim$ and $\textrm{String}_{\infty}\left(p\right)/\sim$
are in bijective correspondences with $\mathbb{N}_{0}$ and $\mathbb{Z}_{p}$,
respectively, by way of the map $\textrm{Dig}\textrm{Sum}_{p}$.

\vphantom{}

II. Two finite strings $\mathbf{i}$ and $\mathbf{j}$ satisfy $\mathbf{i}\sim\mathbf{j}$
if and only if either:

\vphantom{}

a) $\mathbf{i}$ can be obtained from $\mathbf{j}$ by adding finitely
many $0$s to the right of $\mathbf{j}$.

\vphantom{}

b) $\mathbf{j}$ can be obtained from $\mathbf{i}$ by adding finitely
many $0$s to the right of $\mathbf{i}$.

\vphantom{}

III. A finite string $\mathbf{i}$ and an infinite string $\mathbf{j}$
satisfy $\mathbf{i}\sim\mathbf{j}$ if and only if $\mathbf{j}$ contains
finitely many non-zero entries. 
\end{prop}
Proof: Immediate from the definitions.

Q.E.D.

\vphantom{}

The principal utility of this formalism is the way in which composition
sequences $H_{\mathbf{j}}$ interact with concatenations of strings. 
\begin{defn}[\textbf{Concatenation}]
We introduce the \textbf{concatenation operation}\index{concatenation!operation}
$\wedge:\textrm{String}\left(p\right)\times\textrm{String}_{\infty}\left(p\right)\rightarrow\textrm{String}_{\infty}\left(p\right)$,
defined by: 
\begin{equation}
\mathbf{i}\wedge\mathbf{j}=\left(i_{1},\ldots,i_{\left|\mathbf{i}\right|}\right)\wedge\left(j_{1},\ldots,j_{\left|\mathbf{j}\right|}\right)\overset{\textrm{def}}{=}\left(i_{1},\ldots,i_{\left|\mathbf{i}\right|},j_{1},\ldots,j_{\left|\mathbf{j}\right|}\right)\label{eq:Definition of Concatenation}
\end{equation}
with the definition being modified in the obvious way to when $\mathbf{j}$
is of infinite length.

\vphantom{}

Additionally, for any integer $m\geq1$ and any finite string $\mathbf{j}$,
we write $\mathbf{j}^{\wedge m}$ to denote the concatenation of $m$
copies of $\mathbf{j}$: 
\begin{equation}
\mathbf{j}^{\wedge m}\overset{\textrm{def}}{=}\left(\underbrace{j_{1},\ldots,j_{\left|\mathbf{j}\right|},j_{1},\ldots,j_{\left|\mathbf{j}\right|},\ldots,j_{1},\ldots,j_{\left|\mathbf{j}\right|}}_{m\textrm{ times}}\right)\label{eq:Definition of concatenation exponentiation}
\end{equation}
\end{defn}
\vphantom{}

In the next subsection, the proposition below will be the foundation
for various useful functional equations: 
\begin{prop}
\label{prop:H string formalism}\ 
\begin{equation}
H_{\mathbf{i}\wedge\mathbf{j}}\left(x\right)=H_{\mathbf{i}}\left(H_{\mathbf{j}}\left(x\right)\right),\textrm{ }\forall x\in\mathbb{R},\textrm{ }\forall\mathbf{i},\mathbf{j}\in\textrm{String}\left(p\right)\label{eq:H string formalism}
\end{equation}
\end{prop}
Proof: 
\begin{equation}
H_{\mathbf{i}\wedge\mathbf{j}}\left(x\right)=H_{i_{1},\ldots,i_{\left|\mathbf{i}\right|},j_{1},\ldots,j_{\left|\mathbf{j}\right|}}\left(x\right)=\left(H_{i_{1},\ldots,i_{\left|\mathbf{i}\right|}}\circ H_{j_{1},\ldots,j_{\left|\mathbf{j}\right|}}\right)\left(x\right)=H_{\mathbf{i}}\left(H_{\mathbf{j}}\left(x\right)\right)
\end{equation}

Q.E.D.

\vphantom{}

In working with strings and the ($p$-adic) integers they represent,
the following functions will be useful to help bridge string-world
and number-world. 
\begin{defn}[$\lambda_{p}$ \textbf{and} $\#_{p:0},\#_{p:1},\ldots,\#_{p:p-1}$]
We write $\lambda_{p}:\mathbb{N}_{0}\rightarrow\mathbb{N}_{0}$ to
denote the function:\nomenclature{$\lambda_{p}\left(n\right)$}{$\overset{\textrm{def}}{=}\left\lceil \log_{p}\left(n+1\right)\right\rceil$ \nopageref }
\begin{equation}
\lambda_{p}\left(n\right)\overset{\textrm{def}}{=}\left\lceil \log_{p}\left(n+1\right)\right\rceil ,\textrm{ }\forall n\in\mathbb{N}_{0}\label{eq:definition of lambda rho}
\end{equation}
which gives the number of $p$-adic digits of $n$. Every $n\in\mathbb{N}_{0}$
can be uniquely written as: 
\begin{equation}
n=c_{0}+c_{1}p+\cdots+c_{\lambda_{p}\left(n\right)-1}p^{\lambda_{p}\left(n\right)-1}
\end{equation}
where the $c_{j}$s are integers in $\left\{ 0,\ldots,p-1\right\} $.
Additionally, note that: 
\begin{equation}
\left\lceil \log_{p}\left(n+1\right)\right\rceil =\left\lfloor \log_{p}n\right\rfloor +1,\textrm{ }\forall n\in\mathbb{N}_{1}\label{eq:Floor and Ceiling expressions for lambda_rho}
\end{equation}
and that $\lambda_{p}\left(n\right)\leq\left|\mathbf{j}\right|$ is
satisfied for any string $\mathbf{j}$ representing $n$.

Next, for each $n\geq1$ and each $k\in\left\{ 0,\ldots,p-1\right\} $,
we write $\#_{p:k}\left(n\right)$\nomenclature{$\#_{p:k}\left(n\right)$}{the number of $k$s in the $p$-adic digits of $n$ }
to denote the number of $k$s present in the $p$-adic expansion of
$n$. We define $\#_{p:k}\left(0\right)$ to be $0$ for all $k$.
In a minor abuse of notation, we also use $\#_{p:k}$ to denote the
number of $k$s which occur in a given string: 
\begin{equation}
\#_{p:k}\left(\mathbf{j}\right)\overset{\textrm{def}}{=}\textrm{number of }k\textrm{s in }\mathbf{j},\textrm{ }\forall k\in\mathbb{Z}/p\mathbb{Z},\textrm{ }\forall\mathbf{j}\in\textrm{String}\left(p\right)\label{eq:Definition of number of ks in rho adic digits of bold j}
\end{equation}
\end{defn}
\begin{rem}
$\lambda_{2}\left(n\right)$ is what computer scientists called the
number of \textbf{bits }in $n$. 
\end{rem}
\begin{example}
For: 
\begin{align*}
n & =1\cdot6^{0}+2\cdot6^{1}+0\cdot6^{2}+5\cdot6^{3}+0\cdot6^{4}+2\cdot6^{5}\\
 & =1+2\cdot6+5\cdot6^{3}+2\cdot6^{5}\\
 & =16645
\end{align*}
we have: 
\begin{align*}
\#_{6:0}\left(16645\right) & =2\\
\#_{6:1}\left(16645\right) & =1\\
\#_{6:2}\left(16645\right) & =2\\
\#_{6:3}\left(16645\right) & =0\\
\#_{6:4}\left(16645\right) & =0\\
\#_{6:5}\left(16645\right) & =1
\end{align*}
\end{example}
\vphantom{}

All important functions in this dissertation satisfy equally important
functional equations. Here are the functional equations for $\lambda_{p}$
and the $\#_{p:k}$s. 
\begin{prop}[\textbf{Functional Equations for $\lambda_{p}$ and the $\#_{p:k}$s}]
\ 

\vphantom{}

I. For $k\in\left\{ 1,\ldots,p-1\right\} $:
\begin{align}
\#_{p:k}\left(p^{n}a+b\right) & =\#_{p:k}\left(a\right)+\#_{p:k}\left(b\right),\textrm{ }\forall a\in\mathbb{N}_{0},\textrm{ }\forall n\in\mathbb{N}_{1},\textrm{ }\forall b\in\left\{ 0,\ldots,p^{n}-1\right\} \label{eq:number-symbol functional equations}
\end{align}
For $k=0$, we have:
\begin{equation}
\#_{p:0}\left(p^{n}a+b\right)=\#_{p:0}\left(a\right)+\#_{p:0}\left(b\right)+n-\lambda_{p}\left(b\right),\textrm{ }\forall a\in\mathbb{N}_{1},\textrm{ }\forall n\in\mathbb{N}_{1},\textrm{ }\forall b\in\left\{ 0,\ldots,p^{n}-1\right\} \label{eq:number-symbol functional equations, k is 0}
\end{equation}

\vphantom{}

II. 
\begin{align}
\lambda_{p}\left(p^{n}a+b\right) & =\lambda_{p}\left(a\right)+n,\textrm{ }\forall a\in\mathbb{N}_{1},\textrm{ }\forall n\in\mathbb{N}_{1},\textrm{ }\forall b\in\left\{ 0,\ldots,p^{n}-1\right\} \label{eq:lambda functional equations}
\end{align}

\vphantom{}

III. 
\begin{equation}
\sum_{k=0}^{p-1}\#_{p:k}\left(n\right)=\lambda_{p}\left(n\right)\label{eq:Relation between lambda_rho and the number of rho adic digits}
\end{equation}
\end{prop}
Proof: A straightforward computation. (\ref{eq:number-symbol functional equations, k is 0})
follows from using (\ref{eq:Relation between lambda_rho and the number of rho adic digits})
to write:
\[
\#_{p:0}\left(n\right)=\lambda_{p}\left(n\right)-\sum_{k=1}^{p-1}\#_{p:k}\left(n\right)
\]
and then applying (\ref{eq:number-symbol functional equations}) and
(\ref{eq:lambda functional equations}) to compute $\#_{p:0}\left(p^{n}a+b\right)$.

Q.E.D.

\subsection{\label{subsec:2.2.2 Construction-of}Construction of $\chi_{H}$}

Since the $H_{j}$s are affine linear maps, so is any composition
sequence $H_{\mathbf{j}}$. As such, for any $\mathbf{j}\in\textrm{String}\left(p\right)$
we can write:

\begin{equation}
H_{\mathbf{j}}\left(x\right)=H_{\mathbf{j}}^{\prime}\left(0\right)x+H_{\mathbf{j}}\left(0\right),\textrm{ }\forall\mathbf{j}\in\textrm{String}\left(p\right)\label{eq:ax+b formula for h_bold_j}
\end{equation}
If $\mathbf{j}$ is a string for which $H_{\mathbf{j}}\left(x\right)=x$,
this becomes: 
\begin{equation}
x=H_{\mathbf{j}}^{\prime}\left(0\right)x+H_{\mathbf{j}}\left(0\right)\label{eq:affine formula for little x}
\end{equation}
from which we obtain: 
\begin{equation}
x=\frac{H_{\mathbf{j}}\left(0\right)}{1-H_{\mathbf{j}}^{\prime}\left(0\right)}\label{eq:Formula for little x in terms of bold j}
\end{equation}
Note that the map which sends a string $\mathbf{j}$ to the rational
number $x$ satisfying $H_{\mathbf{j}}\left(x\right)=x$ can then
be written as: 
\begin{equation}
\mathbf{j}\mapsto\frac{H_{\mathbf{j}}\left(0\right)}{1-H_{\mathbf{j}}^{\prime}\left(0\right)}\label{eq:Formula for X}
\end{equation}

While we could analyze this map directly, the resultant function is
ill-behaved because of the denominator term. Consequently, instead
of studying (\ref{eq:Formula for X}) directly, we will focus on the
map $\mathbf{j}\mapsto H_{\mathbf{j}}\left(0\right)$, which we denote
by $\chi_{H}$. 
\begin{defn}[$\chi_{H}$]
\label{def:Chi_H on N_0 in strings}Let $H$ be a $p$-Hydra map
that fixes $0$. Then, for any $\mathbf{j}\in\textrm{String}\left(p\right)$,
we write: 
\begin{equation}
\chi_{H}\left(\mathbf{j}\right)\overset{\textrm{def}}{=}H_{\mathbf{j}}\left(0\right)\label{eq:Definition of Chi_H of bold j}
\end{equation}
Identifying $H_{\varnothing}$ (the composition sequence associated
to the empty string) with the identity map, we have that $\chi_{H}\left(\varnothing\right)=0$. 
\end{defn}
\begin{prop}
\label{prop:Chi_H is well defined mod twiddle}$\chi_{H}$ is well-defined
on $\textrm{String}\left(p\right)/\sim$. That is, $\chi_{H}\left(\mathbf{i}\right)=\chi_{H}\left(\mathbf{j}\right)$
for all $\mathbf{i},\mathbf{j}\in\textrm{String}\left(p\right)$ for
which $\mathbf{i}\sim\mathbf{j}$. 
\end{prop}
Proof: Let $\mathbf{j}\in\textrm{String}\left(p\right)$ be any non-empty
string. Then, by \textbf{Proposition \ref{prop:H string formalism}}:
\begin{equation}
\chi_{H}\left(\mathbf{j}\wedge0\right)=H_{\mathbf{j}\wedge0}\left(0\right)=H_{\mathbf{j}}\left(H_{0}\left(0\right)\right)=H_{\mathbf{j}}\left(0\right)=\chi_{H}\left(\mathbf{j}\right)
\end{equation}
since $H_{0}\left(0\right)=0$. \textbf{Proposition \ref{prop:string number equivalence}},
we know that given any two equivalent finite strings $\mathbf{i}\sim\mathbf{j}$,
we can obtain one of the two by concatenating finitely many $0$s
to the right of the other. Without loss of generality, suppose that
$\mathbf{i}=\mathbf{j}\wedge0^{\wedge n}$, where $0^{\wedge n}$
is a string of $n$ consecutive $0$s. Then: 
\begin{equation}
\chi_{H}\left(\mathbf{i}\right)=\chi_{H}\left(\mathbf{j}\wedge0^{\wedge n}\right)=\chi_{H}\left(\mathbf{j}\wedge0^{\wedge\left(n-1\right)}\right)=\cdots=\chi_{H}\left(\mathbf{j}\right)
\end{equation}
Hence, $\chi_{H}$ is well defined on $\textrm{String}\left(p\right)/\sim$.

Q.E.D. 
\begin{lem}[$\chi_{H}$ \textbf{on} $\mathbb{N}_{0}$]
\index{chi_{H}@$\chi_{H}$}\label{lem:Chi_H on N_0}We can realize
$\chi_{H}$ as a function $\mathbb{N}_{0}\rightarrow\mathbb{Q}$ by
defining: 
\begin{equation}
\chi_{H}\left(n\right)\overset{\textrm{def}}{=}\chi_{H}\left(\mathbf{j}\right)=H_{\mathbf{j}}\left(0\right)\label{eq:Definition of Chi_H of n}
\end{equation}
where $\mathbf{j}\in\textrm{String}\left(p\right)$ is any string
representing $n$; we define $\chi_{H}\left(0\right)$ to be $\chi_{H}\left(\varnothing\right)$
(the empty string), which is just $0$. 
\end{lem}
Proof: By \textbf{Proposition \ref{prop:string number equivalence}},
$\textrm{String}\left(p\right)/\sim$ is in a bijection with $\mathbb{N}_{0}$.
\textbf{Proposition \ref{prop:Chi_H is well defined mod twiddle}
}tells us that $\chi_{H}$ is well-defined on $\textrm{String}\left(p\right)/\sim$,
so, by using the aforementioned bijection to identify $\mathbb{N}_{0}$
with $\textrm{String}\left(p\right)/\sim$, the rule $\chi_{H}\left(n\right)\overset{\textrm{def}}{=}\chi_{H}\left(\mathbf{j}\right)$
is well-defined, since it is independent of which $\mathbf{j}\in\textrm{String}\left(p\right)$
we choose to represent $n$.

Q.E.D. 
\begin{defn}[$M_{H}$]
\index{$M_{H}$}\label{def:M_H on N_0 in strings}Let $H$ be a $p$-Hydra
map. Then, we define $M_{H}:\mathbb{N}_{0}\rightarrow\mathbb{Q}$
by: 
\begin{equation}
M_{H}\left(n\right)\overset{\textrm{def}}{=}M_{H}\left(\mathbf{j}\right)\overset{\textrm{def}}{=}H_{\mathbf{j}}^{\prime}\left(0\right),\textrm{ }\forall n\geq1\label{eq:Definition of M_H}
\end{equation}
where $\mathbf{j}\in\textrm{String}\left(p\right)$ is the shortest
element of $\textrm{String}\left(p\right)$ representing $n$. We
define $M_{H}\left(0\right)$ to be $1$.
\end{defn}
\begin{example}
\emph{Unlike} $\chi_{H}$, observe that $M_{H}$ is \emph{not well-defined
over} $\textrm{String}\left(p\right)/\sim$ there can be $H$ and
distinct strings $\mathbf{i},\mathbf{j}\in\textrm{String}\left(p\right)$
so that $\mathbf{i}\sim\mathbf{j}$, yet $M_{H}\left(\mathbf{i}\right)\neq M_{H}\left(\mathbf{j}\right)$.

For example, $H=T_{3}$, the strings $\left(0,1\right)$ and $\left(0,1,0\right)$
both represent the integer $2$. Nevertheless: 
\begin{align*}
H_{0,1}\left(x\right) & =H_{0}\left(H_{1}\left(x\right)\right)=\frac{3}{4}x+\frac{1}{4}\\
H_{0,1,0}\left(x\right) & =H_{0}\left(H_{1}\left(H_{0}\left(x\right)\right)\right)=\frac{3}{8}x+\frac{1}{4}
\end{align*}
Thus, $\left(0,1\right)\sim\left(0,1,0\right)$, but $M_{T_{3}}\left(0,1\right)=3/4$,
while $M_{T_{3}}\left(0,1,0\right)=3/8$.

With these definitions in place, we can now use the string formalism
to establish vital string-based functional equation identities for
$\chi_{H}$ and $M_{H}$. I will refer to these identities as \textbf{concatenation
identities}\index{concatenation!identities}. 
\end{example}
\begin{prop}[\textbf{An expression for $H_{\mathbf{j}}$}]
\label{prop:H_boldj in terms of M_H and Chi_H}Let $H$ be a $p$-Hydra
map. Then: 
\begin{equation}
H_{\mathbf{j}}\left(x\right)=M_{H}\left(\mathbf{j}\right)x+\chi_{H}\left(\mathbf{j}\right),\textrm{ }\forall\mathbf{j}\in\textrm{String}\left(p\right),\textrm{ }\forall x\in\mathbb{Q}\label{eq:Formula for composition sequences of H}
\end{equation}
\end{prop}
Proof: Since the $H_{j}$s are maps of the form $ax+b$, any composition
sequence of the $H_{j}$s will also be of that form. Consequently,
the venerable slope-intercept formula yields: 
\[
H_{\mathbf{j}}\left(x\right)=H_{\mathbf{j}}^{\prime}\left(0\right)x+H_{\mathbf{j}}\left(0\right)\overset{\textrm{def}}{=}M_{H}\left(\mathbf{j}\right)x+\chi_{H}\left(\mathbf{j}\right)
\]

Q.E.D.

\vphantom{}

The lemmata below give the concatenation identities for $\chi_{H}$
and $M_{H}$. 
\begin{lem}[$\chi_{H}$ \textbf{Concatenation Identity}]
\label{lem:Chi_H concatenation identity}\ 

\begin{equation}
\chi_{H}\left(\mathbf{i}\wedge\mathbf{j}\right)=H_{\mathbf{i}}\left(\chi_{H}\left(\mathbf{j}\right)\right),\textrm{ }\forall\mathbf{i},\mathbf{j}\in\textrm{String}\left(p\right)\label{eq:Chi_H concatenation identity}
\end{equation}
That is to say, for $\mathbf{i}=\left(i_{1},i_{2},\ldots,i_{m}\right)\in\textrm{String}\left(p\right)$,
where $m\geq1$, we have: 
\begin{equation}
\chi_{H}\left(i_{1}+i_{2}p+\cdots+i_{m}p^{m-1}+p^{m}n\right)=H_{\mathbf{i}}\left(\chi_{H}\left(n\right)\right),\textrm{ }\forall n\in\mathbb{N}_{0}\label{eq:Chi_H concatenation identity, alternate version}
\end{equation}
\end{lem}
Proof: Letting $\mathbf{i}$, $\mathbf{j}$, and $x$ be arbitrary,
we have that: 
\begin{align*}
H_{\mathbf{i}\wedge\mathbf{j}}\left(x\right) & =H_{\mathbf{i}\wedge\mathbf{j}}^{\prime}\left(0\right)x+H_{\mathbf{i}\wedge\mathbf{j}}\left(0\right)\\
 & =H_{\mathbf{i}\wedge\mathbf{j}}^{\prime}\left(0\right)x+H_{\mathbf{i}}\left(H_{\mathbf{j}}\left(0\right)\right)
\end{align*}
Setting $x=0$ yields: 
\[
\underbrace{H_{\mathbf{i}\wedge\mathbf{j}}\left(0\right)}_{\chi_{H}\left(\mathbf{i}\wedge\mathbf{j}\right)}=\underbrace{H_{\mathbf{i}}\left(H_{\mathbf{j}}\left(0\right)\right)}_{H_{\mathbf{j}}\left(\chi_{H}\left(\mathbf{j}\right)\right)}
\]

Q.E.D. 
\begin{prop}[$M_{H}$ \textbf{Concatenation Identity}]
\label{prop:M_H concatenation identity}\ 
\begin{equation}
M_{H}\left(\mathbf{i}\wedge\mathbf{j}\right)=M_{H}\left(\mathbf{i}\right)M_{H}\left(\mathbf{j}\right),\textrm{ }\forall\mathbf{i},\mathbf{j}\in\textrm{String}\left(p\right),\textrm{ }\left|\mathbf{i}\right|,\left|\mathbf{j}\right|\geq1\label{eq:M_H concatenation identity}
\end{equation}
That is to say, for $\mathbf{i}=\left(i_{1},i_{2},\ldots,i_{m}\right)\in\textrm{String}\left(p\right)$,
where $m\geq1$, we have: 
\begin{equation}
M_{H}\left(i_{1}+i_{2}p+\cdots+i_{m}p^{m-1}+p^{m}n\right)=M_{H}\left(\mathbf{i}\right)M_{H}\left(n\right),\textrm{ }\forall n\in\mathbb{N}_{0}\label{eq:Inductive identity for M_H}
\end{equation}
\end{prop}
Proof: We have: 
\begin{equation}
H_{\mathbf{i}\wedge\mathbf{j}}\left(x\right)=H_{\mathbf{i}}\left(H_{\mathbf{j}}\left(x\right)\right)
\end{equation}
Differentiating both sides with respect to $x$ and applying the Chain
rule gives: 
\[
H_{\mathbf{i}\wedge\mathbf{j}}^{\prime}\left(x\right)=H_{\mathbf{i}}^{\prime}\left(H_{\mathbf{j}}\left(x\right)\right)H_{\mathbf{j}}^{\prime}\left(x\right)
\]
Now, set $x=0$: 
\begin{equation}
\underbrace{H_{\mathbf{i}\wedge\mathbf{j}}^{\prime}\left(0\right)}_{M_{H}\left(\mathbf{i}\wedge\mathbf{j}\right)}=H_{\mathbf{i}}^{\prime}\left(H_{\mathbf{j}}\left(0\right)\right)\underbrace{H_{\mathbf{j}}^{\prime}\left(0\right)}_{M_{H}\left(\mathbf{j}\right)}
\end{equation}
Since $H_{\mathbf{i}}\left(x\right)$ is an affine linear map of the
form $H_{\mathbf{i}}^{\prime}\left(0\right)x+H_{\mathbf{i}}\left(0\right)$,
its derivative is the constant function $H_{\mathbf{i}}^{\prime}\left(y\right)=H_{\mathbf{i}}^{\prime}\left(0\right)$.
Thus: 
\begin{equation}
M_{H}\left(\mathbf{i}\wedge\mathbf{j}\right)=H_{\mathbf{i}}^{\prime}\left(H_{\mathbf{j}}\left(0\right)\right)M_{H}\left(\mathbf{j}\right)=H_{\mathbf{i}}^{\prime}\left(0\right)M_{H}\left(\mathbf{j}\right)=M_{H}\left(\mathbf{i}\right)M_{H}\left(\mathbf{j}\right)
\end{equation}

Q.E.D.

\vphantom{}

We build a bridge between string-world and number-world by restating
the concatenation identities in terms of systems of functional equations
involving integer inputs. 
\begin{prop}[\textbf{Functional Equations for} $M_{H}$]
\label{prop:M_H functional equation}\index{functional equation!$M_{H}$}\ 
\end{prop}
\begin{equation}
M_{H}\left(pn+j\right)=\frac{\mu_{j}}{p}M_{H}\left(n\right),\textrm{ }\forall n\geq0\textrm{ \& }\forall j\in\mathbb{Z}/p\mathbb{Z}:pn+j\neq0\label{eq:M_H functional equations}
\end{equation}

\begin{rem}
To be clear, this functional equation \emph{does not hold} when $pn+j=0$. 
\end{rem}
Proof: Let $n\geq1$, and let $\mathbf{j}=\left(j_{1},\ldots,j_{m}\right)$
be the shortest string in $\textrm{String}\left(p\right)$ which represents
$n$. Then, observe that for $\mathbf{k}=\left(j,j_{1},j_{2},\ldots,j_{m}\right)$,
where $j\in\mathbb{Z}/p\mathbb{Z}$ is arbitrary, we have that: 
\[
\mathbf{k}=j\wedge\mathbf{j}\sim\left(j+pn\right)
\]
\textbf{Proposition \ref{prop:M_H concatenation identity}} lets us
write: 
\[
M_{H}\left(pn+j\right)=M_{H}\left(j\wedge\mathbf{j}\right)=M_{H}\left(j\right)M_{H}\left(\mathbf{j}\right)=\frac{\mu_{j}}{p}M_{H}\left(n\right)
\]
When $n=0$, these equalities hold for all $j\in\left\{ 1,\ldots,p-1\right\} $,
seeing as: 
\[
\left(j+p\cdot0\right)\sim j\wedge\varnothing
\]
As for the one exceptional case\textemdash $n=j=0$\textemdash note
that we obtain: 
\[
M_{H}\left(p\cdot0+0\right)=M_{H}\left(0\right)\overset{\textrm{def}}{=}1
\]
but: 
\[
\frac{\mu_{0}}{p}M_{H}\left(0\right)\overset{\textrm{def}}{=}\frac{\mu_{0}}{p}\times1=\frac{\mu_{0}}{p}
\]
and, as we saw, $\mu_{0}/p\neq1$.

Q.E.D.

\vphantom{}

Using the concatenation identity for $\chi_{H}$, we can establish
the first of two characterizations for it in terms of functional equations. 
\begin{lem}[\textbf{Functional Equations for} $\chi_{H}$]
\label{lem:Chi_H functional equation on N_0 and uniqueness}$\chi_{H}$
is the unique rational-valued function on $\mathbb{N}_{0}$ satisfying
the system of functional equations:\index{chi{H}@$\chi_{H}$!functional equation}\index{functional equation!chi_{H}@$\chi_{H}$}
\begin{equation}
\chi_{H}\left(pn+j\right)=\frac{a_{j}\chi_{H}\left(n\right)+b_{j}}{d_{j}},\textrm{ }\forall n\in\mathbb{N}_{0},\textrm{ }\forall\mathbb{Z}/p\mathbb{Z}\label{eq:Chi_H functional equations}
\end{equation}
\end{lem}
\begin{rem}
These functional equations can be written more compactly as: 
\begin{equation}
\chi_{H}\left(pn+j\right)=H_{j}\left(\chi_{H}\left(n\right)\right),\textrm{ }\forall n\in\mathbb{N}_{0},\textrm{ }\forall j\in\mathbb{Z}/p\mathbb{Z}\label{eq:Alternate statement of Chi_H functional equations}
\end{equation}
\end{rem}
Proof:

I. Let $\mathbf{i}\sim n$ and let $j\in\mathbb{Z}/p\mathbb{Z}$ be
arbitrary. Then $pn+j\sim j\wedge\mathbf{i}$, and hence, by $\chi_{H}$'s
concatenation identity (\textbf{Lemma \ref{lem:Chi_H concatenation identity}}):
\begin{equation}
\chi_{H}\left(pn+j\right)=\chi_{H}\left(j\wedge\mathbf{i}\right)=H_{j}\left(\chi_{H}\left(\mathbf{i}\right)\right)=H_{j}\left(\chi_{H}\left(n\right)\right)
\end{equation}
Thus, $\chi_{H}$ is a solution of the given system of functional
equations. The reason why we need not exclude the case where $n=j=0$
(as we had to do with $M_{H}$) is because $H_{0}\left(\chi_{H}\left(0\right)\right)=H_{0}\left(0\right)=0$.

\vphantom{}

II. On the other hand, let $f:\mathbb{N}_{0}\rightarrow\mathbb{Q}$
be any function so that: 
\begin{equation}
f\left(pn+j\right)=H_{j}\left(f\left(n\right)\right),\textrm{ }\forall n\in\mathbb{N}_{0},\textrm{ }\forall j\in\mathbb{Z}/p\mathbb{Z}
\end{equation}
Setting $n=j=0$ gives:
\begin{align*}
f\left(0\right) & =H_{0}\left(f\left(0\right)\right)\\
\left(H_{0}\left(x\right)=H_{0}^{\prime}\left(0\right)x+\underbrace{H_{0}\left(0\right)}_{0}\right); & =H_{0}^{\prime}\left(0\right)f\left(0\right)
\end{align*}
Since $H$ is a $p$-Hydra map: 
\[
H_{0}^{\prime}\left(0\right)=\frac{\mu_{0}}{p}\notin\left\{ 0,1\right\} 
\]
Thus, $f\left(0\right)=H_{0}^{\prime}\left(0\right)f\left(0\right)$
forces $f\left(0\right)=0$. Plugging in $n=0$ then leaves us with:
\[
f\left(j\right)=H_{j}\left(f\left(0\right)\right)=H_{j}\left(0\right),\textrm{ }\forall j\in\mathbb{Z}/p\mathbb{Z}
\]
Writing $n$ $p$-adically as: 
\begin{equation}
n=j_{1}+j_{2}p+\cdots+j_{L}p^{L-1}
\end{equation}
the identity $f\left(pn+j\right)=H_{j}\left(f\left(n\right)\right)$
is then equivalent to: 
\begin{align*}
f\left(j+j_{1}p+j_{2}p^{2}+\cdots+j_{L}p^{L}\right) & =H_{j}\left(f\left(j_{1}+j_{2}p+\cdots+j_{L}p^{L-1}\right)\right)\\
 & =H_{j}\left(H_{j_{1}}\left(f\left(j_{2}+j_{3}p+\cdots+j_{L}p^{L-2}\right)\right)\right)\\
 & \vdots\\
 & =\left(H_{j}\circ H_{j_{1}}\circ\cdots\circ H_{j_{L}}\right)\left(f\left(0\right)\right)\\
 & =\left(H_{j}\circ H_{j_{1}}\circ\cdots\circ H_{j_{L}}\right)\left(0\right)\\
 & =H_{j,j_{1},\ldots,j_{L}}\left(0\right)
\end{align*}
So, for any string $\mathbf{j}^{\prime}=\left(j,j_{1},\ldots,j_{L}\right)$,
we have: 
\[
f\left(\mathbf{j}^{\prime}\right)=H_{\mathbf{j}^{\prime}}\left(0\right)\overset{\textrm{def}}{=}\chi_{H}\left(\mathbf{j}^{\prime}\right)
\]
where, note: $\mathbf{j}^{\prime}=j\wedge\mathbf{j}$ and $\mathbf{j}^{\prime}\sim n$.
In other words, if $f$ solves the given system of functional equations,
it is in fact equal to $\chi_{H}$ at every finite string, and hence,
at every non-negative integer. This shows the system's solutions are
unique.

Q.E.D.

\vphantom{}

Next we compute explicit formulae for $M_{H}$ and $\chi_{H}$ in
terms of the numbers present in $H$. These will be needed in to make
the $q$-adic estimates needed to establish the extension/interpolation
of $\chi_{H}$ from a function $\mathbb{N}_{0}\rightarrow\mathbb{Q}$
to a function $\mathbb{Z}_{p}\rightarrow\mathbb{Z}_{q}$, for an appropriate
choice of a prime $q$. 
\begin{prop}
\label{prop:Explicit Formulas for M_H}$M_{H}\left(n\right)$ and
$M_{H}\left(\mathbf{j}\right)$ can be explicitly given in terms of
the constants associated to $H$ by the formulae: 
\begin{equation}
M_{H}\left(n\right)=\frac{1}{p^{\lambda_{p}\left(n\right)}}\prod_{k=0}^{p-1}\mu_{k}^{\#_{p:k}\left(n\right)}=\prod_{k=0}^{p-1}\frac{a_{k}^{\#_{p:k}\left(n\right)}}{d_{k}^{\#_{p:k}\left(n\right)}}\label{eq:Formula for M_H of n}
\end{equation}
and: 
\begin{equation}
M_{H}\left(\mathbf{j}\right)=\prod_{k=1}^{\left|\mathbf{j}\right|}\frac{a_{j_{k}}}{d_{j_{k}}}=\frac{1}{p^{\left|\mathbf{j}\right|}}\prod_{k=1}^{\left|\mathbf{j}\right|}\mu_{j_{k}}=\frac{1}{p^{\left|\mathbf{j}\right|}}\prod_{k=0}^{p-1}\mu_{k}^{\#_{p:k}\left(\mathbf{j}\right)}=\prod_{k=0}^{p-1}\frac{a_{k}^{\#_{p:k}\left(\mathbf{j}\right)}}{d_{k}^{\#_{p:k}\left(\mathbf{j}\right)}}\label{eq:formula for M_H of bold-j}
\end{equation}
respectively. We adopt the convention that the $k$-product in \emph{(\ref{eq:formula for M_H of bold-j})}
is defined to be $1$ when $\left|\mathbf{j}\right|=0$.

Finally, for any $\mathfrak{z}\in\mathbb{Z}_{p}$ and any $N\geq1$,
we have the formula: 
\begin{equation}
M_{H}\left(\left[\mathfrak{z}\right]_{p^{N}}\right)=\frac{\prod_{j=0}^{p-1}\mu_{j}^{\#_{p:j}\left(\left[\mathfrak{z}\right]_{p^{N}}\right)}}{p^{\lambda_{p}\left(\left[\mathfrak{z}\right]_{p^{N}}\right)}}\label{eq:M_H of z mod rho to the N}
\end{equation}
obtained by setting $n=\left[\mathfrak{z}\right]_{p^{N}}$ in \emph{(\ref{eq:Formula for M_H of n})}. 
\end{prop}
Proof: The proof is purely computational. I omit it, seeing as it
occurs in the proof of \textbf{Proposition \ref{prop:Explicit formula for Chi_H of bold j}},
given below.

Q.E.D. 
\begin{prop}
\label{prop:Explicit formula for Chi_H of bold j} 
\begin{equation}
\chi_{H}\left(\mathbf{j}\right)=\sum_{m=1}^{\left|\mathbf{j}\right|}\frac{b_{j_{m}}}{d_{j_{m}}}\prod_{k=1}^{m-1}\frac{a_{j_{k}}}{d_{j_{k}}}=\sum_{m=1}^{\left|\mathbf{j}\right|}\frac{b_{j_{m}}}{a_{j_{m}}}\left(\prod_{k=1}^{m}\mu_{j_{k}}\right)p^{-m},\textrm{ }\forall\mathbf{j}\in\textrm{String}\left(p\right)\label{eq:Formula for Chi_H in terms of bold-j}
\end{equation}
where the $k$-product is defined to be $1$ when $m=1$. 
\end{prop}
Proof: Let $\mathbf{j}=\left(j_{1},\ldots,j_{\left|\mathbf{j}\right|}\right)\in\textrm{String}\left(p\right)$
be arbitrary. Since $\chi_{H}\left(\mathbf{j}\right)=\chi_{H}\left(\mathbf{i}\right)$
for any $\mathbf{i}\in\textrm{String}\left(p\right)$ for which $\mathbf{i}\sim\mathbf{j}$,
we can assume without loss of generality that $\mathbf{j}$'s right-most
entry is non-zero; that is: $j_{\left|\mathbf{j}\right|}\neq0$. The
proof follows by examining: 
\begin{equation}
H_{\mathbf{j}}\left(x\right)=M_{H}\left(\mathbf{j}\right)x+\chi_{H}\left(\mathbf{j}\right),\textrm{ }\forall x\in\mathbb{R}
\end{equation}
and carefully writing out the composition sequence in full: 
\begin{align*}
H_{\mathbf{j}}\left(x\right) & =H_{\left(j_{1},\ldots,j_{\left|\mathbf{j}\right|}\right)}\left(x\right)\\
 & =\frac{a_{j_{1}}\frac{a_{j_{2}}\left(\cdots\right)+b_{j_{2}}}{d_{j_{2}}}+b_{j_{1}}}{d_{j_{1}}}\\
 & =\overbrace{\left(\prod_{k=1}^{\left|\mathbf{j}\right|}\frac{a_{j_{k}}}{d_{j_{k}}}\right)}^{M_{H}\left(\mathbf{j}\right)}x+\overbrace{\frac{b_{j_{1}}}{d_{j_{1}}}+\frac{a_{j_{1}}}{d_{j_{1}}}\frac{b_{j_{2}}}{d_{j_{2}}}+\frac{a_{j_{2}}}{d_{j_{2}}}\frac{a_{j_{1}}}{d_{j_{1}}}\frac{b_{j_{3}}}{d_{j_{3}}}+\cdots+\left(\prod_{k=1}^{\left|\mathbf{j}\right|-1}\frac{a_{j_{k}}}{d_{j_{k}}}\right)\frac{b_{j_{\left|\mathbf{j}\right|}}}{d_{j_{\left|\mathbf{j}\right|}}}}^{\chi_{H}\left(\mathbf{j}\right)}\\
 & =M_{H}\left(\mathbf{j}\right)x+\underbrace{\sum_{m=1}^{\left|\mathbf{j}\right|}\frac{b_{j_{m}}}{d_{j_{m}}}\prod_{k=1}^{m-1}\frac{a_{j_{k}}}{d_{j_{k}}}}_{\chi_{H}\left(\mathbf{j}\right)}
\end{align*}
where we use the convention that $\prod_{k=1}^{m-1}\frac{a_{j_{k}}}{d_{j_{k}}}=1$
when $m=1$. This proves: 
\begin{equation}
\chi_{H}\left(\mathbf{j}\right)=\sum_{m=1}^{\left|\mathbf{j}\right|}\frac{b_{j_{m}}}{d_{j_{m}}}\prod_{k=1}^{m-1}\frac{a_{j_{k}}}{d_{j_{k}}}
\end{equation}

Next, using $\mu_{j}=\frac{pa_{j}}{d_{j}}$, for $m>1$, the $m$th
term $\frac{b_{j_{m}}}{d_{j_{m}}}\prod_{k=1}^{m-1}\frac{a_{j_{k}}}{d_{j_{k}}}$
can be written as: 
\begin{align*}
\frac{b_{j_{m}}}{d_{j_{m}}}\prod_{k=1}^{m-1}\frac{a_{j_{k}}}{d_{j_{k}}} & =\frac{1}{p}\cdot\frac{b_{j_{m}}}{a_{j_{m}}}\cdot\frac{pa_{j_{m}}}{d_{j_{m}}}\cdot\prod_{k=1}^{m-1}\left(\frac{pa_{j_{k}}}{d_{j_{k}}}\cdot\frac{1}{p}\right)\\
 & =\frac{1}{p}\frac{b_{j_{m}}}{a_{j_{m}}}\frac{\mu_{j_{m}}}{p^{m-1}}\prod_{k=1}^{m-1}\mu_{j_{k}}\\
 & =\frac{b_{j_{m}}}{a_{j_{m}}}\left(\prod_{k=1}^{m}\mu_{j_{k}}\right)p^{-m}
\end{align*}
On the other hand, when $m=1$, by our $k$-product convention, we
have: 
\begin{equation}
\frac{b_{j_{1}}}{d_{j_{1}}}\underbrace{\prod_{k=1}^{1-1}\frac{a_{j_{k}}}{d_{j_{k}}}}_{1}=\frac{b_{j_{1}}}{d_{j_{1}}}
\end{equation}
Thus: 
\begin{align*}
\chi_{H}\left(\mathbf{j}\right) & =\frac{b_{j_{1}}}{d_{j_{1}}}+\sum_{m=2}^{\left|\mathbf{j}\right|}\frac{b_{j_{m}}}{a_{j_{m}}}\left(\prod_{k=1}^{m}\mu_{j_{k}}\right)p^{-m}\\
 & =\frac{1}{p}\frac{b_{j_{1}}}{a_{j_{1}}}\frac{pa_{j_{1}}}{d_{j_{1}}}+\sum_{m=2}^{\left|\mathbf{j}\right|}\frac{b_{j_{m}}}{a_{j_{m}}}\left(\prod_{k=1}^{m}\mu_{j_{k}}\right)p^{-m}\\
 & =\frac{b_{j_{1}}}{a_{j_{1}}}\cdot\mu_{j_{1}}\cdot p^{-1}+\sum_{m=2}^{\left|\mathbf{j}\right|}\frac{b_{j_{m}}}{a_{j_{m}}}\left(\prod_{k=1}^{m}\mu_{j_{k}}\right)p^{-m}\\
 & =\sum_{m=1}^{\left|\mathbf{j}\right|}\frac{b_{j_{m}}}{a_{j_{m}}}\left(\prod_{k=1}^{m}\mu_{j_{k}}\right)p^{-m}
\end{align*}
as desired.

Q.E.D.

\vphantom{}

To interpolate $\chi_{H}$ from a rational-valued function on $\mathbb{N}_{0}$
to a $q$-adic-valued function on $\mathbb{Z}_{p}$, we need to define
some qualitative properties of Hydra maps so as to distinguish those
cases where this interpolation will actually exist.
\begin{defn}[\textbf{Qualitative Terminology for $p$-Hydra maps}]
\label{def:Qualitative conditions on a p-Hydra map}Let $H$ be a
$p$-Hydra map. For the definitions below, we DO NOT require $p$
to be prime.

\vphantom{}

I. We say $H$ is \textbf{simple }if $d_{j}=p$ for all $j$.

\vphantom{}

II. We say $H$ is \textbf{semi-simple }if $\gcd\left(a_{j},d_{k}\right)=1$
for all $j,k\in\left\{ 0,\ldots,p-1\right\} $.

\vphantom{}

III. We say $H$ is \textbf{degenerate }if $a_{j}=1$ for some non-zero
$j$. If $a_{j}\neq1$ for all non-zero $j$, we say $H$ is \textbf{non-degenerate}.

\vphantom{}

IV\footnote{This definition will not be relevant until Chapter 3.}.
We say $H$ is \textbf{contracting }if $a_{0}/d_{0}$ (a.k.a. $\mu_{0}/p$)
is $<1$, and say $H$ is \textbf{expanding }if $a_{0}/d_{0}$ (a.k.a.
$\mu_{0}/p$) is $>1$.

\vphantom{}

V. We say $H$ is \textbf{monogenic}\footnote{This term is based on the Greek for ``one kind/species''.}
whenever: 
\begin{equation}
\gcd\left(\left\{ a_{j}:j\in\left\{ 0,\ldots,p-1\right\} \textrm{ \& }a_{j}\neq1\right\} \right)\geq2\label{eq:Monogenic gcd}
\end{equation}
When this $\textrm{gcd}$ is $1$, we say $H$ is \textbf{polygenic}.
If $H$ is monogenic, we then write $q_{H}$\nomenclature{$q_{H}$}{ }
to denote the smallest prime divisor of (\ref{eq:Monogenic gcd}).
For brevity, sometimes we will drop the $H$ and just write this as
$q$ instead of $q_{H}$.

\vphantom{}

VI. We say $H$ is \textbf{basic }if $H$ is simple, non-degenerate,
and monogenic.

\vphantom{}

VII. We say $H$ is \textbf{semi-basic }if $H$ is semi-simple, non-degenerate,
and monogenic.
\end{defn}
\begin{rem}
If $H$ is simple, note that $\mu_{j}=a_{j}$ for all $j$.
\end{rem}
\begin{rem}
It is worth noting that all of the results in this section hold when
$q_{H}$ is \emph{any }prime divisor of (\ref{eq:Monogenic gcd}).
\end{rem}
\vphantom{}

Next, a minor technical result which shows that our qualitative definitions
are actually useful. 
\begin{prop}[\textbf{Co-primality of $d_{j}$ and $q_{H}$}]
\label{prop:co-primality of d_j and q_H}Let $H$ be a semi-basic
$p$-Hydra map. Then, $\gcd\left(d_{j},q_{H}\right)=1$ for all $j\in\left\{ 0,\ldots,p-1\right\} $. 
\end{prop}
Proof: By way of contradiction, suppose there is a $k\in\left\{ 0,\ldots,p-1\right\} $
so that $\gcd\left(d_{k},q_{H}\right)>1$. Since $H$ is semi-basic,
$H$ is monogenic and non-degenerate, and so $q_{H}$ is $\geq2$
and, moreover, $a_{j}\neq1$ for any $j\in\left\{ 1,\ldots,p-1\right\} $.
Thus, there is a $k$ so that $d_{k}$ and $q_{H}$ have a non-trivial
common factor, and there is a $j$ so that $q_{H}$ and $a_{j}$ have
a non-trivial common factor (namely, $q_{H}$). Since every factor
of $q_{H}$ divides $a_{j}$, the non-trivial common factor of $d_{k}$
and $q_{H}$ must also divide $a_{j}$. This forces $\gcd\left(a_{j},d_{k}\right)>1$.
However, as a semi-basic $p$-Hydra map, $H$ is semi-simple, and
so $\gcd\left(a_{j},d_{k}\right)=1$ for all $j,k\in\left\{ 0,\ldots,p-1\right\} $\textemdash and
that is a contradiction!

Thus, it must be that $\gcd\left(d_{k},q_{H}\right)=1$ for all $k\in\left\{ 0,\ldots,p-1\right\} $.

Q.E.D.

\vphantom{}

While our next result is also technical, it is not minor in the least.
It demonstrates that $M_{H}\left(\mathbf{j}\right)$ will be $q_{H}$-adically
small whenever $\mathbf{j}$ has many non-zero entries. This is the
essential ingredient for the proof of the existence of $\chi_{H}$'s
$p$-adic interpolation.
\begin{prop}[\textbf{$q_{H}$-adic behavior of $M_{H}\left(\mathbf{j}\right)$
as $\left|\mathbf{j}\right|\rightarrow\infty$}]
\label{prop:q-adic behavior of M_H of j as the number of non-zero digits tends to infinity}Let
$H$ be a semi-basic $p$-Hydra map which fixes $0$. Then: 
\begin{equation}
\left|M_{H}\left(\mathbf{j}\right)\right|_{q_{H}}\leq q_{H}^{-\sum_{k=1}^{p-1}\#_{p:k}\left(\mathbf{j}\right)}\label{eq:M_H of bold j has q-adic valuation at least as large as the number of non-zero digits of bold-j}
\end{equation}
In particular, for any sequence of strings $\left\{ \mathbf{j}_{n}\right\} _{n\geq1}$
in $\textrm{String}\left(p\right)$ so that $\lim_{n\rightarrow\infty}\left|\mathbf{j}_{n}\right|=\infty$,
we have that $\lim_{n\rightarrow\infty}\left|M_{H}\left(\mathbf{j}_{n}\right)\right|_{q_{H}}=0$
whenever the number of non-zero entries in $\mathbf{j}_{n}$ tends
to $\infty$ as $n\rightarrow\infty$. 
\end{prop}
Proof: Let $\mathbf{j}$ be a non-zero element of $\textrm{String}\left(p\right)$.
Then, by \textbf{Proposition \ref{prop:Explicit Formulas for M_H}}:
\[
M_{H}\left(\mathbf{j}\right)=\prod_{\ell=1}^{\left|\mathbf{j}\right|}\frac{a_{j_{\ell}}}{d_{j_{\ell}}}
\]
Since $H$ is semi-basic, \textbf{Proposition \ref{prop:co-primality of d_j and q_H}}
tells us that every $d_{j}$ is co-prime to $q_{H}$, and so $\left|d_{j}\right|_{q_{H}}=1$
for any $j$. On the other hand, the non-degeneracy and monogenicity
of $H$ tells us that $a_{j}$ is a non-zero integer multiple of $q_{H}$
for all $j\in\left\{ 1,\ldots,p-1\right\} $; thus, $\left|a_{j}\right|_{q_{H}}\leq1/q_{H}$
for every $\ell$ for which $j_{\ell}\neq0$. Taking $q_{H}$-adic
absolute values of $M_{H}\left(\mathbf{j}\right)$ then gives: 
\begin{equation}
\left|M_{H}\left(\mathbf{j}\right)\right|_{q_{H}}=\prod_{\ell=1}^{\left|\mathbf{j}\right|}\left|\frac{a_{j_{\ell}}}{d_{j_{\ell}}}\right|_{q_{H}}=\prod_{\ell:a_{j_{\ell}}\neq0}\left|a_{j_{\ell}}\right|_{q_{H}}\leq q_{H}^{-\left|\left\{ \ell\in\left\{ 1,\ldots,\left|\mathbf{j}\right|\right\} :j_{\ell}\neq0\right\} \right|}
\end{equation}
Here, $\left|\left\{ \ell\in\left\{ 1,\ldots,\left|\mathbf{j}\right|\right\} :j_{\ell}\neq0\right\} \right|$
is precisely the number of non-zero entries of $\mathbf{j}$, which
is: 
\begin{equation}
\left|\left\{ \ell\in\left\{ 1,\ldots,\left|\mathbf{j}\right|\right\} :j_{\ell}\neq0\right\} \right|=\sum_{k=1}^{p-1}\#_{p:k}\left(\mathbf{j}\right)
\end{equation}
For a sequence of $\mathbf{j}_{n}$s tending to $\infty$ in length,
this also shows that $\left|M_{H}\left(\mathbf{j}_{n}\right)\right|_{q_{H}}\rightarrow0$
as $n\rightarrow\infty$, provided that the number of non-zero entries
in $\mathbf{j}_{n}$ also tends to $\infty$ as $n\rightarrow\infty$.

Q.E.D.

\vphantom{}

Now the main result of this subsection: the $\left(p,q\right)$-adic
characterization of $\chi_{H}$. 
\begin{lem}[\textbf{$\left(p,q\right)$-adic Characterization of }$\chi_{H}$]
\label{lem:Unique rising continuation and p-adic functional equation of Chi_H}Let
$H$ be a semi-basic $p$-Hydra map which fixes $0$. Then, the limit:
\begin{equation}
\chi_{H}\left(\mathfrak{z}\right)\overset{\mathbb{Z}_{q_{H}}}{=}\lim_{n\rightarrow\infty}\chi_{H}\left(\left[\mathfrak{z}\right]_{p^{n}}\right)\label{eq:Rising Continuity Formula for Chi_H}
\end{equation}
exists for all $\mathfrak{z}\in\mathbb{Z}_{p}$, and thereby defines
an interpolation of $\chi_{H}$ to a function $\chi_{H}:\mathbb{Z}_{p}\rightarrow\mathbb{Z}_{q_{H}}$.
Moreover:

\vphantom{}

I.

\begin{equation}
\chi_{H}\left(p\mathfrak{z}+j\right)=\frac{a_{j}\chi_{H}\left(\mathfrak{z}\right)+b_{j}}{d_{j}},\textrm{ }\forall\mathfrak{z}\in\mathbb{Z}_{p},\textrm{ }\forall j\in\mathbb{Z}/p\mathbb{Z}\label{eq:Functional Equations for Chi_H over the rho-adics}
\end{equation}

\vphantom{}

II. The interpolation $\chi_{H}:\mathbb{Z}_{p}\rightarrow\mathbb{Z}_{q_{H}}$
defined by the limit \emph{(\ref{eq:Rising Continuity Formula for Chi_H})}
is the \textbf{unique} function $f:\mathbb{Z}_{p}\rightarrow\mathbb{Z}_{q_{H}}$
satisfying the functional equations:

\emph{
\begin{equation}
f\left(p\mathfrak{z}+j\right)=\frac{a_{j}f\left(\mathfrak{z}\right)+b_{j}}{d_{j}},\textrm{ }\forall\mathfrak{z}\in\mathbb{Z}_{p},\textrm{ }\forall j\in\mathbb{Z}/p\mathbb{Z}\label{eq:unique p,q-adic functional equation of Chi_H}
\end{equation}
}along with the \textbf{rising-continuity}\footnote{As discussed in Section \ref{sec:3.2 Rising-Continuous-Functions},
I say a function is \textbf{rising-continuous }whenever it satisfies
the limit condition (\ref{eq:unique p,q-adic rising continuity of Chi_H}).
(\ref{eq:Rising Continuity Formula for Chi_H}) shows that $\chi_{H}$
is rising-continuous.} \textbf{condition}:\textbf{ } 
\begin{equation}
f\left(\mathfrak{z}\right)\overset{\mathbb{Z}_{q_{H}}}{=}\lim_{n\rightarrow\infty}f\left(\left[\mathfrak{z}\right]_{p^{n}}\right),\textrm{ }\forall\mathfrak{z}\in\mathbb{Z}_{p}\label{eq:unique p,q-adic rising continuity of Chi_H}
\end{equation}
 
\end{lem}
Proof: First, we show the existence of the limit (\ref{eq:Rising Continuity Formula for Chi_H}).
If $\mathfrak{z}\in\mathbb{N}_{0}\cap\mathbb{Z}_{p}$, then $\left[\mathfrak{z}\right]_{p^{N}}=\mathfrak{z}$
for all sufficiently large $N$, showing that the limit (\ref{eq:Rising Continuity Formula for Chi_H})
exists in that case. So, suppose $\mathfrak{z}\in\mathbb{Z}_{p}^{\prime}$.
Next, let $\mathbf{j}_{n}=\left(j_{1},\ldots,j_{n}\right)$ be the
shortest string representing $\left[\mathfrak{z}\right]_{p^{n}}$;
that is: 
\begin{equation}
\left[\mathfrak{z}\right]_{p^{n}}=\sum_{k=1}^{n}j_{k}p^{k-1}
\end{equation}
By \textbf{Lemma \ref{eq:Definition of Chi_H of n}}, (\ref{eq:Rising Continuity Formula for Chi_H})
can be written as:
\begin{equation}
\chi_{H}\left(\mathfrak{z}\right)=\lim_{n\rightarrow\infty}\chi_{H}\left(\mathbf{j}_{n}\right)
\end{equation}
Using \textbf{Proposition \ref{prop:Explicit formula for Chi_H of bold j}},
we can write: 
\begin{equation}
\chi_{H}\left(\mathbf{j}_{n}\right)=H_{\mathbf{j}_{n}}\left(0\right)=\sum_{m=1}^{\left|\mathbf{j}_{n}\right|}\frac{b_{j_{m}}}{d_{j_{m}}}\prod_{k=1}^{m-1}\frac{a_{j_{k}}}{d_{j_{k}}}
\end{equation}
where the product is defined to be $1$ when $m=1$. Then, using \textbf{Proposition
\ref{prop:Explicit Formulas for M_H}}, we get: 
\begin{equation}
\prod_{k=1}^{m-1}\frac{a_{j_{k}}}{d_{j_{k}}}=M_{H}\left(\mathbf{j}_{m-1}\right)
\end{equation}
where, again, the product is defined to be $1$ when $m=1$. So, our
previous equation for $\chi_{H}\left(\mathbf{j}_{n}\right)$ can be
written as: 
\begin{equation}
\chi_{H}\left(\mathbf{j}_{n}\right)=\sum_{m=1}^{\left|\mathbf{j}_{n}\right|}\frac{b_{j_{m}}}{d_{j_{m}}}M_{H}\left(\mathbf{j}_{m-1}\right)
\end{equation}
Taking limits yields: 
\begin{equation}
\lim_{n\rightarrow\infty}\chi_{H}\left(\left[\mathfrak{z}\right]_{p^{n}}\right)=\lim_{n\rightarrow\infty}\chi_{H}\left(\mathbf{j}_{n}\right)=\lim_{n\rightarrow\infty}\sum_{m=1}^{\infty}\frac{b_{j_{m}}}{d_{j_{m}}}M_{H}\left(\mathbf{j}_{m-1}\right)\label{eq:Formal rising limit of Chi_H}
\end{equation}
Thanks to the ultrametric topology of $\mathbb{Z}_{q_{H}}$, the series
on the right will converge in $\mathbb{Z}_{q_{H}}$ if and only if
is $m$th term tends to $0$ $q_{H}$-adically: 
\begin{equation}
\lim_{m\rightarrow\infty}\left|\frac{b_{j_{m}}}{d_{j_{m}}}M_{H}\left(\mathbf{j}_{m-1}\right)\right|_{q_{H}}=0
\end{equation}
Because the $b_{j}$s and $d_{j}$s belong to a finite set, with the
$d_{j}$s being all non-zero: 
\begin{equation}
\sup_{j\in\left\{ 0,\ldots,p-1\right\} }\left|\frac{b_{j}}{d_{j}}\right|_{q_{H}}<\infty
\end{equation}
and so: 
\begin{equation}
\left|\frac{b_{j_{m}}}{d_{j_{m}}}M_{H}\left(\mathbf{j}_{m-1}\right)\right|_{q_{H}}\ll\left|M_{H}\left(\mathbf{j}_{m-1}\right)\right|_{q_{H}}
\end{equation}
Because $\mathfrak{z}\in\mathbb{Z}_{p}^{\prime}$, $\mathfrak{z}$
has infinitely many non-zero digits. As such, both the length and
number of non-zero entries of $\mathbf{j}_{m-1}$ tend to infinity
as $m\rightarrow\infty$. By\textbf{ Proposition \ref{prop:q-adic behavior of M_H of j as the number of non-zero digits tends to infinity}},
this forces $\left|M_{H}\left(\mathbf{j}_{m-1}\right)\right|_{q_{H}}\rightarrow0$
as $m\rightarrow\infty$. So, the $m$th term of our series decays
to zero $q_{H}$-adically, which then guarantees the convergence of
$\lim_{n\rightarrow\infty}\chi_{H}\left(\left[\mathfrak{z}\right]_{p^{n}}\right)$
in $\mathbb{Z}_{q_{H}}$.

\vphantom{}

To prove (I), by \textbf{Lemma \ref{lem:Chi_H functional equation on N_0 and uniqueness}},
we know that $\chi_{H}$ satisfies the functional equations (\ref{eq:Functional Equations for Chi_H over the rho-adics})
for $\mathfrak{z}\in\mathbb{N}_{0}$. Taking limits of the functional
equations using (\ref{eq:Rising Continuity Formula for Chi_H}) then
proves the equations hold for all $\mathfrak{z}\in\mathbb{Z}_{p}$.

\vphantom{}

For (II), we have just shown that $\chi_{H}$ satisfies the interpolation
condition and the functional equations. So, conversely, suppose $f:\mathbb{Z}_{p}\rightarrow\mathbb{Z}_{q_{H}}$
satisfies (\ref{eq:unique p,q-adic functional equation of Chi_H})
and (\ref{eq:unique p,q-adic rising continuity of Chi_H}). Then,
by \textbf{Lemma \ref{lem:Chi_H functional equation on N_0 and uniqueness}},
the restriction of $f$ to $\mathbb{N}_{0}$ must be equal to $\chi_{H}$.
(\ref{eq:Rising Continuity Formula for Chi_H}) shows that the values
$f\left(\mathfrak{z}\right)\overset{\mathbb{Z}_{q_{H}}}{=}\lim_{n\rightarrow\infty}f\left(\left[\mathfrak{z}\right]_{p^{n}}\right)$
obtained by taking limits necessarily forces $f\left(\mathfrak{z}\right)=\chi_{H}\left(\mathfrak{z}\right)$
for all $\mathfrak{z}\in\mathbb{Z}_{p}$.

Q.E.D.

\subsection{\label{subsec:2.2.3 The-Correspondence-Principle}The Correspondence
Principle}

THROUGHOUT THIS SUBSECTION, WE ASSUME $H$ IS INTEGRAL.

\vphantom{}

Our goal\textemdash the \textbf{Correspondence Principle}\textemdash relates
the periodic points of $H$ in $\mathbb{Z}$ to the set $\mathbb{Z}\cap\chi_{H}\left(\mathbb{Z}_{p}\right)$,
where we view both $\mathbb{Z}$ and $\chi_{H}\left(\mathbb{Z}_{p}\right)$
as being embedded in $\mathbb{Z}_{q_{H}}$. The distinction between
basic and semi-basic $p$-Hydra maps will play a key role here, as
there is an additional property (\textbf{propriety}, to be defined
below) satisfied by semi-basic integral $p$-Hydra maps which is needed
to ensure that the set $\mathbb{Z}\cap\chi_{H}\left(\mathbb{Z}_{p}\right)$
fully characterizes the periodic points of $H$.

Before proceeding, we will also need a certain self-map of $\mathbb{Z}_{p}$
which I denote by $B_{p}$, so as to formulate one more functional
equation satisfied by $\chi_{H}$, one directly implicated in the
Correspondence Principle. First, however, a useful\textemdash and
motivating\textemdash observation. 
\begin{prop}
\label{prop:Concatenation exponentiation}Let $p$ be an integer $\geq2$,
let $n\in\mathbb{N}_{1}$, and let $\mathbf{j}\in\textrm{String}\left(p\right)$
be the shortest string representing $n$. Then, for all $m\in\mathbb{N}_{1}$:
\begin{equation}
\mathbf{j}^{\wedge m}\sim n\frac{1-p^{m\lambda_{p}\left(n\right)}}{1-p^{\lambda_{p}\left(n\right)}}\label{eq:Proposition 1.2.10}
\end{equation}
That is to say, the quantity on the right is the integer whose sequence
of $p$-adic digits consists of $m$ concatenated copies of the $p$-adic
digits of $n$ (or, equivalently, the entries of $\mathbf{j}$). 
\end{prop}
Proof: Since: 
\[
n=j_{1}+j_{2}p+\cdots+j_{\lambda_{p}\left(n\right)-1}p^{\lambda_{p}\left(n\right)-1}
\]
we have that: 
\begin{align*}
\mathbf{j}^{\wedge m} & \sim j_{1}+j_{2}p+\cdots+j_{\lambda_{p}\left(n\right)-1}p^{\lambda_{p}\left(n\right)-1}\\
 & +p^{\lambda_{p}\left(n\right)}\left(j_{1}+j_{2}p+\cdots+j_{\lambda_{p}\left(n\right)-1}p^{\lambda_{p}\left(n\right)-1}\right)\\
 & +\cdots\\
 & +p^{\left(m-1\right)\lambda_{p}\left(n\right)}\left(j_{1}+j_{2}p+\cdots+j_{\lambda_{p}\left(n\right)-1}p^{\lambda_{p}\left(n\right)-1}\right)\\
 & =n+np^{\lambda_{p}\left(n\right)}+np^{2\lambda_{p}\left(n\right)}+\cdots+np^{\left(m-1\right)\lambda_{p}\left(n\right)}\\
 & =n\frac{1-p^{m\lambda_{p}\left(n\right)}}{1-p^{\lambda_{p}\left(n\right)}}
\end{align*}
as desired.

Q.E.D. 
\begin{defn}[$B_{p}$]
\nomenclature{$B_{p}$}{ }Let $p$ be an integer $\geq2$. Then,
we define $B_{p}:\mathbb{N}_{0}\rightarrow\mathbb{Q}\cap\mathbb{Z}_{p}$
by: 
\begin{equation}
B_{p}\left(n\right)\overset{\textrm{def}}{=}\begin{cases}
0 & \textrm{if }n=0\\
\frac{n}{1-p^{\lambda_{p}\left(n\right)}} & \textrm{if }n\geq1
\end{cases}\label{eq:Definition of B projection function}
\end{equation}
\end{defn}
\begin{rem}
$B_{p}$ has a simple interpretation in terms of $p$-adic expansions:
it sends $n$ to the $p$-adic integer whose sequence of $p$-adic
digits consists of infinitely many concatenated copies of the sequence
of $p$-adic digits of $n$: 
\begin{equation}
B_{p}\left(n\right)\overset{\mathbb{Z}_{p}}{=}\lim_{m\rightarrow\infty}n\frac{1-p^{m\lambda_{p}\left(n\right)}}{1-p^{\lambda_{p}\left(n\right)}}\overset{\mathbb{Z}_{p}}{=}\frac{n}{1-p^{\lambda_{p}\left(n\right)}}
\end{equation}
In particular, since the sequence of $p$-adic digits of $B_{p}\left(n\right)$
is, by construction, periodic, the geometric series formula in $\mathbb{Z}_{p}$
guarantees the $p$-adic integer $B_{p}\left(n\right)$ is in fact
an element of $\mathbb{Q}$. 
\end{rem}
\begin{rem}
Note that $B_{p}$ extends to a function on $\mathbb{Z}_{p}$, one
whose restriction to $\mathbb{Z}_{p}^{\prime}$ is identity map. 
\end{rem}
\vphantom{}

Using $\chi_{H}$ and $B_{p}$, we can compactly and elegantly express
the diophantine equation (\ref{eq:The Bohm-Sontacchi Criterion})
of the Böhm-Sontacchi Criterion in terms of a functional equation.
The \textbf{Correspondence Principle }emerges almost immediately from
this functional equation. 
\begin{lem}[\textbf{Functional Equation for }$\chi_{H}\circ B_{p}$]
\label{lem:Chi_H o B_p functional equation}Let\index{functional equation!chi_{H}circ B_{p}@$\chi_{H}\circ B_{p}$}
$H$ be semi-basic. Then:
\begin{equation}
\chi_{H}\left(B_{p}\left(n\right)\right)\overset{\mathbb{Z}_{qH}}{=}\frac{\chi_{H}\left(n\right)}{1-M_{H}\left(n\right)},\textrm{ }\forall n\in\mathbb{N}_{1}\label{eq:Chi_H B functional equation}
\end{equation}
\end{lem}
Proof: Let $n\in\mathbb{N}_{1}$, and let $\mathbf{j}\in\textrm{String}\left(p\right)$
be the shortest string representing $n$. Now, for all $m\geq1$ and
all $k\in\mathbb{N}_{1}$: 
\begin{equation}
H_{\mathbf{j}^{\wedge m}}\left(k\right)=H_{\mathbf{j}}^{\circ m}\left(k\right)
\end{equation}
Since: 
\[
H_{\mathbf{j}}\left(k\right)=M_{H}\left(\mathbf{j}\right)k+\chi_{H}\left(\mathbf{j}\right)
\]
the geometric series formula gives us: 
\begin{equation}
H_{\mathbf{j}^{\wedge m}}\left(k\right)=H_{\mathbf{j}}^{\circ m}\left(k\right)=\left(M_{H}\left(\mathbf{j}\right)\right)^{m}k+\frac{1-\left(M_{H}\left(\mathbf{j}\right)\right)^{m}}{1-M_{H}\left(\mathbf{j}\right)}\chi_{H}\left(\mathbf{j}\right)
\end{equation}
Since we also have: 
\begin{equation}
H_{\mathbf{j}^{\wedge m}}\left(k\right)=M_{H}\left(\mathbf{j}^{\wedge m}\right)k+\chi_{H}\left(\mathbf{j}^{\wedge m}\right)
\end{equation}
this then yields: 
\begin{equation}
M_{H}\left(\mathbf{j}^{\wedge m}\right)k+\chi_{H}\left(\mathbf{j}^{\wedge m}\right)=\left(M_{H}\left(\mathbf{j}\right)\right)^{m}k+\frac{1-\left(M_{H}\left(\mathbf{j}\right)\right)^{m}}{1-M_{H}\left(\mathbf{j}\right)}\chi_{H}\left(\mathbf{j}\right)
\end{equation}
By \textbf{Proposition \ref{prop:M_H concatenation identity}}, we
see that $M_{H}\left(\mathbf{j}^{\wedge m}\right)k=\left(M_{H}\left(\mathbf{j}\right)\right)^{m}k$.
Cancelling these terms from both sides leaves us with: 
\begin{equation}
\chi_{H}\left(\mathbf{j}^{\wedge m}\right)=\frac{1-\left(M_{H}\left(\mathbf{j}\right)\right)^{m}}{1-M_{H}\left(\mathbf{j}\right)}\chi_{H}\left(\mathbf{j}\right)=\frac{1-\left(M_{H}\left(\mathbf{j}\right)\right)^{m}}{1-M_{H}\left(n\right)}\chi_{H}\left(n\right)\label{eq:Chi_H B functional equation, ready to take limits}
\end{equation}

Now, by \textbf{Proposition \ref{prop:Concatenation exponentiation}},
we have that: 
\begin{equation}
\chi_{H}\left(\mathbf{j}^{\wedge m}\right)=\chi_{H}\left(n\frac{1-p^{m\lambda_{p}\left(n\right)}}{1-p^{\lambda_{p}\left(n\right)}}\right)
\end{equation}
where the equality is of rational numbers in $\mathbb{R}$. Since
$n\in\mathbb{N}_{1}$, we have that: 
\begin{equation}
B_{p}\left(n\right)=\frac{n}{1-p^{\lambda_{p}\left(n\right)}}
\end{equation}
is a $p$-adic integer. Moreover, as is immediate from the proof of
\textbf{Proposition \ref{prop:Concatenation exponentiation}}, the
projection of this $p$-adic integer modulo $p^{m}$ is: 
\begin{equation}
\left[B_{p}\left(n\right)\right]_{p^{m}}=n\frac{1-p^{m\lambda_{p}\left(n\right)}}{1-p^{\lambda_{p}\left(n\right)}}
\end{equation}
which is, of course, exactly the rational integer represented by the
string $\mathbf{j}^{\wedge m}$. In other words: 
\begin{equation}
\frac{1-\left(M_{H}\left(\mathbf{j}\right)\right)^{m}}{1-M_{H}\left(n\right)}\chi_{H}\left(n\right)=\chi_{H}\left(\mathbf{j}^{\wedge m}\right)=\chi_{H}\left(n\frac{1-p^{m\lambda_{p}\left(n\right)}}{1-p^{\lambda_{p}\left(n\right)}}\right)=\chi_{H}\left(\left[B_{p}\left(n\right)\right]_{p^{m}}\right)
\end{equation}
By \textbf{Lemma \ref{lem:Unique rising continuation and p-adic functional equation of Chi_H}},
we have that: 
\begin{equation}
\chi_{H}\left(B_{p}\left(n\right)\right)\overset{\mathbb{Z}_{q}}{=}\lim_{m\rightarrow\infty}\chi_{H}\left(\left[B_{p}\left(n\right)\right]_{p^{m}}\right)=\lim_{m\rightarrow\infty}\frac{1-\left(M_{H}\left(\mathbf{j}\right)\right)^{m}}{1-M_{H}\left(n\right)}\chi_{H}\left(n\right)
\end{equation}
Finally, because $H$ is semi-basic, it follows by \textbf{Proposition
\ref{prop:q-adic behavior of M_H of j as the number of non-zero digits tends to infinity}
}that $\left|M_{H}\left(\mathbf{j}\right)\right|_{q}<1$ (since $\mathbf{j}$
has a non-zero entry), and hence, that $\left(M_{H}\left(\mathbf{j}\right)\right)^{m}$
tends to $0$ in $\mathbb{Z}_{q}$ as $m\rightarrow\infty$. Thus:
\begin{equation}
\chi_{H}\left(B_{p}\left(n\right)\right)\overset{\mathbb{Z}_{q}}{=}\lim_{m\rightarrow\infty}\frac{1-\left(M_{H}\left(\mathbf{j}\right)\right)^{m}}{1-M_{H}\left(n\right)}\chi_{H}\left(n\right)\overset{\mathbb{Z}_{q}}{=}\frac{\chi_{H}\left(n\right)}{1-M_{H}\left(n\right)}
\end{equation}

Q.E.D. 
\begin{rem}
As seen in the above proof, (\ref{eq:Chi_H B functional equation})
is really just the geometric series summation formula, convergent
in $\mathbb{Z}_{q}$ whenever $\left|M_{H}\left(n\right)\right|_{q_{H}}<1$.
In computing $\chi_{H}\left(B_{p}\left(n\right)\right)$ by using
(\ref{eq:Rising Continuity Formula for Chi_H}) with $\mathfrak{z}=B_{p}\left(n\right)$
and $n\geq1$, the series in (\ref{eq:Formal rising limit of Chi_H})
will reduce to a geometric series, one which converges in $\mathbb{Z}_{q}$
to $\chi_{H}\left(n\right)/\left(1-M_{H}\left(n\right)\right)$. Additionally,
for any $n\geq1$ for which the \emph{archimedean }absolute value
$\left|M_{H}\left(n\right)\right|$ is less than $1$, the \emph{universality}
of the geometric series formula: 
\begin{equation}
\sum_{n=0}^{\infty}x^{n}=\frac{1}{1-x}
\end{equation}
guarantees that the limit to which the series in (\ref{eq:Formal rising limit of Chi_H})
converges to in $\mathbb{R}$ in this case is the same as the limit
to which it converges in $\mathbb{Z}_{q_{H}}$ \cite{Conrad on p-adic series}.
This is worth mentioning because, as we shall see, the values of $n\geq1$
for which $\chi_{H}\left(n\right)/\left(1-M_{H}\left(n\right)\right)$
is a periodic point of $H$ in $\mathbb{N}_{1}$ are exactly those
values of $n\geq1$ for which the non-archimedean absolute value $\left|M_{H}\left(n\right)\right|_{q_{H}}$
is $<1$. 
\end{rem}
\vphantom{}

Now we introduce the final bit of qualitative controls on $H$'s behavior
which we will needed to prove the \textbf{Correspondence Principle}.
First, however, a proposition about $H$'s extendibility to the $p$-adics.
\begin{prop}
\label{prop:p-adic extension of H}Let $H:\mathbb{Z}\rightarrow\mathbb{Z}$
be a $p$-Hydra map. Then, $H$ admits an extension to a continuous
map $\mathbb{Z}_{p}\rightarrow\mathbb{Z}_{p}$ defined by:
\begin{equation}
H\left(\mathfrak{z}\right)\overset{\textrm{def}}{=}\sum_{j=0}^{p-1}\left[\mathfrak{z}\overset{p}{\equiv}j\right]\frac{a_{j}\mathfrak{z}+b_{j}}{d_{j}},\textrm{ }\forall\mathfrak{z}\in\mathbb{Z}_{p}\label{eq:p-adic extension of H}
\end{equation}
meaning that, $H\left(\mathfrak{z}\right)=H_{j}\left(\mathfrak{z}\right)$
for all $j\in\left\{ 0,\ldots,p-1\right\} $ and all $\mathfrak{z}\in j+p\mathbb{Z}_{p}$.
Moreover, the function $f:\mathbb{Z}_{p}\rightarrow\mathbb{Z}_{p}$
defined by:
\begin{equation}
f\left(\mathfrak{z}\right)=\sum_{j=0}^{p-1}\left[\mathfrak{z}\overset{p}{\equiv}j\right]\frac{a_{j}\mathfrak{z}+b_{j}}{d_{j}},\textrm{ }\forall\mathfrak{z}\in\mathbb{Z}_{p}
\end{equation}
is the unique continuous function $\mathbb{Z}_{p}\rightarrow\mathbb{Z}_{p}$
whose restriction to $\mathbb{Z}$ is equal to $H$.
\end{prop}
Proof: Let $f\left(\mathfrak{z}\right)=\sum_{j=0}^{p-1}\left[\mathfrak{z}\overset{p}{\equiv}j\right]\frac{a_{j}\mathfrak{z}+b_{j}}{d_{j}}$;
continuity is immediate, as is the fact that the restriction of $f$
to $\mathbb{Z}$ is equal to $H$. As for uniqueness, if $f:\mathbb{Z}_{p}\rightarrow\mathbb{Z}_{p}$
is any continuous function whose restriction on $\mathbb{Z}$ is equal
to $H$, the density of $\mathbb{Z}$ in $\mathbb{Z}_{p}$ coupled
with the continuity of $f$ then forces $f\left(\mathfrak{z}\right)=\sum_{j=0}^{p-1}\left[\mathfrak{z}\overset{p}{\equiv}j\right]\frac{a_{j}\mathfrak{z}+b_{j}}{d_{j}}$
for all $\mathfrak{z}\in\mathbb{Z}_{p}$. Hence, the extension of
$H$ to the $p$-adics is unique.

Q.E.D.

\vphantom{}

As mentioned in \textbf{Subsection \ref{subsec:2.1.2 It's-Probably-True}},
the $p$-adic extension of a Hydra map given by \textbf{Proposition
\ref{prop:p-adic extension of H}} has been studied in conjunction
with investigations of the Collatz Conjecture and the $qx+1$ maps,
however\textemdash paradoxically\textemdash doing so seems to lose
most information about the map, viz. \textbf{Theorem \ref{thm:shift map}},
which shows that the $2$-adic extension of the Shortened $qx+1$
map is equivalent to the shift map. For us, however, the utility of
the $p$-adic extension of $H$ is in keeping track of what happens
if we apply a ``wrong'' choice of branch to a given rational number,
as described below.
\begin{defn}[\textbf{Wrong Values and Propriety}]
\label{def:1D wrong values and properiety}\index{Hydra map!wrong value}\index{Hydra map!proper}
Let $H$ be a $p$-Hydra map.

\vphantom{}

I. We say a $p$-adic number $\mathfrak{y}\in\mathbb{Q}_{p}$ is a
\textbf{wrong value for $H$ }whenever there is a $\mathbf{j}\in\textrm{String}\left(p\right)$
and a $\mathfrak{z}\in\mathbb{Z}_{p}$ so that $\mathfrak{y}=H_{\mathbf{j}}\left(\mathfrak{z}\right)$
and $H_{\mathbf{j}}\left(\mathfrak{z}\right)\neq H^{\circ\left|\mathbf{j}\right|}\left(\mathfrak{z}\right)$.
We then call $\mathfrak{z}$ a \textbf{seed }of $\mathfrak{y}$.

\vphantom{}

II. We say $H$ is \textbf{proper }if $\left|H_{j}\left(\mathfrak{z}\right)\right|_{p}>1$
holds for any $\mathfrak{z}\in\mathbb{Z}_{p}$ and any $j\in\mathbb{Z}/p\mathbb{Z}$
so that $\left[\mathfrak{z}\right]_{p}\neq j$.
\end{defn}
\begin{rem}
If $\mathfrak{y}$ is a wrong value of $H$ with seed $\mathfrak{z}$
and string $\mathbf{j}$ so that $\mathfrak{y}=H_{\mathbf{j}}\left(\mathfrak{z}\right)$,
note that for any $\mathbf{i}\in\textrm{String}\left(p\right)$, $H_{\mathbf{i}}\left(\mathfrak{y}\right)=H_{\mathbf{i}}\left(H_{\mathbf{j}}\left(\mathfrak{z}\right)\right)=H_{\mathbf{i}\wedge\mathbf{j}}\left(\mathfrak{z}\right)$
will \emph{also }be a wrong value of $H$ with seed $\mathfrak{z}$.
Thus, branches of $H$ will always map wrong values to wrong values.
\end{rem}
\begin{rem}
Propriety is a stronger version of the integrality condition. If $H$
is an integral Hydra map, then $H_{j}\left(n\right)\notin\mathbb{Z}$
for any $n\in\mathbb{Z}$ and $j\in\mathbb{Z}/p\mathbb{Z}$ with $\left[n\right]_{p}\neq j$.
On the other hand, if $H$ is proper, then $H_{j}\left(\mathfrak{z}\right)\notin\mathbb{Z}_{p}$
for any $\mathfrak{z}\in\mathbb{Z}_{p}$ and $j\in\mathbb{Z}/p\mathbb{Z}$
with $\left[\mathfrak{z}\right]_{p}\neq j$.
\end{rem}
\vphantom{}

Like with the qualitative definitions from the previous subsection,
we some technical results to show that our new definitions are useful.
\begin{prop}
\label{prop:Q_p / Z_p prop}Let $H$ be any $p$-Hydra map. Then,
$H_{j}\left(\mathbb{Q}_{p}\backslash\mathbb{Z}_{p}\right)\subseteq\mathbb{Q}_{p}\backslash\mathbb{Z}_{p}$
for all $j\in\mathbb{Z}/p\mathbb{Z}$.
\end{prop}
\begin{rem}
Here, we are viewing the $H_{j}$s as functions $\mathbb{Q}_{p}\rightarrow\mathbb{Q}_{p}$.
\end{rem}
Proof: Let $\mathfrak{y}\in\mathbb{Q}_{p}\backslash\mathbb{Z}_{p}$
and $j\in\mathbb{Z}/p\mathbb{Z}$ be arbitrary. Note that as a $p$-Hydra
map, we have $p\mid d_{j}$ and $\gcd\left(a_{j},d_{j}\right)=1$.
Consequently, $p$ does not divide $a_{j}$. and so, $\left|a_{j}\right|_{p}=1$.
Thus:
\begin{equation}
\left|a_{j}\mathfrak{y}\right|_{p}=\left|\mathfrak{y}\right|_{p}>1
\end{equation}
So, $a_{j}\mathfrak{y}\in\mathbb{Q}_{p}\backslash\mathbb{Z}_{p}$.
Since $b_{j}$ is a rational integer, the $p$-adic ultrametric inequality
guarantees that $a_{j}\mathfrak{y}+b_{j}\in\mathbb{Q}_{p}\backslash\mathbb{Z}_{p}$.
Finally, $p\mid d_{j}$ implies $\left|d_{j}\right|_{p}<1$, and so:
\begin{equation}
\left|H_{j}\left(\mathfrak{y}\right)\right|_{p}=\left|\frac{a_{j}\mathfrak{y}+b_{j}}{d_{j}}\right|_{p}\geq\left|a_{j}\mathfrak{y}+b_{j}\right|_{p}>1
\end{equation}
which shows that $H_{j}\left(\mathfrak{y}\right)\in\mathbb{Q}_{p}\backslash\mathbb{Z}_{p}$.

Q.E.D.

\vphantom{}

The lemmata given below are the key ingredients in our proofs of the
\textbf{Correspondence Principle}. These will utilize \textbf{Proposition
\ref{prop:Q_p / Z_p prop}}.
\begin{lem}
\label{lem:integrality lemma}Let $H$ be a semi-basic $p$-Hydra
map, where $p$ is prime. Then, $H$ is proper if and only if $H$
is integral.
\end{lem}
Proof: Let $H$ be semi-basic.

I. (Proper implies integral) Suppose $H$ is proper, and let $n\in\mathbb{Z}$
be arbitrary. Then, clearly, when $j=\left[n\right]_{p}$, we have
that $H_{j}\left(n\right)\in\mathbb{Z}$. So, suppose $j\in\mathbb{Z}/p\mathbb{Z}$
satisfies $j\neq\left[n\right]_{p}$. Since $H$ is proper, the fact
that $n\in\mathbb{Z}\subseteq\mathbb{Z}_{p}$ and $j\neq\left[n\right]_{p}$
tells us that $\left|H_{j}\left(n\right)\right|_{p}>1$. Hence, $H_{j}\left(n\right)$
is not a $p$-adic integer, and thus, cannot be a rational integer
either. This proves $H$ is integral.

\vphantom{}

II. (Integral implies proper) Suppose $H$ is integral, and\textemdash by
way of contradiction\textemdash suppose $H$ is \emph{not}\textbf{
}proper. Then, there is a $\mathfrak{z}\in\mathbb{Z}_{p}$ and a $j\in\mathbb{Z}/p\mathbb{Z}$
with $j\neq\left[\mathfrak{z}\right]_{p}$ so that $\left|H_{j}\left(\mathfrak{z}\right)\right|_{p}\leq1$.

Now, writing $\mathfrak{z}=\sum_{n=0}^{\infty}c_{n}p^{n}$:
\begin{align*}
H_{j}\left(\mathfrak{z}\right) & =\frac{a_{j}\sum_{n=0}^{\infty}c_{n}p^{n}+b_{j}}{d_{j}}\\
 & =\frac{a_{j}c_{0}+b_{j}}{d_{j}}+\frac{1}{d_{j}}\sum_{n=1}^{\infty}c_{n}p^{n}\\
\left(c_{0}=\left[\mathfrak{z}\right]_{p}\right); & =H_{j}\left(\left[\mathfrak{z}\right]_{p}\right)+\frac{p}{d_{j}}\underbrace{\sum_{n=1}^{\infty}c_{n}p^{n-1}}_{\textrm{call this }\mathfrak{y}}
\end{align*}
Because $H$ is a $p$-Hydra map, $d_{j}$ must be a divisor of $p$.
The primality of $p$ then forces $d_{j}$ to be either $1$ or $p$.
In either case, we have that $p\mathfrak{y}/d_{j}$ is an element
of $\mathbb{Z}_{p}$. Our contradictory assumption $\left|H_{j}\left(\mathfrak{z}\right)\right|_{p}\leq1$
tells us that $H_{j}\left(\mathfrak{z}\right)$ is also in $\mathbb{Z}_{p}$,
and so:
\begin{equation}
H_{j}\left(\left[\mathfrak{z}\right]_{p}\right)=H_{j}\left(\mathfrak{z}\right)-\frac{p\mathfrak{y}}{d_{j}}\in\mathbb{Z}_{p}
\end{equation}
Since $H$ was given to be integral, $j\neq\left[\mathfrak{z}\right]_{p}$
implies $H_{j}\left(\left[\mathfrak{z}\right]_{p}\right)\notin\mathbb{Z}$.
As such, $H_{j}\left(\left[\mathfrak{z}\right]_{p}\right)$ is a $p$-adic
integer which is \emph{not }a rational integer. Since $H_{j}\left(\left[\mathfrak{z}\right]_{p}\right)$
is a non-integer rational number which is a $p$-adic integer, the
denominator of $H_{j}\left(\left[\mathfrak{z}\right]_{p}\right)=\frac{a_{j}\left[\mathfrak{z}\right]_{p}+b_{j}}{d_{j}}$
must be divisible by some prime $q\neq p$, and hence, $q\mid d_{j}$.
However, we saw that $p$ being prime forced $d_{j}\in\left\{ 1,p\right\} $\textemdash this
is impossible!

Thus, it must be that $H$ is proper.

Q.E.D.
\begin{lem}
\label{lem:wrong values lemma}Let $H$ be a proper $p$-Hydra map.
All wrong values of $H$ are elements of $\mathbb{Q}_{p}\backslash\mathbb{Z}_{p}$.
\end{lem}
Proof: Let $H$ be a proper $p$-Hydra map, let $\mathfrak{z}\in\mathbb{Z}_{p}$,
and let $i\in\mathbb{Z}/p\mathbb{Z}$ be such that $\left[\mathfrak{z}\right]_{p}\neq i$.
Then, by definition of properness, the quantity: 
\begin{equation}
H_{i}\left(\mathfrak{z}\right)=\frac{a_{i}\mathfrak{z}+b_{i}}{d_{i}}
\end{equation}
has $p$-adic absolute value $>1$. By \textbf{Proposition \ref{prop:Q_p / Z_p prop}},
this then forces $H_{\mathbf{j}}\left(H_{i}\left(\mathfrak{z}\right)\right)$
to be an element of $\mathbb{Q}_{p}\backslash\mathbb{Z}_{p}$ for
all $\mathbf{j}\in\textrm{String}\left(p\right)$. Since every wrong
value with seed $\mathfrak{z}$ is of the form $H_{\mathbf{j}}\left(H_{i}\left(\mathfrak{z}\right)\right)$
for some $\mathbf{j}\in\textrm{String}\left(p\right)$, some $\mathfrak{z}\in\mathbb{Z}_{p}$,
and some $i\in\mathbb{Z}/p\mathbb{Z}$ for which $\left[\mathfrak{z}\right]_{p}\neq i$,
this shows that every wrong value of $H$ is in $\mathbb{Q}_{p}\backslash\mathbb{Z}_{p}$.
So, $H$ is proper.

Q.E.D.
\begin{lem}
\label{lem:properness lemma}Let $H$ be a proper $p$-Hydra map,
let $\mathfrak{z}\in\mathbb{Z}_{p}$, and let $\mathbf{j}\in\textrm{String}\left(p\right)$.
If $H_{\mathbf{j}}\left(\mathfrak{z}\right)=\mathfrak{z}$, then $H^{\circ\left|\mathbf{j}\right|}\left(\mathfrak{z}\right)=\mathfrak{z}$.
\end{lem}
Proof: Let $H$, $\mathfrak{z}$, and $\mathbf{j}$ be as given. By
way of contradiction, suppose $H^{\circ\left|\mathbf{j}\right|}\left(\mathfrak{z}\right)\neq\mathfrak{z}$.
But then $\mathfrak{z}=H_{\mathbf{j}}\left(\mathfrak{z}\right)$ implies
$H^{\circ\left|\mathbf{j}\right|}\left(\mathfrak{z}\right)\neq H_{\mathbf{j}}\left(\mathfrak{z}\right)$.
Hence, $H_{\mathbf{j}}\left(\mathfrak{z}\right)$ is a wrong value
of $H$ with seed $\mathfrak{z}$. \textbf{Lemma} \ref{lem:wrong values lemma}
then forces $H_{\mathbf{j}}\left(\mathfrak{z}\right)\in\mathbb{Q}_{p}\backslash\mathbb{Z}_{p}$.
However, $H_{\mathbf{j}}\left(\mathfrak{z}\right)=\mathfrak{z}$,
and $\mathfrak{z}$ was given to be in $\mathbb{Z}_{p}$. This is
impossible!

Consequently, $H_{\mathbf{j}}\left(\mathfrak{z}\right)=\mathfrak{z}$
implies $H^{\circ\left|\mathbf{j}\right|}\left(\mathfrak{z}\right)=\mathfrak{z}$.

Q.E.D.

\vphantom{}

Next, we need an essentially trivial observation: if $H_{\mathbf{j}}\left(n\right)=n$,
$H\left(0\right)=0$, and $\left|\mathbf{j}\right|\geq2$, then the
$p$-ity vector $\mathbf{j}$ must contain \emph{at least one} non-zero
entry. 
\begin{prop}
\label{prop:essentially trivial observation}Let $H$ be a non-degenerate
$p$-Hydra map which fixes $0$, let $n\in\mathbb{Z}$ be a periodic
point of $H$, and let $\Omega$ be the cycle of $H$ which contains
$n$. Letting $\mathbf{j}=\left(j_{1},\ldots,j_{\left|\Omega\right|}\right)\in\left(\mathbb{Z}/p\mathbb{Z}\right)^{\left|\Omega\right|}$
be the unique string so that: 
\begin{equation}
H_{\mathbf{j}}\left(n\right)=H^{\circ\left|\Omega\right|}\left(n\right)=n
\end{equation}
we have that either $\left|\Omega\right|=1$ (i.e., $n$ is a fixed
point of $H$), or, there is an $\ell\in\left\{ 1,\ldots,\left|\Omega\right|\right\} $
so that $j_{\ell}\neq0$. 
\end{prop}
Proof: Let $H$, $n$, and $\mathbf{j}$ be as given, and suppose
$\left|\Omega\right|\geq2$. Since $H$ is non-degenerate, $a_{j}\neq1$
for any $j\in\left\{ 1,\ldots,p-1\right\} $. Now, \emph{by way of
contradiction}, suppose that $j_{\ell}=0$ for every $\ell$.

Since $\mathbf{j}$ contains only $0$s, we can write: 
\begin{equation}
n=H_{\mathbf{j}}\left(n\right)=H_{0}^{\circ\left|\Omega\right|}\left(n\right)
\end{equation}
Because $H_{0}$ is an affine linear map of the form $H_{0}\left(x\right)=\frac{a_{0}x+b_{0}}{d_{0}}$
where $b_{0}=0$, $x=0$ is the one and only fixed point of $H_{0}$
in $\mathbb{R}$. In fact, the same is true of $H_{0}^{\circ m}\left(x\right)$
for all $m\geq1$. Thus, $n=H_{0}^{\circ\left|\Omega\right|}\left(n\right)$
forces $n=0$. So, $H\left(n\right)=H_{0}\left(n\right)=n$, which
is to say, the cycle $\Omega$ to which $n$ belongs contains only
one element: $n$ itself. This contradicts the assumption that $\left|\Omega\right|\geq2$.

Q.E.D.

\vphantom{}

Now, at last, we can prove the \textbf{Correspondence Principle}.
This is done in four different ways, the last of which (\textbf{Corollary
\ref{cor:CP v4}}) is the simplest and most elegant. The first version\textemdash \textbf{Theorem
\ref{thm:CP v1}}, proved below\textemdash is more general than \textbf{Corollary
\ref{cor:CP v4}}, which is proved by using \textbf{Theorem \ref{thm:CP v1}
}as a stepping stone. Whereas \textbf{Corollary \ref{cor:CP v4}}
establishes a correspondence between the values attained by $\chi_{H}$
and the periodic \emph{points }of $H$, \textbf{Theorem \ref{thm:CP v1}
}establishes a slightly weaker correspondence between the values attained
by $\chi_{H}$ and the \emph{cycles }of $H$. \textbf{Corollary \ref{cor:CP v4}}
refines this correspondence to one with individual periodic points
of $\chi_{H}$. The other two versions of the Correspondence Principle
are closer to the Böhm-Sontacchi Criterion, with \textbf{Corollary
\ref{cor:CP v3}} providing Böhm-Sontacchi-style diophantine equations
for the non-zero periodic points of $H$. In that respect, \textbf{Corollary
\ref{cor:CP v3} }is not entirely new; a similar diophantine equation
characterization of periodic points for certain general Collatz-type
maps on $\mathbb{Z}$ can be found in Matthews' slides \cite{Matthews' slides}.
The main innovations of my approach is the reformulation of these
diophantine equations in terms of the function $\chi_{H}$, along
with the tantalizing half-complete characterization $\chi_{H}$ provides
for the divergent trajectories of $H$ (\textbf{Theorem \ref{thm:Divergent trajectories come from irrational z}}). 
\begin{thm}[\textbf{Correspondence Principle, Ver. 1}]
\label{thm:CP v1}Let\index{Correspondence Principle}$H$ be semi-basic.
Then:

\vphantom{}

I. Let $\Omega$ be any cycle of $H$ in $\mathbb{Z}$, with $\left|\Omega\right|\geq2$.
Then, there exist $x\in\Omega$ and $n\in\mathbb{N}_{1}$ so that:
\begin{equation}
\chi_{H}\left(B_{p}\left(n\right)\right)\overset{\mathbb{Z}_{q_{H}}}{=}x
\end{equation}
That is, there is an $n$ so that the infinite series defining $\chi_{H}\left(B_{p}\left(n\right)\right)$
converges\footnote{This series is obtained by evaluating $\chi_{H}\left(\left[B_{p}\left(n\right)\right]_{p^{m}}\right)$
using \textbf{Proposition \ref{prop:Explicit formula for Chi_H of bold j}},
and then taking the limit in $\mathbb{Z}_{q}$ as $m\rightarrow\infty$.} in $\mathbb{Z}_{q_{H}}$ to $x$.

\vphantom{}

II. Suppose also that $H$ is integral, and let $n\in\mathbb{N}_{1}$.
If the rational number $\chi_{H}\left(B_{p}\left(n\right)\right)$
is in $\mathbb{Z}_{p}$, then $\chi_{H}\left(B_{p}\left(n\right)\right)$
is a periodic point of $H$ in $\mathbb{Z}_{p}$. In particular, if
$\chi_{H}\left(B_{p}\left(n\right)\right)$ is in $\mathbb{Z}$, then
$\chi_{H}\left(B_{p}\left(n\right)\right)$ is a periodic point of
$H$ in $\mathbb{Z}$.
\end{thm}
Proof:

I. Let $\Omega$ be a cycle of $H$ in $\mathbb{Z}$, with $\left|\Omega\right|\geq2$.
Given any periodic point $x\in\Omega$ of $H$, there is going to
be a string $\mathbf{j}\in\textrm{String}\left(p\right)$ of length
$\left|\Omega\right|$ so that $H_{\mathbf{j}}\left(x\right)=x$.
In particular, since $\left|\Omega\right|\geq2$,\textbf{ Proposition
\ref{prop:essentially trivial observation}} tells us that $\mathbf{j}$
contains \emph{at least one} non-zero entry. A moment's thought shows
that for any $x^{\prime}\in\Omega$ \emph{other }than $x$, the entries
of the string $\mathbf{i}$ for which $H_{\mathbf{i}}\left(x^{\prime}\right)=x^{\prime}$
must be a cyclic permutation of the entries of $\mathbf{j}$.

For example, for the Shortened Collatz Map and the cycle $\left\{ 1,2\right\} $,
our branches are: 
\begin{align*}
H_{0}\left(x\right) & =\frac{x}{2}\\
H_{1}\left(x\right) & =\frac{3x+1}{2}
\end{align*}
The composition sequence which sends $1$ back to itself first applies
$H_{1}$ to $1$, followed by $H_{0}$: 
\begin{equation}
H_{0}\left(H_{1}\left(1\right)\right)=H_{0}\left(2\right)=1
\end{equation}
On the other hand, the composition sequence which sends $2$ back
to itself first applies $H_{0}$ to $2$, followed by $H_{1}$: 
\begin{equation}
H_{1}\left(H_{0}\left(2\right)\right)=H_{1}\left(1\right)=2
\end{equation}
and the strings $\left(0,1\right)$ and $\left(1,0\right)$ are cyclic
permutations of one another.

In this way, note that there must then exist an $x\in\Omega$ and
a $\mathbf{j}$ (for which $H_{\mathbf{j}}\left(x\right)=x$) so that
\emph{the right-most entry of $\mathbf{j}$ is non-zero}. From this
point forward, we will fix $x$ and $\mathbf{j}$ so that this condition
is satisfied.

That being done, the next observation to make is that because $H_{\mathbf{j}}$
sends $x$ to $x$, further applications of $H_{\mathbf{j}}$ will
keep sending $x$ back to $x$. Expressing this in terms of concatenation
yields: 
\begin{equation}
x=H_{\mathbf{j}}^{\circ m}\left(x\right)=\underbrace{\left(H_{\mathbf{j}}\circ\cdots\circ H_{\mathbf{j}}\right)}_{m\textrm{ times}}\left(x\right)=H_{\mathbf{j}^{\wedge m}}\left(x\right)
\end{equation}
Writing the affine map $H_{\mathbf{j}^{\wedge m}}$ in affine form
gives us: 
\begin{equation}
x=H_{\mathbf{j}^{\wedge m}}\left(x\right)=M_{H}\left(\mathbf{j}^{\wedge m}\right)x+\chi_{H}\left(\mathbf{j}^{\wedge m}\right)=\left(M_{H}\left(\mathbf{j}\right)\right)^{m}x+\chi_{H}\left(\mathbf{j}^{\wedge m}\right)
\end{equation}
where the right-most equality follows from $M_{H}$'s concatenation
identity (\textbf{Proposition \ref{prop:M_H concatenation identity}}).
Since $\mathbf{j}$ the right-most entry of $\mathbf{j}$ is non-zero,
\textbf{Proposition \ref{prop:q-adic behavior of M_H of j as the number of non-zero digits tends to infinity}}
tells us that $\left|M_{H}\left(\mathbf{j}\right)\right|_{q_{H}}<1$.
As such, $\left|\left(M_{H}\left(\mathbf{j}\right)\right)^{m}\right|_{q}\rightarrow0$
as $m\rightarrow\infty$. So, taking limits as $m\rightarrow\infty$,
we get: 
\begin{equation}
x\overset{\mathbb{Z}_{q_{H}}}{=}\lim_{m\rightarrow\infty}\left(\left(M_{H}\left(\mathbf{j}\right)\right)^{m}x+\chi_{H}\left(\mathbf{j}^{\wedge m}\right)\right)\overset{\mathbb{Z}_{q_{H}}}{=}\lim_{m\rightarrow\infty}\chi_{H}\left(\mathbf{j}^{\wedge m}\right)
\end{equation}

Now, let $n$ be the integer represented by $\mathbf{j}$. Since $\mathbf{j}$'s
right-most entry is non-zero, $\mathbf{j}$ is then the \emph{shortest}
string representing $n$; moreover, $n$ is non-zero. In particular,
we have that $\lambda_{p}\left(n\right)=\left|\mathbf{j}\right|>0$,
and so, using \textbf{Proposition \ref{prop:Concatenation exponentiation}},\textbf{
}we can write: 
\begin{equation}
\mathbf{j}^{\wedge m}\sim n\frac{1-p^{m\lambda_{p}\left(n\right)}}{1-p^{\lambda_{p}\left(n\right)}},\textrm{ }\forall m\geq1
\end{equation}
So, like in the proof of \textbf{Lemma \ref{lem:Chi_H o B_p functional equation}},
we then have: 
\begin{align*}
x & \overset{\mathbb{Z}_{q_{H}}}{=}\lim_{m\rightarrow\infty}\chi_{H}\left(\mathbf{j}^{\wedge m}\right)\\
 & \overset{\mathbb{Z}_{q_{H}}}{=}\lim_{m\rightarrow\infty}\chi_{H}\left(n\frac{1-p^{m\lambda_{p}\left(n\right)}}{1-p^{\lambda_{p}\left(n\right)}}\right)\\
 & \overset{\mathbb{Z}_{q_{H}}}{=}\chi_{H}\left(\frac{n}{1-p^{\lambda_{p}\left(n\right)}}\right)\\
 & =\chi_{H}\left(B_{p}\left(n\right)\right)
\end{align*}
This proves the existence of the desired $x$ and $n$.

\vphantom{}

II. Suppose $H$ is integral. First, before we assume the hypotheses
of (II), let us do a brief computation. 
\begin{claim}
\label{claim:2.2}Let $n$ be any integer $\geq1$, and let $\mathbf{j}\in\textrm{String}\left(p\right)$
be the shortest string representing $n$. Then, $H_{\mathbf{j}}$
is a continuous map $\mathbb{Z}_{q_{H}}\rightarrow\mathbb{Z}_{q_{H}}$,
and, moreover, the $q_{H}$-adic integer $\chi_{H}\left(B_{p}\left(n\right)\right)$
is fixed by $H_{\mathbf{j}}$.

Proof of claim: Let $n$ and $\mathbf{j}$ be as given. Using $\chi_{H}$'s
concatenation identity (\textbf{Lemma \ref{lem:Chi_H concatenation identity}}),
we can write: 
\begin{equation}
\chi_{H}\left(\mathbf{j}^{\wedge k}\right)=H_{\mathbf{j}}\left(\chi_{H}\left(\mathbf{j}^{\wedge\left(k-1\right)}\right)\right)=\cdots=H_{\mathbf{j}^{\wedge\left(k-1\right)}}\left(\chi_{H}\left(\mathbf{j}\right)\right)\label{eq:Concatenation identity for Chi_H}
\end{equation}
By \textbf{Proposition \ref{prop:Concatenation exponentiation}},
we know that $B_{p}\left(n\right)$ represents $\lim_{k\rightarrow\infty}\mathbf{j}^{\wedge k}$
(i.e. $\textrm{DigSum}_{p}\left(\lim_{k\rightarrow\infty}\mathbf{j}^{\wedge k}\right)=B_{p}\left(n\right)$).
In particular, we have that: 
\begin{equation}
\left[B_{p}\left(n\right)\right]_{p^{k}}=\mathbf{j}^{\wedge k}=n\frac{1-p^{k\lambda_{p}\left(n\right)}}{1-p^{\lambda_{p}\left(n\right)}}
\end{equation}
Letting $k\rightarrow\infty$, the limit (\ref{eq:Rising Continuity Formula for Chi_H})
tells us that: 
\begin{equation}
\lim_{k\rightarrow\infty}H_{\mathbf{j}^{\wedge\left(k-1\right)}}\left(\chi_{H}\left(n\right)\right)\overset{\mathbb{Z}_{q_{H}}}{=}\chi_{H}\left(\frac{n}{1-p^{\lambda_{p}\left(n\right)}}\right)=\chi_{H}\left(B_{p}\left(n\right)\right)\label{eq:Iterating H_bold-j on Chi_H}
\end{equation}
The heuristic here is that: 
\begin{equation}
H_{\mathbf{j}}\left(\lim_{k\rightarrow\infty}H_{\mathbf{j}^{\wedge\left(k-1\right)}}\right)=\lim_{k\rightarrow\infty}H_{\mathbf{j}^{\wedge\left(k-1\right)}}
\end{equation}
and, hence, that $\chi_{H}\left(B_{p}\left(n\right)\right)$ will
be fixed by $H_{\mathbf{j}}$. To make this rigorous, note that $\gcd\left(a_{j},p\right)=1$
for all $j$, because $H$ is semi-basic. As such, for any element
of $\textrm{String}\left(p\right)$\textemdash such as our $\mathbf{j}$\textemdash the
quantities $H_{\mathbf{j}}^{\prime}\left(0\right)$ and $H_{\mathbf{j}}\left(0\right)$
are then rational numbers which lie in $\mathbb{Z}_{q_{H}}$. This
guarantees that the function $\mathfrak{z}\mapsto H_{\mathbf{j}}^{\prime}\left(0\right)\mathfrak{z}+H_{\mathbf{j}}\left(0\right)$
(which is, of course $H_{\mathbf{j}}\left(\mathfrak{z}\right)$) is
then a well-defined \emph{continuous} map $\mathbb{Z}_{q_{H}}\rightarrow\mathbb{Z}_{q_{H}}$.
Applying $H_{\mathbf{j}}$ to $\chi_{H}\left(B_{p}\left(n\right)\right)$,
we obtain: 
\begin{align*}
H_{\mathbf{j}}\left(\chi_{H}\left(B_{p}\left(n\right)\right)\right) & \overset{\mathbb{Z}_{q_{H}}}{=}H_{\mathbf{j}}\left(\lim_{k\rightarrow\infty}\chi_{H}\left(\mathbf{j}^{\wedge\left(k-1\right)}\right)\right)\\
 & \overset{\mathbb{Z}_{q_{H}}}{=}\lim_{k\rightarrow\infty}H_{\mathbf{j}}\left(\chi_{H}\left(\mathbf{j}^{\wedge\left(k-1\right)}\right)\right)\\
 & \overset{\mathbb{Z}_{q_{H}}}{=}\lim_{k\rightarrow\infty}H_{\mathbf{j}}\left(H_{\mathbf{j}^{\wedge\left(k-1\right)}}\left(0\right)\right)\\
 & \overset{\mathbb{Z}_{q_{H}}}{=}\lim_{k\rightarrow\infty}H_{\mathbf{j}^{\wedge k}}\left(0\right)\\
 & \overset{\mathbb{Z}_{q_{H}}}{=}\lim_{k\rightarrow\infty}\chi_{H}\left(\mathbf{j}^{\wedge k}\right)\\
 & \overset{\mathbb{Z}_{q_{H}}}{=}\chi_{H}\left(B_{p}\left(n\right)\right)
\end{align*}
The interchange of $H_{\mathbf{j}}$ and $\lim_{k\rightarrow\infty}$
on the second line is justified by the continuity of $H_{\mathbf{j}}$
on $\mathbb{Z}_{q_{H}}$. This proves the claim.
\end{claim}
\vphantom{}

Now, let us actually \emph{assume} that $n$ satisfies the hypotheses
of (II): suppose\emph{ }$\chi_{H}\left(B_{p}\left(n\right)\right)\in\mathbb{Z}_{p}$.
By \textbf{Claim \ref{claim:2.2}}, we know that $H_{\mathbf{j}}\left(\chi_{H}\left(B_{p}\left(n\right)\right)\right)=\chi_{H}\left(B_{p}\left(n\right)\right)$.
Because $H$ was given to be integral, the primality of $p$ and $H$'s
semi-basicness guarantee that $H$ is proper (\textbf{Lemma \ref{lem:integrality lemma}}),
and so, we can apply \textbf{Lemma \ref{lem:properness lemma}} to
conclude that: 
\begin{equation}
H^{\circ\left|\mathbf{j}\right|}\left(\chi_{H}\left(B_{p}\left(n\right)\right)\right)=\chi_{H}\left(B_{p}\left(n\right)\right)
\end{equation}
where $\left|\mathbf{j}\right|=\lambda_{p}\left(n\right)\geq1$. So,
$\chi_{H}\left(B_{p}\left(n\right)\right)$ is a periodic point of
$H$ in $\mathbb{Z}_{p}$. In particular, if $\chi_{H}\left(B_{p}\left(n\right)\right)$
is in $\mathbb{Z}$, then every element of the cycle generated by
applying the $p$-adic extension of $H$ to $\chi_{H}\left(B_{p}\left(n\right)\right)$
is an element of $\mathbb{Z}$, and $\chi_{H}\left(B_{p}\left(n\right)\right)$
is then a periodic point of $H$ in $\mathbb{Z}$.

Q.E.D.

\vphantom{}

Using the functional equation in \textbf{Lemma \ref{lem:Chi_H o B_p functional equation}},
we can restate the Correspondence Principle in terms of $\chi_{H}$
and $M_{H}$. 
\begin{cor}[\textbf{Correspondence Principle, Ver. 2}]
\label{cor:CP v2}\index{Correspondence Principle} Suppose that
$H$ is semi-basic.

\vphantom{}

I. Let $\Omega\subseteq\mathbb{Z}$ be any cycle of $H$. Then, viewing
$\Omega\subseteq\mathbb{Z}$ as a subset of $\mathbb{Z}_{q_{H}}$,
the intersection $\chi_{H}\left(\mathbb{Z}_{p}\right)\cap\Omega$
is non-empty. Moreover, for every $x\in\chi_{H}\left(\mathbb{Z}_{p}\right)\cap\Omega$,
there is an $n\in\mathbb{N}_{1}$ so that: 
\begin{equation}
x=\frac{\chi_{H}\left(n\right)}{1-M_{H}\left(n\right)}
\end{equation}
\vphantom{}

II. Suppose in addition that $H$ is integral, and let $n\in\mathbb{N}_{1}$.
If the quantity $x$ given by: 
\begin{equation}
x=\frac{\chi_{H}\left(n\right)}{1-M_{H}\left(n\right)}
\end{equation}
is a $p$-adic integer, then $x$ is a periodic point of $H$ in $\mathbb{Z}_{p}$;
if $x$ is in $\mathbb{Z}$, then $x$ is a periodic point of $H$
in $\mathbb{Z}$. Moreover, if $x\in\mathbb{Z}$, then $x$ is positive
if and only if $M_{H}\left(n\right)<1$, and $x$ is negative if and
only if $M_{H}\left(n\right)>1$. 
\end{cor}
Proof: Re-write the results of \textbf{Theorem \ref{thm:CP v1}} using
\textbf{Lemma \ref{lem:Chi_H o B_p functional equation}}. The positivity/negativity
of $x$ stipulated in (II) follows by noting that $\chi_{H}\left(n\right)$
and $M_{H}\left(n\right)$ are positive rational numbers for all $n\in\mathbb{N}_{1}$.

Q.E.D.

\vphantom{}

Next, we have Böhm-Sontacchi-style diophantine equations characterizing
of $H$'s periodic points.
\begin{cor}[\textbf{Correspondence Principle, Ver. 3}]
\label{cor:CP v3}Let $H$ be a semi-basic $p$-Hydra map. \index{Böhm-Sontacchi criterion}\index{diophantine equation}Then:

\vphantom{}

I. Let $\Omega$ be a cycle of $H$ in $\mathbb{Z}$ containing at
least two elements. Then, there is an $x\in\Omega$, and a $\mathbf{j}\in\textrm{String}\left(p\right)$
of length $\geq2$ so that: 
\begin{equation}
x=\frac{\sum_{n=1}^{\left|\mathbf{j}\right|}\frac{b_{j_{n}}}{a_{j_{n}}}\left(\prod_{k=1}^{n}\mu_{j_{k}}\right)p^{\left|\mathbf{j}\right|-n}}{p^{\left|\mathbf{j}\right|}-\prod_{k=1}^{\left|\mathbf{j}\right|}\mu_{j_{k}}}\label{eq:Generalized Bohm-Sontacchi criterion-1}
\end{equation}

\vphantom{}

II. Suppose $H$ is integral, and let $x\left(\mathbf{j}\right)=x\left(j_{1},\ldots,j_{N}\right)$
denote the quantity: 
\begin{equation}
x\left(\mathbf{j}\right)\overset{\textrm{def}}{=}\frac{\sum_{n=1}^{\left|\mathbf{j}\right|}\frac{b_{j_{n}}}{a_{j_{n}}}\left(\prod_{k=1}^{n}\mu_{j_{k}}\right)p^{\left|\mathbf{j}\right|-n}}{p^{\left|\mathbf{j}\right|}-\prod_{k=1}^{\left|\mathbf{j}\right|}\mu_{j_{k}}},\textrm{ }\forall\mathbf{j}\in\textrm{String}\left(p\right)\label{eq:definition of x of the js (bold version)}
\end{equation}
If $\mathbf{j}\in\textrm{String}\left(p\right)$ has $\left|\mathbf{j}\right|\geq2$
and makes $x\left(\mathbf{j}\right)\in\mathbb{Z}$, then $x\left(\mathbf{j}\right)$
is a periodic point of $H$. 
\end{cor}
Proof: Use \textbf{Propositions \ref{prop:Explicit Formulas for M_H}}
and \textbf{\ref{prop:Explicit formula for Chi_H of bold j}} on (I)
and (II) of \textbf{Corollary \ref{cor:CP v2}}.

Q.E.D.

\vphantom{}

Finally, we can state and prove the most elegant version of the Correspondence
Principle for integral semi-basic Hydra maps. This is also the most
\emph{powerful }version of the Correspondence Principle, because it
will allow us to establish a connection between $\chi_{H}$ and $H$'s
\textbf{\emph{divergent points}}.
\begin{cor}[\textbf{Correspondence Principle, Ver. 4}]
\label{cor:CP v4} \index{Correspondence Principle}Let $H$ be an
integral semi-basic $p$-Hydra map. Then, the set of all non-zero
periodic points of $H$ in $\mathbb{Z}$ is equal to $\mathbb{Z}\cap\chi_{H}\left(\mathbb{Q}\cap\mathbb{Z}_{p}^{\prime}\right)$,
where $\chi_{H}\left(\mathbb{Q}\cap\mathbb{Z}_{p}^{\prime}\right)$
is viewed as a subset of $\mathbb{Z}_{q_{H}}$.
\end{cor}
\begin{rem}
The implication ``If $x\in\mathbb{Z}\backslash\left\{ 0\right\} $
is a periodic point, then $x\in\chi_{H}\left(\mathbb{Q}\cap\mathbb{Z}_{p}^{\prime}\right)$''
\emph{does not }require $H$ to be integral.
\end{rem}
Proof: First, note that $p$ is prime and since $H$ is integral and
semi-basic, \textbf{Lemma \ref{lem:integrality lemma}} tell us that
$H$ is proper.

I. Let $x$ be a non-zero periodic point of $H$, and let $\Omega$
be the unique cycle of $H$ in $\mathbb{Z}$ which contains $x$.
By Version 1 of the Correspondence Principle (\textbf{Theorem \ref{thm:CP v1}}),
there exists a $y\in\Omega$ and a $\mathfrak{z}=B_{p}\left(n\right)\subset\mathbb{Z}_{p}$
(for some $n\in\mathbb{N}_{1}$) so that $\chi_{H}\left(\mathfrak{z}\right)=y$.
Since $y$ is in $\Omega$, there is an $k\geq1$ so that $x=H^{\circ k}\left(y\right)$.
In particular, there is a \emph{unique} length $k$ string $\mathbf{i}\in\textrm{String}\left(p\right)$
so that $H_{\mathbf{i}}\left(y\right)=H^{\circ k}\left(y\right)=x$.

Now, let $\mathbf{j}\in\textrm{String}_{\infty}\left(p\right)$ represent
$\mathfrak{z}$; note that $\mathbf{j}$ is infinite and that its
entries are periodic. Using $\chi_{H}$'s concatenation identity (\textbf{Lemma
\ref{lem:Chi_H concatenation identity}}), we can write: 
\begin{equation}
x=H_{\mathbf{i}}\left(y\right)=H_{\mathbf{i}}\left(\chi_{H}\left(\mathfrak{z}\right)\right)=\chi_{H}\left(\mathbf{i}\wedge\mathbf{j}\right)
\end{equation}
Next, let $\mathfrak{x}$ denote the $p$-adic integer $\mathfrak{x}$
whose sequence of $p$-adic digits is $\mathbf{i}\wedge\mathbf{j}$;
note that $\mathfrak{x}$ is then \emph{not} an element of $\mathbb{N}_{0}$.
By the above, we have that $\chi_{H}\left(\mathfrak{x}\right)=x$,
and hence that $x\in\mathbb{Z}\cap\chi_{H}\left(\mathbb{Z}_{p}\right)$.

Finally, since $\mathfrak{z}=B_{p}\left(n\right)$, its $p$-adic
digits are periodic, which forces $\mathfrak{z}$ to be an element
of $\mathbb{Q}\cap\mathbb{Z}_{p}^{\prime}$. Indeed, $\mathfrak{z}$
is not in $\mathbb{N}_{0}$ because $n\geq1$, and $B_{p}\left(n\right)\in\mathbb{N}_{0}$
if and only if $n=0$. So, letting $m$ be the rational integer represented
by length-$k$ string $\mathbf{i}$, we have: 
\begin{equation}
\mathfrak{x}\sim\mathbf{i}\wedge\mathbf{j}\sim m+p^{\lambda_{p}\left(m\right)}\mathfrak{z}
\end{equation}
This shows that $\mathfrak{x}\in\mathbb{Q}\cap\mathbb{Z}_{p}^{\prime}$,
and hence, that $x=\chi_{H}\left(\mathfrak{x}\right)\in\mathbb{Z}\cap\chi_{H}\left(\mathbb{Q}\cap\mathbb{Z}_{p}^{\prime}\right)$.

\vphantom{}

II. Suppose $x\in\mathbb{Z}\cap\chi_{H}\left(\mathbb{Q}\cap\mathbb{Z}_{p}^{\prime}\right)$,
with $x=\chi_{H}\left(\mathfrak{z}\right)$ for some $\mathfrak{z}\in\mathbb{Q}\cap\mathbb{Z}_{p}^{\prime}$.
As a rational number which is both a $p$-adic integer \emph{and}
not an element of $\mathbb{N}_{0}$, the $p$-adic digits of $\mathfrak{z}$
are \emph{eventually} periodic. As such, there are integers $m$ and
$n$ (with $n\neq0$) so that: 
\begin{equation}
\mathfrak{z}=m+p^{\lambda_{p}\left(m\right)}B_{p}\left(n\right)
\end{equation}
Here, $n$'s $p$-adic digits generate the periodic part of $\mathfrak{z}$'s
digits, while $m$'s $p$-adic digits are the finite-length sequence
in $\mathfrak{z}$'s digits that occurs before the periodicity sets
in. Now, let $\mathbf{i}$ be the finite string representing $m$,
and let $\mathbf{j}$ be the infinite string representing $B_{p}\left(n\right)$.
Then, $\mathfrak{z}=\mathbf{i}\wedge\mathbf{j}$. So, by \textbf{Lemmata
\ref{lem:Chi_H concatenation identity}} and \textbf{\ref{lem:Chi_H o B_p functional equation}}
(the concatenation identity and $\chi_{H}\circ B_{p}$ functional
equation, respectively): 
\begin{equation}
x=\chi_{H}\left(\mathbf{i}\wedge\mathbf{j}\right)=H_{\mathbf{i}}\left(\chi_{H}\left(\mathbf{j}\right)\right)=H_{\mathbf{i}}\left(\chi_{H}\left(B_{p}\left(n\right)\right)\right)=H_{\mathbf{i}}\left(\underbrace{\frac{\chi_{H}\left(n\right)}{1-M_{H}\left(n\right)}}_{\textrm{call this }y}\right)
\end{equation}
where $y\overset{\textrm{def}}{=}\frac{\chi_{H}\left(n\right)}{1-M_{H}\left(n\right)}$
is a rational number.
\begin{claim}
$\left|y\right|_{p}\leq1$.

Proof of claim: First, since $y$ is a rational number, it lies in
$\mathbb{Q}_{p}$. So, by way of contradiction, suppose $\left|y\right|_{p}>1$.
By \textbf{Lemma \ref{lem:wrong values lemma}}, because $H$ is proper,
every branch specified by $\mathbf{i}$ maps $\mathbb{Q}_{p}\backslash\mathbb{Z}_{p}$
into $\mathbb{Q}_{p}\backslash\mathbb{Z}_{p}$. Consequently, $\left|H_{\mathbf{i}}\left(y\right)\right|_{p}>1$.
However, $H_{\mathbf{i}}\left(y\right)=x$, and $x$ is a rational
integer; hence, $1<\left|H_{\mathbf{i}}\left(y\right)\right|_{p}=\left|x\right|_{p}\leq1$.
This is impossible! So, it must be that $\left|y\right|_{p}\leq1$.
This proves the claim.
\end{claim}
\begin{claim}
\label{claim:2.3}$x=H^{\circ\left|\mathbf{i}\right|}\left(y\right)$

Proof of claim: Suppose the equality failed. Then $H_{\mathbf{i}}\left(y\right)=x\neq H^{\circ\left|\mathbf{i}\right|}\left(y\right)$,
and so $x=H_{\mathbf{i}}\left(y\right)$ is a wrong value of $H$
with seed $y$. Because $H$ is proper,\textbf{ Lemma \ref{lem:wrong values lemma}}
forces $\left|x\right|_{p}=\left|H_{\mathbf{i}}\left(y\right)\right|_{p}>1$.
However, $\left|x\right|_{p}\leq1$. This is just as impossible as
it was in the previous paragraph. This proves the claim.
\end{claim}
\vphantom{}

Finally, let $\mathbf{v}$ be the shortest string representing $n$,
so that $\mathbf{j}$ (the string representing $B_{p}\left(n\right)$)
is obtained by concatenating infinitely many copies of $\mathbf{v}$.
Because $q_{H}$ is co-prime to all the $d_{j}$s, \emph{note that
$H_{\mathbf{v}}$ is continuous on $\mathbb{Z}_{q_{H}}$}. As such:
\begin{equation}
\chi_{H}\left(B_{p}\left(n\right)\right)\overset{\mathbb{Z}_{q_{H}}}{=}\lim_{k\rightarrow\infty}\chi_{H}\left(\mathbf{v}^{\wedge k}\right)
\end{equation}
implies:
\begin{align*}
H_{\mathbf{v}}\left(\chi_{H}\left(B_{p}\left(n\right)\right)\right) & \overset{\mathbb{Z}_{q_{H}}}{=}\lim_{k\rightarrow\infty}H_{\mathbf{v}}\left(\chi_{H}\left(\mathbf{v}^{\wedge k}\right)\right)\\
 & \overset{\mathbb{Z}_{q_{H}}}{=}\lim_{k\rightarrow\infty}\chi_{H}\left(\mathbf{v}^{\wedge\left(k+1\right)}\right)\\
 & \overset{\mathbb{Z}_{q_{H}}}{=}\chi_{H}\left(B_{p}\left(n\right)\right)
\end{align*}
Hence, $H_{\mathbf{v}}\left(y\right)=y$. Since $\left|y\right|_{p}\leq1$,
the propriety of $H$ allows us to apply \textbf{Lemma \ref{lem:properness lemma}}
to conclude $H_{\mathbf{v}}\left(y\right)=H^{\circ\left|\mathbf{v}\right|}\left(y\right)$.
Thus, $y$ is a periodic point of $H:\mathbb{Z}_{p}\rightarrow\mathbb{Z}_{p}$.

By \textbf{Claim \ref{claim:2.3}}, $H$ iterates $y$ to $x$. Since
$y$ is a periodic point of $H$ in $\mathbb{Z}_{p}$, this forces
$x$ and $y$ to belong to the same cycle of $H$ in $\mathbb{Z}_{p}$,
with $x$ being a periodic point of $H$. As such, just as $H$ iterates
$y$ to $x$, so too does $H$ iterate $x$ to $y$. Likewise, since
$x$ is a rational integer, so too is $y$. Thus, $x$ belongs to
a cycle of $H$ in $\mathbb{Z}$, as desired.

Q.E.D.
\begin{example}
To illustrate the Correspondence Principle in action, observe that
the cycle $1\rightarrow2\rightarrow1$ of the $2$-Hydra map $H=T_{3}$
applies the even branch ($H_{0}$) second and the odd branch ($H_{1}$)
first. Thus, the string $\mathbf{j}$ such that $H_{\mathbf{j}}\left(1\right)=1$
is $\mathbf{j}=\left(0,1\right)$: 
\begin{equation}
1=H_{0}\left(H_{1}\left(1\right)\right)=H_{0}\left(2\right)=1
\end{equation}
The integer $n$ represented by $\mathbf{j}$ is $n=0\cdot2^{0}+1\cdot2^{1}=2$.
Thus: 
\begin{equation}
\centerdot01010101\ldots\overset{\mathbb{Z}_{2}}{=}B_{2}\left(2\right)=\frac{2}{1-2^{\lambda_{2}\left(2\right)}}=\frac{2}{1-4}=-\frac{2}{3}
\end{equation}
and so: 
\begin{equation}
\chi_{3}\left(-\frac{2}{3}\right)=\chi_{3}\left(B_{2}\left(2\right)\right)=\frac{\chi_{3}\left(2\right)}{1-M_{3}\left(2\right)}=\frac{\frac{1}{4}}{1-\frac{3}{4}}=\frac{1}{4-3}=1
\end{equation}
where: 
\begin{equation}
\chi_{3}\left(2\right)=\frac{1}{2}\chi_{3}\left(1\right)=\frac{1}{2}\left(\frac{3\chi_{3}\left(0\right)+1}{2}\right)=\frac{1}{2}\cdot\frac{0+1}{2}=\frac{1}{4}
\end{equation}
and: 
\begin{equation}
M_{3}\left(2\right)=\frac{3^{\#_{1}\left(2\right)}}{2^{\lambda_{2}\left(2\right)}}=\frac{3^{1}}{2^{2}}=\frac{3}{4}
\end{equation}
\end{example}
\vphantom{}

We end our exploration of the Correspondence Principle with a question
mark. As infamous as Collatz-type problems' difficulty might be, a
combination of familiarization with the available literature and personal
exploration of the problems themselves suggests that, as difficult
as the study of these maps' periodic points might be, the question
of their divergent trajectories could very well be an order of magnitude
more difficult. The diophantine equation characterizations of periodic
points given in \textbf{Corollary \ref{cor:CP v3} }show that, at
the bare minimum, the question of periodic points can has an \emph{expressible
}non-trivial reformulation in terms of an equation of finitely many
variables. Yes, there might be infinitely many possible equations
to consider, and solving any single one of them\textemdash let alone
all of them\textemdash is difficult\textemdash maybe even extraordinarily
difficulty. But, at least there is something we can \emph{attempt
}to confront. The finitude of any cycle of $H$ guarantees at least
this much.

For divergent trajectories, on the other hand, the prospects are bleak
indeed. The only finite aspect of their existence that come to mind
are the fact that they are bounded from below (if divergent to $+\infty$)
or above (if divergent to $-\infty$). I am not certain if it has
been shown, given an $n\in\mathbb{Z}$ belonging to a divergent trajectory
of $T_{3}$, that there are iterates $T_{3}^{\circ k}\left(n\right)$
which are divisible by arbitrarily large powers of $2$; that is:
\begin{equation}
\sup_{k\geq0}v_{2}\left(T_{3}^{\circ k}\left(n\right)\right)=\infty
\end{equation}
where $v_{2}$ is the $2$-adic valuation. Even with a positive or
negative answer to this question, it is not at all clear how it might
be useful for tackling the question of divergent trajectories, assuming
it would be useful at all. The extra difficulty of divergent trajectories
is all the more startling considering the likes of the shortened $5x+1$
map, $T_{5}$. We know that the set of divergent trajectories of $T_{5}$
in $\mathbb{N}_{1}$ has density $1$, yet not a single positive integer
has been proven to belong to a divergent trajectory! \cite{Lagarias-Kontorovich Paper}

Late in the writing of this dissertation (March 2022), I was surprised
to realize that, for an almost trivial reason, the Correspondence
Principle implies a connection between $\chi_{H}$ and $H$'s divergent
trajectories, and a particularly beautiful one at that. \textbf{Corollary
\ref{cor:CP v4} }shows that the periodic points are controlled by
values in $\mathbb{Z}$ attained by $\chi_{H}$ for \emph{rational
}$p$-adic inputs; that is, those $p$-adic integers whose $p$-adic
digit sequences are eventually periodic. But... what about \emph{irrational
}$p$-adic integers?\textemdash that is, $\mathbb{Z}_{p}\backslash\mathbb{Q}$,
the set of $\mathfrak{z}\in\mathbb{Z}_{p}$ whose sequence of $p$-adic
digits are not eventually periodic.

Modulo some simple conditions on $H$, we can prove that irrational
$\mathfrak{z}$s which make $\chi_{H}\left(\mathfrak{z}\right)$ an
integer then make $\chi_{H}\left(\mathfrak{z}\right)$ an element
of a divergent trajectory; this is \textbf{Theorem \ref{thm:Divergent trajectories come from irrational z}}.
\begin{prop}
\label{prop:Preparing for application to divergent trajectories}Let
$H$ be a semi-basic $p$-Hydra map, and suppose that $\left|H_{j}\left(0\right)\right|_{q_{H}}=1$
for all $j\in\left\{ 1,\ldots,p-1\right\} $. Then $\chi_{H}\left(\mathfrak{z}\right)\neq0$
for any $\mathfrak{z}\in\mathbb{Z}_{p}\backslash\mathbb{Q}$. 
\end{prop}
Proof: Let $H$ as given. By way of contradiction, suppose that $\chi_{H}\left(\mathfrak{z}\right)=0$
for some $\mathfrak{z}\in\mathbb{Z}_{p}\backslash\mathbb{Q}$. Now,
let $j\in\left\{ 0,\ldots,p-1\right\} $ be the first $p$-adic digit
of $\mathfrak{z}$; that is $j=\left[\mathfrak{z}\right]_{p}$. Then,
letting $\mathfrak{z}^{\prime}$ denote the $p$-adic integer $\left(\mathfrak{z}-j\right)/p$,
we can write: 
\begin{equation}
0=\chi_{H}\left(\mathfrak{z}\right)=\chi_{H}\left(p\mathfrak{z}^{\prime}+j\right)=H_{j}\left(\chi_{H}\left(\mathfrak{z}^{\prime}\right)\right)
\end{equation}
Next, suppose $j=0$. Then, since $H\left(0\right)=0$:
\begin{equation}
0=H_{0}\left(\chi_{H}\left(\mathfrak{z}^{\prime}\right)\right)=\frac{\mu_{0}}{p}\chi_{H}\left(\mathfrak{z}^{\prime}\right)
\end{equation}
This forces $\chi_{H}\left(\mathfrak{z}^{\prime}\right)=0$, seeing
as $\mu_{0}\neq0$. In this way, if the first $n$ $p$-adic digits
of $\mathfrak{z}$ are all $0$, we can remove those digits from $\mathfrak{z}$
to obtain a $p$-adic integer $\mathfrak{z}^{\left(n\right)}\overset{\textrm{def}}{=}p^{-n}\mathfrak{z}$
with the property that $\left|\mathfrak{z}^{\left(n\right)}\right|_{p}=1$
(i.e., $\left[\mathfrak{z}^{\left(n\right)}\right]_{p}\neq0$). If
\emph{all }of the digits of $\mathfrak{z}$ are $0$, this then makes
$\mathfrak{z}=0$, which would contradict the given that $\mathfrak{z}$
was an element of $\mathbb{Z}_{p}\backslash\mathbb{Q}$.

So, without loss of generality, we can assume that $j\neq0$. Hence:
\begin{align*}
0 & =H_{j}\left(\chi_{H}\left(\mathfrak{z}^{\prime}\right)\right)\\
 & \Updownarrow\\
\chi_{H}\left(\mathfrak{z}^{\prime}\right) & =-\frac{pH_{j}\left(0\right)}{\mu_{j}}
\end{align*}
Since $H$ is semi-basic, the fact that $j\neq0$ tells us that $\left|\mu_{j}\right|_{q_{H}}<1$
and $\left|p\right|_{q_{H}}=1$. This means $\left|H_{j}\left(0\right)\right|_{q_{H}}<1$;
else we would have that $\left|\chi_{H}\left(\mathfrak{z}^{\prime}\right)\right|_{q_{H}}>1$,
which is impossible seeing as $\chi_{H}\left(\mathbb{Z}_{p}\right)\subseteq\mathbb{Z}_{q_{H}}$
and every $q_{H}$-adic integer has\textbf{ }$q_{H}$-adic absolute
value $\leq1$. So, $\left|H_{j}\left(0\right)\right|_{q_{H}}<1$\textemdash but,
we were given that $\left|H_{j}\left(0\right)\right|_{q_{H}}=1$ for
all $j\in\left\{ 1,\ldots,p-1\right\} $. There's our first contradiction.

So, $\mathfrak{z}$ has no non-zero $p$-adic digits, which forces
$\mathfrak{z}$ to be $0$. But, once again, $\mathfrak{z}$ was given
to be an element of $\mathbb{Z}_{p}\backslash\mathbb{Q}$. There's
our second and final contradiction.

So, $\chi_{H}\left(\mathfrak{z}\right)$ must be non-zero.

Q.E.D.
\begin{thm}
\label{thm:Divergent trajectories come from irrational z}\index{Hydra map!divergent trajectories}Let
$H$ be a proper, integral, contracting, semi-basic $p$-Hydra map
fixing $0$ so that $\left|H_{j}\left(0\right)\right|_{q_{H}}=1$
for all $j\in\left\{ 1,\ldots,p-1\right\} $.

Let $\mathfrak{z}\in\mathbb{Z}_{p}\backslash\mathbb{Q}$ be such that
$\chi_{H}\left(\mathfrak{z}\right)\in\mathbb{Z}$. Then $\chi_{H}\left(\mathfrak{z}\right)$
belongs to a divergent trajectory of $H$. 
\end{thm}
Proof: Let $H$ and $\mathfrak{z}$ be as given. By \textbf{Proposition
\ref{prop:Preparing for application to divergent trajectories}},
$\chi_{H}\left(\mathfrak{z}\right)\neq0$.

Now, by way of contradiction, suppose that $\chi_{H}\left(\mathfrak{z}\right)$
did not belong to a divergent trajectory of $H$. By the basic theory
of dynamical systems on $\mathbb{Z}$ (\textbf{Theorem \ref{thm:orbit classes partition domain}}),
every element of $\mathbb{Z}$ belongs to either a divergent trajectory
of $H$, or to an orbit class of $H$ which contains a cycle. Thus,
it must be that $\chi_{H}\left(\mathfrak{z}\right)$ is pre-periodic.
As such, there is an $n\geq0$ so that $H^{\circ n}\left(\chi_{H}\left(\mathfrak{z}\right)\right)$
is a periodic point of $H$.

Using $\chi_{H}$'s functional equations from \textbf{Lemma \ref{lem:Unique rising continuation and p-adic functional equation of Chi_H}},
it then follows that $H^{\circ n}\left(\chi_{H}\left(\mathfrak{z}\right)\right)=\chi_{H}\left(\mathfrak{y}\right)$
where $\mathfrak{y}=m+p^{\lambda_{p}\left(m\right)}\mathfrak{z}$,
where $m$ is the unique non-negative integer whose sequence of $p$-adic
digits corresponds to the composition sequence of branches of $H$
brought to bear when we apply $H^{\circ n}$ to send $\chi_{H}\left(\mathfrak{z}\right)$
to $\chi_{H}\left(\mathfrak{y}\right)$. Moreover, since $\chi_{H}\left(\mathfrak{z}\right)\neq0$,
the fact that $\left\{ 0\right\} $ is an isolated cycle of $H$ guarantees
that the periodic point $\chi_{H}\left(\mathfrak{y}\right)$ is non-zero.
Since $H$ is proper, we can apply Version 4 of \textbf{Correspondence
Principle }(\textbf{Corollary \ref{cor:CP v4}}). This tells us $\mathfrak{y}$
is an element of $\mathbb{Q}\cap\mathbb{Z}_{p}^{\prime}$. Since $\mathfrak{y}$
is rational, its $p$-adic digits are eventually periodic. However,
$\mathfrak{y}=m+p^{\lambda_{p}\left(m\right)}\mathfrak{z}$, where
$\mathfrak{z}$ is ``irrational'' ($\mathfrak{z}\in\mathbb{Z}_{p}\backslash\mathbb{Q}$)\textemdash the
$p$-adic digits of $\mathfrak{z}$ are aperiodic, which means the
same is true of $\mathfrak{y}$. This is a contradiction!

So, our assumption that $\chi_{H}\left(\mathfrak{z}\right)$ did not
belong to a divergent trajectory of $H$ forces the digits of $\mathfrak{y}$
(which are not periodic) to be eventually periodic\textemdash a clear
impossibility. This proves that $\chi_{H}\left(\mathfrak{z}\right)$
must be a divergent point of $H$.

Q.E.D.

\vphantom{}

This theorem suggests the following: 
\begin{conjecture}[\textbf{A Correspondence Principle for Divergent Points?}]
\label{conj:correspondence theorem for divergent trajectories}Provided
that $H$ satisfies certain prerequisites such as the hypotheses of
\textbf{Theorem \ref{thm:Divergent trajectories come from irrational z}},
$x\in\mathbb{Z}$ belongs to a divergent trajectory under $H$ if
and only if there is a $\mathfrak{z}\in\mathbb{Z}_{p}\backslash\mathbb{Q}$
so that $\chi_{H}\left(\mathfrak{z}\right)\in\mathbb{Z}$. 
\end{conjecture}
\vphantom{}

The difficulty of this conjecture is that, unlike the Correspondence
Principle for periodic points of $H$, I have yet to find \emph{constructive
}method of producing a $\mathfrak{z}\in\mathbb{Z}_{p}\backslash\mathbb{Q}$
for which $\chi_{H}\left(\mathfrak{z}\right)$ is in $\mathbb{Z}$.

Finally, we have:
\begin{conjecture}[\textbf{Divergence Conjectures for $qx+1$}]
Let $\chi_{q}\left(\mathfrak{z}\right)$ denote the numen of the
shortened $qx+1$ map\index{$qx+1$ map}, where $q$ is an odd prime.
Then, there exists a $\mathfrak{z}\in\mathbb{Z}_{2}\backslash\mathbb{Q}$
so that $\chi_{q}\left(\mathfrak{z}\right)\in\mathbb{Z}$ if and only
if $q=3$.
\end{conjecture}

\subsection{\label{subsec:2.2.4 Other-Avenues}Other Avenues}

\subsubsection{\label{subsec:Relaxing-the-Requirements}Relaxing the Requirements
of Monogenicity and Non-Degeneracy}

With regard to the convergence of $\chi_{H}$ over the $p$-adics,
the key expression is (\ref{eq:Formal rising limit of Chi_H}): 
\begin{equation}
\lim_{m\rightarrow\infty}\chi_{H}\left(\left[\mathfrak{z}\right]_{p^{m}}\right)\overset{\mathbb{Z}_{q}}{=}\sum_{\ell=1}^{\infty}\frac{b_{j_{\ell}\left(\mathfrak{z}\right)}}{d_{j_{\ell}\left(\mathfrak{z}\right)}}M_{H}\left(\mathbf{j}_{\ell-1}\left(\mathfrak{z}\right)\right)
\end{equation}
where, recall, we write the $p$-adic integer $\mathfrak{z}$ as:
\begin{equation}
\mathfrak{z}=j_{1}\left(\mathfrak{z}\right)+j_{2}\left(\mathfrak{z}\right)p+j_{3}\left(\mathfrak{z}\right)p^{2}+\cdots
\end{equation}
and write $\mathbf{j}_{\ell-1}\left(\mathfrak{z}\right)=\left(j_{1}\left(\mathfrak{z}\right),\ldots,j_{\ell-1}\left(\mathfrak{z}\right)\right)\in\left(\mathbb{Z}/p\mathbb{Z}\right)^{\ell-1}$,
with the convention that $M_{H}\left(\mathbf{j}_{0}\left(\mathfrak{z}\right)\right)=1$.
By the concatenation identity \textbf{Proposition \ref{prop:M_H concatenation identity}},
to obtain $M_{H}\left(\mathbf{j}_{\ell}\left(\mathfrak{z}\right)\right)$,
we multiply $M_{H}\left(\mathbf{j}_{\ell-1}\left(\mathfrak{z}\right)\right)$
by a factor of $a_{j_{\ell}\left(\mathfrak{z}\right)}/d_{j_{\ell}\left(\mathfrak{z}\right)}$.
The various conditions contained in the notion of semi-basicness were
designed so as to ensure that $\left|a_{j}/d_{j}\right|_{q_{H}}\geq1$
occurs if and only if $j=0$. This $q_{H}$-adic condition, in turn,
is what guarantees that (\ref{eq:Formal rising limit of Chi_H}) will
converge in $\mathbb{Z}_{q_{H}}$ for any $\mathfrak{z}\in\mathbb{Z}_{p}^{\prime}$\textemdash that
is, any $\mathfrak{z}$ with infinitely many non-zero $p$-adic digits.

For the sake of clarity, rather than continue to work with the encumbrance
of our indexing notation and an arbitrary $p$-Hydra map, we will
illustrate the obstacles in relaxing the various conditions (monogenicity,
(semi)-simplicity, non-degeneracy), by way of an example. 
\begin{example}[\textbf{A polygenic Hydra map}]
\label{exa:Polygenic example, part 1}Consider the $3$-Hydra map
$H:\mathbb{Z}\rightarrow\mathbb{Z}$ defined by: 
\begin{equation}
H\left(n\right)\overset{\textrm{def}}{=}\begin{cases}
\frac{an}{3} & \textrm{if }n\overset{3}{\equiv}0\\
\frac{bn+b^{\prime}}{3} & \textrm{if }n\overset{3}{\equiv}1\\
\frac{cn+2c^{\prime}}{3} & \textrm{if }n\overset{3}{\equiv}2
\end{cases}\label{eq:Toy model for a polygenic 3-Hydra map}
\end{equation}
where $a$, $b$, and $c$ are positive integers co-prime to $3$,
and where $b^{\prime}$ and $c^{\prime}$ are positive integers so
that $b+b^{\prime}\overset{3}{\equiv}0$ and $c+c^{\prime}\overset{3}{\equiv}0$,
respectively; $b^{\prime}$ and $c^{\prime}$ guarantee that each
of the branches of $H$ outputs a non-negative integer. For this $H$,
we have that: 
\begin{equation}
M_{H}\left(\mathbf{j}\right)=\frac{a^{\#_{3:0}\left(\mathbf{j}\right)}\times b^{\#_{3:1}\left(\mathbf{j}\right)}\times c^{\#_{3:2}\left(\mathbf{j}\right)}}{3^{\left|\mathbf{j}\right|}}
\end{equation}
In order for (\ref{eq:Formal rising limit of Chi_H}) to converge,
we need $M_{H}\left(\mathbf{j}\right)$ to become small in some sense
as $\left|\mathbf{j}\right|\rightarrow\infty$ with infinitely many
non-zero digits. As such, the value of $a$ is irrelevant to our convergence
issue; only the constants $b$ and $c$, associated\textemdash respectively\textemdash to
the branches $H_{1}$ and $H_{2}$ are of import, because any $\mathbf{j}$
with only finitely many non-zero entries reduces to an integer.

In the monogenic case dealt with in this dissertation, we required
both $b$ and $c$ to be multiples of some integer $q_{H}\geq2$.
\textbf{Proposition \ref{prop:q-adic behavior of M_H of j as the number of non-zero digits tends to infinity}}
showed that this condition guaranteed the $q_{H}$-adic decay of $M_{H}\left(\mathbf{j}\right)$
to $0$ as the number of non-zero digits of $\mathbf{j}$ increased
to infinity. On the other hand, if this were the \index{Hydra map!polygenic}\textbf{polygenic
}case, then $b$ and $c$ would be co-prime to one another. Note that
if one or both of $b$ and $c$ are equal to $1$, then we have a
degenerate case, because multiplication by $1$ does not decrease
$p$-adic size for any prime $p$. So, we can suppose that $b$ and
$c$ are co-prime. For simplicity, let $b$ and $c$ be distinct primes,
with neither $b$ nor $c$ being $3$. Then, there are several possibilities
for $M_{H}\left(\mathbf{j}\right)$ as $\mathbf{j}\rightarrow\infty$
with $\mathbf{j}$ having infinitely many non-zero entries: 
\end{example}
\begin{enumerate}
\item If $\mathbf{j}$ has infinitely many $1$s (this corresponds to multiplications
by $b$), then $M_{H}\left(\mathbf{j}\right)$ tends to $0$ in $b$-adic
magnitude. 
\item If $\mathbf{j}$ has infinitely many $2$s (this corresponds to multiplications
by $c$), then $M_{H}\left(\mathbf{j}\right)$ tends to $0$ in $c$-adic
magnitude. 
\item If only one of $1$ or $2$ occurs infinitely many times in $\mathbf{j}$,
then $M_{H}\left(\mathbf{j}\right)$ will tend to $0$ in either $b$-adic
or $c$-adic magnitude, \emph{but not in both}\textemdash \emph{that}
occurs if and only if $\mathbf{j}$ contains both infinitely many
$1$s \emph{and} infinitely many $2$s.
\end{enumerate}
\begin{rem}
In hindsight, I believe my notion of a \textbf{frame} (see Subsection
\ref{subsec:3.3.3 Frames}) might be suited for the polygenic case.
Considering \textbf{Example \ref{exa:Polygenic example, part 1}},
we could assign the $b$-adic and $c$-adic convergence, respectively,
to the sets of $3$-adic integers whose representative strings $\mathbf{j}$
have infinitely many $1$s and $2$s, respectively. This would then
allow for a meaningful definition of $\chi_{H}$, and thereby potentially
open the door to applying my methods to polygenic Hydra maps.
\end{rem}

\subsubsection{\label{subsec:Baker,-Catalan,-and}Baker, Catalan, and Collatz}

Any proof of the Weak Collatz Conjecture\footnote{Recall, this is the assertion that the only periodic points of the
Collatz map in the positive integers are $1$, $2$, and $4$.} will necessarily entail a significant advancement in\index{transcendental number theory}
transcendental number theory. A proof of the Weak Collatz Conjecture
would necessarily yield a proof of \textbf{Baker's Theorem} far simpler
than any currently known method \cite{Tao Blog}. \index{Baker's Theorem}Baker's
Theorem, recall, concerns lower bounds on the absolute values of \textbf{linear
forms of logarithms}, which are expressions of the form: 
\begin{equation}
\beta_{1}\ln\alpha_{1}+\cdots+\beta_{N}\ln\alpha_{N}\label{eq:linear form in logarithm}
\end{equation}
where the $\beta_{n}$s and $\alpha_{n}$s are complex algebraic numbers
(with all of the $\alpha_{n}$s being non-zero). Most proofs of Baker's
Theorem employ a variant of what is known as Baker's Method, in which
one constructs of an analytic function (called an \textbf{auxiliary
function})\textbf{ }with a large number of zeroes of specified degree
so as to obtain contradictions on the assumption that (\ref{eq:linear form in logarithm})
is small in absolute value. The import of Baker's Theorem is that
it allows one to obtain lower bounds on expressions such as $\left|2^{m}-3^{n}\right|$,
where $m$ and $n$ are positive integers, and it was for applications
such as these in conjunction with the study of Diophantine equations
that Baker earned the Fields Medal. See \cite{Baker's Transcendental Number Theory}
for a comprehensive account of the subject; also, parts of \cite{Cohen Number Theory}.

With $\chi_{H}$ and the \textbf{Correspondence Principle }at our
disposal, we can get a glimpse at the kinds of advancements in transcendental
number theory that might be needed in order to resolve the Weak Collatz
Conjecture. To begin, it is instructive to consider the following
table of values for $\chi_{H}$ and related functions in the case
of the Shortened $qx+1$ maps $T_{q}$, where $q$ is an odd prime.
For brevity, we write $\chi_{q}$ to denote $\chi_{T_{q}}:\mathbb{Z}_{2}\rightarrow\mathbb{Z}_{q}$. 
\begin{center}
\begin{table}
\begin{centering}
\begin{tabular}{|c|c|c|c|c|c|c|}
\hline 
$n$ & $\#_{1}\left(n\right)$ & $\lambda_{2}\left(n\right)$ & $\chi_{q}\left(n\right)$ & $\chi_{q}\left(B_{2}\left(n\right)\right)$ & $\chi_{3}\left(B_{2}\left(n\right)\right)$ & $\chi_{5}\left(B_{2}\left(n\right)\right)$\tabularnewline
\hline 
\hline 
$0$ & $0$ & $0$ & $0$ & $0$ & $0$ & $0$\tabularnewline
\hline 
$1$ & $1$ & $1$ & $\frac{1}{2}$ & $\frac{1}{2-q}$ & $-1$ & $-\frac{1}{3}$\tabularnewline
\hline 
$2$ & $1$ & $2$ & $\frac{1}{4}$ & $\frac{1}{4-q}$ & $1$ & $-1$\tabularnewline
\hline 
$3$ & $2$ & $2$ & $\frac{2+q}{4}$ & $\frac{2+q}{4-q^{2}}$ & $-1$ & $-\frac{1}{3}$\tabularnewline
\hline 
$4$ & $1$ & $3$ & $\frac{1}{8}$ & $\frac{1}{8-q}$ & $\frac{1}{5}$ & $\frac{1}{3}$\tabularnewline
\hline 
$5$ & $2$ & $3$ & $\frac{4+q}{8}$ & $\frac{4+q}{8-q^{2}}$ & $-7$ & $-\frac{9}{17}$\tabularnewline
\hline 
$6$ & $2$ & $3$ & $\frac{2+q}{8}$ & $\frac{2+q}{8-q^{2}}$ & $-5$ & $-\frac{7}{17}$\tabularnewline
\hline 
$7$ & $3$ & $3$ & $\frac{4+2q+q^{2}}{8}$ & $\frac{4+2q+q^{2}}{8-q^{3}}$ & $-1$ & $-\frac{1}{3}$\tabularnewline
\hline 
$8$ & $1$ & $4$ & $\frac{1}{16}$ & $\frac{1}{16-q}$ & $\frac{1}{13}$ & $\frac{1}{11}$\tabularnewline
\hline 
$9$ & $2$ & $4$ & $\frac{8+q}{16}$ & $\frac{8+q}{16-q^{2}}$ & $\frac{11}{7}$ & $-\frac{13}{9}$\tabularnewline
\hline 
$10$ & $2$ & $4$ & $\frac{4+q}{16}$ & $\frac{4+q}{16-q^{2}}$ & $1$ & $-1$\tabularnewline
\hline 
$11$ & $3$ & $4$ & $\frac{8+4q+q^{2}}{16}$ & $\frac{8+4q+q^{2}}{16-q^{3}}$ & $-\frac{29}{11}$ & $-\frac{53}{109}$\tabularnewline
\hline 
$12$ & $2$ & $4$ & $\frac{2+q}{16}$ & $\frac{2+q}{16-q^{2}}$ & $\frac{5}{7}$ & $-\frac{7}{9}$\tabularnewline
\hline 
$13$ & $3$ & $4$ & $\frac{8+2q+q^{2}}{16}$ & $\frac{8+2q+q^{2}}{16-q^{3}}$ & $-\frac{23}{11}$ & $-\frac{43}{109}$\tabularnewline
\hline 
$14$ & $3$ & $4$ & $\frac{4+2q+q^{2}}{16}$ & $\frac{4+2q+q^{2}}{16-q^{3}}$ & $\frac{19}{7}$ & $-\frac{39}{109}$\tabularnewline
\hline 
$15$ & $4$ & $4$ & $\frac{8+4q+2q^{2}+q^{3}}{16}$ & $\frac{8+4q+2q^{2}+q^{3}}{16-q^{4}}$ & $-1$ & $-\frac{1}{3}$\tabularnewline
\hline 
\end{tabular}
\par\end{centering}
\caption{Values of $\chi_{q}\left(n\right)$ and related functions}
\end{table}
\par\end{center}
\begin{example}
\index{Correspondence Principle}As per the Correspondence Principle,\textbf{
}note that the integer values attained by $\chi_{3}\left(B_{2}\left(n\right)\right)$
and $\chi_{5}\left(B_{2}\left(n\right)\right)$ are all periodic points
of the maps $T_{3}$ and $T_{5}$, respectively; this includes fixed
points at negative integers, as well. Examining $\chi_{q}\left(B_{2}\left(n\right)\right)$,
we see certain patterns, such as:
\begin{equation}
\chi_{q}\left(B_{2}\left(2^{n}-1\right)\right)=\frac{1}{2-q},\textrm{ }\forall n\in\mathbb{N}_{1}
\end{equation}
More significantly, $\chi_{3}\left(B_{2}\left(n\right)\right)$ appears
to be more likely to be positive than negative, whereas the opposite
appears to hold true for $q=5$ (and, heuristically, for all $q\geq5$).
Of special interest, however, is: 
\begin{equation}
\chi_{q}\left(B_{2}\left(10\right)\right)=\frac{4+q}{16-q^{2}}=\frac{1}{4-q}
\end{equation}
\end{example}
By Version 1 of the \textbf{Correspondence} \textbf{Principle} (\textbf{Theorem
\ref{thm:CP v1}}), every cycle $\Omega\subseteq\mathbb{Z}$ of $T_{q}$
with $\left|\Omega\right|\geq2$ contains an integer $x$ of the form:
\begin{equation}
x=\chi_{q}\left(B_{2}\left(n\right)\right)=\frac{\chi_{q}\left(n\right)}{1-\frac{q^{\#_{2:1}\left(n\right)}}{2^{\lambda_{2}\left(n\right)}}}=\frac{2^{\lambda_{2}\left(n\right)}\chi_{q}\left(n\right)}{2^{\lambda_{2}\left(n\right)}-q^{\#_{2:1}\left(n\right)}}\label{eq:Rational expression of odd integer periodic points of ax+1}
\end{equation}
for some $n\in\mathbb{N}_{1}$. In fact, every periodic point $x$
of $T_{q}$ in the odd integers can be written in this form. As such
$\left|2^{\lambda_{2}\left(n\right)}-q^{\#_{2:1}\left(n\right)}\right|=1$
is a \emph{sufficient condition }for $\chi_{q}\left(B_{2}\left(n\right)\right)$
to be a periodic point of $T_{q}$. However, as the $n=10$ case shows,
this is not\emph{ }a \emph{necessary }condition: there can be values
of $n$ where $2^{\lambda_{2}\left(n\right)}-q^{\#_{2:1}\left(n\right)}$
is large in archimedean absolute value, and yet nevertheless divides
the numerator on the right-hand side of (\ref{eq:Rational expression of odd integer periodic points of ax+1}),
thereby reducing $\chi_{q}\left(B_{2}\left(n\right)\right)$ to an
integer. In fact, thanks to P. Mih\u{a}ilescu \index{Mihu{a}ilescu, Preda@Mih\u{a}ilescu, Preda}'s
resolution of \index{Catalan's Conjecture}\textbf{ Catalan's Conjecture},
it would seem that Baker's Method-style estimates on the archimedean\emph{
}size of $2^{\lambda_{2}\left(n\right)}-q^{\#_{2:1}\left(n\right)}$
will be of little use in understanding $\chi_{q}\left(B_{2}\left(n\right)\right)$: 
\begin{thm}[\textbf{Mih\u{a}ilescu's Theorem}\footnote{Presented in \cite{Cohen Number Theory}.}]
The only choice of $x,y\in\mathbb{N}_{1}$ and $m,n\in\mathbb{N}_{2}$
for which: 
\begin{equation}
x^{m}-y^{n}=1\label{eq:Mihailescu's Theorem}
\end{equation}
are $x=3$, $m=2$, $y=2$, $n=3$ (that is, $3^{2}-2^{3}=1$). 
\end{thm}
\vphantom{}

With Mih\u{a}ilescu's Theorem, it is easy to see that for any odd
integer $q\geq3$, $\left|q^{\#_{2:1}\left(n\right)}-2^{\lambda_{2}\left(n\right)}\right|$
will never be equal to $1$ for any $n\geq8$, because the exponent
of $2$ (that is, $\lambda_{2}\left(n\right)$) will be $\geq4$ for
all $n\geq8$ (any such $n$ has $4$ or more binary digits). Consequently,
for any odd prime $q$, if $n\geq8$ makes $\chi_{q}\left(B_{2}\left(n\right)\right)$
into a rational integer (and hence, a periodic point of $T_{q}$),
it \emph{must }be that the numerator $2^{\lambda_{2}\left(n\right)}\chi_{q}\left(n\right)$
in (\ref{eq:Rational expression of odd integer periodic points of ax+1})
is a multiple of $2^{\lambda_{2}\left(n\right)}-q^{\#_{2:1}\left(n\right)}$.

The set of $p$-Hydra maps for which conclusions of this sort hold
can be enlarged in several ways. First, in general, note that we have:
\begin{equation}
\chi_{H}\left(B_{p}\left(n\right)\right)=\frac{\chi_{H}\left(n\right)}{1-M_{H}\left(n\right)}=\frac{p^{\lambda_{p}\left(n\right)}\chi_{H}\left(n\right)}{p^{\lambda_{p}\left(n\right)}\left(1-M_{H}\left(n\right)\right)}=\frac{p^{\lambda_{p}\left(n\right)}\chi_{H}\left(n\right)}{p^{\lambda_{p}\left(n\right)}-\prod_{j=0}^{p-1}\mu_{j}^{\#_{p:j}\left(n\right)}}\label{eq:Chi_H o B functional equation with fraction in simplest form}
\end{equation}
If $\mu_{0}=1$ and there is an integer $\mu\geq2$ so that $\mu_{1},\ldots,\mu_{p-1}$
are all positive integer powers of $\mu$ (i.e., for each $j\in\left\{ 1,\ldots,p-1\right\} $
there is an $r_{j}\in\mathbb{N}_{1}$ so that $\mu_{j}=\mu^{r_{j}}$)
then, when $n\geq1$, the denominator of the right-hand side takes
the form $p^{a}-\mu^{b}$ for some positive integers $a$ and $b$
depending solely on $n$. By Mih\u{a}ilescu's Theorem, the only way
to have $\left|p^{a}-\mu^{b}\right|=1$ is for either $p^{a}-\mu^{b}=3^{2}-2^{3}$
or $p^{a}-\mu^{b}=2^{3}-3^{2}$. As such, for this case, once $n$
is large enough so that both $a$ and $b$ are $\geq3$, the only
way for $\chi_{H}\left(B_{p}\left(n\right)\right)$ to be a rational
integer for such an $n$ is if $p^{\lambda_{p}\left(n\right)}-\prod_{j=0}^{p-1}\mu_{j}^{\#_{p:j}\left(n\right)}$
is a divisor of the integer $p^{\lambda_{p}\left(n\right)}\chi_{H}\left(n\right)$
of magnitude $\geq2$.

In the most general case, if we allow the $\mu_{j}$ to take on different
values (while still requiring $H$ to be basic so that the Correspondence
Principle applies), in order to replicate the application of Mih\u{a}ilescu's
Theorem to the $p^{a}-\mu^{b}$ case, examining the the denominator
term $p^{\lambda_{p}\left(n\right)}-\prod_{j=0}^{p-1}\mu_{j}^{\#_{p:j}\left(n\right)}$
from (\ref{eq:Chi_H o B functional equation with fraction in simplest form})
and ``de-parameterizing'' it, the case where the denominator reduces
to $1$ or $-1$ can be represented by a Diophantine equation\index{diophantine equation}
of the form: 
\begin{equation}
\left|x^{m}-y_{1}^{n_{1}}y_{2}^{n_{2}}\cdots y_{J}^{n_{J}}\right|=1\label{eq:Generalized Catalan Diophantine Equation}
\end{equation}
where $J$ is a fixed integer $\geq2$, and $x$, $y_{1},\ldots,y_{J}$,
$m$, and $n_{1},\ldots,n_{J}$ are positive integer variables. For
any fixed choice of $m$ and the $n_{j}$s, \textbf{Faltings' Theorem}\index{Faltings' Theorem}\textbf{
}(the resolution of \textbf{Mordell's Conjecture}) guarantees there
are only finitely many rational numbers $x,y_{1},\ldots,y_{J}$ for
which (\ref{eq:Generalized Catalan Diophantine Equation}) holds true.
On the other hand, for any fixed choice of $x,y_{1},\ldots,y_{J}$
(where we identify $x$ with $p$ and the $y_{j}$s with the $\mu_{j}$s),
if it can be shown that (\ref{eq:Generalized Catalan Diophantine Equation})
holds for only finitely many choices of $m$ and the $n_{j}$s, it
will follow that the corresponding case for $\left|p^{\lambda_{p}\left(n\right)}-\prod_{j=0}^{p-1}\mu_{j}^{\#_{p:j}\left(n\right)}\right|$
is not equal to $1$ for all sufficiently large positive integers
$n$. In that situation, any cycles of $H$ not accounted for by values
of $n$ for which $\left|p^{\lambda_{p}\left(n\right)}-\prod_{j=0}^{p-1}\mu_{j}^{\#_{p:j}\left(n\right)}\right|=1$
must occur as a result of $\left|p^{\lambda_{p}\left(n\right)}-\prod_{j=0}^{p-1}\mu_{j}^{\#_{p:j}\left(n\right)}\right|$
being an integer $\geq2$ which divides $p^{\lambda_{p}\left(n\right)}\chi_{H}\left(n\right)$.

In light of this, rather than the archimedean/euclidean \emph{size
}of $p^{\lambda_{p}\left(n\right)}-\prod_{j=0}^{p-1}\mu_{j}^{\#_{p:j}\left(n\right)}$,
it appears we must study its \emph{multiplicative }structure, as well
as that of $p^{\lambda_{p}\left(n\right)}\chi_{H}\left(n\right)$\textemdash which
is to say, the set of these integers' prime divisors, and how this
set depends on $n$. In other words, we ought to study the $p$-adic
absolute values of $\left|\chi_{H}\left(n\right)\right|_{p}$ and
$\left|1-M_{H}\left(n\right)\right|_{p}$ for various values of $n\geq1$
and primes $p$. It may also be of interest to use $p$-adic methods
(Fourier series, Mellin transform) to study $\left|\chi_{H}\left(\mathfrak{z}\right)\right|_{p}$
as a real-valued function of a $p$-adic integer variable for varying
values of $p$.

\subsubsection{\label{subsec:Connections-to-Tao}Connections to Tao (2019)}

The origins of the present paper lie in work the author did in late
2019 and early 2020. Around the same time, in late 2019, Tao published
a paper on the Collatz Conjecture that applied probabilistic techniques
in a novel way, by use of what Tao\index{Tao, Terence} called ``\index{Syracuse Random Variables}Syracuse
Random Variables'' \cite{Tao Probability paper}. Despite the many
substantial differences between these two approaches (Tao's is probabilistic,
ours is not), they are linked at a fundamental level, thanks to the
function $\chi_{3}$\textemdash our shorthand for $\chi_{T_{3}}$,
the $\chi_{H}$ associated to the Shortened Collatz map. Given that
both papers have their genesis in an examination of the behavior of
different combinations of the iterates of branches of $T_{3}$, it
is not entirely surprising that both approaches led to $\chi_{3}$.

Tao's approach involves constructing his Syracuse Random Variables
and then comparing them to a set-up involving tuples of geometric
random variables, the comparison in question being an estimate on
the ``distance'' between the Syracuse Random Variables and the geometric
model, as measured by the total variation norm for discrete random
variables. To attain these results, the central challenge Tao overcomes
is obtaining explicit estimates for the decay of the characteristic
function of the Syracuse Random Variables. In our terminology, Tao
establishes decay estimates for the archimedean absolute value of
the function $\varphi_{3}:\hat{\mathbb{Z}}_{3}\rightarrow\mathbb{C}$
defined by:
\begin{equation}
\varphi_{3}\left(t\right)\overset{\textrm{def}}{=}\int_{\mathbb{Z}_{2}}e^{-2\pi i\left\{ t\chi_{3}\left(\mathfrak{z}\right)\right\} _{3}}d\mathfrak{z},\textrm{ }\forall t\in\hat{\mathbb{Z}}_{3}\label{eq:Tao's Characteristic Function}
\end{equation}
where $\left\{ \cdot\right\} _{3}$ is the $3$-adic fractional part,
$d\mathfrak{z}$ is the Haar probability measure on $\mathbb{Z}_{2}$,
and $\hat{\mathbb{Z}}_{3}=\mathbb{Z}\left[\frac{1}{3}\right]/\mathbb{Z}=\mathbb{Q}_{3}/\mathbb{Z}_{3}$
is the Pontryagin dual of $\mathbb{Z}_{3}$, identified here with
the set of all rational numbers in $\left[0,1\right)$ whose denominators
are non-negative integer powers of the prime number $3$. Tao's decay
estimate is given in \textbf{Proposition 1.17 }of his paper, where
the above integral appears in the form of an expected value \cite{Tao Probability paper}.

With this in mind, a natural next step for furthering both this paper
and Tao's would be to study the ``characteristic function\index{characteristic function}''
of $\chi_{H}$ for an arbitrary semi-basic $p$-Hydra maps $H$; this
is the continuous function $\varphi_{H}:\hat{\mathbb{Z}}_{q_{H}}\rightarrow\mathbb{C}$
defined by: 
\begin{equation}
\varphi_{H}\left(t\right)\overset{\textrm{def}}{=}\int_{\mathbb{Z}_{p}}e^{-2\pi i\left\{ t\chi_{H}\left(\mathfrak{z}\right)\right\} _{q_{H}}}d\mathfrak{z},\textrm{ }\forall t\in\hat{\mathbb{Z}}_{q_{H}}\label{eq:Characteristic function of Chi_H}
\end{equation}
where $\left\{ \cdot\right\} _{q_{H}}$ is the $q_{H}$-adic fractional
part, $d\mathfrak{z}$ is the Haar probability measure on $\mathbb{Z}_{p}$,
and $\hat{\mathbb{Z}}_{q_{H}}=\mathbb{Z}\left[\frac{1}{q_{H}}\right]/\mathbb{Z}=\mathbb{Q}_{q_{H}}/\mathbb{Z}_{q_{H}}$
is the Pontryagin dual of $\mathbb{Z}_{q_{H}}$, identified here with
the set of all rational numbers in $\left[0,1\right)$ whose denominators
are non-negative integer powers of $q_{H}$.

We can establish functional equations\index{functional equation}
can be established for $\varphi_{H}$. 
\begin{prop}
Let $H$ be a semi-basic $p$-Hydra map which fixes $0$, and write
$q=q_{H}$. Then, the characteristic function $\varphi_{H}:\hat{\mathbb{Z}}_{q}\rightarrow\mathbb{C}$
defined by \emph{(\ref{eq:Characteristic function of Chi_H})} satisfies
the functional equation: 
\begin{equation}
\varphi_{H}\left(t\right)=\frac{1}{p}\sum_{j=0}^{p-1}e^{-2\pi i\left\{ \frac{b_{j}t}{d_{j}}\right\} _{q}}\varphi_{H}\left(\frac{a_{j}t}{d_{j}}\right),\textrm{ }\forall t\in\hat{\mathbb{Z}}_{q}\label{eq:Phi_H functional equation}
\end{equation}
\end{prop}
\begin{rem}
Since $H$ fixes $0$ and is semi-basic, \textbf{Proposition \ref{prop:co-primality of d_j and q_H}}
tells us that all of the $d_{j}$s are co-prime to $q$. Thus, the
multiplication-by-$d_{j}$ map $t\mapsto d_{j}t$ is a group automorphism
of $\hat{\mathbb{Z}}_{q}$, and as such, $b_{j}t/d_{j}$ and $a_{j}t/d_{j}$
denote the images of $b_{j}t$ and $a_{j}t$, respectively, under
the inverse of this automorphism. One can remove the denominator terms
if desired by writing: 
\begin{align*}
d & \overset{\textrm{def}}{=}\textrm{lcm}\left(d_{0},\ldots,d_{p-1}\right)\\
d_{j}^{\prime} & \overset{\textrm{def}}{=}\frac{d}{d_{j}},\textrm{ }\forall j\in\left\{ 0,\ldots,p-1\right\} \in\mathbb{Z}\\
\alpha_{j} & \overset{\textrm{def}}{=}a_{j}d_{j}^{\prime}\in\mathbb{Z}\\
\beta_{j} & \overset{\textrm{def}}{=}b_{j}d_{j}^{\prime}\in\mathbb{Z}
\end{align*}
and then replacing $t$ with $dt$ in (\ref{eq:Phi_H functional equation})
to obtain: 
\begin{equation}
\varphi_{H}\left(dt\right)=\frac{1}{p}\sum_{j=0}^{p-1}e^{-2\pi i\beta_{j}t}\varphi_{H}\left(\alpha_{j}t\right),\textrm{ }\forall t\in\hat{\mathbb{Z}}_{q}
\end{equation}
\end{rem}
Proof: The proof is a straight-forward computation. Details regarding
Fourier analysis on $p$-adic rings can be found in \cite{Automorphic Representations}.
\begin{align*}
\varphi_{H}\left(t\right) & =\int_{\mathbb{Z}_{p}}e^{-2\pi i\left\{ t\chi_{H}\left(\mathfrak{z}\right)\right\} _{q}}d\mathfrak{z}\\
 & =\sum_{j=0}^{p-1}\int_{j+p\mathbb{Z}_{p}}e^{-2\pi i\left\{ t\chi_{H}\left(\mathfrak{z}\right)\right\} _{q}}d\mathfrak{z}\\
\left(\mathfrak{y}=\frac{\mathfrak{z}-j}{p}\right); & =\frac{1}{p}\sum_{j=0}^{p-1}\int_{\mathbb{Z}_{p}}e^{-2\pi i\left\{ t\chi_{H}\left(p\mathfrak{y}+j\right)\right\} _{q}}d\mathfrak{y}\\
 & =\frac{1}{p}\sum_{j=0}^{p-1}\int_{\mathbb{Z}_{p}}e^{-2\pi i\left\{ t\frac{a_{j}\chi_{H}\left(\mathfrak{y}\right)+b_{j}}{d_{j}}\right\} _{q}}d\mathfrak{y}\\
 & =\frac{1}{p}\sum_{j=0}^{p-1}e^{-2\pi i\left\{ \frac{b_{j}t}{d_{j}}\right\} _{q}}\int_{\mathbb{Z}_{p}}e^{-2\pi i\left\{ \frac{a_{j}t}{d_{j}}\chi_{H}\left(\mathfrak{y}\right)\right\} _{q}}d\mathfrak{y}\\
 & =\frac{1}{p}\sum_{j=0}^{p-1}e^{-2\pi i\left\{ \frac{b_{j}t}{d_{j}}\right\} _{q}}\varphi_{H}\left(\frac{a_{j}t}{d_{j}}\right)
\end{align*}

Q.E.D.

\vphantom{}

The Fourier analysis here naturally leads to probabilistic notions.
For any measurable function $f:\mathbb{Z}_{p}\rightarrow\mathbb{Z}_{q}$,
any $n\in\mathbb{N}_{0}$ and any $m\in\mathbb{Z}$, we write:

\begin{equation}
\textrm{P}\left(f\overset{q^{n}}{\equiv}m\right)\overset{\textrm{def}}{=}\int_{\mathbb{Z}_{p}}\left[f\left(\mathfrak{z}\right)\overset{q^{n}}{\equiv}m\right]d\mathfrak{z}\label{eq:Definition of the probability that f is congruent to k mod q to the n}
\end{equation}
where $d\mathfrak{z}$ is the Haar probability measure on $\mathbb{Z}_{p}$
and $\left[f\left(\mathfrak{z}\right)\overset{q^{n}}{\equiv}m\right]$
is $1$ if $f\left(\mathfrak{z}\right)$ is congruent to $m$ mod
$q^{n}$ and is $0$ otherwise. We adopt the convention that $\mathfrak{x}\overset{1}{\equiv}\mathfrak{y}$
is true for any $q$-adic integers $\mathfrak{x}$ and $\mathfrak{y}$,
and hence, that: 
\begin{equation}
\textrm{P}\left(f\overset{q^{0}}{\equiv}m\right)=\textrm{P}\left(f\overset{1}{\equiv}m\right)=\int_{\mathbb{Z}_{p}}\left[f\left(\mathfrak{z}\right)\overset{1}{\equiv}m\right]d\mathfrak{z}=\int_{\mathbb{Z}_{p}}d\mathfrak{z}=1
\end{equation}
for any $f$ and any $m$.

Fourier relations follow from the identity: 
\begin{equation}
\left[\mathfrak{x}\overset{q^{n}}{\equiv}\mathfrak{y}\right]=\frac{1}{q^{n}}\sum_{k=0}^{q^{n}-1}e^{2\pi i\left\{ k\frac{\mathfrak{x}-\mathfrak{y}}{q^{n}}\right\} _{q}},\textrm{ }\forall n\in\mathbb{N}_{0},\textrm{ }\forall\mathfrak{x},\mathfrak{y}\in\mathbb{Z}_{q}
\end{equation}
This gives us the typical Fourier relations between probabilities
and characteristic functions: 
\begin{align}
\textrm{P}\left(\chi_{H}\overset{q^{n}}{\equiv}m\right) & =\frac{1}{q^{n}}\sum_{k=0}^{q^{n}-1}e^{-\frac{2\pi ikm}{q^{n}}}\varphi_{H}\left(\frac{k}{q^{n}}\right)\\
\varphi_{H}\left(\frac{m}{q^{n}}\right) & =\sum_{k=0}^{q^{n}-1}e^{\frac{2\pi imk}{q^{n}}}\textrm{P}\left(\chi_{H}\overset{q^{n}}{\equiv}k\right)
\end{align}
and the Parseval Identity: 
\begin{equation}
\frac{1}{q^{n}}\sum_{k=0}^{q^{n}-1}\left|\varphi_{H}\left(\frac{k}{q^{n}}\right)\right|^{2}=\sum_{k=0}^{q^{n}-1}\left(\textrm{P}\left(\chi_{H}\overset{q^{n}}{\equiv}k\right)\right)^{2}
\end{equation}
Under certain circumstances (such as the case of the $T_{a}$ maps),
one can use (\ref{eq:Phi_H functional equation}) to recursively solve
for $\varphi_{H}$\textemdash or for $\textrm{P}\left(\chi_{H}\overset{q^{n}}{\equiv}k\right)$,
after performing a discrete Fourier transform. Doing so for $H=T_{3}$
yields the recursive formula given by Tao in \textbf{Lemma 1.12 }of
his paper \cite{Tao Probability paper}. That being said, it remains
to be seen whether or not recursive formulae can be derived for these
probabilities and characteristic functions in the case of a general
$p$-Hydra map. It may also be of interest to study the expressions
$\varphi_{H,p}:\hat{\mathbb{Z}}_{p}\rightarrow\mathbb{C}$ given by:
\begin{equation}
\varphi_{H,p}\left(t\right)\overset{\textrm{def}}{=}\int_{\mathbb{Z}_{p}}e^{-2\pi i\left\{ t\chi_{H}\left(\mathfrak{z}\right)\right\} _{p}}d\mathfrak{z}\overset{\textrm{def}}{=}\lim_{N\rightarrow\infty}\frac{1}{p^{N}}\sum_{n=0}^{p^{N}-1}e^{-2\pi i\left\{ t\chi_{H}\left(n\right)\right\} _{p}},\textrm{ }\forall t\in\hat{\mathbb{Z}}_{p}\label{eq:p-adic characteristic function for Chi_H}
\end{equation}
where $p$ is an arbitrary prime. Doing so, should, conceivably, advance
our understanding of the divisibility of $\chi_{H}\left(n\right)$
by $p$ as $n$ varies.

\subsubsection{\label{subsec:Dirichlet-Series-and}Dirichlet Series and Complex-Analytic
Methods}

Like in most situations, given the Correspondence Principle and the
associated identity from \textbf{Lemma \ref{lem:Chi_H o B_p functional equation}}:

\begin{equation}
\chi_{H}\left(B_{p}\left(n\right)\right)=\frac{\chi_{H}\left(n\right)}{1-M_{H}\left(n\right)}
\end{equation}
one asks: \emph{what can be done with this?} Prior to stumbling upon
the $\left(p,q\right)$-adic approach, my intention was to use classical,
complex-analytic techniques of analytic number theory, such as Dirichlet
Series\index{Dirichlet series}, Mellin transforms, contour integration,
and the Residue Theorem. According to the Correspondence Principle,
classifying the periodic points of $H$ in $\mathbb{Z}$ is just a
matter of finding those integers $x\in\mathbb{Z}$ and $n\geq0$ so
that $x=\chi_{H}\left(B_{p}\left(n\right)\right)$; which is to say:
\begin{equation}
\left(1-M_{H}\left(n\right)\right)x-\chi_{H}\left(n\right)=0\label{eq:nth term of Dirichlet series}
\end{equation}
Dividing everything by $\left(n+1\right)^{s}$ and summing over $n\geq0$
gives: 
\begin{equation}
\left(\zeta\left(s\right)-\sum_{n=0}^{\infty}\frac{M_{H}\left(n\right)}{\left(n+1\right)^{s}}\right)x-\sum_{n=0}^{\infty}\frac{\chi_{H}\left(n\right)}{\left(n+1\right)^{s}}=0\label{eq:Dirichlet series}
\end{equation}
for all complex $s$ for which the Dirichlet series converge. By exploiting
$M_{H}$ and $\chi_{H}$'s functional equations, the right-hand side
of (\ref{eq:nth term of Dirichlet series}) can be written in terms
of a contour integral of the right-hand side of (\ref{eq:Dirichlet series})
using a generalized version of \index{Perron's Formula}\textbf{Perron's
Formula} (see\footnote{The series of papers by Flajolet et. al. on applications of the Mellin
transform in analytic combinatorics (starting with \cite{Flajolet - Mellin Transforms})
provide a comprehensive and illustrative exposition of the Mellin
transform\index{Mellin transform!complex} and its many, \emph{many
}applications to combinatorics, number theory, algorithm analysis,
and more.} \cite{Flajolet - Digital sums}). As usual, the hope was to shift
the contour of integration by first analytically continuing (\ref{eq:Dirichlet series})
to the left half-plane, and then using the Residue Theorem to obtain
to obtain a more enlightening expression for (\ref{eq:nth term of Dirichlet series}).
Unfortunately, while the Dirichlet series in (\ref{eq:Dirichlet series})
\emph{do }admit analytic continuations to meromorphic functions on
$\mathbb{C}$\textemdash where they have half-lattices of poles to
the left of their abscissae of absolute convergence\textemdash these
continuations exhibit hyper-exponential growth as $\textrm{Re}\left(s\right)\rightarrow-\infty$.
This appears to greatly limit the usefulness of this complex analytic
approach. I have unpublished work in this vein; a messy rough draft
of which (containing much overlap with this dissertation) can be found
on arXiv (see \cite{Mellin transform paper}). I intend to pursue
these issues more delicately and comprehensively at some later date.
That being said, I briefly revisit this approach at the end of Section
\ref{sec:4.1 Preparatory-Work--}, so the curious reader can turn
there for a first taste of that approach.

\chapter{\label{chap:3 Methods-of--adic}Methods of $\left(p,q\right)$-adic
Analysis}

\includegraphics[scale=0.45]{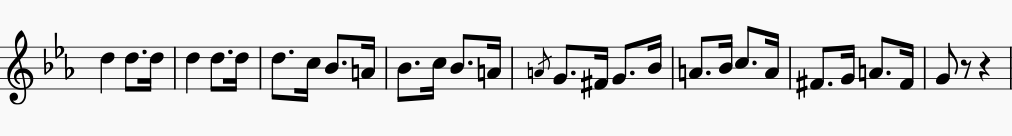}
\begin{quote}
\vphantom{}
\begin{flushright}
$\ldots$many people feel that functions $\mathbb{Z}_{p}\rightarrow\mathbb{Q}_{p}$
are more interesting than functions $\mathbb{Z}_{p}\rightarrow\mathbb{Q}_{q}$,
which is understandable. 
\par\end{flushright}
\begin{flushright}
\textemdash W. M. Schikhof\footnote{From \emph{Ultrametric Calculus} (\cite{Ultrametric Calculus}), on
the bottom of page 97.} 
\par\end{flushright}
\begin{flushright}
People will want to know: what is the \emph{point} of $\left(p,q\right)$-adic
functions? You need a good answer. 
\par\end{flushright}
\begin{flushright}
\textemdash K. Conrad\footnote{<https://mathoverflow.net/questions/409806/in-need-of-help-with-parsing-non-archimedean-function-theory>} 
\par\end{flushright}

\end{quote}
\vphantom{}

The purpose of this chapter is to present and develop an analytical
theory capable of dealing with functions like $\chi_{H}$. Section
\ref{sec:3.1 A-Survey-of} presents the ``classical'' theory, much
of which I independently re-discovered before learning that Schikhof\index{Schikhof, W. M.}
had already done it\footnote{One of the dangers of working in isolation as I did is the high probability
of ``re-inventing the wheel'', as one of my professors (Nicolai
Haydn) put it.} all\textemdash and more\textemdash in his own PhD dissertation, back
in 1967. Prototypically, the ``classical'' material concerns the
study of $C\left(\mathbb{Z}_{p},\mathbb{Q}_{q}\right)$\textemdash the
space of continuous functions from $\mathbb{Z}_{p}$ to $\mathbb{Q}_{q}$\textemdash and
the related interplay between functional analysis and Fourier analysis
in that setting. Section \ref{sec:3.1 A-Survey-of} explains these
theoretical details, as well as providing a primer on how to do $\left(p,q\right)$-adic
computations.

The other two sections of this chapter, however, are new; they extend
the ``classical'' material to include functions like $\chi_{H}$.
Section \ref{sec:3.2 Rising-Continuous-Functions} introduces the
notion of \textbf{rising-continuity}, which provides a natural extension
of the $C\left(\mathbb{Z}_{p},\mathbb{Q}_{q}\right)$ theory. Section
\ref{sec:3.3 quasi-integrability}\textemdash the largest of the three\textemdash is
the most important section of Chapter 3, containing the most significant
innovations, as well as the greatest complexities.

The central focus of Section \ref{sec:3.3 quasi-integrability} is
the phenomenon I call \textbf{quasi-integrability}. Classically, the
range of $\left(p,q\right)$-adic which can be meaningfully integrated
is astonishingly narrow, consisting solely of continuous functions.
By re-interpreting discontinuous functions like $\chi_{H}$ as \emph{measures},
we can meaningfully integrate them, as well as integrate their product
with any continuous $\left(p,q\right)$-adic function; this is done
in Subsection \ref{subsec:3.3.5 Quasi-Integrability}; the quasi-integrable
functions are precisely those $\left(p,q\right)$-adic functions for
which this method applies. However, in formulating quasi-integrability,
we will have to deal with an extremely unorthodox situation: infinite
series of a $p$-adic integer variable $\mathfrak{z}$ which converge
at every $\mathfrak{z}\in\mathbb{Z}_{p}$, but for which the topology
used to sum the series is allowed to vary from point to point. I introduce
the concept of a \textbf{frame }in order to bring order and rigor
to this risky, messy business (see Subsection \ref{subsec:3.3.3 Frames}).

Subsections \ref{subsec:3.3.1 Heuristics-and-Motivations} and \ref{subsec:3.3.2 The--adic-Dirichlet}
contain motivational examples and a connection to the Dirichlet kernel
of classical Fourier analysis, respectively, and thereby serve to
set the stage for frames and quasi-integrability. Subsection \ref{subsec:3.3.4 Toward-a-Taxonomy}
introduces some tedious terminology for describing noteworthy types
of $\left(p,q\right)$-adic measures. Its most important features
are the \textbf{Fourier Resummation Lemmata}, which will be used extensively
in our analysis of $\chi_{H}$ in Chapter 4. Subsection \ref{subsec:3.3.6 L^1 Convergence}
briefly introduces another hitherto-ignored area of non-archimedean
analysis: the Banach space of functions $\chi:\mathbb{Z}_{p}\rightarrow\mathbb{C}_{q}$
for which the real-valued function $\mathfrak{z}\in\mathbb{Z}_{p}\mapsto\left|\chi\left(\mathfrak{z}\right)\right|_{q}\in\mathbb{R}$
is integrable with respect to the real-valued Haar probability measure
on $\mathbb{Z}_{p}$. As will be shown in Subsection \ref{subsec:4.3.2 A-Wisp-of}
of Chapter 4, these spaces provide an approachable setting for studying
$\chi_{H}$. Whether or not significant progress can be made in that
setting remains to be seen.

The final subsection of \ref{sec:3.3 quasi-integrability}, Subsection
\ref{subsec:3.3.7 -adic-Wiener-Tauberian}, contains the second of
the this dissertation's three main results, a $\left(p,q\right)$-adic
generalization of \textbf{Wiener's Tauberian Theorem}. I prove two
versions, one for continuous $\left(p,q\right)$-adic functions, and
another for $\left(p,q\right)$-adic measures. The latter version
will be the one used in Chapter 4\textemdash at the end of Subsection
\ref{subsec:4.2.2}\textemdash to establish the \textbf{Tauberian
Spectral Theorem }for $p$-Hydra maps, my dissertation's third main
result. A \emph{sub}-subsection at the end of \ref{subsec:3.3.7 -adic-Wiener-Tauberian}
shows how the study of $\chi_{H}$ can be $q$-adically approximated
by the study of the eigenvalues of a family of increasingly large
Hermitian matrices.

\section{\label{sec:3.1 A-Survey-of}A Survey of the ``Classical'' theory}

Due to the somewhat peculiar position currently occupied by $\left(p,q\right)$-adic
analysis\index{$p,q$-adic!analysis}\textemdash too exotic to be thoroughly
explored, yet too specific to merit consideration at length by specialists
in non-archimedean analysis\textemdash I have taken the liberty of
providing a brief essay on the historical and theoretical context
of the subject and its various flavors; this can be found in Subsection
\ref{subsec:3.1.1 Some-Historical-and}.

Before that, however, I think it will be helpful to point out the
most important (and bizarre) results of $\left(p,q\right)$-adic analysis,
particularly for the benefit of readers who are used to working with
real- or complex-valued functions. The main features are: 
\begin{itemize}
\item The existence of a translation-invariant $\mathbb{C}_{q}$-valued
``measure'' on $\mathbb{Z}_{p}$, unique up to a choice of a normalization
constant, and compatible with an algebra of sets, albeit not the Borel
$\sigma$-algebra. \index{Borel!sigma-algebra@$\sigma$-algebra}Moreover,
this measure (the \textbf{$\left(p,q\right)$-adic Haar measure})
is a continuous linear functional on $C\left(\mathbb{Z}_{p},\mathbb{Q}_{q}\right)$.
This is very good news for us, because it means we can integrate and
do Fourier analysis in the $\left(p,q\right)$-adic setting. 
\item The lack of a concept of ``almost everywhere\emph{''}. Although
a substantive theory of $\left(p,q\right)$-adic measures exists,
it has no ``almost everywhere\emph{''}.\emph{ }In fact, unless the
normalization constant of the $\left(p,q\right)$-adic Haar measure
is chosen to be $0$, \emph{the only measurable set of measure $0$
is the empty set!} 
\item Because there is no way to take a difference quotient for a $\left(p,q\right)$-adic
function (how do you divide a $q$-adic number by a $p$-adic number?),
there is no differential calculus\footnote{However, my methods allow for us to define a $\left(p,q\right)$-adic
Mellin transform (see \textbf{Remark \ref{rem:pq adic mellin transform}}
on page \pageref{rem:pq adic mellin transform}). Since the Mellin
transform is used to formulate differentiation in the distributional
sense for functions $f:\mathbb{Q}_{p}\rightarrow\mathbb{C}$, it seems
quite likely that that same approach, done in the $\left(p,q\right)$-adic
context, can be used to give $\left(p,q\right)$-adic analysis a concept
of distributional derivatives, and thereby open the door to analytical
investigation of $\left(p,q\right)$-adic differential equations.} in $\left(p,q\right)$-adic analysis. There are also no such things
as polynomials, rational functions, or analytic functions. 
\item \emph{A function is integrable if and only if it is continuous!} This
equivalence is intimately connected to the absence of non-empty sets
of measure zero. 
\item We will have to leave the triangle inequality for integrals to wait
for us outside by the door; it has no place in $\left(p,q\right)$-adic
analysis. In this subject, there is generally no meaningful relationship
between $\left|\int f\right|_{q}$ and $\int\left|f\right|_{q}$.
As much as we would want to write $\left|\int f\right|_{q}\leq\int\left|f\right|_{q}$,
it simply isn't justified. \textbf{Example \ref{exa:triangle inequality failure}
}on page \pageref{exa:triangle inequality failure} provides the unfortunate
details. 
\end{itemize}
A first encounter with the theory of $\left(p,q\right)$-adic Fourier
analysis is a surreal experience. Indeed, in the ``classical'' $\left(p,q\right)$-adic
theory, given a $\left(p,q\right)$-adic function $f:\mathbb{Z}_{p}\rightarrow\mathbb{C}_{q}$
(where, of course, $p\neq q$),\emph{ the following are equivalent}: 
\begin{enumerate}
\item $f$ possesses a well-defined Fourier transform $\hat{f}$. 
\item $f$ is integrable. 
\item $f$ is continuous. 
\item The Fourier series representation of $f$ (using $\hat{f}$) is \emph{absolutely
convergent }uniformly over $\mathbb{Z}_{p}$. 
\end{enumerate}
These three wonders\textemdash especially the spectacular equivalence
of (1) through (4)\textemdash are, among some other factors, the principal
reasons why $\left(p,q\right)$-adic analysis has lain neglected for
the over-half-a-century that has elapsed since its inception. The
subject's rigid inflexibility and astonishing lack of subtlety goes
a long way toward explaining why it has been consigned to eek out
a lonely existence in the cabinet of mathematical curiosities. Of
course, part of the purpose of this dissertation is to argue that
this consignment was premature.

Lastly, before freeing the reader to pursue the exposition I have
prepared for them, let me first give an exposition for that exposition,
which is somewhat non-standard, to the extent that there even \emph{is
}a ``standard'' pedagogical approach to $\left(p,q\right)$-adic
analysis.

One of the main features of non-archimedean analysis as a subject
is that the available integration theories will depend on the fields
used for the functions' domains and codomains. Of these, the one we
will use is what Khrennikov calls the \textbf{Monna-Springer integral}
\cite{Quantum Paradoxes}. Schikhof in an appendix of \emph{Ultrametric
Analysis} \cite{Ultrametric Calculus} and van Rooij in a chapter
of \emph{Non-Archimedean Functional Analysis} \cite{van Rooij - Non-Archmedean Functional Analysis}
also give expositions of this theory, though not by that name. For
the probability-minded, Khrennikov gives a second, full-throated presentation
of the integration theory in his article \cite{Probabilities taking values in non-archimedean fields}
on non-archimedean probability theories\index{non-archimedean!probability theory}.

As anyone who has gone through a standard course on integration and
measures can attest, presentations of measure theory and measure-theoretic
integration often lean toward the abstract. The same is true for the
non-archimedean case. In certain respect, the non-archimedean abstraction
outdoes the archimedean. These presentations typically begin with
one of the following: 
\begin{itemize}
\item An explanation of how the set-based approach to measure theory utilized
in most university courses is not compatible with the non\textendash archimedean
case. Instead, they turn to the functional analysis approach advocated
by the Bourbaki group, defining measures as linear functionals satisfying
a certain norm condition. 
\item To get around the fact that Borel sets aren't well-suited to form
a set-based theory of measures, the presentation takes as a primitive
notion the existence of a collection of sets satisfying most of the
properties you would expect, along with a few extra technicalities
needed to avoid disaster. 
\end{itemize}
Although we will \emph{eventually }get to this sort of material (see
Subsection \ref{subsec:3.1.6 Monna-Springer-Integration}), because
of the paramount importance of explicit, specific, \emph{concrete}
computations for this dissertation, I have weighted my presentation
toward such practical matters. With the exception of \ref{subsec:3.1.6 Monna-Springer-Integration},
I have tried to keep the level of generality to an absolute minimum,
so as to allow the basic mechanics to shine through as clearly as
possible.

After the historical essay in Subsection \ref{subsec:3.1.1 Some-Historical-and},
the exposition begins in Subsection \ref{subsec:3.1.2 Banach-Spaces-over}
with a sliver of the more general theory of Banach spaces over a non-archimedean
field. Subsection \ref{subsec:3.1.3 The-van-der} introduces the \textbf{van
der Put basis} for spaces of functions on $\mathbb{Z}_{p}$, which
will be used throughout this section to explicitly construct spaces
like $C\left(\mathbb{Z}_{p},\mathbb{Q}_{q}\right)$ and compute integrals
and Fourier transforms.

After this, we explore the $\left(p,q\right)$-adic Fourier transform
of a continuous $\left(p,q\right)$-adic function (Subsection \ref{subsec:3.1.4. The--adic-Fourier}),
and then use it to construct $\left(p,q\right)$-adic measures. The
initial viewpoint will be of a measure as an object defined by its
\textbf{Fourier-Stieltjes transform} (Subsection \ref{subsec:3.1.5-adic-Integration-=00003D000026}).
This view of $\left(p,q\right)$-adic measures will be very important
in Section \ref{sec:3.3 quasi-integrability} and throughout Chapters
\ref{chap:4 A-Study-of} and \ref{chap:6 A-Study-of}. In the process,
we will see how to compute continuous $\left(p,q\right)$-adic functions'
Fourier transforms\textemdash two different formulae will be given\textemdash along
with how to integrate functions with respect to $\left(p,q\right)$-adic
measures, convolutions, and the basic\textemdash but fundamental\textemdash change
of variable formula for affine transformations ($\mathfrak{y}=\mathfrak{a}\mathfrak{z}+\mathfrak{b}$)
in integrals. We will constantly draw from this practical information,
knowledge of which is essential for further independent study of $\chi_{H}$
and related matters,

We conclude our survey by returning to the abstract with Subsection
\ref{subsec:3.1.6 Monna-Springer-Integration}, giving an exposition
of Monna-Springer integration theory. This provides the ``correct''
generalization of the explicit, concrete methods of $\left(p,q\right)$-adic
integration dealt with in the earlier subsections. As a pedagogical
approach, I feel that grounding the Monna-Springer integral in the
concrete case of $\left(p,q\right)$-adic Fourier series and $\left(p,q\right)$-adic
Fourier-Stieltjes transforms makes the abstract version less intimidating.
That being said, my inclusion of the Monna-Springer theory is primarily
for completeness' sake. $\left(p,q\right)$-adic analysis is extraordinarily
pathological, even by the bizarre standards of non-archimedean analysis.
Familiar results such as the Dominated Convergence Theorem or Hölder's
Inequality reduce to trivialities in the $\left(p,q\right)$-adic
setting.

\subsection{\label{subsec:3.1.1 Some-Historical-and}Some Historical and Theoretical
Context}

In the title of Section \ref{sec:3.1 A-Survey-of}, the word ``classical''
is enclosed in scare quotes. This seems like a fair compromise to
me. On the one hand, in $\left(p,q\right)$-adic analysis we have
a nearly-sixty-year-old subject ripe for generalization; on the other,
to the extent that a well-articulated theory of $\left(p,q\right)$-adic
analysis exists, it is too esoteric and (understandably) under-appreciated
to truly merit the status of a ``classic''. Until now, the $\left(p,q\right)$-adic
analysis chronicled in this Chapter has smoldered in double disregard.
Why, indeed, would anyone in their right mind want to investigate
a subject which is too exotic and inert for the general community,
but also too \emph{insufficiently }general to have earned clear place
in the work of the exotic specialists? Somewhat surprisingly, the
same could be said for the $p$-adic numbers themselves.

Kurt Hensel first published his discovery of $p$-adic numbers in
1897 \cite{Hensel-s original article,Gouvea's p-adic number history slides,Journey throughout the history of p-adic numbers}.
At the time, they did not make much of a splash. The eventual inclusion
of $p$-adic numbers in the mathematical canon came about purely by
serendipity, with Helmut Hasse in the role of the inadvertent midwife.
It would make for an interesting scene in a Wes Anderson film: in
1913, Hensel publishes a book on the $p$-adics; sometime between
then and 1920, a copy of said book managed to find way to a certain
antique shop in Göttingen. In 1920, Hasse\textemdash a then-student
at the University of Göttingen\textemdash stumbled across said book
in said antique shop while wandering about town as graduate students
in Göttingen are wont to do when they feel a hankering for a breath
of fresh air. Through this chain of events, Hasse chose to transfer
to the University of Marburg, at which Hensel taught, with the express
purpose of learning everything he could about Hensel's marvelous,
star-crossed innovation. The result? In 1922, Hasse publishes his
doctoral thesis, in which he establishes what is now called the \textbf{Hasse-Minkowski
Theorem}, otherwise known as the \textbf{local-global principle} \textbf{for
quadratic forms }\cite{Gouvea's p-adic number history slides,On the origins of p-adic analysis}.

For the uninitiated, the local-global principle is the aspiration
that one can find construct a tuple of rational numbers which solve
a given polynomial equation in several variables by solving the equation
in $\mathbb{R}$ as well as in the $p$-adics for all primes $p$,
and then using the Chinese Remainder Theorem to stitch together the
real and $p$-adic solutions into a solution over $\mathbb{Q}$. (Conrad's
online notes \cite{Conrad - Local Global Principle} give an excellent
demonstration of this method in practice.) Fascinatingly, although
the Hasse-Minkowski Theorem assures us this approach will work in
polynomials for which each variable has an exponent of at most $2$,
the local-global principle \emph{does not }apply for polynomial equations
of arbitrary degree. A particularly well-known counterexample is due
to Selmer \cite{Conrad - Local Global Principle,Local-Global Principle Failure}:
\begin{thm}[\textbf{\textit{Selmer's Counterexample}}]
The cubic diophantine equation:
\begin{equation}
3x^{3}+4y^{3}+5z^{3}=0
\end{equation}
has solutions $\left(x,y,z\right)\neq\left(0,0,0\right)$ over $\mathbb{R}$
and over $\mathbb{Q}_{p}$ for every prime $p$, yet its only solution
over $\mathbb{Q}$ is $\left(x,y,z\right)=\left(0,0,0\right)$.
\end{thm}
\vphantom{}

Understanding the conditions where it does or does not apply is still
an active area of research in arithmetic geometry. Indeed, as Conrad
points out, the famous \textbf{Birch and Swinnerton-Dyer Conjecture
}concerns, in part, a relation between the rational points of an elliptic
curve $E$ over $\mathbb{Q}$ and the points of $E$ over $\mathbb{R}$
and the $p$-adics \cite{Conrad - Local Global Principle}.

The local-global philosophy has since had an explosive influence in
mathematics, above all in algebraic number theory and adjacent fields,
where it opened the door to the introduction of $p$-adic and adèlic
approaches. A central protagonist of this story is the late, great
John Tate (1925 \textendash{} 2019)\index{Tate, John}. In his epochal
thesis of 1950 (``\emph{Fourier analysis in number fields and Hecke's
zeta functions}''), Tate revolutionized algebraic number theory by
demonstrating how the \textbf{Euler Products} for $L$-functions such
as the Riemann Zeta Function: 
\begin{equation}
\zeta\left(s\right)\overset{\mathbb{C}}{=}\prod_{p\in\mathbb{P}}\frac{1}{1-p^{-s}},\textrm{ }\forall\textrm{ }\textrm{Re}\left(s\right)>1\label{eq:Euler Product for the RZF}
\end{equation}
(where $\mathbb{P}$ denotes the set of all prime numbers) could be
understood in the local-global spirit as the product of so-called
\textbf{local factors} associated to each prime $p$ \cite{Tate's thesis}.
This recontextualization occurred by way of Mellin-transform type
integrals over the $p$-adic numbers.

Letting $d\mathfrak{z}_{p}$ denote the real-valued Haar probability
measure on $\mathbb{Z}_{p}$, the aforementioned local factors of
(\ref{eq:Euler Product for the RZF}) can be expressed by the integral:
\begin{equation}
\frac{1}{1-p^{-s}}\overset{\mathbb{C}}{=}\frac{p}{p-1}\int_{\mathbb{Z}_{p}}\left|\mathfrak{z}_{p}\right|_{p}^{s-1}d\mathfrak{z}_{p},\textrm{ }\forall\textrm{Re}\left(s\right)>1
\end{equation}
This can also be re-written as: 
\begin{equation}
\frac{1}{1-p^{-s}}\overset{\mathbb{C}}{=}\int_{\mathbb{Z}_{p}}\left|\mathfrak{z}_{p}\right|_{p}^{s-1}d^{\times}\mathfrak{z}_{p},\textrm{ }\forall\textrm{Re}\left(s\right)>1\label{eq:RZF local factors}
\end{equation}
where: 
\begin{equation}
d^{\times}\mathfrak{z}_{p}\overset{\textrm{def}}{=}\frac{p}{p-1}d\mathfrak{z}_{p}
\end{equation}
is the real-valued Haar probability measure for the group $\left(\mathbb{Z}_{p}^{\times},\times\right)$
of multiplicatively invertible elements of $\mathbb{Z}_{p}$. Taking
the product of the local factors (\ref{eq:RZF local factors}) over
all primes $p$, the tensor product transforms the product of the
integrals into a single integral over the profinite integers, $\check{\mathbb{Z}}$\nomenclature{$\check{\mathbb{Z}}$}{the ring of profinite integers, $\prod_{p\in\mathbb{P}}\mathbb{Z}_{p}$ \nopageref}:
\begin{equation}
\zeta\left(s\right)\overset{\mathbb{C}}{=}\int_{\check{\mathbb{Z}}}\left|\mathfrak{z}\right|^{s}d^{\times}\mathfrak{z}\overset{\textrm{def}}{=}\prod_{p\in\mathbb{P}}\int_{\mathbb{Z}_{p}}\left|\mathfrak{z}_{p}\right|_{p}^{s-1}d^{\times}\mathfrak{z}_{p},\textrm{ }\forall\textrm{Re}\left(s\right)>1
\end{equation}
Here:
\begin{equation}
\check{\mathbb{Z}}\overset{\textrm{def}}{=}\prod_{p\in\mathbb{P}}\mathbb{Z}_{p}\label{eq:Definition of the Profinite Integers}
\end{equation}
and so, every $\mathfrak{z}\in\check{\mathbb{Z}}$ is a $\mathbb{P}$-tuple
$\left(\mathfrak{z}_{2},\mathfrak{z}_{3},\mathfrak{z}_{5},\ldots\right)$
of $p$-adic integers, one from each $p$. $d^{\times}\mathfrak{z}$,
meanwhile, is the Haar probability measure for the multiplicative
group of units of $\check{\mathbb{Z}}$. The significance of this
new perspective and its impact on subsequent developments in number
theory cannot be overstated. Indeed, in hindsight, it is hard to believe
that anyone could have ever doubted the value of the $p$-adics to
mathematics.

But the story does not end here. If anything, Tate's work is like
an intermezzo, bridging two worlds by transmuting the one into the
other. Starting in the mid-twentieth century, mathematical analysis
began to pick up where Hensel left off \cite{Gouvea's p-adic number history slides,Journey throughout the history of p-adic numbers},
exploring the $p$-adic numbers' ultrametric structure and the implications
it would have for analysis and its theories in those settings\textemdash continuity,
differentiability, integrability, and so on. Although Hensel himself
had taken the first steps toward investigating the properties of functions
of his $p$-adic numbers, a general programme for analysis on valued
fields other than $\mathbb{R}$ or $\mathbb{C}$\textemdash independent
of any number-theoretic investigations was first posed by A. F. Monna\index{Monna, A. F.}
in a series of papers in 1943 \cite{van Rooij - Non-Archmedean Functional Analysis}.
From this grew the subject that would come to be known as \emph{non-archimedean
analysis}.

In this sense, non-archimedean analysis\textemdash as opposed to the
more number-theoretically loaded term ``$p$-adic analysis''\textemdash has
its origins in the activities of graduate students, this time Netherlands
in the 1960s, when W. M. Schikhof, A. C. M. van Rooij\index{van Rooij, A. C. M.},
M. van der Put\index{van der Put, Marius}, and J. van Tiel established
a so-called ``$p$-adics and fine dining'' society, intent on exploring
the possibility of analysis\textemdash particularly \emph{functional}
analysis\textemdash over non-archimedean valued fields. The subject
is also sometimes called \emph{ultrametric analysis}\footnote{In this vein, one of the subject's principal introductory texts is
called \emph{Ultrametric Calculus} \cite{Ultrametric Calculus}.}, to emphasize the fundamental distinction between itself and classical
analysis.

In the grand scheme of things, non-archimedean valued fields aren't
as wild as one might think. As was proved by Ostrowski in 1915, the
only valued fields of characteristic zero obtainable as metric completions
of $\mathbb{Q}$ are $\mathbb{R}$, $\mathbb{C}$, and the $p$-adics.
A necessary and immediate consequence of this straight-forward classification
is that non-archimedean analysis comes in four different flavors depending
on the nature of the fields we use for our domains and co-domains.

Let $\mathbb{F}$ and $K$ be non-archimedean valued fields, complete
with respect to their absolute values. Moreover, let the residue fields
of $\mathbb{F}$ and $K$ have different characteristic; say, let
$\mathbb{F}$ be a $p$-adic field, and let $K$ be a $q$-adic field.
By combinatorics, note that there are really only four\emph{ }different
types of non-archimedean functions we can consider: functions from
an archimedean field to a non-archimedean field ($\mathbb{C}\rightarrow\mathbb{F}$);
from a non-archimedean field to an archimedean field ($\mathbb{F}\rightarrow\mathbb{C}$);
from \emph{and} to the same non-archimedean field ($\mathbb{F}\rightarrow\mathbb{F}$);
and those between two archimedean fields with non-equal residue fields
($\mathbb{F}\rightarrow K$). In order, I refer to these as \textbf{$\left(\infty,p\right)$-adic
analysis} ($\mathbb{C}\rightarrow\mathbb{F}$), \textbf{$\left(p,\infty\right)$-adic
analysis} ($\mathbb{F}\rightarrow\mathbb{C}$), \textbf{$\left(p,p\right)$-adic
analysis} ($\mathbb{F}\rightarrow\mathbb{F}$), and \textbf{$\left(p,q\right)$-adic
analysis} ($\mathbb{F}\rightarrow K$). In this terminology, real
and complex analysis would fall under the label \textbf{$\left(\infty,\infty\right)$-adic
analysis}.

The first flavor\textemdash $\left(\infty,p\right)$-adic analysis
($\mathbb{C}\rightarrow\mathbb{F}$)\textemdash I must confess, I
know little about, and have seen little, if any, work done with it.
Hopefully, this will soon be rectified. Each of the remaining three,
however, merit deeper discussion.

\index{analysis!left(p,inftyright)-adic@$\left(p,\infty\right)$-adic}The
second flavor of non-archimedean analysis\textemdash $\left(p,\infty\right)$-adic
analysis ($\mathbb{F}\rightarrow\mathbb{C}$) is the \emph{least}
exotic of the three. Provided that $\mathbb{F}$ is locally compact
(so, an at most finite-degree extension of $\mathbb{Q}_{p}$), this
theory is addressed by modern abstract harmonic analysis, specifically
through \textbf{Pontryagin duality} and the theory of Fourier analysis
on a locally compact abelian group. Tate's Thesis is of this flavor.
In addition to Tate's legacy, the past thirty years have seen a not-insignificant
enthusiasm in $\left(p,\infty\right)$-adic analysis from theoretical
physics, of all places \cite{First 30 years of p-adic mathematical physics,Quantum Paradoxes,p-adic space-time}
(with Volovich's 1987 paper \cite{p-adic space-time} being the ``classic''
of this burgeoning subject). This also includes the adèlic variant
of $\left(p,\infty\right)$-adic analysis, where $\mathbb{F}$ is
instead replaced with $\mathbb{A}_{\mathbb{Q}}$, the \textbf{adèle
ring }of $\mathbb{Q}$; see \cite{Adelic Harmonic Oscillator}, for
instance, for an \emph{adèlic} harmonic oscillator. This body of work
is generally known by the name $p$-adic (mathematical) physics and
$p$-adic quantum mechanics, and is as mathematically intriguing as
it is conceptually audacious. The hope in the depths of this Pandora's
box is the conviction that non-archimedean mathematics might provide
a better model for reality at its smallest scales. On the side of
pure mathematics, $\left(p,\infty\right)$-adic analysis is a vital
tool in representation theory (see, for instance \cite{Automorphic Representations}).
In hindsight, this is not surprising: the abstract harmonic analysis
used to formulate Fourier theory over locally compact abelian groups
must appeal to representation theory in order to formulate Fourier
analysis over an arbitrary compact group.

Of the many flavors of non-archimedean analysis, I feel that only
the third\textemdash \index{analysis!left(p,pright)-adic@$\left(p,p\right)$-adic}$\left(p,p\right)$-adic
anlaysis ($\mathbb{F}\rightarrow\mathbb{F}$)\textemdash truly deserves
to be called ``\index{analysis!$p$-adic}$p$-adic analysis''. Of
course, even this flavor has its own variations. One can study functions
$\mathbb{Z}_{p}\rightarrow\mathbb{F}$, as well as functions $\mathbb{F}\rightarrow\mathbb{F}^{\prime}$,
where $\mathbb{F}^{\prime}$ is a field extension of $\mathbb{F}$,
possibly of infinite degree. The utility, flexibility, adaptability,
and number-theoretic import of $\left(p,p\right)$-adic analysis have
all contributed to the high esteem and widespread popularity this
subject now enjoys. $\left(p,p\right)$-adic analysis can be studied
in its own right for a single, fixed prime $p$, or it can be used
in conjunction with the local-global philosophy, wherein one does
$\left(p,p\right)$-adic analysis for many\textemdash even infinitely
many\textemdash different values of $p$ simultaneously. This synergy
is on full display in one of $\left(p,p\right)$-adic analysis most
famous early triumphs: Bernard Dwork\index{Dwork, Bernard}'s proof
(\cite{Dwork Zeta Rationality}) of the rationality of the zeta function
for a finite field\textemdash the first of the three parts of the
\textbf{Weil Conjectures }to be proven true \cite{Dwork Zeta Rationality,Koblitz's other book,Robert's Book,p-adic proof for =00003D0003C0}.

Much like the local-global principle itself, Dwork's application of
$p$-adic analysis is a stratospheric elevation of the Fundamental
Theorem of Arithmetic: every integer $\geq2$ is uniquely factorable
as a product of primes. This is evident in the so-called \textbf{Product
Formula} of algebraic number theory, applied to $\mathbb{Q}$:

\begin{equation}
\left|a\right|_{\infty}\cdot\prod_{p\in\mathbb{P}}\left|a\right|_{p}=1,\textrm{ }\forall a\in\mathbb{Q}\backslash\left\{ 0\right\} \label{eq:The Product Formula}
\end{equation}
where $\left|a\right|_{p}$ is the $p$-adic absolute value of $a$,
and $\left|a\right|_{\infty}$ is the usual absolute value (on $\mathbb{R}$)
of $a$. Indeed, (\ref{eq:The Product Formula}) is but a restatement
of the Fundamental Theorem of Arithmetic.

The key tool in Dwork's proof was the following theorem: 
\begin{thm}[\textbf{Dwork}\footnote{Cited in \cite{p-adic proof for =00003D0003C0}.}]
Let $f\left(z\right)=\sum_{n=0}^{\infty}a_{n}z^{n}$ be a power series
in the complex variable $\mathbb{C}$ with coefficients $\left\{ a_{n}\right\} _{n\geq0}$
in $\mathbb{Q}$. Let $S$ be a set of finitely many places\footnote{For the uninitiated, a ``place'' is a fancy term for a prime number.
In algebraic number theory, primes\textemdash and, more generally,
prime \emph{ideals}\textemdash of a given field or ring are associated
with the valuations they induce. Due to a technicality, much like
a Haar measure on a group, valuations are not unique. The absolute
value $\mathfrak{z}\in\mathbb{Z}_{2}\mapsto3^{-v_{2}\left(\mathfrak{z}\right)}$,
for example, defines a topology on $\mathbb{Z}_{2}$ which is equivalent
to that of the standard $2$-adic absolute value. The valuation associated
to this modified absolute value would be $\mathfrak{z}\mapsto v_{2}\left(\mathfrak{z}\right)\log_{2}3$.
A \textbf{place}, then, is an equivalence class of valuations associated
to a given prime ideal, so that any two valuations in the class induce
equivalent topologies on the ring or field under consideration.} of $\mathbb{Q}$ which contains the archimedean complex place\footnote{Meaning the valuation corresponding to the ordinary absolute value}.
If $S$ satisfies:

\vphantom{}

I. For any\footnote{Here, $p$ represents \emph{both} a prime number and the associated
non-archimedean absolute value.} $p\notin S$, $\left|a_{n}\right|_{p}\leq1$ for all $n\geq0$ (i.e.,
$a_{n}\in\mathbb{Z}_{p}$);

\vphantom{}

II. For any\footnote{Here, $v$ is a absolute value on $\mathbb{Q}$, and $\mathbb{C}_{v}$
is the algebraic closure of the $v$-adic completion of $\mathbb{Q}$.} $p\in S$, there is a disk $D_{p}\subseteq\mathbb{C}_{p}$ of radius
$R_{p}$ so that $f\left(z\right)$ extends to a meromorphic function
$D_{p}\rightarrow\mathbb{C}_{p}$ and that $\prod_{p\in S}R_{p}>1$;

\vphantom{}

Then $f\left(z\right)$ is a rational function. 
\end{thm}
\vphantom{}

This is a $p$-adic generalization of an observation, first made by
Borel, that if a power series $\sum_{n=0}^{\infty}a_{n}z^{n}$ with
\emph{integer} coefficients has a radius of convergence strictly greater
than $1$, the power series is necessarily\footnote{Proof: let $R$ be the radius of convergence, and suppose $R>1$.
Then, there is an $r\in\left(1,R\right)$ for which $f\left(z\right)=\sum_{n=0}^{\infty}a_{n}z^{n}$
converges uniformly for $\left|z\right|=r$, and hence, by the standard
Cauchy estimate: 
\[
\left|a_{n}\right|=\left|\frac{f^{\left(n\right)}\left(0\right)}{n!}\right|\leq\frac{1}{r^{n}}\sup_{\left|z\right|=r}\left|f\left(z\right)\right|
\]
Since $r>1$, the upper bound will be $<1$ for all sufficiently large
$n$. However, if the $a_{n}$s are all integers, then any $n$ for
which $\left|a_{n}\right|<1$ necessarily forces $a_{n}=0$. Hence
$f$ is a polynomial.} a polynomial. It is a testament to the power of Dwork's Theorem\textemdash via
a generalization due Bertrandias\textemdash that it can be used to
prove the \index{Lindemann-Weierstrass Theorem}\textbf{Lindemann-Weierstrass
Theorem}\textemdash the Lindemann-Weierstrass Theorem being the result
that proved $\pi$\index{transcendence of pi@transcendence of $\pi$}
to be a transcendental number, thereby resolving the millennia-old
problem of squaring the circle\textemdash the answer being a resounding
\emph{no}. See \cite{p-adic proof for =00003D0003C0} for a lovely
exposition of this argument.

The fourth, most general flavor of non-archimedean analysis concerns
functions between arbitrary non-archimedean valued fields. The subject
has gone even further, studying functions $\Omega\rightarrow\mathbb{F}$,
where $\Omega$ is any appropriately well-behaved Hausdorff space.
Of the different flavors, I feel this is the one most deserving of
the title ``\textbf{non-archimedean analysis}\index{analysis!non-archimedean}''.\textbf{
}A good deal of the work done in non-archimedean analysis is in the
vein of so-called ``soft analysis'' variety, using the abstract
and often heavily algebraified language common to modern real or complex
functional analysis\textemdash ``maximal ideals'', ``convex sets'',
``nets'', and the like\textemdash to investigate exotic generalities.
Aside from generality's sake, part of the reason for the lack of a
significant distinction in the literature between non-archimedean
analysis and non-archimedean \emph{functional }analysis is one the
central difficulties of non-archimedean analysis: \emph{the absence
of a ``nice'' analogue of differentiation in the classical sense.}

In the naïveté of my youth, I dared to dream that calculus in multiple
variables\textemdash once I learned it\textemdash would be every bit
as elegant and wondrous as its single-variable counterparts. Inevitably,
I was disabused of this notion. When it comes to differential calculus,
single-variable $\left(p,\infty\right)$-adic analysis (ex: $\mathbb{Z}_{p}\rightarrow\mathbb{C}$)
and non-archimedean analysis (such as $\mathbb{Z}_{p}\rightarrow\mathbb{C}_{q}$)
manage to surpass multi-variable real analysis in theoretical complexity
without needing to draw on even half as much indecency of notation.
The main problem in these settings is that there is no way to define
a \emph{difference quotient }for functions. The topologies and algebraic
operations of functions' inputs and outputs are \emph{completely}
different. You can't divide a complex number by a $p$-adic integer,
and so on and so forth.

On the other hand, the theory of integration rests on much firmer
ground. Thanks to Haar measures and Pontryagin Duality, integration
in $\left(p,\infty\right)$-adic analysis falls squarely under the
jurisdiction of standard measure theory and abstract harmonic analysis.
In fact, the parallels between $\left(p,\infty\right)$-adic Fourier
analysis and classical Fourier analysis on tori or euclidean space
are so nice that, in 1988, Vasilii Vladimirov was able to use Fourier
analysis to introduce $\left(p,\infty\right)$-adic differentiation
in the sense of distributions. This was done by co-opting the classic
\emph{complex-valued }$p$-adic Mellin integral \cite{Vladimirov - the big paper about complex-valued distributions over the p-adics}:
\begin{equation}
\int_{\mathbb{Q}_{p}}\left|\mathfrak{z}\right|_{p}^{s-1}f\left(\mathfrak{z}\right)d\mathfrak{z}
\end{equation}
where $d\mathfrak{z}$ is the \emph{real-valued} Haar measure on $\mathbb{Q}_{p}$,
normalized to be a probability measure on $\mathbb{Z}_{p}$ and where
$s$ is a complex variable. Here, $f$ is a complex-valued function
on $\mathbb{Q}_{p}$, and\textemdash for safety\textemdash we assume
$f$ is compactly supported. The fruit of this approach is the \index{fractional differentiation}\index{Vladimirov operator}\textbf{
Vladimirov (Fractional) Differentiation / Integration operator} of
order $\alpha$ (where $\alpha$ is a real number $>0$). For a function
$f:\mathbb{Q}_{p}\rightarrow\mathbb{C}$, this is given by the integral:
\begin{equation}
D^{\alpha}\left\{ f\right\} \left(\mathfrak{z}\right)\overset{\textrm{def}}{=}\frac{1}{\Gamma_{p}\left(-\alpha\right)}\int_{\mathbb{Q}_{p}}\frac{f\left(\mathfrak{z}\right)-f\left(\mathfrak{y}\right)}{\left|\mathfrak{z}-\mathfrak{y}\right|_{p}^{1+\alpha}}d\mathfrak{y}\label{eq:Definition of the Vladimirov Fractional Differentiation Operator}
\end{equation}
provided the integral exists. Here, $\Gamma_{p}$ is the physicists'
notation for: 
\begin{equation}
\Gamma_{p}\left(-\alpha\right)\overset{\textrm{def}}{=}\frac{p^{\alpha}-1}{1-p^{-1-\alpha}}\label{eq:Physicist's gamma_p}
\end{equation}
This operator can be extended to negative $\alpha$ by being rewritten
so as to produce the $\left(p,\infty\right)$-adic equivalent of the
usual, Fourier-transform-based definition of fractional derivatives
in fractional-order Sobolev spaces.

As for \emph{non-archimedean} analysis (including $\left(p,q\right)$),
there is a well developed theory of integration, albeit with a strong
functional analytic flavor. Integrals and measures are formulated
not as $K$-valued functions of sets, but continuous, $K$-valued
linear functionals on the Banach space of continuous functions $X\rightarrow K$,
where $X$ is an ultrametric space and $K$ is a metrically complete
non-archimedean field. An exposition of this theory is given in Subsection
\ref{subsec:3.1.6 Monna-Springer-Integration}.

To return to physics for a moment, one of the principal \emph{mathematical
}reasons for interest in $p$-adic analysis among quantum physicists
has been the hope of realize a theory of probability based on axioms
other than Kolmogorov's. The \emph{idée fixe }of this approach is
that sequences of rational numbers (particularly those representing
frequencies of events) which diverge in $\mathbb{R}$ can be convergent
in $\mathbb{Q}_{p}$ for an appropriate choice of $p$. The hope is
that this will allow for a resurrection of the pre-Kolmogorov frequency-based
probability theory, with the aim of using the frequency interpretation
to assign probabilistic meaning to sequences of rational numbers which
converge in a $p$-adic topology. Andrei Khrennikov and his associates
have spearheaded these efforts; \cite{Measure-theoretic approach to p-adic probability theory}
is an excellent (though dense) exposition of such a non-archimedean
probability theory, one which strives to establish a more classical
measure-based notion of probability. Khrennikov's eclectic monograph
\cite{Quantum Paradoxes} also makes for a fascinating read, containing
much more information, and\textemdash most importantly\textemdash worked
examples from both physics and probability.

Now for the bad news. Both $\left(p,\infty\right)$-adic and non-archimedean
analysis are effectively incompatible with any notion of analytic
functions; there are no power series, nor rational functions. $\left(p,p\right)$-adic
analysis \emph{does not }suffer this limitation, which is another
reason for its longstanding success. Even so, this theory of $p$-adic
analytic functions still managed to get itself bewitched when reality
was in the drafting stage: the non-archimedean structure of the $p$-adics
makes even the most foundational principles of classical analytic
function theory completely untenable in the $p$-adic setting.

While one \emph{can} define the derivative of a function, say, $f:\mathbb{Z}_{p}\rightarrow\mathbb{Q}_{p}$
as the limit of difference quotient: 
\begin{equation}
f^{\prime}\left(\mathfrak{z}\right)=\lim_{\mathfrak{y}\rightarrow\mathfrak{z}}\frac{f\left(\mathfrak{y}\right)-f\left(\mathfrak{z}\right)}{\mathfrak{y}-\mathfrak{z}}\label{eq:Naive differentiability}
\end{equation}
doing so courts disaster. Because there is sense in the world, this
notion of ``derivative\index{derivative}'' \emph{does }work for
analytic functions. However, it fails \emph{miserably }for everything
else. Yes, this kind of differentiability still implies continuity\textemdash that
much, we are allowed to keep. Nevertheless, there are simple examples
of \emph{differentiable functions $f:\mathbb{Z}_{p}\rightarrow\mathbb{Q}_{p}$
}(in the sense that (\ref{eq:Naive differentiability}) exists at
every point)\emph{ which are injective, but whose derivatives are
everywhere zero!}
\begin{example}
\label{exa:p-adic differentiation is crazy}Given\footnote{As a minor spoiler, this example proves to be useful to one of my
ideas for studying $\chi_{H}$, because \emph{of course }it would.} any $\mathfrak{z}\in\mathbb{Z}_{p}$, write $\mathfrak{z}$ as $\sum_{n=0}^{\infty}\mathfrak{z}_{n}p^{n}$,
where, $\mathfrak{z}_{n}\in\left\{ 0,\ldots,p-1\right\} $ is the
$n$th $p$-adic digit of $\mathfrak{z}$, and consider the function\footnote{Amusingly, this function will turn out to be relevant to studying
Collatz and the like; see the lead-up to \textbf{Example \ref{exa:L^1 method example}}
(page \pageref{exa:L^1 method example}) from Subsection \ref{subsec:4.3.4 Archimedean-Estimates}
in Chapter 4 for the details.} $\phi:\mathbb{Z}_{p}\rightarrow\mathbb{Q}_{p}$ defined by: 
\begin{equation}
\phi\left(\mathfrak{z}\right)=\sum_{n=0}^{\infty}\mathfrak{z}_{n}p^{2n}\label{eq:Phi for the L^1 method}
\end{equation}
Because $\left|\mathfrak{z}-\mathfrak{y}\right|_{p}=1/p^{n+1}$, where
$n$ is the largest integer $\geq0$ for which $\mathfrak{z}_{n}=\mathfrak{y}_{n}$
(with the convention that $n=-1$ if no such $n$ exists), we have
that: 
\begin{equation}
\left|\phi\left(\mathfrak{z}\right)-\phi\left(\mathfrak{y}\right)\right|_{p}=\frac{1}{p^{2\left(n+1\right)}}=\left|\mathfrak{z}-\mathfrak{y}\right|_{p}^{2}
\end{equation}
Thus, $\phi$ is continuous, injective, and differentiable, with:
\begin{equation}
\left|\phi^{\prime}\left(\mathfrak{z}\right)\right|_{p}=\lim_{\mathfrak{y}\rightarrow\mathfrak{z}}\left|\frac{\phi\left(\mathfrak{y}\right)-\phi\left(\mathfrak{z}\right)}{\mathfrak{y}-\mathfrak{z}}\right|_{p}=\lim_{\mathfrak{y}\rightarrow\mathfrak{z}}\left|\mathfrak{y}-\mathfrak{z}\right|_{p}=0,\textrm{ }\forall\mathfrak{z}\in\mathbb{Z}_{p}
\end{equation}
which is \emph{not} what one would like to see for an injective function!
The notion of \textbf{strict differentiability }(see \cite{Robert's Book,Ultrametric Calculus})
is needed to get a non-degenerate notion of differentiability for
non-analytic functions. In short, this is the requirement that we
take the limit of (\ref{eq:Naive differentiability}) as\emph{ both
}$\mathfrak{y}$ and $\mathfrak{z}$ to converge to a single point
$\mathfrak{z}_{0}$, and that, regardless of path chosen, the limit
exists; this then shows that $f$ is ``strictly differentiable''
at $\mathfrak{z}_{0}$.
\end{example}
\vphantom{}

To continue our tour of the doom and the gloom, the $p$-adic Mean
Value Theorem is but a shadow of its real-analytic self; it might
as well not exist at all. Unsurprisingly, $p$-adic analysis suffers
when it comes to integration. Without a dependable Mean Value Theorem,
the relationship between integration and anti-differentiation familiar
to us from real and complex analysis fails to hold in the $\left(p,p\right)$-adic
context. But the story only gets worse from there.

Surely, one of the most beautiful syzygies in modern analysis is the
intermingling of measure theory and functional analysis. This is made
all the more elegant when Haar measures and harmonic analysis are
brought into the mix. Unfortunately, in $\left(p,p\right)$-adic waters,
this elegance bubbles, burns, and crumbles into nothingness. Basic
considerations of algebra and arithmetic show that any $\left(p,p\right)$-adic
``Haar measure'' (translation-invariant $\mathbb{F}$-valued function
of ``nice'' sets in $\mathbb{Z}_{p}$) must assign the rational
number $1/p^{n}$ as the measure of sets of the form $k+p^{n}\mathbb{Z}_{p}$,
where $k\in\mathbb{Z}$ and $n\in\mathbb{N}_{0}$. Because the $p$-adic
absolute value of $1/p^{n}$ tends to $\infty$ as $n\rightarrow\infty$,
the only \emph{continuous}, translation-invariant $\mathbb{Q}_{p}$-valued
function of sets on $\mathbb{Z}_{p}$ is the constant function $0$\textemdash the
zero measure. This same result holds for any ring extension of $\mathbb{Z}_{p}$,
as well as for $\mathbb{Q}_{p}$ and any field extension thereof.

Fortunately, there ways of getting around this sorry state of affairs,
even if the result ends up being completely different than classical
integration theory. The first of these is the \textbf{Volkenborn integral},
which\index{integral!Volkenborn} generalizes the classical notion
of an integral as the limit of a Riemann sum\index{Riemann sum}.
The Volkenborn integral\footnote{Robert covers Volkenborn integration in Section 5 of his chapter on
Differentiation in \cite{Robert's Book}, starting at page 263. One
can also turn to Borm's dissertation \cite{pp adic Fourier theory}.} of a function $f:\mathbb{Z}_{p}\rightarrow\mathbb{C}_{q}$ is defined
by the limit: 
\begin{equation}
\int_{\mathbb{Z}_{p}}f\left(\mathfrak{z}\right)d_{\textrm{Volk}}\mathfrak{z}\overset{\textrm{def}}{=}\lim_{N\rightarrow\infty}\frac{1}{p^{N}}\sum_{n=0}^{p^{N}-1}f\left(n\right)\label{eq:definition of the volkenborn integral}
\end{equation}
which\textemdash the reader should note\textemdash is \emph{not }translation
invariant (see \cite{Robert's Book})! The Volkenborn integral turns
out to be intimately connected with the finite differences and the
indefinite summation operator. In all likelihood, the most noteworthy
property of the Volkenborn integral is that, for any prime $p$ and
any integer $n\geq0$, the integral of the function $\mathfrak{z}^{n}$
is equal to the $n$th Bernoulli number (with $B_{1}=-1/2$).

In practice, however, the Volkenborn integral is neglected in favor
for the Amice-Iwasawa\index{Amice, Yvette} approach, otherwise known
as the theory of \index{$p$-adic!distribution}\textbf{$p$-adic distributions}\footnote{Note, these should \textbf{\emph{not}}\emph{ }be confused with the
distributions / ``generalized functions'' dealt with by Vladimirov
in \cite{Vladimirov - the big paper about complex-valued distributions over the p-adics}.
Amice-Iwasawa ``$p$-adic distributions'' are linear functionals
taking values in $\mathbb{Q}_{p}$ or a finite degree field extension
thereof. Vladimirov's ``$p$-adic distributions'' are linear functionals
taking values in $\mathbb{C}$.} (\cite{p-adic L-functions paper,Iwasawa}). Let $\mathbb{F}$ be
a metrically complete field extension of $\mathbb{Q}_{p}$. The idea
is to exploit the fact that the space $C\left(\mathbb{Z}_{p},\mathbb{F}\right)$
has as a basis the binomial coefficient polynomials ($\left\{ \binom{\mathfrak{z}}{n}\right\} _{n\geq0}$):
\begin{align*}
\binom{\mathfrak{z}}{0} & =1\\
\binom{\mathfrak{z}}{1} & =\mathfrak{z}\\
\binom{\mathfrak{z}}{2} & =\frac{\mathfrak{z}\left(\mathfrak{z}-1\right)}{2}\\
 & \vdots
\end{align*}
That these functions form a basis of $C\left(\mathbb{Z}_{p},\mathbb{F}\right)$
is a classic result of $p$-adic analysis first proved by Kurt Mahler
in 1958\footnote{\cite{Mahler Series} is Mahler's original paper, but his result is
now standard material in any worthwhile course on $p$-adc analysis
\cite{Mahler,Robert's Book,Ultrametric Calculus}.}. This basis is now known as the\textbf{ Mahler basis}, in his honor.
From this, one can define $\mathbb{F}$-valued measures and distributions
by specifying what they do to $\binom{\mathfrak{z}}{n}$ for each
$n$, and then extending by linearity. This construction is of considerable
importance in modern algebraic number theory, where it is one of several
ways of defining $p$-adic $L$ functions\index{$p$-adic!$L$-function};
the fact that the several different methods of constructing $p$-adic
$L$ functions are all equivalent was a major landmark of twentieth
century number theory \cite{p-adic L-functions paper,Iwasawa}. So,
it's no wonder the Volkenborn integral has gotten the proverbial cold
shoulder.

Despite all this, the particular sub-discipline this dissertation
focuses on remains essentially untouched.\textbf{ }In terms of the
four flavor classification given above, \textbf{$\left(p,q\right)$-adic
analysis}\index{analysis!left(p,qright)-adic@$\left(p,q\right)$-adic}\textemdash the
study of functions from $\mathbb{Z}_{p}$ (or a metrically complete
field extension thereof) to $\mathbb{Z}_{q}$ (or a metrically complete
field extension thereof), where $p$ and $q$ are distinct primes,
is technically a form of non-archimedean analysis. Even so, this dissertation
stands in stark contrast to much of the established literature in
non-archimedean analysis not only by its unusual content, but also
by virtue of its concreteness. One generally needs to turn to works
of mathematical physics like \cite{Quantum Paradoxes} to get comparable
levels of computational depth. Also significant, unlike several works
going by the title ``non-archimedean analysis'' (or some variation
thereof), my dissertation has the distinction of actually \emph{being
}a work of mathematical analysis, as opposed to algebra going by a
false name\footnote{A particularly egregious example of this is \cite{Bosch lying title},
a text on $p$-adic algebraic geometry and Tate's ``rigid analytic
geometry'' which has the \emph{nerve} to call itself \emph{Non-archimedean
analysis}, despite not containing so much as a single line of bonafide
analysis. Other examples of such wolves-in-sheep's-clothing are \cite{Schneider awful book,More Schneider lies},
both due to Peter Schneider.}. And, by mathematical analysis, I mean \emph{hard} analysis: $\epsilon$s
and $\delta$s, approximations, asymptotics, and the like, and all
the detailed computations entailed.

The past twenty or so years have seen many new strides in $p$-adic
and non-archimedean analysis, particularly in the vein of generalizing
as many of classical analysis' most useful tools, principally to the
setting of $\left(p,\infty\right)$-adic analysis, but also to $\left(p,p\right)$
as well. The tools in question include, asymptotic analysis, Tauberian
theorems, distributions, oscillatory integrals (a $p$-adic van der
Corput lemma), , and many others; see, \cite{p-adic Tauberian}, \cite{Vladimirov - the big paper about complex-valued distributions over the p-adics,Volo - p-adic Wiener},
\cite{p-adic van der Corput lemma,Real and p-Adic Oscillatory integrals}).
However, because $\left(p,q\right)$-adic analysis works with functions
taking values in non-archimedean fields, the vast majority of these
tools and their associated proofs do not extend to the $\left(p,q\right)$-adic
case\textemdash at least, not yet. And it is precisely there where
we will begin our work.

\subsection{\label{subsec:3.1.2 Banach-Spaces-over}Banach Spaces over a Non-Archimedean
Field}

THROUGHOUT THIS SUBSECTION, $K$ IS A COMPLETE NON-ARCHIMEDEAN VALUED
FIELD.

\vphantom{}

In classical analysis, Banach spaces frequently trot onto the stage
when one moves from asking ``biographical'' questions about specific
functions, specific equations, and the like, to ``sociological''
ones\textemdash those regarding whole classes of functions. Much like
with real or complex analysis, there is quite a lot one can do in
$\left(p,p\right)$-adic analysis without needing to appeal to the
general theories of Banach spaces and functional analysis. Power series
and Mahler's basis for continuous $\left(p,p\right)$-adic functions
give $p$-adic analysis enough machinery to make concreteness just
as viable as in the real or complex case. In $\left(p,q\right)$-adic
analysis, however, the lack of any notions of differentiation, power
series, and the non-applicability of the Mahler Basis make for an
entirely different set of circumstances. Rather than drawing from
the broader perspectives offered by functional analysis as a means
of enriching how we approach what we already know, Banach spaces are
foundational to non-archimedean analysis precisely because they can
be used to determine what is or is not actually \emph{possible} in
the subject.

As a simple example, consider Hilbert spaces and orthogonal functions.
Investigations into Fourier series and differential equations in the
nineteenth century motivated subsequent explorations of function spaces
as a whole, especially in the abstract. The Cauchy-Schwarz Inequality
came before inner product spaces. In non-archimedean analysis, however,
it can be shown that the only inner product spaces over non-archimedean
fields which are Hilbert spaces in the classical sense\footnote{Specifically, their norm is induced by a bilinear form, and orthogonal
projections exist for every closed subspace.} are necessarily finite-dimensional (\textbf{Theorem 4 }from \cite{Schikhof Banach Space Paper}).
Sobering information such as this gives non-archimedean analysis an
idea of what they should or should not set their hopes on.

Most of the exposition given here is taken from Schikhof's excellent
article on Banach spaces over a non-archimedean field \cite{Schikhof Banach Space Paper},
which\textemdash unlike van Rooij's book \cite{van Rooij - Non-Archmedean Functional Analysis}\textemdash is
still in print.

We begin with the basic definitions. 
\begin{defn}[S\textbf{emi-norms, normed vector spaces, etc.} \cite{Ultrametric Calculus}]
\label{def:seminorms etc.}Let $E$ be a vector space over a non-archimedean
valued field $K$ (a.k.a., $E$ is a \textbf{$K$-vector space} or
\textbf{$K$-linear space}). A\index{non-archimedean!semi-norm} \textbf{(non-archimedean)
semi-norm }on $E$ is a function $\rho:E\rightarrow\mathbb{R}$ such
that for all $x,y\in E$ and all $\mathfrak{a}\in K$:

\vphantom{}

I. $\rho\left(x\right)\geq0$;

\vphantom{}

II. $\rho\left(\mathfrak{a}x\right)=\left|\mathfrak{a}\right|_{K}\rho\left(x\right)$;

\vphantom{}

III. $\rho\left(x+y\right)\leq\max\left\{ \rho\left(x\right),\rho\left(y\right)\right\} $.

\vphantom{}

Additionally, $\rho$ is said to be a \textbf{(non-archimedean) norm
}over\index{non-archimedean!norm} $E$ whenever it satisfies the
additional condition:

\vphantom{}

IV. $\rho\left(x\right)=0$ if and only if $x=0$.

\vphantom{}

Next, given a non-archimedean norm $\left\Vert \cdot\right\Vert $,
the pair $\left(E,\left\Vert \cdot\right\Vert \right)$ is called
a \textbf{(non-archimedean) normed vector space}; this is then an
ultrametric space with respect to the metric induced by the norm.
We say $\left(E,\left\Vert \cdot\right\Vert \right)$ is a \textbf{(non-archimedean)
Banach space }when\index{non-archimedean!Banach space}\index{Banach algebra!non-archimedean}
it is complete as a metric space.

Finally, a \textbf{(non-archimedean) Banach algebra }is\index{non-archimedean!Banach algebra}
a Banach space $\left(E,\left\Vert \cdot\right\Vert \right)$ with
a multiplication operation $m:E\times E\rightarrow E$ so that $\left\Vert m\left(x,y\right)\right\Vert \leq\left\Vert x\right\Vert \left\Vert y\right\Vert $
for all $x,y\in E$. Also, a $K$-Banach space (resp. algebra) is
a Banach space (resp. algebra) over the field $K$. 
\end{defn}
\begin{prop}
\label{prop:series characterization of a Banach space}Let $E$ be
a non-archimedean normed vector space. Then, $E$ is complete if and
only if, for every sequence $\left\{ x_{n}\right\} _{n\geq0}$ in
$E$ for which $\lim_{n\rightarrow\infty}\left\Vert x_{n}\right\Vert \overset{\mathbb{R}}{=}0$,
the series $\sum_{n=0}^{\infty}x_{n}$ converges to a limit in $E$. 
\end{prop}
Proof:

I. Suppose $E$ is complete. Then, $\left(E,\left\Vert \cdot\right\Vert \right)$
is a complete ultrametric space, and as such, an infinite series $\sum_{n=0}^{\infty}x_{n}$
in $E$ converges if and only if $\lim_{n\rightarrow\infty}\left\Vert x_{n}\right\Vert \overset{\mathbb{R}}{=}0$.

\vphantom{}

II. Conversely, suppose an infinite series $\sum_{n=0}^{\infty}x_{n}$
in $E$ converges if and only if $\lim_{n\rightarrow\infty}\left\Vert x_{n}\right\Vert \overset{\mathbb{R}}{=}0$.
Then, letting $\left\{ x_{n}\right\} _{n\geq0}$ be a Cauchy sequence
in $E$, we can choose a sequence $n_{1}<n_{2}<\cdots$ so that, for
any $j$, $\left\Vert x_{m}-x_{n}\right\Vert <2^{-j}$ holds for all
$m,n\geq n_{j}$. Setting $y_{1}=x_{n_{1}}$ and $y_{j}=x_{n_{j}}-x_{n_{j-1}}$
gives: 
\begin{equation}
\lim_{J\rightarrow\infty}x_{n_{J}}\overset{E}{=}\lim_{J\rightarrow\infty}\sum_{j=1}^{J}y_{j}
\end{equation}
By construction, the Cauchyness of the $x_{n}$s guarantees that $\left\Vert y_{j}\right\Vert \rightarrow0$
as $j\rightarrow\infty$. As such, our assumption ``every infinite
series in $E$ converges if and only if its terms tend to zero in
$E$-norm'' forces $\lim_{J\rightarrow\infty}\sum_{j=1}^{J}y_{j}$
(and hence, $\lim_{J\rightarrow\infty}x_{n_{J}}$) to converge to
a limit in $E$. Since the $x_{n}$s are Cauchy, the existence of
a convergent subsequence then forces the entire sequence to converge
in $E$, with: 
\begin{equation}
\lim_{n\rightarrow\infty}x_{n}\overset{E}{=}\sum_{j=1}^{\infty}y_{j}
\end{equation}
Because the $x_{n}$s were arbitrary, this shows that $E$ is complete.

Q.E.D.

\vphantom{}

Next, we give the fundamental examples of non-archimedean Banach spaces.
Throughout, we will suppose that $K$ is a complete non-archimedean
valued field.

We address the finite-dimensional cases first. 
\begin{defn}
Given a vector space $V$ over a non-archimedean field $K$, two non-archimedean
norms $\left\Vert \cdot\right\Vert $ and $\left\Vert \cdot\right\Vert ^{\prime}$
defined on $V$ are said to be \textbf{equivalent }whenever there
are real constants $0<c\leq C<\infty$ so that: 
\begin{equation}
c\left\Vert \mathbf{x}\right\Vert \leq\left\Vert \mathbf{x}\right\Vert ^{\prime}\leq C\left\Vert \mathbf{x}\right\Vert 
\end{equation}
\end{defn}
\cite{Robert's Book} gives the following basic results: 
\begin{thm}
Let $V$ and $W$ be normed vector spaces over $\mathbb{Q}_{p}$.

\vphantom{}

I. If $V$ is finite-dimensional, then all non-archimedean norms on
$V$ are equivalent.

\vphantom{}

II. If $V$ is locally compact, then $V$ is finite dimensional.

\vphantom{}

III. If $V$ is locally compact, then $V$ has the Heine-Borel property\textemdash sets
in $V$ are compact if and only if they are closed and bounded.

\vphantom{}

IV. If $V$ and $W$ are finite-dimensional, then every linear map
$L:V\rightarrow W$ is continuous. 
\end{thm}
\begin{defn}
\ 

I. In light of the above, we shall use only a single norm on finite-dimensional
vector spaces over $\mathbb{Q}_{p}$: the $\ell^{\infty}$-norm (written
$\left\Vert \cdot\right\Vert _{p}$), which outputs the maximum of
the $p$-adic absolute values of the entries of an element of the
vector space.

\vphantom{}

II. Given integers $\rho,c\geq1$, we write \nomenclature{$K^{\rho,c}$}{$\rho\times c$ matrices with entries in $K$ \nopageref}$K^{\rho,c}$
to denote the set of all $\rho\times c$ matrices with entries in
$K$. We make this a non-archimedean normed vector space by equipping
it with the $\infty$-norm, which, for any $\mathbf{A}\in K^{\rho,c}$
with entries $\left\{ a_{i,j}\right\} _{1\leq i\leq\rho,1,j\leq c}$,
is given by: \nomenclature{$\left\Vert \mathbf{A}\right\Vert _{K}$}{max. of the $K$-adic absolute values of a matrix's entries \nopageref}\nomenclature{$\left\Vert \mathbf{a}\right\Vert _{K}$}{max. of the $K$-adic absolute values of a vector's entries \nopageref}
\begin{equation}
\left\Vert \mathbf{A}\right\Vert _{K}\overset{\textrm{def}}{=}\max_{1\leq i\leq\rho,1\leq j\leq c}\left|a_{i,j}\right|_{K}\label{eq:Definition of the non-archimedean matrix norm}
\end{equation}
When $c=1$, we identify $K^{\rho,1}$ with the set of all $\rho\times1$
column vectors.

\vphantom{}

III. Given $d\geq2$, we write $\textrm{GL}_{d}\left(K\right)$\nomenclature{$\textrm{GL}_{d}\left(K\right)$}{set of invertible $d\times d$ matrices with entries in $K$ \nopageref}
to denote the set of all invertible $d\times d$ matrices with entries
in $K$, made into a non-abelian group by matrix multiplication.

\vphantom{}

IV. Because we will need it for parts of Chapter 5, we define \nomenclature{$\left\Vert \mathbf{x}\right\Vert _{\infty}$}{$\overset{\textrm{def}}{=}\max\left\{ \left|x_{1}\right|,\ldots,\left|x_{d}\right|\right\}$}$\left\Vert \cdot\right\Vert _{\infty}$
on $\mathbb{C}^{d}$ (or any subset thereof) by:
\begin{equation}
\left\Vert \mathbf{x}\right\Vert _{\infty}\overset{\textrm{def}}{=}\max\left\{ \left|x_{1}\right|,\ldots,\left|x_{d}\right|\right\} \label{eq:Definition of infinity norm}
\end{equation}
where $\mathbf{x}=\left(x_{1},\ldots,x_{d}\right)$. Here, the absolute
values are those on $\mathbb{C}$.
\end{defn}
\begin{prop}
$K^{\rho,c}$ is a non-archimedean Banach space over $K$. If $\rho=c$,
it is also a Banach algebra, with the multiplication operation being
matrix multiplication. 
\end{prop}
We will not mention matrix spaces until we get to Chapters 5 and 6,
where we they will be ubiquitous.

Next, we introduce the fundamental examples of \emph{infinite}-dimensional
non-archimedean Banach spaces. 
\begin{defn}
Let $X$ be a set and let $K$ be a metrically-complete non-archimedean
valued field.

\vphantom{}

I. A function $f:X\rightarrow K$ is said to be \textbf{bounded }if:
\begin{equation}
\sup_{x\in X}\left|f\left(x\right)\right|_{K}<\infty
\end{equation}
We get a norm \nomenclature{$\left\Vert f\right\Vert _{X,K}$}{$\sup_{x\in X}\left|f\left(x\right)\right|_{K}$}
out of this by defining: 
\begin{equation}
\left\Vert f\right\Vert _{X,K}\overset{\textrm{def}}{=}\sup_{x\in X}\left|f\left(x\right)\right|_{K}\label{eq:Definition of X,K norm}
\end{equation}
We call this the $K$\textbf{-supremum norm }\index{norm!$K$ supremum}
on $X$. We write \nomenclature{$B\left(X,K\right)$}{Bounded $K$-valued functions on $X$}$B\left(X,K\right)$
to denote the set of all bounded $K$-valued functions on $X$. This
is a non-archimedean Banach space under $\left\Vert \cdot\right\Vert _{X,K}$.

When $X=\mathbb{Z}_{p}$, we call the norm $\left\Vert \cdot\right\Vert _{X,K}$
the \textbf{$\left(p,K\right)$-adic norm}\index{norm!left(p,Kright)-adic@$\left(p,K\right)$-adic},
and denote it by $\left\Vert \cdot\right\Vert _{p,K}$\nomenclature{$\left\Vert \cdot\right\Vert _{p,K}$}{ }.

When $K$ is a $q$-adic field, we call this the \index{norm!left(p,qright)-adic@$\left(p,q\right)$-adic}\index{$p,q$-adic!norm}\textbf{$\left(p,q\right)$-adic
norm}, and denote it by \nomenclature{$\left\Vert \cdot\right\Vert _{p,q}$}{ }$\left\Vert \cdot\right\Vert _{p,q}$.

In a minor abuse of notation, we also use $\left\Vert \cdot\right\Vert _{p,K}$
and $\left\Vert \cdot\right\Vert _{p,q}$ to denote the $K$ supremum
norm on $X$ when $X$ is $\hat{\mathbb{Z}}_{p}$.

\vphantom{}

II. If $X$ is a compact Hausdorff space, we write \nomenclature{$C\left(X,K\right)$}{set of continuous $K$-valued functions on $X$}$C\left(X,K\right)$
to denote the set of all continuous functions $f:X\rightarrow K$.
This is a non-archimedean Banach space\footnote{This construction is of vital importance to us, because it allows
us to consider continuous $K$-valued functions defined on a compact
ultrametric space $X$ like $\mathbb{Z}_{p}$, even if $K$'s residue
field has characteristic $q\neq p$.} under $\left\Vert \cdot\right\Vert _{X,K}$.

\vphantom{}

III. We write $\text{\ensuremath{\ell^{\infty}\left(K\right)}}$ (or
$\ell^{\infty}$, if $K$ is not in question) to denote the set of
all bounded sequences in $K$. This is a non-archimedean Banach space
under the norm: 
\begin{equation}
\left\Vert \left(a_{1},a_{2},\ldots\right)\right\Vert _{K}\overset{\textrm{def}}{=}\sup_{n\geq1}\left|a_{n}\right|_{K}\label{eq:K sequence norm}
\end{equation}

\vphantom{}

IV. We write \nomenclature{$c_{0}\left(K\right)$}{set of $K$-valued sequences converging to $0$}$c_{0}\left(K\right)$
(or $c_{0}$, if $K$ is not in question) to denote the set of all
sequences in $K$ which converge to $0$ in $K$'s absolute value.
This is a Banach space under (\ref{eq:K sequence norm}) and is a
closed subspace of $\ell^{\infty}$.

\vphantom{}

V. We write\footnote{I've adapted this notation from what van Rooij calls ``$c_{0}\left(X\right)$''.
\cite{van Rooij - Non-Archmedean Functional Analysis}} $c_{0}\left(X,K\right)$ to denote the set of all $f\in B\left(X,K\right)$
so that, for every $\epsilon>0$, there are only finitely many $x\in X$
for which $\left|f\left(x\right)\right|_{K}\geq\epsilon$. This is
a Banach space with respect to $\left\Vert \cdot\right\Vert _{X,K}$,
being a closed subspace of $B\left(X,K\right)$.

\vphantom{}

VI. \nomenclature{$\text{\ensuremath{\ell^{1}\left(K\right)}}$}{set of absolutely summable $K$-valued sequences}We
write $\text{\ensuremath{\ell^{1}\left(K\right)}}$ (or $\ell^{1}$,
if $K$ is not in question) to denote the set of all sequences $\mathbf{c}:\mathbb{N}_{0}\rightarrow K$
so that: 
\[
\sum_{n=0}^{\infty}\left|\mathbf{c}\left(n\right)\right|_{K}<\infty
\]
We write: 
\begin{equation}
\left\Vert \mathbf{c}\right\Vert _{1}\overset{\textrm{def}}{=}\sum_{n=0}^{\infty}\left|\mathbf{c}\left(n\right)\right|_{K}\label{eq:Definition of ell-1 norm}
\end{equation}
$\ell^{1}$ is an \emph{archimedean }$K$-Banach space under $\left\Vert \cdot\right\Vert _{1}$\textemdash that
is, $\left\Vert \mathbf{a}+\mathbf{b}\right\Vert _{1}\leq\left\Vert \mathbf{a}\right\Vert _{1}+\left\Vert \mathbf{b}\right\Vert _{1}$,
however, the strong triangle inequality does not hold on it; that
is, there exist $\mathbf{a},\mathbf{b}$ for which $\left\Vert \mathbf{a}+\mathbf{b}\right\Vert _{1}>\max\left\{ \left\Vert \mathbf{a}\right\Vert _{1},\left\Vert \mathbf{b}\right\Vert _{1}\right\} $.
This also goes to show that not every Banach space over a non-archimedean
field is, itself, non-archimedean. 
\end{defn}
\begin{rem}
As a rule, when writing expressions involving $\left\Vert \cdot\right\Vert $,
I will use a single subscript (ex: $\left\Vert \cdot\right\Vert _{K}$)
whenever only the absolute value of $K$ is being used. Expressions
with \emph{two} subscripts (ex: $\left\Vert \cdot\right\Vert _{p,K}$,
$\left\Vert \cdot\right\Vert _{p,q}$), meanwhile, indicate that we
are using the absolute value indicated by the subscript on the right
as well as taking a supremum of an input variable with respect to
the absolute value indicated by the subscript on the left. 
\end{rem}
\vphantom{}

Although $\ell^{1}$ and its continuous counterpart $L^{1}$ are of
fundamental importance in real and complex analysis, their non-archimedean
analogues are not nearly as distinguished. The first hint that something
is amiss comes from the fact that $\ell^{1}$ is archimedean, rather
than non-archimedean.\emph{ Traditionally}\footnote{Schikhof himself says as much at the top of page 6 in \cite{Schikhof Banach Space Paper},
writing ``Absolute summability plays no a role in {[}non-archimedean{]}
analysis.''}, $\ell^{1}$ is of no interest in non-archimedean analysts. Unsurprisingly,
it turns out that $\ell^{1}$\textemdash specifically, its continuous
analogue, $L^{1}$) is pertinent to the study of $\chi_{H}$. This
is the subject of Subsection \ref{subsec:3.3.6 L^1 Convergence}.

Next, we have linear operators and functionals. 
\begin{defn}
Let $E$ and $F$ be non-archimedean Banach spaces, both over the
field $K$. Then one writes $\mathcal{L}\left(E,F\right)$ to denote
the set of all continuous linear operators from $E$ to $F$. We equip
$\mathcal{L}\left(E,F\right)$ with the \textbf{operator norm}:

\begin{equation}
\left\Vert T\right\Vert \overset{\textrm{def}}{=}\sup_{x\in E}\frac{\left\Vert T\left(x\right)\right\Vert _{F}}{\left\Vert x\right\Vert _{E}}\label{eq:Definition of operator norm}
\end{equation}
\end{defn}
\begin{prop}
$\mathcal{L}\left(E,F\right)$ is a non-archimedean Banach space over
$K$ with respect to the operator norm. 
\end{prop}
\begin{prop}
A linear operator $T:E\rightarrow F$ between two non-archimedean
Banach spaces $E$ and $F$ over the field $K$ is continuous if and
only if it has finite operator norm. Consequently, \textbf{bounded
linear operators }and \textbf{continuous linear operators }are one
and the same, just like in the classical case. 
\end{prop}
\begin{defn}
Given a non-archimedean $K$-Banach space $E$, the \textbf{dual}\index{Banach space!dual of}
of $E$ (denoted $E^{\prime}$) is, as usual, the $K$-Banach space
of all continuous linear functionals from $E$ to $K$. Given an $f\in E^{\prime}$,
we write $\left\Vert f\right\Vert _{E^{\prime}}$ to denote the norm
of $f$; this is defined by the usual variant of (\ref{eq:Definition of operator norm}):
\begin{equation}
\left\Vert f\right\Vert _{E^{\prime}}\overset{\textrm{def}}{=}\sup_{x\in E}\frac{\left|f\left(x\right)\right|_{K}}{\left\Vert x\right\Vert _{E}}\label{eq:Definition of the norm of a linear functional}
\end{equation}
\end{defn}
\begin{rem}
For the case of the space $C\left(\mathbb{Z}_{p},K\right)$, where
$K$ is a $q$-adic field, where $p$ and $q$ are distinct primes,
note that its dual $C\left(\mathbb{Z}_{p},K\right)^{\prime}$ \nomenclature{$C\left(\mathbb{Z}_{p},K\right)^{\prime}$}{set of continuous $K$-valued linear functionals on $C\left(\mathbb{Z}_{p},K\right)$}
will consist of all continuous $K$-valued linear functionals on $C\left(\mathbb{Z}_{p},K\right)$.
\end{rem}
\begin{fact}
\emph{Importantly, because of their proofs' dependence on the }\textbf{\emph{Baire
Category Theorem}}\emph{,\index{Uniform Boundedness Principle} the
}\textbf{\emph{Uniform Boundedness Principle}}\emph{, the }\textbf{\emph{Closed
Graph Theorem}}\emph{, the }\textbf{\emph{Open Mapping Theorem}}\emph{,
and the }\textbf{\emph{Banach-Steinhaus Theorem}}\emph{ }all apply
to non-archimedean Banach spaces\emph{ \cite{Schikhof Banach Space Paper}.
The Hahn-Banach Theorem, however, is more subtle, as is discussed
below.} 
\end{fact}
\vphantom{}

As we now move to discuss notions of duality, the notion of \index{spherically complete}\textbf{spherical
completeness}\footnote{Recall, a valued field $K$ (or, more generally, any metric space)
is said to be \textbf{spherically complete} whenever any sequence
of nested non-empty balls $B_{1}\supseteq B_{2}\supseteq\cdots$ in
$K$ has a non-empty intersection. $\mathbb{Q}_{p}$ and any finite
extension thereof are spherically complete; $\mathbb{C}_{p}$, however,
is not.} will raise its head in surprising ways. The first instance of this
is in the non-archimedean Hahn-Banach Theorem. First, however, another
definition, to deal with the fact that non-archimedean fields need
not be separable as topological spaces. 
\begin{defn}
A $K$-Banach space $E$ is said\index{Banach space!of countable type}
to be \textbf{of countable type }if there is a countable set whose
$K$-span is dense in $E$. 
\end{defn}
\begin{thm}[\textbf{Non-Archimedean Hahn-Banach Theorems}\footnote{Theorems 9 and 10 in \cite{Schikhof Banach Space Paper}.}\index{Hahn-Banach Theorem}]
Let $E$ be a Banach space over a complete non-archimedean valued
field $K$, let $D$ be a subspace of $E$. Let $D^{\prime}$ be the
continuous dual of $D$ (continuous linear functionals $D\rightarrow K$)
and let $f\in D^{\prime}$.

\vphantom{}

I. If $K$ is \emph{spherically complete} (for instance, if $K$ is
$\mathbb{Q}_{q}$ or a finite, metrically-complete extension thereof),
there then exists a continuous linear functional $\overline{f}\in E^{\prime}$
whose restriction to $D$ is equal to $f$, and so that $\left\Vert \overline{f}\right\Vert _{E^{\prime}}=\left\Vert f\right\Vert _{D^{\prime}}$.

\vphantom{}

II. If $E$ is of countable type, then\textemdash even if $K$ is
\textbf{\emph{not }}spherically complete\textemdash for any $\epsilon>0$,
there exists a continuous linear functional $\overline{f}:E\rightarrow K$
whose restriction to $D$ is equal to $f$. Moreover, $\left\Vert \overline{f}\right\Vert _{E^{\prime}}\leq\left(1+\epsilon\right)\left\Vert f\right\Vert _{D^{\prime}}$. 
\end{thm}
\begin{rem}
The requirement that $E$ be of countable type \emph{cannot} be dropped,
nor can $\epsilon$ ever be set to $0$. There exists a two-dimensional
non-archimedean Banach space $E$ over a non-spherically-complete
$K$ and a subspace $D\subseteq E$ and an $f\in D^{\prime}$ such
that no extension $\overline{f}\in E^{\prime}$ exists satisfying
$\left\Vert \overline{f}\right\Vert _{E^{\prime}}\leq\left\Vert f\right\Vert _{D^{\prime}}$
\cite{Schikhof Banach Space Paper}. This counterexample can be found
on page 68 of van Rooij's book \cite{van Rooij - Non-Archmedean Functional Analysis}. 
\end{rem}
\vphantom{}

Spherical completeness is responsible for non-archimedean Banach spaces'
most significant deviations from the classical theory, which we now
chronicle. The names I give these theorems were taken from \cite{van Rooij - Non-Archmedean Functional Analysis}. 
\begin{thm}[\textbf{Fleischer's Theorem}\footnote{See \cite{Schikhof Banach Space Paper,van Rooij - Non-Archmedean Functional Analysis}.
This theorem is given as \textbf{Theorem}}]
Let $K$ be \textbf{spherically complete}. Then, a $K$-Banach space
is reflexive\index{Banach space!reflexive} if and only if it is finite-dimensional. 
\end{thm}
\begin{thm}[\textbf{van der Put's Theorem}\footnote{See \cite{Schikhof Banach Space Paper,van Rooij - Non-Archmedean Functional Analysis}.
This theorem is a specific part (iii) of \textbf{Theorem 4.21 }on
page 118 of \cite{van Rooij - Non-Archmedean Functional Analysis}.}]
\label{thm:c_o and ell_infinit are each others duals in a spherically incomplete NA field}Let
$K$ be\textbf{\emph{ spherically incomplete}}. Then, for any countable\footnote{To give the reader a taste of the level of generality of van Rooij's
text, while I state the theorem for countable $X$, the formulation
given in \cite{van Rooij - Non-Archmedean Functional Analysis} states
that it holds for any ``small'' set $X$, a definition that involves
working with cardinal numbers and arbitrary $\sigma$-additive measures
on Boolean rings. On page 31 of the same, van Rooij says that ``obviously'',
a set with cardinality $\aleph_{0}$ (that is, a countable set) is
``small''. I choose to believe him, from whence I obtained the version
of the theorem given here.} set $X$, the $K$-Banach spaces $c_{0}\left(X,K\right)$ and $B\left(X,K\right)$
are reflexive\index{Banach space!reflexive}, being one another's
duals in a natural way. 
\end{thm}
\begin{cor}
If $K$ is \textbf{\emph{spherically incomplete}}, then \textbf{every}
$K$-Banach space of countable type is reflexive. \cite{van Rooij - Non-Archmedean Functional Analysis} 
\end{cor}

\subsection{\label{subsec:3.1.3 The-van-der}The van der Put Basis}

IN THIS SUBSECTION, $p$ IS A PRIME. $K$ IS A METRICALLY COMPLETE
VALUED FIELD (NOT NECESSARILY NON-ARCHIMEDEAN) OF CHARACTERISTIC ZERO.

\vphantom{}

A principal reason for the ``simplicity'' of $\left(p,q\right)$-adic
analysis in comparison with more general theories of non-archimedean
analysis is that, much like in the $\left(p,p\right)$-adic case,
the Banach space of continuous functions $\left(p,q\right)$-adic
functions admits a countable basis. This is \textbf{the van der Put
basis}. This basis actually works for the space of continuous functions
$f:\mathbb{Z}_{p}\rightarrow K$ where $K$ is \emph{any }metrically
complete non-archimedean field, unlike the Mahler basis, which only
works when $K$ is $\mathbb{Q}_{p}$ or an extension thereof. Because
of the van der Put basis, we will be able to explicitly compute integrals
and Fourier transforms, which will be of the utmost importance for
our analyses of $\chi_{H}$. 
\begin{defn}[\textbf{The van der Put Basis} \cite{Robert's Book,Ultrametric Calculus}]
\label{def:vdP basis, n_minus, vpD coefficients}\ 

\vphantom{}

I. We call the indicator functions $\left[\mathfrak{z}\overset{p^{\lambda_{p}\left(j\right)}}{\equiv}j\right]$
the \index{van der Put!basis}\textbf{van der Put (basis) functions
}and refer to the set $\left\{ \left[\mathfrak{z}\overset{p^{\lambda_{p}\left(j\right)}}{\equiv}j\right]\right\} _{j\in\mathbb{N}_{0}}$
as the \textbf{van der Put basis}.

\vphantom{}

II. Given an integer $n\geq1$, we can express $n$ $p$-adically
as: 
\begin{equation}
n=\sum_{k=0}^{\lambda_{p}\left(n\right)-1}n_{k}p^{k}
\end{equation}
where the $n_{k}$s are the $p$-adic digits of $n$. Following Robert
and Schikhof (\cite{Robert's Book,Ultrametric Calculus}), we write
\nomenclature{$n_{-}$}{ }$n_{-}$ to denote the integer obtained
by deleting the $k=\lambda_{p}\left(n\right)-1$ term from the above
series; equivalently, $n_{-}$ is the integer obtained by deleting
the right-most non-zero digit in $n$'s $p$-adic representation of
$n$. That is: 
\begin{equation}
n_{-}=n-n_{\lambda_{p}\left(n\right)-1}p^{\lambda_{p}\left(n\right)-1}\label{eq:Definition of n minus}
\end{equation}
We can also write this as the value of $n$ modulo $p^{\lambda_{p}\left(n\right)-1}$:
\begin{equation}
n_{-}=\left[n\right]_{p^{\lambda_{p}\left(n\right)-1}}
\end{equation}

\vphantom{}

III. Let $f:\mathbb{Z}_{p}\rightarrow K$ be any function. Then, for
all $n\in\mathbb{N}_{0}$, we define the constants $\left\{ c_{n}\left(f\right)\right\} _{n\geq0}\in K$
by: 
\begin{equation}
c_{n}\left(f\right)\overset{\textrm{def}}{=}\begin{cases}
f\left(0\right) & \textrm{if }n=0\\
f\left(n\right)-f\left(n_{-}\right) & \textrm{if }n\geq1
\end{cases},\textrm{ }\forall n\in\mathbb{N}_{0}\label{eq:Def of c_n of f}
\end{equation}
We call \nomenclature{$c_{n}\left(f\right)$}{$n$th van der Put coefficient of $f$}$c_{n}\left(f\right)$
the $n$th \index{van der Put!coefficients}\textbf{van der Put coefficient
}of $f$.

\vphantom{}

IV. We write \nomenclature{$\textrm{vdP}\left(\mathbb{Z}_{p},K\right)$}{the set of formal van der Put series with coefficients in $K$}$\textrm{vdP}\left(\mathbb{Z}_{p},K\right)$
to denote the $K$-vector space of all \textbf{formal van der Put
series}\index{van der Put!series}. The elements of $\textrm{vdP}\left(\mathbb{Z}_{p},K\right)$
are formal sums: 
\begin{equation}
\sum_{n=0}^{\infty}\mathfrak{a}_{n}\left[\mathfrak{z}\overset{p^{\lambda_{p}\left(n\right)}}{\equiv}n\right]
\end{equation}
where the $\mathfrak{a}_{n}$s are elements of $K$.

The \textbf{domain of convergence }of a formal van der Put series
is defined as the set of all $\mathfrak{z}\in\mathbb{Z}_{p}$ for
which the series converges in $K$. We call $\textrm{vdP}\left(\mathbb{Z}_{p},\mathbb{C}_{q}\right)$
the space of \textbf{formal $\left(p,q\right)$-adic van der Put series}.

\vphantom{}

V. Let \nomenclature{$F\left(\mathbb{Z}_{p},K\right)$}{$K$-valued functions on $\mathbb{Z}_{p}$}$F\left(\mathbb{Z}_{p},K\right)$
denote the $K$-linear space of all $K$-valued functions on $\mathbb{Z}_{p}$.
Then, we define the linear operator \nomenclature{$S_{p}$}{van-der-Put-series-creating operator}$S_{p}:F\left(\mathbb{Z}_{p},K\right)\rightarrow\textrm{vdP}\left(\mathbb{Z}_{p},K\right)$
by: 
\begin{equation}
S_{p}\left\{ f\right\} \left(\mathfrak{z}\right)\overset{\textrm{def}}{=}\sum_{n=0}^{\infty}c_{n}\left(f\right)\left[\mathfrak{z}\overset{p^{\lambda_{p}\left(n\right)}}{\equiv}n\right],\textrm{ }\forall f\in F\left(\mathbb{Z}_{p},K\right)\label{eq:Definition of S_p of f}
\end{equation}
We also define partial sum operators \nomenclature{$S_{p:N}$}{$N$th partial van-der-Put-series-creating operator}$S_{p:N}:F\left(\mathbb{Z}_{p},K\right)\rightarrow C\left(\mathbb{Z}_{p},K\right)$
by: 
\begin{equation}
S_{p:N}\left\{ f\right\} \left(\mathfrak{z}\right)\overset{\textrm{def}}{=}\sum_{n=0}^{p^{N}-1}c_{n}\left(f\right)\left[\mathfrak{z}\overset{p^{\lambda_{p}\left(n\right)}}{\equiv}n\right],\textrm{ }\forall f\in F\left(\mathbb{Z}_{p},K\right)\label{eq:Definition of S_p N of f}
\end{equation}
\end{defn}
\begin{rem}
Fixing $\mathfrak{z}\in\mathbb{N}_{0}$, observe that $\mathfrak{z}\overset{p^{\lambda_{p}\left(n\right)}}{\equiv}n$
can hold true for at most finitely many $n$. As such, every formal
van der Put series necessarily converges in $K$ at every $\mathfrak{z}\in\mathbb{N}_{0}$. 
\end{rem}
\begin{example}[\textbf{$\lambda$-Decomposition}]
Given a sum of the form with coefficients in $\mathbb{C}_{q}$: 
\begin{equation}
\sum_{n=0}^{\infty}\mathfrak{a}_{n}\left[\mathfrak{z}\overset{p^{\lambda_{p}\left(n\right)}}{\equiv}n\right]
\end{equation}
we can decompose it by splitting the sum into distinct pieces, with
$\lambda_{p}\left(n\right)$ taking a single value over each individual
piece. In particular, noting that for any $m\in\mathbb{N}_{1}$: 
\begin{equation}
\lambda_{p}\left(n\right)=m\Leftrightarrow p^{m-1}\leq n\leq p^{m}-1
\end{equation}
we can write: 
\begin{equation}
\sum_{n=0}^{\infty}\mathfrak{a}_{n}\left[\mathfrak{z}\overset{p^{\lambda_{p}\left(n\right)}}{\equiv}n\right]\overset{\mathbb{C}_{q}}{=}\mathfrak{a}_{0}+\sum_{m=1}^{\infty}\sum_{n=p^{m-1}}^{p^{m}-1}\mathfrak{a}_{n}\left[\mathfrak{z}\overset{p^{m}}{\equiv}n\right]\label{eq:Lambda Decomposition}
\end{equation}
We will use this decomposition frequently enough for it to be worth
having its own name: \emph{the Lambda decomposition}, or \emph{$\lambda$-decomposition\index{lambda-decomposition@$\lambda$-decomposition}},
for short. 
\end{example}
\vphantom{}

For us, the chief utility of van der Put series is the fact that they
tell us how to compute a type of function-limit which we have seen
to great effect in Chapter 2. The next proposition gives the details. 
\begin{prop}[\textbf{\textit{van der Put Identity}}]
\label{prop:vdP identity}Let $f\in B\left(\mathbb{Z}_{p},K\right)$.
Then, for any $\mathfrak{z}\in\mathbb{Z}_{p}$ which is in the domain
of convergence of $S_{p}\left\{ f\right\} $: 
\begin{equation}
S_{p}\left\{ f\right\} \left(\mathfrak{z}\right)\overset{K}{=}\lim_{k\rightarrow\infty}f\left(\left[\mathfrak{z}\right]_{p^{k}}\right)\label{eq:van der Put identity}
\end{equation}
We call this the \textbf{van der Put identity}\index{van der Put!identity}.
We also have: 
\begin{equation}
\sum_{n=0}^{p^{N}-1}c_{n}\left(f\right)\left[\mathfrak{z}\overset{p^{\lambda_{p}\left(n\right)}}{\equiv}n\right]\overset{\mathbb{F}}{=}f\left(\left[\mathfrak{z}\right]_{p^{N}}\right),\textrm{ }\forall\mathfrak{z}\in\mathbb{Z}_{p}\label{eq:truncated van der Put identity}
\end{equation}
which we call the (\textbf{$N$th}) \textbf{truncated van der Put
identity}. 
\end{prop}
\begin{rem}
This is part of Exercise 62.B on page 192 of \cite{Ultrametric Calculus}. 
\end{rem}
Proof: Performing a $\lambda$-decomposition yields:

\begin{equation}
S_{p}\left\{ f\right\} \left(\mathfrak{z}\right)=\sum_{n=0}^{\infty}c_{n}\left(f\right)\left[\mathfrak{z}\overset{p^{\lambda_{p}\left(n\right)}}{\equiv}n\right]\overset{K}{=}c_{0}\left(f\right)+\sum_{k=1}^{\infty}\sum_{n=p^{k-1}}^{p^{k}-1}c_{n}\left(f\right)\left[\mathfrak{z}\overset{p^{k}}{\equiv}n\right]
\end{equation}
Fixing $\mathfrak{z}$ and letting $k$ be arbitrary, we note that
there is at most one value of $n\in\left\{ p^{k-1},\ldots,p^{k}-1\right\} $
which solves the congruence $\mathfrak{z}\overset{p^{k}}{\equiv}n$;
that value is $n=\left[\mathfrak{z}\right]_{p^{k}}$. Thus: 
\begin{equation}
\sum_{n=p^{k-1}}^{p^{k}-1}c_{n}\left(f\right)\left[\mathfrak{z}\overset{p^{k}}{\equiv}n\right]\overset{K}{=}c_{\left[\mathfrak{z}\right]_{p^{k}}}\left(f\right)\left[\lambda_{p}\left(\left[\mathfrak{z}\right]_{p^{k}}\right)=k\right]\label{eq:Inner term of vdP lambda decomposition}
\end{equation}
where the Iverson bracket on the right indicates that the right-hand
side is $0$ whenever $\lambda_{p}\left(\left[\mathfrak{z}\right]_{p^{k}}\right)$
is not equal to $k$. This failure of equality occurs precisely when
there is no $n\in\left\{ p^{k-1},\ldots,p^{k}-1\right\} $ which solves
the congruence $\mathfrak{z}\overset{p^{k}}{\equiv}n$.

Next, writing $\mathfrak{z}$ as $\mathfrak{z}=\sum_{j=0}^{\infty}\mathfrak{z}_{j}p^{j}$,
where the $\mathfrak{z}_{j}$s are the $p$-adic digits of $\mathfrak{z}$,
we have that: 
\begin{equation}
c_{\left[\mathfrak{z}\right]_{p^{k}}}\left(f\right)\overset{K}{=}f\left(\left[\mathfrak{z}\right]_{p^{k}}\right)-f\left(\left[\mathfrak{z}\right]_{p^{k}}-\mathfrak{z}_{\lambda_{p}\left(\left[\mathfrak{z}\right]_{p^{k}}\right)-1}p^{\lambda_{p}\left(\left[\mathfrak{z}\right]_{p^{k}}\right)-1}\right)
\end{equation}
Note that $\left[\mathfrak{z}\right]_{p^{k}}$ can have \emph{at most}
$k$ $p$-adic digits. In particular, $\lambda_{p}\left(\left[\mathfrak{z}\right]_{p^{k}}\right)=j$
if and only if the $\mathfrak{z}_{j-1}$ is the right-most non-zero
$p$-adic digit of $\mathfrak{z}$. Consequently, letting $j\in\left\{ 0,\ldots,k\right\} $
denote the integer $\lambda_{p}\left(\left[\mathfrak{z}\right]_{p^{k}}\right)$,
we have that $\left[\mathfrak{z}\right]_{p^{k}}=\left[\mathfrak{z}\right]_{p^{j}}$.
We can then write: 
\begin{align*}
c_{\left[\mathfrak{z}\right]_{p^{k}}}\left(f\right) & \overset{K}{=}f\left(\left[\mathfrak{z}\right]_{p^{j}}\right)-f\left(\left[\mathfrak{z}\right]_{p^{j}}-\mathfrak{z}_{j-1}p^{j-1}\right)\\
 & =f\left(\left[\mathfrak{z}\right]_{p^{j}}\right)-f\left(\sum_{i=0}^{j-1}\mathfrak{z}_{i}p^{i}-\mathfrak{z}_{j-1}p^{j-1}\right)\\
 & =f\left(\left[\mathfrak{z}\right]_{p^{j}}\right)-f\left(\sum_{i=0}^{j-2}\mathfrak{z}_{i}p^{i}\right)\\
 & =f\left(\left[\mathfrak{z}\right]_{p^{j}}\right)-f\left(\left[\mathfrak{z}\right]_{p^{j-1}}\right)
\end{align*}
That is to say: 
\begin{equation}
c_{\left[\mathfrak{z}\right]_{p^{k}}}\left(\chi\right)\overset{K}{=}\chi\left(\left[\mathfrak{z}\right]_{p^{\lambda_{p}\left(\left[\mathfrak{z}\right]_{p^{k}}\right)}}\right)-\chi\left(\left[\mathfrak{z}\right]_{p^{\lambda_{p}\left(\left[\mathfrak{z}\right]_{p^{k}}\right)-1}}\right)
\end{equation}
for all $\mathfrak{z}\in\mathbb{Z}_{p}\backslash\left\{ 0\right\} $,
and all $k$ large enough so that $\lambda_{p}\left(\left[\mathfrak{z}\right]_{p^{k}}\right)\geq1$;
that is, all $k>v_{p}\left(\mathfrak{z}\right)$.

So, we have: 
\begin{equation}
c_{\left[\mathfrak{z}\right]_{p^{k}}}\left(f\right)\left[\lambda_{p}\left(\left[\mathfrak{z}\right]_{p^{k}}\right)=k\right]\overset{K}{=}\left(f\left(\left[\mathfrak{z}\right]_{p^{k}}\right)-f\left(\left[\mathfrak{z}\right]_{p^{k-1}}\right)\right)\left[\lambda_{p}\left(\left[\mathfrak{z}\right]_{p^{k}}\right)=k\right]\label{eq:Iverson Bracket Check for vdP identity}
\end{equation}

\begin{claim}
The Iverson bracket on the right-hand side of (\ref{eq:Iverson Bracket Check for vdP identity})
can be removed.

Proof of claim: If $\lambda_{p}\left(\left[\mathfrak{z}\right]_{p^{k}}\right)=k$,
then the Iverson bracket is $1$, and it disappears on its own without
causing us any trouble. So, suppose $\lambda_{p}\left(\left[\mathfrak{z}\right]_{p^{k}}\right)\neq k$.
Then the right-most digit of $\left[\mathfrak{z}\right]_{p^{k}}$\textemdash i.e.,
$\mathfrak{z}_{k-1}$\textemdash is $0$. But then, $\left[\mathfrak{z}\right]_{p^{k}}=\left[\mathfrak{z}\right]_{p^{k-1}}$,
which makes the right-hand side of (\ref{eq:Iverson Bracket Check for vdP identity})
anyway, because $f\left(\left[\mathfrak{z}\right]_{p^{k}}\right)-f\left(\left[\mathfrak{z}\right]_{p^{k-1}}\right)$
is then equal to $0$. So, \emph{any case} which would cause the Iverson
bracket to vanish causes $f\left(\left[\mathfrak{z}\right]_{p^{k}}\right)-f\left(\left[\mathfrak{z}\right]_{p^{k-1}}\right)$
to vanish. This tells us that the Iverson bracket in (\ref{eq:Iverson Bracket Check for vdP identity})
isn't actually needed. As such, we are justified in writing: 
\[
c_{\left[\mathfrak{z}\right]_{p^{k}}}\left(f\right)\left[\lambda_{p}\left(\left[\mathfrak{z}\right]_{p^{k}}\right)=k\right]\overset{K}{=}f\left(\left[\mathfrak{z}\right]_{p^{k}}\right)-f\left(\left[\mathfrak{z}\right]_{p^{k-1}}\right)
\]
This proves the claim. 
\end{claim}
\vphantom{}

That being done, we can write: 
\begin{align*}
\sum_{n=0}^{\infty}c_{n}\left(f\right)\left[\mathfrak{z}\overset{p^{\lambda_{p}\left(n\right)}}{\equiv}n\right] & \overset{K}{=}c_{0}\left(f\right)+\sum_{k=1}^{\infty}\sum_{n=p^{k-1}}^{p^{k}-1}c_{n}\left(f\right)\left[\mathfrak{z}\overset{p^{k}}{\equiv}n\right]\\
 & =c_{0}\left(f\right)+\sum_{k=1}^{\infty}c_{\left[\mathfrak{z}\right]_{p^{k}}}\left(f\right)\left[\lambda_{p}\left(\left[\mathfrak{z}\right]_{p^{k}}\right)=k\right]\\
 & =c_{0}\left(f\right)+\underbrace{\sum_{k=1}^{\infty}\left(f\left(\left[\mathfrak{z}\right]_{p^{k}}\right)-f\left(\left[\mathfrak{z}\right]_{p^{k-1}}\right)\right)}_{\textrm{telescoping}}\\
\left(c_{0}\left(f\right)=f\left(\left[\mathfrak{z}\right]_{p^{0}}\right)=f\left(0\right)\right); & \overset{K}{=}f\left(0\right)+\lim_{k\rightarrow\infty}f\left(\left[\mathfrak{z}\right]_{p^{k}}\right)-f\left(0\right)\\
 & \overset{K}{=}\lim_{k\rightarrow\infty}f\left(\left[\mathfrak{z}\right]_{p^{k}}\right)
\end{align*}
which is the van der Put identity. If we instead only sum $n$ from
$0$ up to $p^{N}-1$, the upper limit of the $k$-sum in the above
is changed from $\infty$ to $N$, which yields (\ref{eq:truncated van der Put identity}).

Q.E.D.

\vphantom{}

With the van der Put basis, we can give a satisfying characterization
of the $\left(p,K\right)$-adic continuity of a function in terms
of the decay of its van der Put coefficients and the uniform convergence
of its van der Put series. 
\begin{thm}[\textbf{van der Put Basis Theorem}\footnote{\cite{Robert's Book} gives this theorem on pages 182 and 183, however,
it is stated there only for the case where $K$ is an extension of
$\mathbb{Q}_{p}$.}]
\label{thm:vdP basis theorem}Let $K$ be non-archimedean, and let
$f:\mathbb{Z}_{p}\rightarrow K$ be any function. Then, the following
are equivalent:

\vphantom{}

I. $f$ is continuous.

\vphantom{}

II. $\lim_{n\rightarrow\infty}\left|c_{n}\left(f\right)\right|_{K}=0$.

\vphantom{}

III. $S_{p:N}\left\{ f\right\} $ converges uniformly to $f$ in $\left(p,K\right)$-adic
norm; that is: 
\begin{equation}
\lim_{N\rightarrow\infty}\sup_{\mathfrak{z}\in\mathbb{Z}_{p}}\left|f\left(\mathfrak{z}\right)-\sum_{n=0}^{p^{N}-1}c_{n}\left(f\right)\left[\mathfrak{z}\overset{p^{\lambda_{p}\left(n\right)}}{\equiv}n\right]\right|_{K}=0
\end{equation}
\end{thm}
Proof:

\textbullet{} ($\textrm{(I)}\Rightarrow\textrm{(II)}$) Suppose $f$
is continuous. Since $n_{-}=n-n_{\lambda_{p}\left(n\right)-1}p^{\lambda_{p}\left(n\right)-1}$
for all $n\geq1$, we have that: 
\begin{equation}
\left|n-n_{-}\right|_{p}=\left|n_{\lambda_{p}\left(n\right)-1}p^{\lambda_{p}\left(n\right)-1}\right|_{p}=\frac{1}{p^{\lambda_{p}\left(n\right)-1}}
\end{equation}
Since this tends to $0$ as $n\rightarrow\infty$, the continuity
of $f$ guarantees that $\left|c_{n}\left(f\right)\right|_{K}=\left|f\left(n\right)-f\left(n_{-}\right)\right|_{K}$
tends to $0$ as $n\rightarrow\infty$.

\vphantom{}

\textbullet{} ($\textrm{(II)}\Rightarrow\textrm{(III)}$) Suppose
$\lim_{n\rightarrow\infty}\left|c_{n}\left(f\right)\right|_{K}=0$.
Then, since $K$ is non-archimedean, the absolute value of the difference
between $S_{p}\left\{ f\right\} $ and its $N$th partial sum satisfies:
\begin{align*}
\left|S_{p}\left\{ f\right\} \left(\mathfrak{z}\right)-\sum_{n=0}^{p^{N}-1}c_{n}\left(f\right)\left[\mathfrak{z}\overset{p^{\lambda_{p}\left(n\right)}}{\equiv}n\right]\right|_{K} & =\left|\sum_{n=p^{N}}^{\infty}c_{n}\left(f\right)\left[\mathfrak{z}\overset{p^{\lambda_{p}\left(n\right)}}{\equiv}n\right]\right|_{K}\\
\left(\textrm{ultrametric ineq.}\right); & \leq\sup_{n\geq p^{N}}\left|c_{n}\left(f\right)\left[\mathfrak{z}\overset{p^{\lambda_{p}\left(n\right)}}{\equiv}n\right]\right|_{K}\\
\left(\left[\mathfrak{z}\overset{p^{\lambda_{p}\left(n\right)}}{\equiv}n\right]\in\left\{ 0,1\right\} ,\textrm{ }\forall n,\mathfrak{z}\right); & \leq\sup_{n\geq p^{N}}\left|c_{n}\left(f\right)\right|_{K}\\
 & \leq\sup_{n\geq N}\left|c_{n}\left(f\right)\right|_{K}
\end{align*}
Since $\left|c_{n}\left(f\right)\right|_{K}$ is given to tend to
$0$ as $N\rightarrow\infty$, we then have that $\sup_{n\geq N}\left|c_{n}\left(f\right)\right|_{K}\rightarrow0$
as $N\rightarrow\infty$, and hence, that $S_{p}\left\{ f\right\} $
converges uniformly over $\mathbb{Z}_{p}$. The van der Put identity
(\ref{eq:van der Put identity}) shows that $f$ is the point-wise
limit of $S_{p:N}\left\{ f\right\} $, which in turn shows that $S_{p}\left\{ f\right\} $
converges point-wise to $f$. Since $S_{p}\left\{ f\right\} $ converges
uniformly, this proves that $S_{p:N}\left\{ f\right\} $ converges
uniformly to $f$.

\vphantom{}

\textbullet{} ($\textrm{(III)}\Rightarrow\textrm{(I)}$). Suppose
$S_{p}\left\{ f\right\} $ converges uniformly to $f$. Since $S_{p}\left\{ f\right\} $
is the limit of a sequence of continuous functions ($\left\{ S_{p:N}\left\{ f\right\} \right\} _{N\geq0}$)
which converge uniformly, its limit is necessarily continuous. Thus,
$f$ is continuous.

Q.E.D.

\vphantom{}

As an aside, we note the following: 
\begin{prop}
\label{prop:Convergence of real-valued vdP series}Let $f\in C\left(\mathbb{Z}_{p},K\right)$.
Then, the van der Put series for the continuous real-valued function
$\mathfrak{z}\in\mathbb{Z}_{p}\mapsto\left|f\left(\mathfrak{z}\right)\right|_{q}\in\mathbb{R}$
is uniformly convergent. 
\end{prop}
Proof: Since the $q$-adic absolute value is a continuous function
from $K$ to $\mathbb{R}$, the continuity of $f:\mathbb{Z}_{p}\rightarrow K$
guarantees that $\left|f\left(\mathfrak{z}\right)\right|_{q}$ is
a continuous\textemdash in fact, \emph{uniformly }continuous\textemdash real-valued
function on $\mathbb{Z}_{p}$.

The van der Put coefficients of $\left|f\right|_{q}$ are then: 
\begin{equation}
c_{n}\left(\left|f\right|_{q}\right)=\begin{cases}
\left|f\left(0\right)\right|_{q} & \textrm{if }n=0\\
\left|f\left(n\right)\right|_{q}-\left|f\left(n_{-}\right)\right|_{q} & \textrm{if }n\geq1
\end{cases}
\end{equation}
As such, (\ref{eq:van der Put identity}) yields: 
\begin{equation}
\sum_{n=0}^{\infty}c_{n}\left(\left|f\right|_{q}\right)\left[\mathfrak{z}\overset{p^{\lambda_{p}\left(n\right)}}{\equiv}n\right]\overset{\mathbb{R}}{=}\lim_{n\rightarrow\infty}\left|f\left(\left[\mathfrak{z}\right]_{p^{n}}\right)\right|_{q}
\end{equation}
Since $\left|f\right|_{q}$ is continuous, we see that the van der
Put series for $\left|f\right|_{q}$ then converges point-wise to
$\left|f\right|_{q}$. Thus, we need only upgrade the convergence
from point-wise to uniform.

To do this, we show that the sequence $\left\{ \left|f\left(\left[\mathfrak{z}\right]_{p^{n}}\right)\right|_{q}\right\} _{n\geq1}$
is uniformly Cauchy. Let $m$ and $n$ be arbitrary positive integers.
Then, by the reverse triangle inequality for the $q$-adic absolute
value: 
\begin{equation}
\left|\left|f\left(\left[\mathfrak{z}\right]_{p^{m}}\right)\right|_{q}-\left|f\left(\left[\mathfrak{z}\right]_{p^{n}}\right)\right|_{q}\right|\leq\left|f\left(\left[\mathfrak{z}\right]_{p^{m}}\right)-f\left(\left[\mathfrak{z}\right]_{p^{n}}\right)\right|_{q}
\end{equation}
Since the van der Put series for $f$ converges uniformly in $K$
over $\mathbb{Z}_{p}$, the $f\left(\left[\mathfrak{z}\right]_{p^{n}}\right)$s
are therefore uniformly Cauchy, and hence, so too are the $\left|f\left(\left[\mathfrak{z}\right]_{p^{n}}\right)\right|_{q}$s.
This establishes the desired uniform convergence.

Q.E.D.

\vphantom{}

With the van der Put basis, we can completely characterize $C\left(\mathbb{Z}_{p},K\right)$
when $K$ is non-archimedean. 
\begin{thm}
\label{thm:C(Z_p,K) is iso to c_0 K}If $K$ is non-archimedean, then
the Banach space $C\left(\mathbb{Z}_{p},K\right)$ is isometrically
isomorphic to $c_{0}\left(K\right)$ (the space of sequences in $K$
that converge to $0$ in $K$). 
\end{thm}
Proof: Let $L:C\left(\mathbb{Z}_{p},K\right)\rightarrow c_{0}\left(K\right)$
be the map which sends every $f\in C\left(\mathbb{Z}_{p},K\right)$
to the sequence $\mathbf{c}\left(f\right)=\left(c_{0}\left(f\right),c_{1}\left(f\right),\ldots\right)$
whose entries are the van der Put coefficients of $f$. Since: 
\begin{equation}
c_{n}\left(\alpha f+\beta g\right)=\alpha f\left(n\right)+\beta g\left(n\right)-\left(\alpha f\left(n_{-}\right)+\beta g\left(n_{-}\right)\right)=\alpha c_{n}\left(f\right)+\beta c_{n}\left(g\right)
\end{equation}
for all $\alpha,\beta\in K$ and all $f,g\in C\left(\mathbb{Z}_{p},K\right)$,
we have that $L$ is linear. By the \textbf{van der Put Basis Theorem}
(\textbf{Theorem \ref{thm:vdP basis theorem}}), for any $\mathbf{c}=\left(c_{0},c_{1},\ldots\right)$
in $c_{0}\left(K\right)$, the van der Put series $\sum_{n=0}^{\infty}c_{n}\left[\mathfrak{z}\overset{p^{\lambda_{p}\left(n\right)}}{\equiv}n\right]$
is an element of $C\left(\mathbb{Z}_{p},K\right)$ which $L$ sends
to $\mathbf{c}$; thus, $L$ is surjective. All that remains is to
establish norm-preservation; injectivity will then follow automatically.

Since the norm on $c_{0}\left(K\right)$ is the supremum norm: 
\begin{equation}
\left\Vert \mathbf{c}\right\Vert _{K}=\sup_{n\geq0}\left|c_{n}\right|_{K},\textrm{ }\forall\mathbf{c}\in c_{0}\left(K\right)
\end{equation}
observe that: 
\begin{equation}
\left\Vert L\left\{ f\right\} \right\Vert _{p,K}=\left\Vert \mathbf{c}\left(f\right)\right\Vert _{K}=\sup_{n\geq0}\left|c_{n}\left(f\right)\right|_{K}\geq\underbrace{\left|\sum_{n=0}^{\infty}c_{n}\left(f\right)\left[\mathfrak{z}\overset{p^{\lambda_{p}\left(n\right)}}{\equiv}n\right]\right|_{K}}_{\left|f\left(\mathfrak{z}\right)\right|_{K}}
\end{equation}
for all $f\in C\left(\mathbb{Z}_{p},K\right)$ and all $\mathfrak{z}\in\mathbb{Z}_{p}$,
and hence: 
\[
\left\Vert L\left\{ f\right\} \right\Vert _{p,K}\geq\sup_{\mathfrak{z}\in\mathbb{Z}_{p}}\left|f\left(\mathfrak{z}\right)\right|_{K}=\left\Vert f\right\Vert _{p,K},\textrm{ }\forall f\in C\left(\mathbb{Z}_{p},K\right)
\]

Finally, suppose there is an $f\in C\left(\mathbb{Z}_{p},K\right)$
so that the lower bound is strict: 
\begin{equation}
\left\Vert L\left\{ f\right\} \right\Vert _{p,K}>\left\Vert f\right\Vert _{p,K}
\end{equation}
Since $\lim_{n\rightarrow\infty}\left|c_{n}\left(f\right)\right|_{K}=0$
the principles of ultrametric analysis (see page \pageref{fact:Principles of Ultrametric Analysis})
guarantee the existence of a $\nu\geq0$ so that $\sup_{n\geq0}\left|c_{n}\left(f\right)\right|_{K}=\left|c_{\nu}\left(f\right)\right|_{K}$,
and hence, so that: 
\begin{equation}
\left\Vert L\left\{ f\right\} \right\Vert _{p,K}=\left|c_{\nu}\left(f\right)\right|_{K}=\left|f\left(\nu\right)-f\left(\nu_{-}\right)\right|_{K}
\end{equation}
Thus, our assumption $\left\Vert L\left\{ f\right\} \right\Vert _{p,K}>\left\Vert f\right\Vert _{p,K}$
forces: 
\begin{equation}
\sup_{\mathfrak{z}\in\mathbb{Z}_{p}}\left|f\left(\mathfrak{z}\right)\right|_{K}<\left|f\left(\nu\right)-f\left(\nu_{-}\right)\right|_{K}
\end{equation}
In particular, we have that the right-hand side is strictly greater
than both $\left|f\left(\nu\right)\right|_{K}$ and $\left|f\left(\nu_{-}\right)\right|_{K}$;
that is: 
\begin{equation}
\max\left\{ \left|f\left(\nu\right)\right|_{K},\left|f\left(\nu_{-}\right)\right|_{K}\right\} <\left|f\left(\nu\right)-f\left(\nu_{-}\right)\right|_{K}
\end{equation}
However, the strong triangle inequality allows us to write: 
\begin{equation}
\max\left\{ \left|f\left(\nu\right)\right|_{K},\left|f\left(\nu_{-}\right)\right|_{K}\right\} <\left|f\left(\nu\right)-f\left(\nu_{-}\right)\right|_{K}\leq\max\left\{ \left|f\left(\nu\right)\right|_{K},\left|f\left(\nu_{-}\right)\right|_{K}\right\} 
\end{equation}
which is impossible.

Thus, no such $f$ exists, and it must be that $\left\Vert L\left\{ f\right\} \right\Vert _{p,K}=\left\Vert f\right\Vert _{p,K}$,
for all $f\in C\left(\mathbb{Z}_{p},K\right)$, and $L$ is therefore
an isometry.

Q.E.D. 
\begin{cor}
Every $f\in C\left(\mathbb{Z}_{p},K\right)$ is uniquely representable
as a van der Put series. 
\end{cor}
Proof: Because map $L$ from the proof of \ref{thm:C(Z_p,K) is iso to c_0 K}
is an isometry, it is necessarily injective.

Q.E.D. 
\begin{cor}
If $K$ is spherically incomplete, then $C\left(\mathbb{Z}_{p},K\right)$
is a reflexive Banach space, and its dual, is then isometrically isomorphic
to $\ell^{\infty}\left(K\right)$. 
\end{cor}
Proof: Since \textbf{Theorem \ref{thm:C(Z_p,K) is iso to c_0 K}}
shows that $C\left(\mathbb{Z}_{p},K\right)$ is isometrically isomorphic
to $c_{0}\left(K\right)$, by \textbf{van der Put's Theorem} (\textbf{Theorem
\ref{thm:c_o and ell_infinit are each others duals in a spherically incomplete NA field}},
page \pageref{thm:c_o and ell_infinit are each others duals in a spherically incomplete NA field}),
the spherical incompleteness of $K$ then guarantees that the Banach
space $C\left(\mathbb{Z}_{p},K\right)$ is reflexive, and that its
dual is isometrically isomorphic to $\ell^{\infty}\left(K\right)$.

Q.E.D. 
\begin{thm}
The span of the van der Put basis is dense in $B\left(\mathbb{Z}_{p},K\right)$,
the space of bounded functions $f:\mathbb{Z}_{p}\rightarrow K$. 
\end{thm}
Proof: The span (finite $K$-linear combinations) of the indicator
functions: 
\[
\left\{ \left[\mathfrak{z}\overset{p^{n}}{\equiv}k\right]:n\geq0,k\in\left\{ 0,\ldots,p^{n}-1\right\} \right\} 
\]
is precisely the set of locally constant functions $\mathbb{Z}_{p}\rightarrow\mathbb{C}_{q}$.
Since each of these indicator functions is continuous, it can be approximated
in supremum norm to arbitrary accuracy by the span of the van der
Put basis functions. Since the locally constant functions are dense
in $B\left(\mathbb{Z}_{p},K\right)$, we have that the $\left\{ \left[\mathfrak{z}\overset{p^{n}}{\equiv}k\right]:n\geq0,k\in\left\{ 0,\ldots,p^{n}-1\right\} \right\} $s
are dense in $B\left(\mathbb{Z}_{p},K\right)$. This proves the van
der Put basis is dense in $B\left(\mathbb{Z}_{p},K\right)$.

Q.E.D.

\subsection{\label{subsec:3.1.4. The--adic-Fourier}The $\left(p,q\right)$-adic
Fourier Transform\index{Fourier!transform}}
\begin{defn}
A \textbf{$\left(p,q\right)$-adic Fourier series}\index{Fourier!series}\textbf{
}is a sum of the form: 
\begin{equation}
f\left(\mathfrak{z}\right)\overset{\mathbb{C}_{q}}{=}\sum_{t\in\hat{\mathbb{Z}}_{p}}\hat{f}\left(t\right)e^{2\pi i\left\{ t\mathfrak{z}\right\} _{p}}\label{eq:General Form of a (p,q)-adic Fourier Series}
\end{equation}
where $\hat{f}$ is any function $\hat{\mathbb{Z}}_{p}\rightarrow\mathbb{C}_{q}$,
and where $\mathfrak{z}$ is a $p$-adic variable (in $\mathbb{Z}_{p}$,
or, more generally, $\mathbb{Q}_{p}$). We say the Fourier series
is\textbf{ convergent at $\mathfrak{z}_{0}\in\mathbb{Q}_{p}$ }whenever
the sum $f\left(\mathfrak{z}_{0}\right)$ converges in $\mathbb{C}_{q}$.
As per our abuse of notation, we interpret $e^{2\pi i\left\{ t\mathfrak{z}_{0}\right\} _{p}}$
as a function of $t\in\hat{\mathbb{Z}}_{p}$ which, for each $t$,
outputs a certain $p$-power root of unity in $\mathbb{C}_{q}$ in
accordance with a fixed embedding of $\overline{\mathbb{Q}}$ in both
$\mathbb{C}_{q}$ and $\mathbb{C}$. Additionally, we say the Fourier
series \textbf{converges} \textbf{(point-wise) }on a set $U\subseteq\mathbb{Q}_{p}$
if it converges at each point in $U$. 
\end{defn}
\vphantom{}

For us, the most important Fourier series is the one for indicator
functions of sets of the form $\mathfrak{a}+p^{n}\mathbb{Z}_{p}$,
where $\mathfrak{a}\in\mathbb{Z}_{p}$ and $n\in\mathbb{N}_{0}$ 
\begin{prop}
\label{prop:Indicator for p-adic coset Fourier identity}Fix $\mathfrak{a}\in\mathbb{Z}_{p}$
and $n\in\mathbb{Z}_{0}$. Then, Fourier series of the indicator function
$\left[\mathfrak{z}\overset{p^{n}}{\equiv}\mathfrak{a}\right]$ is:
\begin{equation}
\left[\mathfrak{z}\overset{p^{n}}{\equiv}\mathfrak{a}\right]\overset{\overline{\mathbb{Q}}}{=}\frac{1}{p^{n}}\sum_{k=0}^{p^{n}-1}e^{2\pi i\left\{ k\frac{\mathfrak{z}-\mathfrak{a}}{p^{n}}\right\} _{p}},\textrm{ }\forall\mathfrak{z}\in\mathbb{Z}_{p}\label{eq:Fourier Identity for indicator functions}
\end{equation}
\end{prop}
\begin{rem}
This can also be written as:
\begin{equation}
\left[\mathfrak{z}\overset{p^{n}}{\equiv}\mathfrak{a}\right]\overset{\overline{\mathbb{Q}}}{=}\frac{1}{p^{n}}\sum_{\left|t\right|_{p}\leq p^{n}}e^{2\pi i\left\{ t\left(\mathfrak{z}-\mathfrak{a}\right)\right\} _{p}},\textrm{ }\forall\mathfrak{z}\in\mathbb{Z}_{p}
\end{equation}
\end{rem}
\begin{rem}
As indicated by the $\overline{\mathbb{Q}}$ over the equality, the
above identity is valid in $\overline{\mathbb{Q}}$, and hence, in
any field extension of $\overline{\mathbb{Q}}$, such as $\mathbb{C}$,
or $\mathbb{C}_{q}$ for any prime $q$. 
\end{rem}
Proof: 
\[
\frac{1}{p^{n}}\sum_{k=0}^{p^{n}-1}e^{2\pi i\left\{ k\frac{\mathfrak{z}-\mathfrak{a}}{p^{n}}\right\} _{p}}=\frac{1}{p^{n}}\sum_{k=0}^{p^{n}-1}e^{2\pi ik\left[\mathfrak{z}-\mathfrak{a}\right]_{p^{n}}/p^{n}}=\begin{cases}
1 & \textrm{if }\left[\mathfrak{z}-\mathfrak{a}\right]_{p^{n}}=0\\
0 & \textrm{else}
\end{cases}
\]
The formula follows from the fact that $\left[\mathfrak{z}-\mathfrak{a}\right]_{p^{n}}=0$
is equivalent to $\mathfrak{z}\overset{p^{n}}{\equiv}\mathfrak{a}$.

Q.E.D. 
\begin{prop}
A $\left(p,q\right)$-adic Fourier series $\sum_{t\in\hat{\mathbb{Z}}_{p}}\hat{f}\left(t\right)e^{2\pi i\left\{ t\mathfrak{z}\right\} _{p}}$
converges in $\mathbb{C}_{q}$ uniformly with respect to $\mathfrak{z}\in\mathbb{Q}_{p}$
if and only if $\hat{f}\in c_{0}\left(\hat{\mathbb{Z}}_{p},\mathbb{C}_{q}\right)$;
that is, if and only if: 
\begin{equation}
\lim_{n\rightarrow\infty}\max_{\left|t\right|_{p}=p^{n}}\left|\hat{f}\left(t\right)\right|_{q}\overset{\mathbb{R}}{=}0
\end{equation}
\end{prop}
Proof: Consider an arbitrary $\left(p,q\right)$-adic Fourier series
$\sum_{t\in\hat{\mathbb{Z}}_{p}}\hat{f}\left(t\right)e^{2\pi i\left\{ t\mathfrak{z}\right\} _{p}}$.
Due to the ultrametric topology of $\mathbb{C}_{q}$, the series will
converge at any given $\mathfrak{z}\in\mathbb{Q}_{p}$ if and only
if: 
\begin{equation}
\lim_{n\rightarrow\infty}\max_{\left|t\right|_{p}=p^{n}}\left|\hat{f}\left(t\right)e^{2\pi i\left\{ t\mathfrak{z}\right\} _{p}}\right|_{q}\overset{\mathbb{R}}{=}0
\end{equation}
In our abuse of notation, we have that $e^{2\pi i\left\{ t\mathfrak{z}\right\} _{p}}$
is a $p$-power root of unity in $\mathbb{C}_{q}$ for all $t\in\hat{\mathbb{Z}}_{p}$
and all $\mathfrak{z}\in\mathbb{Q}_{p}$; consequently, $\left|e^{2\pi i\left\{ t\mathfrak{z}\right\} _{p}}\right|_{q}=1$
for all $t$ and $\mathfrak{z}$. Hence: 
\begin{equation}
\lim_{n\rightarrow\infty}\max_{\left|t\right|_{p}=p^{n}}\left|\hat{f}\left(t\right)e^{2\pi i\left\{ t\mathfrak{z}\right\} _{p}}\right|_{q}\overset{\mathbb{R}}{=}\lim_{n\rightarrow\infty}\max_{\left|t\right|_{p}=p^{n}}\left|\hat{f}\left(t\right)\right|_{q},\textrm{ }\forall\mathfrak{z}\in\mathbb{Q}_{p}
\end{equation}

Q.E.D.

\vphantom{}

Before we prove more on convergence, existence, or uniqueness, let
us give the formal definition for what the Fourier coefficients of
$f:\mathbb{Z}_{p}\rightarrow\mathbb{C}_{q}$ \emph{should} be, if
they happened to exist. 
\begin{defn}
\label{def:pq adic Fourier coefficients}For a function $f:\mathbb{Z}_{p}\rightarrow\mathbb{C}_{q}$,
for each $t\in\hat{\mathbb{Z}}_{p}$, we define the $t$th $\left(p,q\right)$\textbf{-adic
Fourier coefficient }of\index{Fourier!coefficients} $f$, denoted
$\hat{f}\left(t\right)$, by the rule: 
\begin{equation}
\hat{f}\left(t\right)\overset{\mathbb{C}_{q}}{=}\sum_{n=\frac{1}{p}\left|t\right|_{p}}^{\infty}\frac{c_{n}\left(f\right)}{p^{\lambda_{p}\left(n\right)}}e^{-2n\pi it}\label{eq:Definition of (p,q)-adic Fourier Coefficients}
\end{equation}
provided that the sum on the right is convergent in $\mathbb{C}_{q}$. 
\end{defn}
\vphantom{}

The existence criteria for $\hat{f}\left(t\right)$ are extremely
simple: 
\begin{lem}
Let $f:\mathbb{Z}_{p}\rightarrow\mathbb{C}_{q}$ be a function. Then,
the series defining $\hat{f}\left(t\right)$ will converge in $\mathbb{C}_{q}$
uniformly with respect to $t\in\hat{\mathbb{Z}}_{p}$ if and only
if: 
\begin{equation}
\lim_{n\rightarrow\infty}\left|\frac{c_{n}\left(f\right)}{p^{\lambda_{p}\left(n\right)}}\right|_{q}\overset{\mathbb{R}}{=}0
\end{equation}
In particular, the uniform convergence of $\hat{f}\left(t\right)$
for all values of $t$ will occur if and only if $f$ is continuous. 
\end{lem}
Proof: $\mathbb{C}_{q}$'s non-archimedean topology guarantees that:
\begin{equation}
\hat{f}\left(t\right)\overset{\mathbb{C}_{q}}{=}\sum_{n=\frac{1}{p}\left|t\right|_{p}}^{\infty}\frac{c_{n}\left(f\right)}{p^{\lambda_{p}\left(n\right)}}e^{-2n\pi it}
\end{equation}
converges if and only if: 
\begin{equation}
\lim_{n\rightarrow\infty}\left|\frac{c_{n}\left(f\right)}{p^{\lambda_{p}\left(n\right)}}e^{-2n\pi it}\right|_{q}\overset{\mathbb{R}}{=}\lim_{n\rightarrow\infty}\left|\frac{c_{n}\left(f\right)}{p^{\lambda_{p}\left(n\right)}}\right|_{q}\overset{\mathbb{R}}{=}0
\end{equation}

Q.E.D. 
\begin{defn}
We write \nomenclature{$c_{0}\left(\hat{\mathbb{Z}}_{p},\mathbb{C}_{q}\right)$}{set of functions $\hat{f}:\hat{\mathbb{Z}}_{p}\rightarrow\mathbb{C}_{q}$ for which $\lim_{n\rightarrow\infty}\max_{\left|t\right|_{p}=p^{n}}\left|\hat{f}\left(t\right)\right|_{q}\overset{\mathbb{R}}{=}0$}$c_{0}\left(\hat{\mathbb{Z}}_{p},\mathbb{C}_{q}\right)$
to denote the $\mathbb{C}_{q}$-vector space of all functions $\hat{f}:\hat{\mathbb{Z}}_{p}\rightarrow\mathbb{C}_{q}$
for which: 
\begin{equation}
\lim_{n\rightarrow\infty}\max_{\left|t\right|_{p}=p^{n}}\left|\hat{f}\left(t\right)\right|_{q}\overset{\mathbb{R}}{=}0
\end{equation}
We make $c_{0}\left(\hat{\mathbb{Z}}_{p},\mathbb{C}_{q}\right)$ into
a non-archimedean Banach space by equipping it with the norm: 
\begin{equation}
\left\Vert \hat{f}\right\Vert _{p,q}\overset{\textrm{def}}{=}\sup_{t\in\hat{\mathbb{Z}}_{p}}\left|\hat{f}\left(t\right)\right|_{q}
\end{equation}
\end{defn}
\begin{thm}
\label{thm:formula for Fourier series}Let $f:\mathbb{Z}_{p}\rightarrow\mathbb{C}_{q}$
be a function. Then:

\vphantom{}

I. $f$ admits a uniformly convergent $\left(p,q\right)$-adic Fourier
series if and only if $\hat{f}\in c_{0}\left(\hat{\mathbb{Z}}_{p},\mathbb{C}_{q}\right)$.
In particular, $f$ is continuous if and only if $\hat{f}\in c_{0}\left(\hat{\mathbb{Z}}_{p},\mathbb{C}_{q}\right)$.

\vphantom{}

II. If $f$ admits a uniformly convergent $\left(p,q\right)$-adic
Fourier series, the series is unique, and is given by: 
\begin{equation}
f\left(\mathfrak{z}\right)\overset{\mathbb{C}_{q}}{=}\sum_{t\in\hat{\mathbb{Z}}_{p}}\hat{f}\left(t\right)e^{2\pi i\left\{ t\mathfrak{z}\right\} _{p}},\textrm{ }\forall\mathfrak{z}\in\mathbb{Z}_{p}
\end{equation}
\end{thm}
Proof: Let $f:\mathbb{Z}_{p}\rightarrow\mathbb{C}_{q}$ be a function.
For the moment, let us proceed formally. Using (\ref{eq:Fourier Identity for indicator functions}),
we can re-write the crossed van der Put series of $f:\mathbb{Z}_{p}\rightarrow\mathbb{C}_{q}$
as a Fourier series: 
\begin{align*}
S_{p}\left\{ f\right\} \left(\mathfrak{z}\right) & \overset{\mathbb{C}_{q}}{=}c_{0}\left(f\right)+\sum_{n=1}^{\infty}c_{n}\left(f\right)\left[\mathfrak{z}\overset{p^{\lambda_{p}\left(n\right)}}{\equiv}n\right]\\
 & =c_{0}\left(f\right)+\sum_{n=1}^{\infty}c_{n}\left(f\right)\frac{1}{p^{\lambda_{p}\left(n\right)}}\sum_{k=0}^{p^{\lambda_{p}\left(n\right)}-1}\exp\left(2k\pi i\left\{ \frac{\mathfrak{z}-n}{p^{\lambda_{p}\left(n\right)}}\right\} _{p}\right)\\
 & =c_{0}\left(f\right)+\sum_{n=1}^{\infty}\frac{c_{n}\left(f\right)}{p^{\lambda_{p}\left(n\right)}}\sum_{k=0}^{p^{\lambda_{p}\left(n\right)}-1}\exp\left(-2n\pi i\frac{k}{p^{\lambda_{p}\left(n\right)}}\right)\exp\left(2\pi i\left\{ \frac{k\mathfrak{z}}{p^{\lambda_{p}\left(n\right)}}\right\} _{p}\right)
\end{align*}
Now, observe that for any function $g$: 
\begin{align}
\sum_{k=0}^{p^{m}-1}g\left(\frac{k}{p^{m}}\right) & =\sum_{k=0}^{p^{m-1}-1}g\left(\frac{k}{p^{m-1}}\right)+\sum_{k\in\left(\mathbb{Z}/p^{m}\mathbb{Z}\right)^{\times}}g\left(\frac{k}{p^{m}}\right)\label{eq:Prufer-group sum decomposition - inductive step}
\end{align}
As such, by induction: 
\begin{align}
\sum_{k=0}^{p^{m}-1}g\left(\frac{k}{p^{m}}\right) & =g\left(0\right)+\sum_{j=1}^{m}\sum_{k\in\left(\mathbb{Z}/p^{j}\mathbb{Z}\right)^{\times}}g\left(\frac{k}{p^{j}}\right),\textrm{ }\forall m\in\mathbb{N}_{1}\label{eq:Prufer-group sum decomposition}
\end{align}
The sum on the right is taken over all irreducible non-integer fractions
whose denominator is a divisor of $p^{m}$. Adding $0$ into the mix,
we see that the sum on the right is then exactly the same as summing
$g\left(t\right)$ over every element of $\hat{\mathbb{Z}}_{p}$ whose
$p$-adic magnitude is $\leq p^{m}$. Hence: 
\begin{equation}
\sum_{k=0}^{p^{m}-1}g\left(\frac{k}{p^{m}}\right)=\sum_{\left|t\right|_{p}\leq p^{m}}g\left(t\right)\label{eq:Prufer-group sum decomposition - Re-indexing}
\end{equation}
and so: 
\[
\sum_{k=0}^{p^{\lambda_{p}\left(n\right)}-1}\exp\left(-2n\pi i\frac{k}{p^{\lambda_{p}\left(n\right)}}\right)\exp\left(2\pi i\left\{ \frac{k\mathfrak{z}}{p^{\lambda_{p}\left(n\right)}}\right\} _{p}\right)=\sum_{\left|t\right|_{p}\leq p^{\lambda_{p}\left(n\right)}}e^{-2n\pi it}e^{2\pi i\left\{ t\mathfrak{z}\right\} _{p}}
\]
and so: 
\begin{align*}
S_{p}\left\{ f\right\} \left(\mathfrak{z}\right) & \overset{\mathbb{C}_{q}}{=}c_{0}\left(f\right)+\sum_{n=1}^{\infty}\frac{c_{n}\left(f\right)}{p^{\lambda_{p}\left(n\right)}}\sum_{\left|t\right|_{p}\leq p^{\lambda_{p}\left(n\right)}}e^{-2n\pi it}e^{2\pi i\left\{ t\mathfrak{z}\right\} _{p}}\\
 & =c_{0}\left(f\right)+\sum_{n=1}^{\infty}\frac{c_{n}\left(f\right)}{p^{\lambda_{p}\left(n\right)}}\left(1+\sum_{0<\left|t\right|_{p}\leq p^{\lambda_{p}\left(n\right)}}e^{-2n\pi it}e^{2\pi i\left\{ t\mathfrak{z}\right\} _{p}}\right)\\
 & =\underbrace{c_{0}\left(f\right)+\sum_{n=1}^{\infty}\frac{c_{n}\left(f\right)}{p^{\lambda_{p}\left(n\right)}}}_{\sum_{n=0}^{\infty}\frac{c_{n}\left(f\right)}{p^{\lambda_{p}\left(n\right)}}}+\sum_{n=1}^{\infty}\frac{c_{n}\left(f\right)}{p^{\lambda_{p}\left(n\right)}}\sum_{0<\left|t\right|_{p}\leq p^{\lambda_{p}\left(n\right)}}e^{-2n\pi it}e^{2\pi i\left\{ t\mathfrak{z}\right\} _{p}}
\end{align*}

Next, observe that a given $t\in\hat{\mathbb{Z}}_{p}\backslash\left\{ 0\right\} $
will occur in the innermost summand if and only if: 
\begin{align*}
p^{\lambda_{p}\left(n\right)} & \geq\left|t\right|_{p}\\
 & \Updownarrow\\
n & \geq p^{-\nu_{p}\left(t\right)-1}=\frac{1}{p}\left|t\right|_{p}
\end{align*}
So, re-indexing in terms of $t$, we have that:

\begin{align*}
S_{p}\left\{ f\right\} \left(\mathfrak{z}\right) & \overset{\mathbb{C}_{q}}{=}\sum_{n=0}^{\infty}\frac{c_{n}\left(f\right)}{p^{\lambda_{p}\left(n\right)}}+\sum_{t\in\hat{\mathbb{Z}}_{p}\backslash\left\{ 0\right\} }\sum_{n=\frac{1}{p}\left|t\right|_{p}}^{\infty}\frac{c_{n}\left(f\right)}{p^{\lambda_{p}\left(n\right)}}e^{-2n\pi it}e^{2\pi i\left\{ t\mathfrak{z}\right\} _{p}}\\
 & =\sum_{t\in\hat{\mathbb{Z}}_{p}}\underbrace{\sum_{n=\frac{1}{p}\left|t\right|_{p}}^{\infty}\frac{c_{n}\left(f\right)}{p^{\lambda_{p}\left(n\right)}}e^{-2n\pi it}}_{\hat{f}\left(t\right)}e^{2\pi i\left\{ t\mathfrak{z}\right\} _{p}}\\
 & =\sum_{t\in\hat{\mathbb{Z}}_{p}}\hat{f}\left(t\right)e^{2\pi i\left\{ t\mathfrak{z}\right\} _{p}}
\end{align*}

As such, if $\hat{f}\in c_{0}\left(\hat{\mathbb{Z}}_{p},\mathbb{C}_{q}\right)$,
then: 
\[
S_{p}\left\{ f\right\} \left(\mathfrak{z}\right)=\sum_{t\in\hat{\mathbb{Z}}_{p}}\hat{f}\left(t\right)e^{2\pi i\left\{ t\mathfrak{z}\right\} _{p}}
\]
converges in $\mathbb{C}_{q}$ uniformly over $\mathbb{Z}_{p}$ and
is equal to $f$, which is necessarily uniformly continuous. Conversely,
if $f$ is continuous (which then makes $f$ uniformly continuous,
since its domain is compact), then all of the above formal computations
count as grouping and re-ordering of the uniformly convergent series:
\begin{equation}
f\left(\mathfrak{z}\right)\overset{\mathbb{C}_{q}}{=}\sum_{n=0}^{\infty}c_{n}\left(f\right)\left[\mathfrak{z}\overset{p^{\lambda_{p}\left(n\right)}}{\equiv}n\right]
\end{equation}
This proves that: 
\begin{equation}
S_{p}\left\{ f\right\} \left(\mathfrak{z}\right)=\sum_{t\in\hat{\mathbb{Z}}_{p}}\hat{f}\left(t\right)e^{2\pi i\left\{ t\mathfrak{z}\right\} _{p}}
\end{equation}
converges in $\mathbb{C}_{q}$ uniformly over $\mathbb{Z}_{p}$, which
then forces $\hat{f}\in c_{0}\left(\hat{\mathbb{Z}}_{p},\mathbb{C}_{q}\right)$.
Uniqueness follows from the uniqueness of the van der Put coefficients
of a continuous $\left(p,q\right)$-adic function.

Q.E.D. 
\begin{cor}
\label{cor:pq adic Fourier transform is an isometric isomorphism}The
$\left(p,q\right)$-adic Fourier transform: 
\begin{equation}
f\in C\left(\mathbb{Z}_{p},\mathbb{C}_{q}\right)\mapsto\hat{f}\in c_{0}\left(\hat{\mathbb{Z}}_{p},\mathbb{C}_{q}\right)
\end{equation}
is an isometric isomorphism of the Banach spaces $C\left(\mathbb{Z}_{p},\mathbb{C}_{q}\right)$
and $c_{0}\left(\hat{\mathbb{Z}}_{p},\mathbb{C}_{q}\right)$. 
\end{cor}
Proof: \textbf{Theorem \ref{thm:formula for Fourier series} }shows
that the Fourier transform is a bijection between $C\left(\mathbb{Z}_{p},\mathbb{C}_{q}\right)$
and $c_{0}\left(\hat{\mathbb{Z}}_{p},\mathbb{C}_{q}\right)$. The
linearity of the Fourier transform then shows that it is an isomorphism
of vector space (a.k.a., linear spaces). As for continuity, let $\hat{f}\in c_{0}\left(\hat{\mathbb{Z}}_{p},\mathbb{C}_{q}\right)$
be arbitrary. Then, we can write: 
\begin{align*}
\left\Vert \hat{f}\right\Vert _{p,q} & =\sup_{t\in\hat{\mathbb{Z}}_{p}}\left|\sum_{n=\frac{1}{p}\left|t\right|_{p}}^{\infty}\frac{c_{n}\left(f\right)}{p^{\lambda_{p}\left(n\right)}}e^{-2n\pi it}\right|_{q}\\
 & \leq\sup_{t\in\hat{\mathbb{Z}}_{p}}\sup_{n\geq\frac{1}{p}\left|t\right|_{p}}\left|\frac{c_{n}\left(f\right)}{p^{\lambda_{p}\left(n\right)}}e^{-2n\pi it}\right|_{q}\\
 & =\sup_{t\in\hat{\mathbb{Z}}_{p}}\sup_{n\geq\frac{1}{p}\left|t\right|_{p}}\left|c_{n}\left(f\right)\right|_{q}\\
 & =\sup_{n\geq0}\left|c_{n}\left(f\right)\right|_{q}
\end{align*}
For $n\geq1$, $c_{n}\left(f\right)=f\left(n\right)-f\left(n_{-}\right)$.
As such, the ultrametric inequality implies that: 
\[
\left\Vert \hat{f}\right\Vert _{p,q}\leq\sup_{n\geq0}\left|c_{n}\left(f\right)\right|_{q}\leq\sup_{n\geq0}\left|f\left(n\right)\right|_{q}\leq\sup_{\mathfrak{z}\in\mathbb{Z}_{p}}\left|f\left(\mathfrak{z}\right)\right|_{q}=\left\Vert f\right\Vert _{p,q}
\]
On the other hand, if $f\in C\left(\mathbb{Z}_{p},\mathbb{C}_{q}\right)$,
then: 
\begin{align*}
\left\Vert f\right\Vert _{p,q} & =\sup_{\mathfrak{z}\in\mathbb{Z}_{p}}\left|\sum_{t\in\hat{\mathbb{Z}}_{p}}\hat{f}\left(t\right)e^{2\pi it\left\{ t\mathfrak{z}\right\} _{p}}\right|_{q}\\
 & \leq\sup_{\mathfrak{z}\in\mathbb{Z}_{p}}\sup_{t\in\hat{\mathbb{Z}}_{p}}\left|\hat{f}\left(t\right)\right|_{q}\\
 & \leq\sup_{t\in\hat{\mathbb{Z}}_{p}}\left|\hat{f}\left(t\right)\right|_{q}\\
 & =\left\Vert \hat{f}\right\Vert _{p,q}
\end{align*}
Since the inequality holds in both directions, this forces $\left\Vert f\right\Vert _{p,q}=\left\Vert \hat{f}\right\Vert _{p,q}$.
Thus, the $\left(p,q\right)$-adic Fourier transform is an isometric
isomorphism, as desired.

Q.E.D.

\subsection{\label{subsec:3.1.5-adic-Integration-=00003D000026}$\left(p,q\right)$-adic
Integration and the Fourier-Stieltjes Transform}

Because we are focusing on the concrete case of functions on $\mathbb{Z}_{p}$
taking values in a non-archimedean field $K$ with residue field of
characteristic $q\neq p$, the van der Put basis for $C\left(\mathbb{Z}_{p},K\right)$
allows us to explicitly define and compute integrals and measures
in terms of their effects on this basis, and\textemdash equivalently\textemdash in
terms of the $\left(p,q\right)$-adic characters $e^{2\pi i\left\{ t\mathfrak{z}\right\} _{p}}$. 
\begin{defn}
A\textbf{ $\left(p,q\right)$-adic measure}\index{measure!left(p,qright)-adic@$\left(p,q\right)$-adic}\index{measure}\index{$p,q$-adic!measure}
$d\mu$ is a continuous $K$-valued linear functional $C\left(\mathbb{Z}_{p},K\right)$;
that is, $d\mu\in C\left(\mathbb{Z}_{p},K\right)^{\prime}$. Given
a function $f\in C\left(\mathbb{Z}_{p},K\right)$, we write: 
\begin{equation}
\int_{\mathbb{Z}_{p}}f\left(\mathfrak{z}\right)d\mu\left(\mathfrak{z}\right)
\end{equation}
to denote the image of $f$ under the linear functional $d\mu$. We
call this the \textbf{integral }of $f$ with respect to $d\mu$.

\vphantom{}

Like one would expect, we can multiply measures by continuous functions.
Given $g\in C\left(\mathbb{Z}_{p},K\right)$ and a measure $d\mu$,
we define the measure $d\nu=gd\mu$ by the rule: 
\begin{equation}
\nu\left(f\right)\overset{K}{=}\int_{\mathbb{Z}_{p}}f\left(\mathfrak{z}\right)d\nu\left(\mathfrak{z}\right)\overset{\textrm{def}}{=}\int_{\mathbb{Z}_{p}}f\left(\mathfrak{z}\right)g\left(\mathfrak{z}\right)d\mu\left(\mathfrak{z}\right)\overset{K}{=}\mu\left(f\times g\right)
\end{equation}
This definition makes sense since the continuity of $f$ and $g$
guarantees the continuity of their product.

\vphantom{}

Finally, we say a measure $d\mu$ is\index{translation invariance}
\textbf{translation invariant }whenever: 
\begin{equation}
\int_{\mathbb{Z}_{p}}f\left(\mathfrak{z}+\mathfrak{a}\right)d\mu\left(\mathfrak{z}\right)\overset{K}{=}\int_{\mathbb{Z}_{p}}f\left(\mathfrak{z}\right)d\mu\left(\mathfrak{z}\right),\textrm{ }\forall f\in C\left(\mathbb{Z}_{p},K\right),\textrm{ }\forall\mathfrak{a}\in\mathbb{Z}_{p}
\end{equation}
Note that the zero-measure (which sends every function to $0$) is
translation-invariant, so this notion is not vacuous. 
\end{defn}
\begin{prop}[\textbf{Integral-Series Interchange}]
\label{prop:pq-adic integral interchange rule}Let $d\mu$ be a $\left(p,q\right)$-adic
measure, and let $\left\{ F_{N}\right\} _{N\geq1}$ be a convergent
sequence in $C\left(\mathbb{Z}_{p},\mathbb{C}_{q}\right)$. Then:
\begin{equation}
\lim_{N\rightarrow\infty}\int_{\mathbb{Z}_{p}}F_{N}\left(\mathfrak{z}\right)d\mu\left(\mathfrak{z}\right)\overset{\mathbb{C}_{q}}{=}\int_{\mathbb{Z}_{p}}\left(\lim_{N\rightarrow\infty}F_{N}\left(\mathfrak{z}\right)\right)d\mu\left(\mathfrak{z}\right)\label{eq:pq adic integral interchange rule}
\end{equation}
Equivalently, limits and integrals can be interchanged whenever the
limit converges uniformly. 
\end{prop}
Proof: This follows immediately from the definition of $\left(p,q\right)$-adic
measures as continuous linear functionals on $C\left(\mathbb{Z}_{p},\mathbb{C}_{q}\right)$.

Q.E.D. 
\begin{prop}[\textbf{Integration term-by-term}]
\label{prop:term-by-term evaluation of integrals}Let $d\mu$ be
a $\left(p,q\right)$-adic measure, and let $f\in C\left(\mathbb{Z}_{p},K\right)$.

I. The integral of $f$ can be evaluated term-by term using its van
der Put series \index{van der Put!series}: 
\begin{equation}
\int_{\mathbb{Z}_{p}}f\left(\mathfrak{z}\right)d\mu\left(\mathfrak{z}\right)\overset{K}{=}\sum_{n=0}^{\infty}c_{n}\left(f\right)\int_{\mathbb{Z}_{p}}\left[\mathfrak{z}\overset{p^{\lambda_{p}\left(n\right)}}{\equiv}n\right]d\mu\left(\mathfrak{z}\right)\label{eq:Integration in terms of vdP basis}
\end{equation}

II. The integral of $f$ can be evaluated term-by term using its Fourier
series: 
\begin{equation}
\int_{\mathbb{Z}_{p}}f\left(\mathfrak{z}\right)d\mu\left(\mathfrak{z}\right)\overset{K}{=}\sum_{t\in\hat{\mathbb{Z}}_{p}}\hat{f}\left(t\right)\int_{\mathbb{Z}_{p}}e^{2\pi i\left\{ t\mathfrak{z}\right\} _{p}}d\mu\left(\mathfrak{z}\right),\forall f\in C\left(\mathbb{Z}_{p},K\right)\label{eq:Integration in terms of Fourier basis}
\end{equation}
\end{prop}
Proof: Let $f$ and $d\mu$ be as given. By \textbf{Theorems \ref{thm:vdP basis theorem}}
and \textbf{\ref{thm:formula for Fourier series}}, the continuity
of $f$ guarantees the uniform convergence of the van der Put and
Fourier series of $f$, respectively. \textbf{Proposition \ref{prop:pq-adic integral interchange rule}}
then justifies the interchanges of limits and integrals in (I) and
(II).

Q.E.D.

\vphantom{}

Because of the van der Put basis, every $\left(p,q\right)$-adic measure
is uniquely determined by what it does to the $n$th van der Put basis
function. On the other hand, since every $f\in C\left(\mathbb{Z}_{p},K\right)$
can also be represented by a uniformly continuous $\left(p,q\right)$-adic
Fourier series, a measure is therefore uniquely determined by what
it does to the $\left(p,q\right)$-adic characters $e^{-2\pi i\left\{ t\mathfrak{z}\right\} _{p}}$.
This leads to Fourier-Stieltjes transform of a $\left(p,q\right)$-adic
measure. 
\begin{defn}
Let $d\mu$ be a $\left(p,q\right)$-adic measure. Then, the \textbf{Fourier-Stieltjes
transform}\index{Fourier-Stieltjes transform} of $d\mu$, denoted
$\hat{\mu}$, is the function $\hat{\mu}:\hat{\mathbb{Z}}_{p}\rightarrow\mathbb{C}_{q}$
defined by: 
\begin{equation}
\hat{\mu}\left(t\right)\overset{\mathbb{C}_{q}}{=}\int_{\mathbb{Z}_{p}}e^{-2\pi i\left\{ t\mathfrak{z}\right\} _{p}}d\mu\left(\mathfrak{z}\right),\textrm{ }\forall t\in\hat{\mathbb{Z}}_{p}\label{eq:Fourier-Stieltjes transform of a measure}
\end{equation}
where, for each $t$, $\int_{\mathbb{Z}_{p}}e^{-2\pi i\left\{ t\mathfrak{z}\right\} _{p}}d\mu\left(\mathfrak{z}\right)$
denotes the image of the function $e^{-2\pi i\left\{ t\mathfrak{z}\right\} _{p}}=e^{2\pi i\left\{ \left(-t\right)\mathfrak{z}\right\} _{p}}$
under $d\mu$. 
\end{defn}
\begin{rem}
\label{rem:parseval-plancherel formula for integration against measures}In
this notation, (\ref{eq:Integration in terms of Fourier basis}) can
be written as: 
\begin{equation}
\int_{\mathbb{Z}_{p}}f\left(\mathfrak{z}\right)d\mu\left(\mathfrak{z}\right)\overset{K}{=}\sum_{t\in\hat{\mathbb{Z}}_{p}}\hat{f}\left(t\right)\hat{\mu}\left(-t\right)\overset{K}{=}\sum_{t\in\hat{\mathbb{Z}}_{p}}\hat{f}\left(-t\right)\hat{\mu}\left(t\right)\label{eq:Integration in terms of Fourier-Stieltjes coefficients}
\end{equation}
This formula will be \emph{vital }for studying $\chi_{H}$. 
\end{rem}
\begin{thm}
\label{thm:FST is an iso from measures to ell infinity}The Fourier-Stieltjes
transform defines an isometric isomorphism from the Banach space $C\left(\mathbb{Z}_{p},\mathbb{C}_{q}\right)^{\prime}$
onto the Banach space $B\left(\hat{\mathbb{Z}}_{p},\mathbb{C}_{q}\right)$
of bounded $\mathbb{C}_{q}$-valued functions on $\hat{\mathbb{Z}}_{p}$. 
\end{thm}
\begin{rem}
The more general version of this result is so significant to non-archimedean
functional analysis that van Rooij put it on the cover of his book
\cite{van Rooij - Non-Archmedean Functional Analysis}! 
\end{rem}
Proof: By \textbf{Corollary \ref{cor:pq adic Fourier transform is an isometric isomorphism}},
Fourier transform is an isometric isomorphism of the Banach spaces
$C\left(\mathbb{Z}_{p},\mathbb{C}_{q}\right)$ and $c_{0}\left(\hat{\mathbb{Z}}_{p},\mathbb{C}_{q}\right)$.
By \textbf{Theorem \ref{thm:c_o and ell_infinit are each others duals in a spherically incomplete NA field}},
since $\hat{\mathbb{Z}}_{p}$ is countable and $\mathbb{C}_{q}$ is
spherically incomplete, $c_{0}\left(\hat{\mathbb{Z}}_{p},\mathbb{C}_{q}\right)$'s
dual is then $B\left(\hat{\mathbb{Z}}_{p},\mathbb{C}_{q}\right)$.
The isometric isomorphism the Fourier transform between the base spaces
$C\left(\mathbb{Z}_{p},\mathbb{C}_{q}\right)$ and $c_{0}\left(\hat{\mathbb{Z}}_{p},\mathbb{C}_{q}\right)$
then indices an isometric isomorphism between their duals, as desired.

Q.E.D. 
\begin{prop}
\label{prop:Fourier transform of haar measure}Let $d\mu,d\nu\in C\left(\mathbb{Z}_{p},K\right)^{\prime}$
be the measures defined by: 
\begin{equation}
\int_{\mathbb{Z}_{p}}\left[\mathfrak{z}\overset{p^{n}}{\equiv}k\right]d\mu\left(\mathfrak{z}\right)=\frac{1}{p^{n}},\textrm{ }\forall n\in\mathbb{N}_{0},\textrm{ }\forall k\in\mathbb{Z}
\end{equation}
and: 
\begin{equation}
\int_{\mathbb{Z}_{p}}e^{-2\pi i\left\{ t\mathfrak{z}\right\} _{p}}d\nu\left(\mathfrak{z}\right)=\mathbf{1}_{0}\left(t\right),\textrm{ }\forall t\in\hat{\mathbb{Z}}_{p}
\end{equation}
Then, $d\mu=d\nu$. 
\end{prop}
Proof: Direct computation using the formula: 
\begin{equation}
\left[\mathfrak{z}\overset{p^{n}}{\equiv}k\right]=\frac{1}{p^{n}}\sum_{\left|t\right|_{p}\leq p^{n}}e^{-2\pi i\left\{ t\left(\mathfrak{z}-k\right)\right\} _{p}}
\end{equation}
from \textbf{Proposition \ref{prop:Indicator for p-adic coset Fourier identity}}.

Q.E.D. 
\begin{defn}
Following \textbf{Proposition \ref{prop:Fourier transform of haar measure}},
I write $d\mathfrak{z}$ to denote the $\left(p,q\right)$-adic measure
defined by: 
\begin{equation}
\int_{\mathbb{Z}_{p}}\left[\mathfrak{z}\overset{p^{n}}{\equiv}k\right]d\mathfrak{z}\overset{\textrm{def}}{=}\frac{1}{p^{n}},\textrm{ }\forall n\in\mathbb{N}_{0},\textrm{ }\forall k\in\mathbb{Z}\label{eq:Definition of (p,q)-adic Haar probability measure of k + p^n Z_p}
\end{equation}
and: 
\begin{equation}
\int_{\mathbb{Z}_{p}}e^{-2\pi i\left\{ t\mathfrak{z}\right\} _{p}}d\mathfrak{z}\overset{\textrm{def}}{=}\mathbf{1}_{0}\left(t\right),\textrm{ }\forall t\in\hat{\mathbb{Z}}_{p}\label{eq:Definition of the Fourier-Stieltjes transform of the (p,q)-adic Haar probability measure}
\end{equation}
I call $d\mathfrak{z}$ the \textbf{$\left(p,q\right)$-adic Haar
(probability) measure}\index{measure!Haar}. Note that this is the
\emph{unique} measure satisfying the above two identities. 
\end{defn}
\begin{thm}
The set of \index{translation invariance}translation-invariant $\left(p,q\right)$-adic
measures is the one-dimensional subspace of $B\left(\hat{\mathbb{Z}}_{p},K\right)$
spanned by $d\mathfrak{z}$. Additionally, $d\mathfrak{z}$ is the
unique translation-invariant $\left(p,q\right)$-adic measure satisfying
the unit normalization condition: 
\begin{equation}
\int_{\mathbb{Z}_{p}}d\mathfrak{z}\overset{K}{=}1
\end{equation}
\end{thm}
Proof: First, let $d\mu=\mathfrak{c}d\mathfrak{z}$ for some $\mathfrak{c}\in K$.
Then, letting $f\in C\left(\mathbb{Z}_{p},\mathbb{C}_{q}\right)$
and $\mathfrak{a}\in\mathbb{Z}_{p}$ be arbitrary: 
\begin{align*}
\int_{\mathbb{Z}_{p}}f\left(\mathfrak{z}+\mathfrak{a}\right)d\mu\left(\mathfrak{z}\right) & =\mathfrak{c}\int_{\mathbb{Z}_{p}}\sum_{t\in\hat{\mathbb{Z}}_{p}}\hat{f}\left(t\right)e^{2\pi i\left\{ t\left(\mathfrak{z}+\mathfrak{a}\right)\right\} _{p}}d\mathfrak{z}\\
\left(\textrm{\textbf{Proposition }\textbf{\ref{prop:pq-adic integral interchange rule}}}\right); & =\mathfrak{c}\sum_{t\in\hat{\mathbb{Z}}_{p}}\hat{f}\left(t\right)e^{2\pi i\left\{ t\mathfrak{a}\right\} _{p}}\int_{\mathbb{Z}_{p}}e^{2\pi i\left\{ t\mathfrak{z}\right\} _{p}}d\mathfrak{z}\\
 & =\mathfrak{c}\sum_{t\in\hat{\mathbb{Z}}_{p}}\hat{f}\left(t\right)e^{2\pi i\left\{ t\mathfrak{a}\right\} _{p}}\mathbf{1}_{0}\left(-t\right)\\
 & =\mathfrak{c}\hat{f}\left(0\right)\\
 & =\mathfrak{c}\int_{\mathbb{Z}_{p}}f\left(\mathfrak{z}\right)d\mathfrak{z}\\
 & =\int_{\mathbb{Z}_{p}}f\left(\mathfrak{z}\right)d\mu\left(\mathfrak{z}\right)
\end{align*}
So, every scalar multiple of $d\mathfrak{z}$ is translation-invariant.

Next, let $d\mu$ be translation-invariant. Then, for all $t\in\hat{\mathbb{Z}}_{p}$:
\begin{equation}
\hat{\mu}\left(t\right)=\int_{\mathbb{Z}_{p}}e^{-2\pi i\left\{ t\mathfrak{z}\right\} _{p}}d\mu\left(\mathfrak{z}\right)=\int_{\mathbb{Z}_{p}}e^{-2\pi i\left\{ t\left(\mathfrak{z}+1\right)\right\} _{p}}d\mu\left(\mathfrak{z}\right)=e^{-2\pi it}\hat{\mu}\left(t\right)
\end{equation}
Since $\hat{\mu}\left(t\right)\in\mathbb{C}_{q}$ for all $t$, the
above forces $\hat{\mu}\left(t\right)=0$ for all $t\in\hat{\mathbb{Z}}_{p}\backslash\left\{ 0\right\} $.
As such: 
\[
\hat{\mu}\left(t\right)=\hat{\mu}\left(0\right)\mathbf{1}_{0}\left(t\right)
\]
which shows $d\mu\left(\mathfrak{z}\right)=\hat{\mu}\left(0\right)d\mathfrak{z}$.
Thus, $d\mu$ is in the span of $d\mathfrak{z}$. Hence, every translation-invariant
measure is a scalar multiple of $d\mathfrak{z}$.

Finally, since $\int_{\mathbb{Z}_{p}}d\mu\left(\mathfrak{z}\right)=\hat{\mu}\left(0\right)$,
we have that $d\mathfrak{z}$ itself is the unique translation-invariant
$\left(p,q\right)$-adic measure $d\mu\left(\mathfrak{z}\right)$
satisfying $\int_{\mathbb{Z}_{p}}d\mu\left(\mathfrak{z}\right)=1$.

Q.E.D. 
\begin{lem}[\textbf{Integral Formula for the Fourier Transform}]
For\index{Fourier!transform} any $f\in C\left(\mathbb{Z}_{p},\mathbb{C}_{q}\right)$,
its Fourier transform $\hat{f}:\hat{\mathbb{Z}}_{p}\rightarrow\mathbb{C}_{q}$
is given by the integral formula: 
\begin{equation}
\hat{f}\left(t\right)\overset{\mathbb{C}_{q}}{=}\int_{\mathbb{Z}_{p}}f\left(\mathfrak{z}\right)e^{-2\pi i\left\{ t\mathfrak{z}\right\} _{p}}d\mathfrak{z},\textrm{ }\forall t\in\hat{\mathbb{Z}}_{p}\label{eq:Integral Formula for the Fourier transform}
\end{equation}
\end{lem}
Proof: Fix an arbitrary $t\in\hat{\mathbb{Z}}_{p}$. Then, by definition,
(\ref{eq:Integral Formula for the Fourier transform}) is the image
of the continuous $\left(p,q\right)$-adic function $f\left(\mathfrak{z}\right)e^{-2\pi i\left\{ t\mathfrak{z}\right\} _{p}}$
under the $\left(p,q\right)$-adic Haar measure. Expressing $f$ as
a Fourier series, we have that: 
\[
e^{-2\pi i\left\{ t\mathfrak{z}\right\} _{p}}f\left(\mathfrak{z}\right)\overset{\mathbb{C}_{q}}{=}e^{-2\pi i\left\{ t\mathfrak{z}\right\} _{p}}\sum_{s\in\hat{\mathbb{Z}}_{p}}\hat{f}\left(s\right)e^{2\pi i\left\{ s\mathfrak{z}\right\} _{p}}=\sum_{s\in\hat{\mathbb{Z}}_{p}}\hat{f}\left(s\right)e^{2\pi i\left\{ \left(s-t\right)\mathfrak{z}\right\} _{p}}
\]
for all $\mathfrak{z}\in\mathbb{Z}_{p}$. Thus, by linearity, since
the Fourier series is uniformly convergent in $\mathbb{C}_{q}$, \textbf{Proposition
\ref{prop:term-by-term evaluation of integrals}} justifies integrating
term-by-term: 
\begin{align*}
\int_{\mathbb{Z}_{p}}f\left(\mathfrak{z}\right)e^{-2\pi i\left\{ t\mathfrak{z}\right\} _{p}}d\mathfrak{z} & =\sum_{s\in\hat{\mathbb{Z}}_{p}}\hat{f}\left(s\right)\int_{\mathbb{Z}_{p}}e^{2\pi i\left\{ \left(s-t\right)\mathfrak{z}\right\} _{p}}d\mathfrak{z}\\
 & =\sum_{s\in\hat{\mathbb{Z}}_{p}}\hat{f}\left(s\right)\underbrace{\mathbf{1}_{0}\left(s-t\right)}_{1\textrm{iff }s=t}\\
 & =\hat{f}\left(t\right)
\end{align*}
which gives the desired result.

Q.E.D.

\vphantom{}

With this integral formula, we then get the usual integral formulas
for convolution\index{convolution!of left(p,qright)-adic functions@of $\left(p,q\right)$-adic functions},
and they behave exactly as we would expect them to behave. 
\begin{defn}[\textbf{Convolution}]
We define the \textbf{convolution }of functions on $\mathbb{Z}_{p}$
and $\hat{\mathbb{Z}}_{p}$, respectively, by: 
\begin{equation}
\left(f*g\right)\left(\mathfrak{z}\right)\overset{\mathbb{C}_{q}}{=}\int_{\mathbb{Z}_{p}}f\left(\mathfrak{z}-\mathfrak{y}\right)g\left(\mathfrak{y}\right)d\mathfrak{y},\textrm{ }\forall\mathfrak{z}\in\mathbb{Z}_{p},\forall f,g\in C\left(\mathbb{Z}_{p},\mathbb{C}_{q}\right)\label{eq:Continuous Convolution Definition}
\end{equation}
\begin{equation}
\left(\hat{f}*\hat{g}\right)\left(t\right)\overset{\mathbb{C}_{q}}{=}\sum_{\tau\in\hat{\mathbb{Z}}_{p}}\hat{f}\left(t-\tau\right)\hat{g}\left(\tau\right),\textrm{ }\forall t\in\hat{\mathbb{Z}}_{p},\forall\hat{f},\hat{g}\in c_{0}\left(\hat{\mathbb{Z}}_{p},\mathbb{C}_{q}\right)\label{eq:Discrete Convolution Definition}
\end{equation}
\end{defn}
\begin{thm}[\textbf{The} \textbf{Convolution Theorem}]
\label{thm:Convolution Theorem}\index{convolution!theorem}For $f,g\in C\left(\mathbb{Z}_{p},\mathbb{C}_{q}\right)$:
\begin{align}
\widehat{\left(f*g\right)}\left(t\right) & =\hat{f}\left(t\right)\hat{g}\left(t\right)\label{eq:Convolution Theorem 1}\\
\widehat{\left(fg\right)}\left(t\right) & =\left(\hat{f}*\hat{g}\right)\left(t\right)\label{eq:Convolution Theorem 2}
\end{align}
\end{thm}
\vphantom{}

We can also take convolutions of measures: 
\begin{defn}
Let $d\mu,d\nu\in C\left(\mathbb{Z}_{p},\mathbb{C}_{q}\right)^{\prime}$
be $\left(p,q\right)$-adic measures. Then, the \textbf{convolution
}of $d\mu$ and $d\nu$, denoted $d\mu*d\nu$, is the measure defined
by the Fourier-Stieltjes transform: 
\begin{equation}
\widehat{\left(d\mu*d\nu\right)}\left(t\right)\overset{\textrm{def}}{=}\hat{\mu}\left(t\right)\hat{\nu}\left(t\right)\label{eq:Definition of the convolution of two pq adic measures}
\end{equation}
That is: 
\begin{equation}
\int_{\mathbb{Z}_{p}}f\left(\mathfrak{z}\right)\left(d\mu*d\nu\right)\left(\mathfrak{z}\right)\overset{\mathbb{C}_{q}}{=}\sum_{t\in\hat{\mathbb{Z}}_{p}}\hat{f}\left(-t\right)\hat{\mu}\left(t\right)\hat{\nu}\left(t\right),\textrm{ }\forall f\in C\left(\mathbb{Z}_{p},\mathbb{C}_{q}\right)\label{eq:Definition of the action of the convolution of pq adic measures on a function}
\end{equation}
\end{defn}
\vphantom{}

Next up, \textbf{Parseval-Plancherel Identity}:\index{Parseval-Plancherel Identity} 
\begin{thm}[\textbf{The Parseval-Plancherel Identity}]
\label{thm:Parseval-Plancherel Identity}Let $f,g\in C\left(\mathbb{Z}_{p},\mathbb{C}_{q}\right)$.
Then: 
\end{thm}
\begin{equation}
\int_{\mathbb{Z}_{p}}f\left(\mathfrak{z}\right)g\left(\mathfrak{z}\right)d\mathfrak{z}\overset{\mathbb{C}_{q}}{=}\sum_{t\in\hat{\mathbb{Z}}_{p}}\hat{f}\left(t\right)\hat{g}\left(-t\right),\textrm{ }\forall f,g\in C\left(\mathbb{Z}_{p},\mathbb{C}_{q}\right)\label{eq:Parseval-Plancherel Identity}
\end{equation}

\vphantom{}

In all of these formulae, convergence, rearrangements, and interchanges
are justified by the $q$-adic decay of $\hat{f}$ and $\hat{g}$
which is guaranteed to occur thanks to the continuity of $f$ and
$g$. 
\begin{prop}
Let $g\in C\left(\mathbb{Z}_{p},\mathbb{C}_{q}\right)$ and let $d\mu\in C\left(\mathbb{Z}_{p},\mathbb{C}_{q}\right)^{\prime}$.
Viewing $g$ as the measure: 
\[
f\mapsto\int_{\mathbb{Z}_{p}}f\left(\mathfrak{z}\right)g\left(\mathfrak{z}\right)d\mathfrak{z}
\]
we can then define the convolution of $g$ and $d\mu$ as the measure
with the $\left(p,q\right)$-adic Fourier-Stieltjes transform: 
\[
\widehat{\left(g*d\mu\right)}\left(t\right)=\sum_{t\in\hat{\mathbb{Z}}_{p}}\hat{g}\left(t\right)\hat{\mu}\left(t\right)
\]
The measure is $f*d\mu$ is then absolutely continuous with respect
to the $\left(p,q\right)$-adic Haar measure $d\mathfrak{z}$, meaning
that there is a continuous $\left(p,q\right)$-adic function $g\in C\left(\mathbb{Z}_{p},\mathbb{C}_{q}\right)$
so that: 
\[
\left(f*d\mu\right)\left(\mathfrak{z}\right)=g\left(\mathfrak{z}\right)d\mathfrak{z}
\]
In fact, we have that $\hat{g}\left(t\right)=\hat{f}\left(t\right)\hat{\mu}\left(t\right)$,
and hence, we can\footnote{That is, the convolution of a continuous $\left(p,q\right)$-adic
function and a $\left(p,q\right)$-adic measure is a continuous function.} view $f*d\mu$ as a continuous function $\mathbb{Z}_{p}\rightarrow\mathbb{C}_{q}$. 
\end{prop}
Proof: Since $f$ is continuous, $\left|\hat{f}\left(t\right)\right|_{q}\rightarrow0$
as $\left|t\right|_{p}\rightarrow\infty$. Since $d\mu$ is a measure,
$\left\Vert \hat{\mu}\right\Vert _{p,q}=\sup_{t\in\hat{\mathbb{Z}}_{p}}\left|\hat{\mu}\left(t\right)\right|_{q}<\infty$,
hence, defining $\hat{g}\left(t\right)\overset{\textrm{def}}{=}\hat{f}\left(t\right)\hat{\mu}\left(t\right)$,
we have that $\hat{g}\in c_{0}\left(\hat{\mathbb{Z}}_{p},\mathbb{C}_{q}\right)$,
which shows that $\hat{g}$ is the Fourier transform of the continuous
$\left(p,q\right)$-adic function: 
\[
g\left(\mathfrak{z}\right)\overset{\textrm{def}}{=}\sum_{t\in\hat{\mathbb{Z}}_{p}}\hat{g}\left(t\right)e^{2\pi i\left\{ t\mathfrak{z}\right\} _{p}},\textrm{ }\forall\mathfrak{z}\in\mathbb{Z}_{p}
\]
Furthermore, for any $h\in C\left(\mathbb{Z}_{p},\mathbb{C}_{q}\right)$:
\begin{align*}
\int_{\mathbb{Z}_{p}}h\left(\mathfrak{z}\right)\left(f*d\mu\right)\left(\mathfrak{z}\right) & \overset{\mathbb{C}_{q}}{=}\sum_{t\in\hat{\mathbb{Z}}_{p}}\hat{h}\left(-t\right)\hat{f}\left(t\right)\hat{\mu}\left(t\right)\\
 & \overset{\mathbb{C}_{q}}{=}\sum_{t\in\hat{\mathbb{Z}}_{p}}\hat{h}\left(-t\right)\hat{g}\left(t\right)\\
 & \overset{\mathbb{C}_{q}}{=}\int_{\mathbb{Z}_{p}}h\left(\mathfrak{z}\right)g\left(\mathfrak{z}\right)d\mathfrak{z}
\end{align*}
Hence, $\left(f*d\mu\right)\left(\mathfrak{z}\right)$ and $g\left(\mathfrak{z}\right)d\mathfrak{z}$
are identical as measures.

Q.E.D.

\vphantom{}

For doing computations with $\left(p,q\right)$-adic integrals, the
following simple change-of-variables\index{change of variables!affine substitutions}
formula will be of the utmost import: 
\begin{lem}[\textbf{Change of Variables \textendash{} Affine substitutions}]
\label{lem:Affine substitution change of variable formula}Let $f\in C\left(\mathbb{Z}_{p},K\right)$,
where $K$ is a complete ring extension of $\mathbb{Z}_{q}$ which
is itself contained in $\mathbb{C}_{q}$. Then: 
\begin{equation}
\int_{\mathfrak{a}\mathbb{Z}_{p}}f\left(\mathfrak{z}\right)d\mathfrak{z}=\int_{p^{v_{p}\left(\mathfrak{a}\right)}\mathbb{Z}_{p}}f\left(\mathfrak{z}\right)d\mathfrak{z},\textrm{ }\forall\mathfrak{a}\in\mathbb{Z}_{p}\backslash\left\{ 0\right\} \label{eq:Multiplying domain of integration by a}
\end{equation}
\begin{equation}
\int_{\mathbb{Z}_{p}}f\left(\mathfrak{a}\mathfrak{z}+\mathfrak{b}\right)d\mathfrak{z}\overset{K}{=}\frac{1}{\left|\mathfrak{a}\right|_{p}}\int_{\mathfrak{a}\mathbb{Z}_{p}+\mathfrak{b}}f\left(\mathfrak{z}\right)d\mathfrak{z},\textrm{ }\forall\mathfrak{a}\in\mathbb{Z}_{p}\backslash\left\{ 0\right\} ,\forall\mathfrak{b}\in\mathbb{Z}_{p}\label{eq:Change of variables formula - affine linear substitutions}
\end{equation}
\end{lem}
\begin{rem}
When working with these integrals, a key identity is: 
\[
\int_{p^{n}\mathbb{Z}_{p}+k}f\left(\mathfrak{z}\right)d\mu\left(\mathfrak{z}\right)=\int_{\mathbb{Z}_{p}}\left[\mathfrak{z}\overset{p^{n}}{\equiv}k\right]f\left(\mathfrak{z}\right)d\mu\left(\mathfrak{z}\right)
\]
which holds for all $f\in C\left(\mathbb{Z}_{p},\mathbb{C}_{q}\right)$,
all $d\mu\in C\left(\mathbb{Z}_{p},\mathbb{C}_{q}\right)$, and all
$n,k\in\mathbb{Z}$ with $n\geq0$. 
\end{rem}
Proof: (\ref{eq:Multiplying domain of integration by a}) follows
from the fact that $\mathfrak{a}\neq0$ implies there are (unique)
$m\in\mathbb{N}_{0}$ and $\mathfrak{u}\in\mathbb{Z}_{p}^{\times}$
so that $\mathfrak{a}=p^{m}\mathfrak{u}$. Since multiplication by
$\mathfrak{u}$ is a measure-preserving automorphism of the group
$\left(\mathbb{Z}_{p},+\right)$, $\mathfrak{a}\mathbb{Z}_{p}=p^{m}\mathfrak{u}\mathbb{Z}_{p}=p^{m}\mathbb{Z}_{p}$.
Here, $m=v_{p}\left(\mathfrak{z}\right)$.

For (\ref{eq:Change of variables formula - affine linear substitutions}),
by the translation-invariance of $d\mathfrak{z}$, it suffices to
prove: 
\[
\int_{\mathbb{Z}_{p}}f\left(\mathfrak{a}\mathfrak{z}\right)d\mathfrak{z}\overset{K}{=}\frac{1}{\left|\mathfrak{a}\right|_{p}}\int_{\mathfrak{a}\mathbb{Z}_{p}}f\left(\mathfrak{z}\right)d\mathfrak{z}
\]
In fact, using the van der Put basis, we need only verify the formula
for the indicator functions $f\left(\mathfrak{z}\right)=\left[\mathfrak{z}\overset{p^{n}}{\equiv}k\right]$.
Because the integral of such an $f$ is rational-valued, note that
the computation will be the same for the case of the \emph{real-valued
}Haar probability measure $d\mathfrak{z}$ on $\mathbb{Z}_{p}$. For
any field extension $\mathbb{F}$ of $\mathbb{Q}$, any translation-invariant
linear functional on the space of $\mathbb{Q}$-valued functions on
$\mathbb{Z}_{p}$ normalized to send the constant function $1$ to
the number $1$ necessarily sends the functions $\left[\mathfrak{z}\overset{p^{n}}{\equiv}k\right]$
to the number $1/p^{n}$. Thus, the proof for the real-valued case,
such as can be found in \cite{Bell - Harmonic Analysis on the p-adics}
or \cite{Automorphic Representations} applies, and we are done.

Q.E.D. 
\begin{rem}
(\ref{eq:Change of variables formula - affine linear substitutions})
is \emph{also }valid in the case where $f:\mathbb{Z}_{p}\rightarrow\mathbb{C}$
($\mathbb{C}$, not $\mathbb{C}_{q}$!) is integrable with respect
to the \emph{real-valued }Haar probability measure on $\mathbb{Z}_{p}$. 
\end{rem}
\begin{rem}
Seeing as we write $\left[\mathfrak{z}\overset{p^{n}}{\equiv}\mathfrak{b}\right]$
to denote the indicator function for the set $\mathfrak{b}+p^{n}\mathbb{Z}_{p}$,
we then have the identities:

\begin{equation}
\int_{\mathbb{Z}_{p}}\left[\mathfrak{z}\overset{p^{n}}{\equiv}\mathfrak{b}\right]f\left(\mathfrak{z}\right)d\mathfrak{z}=\int_{\mathfrak{b}+p^{n}\mathbb{Z}_{p}}f\left(\mathfrak{z}\right)d\mathfrak{z}=\frac{1}{p^{n}}\int_{\mathbb{Z}_{p}}f\left(p^{n}\mathfrak{y}+\mathfrak{b}\right)d\mathfrak{y}\label{eq:change of variable trifecta}
\end{equation}
\end{rem}
\vphantom{}

Lastly, we have the basic integral inequalities. Alas, the triangle
inequality is not among them. 
\begin{prop}[\textbf{Integral Inequalities}]
\index{triangle inequality!left(p,qright)-adic@$\left(p,q\right)$-adic}Let
$f\in C\left(\mathbb{Z}_{p},\mathbb{C}_{q}\right)$. Then:

\vphantom{}

I. 
\begin{equation}
\left|\int_{\mathbb{Z}_{p}}f\left(\mathfrak{z}\right)d\mathfrak{z}\right|_{q}\leq\sup_{\mathfrak{z}\in\mathbb{Z}_{p}}\left|f\left(\mathfrak{z}\right)\right|_{q}=\left\Vert f\right\Vert _{p,q}\label{eq:Triangle Inequality for the (p,q)-adic Haar measure}
\end{equation}

\vphantom{}

II. For any measure $d\mu\in C\left(\mathbb{Z}_{p},\mathbb{C}_{q}\right)^{\prime}$:
\begin{equation}
\left|\int_{\mathbb{Z}_{p}}f\left(\mathfrak{z}\right)d\mu\left(\mathfrak{z}\right)\right|_{q}\leq\sup_{t\in\hat{\mathbb{Z}}_{p}}\left|\hat{f}\left(-t\right)\hat{\mu}\left(t\right)\right|_{q}\label{eq:Triangle inequality for an arbitrary (p,q)-adic measure (with Fourier)}
\end{equation}

\vphantom{}

III.\index{triangle inequality!van der Put series} 
\begin{equation}
\int_{\mathbb{Z}_{p}}\left|f\left(\mathfrak{z}\right)\right|_{q}d\mathfrak{z}\leq\sum_{n=0}^{\infty}\frac{\left|c_{n}\left(f\right)\right|_{q}}{p^{\lambda_{p}\left(n\right)}}\label{eq:Triangle inequality for Integral of q-adic absolute value}
\end{equation}
\end{prop}
\begin{rem}
As shown in \textbf{Example} \textbf{\ref{exa:triangle inequality failure}}
(page \pageref{exa:triangle inequality failure}), unlike in real
or complex analysis, there is not necessarily any relationship between
$\left|\int f\right|_{q}$ and $\int\left|f\right|_{q}$. 
\end{rem}
Proof: (I) is the consequence of the Theorem from page 281 of \emph{Ultrametric
Calculus }(\cite{Ultrametric Calculus}) and the Exercise that occurs
immediately after it. (II), meanwhile, is an immediate consequence
of (\ref{eq:Integration in terms of Fourier-Stieltjes coefficients}).
(III) is the only one that requires an argument.

For (III), let $f\in C\left(\mathbb{Z}_{p},\mathbb{C}_{q}\right)$.
Then: 
\begin{equation}
\int_{\mathbb{Z}_{p}}f\left(\mathfrak{z}\right)d\mathfrak{z}\overset{K}{=}\sum_{n=0}^{\infty}c_{n}\left(f\right)\int_{\mathbb{Z}_{p}}\left[\mathfrak{z}\overset{p^{\lambda_{p}\left(n\right)}}{\equiv}n\right]d\mathfrak{z}=\sum_{n=0}^{\infty}\frac{c_{n}\left(f\right)}{p^{\lambda_{p}\left(n\right)}}
\end{equation}
where the interchange of sum and integral/linear functional is valid
due to the convergence of the series, guaranteed by the $q$-adic
decay of the $c_{n}\left(f\right)$s. Now, since $f$ is continuous,
$\left|f\left(\mathfrak{z}\right)\right|_{q}$ is a continuous function
from $\mathbb{Z}_{p}$ to $\mathbb{R}$\textemdash uniformly continuous,
in fact, by the compactness of $\mathbb{Z}_{p}$. Thus, by \textbf{Proposition
\ref{prop:Convergence of real-valued vdP series}}, $\left|f\left(\mathfrak{z}\right)\right|_{q}$'s
van der Put series: 
\begin{equation}
\sum_{n=0}^{\infty}c_{n}\left(\left|f\right|_{q}\right)\left[\mathfrak{z}\overset{p^{\lambda_{p}\left(n\right)}}{\equiv}n\right]=\left|f\left(0\right)\right|_{q}+\sum_{n=0}^{\infty}\underbrace{\left(\left|f\left(n\right)\right|_{q}-\left|f\left(n_{-}\right)\right|_{q}\right)}_{c_{n}\left(\left|f\right|_{q}\right)}\left[\mathfrak{z}\overset{p^{\lambda_{p}\left(n\right)}}{\equiv}n\right]
\end{equation}
is uniformly convergent in $\mathbb{R}$. As such, when integrating
it with respect to the real-valued $p$-adic Haar measure, we may
interchange integration and summation: 
\begin{equation}
\int_{\mathbb{Z}_{p}}\left|f\left(\mathfrak{z}\right)\right|_{q}d\mathfrak{z}=\sum_{n=0}^{\infty}c_{n}\left(\left|f\right|_{q}\right)\int_{\mathbb{Z}_{p}}\left[\mathfrak{z}\overset{p^{\lambda_{p}\left(n\right)}}{\equiv}n\right]d\mathfrak{z}=\sum_{n=0}^{\infty}\frac{c_{n}\left(\left|f\right|_{q}\right)}{p^{\lambda_{p}\left(n\right)}}
\end{equation}
and the infinite series on the right is then necessarily convergent,
since $f$'s uniformly continuity guarantees its integrability. Here:
\begin{equation}
\sum_{n=0}^{\infty}\frac{c_{n}\left(\left|f\right|_{q}\right)}{p^{\lambda_{p}\left(n\right)}}=\sum_{n=0}^{\infty}\frac{\left|f\left(n\right)\right|_{q}-\left|f\left(n_{-}\right)\right|_{q}}{p^{\lambda_{p}\left(n\right)}}
\end{equation}
Hence: 
\begin{align*}
\int_{\mathbb{Z}_{p}}\left|f\left(\mathfrak{z}\right)\right|_{q}d\mathfrak{z} & =\left|\int_{\mathbb{Z}_{p}}\left|f\left(\mathfrak{z}\right)\right|_{q}d\mathfrak{z}\right|\\
 & \leq\sum_{n=0}^{\infty}\frac{\left|\left|f\left(n\right)\right|_{q}-\left|f\left(n_{-}\right)\right|_{q}\right|}{p^{\lambda_{p}\left(n\right)}}\\
\left(\textrm{reverse }\Delta\textrm{-ineq.}\right); & \leq\sum_{n=0}^{\infty}\frac{\left|f\left(n\right)-f\left(n_{-}\right)\right|_{q}}{p^{\lambda_{p}\left(n\right)}}\\
 & =\sum_{n=0}^{\infty}\frac{\left|c_{n}\left(f\right)\right|_{q}}{p^{\lambda_{p}\left(n\right)}}
\end{align*}
as desired.

Q.E.D.

\vphantom{}

The most important practical consequence of this approach to $\left(p,q\right)$-adic
integration is that it trivializes most of the basic concerns of classical
real analysis. Because of our definition of $\left(p,q\right)$-adic
measures as continuous linear functionals on the space of continuous
$\left(p,q\right)$-adic functions, the equivalence of continuity
with representability by a uniformly convergent Fourier series along
with the action of measures on functions in terms of their Fourier
series expansions shows that all continuous $\left(p,q\right)$-adic
functions are integrable with respect to the $\left(p,q\right)$-adic
Haar probability measure. Indeed, the following argument illustrates
this quite nicely. 
\begin{example}
\label{exa:end of 3.1.5. example}Let $f:\mathbb{Z}_{p}\rightarrow\mathbb{C}_{q}$
be given by the van der Put series: 
\begin{equation}
f\left(\mathfrak{z}\right)\overset{\mathbb{C}_{q}}{=}\sum_{n=0}^{\infty}\mathfrak{a}_{n}\left[\mathfrak{z}\overset{p^{\lambda_{p}\left(n\right)}}{\equiv}n\right]
\end{equation}
where the series is assumed to converge merely point-wise. Letting:
\begin{equation}
f_{N}\left(\mathfrak{z}\right)\overset{\textrm{def}}{=}\sum_{n=0}^{p^{N}-1}\mathfrak{a}_{n}\left[\mathfrak{z}\overset{p^{\lambda_{p}\left(n\right)}}{\equiv}n\right]
\end{equation}
observe that, like in the classic real-analytical Lebesgue theory,
the $f_{N}$s are step functions via (\ref{eq:truncated van der Put identity}):
\begin{equation}
f_{N}\left(\mathfrak{z}\right)=f\left(\left[\mathfrak{z}\right]_{p^{N}}\right)=\sum_{n=0}^{p^{N}-1}f\left(n\right)\left[\mathfrak{z}\overset{p^{N}}{\equiv}n\right]
\end{equation}
which converge to $f$. Similar to what is done in the theory of Lebesgue
integration of real-valued functions on measure spaces, in Subsection
\ref{subsec:3.1.6 Monna-Springer-Integration}, we will construct
the class of integrable functions by considering limits of sequences
of step functions. To that end, observe that: 
\begin{equation}
\int_{\mathbb{Z}_{p}}f_{N}\left(\mathfrak{z}\right)d\mathfrak{z}\overset{\mathbb{C}_{q}}{=}\sum_{n=0}^{p^{N}-1}\mathfrak{a}_{n}\int_{\mathbb{Z}_{p}}\left[\mathfrak{z}\overset{p^{\lambda_{p}\left(n\right)}}{\equiv}n\right]d\mathfrak{z}=\sum_{n=0}^{p^{N}-1}\frac{\mathfrak{a}_{n}}{p^{\lambda_{p}\left(n\right)}}
\end{equation}
Since $p$ and $q$ are co-prime, $\left|\mathfrak{a}_{n}p^{-\lambda_{p}\left(n\right)}\right|_{q}=\left|\mathfrak{a}_{n}\right|_{q}$,
which shows that: 
\begin{equation}
\lim_{N\rightarrow\infty}\int_{\mathbb{Z}_{p}}f_{N}\left(\mathfrak{z}\right)d\mathfrak{z}
\end{equation}
exists in $\mathbb{C}_{q}$ if and only if $\lim_{n\rightarrow\infty}\left|\mathfrak{a}_{n}\right|_{q}\overset{\mathbb{R}}{=}0$.
However, as we saw in \textbf{Theorem \ref{thm:vdP basis theorem}},
$\left|\mathfrak{a}_{n}\right|_{q}\rightarrow0$ if and only if $f$
is continuous! Schikhof leaves as an exercise in \emph{Ultrametric
Calculus}' discussion of the van der Put basis the fact that $f\left(\left[\mathfrak{z}\right]_{p^{N}}\right)$
is the best possible approximation of $f$ in terms of the simple
functions $\left[\mathfrak{z}\overset{p^{N}}{\equiv}n\right]$ for
$n\in\left\{ 0,\ldots,p^{N}-1\right\} $ \cite{Ultrametric Calculus}.
As such, if we are going to use step-function-approximation to construct
integrable functions, this analysis of ours tells us that only the
continuous functions are going to be integrable. 
\end{example}
\vphantom{}

Before we conclude, we need to cover one last special case where integration
becomes quite simple. 
\begin{defn}
Let $\mathbb{F}$ be any field, let $f:\mathbb{Z}_{p}\rightarrow\mathbb{F}$,
and let $n\in\mathbb{N}_{0}$. We say $f$ is \textbf{constant over
inputs modulo $p^{n}$} whenever $f\left(\mathfrak{z}\right)=f\left(\left[\mathfrak{z}\right]_{p^{n}}\right)$
for all $\mathfrak{z}\in\mathbb{Z}_{p}$. 
\end{defn}
\begin{prop}
Let $f:\mathbb{Z}_{p}\rightarrow\mathbb{F}$ be constant over inputs
modulo $p^{N}$. Then: 
\begin{equation}
f\left(\mathfrak{z}\right)=\sum_{n=0}^{p^{N}-1}f\left(n\right)\left[\mathfrak{z}\overset{p^{N}}{\equiv}n\right]\label{eq:representation of functions constant over inputs mod p^N}
\end{equation}
\end{prop}
Proof: Immediate.

Q.E.D. 
\begin{prop}
\label{prop:How to integrate locally constant functions}Let $\mathbb{F}$
be a metrically complete valued field, and let $d\mu$ be a translation
invariant $\mathbb{F}$-valued probability measure on $\mathbb{Z}_{p}$;
that is: 
\begin{equation}
\int_{\mathbb{Z}_{p}}d\mu\left(\mathfrak{z}\right)\overset{\mathbb{F}}{=}1
\end{equation}
and: 
\begin{equation}
\int_{\mathbb{Z}_{p}}f\left(\mathfrak{z}\right)d\mu\left(\mathfrak{z}\right)=\int_{\mathbb{Z}_{p}}f\left(\mathfrak{z}+\mathfrak{a}\right)d\mu\left(\mathfrak{z}\right)\in\mathbb{F},\textrm{ }\forall\mathfrak{a}\in\mathbb{Z}_{p},\textrm{ }\forall f\in C\left(\mathbb{Z}_{p},\mathbb{F}\right)
\end{equation}
If $f$ is constant over inputs modulo $p^{N}$, then: 
\begin{equation}
\int_{\mathbb{Z}_{p}}f\left(\mathfrak{z}\right)d\mu=\frac{1}{p^{N}}\sum_{n=0}^{p^{N}-1}f\left(n\right)
\end{equation}
\end{prop}
Proof: By the translation invariance of the $\mathbb{F}$-valued probability
measure $d\mu$, it follows that: 
\begin{align*}
1 & =\int_{\mathbb{Z}_{p}}d\mu\left(\mathfrak{z}\right)\\
 & =\int_{\mathbb{Z}_{p}}\sum_{n=0}^{p^{N}-1}\left[\mathfrak{z}\overset{p^{N}}{\equiv}n\right]d\mu\left(\mathfrak{z}\right)\\
 & =\sum_{n=0}^{p^{N}-1}\int_{\mathbb{Z}_{p}}\left[\mathfrak{z}\overset{p^{N}}{\equiv}0\right]d\mu\left(\mathfrak{z}\right)\\
 & =p^{N}\int_{\mathbb{Z}_{p}}\left[\mathfrak{z}\overset{p^{N}}{\equiv}0\right]d\mu\left(\mathfrak{z}\right)
\end{align*}
Hence: 
\[
\int_{\mathbb{Z}_{p}}\left[\mathfrak{z}\overset{p^{N}}{\equiv}0\right]d\mu\left(\mathfrak{z}\right)\overset{\mathbb{F}}{=}\frac{1}{p^{N}}
\]
and so, by translation invariance: 
\[
\int_{\mathbb{Z}_{p}}\left[\mathfrak{z}\overset{p^{N}}{\equiv}n\right]d\mu\left(\mathfrak{z}\right)\overset{\mathbb{F}}{=}\frac{1}{p^{N}},\textrm{ }\forall n,N
\]
By linearity, integrating the formula (\ref{eq:representation of functions constant over inputs mod p^N})
yields the desired formula for the integral of $f\left(\mathfrak{z}\right)d\mu\left(\mathfrak{z}\right)$.

Q.E.D.

\vphantom{}

Note that this method can be used to integrate both $\left(p,q\right)$-adic
functions \emph{and }$\left(p,\infty\right)$-adic functions, seeing
as (provided that $p\neq q$), both those cases admit translation
invariant probability measures on $\mathbb{Z}_{p}$.

\subsection{\label{subsec:3.1.6 Monna-Springer-Integration}Monna-Springer Integration}

FOR THIS SUBSECTION, $K$ IS A COMPLETE ULTRAMETRIC FIELD, $p$ AND
$q$ ARE DISTINCT PRIMES, AND $X$ IS AN ARBITRARY ULTRAMETRIC SPACE.

\vphantom{}

The simplicity, directness, practicality and technical clarity of
the basis-based approach to $\left(p,q\right)$-adic integration we
have covered so far all speak their usefulness in doing $\left(p,q\right)$-adic
analysis. Indeed, it will completely suffice for all of our investigations
of $\chi_{H}$. That this approach works out so nicely is all thanks
to the translation-invariance of the $\left(p,q\right)$-adic Haar
measure. At the same time, note that we have failed to mention many
of the issues central to most any theory of integration: measurable
sets, measurable functions, and the interchanging of limits and integrals.
To do these topics justice, we need Monna-Springer integration theory:
the method for defining measures and integration of $K$-valued functions
on an arbitrary ultrametric space $X$. Presenting Monna-Springer
theory means bombarding the reader with many definitions and almost
just as many infima and suprema. As such, to maximize clarity, I have
supplemented the standard exposition with examples from the $\left(p,q\right)$-case
and additional discussion.

To begin, it is worth reflecting on the notion of measure defined
in Subsection \ref{subsec:3.1.5-adic-Integration-=00003D000026}.
Because we defined measure as continuous linear functionals on the
space of continuous $\left(p,q\right)$-adic functions, we will be
able to say a set in $\mathbb{Z}_{p}$ is $\left(p,q\right)$-adically
measurable'' precisely when the indicator function for that set is
$\left(p,q\right)$-adically continuous. In classical integration
theory\textemdash particularly when done in the Bourbaki-style functional
analysis approach\textemdash indicator functions provide a bridge
between the world of functions and the world of sets. Indicator functions
for measurable subsets of $\mathbb{R}$ can be approximated in $L^{1}$
to arbitrary accuracy by piece-wise continuous functions.

What makes non-archimedean integration theory so different from its
archimedean counterpart are the restrictions which $\mathbb{C}_{q}$'s
ultrametric structure imposes on continuous functions. Because of
the equivalence of integrable functions and continuous functions,
the only sets in $\mathbb{Z}_{p}$ which we can call $\left(p,q\right)$-adically
measurable are those whose indicator functions are $\left(p,q\right)$-adically
continuous. Since indicator functions $\left[\mathfrak{z}\overset{p^{n}}{\equiv}k\right]$
are continuous, it follows that any finite union, intersection, or
complement of sets of the form $k+p^{n}\mathbb{Z}_{p}$ ``$\left(p,q\right)$-adically
measurable''. Note that all of these sets are compact clopen sets.
On the other hand, single points will \emph{not }be measurable. 
\begin{defn}
Fix $\mathfrak{a}\in\mathbb{Z}_{p}$. Then, define the function \nomenclature{$\mathbf{1}_{\mathfrak{a}}$}{ }$\mathbf{1}_{\mathfrak{a}}:\mathbb{Z}_{p}\rightarrow K$
by: 
\begin{equation}
\mathbf{1}_{\mathfrak{a}}\left(\mathfrak{z}\right)=\begin{cases}
1 & \textrm{if }\mathfrak{z}=\mathfrak{a}\\
0 & \textrm{else}
\end{cases}\label{eq:Definition of a one-point function}
\end{equation}
I call any non-zero scalar multiple of (\ref{eq:Definition of a one-point function})
a \textbf{one-point function} supported at $\mathfrak{a}$. More generally,
I call any finite linear combination of the form: 
\begin{equation}
\sum_{n=1}^{N}\mathfrak{c}_{n}\mathbf{1}_{\mathfrak{a}_{n}}\left(\mathfrak{z}\right)
\end{equation}
an\textbf{ $N$-point function}; here, the $\mathfrak{a}_{n}$s in
$\mathbb{Z}_{p}$ and the $\mathfrak{c}_{n}$s in $K\backslash\left\{ 0\right\} $. 
\end{defn}
\begin{rem}
It is a simple exercise to prove the formulae: 
\begin{equation}
\mathbf{1}_{0}\left(\mathfrak{z}\right)=1-\sum_{n=0}^{\infty}\sum_{k=1}^{p-1}\left[\mathfrak{z}\overset{p^{n+1}}{\equiv}kp^{n}\right]\label{eq:van der Point series for one point function at 0}
\end{equation}
\begin{equation}
\mathbf{1}_{\mathfrak{a}}\left(\mathfrak{z}\right)=1-\sum_{n=0}^{\infty}\sum_{k=1}^{p-1}\left[\mathfrak{z}-\mathfrak{a}\overset{p^{n+1}}{\equiv}kp^{n}\right]\label{eq:vdP series for a one-point function}
\end{equation}
\end{rem}
\begin{example}
The indicator function of a single point in $\mathbb{Z}_{p}$ is not
integrable with respect to the $\left(p,q\right)$-adically Haar probability
measure.

To see this, note that (\ref{eq:van der Point series for one point function at 0})
is actually a van der Put series expression for $\mathbf{1}_{0}\left(\mathfrak{z}\right)$,
and the non-zero van der Put coefficients of that series are either
$1$ or $-1$. By \textbf{Example \ref{exa:end of 3.1.5. example}}
(see page \pageref{exa:end of 3.1.5. example}), since these coefficients
do not tend to zero, the integrals of the step-function-approximations
of $\mathbf{1}_{0}\left(\mathfrak{z}\right)$ do not converge $q$-adically
to a limit. For the case of a general one-point function, the translation
invariance of $d\mathfrak{z}$ allows us to extend the argument for
$\mathbf{1}_{0}$ to $\mathbf{1}_{\mathfrak{a}}$ for any $\mathfrak{a}\in\mathbb{Z}_{p}$. 
\end{example}
\vphantom{}

Seeing as we will \emph{define }the most general class of $\left(p,q\right)$-adically
Haar-measurable subsets of $\mathbb{Z}_{p}$ as being precisely those
sets whose indicator functions are the limit of a sequence of $\left(p,q\right)$-adically
integrable functions, the failure of one-point functions to be integrable
tells us that a set consisting of one point (or, more generally, any
set consisting of finitely many points) \emph{is not going to be $\left(p,q\right)$-adically
Haar integrable!} This is a \emph{fundamental} contrast with the real
case, where the vanishing measure of a single point tells us that
integrability is unaffected by changing functions' values at finitely
many (or even countably many) points. In $\left(p,q\right)$-adic
analysis, therefore, there are no notions of ``almost everywhere''
or ``null sets''. From a more measure-theoretic standpoint, this
also shows that it is too much to expect an arbitrary Borel set\index{Borel!set}
to be $\left(p,q\right)$-adically measurable. In light of this, Monna-Springer
integration instead begins by working with a slightly smaller, more
manageable algebra of sets: the compact clopen sets. 
\begin{defn}
We write $\Omega\left(X\right)$ to denote the collection of all compact
open subsets of $X$ (equivalently, compact clopen sets). We make
$\Omega\left(X\right)$ an algebra of sets by equipping it with the
operations of unions, intersections, and complements. 
\end{defn}
\vphantom{}

In more abstract treatments of non-archimedean integration (such as
\cite{van Rooij - Non-Archmedean Functional Analysis,Measure-theoretic approach to p-adic probability theory,van Rooij and Schikhof "Non-archimedean integration theory"}),
one begins by fixing a ring of sets $\mathcal{R}$ to use to construct
simple functions, which then give rise to a notion of measure. Here,
we have chosen $\mathcal{R}$ to be $\Omega\left(X\right)$. Taking
limits of simple functions in $L^{1}$ with respect to a given measure
$d\mu$ then enlarges $\mathcal{R}$ to an algebra of sets, denoted
$\mathcal{R}_{\mu}$, which is the maximal extension of $\mathcal{R}$
containing subsets of $X$ which admit a meaningful notion of $\mu$-measure.
This is, in essence, what we will do here, albeit in slightly more
concrete terms. For now, though, let us continue through the definitions.

Before proceeding, I should also mention that my exposition here is
a hybrid of the treatments given by Schikhof in one of the Appendices
of \cite{Ultrametric Calculus} and the chapter on Monna-Springer
integration theory given by Khrennikov in his book \cite{Quantum Paradoxes}
and in the paper \cite{Measure-theoretic approach to p-adic probability theory}. 
\begin{defn}[\textbf{Integrals and Measures}\footnote{Taken from \cite{Ultrametric Calculus}.}]
\ 

\vphantom{}

I. Elements of the dual space $C\left(X,K\right)^{\prime}$ are \index{integral!Monna-Springer}
called \textbf{integrals }on $C\left(X,K\right)$.

\vphantom{}

II. A \index{measure}\textbf{ measure}\index{measure!non-archimedean}
is a function $\mu:\Omega\left(X\right)\rightarrow K$ which satisfies
the conditions:

\vphantom{}

i. (Additivity): $\mu\left(U\cup V\right)=\mu\left(U\right)+\mu\left(V\right)$
for all $U,V\in\Omega\left(X\right)$ with $U\cap V=\varnothing$.

\vphantom{}

ii. (Boundedness): the real number $\left\Vert \mu\right\Vert $ (the\index{measure!total variation}
\textbf{total variation }of $\mu$) defined by: 
\begin{equation}
\left\Vert \mu\right\Vert \overset{\textrm{def}}{=}\sup\left\{ \left|\mu\left(U\right)\right|_{K}:U\in\Omega\left(X\right)\right\} \label{eq:Definition of the norm of a measure}
\end{equation}
is finite.

\vphantom{}

III. We write $M\left(X,K\right)$ to denote the vector space of $K$-valued
measures on $X$. Equipping this vector space with (\ref{eq:Definition of the norm of a measure})
makes it into a non-archimedean normed vector space. 
\end{defn}
\vphantom{}

Even though our starting algebra of measurable sets ($\Omega\left(X\right)$)
is smaller than its classical counterpart, there is still a correspondence
between measures and continuous linear functionals \cite{Ultrametric Calculus}. 
\begin{thm}
For each $\varphi\in C\left(X,K\right)^{\prime}$ (which sends a function
$f$ to the scalar $\varphi\left(f\right)$) the map $\mu_{\varphi}:\Omega\left(X\right)\rightarrow K$
defined by: 
\[
\mu_{\varphi}\left(U\right)\overset{\textrm{def}}{=}\varphi\left(\mathbf{1}_{U}\right),\textrm{ }\forall U\in\Omega\left(X\right)
\]
(where $\mathbf{1}_{U}$ is the indicator function of $U$) is a measure.
Moreover, the map: 
\[
\varphi\in C\left(X,K\right)^{\prime}\mapsto\mu_{\varphi}\in M\left(X,K\right)
\]
is a $K$-linear isometry of $C\left(X,K\right)^{\prime}$ onto $M\left(X,K\right)$.
We also have that: 
\begin{equation}
\left|\varphi\left(f\right)\right|_{K}\leq\left\Vert f\right\Vert _{X,K}\left\Vert \mu_{\varphi}\right\Vert ,\textrm{ }\forall f\in C\left(X,K\right),\forall\varphi\in C\left(X,K\right)^{\prime}\label{eq:Absolute Value of the output of a functional in terms of the associated measure}
\end{equation}
\end{thm}
\begin{notation}
In light of the above theorem, we shall now adopt the following notations
for $K$-valued measures and elements of $C\left(X,K\right)^{\prime}$:

\vphantom{}

I. Elements of $C\left(X,K\right)^{\prime}$ will all be written with
a $d$ in front of them: $d\mathfrak{z}$, $d\mu$, $d\mu\left(\mathfrak{z}\right)$,
$dA_{3}$, etc.

\vphantom{}

II. Given $d\mu\in C\left(X,K\right)^{\prime}$ and $f\in C\left(X,K\right)$,
we denote the image of $f$ under $d\mu$ by: 
\[
\int_{X}f\left(\mathfrak{z}\right)d\mu\left(\mathfrak{z}\right)
\]

\vphantom{}

III. Given $d\mu\in C\left(\mathbb{Z}_{p},\mathbb{C}_{q}\right)^{\prime}$
and $U\in\Omega\left(\mathbb{Z}_{p}\right)$, we will denote the $\mu$-measure
of $U$ by the symbols $\mu\left(U\right)$ and $\int_{U}d\mu$.

\vphantom{}

IV. Given any function $g\in C\left(X,K\right)$, the operation of
point-wise multiplication: 
\[
f\in C\left(X,K\right)\mapsto fg\in C\left(X,K\right)
\]
defines a continuous linear operator on $C\left(X,K\right)$. As such,
we can multiply elements of $C\left(X,K\right)^{\prime}$ by continuous
functions to obtain new elements of $C\left(X,K\right)^{\prime}$.
Given $d\mu\in C\left(X,K\right)^{\prime}$, we write $g\left(\mathfrak{z}\right)d\mu$
and $g\left(\mathfrak{z}\right)d\mu\left(\mathfrak{z}\right)$ to
denote the measure which first multiplies a function $f$ by $g$
and then integrates the product against $d\mu$: 
\begin{equation}
\int_{\mathbb{Z}_{p}}f\left(\mathfrak{z}\right)\left(g\left(\mathfrak{z}\right)d\mu\left(\mathfrak{z}\right)\right)\overset{\textrm{def}}{=}\int_{\mathbb{Z}_{p}}\left(f\left(\mathfrak{z}\right)g\left(\mathfrak{z}\right)\right)d\mu\left(\mathfrak{z}\right)\label{eq:Definition of the product of a measure by a continuous function}
\end{equation}

\vphantom{}

V. Given any $f\in C\left(X,K\right)$, any $d\mu\in C\left(X,K\right)^{\prime}$,
and any $U\in\Omega\left(X\right)$, we write: 
\[
\int_{U}f\left(\mathfrak{z}\right)d\mu\left(\mathfrak{z}\right)
\]
to denote: 
\[
\int_{\mathbb{Z}_{p}}\left(f\left(\mathfrak{z}\right)\mathbf{1}_{U}\left(\mathfrak{z}\right)\right)d\mu\left(\mathfrak{z}\right)
\]
\end{notation}
\vphantom{}

Like with the real case, there is also a Riemann-sum\index{Riemann sum}-type
formulation of the integrability of a function with respect to a given
measure \cite{Ultrametric Calculus}: 
\begin{prop}[\textbf{$\left(p,q\right)$-adic integration via Riemann sums}]
Let $d\mu\in C\left(\mathbb{Z}_{p},\mathbb{C}_{q}\right)^{\prime}$,
and let $f\in C\left(\mathbb{Z}_{p},\mathbb{C}_{q}\right)$. Then,
for any $\epsilon>0$, there exists a $\delta>0$, such that for each
partition of $\mathbb{Z}_{p}$ into pair-wise disjoint balls $B_{1},\ldots,B_{N}\in\Omega\left(\mathbb{Z}_{p}\right)$
, all of radius $<\delta$, and any choice of $\mathfrak{a}_{1}\in B_{1}$,
$\mathfrak{a}_{2}\in B_{2}$,$\ldots$: 
\begin{equation}
\left|\int_{\mathbb{Z}_{p}}f\left(\mathfrak{z}\right)d\mu\left(\mathfrak{z}\right)-\sum_{n=1}^{N}f\left(\mathfrak{a}_{n}\right)\mu\left(B_{j}\right)\right|_{q}<\epsilon\label{eq:Riemann sum criterion for (p,q)-adic integrability}
\end{equation}
\end{prop}
\begin{rem}
The pair-wise disjoint balls of radius $<\delta$ must be replaced
by pair-wise disjoint open sets of diameter $<\delta$ in the case
where $\mathbb{Z}_{p}$ is replaced by a compact ultrametric space
$X$ which does not possess the Heine-Borel property. If $X$ possesses
the Heine-Borel property, then the open balls may be used instead
of the more general open sets. 
\end{rem}
\begin{thm}[\textbf{Fubini's Theorem}\index{Fubini's Theorem} \cite{Ultrametric Calculus}]
Let $d\mu,d\nu$ be $K$-valued measures on $C\left(X,K\right)$
and $C\left(Y,K\right)$, respectively, where $X$ and $Y$ are compact
ultrametric spaces. Then, for any continuous $f:X\times Y\rightarrow K$:
\begin{equation}
\int_{X}\int_{Y}f\left(x,y\right)d\mu\left(x\right)d\nu\left(y\right)=\int_{X}\int_{Y}f\left(x,y\right)d\nu\left(y\right)d\mu\left(x\right)\label{eq:Fubini's Theorem}
\end{equation}
\end{thm}
\begin{example}[\textbf{A $\left(p,p\right)$-adic measure of great importance}\footnote{Given as part of an exercise on page 278 in the Appendix of \cite{Ultrametric Calculus}.}]

Let $\mathbb{C}_{p}^{+}$ denote the set: 
\begin{equation}
\mathbb{C}_{p}^{+}\overset{\textrm{def}}{=}\left\{ \mathfrak{s}\in\mathbb{C}_{p}:\left|1-\mathfrak{s}\right|_{p}<1\right\} 
\end{equation}
Schikhof calls such $\mathfrak{s}$ ``positive'' \cite{Ultrametric Calculus},
hence the superscript $+$. For any $\mathfrak{s}\in\mathbb{C}_{p}\backslash\mathbb{C}_{p}^{+}$,
define $d\mu_{\mathfrak{s}}\in M\left(\mathbb{Z}_{p},\mathbb{C}_{p}\right)$
by: 
\begin{equation}
\int_{\mathbb{Z}_{p}}\left[\mathfrak{z}\overset{p^{n}}{\equiv}k\right]d\mu_{\mathfrak{s}}\left(\mathfrak{z}\right)\overset{\mathbb{C}_{p}}{=}\frac{\mathfrak{s}^{k}}{1-\mathfrak{s}^{p^{n}}},\textrm{ }\forall n\in\mathbb{N}_{0},\textrm{ }\forall k\in\mathbb{Z}/p^{n}\mathbb{Z}
\end{equation}
This measure satisfies $\left\Vert \mu_{\mathfrak{s}}\right\Vert <1$,
and is used by Koblitz (among others) \cite{Koblitz's other book}
to define the \textbf{$p$-adic zeta function}\index{$p$-adic!Zeta function}.
A modification of this approach is one of the methods used to construct
more general $p$-adic $L$-functions\footnote{One of the major accomplishments of $p$-adic analysis ($\left(p,p\right)$,
not $\left(p,q\right)$) was to confirm that the measure-based approach
to constructing $p$-adic $L$-functions is, in fact, the same as
the $p$-adic $L$-functions produced by interpolating the likes of
Kummer's famous congruences for the Bernoulli numbers \cite{p-adic L-functions paper}.}\index{$p$-adic!$L$-function}\cite{p-adic L-functions paper}. Note
that because $d\mu_{\mathfrak{s}}$ is a $\left(p,p\right)$-adic
measure, it has no hope of being translation invariant, seeing as
the only translation-invariant $\left(p,p\right)$-adic linear functional
is the zero map. 
\end{example}
\vphantom{}

It is at this point, however, where we begin to step away from Subsections
\ref{subsec:3.1.4. The--adic-Fourier} and \ref{subsec:3.1.5-adic-Integration-=00003D000026}
and make our way toward abstraction. The standard construction of
a measure-theoretic integral is to first define it for the indicator
function of a measurable set, and then for finite linear combinations
thereof, and finally taking $\sup$s for the general case. Unfortunately,
this does not work well, both due to the shortage of measurable sets,
and\textemdash more gallingly\textemdash due to a terrible betrayal
by one of analysts' closest companions: the triangle inequality for
integrals. 
\begin{example}[\textbf{Failure of the triangle inequality for $\left(p,q\right)$-adic
integrals}]
\textbf{\label{exa:triangle inequality failure}}\index{triangle inequality!failure of}Letting
$d\mathfrak{z}$ denote the $\left(p,q\right)$-adic Haar measure
(which sends the indicator function $\left[\mathfrak{z}\overset{p^{n}}{\equiv}j\right]$
to the scalar $\frac{1}{p^{n}}\in\mathbb{Q}_{q}$), note that for
any $n\geq1$, any distinct $j,k\in\left\{ 0,\ldots,p^{n}-1\right\} $,
and any $\mathfrak{a},\mathfrak{b}\in\mathbb{Q}_{q}$: 
\begin{equation}
\left|\mathfrak{a}\left[\mathfrak{z}\overset{p^{n}}{\equiv}j\right]+\mathfrak{b}\left[\mathfrak{z}\overset{p^{n}}{\equiv}k\right]\right|_{q}\overset{\mathbb{R}}{=}\left|\mathfrak{a}\right|_{q}\left[\mathfrak{z}\overset{p^{n}}{\equiv}j\right]+\left|\mathfrak{b}\right|_{q}\left[\mathfrak{z}\overset{p^{n}}{\equiv}k\right],\textrm{ }\forall\mathfrak{z}\in\mathbb{Z}_{p}
\end{equation}
because, for any $\mathfrak{z}$, at most one of $\left[\mathfrak{z}\overset{p^{n}}{\equiv}j\right]$
and $\left[\mathfrak{z}\overset{p^{n}}{\equiv}k\right]$ can be equal
to $1$. So, letting: 
\begin{equation}
f\left(\mathfrak{z}\right)=\mathfrak{a}\left[\mathfrak{z}\overset{p^{n}}{\equiv}j\right]+\mathfrak{b}\left[\mathfrak{z}\overset{p^{n}}{\equiv}k\right]
\end{equation}
we have that: 
\begin{equation}
\left|\int_{\mathbb{Z}_{p}}f\left(\mathfrak{z}\right)d\mathfrak{z}\right|_{q}\overset{\mathbb{R}}{=}\left|\frac{\mathfrak{a}+\mathfrak{b}}{p^{n}}\right|_{q}
\end{equation}
and: 
\begin{equation}
\int_{\mathbb{Z}_{p}}\left|f\left(\mathfrak{z}\right)\right|_{q}d\mathfrak{z}\overset{\mathbb{R}}{=}\int_{\mathbb{Z}_{p}}\left(\left|\mathfrak{a}\right|_{q}\left[\mathfrak{z}\overset{p^{n}}{\equiv}j\right]+\left|\mathfrak{b}\right|_{q}\left[\mathfrak{z}\overset{p^{n}}{\equiv}k\right]\right)d\mathfrak{z}=\frac{\left|\mathfrak{a}\right|_{q}+\left|\mathfrak{b}\right|_{q}}{p^{n}}
\end{equation}
where the integral in the middle on the second line is with respect
to the standard \emph{real-valued} Haar probability measure on $\mathbb{Z}_{p}$.
If we choose $\mathfrak{a}=1$ and $\mathfrak{b}=q^{m}-1$, where
$m\in\mathbb{N}_{1}$, then: 
\begin{equation}
\left|\int_{\mathbb{Z}_{p}}f\left(\mathfrak{z}\right)d\mathfrak{z}\right|_{q}\overset{\mathbb{R}}{=}\left|\frac{q^{m}}{p^{n}}\right|_{q}=\frac{1}{q^{m}}
\end{equation}
\begin{equation}
\int_{\mathbb{Z}_{p}}\left|f\left(\mathfrak{z}\right)\right|_{q}d\mathfrak{z}\overset{\mathbb{R}}{=}\frac{\left|1\right|_{q}+\left|q^{m}-1\right|_{q}}{p^{n}}=\frac{2}{p^{n}}
\end{equation}
Because of this, the symbol $\sim$ in the expression: 
\begin{equation}
\left|\int_{\mathbb{Z}_{p}}f\left(\mathfrak{z}\right)d\mathfrak{z}\right|_{q}\sim\int_{\mathbb{Z}_{p}}\left|f\left(\mathfrak{z}\right)\right|_{q}d\mathfrak{z}
\end{equation}
can be arranged to be any of the symbols $\leq$, $\geq$, $<$, $>$,
or $=$ by making an appropriate choice of the values of $m$, $n$,
$p$, and $q$. A very distressing outcome, indeed.
\end{example}
\vphantom{}

The \emph{point} of Monna-Springer theory, one could argue, is to
find a work-around for this depressing state of affairs. That work-around
takes the form of an observation. What we would like to do is to write:
\begin{equation}
\left|\int f\right|\leq\int\left|f\right|
\end{equation}
However, we can't. Instead, we take out the middle man, which, in
this case, is the absolute value on the left-hand side. Rather than
view the triangle inequality as an operation we apply to an integral
already in our possession, we view it as an operation we apply directly
to the function $f$ we intend to integrate: 
\begin{equation}
f\mapsto\int\left|f\right|
\end{equation}
where, for the sake of discussion, we view everything as happening
on $\mathbb{R}$; $f:\mathbb{R}\rightarrow\mathbb{R}$, integration
is with respect to the Lebesgue measure, etc. Because the vanishing
of the integral only forces $f$ to be $0$ \emph{almost }everywhere,
the above map is not a true norm, but a semi-norm. Classically, remember,
$L^{1}$-norm only becomes a norm after we mod out by functions which
are zero almost everywhere.

More generally, given\textemdash say\textemdash a Radon measure $d\mu$,
the map $f\mapsto\int\left|f\right|d\mu$ defines a semi-norm which
we can use to compute the $\mu$-measure of a set $E$ by applying
it to continuous functions which approximate $E$'s indicator function.
While we will not be able to realize this semi-norm as the integral
of $\left|f\right|$ $d\mu$ in the non-archimedean case, we can still
obtain a semi-norm by working from the ground up. Given a measure
$\mu\in C\left(X,K\right)^{\prime}$, we consider the semi-norm:

\begin{equation}
f\in C\left(X,Y\right)\mapsto\sup_{g\in C\left(X,K\right)\backslash\left\{ 0\right\} }\frac{1}{\left\Vert g\right\Vert _{X,K}}\left|\int_{X}f\left(\mathfrak{z}\right)g\left(\mathfrak{z}\right)d\mu\left(\mathfrak{z}\right)\right|_{K}\in\mathbb{R}
\end{equation}
In the appendices of \cite{Ultrametric Calculus}, Schikhof takes
this formula as the definition for the semi-norm $N_{\mu}$ induced
by the measure $\mu$, which he then defines as a function on $\Omega\left(X\right)$
by defining $N_{\mu}\left(U\right)$ to be $N_{\mu}\left(\mathbf{1}_{U}\right)$,
where $\mathbf{1}_{U}$ is the indicator function of $U$. Taking
the infimum over all $U\in\Omega\left(X\right)$ containing a given
point $\mathfrak{z}$ then allows for $N_{\mu}$ to be extended to
a real-valued function on $X$ by way of the formula: 
\begin{equation}
\mathfrak{z}\in X\mapsto\inf_{U\in\Omega\left(X\right):\mathfrak{z}\in U}\sup_{g\in C\left(X,K\right)\backslash\left\{ 0\right\} }\frac{\left|\int_{X}\mathbf{1}_{U}\left(\mathfrak{y}\right)g\left(\mathfrak{y}\right)d\mu\left(\mathfrak{y}\right)\right|_{K}}{\left\Vert g\right\Vert _{X,K}}\in\mathbb{R}
\end{equation}

While this approach is arguably somewhat more intuitive, the above
definitions are not very pleasing to look at, and using them is a
hassle. Instead, I follow Khrennikov et. al. \cite{Measure-theoretic approach to p-adic probability theory}
in giving the following formulas as the definition of $N_{\mu}$: 
\begin{defn}
\label{def:Khrennikov N_mu definition}Let\footnote{Taken from \cite{Measure-theoretic approach to p-adic probability theory}.}
$d\mu\in C\left(X,K\right)^{\prime}$.

\vphantom{}

I. For any $U\in\Omega\left(X\right)$, we define $N_{\mu}\left(U\right)$
by: 
\begin{equation}
N_{\mu}\left(U\right)=\sup_{V\in\Omega\left(X\right):V\subset U}\left|\mu\left(V\right)\right|_{K},\textrm{ }\forall U\in\Omega\left(X\right)\label{eq:mu seminorm of a compact clopen set}
\end{equation}

\vphantom{}

II. For any $\mathfrak{z}\in X$, we define $N_{\mu}\left(\mathfrak{z}\right)$
by: 
\begin{equation}
N_{\mu}\left(\mathfrak{z}\right)\overset{\textrm{def}}{=}\inf_{U\in\Omega\left(X\right):\mathfrak{z}\in U}N_{\mu}\left(\mathbf{1}_{U}\right)\label{eq:Definition of N_mu at a point}
\end{equation}

\vphantom{}

III. Given any $f:X\rightarrow K$ (not necessarily continuous), we
define 
\begin{equation}
N_{\mu}\left(f\right)\overset{\textrm{def}}{=}\sup_{\mathfrak{z}\in X}N_{\mu}\left(\mathfrak{z}\right)\left|f\left(\mathfrak{z}\right)\right|_{K}\label{eq:Definition of N_mu for an arbitrary f}
\end{equation}
\end{defn}
\vphantom{}

In \cite{Ultrametric Calculus}, Schikhof proves that the formulas
of \textbf{Definition \ref{def:Khrennikov N_mu definition} }as a
theorem. Instead, following Khrennikov et. al. in \cite{Measure-theoretic approach to p-adic probability theory},
I state Schikhof's definitions of $N_{\mu}$ as a theorem (\textbf{Theorem
\ref{thm:Schikhof's N_mu definition}}), deducible from \textbf{Definition
\ref{def:Khrennikov N_mu definition}}. 
\begin{thm}
\label{thm:Schikhof's N_mu definition}\ 

\vphantom{}

I. For all $f\in C\left(X,K\right)$: 
\begin{equation}
N_{\mu}\left(f\right)=\sup_{g\in C\left(X,K\right)\backslash\left\{ 0\right\} }\frac{1}{\left\Vert g\right\Vert _{X,K}}\left|\int_{X}f\left(\mathfrak{z}\right)g\left(\mathfrak{z}\right)d\mu\left(\mathfrak{z}\right)\right|_{K}\label{eq:Schikhof's definition of N_mu}
\end{equation}

\vphantom{}

II. 
\begin{equation}
N_{\mu}\left(\mathfrak{z}\right)=\inf_{U\in\Omega\left(X\right):\mathfrak{z}\in U}\sup_{g\in C\left(X,K\right)\backslash\left\{ 0\right\} }\frac{\left|\int_{U}g\left(\mathfrak{y}\right)d\mu\left(\mathfrak{y}\right)\right|_{K}}{\left\Vert g\right\Vert _{X,K}}
\end{equation}
\end{thm}
\begin{rem}
The notation Schikhof uses in \cite{Ultrametric Calculus} for the
semi-norms is slightly different than the one given here, which I
have adopted from Khrennikov's exposition of the Monna-Springer integral
in \cite{Quantum Paradoxes}, which is simpler and more consistent. 
\end{rem}
\vphantom{}

In order to be able to extend the notion of $\mu$-integrable functions
from $C\left(X,K\right)$ to a large space, we need a non-archimedean
version of the triangle inequality for integrals. As we saw above,
however, this doesn't really exist in the non-archimedean context.
This is where $N_{\mu}$ comes in; it is our replacement for what
where $\int\left|f\right|d\mu$ would be in the archimedean case. 
\begin{prop}[\textbf{Monna-Springer Triangle Inequality} \cite{Ultrametric Calculus,Quantum Paradoxes}]
Let $X$ be a compact ultrametric space, and let $d\mu\in C\left(X,K\right)^{\prime}$.
Then: 
\begin{equation}
\left|\int_{X}f\left(\mathfrak{z}\right)d\mu\left(\mathfrak{z}\right)\right|_{K}\leq N_{\mu}\left(f\right)\leq\left\Vert \mu\right\Vert \left\Vert f\right\Vert _{X,K}\label{eq:Non-archimedean triangle inequality}
\end{equation}
\end{prop}
\begin{rem}
Unfortunately it is often the case\textemdash especially so for $\left(p,q\right)$-adic
analysis\textemdash that this estimate is little better than bounding
an integral by the product of the $\sup$ of the integrand and the
measure of the domain of integration. 
\end{rem}
\vphantom{}

Having done all this, the idea is to use $N_{\mu}$ to define what
it means for an arbitrary function $f:X\rightarrow K$ to be $\mu$-integrable.
This is done in the way one would expect: identifying $f,g\in C\left(X,K\right)$
whenever $N_{\mu}\left(f-g\right)=0$, and then modding $C\left(X,K\right)$
out by this equivalence relation, so as to upgrade $N_{\mu}$ from
a semi-norm to a fully-fledged norm. We then take the completion of
$C\left(X,K\right)$ with respect to this norm to obtain our space
of $\mu$-integrable functions. \cite{Quantum Paradoxes,Measure-theoretic approach to p-adic probability theory}
and \cite{Ultrametric Calculus} give us the definitions that will
guide us through this process. 
\begin{defn}
\ 

\vphantom{}

I. Because the indicator function for any set $S\subseteq X$ is then
a function $X\rightarrow K$, we can use (\ref{eq:Definition of N_mu for an arbitrary f})
to define $N_{\mu}$ for an arbitrary $S\subseteq X$: 
\begin{equation}
N_{\mu}\left(S\right)\overset{\textrm{def}}{=}\inf_{U\in\Omega\left(X\right):S\subseteq U}N_{\mu}\left(U\right)=\inf_{U\in\Omega\left(X\right):S\subseteq U}\sup_{V\in\Omega\left(X\right):V\subset U}\left|\mu\left(V\right)\right|_{K}\label{eq:Definition of N_mu of an arbitrary set}
\end{equation}

\vphantom{}

II.\textbf{ }A function $f:X\rightarrow K$ (not necessarily continuous)
is said to be \textbf{$\mu$-integrable }whenever there is a sequence
$\left\{ f_{n}\right\} _{n\geq1}$ in $C\left(X,K\right)$ so that
$\lim_{n\rightarrow\infty}N_{\mu}\left(f-f_{n}\right)=0$, where we
compute $N_{\mu}\left(f-f_{n}\right)$ using (\ref{eq:Definition of N_mu for an arbitrary f}).
We then write $\mathcal{L}_{\mu}^{1}\left(X,K\right)$ to denote the
vector space of all $\mu$-integrable functions from $X$ to $K$.

\vphantom{}

III. We say a set $S\subseteq X$ is \textbf{$\mu$-measurable }whenever
$N_{\mu}\left(S\right)$ is finite. Equivalently, $S$ is $\mu$-measurable
whenever there is a sequence $\left\{ f_{n}\right\} _{n\geq1}$ in
$C\left(X,K\right)$ so that $\lim_{n\rightarrow\infty}N_{\mu}\left(\mathbf{1}_{S}-f_{n}\right)=0$.
We write $\mathcal{R}_{\mu}$ to denote the collection of all $\mu$-measurable
subsets of $X$.

\vphantom{}

IV. A \textbf{non-archimedean measure space}\index{non-archimedean!measure space}\textbf{
(over $K$) }is a triple $\left(X,d\mu,\mathcal{R}_{\mu}\right)$,
where $X$ is an ultrametric space, $d\mu\in C\left(X,K\right)^{\prime}$,
and $\mathcal{R}_{\mu}$ is the collection of all $\mu$-measurable
subsets of $X$.

\vphantom{}

V. Given two non-archimedean measure spaces $\left(X_{1},d\mu_{1},\mathcal{R}_{\mu_{1}}\right)$
and $\left(X_{2},d\mu_{2},\mathcal{R}_{\mu_{2}}\right)$ with we say
a function $\phi:X_{1}\rightarrow X_{2}$ is \textbf{$\mu$-measurable}
whenever $\phi^{-1}\left(V\right)\in\mathcal{R}_{\mu_{1}}$ for all
$V\in\mathcal{R}_{\mu_{2}}$.

\vphantom{}

VI. We define $E_{\mu}$ as the set: 
\begin{equation}
E_{\mu}\overset{\textrm{def}}{=}\left\{ f:X\rightarrow K\textrm{ \& }N_{\mu}\left(f\right)<\infty\right\} \label{eq:Definition of E_mu}
\end{equation}
Given any $f\in E_{\mu}$, we say $f$ is \textbf{$\mu$-negligible
}whenever $N_{\mu}\left(f\right)=0$. Likewise, given any set $S\subseteq X$,
we say $S$ is \textbf{$\mu$-negligible} whenever $N_{\mu}\left(S\right)=0$.

\vphantom{}

VII. Given any real number $\alpha>0$, we define: 
\begin{equation}
X_{\mu:\alpha}\overset{\textrm{def}}{=}\left\{ \mathfrak{z}\in X:N_{\mu}\left(\mathfrak{z}\right)\geq\alpha\right\} \label{eq:Definition of X mu alpha}
\end{equation}
\begin{equation}
X_{\mu:+}\overset{\textrm{def}}{=}\bigcup_{\alpha>0}X_{\mu:\alpha}\label{eq:Definition of X mu plus}
\end{equation}
\begin{equation}
X_{\mu:0}\overset{\textrm{def}}{=}\left\{ \mathfrak{z}\in X:N_{\mu}\left(\mathfrak{z}\right)=0\right\} \label{eq:Definition of X mu 0}
\end{equation}
Note that $X=X_{\mu:0}\cup X_{\mu:+}$. We omit $\mu$ and just write
$X_{\alpha}$, $X_{+}$, and $X_{0}$ whenever there is no confusion
as to the measure $\mu$ we happen to be using.

\vphantom{}

VIII. We write $L_{\mu}^{1}\left(X,K\right)$ to denote the space
of equivalence classes of $\mathcal{L}_{\mu}^{1}\left(X,K\right)$
under the relation: 
\[
f\sim g\Leftrightarrow f-g\textrm{ is }\mu\textrm{-negligible}
\]
\end{defn}
\vphantom{}

Using these definitions, we can finally construct an honest-to-goodness
space of integrable non-archimedean functions \cite{Ultrametric Calculus}. 
\begin{thm}
Let $d\mu\in C\left(X,K\right)^{\prime}$. Then:

\vphantom{}

I. $\mathcal{L}_{\mu}^{1}\left(X,K\right)$ contains $C\left(X,K\right)$,
as well as all $\mu$-negligible functions from $X$ to $K$.

\vphantom{}

II. $d\mu$ can be extended from $C\left(X,K\right)$ to $\mathcal{L}_{\mu}^{1}\left(X,K\right)$;
denote this extension by $\overline{d\mu}$. Then: 
\[
\left|\int_{X}f\left(\mathfrak{z}\right)\overline{d\mu}\left(\mathfrak{z}\right)\right|_{K}\leq N_{\mu}\left(f\right),\textrm{ }\forall f\in\mathcal{L}_{\mu}^{1}\left(X,K\right)
\]

\vphantom{}

III. $L_{\mu}^{1}\left(X,K\right)$ is a Banach space over $K$ with
respect to the norm induced by $N_{\mu}$. 
\end{thm}
\vphantom{}

In $\left(p,q\right)$-adic case, however, these constructions turn
out to be overkill. 
\begin{thm}
Let $d\mu$ be the $\left(p,q\right)$-adic Haar probability measure
$d\mathfrak{z}$. Then $N_{\mu}\left(\mathfrak{z}\right)=1$ for all
$\mathfrak{z}\in\mathbb{Z}_{p}$. Equivalently, for $X=\mathbb{Z}_{p}$,
we have\footnote{This is given as an exercise on Page 281 of \cite{Ultrametric Calculus}.}:

\begin{align*}
X_{0} & =\varnothing\\
X_{+} & =X\\
X_{\alpha} & =\begin{cases}
\varnothing & \textrm{if }0<\alpha<1\\
X & \textrm{if }\alpha=1\\
\varnothing & \alpha>1
\end{cases}
\end{align*}
\end{thm}
\vphantom{}

Put in words, the only set which is negligible with respect to $d\mathfrak{z}$
is the \emph{empty set!} This is the reason why integrability and
continuity are synonymous in $\left(p,q\right)$-adic analysis. 
\begin{cor}[\textbf{The Fundamental Theorem of $\left(p,q\right)$-adic Analysis}\footnote{The name is my own.}]
Letting $d\mathfrak{z}$ denote the $\left(p,q\right)$-adic Haar
measure, we have the equality (not just an isomorphism, but an\emph{
}\textbf{equality}!): 
\begin{equation}
L_{d\mathfrak{z}}^{1}\left(\mathbb{Z}_{p},\mathbb{C}_{q}\right)=C\left(\mathbb{Z}_{p},\mathbb{C}_{q}\right)\label{eq:in (p,q)-adic analysis, integrable equals continuous}
\end{equation}
Consequently: 
\begin{equation}
\left|\int_{\mathbb{Z}_{p}}f\left(\mathfrak{z}\right)d\mathfrak{z}\right|_{q}\leq\left\Vert f\right\Vert _{p,q},\textrm{ }\forall f\in C\left(\mathbb{Z}_{p},\mathbb{C}_{q}\right)\label{eq:(p,q)-adic triangle inequality}
\end{equation}
\end{cor}
Proof: Since $N_{d\mathfrak{z}}=1$, $N_{d\mathfrak{z}}$ is then
the $\infty$-norm on $C\left(\mathbb{Z}_{p},\mathbb{C}_{q}\right)$.
Hence: 
\[
\mathcal{L}_{d\mathfrak{z}}^{1}\left(\mathbb{Z}_{p},\mathbb{C}_{q}\right)=L_{d\mathfrak{z}}^{1}\left(\mathbb{Z}_{p},\mathbb{C}_{q}\right)=C\left(\mathbb{Z}_{p},\mathbb{C}_{q}\right)
\]

Q.E.D.

\vphantom{}

With a little more work, one can prove many of the basic results familiar
to us from classical analysis. While these hold for Monna-Springer
theory in general, the simplicity of the $\left(p,q\right)$-adic
case makes these theorems significantly less useful than their classical
counterparts. I include them primarily for completeness' sake. 
\begin{thm}[\textbf{Non-Archimedean Dominated Convergence Theorem}\footnote{Given in Exercise 7F on page 263 of \cite{van Rooij - Non-Archmedean Functional Analysis}.}]
Let\index{Dominated Convergence Theorem} $X$ be a locally compact
ultrametric space. Let $g\in L_{\mu}^{1}\left(X,K\right)$ and let
$\left\{ f_{n}\right\} _{n\geq0}$ be a sequence in $L_{\mu}^{1}\left(X,K\right)$
such that:

I. For every set $E\in\mathcal{R}_{\mu}$, $f_{n}$ converges uniformly
to $f$ on $E$;

II. $\left|f_{n}\left(\mathfrak{z}\right)\right|_{K}\leq\left|g\left(\mathfrak{z}\right)\right|_{K}$
for all $\mathfrak{z}\in X$ and all $n$.

Then, $f\in L_{\mu}^{1}\left(X,K\right)$, $\lim_{n\rightarrow\infty}N_{\mu}\left(f_{n}-f\right)=0$
and: 
\begin{equation}
\lim_{n\rightarrow\infty}\int_{X}f_{n}\left(\mathfrak{z}\right)d\mu\left(\mathfrak{z}\right)\overset{K}{=}\int_{X}\lim_{n\rightarrow\infty}f_{n}\left(\mathfrak{z}\right)d\mu\left(\mathfrak{z}\right)\overset{K}{=}\int_{X}f\left(\mathfrak{z}\right)d\mu\left(\mathfrak{z}\right)\label{eq:Non-Archimedean Dominated Convergence Theorem}
\end{equation}
\end{thm}
\begin{rem}
For $\left(p,q\right)$-analysis, the Dominated Convergence Theorem
is simply the statement that: 
\begin{equation}
\lim_{n\rightarrow\infty}\int f_{n}=\int\lim_{n\rightarrow\infty}f_{n}
\end{equation}
occurs if and only if $f_{n}\rightarrow f$ uniformly on $\mathbb{Z}_{p}$:
\begin{equation}
\lim_{n\rightarrow\infty}\sup_{\mathfrak{z}\in\mathbb{Z}_{p}}\left|f\left(\mathfrak{z}\right)-f_{n}\left(\mathfrak{z}\right)\right|_{q}=0
\end{equation}
\end{rem}
\begin{thm}[\textbf{Non-Archimedean Hölder's Inequality} \cite{Quantum Paradoxes}]
Let $f,g\in C\left(X,K\right)$, and let $d\mu\in C\left(X,K\right)^{\prime}$
Then: 
\begin{equation}
N_{\mu}\left(fg\right)\leq N_{\mu}\left(f\right)\left\Vert g\right\Vert _{X,K}\label{eq:(p,q)-adic Holder inequality}
\end{equation}
\end{thm}
\begin{rem}
For\index{Hölder's Inequality} $\left(p,q\right)$-analysis, Hölder's
Inequality is extremely crude, being the statement that $\left|\int fg\right|_{q}\leq\left\Vert f\right\Vert _{X,K}\left\Vert g\right\Vert _{X,K}$. 
\end{rem}
\vphantom{}

Finally\textemdash although we shall not use it beyond the affine
case detailed in \textbf{Lemma \ref{lem:Affine substitution change of variable formula}},
we also have a general change of variables formula. First, however,
the necessary definitional lead-up. 
\begin{defn}
\label{def:pullback measure}Let $\left(X,d\mu,\mathcal{R}_{\mu}\right)$
and $\left(Y,d\nu,\mathcal{R}_{\nu}\right)$ be non-archimedean measure
spaces over $K$. For any measurable $\phi:X\rightarrow Y$, we define
the measure $d\mu_{\phi}:\mathcal{R}_{\nu}\rightarrow K$ by: 
\begin{equation}
\mu_{\phi}\left(V\right)\overset{\textrm{def}}{=}\mu\left(\phi^{-1}\left(V\right)\right),\textrm{ }\forall V\in\mathcal{R}_{\nu}\label{eq:Definition of change-of-variables measure}
\end{equation}
\end{defn}
\begin{lem}
Let $\left(X,d\mu,\mathcal{R}_{\mu}\right)$ and $\left(Y,d\nu,\mathcal{R}_{\nu}\right)$
be non-archimedean measure spaces over $K$. For any measurable $\phi:X\rightarrow Y$,
and for every $\mathcal{R}_{\nu}$ continuous $f:Y\rightarrow K$:
\begin{equation}
N_{\mu_{\phi}}\left(f\right)\leq N_{\mu}\left(f\circ\phi\right)
\end{equation}
\end{lem}
\begin{thm}[\textbf{Change of Variables \textendash{} Monna-Springer Integration
}\cite{Measure-theoretic approach to p-adic probability theory}]
Let\index{change of variables} $\left(X,d\mu,\mathcal{R}_{\mu}\right)$
and $\left(Y,d\nu,\mathcal{R}_{\nu}\right)$ be non-archimedean measure
spaces over $K$. For any measurable $\phi:X\rightarrow Y$, and for
every $\mathcal{R}_{\nu}$ continuous\footnote{As remarked in \cite{Measure-theoretic approach to p-adic probability theory},
it would be desirable if the hypothesis of continuity could be weakened.} $f:Y\rightarrow K$ so that $f\circ\phi\in L_{\mu}^{1}\left(X,K\right)$,
the function $f$ is in $L_{\nu}^{1}\left(Y,K\right)$ and: 
\begin{equation}
\int_{X}f\left(\phi\left(\mathfrak{z}\right)\right)d\mu\left(\mathfrak{z}\right)=\int_{Y}f\left(\mathfrak{y}\right)d\mu_{\phi}\left(y\right)\label{eq:Monna-Springer Integral Change-of-Variables Formula}
\end{equation}
\end{thm}
\newpage{}

\section{\label{sec:3.2 Rising-Continuous-Functions}Rising-Continuous Functions}

IN THIS SECTION, $p$ A PRIME AND $K$ IS A COMPLETE NON-ARCHIMEDEAN
VALUED FIELD.

\vphantom{}

Given the somewhat unusual method we used to construct $\chi_{H}$\textemdash strings
and all\textemdash you would think that the $\chi_{H}$s would be
unusual functions, located out in the left field of ``mainstream''
$\left(p,q\right)$-adic function theory. Surprisingly, this is not
the case. The key property here is the limit: 
\begin{equation}
\chi_{H}\left(\mathfrak{z}\right)\overset{\mathbb{Z}_{q_{H}}}{=}\lim_{n\rightarrow\infty}\chi_{H}\left(\left[\mathfrak{z}\right]_{p^{n}}\right)
\end{equation}
\textbf{Rising-continuous functions} are precisely those functions
which satisfy this limit. In this section, we shall get to know these
functions as a whole and see how they are easily obtained from the
smaller class of continuous functions.

\subsection{\label{subsec:3.2.1 -adic-Interpolation-of}$\left(p,q\right)$-adic
Interpolation of Functions on $\mathbb{N}_{0}$}

It is not an understatement to say that the entire theory of rising-continuous
functions emerges (or, should I say, ``\emph{rises}''?) out of the
van der Put identity ((\ref{eq:van der Put identity}) from \textbf{Proposition
\ref{prop:vdP identity} }on page \pageref{prop:vdP identity}): 
\begin{equation}
\sum_{n=0}^{\infty}c_{n}\left(\chi\right)\left[\mathfrak{z}\overset{p^{\lambda_{p}\left(n\right)}}{\equiv}n\right]\overset{K}{=}\lim_{k\rightarrow\infty}\chi\left(\left[\mathfrak{z}\right]_{p^{k}}\right),\textrm{ }\forall\mathfrak{z}\in\mathbb{Z}_{p}
\end{equation}
We make two observations: 
\begin{enumerate}
\item The field $K$ can be \emph{any }metrically complete valued field\textemdash archimedean
or non-archimedean. 
\item $\chi$ can be\emph{ any }function $\chi:\mathbb{Z}_{p}\rightarrow K$. 
\end{enumerate}
(1) will be arguably even more important than (2), because it highlights
what will quickly become a cornerstone of our approach: the ability
to allow the field (and hence, topology)\emph{ }used to sum our series
\emph{vary depending on $\mathfrak{z}$}. First, however, some terminology. 
\begin{defn}
Let $p$ be an integer $\geq2$.

\vphantom{}

I. Recall, we write \nomenclature{$\mathbb{Z}_{p}^{\prime}$}{$\overset{\textrm{def}}{=}\mathbb{Z}_{p}\backslash\left\{ 0,1,2,3,\ldots\right\}$ \nopageref}$\mathbb{Z}_{p}^{\prime}$
to denote the set: 
\begin{equation}
\mathbb{Z}_{p}^{\prime}\overset{\textrm{def}}{=}\mathbb{Z}_{p}\backslash\mathbb{N}_{0}\label{eq:Definition of Z_p prime}
\end{equation}

\vphantom{}

II. A sequence $\left\{ \mathfrak{z}_{n}\right\} _{n\geq0}$ in $\mathbb{Z}_{p}$
is said to be \textbf{($p$-adically) rising }if, as $n\rightarrow\infty$,
the number of non-zero $p$-adic digits in $\mathfrak{z}_{n}$ tends
to $\infty$. 
\end{defn}
\begin{rem}
$\mathbb{Z}_{p}^{\prime}$ is neither open nor closed, and has an
empty interior. 
\end{rem}
\begin{rem}
Note that for any rising $\left\{ \mathfrak{z}_{n}\right\} _{n\geq0}\subseteq\mathbb{Z}_{p}$,
the number of non-zero $p$-adic digits in $\mathfrak{z}_{n}$s tends
to $\infty$ as $n\rightarrow\infty$; that is, $\lim_{n\rightarrow\infty}\sum_{j=1}^{p-1}\#_{p:j}\left(\mathfrak{z}_{n}\right)\overset{\mathbb{R}}{=}\infty$.
However, the converse of this is not true: there are sequences $\left\{ \mathfrak{z}_{n}\right\} _{n\geq0}\subseteq\mathbb{Z}_{p}$
for which $\lim_{n\rightarrow\infty}\sum_{j=1}^{p-1}\#_{p:j}\left(\mathfrak{z}_{n}\right)\overset{\mathbb{R}}{=}\infty$
but $\lim_{n\rightarrow\infty}\mathfrak{z}_{n}\in\mathbb{N}_{0}$. 
\end{rem}
\begin{example}
For any $p$, let $\mathfrak{z}_{n}$'s sequence of $p$-adic digits
be $n$ consecutive $0$s followed by infinitely many $1$s; then
each $\mathfrak{z}_{n}$ has infinitely many non-zero $p$-adic digits,
but the $\mathfrak{z}_{n}$s converge $p$-adically to $0$. 
\end{example}
\vphantom{}

In ``classical'' non-archimedean analysis, for a function $\chi$,
$\chi$'s continuity is equivalent to the $K$-adic decay of its van
der Put coefficients to zero. Rising-continuous functions, in contrast,
arise when we relax this requirement of decay. Because we still need
for there to be enough regularity to the van der Put coefficients
of $\chi$ in order for its van der Put series to converge point-wise,
it turns out the ``correct'' definition for rising-continuous functions
is the limit condition we established for $\chi_{H}$ in our proof
of \textbf{Lemma \ref{lem:Unique rising continuation and p-adic functional equation of Chi_H}}
(page \pageref{lem:Unique rising continuation and p-adic functional equation of Chi_H}).
\begin{defn}
Let $K$ be a metrically complete non-archimedean valued field. A
function $\chi:\mathbb{Z}_{p}\rightarrow K$ is said to be \textbf{($\left(p,K\right)$-adically)}
\textbf{rising-continuous}\index{rising-continuous!function}\index{rising-continuous!left(p,Kright)-adically@$\left(p,K\right)$-adically}\textbf{
}whenever: 
\begin{equation}
\lim_{n\rightarrow\infty}\chi\left(\left[\mathfrak{z}\right]_{p^{n}}\right)\overset{K}{=}\chi\left(\mathfrak{z}\right),\textrm{ }\forall\mathfrak{z}\in\mathbb{Z}_{p}\label{eq:Definition of a rising-continuous function}
\end{equation}

We write $\tilde{C}\left(\mathbb{Z}_{p},K\right)$ to denote the $K$-linear
space\footnote{Note that this is a subspace of $B\left(\mathbb{Z}_{p},K\right)$.}
of all rising-continuous functions $\chi:\mathbb{Z}_{p}\rightarrow K$.
We call elements of \nomenclature{$\tilde{C}\left(\mathbb{Z}_{p},K\right)$}{set of $K$-valued rising-continuous functions on $\mathbb{Z}_{p}$}$\tilde{C}\left(\mathbb{Z}_{p},K\right)$
\textbf{rising-continuous functions}; or, more pedantically, \textbf{$\left(p,K\right)$-adic
rising-continuous functions}, or $\left(p,q\right)$\textbf{-adic
rising-continuous functions}, when $K$ is a $q$-adic field.
\end{defn}
\begin{rem}
The convergence in (\ref{eq:Definition of a rising-continuous function})
only needs to occur \emph{point-wise}.
\end{rem}
\begin{example}
Every continuous $f:\mathbb{Z}_{p}\rightarrow K$ is also rising-continuous,
but not every rising-continuous function is continuous. As an example,
consider: 
\begin{equation}
f\left(\mathfrak{z}\right)=q^{\textrm{digits}_{p}\left(\mathfrak{z}\right)}\label{eq:Example of a discontinuous rising-continuous function}
\end{equation}
where $\textrm{digits}_{p}:\mathbb{Z}_{p}\rightarrow\mathbb{N}_{0}\cup\left\{ +\infty\right\} $
outputs the number of non-zero $p$-adic digits in $\mathfrak{z}$.
In particular, $\textrm{digits}_{p}\left(\mathfrak{z}\right)=0$ if
and only if $\mathfrak{z}=0$, $\textrm{digits}_{p}\left(\mathfrak{z}\right)=+\infty$
if and only if $\mathfrak{z}\in\mathbb{Z}_{p}^{\prime}$. As defined,
$f\left(\mathfrak{z}\right)$ then satisfies: 
\begin{equation}
f\left(\mathfrak{z}\right)\overset{K}{=}\begin{cases}
q^{\textrm{digits}_{p}\left(\mathfrak{z}\right)} & \textrm{if }\mathfrak{z}\in\mathbb{N}_{0}\\
0 & \textrm{if }\mathfrak{z}\in\mathbb{Z}_{p}^{\prime}
\end{cases}
\end{equation}
Since $\mathbb{Z}_{p}^{\prime}$ is dense in $\mathbb{Z}_{p}$, if
$f$ \emph{were} continuous, the fact that it vanishes on $\mathbb{Z}_{p}$
would force it to be identically zero, which is not the case. Thus,
$f$ \emph{cannot} be continuous.
\end{example}
\begin{example}
Somewhat unfortunately, even if $\chi$ is rising continuous, we cannot
guarantee that $\chi\left(\mathfrak{z}_{n}\right)$ converges to $\chi\left(\mathfrak{z}\right)$
for every $p$-adically rising sequence $\left\{ \mathfrak{z}_{n}\right\} _{n\geq1}$
with limit $\mathfrak{z}_{n}\rightarrow\mathfrak{z}$. $\chi_{q}$\textemdash the
numen of the Shortened $qx+1$ map\textemdash gives us a simple example
of this. As we shall soon prove (see \textbf{Theorem \ref{thm:rising-continuability of Generic H-type functional equations}}
on page \pageref{thm:rising-continuability of Generic H-type functional equations}),
$\chi_{q}:\mathbb{Z}_{2}\rightarrow\mathbb{Z}_{q}$ is rising-continuous.
Now, consider the sequence: 
\begin{equation}
\mathfrak{z}_{n}\overset{\textrm{def}}{=}-2^{n}=\centerdot_{2}\underbrace{0\ldots0}_{n}\overline{1}\ldots
\end{equation}
where $\overline{1}$ indicates that all the remaining $2$-adic digits
of $\mathfrak{z}_{n}$ are $1$s; the first $n$ $2$-adic digits
of $\mathfrak{z}_{n}$ are $0$. This is a rising sequence, yet it
converges $2$-adically to $0$. Moreover, observe that:
\begin{equation}
\chi_{q}\left(-2^{n}\right)=\frac{1}{2^{n}}\chi_{q}\left(-1\right)=\frac{1}{2^{n}}\chi_{q}\left(B_{2}\left(1\right)\right)=\frac{1}{2^{n}}\frac{\chi_{q}\left(1\right)}{1-M_{q}\left(1\right)}=\frac{1}{2^{n}}\frac{1}{2-q}
\end{equation}
Thus, $\chi_{q}\left(-2^{n}\right)$ does not even converge $q$-adically
to a limit as $n\rightarrow\infty$, even though $\left\{ -2^{n}\right\} _{n\geq1}$
is a $2$-adically rising sequence which converges $2$-adically to
$0$.
\end{example}
\begin{prop}
\label{prop:vdP criterion for rising continuity}Let $\chi:\mathbb{Z}_{p}\rightarrow\mathbb{C}_{q}$
be any function. Then, the van der Put series of $\chi$ (that is,
$S_{p}\left\{ \chi\right\} $) converges at $\mathfrak{z}\in\mathbb{Z}_{p}$
if and only if: 
\begin{equation}
\lim_{n\rightarrow\infty}c_{\left[\mathfrak{z}\right]_{p^{n}}}\left(\chi\right)\left[\lambda_{p}\left(\left[\mathfrak{z}\right]_{p^{n}}\right)=n\right]\overset{\mathbb{C}_{q}}{=}0\label{eq:vdP criterion for rising-continuity}
\end{equation}
where $c_{\left[\mathfrak{z}\right]_{p^{n}}}\left(\chi\right)$ is
the $\left[\mathfrak{z}\right]_{p^{n}}$th van der Put coefficient
of $\chi$. 
\end{prop}
Proof: Using (\ref{eq:Inner term of vdP lambda decomposition}), we
can write the van der Put series for $\chi$ as: 
\begin{equation}
S_{p}\left\{ \chi\right\} \left(\mathfrak{z}\right)\overset{\mathbb{C}_{q}}{=}c_{0}\left(\chi\right)+\sum_{n=1}^{\infty}c_{\left[\mathfrak{z}\right]_{p^{n}}}\left(\chi\right)\left[\lambda_{p}\left(\left[\mathfrak{z}\right]_{p^{n}}\right)=n\right],\textrm{ }\forall\mathfrak{z}\in\mathbb{Z}_{p}\label{eq:Lambda-decomposed Chi}
\end{equation}
The ultrametric properties of $\mathbb{C}_{q}$ tell us that the $q$-adic
convergence of this series at any given $\mathfrak{z}\in\mathbb{Z}_{p}$
is equivalent to: 
\begin{equation}
\lim_{k\rightarrow\infty}c_{\left[\mathfrak{z}\right]_{p^{k}}}\left(\chi\right)\left[\lambda_{p}\left(\left[\mathfrak{z}\right]_{p^{k}}\right)=k\right]\overset{\mathbb{C}_{q}}{=}0\label{eq:vdP coefficient decay for lemma}
\end{equation}

Q.E.D. 
\begin{thm}
\label{thm:S_p}The operator $S_{p}$ which sends a function to its
formal van der Put series is a isomorphism from the $K$-linear space
$\tilde{C}\left(\mathbb{Z}_{p},K\right)$ onto the subspace of $\textrm{vdP}\left(\mathbb{Z}_{p},K\right)$
consisting of all formal van der Put series\index{van der Put!series}
which converge at every $\mathfrak{z}\in\mathbb{Z}_{p}$. In particular,
we have that for every $\chi\in\tilde{C}\left(\mathbb{Z}_{p},K\right)$:

\vphantom{}

I. $\chi=S_{p}\left\{ \chi\right\} $;

\vphantom{}

II. $\chi$ is uniquely represented by its van der Put series: 
\begin{equation}
\chi\left(\mathfrak{z}\right)\overset{K}{=}\sum_{n=0}^{\infty}c_{n}\left(\chi\right)\left[\mathfrak{z}\overset{p^{\lambda_{p}\left(n\right)}}{\equiv}n\right],\textrm{ }\forall\mathfrak{z}\in\mathbb{Z}_{p}\label{eq:Chi vdP series}
\end{equation}
where the convergence is point-wise. 
\end{thm}
Proof: Let $\chi\in\tilde{C}\left(\mathbb{Z}_{p},K\right)$ be arbitrary,
by the truncated van der Put identity (\ref{eq:truncated van der Put identity}),
$\chi\left(\left[\mathfrak{z}\right]_{p^{N}}\right)=S_{p:N}\left\{ \chi\right\} \left(\mathfrak{z}\right)$.
Thus, the rising-continuity of $\chi$ guarantees that $S_{p:N}\left\{ \chi\right\} \left(\mathfrak{z}\right)$
converges in $K$ to $\chi\left(\mathfrak{z}\right)$ as $N\rightarrow\infty$
for each $\mathfrak{z}\in\mathbb{Z}_{p}$. By (\ref{eq:vdP criterion for rising-continuity}),
this then implies that the van der Put coefficients of $\chi$ satisfy
the conditions of \textbf{Proposition \ref{prop:vdP criterion for rising continuity}}
for all $\mathfrak{z}\in\mathbb{Z}_{p}$, which then shows that the
van der Put series $S_{p}\left\{ \chi\right\} \left(\mathfrak{z}\right)$
converges at every $\mathfrak{z}\in\mathbb{Z}_{p}$, and\textemdash moreover\textemdash that
it converges to $\chi\left(\mathfrak{z}\right)$. This proves (I).

As for (II), the uniqueness specified therein is equivalent to demonstrating
that $S_{p}$ is an isomorphism in the manner described above. The
proof of this is like so: 
\begin{itemize}
\item (Surjectivity) Let $V\left(\mathfrak{z}\right)$ be any formal van
der Put series which converges $q$-adically at every $\mathfrak{z}\in\mathbb{Z}_{p}$.
Then, letting: 
\begin{equation}
\chi\left(\mathfrak{z}\right)\overset{\textrm{def}}{=}\lim_{N\rightarrow\infty}S_{p:N}\left\{ V\right\} \left(\mathfrak{z}\right)
\end{equation}
we have that $V\left(m\right)=\chi\left(m\right)$, and hence, $V\left(\mathfrak{z}\right)=S_{p}\left\{ \chi\right\} $.
Thus, $S_{p:N}\left\{ \chi\right\} \left(\mathfrak{z}\right)=\chi\left(\left[\mathfrak{z}\right]_{p^{N}}\right)$,
and so $\chi\left(\mathfrak{z}\right)\overset{\textrm{def}}{=}\lim_{N\rightarrow\infty}S_{p:N}\left\{ V\right\} \left(\mathfrak{z}\right)=\chi\left(\left[\mathfrak{z}\right]_{p^{N}}\right)$
shows that $\chi$ is rising-continuous. This shows that $V=S_{p}\left\{ \chi\right\} $,
and thus, that $S_{p}$ is surjective. 
\item (Injectivity) Let $\chi_{1},\chi_{2}\in\tilde{C}\left(\mathbb{Z}_{p},K\right)$
and suppose $S_{p}\left\{ \chi_{1}\right\} =S_{p}\left\{ \chi_{2}\right\} $.
Then, by (I): 
\begin{equation}
\chi_{1}\left(\mathfrak{z}\right)\overset{\textrm{(I)}}{=}S_{p}\left\{ \chi_{1}\right\} \left(\mathfrak{z}\right)=S_{p}\left\{ \chi_{2}\right\} \left(\mathfrak{z}\right)\overset{\textrm{(I)}}{=}\chi_{2}\left(\mathfrak{z}\right),\textrm{ }\forall\mathfrak{z}\in\mathbb{Z}_{p}
\end{equation}
which proves $\chi_{1}=\chi_{2}$, which proves the injectivity of
$S_{p}$. 
\end{itemize}
Thus, $S_{p}$ is an isomorphism.

Q.E.D.

\vphantom{}

While this shows that rising-continuous functions are equivalent to
those van der Put series which converge point-wise at every $\mathfrak{z}\in\mathbb{Z}_{p}$,
they are also more than that. As suggested by the name, rising-continuous
functions naturally occur when we wish to take a function on $\mathbb{N}_{0}$
and interpolate it to one on $\mathbb{Z}_{p}$. This will provide
us with an alternative characterization of rising-continuous functions
as being precisely those functions on $\mathbb{Z}_{p}$ whose restrictions
to $\mathbb{N}_{0}$ admit interpolations to functions on $\mathbb{Z}_{p}$.
I call the process underlying \textbf{rising-continuation}. 
\begin{defn}
\label{def:rising-continuation}Let $\mathbb{F}$ be $\mathbb{Q}$
or a field extension thereof, and let $\chi:\mathbb{N}_{0}\rightarrow\mathbb{F}$
be a function. Letting $p,q$ be integers $\geq2$, with $q$ prime,
we say $\chi$ has (or ``admits'') a \textbf{$\left(p,q\right)$-adic
}\index{rising-continuation}\textbf{rising-continuation} \textbf{(to
$K$)} whenever there is a metrically complete $q$-adic field extension
$K$ of $\mathbb{F}$ and a rising-continuous function $\chi^{\prime}:\mathbb{Z}_{p}\rightarrow K$
so that $\chi^{\prime}\left(n\right)=\chi\left(n\right)$ for all
$n\in\mathbb{N}_{0}$. We call any $\chi^{\prime}$ satisfying this
property a \textbf{($\left(p,q\right)$-adic} or \textbf{$\left(p,K\right)$-adic)}\emph{
}\textbf{rising-continuation }of $\chi$ (to $K$). 
\end{defn}
\begin{prop}
\label{prop:rising-continuation admission}Let $\chi:\mathbb{N}_{0}\rightarrow\mathbb{F}$
be a function admitting a $\left(p,q\right)$-adic rising-continuation
to $K$. Then:

\vphantom{}

I. The rising-continuation of $\chi$ is unique. As such, we will
write $\chi^{\prime}$\nomenclature{$\chi^{\prime}$}{the rising-contination of $\chi:\mathbb{N}_{0}\rightarrow\mathbb{F}$}
to denote the rising continuation of $\chi$.

\vphantom{}

II. 
\begin{equation}
\chi^{\prime}\left(\mathfrak{z}\right)\overset{K}{=}\lim_{j\rightarrow\infty}\chi\left(\left[\mathfrak{z}\right]_{p^{j}}\right)=\sum_{n=0}^{\infty}c_{n}\left(\chi\right)\left[\mathfrak{z}\overset{p^{\lambda_{p}\left(n\right)}}{\equiv}n\right],\textrm{ }\forall\mathfrak{z}\in\mathbb{Z}_{p}\label{eq:Rising continuation limit formula}
\end{equation}
\end{prop}
Proof:

I. Suppose $\chi$ admits two (possibly distinct) rising-continuations,
$\chi^{\prime}$ and $\chi^{\prime\prime}$. To see that $\chi^{\prime}$
and $\chi^{\prime\prime}$ must be the same, we note that since the
restrictions of both $\chi^{\prime}$ and $\chi^{\prime\prime}$ to
$\mathbb{N}_{0}$ are, by definition, equal to $\chi$, $\chi^{\prime}\left(n\right)=\chi^{\prime\prime}\left(n\right)$
for all $n\in\mathbb{N}_{0}$.

So, let $\mathfrak{z}$ be an arbitrary element of $\mathbb{Z}_{p}^{\prime}$;
necessarily, $\mathfrak{z}$ has infinitely many non-zero $p$-adic
digits. Consequently, $\left\{ \left[\mathfrak{z}\right]_{p^{j}}\right\} _{j\geq1}$
is a rising sequence of non-negative integers converging to $\mathfrak{z}$.
As such, using rising-continuity of $\chi^{\prime}$ and $\chi^{\prime\prime}$,
we can write: 
\[
\lim_{j\rightarrow\infty}\chi^{\prime}\left(\left[\mathfrak{z}\right]_{p^{j}}\right)\overset{K}{=}\chi^{\prime}\left(\mathfrak{z}\right)
\]
\[
\lim_{j\rightarrow\infty}\chi^{\prime\prime}\left(\left[\mathfrak{z}\right]_{p^{j}}\right)\overset{K}{=}\chi^{\prime\prime}\left(\mathfrak{z}\right)
\]
Since the $\left[\mathfrak{z}\right]_{p^{j}}$s are integers, this
yields: 
\[
\chi^{\prime}\left(\mathfrak{z}\right)=\lim_{j\rightarrow\infty}\chi^{\prime}\left(\left[\mathfrak{z}\right]_{p^{j}}\right)=\lim_{j\rightarrow\infty}\chi\left(\left[\mathfrak{z}\right]_{p^{j}}\right)=\lim_{j\rightarrow\infty}\chi^{\prime\prime}\left(\left[\mathfrak{z}\right]_{p^{j}}\right)=\chi^{\prime\prime}\left(\mathfrak{z}\right)
\]
Since $\mathfrak{z}$ was arbitrary, we have that $\chi^{\prime}\left(\mathfrak{z}\right)=\chi^{\prime\prime}\left(\mathfrak{z}\right)$
for all $\mathfrak{z}\in\mathbb{Z}_{p}^{\prime}$. Thus, $\chi^{\prime}$
and $\chi^{\prime\prime}$ are equal to one another on both $\mathbb{Z}_{p}^{\prime}$
and $\mathbb{N}_{0}$, and hence, on all of $\mathbb{Z}_{p}$. This
shows that $\chi^{\prime}$ and $\chi^{\prime\prime}$ are in fact
the same function, and therefore proves the uniqueness of $\chi$'s
rising-continuation.

\vphantom{}

II. As a rising-continuous function, $\chi^{\prime}\left(\mathfrak{z}\right)$
is uniquely determined by its values on $\mathbb{N}_{0}$. Since $\chi^{\prime}\mid_{\mathbb{N}_{0}}=\chi$,
$\chi^{\prime}$ and $\chi$ then have the same van der Put coefficients.
Using the van der Put identity (\ref{eq:van der Put identity}) then
gives (\ref{eq:Rising continuation limit formula}).

Q.E.D.

\vphantom{}

Now we have our characterization of rising-continuous functions in
terms of rising-continuability. 
\begin{thm}
\label{thm:characterization of rising-continuability}Let $\mathbb{F}$
be $\mathbb{Q}$ or a field extension thereof, let $\chi:\mathbb{N}_{0}\rightarrow\mathbb{F}$
be a function, let $p,q$ be integers $\geq2$, with $q$ prime, and
let $K$ be a metrically complete $q$-adic field extension of $\mathbb{F}$.
Then, the following are equivalent:

\vphantom{}

I. $\chi$ admits a $\left(p,q\right)$-adic rising-continuation to
$K$.

\vphantom{}

II. For each $\mathfrak{z}\in\mathbb{Z}_{p}$, $\chi\left(\left[\mathfrak{z}\right]_{p^{n}}\right)$
converges to a limit in $K$ as $n\rightarrow\infty$. 
\end{thm}
Proof:

i. Suppose (I) holds. Then, by \textbf{Proposition \ref{prop:rising-continuation admission}},
(\ref{eq:Rising continuation limit formula}) holds, which shows that
(II) is true.

\vphantom{}

ii. Conversely, suppose (II) holds. Then, by the van der Put identity
(\textbf{Proposition \ref{prop:vdP identity}}), $S_{p}\left\{ \chi\right\} $
is a rising-continuous function whose restriction to $\mathbb{N}_{0}$
is equal to $\chi$, which means that $S_{p}\left\{ \chi\right\} $
is the rising-continuation of $\chi$, and thus, that $\chi$ is rising-continuable.

So, (II) and (I) are equivalent.

Q.E.D.

\vphantom{}

As a consequence of this, we then have that $K$-linear space (under
point-wise addition) of all $\chi:\mathbb{N}_{0}\rightarrow\mathbb{F}$
which admit $\left(p,q\right)$-adic rising-continuations to $K$
is then isomorphic to $\tilde{C}\left(\mathbb{Z}_{p},K\right)$, with
the isomorphism being the act of rising-continuation: $\chi\mapsto\chi^{\prime}$.
As such, \textbf{\emph{from now on (for the most part) we will identify
$\chi^{\prime}$ and $\chi=\chi^{\prime}\mid_{\mathbb{N}_{0}}$ and
treat them as one and the same.}} The exceptions to this convention
will be those occasions where we have a $\left(p,q\right)$-adically
continuable function $\chi:\mathbb{N}_{0}\rightarrow\mathbb{F}\subseteq K$
and a function $\eta:\mathbb{Z}_{p}\rightarrow K$, and we wish to
show that $\mu$ is in fact equal $\chi^{\prime}$. Using $\chi^{\prime}$
in this context will then allow us to distinguish between $\mu$ (our
candidate for $\chi^{\prime}$) and the actual ``correct'' rising-continuation
of $\chi$.

Because many of our rising-continuous functions will emerge as interpolations
of solutions of systems of functional equations on $\mathbb{N}_{0}$,
our next theorem will be quite the time-saver. 
\begin{thm}
\label{thm:rising-continuability of Generic H-type functional equations}Let
$H$ be a semi-basic $p$-Hydra map which fixes $0$, and consider
the system of functional equations\index{functional equation}: 
\begin{equation}
f\left(pn+j\right)=\frac{\mu_{j}}{p}f\left(n\right)+c_{j},\textrm{ }\forall j\in\left\{ 0,\ldots,p-1\right\} ,\textrm{ }\forall n\in\mathbb{N}_{0}\label{eq:Generic H-type functional equations}
\end{equation}
where $\mu_{j}/p=H_{j}^{\prime}\left(0\right)$. Then:

\vphantom{}

I. There is a unique function $\chi:\mathbb{N}_{0}\rightarrow\overline{\mathbb{Q}}$
such that $f=\chi$ is a solution of \emph{(\ref{eq:Generic H-type functional equations})}.

\vphantom{}

II. The solution $\chi$ \emph{(\ref{eq:Generic H-type functional equations})}
is rising-continuable\index{rising-continuation} to a function $\chi:\mathbb{Z}_{p}\rightarrow\mathbb{C}_{q_{H}}$.
Moreover, this continuation satisfies: 
\begin{equation}
\chi\left(p\mathfrak{z}+j\right)=\frac{\mu_{j}}{p}\chi\left(\mathfrak{z}\right)+c_{j},\textrm{ }\forall j\in\left\{ 0,\ldots,p-1\right\} ,\textrm{ }\forall\mathfrak{z}\in\mathbb{Z}_{p}\label{eq:Rising-continuation Generic H-type functional equations}
\end{equation}

\vphantom{}

III. The function $\chi:\mathbb{Z}_{p}\rightarrow\mathbb{C}_{q_{H}}$
described in \emph{(III)} is the unique rising-continuous function
$\mathbb{Z}_{p}\rightarrow\mathbb{C}_{q_{H}}$ satisfying \emph{(\ref{eq:Rising-continuation Generic H-type functional equations})}. 
\end{thm}
Proof:

I. Let $f:\mathbb{N}_{0}\rightarrow\overline{\mathbb{Q}}$ be any
solution of (\ref{eq:Generic H-type functional equations}). Setting
$n=j=0$ yields: 
\begin{align*}
f\left(0\right) & =\frac{\mu_{0}}{p}f\left(0\right)+c_{0}\\
 & \Updownarrow\\
f\left(0\right) & =\frac{c_{0}}{1-\frac{\mu_{0}}{p}}
\end{align*}
Since $H$ is semi-basic, $\mu_{0}/p\neq1$, so $f\left(0\right)$
is well-defined. Then, we have that: 
\begin{equation}
f\left(j\right)=\frac{\mu_{j}}{p}f\left(0\right)+c_{j},\textrm{ }\forall j\in\left\{ 0,\ldots,p-1\right\} 
\end{equation}
and, more generally: 
\begin{equation}
f\left(m\right)=\frac{\mu_{\left[m\right]_{p}}}{p}f\left(\frac{m-\left[m\right]_{p}}{p}\right)+c_{\left[m\right]_{p}},\textrm{ }\forall m\in\mathbb{N}_{0}\label{eq:f,m, digit shifting}
\end{equation}
Since the map $m\mapsto\frac{m-\left[m\right]_{p}}{p}$ sends the
non-negative integer: 
\[
m=m_{0}+m_{1}p+\cdots+m_{\lambda_{p}\left(m\right)-1}p^{\lambda_{p}\left(m\right)-1}
\]
to the integer: 
\[
m=m_{1}+m_{2}p+\cdots+m_{\lambda_{p}\left(m\right)-1}p^{\lambda_{p}\left(m\right)-2}
\]
it follows that $m\mapsto\frac{m-\left[m\right]_{p}}{p}$ eventually
iterates every $m\in\mathbb{N}_{0}$ to $0$, and hence (\ref{eq:f,m, digit shifting})
implies that, for every $m\in\mathbb{N}_{0}$, $f\left(m\right)$
is entirely determined by $f\left(0\right)$ and the $c_{j}$s. Since
$f\left(0\right)$ is uniquely determined by $c_{0}$ and $H$, the
equation (\ref{eq:Generic H-type functional equations}) possesses
exactly one solution. Let us denote this solution by $\chi$.

\vphantom{}

II. Because $H$ is semi-basic, we can repeat for $\chi$ the argument
we used to prove the existence of $\chi_{H}$'s rising-continuation
in \textbf{Lemma \ref{lem:Unique rising continuation and p-adic functional equation of Chi_H}}
(page \pageref{lem:Unique rising continuation and p-adic functional equation of Chi_H}):
since any $\mathfrak{z}\in\mathbb{Z}_{p}^{\prime}$ has infinitely
many non-zero $p$-adic digits, the product of $\mu_{j}/p$ taken
over all the digits of such a $\mathfrak{z}$ will converge $q_{H}$-adically
to zero. Then, seeing as (\ref{eq:Generic H-type functional equations})
shows that $\chi\left(\mathfrak{z}\right)$ will be a sum of the form:
\begin{equation}
\beta_{0}+\alpha_{1}\beta_{1}+\alpha_{1}\alpha_{2}\beta_{2}+\alpha_{1}\alpha_{2}\alpha_{3}\beta_{3}+\cdots
\end{equation}
(where for each $n$, $\beta_{n}$ is one of the $c_{j}$s and $\alpha_{n}$
is $\mu_{j_{n}}/p$, where $j_{n}$ is the $n$th $p$-adic digit
of $\mathfrak{z}$) this sum will converge in $\mathbb{C}_{q_{H}}$
for all $\mathfrak{z}\in\mathbb{Z}_{p}^{\prime}$ because of the estimate
guaranteed by $H$'s semi-basicness: 
\begin{equation}
\left|\frac{\mu_{j}}{p}\right|_{q_{H}}<1,\textrm{ }\forall j\in\left\{ 1,\ldots,p-1\right\} 
\end{equation}
This proves $\chi$ is rising-continuable.

\vphantom{}

III. Because $\chi$ admits a rising continuation, its continuation
at any $\mathfrak{z}\in\mathbb{Z}_{p}$ is given by the value of $\chi$'s
van der Put series at $\mathfrak{z}$. As such: 
\begin{equation}
\chi\left(\mathfrak{z}\right)\overset{\mathbb{C}_{q_{H}}}{=}S_{p}\left\{ \chi\right\} \left(\mathfrak{z}\right)\overset{\mathbb{C}_{q_{H}}}{=}\lim_{N\rightarrow\infty}\chi\left(\left[\mathfrak{z}\right]_{p^{N}}\right)
\end{equation}
Since $\chi$ satisfies the functional equations (\ref{eq:Generic H-type functional equations}),
we can write: 
\begin{equation}
\chi\left(p\left[\mathfrak{z}\right]_{p^{N}}+j\right)\overset{\overline{\mathbb{Q}}}{=}\frac{\mu_{j}}{p}\chi\left(\left[\mathfrak{z}\right]_{p^{N}}\right)+c_{j}
\end{equation}
for all $\mathfrak{z}\in\mathbb{Z}_{p}$, all $j\in\left\{ 0,\ldots,p-1\right\} $,
and all $N\geq0$. Letting $\mathfrak{z}\in\mathbb{Z}_{p}^{\prime}$
be arbitrary, we let $N\rightarrow\infty$. This gives us: 
\begin{equation}
\chi\left(p\mathfrak{z}+j\right)\overset{\mathbb{C}_{q_{H}}}{=}\frac{\mu_{j}}{p}\chi\left(\mathfrak{z}\right)+c_{j},\textrm{ }\forall j,\textrm{ }\forall\mathfrak{z}\in\mathbb{Z}_{p}^{\prime}
\end{equation}
Note that these identities automatically hold for $\mathfrak{z}\in\mathbb{N}_{0}$
seeing as those particular cases were governed by (\ref{eq:Generic H-type functional equations}).
The uniqueness of $\chi$ as a solution of (\ref{eq:Rising-continuation Generic H-type functional equations})
follows from the fact that any $f:\mathbb{Z}_{p}\rightarrow\mathbb{C}_{q_{H}}$
satisfying (\ref{eq:Rising-continuation Generic H-type functional equations})
has a restriction to $\mathbb{N}_{0}$ which satisfies (\ref{eq:Generic H-type functional equations}).
This then forces $f=\chi$.

Q.E.D.

\vphantom{}

A slightly more general version of this type of argument (which will
particularly useful in Chapters 4 and 6) is as follows: 
\begin{lem}
\label{lem:rising-continuations preserve functional equations}Fix
integers $p,q\geq2$, let $K$ be a metrically complete $q$-adic
field, and let $\Phi_{j}:\mathbb{Z}_{p}\times K\rightarrow K$ be
continuous for $j\in\left\{ 0,\ldots,p-1\right\} $. If $\chi:\mathbb{N}_{0}\rightarrow K$
is rising-continuable to an element $\chi\in\tilde{C}\left(\mathbb{Z}_{p},K\right)$,
and if: 
\[
\chi\left(pn+j\right)=\Phi_{j}\left(n,\chi\left(n\right)\right),\textrm{ }\forall n\in\mathbb{N}_{0},\textrm{ }\forall j\in\left\{ 0,\ldots,p-1\right\} 
\]
then: 
\[
\chi\left(p\mathfrak{z}+j\right)=\Phi_{j}\left(\mathfrak{z},\chi\left(\mathfrak{z}\right)\right),\textrm{ }\forall\mathfrak{z}\in\mathbb{Z}_{p},\textrm{ }\forall j\in\left\{ 0,\ldots,p-1\right\} 
\]
\end{lem}
Proof: Let everything be as given. Then, since $\chi$ is rising continuous,
we have that: 
\[
\chi\left(\mathfrak{z}\right)\overset{K}{=}\lim_{N\rightarrow\infty}\chi\left(\left[\mathfrak{z}\right]_{p^{N}}\right),\textrm{ }\forall\mathfrak{z}\in\mathbb{Z}_{p}
\]
Since: 
\[
\chi\left(p\left[\mathfrak{z}\right]_{p^{N}}+j\right)=\Phi_{j}\left(\left[\mathfrak{z}\right]_{p^{N}},\chi\left(\left[\mathfrak{z}\right]_{p^{N}}\right)\right)
\]
holds true for all $\mathfrak{z}\in\mathbb{Z}_{p}$ and all $N\geq0$,
the rising-continuity of $\chi$ and the continuity of $\Phi_{j}$
then guarantee that: 
\begin{align*}
\chi\left(p\mathfrak{z}+j\right) & \overset{K}{=}\lim_{N\rightarrow\infty}\chi\left(p\left[\mathfrak{z}\right]_{p^{N}}+j\right)\\
 & \overset{K}{=}\lim_{N\rightarrow\infty}\Phi_{j}\left(\left[\mathfrak{z}\right]_{p^{N}},\chi\left(\left[\mathfrak{z}\right]_{p^{N}}\right)\right)\\
 & \overset{K}{=}\Phi_{j}\left(\mathfrak{z},\chi\left(\mathfrak{z}\right)\right)
\end{align*}
as desired.

Q.E.D.

\subsection{\label{subsec:3.2.2 Truncations-=00003D000026-The}Truncations and
The Banach Algebra of Rising-Continuous Functions}

RECALL THAT WE WRITE $\left\Vert \cdot\right\Vert _{p,q}$ TO DENOTE
THE $\left(p,q\right)$-adic SUPREMUM NORM.

\vphantom{}

In this section, we will examine the structure of the vector space
$\tilde{C}\left(\mathbb{Z}_{p},K\right)$ of rising-continuous functions
$\chi:\mathbb{Z}_{p}\rightarrow K$. We will demonstrate that $\tilde{C}\left(\mathbb{Z}_{p},K\right)$
is a Banach algebra over $K$\textemdash with the usual ``point-wise
multiplication of functions'' as its multiplication operation. Not
only that, we will also see that $\tilde{C}\left(\mathbb{Z}_{p},K\right)$
extends $C\left(\mathbb{Z}_{p},K\right)$, containing it as a proper
sub-algebra. Before we can even begin our discussion of Banach algebras,
however, we need to introduce a simple construction which will be
of extraordinary importance in our analyses of $\chi_{H}$ in Chapter
4 and beyond.

To begin, as we saw in Subsection \ref{subsec:3.1.4. The--adic-Fourier},
the Fourier Transform of a continuous $\left(p,q\right)$-adic function
$f$ can be computed using $f$'s van der Put coefficients via equation
(\ref{eq:Definition of (p,q)-adic Fourier Coefficients}) from \textbf{Definition
\ref{def:pq adic Fourier coefficients}} 
\begin{equation}
\hat{f}\left(t\right)\overset{\mathbb{C}_{q}}{=}\sum_{n=\frac{1}{p}\left|t\right|_{p}}^{\infty}\frac{c_{n}\left(f\right)}{p^{\lambda_{p}\left(n\right)}}e^{-2n\pi it},\textrm{ }\forall t\in\hat{\mathbb{Z}}_{p}
\end{equation}
Because the convergence of this infinite series \emph{requires }$\lim_{n\rightarrow\infty}\left|c_{n}\left(f\right)\right|_{q}=0$,
this formula for $\hat{f}$ is not compatible with general rising-continuous
functions; the van der Put coefficients need not converge $q$-adically
to $0$. As will be shown in Chapter 4, we can get around the non-convergence
of (\ref{eq:Definition of (p,q)-adic Fourier Coefficients}) for an
arbitrary rising-continuous function $\chi$ by replacing $\chi$
with a locally constant approximation. I call these approximations
\textbf{truncations}, and\textemdash as locally constant functions\textemdash they
have the highly desirable property of $\left(p,q\right)$-adic continuity. 
\begin{defn}[\textbf{$N$th truncations}]
\label{def:Nth truncation}For any $\chi\in B\left(\mathbb{Z}_{p},\mathbb{C}_{q}\right)$
and any $N\in\mathbb{N}_{0}$, the \index{$N$th truncation}\textbf{$N$th
truncation of $\chi$}, denoted $\chi_{N}$, is the function $\chi_{N}:\mathbb{Z}_{p}\rightarrow\mathbb{C}_{q}$
defined by: 
\begin{equation}
\chi_{N}\left(\mathfrak{z}\right)\overset{\textrm{def}}{=}\sum_{n=0}^{p^{N}-1}\chi\left(n\right)\left[\mathfrak{z}\overset{p^{N}}{\equiv}n\right],\textrm{ }\forall\mathfrak{z}\in\mathbb{Z}_{p}\label{eq:Definition of Nth truncation}
\end{equation}
We also extend this notation to negative $N$ by defining $\chi_{N}\left(\mathfrak{z}\right)$
to be identically zero whenever $N<0$. 
\end{defn}
\begin{rem}
For any $N\in\mathbb{N}_{0}$: 
\begin{equation}
\chi_{N}\left(\mathfrak{z}\right)=\chi\left(\left[\mathfrak{z}\right]_{p^{N}}\right),\textrm{ }\forall\mathfrak{z}\in\mathbb{Z}_{p}\label{eq:Nth truncation of Chi in terms of Chi}
\end{equation}
Combining (\ref{eq:Nth truncation of Chi in terms of Chi}) and (\ref{eq:truncated van der Put identity}),
we then have that: 
\begin{equation}
\chi_{N}\left(\mathfrak{z}\right)=\sum_{n=0}^{p^{N}-1}\chi\left(n\right)\left[\mathfrak{z}\overset{p^{N}}{\equiv}n\right]=\sum_{n=0}^{p^{N}-1}c_{n}\left(\chi\right)\left[\mathfrak{z}\overset{p^{\lambda_{p}\left(n\right)}}{\equiv}n\right]=\chi\left(\left[\mathfrak{z}\right]_{p^{N}}\right)\label{eq:Nth truncation and truncated van-der-Put identity compared}
\end{equation}
Moreover, note that $\chi_{N}\in C\left(\mathbb{Z}_{p},\mathbb{C}_{q}\right)$
for all $N\in\mathbb{N}_{0}$. 
\end{rem}
\begin{rem}
Although Schikhof did not make any significant investigations with
truncations, he was aware of there existence. There is an exercise
in \cite{Ultrametric Calculus} in the section on the van der Put
basis which asks the reader to show that the $N$th truncation of
$\chi$ is the best possible $\left(p,q\right)$-adic approximation
of $\chi$ which is constant with respect to inputs modulo $p^{N}$. 
\end{rem}
\begin{prop}
\label{prop:Truncations converge pointwise iff rising-continuous}Let
$\chi:\mathbb{Z}_{p}\rightarrow K$ be\textbf{ any}\emph{ }function.
Then, $\chi\in\tilde{C}\left(\mathbb{Z}_{p},K\right)$ if and only
if $\chi_{N}$ converges in $K$ to $\chi$ point-wise on $\mathbb{Z}_{p}$
as $N\rightarrow\infty$.
\end{prop}
Proof: Since $\chi_{N}\left(\mathfrak{z}\right)=\chi\left(\left[\mathfrak{z}\right]_{p^{N}}\right)$,
this proposition is just a restatement of the definition of what it
means for $\chi$ to be rising continuous ($\chi\left(\left[\mathfrak{z}\right]_{p^{N}}\right)$
converges to $\chi\left(\mathfrak{z}\right)$ point-wise everywhere).

Q.E.D. 
\begin{prop}
\label{prop:Unif. convergence of truncation equals continuity}Let
$f:\mathbb{Z}_{p}\rightarrow K$ be \textbf{any} function. Then, $f$
is continuous if and only if $f_{N}$ converges in $K$ to $f$ uniformly
on $\mathbb{Z}_{p}$ as $N\rightarrow\infty$. 
\end{prop}
Proof:

I. If the $f_{N}$s converge uniformly to $f$, $f$ is continuous,
seeing as the $f_{N}$ are locally constant\textemdash and hence,
continuous\textemdash and seeing as how the uniform limit of continuous
functions is continuous.

\vphantom{}

II. Conversely, if $f$ is continuous then, since $\mathbb{Z}_{p}$
is compact, $f$ is uniformly continuous. So, letting $\epsilon>0$,
pick $\delta$ so that $\left|f\left(\mathfrak{z}\right)-f\left(\mathfrak{y}\right)\right|_{q}<\epsilon$
for all $\mathfrak{z},\mathfrak{y}\in\mathbb{Z}_{p}$ with $\left|\mathfrak{z}-\mathfrak{y}\right|_{p}<\delta$.
Note that: 
\[
\left|\mathfrak{z}-\left[\mathfrak{z}\right]_{p^{N}}\right|_{p}\leq p^{-N},\textrm{ }\forall\mathfrak{z}\in\mathbb{Z}_{p},\textrm{ }\forall N\in\mathbb{N}_{0}
\]
So, choose $N$ large enough so that $p^{-N}<\delta$. Then $\left|\mathfrak{z}-\left[\mathfrak{z}\right]_{p^{N}}\right|_{p}<\delta$
for all $\mathfrak{z}\in\mathbb{Z}_{p}$, and so: 
\[
\sup_{\mathfrak{z}\in\mathbb{Z}_{p}}\left|f\left(\mathfrak{z}\right)-f_{N}\left(\mathfrak{z}\right)\right|_{q}=\sup_{\mathfrak{z}\in\mathbb{Z}_{p}}\left|f\left(\mathfrak{z}\right)-f\left(\left[\mathfrak{z}\right]_{p^{N}}\right)\right|_{q}<\epsilon
\]
which proves that the $f_{N}$s converge uniformly to $f$.

Q.E.D. 
\begin{rem}
As is proved in \textbf{Theorem \ref{thm:vdP basis theorem}} (see
(\ref{eq:Fourier transform of Nth truncation in terms of vdP coefficients})
on page \pageref{eq:Fourier transform of Nth truncation in terms of vdP coefficients}),
the formula (\ref{eq:Definition of (p,q)-adic Fourier Coefficients})
can be applied to: 
\begin{equation}
f_{N}\left(\mathfrak{z}\right)=f\left(\left[\mathfrak{z}\right]_{p^{N}}\right)=\sum_{n=0}^{p^{N}-1}c_{n}\left(f\right)\left[\mathfrak{z}\overset{p^{\lambda_{p}\left(n\right)}}{\equiv}n\right]
\end{equation}
whereupon we obtain: 
\begin{equation}
\hat{f}_{N}\left(t\right)=\sum_{n=\frac{\left|t\right|_{p}}{p}}^{p^{N}-1}\frac{c_{n}\left(f\right)}{p^{\lambda_{p}\left(n\right)}}e^{-2\pi int}
\end{equation}
where $\hat{f}_{N}$ is the Fourier transform of the $N$th truncation
of $f$. While this formula for $\hat{f}_{N}$ can be used to perform
the kind of Fourier analysis we shall do in Chapters 4 and 6, it turns
out to be much easier to work with the $N$th truncation directly
instead of using the van der Put coefficients in this manner, primarily
because the formula for $c_{n}\left(\chi_{H}\right)$ becomes unwieldy\footnote{See\textbf{ Proposition \ref{prop:van der Put series for Chi_H}}
on page \pageref{prop:van der Put series for Chi_H} for the details.} when $H$ is a $p$-Hydra map for $p\geq3$. Nevertheless, the van
der Put coefficient expression for $\hat{f}_{N}\left(t\right)$ is
of interest in its own right, and it may be worth investigating how
the above formula might be inverted so as to produce expressions for
$c_{n}\left(f\right)$ in terms of $\hat{f}_{N}\left(t\right)$. That
being said, both of these realizations of a function's $N$th truncation
will be of use to us, both now and in the future. 
\end{rem}
\vphantom{}

Now we can begin our approach of the Banach algebra of rising-continuous\index{rising-continuous!Banach algebra of}
functions. 
\begin{defn}
For any integer $n\geq0$, we define the operator \nomenclature{$\nabla_{p^{n}}$}{ }$\nabla_{p^{n}}:B\left(\mathbb{Z}_{p},\mathbb{C}_{q}\right)\rightarrow B\left(\mathbb{Z}_{p},\mathbb{C}_{q}\right)$
by: 
\begin{equation}
\nabla_{p^{n}}\left\{ \chi\right\} \left(\mathfrak{z}\right)\overset{\textrm{def}}{=}\begin{cases}
\chi\left(\left[\mathfrak{z}\right]_{p^{n}}\right)-\chi\left(\left[\mathfrak{z}\right]_{p^{n-1}}\right) & \textrm{if }n\geq1\\
\chi\left(0\right) & \textrm{if }n=0
\end{cases},\textrm{ }\forall\mathfrak{z}\in\mathbb{Z}_{p},\textrm{ }\forall\chi\in B\left(\mathbb{Z}_{p},\mathbb{C}_{q}\right)\label{eq:Definition of Del p^n of Chi}
\end{equation}
We then define the \textbf{$\nabla$-norm }(``del norm'') $\left\Vert \cdot\right\Vert _{\nabla}:B\left(\mathbb{Z}_{p},\mathbb{C}_{q}\right)\rightarrow\left[0,\infty\right)$
by: 
\begin{equation}
\left\Vert \chi\right\Vert _{\nabla}\overset{\textrm{def}}{=}\sup_{n\geq0}\left\Vert \nabla_{p^{n}}\left\{ \chi\right\} \right\Vert _{p,q}=\sup_{n\geq0}\left(\sup_{\mathfrak{z}\in\mathbb{Z}_{p}}\left|\nabla_{p^{n}}\left\{ \chi\right\} \left(\mathfrak{z}\right)\right|_{q}\right)\label{eq:Defintion of Del-norm}
\end{equation}
\end{defn}
\begin{prop}
Let $\chi\in B\left(\mathbb{Z}_{p},\mathbb{C}_{q}\right)$. Then:

\begin{equation}
S_{p:N}\left\{ \chi\right\} \left(\mathfrak{z}\right)=\sum_{n=0}^{N}\nabla_{p^{n}}\left\{ \chi\right\} \left(\mathfrak{z}\right),\textrm{ }\forall\mathfrak{z}\in\mathbb{Z}_{p},\textrm{ }\forall N\in\mathbb{N}_{0}\label{eq:Partial vdP series in terms of Del}
\end{equation}
Moreover, we have the formal identity: 
\begin{equation}
S_{p}\left\{ \chi\right\} \left(\mathfrak{z}\right)=\sum_{n=0}^{\infty}\nabla_{p^{n}}\left\{ \chi\right\} \left(\mathfrak{z}\right)\label{eq:vdP series in terms of Del}
\end{equation}
This identity holds for any $\mathfrak{z}\in\mathbb{Z}_{p}$ for which
either the left- or right-hand side converges in $\mathbb{C}_{q}$.
In particular, (\ref{eq:vdP series in terms of Del}) holds in $\mathbb{C}_{q}$
for all $\mathfrak{z}\in\mathbb{Z}_{p}$ whenever $\chi\in\tilde{C}\left(\mathbb{Z}_{p},\mathbb{C}_{q}\right)$. 
\end{prop}
Proof: Sum the telescoping series.

Q.E.D. 
\begin{prop}
Let $\chi,\eta\in B\left(\mathbb{Z}_{p},\mathbb{C}_{q}\right)$. Then,
for all $n\in\mathbb{N}_{0}$ and all $\mathfrak{z}\in\mathbb{Z}_{p}$:
\begin{equation}
\nabla_{p^{n}}\left\{ \chi\cdot\eta\right\} \left(\mathfrak{z}\right)=\chi\left(\left[\mathfrak{z}\right]_{p^{n-1}}\right)\nabla_{p^{n}}\left\{ \eta\right\} \left(\mathfrak{z}\right)+\eta\left(\left[\mathfrak{z}\right]_{p^{n}}\right)\nabla_{p^{n}}\left\{ \chi\right\} \left(\mathfrak{z}\right)\label{eq:Quasi-derivation identity for del}
\end{equation}
using the abuse of notation $\chi\left(\left[\mathfrak{z}\right]_{p^{0-1}}\right)$
to denote the zero function in the $n=0$ case. 
\end{prop}
Proof: Using truncation notation: 
\begin{align*}
\chi_{n-1}\left(\mathfrak{z}\right)\nabla_{p^{n}}\left\{ \eta\right\} \left(\mathfrak{z}\right)+\eta_{n}\left(\mathfrak{z}\right)\nabla_{p^{n}}\left\{ \chi\right\} \left(\mathfrak{z}\right) & =\chi_{n-1}\left(\mathfrak{z}\right)\left(\eta_{n}\left(\mathfrak{z}\right)-\eta_{n-1}\left(\mathfrak{z}\right)\right)\\
 & +\eta_{n}\left(\mathfrak{z}\right)\left(\chi_{n}\left(\mathfrak{z}\right)-\chi_{n-1}\left(\mathfrak{z}\right)\right)\\
 & =\bcancel{\chi_{n-1}\left(\mathfrak{z}\right)\eta_{n}\left(\mathfrak{z}\right)}-\chi_{n-1}\left(\mathfrak{z}\right)\eta_{n-1}\left(\mathfrak{z}\right)\\
 & +\eta_{n}\left(\mathfrak{z}\right)\chi_{n}\left(\mathfrak{z}\right)-\cancel{\eta_{n}\left(\mathfrak{z}\right)\chi_{n-1}\left(\mathfrak{z}\right)}\\
 & =\chi_{n}\left(\mathfrak{z}\right)\eta_{n}\left(\mathfrak{z}\right)-\chi_{n-1}\left(\mathfrak{z}\right)\eta_{n-1}\left(\mathfrak{z}\right)\\
 & =\nabla_{p^{n}}\left\{ \chi\cdot\eta\right\} \left(\mathfrak{z}\right)
\end{align*}

Q.E.D. 
\begin{prop}
Let $\chi\in B\left(\mathbb{Z}_{p},K\right)$. Then, $S_{p}\left\{ \chi\right\} $
is identically zero whenever $\chi$ vanishes on $\mathbb{N}_{0}$. 
\end{prop}
Proof: If $\chi\left(n\right)=0$ for all $n\in\mathbb{N}_{0}$, then
$c_{n}\left(\chi\right)=0$ for all $n\in\mathbb{N}_{0}$, and hence
$S_{p}\left\{ \chi\right\} $ is identically zero.

Q.E.D.

\vphantom{}

The moral of this proposition is that the operator $S_{p}$ behaves
like a projection operator on $B\left(\mathbb{Z}_{p},K\right)$. For
example, given any rising-continuous function $\chi$ and any function
$f\in B\left(\mathbb{Z}_{p},K\right)$ with $f\left(n\right)=0$ for
all $n\in\mathbb{N}_{0}$, we have that: 
\begin{equation}
S_{p}\left\{ \chi+f\right\} =S_{p}\left\{ \chi\right\} =\chi
\end{equation}
The reason why we say $S_{p}$ behaves merely ``like'' a projection
operator on $B\left(\mathbb{Z}_{p},K\right)$\textemdash rather than
proving that $S_{p}$ actually \emph{is }a projection operator on
$B\left(\mathbb{Z}_{p},K\right)$\textemdash is because there exist
$f\in B\left(\mathbb{Z}_{p},K\right)$ for which $S_{p}\left\{ f\right\} \left(\mathfrak{z}\right)$
fails to converge at \emph{any} $\mathfrak{z}\in\mathbb{Z}_{p}$.
Case in point:
\begin{example}
Let $f\in B\left(\mathbb{Z}_{p},K\right)$ be defined by: 
\begin{equation}
f\left(\mathfrak{z}\right)\overset{\textrm{def}}{=}\begin{cases}
\mathfrak{z} & \textrm{if }\mathfrak{z}\in\mathbb{N}_{0}\\
0 & \textrm{else}
\end{cases}\label{eq:Example of a bounded (p,q)-adic function whose van der Put series is divergent.}
\end{equation}
Then, we have the formal identity: 
\begin{align*}
S_{p}\left\{ f\right\} \left(\mathfrak{z}\right) & =\sum_{n=1}^{\infty}n\left[\mathfrak{z}\overset{p^{\lambda_{p}\left(n\right)}}{\equiv}n\right]\\
 & =\sum_{n=1}^{\infty}\sum_{m=p^{n-1}}^{p^{n}-1}m\left[\mathfrak{z}\overset{p^{n}}{\equiv}m\right]\\
 & =\sum_{n=1}^{\infty}\left(\sum_{m=0}^{p^{n}-1}m\left[\mathfrak{z}\overset{p^{n}}{\equiv}m\right]-\sum_{m=0}^{p^{n-1}-1}m\left[\mathfrak{z}\overset{p^{n}}{\equiv}m\right]\right)
\end{align*}
Here: 
\[
\sum_{m=0}^{p^{n}-1}m\left[\mathfrak{z}\overset{p^{n}}{\equiv}m\right]=\left[\mathfrak{z}\right]_{p^{n}}
\]
because $m=\left[\mathfrak{z}\right]_{p^{n}}$ is the unique integer
in $\left\{ 0,\ldots,p^{n}-1\right\} $ which is congruent to $\mathfrak{z}$
mod $p^{n}$. On the other hand, for the $m$-sum with $p^{n-1}-1$
as an upper bound, we end up with: 
\[
\sum_{m=0}^{p^{n-1}-1}m\left[\mathfrak{z}\overset{p^{n}}{\equiv}m\right]=\left[\mathfrak{z}\right]_{p^{n}}\left[\left[\mathfrak{z}\right]_{p^{n}}<p^{n-1}\right]
\]
and so: 
\begin{align*}
S_{p}\left\{ f\right\} \left(\mathfrak{z}\right) & =\sum_{n=1}^{\infty}\left[\mathfrak{z}\right]_{p^{n}}\left(1-\left[\left[\mathfrak{z}\right]_{p^{n}}<p^{n-1}\right]\right)\\
 & =\sum_{n=1}^{\infty}\left[\mathfrak{z}\right]_{p^{n}}\left[\left[\mathfrak{z}\right]_{p^{n}}\geq p^{n-1}\right]
\end{align*}
Consequently: 
\[
S_{p}\left\{ f\right\} \left(-1\right)=\sum_{n=1}^{\infty}\left(p^{n}-1\right)\left[p^{n}-1\geq p^{n-1}\right]=\sum_{n=1}^{\infty}\left(p^{n}-1\right)
\]
which does not converge in $\mathbb{C}_{q}$ because $\left|p^{n}-1\right|_{q}$
does not tend to $0$ in $\mathbb{R}$ as $n\rightarrow\infty$. 
\end{example}
\vphantom{}

The above example also shows that $S_{p}$ is as much a continuation-creating
operator as it is a projection operator. Bounded functions which do
not behave well under limits of sequences in $\mathbb{Z}_{p}$ will
have ill-behaved images under $S_{p}$. This motivates the following
definition: 
\begin{defn}
We write \nomenclature{$\tilde{B}\left(\mathbb{Z}_{p},\mathbb{C}_{q}\right)$}{set of $f\in B\left(\mathbb{Z}_{p},\mathbb{C}_{q}\right)$ so that $S_{p}\left\{ f\right\} \in B\left(\mathbb{Z}_{p},\mathbb{C}_{q}\right)$.}$\tilde{B}\left(\mathbb{Z}_{p},\mathbb{C}_{q}\right)$
to denote the set of all $f\in B\left(\mathbb{Z}_{p},\mathbb{C}_{q}\right)$
so that $S_{p}\left\{ f\right\} \in B\left(\mathbb{Z}_{p},\mathbb{C}_{q}\right)$
(i.e., $S_{p}\left\{ f\right\} $ converges point-wise in $\mathbb{C}_{q}$). 
\end{defn}
\begin{rem}
Note that $\tilde{B}\left(\mathbb{Z}_{p},\mathbb{C}_{q}\right)\neq\tilde{C}\left(\mathbb{Z}_{p},\mathbb{C}_{q}\right)$.
The only function in $\tilde{C}\left(\mathbb{Z}_{p},\mathbb{C}_{q}\right)$
which vanishes on $\mathbb{N}_{0}$ is the constant function $0$.
Consequently, any function in $B\left(\mathbb{Z}_{p},\mathbb{C}_{q}\right)$
which is supported on $\mathbb{Z}_{p}^{\prime}$ is an element of
$\tilde{B}\left(\mathbb{Z}_{p},\mathbb{C}_{q}\right)$ which is not
in $\tilde{C}\left(\mathbb{Z}_{p},\mathbb{C}_{q}\right)$. That such
a function is in $\tilde{B}\left(\mathbb{Z}_{p},\mathbb{C}_{q}\right)$
is, of course, because its image under $S_{p}$ is identically $0$. 
\end{rem}
\begin{prop}
The map: 
\[
\left(\chi,\eta\right)\in B\left(\mathbb{Z}_{p},\mathbb{C}_{q}\right)\times B\left(\mathbb{Z}_{p},\mathbb{C}_{q}\right)\mapsto\left\Vert \chi-\eta\right\Vert _{\nabla}\in\left[0,\infty\right)
\]
defines a non-archimedean semi-norm on $B\left(\mathbb{Z}_{p},\mathbb{C}_{q}\right)$
and $\tilde{B}\left(\mathbb{Z}_{p},\mathbb{C}_{q}\right)$.
\end{prop}
Proof:

I. (Positivity) Observe that $\left\Vert \chi\right\Vert _{\nabla}<\infty$
for all $\chi\in B\left(\mathbb{Z}_{p},\mathbb{C}_{q}\right)$. However,
$\left\Vert \cdot\right\Vert _{\nabla}$ is \emph{not }positive-definite;
let $\chi\in B\left(\mathbb{Z}_{p},\mathbb{C}_{q}\right)$ be any
function supported on $\mathbb{Z}_{p}^{\prime}$. Then, $\chi$ vanishes
on $\mathbb{N}_{0}$, and so, $\chi\left(\left[\mathfrak{z}\right]_{p^{n}}\right)=0$
for all $\mathfrak{z}\in\mathbb{Z}_{p}$ and all $n\in\mathbb{N}_{0}$.
This forces $\left\Vert \chi\right\Vert _{\nabla}$ to be zero; however,
$\chi$ is not identically zero. So, $\left\Vert \cdot\right\Vert _{\nabla}$
cannot be a norm on $B\left(\mathbb{Z}_{p},\mathbb{C}_{q}\right)$

\vphantom{}

II. Since $\nabla_{p^{n}}$ is linear for all $n\geq0$, for any $\mathfrak{a}\in\mathbb{C}_{q}$,
we have that: 
\[
\left\Vert \mathfrak{a}\chi\right\Vert _{\nabla}=\sup_{n\geq0}\left\Vert \mathfrak{a}\nabla_{p^{n}}\left\{ \chi\right\} \right\Vert _{p,q}=\left|\mathfrak{a}\right|_{q}\sup_{n\geq0}\left\Vert \nabla_{p^{n}}\left\{ \chi\right\} \right\Vert _{p,q}=\left|\mathfrak{a}\right|_{q}\left\Vert \chi\right\Vert _{\nabla}
\]
for all $\chi\in B\left(\mathbb{Z}_{p},\mathbb{C}_{q}\right)$. As
for the ultrametric inequality, observing that: 
\[
\left\Vert \chi-\eta\right\Vert _{\nabla}=\sup_{n\geq0}\left\Vert \nabla_{p^{n}}\left\{ \chi-\eta\right\} \right\Vert _{p,q}=\sup_{n\geq0}\left\Vert \nabla_{p^{n}}\left\{ \chi\right\} -\nabla_{p^{n}}\left\{ \eta\right\} \right\Vert _{p,q}
\]
the fact that $d\left(\chi,\eta\right)=\left\Vert \chi-\eta\right\Vert _{p,q}$
defines an ultrametric on $B\left(\mathbb{Z}_{p},\mathbb{C}_{q}\right)$
lets us write: 
\begin{align*}
\left\Vert \chi-\eta\right\Vert _{\nabla} & =\sup_{n\geq0}\left\Vert \nabla_{p^{n}}\left\{ \chi\right\} -\nabla_{p^{n}}\left\{ \eta\right\} \right\Vert _{p,q}\\
 & \leq\sup_{n\geq0}\max\left\{ \left\Vert \nabla_{p^{n}}\left\{ \chi\right\} \right\Vert _{p,q},\left\Vert \nabla_{p^{n}}\left\{ \eta\right\} \right\Vert _{p,q}\right\} \\
 & =\max\left\{ \left\Vert \chi\right\Vert _{\nabla},\left\Vert \eta\right\Vert _{\nabla}\right\} 
\end{align*}
which proves the ultrametric inequality.

Thus, $\left\Vert \cdot\right\Vert _{\nabla}$ is a non-archimedean
semi-norm on $B\left(\mathbb{Z}_{p},\mathbb{C}_{q}\right)$. The arguments
for $\tilde{B}\left(\mathbb{Z}_{p},\mathbb{C}_{q}\right)$ are identical.

Q.E.D.

\vphantom{}

Next, we show that $\nabla$-(semi)norm dominates the $\left(p,q\right)$-adic
supremum norm. 
\begin{prop}
\label{prop:continuous embeds in rising-continuous}If $\chi\in\tilde{C}\left(\mathbb{Z}_{p},\mathbb{C}_{q}\right)$,
then: 
\begin{equation}
\left\Vert \chi\right\Vert _{p,q}\leq\left\Vert \chi\right\Vert _{\nabla}\label{eq:Supremum norm is dominated by Del norm}
\end{equation}
\end{prop}
Proof: If $\chi\in\tilde{C}\left(\mathbb{Z}_{p},\mathbb{C}_{q}\right)$,
then the van der Put series for $\chi$ converges $q$-adically to
$\chi$ point-wise over $\mathbb{Z}_{p}$. As such: 
\begin{align*}
\left|\chi\left(\mathfrak{z}\right)\right|_{q} & =\left|\sum_{n=0}^{\infty}\nabla_{p^{n}}\left\{ \chi\right\} \left(\mathfrak{z}\right)\right|_{q}\\
 & \leq\sup_{n\geq0}\left|\nabla_{p^{n}}\left\{ \chi\right\} \left(\mathfrak{z}\right)\right|_{q}
\end{align*}
Taking suprema over $\mathfrak{z}\in\mathbb{Z}_{p}$ gives: 
\begin{align*}
\left\Vert \chi\right\Vert _{p,q} & =\sup_{\mathfrak{z}\in\mathbb{Z}_{p}}\left|\chi\left(\mathfrak{z}\right)\right|_{q}\\
 & \leq\sup_{\mathfrak{z}\in\mathbb{Z}_{p}}\sup_{n\geq0}\left|\nabla_{p^{n}}\left\{ \chi\right\} \left(\mathfrak{z}\right)\right|_{q}\\
\left(\sup_{m}\sup_{n}x_{m,n}=\sup_{n}\sup_{m}x_{m,n}\right); & \leq\sup_{n\geq0}\sup_{\mathfrak{z}\in\mathbb{Z}_{p}}\left|\nabla_{p^{n}}\left\{ \chi\right\} \left(\mathfrak{z}\right)\right|_{q}\\
 & =\left\Vert \chi\right\Vert _{\nabla}
\end{align*}

Q.E.D. 
\begin{lem}
$\left(\tilde{C}\left(\mathbb{Z}_{p},\mathbb{C}_{q}\right),\left\Vert \cdot\right\Vert _{\nabla}\right)$
is a non-archimedean normed linear space. In particular, we can identify
$\tilde{C}\left(\mathbb{Z}_{p},\mathbb{C}_{q}\right)$ with the image
of $\tilde{B}\left(\mathbb{Z}_{p},\mathbb{C}_{q}\right)$ under $S_{p}$.
Moreover, $S_{p}$ is then an isometry of $\tilde{C}\left(\mathbb{Z}_{p},\mathbb{C}_{q}\right)$. 
\end{lem}
Proof: By the First Isomorphism Theorem from abstract algebra, $S_{p}$
maps the $\mathbb{C}_{q}$-linear space $\tilde{B}\left(\mathbb{Z}_{p},\mathbb{C}_{q}\right)$
isomorphically onto $\tilde{B}\left(\mathbb{Z}_{p},\mathbb{C}_{q}\right)/\textrm{Ker}\left(S_{p}\right)$.
Since $\left\Vert \cdot\right\Vert _{\nabla}$ is a semi-norm on $\tilde{B}\left(\mathbb{Z}_{p},\mathbb{C}_{q}\right)$,
and because the set of non-zero $\chi\in\tilde{B}\left(\mathbb{Z}_{p},\mathbb{C}_{q}\right)$
possessing $\left\Vert \chi\right\Vert _{\nabla}=0$ is equal to $\textrm{Ker}\left(S_{p}\right)$,
the non-archimedean semi-norm $\left\Vert \cdot\right\Vert _{\nabla}$
is then a norm on the quotient space $\tilde{B}\left(\mathbb{Z}_{p},\mathbb{C}_{q}\right)/\textrm{Ker}\left(S_{p}\right)$.

To conclude, we need only show that $\tilde{B}\left(\mathbb{Z}_{p},\mathbb{C}_{q}\right)/\textrm{Ker}\left(S_{p}\right)$
can be identified with $\tilde{C}\left(\mathbb{Z}_{p},\mathbb{C}_{q}\right)$.
To do so, let $\chi\in\tilde{C}\left(\mathbb{Z}_{p},\mathbb{C}_{q}\right)$
be arbitrary. Then, $\chi$ is an element of $B\left(\mathbb{Z}_{p},\mathbb{C}_{q}\right)$.
As we saw in the previous subsection (viz., \textbf{Theorem \ref{thm:S_p}}),
$S_{p}$ is the identity map on $\tilde{C}\left(\mathbb{Z}_{p},\mathbb{C}_{q}\right)$.
Consequently, $S_{p}\left\{ \chi\right\} =\chi$, and so, $\chi\in\tilde{B}\left(\mathbb{Z}_{p},\mathbb{C}_{q}\right)$.

For the other direction, let $\chi$ and $\eta$ be two elements of
$\tilde{B}\left(\mathbb{Z}_{p},\mathbb{C}_{q}\right)$ which belong
to the same equivalence class in $\tilde{B}\left(\mathbb{Z}_{p},\mathbb{C}_{q}\right)/\textrm{Ker}\left(S_{p}\right)$.
Then, $S_{p}\left\{ \chi-\eta\right\} $ is the zero function. By
the linearity of $S_{p}$, this means $S_{p}\left\{ \chi\right\} =S_{p}\left\{ \eta\right\} $.
So, $S_{p}\left\{ \chi\right\} $ and $S_{p}\left\{ \eta\right\} $
have the same van der Put series, and that series defines an element
of $B\left(\mathbb{Z}_{p},\mathbb{C}_{q}\right)$. Since this series
converges everywhere point-wise, the van der Put identity (\textbf{Proposition
\ref{prop:vdP identity}}) shows it is an element of $\tilde{C}\left(\mathbb{Z}_{p},\mathbb{C}_{q}\right)$.

Thus, the two spaces are equivalent.

Q.E.D. 
\begin{thm}
$\left(\tilde{C}\left(\mathbb{Z}_{p},\mathbb{C}_{q}\right),\left\Vert \cdot\right\Vert _{\nabla}\right)$
is a non-archimedean Banach algebra over $\mathbb{C}_{q}$, with point-wise
multiplication of functions as the algebra multiplication identity.
Additionally, $\left(\tilde{C}\left(\mathbb{Z}_{p},\mathbb{C}_{q}\right),\left\Vert \cdot\right\Vert _{\nabla}\right)$
contains $\left(C\left(\mathbb{Z}_{p},\mathbb{C}_{q}\right),\left\Vert \cdot\right\Vert _{\nabla}\right)$
as a sub-algebra. 
\end{thm}
Proof: To see that $\left(\tilde{C}\left(\mathbb{Z}_{p},\mathbb{C}_{q}\right),\left\Vert \cdot\right\Vert _{\nabla}\right)$
is a normed algebra, let $\chi,\eta\in\tilde{C}\left(\mathbb{Z}_{p},\mathbb{C}_{q}\right)$.
Using (\ref{eq:Quasi-derivation identity for del}) and truncation
notation, we have that:

\begin{align*}
\left\Vert \chi\cdot\eta\right\Vert _{\nabla} & =\sup_{n\geq0}\left\Vert \nabla_{p^{n}}\left\{ \chi\cdot\eta\right\} \right\Vert _{p,q}\\
 & =\sup_{n\geq0}\left\Vert \chi_{n-1}\cdot\nabla_{p^{n}}\left\{ \eta\right\} +\eta_{n}\cdot\nabla_{p^{n}}\left\{ \chi\right\} \right\Vert _{p,q}\\
\left(\left\Vert \cdot\right\Vert _{p,q}\textrm{ is non-arch.}\right); & \leq\max\left\{ \sup_{m\geq0}\left\Vert \chi_{m-1}\cdot\nabla_{p^{m}}\left\{ \eta\right\} \right\Vert _{p,q},\sup_{n\geq0}\left\Vert \eta_{n}\cdot\nabla_{p^{n}}\left\{ \chi\right\} \right\Vert _{p,q}\right\} 
\end{align*}
For each $m$ and $n$, the function $\chi_{m-1}$, $\nabla_{p^{m}}\left\{ \eta\right\} $,
$\eta_{n}$, and $\nabla_{p^{n}}\left\{ \chi\right\} $ are continuous
functions. Then, because $\left(C\left(\mathbb{Z}_{p},\mathbb{C}_{q}\right),\left\Vert \cdot\right\Vert _{p,q}\right)$
is a Banach algebra, we can write: 
\begin{align*}
\left\Vert \chi\cdot\eta\right\Vert _{\nabla} & \leq\max\left\{ \sup_{m\geq0}\left(\left\Vert \chi_{m-1}\right\Vert _{p,q}\left\Vert \nabla_{p^{m}}\left\{ \eta\right\} \right\Vert _{p,q}\right),\sup_{n\geq0}\left(\left\Vert \eta_{n}\right\Vert _{p,q}\left\Vert \nabla_{p^{n}}\left\{ \chi\right\} \right\Vert _{p,q}\right)\right\} 
\end{align*}
Since $\chi_{m-1}$ is rising-continuous, \textbf{Proposition \ref{prop:continuous embeds in rising-continuous}}
tells us that $\left\Vert \chi_{m-1}\right\Vert _{p,q}\leq\left\Vert \chi_{m-1}\right\Vert _{\nabla}$.
Moreover, because $\chi_{m-1}$ is the restriction of $\chi$ to inputs
mod $p^{m-1}$, we have:
\begin{align*}
\left\Vert \chi_{m-1}\right\Vert _{\nabla} & =\sup_{k\geq0}\sup_{\mathfrak{z}\in\mathbb{Z}_{p}}\left|\nabla_{p^{k}}\left\{ \chi_{m-1}\right\} \left(\mathfrak{z}\right)\right|_{q}\\
 & \leq\sup_{k\geq0}\sup_{\mathfrak{z}\in\mathbb{Z}_{p}}\left|\chi_{m-1}\left(\left[\mathfrak{z}\right]_{p^{k}}\right)-\chi_{m-1}\left(\left[\mathfrak{z}\right]_{p^{k-1}}\right)\right|_{q}\\
 & \leq\sup_{k\geq0}\sup_{\mathfrak{z}\in\mathbb{Z}_{p}}\left|\chi\left(\left[\left[\mathfrak{z}\right]_{p^{k}}\right]_{p^{m-1}}\right)-\chi_{m-1}\left(\left[\left[\mathfrak{z}\right]_{p^{k-1}}\right]_{p^{m-1}}\right)\right|_{q}
\end{align*}
Here: 
\begin{equation}
\left[\left[\mathfrak{z}\right]_{p^{k}}\right]_{p^{m-1}}=\left[\mathfrak{z}\right]_{p^{\min\left\{ k,m-1\right\} }}
\end{equation}
and so: 
\begin{align*}
\left\Vert \chi_{m-1}\right\Vert _{\nabla} & =\sup_{k\geq0}\sup_{\mathfrak{z}\in\mathbb{Z}_{p}}\left|\nabla_{p^{k}}\left\{ \chi_{m-1}\right\} \left(\mathfrak{z}\right)\right|_{q}\\
 & \leq\sup_{k\geq0}\sup_{\mathfrak{z}\in\mathbb{Z}_{p}}\left|\chi_{m-1}\left(\left[\mathfrak{z}\right]_{p^{k}}\right)-\chi_{m-1}\left(\left[\mathfrak{z}\right]_{p^{k-1}}\right)\right|_{q}\\
 & \leq\sup_{k\geq0}\sup_{\mathfrak{z}\in\mathbb{Z}_{p}}\underbrace{\left|\chi\left(\left[\mathfrak{z}\right]_{p^{\min\left\{ k,m-1\right\} }}\right)-\chi\left(\left[\mathfrak{z}\right]_{p^{\min\left\{ k-1,m-1\right\} }}\right)\right|_{q}}_{0\textrm{ when }k\geq m}\\
 & \leq\sup_{0\leq k<m}\sup_{\mathfrak{z}\in\mathbb{Z}_{p}}\left|\chi\left(\left[\mathfrak{z}\right]_{p^{k}}\right)-\chi\left(\left[\mathfrak{z}\right]_{p^{k-1}}\right)\right|_{q}\\
 & \leq\sup_{k\geq0}\sup_{\mathfrak{z}\in\mathbb{Z}_{p}}\left|\chi\left(\left[\mathfrak{z}\right]_{p^{k}}\right)-\chi\left(\left[\mathfrak{z}\right]_{p^{k-1}}\right)\right|_{q}\\
 & =\left\Vert \chi\right\Vert _{\nabla}
\end{align*}
Applying this argument to the $\eta_{n}$ yields: 
\begin{equation}
\left\Vert \eta_{n}\right\Vert _{p,q}\leq\left\Vert \eta_{n}\right\Vert _{\nabla}\leq\left\Vert \eta\right\Vert _{\nabla}
\end{equation}
and hence: 
\begin{align*}
\left\Vert \chi\cdot\eta\right\Vert _{\nabla} & \leq\max\left\{ \sup_{m\geq0}\left(\left\Vert \chi_{m-1}\right\Vert _{p,q}\left\Vert \nabla_{p^{m}}\left\{ \eta\right\} \right\Vert _{p,q}\right),\sup_{n\geq0}\left(\left\Vert \eta_{n}\right\Vert _{p,q}\left\Vert \nabla_{p^{n}}\left\{ \chi\right\} \right\Vert _{p,q}\right)\right\} \\
 & \leq\max\left\{ \left\Vert \chi\right\Vert _{\nabla}\sup_{m\geq0}\left\Vert \nabla_{p^{m}}\left\{ \eta\right\} \right\Vert _{p,q},\left\Vert \eta\right\Vert _{\nabla}\sup_{n\geq0}\left\Vert \nabla_{p^{n}}\left\{ \chi\right\} \right\Vert _{p,q}\right\} \\
 & =\max\left\{ \left\Vert \chi\right\Vert _{\nabla}\left\Vert \eta\right\Vert _{\nabla},\left\Vert \eta\right\Vert _{\nabla}\left\Vert \chi\right\Vert _{\nabla}\right\} \\
 & =\left\Vert \chi\right\Vert _{\nabla}\left\Vert \eta\right\Vert _{\nabla}
\end{align*}
Thus, $\left(\tilde{C}\left(\mathbb{Z}_{p},\mathbb{C}_{q}\right),\left\Vert \cdot\right\Vert _{\nabla}\right)$
is a normed algebra. The inequality from \textbf{Proposition \ref{prop:continuous embeds in rising-continuous}}
then shows that $\left(\tilde{C}\left(\mathbb{Z}_{p},\mathbb{C}_{q}\right),\left\Vert \cdot\right\Vert _{\nabla}\right)$
contains $\left(C\left(\mathbb{Z}_{p},\mathbb{C}_{q}\right),\left\Vert \cdot\right\Vert _{\nabla}\right)$
as a sub-algebra.

To conclude, we need only establish the completeness of $\left(\tilde{C}\left(\mathbb{Z}_{p},\mathbb{C}_{q}\right),\left\Vert \cdot\right\Vert _{\nabla}\right)$.
We do this using \textbf{Proposition \ref{prop:series characterization of a Banach space}}:
it suffices to prove that every sequence in $\tilde{C}\left(\mathbb{Z}_{p},\mathbb{C}_{q}\right)$
which tends to $0$ in $\nabla$-norm is summable to an element of
$\tilde{C}\left(\mathbb{Z}_{p},\mathbb{C}_{q}\right)$.

So, let $\left\{ \chi_{n}\right\} _{n\geq0}$ be a sequence\footnote{This is \emph{not }an instance truncation notation. I apologize for
this abuse of notation.} in $\tilde{C}\left(\mathbb{Z}_{p},\mathbb{C}_{q}\right)$ which tends
to $0$ in $\nabla$-norm. Then, by \textbf{Proposition \ref{prop:series characterization of a Banach space}},
upon letting $M$ and $N$ be arbitrary positive integers with $M\leq N$,
we can write: 
\begin{equation}
\left\Vert \sum_{n=M}^{N}\chi_{n}\left(\mathfrak{z}\right)\right\Vert _{p,q}\leq\left\Vert \sum_{n=M}^{N}\chi_{n}\left(\mathfrak{z}\right)\right\Vert _{\nabla}\leq\max_{n\geq M}\left\Vert \chi_{n}\right\Vert _{\nabla}
\end{equation}
Here, the upper bound tends to $0$ as $M\rightarrow\infty$. Consequently
the series $\sum_{n=0}^{\infty}\chi_{n}\left(\mathfrak{z}\right)$
converges $q$-adically to a sum $X\left(\mathfrak{z}\right)$, and,
moreover, this convergence is \emph{uniform} in $\mathfrak{z}\in\mathbb{Z}_{p}$.
Since the $\chi_{n}$s are rising-continuous, we have: 
\begin{equation}
S_{p}\left\{ \sum_{n=0}^{N}\chi_{n}\right\} \left(\mathfrak{z}\right)=\sum_{n=0}^{N}S_{p}\left\{ \chi_{n}\right\} \left(\mathfrak{z}\right)=\sum_{n=0}^{N}\chi_{n}\left(\mathfrak{z}\right)
\end{equation}
The uniform convergence of the $n$-sum justifies interchanging limits
with $S_{p}$:
\[
X\left(\mathfrak{z}\right)=\lim_{N\rightarrow\infty}\sum_{n=0}^{N}\chi_{n}=\lim_{N\rightarrow\infty}S_{p}\left\{ \sum_{n=0}^{N}\chi_{n}\right\} \left(\mathfrak{z}\right)=S_{p}\left\{ \lim_{N\rightarrow\infty}\sum_{n=0}^{N}\chi_{n}\right\} \left(\mathfrak{z}\right)=S_{p}\left\{ X\right\} \left(\mathfrak{z}\right)
\]
Consequently, $X=S_{p}\left\{ X\right\} $. This proves that $X$
is rising-continuous. So, $\sum_{n=0}^{\infty}\chi_{n}\left(\mathfrak{z}\right)$
converges in $\tilde{C}\left(\mathbb{Z}_{p},\mathbb{C}_{q}\right)$,
establishing the completeness of $\left(\tilde{C}\left(\mathbb{Z}_{p},\mathbb{C}_{q}\right),\left\Vert \cdot\right\Vert _{\nabla}\right)$.

Q.E.D. 
\begin{problem}
Although we shall not explore it here, I think it would be interesting
to investigate the behavior of $\tilde{C}\left(\mathbb{Z}_{p},\mathbb{C}_{q}\right)$
under linear operators $T_{\phi}:\tilde{C}\left(\mathbb{Z}_{p},\mathbb{C}_{q}\right)\rightarrow B\left(\mathbb{Z}_{p},\mathbb{C}_{q}\right)$
of the form: 
\begin{equation}
T_{\phi}\left\{ \chi\right\} \left(\mathfrak{z}\right)\overset{\textrm{def}}{=}\chi\left(\phi\left(\mathfrak{z}\right)\right),\textrm{ }\forall\chi\in\tilde{C}\left(\mathbb{Z}_{p},\mathbb{C}_{q}\right),\textrm{ }\forall\mathfrak{z}\in\mathbb{Z}_{p}\label{eq:Definition of T_phi of Chi}
\end{equation}
where $\phi:\mathbb{Z}_{p}\rightarrow\mathbb{Z}_{p}$. The simplest
$T_{\phi}$s are translation operator\index{translation operator}s:
\begin{equation}
\tau_{\mathfrak{a}}\left\{ \chi\right\} \left(\mathfrak{z}\right)\overset{\textrm{def}}{=}\chi\left(\mathfrak{z}+\mathfrak{a}\right)\label{eq:Definition of the translation operator}
\end{equation}
where $\mathfrak{a}$ is a fixed element of $\mathbb{Z}_{p}$. The
van der Put series (\ref{eq:van der Point series for one point function at 0})
given for $\mathbf{1}_{0}\left(\mathfrak{z}\right)$ in Subsection
\ref{subsec:3.1.6 Monna-Springer-Integration} shows that $\mathbf{1}_{0}\left(\mathfrak{z}\right)\in\tilde{C}\left(\mathbb{Z}_{p},\mathbb{C}_{q}\right)$,
and that $\tau_{n}\left\{ \mathbf{1}_{0}\right\} =\mathbf{1}_{n}$
is in $\tilde{C}\left(\mathbb{Z}_{p},\mathbb{C}_{q}\right)$ for all
$n\geq0$. However, observe that for any $\mathfrak{a}\in\mathbb{Z}_{p}^{\prime}$,
$\mathbf{1}_{\mathfrak{a}}\left(\mathfrak{z}\right)$ vanishes on
$\mathbb{N}_{0}$, and as a consequence, $S_{p}\left\{ \mathbf{1}_{\mathfrak{a}}\right\} $
is identically $0$; that is to say, $\mathbf{1}_{\mathfrak{a}}\left(\mathfrak{z}\right)$
is rising-continuous \emph{if and only if }$\mathfrak{a}\in\mathbb{N}_{0}$.
This shows that $\tilde{C}\left(\mathbb{Z}_{p},\mathbb{C}_{q}\right)$
is not closed under translations.

Aside from being a curiosity in its own right, studying operators
like (\ref{eq:Definition of T_phi of Chi}) is also relevant to our
investigations of Hydra maps. Choosing $\phi=B_{p}$, where: 
\begin{equation}
B_{p}\left(\mathfrak{z}\right)\overset{\mathbb{Z}_{p}}{=}\begin{cases}
0 & \textrm{if }\mathfrak{z}=0\\
\frac{\mathfrak{z}}{1-p^{\lambda_{p}\left(\mathfrak{z}\right)}} & \textrm{if }\mathfrak{z}\in\mathbb{N}_{1}\\
\mathfrak{z} & \textrm{if }\mathfrak{z}\in\mathbb{Z}_{p}^{\prime}
\end{cases}
\end{equation}
is $B_{p}$ as defined as (\ref{eq:Definition of B projection function}),
observe that the functional equation (\ref{eq:Chi_H B functional equation})
can be written as: 
\[
T_{B_{p}}\left\{ \chi_{H}\right\} \left(\mathfrak{z}\right)=\chi_{H}\left(B_{p}\left(\mathfrak{z}\right)\right)=\begin{cases}
\chi_{H}\left(0\right) & \textrm{if }\mathfrak{z}=0\\
\frac{\chi_{H}\left(\mathfrak{z}\right)}{1-M_{H}\left(\mathfrak{z}\right)} & \textrm{else}
\end{cases}
\]
where, by the arguments shown in Chapter 2, $M_{H}\left(\mathfrak{z}\right)\overset{\mathbb{Z}_{q}}{=}0$
for all $\mathfrak{z}\in\mathbb{Z}_{p}^{\prime}$ ( every such $\mathfrak{z}\in\mathbb{Z}_{p}^{\prime}$
has infinitely many non-zero $p$-adic digits). As such, understanding
the action of the operator $T_{B_{p}}$ may be of use to understanding
the periodic points of $H$. 
\end{problem}
\vphantom{}

We now move to investigate the units of $\tilde{C}\left(\mathbb{Z}_{p},\mathbb{C}_{q}\right)$. 
\begin{prop}[\textbf{Unit criteria for} $\tilde{C}\left(\mathbb{Z}_{p},\mathbb{C}_{q}\right)$]
\label{prop:criteria for being a rising-continuous unit}\ 

\vphantom{}

I. Let $\chi\in\tilde{C}\left(\mathbb{Z}_{p},K\right)$. Then, $\chi$
is a unit of $\tilde{C}\left(\mathbb{Z}_{p},K\right)$ whenever: 
\begin{equation}
\inf_{n\in\mathbb{N}_{0}}\left|\chi\left(n\right)\right|_{q}>0\label{eq:Sufficient condition for a rising continuous function to have a reciprocal}
\end{equation}

\vphantom{}

II. Let $\chi\in C\left(\mathbb{Z}_{p},K\right)$. Then, $\chi$ is
a unit of $\tilde{C}\left(\mathbb{Z}_{p},K\right)$ if and only if
$\chi$ has no zeroes. 
\end{prop}
Proof:

I. Let $\chi\in\tilde{C}\left(\mathbb{Z}_{p},K\right)$.

\begin{align*}
\left\Vert \frac{1}{\chi}\right\Vert _{\nabla} & =\sup_{n\geq0}\left\Vert \nabla_{p^{n}}\left\{ \frac{1}{\chi}\right\} \right\Vert _{\nabla}\\
 & =\max\left\{ \left|\frac{1}{\chi\left(0\right)}\right|_{q},\sup_{n\geq1}\sup_{\mathfrak{z}\in\mathbb{Z}_{p}}\left|\frac{1}{\chi\left(\left[\mathfrak{z}\right]_{p^{n}}\right)}-\frac{1}{\chi\left(\left[\mathfrak{z}\right]_{p^{n-1}}\right)}\right|_{q}\right\} \\
 & =\max\left\{ \left|\frac{1}{\chi\left(0\right)}\right|_{q},\sup_{n\geq1}\sup_{\mathfrak{z}\in\mathbb{Z}_{p}}\left|\frac{\nabla_{p^{n}}\left\{ \chi\right\} \left(\mathfrak{z}\right)}{\chi\left(\left[\mathfrak{z}\right]_{p^{n}}\right)\chi\left(\left[\mathfrak{z}\right]_{p^{n-1}}\right)}\right|_{q}\right\} \\
 & \leq\max\left\{ \left|\frac{1}{\chi\left(0\right)}\right|_{q},\sup_{m\geq0}\sup_{\mathfrak{z}\in\mathbb{Z}_{p}}\left|\nabla_{p^{m}}\left\{ \chi\right\} \left(\mathfrak{z}\right)\right|\sup_{k\geq1}\sup_{\mathfrak{z}\in\mathbb{Z}_{p}}\left|\chi\left(\left[\mathfrak{z}\right]_{p^{k}}\right)\right|_{q}^{-1}\cdot\sup_{n\geq1}\sup_{\mathfrak{z}\in\mathbb{Z}_{p}}\left|f\left(\left[\mathfrak{z}\right]_{p^{n-1}}\right)\right|_{q}^{-1}\right\} \\
 & \leq\max\left\{ \left|\frac{1}{\chi\left(0\right)}\right|_{q},\frac{\left\Vert \chi\right\Vert _{\nabla}}{\inf_{n\geq1}\left|\chi\left(n\right)\right|_{q}^{2}}\right\} 
\end{align*}
which will be finite whenever $\inf_{n\in\mathbb{N}_{0}}\left|\chi\left(n\right)\right|_{q}>0$.

\vphantom{}

II. If $\chi\in C\left(\mathbb{Z}_{p},K\right)$ has no zeroes, then
$1/\chi\in C\left(\mathbb{Z}_{p},K\right)$. Since $C\left(\mathbb{Z}_{p},K\right)\subset\tilde{C}\left(\mathbb{Z}_{p},K\right)$,
we then have that $1/\chi\in\tilde{C}\left(\mathbb{Z}_{p},K\right)$.
Conversely, if $1/\chi\in C\left(\mathbb{Z}_{p},K\right)$ has a zero,
then $1/\chi\notin C\left(\mathbb{Z}_{p},K\right)$. If $1/\chi$
were in $\tilde{C}\left(\mathbb{Z}_{p},K\right)$, then since the
norm on $\tilde{C}\left(\mathbb{Z}_{p},K\right)$ dominates the norm
on $C\left(\mathbb{Z}_{p},K\right)$, $1/\chi\in\tilde{C}\left(\mathbb{Z}_{p},K\right)$
force $1/\chi$ to be an element of $C\left(\mathbb{Z}_{p},K\right)$,
which would be a contradiction.

Q.E.D. 
\begin{lem}
Let $\chi\in\tilde{C}\left(\mathbb{Z}_{p},K\right)$. Then, $1/\chi\in\tilde{C}\left(\mathbb{Z}_{p},K\right)$
if and only if $\chi_{N}^{-1}=1/\chi_{N}$ converges point-wise in
$K$ as $N\rightarrow\infty$, in which case, the point-wise limit
is equal to $1/\chi$. 
\end{lem}
Proof: Since $\chi_{N}^{-1}$ is the $N$th truncation of $1/\chi_{N}$,
we have that $1/\chi\in\tilde{C}\left(\mathbb{Z}_{p},K\right)$ if
and only if $\chi_{N}^{-1}$ converges point-wise to a limit, in which
case that limit is necessarily $1/\chi$.

Q.E.D.

\vphantom{}

This shows that we can simplify the study of units of $\tilde{C}\left(\mathbb{Z}_{p},K\right)$
by considering their truncations, which is very nice, seeing as the
$\left(p,q\right)$-adic continuity of an $N$th truncation allows
us to employ the techniques $\left(p,q\right)$-adic Fourier analysis.

Our next inversion criterion takes the form of a characterization
of the rate at which the truncations of a rising continuous function
converge in the limit, and has the appeal of being both necessary
\emph{and }sufficient.
\begin{lem}[\textbf{The Square Root Lemma}\index{Square Root Lemma}]
\label{lem:square root lemma}Let $\chi\in\tilde{C}\left(\mathbb{Z}_{p},K\right)$,
let $\mathfrak{z}\in\mathbb{Z}_{p}^{\prime}$, and let $\mathfrak{c}\in K\backslash\chi\left(\mathbb{N}_{0}\right)$.
Then, $\chi\left(\mathfrak{z}\right)=\mathfrak{c}$ if and only if:
\begin{equation}
\liminf_{n\rightarrow\infty}\frac{\left|\chi\left(\left[\mathfrak{z}\right]_{p^{n}}\right)-\mathfrak{c}\right|_{q}}{\left|\nabla_{p^{n}}\left\{ \chi\right\} \left(\mathfrak{z}\right)\right|_{q}^{1/2}}<\infty\label{eq:Square Root Lemma}
\end{equation}
\end{lem}
\begin{rem}
Since this is a biconditional statement, it is equivalent to its inverse,
which is: \emph{
\begin{equation}
\chi\left(\mathfrak{z}\right)\neq\mathfrak{c}\Leftrightarrow\lim_{n\rightarrow\infty}\frac{\left|\chi\left(\left[\mathfrak{z}\right]_{p^{n}}\right)-\mathfrak{c}\right|_{q}}{\left|\nabla_{p^{n}}\left\{ \chi\right\} \left(\mathfrak{z}\right)\right|_{q}^{1/2}}=\infty\label{eq:Square Root Lemma - Inverse}
\end{equation}
}This is because the $\liminf$
\[
\liminf_{n\rightarrow\infty}\frac{\left|\chi\left(\left[\mathfrak{z}\right]_{p^{n}}\right)-\mathfrak{c}\right|_{q}}{\left|\nabla_{p^{n}}\left\{ \chi\right\} \left(\mathfrak{z}\right)\right|_{q}^{1/2}}=\infty
\]
forces the true limit to exist and be equal to $\infty$. 
\end{rem}
Proof: Let $\chi$, $\mathfrak{z}$, and $\mathfrak{c}$ be as given.
We will prove (\ref{eq:Square Root Lemma - Inverse}).

I. Suppose $\chi\left(\mathfrak{z}\right)\neq\mathfrak{c}$. Then,
it must be that $\liminf_{n\rightarrow\infty}\left|\chi\left(\left[\mathfrak{z}\right]_{p^{n}}\right)-\mathfrak{c}\right|_{q}>0$.
On the other hand, since $\chi\in\tilde{C}\left(\mathbb{Z}_{p},K\right)$,
the denominator of (\ref{eq:Square Root Lemma - Inverse}) tends to
$0$. This shows (\ref{eq:Square Root Lemma - Inverse}) holds whenever
$\chi\left(\mathfrak{z}\right)\neq\mathfrak{c}$.

\vphantom{}

II. Suppose: 
\[
\lim_{n\rightarrow\infty}\frac{\left|\chi\left(\left[\mathfrak{z}\right]_{p^{n}}\right)-\mathfrak{c}\right|_{q}}{\left|\nabla_{p^{n}}\left\{ \chi\right\} \left(\mathfrak{z}\right)\right|_{q}^{1/2}}=\infty
\]
Now, letting $\mathfrak{y}\in\mathbb{Z}_{p}^{\prime}$ be arbitrary,
observe that we have the formal identity: 
\begin{align*}
\frac{1}{\chi\left(\mathfrak{y}\right)-\mathfrak{c}} & =\frac{1}{\chi\left(0\right)-\mathfrak{c}}+\sum_{n=1}^{\infty}\left(\frac{1}{\chi\left(\left[\mathfrak{y}\right]_{p^{n}}\right)-\mathfrak{c}}-\frac{1}{\chi\left(\left[\mathfrak{y}\right]_{p^{n-1}}\right)-\mathfrak{c}}\right)\\
 & =\frac{1}{\chi\left(0\right)-\mathfrak{c}}-\sum_{n=1}^{\infty}\frac{\overbrace{\chi\left(\left[\mathfrak{y}\right]_{p^{n}}\right)-\chi\left(\left[\mathfrak{y}\right]_{p^{n-1}}\right)}^{\nabla_{p^{n}}\left\{ \chi\right\} \left(\mathfrak{y}\right)}}{\left(\chi\left(\left[\mathfrak{y}\right]_{p^{n}}\right)-\mathfrak{c}\right)\left(\chi\left(\left[\mathfrak{y}\right]_{p^{n-1}}\right)-\mathfrak{c}\right)}
\end{align*}
Note that we needn't worry about $1/\left(\chi\left(0\right)-\mathfrak{c}\right)$,
since $\mathfrak{c}$ was given to not be of the form $\mathfrak{c}=\chi\left(n\right)$
for any $n\in\mathbb{N}_{0}$.

Now, observe that $\chi\left(\mathfrak{z}\right)\neq\mathfrak{c}$
if and only if the above series converges in $K$ at $\mathfrak{y}=\mathfrak{z}$.
This, in turn, occurs if and only if: 
\[
\lim_{n\rightarrow\infty}\left|\frac{\nabla_{p^{n}}\left\{ \chi\right\} \left(\mathfrak{z}\right)}{\left(\chi\left(\left[\mathfrak{z}\right]_{p^{n}}\right)-\mathfrak{c}\right)\left(\chi\left(\left[\mathfrak{z}\right]_{p^{n-1}}\right)-\mathfrak{c}\right)}\right|_{q}=0
\]
and hence, if and only if: 
\begin{equation}
\lim_{n\rightarrow\infty}\left|\frac{\left(\chi\left(\left[\mathfrak{z}\right]_{p^{n}}\right)-\mathfrak{c}\right)\left(\chi\left(\left[\mathfrak{z}\right]_{p^{n-1}}\right)-\mathfrak{c}\right)}{\nabla_{p^{n}}\left\{ \chi\right\} \left(\mathfrak{z}\right)}\right|_{q}=\infty\label{eq:Desired Limit}
\end{equation}
So, in order to show that (\ref{eq:Square Root Lemma - Inverse})
implies $\chi\left(\mathfrak{z}\right)\neq\mathfrak{c}$, we need
only show that (\ref{eq:Square Root Lemma - Inverse}) implies (\ref{eq:Desired Limit}).

Here, note that: 
\[
\chi\left(\left[\mathfrak{z}\right]_{p^{n-1}}\right)=\chi\left(\left[\mathfrak{z}\right]_{p^{n}}\right)-\nabla_{p^{n}}\left\{ \chi\right\} \left(\mathfrak{z}\right)
\]
and so: 
\begin{align*}
\left|\frac{\left(\chi\left(\left[\mathfrak{z}\right]_{p^{n}}\right)-\mathfrak{c}\right)\left(\chi\left(\left[\mathfrak{z}\right]_{p^{n-1}}\right)-\mathfrak{c}\right)}{\nabla_{p^{n}}\left\{ \chi\right\} \left(\mathfrak{z}\right)}\right|_{q} & =\left|\frac{\left(\chi\left(\left[\mathfrak{z}\right]_{p^{n}}\right)-\mathfrak{c}\right)\left(\chi\left(\left[\mathfrak{z}\right]_{p^{n}}\right)-\mathfrak{c}-\nabla_{p^{n}}\left\{ \chi\right\} \left(\mathfrak{z}\right)\right)}{\nabla_{p^{n}}\left\{ \chi\right\} \left(\mathfrak{z}\right)}\right|_{q}\\
 & =\frac{\left|\chi\left(\left[\mathfrak{z}\right]_{p^{n}}\right)-\mathfrak{c}\right|_{q}}{\left|\nabla_{p^{n}}\left\{ \chi\right\} \left(\mathfrak{z}\right)\right|_{q}^{1/2}}\frac{\left|\chi\left(\left[\mathfrak{z}\right]_{p^{n}}\right)-\mathfrak{c}-\nabla_{p^{n}}\left\{ \chi\right\} \left(\mathfrak{z}\right)\right|_{q}}{\left|\nabla_{p^{n}}\left\{ \chi\right\} \left(\mathfrak{z}\right)\right|_{q}^{1/2}}\\
 & \geq\frac{\left|\chi\left(\left[\mathfrak{z}\right]_{p^{n}}\right)-\mathfrak{c}\right|_{q}}{\left|\nabla_{p^{n}}\left\{ \chi\right\} \left(\mathfrak{z}\right)\right|_{q}^{1/2}}\left|\frac{\left|\chi\left(\left[\mathfrak{z}\right]_{p^{n}}\right)-\mathfrak{c}\right|_{q}}{\left|\nabla_{p^{n}}\left\{ \chi\right\} \left(\mathfrak{z}\right)\right|_{q}^{1/2}}-\left|\nabla_{p^{n}}\left\{ \chi\right\} \left(\mathfrak{z}\right)\right|_{q}^{1/2}\right|
\end{align*}

Now, since $\chi\in\tilde{C}\left(\mathbb{Z}_{p},K\right)$, we have
$\left|\nabla_{p^{n}}\left\{ \chi\right\} \left(\mathfrak{z}\right)\right|_{q}^{1/2}\rightarrow0$
as $n\rightarrow\infty$. On the other hand, for this part of the
proof, we assumed that: 
\begin{equation}
\lim_{n\rightarrow\infty}\frac{\left|\chi\left(\left[\mathfrak{z}\right]_{p^{n}}\right)-\mathfrak{c}\right|_{q}}{\left|\nabla_{p^{n}}\left\{ \chi\right\} \left(\mathfrak{z}\right)\right|_{q}^{1/2}}=\infty
\end{equation}
Consequently, since $\mathfrak{z}$ is fixed: 
\begin{align*}
\lim_{n\rightarrow\infty}\left|\frac{\left(\chi\left(\left[\mathfrak{z}\right]_{p^{n}}\right)-\mathfrak{c}\right)\left(\chi\left(\left[\mathfrak{z}\right]_{p^{n-1}}\right)-\mathfrak{c}\right)}{\nabla_{p^{n}}\left\{ \chi\right\} \left(\mathfrak{z}\right)}\right|_{q} & \geq\lim_{n\rightarrow\infty}\frac{\left|\chi\left(\left[\mathfrak{z}\right]_{p^{n}}\right)-\mathfrak{c}\right|_{q}}{\left|\nabla_{p^{n}}\left\{ \chi\right\} \left(\mathfrak{z}\right)\right|_{q}^{1/2}}\left|\frac{\left|\chi\left(\left[\mathfrak{z}\right]_{p^{n}}\right)-\mathfrak{c}\right|_{q}}{\left|\nabla_{p^{n}}\left\{ \chi\right\} \left(\mathfrak{z}\right)\right|_{q}^{1/2}}-\left|\nabla_{p^{n}}\left\{ \chi\right\} \left(\mathfrak{z}\right)\right|_{q}^{1/2}\right|\\
 & =\infty\cdot\left|\infty-0\right|\\
 & =\infty
\end{align*}
which is (\ref{eq:Desired Limit}).

Q.E.D.

\vphantom{}

This time, the moral is that $\mathfrak{c}\in K\backslash\chi\left(\mathbb{N}_{0}\right)$
is in the image of $\chi$ if and only if there is a subsequence of
$n$s along which $\left|\chi\left(\left[\mathfrak{z}\right]_{p^{n}}\right)-\mathfrak{c}\right|_{q}$
tends to $0$ \emph{at least as quickly as} $\left|\nabla_{p^{n}}\left\{ \chi\right\} \left(\mathfrak{z}\right)\right|_{q}^{1/2}$. 
\begin{prop}
Let $\chi\in\tilde{C}\left(\mathbb{Z}_{p},K\right)$, let $\mathfrak{z}\in\mathbb{Z}_{p}^{\prime}$,
and let $\mathfrak{c}\in K\backslash\chi\left(\mathbb{N}_{0}\right)$.
If: 
\[
\liminf_{n\rightarrow\infty}\frac{\left|\chi\left(\left[\mathfrak{z}\right]_{p^{n}}\right)-\mathfrak{c}\right|_{q}}{\left|\nabla_{p^{n}}\left\{ \chi\right\} \left(\mathfrak{z}\right)\right|_{q}^{1/2}}<\infty
\]
then: 
\[
\liminf_{n\rightarrow\infty}\frac{\left|\chi\left(\left[\mathfrak{z}\right]_{p^{n}}\right)-\mathfrak{c}\right|_{q}}{\left|\nabla_{p^{n}}\left\{ f\right\} \left(\mathfrak{z}\right)\right|_{q}^{1/2}}>0
\]
occurs if and only if, for any strictly increasing sequence $\left\{ n_{j}\right\} _{j\geq1}\subseteq\mathbb{N}_{0}$
so that: 
\[
\lim_{j\rightarrow\infty}\frac{\left|\chi\left(\left[\mathfrak{z}\right]_{p^{n_{j}}}\right)-\mathfrak{c}\right|_{q}}{\left|\nabla_{p^{n_{j}}}\left\{ \chi\right\} \left(\mathfrak{z}\right)\right|_{q}^{1/2}}=\liminf_{n\rightarrow\infty}\frac{\left|\chi\left(\left[\mathfrak{z}\right]_{p^{n}}\right)-\mathfrak{c}\right|_{q}}{\left|\nabla_{p^{n}}\left\{ \chi\right\} \left(\mathfrak{z}\right)\right|_{q}^{1/2}}
\]
the congruences: 
\begin{equation}
v_{q}\left(\nabla_{p^{n_{j+1}}}\left\{ \chi\right\} \left(\mathfrak{z}\right)\right)\overset{2}{\equiv}v_{q}\left(\nabla_{p^{n_{j}}}\left\{ \chi\right\} \left(\mathfrak{z}\right)\right)
\end{equation}
hold for all sufficiently large $j$; that is, the parity of the $q$-adic
valuation of $\nabla_{p^{n_{j}}}\left\{ \chi\right\} \left(\mathfrak{z}\right)$
becomes constant for all sufficiently large $j$. 
\end{prop}
Proof: Suppose: 
\[
\liminf_{n\rightarrow\infty}\frac{\left|\chi\left(\left[\mathfrak{z}\right]_{p^{n}}\right)-\mathfrak{c}\right|_{q}}{\left|\nabla_{p^{n}}\left\{ \chi\right\} \left(\mathfrak{z}\right)\right|_{q}^{1/2}}>0
\]
In particular, let us write the value of this $\liminf$ as $q^{C}$,
where $C$ is a positive real number. Letting the $n_{j}$s be any
sequence as described above, by way of contradiction, suppose there
are two disjoint subsequences of $j$s so that $v_{q}\left(\nabla_{p^{n_{j}}}\left\{ \chi\right\} \left(\mathfrak{z}\right)\right)$
is always even for all $j$ in the first subsequence and is always
odd for all $j$ in the second subsequence. Now, for \emph{any} $j$,
let $v_{j}$ be the positive integer so that $\left|\chi\left(\left[\mathfrak{z}\right]_{p^{n_{j}}}\right)-\mathfrak{c}\right|_{q}=q^{-v_{j}}$.
Then: 
\begin{align*}
q^{C}=\lim_{j\rightarrow\infty}q^{\frac{1}{2}v_{q}\left(\nabla_{p^{n_{j}}}\left\{ \chi\right\} \left(\mathfrak{z}\right)\right)}\left|\chi\left(\left[\mathfrak{z}\right]_{p^{n_{j}}}\right)-\mathfrak{c}\right|_{q} & =\lim_{j\rightarrow\infty}q^{\frac{1}{2}v_{q}\left(\nabla_{p^{n_{j}}}\left\{ \chi\right\} \left(\mathfrak{z}\right)\right)-v_{j}}
\end{align*}
and hence: 
\begin{equation}
C=\lim_{j\rightarrow\infty}\left(\frac{1}{2}v_{q}\left(\nabla_{p^{n_{j}}}\left\{ \chi\right\} \left(\mathfrak{z}\right)\right)-v_{j}\right)\label{eq:Limit formula for K}
\end{equation}
for \emph{both }subsequences of $j$s. For the $j$s in the first
subsequence (where $v_{q}\left(\nabla_{p^{n_{j}}}\left\{ \chi\right\} \left(\mathfrak{z}\right)\right)$
is always even), the limit on the right hand side is of a subsequence
in the set $\mathbb{Z}$. For the $j$s in the second subsequence
(where $v_{q}\left(\nabla_{p^{n_{j}}}\left\{ \chi\right\} \left(\mathfrak{z}\right)\right)$
is always odd), however, the limit on the right hand side is that
of a sequence in $\mathbb{Z}+\frac{1}{2}$. But this is impossible:
it forces $C$ to be an element of both $\mathbb{Z}$ and $\mathbb{Z}+\frac{1}{2}$,
and those two sets are disjoint.

Thus, the parity of $v_{q}\left(\nabla_{p^{n_{j}}}\left\{ \chi\right\} \left(\mathfrak{z}\right)\right)$
must become constant for all sufficiently large $j$ whenever $C>-\infty$.

Q.E.D.

\vphantom{}

Given the criteria described in these last few results, it seems natural
to introduce the following definitions: 
\begin{defn}
Let $\chi\in\tilde{C}\left(\mathbb{Z}_{p},K\right)$. A pair $\left(\mathfrak{z},\mathfrak{c}\right)$,
where $\mathfrak{z}\in\mathbb{Z}_{p}^{\prime}$ and $\mathfrak{c}\in\chi\left(\mathbb{Z}_{p}\right)\backslash\chi\left(\mathbb{N}_{0}\right)$
is said to be a \index{quick approach}\textbf{quick approach of $\chi$
/ quickly approached by $\chi$ }whenever: 
\[
\liminf_{n\rightarrow\infty}\frac{\left|\chi\left(\left[\mathfrak{z}\right]_{p^{n}}\right)-\mathfrak{c}\right|_{q}}{\left|\nabla_{p^{n}}\left\{ \chi\right\} \left(\mathfrak{z}\right)\right|_{q}^{1/2}}=0
\]
We say that the pair is a\textbf{ slow approach of $\chi$ / slowly
approached by $\chi$ }whenever the $\liminf$ is finite and non-zero. 
\end{defn}
\begin{question}
Let $\chi\in\tilde{C}\left(\mathbb{Z}_{p},K\right)$. What can be
said about the speed of approach of $\chi$ to $\mathfrak{c}$ at
any given $\mathfrak{z}\in\mathbb{Z}_{p}^{\prime}$? Can there be
infinitely many (or at most finitely many\textemdash or none?) points
of quick approach? Slow approach? Are there other properties we can
use to characterize points of quick or slow approach?\footnote{In Chapter 4, we will show that, for every semi-basic, contracting
$2$-Hydra map which fixes $0$, every point in $\chi_{H}\left(\mathbb{Z}_{p}\right)\backslash\chi_{H}\left(\mathbb{N}_{0}\right)$
is a point of quick approach $\chi_{H}$. The more general case of
a $p$-Hydra map likely also holds, but I did not have sufficient
time to investigate it fully.} 
\end{question}
\begin{rem}
To share a bit of my research process, although I discovered the notion
of quasi-integrability back in the second half of 2020, it wasn't
until the 2021 holiday season that I began to unify my $\left(p,q\right)$-adic
work in earnest. With the Correspondence Principle on my mind, I sought
to find ways to characterize the image of a $\left(p,q\right)$-adic
function. In particular, I found myself gravitating toward the question:
\emph{given a rising-continuous function $\chi:\mathbb{Z}_{p}\rightarrow\mathbb{Z}_{q}$,
when is $1/\chi$ rising-continuous? }Given a $\chi_{H}$, investigating
the singular behavior of $1/\left(\chi_{H}\left(\mathfrak{z}\right)-x\right)$
where $x\in\mathbb{Z}\backslash\left\{ 0\right\} $ would seem to
be a natural way to exploit the Correspondence Principle, because
an $x\in\mathbb{Z}\backslash\left\{ 0\right\} $ would be a periodic
point of $H$ if and only if $1/\left(\chi_{H}\left(\mathfrak{z}\right)-x\right)$
had a singularity for some $\mathfrak{z}\in\mathbb{Z}_{p}^{\prime}$.
I drew inspiration from W. Cherry's notes on non-archimedean function
theory \cite{Cherry non-archimedean function theory notes}, and the
existence of a $p$-adic Nevanlinna theory\textemdash Nevanlinna theory
being one of analytic function theorists' most powerful tools for
keeping track of the poles of meromorphic functions.

Due to the (present) lack of a meaningful notion of an \emph{analytic
}$\left(p,q\right)$-adic function, however, I had to go down another
route to investigate reciprocals of rising-continuous functions. Things
began to pick up speed once I realized that the set of rising-continuous
functions could be realized as a non-archimedean Banach algebra under
point-wise multiplication. In that context, exploring the Correspondence
Principle becomes extremely natural: it is the quest to understand
those values of $x\in\mathbb{Z}\backslash\left\{ 0\right\} $ for
which $\chi_{H}\left(\mathfrak{z}\right)-x$ is not a unit of the
Banach algebra of rising-continuous functions. 
\end{rem}

\subsubsection{A Very Brief Foray Into Berkovitch Space\index{Berkovitch space}}

To give yet another reason why $\left(p,q\right)$-adic analysis has
been considered ``uninteresting'' by the mathematical community
at large: their zeroes do not enjoy any special properties. Case in
point: 
\begin{prop}
\label{prop:(p,q)-adic functions are "uninteresting"}Let $A$ be
any closed subset of $\mathbb{Z}_{p}$. Then, there is an $f\in C\left(\mathbb{Z}_{p},K\right)$
such that $f\left(\mathfrak{z}\right)=0$ if and only if $\mathfrak{z}\in A$. 
\end{prop}
Proof: Let $A$ be a closed subset of $\mathbb{Z}_{p}$. For any $\mathfrak{z}\in\mathbb{Z}_{p}\backslash A$,
the Hausdorff distance: 
\begin{equation}
d\left(\mathfrak{z},A\right)\overset{\textrm{def}}{=}\inf_{\mathfrak{a}\in A}\left|\mathfrak{z}-\mathfrak{a}\right|_{p}\label{eq:Definition of the distance of a point in Z_p from a set in Z_p}
\end{equation}
of $\mathfrak{z}$ form $A$ will be of the form $\frac{1}{p^{n}}$
for some $n\in\mathbb{N}_{0}$. Consequently, we can define a function
$g:\mathbb{Z}_{p}\backslash A\rightarrow\mathbb{N}_{0}$ so that:
\begin{equation}
d\left(x,A\right)=p^{-g\left(x\right)},\textrm{ }\forall x\in\mathbb{Z}_{p}\backslash A
\end{equation}
Then, the function: 
\begin{equation}
f\left(x\right)\overset{\textrm{def}}{=}\left[x\notin A\right]q^{g\left(x\right)}\label{eq:Function with zeroes on an arbitrary closed subset of Z_p}
\end{equation}
will be continuous, and will have $A$ as its vanishing set.

Q.E.D.

\vphantom{}

In this short sub-subsection, we deal with some topological properties
of rising-continuous functions which can be proven using \textbf{Berkovitch
spaces}. For the uninitiated\textemdash such as myself\textemdash Berkovitch
space is the capstone of a series of turbulent efforts by algebraic
and arithmetic geometers in the second half of the twentieth century
which had as their \emph{idée-fixe }the hope of circumventing the
Achilles' heel of the native theory of $\left(p,p\right)$-adic\footnote{i.e., $f:\mathbb{C}_{p}\rightarrow\mathbb{C}_{p}$.}
analytic functions: the abysmal failure of Weierstrass' germ-based
notions of analytic continuation\index{analytic continuation}. Because
of the $p$-adics' ultrametric structure, every $\mathfrak{z}$ in
the disk $\left\{ \mathfrak{z}\in\mathbb{C}_{p}:\left|\mathfrak{z}-\mathfrak{a}\right|_{p}<r\right\} $
is, technically, at the center of said disk; here $\mathfrak{a}\in\mathbb{C}_{p}$
and $r>0$. Along with the $p$-adics' stark convergence properties
(viz. the principles of ultrametric analysis on \pageref{fact:Principles of Ultrametric Analysis})
this strange property of non-archimedean disks leads to the following
unexpected behaviors in a $p$-adic power series $f\left(\mathfrak{z}\right)=\sum_{n=0}^{\infty}\mathfrak{a}_{n}\mathfrak{z}^{n}$.
Here, $r$ denotes $f$'s radius of convergence $r>0$.
\begin{enumerate}
\item Either $f\left(\mathfrak{z}\right)$ converges at every $\mathfrak{z}\in\mathbb{C}_{p}$
with $\left|\mathfrak{z}\right|_{p}=r$, or diverges at every such
$\mathfrak{z}$. 
\item Developing $f\left(\mathfrak{z}\right)$ into a power series at any
$\mathfrak{z}_{0}\in\mathbb{C}_{p}$ with $\left|\mathfrak{z}_{0}\right|_{p}<r$
results in a power series whose radius of convergence is equal to
$r$. So, no analytic continuation is possible purely by way of rejiggering
$f$'s power series. 
\end{enumerate}
The construction of Berkovitch space provides a larger space containing
$\mathbb{C}_{p}$ which is Hausdorff, and\textemdash crucially, unlike
the totally disconnected space $\mathbb{C}_{p}$\textemdash is actually
arc-wise connected.\footnote{Let $X$ be a topological space. Two points $x,y\in X$ are said to
be \textbf{topologically distinguishable }if there exists an open
set $U\subseteq X$ which contains one and only one of $x$ or $y$.
$X$ is then said to be \textbf{arc-wise connected }whenever, for
any topologically distinguishable points $x,y\in X$, there is an
embedding $\alpha:\left[0,1\right]\rightarrow X$ (continuous and
injective) such that $\alpha\left(0\right)=x$ and $\alpha\left(1\right)=y$;
such an $\alpha$ is called an \textbf{arc}.} The idea\footnote{I base my exposition here on what is given in \cite{The Berkovitch Space Paper}.}
(a minor bit of trepanation) is to replace points in $\mathbb{C}_{p}$
with \index{semi-norm}\emph{semi-norms}. This procedure holds, in
fact, for an arbitrary Banach algebra over a non-archimedean field.
Before going forward, let me mention that I include the following
definitions (taken from \cite{The Berkovitch Space Paper}) primarily
for the reader's edification, provided the reader is the sort of person
who finds stratospheric abstractions edifying. Meanwhile, the rest
of us mere mortals can safely skip over to the proofs at the end of
this sub-subsection. 
\begin{defn}
Let $K$ be an algebraically closed non-archimedean field, complete
with respect to its ultrametric. Let $A$ be a commutative $K$-Banach
algebra with identity $1_{A}$, and with norm $\left\Vert \cdot\right\Vert $.
A function $\rho:A\times A\rightarrow\left[0,\infty\right)$ is an
\textbf{algebra semi-norm }if\index{Banach algebra!semi-norm}, in
addition to being a semi-norm on $A$ when $A$ is viewed as a $K$-Banach
space, we have that:

\vphantom{}

I. $\rho\left(1_{A}\right)=1$;

\vphantom{}

II. $\rho\left(f\cdot g\right)\leq\rho\left(f\right)\cdot\rho\left(g\right)$
for all $f,g\in A$.

$\rho$ is said to be \textbf{multiplicative} whenever\index{semi-norm!multiplicative}
(II) holds with equality for all $f,g\in A$.

The \textbf{multiplicative spectrum}\index{Banach algebra!multiplicative spectrum}\textbf{
of $A$}, denoted $\mathcal{M}\left(A,\left\Vert \cdot\right\Vert \right)$,
is the set of multiplicative semi-norms on $A$ which are continuous
with respect to the topology defined by $\left\Vert \cdot\right\Vert $;
i.e., for each $\rho\in\mathcal{M}\left(A,\left\Vert \cdot\right\Vert \right)$,
there is a real constant $K\geq0$ so that: 
\begin{equation}
\rho\left(f\right)\leq K\left\Vert f\right\Vert ,\textrm{ }\forall f\in A\label{eq:multiplicative semi-norm thing for Berkovitch space}
\end{equation}
\end{defn}
\vphantom{}

Using this bizarre setting, algebraic geometers then construct affine
spaces over $K$. 
\begin{defn}
We write $\mathbb{A}_{K}^{1}$ to denote the space of all multiplicative
semi-norms on the space $K\left[X\right]$ of formal polynomials in
the indeterminate $X$ with coefficients in $K$. $\mathbb{A}_{K}^{1}$
is called the \textbf{affine line over $K$}\index{affine line over K@affine line over \textbf{$K$}}. 
\end{defn}
\begin{fact}
Both $\mathcal{M}\left(A,\left\Vert \cdot\right\Vert \right)$ and
$\mathbb{\mathbb{A}}_{K}^{1}$ are made into topological spaces with
the topology of point-wise convergence. 
\end{fact}
\begin{thm}
$\mathcal{M}\left(A,\left\Vert \cdot\right\Vert \right)$ is Hausdorff
and compact. $\mathbb{A}_{K}^{1}$ is Hausdorff, locally compact,
and arc-wise connected \cite{The Berkovitch Space Paper}. 
\end{thm}
\begin{defn}
For each $t\in A$, write $\theta_{t}$ to denote the $K$-algbra
homomorphism $\theta_{t}:K\left[X\right]\longrightarrow A$ defined
by $\theta_{t}\left(X\right)=t$.

For each $t\in A$, the \index{Gelfand transform}\textbf{Gelfand
transform of $t$} is then defined to be the map $\mathcal{G}_{t}:\mathcal{M}\left(A,\left\Vert \cdot\right\Vert \right)\rightarrow\mathbb{A}_{K}^{1}$
given by: 
\begin{equation}
\mathcal{G}_{t}\left\{ \rho\right\} \overset{\textrm{def}}{=}\rho\circ\theta_{t}\label{Gelfand transform of t}
\end{equation}
$\mathcal{G}_{t}\left\{ \rho\right\} $ is then a multiplicative semi-norm
on $K\left[X\right]$ for all $\rho$ and $t$. The map $\mathcal{G}$
which sends $t\in A$ to the map $\mathcal{G}_{t}$ is called the
\textbf{Gelfand transform}. 
\end{defn}
\begin{thm}
$\mathcal{G}_{t}$ is continuous for every $t\in A$. 
\end{thm}
\vphantom{}

The point of this, apparently, is that $\mathcal{G}$ allows any $t\in A$
to be interpreted as the continuous function $\mathcal{G}_{t}$. Now
for some more definitions: 
\begin{defn}
Given $t\in A$:

\vphantom{}

I. We write $\textrm{s}\left(t\right)$ to denote the \textbf{scalar
spectrum }of $t$: the set of $\lambda\in K$ for which $t-\lambda1_{A}$
is not invertible in $A$.

\vphantom{}

II. We write $\sigma\left(t\right)$ to denote the \textbf{spectrum
}of $t$: the image of $\mathcal{M}\left(A,\left\Vert \cdot\right\Vert \right)$
under $\mathcal{G}_{t}$. 
\end{defn}
\vphantom{}

Since $K$ is a Banach algebra over itself (with $\left|\cdot\right|_{K}$
as the norm), $\mathcal{M}\left(K,\left|\cdot\right|_{K}\right)$
consists of a single semi-norm: the absolute value on $K$. This allows
$K$ to be embedded in $\mathbb{A}_{K}^{1}$ by way of the Gelfand
transform $\mathcal{G}$, identifying any $\mathfrak{a}\in K$ with
$\mathcal{G}_{\mathfrak{a}}$. Since $\mathcal{M}\left(K,\left|\cdot\right|_{K}\right)$
contains only the single element $\left|\cdot\right|_{K}$, the image
of $\mathcal{G}_{\mathfrak{a}}$ ($\sigma\left(\mathfrak{a}\right)$\textemdash the
spectrum of $\mathfrak{a}$) is the semi-norm which sends $P\in K\left[X\right]$
to $\left|\theta_{\mathfrak{a}}\left(P\right)\right|_{K}$. That is,
for any $\mathfrak{a}\in K$, we we identify with $\mathfrak{a}$
the semi-norm $P\mapsto\left|\theta_{\mathfrak{a}}\left(P\right)\right|_{K}$
on $K\left[X\right]$. 
\begin{defn}
The \textbf{closed disk }in $\mathbb{A}_{K}^{1}$ centered at $\mathfrak{a}\in K$
of radius $r>0$ is defined as the closure in $\mathbb{A}_{K}^{1}$
(with respect to point-wise convergence) of the set of semi-norms
$P\in K\left[X\right]\mapsto\left|\theta_{\mathfrak{z}}\left(P\right)\right|_{K}$
as $\mathfrak{z}$ varies over all elements of $K$ with $\left|\mathfrak{z}-\mathfrak{a}\right|_{K}\leq r$. 
\end{defn}
\begin{thm}[\textbf{The First Spectral Radius Formula }\cite{The Berkovitch Space Paper}]
Let $A$ be a commutative $K$-Banach algebra with identity $1_{A}$,
and with norm $\left\Vert \cdot\right\Vert $, and let $t\in A$.
Then:

\vphantom{}

I. $\sigma\left(t\right)$ is non-empty and compact.

\vphantom{}

II. The smallest closed disk in $\mathbb{A}_{K}^{1}$ centered at
$0$ which contains $\sigma\left(t\right)$ has radius: 
\begin{equation}
\lim_{n\rightarrow\infty}\left\Vert t^{n}\right\Vert ^{1/n}
\end{equation}

\vphantom{}

III. $K\cap\sigma\left(t\right)=s\left(t\right)$, where by $K$,
we mean the copy of $K$ embedded in $\mathbb{A}_{K}^{1}$ as described
above.
\end{thm}
\vphantom{}

We now recall some basic definitions from topology: 
\begin{defn}
Let $X$ and $Y$ be topological spaces. A map $M:X\rightarrow Y$
is said to be:

\vphantom{}

I.\textbf{ Proper},\textbf{ }if $M^{-1}\left(C\right)$ is compact
in $X$ for all compact sets $C\subseteq Y$.

\vphantom{}

II. \textbf{Closed}, if $M\left(S\right)$ is closed in $Y$ for all
closed sets $S\subseteq X$. 
\end{defn}
\begin{fact}
Let $X$ and $Y$ be topological spaces, and let $M:X\rightarrow Y$
be continuous. If $X$ is compact and $Y$ is Hausdorff, then $M$
is both proper and closed. 
\end{fact}
\begin{prop}
$s\left(t\right)$ is compact in $K$. 
\end{prop}
Proof: Let $t\in A$ be arbitrary. Since $\mathcal{G}_{t}:\mathcal{M}\left(A,\left\Vert \cdot\right\Vert \right)\rightarrow\mathbb{A}_{K}^{1}$
is continuous, since $\mathcal{M}\left(A,\left\Vert \cdot\right\Vert \right)$
is compact, and since $\mathbb{A}_{K}^{1}$ is Hausdorff, it then
follows that $\mathcal{G}_{t}$ is both proper and closed. Consequently,
since $K$ is closed, so is its copy in $\mathbb{A}_{K}^{1}$. Since
$\sigma\left(t\right)$ is compact in $\mathbb{A}_{K}^{1}$, this
means that $s\left(t\right)=K\cap\sigma\left(t\right)$ is compact
in $\mathbb{A}_{K}^{1}\cap K$, and hence in $K$.

Q.E.D.

\vphantom{}

The non-algebraic reader can now resume paying attention. 
\begin{lem}
\label{lem:compactness of the image of a rising-continuous function}Let
$\chi\in\tilde{C}\left(\mathbb{Z}_{p},K\right)$. Then, $\chi\left(\mathbb{Z}_{p}\right)$
is a compact subset of $K$. 
\end{lem}
Proof: First, we embed $\chi$ in $\tilde{C}\left(\mathbb{Z}_{p},\mathbb{C}_{q}\right)$.
Since $\tilde{C}\left(\mathbb{Z}_{p},\mathbb{C}_{q}\right)$ is a
commutative Banach algebra over the algebraically closed, metrically
complete non-archimedean field $\mathbb{C}_{q}$, it follows by the
previous proposition that the scalar spectrum of $\chi$ is compact
in $\mathbb{C}_{q}$. Since the scalar spectrum is precisely the set
of $\mathfrak{c}\in\mathbb{C}_{q}$ for which $\frac{1}{\chi\left(\mathfrak{z}\right)-\mathfrak{c}}\in\tilde{C}\left(\mathbb{Z}_{p},\mathbb{C}_{q}\right)$,
it then follows that the scalar spectrum of $\chi$\textemdash a compact
set in $\mathbb{C}_{q}$\textemdash is equal to $\chi\left(\mathbb{Z}_{p}\right)$,
the image of $\chi$ over $\mathbb{Z}_{p}$.

Since $K$ is closed in $\mathbb{C}_{q}$, the set $\chi\left(\mathbb{Z}_{p}\right)\cap K$
(an intersection of a compact set and a closed set) is compact in
both $\mathbb{C}_{q}$ and $K$. Since, by definition of $K$, $\chi\left(\mathbb{Z}_{p}\right)=\chi\left(\mathbb{Z}_{p}\right)\cap K$,
we have that $\chi\left(\mathbb{Z}_{p}\right)$ is compact in $K$.

Q.E.D.

\vphantom{}

Before we get to our next results, note the following: 
\begin{fact}
Let $U$ be any clopen subset of $\mathbb{Z}_{p}$. Then, the indicator
function for $U$ is continuous as a function from $\mathbb{Z}_{p}$
to any topological ring (where the ``$1$'' taken by the function
over $U$ is the multiplicative identity element of the ring). 
\end{fact}
\begin{thm}
Let $\chi\in\tilde{C}\left(\mathbb{Z}_{p},K\right)$. Then, $\chi\left(U\right)$
is a compact subset of $K$ for all clopen sets $U\subseteq\mathbb{Z}_{p}$. 
\end{thm}
Proof: Let $U$ be any clopen subset of $\mathbb{Z}_{p}$ with non-empty
complement. Then, let $V\overset{\textrm{def}}{=}\mathbb{Z}_{p}\backslash U$
be the complement of $U$, and let $\mathfrak{a}\in\chi\left(U\right)$
be any value attained by $\chi$ on $U$. Then, define: 
\begin{equation}
f\left(\mathfrak{z}\right)\overset{\textrm{def}}{=}\left[\mathfrak{x}\in U\right]\chi\left(\mathfrak{z}\right)+\mathfrak{a}\left[\mathfrak{x}\in V\right]\label{eq:Berkovitch space f construction}
\end{equation}
Since both $V$ and $U$ are clopen sets, their indicator functions
$\left[\mathfrak{x}\in U\right]$ and $\left[\mathfrak{x}\in V\right]$
are continuous functions from $\mathbb{Z}_{p}$ to $K$. Since $\tilde{C}\left(\mathbb{Z}_{p},K\right)$
is an algebra, we have that $f\in\tilde{C}\left(\mathbb{Z}_{p},K\right)$.
As such, by \textbf{Lemma \ref{lem:compactness of the image of a rising-continuous function}},
$f\left(\mathbb{Z}_{p}\right)$ is compact in $K$. Observing that:
\begin{align*}
f\left(\mathbb{Z}_{p}\right) & =f\left(U\right)\cup f\left(V\right)\\
 & =\chi\left(U\right)\cup\left\{ \mathfrak{a}\right\} \\
\left(\mathfrak{a}\in\chi\left(U\right)\right); & =\chi\left(U\right)
\end{align*}
we then have that $\chi\left(U\right)$ is compact in $K$, as desired.

Q.E.D. 
\begin{thm}
Let $\chi\in\tilde{C}\left(\mathbb{Z}_{p},K\right)$, and suppose
there is a set $S\subseteq\mathbb{Z}_{p}$ and a real constant $v\geq0$
so that $\left|\chi\left(\mathfrak{z}\right)\right|_{q}>q^{-v}$ for
all $\mathfrak{z}\in S$. Then, $\chi\left(U\right)$ is compact in
$K$ for any set $U\subseteq S$ which is closed in $\mathbb{Z}_{p}$. 
\end{thm}
\begin{rem}
In other words, if $\chi$ is bounded away from zero on some set $S$,
then the restriction of $\chi$ to $S$ is a proper map. 
\end{rem}
Proof: Since $U$ is closed in $\mathbb{Z}_{p}$, as we saw in \textbf{Proposition
\ref{prop:(p,q)-adic functions are "uninteresting"}}, there is a
function $g:\mathbb{Z}_{p}\rightarrow\mathbb{N}_{0}$ so that $h:\mathbb{Z}_{p}\rightarrow K$
defined by: 
\begin{equation}
h\left(\mathfrak{z}\right)\overset{\textrm{def}}{=}\left[\mathfrak{z}\notin U\right]q^{g\left(\mathfrak{z}\right)}\label{eq:h construction for Berkovitch space}
\end{equation}
is continuous and vanishes if and only if $\mathfrak{z}\in U$. Then,
letting: 
\begin{equation}
f\left(\mathfrak{z}\right)\overset{\textrm{def}}{=}\left[\mathfrak{x}\in U\right]\chi\left(\mathfrak{z}\right)+q^{v}h\left(\mathfrak{z}\right)\label{eq:Second f construction for Berkovitch space}
\end{equation}
we have that $f\in\tilde{C}\left(\mathbb{Z}_{p},K\right)$, and that:
\begin{equation}
f\left(\mathbb{Z}_{p}\right)=\chi\left(U\right)\cup q^{v+g\left(V\right)}
\end{equation}
is compact in $K$, where $q^{v+g\left(V\right)}=\left\{ q^{v+n}:n\in g\left(V\right)\right\} $.
Since $\left|\chi\left(\mathfrak{z}\right)\right|_{q}>q^{-v}$ for
all $\mathfrak{z}\in S$, we have; 
\begin{equation}
\chi\left(U\right)\subseteq\chi\left(S\right)\subseteq K\backslash q^{v}\mathbb{Z}_{q}
\end{equation}
As such: 
\begin{align*}
\left(K\backslash q^{v}\mathbb{Z}_{q}\right)\cap f\left(\mathbb{Z}_{p}\right) & =\left(\chi\left(U\right)\cup q^{v+g\left(V\right)}\right)\cap\left(K\backslash q^{v}\mathbb{Z}_{q}\right)\\
 & =\left(\chi\left(U\right)\cap\left(K\backslash q^{v}\mathbb{Z}_{q}\right)\right)\cup\left(q^{v+g\left(V\right)}\cap K\cap\left(q^{v}\mathbb{Z}_{q}\right)^{c}\right)\\
 & =\underbrace{\left(\chi\left(U\right)\cap\left(K\backslash q^{v}\mathbb{Z}_{q}\right)\right)}_{\chi\left(U\right)}\cup\underbrace{\left(q^{v+g\left(V\right)}\cap K\cap\left(q^{v}\mathbb{Z}_{q}\right)^{c}\right)}_{\varnothing}\\
\left(q^{v+g\left(V\right)}\cap\left(q^{v}\mathbb{Z}_{q}\right)^{c}=\varnothing\right); & =\chi\left(U\right)
\end{align*}
Since $K\backslash q^{v}\mathbb{Z}_{q}$ is closed and $f\left(\mathbb{Z}_{p}\right)$
is compact, $\left(K\backslash q^{v}\mathbb{Z}_{q}\right)\cap f\left(\mathbb{Z}_{p}\right)=\chi\left(U\right)$
is compact, as desired.

Q.E.D. 
\begin{lem}
$f\in C\left(\mathbb{Z}_{p},K\right)$ is a unit if and only if $f$
has no zeroes. 
\end{lem}
Proof:

I. If $f$ is a unit, then $\frac{1}{f}$ is continuous, and hence,
$f$ can have no zeroes.

\vphantom{}

II. Conversely, suppose $f$ has no zeroes. Since $f$ is continuous,
its image $f\left(\mathbb{Z}_{p}\right)$ is a compact subset of $K$
(even if $K$ itself is \emph{not} locally compact). Since $K$ is
a topological field, the reciprocation map $\iota\left(\mathfrak{y}\right)\overset{\textrm{def}}{=}\frac{1}{\mathfrak{y}}$
is a self-homeomorphism of $K\backslash\left\{ 0\right\} $. Here,
we observe that the pre-image of a set $Y$ under $1/f$ is then equal
to $f^{-1}\left(\iota^{-1}\left(Y\right)\right)$.

So, let $Y$ be an arbitrary closed subset of $K$. Then, the pre-image
of $Y$ under $1/f$ is equal to: 
\begin{equation}
f^{-1}\left(\iota^{-1}\left(Y\right)\right)=f^{-1}\left(\iota^{-1}\left(Y\right)\cap f\left(\mathbb{Z}_{p}\right)\right)
\end{equation}
Since $f\left(\mathbb{Z}_{p}\right)$ is compact and bounded away
from zero: 
\begin{equation}
\iota^{-1}\left(Y\right)\cap f\left(\mathbb{Z}_{p}\right)
\end{equation}
is closed and bounded away from zero, and so, by the continuity of
$f$, $f^{-1}\left(\iota^{-1}\left(Y\right)\cap f\left(\mathbb{Z}_{p}\right)\right)$
is closed. Thus, the pre-image of every closed subset of $K$ under
$1/f$ is closed; this proves that $1/f$ is continuous whenever $f$
has no zeroes.

Q.E.D. 
\begin{lem}
\label{lem:3.17}$\chi\in\tilde{C}\left(\mathbb{Z}_{p},K\right)$
is a unit of $\tilde{C}\left(\mathbb{Z}_{p},K\right)$ whenever the
following conditions hold:

\vphantom{}

I. $\chi$ has no zeroes.

\vphantom{}

II. For each $\mathfrak{z}\in\mathbb{Z}_{p}$: 
\begin{equation}
\sup_{n\geq1}\left|\frac{\chi\left(\left[\mathfrak{z}\right]_{p^{n}}\right)}{\chi\left(\left[\mathfrak{z}\right]_{p^{n-1}}\right)}\right|_{q}<\infty
\end{equation}
Note that this bound need not be uniform in $\mathfrak{z}$. 
\end{lem}
Proof: Let $\chi\in\tilde{C}\left(\mathbb{Z}_{p},K\right)$, and suppose
$0\notin\chi\left(\mathbb{Z}_{p}\right)$. Then, by the van der Put
identity, we have that $\frac{1}{\chi}$ will be rising-continuous
if and only if: 
\[
\lim_{n\rightarrow\infty}\left|\frac{1}{\chi\left(\left[\mathfrak{z}\right]_{p^{n}}\right)}-\frac{1}{\chi\left(\left[\mathfrak{z}\right]_{p^{n-1}}\right)}\right|_{q}=0
\]
for all $\mathfrak{z}\in\mathbb{Z}_{p}$, which is the same as: 
\[
\lim_{n\rightarrow\infty}\left|\frac{\nabla_{p^{n}}\left\{ \chi\right\} \left(\mathfrak{z}\right)}{\chi\left(\left[\mathfrak{z}\right]_{p^{n}}\right)\chi\left(\left[\mathfrak{z}\right]_{p^{n-1}}\right)}\right|_{q}=0
\]

By the \textbf{Square-Root Lemma} (\textbf{Lemma \ref{lem:square root lemma}}),
since $0$ is not in the image of $\chi$, it must be that: 
\begin{equation}
\liminf_{n\rightarrow\infty}\frac{\left|\chi\left(\left[\mathfrak{z}\right]_{p^{n}}\right)\right|_{q}}{\left|\nabla_{p^{n}}\left\{ \chi\right\} \left(\mathfrak{z}\right)\right|_{q}^{1/2}}=\infty
\end{equation}
which then forces: 
\begin{equation}
\lim_{n\rightarrow\infty}\frac{\left|\chi\left(\left[\mathfrak{z}\right]_{p^{n}}\right)\right|_{q}}{\left|\nabla_{p^{n}}\left\{ \chi\right\} \left(\mathfrak{z}\right)\right|_{q}^{1/2}}=\infty
\end{equation}
which forces: 
\begin{equation}
\lim_{n\rightarrow\infty}\frac{\left|\nabla_{p^{n}}\left\{ \chi\right\} \left(\mathfrak{z}\right)\right|_{q}^{1/2}}{\left|\chi\left(\left[\mathfrak{z}\right]_{p^{n}}\right)\right|_{q}}=0
\end{equation}
Now, the condition required for $1/\chi$ to be rising-continuous
is: 
\begin{equation}
\lim_{n\rightarrow\infty}\left|\frac{\nabla_{p^{n}}\left\{ \chi\right\} \left(\mathfrak{z}\right)}{\chi\left(\left[\mathfrak{z}\right]_{p^{n}}\right)\chi\left(\left[\mathfrak{z}\right]_{p^{n-1}}\right)}\right|_{q}=0
\end{equation}
This can be re-written as: 
\[
\lim_{n\rightarrow\infty}\frac{\left|\nabla_{p^{n}}\left\{ \chi\right\} \left(\mathfrak{z}\right)\right|_{q}^{1/2}}{\left|\chi\left(\left[\mathfrak{z}\right]_{p^{n}}\right)\right|_{q}}\cdot\frac{\left|\chi\left(\left[\mathfrak{z}\right]_{p^{n}}\right)\right|_{q}}{\left|\chi\left(\left[\mathfrak{z}\right]_{p^{n-1}}\right)\right|_{q}}\cdot\frac{\left|\nabla_{p^{n}}\left\{ \chi\right\} \left(\mathfrak{z}\right)\right|_{q}^{1/2}}{\left|\chi\left(\left[\mathfrak{z}\right]_{p^{n}}\right)\right|_{q}}=0
\]
The assumptions on $\chi$, meanwhile, guarantee both: 
\begin{equation}
\frac{\left|\nabla_{p^{n}}\left\{ \chi\right\} \left(\mathfrak{z}\right)\right|_{q}^{1/2}}{\left|\chi\left(\left[\mathfrak{z}\right]_{p^{n}}\right)\right|_{q}}\rightarrow0
\end{equation}
and: 
\begin{equation}
\sup_{n\geq1}\left|\frac{\chi\left(\left[\mathfrak{z}\right]_{p^{n}}\right)}{\chi\left(\left[\mathfrak{z}\right]_{p^{n-1}}\right)}\right|_{q}<\infty
\end{equation}
Consequently, the above limit is zero for each $\mathfrak{z}\in\mathbb{Z}_{p}$.
This proves $1/\chi$ is rising-continuous.

Q.E.D. 
\begin{thm}
Let $\tilde{C}\left(\mathbb{Z}_{p},\mathbb{T}_{q}\right)$ denote
the set of all $\chi\in\tilde{C}\left(\mathbb{Z}_{p},\mathbb{C}_{q}\right)$
for which $\left|\chi\left(\mathfrak{z}\right)\right|_{q}=1$ for
all $\mathfrak{z}\in\mathbb{Z}_{p}$, where $\mathbb{T}_{q}$ is the
set of all $q$-adic complex numbers with a $q$-adic absolute value
of $1$. Then, $\tilde{C}\left(\mathbb{Z}_{p},\mathbb{T}_{q}\right)$
is an abelian group under the operation of point-wise multiplication;
the identity element is the constant function $1$. 
\end{thm}
Proof: Every $\chi\in\tilde{C}\left(\mathbb{Z}_{p},\mathbb{T}_{q}\right)$
satisfies the conditions of \textbf{Lemma \ref{lem:3.17}}, because
$\left|\chi\left(\mathfrak{z}\right)\right|_{q}=\left|\chi\left(\left[\mathfrak{z}\right]_{p^{n}}\right)\right|_{q}=1$
for all $\mathfrak{z}\in\mathbb{Z}_{p}$ and all $n\in\mathbb{N}_{0}$.

Q.E.D.

\newpage{} 

\section{\label{sec:3.3 quasi-integrability}$\left(p,q\right)$-adic Measures
and Quasi-Integrability}

THROUGHOUT THIS SECTION $p$ AND $q$ DENOTE DISTINCT PRIME NUMBERS.

\vphantom{}

Despite the pronouncements from Schikhof and Konrad used in this chapter's
epigraphs, there is more to $\left(p,q\right)$-adic integration than
meets the eye. Consider a function $\hat{\mu}:\hat{\mathbb{Z}}_{p}\rightarrow\overline{\mathbb{Q}}$.
Viewing $\overline{\mathbb{Q}}$ as being embedded in $\mathbb{C}_{q}$,
observe that if $\sup_{t\in\hat{\mathbb{Z}}_{p}}\left|\hat{\mu}\left(t\right)\right|_{q}<\infty$,
we can then make sense of $\hat{\mu}$ as the Fourier-Stieltjes transform
of the measure $d\mu\in C\left(\mathbb{Z}_{p},\mathbb{C}_{q}\right)^{\prime}$
defined by the formula: 
\begin{equation}
\int_{\mathbb{Z}_{p}}f\left(\mathfrak{z}\right)d\mu\left(\mathfrak{z}\right)=\sum_{t\in\hat{\mathbb{Z}}_{p}}\hat{f}\left(-t\right)\hat{\mu}\left(t\right),\textrm{ }\forall f\in C\left(\mathbb{Z}_{p},\mathbb{C}_{q}\right)
\end{equation}
Now, for any $g\in C\left(\mathbb{Z}_{p},\mathbb{C}_{q}\right)$,
the $\left(p,q\right)$-adic measure $d\mu\left(\mathfrak{z}\right)=g\left(\mathfrak{z}\right)d\mathfrak{z}$
satisfies: 
\begin{equation}
g\left(\mathfrak{z}\right)\overset{\mathbb{C}_{q}}{=}\sum_{t\in\hat{\mathbb{Z}}_{p}}\hat{\mu}\left(t\right)e^{2\pi i\left\{ t\mathfrak{z}\right\} _{p}},\textrm{ }\forall\mathfrak{z}\in\mathbb{Z}_{p}\label{eq:g in terms of mu hat}
\end{equation}
because the continuity of $g$ guarantees that $\hat{\mu}\in c_{0}\left(\hat{\mathbb{Z}}_{p},\mathbb{C}_{q}\right)$.
For a measure like $d\mu\left(\mathfrak{z}\right)=g\left(\mathfrak{z}\right)d\mathfrak{z}$,
the Fourier series on the right-hand side of (\ref{eq:g in terms of mu hat})
makes perfect sense. On the other hand, for an arbitrary bounded function
$\hat{\mu}:\hat{\mathbb{Z}}_{p}\rightarrow\mathbb{C}_{q}$, there
is no guarantee that the right-hand side of (\ref{eq:g in terms of mu hat})\textemdash viewed
as the limit\footnote{I will, at times, refer to (\ref{eq:The Limit of Interest}) as the
``limit of interest''.} of the partial sums: 
\begin{equation}
\lim_{N\rightarrow\infty}\sum_{\left|t\right|_{p}\leq p^{N}}\hat{\mu}\left(t\right)e^{2\pi i\left\{ t\mathfrak{z}\right\} _{p}}\label{eq:The Limit of Interest}
\end{equation}
\textemdash will even exist for any given $\mathfrak{z}\in\mathbb{Z}_{p}$.

In this section, we will investigate the classes of $\left(p,q\right)$-adic
measures $d\mu$ so that (\ref{eq:The Limit of Interest}) makes sense
for some\textemdash or even all\textemdash values of $\mathfrak{z}\in\mathbb{Z}_{p}$\textemdash even
if $\hat{\mu}\notin c_{0}\left(\hat{\mathbb{Z}}_{p},\mathbb{C}_{q}\right)$.
Because measures are, by definition, ``integrable'', we can enlarge
the tent of $\left(p,q\right)$-adically ``integrable'' functions
by including those functions which just so happen to be given by (\ref{eq:The Limit of Interest})
for some bounded $\hat{\mu}$. I call these \textbf{quasi-integrable
functions}.\textbf{ }Given a quasi-integrable function\index{quasi-integrability}
$\chi:\mathbb{Z}_{p}\rightarrow\mathbb{C}_{q}$, upon identifying
$\chi$ with the measure $d\mu$, integration against $\chi$ can
then be defined as integration against $d\mu$ by way of the formula:
\begin{equation}
\int_{\mathbb{Z}_{p}}f\left(\mathfrak{z}\right)\chi\left(\mathfrak{z}\right)d\mathfrak{z}\overset{\mathbb{C}_{q}}{=}\sum_{t\in\hat{\mathbb{Z}}_{p}}\hat{f}\left(-t\right)\hat{\mu}\left(t\right),\textrm{ }\forall f\in C\left(\mathbb{Z}_{p},\mathbb{C}_{q}\right)
\end{equation}

I have two points to make before we begin. First, \emph{the reader
should be aware that} \emph{this chapter contains a large number of
qualitative definitions}. As will be explained in detail throughout,
this is primarily because, at present, there does not appear to be
a way to get desirable conclusions properties out of generic $\left(p,q\right)$-adic
functions without relying on explicit formulae and computations thereof,
such as (\ref{eq:Definition of right-ended p-adic structural equations}),
(\ref{eq:Radial-Magnitude Fourier Resummation Lemma - p-adically distributed case}),
and (\ref{eq:reciprocal of the p-adic absolute value of z is quasi-integrable}),
and for these explicit computations to be feasible, we need our functions
to come pre-equipped with at least \emph{some }concrete algebraic
structure. By and by, I will indicate where I have been forced to
complicate definitions for the sake of subsequent proof. It is obviously
of interest to see how much these definitions can be streamlined and
simplified.

Secondly and finally, I have to say that there are still many basic
technical considerations regarding quasi-integrability and quasi-integrable
functions which I have not been able to develop to my satisfaction
due to a lack of time and/or cleverness. The relevant questions are
stated and discussed on pages \pageref{que:3.3}, \pageref{que:3.2},
pages \pageref{que:3.4} through \pageref{que:3.6}, and page \pageref{que:3.7}.

\subsection{\label{subsec:3.3.1 Heuristics-and-Motivations}Heuristics and Motivations}

To the extent that I have discussed my ideas with other mathematicians
who did not already know me personally, the reactions I have gotten
from bringing up frames seems to indicate that my ideas are considered
``difficult'', maybe even radical. But they were difficult for me
long before they were difficult for anyone else. It took over a year
before I finally understood what was going on. My discoveries of quasi-integrability
and the necessity of frames occurred completely by chance. In performing
a certain lengthy computation, I seemed to get a different answer
every time. It was only after breaking up the computations into multiple
steps\textemdash performing each separately and in greater generality\textemdash that
I realized what was going on.

To that end, rather than pull a Bourbaki, I think it will be best
if we first consider the main examples that I struggled with. It's
really the natural thing to do. Pure mathematics \emph{is} a natural
science, after all. Its empirical data are the patterns we notice
and the observations we happen to make. The hypotheses are our latest
ideas for a problem, and the experiments are our efforts to see if
those ideas happen to bear fruit. In this respect, frames and quasi-integrability
should not be seen as platonic ideals, but as my effort to describe
and speak clearly about an unusual phenomenon. The examples below
are the context in which the reader should first try to understand
my ideas, because this was the context in which all of my work was
originally done.

Instead of saving the best for last, let's go right to heart of the
matter and consider the example that caused me the most trouble, by
far. First, however, I must introduce the main antagonist to feature
in our studies of the shortened $qx+1$ maps.
\begin{defn}[$\hat{A}_{q}\left(t\right)$]
Let \index{hat{A}{q}@$\hat{A}_{q}$}$q$ be an odd prime. \nomenclature{$\hat{A}_{q}\left(t\right)$}{ }I
define $\hat{A}_{q}:\hat{\mathbb{Z}}_{2}\rightarrow\overline{\mathbb{Q}}$
by: 
\begin{equation}
\hat{A}_{q}\left(t\right)\overset{\textrm{def}}{=}\begin{cases}
1 & \textrm{if }t=0\\
\prod_{n=0}^{-v_{2}\left(t\right)-1}\frac{1+qe^{-2\pi i\left(2^{n}t\right)}}{4} & \textrm{else}
\end{cases},\textrm{ }\forall t\in\hat{\mathbb{Z}}_{2}\label{eq:Definition of A_q hat}
\end{equation}
Noting that:
\begin{equation}
\sup_{t\in\hat{\mathbb{Z}}_{2}}\left|\hat{A}_{q}\left(t\right)\right|_{q}=1
\end{equation}
it follows that we can define a $\left(2,q\right)$-adic measure \nomenclature{$dA_{q}$}{ }$dA_{q}\in C\left(\mathbb{Z}_{2},\mathbb{C}_{q}\right)^{\prime}$
by way of the formula: 
\begin{equation}
\int_{\mathbb{Z}_{2}}f\left(\mathfrak{z}\right)dA_{q}\left(\mathfrak{z}\right)\overset{\mathbb{C}_{q}}{=}\sum_{t\in\hat{\mathbb{Z}}_{q}}\hat{f}\left(-t\right)\hat{A}_{q}\left(t\right),\textrm{ }\forall f\in C\left(\mathbb{Z}_{2},\mathbb{C}_{q}\right)
\end{equation}
As defined, $dA_{q}$ has $\hat{A}_{q}$ as its Fourier-Stieltjes
transform\index{$dA_{q}$!Fourier-Stieltjes transf.}.
\end{defn}
\vphantom{}

The $\hat{A}_{q}$s will be of \emph{immen}se importance in Chapter
4, where we will see how the behavior of $dA_{q}$ is for $q=3$ is
drastically different from all other odd primes. For now, though,
let's focus on the $q=3$ case\textemdash the Collatz Case\textemdash seeing
as it provides us an example of a particularly pathological $\left(p,q\right)$-adic
measure\textemdash what I call a \textbf{degenerate measure}\textemdash as
well as a motivation for the topology-mixing inspiration that birthed
my concept of a frame.
\begin{prop}
\label{prop:Nth partial sum of Fourier series of A_3 hat}\index{hat{A}{3}@$\hat{A}_{3}$}
\begin{equation}
\sum_{\left|t\right|_{2}\leq2^{N}}\hat{A}_{3}\left(t\right)e^{2\pi i\left\{ t\mathfrak{z}\right\} _{2}}=\frac{3^{\#_{1}\left(\left[\mathfrak{z}\right]_{2^{N}}\right)}}{2^{N}},\textrm{ }\forall\mathfrak{z}\in\mathbb{Z}_{2},\textrm{ }\forall N\geq0\label{eq:Convolution of dA_3 with D_N}
\end{equation}
where $\#_{1}\left(m\right)\overset{\textrm{def}}{=}\#_{2:1}\left(m\right)$
is the number of $1$s in the binary expansion of $m$. 
\end{prop}
Proof: As will be proven in \textbf{Proposition \ref{prop:Generating function identities}}
(page \pageref{prop:Generating function identities}), we have a generating
function identity: 
\begin{equation}
\prod_{m=0}^{n-1}\left(1+az^{2^{m}}\right)=\sum_{m=0}^{2^{n}-1}a^{\#_{1}\left(m\right)}z^{m}\label{eq:Product-to-sum identity for number of 1s digits}
\end{equation}
which holds for any field $\mathbb{F}$ of characteristic $0$ and
any $a,z\in\mathbb{F}$, and any $n\geq1$.

Consequently: 
\begin{equation}
\prod_{n=0}^{-v_{2}\left(t\right)-1}\frac{1+3e^{-2\pi i\left(2^{n}t\right)}}{4}=\frac{1}{\left|t\right|_{2}^{2}}\sum_{m=0}^{\left|t\right|_{2}-1}3^{\#_{1}\left(m\right)}e^{-2\pi imt}
\end{equation}
As such, for any $n\geq1$: 
\begin{align*}
\sum_{\left|t\right|_{2}=2^{n}}\hat{A}_{3}\left(t\right)e^{2\pi i\left\{ t\mathfrak{z}\right\} _{2}} & =\sum_{\left|t\right|_{2}=2^{n}}\left(\frac{1}{\left|t\right|_{2}^{2}}\sum_{m=0}^{\left|t\right|_{2}-1}3^{\#_{1}\left(m\right)}e^{-2\pi imt}\right)e^{2\pi i\left\{ t\mathfrak{z}\right\} _{2}}\\
 & =\frac{1}{4^{n}}\sum_{m=0}^{2^{n}-1}3^{\#_{1}\left(m\right)}\sum_{\left|t\right|_{2}=2^{n}}e^{2\pi i\left\{ t\left(\mathfrak{z}-m\right)\right\} _{2}}\\
 & =\frac{1}{4^{n}}\sum_{m=0}^{2^{n}-1}3^{\#_{1}\left(m\right)}\left(2^{n}\left[\mathfrak{z}\overset{2^{n}}{\equiv}m\right]-2^{n-1}\left[\mathfrak{z}\overset{2^{n-1}}{\equiv}m\right]\right)\\
 & =\frac{1}{2^{n}}\sum_{m=0}^{2^{n}-1}3^{\#_{1}\left(m\right)}\left[\mathfrak{z}\overset{2^{n}}{\equiv}m\right]-\frac{1}{2^{n+1}}\sum_{m=0}^{2^{n}-1}3^{\#_{1}\left(m\right)}\left[\mathfrak{z}\overset{2^{n-1}}{\equiv}m\right]
\end{align*}
Here, note that $m=\left[\mathfrak{z}\right]_{2^{n}}$ is the unique
integer $m\in\left\{ 0,\ldots,2^{n}-1\right\} $ for which $\mathfrak{z}\overset{2^{n}}{\equiv}m$
holds true. Thus: 
\begin{equation}
\sum_{m=0}^{2^{n}-1}3^{\#_{1}\left(m\right)}\left[\mathfrak{z}\overset{2^{n}}{\equiv}m\right]=3^{\#_{1}\left(\left[\mathfrak{z}\right]_{2^{n}}\right)}
\end{equation}
Likewise: 
\begin{align*}
\sum_{m=0}^{2^{n}-1}3^{\#_{1}\left(m\right)}\left[\mathfrak{z}\overset{2^{n-1}}{\equiv}m\right] & =\sum_{m=0}^{2^{n-1}-1}3^{\#_{1}\left(m\right)}\left[\mathfrak{z}\overset{2^{n-1}}{\equiv}m\right]+\sum_{m=2^{n-1}}^{2^{n}-1}3^{\#_{1}\left(m\right)}\left[\mathfrak{z}\overset{2^{n-1}}{\equiv}m\right]\\
 & =3^{\#_{1}\left(\left[\mathfrak{z}\right]_{2^{n-1}}\right)}+\sum_{m=0}^{2^{n-1}-1}3^{\#_{1}\left(m+2^{n-1}\right)}\left[\mathfrak{z}\overset{2^{n-1}}{\equiv}m+2^{n-1}\right]\\
 & =3^{\#_{1}\left(\left[\mathfrak{z}\right]_{2^{n-1}}\right)}+\sum_{m=0}^{2^{n-1}-1}3^{1+\#_{1}\left(m\right)}\left[\mathfrak{z}\overset{2^{n-1}}{\equiv}m\right]\\
 & =3^{\#_{1}\left(\left[\mathfrak{z}\right]_{2^{n-1}}\right)}+3\cdot3^{\#_{1}\left(\left[\mathfrak{z}\right]_{2^{n-1}}\right)}\\
 & =4\cdot3^{\#_{1}\left(\left[\mathfrak{z}\right]_{2^{n-1}}\right)}
\end{align*}
and so: 
\[
\sum_{\left|t\right|_{2}=2^{n}}\hat{A}_{3}\left(t\right)e^{2\pi i\left\{ t\mathfrak{z}\right\} _{2}}=\frac{3^{\#_{1}\left(\left[\mathfrak{z}\right]_{2^{n}}\right)}}{2^{n}}-\frac{4\cdot3^{\#_{1}\left(\left[\mathfrak{z}\right]_{2^{n-1}}\right)}}{2^{n+1}}=\frac{3^{\#_{1}\left(\left[\mathfrak{z}\right]_{2^{n}}\right)}}{2^{n}}-\frac{3^{\#_{1}\left(\left[\mathfrak{z}\right]_{2^{n-1}}\right)}}{2^{n-1}}
\]
Consequently: 
\begin{align*}
\sum_{\left|t\right|_{2}\leq2^{N}}\hat{A}_{3}\left(t\right)e^{2\pi i\left\{ t\mathfrak{z}\right\} _{2}} & =\hat{A}_{3}\left(0\right)+\sum_{n=1}^{N}\sum_{\left|t\right|_{2}=2^{n}}\hat{A}_{3}\left(t\right)e^{2\pi i\left\{ t\mathfrak{z}\right\} _{2}}\\
 & =1+\sum_{n=1}^{N}\left(\frac{3^{\#_{1}\left(\left[\mathfrak{z}\right]_{2^{n}}\right)}}{2^{n}}-\frac{3^{\#_{1}\left(\left[\mathfrak{z}\right]_{2^{n-1}}\right)}}{2^{n-1}}\right)\\
\left(\textrm{telescoping series}\right); & =1+\frac{3^{\#_{1}\left(\left[\mathfrak{z}\right]_{2^{N}}\right)}}{2^{N}}-\underbrace{3^{\#_{1}\left(\left[\mathfrak{z}\right]_{2^{0}}\right)}}_{=3^{0}=1}\\
\left(\left[\mathfrak{z}\right]_{2^{0}}=0\right); & =\frac{3^{\#_{1}\left(\left[\mathfrak{z}\right]_{2^{N}}\right)}}{2^{N}}
\end{align*}

Q.E.D.

\vphantom{}

With this formula, the reader will see exactly what I when I say the
measure $dA_{3}$ is ``\index{$dA_{3}$!degeneracy}degenerate''.
\begin{prop}
\label{prop:dA_3 is a degenerate measure}We have: 
\begin{equation}
\lim_{N\rightarrow\infty}\sum_{\left|t\right|_{2}\leq2^{N}}\hat{A}_{3}\left(t\right)e^{2\pi i\left\{ t\mathfrak{z}\right\} _{2}}\overset{\mathbb{C}_{3}}{=}0,\textrm{ }\forall\mathfrak{z}\in\mathbb{Z}_{2}^{\prime}\label{eq:Degeneracy of dA_3 on Z_2 prime}
\end{equation}
and: 
\begin{equation}
\lim_{N\rightarrow\infty}\sum_{\left|t\right|_{2}\leq2^{N}}\hat{A}_{3}\left(t\right)e^{2\pi i\left\{ t\mathfrak{z}\right\} _{2}}\overset{\mathbb{C}}{=}0,\textrm{ }\forall\mathfrak{z}\in\mathbb{N}_{0}\label{eq:Degeneracy of dA_3 on N_0}
\end{equation}
\end{prop}
\begin{rem}
That is to say, the upper limit occurs in the topology of $\mathbb{C}_{3}$
while the lower limit occurs in the topology of $\mathbb{C}$. The
reason this works is because, for each $N$, the Fourier series is
a sum of finitely many elements of $\overline{\mathbb{Q}}$, which
we have embedded in both $\mathbb{C}_{3}$ and $\mathbb{C}$, and
that this sum is invariant under the action of $\textrm{Gal}\left(\overline{\mathbb{Q}}/\mathbb{Q}\right)$.
\end{rem}
Proof: By (\ref{eq:Convolution of dA_3 with D_N}), we have: 
\begin{equation}
\lim_{N\rightarrow\infty}\sum_{\left|t\right|_{2}\leq2^{N}}\hat{A}_{3}\left(t\right)e^{2\pi i\left\{ t\mathfrak{z}\right\} _{2}}=\lim_{N\rightarrow\infty}\frac{3^{\#_{1}\left(\left[\mathfrak{z}\right]_{2^{N}}\right)}}{2^{N}}
\end{equation}
When $\mathfrak{z}\in\mathbb{Z}_{2}^{\prime}$, the number $\#_{1}\left(\left[\mathfrak{z}\right]_{2^{N}}\right)$
tends to $\infty$ as $N\rightarrow\infty$, thereby guaranteeing
that the above limit converges to $0$ in $\mathbb{C}_{3}$. On the
other hand, when $\mathfrak{z}\in\mathbb{N}_{0}$, $\#_{1}\left(\left[\mathfrak{z}\right]_{2^{N}}\right)=\#_{1}\left(\mathfrak{z}\right)$
for all $N\geq\lambda_{2}\left(\mathfrak{z}\right)$, and so, the
above limit then tends to $0$ in $\mathbb{C}$.

Q.E.D. 
\begin{prop}
We have the integral formula: 
\begin{equation}
\int_{\mathbb{Z}_{2}}f\left(\mathfrak{z}\right)dA_{3}\left(\mathfrak{z}\right)\overset{\mathbb{C}_{3}}{=}\lim_{N\rightarrow\infty}\frac{1}{4^{N}}\sum_{n=0}^{2^{N}-1}f\left(n\right)3^{\#_{1}\left(n\right)},\textrm{ }\forall f\in C\left(\mathbb{Z}_{2},\mathbb{C}_{3}\right)\label{eq:Riemann sum formula for the dA_3 integral of a continuous (2,3)-adic function}
\end{equation}
\end{prop}
Proof: Since: 
\begin{equation}
\left[\mathfrak{z}\overset{2^{N}}{\equiv}n\right]=\frac{1}{2^{N}}\sum_{\left|t\right|_{2}\leq2^{N}}e^{2\pi i\left\{ t\left(\mathfrak{z}-n\right)\right\} _{2}}=\frac{1}{2^{N}}\sum_{\left|t\right|_{2}\leq2^{N}}e^{-2\pi i\left\{ tn\right\} _{2}}e^{2\pi i\left\{ t\mathfrak{z}\right\} _{2}}
\end{equation}
we have: 
\begin{equation}
\int_{\mathbb{Z}_{2}}\left[\mathfrak{z}\overset{2^{N}}{\equiv}n\right]dA_{3}\left(\mathfrak{z}\right)=\frac{1}{2^{N}}\sum_{\left|t\right|_{2}\leq2^{N}}\hat{A}_{3}\left(t\right)e^{2\pi i\left\{ tn\right\} _{2}}=\frac{3^{\#_{1}\left(n\right)}}{4^{N}}
\end{equation}
for all $N\geq0$ and $n\in\left\{ 0,\ldots,2^{N}-1\right\} $. So,
letting $f\in C\left(\mathbb{Z}_{2},\mathbb{C}_{3}\right)$ be arbitrary,
taking $N$th truncations yields: 
\begin{equation}
\int_{\mathbb{Z}_{2}}f_{N}\left(\mathfrak{z}\right)dA_{3}\left(\mathfrak{z}\right)=\sum_{n=0}^{2^{N}-1}f\left(n\right)\int_{\mathbb{Z}_{2}}\left[\mathfrak{z}\overset{2^{N}}{\equiv}n\right]dA_{3}\left(\mathfrak{z}\right)=\frac{1}{4^{N}}\sum_{n=0}^{2^{N}-1}f\left(n\right)3^{\#_{1}\left(n\right)}
\end{equation}
Since $f$ is continuous, the $f_{N}$ converge uniformly to $f$
(\textbf{Proposition \ref{prop:Unif. convergence of truncation equals continuity}}).
Since $dA_{3}\in C\left(\mathbb{Z}_{2},\mathbb{C}_{3}\right)^{\prime}$,
this guarantees: 
\begin{equation}
\int_{\mathbb{Z}_{2}}f\left(\mathfrak{z}\right)dA_{3}\left(\mathfrak{z}\right)\overset{\mathbb{C}_{3}}{=}\lim_{N\rightarrow\infty}\int_{\mathbb{Z}_{2}}f_{N}\left(\mathfrak{z}\right)dA_{3}\left(\mathfrak{z}\right)=\lim_{N\rightarrow\infty}\frac{1}{4^{N}}\sum_{n=0}^{2^{N}-1}f\left(n\right)3^{\#_{1}\left(n\right)}
\end{equation}

Q.E.D.

\vphantom{}

Thus, for $d\mu=dA_{3}$, our limit of interest (\ref{eq:The Limit of Interest})
converges to $0$, yet $dA_{3}$ is not the zero measure: this is
$dA_{3}$'s ``degeneracy''. $dA_{3}$ also places us in the extremely
unusual position of resorting to \emph{different topologies} for different
inputs in order to guarantee the point-wise existence of our limit
for every $\mathfrak{z}\in\mathbb{Z}_{2}$. The purpose of Section
\ref{sec:3.3 quasi-integrability} is to demonstrate that this procedure
can be done consistently, and, moreover, in a way that leads to useful
results.

Contrary to what Monna-Springer theory would have us believe, our
next examples show that there are discontinuous\textemdash even \emph{singular!}\textemdash $\left(p,q\right)$-adic
functions which can be meaningfully integrated. 
\begin{prop}
\label{prop:sum of v_p}Let $p$ be an integer $\geq2$. Then, for
each $\mathfrak{z}\in\mathbb{Z}_{p}\backslash\left\{ 0\right\} $:
\begin{equation}
\sum_{0<\left|t\right|_{p}\leq p^{N}}v_{p}\left(t\right)e^{2\pi i\left\{ t\mathfrak{z}\right\} _{p}}\overset{\overline{\mathbb{Q}}}{=}\frac{p\left|\mathfrak{z}\right|_{p}^{-1}-1}{p-1},\textrm{ }\forall N>v_{p}\left(\mathfrak{z}\right)\label{eq:Fourier sum of v_p of t}
\end{equation}
Here, the use of $\overset{\overline{\mathbb{Q}}}{=}$ is to indicate
that the equality is one of elements of $\overline{\mathbb{Q}}$. 
\end{prop}
Proof: Fixing $\mathfrak{z}\in\mathbb{Z}_{p}\backslash\left\{ 0\right\} $,
we have: 
\begin{align*}
\sum_{0<\left|t\right|_{p}\leq p^{N}}v_{p}\left(t\right)e^{2\pi i\left\{ t\mathfrak{z}\right\} _{p}} & =\sum_{n=1}^{N}\sum_{\left|t\right|_{p}=p^{n}}\left(-n\right)e^{2\pi i\left\{ t\mathfrak{z}\right\} _{p}}\\
 & =\sum_{n=1}^{N}\left(-n\right)\left(p^{n}\left[\mathfrak{z}\overset{p^{n}}{\equiv}0\right]-p^{n-1}\left[\mathfrak{z}\overset{p^{n-1}}{\equiv}0\right]\right)
\end{align*}
Simplifying gives: 
\begin{equation}
\sum_{0<\left|t\right|_{p}\leq p^{N}}v_{p}\left(t\right)e^{2\pi i\left\{ t\mathfrak{z}\right\} _{p}}=-Np^{N}\left[\mathfrak{z}\overset{p^{N}}{\equiv}0\right]+\sum_{n=0}^{N-1}p^{n}\left[\mathfrak{z}\overset{p^{n}}{\equiv}0\right]
\end{equation}
Because $\mathfrak{z}$ is non-zero, $\mathfrak{z}=p^{v_{p}\left(\mathfrak{z}\right)}\mathfrak{u}$
for some $\mathfrak{u}\in\mathbb{Z}_{p}^{\times}$, where $v_{p}\left(\mathfrak{z}\right)$
is a non-negative rational integer. As such, the congruence $\mathfrak{z}\overset{p^{n}}{\equiv}0$
will fail to hold for all $n\geq v_{p}\left(\mathfrak{z}\right)+1$.
As such, we can write:: 
\begin{equation}
\sum_{0<\left|t\right|_{p}\leq p^{N}}v_{p}\left(t\right)e^{2\pi i\left\{ t\mathfrak{z}\right\} _{p}}=\sum_{n=0}^{v_{p}\left(\mathfrak{z}\right)}p^{n}\left[\mathfrak{z}\overset{p^{n}}{\equiv}0\right],\textrm{ }\forall N\geq v_{p}\left(\mathfrak{z}\right)+1,\textrm{ }\forall\mathfrak{z}\in\mathbb{Z}_{p}\backslash\left\{ 0\right\} 
\end{equation}
Since $\mathfrak{z}\overset{p^{n}}{\equiv}0$ holds true if and only
if $n\in\left\{ 0,\ldots,v_{p}\left(\mathfrak{z}\right)\right\} $,
the Iverson brackets in the above formula all evaluate to $1$. This
leaves us with a finite geometric series: 
\begin{equation}
\sum_{0<\left|t\right|_{p}\leq p^{N}}v_{p}\left(t\right)e^{2\pi i\left\{ t\mathfrak{z}\right\} _{p}}=\sum_{n=0}^{v_{p}\left(\mathfrak{z}\right)}p^{n}\left[\mathfrak{z}\overset{p^{n}}{\equiv}0\right]=\sum_{n=0}^{v_{p}\left(\mathfrak{z}\right)}p^{n}=\frac{p^{v_{p}\left(\mathfrak{z}\right)+1}-1}{p-1}
\end{equation}
for all $N\geq v_{p}\left(\mathfrak{z}\right)+1$. Noting that $p\left|\mathfrak{z}\right|_{p}^{-1}=p^{v_{p}\left(\mathfrak{z}\right)+1}$
then completes the computation.

Q.E.D.

\vphantom{}

In general, there is no way to define (\ref{eq:Fourier sum of v_p of t})
when $\mathfrak{z}=0$. Indeed: 
\begin{equation}
\sum_{0<\left|t\right|_{p}\leq p^{N}}v_{p}\left(t\right)=\sum_{n=1}^{N}\sum_{\left|t\right|_{p}=p^{n}}\left(-n\right)=-\sum_{n=1}^{N}n\varphi\left(p^{n}\right)=-\left(p-1\right)\sum_{n=1}^{N}np^{n-1}
\end{equation}
where $\varphi$ is the \textbf{Euler Totient Function}. Since $p$
and $q$ are distinct primes, the $q$-adic absolute value of the
$n$th term of the partial sum of the series is: 
\begin{equation}
\left|np^{n-1}\right|_{q}=\left|n\right|_{q}
\end{equation}
which does not tend to $0$ as $n\rightarrow\infty$. Thus, (\ref{eq:Fourier sum of v_p of t})
does not converge $q$-adically when $\mathfrak{z}=0$. On the other
hand, in the topology of $\mathbb{C}$, the series diverges to $\infty$.

Nevertheless, because of the $\hat{\mu}$ defined above, we can write:
\begin{equation}
\frac{p\left|\mathfrak{z}\right|_{p}^{-1}-1}{p-1}d\mathfrak{z}
\end{equation}
to denote the measure defined by: 
\begin{equation}
\int_{\mathbb{Z}_{p}}\frac{p\left|\mathfrak{z}\right|_{p}^{-1}-1}{p-1}f\left(\mathfrak{z}\right)d\mathfrak{z}\overset{\textrm{def}}{=}\sum_{t\in\hat{\mathbb{Z}}_{p}\backslash\left\{ 0\right\} }\hat{f}\left(-t\right)v_{p}\left(t\right),\textrm{ }\forall f\in C\left(\mathbb{Z}_{p},\mathbb{C}_{q}\right)
\end{equation}
In this sense, the $\left(p,q\right)$-adic function: 
\begin{equation}
\frac{p\left|\mathfrak{z}\right|_{p}^{-1}-1}{p-1}
\end{equation}
is ``quasi-integrable'' in that despite having a singularity at
$\mathfrak{z}=0$, we can still use the above integral formula to
define its integral: 
\begin{equation}
\int_{\mathbb{Z}_{p}}\frac{p\left|\mathfrak{z}\right|_{p}^{-1}-1}{p-1}d\mathfrak{z}\overset{\textrm{def}}{=}0
\end{equation}
as the image of the constant function $1$ under the associated measure.

With regard to our notational formalism, it is then natural to ask:
\begin{equation}
\int_{\mathbb{Z}_{p}}\frac{p\left|\mathfrak{z}\right|_{p}^{-1}-1}{p-1}d\mathfrak{z}\overset{?}{=}\frac{p}{p-1}\int_{\mathbb{Z}_{p}}\left|\mathfrak{z}\right|_{p}^{-1}d\mathfrak{z}-\frac{1}{p-1}\underbrace{\int_{\mathbb{Z}_{p}}d\mathfrak{z}}_{1}
\end{equation}
The answer is \emph{yes}. 
\begin{prop}
Let $p$ be an integer $\geq2$. Then, for all $\mathfrak{z}\in\mathbb{Z}_{p}\backslash\left\{ 0\right\} $,
we have: 
\begin{equation}
\left|\mathfrak{z}\right|_{p}^{-1}\overset{\overline{\mathbb{Q}}}{=}\frac{1}{p}+\frac{p-1}{p}\sum_{0<\left|t\right|_{p}\leq p^{N}}v_{p}\left(t\right)e^{2\pi i\left\{ t\mathfrak{z}\right\} _{p}},\textrm{ }\forall N>v_{p}\left(\mathfrak{z}\right)\label{eq:reciprocal of the p-adic absolute value of z is quasi-integrable}
\end{equation}
\end{prop}
Proof: Let: 
\begin{equation}
\frac{1}{p-1}d\mathfrak{z}
\end{equation}
denote the $\left(p,q\right)$-adic measure which is the scalar multiple
of the $\left(p,q\right)$-adic Haar probability measure by $1/\left(p-1\right)$.
Then, by definition, for any $f\in C\left(\mathbb{Z}_{p},\mathbb{C}_{q}\right)$:
\begin{align*}
\int_{\mathbb{Z}_{p}}\frac{p\left|\mathfrak{z}\right|_{p}^{-1}-1}{p-1}f\left(\mathfrak{z}\right)d\mathfrak{z}+\int_{\mathbb{Z}_{p}}\frac{f\left(\mathfrak{z}\right)}{p-1}d\mathfrak{z} & =\sum_{t\in\hat{\mathbb{Z}}_{p}\backslash\left\{ 0\right\} }\hat{f}\left(-t\right)v_{p}\left(t\right)+\frac{\hat{f}\left(0\right)}{p-1}\\
 & =\sum_{t\in\hat{\mathbb{Z}}_{p}}\hat{f}\left(-t\right)\hat{\nu}\left(t\right)
\end{align*}
where $\hat{\nu}:\hat{\mathbb{Z}}_{p}\rightarrow\mathbb{C}_{q}$ is
defined by: 
\begin{equation}
\hat{\nu}\left(t\right)\overset{\textrm{def}}{=}\begin{cases}
\frac{1}{p-1} & \textrm{if }t=0\\
v_{p}\left(t\right) & \textrm{else}
\end{cases},\textrm{ }\forall t\in\hat{\mathbb{Z}}_{p}
\end{equation}
Running through the computation of $\tilde{\nu}_{N}\left(\mathfrak{z}\right)$
yields: 
\begin{equation}
\sum_{\left|t\right|_{p}\leq p^{N}}\nu\left(t\right)e^{2\pi i\left\{ t\mathfrak{z}\right\} _{p}}\overset{\mathbb{Q}}{=}\frac{p\left|\mathfrak{z}\right|_{p}^{-1}}{p-1},\textrm{ }\forall\mathfrak{z}\in\mathbb{Z}_{p}\backslash\left\{ 0\right\} ,\textrm{ }\forall N>v_{p}\left(\mathfrak{z}\right)
\end{equation}
Hence, $\mathfrak{z}\in\mathbb{Z}_{p}\backslash\left\{ 0\right\} $
and $N>v_{p}\left(\mathfrak{z}\right)$ imply: 
\begin{align*}
\left|\mathfrak{z}\right|_{p}^{-1} & \overset{\overline{\mathbb{Q}}}{=}\sum_{\left|t\right|_{p}\leq p^{N}}\frac{p-1}{p}\nu\left(t\right)e^{2\pi i\left\{ t\mathfrak{z}\right\} _{p}}\\
 & \overset{\overline{\mathbb{Q}}}{=}\frac{1}{p}+\frac{p-1}{p}\sum_{0<\left|t\right|_{p}\leq p^{N}}v_{p}\left(t\right)e^{2\pi i\left\{ t\mathfrak{z}\right\} _{p}}
\end{align*}
as desired.

Q.E.D. 
\begin{example}
Letting $\mathfrak{a}\in\mathbb{Z}_{p}$ be arbitrary, for all $\mathfrak{z}\in\mathbb{Z}_{p}\backslash\left\{ 0\right\} $
and all $N>v_{p}\left(\mathfrak{z}\right)$, we have: 
\[
\left|\mathfrak{z}-\mathfrak{a}\right|_{p}^{-1}\overset{\mathbb{Q}}{=}\frac{1}{p}+\frac{p-1}{p}\sum_{0<\left|t\right|_{p}\leq p^{N}}v_{p}\left(t\right)e^{-2\pi i\left\{ t\mathfrak{a}\right\} _{p}}e^{2\pi i\left\{ t\mathfrak{z}\right\} _{p}}
\]
Consequently, letting $\left\{ \mathfrak{a}_{n}\right\} _{n\geq0}$
be any sequence of distinct elements of $\mathbb{Z}_{p}$ and letting
$\left\{ \mathfrak{b}_{n}\right\} _{n\geq0}$ be any sequence in $\mathbb{C}_{q}$
tending to $0$ in $q$-adic magnitude, the function: 
\begin{equation}
\sum_{n=0}^{\infty}\frac{\mathfrak{b}_{n}}{\left|\mathfrak{z}-\mathfrak{a}\right|_{p}}\label{eq:"quasi-integrable" function with infinitely many singularities}
\end{equation}
is ``quasi-integrable'', because we can get a measure out of it
by writing: 
\begin{equation}
\int_{\mathbb{Z}_{p}}f\left(\mathfrak{z}\right)\sum_{n=0}^{\infty}\frac{\mathfrak{b}_{n}}{\left|\mathfrak{z}-\mathfrak{a}\right|_{p}}d\mathfrak{z}
\end{equation}
as:
\begin{equation}
\frac{\hat{f}\left(0\right)}{p}\sum_{n=0}^{\infty}\mathfrak{b}_{n}+\frac{p-1}{p}\sum_{t\neq0}v_{p}\left(t\right)\left(\sum_{n=0}^{\infty}\mathfrak{b}_{n}e^{-2\pi i\left\{ t\mathfrak{a}_{n}\right\} _{p}}\right)\hat{f}\left(-t\right)
\end{equation}
for any $f\in C\left(\mathbb{Z}_{p},\mathbb{C}_{q}\right)$. Moreover,
because of the failure of (\ref{eq:reciprocal of the p-adic absolute value of z is quasi-integrable})
to converge in either $\mathbb{C}_{q}$ or $\mathbb{C}$ for $\mathfrak{z}=0$,
the Fourier series for (\ref{eq:"quasi-integrable" function with infinitely many singularities})
then fails to converge in either $\mathbb{C}_{q}$ or $\mathbb{C}$
for $\mathfrak{z}=\mathfrak{a}_{n}$ for any $n$.

Even then, this is not the most pathological case we might have to
deal with. Once again, thanks to the invariance of (\ref{eq:reciprocal of the p-adic absolute value of z is quasi-integrable})
under the action of $\textrm{Gal}\left(\overline{\mathbb{Q}}/\mathbb{Q}\right)$,
the point-wise convergence of (\ref{eq:reciprocal of the p-adic absolute value of z is quasi-integrable})
on $\mathbb{Z}_{p}\backslash\left\{ 0\right\} $ has the useful property
of occurring in \emph{every} field extension of $\overline{\mathbb{Q}}$:
for each $\mathfrak{z}$, (\ref{eq:reciprocal of the p-adic absolute value of z is quasi-integrable})
is constant for all $N>v_{p}\left(\mathfrak{z}\right)$. On the other
hand, suppose the limit of interest (\ref{eq:The Limit of Interest})
yields a rising-continuous function $\chi\in\tilde{C}\left(\mathbb{Z}_{p},\mathbb{C}_{q}\right)$
which converges in $\mathbb{C}_{q}$ for all $\mathfrak{z}\in\mathbb{Z}_{p}^{\prime}$
and converges in $\mathbb{C}$ for all $\mathfrak{z}\in\mathbb{N}_{0}$.
Then, translating by some $\mathfrak{a}\in\mathbb{Z}_{p}$, we can
shift the sets on which $\chi$ converges $q$-adically by considering
the Fourier series for $\chi\left(\mathfrak{z}+\mathfrak{a}\right)$.
But then, letting $\mathfrak{a}$ and $\mathfrak{b}$ be distinct
$p$-adic integers, how might we make sense of the convergence of
the Fourier series of: 
\begin{equation}
\chi\left(\mathfrak{z}+\mathfrak{a}\right)+\chi\left(\mathfrak{z}+\mathfrak{b}\right)
\end{equation}
In this set-up, there may be a $\mathfrak{z}_{0}\in\mathbb{Z}_{p}$
for which $\chi\left(\mathfrak{z}+\mathfrak{a}\right)$'s Fourier
series converges $q$-adically and $\chi\left(\mathfrak{z}+\mathfrak{b}\right)$'s
series converges in $\mathbb{C}$ and \emph{not }in $\mathbb{C}_{q}$!
Worse yet, there is a very natural reason to \emph{want} to work with
linear combinations of the form: 
\begin{equation}
\sum_{n}\mathfrak{b}_{n}\chi\left(\mathfrak{z}+\mathfrak{a}_{n}\right)
\end{equation}
for $\mathfrak{b}_{n}$s in $\mathbb{C}_{q}$ and $\mathfrak{a}_{n}$s
in $\mathbb{Z}_{p}$; these arise naturally when considering $\left(p,q\right)$-adic
analogues of the Wiener Tauberian Theorem. Results of this type describe
conditions in which (and the extent to which) a quasi-integrable function
$\chi\left(\mathfrak{z}\right)$ will have well-defined reciprocal.
Since we want to determine the values $x\in\mathbb{Z}$ for which
$\chi_{H}\left(\mathfrak{z}\right)-x$ vanishes for some $\mathfrak{z}\in\mathbb{Z}_{p}$,
these issues are, obviously, of the utmost import to us. This particular
problem will be dealt with by defining vector spaces of $\left(p,q\right)$-adically
bounded functions $\hat{\mathbb{Z}}_{p}\rightarrow\overline{\mathbb{Q}}$,
so as to guarantee that we can take linear combinations of translates
of functions on $\mathbb{Z}_{p}$ without having to worry about getting
our topologies in a twist.
\end{example}
\vphantom{}

Whereas the previous two examples were liberating, the next example
will be dour and sobering. It tells us that our newly acquired freedoms
\emph{do} come at a price. 
\begin{example}
Having shown that: 
\begin{equation}
\sum_{0<\left|t\right|_{p}\leq p^{N}}v_{p}\left(t\right)e^{2\pi i\left\{ t\mathfrak{z}\right\} _{p}}\overset{\overline{\mathbb{Q}}}{=}\frac{p\left|\mathfrak{z}\right|_{p}^{-1}-1}{p-1},\textrm{ }\forall\mathfrak{z}\in\mathbb{Z}_{p}\backslash\left\{ 0\right\} ,\textrm{ }\forall N>v_{p}\left(\mathfrak{z}\right)
\end{equation}
upon letting $N\rightarrow\infty$, we can write: 
\begin{equation}
\sum_{t\in\hat{\mathbb{Z}}_{p}\backslash\left\{ 0\right\} }v_{p}\left(t\right)e^{2\pi i\left\{ t\mathfrak{z}\right\} _{p}}\overset{\mathbb{F}}{=}\frac{p\left|\mathfrak{z}\right|_{p}^{-1}-1}{p-1}
\end{equation}
where $\mathbb{F}$ is any valued field extension of $\overline{\mathbb{Q}}$;
note, completeness of $\mathbb{F}$ is \emph{not }required! Now, let
$p=2$. Then, we can write: 
\begin{equation}
\sum_{\left|t\right|_{2}\leq2^{N}}\left(1-\mathbf{1}_{0}\left(t\right)\right)v_{2}\left(t\right)e^{2\pi i\left\{ t\mathfrak{z}\right\} _{2}}\overset{\overline{\mathbb{Q}}}{=}2\left|\mathfrak{z}\right|_{2}^{-1}-1,\textrm{ }\forall N\geq v_{2}\left(\mathfrak{z}\right)
\end{equation}
Here, $1-\mathbf{1}_{0}\left(t\right)$ is $0$ when $t\overset{1}{\equiv}0$
and is $1$ for all other $t$. Now, let us add the Fourier series
generated by $\hat{A}_{3}$. This gives:
\begin{equation}
\sum_{\left|t\right|_{2}\leq2^{N}}\left(\hat{A}_{3}\left(t\right)+\left(1-\mathbf{1}_{0}\left(t\right)\right)v_{2}\left(t\right)\right)e^{2\pi i\left\{ t\mathfrak{z}\right\} _{2}}\overset{\overline{\mathbb{Q}}}{=}2\left|\mathfrak{z}\right|_{2}^{-1}-1+\frac{3^{\#_{1}\left(\left[\mathfrak{z}\right]_{2^{N}}\right)}}{2^{N}}
\end{equation}
As $N\rightarrow\infty$, when $\mathfrak{z}\in\mathbb{Z}_{2}^{\prime}$,
the right-hand side converges to $2\left|\mathfrak{z}\right|_{2}^{-1}-1$
in $\mathbb{C}_{3}$ ; for $\mathbb{\mathfrak{z}}\in\mathbb{N}_{0}$,
meanwhile, the right-hand side converges to $2\left|\mathfrak{z}\right|_{2}^{-1}-1$
in $\mathbb{C}$. The same is true if we remove $\hat{A}_{3}\left(t\right)$.
As such, we have two different formulas: 
\begin{equation}
f\mapsto\sum_{t\in\hat{\mathbb{Z}}_{2}\backslash\left\{ 0\right\} }\hat{f}\left(t\right)v_{2}\left(-t\right)
\end{equation}
\begin{equation}
f\mapsto\sum_{t\in\hat{\mathbb{Z}}_{2}}\hat{f}\left(t\right)\left(\hat{A}_{3}\left(-t\right)+\left(1-\mathbf{1}_{0}\left(-t\right)\right)v_{2}\left(-t\right)\right)
\end{equation}
representing two \emph{entirely different }$\left(2,3\right)$-adic
measures, and yet, \emph{both }of them constitute perfectly reasonable
ways of defining the integral of $\left(2\left|\mathfrak{z}\right|_{2}^{-1}-1\right)f\left(\mathfrak{z}\right)$.
These are:
\begin{equation}
\int_{\mathbb{Z}_{2}}\left(2\left|\mathfrak{z}\right|_{2}^{-1}-1\right)f\left(\mathfrak{z}\right)d\mathfrak{z}\overset{\mathbb{C}_{3}}{=}\sum_{t\in\hat{\mathbb{Z}}_{2}\backslash\left\{ 0\right\} }\hat{f}\left(t\right)v_{2}\left(-t\right)
\end{equation}
and: 
\begin{equation}
\int_{\mathbb{Z}_{2}}\left(2\left|\mathfrak{z}\right|_{2}^{-1}-1\right)f\left(\mathfrak{z}\right)d\mathfrak{z}\overset{\mathbb{C}_{3}}{=}\sum_{t\in\hat{\mathbb{Z}}_{2}}\hat{f}\left(t\right)\left(\hat{A}_{3}\left(-t\right)+\left(1-\mathbf{1}_{0}\left(-t\right)\right)v_{2}\left(-t\right)\right)
\end{equation}
Both of these are valid because the $\hat{\mu}$ we are using against
$\hat{f}$ are $3$-adically bounded. However, because of $dA_{3}$'s
degeneracy, even if $\int fdA_{3}\neq0$, when we try to consider
the Fourier series generated by $\hat{\mu}$, the part generated by
$\hat{A}_{3}$ vanishes into thin air.

So, while quasi-integrability allows us to integrate discontinuous
functions, there will be no canonical choice for these integrals'
values. Indeed, as we will see in Subsection \ref{subsec:3.3.5 Quasi-Integrability},
the integral of a quasi-integrable function will only be well-defined
modulo an arbitrary degenerate measure.
\end{example}
\vphantom{}

Before we begin our study of these matters in earnest, we record the
following result, which\textemdash as described below\textemdash can
be used to define a ``$p,q$-adic Mellin transform\index{Mellin transform!left(p,qright)-adic@$\left(p,q\right)$-adic}\index{$p,q$-adic!Mellin transform}''.
This may be of great interest in expanding the scope of $\left(p,q\right)$-adic
analysis beyond what was previously thought possible.
\begin{lem}
Let $p$ and $q$ be integers $\geq2$ so that $q\mid\left(p-1\right)$.
Then, for any non-zero rational number $r$ and any $\mathfrak{z}\in\mathbb{Z}_{p}\backslash\left\{ 0\right\} $:
\begin{equation}
\left|\mathfrak{z}\right|_{p}^{-r}\overset{\mathbb{C}_{q}}{=}\sum_{\left|t\right|_{p}\leq p^{N}}\left(\frac{1}{p^{r}}\sum_{m=0}^{\infty}\binom{r}{m}\left(p-1\right)^{m}\sum_{\begin{array}{c}
0<\left|s_{1}\right|_{p},\ldots,\left|s_{m}\right|_{p}\leq p^{N}\\
s_{1}+\cdots+s_{m}=t
\end{array}}\prod_{j=1}^{m}v_{p}\left(s_{j}\right)\right)e^{2\pi i\left\{ t\mathfrak{z}\right\} _{p}}\label{eq:-rth power of the p-adic absolute value is sometimes quasi-integrable}
\end{equation}
for all $N>v_{p}\left(\mathfrak{z}\right)$. Here, the condition on
$q\mid\left(p-1\right)$ guarantees that the $m$-series converges
$q$-adically. Also, the sum over $s_{1},\ldots,s_{m}$ is defined
to be $1$ when $m=0$.

\emph{Note}: if $r$ is a positive integer, the infinite series reduces
to the finite sum $\sum_{m=0}^{r}$, and as such, the condition that
$q\mid\left(p-1\right)$ can be dropped.
\end{lem}
Proof: Letting $p$, $q$, $\mathfrak{z}$, and $N$ be as given,
we take (\ref{eq:reciprocal of the p-adic absolute value of z is quasi-integrable})
and raise it to the $r$th power: 
\begin{align*}
\left|\mathfrak{z}\right|_{p}^{-r} & =\left(\frac{1}{p}+\frac{p-1}{p}\sum_{0<\left|t\right|_{p}\leq p^{N}}v_{p}\left(t\right)e^{2\pi i\left\{ t\mathfrak{z}\right\} _{p}}\right)^{r}\\
 & =\frac{1}{p^{r}}\left(1+\left(p-1\right)\sum_{0<\left|t\right|_{p}\leq p^{N}}v_{p}\left(t\right)e^{2\pi i\left\{ t\mathfrak{z}\right\} _{p}}\right)^{r}\\
 & =\frac{1}{p^{r}}\sum_{m=0}^{\infty}\binom{r}{m}\left(p-1\right)^{m}\underbrace{\sum_{0<\left|t_{1}\right|_{p},\ldots,\left|t_{m}\right|_{p}\leq p^{N}}v_{p}\left(t_{1}\right)\cdots v_{p}\left(t_{m}\right)e^{2\pi i\left\{ \left(t_{1}+\cdots+t_{m}\right)\mathfrak{z}\right\} _{p}}}_{1\textrm{ when }m=0}
\end{align*}
Note that the use of the binomial series for $\left(1+\mathfrak{y}\right)^{r}$
is valid here, seeing as: 
\begin{equation}
\left(1+\mathfrak{y}\right)^{r}\overset{\mathbb{C}_{q}}{=}\sum_{m=0}^{\infty}\binom{r}{m}\mathfrak{y}^{m},\textrm{ }\forall\mathfrak{y}\in\mathbb{C}_{q}:\left|\mathfrak{y}\right|_{q}<1
\end{equation}
and that: 
\begin{equation}
\left|\left(p-1\right)\sum_{0<\left|t\right|_{p}\leq p^{N}}v_{p}\left(t\right)e^{2\pi i\left\{ t\mathfrak{z}\right\} _{p}}\right|_{q}\leq\left|p-1\right|_{q}\cdot\underbrace{\max_{0<\left|t\right|_{p}\leq p^{N}}\left|v_{p}\left(t\right)\right|_{q}}_{\leq1}\leq\left|p-1\right|_{q}
\end{equation}
where the divisibility of $p-1$ by $q$ then guarantees $\left|p-1\right|_{q}<1$,
and hence, that the series in $m$ converges uniformly in $t$ and
$\mathfrak{z}$.

We then have: 
\begin{eqnarray*}
 & \sum_{0<\left|t_{1}\right|_{p},\ldots,\left|t_{m}\right|_{p}\leq p^{N}}v_{p}\left(t_{1}\right)\cdots v_{p}\left(t_{m}\right)e^{2\pi i\left\{ \left(t_{1}+\cdots+t_{m}\right)\mathfrak{z}\right\} _{p}}\\
 & =\\
 & \sum_{\left|t\right|_{p}\leq p^{N}}\left(\sum_{\begin{array}{c}
0<\left|s_{1}\right|_{p},\ldots,\left|s_{m}\right|_{p}\leq p^{N}\\
s_{1}+\cdots+s_{m}=t
\end{array}}v_{p}\left(s_{1}\right)\cdots v_{p}\left(s_{m}\right)\right)e^{2\pi i\left\{ t\mathfrak{z}\right\} _{p}}
\end{eqnarray*}
Hence: 
\begin{align*}
\left|\mathfrak{z}\right|_{p}^{-r} & =\sum_{\left|t\right|_{p}\leq p^{N}}\left(\frac{1}{p^{r}}\sum_{m=0}^{\infty}\binom{r}{m}\left(p-1\right)^{m}\underbrace{\sum_{\begin{array}{c}
0<\left|s_{1}\right|_{p},\ldots,\left|s_{m}\right|_{p}\leq p^{N}\\
s_{1}+\cdots+s_{m}=t
\end{array}}\prod_{j=1}^{m}v_{p}\left(s_{j}\right)}_{1\textrm{ when }m=0}\right)e^{2\pi i\left\{ t\mathfrak{z}\right\} _{p}}
\end{align*}
for all $N>v_{p}\left(\mathfrak{z}\right)$, as desired.

Q.E.D. 
\begin{rem}[\textbf{A $\left(p,q\right)$-adic Mellin transform}]
\label{rem:pq adic mellin transform}The above allows us to define
a $\left(p,q\right)$-adic Mellin transform, $\mathscr{M}_{p,q}$,
by: 
\begin{equation}
\mathscr{M}_{p,q}\left\{ f\right\} \left(r\right)\overset{\textrm{def}}{=}\int_{\mathbb{Z}_{p}}\left|\mathfrak{z}\right|_{p}^{r-1}f\left(\mathfrak{z}\right)d\mathfrak{z},\textrm{ }\forall r\in\mathbb{Q},\textrm{ }\forall f\in C\left(\mathbb{Z}_{p},\mathbb{C}_{q}\right)\label{eq:Definition of (p,q)-adic Mellin transform}
\end{equation}
whenever $q\mid\left(p-1\right)$. In particular, we then have: 
\begin{equation}
\int_{\mathbb{Z}_{p}}\left|\mathfrak{z}\right|_{p}^{r}f\left(\mathfrak{z}\right)d\mathfrak{z}=\sum_{t\in\hat{\mathbb{Z}}_{p}}\left(p^{r}\sum_{m=0}^{\infty}\binom{1-r}{m}\left(p-1\right)^{m}\sum_{\begin{array}{c}
\mathbf{s}\in\left(\hat{\mathbb{Z}}_{p}\backslash\left\{ 0\right\} \right)^{m}\\
\Sigma\left(\mathbf{s}\right)=t
\end{array}}v_{p}\left(\mathbf{s}\right)\right)\hat{f}\left(-t\right)\label{eq:Formula for (p,q)-adic Mellin transform}
\end{equation}
where the $m$-sum is $1$ when $=0$, and where: 
\begin{align}
\Sigma\left(\mathbf{s}\right) & \overset{\textrm{def}}{=}\sum_{j=1}^{m}s_{j}\\
v_{p}\left(\mathbf{s}\right) & \overset{\textrm{def}}{=}\prod_{j=1}^{m}v_{p}\left(s_{j}\right)
\end{align}
for all $\mathbf{s}=\left(s_{1},\ldots,s_{m}\right)\in\left(\hat{\mathbb{Z}}_{p}\backslash\left\{ 0\right\} \right)^{m}$.
If we allow for $f$ to take values in the universal $q$-adic field
$\Omega_{q}$ (the spherical completion of $\mathbb{C}_{q}$, see
\cite{Robert's Book} for details), it may be possible to interpret
$\mathscr{M}_{p,q}\left\{ f\right\} \left(r\right)$ for an arbitrary
real number $r$ (which would mean that $\mathscr{M}_{p,q}\left\{ f\right\} :\mathbb{R}\rightarrow\Omega_{q}$),
although care would obviously need to be taken to properly define
and interpret (\ref{eq:Formula for (p,q)-adic Mellin transform})
in this case. Furthermore, if we can show that, given a quasi-integrable
function $\chi$, the product $\left|\mathfrak{z}\right|_{p}^{r-1}\chi\left(\mathfrak{z}\right)$
will be quasi-integrable for $r\in\mathbb{Q}$ under certain circumstances,
we could then define the $\left(p,q\right)$-adic Mellin transform
of $\chi$. 
\end{rem}
\begin{rem}
It is worth noting that this Mellin transform formalism potentially
opens the door for a $\left(p,q\right)$-adic notion of differentiability
in the sense of distributions. This notion of differentiability is
already well-established for the $\left(p,\infty\right)$-adic case
(real/complex-valued functions), having originated a 1988 paper by
V. Vladimirov, \emph{Generalized functions over the field of $p$-adic
numbers} \cite{Vladimirov - the big paper about complex-valued distributions over the p-adics}.
This paper is a comprehensive exposé of everything one needs to know
to work with distributions in the $\left(p,\infty\right)$-adic context.
This method, now known as the \textbf{Vladimirov operator}\index{Vladimirov operator}\textbf{
}or \index{$p$-adic!fractional differentiation}\textbf{$p$-adic
(fractional) differentiation }has since become standard in the areas
of $\left(p,\infty\right)$-adic mathematical theoretical physics
(``$p$-adic mathematical physics'') where it is employed; see,
for instance, the article \cite{First 30 years of p-adic mathematical physics},
which gives a summary of the first thirty years' worth of developments
in the subject.

To borrow from \cite{First 30 years of p-adic mathematical physics}
(although, it should be noted there is a typographical error in their
statement of the definition), for $\alpha\in\mathbb{C}$, the order
$\alpha$ $p$-adic fractional derivative of a function $f:\mathbb{Q}_{p}\rightarrow\mathbb{C}$
is defined\footnote{This definition is slightly more general than the one given in Subsection
\ref{subsec:3.1.1 Some-Historical-and}, which was only valid for
real $\alpha>0$; we re-obtain (\ref{eq:Definition of the Vladimirov Fractional Differentiation Operator})
from (\ref{eq:First formula for the p-adic fractional derivative})
below.} by: 
\begin{equation}
D^{\alpha}\left\{ f\right\} \left(\mathfrak{z}\right)\overset{\textrm{def}}{=}\int_{\mathbb{Q}_{p}}\left|\mathfrak{y}\right|_{p}^{\alpha}\hat{f}\left(-\mathfrak{y}\right)e^{2\pi i\left\{ \mathfrak{y}\mathfrak{z}\right\} _{p}}d\mathfrak{y}\label{eq:First formula for the p-adic fractional derivative}
\end{equation}
where $d\mathfrak{y}$ is the real-valued $p$-adic Haar measure on
$\mathbb{Q}_{p}$, normalized to be a probability measure on $\mathbb{Z}_{p}$,
and where: 
\begin{equation}
\hat{f}\left(\mathfrak{y}\right)\overset{\textrm{def}}{=}\int_{\mathbb{Q}_{p}}f\left(\mathfrak{x}\right)e^{-2\pi i\left\{ \mathfrak{y}\mathfrak{x}\right\} _{p}}d\mathfrak{x},\textrm{ }\forall\mathfrak{y}\in\mathbb{Q}_{p}\label{eq:Fourier transform over Q_p}
\end{equation}
is the Fourier transform of $f$. All of these require $f$ to be
sufficiently well-behaved\textemdash say, compactly supported. When
the support of $f$ lies in $\mathbb{Z}_{p}$, (\ref{eq:Fourier transform over Q_p})
reduces to: 
\[
\hat{f}\left(t\right)\overset{\mathbb{C}}{=}\int_{\mathbb{Z}_{p}}f\left(\mathfrak{x}\right)e^{-2\pi i\left\{ t\mathfrak{x}\right\} _{p}}d\mathfrak{x}
\]
and (\ref{eq:First formula for the p-adic fractional derivative})
becomes: 
\begin{equation}
D^{\alpha}\left\{ f\right\} \left(\mathfrak{z}\right)=\sum_{t\in\hat{\mathbb{Z}}_{p}}\left|t\right|_{p}^{\alpha}\hat{f}\left(-t\right)e^{2\pi i\left\{ t\mathfrak{z}\right\} _{p}}\label{eq:p-adic fractional derivative over Z_p}
\end{equation}
When $\alpha>0$, (\ref{eq:First formula for the p-adic fractional derivative})
can be written as: 
\begin{equation}
D^{\alpha}\left\{ f\right\} \left(\mathfrak{z}\right)=\frac{1}{\Gamma_{p}\left(-\alpha\right)}\int_{\mathbb{Q}_{p}}\frac{f\left(\mathfrak{z}\right)-f\left(\mathfrak{y}\right)}{\left|\mathfrak{z}-\mathfrak{y}\right|_{p}^{\alpha+1}}d\mathfrak{y}\label{eq:Second formula for the p-adic fractional derivative}
\end{equation}
where, recall, $\Gamma_{p}\left(-\alpha\right)$ is the physicist's
notation for the normalization constant: 
\begin{equation}
\Gamma_{p}\left(-\alpha\right)=\frac{p^{\alpha}-1}{1-p^{-1-\alpha}}
\end{equation}

\cite{First 30 years of p-adic mathematical physics} contains no
less than \emph{three-hundred forty two} references, describing innumerable
applications of this fractional derivative, including ``spectral
properties, operators on bounded regions, analogs of elliptic and
parabolic equations, a wave-type equation'', and many more. Given
that my $\left(p,q\right)$-adic implementation of the Mellin transform
(\ref{eq:Definition of (p,q)-adic Mellin transform}) can then be
used to make sense of (\ref{eq:First formula for the p-adic fractional derivative})
for any $f\in C\left(\mathbb{Z}_{p},\mathbb{C}_{q}\right)$ and any
$r\in\mathbb{Q}$, this seems like a research direction worth exploring.
\end{rem}

\subsection{\label{subsec:3.3.2 The--adic-Dirichlet}The $p$-adic Dirichlet
Kernel and $\left(p,q\right)$-adic Fourier Resummation Lemmata}

In this subsection, we will introduce specific types of $\left(p,q\right)$-adic
measures for which (\ref{eq:The Limit of Interest}) is sufficiently
well-behaved for us to construct useful spaces of quasi-integrable
functions. Consequently, this subsection will be low on concepts and
high on computations\textemdash though that's nothing in comparison
with what we'll confront in Chapters 4 and 6\textemdash but, I digress.

A recurring device in Section \ref{sec:3.3 quasi-integrability} is
to start with a $\overline{\mathbb{Q}}$-valued function defined on
$\hat{\mathbb{Z}}_{p}$ and then use it to create functions on $\mathbb{Z}_{p}$
by summing the associated Fourier series. Delightfully, this procedure
is a direct $\left(p,q\right)$-adic analogue of a fundamental construction
of classical Fourier analysis: convolution against the Dirichlet kernel.
\begin{defn}[\textbf{The $p$-adic Dirichlet Kernel}]
\index{Fourier series!$N$th partial sum}\ 

\vphantom{}

I. Let $\mathbb{F}$ be an algebraically closed field extension of
$\mathbb{Q}$, and let $\hat{\mu}:\hat{\mathbb{Z}}_{p}\rightarrow\mathbb{F}$.
For each $N\in\mathbb{N}_{0}$, we then define \nomenclature{$\tilde{\mu}_{N}$}{$N$th partial sum of Fourier series generated by $\hat{\mu}$}$\tilde{\mu}_{N}:\mathbb{Z}_{p}\rightarrow\mathbb{F}$
by: 
\begin{equation}
\tilde{\mu}_{N}\left(\mathfrak{z}\right)\overset{\textrm{def}}{=}\sum_{\left|t\right|_{p}\leq p^{N}}\hat{\mu}\left(t\right)e^{2\pi i\left\{ t\mathfrak{z}\right\} _{p}}\label{eq:Definition of mu_N twiddle}
\end{equation}
where the sum is taken over all $t\in\hat{\mathbb{Z}}_{p}$ satisfying
$\left|t\right|_{p}\leq p^{N}$. Note that $\tilde{\mu}_{N}$ is then
locally constant. As such, for any valued field $K$ extending $\mathbb{F}$,
$\tilde{\mu}_{N}$ will be continuous as a function from $\mathbb{Z}_{p}$
to $K$. We call $\tilde{\mu}_{N}$ the \textbf{$N$th partial sum
of the Fourier series generated by $\hat{\mu}$}.

\vphantom{}

II. For each $N\in\mathbb{N}_{0}$, we define the function $D_{p:N}:\mathbb{Z}_{p}\rightarrow\mathbb{Q}$
by: 
\begin{equation}
D_{p:N}\left(\mathfrak{z}\right)\overset{\textrm{def}}{=}p^{N}\left[\mathfrak{z}\overset{p^{N}}{\equiv}0\right]\label{eq:Definition of the p-adic Dirichlet Kernel}
\end{equation}
We call \nomenclature{$D_{p:N}$}{$N$th $p$-adic Dirichlet kernel}$D_{p:N}$
the $N$th \textbf{$p$-adic Dirichlet kernel}\index{$p$-adic!Dirichlet kernel}.
$\left\{ D_{p:N}\right\} _{N\geq0}$ is the \textbf{family of $p$-adic
Dirichlet kernels}. Note that since each $D_{p:N}$ is locally constant,
each is an element of $C\left(\mathbb{Z}_{p},K\right)$. Moreover,
the Fourier transform of $D_{p:N}$ has finite support, with: 
\begin{equation}
\hat{D}_{p:N}\left(t\right)=\mathbf{1}_{0}\left(p^{N}t\right)\overset{\textrm{def}}{=}\begin{cases}
1 & \textrm{if }\left|t\right|_{p}\leq p^{N}\\
0 & \textrm{if }\left|t\right|_{p}>p^{N}
\end{cases}\label{eq:Fourier Transform of the p-adic Dirichlet Kernel}
\end{equation}
We note here that since the $p$-adic Dirichlet kernels are rational
valued, we can also compute their real/complex valued Fourier transforms
with respect to the real-valued Haar probability measure on $\mathbb{Z}_{p}$,
and\textemdash moreover\textemdash the resulting $\hat{D}_{p:N}\left(t\right)$
will be \emph{the same} as the one given above. In other words, the
$\left(p,q\right)$-adic and $\left(p,\infty\right)$-adic Fourier
transforms of $D_{p:N}$ and the procedures for computing them are
\emph{formally identical}.

I use the term \textbf{$\left(p,\infty\right)$-adic Dirichlet Kernel}
when viewing the $D_{p:N}$s as taking values in $\mathbb{C}$; I
use the term \textbf{$\left(p,K\right)$-adic Dirichlet Kernel }when
viewing the $D_{p:N}$s as taking values in a metrically complete
$q$-adic field $K$; I use the term \textbf{$\left(p,q\right)$-adic
Dirichlet Kernel} for when the specific field $K$ is either not of
concern or is $\mathbb{C}_{q}$. 
\end{defn}
\begin{rem}
Regarding the name ``$p$-adic Dirichlet Kernel'', it is worth noting
that when one is doing $\left(p,\infty\right)$-adic analysis, the
family of $p$-adic Dirichlet Kernels then form an approximate identity,
in the sense that for any \nomenclature{$L^{1}\left(\mathbb{Z}_{p},\mathbb{C}\right)$}{set of absolutely integrable $f:\mathbb{Z}_{p}\rightarrow\mathbb{C}$}$f\in L^{1}\left(\mathbb{Z}_{p},\mathbb{C}\right)$
(complex-valued function on $\mathbb{Z}_{p}$, integrable with respect
to the real-valued Haar probability measure on $\mathbb{Z}_{p}$),
it can be shown that:
\begin{equation}
\lim_{N\rightarrow\infty}\left(D_{p:N}*f\right)\left(\mathfrak{z}\right)\overset{\mathbb{C}}{=}\lim_{N\rightarrow\infty}p^{N}\int_{\mathbb{Z}_{p}}\left[\mathfrak{y}\overset{p^{N}}{\equiv}\mathfrak{z}\right]f\left(\mathfrak{y}\right)d\mathfrak{y}=f\left(\mathfrak{z}\right)\label{eq:p infinity adic Lebesgue differentiation theorem}
\end{equation}
for \emph{almost} every $\mathfrak{z}\in\mathbb{Z}_{p}$; in particular,
the limit holds for all $\mathfrak{z}$ at which $f$ is continuous.
Indeed, (\ref{eq:p infinity adic Lebesgue differentiation theorem})
is an instance of the \textbf{Lebesgue Differentiation Theorem}\index{Lebesgue Differentiation Theorem}\footnote{Like in the archimedean case, the theorem is a consequence of a covering
lemma and estimates for the ($\left(p,\infty\right)$-adic) Hardy-Littlewood
maximal function (see \cite{Geometric Measure Theory}, for instance).}\textbf{ }for real and complex valued functions of a $p$-adic variable.
Our use of the $p$-adic Dirichlet Kernels in the $\left(p,q\right)$-adic
context therefore falls under the topic of summability kernels and
approximate identities\index{approximate identity}. It remains to
be seen if convolving other types of kernels with $\left(p,q\right)$-adic
measures or $\left(p,q\right)$-adic functions will yield interesting
results. 
\end{rem}
\vphantom{}

For now, let us establish the identities of import. 
\begin{prop}[\textbf{$D_{p:N}$ is an approximate identity}]
\label{prop:D_p is an approximate identity}\ 

\vphantom{}

I. Let $f\in C\left(\mathbb{Z}_{p},\mathbb{C}_{q}\right)$. Then,
for all $N\in\mathbb{N}_{0}$: 
\begin{equation}
\left(D_{p:N}*f\right)\left(\mathfrak{z}\right)\overset{\mathbb{C}_{q}}{=}\int_{\mathbb{Z}_{p}}p^{N}\left[\mathfrak{y}\overset{p^{N}}{\equiv}\mathfrak{z}\right]f\left(\mathfrak{y}\right)d\mathfrak{y}=\sum_{\left|t\right|_{p}\leq p^{N}}\hat{f}\left(t\right)e^{2\pi i\left\{ t\mathfrak{z}\right\} _{p}},\textrm{ }\forall\mathfrak{z}\in\mathbb{Z}_{p}
\end{equation}
Moreover, as $N\rightarrow\infty$, $D_{p:N}*f$ converges uniformly
in $\mathbb{C}_{q}$ to $f$: 
\begin{equation}
\lim_{N\rightarrow\infty}\sup_{\mathfrak{z}\in\mathbb{Z}_{p}}\left|f\left(\mathfrak{z}\right)-\left(D_{p:N}*f\right)\left(\mathfrak{z}\right)\right|_{q}\overset{\mathbb{R}}{=}0
\end{equation}
Also, this shows that the evaluation map $f\mapsto f\left(\mathfrak{z}_{0}\right)$
is a continuous linear functional on $C\left(\mathbb{Z}_{p},\mathbb{C}_{q}\right)$.

\vphantom{}

II. Let $d\mu\in C\left(\mathbb{Z}_{p},\mathbb{C}_{q}\right)^{\prime}$.
Then, for all $N\in\mathbb{N}_{0}$: 
\begin{equation}
\left(D_{p:N}*d\mu\right)\left(\mathfrak{z}\right)\overset{\mathbb{C}_{q}}{=}\int_{\mathbb{Z}_{p}}p^{N}\left[\mathfrak{y}\overset{p^{N}}{\equiv}\mathfrak{z}\right]d\mu\left(\mathfrak{y}\right)=\tilde{\mu}_{N}\left(\mathfrak{z}\right),\textrm{ }\forall\mathfrak{z}\in\mathbb{Z}_{p}
\end{equation}
For fixed $N$, convergence occurs uniformly with respect to $\mathfrak{z}$.
Additionally, if $\hat{\mu}$ takes values in $\overline{\mathbb{Q}}$,
the the above equality also holds in $\mathbb{C}$, and, for fixed
$N$, the convergence there is uniform with respect to $\mathfrak{z}$. 
\end{prop}
Proof: The \textbf{Convolution Theorem }(\textbf{Theorem \ref{thm:Convolution Theorem}})
for the\textbf{ }$\left(p,q\right)$-adic Fourier Transform tells
us that: 
\begin{equation}
\left(D_{p:N}*f\right)\left(\mathfrak{z}\right)=\sum_{t\in\hat{\mathbb{Z}}_{p}}\hat{D}_{p:N}\left(t\right)\hat{f}\left(t\right)e^{2\pi i\left\{ t\mathfrak{z}\right\} _{p}}
\end{equation}
and: 
\begin{equation}
\left(D_{p:N}*d\mu\right)\left(\mathfrak{z}\right)=\sum_{t\in\hat{\mathbb{Z}}_{p}}\hat{D}_{p:N}\left(t\right)\hat{\mu}\left(t\right)e^{2\pi i\left\{ t\mathfrak{z}\right\} _{p}}
\end{equation}
Using (\ref{eq:Fourier Transform of the p-adic Dirichlet Kernel})
gives the desired results. The uniform convergence of $D_{p:N}*f$
to $f$ follows from: 
\begin{equation}
\sup_{\mathfrak{z}\in\mathbb{Z}_{p}}\left|f\left(\mathfrak{z}\right)-\left(D_{p:N}*f\right)\left(\mathfrak{z}\right)\right|_{q}=\left|\sum_{\left|t\right|_{p}>p^{N}}\hat{f}\left(t\right)e^{2\pi i\left\{ t\mathfrak{z}\right\} _{p}}\right|_{q}\leq\sup_{\left|t\right|_{p}>p^{N}}\left|\hat{f}\left(t\right)\right|_{q}
\end{equation}
with the upper bound tending to $0$ in $\mathbb{R}$ as $N\rightarrow\infty$,
since $f\in C\left(\mathbb{Z}_{p},\mathbb{C}_{q}\right)\Leftrightarrow\hat{f}\in c_{0}\left(\hat{\mathbb{Z}}_{p},\mathbb{C}_{q}\right)$.

Finally, fixing $\mathfrak{z}_{0}\in\mathbb{Z}_{p}$, the maps $f\mapsto\left(D_{p:N}*f\right)\left(\mathfrak{z}_{0}\right)$
are a family of continuous linear functionals on $C\left(\mathbb{Z}_{p},\mathbb{C}_{q}\right)$
which are indexed by $N$ which converge to the evaluation map $f\mapsto f\left(\mathfrak{z}_{0}\right)$
in supremum norm as $N\rightarrow\infty$. Since the dual space of
$C\left(\mathbb{Z}_{p},\mathbb{C}_{q}\right)$ is a Banach space,
it is closed under such limits. This proves the evaluation map $f\mapsto f\left(\mathfrak{z}_{0}\right)$
is an element of $C\left(\mathbb{Z}_{p},\mathbb{C}_{q}\right)^{\prime}$.

Q.E.D.

\vphantom{}

With respect to the above, the introduction of frames in Subsection
\ref{subsec:3.3.3 Frames} is little more than a backdrop for studying
the point-wise behavior of the limit $\left(D_{p:N}*d\mu\right)\left(\mathfrak{z}\right)$
as $N\rightarrow\infty$ for various $d\mu$. Since we are potentially
going to need to utilize different topologies to make sense of this
limit, we might as well come up with a consistent set of terminology
for describing our choices of topologies, instead of having to repeatedly
specify the topologies \emph{every time }we want to talk about a limit\textemdash this
is what frames are for.

For now, though, let us continue our investigation of $\tilde{\mu}_{N}=D_{p:N}*d\mu$. 
\begin{prop}
\label{prop:Criterion for zero measure in terms of partial sums of Fourier series}Let
$d\mu\in C\left(\mathbb{Z}_{p},\mathbb{C}_{q}\right)^{\prime}$. Then
$d\mu$ is the zero measure if and only if: 
\begin{equation}
\lim_{N\rightarrow\infty}\left\Vert \tilde{\mu}_{N}\right\Vert _{p,q}\overset{\mathbb{R}}{=}0\label{eq:Criterion for a (p,q)-adic measure to be zero}
\end{equation}
where, recall, $\left\Vert \tilde{\mu}_{N}\right\Vert _{p,q}=\sup_{\mathfrak{z}\in\mathbb{Z}_{p}}\left|\tilde{\mu}_{N}\left(\mathfrak{z}\right)\right|_{q}$. 
\end{prop}
Proof:

I. If $d\mu$ is the zero measure, the $\hat{\mu}$ is identically
zero, and hence, so is $\tilde{\mu}_{N}$.

\vphantom{}

II. Conversely, suppose $\lim_{N\rightarrow\infty}\left\Vert \tilde{\mu}_{N}\right\Vert _{p,q}\overset{\mathbb{R}}{=}0$.
Then, let $f\in C\left(\mathbb{Z}_{p},\mathbb{C}_{q}\right)$ be arbitrary.
Since $\tilde{\mu}_{N}$ is continuous for each $N$ and has $\mathbf{1}_{0}\left(p^{N}t\right)\hat{\mu}\left(t\right)$
as its Fourier transform, we can write: 
\begin{equation}
\int_{\mathbb{Z}_{p}}f\left(\mathfrak{z}\right)\tilde{\mu}_{N}\left(\mathfrak{z}\right)d\mathfrak{z}\overset{\mathbb{C}_{q}}{=}\sum_{\left|t\right|_{p}\leq p^{N}}\hat{f}\left(-t\right)\hat{\mu}\left(t\right)
\end{equation}
Then, by the $\left(p,q\right)$-adic triangle inequality (\ref{eq:(p,q)-adic triangle inequality}):
\[
\left|\sum_{\left|t\right|_{p}\leq p^{N}}\hat{f}\left(-t\right)\hat{\mu}\left(t\right)\right|_{q}\overset{\mathbb{R}}{=}\left|\int_{\mathbb{Z}_{p}}f\left(\mathfrak{z}\right)\tilde{\mu}_{N}\left(\mathfrak{z}\right)d\mathfrak{z}\right|_{q}\leq\left\Vert f\cdot\tilde{\mu}_{N}\right\Vert _{p,q}\leq\left\Vert f\right\Vert _{p,q}\left\Vert \tilde{\mu}_{N}\right\Vert _{p,q}
\]
Because the upper bound tends to $0$ in $\mathbb{R}$ as $N\rightarrow\infty$,
and because: 
\begin{equation}
\lim_{N\rightarrow\infty}\sum_{\left|t\right|_{p}\leq p^{N}}\hat{f}\left(-t\right)\hat{\mu}\left(t\right)\overset{\mathbb{C}_{q}}{=}\int_{\mathbb{Z}_{p}}f\left(\mathfrak{z}\right)d\mu\left(\mathfrak{z}\right)
\end{equation}
this shows that: 
\begin{equation}
\left|\int_{\mathbb{Z}_{p}}f\left(\mathfrak{z}\right)d\mu\left(\mathfrak{z}\right)\right|_{q}\overset{\mathbb{R}}{=}\lim_{N\rightarrow\infty}\left|\sum_{\left|t\right|_{p}\leq p^{N}}\hat{f}\left(-t\right)\hat{\mu}\left(t\right)\right|_{q}\leq\left\Vert f\right\Vert _{p,q}\cdot\lim_{N\rightarrow\infty}\left\Vert \tilde{\mu}_{N}\right\Vert _{p,q}=0
\end{equation}
Since $f$ was arbitrary, we conclude that $d\mu$ is the zero measure.

Q.E.D.

\vphantom{}Next, we have the \textbf{Fourier resummation lemmata}
which\index{resummation lemmata} will synergize with our analysis
of quasi-integrable functions. To make these easier to state and remember,
here are some definitions describing the conditions which the lemmata
will impose on $\hat{\mu}$.
\begin{defn}
Consider a field $\mathbb{F}$ of characteristic $0$, a function
$\kappa:\mathbb{N}_{0}\rightarrow\mathbb{F}$, and let $K$ be a metrically
complete valued field extension of $\mathbb{F}$ which contains $\kappa\left(\mathbb{N}_{0}\right)$.

\vphantom{}

I. We say $\kappa$ is \index{tame}\textbf{$\left(p,K\right)$-adically
tame} on a set $X\subseteq\mathbb{Z}_{p}$ whenever $\lim_{n\rightarrow\infty}\left|\kappa\left(\left[\mathfrak{z}\right]_{p^{n}}\right)\right|_{K}=0$
for all $\mathfrak{z}\in X$. We do not mention $X$ when $X=\mathbb{Z}_{p}$.
If $K$ is a $q$-adic field, we say $\kappa$ is $\left(p,q\right)$-adically
tame on $X$; if $K$ is archimedean, we say $\kappa$ is $\left(p,\infty\right)$-adically
tame on $X$.

\vphantom{}

II. For a prime $p$, we say $\kappa$ has \index{$p$-adic!structure}\textbf{$p$-adic
structure}\footnote{It would seem that this property is necessary in order for our computations
to simplify in interesting ways. Life is much more difficult without
it.}\textbf{ }whenever there are constants $a_{0},\ldots,a_{p-1}\in\kappa\left(\mathbb{N}_{0}\right)$
so that: 
\begin{equation}
\kappa\left(pn+j\right)=a_{j}\kappa\left(n\right),\textrm{ }\forall n\geq0,\textrm{ }\forall j\in\left\{ 0,\ldots,p-1\right\} \label{eq:Definition of left-ended p-adic structural equations}
\end{equation}
We call (\ref{eq:Definition of left-ended p-adic structural equations})
the \textbf{left-ended structural equations }of $\kappa$. 
\end{defn}
\begin{rem}
Note that the set of functions $\kappa:\mathbb{N}_{0}\rightarrow\mathbb{F}$
with $p$-adic structure is then a linear space over $\mathbb{F}$,
as well as over any field extension $K$ of $\mathbb{F}$. 
\end{rem}
\begin{rem}
Comparing the above to what we did with $\chi_{H}$, one could distinguish
between (\ref{eq:Definition of left-ended p-adic structural equations})
(which could be called \textbf{linear $p$-adic structure}) and a
more $\chi_{H}$-like structure: 
\begin{equation}
\kappa\left(pn+j\right)=a_{j}\kappa\left(n\right)+b_{j}
\end{equation}
which we would call \textbf{affine $p$-adic structure}. For the present
purposes however, we shall only concern ourselves with linear $p$-adic
structure, so all mentions of ``$p$-adic structure'' in this and
any subsequent part of this dissertation means only the \emph{linear
}variety of $p$-adic structure defined in (\ref{eq:Definition of left-ended p-adic structural equations}),
unless specifically stated otherwise. 
\end{rem}
As defined, functions with $p$-adic structure clearly generalize
the systems of functional equations that characterized $\chi_{H}$
in Chapter 2. Because the functional equations satisfied by these
functions appear to be crucial when it comes to working with quasi-integrable
functions, it is important that we establish their properties as a
general class of functions.

As for the terminology, I use ``left-ended'' to describe the equations
in (\ref{eq:Definition of left-ended p-adic structural equations})
because, given an integer $m=pn+j$ (where $j=\left[m\right]_{p}$),
they show what happens when we pull out the left-most $p$-adic digit
of $m$. Naturally, there is a ``right-ended'' counterpart to (\ref{eq:Definition of left-ended p-adic structural equations}),
obtained by pulling out the right-most $p$-adic digit of $m$. Moreover,
as the next proposition shows, $\kappa$ satisfies one of these two
systems of functional equations if and only if it satisfies both systems. 
\begin{prop}[\textbf{Properties of functions with $p$-adic structure}]
\label{prop:properties of functions with p-adic structure} Let $\mathbb{F}$
be any field of characteristic zero, and consider a function $\kappa:\mathbb{N}_{0}\rightarrow\mathbb{F}$.
Then, $\kappa$ has $p$-adic structure if and only if there are constants
$a_{0}^{\prime},\ldots,a_{p-1}^{\prime}$ so that: 
\begin{equation}
\kappa\left(n+jp^{k}\right)=a_{j}^{\prime}\kappa\left(n\right),\textrm{ }\forall n\geq0,\textrm{ }\forall m\in\left\{ 0,\ldots,p^{n}-1\right\} ,\textrm{ }\forall j\in\left\{ 1,\ldots,p-1\right\} \label{eq:Definition of right-ended p-adic structural equations}
\end{equation}
We call \emph{(\ref{eq:Definition of right-ended p-adic structural equations})}
the \textbf{right-ended structural equations }of $\kappa$. 
\end{prop}
Proof: Fix an integer $m\in\mathbb{N}_{0}$. Then, we can write $m$
uniquely in $p$-adic form as: 
\begin{equation}
m=\sum_{\ell=0}^{\lambda_{p}\left(m\right)-1}m_{\ell}p^{\ell}
\end{equation}
where $m_{\ell}\in\left\{ 0,\ldots,p-1\right\} $ for all $\ell$,
with $m_{\lambda_{p}\left(m\right)-1}\neq0$.

I. Let $\kappa$ satisfy the left-ended structural equations (\ref{eq:Definition of left-ended p-adic structural equations}).
Then: 
\begin{align*}
\kappa\left(m\right) & =\kappa\left(p\left(\sum_{\ell=1}^{\lambda_{p}\left(m\right)-1}m_{\ell}p^{\ell-1}\right)+m_{0}\right)\\
\left(\textrm{use }(\ref{eq:Definition of left-ended p-adic structural equations})\right); & =a_{m_{0}}\kappa\left(\sum_{\ell=1}^{\lambda_{p}\left(m\right)-1}m_{\ell}p^{\ell-1}\right)\\
 & \vdots\\
 & =\left(\prod_{\ell=0}^{\lambda_{p}\left(m\right)-1}a_{m_{\ell}}\right)\kappa\left(0\right)
\end{align*}
Pulling out the right-most $p$-adic digit of $m$ gives: 
\begin{equation}
m=\left(\sum_{\ell=0}^{\lambda_{p}\left(m\right)-2}m_{\ell}p^{\ell-1}\right)+m_{\lambda_{p}\left(m\right)-1}p^{\lambda_{p}\left(m\right)-1}
\end{equation}
Since the above implies: 
\begin{equation}
\kappa\left(\sum_{\ell=0}^{\lambda_{p}\left(m\right)-2}m_{\ell}p^{\ell-1}\right)=\left(\prod_{\ell=0}^{\lambda_{p}\left(m\right)-2}a_{m_{\ell}}\right)\kappa\left(0\right)
\end{equation}
it follows that: 
\begin{align*}
\left(\prod_{\ell=1}^{\lambda_{p}\left(m\right)-1}a_{m_{\ell}}\right)\kappa\left(0\right) & =\kappa\left(m\right)\\
 & =\kappa\left(\sum_{\ell=0}^{\lambda_{p}\left(m\right)-2}m_{\ell}p^{\ell-1}+m_{\lambda_{p}\left(m\right)-1}p^{\lambda_{p}\left(m\right)-1}\right)\\
 & =\left(\prod_{\ell=0}^{\lambda_{p}\left(m\right)-1}a_{m_{\ell}}\right)\kappa\left(0\right)\\
 & =a_{m_{\lambda_{p}\left(m\right)-1}}\left(\prod_{\ell=0}^{\lambda_{p}\left(m\right)-2}a_{m_{\ell}}\right)\kappa\left(0\right)\\
 & =a_{m_{\lambda_{p}\left(m\right)-1}}\kappa\left(\sum_{\ell=0}^{\lambda_{p}\left(m\right)-2}m_{\ell}p^{\ell-1}\right)
\end{align*}
This proves (\ref{eq:Definition of right-ended p-adic structural equations}).

\vphantom{}

II. Supposing instead that $\kappa$ satisfies (\ref{eq:Definition of right-ended p-adic structural equations}),
perform the same argument as for (I), but with $a^{\prime}$s instead
of $a$s, and pulling out $a_{m_{0}}$ instead of $a_{m_{\lambda_{p}\left(m\right)-1}}$.
This yields (\ref{eq:Definition of left-ended p-adic structural equations}).

Q.E.D.
\begin{prop}
\label{prop:value of a p-adic structure at 0}Let $\kappa:\mathbb{N}_{0}\rightarrow\mathbb{F}$
have linear $p$-adic structure, and be non-identically zero. Then
the constants $a_{j}$ and $a_{j}^{\prime}$ from the left- and right-ended
structural equations for $\kappa$ are: 
\begin{align}
a_{j} & =\kappa\left(j\right)\\
a_{j}^{\prime} & =\frac{\kappa\left(j\right)}{\kappa\left(0\right)}
\end{align}
\end{prop}
Proof: Use the structural equations to evaluate $\kappa\left(j\right)$
for $j\in\left\{ 0,\ldots,p-1\right\} $.

Q.E.D. 
\begin{prop}
\label{prop:formula for functions with p-adic structure}Let $\left\{ c_{n}\right\} _{n\geq0}$
be a sequence in $\left\{ 0,\ldots,p-1\right\} $. If $\kappa$ is
non-zero and has linear $p$-adic structure: 
\begin{equation}
\kappa\left(\sum_{m=0}^{N-1}c_{n}p^{m}\right)=\frac{1}{\left(\kappa\left(0\right)\right)^{N}}\prod_{n=0}^{N-1}\kappa\left(c_{n}\right)\label{eq:Kappa action in terms of p-adic digits}
\end{equation}
\end{prop}
Proof: Use $\kappa$'s right-ended structural equations along with
the formula for the $a_{j}^{\prime}$s given in \textbf{Proposition
\ref{prop:value of a p-adic structure at 0}}.

Q.E.D. 
\begin{lem}
\label{lem:structural equations uniquely determine p-adic structured functions}Let
$\mathbb{F}$ be any field of characteristic zero. Then, every $p$-adically
structured $\kappa:\mathbb{N}_{0}\rightarrow\mathbb{F}$ is uniquely
determined by its structural equations. Moreover:

\vphantom{}

I. $\kappa$ is identically zero if and only if $\kappa\left(0\right)=0$.

\vphantom{}

II. There exists an $n\in\mathbb{N}_{0}$ so that $\kappa\left(n\right)\neq0$
if and only if $\kappa\left(0\right)=1$. 
\end{lem}
Proof: We first prove that $\kappa$ is identically zero if and only
if $\kappa\left(0\right)=0$.

To do this, we need only show that $\kappa\left(0\right)=0$ implies
$\kappa$ is identically zero. So, let $x$ be an arbitrary integer
$\geq1$. Observe that every positive integer $x$ can be written
in the form $x_{-}+jp^{n}$, where $n=\lambda_{p}\left(x\right)$,
where $x_{-}$ is the integer obtained by deleting $j$\textemdash the
right-most digit of $x$'s $p$-adic representation. Moreover, by
construction, $x_{-}$ is in $\left\{ 0,\ldots,p^{n}-1\right\} $.
So, since $\kappa$ has linear $p$-adic structure, the right-ended
structural equations yield: 
\begin{equation}
\kappa\left(x\right)=\kappa\left(x_{-}+jp^{n}\right)=c_{j}\kappa\left(x_{-}\right)
\end{equation}
Noting that the map $n\in\mathbb{N}_{0}\mapsto n_{-}\in\mathbb{N}_{0}$
eventually iterates every non-negative integer to $0$, we can write:
\begin{equation}
\kappa\left(x\right)=\left(\prod_{j}c_{j}\right)\kappa\left(0\right)=0
\end{equation}
Here, the product is taken over finitely many $j$s. In particular,
these $j$s are precisely the non-zero $p$-adic digits of $x$. Since
$x$ was an arbitrary integer $\geq1$, and since $\kappa$ was given
to vanish on $\left\{ 0,\ldots,p-1\right\} $, this shows that $\kappa$
is identically zero whenever $\kappa\left(0\right)=0$. Conversely,
if $\kappa$ is identically zero, then, necessarily, $\kappa\left(0\right)=0$.

Next, suppose $\kappa$ is $p$-adically structured and is \emph{not}
identically zero. Then, by the left-ended structural equations (\ref{eq:Definition of left-ended p-adic structural equations}),
we have $\kappa\left(0\right)=a_{0}\kappa\left(0\right)$, where $\kappa\left(0\right)\neq0$.
This forces $a_{0}=1$. Applying \textbf{Proposition \ref{prop:value of a p-adic structure at 0}}
gives $\kappa\left(0\right)=a_{0}$. So, $\kappa\left(0\right)=1$.

Finally, let $\kappa_{1}$ and $\kappa_{2}$ be two $p$-adically
structured functions satisfying the same set of structural equations.
Then, $\kappa_{1}-\kappa_{2}$ is then $p$-adically structured and
satisfies: 
\begin{equation}
\left(\kappa_{1}-\kappa_{2}\right)\left(0\right)=\kappa_{1}\left(0\right)-\kappa_{2}\left(0\right)=1-1=0
\end{equation}
which, as we just saw, then forces $\kappa_{1}-\kappa_{2}$ to be
identically zero. This proves that $\kappa$ is uniquely determined
by its structural equations.

Q.E.D.

\vphantom{}

For us, the most important case will be when $\kappa$ is both $\left(p,K\right)$-adically
tame over $\mathbb{Z}_{p}^{\prime}$ \emph{and }$p$-adically structured. 
\begin{defn}
Let $\mathbb{F}$ be a field of characteristic $0$, and let $K$
be a metrically complete valued field extension of $\mathbb{F}$.
We say $\kappa$ is \textbf{$\left(p,K\right)$-adically regular}
whenever\index{$p,q$-adic!regular function} it has linear $p$-adic
structure and is $\left(p,K\right)$-adically tame over $\mathbb{Z}_{p}^{\prime}$.
In the case $K$ is a $q$-adic field, we call this \textbf{$\left(p,q\right)$-adic
regularity};\textbf{ }when $K$ is $\mathbb{R}$ or $\mathbb{C}$,
we call this \textbf{$\left(p,\infty\right)$-adic regularity} when
$K=\mathbb{C}$. 
\end{defn}
\vphantom{}

This might seem like a lot of definitions, but we are really only
just beginning to explore the possibilities. The next few definitions
describe various important ``standard forms'' for $\left(p,q\right)$-adic
measures in terms of their Fourier-Stieltjes coefficients. 
\begin{defn}
Let $K$ be an algebraically closed field of characteristic zero and
consider a function $\hat{\mu}:\hat{\mathbb{Z}}_{p}\rightarrow K$.

\vphantom{}

I. We say $\hat{\mu}$ is \textbf{radially symmetric }whenever\index{radially symmetric}:
\begin{equation}
\hat{\mu}\left(t\right)=\hat{\mu}\left(\frac{1}{\left|t\right|_{p}}\right),\textrm{ }\forall t\in\hat{\mathbb{Z}}_{p}\backslash\left\{ 0\right\} \label{eq:Definition of a radially symmetric function}
\end{equation}
That is, the value of $\hat{\mu}\left(t\right)$ depends only on the
$p$-adic absolute value of $t$. Similarly, we say a measure $d\mu$
is \textbf{radially symmetric} whenever it has a radially symmetric
Fourier-Stieltjes transform.

\vphantom{}

II. We say $\hat{\mu}$ is \textbf{magnitudinal }whenever\index{magnitudinal}
there is a function $\kappa:\mathbb{N}_{0}\rightarrow K$ satisfying:
\begin{equation}
\hat{\mu}\left(t\right)=\begin{cases}
\kappa\left(0\right) & \textrm{if }t=0\\
\sum_{m=0}^{\left|t\right|_{p}-1}\kappa\left(m\right)e^{-2\pi imt} & \textrm{else}
\end{cases},\textrm{ }\forall t\in\hat{\mathbb{Z}}_{p}\label{eq:Definition of a Magnitudinal function}
\end{equation}
Here, we write $e^{-2\pi imt}$ to denote the $m$th power of the
particular $\left|t\right|_{p}$th root of unity in $K$ indicated
by $e^{-2\pi imt}$. We call an expression of the form of the right-hand
side of (\ref{eq:Definition of a Magnitudinal function}) a \textbf{magnitudinal
series}. Similarly, we say a measure $d\mu$ is \textbf{magnitudinal}
whenever it has a magnitude dependent Fourier-Stieltjes transform.

Additionally, for $q$ prime or $q=\infty$, we say that such a $\hat{\mu}$
and $d\mu$ are \textbf{$\left(p,q\right)$-adically regular and magnitudinal}
whenever $\kappa$ is $\left(p,q\right)$-adically regular. We say
$\hat{\mu}$ and $d\mu$ \textbf{have $p$-adic structure }whenever
$\kappa$ does.

\vphantom{}

III. We say a function $\hat{\mu}:\hat{\mathbb{Z}}_{p}\rightarrow\mathbb{F}$
is \index{radially-magnitudinal}\textbf{radially-magnitudinal }if
it is of the form $\hat{\mu}\left(t\right)=\hat{\nu}\left(t\right)\hat{\eta}\left(t\right)$
for a radially symmetric $\hat{\nu}$ and a magnitude\index{$p,q$-adic!regular magnitudinal measure}
dependent $\hat{\eta}$. A \textbf{$\left(p,q\right)$-adically regular,
radially-magnitudinal }$\hat{\mu}$ is one for which $\hat{\eta}$
is $\left(p,q\right)$-adically magnitudinal. \index{measure!regular, radially-magnitudinal} 
\end{defn}
\begin{lem}[\textbf{Fourier Resummation \textendash{} Radially Symmetric Measures}]
\textbf{\label{lem:Rad-Sym Fourier Resum Lemma}}Let\index{resummation lemmata}
$\hat{\mu}:\hat{\mathbb{Z}}_{p}\rightarrow\overline{\mathbb{Q}}$
be radially symmetric. Then: 
\begin{equation}
\sum_{\left|t\right|_{p}\leq p^{N}}\hat{\mu}\left(t\right)e^{2\pi i\left\{ t\mathfrak{z}\right\} _{p}}\overset{\overline{\mathbb{Q}}}{=}\begin{cases}
\hat{\mu}\left(0\right)+\left(p-1\right)\sum_{n=1}^{N}p^{n-1}\hat{\mu}\left(\frac{1}{p^{n}}\right) & \textrm{for }\mathfrak{z}=0\textrm{ and }N\geq0\\
\sum_{n=0}^{v_{p}\left(\mathfrak{z}\right)}\left(\hat{\mu}\left(p^{-n}\right)-\hat{\mu}\left(p^{-n-1}\right)\right)p^{n} & \textrm{for }\mathfrak{z}\in\mathbb{Z}_{p}^{\prime}\textrm{ \& }N>v_{p}\left(\mathfrak{z}\right)
\end{cases}\label{eq:Radial Fourier Resummation Lemma}
\end{equation}
The lower line can also be written as: 
\begin{equation}
\sum_{n=0}^{v_{p}\left(\mathfrak{z}\right)}\left(\hat{\mu}\left(p^{-n}\right)-\hat{\mu}\left(p^{-n-1}\right)\right)p^{n}=\hat{\mu}\left(0\right)-\left|\mathfrak{z}\right|_{p}^{-1}\hat{\mu}\left(\frac{\left|\mathfrak{z}\right|_{p}}{p}\right)+\left(1-\frac{1}{p}\right)\sum_{n=1}^{v_{p}\left(\mathfrak{z}\right)}\hat{\mu}\left(\frac{1}{p^{n}}\right)p^{n}\label{eq:Radial Fourier Resummation - Lower Line Simplification}
\end{equation}
\end{lem}
Proof: When $\mathfrak{z}=0$: 
\begin{align*}
\sum_{\left|t\right|_{p}\leq p^{N}}\hat{\mu}\left(t\right)e^{2\pi i\left\{ t\mathfrak{z}\right\} _{p}} & =\sum_{\left|t\right|_{p}\leq p^{N}}\hat{\mu}\left(t\right)\\
 & =\hat{\mu}\left(0\right)+\sum_{n=1}^{N}\sum_{\left|t\right|_{p}=p^{n}}\hat{\mu}\left(t\right)\\
\left(\hat{\mu}\textrm{ is radially symm.}\right); & =\hat{\mu}\left(0\right)+\sum_{n=1}^{N}\sum_{\left|t\right|_{p}=p^{n}}\hat{\mu}\left(\frac{1}{p^{n}}\right)\\
\left(\sum_{\left|t\right|_{p}=p^{n}}1=\varphi\left(p^{n}\right)=\left(p-1\right)p^{n-1}\right); & =\hat{\mu}\left(0\right)+\left(p-1\right)\sum_{n=1}^{N}p^{n-1}\hat{\mu}\left(\frac{1}{p^{n}}\right)
\end{align*}

Next, letting $\mathfrak{z}$ be non-zero:

\begin{align*}
\sum_{\left|t\right|_{p}\leq p^{N}}\hat{\mu}\left(t\right)e^{2\pi i\left\{ t\mathfrak{z}\right\} _{p}} & =\hat{\mu}\left(0\right)+\sum_{n=1}^{N}\sum_{\left|t\right|_{p}=p^{n}}\hat{\mu}\left(t\right)e^{2\pi i\left\{ t\mathfrak{z}\right\} _{p}}\\
\left(\hat{\mu}\textrm{ is radially symm.}\right); & =\hat{\mu}\left(0\right)+\sum_{n=1}^{N}\sum_{\left|t\right|_{p}=p^{n}}\hat{\mu}\left(\frac{1}{\left|t\right|_{p}}\right)e^{2\pi i\left\{ t\mathfrak{z}\right\} _{p}}\\
 & =\hat{\mu}\left(0\right)+\sum_{n=1}^{N}\hat{\mu}\left(p^{-n}\right)\sum_{\left|t\right|_{p}=p^{n}}e^{2\pi i\left\{ t\mathfrak{z}\right\} _{p}}\\
 & =\hat{\mu}\left(0\right)+\sum_{n=1}^{N}\hat{\mu}\left(p^{-n}\right)\left(p^{n}\left[\mathfrak{z}\overset{p^{n}}{\equiv}0\right]-p^{n-1}\left[\mathfrak{z}\overset{p^{n-1}}{\equiv}0\right]\right)\\
 & =\hat{\mu}\left(0\right)+\sum_{n=1}^{N}\hat{\mu}\left(p^{-n}\right)p^{n}\left[\mathfrak{z}\overset{p^{n}}{\equiv}0\right]-\sum_{n=1}^{N}\hat{\mu}\left(p^{-n}\right)p^{n-1}\left[\mathfrak{z}\overset{p^{n-1}}{\equiv}0\right]\\
 & =\hat{\mu}\left(0\right)+\sum_{n=1}^{N}\hat{\mu}\left(p^{-n}\right)p^{n}\left[\mathfrak{z}\overset{p^{n}}{\equiv}0\right]-\sum_{n=0}^{N-1}\hat{\mu}\left(p^{-n-1}\right)p^{n}\left[\mathfrak{z}\overset{p^{n}}{\equiv}0\right]\\
\left(\hat{\mu}\left(t+1\right)=\hat{\mu}\left(t\right)\right); & =p^{N}\hat{\mu}\left(p^{-N}\right)\left[\mathfrak{z}\overset{p^{N}}{\equiv}0\right]+\sum_{n=0}^{N-1}p^{n}\left(\hat{\mu}\left(p^{-n}\right)-\hat{\mu}\left(p^{-n-1}\right)\right)\left[\mathfrak{z}\overset{p^{n}}{\equiv}0\right]\\
\left(\mathfrak{z}\overset{p^{n}}{\equiv}0\Leftrightarrow n\leq v_{p}\left(\mathfrak{z}\right)\right); & =p^{N}\hat{\mu}\left(p^{-N}\right)\left[\mathfrak{z}\overset{p^{N}}{\equiv}0\right]+\sum_{n=0}^{\min\left\{ v_{p}\left(\mathfrak{z}\right)+1,N\right\} -1}\left(\hat{\mu}\left(p^{-n}\right)-\hat{\mu}\left(p^{-n-1}\right)\right)p^{n}
\end{align*}

Since $\mathfrak{z}\neq0$, $v_{p}\left(\mathfrak{z}\right)$ is an
integer, and, as such, $\min\left\{ v_{p}\left(\mathfrak{z}\right)+1,N\right\} =v_{p}\left(\mathfrak{z}\right)+1$
for all sufficiently large $N$\textemdash specifically, all $N\geq1+v_{p}\left(\mathfrak{z}\right)$.
For such $N$, we also have that $\mathfrak{z}\overset{p^{N}}{\cancel{\equiv}}0$,
and thus, that: 
\begin{equation}
\sum_{\left|t\right|_{p}\leq p^{N}}\hat{\mu}\left(t\right)e^{2\pi i\left\{ t\mathfrak{z}\right\} _{p}}=\sum_{n=0}^{v_{p}\left(\mathfrak{z}\right)}\left(\hat{\mu}\left(p^{-n}\right)-\hat{\mu}\left(p^{-n-1}\right)\right)p^{n},\textrm{ }\forall\mathfrak{z}\in\mathbb{Z}_{p}^{\prime},\textrm{ }\forall N>v_{p}\left(\mathfrak{z}\right)
\end{equation}
As for (\ref{eq:Radial Fourier Resummation - Lower Line Simplification}),
we can write: 
\begin{align*}
\sum_{n=0}^{v_{p}\left(\mathfrak{z}\right)}\left(\hat{\mu}\left(p^{-n}\right)-\hat{\mu}\left(p^{-n-1}\right)\right)p^{n} & =\sum_{n=0}^{v_{p}\left(\mathfrak{z}\right)}\hat{\mu}\left(\frac{1}{p^{n}}\right)p^{n}-\frac{1}{p}\sum_{n=0}^{v_{p}\left(\mathfrak{z}\right)}\hat{\mu}\left(\frac{1}{p^{n+1}}\right)p^{n+1}\\
 & =\sum_{n=0}^{v_{p}\left(\mathfrak{z}\right)}\hat{\mu}\left(\frac{1}{p^{n}}\right)p^{n}-\frac{1}{p}\sum_{n=1}^{v_{p}\left(\mathfrak{z}\right)+1}\hat{\mu}\left(\frac{1}{p^{n}}\right)p^{n}\\
 & =\hat{\mu}\left(0\right)-p^{v_{p}\left(\mathfrak{z}\right)}\hat{\mu}\left(\frac{p^{-v_{p}\left(\mathfrak{z}\right)}}{p}\right)+\left(1-\frac{1}{p}\right)\sum_{n=1}^{v_{p}\left(\mathfrak{z}\right)}\hat{\mu}\left(\frac{1}{p^{n}}\right)p^{n}\\
 & =\hat{\mu}\left(0\right)-\left|\mathfrak{z}\right|_{p}^{-1}\hat{\mu}\left(\frac{\left|\mathfrak{z}\right|_{p}}{p}\right)+\left(1-\frac{1}{p}\right)\sum_{n=1}^{v_{p}\left(\mathfrak{z}\right)}\hat{\mu}\left(\frac{1}{p^{n}}\right)p^{n}
\end{align*}

Q.E.D. 
\begin{lem}[\textbf{Fourier Resummation \textendash{} Magnitudinal measures}]
 Let \index{resummation lemmata}$\hat{\mu}:\hat{\mathbb{Z}}_{p}\rightarrow\overline{\mathbb{Q}}$
be magnitudinal and not identically $0$. Then, for all $N\in\mathbb{N}_{1}$
and all $\mathfrak{z}\in\mathbb{Z}_{p}$: 
\begin{equation}
\sum_{\left|t\right|_{p}\leq p^{N}}\hat{\mu}\left(t\right)e^{2\pi i\left\{ t\mathfrak{z}\right\} _{p}}\overset{\overline{\mathbb{Q}}}{=}p^{N}\kappa\left(\left[\mathfrak{z}\right]_{p^{N}}\right)-\sum_{n=0}^{N-1}\sum_{j=1}^{p-1}p^{n}\kappa\left(\left[\mathfrak{z}\right]_{p^{n}}+jp^{n}\right)\label{eq:Magnitude Fourier Resummation Lemma}
\end{equation}
If, in addition, $\hat{\mu}$ has $p$-adic structure, then: 
\begin{equation}
\sum_{\left|t\right|_{p}\leq p^{N}}\hat{\mu}\left(t\right)e^{2\pi i\left\{ t\mathfrak{z}\right\} _{p}}\overset{\overline{\mathbb{Q}}}{=}p^{N}\kappa\left(\left[\mathfrak{z}\right]_{p^{N}}\right)-\left(\sum_{j=1}^{p-1}\kappa\left(j\right)\right)\sum_{n=0}^{N-1}p^{n}\kappa\left(\left[\mathfrak{z}\right]_{p^{n}}\right)\label{eq:Magnitudinal Fourier Resummation Lemma - p-adically distributed case}
\end{equation}
\end{lem}
Proof:

\begin{align*}
\sum_{\left|t\right|_{p}\leq p^{N}}\hat{\mu}\left(t\right)e^{2\pi i\left\{ t\mathfrak{z}\right\} _{p}} & =\hat{\mu}\left(0\right)+\sum_{n=1}^{N}\sum_{\left|t\right|_{p}=p^{n}}\hat{\mu}\left(t\right)e^{2\pi i\left\{ t\mathfrak{z}\right\} _{p}}\\
 & =\hat{\mu}\left(0\right)+\sum_{n=1}^{N}\sum_{\left|t\right|_{p}=p^{n}}\left(\sum_{m=0}^{p^{n}-1}\kappa\left(m\right)e^{-2\pi imt}\right)e^{2\pi i\left\{ t\mathfrak{z}\right\} _{p}}\\
 & =\hat{\mu}\left(0\right)+\sum_{n=1}^{N}\sum_{m=0}^{p^{n}-1}\kappa\left(m\right)\sum_{\left|t\right|_{p}=p^{n}}e^{2\pi i\left\{ t\left(\mathfrak{z}-m\right)\right\} _{p}}\\
 & =\hat{\mu}\left(0\right)+\sum_{n=1}^{N}\sum_{m=0}^{p^{n}-1}\kappa\left(m\right)\left(p^{n}\left[\mathfrak{z}\overset{p^{n}}{\equiv}m\right]-p^{n-1}\left[\mathfrak{z}\overset{p^{n-1}}{\equiv}m\right]\right)
\end{align*}
Next, note that as $m$ varies in $\left\{ 0,\ldots,p^{N}-1\right\} $,
the only value of $m$ satisfying $\mathfrak{z}\overset{p^{N}}{\equiv}m$
is $m=\left[\mathfrak{z}\right]_{p^{N}}$. Consequently: 
\begin{equation}
\sum_{m=0}^{p^{n}-1}\kappa\left(m\right)p^{n}\left[\mathfrak{z}\overset{p^{n}}{\equiv}m\right]=p^{n}\kappa\left(\left[\mathfrak{z}\right]_{p^{n}}\right)
\end{equation}
On the other hand: 
\begin{align*}
\sum_{m=0}^{p^{n}-1}\kappa\left(m\right)p^{n-1}\left[\mathfrak{z}\overset{p^{n-1}}{\equiv}m\right] & =\sum_{m=p^{n-1}}^{p^{n}-1}\kappa\left(m\right)p^{n-1}\left[\mathfrak{z}\overset{p^{n-1}}{\equiv}m\right]+\sum_{m=0}^{p^{n-1}-1}\kappa\left(m\right)p^{n-1}\left[\mathfrak{z}\overset{p^{n-1}}{\equiv}m\right]\\
 & =\sum_{m=0}^{p^{n}-p^{n-1}-1}\kappa\left(m+p^{n-1}\right)p^{n-1}\left[\mathfrak{z}\overset{p^{n-1}}{\equiv}m+p^{n-1}\right]+p^{n-1}\kappa\left(\left[\mathfrak{z}\right]_{p^{n-1}}\right)\\
 & =\sum_{m=0}^{\left(p-1\right)p^{n-1}-1}\kappa\left(m+p^{n-1}\right)p^{n-1}\left[\mathfrak{z}\overset{p^{n-1}}{\equiv}m\right]+p^{n-1}\kappa\left(\left[\mathfrak{z}\right]_{p^{n-1}}\right)\\
 & =\sum_{j=1}^{p-1}\sum_{m=\left(j-1\right)p^{n-1}}^{jp^{n-1}-1}\kappa\left(m+p^{n-1}\right)p^{n-1}\left[\mathfrak{z}\overset{p^{n-1}}{\equiv}m\right]+p^{n-1}\kappa\left(\left[\mathfrak{z}\right]_{p^{n-1}}\right)\\
 & =\sum_{j=1}^{p-1}\sum_{m=0}^{p^{n-1}-1}\kappa\left(m+jp^{n-1}\right)p^{n-1}\left[\mathfrak{z}\overset{p^{n-1}}{\equiv}m\right]+p^{n-1}\kappa\left(\left[\mathfrak{z}\right]_{p^{n-1}}\right)\\
 & =\sum_{j=1}^{p-1}p^{n-1}\kappa\left(\left[\mathfrak{z}\right]_{p^{n-1}}+jp^{n-1}\right)+p^{n-1}\kappa\left(\left[\mathfrak{z}\right]_{p^{n-1}}\right)\\
 & =\sum_{j=0}^{p-1}p^{n-1}\kappa\left(\left[\mathfrak{z}\right]_{p^{n-1}}+jp^{n-1}\right)
\end{align*}
Thus: 
\begin{align*}
\sum_{\left|t\right|_{p}\leq p^{N}}\hat{\mu}\left(t\right)e^{2\pi i\left\{ t\mathfrak{z}\right\} _{p}} & =\hat{\mu}\left(0\right)+\sum_{n=1}^{N}\sum_{m=0}^{p^{n}-1}\kappa\left(m\right)\left(p^{n}\left[\mathfrak{z}\overset{p^{n}}{\equiv}m\right]-p^{n-1}\left[\mathfrak{z}\overset{p^{n-1}}{\equiv}m\right]\right)\\
 & =\hat{\mu}\left(0\right)+\sum_{n=1}^{N}\left(p^{n}\kappa\left(\left[\mathfrak{z}\right]_{p^{n}}\right)-\sum_{j=0}^{p-1}p^{n-1}\kappa\left(\left[\mathfrak{z}\right]_{p^{n-1}}+jp^{n-1}\right)\right)\\
 & =\hat{\mu}\left(0\right)+\underbrace{\sum_{n=1}^{N}\left(p^{n}\kappa\left(\left[\mathfrak{z}\right]_{p^{n}}\right)-p^{n-1}\kappa\left(\left[\mathfrak{z}\right]_{p^{n-1}}\right)\right)}_{\textrm{telescoping}}-\sum_{n=1}^{N}\sum_{j=1}^{p-1}p^{n-1}\kappa\left(\left[\mathfrak{z}\right]_{p^{n-1}}+jp^{n-1}\right)\\
 & =\hat{\mu}\left(0\right)-\underbrace{\kappa\left(0\right)}_{\hat{\mu}\left(0\right)}+p^{N}\kappa\left(\left[\mathfrak{z}\right]_{p^{N}}\right)-\sum_{n=0}^{N-1}\sum_{j=1}^{p-1}p^{n}\kappa\left(\left[\mathfrak{z}\right]_{p^{n}}+jp^{n}\right)
\end{align*}
as desired.

Finally, when $\kappa$ has $p$-adic structure, $\kappa$'s right-ended
structural equations let us write: 
\[
\sum_{n=0}^{N-1}\sum_{j=1}^{p-1}p^{n}\kappa\left(\left[\mathfrak{z}\right]_{p^{n}}+jp^{n}\right)=\left(\sum_{j=1}^{p-1}\kappa\left(j\right)\right)\sum_{n=0}^{N-1}p^{n}\kappa\left(\left[\mathfrak{z}\right]_{p^{n}}\right)
\]

Q.E.D. 
\begin{lem}[\textbf{Fourier Resummation \textendash{} Radially-Magnitudinal measures}]
\textbf{\label{lem:Radially-Mag Fourier Resummation Lemma}} Let
$\hat{\mu}:\hat{\mathbb{Z}}_{p}\rightarrow\overline{\mathbb{Q}}$
be a function which is not identically $0$, and suppose that $\hat{\mu}\left(t\right)=\hat{\nu}\left(t\right)\hat{\eta}\left(t\right)$,
where $\hat{\nu},\hat{\eta}:\hat{\mathbb{Z}}_{p}\rightarrow\overline{\mathbb{Q}}$
are radially symmetric and magnitude-dependent, respectively. Then,
for all $N\in\mathbb{N}_{1}$ and all $\mathfrak{z}\in\mathbb{Z}_{p}$:
\begin{align}
\sum_{\left|t\right|_{p}\leq p^{N}}\hat{\mu}\left(t\right)e^{2\pi i\left\{ t\mathfrak{z}\right\} _{p}} & \overset{\overline{\mathbb{Q}}}{=}\hat{\nu}\left(0\right)\kappa\left(0\right)+\sum_{n=1}^{N}\hat{\nu}\left(p^{-n}\right)\left(p^{n}\kappa\left(\left[\mathfrak{z}\right]_{p^{n}}\right)-p^{n-1}\kappa\left(\left[\mathfrak{z}\right]_{p^{n-1}}\right)\right)\label{eq:Radial-Magnitude Fourier Resummation Lemma}\\
 & -\sum_{j=1}^{p-1}\sum_{n=1}^{N}p^{n-1}\hat{\nu}\left(p^{-n}\right)\kappa\left(\left[\mathfrak{z}\right]_{p^{n-1}}+jp^{n-1}\right)\nonumber 
\end{align}
If, in addition, $\kappa$ has $p$-adic structure, then: 
\begin{align}
\sum_{\left|t\right|_{p}\leq p^{N}}\hat{\mu}\left(t\right)e^{2\pi i\left\{ t\mathfrak{z}\right\} _{p}} & \overset{\overline{\mathbb{Q}}}{=}p^{N}\hat{\nu}\left(p^{-N}\right)\kappa\left(\left[\mathfrak{z}\right]_{p^{N}}\right)\label{eq:Radial-Magnitude Fourier Resummation Lemma - p-adically distributed case}\\
 & +\sum_{n=0}^{N-1}\left(\hat{\nu}\left(p^{-n}\right)-\left(\sum_{j=0}^{p-1}\kappa\left(j\right)\right)\hat{\nu}\left(p^{-n-1}\right)\right)\kappa\left(\left[\mathfrak{z}\right]_{p^{n}}\right)p^{n}\nonumber 
\end{align}
\end{lem}
\begin{rem}
In order to make things more compact, we will frequently adopt the
notation: 
\begin{align}
A & \overset{\textrm{def}}{=}\sum_{j=1}^{p-1}\kappa\left(j\right)\label{eq:Definition of Big A}\\
v_{n} & \overset{\textrm{def}}{=}\hat{\nu}\left(\frac{1}{p^{n}}\right)\label{eq:Definition of v_n}
\end{align}
so that (\ref{eq:Radial-Magnitude Fourier Resummation Lemma - p-adically distributed case})
becomes: 
\begin{equation}
\sum_{\left|t\right|_{p}\leq p^{N}}\hat{\mu}\left(t\right)e^{2\pi i\left\{ t\mathfrak{z}\right\} _{p}}\overset{\overline{\mathbb{Q}}}{=}v_{N}\kappa\left(\left[\mathfrak{z}\right]_{p^{N}}\right)p^{N}+\sum_{n=0}^{N-1}\left(v_{n}-\left(1+A\right)v_{n+1}\right)\kappa\left(\left[\mathfrak{z}\right]_{p^{n}}\right)p^{n}\label{eq:Radial-Magnitude p-adic distributed Fourier Resummation Identity, simplified}
\end{equation}
\end{rem}
Proof:

\begin{align*}
\sum_{\left|t\right|_{p}\leq p^{N}}\hat{\mu}\left(t\right)e^{2\pi i\left\{ t\mathfrak{z}\right\} _{p}} & =\hat{\mu}\left(0\right)+\sum_{n=1}^{N}\sum_{\left|t\right|_{p}=p^{n}}\hat{\nu}\left(t\right)\hat{\eta}\left(t\right)e^{2\pi i\left\{ t\mathfrak{z}\right\} _{p}}\\
\left(\hat{\nu}\left(t\right)=\hat{\nu}\left(\frac{1}{\left|t\right|_{p}}\right)\right); & =\hat{\mu}\left(0\right)+\sum_{n=1}^{N}\hat{\nu}\left(p^{-n}\right)\sum_{\left|t\right|_{p}=p^{n}}\hat{\eta}\left(t\right)e^{2\pi i\left\{ t\mathfrak{z}\right\} _{p}}
\end{align*}
From (\ref{eq:Magnitude Fourier Resummation Lemma}), we have that:
\[
\sum_{\left|t\right|_{p}=p^{n}}\hat{\eta}\left(t\right)e^{2\pi i\left\{ t\mathfrak{z}\right\} _{p}}=p^{n}\kappa\left(\left[\mathfrak{z}\right]_{p^{n}}\right)-\sum_{j=0}^{p-1}p^{n-1}\kappa\left(\left[\mathfrak{z}\right]_{p^{n-1}}+jp^{n-1}\right)
\]
and hence: 
\begin{align*}
\sum_{\left|t\right|_{p}\leq p^{N}}\hat{\mu}\left(t\right)e^{2\pi i\left\{ t\mathfrak{z}\right\} _{p}} & =\hat{\mu}\left(0\right)+\sum_{n=1}^{N}\hat{\nu}\left(p^{-n}\right)\sum_{\left|t\right|_{p}=p^{n}}\hat{\eta}\left(t\right)e^{2\pi i\left\{ t\mathfrak{z}\right\} _{p}}\\
 & =\hat{\mu}\left(0\right)+\sum_{n=1}^{N}\hat{\nu}\left(p^{-n}\right)\kappa\left(\left[\mathfrak{z}\right]_{p^{n}}\right)p^{n}\\
 & -\sum_{j=0}^{p-1}\sum_{n=1}^{N}\hat{\nu}\left(p^{-n}\right)\kappa\left(\left[\mathfrak{z}\right]_{p^{n-1}}+jp^{n-1}\right)p^{n-1}
\end{align*}
If, in addition, $\hat{\mu}$ is $p$-adically distributed, then:
\begin{align*}
\sum_{\left|t\right|_{p}\leq p^{N}}\hat{\mu}\left(t\right)e^{2\pi i\left\{ t\mathfrak{z}\right\} _{p}} & =\hat{\mu}\left(0\right)+\sum_{n=1}^{N}\hat{\nu}\left(p^{-n}\right)\kappa\left(\left[\mathfrak{z}\right]_{p^{n}}\right)p^{n}\\
 & -\left(\sum_{j=0}^{p-1}\kappa\left(j\right)\right)\sum_{n=1}^{N}\hat{\nu}\left(p^{-n}\right)\kappa\left(\left[\mathfrak{z}\right]_{p^{n-1}}\right)p^{n-1}\\
\left(\hat{\mu}\left(0\right)=\hat{\nu}\left(0\right)\underbrace{\hat{\eta}\left(0\right)}_{\kappa\left(0\right)}\right); & =\underbrace{\hat{\nu}\left(0\right)\kappa\left(0\right)-\hat{\nu}\left(\frac{1}{p}\right)\kappa\left(0\right)\left(\sum_{j=0}^{p-1}\kappa\left(j\right)\right)}_{n=0\textrm{ term of the series on the bottom line}}+p^{N}\hat{\nu}\left(p^{-N}\right)\kappa\left(\left[\mathfrak{z}\right]_{p^{N}}\right)\\
 & +\sum_{n=1}^{N-1}\left(\hat{\nu}\left(p^{-n}\right)-\left(\sum_{j=0}^{p-1}\kappa\left(j\right)\right)\hat{\nu}\left(p^{-n-1}\right)\right)\kappa\left(\left[\mathfrak{z}\right]_{p^{n}}\right)p^{n}
\end{align*}

Q.E.D. 
\begin{prop}[\textbf{Convolution with $v_{p}$}]
\label{prop:v_p of t times mu hat sum}Let $\hat{\mu}:\hat{\mathbb{Z}}_{p}\rightarrow\mathbb{C}_{q}$
be any function. Then: 
\begin{equation}
\sum_{0<\left|t\right|_{p}\leq p^{N}}v_{p}\left(t\right)\hat{\mu}\left(t\right)e^{2\pi i\left\{ t\mathfrak{z}\right\} _{p}}\overset{\mathbb{C}_{q}}{=}-N\tilde{\mu}_{N}\left(\mathfrak{z}\right)+\sum_{n=0}^{N-1}\tilde{\mu}_{n}\left(\mathfrak{z}\right)\label{eq:Fourier sum of v_p times mu-hat}
\end{equation}
\end{prop}
Proof: 
\begin{align*}
\sum_{0<\left|t\right|_{p}\leq p^{N}}v_{p}\left(t\right)\hat{\mu}\left(t\right)e^{2\pi i\left\{ t\mathfrak{z}\right\} _{p}} & =\sum_{n=1}^{N}\sum_{\left|t\right|_{p}=p^{n}}\left(-n\right)\hat{\mu}\left(t\right)e^{2\pi i\left\{ t\mathfrak{z}\right\} _{p}}\\
 & =-\sum_{n=1}^{N}n\left(\sum_{\left|t\right|_{p}\leq p^{n}}\hat{\mu}\left(t\right)e^{2\pi i\left\{ t\mathfrak{z}\right\} _{p}}-\sum_{\left|t\right|_{p}\leq p^{n-1}}\hat{\mu}\left(t\right)e^{2\pi i\left\{ t\mathfrak{z}\right\} _{p}}\right)\\
 & =-\sum_{n=1}^{N}n\left(\tilde{\mu}_{n}\left(\mathfrak{z}\right)-\tilde{\mu}_{n-1}\left(\mathfrak{z}\right)\right)\\
 & =\sum_{n=1}^{N}n\tilde{\mu}_{n-1}\left(\mathfrak{z}\right)-\sum_{n=1}^{N}n\tilde{\mu}_{n}\left(\mathfrak{z}\right)\\
 & =\sum_{n=1}^{N}\left(\tilde{\mu}_{n-1}\left(\mathfrak{z}\right)+\left(n-1\right)\tilde{\mu}_{n-1}\left(\mathfrak{z}\right)\right)-\sum_{n=1}^{N}n\tilde{\mu}_{n}\left(\mathfrak{z}\right)\\
 & =\sum_{n=1}^{N}\tilde{\mu}_{n-1}\left(\mathfrak{z}\right)+\sum_{n=0}^{N-1}\underbrace{n\tilde{\mu}_{n}\left(\mathfrak{z}\right)}_{0\textrm{ when }n=0}-\sum_{n=1}^{N}n\tilde{\mu}_{n}\left(\mathfrak{z}\right)\\
 & =\sum_{n=0}^{N-1}\tilde{\mu}_{n}\left(\mathfrak{z}\right)-N\tilde{\mu}_{N}\left(\mathfrak{z}\right)
\end{align*}

Q.E.D.

\subsection{\label{subsec:3.3.3 Frames}Frames}
\begin{rem}
In order to be fully comprehensive, this sub-section takes a significant
dive into abstraction. This is meant primarily for readers with a
good background in non-archimedean analysis, as well as anyone who
can handle a significant number of definitions being thrown their
way. For readers more interested in how all of this relates to $\chi_{H}$
and Hydra maps, they can safely skip this section so long as the keep
the following concepts in mind:

Consider a function $f:\mathbb{Z}_{p}\rightarrow\mathbb{C}_{q}$.
A sequence of functions $\left\{ f_{n}\right\} _{n\geq1}$ on $\mathbb{Z}_{p}$
is said to \textbf{converge} \textbf{to $f$ with respect to the standard
$\left(p,q\right)$-adic frame }if\index{frame!standard left(p,qright)-adic@standard $\left(p,q\right)$-adic}:

\vphantom{}

I. For all $n$, $f_{n}\left(\mathfrak{z}\right)\in\overline{\mathbb{Q}}$
for all $\mathfrak{z}\in\mathbb{Z}_{p}$.

\vphantom{}

II. For all $\mathfrak{z}\in\mathbb{N}_{0}$, $f\left(\mathfrak{z}\right)\in\overline{\mathbb{Q}}$
and $\lim_{n\rightarrow\infty}f_{n}\left(\mathfrak{z}\right)\overset{\mathbb{C}}{=}f\left(\mathfrak{z}\right)$
(meaning the convergence is in the topology of $\mathbb{C}$).

\vphantom{}

III. For all $\mathfrak{z}\in\mathbb{Z}_{p}^{\prime}$, $f\left(\mathfrak{z}\right)\in\mathbb{C}_{q}$
and $\lim_{n\rightarrow\infty}f_{n}\left(\mathfrak{z}\right)\overset{\mathbb{C}_{q}}{=}f\left(\mathfrak{z}\right)$
(meaning the convergence is in the topology of $\mathbb{C}_{q}$).

\vphantom{}

We write $\mathcal{F}_{p,q}$ to denote the standard $\left(p,q\right)$-adic
frame and write $\lim_{n\rightarrow\infty}f_{n}\left(\mathfrak{z}\right)\overset{\mathcal{F}_{p,q}}{=}f\left(\mathfrak{z}\right)$
to denote convergence with respect to the standard frame. Also, if
the $f_{n}$s are the partial sums of an infinite series, we say that
$\lim_{n\rightarrow\infty}f_{n}\left(\mathfrak{z}\right)$ is an \textbf{$\mathcal{F}_{p,q}$-series
}for/of\index{mathcal{F}-@$\mathcal{F}$-!series} $f$.
\end{rem}
\vphantom{}

HERE BEGINS THE COMPREHENSIVE ACCOUNT OF FRAMES

\vphantom{}

To begin, let us recall the modern notion of what a ``function''
actually \emph{is}. Nowadays a ``function'' is merely a rule which
to each element of a designated input set $A$ associates a single
element of a designated output set, $B$. The first appearance of
this type of understanding is generally attributed to Dirichlet. Prior
to his work (and, honestly, for a good time after that), functions
were often synonymous with what we now call analytic functions\textemdash those
which are expressible as a power series, or some combination of functions
of that type. This synonymy is understandable: power series are wonderful.
Not only do they give us something concrete to manipulate, they even
allow us to compute functions' values to arbitrary accuracy\textemdash{}
provided that the series converges, of course.

For example, suppose we wish to use the geometric series formula:
\begin{equation}
\sum_{n=0}^{\infty}z^{n}=\frac{1}{1-z}
\end{equation}
to compute the value of the function $1/\left(1-z\right)$. As we
all know, for any real or complex number $z$ with $\left|z\right|<1$,
the series will converge in $\mathbb{R}$ or $\mathbb{C}$ to the
right hand side. However, \emph{what if we plug in $z=3/2$}? There,
the topologies of $\mathbb{R}$ and $\mathbb{C}$ do not help us.
While we \emph{could} use analytic continuation\textemdash rejiggering
the series until we get something which converges at $3/2$\textemdash that
particular method doesn't make much sense in non-archimedean analysis,
let alone $\left(p,q\right)$-adic analysis, so that route isn't available
to us. Instead, we can observe that our obstacle isn't the geometric
series formula itself, but rather the topologies of $\mathbb{R}$
and $\mathbb{C}$. As long as we stubbornly insist on living in those
archimedean universes and no others, we cannot sum the series as written
for any $\left|z\right|\geq1$. However, just like with speculative
fiction, there is no reason we must confine ourselves to a single
universe.

If we go through the magic portal to the world of the $3$-adic topology,
we can easily compute the true value of $1/\left(1-z\right)$ at $z=3/2$
using the power series formula. There: 
\begin{equation}
\sum_{n=0}^{\infty}\left(\frac{3}{2}\right)^{n}\overset{\mathbb{Z}_{3}}{=}\frac{1}{1-\frac{3}{2}}
\end{equation}
is a perfectly sensible statement, both ordinary and rigorously meaningful.
Most importantly, even though we have \emph{summed }the series in
an otherworldly topology, the answer ends up being an ordinary real
number, one which agrees with the value of $1/\left(1-z\right)$ at
$z=3/2$ in $\mathbb{R}$ or $\mathbb{C}$. More generally, for any
rational number $p/q$ where $p$ and $q$ are co-prime integers with
$q\neq0$, we can go about computing $1/\left(1-z\right)$ at $z=p/q$
by summing the series: 
\begin{equation}
\sum_{n=0}^{\infty}\left(\frac{p}{q}\right)^{n}
\end{equation}
in $\mathbb{Z}_{p}$. That this all works is thanks to the universality
of the geometric series universality \index{geometric series universality}(see
page \pageref{fact:Geometric series universality}). In this way,
we can get more mileage out of our series formulae. Even though the
geometric series converges in the topology of $\mathbb{C}$ only for
complex numbers $\left|z\right|<1$, we can make the formula work
at every complex algebraic number $\alpha$ of absolute value $>1$
by summing: 
\begin{equation}
\sum_{n=0}^{\infty}\alpha^{n}
\end{equation}
in the topology of the $\alpha$-adic completion of the field $\mathbb{Q\left(\alpha\right)}$.
This then gives us the correct value of the complex-valued function
$1/\left(1-z\right)$ at $z=\alpha$. The only unorthodoxy is the
route we took to get there. Nevertheless, even this is fully in line
with the modern notion of functions: the output a function assigns
to a given input should be defined \emph{independently} of any method
(if any) we use to explicitly compute it. Treating $1/\left(1-z\right)$
as a function from $\overline{\mathbb{Q}}$ to $\overline{\mathbb{Q}}$,
we can use the geometric series to compute the value of $1/\left(1-z\right)$
at any $\alpha\in\overline{\mathbb{Q}}$ by choosing a topological
universe (the space of $\alpha$-adic numbers) where the series formula
happens to make sense.

This discussion of $1/\left(1-z\right)$ contains in a nutshell all
the key concepts that come into play when defining frames. In this
dissertation, the $\left(p,q\right)$-adic functions like $\chi_{H}$
we desire to study \emph{will not }possess series representations
which converge in a single topology at every point of the function's
domain. However, if we allow ourselves to vary from point to point
the specific topology we use to sum these series, we can arrange things
so that the series is meaningful at every possible input value. That
this all \emph{works} is thanks to the observation that the values
we get by varying the topology in this way are the correct outputs
that the function assigns to its inputs, even if\textemdash as stated
above\textemdash the route we took to compute them happened to be
unusual. So, when working in a certain topology $\mathcal{T}$, instead
of admitting defeat and giving up hope when our series become non-convergent\textemdash or,
worse, \emph{divergent}\textemdash in $\mathcal{T}$'s topology, we
will instead hop over to a \emph{different} topological universe $\mathcal{T}^{\prime}$
in which the series \emph{is} convergent.

It almost goes without saying that world-hopping in this way is fraught
with pitfalls. As Hensel's famous blunder (see Subsection \ref{subsec:1.3.5 Hensel's-Infamous-Blunder})
perfectly illustrates, just because there exists \emph{some }topology
in which a series converges, it does not mean that the sum of the
series has meaning from the perspective of a different topology\textemdash the
one exception, of course, being the universality of the geometric
series. Though cumbersome, the abstractions presented in this subsection
are a necessary security measure. We need them to ensure that the
series and functions we work with will be compatible with the itinerary
of world-hopping we will use to make sense of them.

As alluded to in Subsection \ref{subsec:3.3.1 Heuristics-and-Motivations},
the general idea behind a quasi-integrable function $\chi$ on $\mathbb{Z}_{p}$
is that we can express $\chi$ as the limit of interest (\ref{eq:The Limit of Interest})
for some $\hat{\mu}:\hat{\mathbb{Z}}_{p}\rightarrow\overline{\mathbb{Q}}$.
Moreover, as demonstrated by $\hat{A}_{3}$ and the discussion of
the geometric series above, given a function $\chi:\mathbb{Z}_{p}\rightarrow\mathbb{K}$,
it will be too restrictive to mandate that the partial sums of our
Fourier series: 
\begin{equation}
\sum_{\left|t\right|_{p}\leq p^{N}}\hat{\chi}\left(t\right)e^{2\pi i\left\{ t\mathfrak{z}\right\} _{p}}
\end{equation}
converge point-wise on $\mathbb{Z}_{p}$ with respect to a \emph{single
}topology. Instead, we will only require that at every $\mathfrak{z}\in\mathbb{Z}_{p}$
for which the function $\chi\left(\mathfrak{z}\right)$ takes a finite
value, either the sequence of partial sums in the limit of interest
(\ref{eq:The Limit of Interest}) equals its limit for all sufficiently
large $N$ (as was the case in \textbf{Proposition \ref{prop:sum of v_p}}),
or, that there exists\emph{ }a valued field $K_{\mathfrak{z}}$ with
absolute value $\left|\cdot\right|_{K_{\mathfrak{z}}}$ for which:
\begin{equation}
\lim_{N\rightarrow\infty}\left|\chi\left(\mathfrak{z}\right)-\sum_{\left|t\right|_{p}\leq p^{N}}\hat{\chi}\left(t\right)e^{2\pi i\left\{ t\mathfrak{z}\right\} _{p}}\right|_{K_{\mathfrak{z}}}=0\label{eq:Quasi-integrability in a nutshell}
\end{equation}
With this procedure in mind, frames will be formulated as a collection
of pairs $\left(\mathfrak{z},K_{\mathfrak{z}}\right)$, where $\mathfrak{z}$
is a point in $\mathbb{Z}_{p}$ and $K_{\mathfrak{z}}$ is a field
so that (\ref{eq:The Limit of Interest}) converges at $\mathfrak{z}$
in the topology of $K_{\mathfrak{z}}$.

That being said, we still have to have \emph{some} limitations. In
applying $K_{\mathfrak{z}}$-absolute values to an equation like (\ref{eq:Quasi-integrability in a nutshell}),
it is implicitly required that all of the terms the sum are elements
of $K_{\mathfrak{z}}$. To surmount this, note that for any $t$ and
$\mathfrak{z}$, $e^{2\pi i\left\{ t\mathfrak{z}\right\} _{p}}$ is
some root of unity. Roots of unity exist abstractly in any algebraically
closed field. So, we can fit these complex exponentials inside all
of the $K_{\mathfrak{z}}$s by requiring the $K_{\mathfrak{z}}$s
to be \textbf{algebraically closed}.

There is also the issue of $\chi$ itself. Suppose that there are
$\mathfrak{a},\mathfrak{b}\in\mathbb{Z}_{p}$ so that $K_{\mathfrak{a}}=\mathbb{C}_{q}$
and $K_{\mathfrak{b}}=\mathbb{C}_{r}$, for distinct primes $p,q,r$?
Unless it just so happens that $\chi\left(\mathfrak{a}\right)$ and
$\chi\left(\mathfrak{b}\right)$ both lie in $\overline{\mathbb{Q}}$,
we will not be able to deal with all of $\chi$ by forcing $\chi$'s
co-domain to be a single field. So, we will need to be flexible. Just
as we allow for the topology of convergence to vary from point to
point, so too will we allow for the co-domain of our functions to
be more than just one field. This then suggests the following set
up: 
\begin{itemize}
\item A ``big'' space in which the limits of our (Fourier) series shall
live. While for us, this will end up always being in $\mathbb{C}_{q}$,
in general, it can be more than a single field\footnote{Originally, I required the ``big'' space to be a single field, but,
upon contemplating what would need to be done in order to make frames
applicable to the study of $\chi_{H}$ in the case of a polygenic
Hydra map, I realized I would need to allow for $\chi_{H}$ to take
values in more than one field.}. 
\item A choice of topologies of convergence $\left(\mathfrak{z},K_{\mathfrak{z}}\right)$
for every point in $\mathbb{Z}_{p}$. These shall be topological universes
we use to sum our Fourier series so as to arrive at elements of the
big space. 
\item We should require both the big space and the $K_{\mathfrak{z}}$s,
to contain $\overline{\mathbb{Q}}$, and will require the coefficients
of our Fourier series to be contained \emph{in }$\overline{\mathbb{Q}}$.
This way, we will be able to compute partial sums of Fourier series
in all of the $K_{\mathfrak{z}}$s simultaneously by working in $\overline{\mathbb{Q}}$,
the space they share in common. 
\end{itemize}
To make this work, we will first start with a choice of the pairs
$\left(\mathfrak{z},K_{\mathfrak{z}}\right)$. This will get us a
``$p$-adic frame'', denoted $\mathcal{F}$. We will then define
$I\left(\mathcal{F}\right)$\textemdash the \textbf{image }of $\mathcal{F}$\textemdash to
denote the big space; this will just be the union of all the $K_{\mathfrak{z}}$s.
Lastly, we will have $C\left(\mathcal{F}\right)$, the $\overline{\mathbb{Q}}$-linear
space of all the functions which are compatible with this set up;
the archetypical elements of these spaces will be partial sums of
Fourier series generated by functions $\hat{\mathbb{Z}}_{p}\rightarrow\overline{\mathbb{Q}}$. 
\begin{defn}[\textbf{Frames}]
\label{def:Frames}A \textbf{$p$-adic quasi-integrability frame
}\index{frame}(or just ``frame'', for short) $\mathcal{F}$ consists
of the following pieces of data:

\vphantom{}

I. A non-empty set \nomenclature{$U_{\mathcal{F}}$}{domain of $\mathcal{F}$}$U_{\mathcal{F}}\subseteq\mathbb{Z}_{p}$,
called the \textbf{domain }of $\mathcal{F}$.\index{frame!domain}

\vphantom{}

II. For every $\mathfrak{z}\in U_{\mathcal{F}}$, a topological field
$K_{\mathfrak{z}}$ \nomenclature{$K_{\mathfrak{z}}$}{ } with $K_{\mathfrak{z}}\supseteq\overline{\mathbb{Q}}$.
We allow for $K_{\mathfrak{z}}=\overline{\mathbb{Q}}$. Additionally,
for any $\mathfrak{z}\in U_{\mathcal{F}}$, we require the topology
of $K_{\mathfrak{z}}$ to be either the discrete topology\emph{ or}
the topology induced by an absolute value, denoted $\left|\cdot\right|_{K_{\mathfrak{z}}}$\nomenclature{$\left|\cdot\right|_{K_{\mathfrak{z}}}$}{ },
so that $\left(K_{\mathfrak{z}},\left|\cdot\right|_{\mathfrak{z}}\right)$
is then a metrically complete valued field.
\end{defn}
\vphantom{}

The two most important objects associated with a given frame are its
image and the space of compatible functions:
\begin{defn}[\textbf{Image and Compatible Functions}]
\label{def:Frame terminology}Let $\mathcal{F}$ be a $p$-adic frame.

\vphantom{}

I. The \textbf{image }of $\mathcal{F}$, denoted $I\left(\mathcal{F}\right)$,
is the \emph{set }defined by: 
\begin{equation}
I\left(\mathcal{F}\right)\overset{\textrm{def}}{=}\bigcup_{\mathfrak{z}\in U_{\mathcal{F}}}K_{\mathfrak{z}}\label{eq:The image of a frame}
\end{equation}
\nomenclature{$I\left(\mathcal{F}\right)$}{The image of $\mathcal{F}$}\index{frame!image}

\vphantom{}

II. A function $\chi:U_{\mathcal{F}}\rightarrow I\left(\mathcal{F}\right)$
is said to be \textbf{$\mathcal{F}$-compatible} / \textbf{compatible
}(\textbf{with $\mathcal{F}$}) whenever \index{mathcal{F}-@$\mathcal{F}$-!compatible}\index{frame!compatible functions}
$\chi\left(\mathfrak{z}\right)\in K_{\mathfrak{z}}$ for all $\mathfrak{z}\in U_{\mathcal{F}}$.
I write $C\left(\mathcal{F}\right)$ \nomenclature{$C\left(\mathcal{F}\right)$}{set of $\mathcal{F}$-compatible functions}
to denote the set of all $\mathcal{F}$-compatible functions. 
\end{defn}
\begin{rem}
Note that $C\left(\mathcal{F}\right)$ is a linear space over $\overline{\mathbb{Q}}$
with respect to point-wise addition of functions and scalar multiplication. 
\end{rem}
\vphantom{}

Next, we have some terminology which is of use when working with a
frame. 
\begin{defn}[\textbf{Frame Terminology}]
Let $\mathcal{F}$ be a $p$-adic frame.

\vphantom{}

I. I call $\mathbb{Z}_{p}\backslash U_{\mathcal{F}}$ the \textbf{set
of singularities }of the frame. I say that $\mathcal{F}$ is \textbf{non-singular
}whenever $U_{\mathcal{F}}=\mathbb{Z}_{p}$.\index{frame!non-singular}

\vphantom{}

II. I write $\textrm{dis}$ to denote the discrete topology. I write
$\infty$ to denote the topology of $\mathbb{C}$ (induced by $\left|\cdot\right|$),
and I write $p$ to denote the $p$-adic topology (that of $\mathbb{C}_{p}$,
induced by $\left|\cdot\right|_{p}$). I write $\textrm{non}$ to
denote an arbitrary non-archimedean topology (induced by some non-archimedean
absolute value on $\overline{\mathbb{Q}}$).

\vphantom{}

III. Given a topology $\tau$ (so, $\tau$ could be $\textrm{dis}$,
$\infty$, $\textrm{non}$, or $p$), I write\footnote{So, we can write $U_{\textrm{dis}}\left(\mathcal{F}\right)$, $U_{\infty}\left(\mathcal{F}\right)$,
$U_{\textrm{non}}\left(\mathcal{F}\right)$, or $U_{p}\left(\mathcal{F}\right)$.} $U_{\tau}\left(\mathcal{F}\right)$ to denote the set of $\mathfrak{z}\in U_{\mathcal{F}}$
so that $K_{\mathfrak{z}}$ has been equipped with $\tau$. I call
the $U_{\tau}\left(\mathcal{F}\right)$\textbf{ $\tau$-convergence
domain }or \textbf{domain of $\tau$ convergence }of $\mathcal{F}$.
In this way, we can speak of the \textbf{discrete convergence domain},
the \textbf{archimedean convergence domain}, the \textbf{non-archimedean
convergence domain}, and the \textbf{$p$-adic convergence domain}
of a given frame $\mathcal{F}$.

\vphantom{}

IV. I say $\mathcal{F}$ is \textbf{simple} \index{frame!simple}if
either $U_{\textrm{non}}\left(\mathcal{F}\right)$ is empty or there
exists a single metrically complete non-archimedean field extension
$K$ of $\overline{\mathbb{Q}}$ so that $K=K_{\mathfrak{z}}$ for
all $\mathfrak{z}\in U_{\textrm{non}}\left(\mathcal{F}\right)$. (That
is to say, for a simple frame, we use at most one non-archimedean
topology.)

\vphantom{}

V. I say $\mathcal{F}$ is \textbf{proper }\index{frame!proper}proper
whenever $K_{\mathfrak{z}}$ is not a $p$-adic field for any $\mathfrak{z}\in U_{\textrm{non}}\left(\mathcal{F}\right)$. 
\end{defn}
\vphantom{}

UNLESS STATED OTHERWISE, ALL FRAMES ARE ASSUMED TO BE PROPER.

\vphantom{}

Only one frame will ever be used in this dissertation. It is described
below. 
\begin{defn}
The \textbf{standard ($\left(p,q\right)$-adic) frame}\index{frame!standard left(p,qright)-adic@standard $\left(p,q\right)$-adic},
denoted $\mathcal{F}_{p,q}$, is the frame for which the topology
of $\mathbb{C}$ is associated to $\mathbb{N}_{0}$ and the topology
of $\mathbb{C}_{q}$ is associated to $\mathbb{Z}_{p}^{\prime}$. 
\end{defn}
\vphantom{}

With frames, we can avoid circumlocution when discussing topologies
of convergence for $\left(p,q\right)$-adic functions and sequences
thereof. 
\begin{defn}[\textbf{$\mathcal{F}$-convergence}]
Given a frame a $\mathcal{F}$ and a function $\chi\in C\left(\mathcal{F}\right)$,
we say a sequence $\left\{ \chi_{n}\right\} _{n\geq1}$ in $C\left(\mathcal{F}\right)$
\textbf{converges to $\chi$ over $\mathcal{F}$ }(or \textbf{is $\mathcal{F}$-convergent
to $\chi$}) whenever, for each $\mathfrak{z}\in U_{\mathcal{F}}$,
we have: 
\begin{equation}
\lim_{n\rightarrow\infty}\chi_{n}\left(\mathfrak{z}\right)\overset{K_{\mathfrak{z}}}{=}\chi\left(\mathfrak{z}\right)
\end{equation}
Note that if $K_{\mathfrak{z}}$ has the discrete topology, the limit
implies that $\chi_{n}\left(\mathfrak{z}\right)=\chi\left(\mathfrak{z}\right)$
for all sufficiently large $n$.\index{mathcal{F}-@$\mathcal{F}$-!convergence}\index{frame!convergence}

The convergence is point-wise with respect to $\mathfrak{z}$. We
then call $\chi$ the \index{mathcal{F}-@$\mathcal{F}$-!limit}\textbf{$\mathcal{F}$-limit
of the $\chi_{n}$s}. \index{frame!limit}More generally, we say $\left\{ \chi_{n}\right\} _{n\geq1}$
is \textbf{$\mathcal{F}$-convergent / converges over $\mathcal{F}$}
whenever there is a function $\chi\in C\left(\mathcal{F}\right)$
so that the $\chi_{n}$s are $\mathcal{F}$-convergent to $\chi$.
In symbols, we denote $\mathcal{F}$-convergence by: 
\begin{equation}
\lim_{n\rightarrow\infty}\chi_{n}\left(\mathfrak{z}\right)\overset{\mathcal{F}}{=}\chi\left(\mathfrak{z}\right),\textrm{ }\forall\mathfrak{z}\in U_{\mathcal{F}}\label{eq:Definition-in-symbols of F-convergence}
\end{equation}
or simply: 
\begin{equation}
\lim_{n\rightarrow\infty}\chi_{n}\overset{\mathcal{F}}{=}\chi\label{eq:Simplified version of expressing f as the F-limit of f_ns}
\end{equation}
\end{defn}
\begin{rem}
With respect to the standard $\left(p,q\right)$-adic frame, $\mathcal{F}_{p,q}$-convergence
means $\chi_{n}\left(\mathfrak{z}\right)$ converges to $\chi\left(\mathfrak{z}\right)$
in the topology of $\mathbb{C}_{q}$ for all $\mathfrak{z}\in\mathbb{Z}_{p}^{\prime}$
and converges to $\chi\left(\mathfrak{z}\right)$ in the topology
of $\mathbb{C}$ for all $\mathfrak{z}\in\mathbb{N}_{0}$. 
\end{rem}
\begin{rem}
In an abuse of notation, we will sometimes write $\mathcal{F}$ convergence
as: 
\begin{equation}
\lim_{n\rightarrow\infty}\left|\chi_{n}\left(\mathfrak{z}\right)-\chi\left(\mathfrak{z}\right)\right|_{K_{\mathfrak{z}}},\textrm{ }\forall\mathfrak{z}\in U_{\mathcal{F}}
\end{equation}
The abuse here is that the absolute value $\left|\cdot\right|_{K_{\mathfrak{z}}}$
will not exist if $\mathfrak{z}\in U_{\textrm{dis}}\left(\mathcal{F}\right)$.
As such, for any $\mathfrak{z}\in U_{\textrm{dis}}\left(\mathcal{F}\right)$,
we define: 
\begin{equation}
\lim_{n\rightarrow\infty}\left|\chi_{n}\left(\mathfrak{z}\right)-\chi\left(\mathfrak{z}\right)\right|_{K_{\mathfrak{z}}}=0
\end{equation}
to mean that, for all sufficiently large $n$, the equality $\chi_{n}\left(\mathfrak{z}\right)=\chi_{n+1}\left(\mathfrak{z}\right)$
holds in $K_{\mathfrak{z}}$. 
\end{rem}
\begin{prop}
Let $\mathcal{F}$ be a $p$-adic frame, and let $\hat{\mu}:\hat{\mathbb{Z}}_{p}\rightarrow\overline{\mathbb{Q}}$.
Then, $\tilde{\mu}_{N}\left(\mathfrak{z}\right)$\textemdash the $N$th
partial sum of the Fourier series generated by $\hat{\mu}$\textemdash is
an element of $C\left(\mathcal{F}\right)$. 
\end{prop}
Proof: Because $\hat{\mu}$ takes values in $\overline{\mathbb{Q}}$,
$\tilde{\mu}_{N}$ is a function from $\mathbb{Z}_{p}$ to $\overline{\mathbb{Q}}$.
As such, restricting $\tilde{\mu}_{N}$ to $U_{\mathcal{F}}$ then
gives us a function $U_{\mathcal{F}}\rightarrow I\left(\mathcal{F}\right)$,
because $\overline{\mathbb{Q}}\subseteq I\left(\mathcal{F}\right)$.

Q.E.D.

\vphantom{}

Now that we have the language of frames at our disposal, we can begin
to define classes of $\left(p,q\right)$-adic measures which will
be of interest to us. 
\begin{defn}
For any prime $q\neq p$, we say $d\mu\in C\left(\mathbb{Z}_{p},\mathbb{C}_{q}\right)^{\prime}$
is\index{measure!algebraic} \textbf{algebraic},\textbf{ }if its Fourier-Stieltjes
transform $\hat{\mu}$ takes values in $\overline{\mathbb{Q}}$. 
\end{defn}
\begin{defn}[\textbf{$\mathcal{F}$-measures}]
Given a $p$-adic frame $\mathcal{F}$, we write \nomenclature{$M\left(\mathcal{F}\right)$}{set of $\mathcal{F}$-measures}$M\left(\mathcal{F}\right)$
to denote the set of all functions $\hat{\mu}:\hat{\mathbb{Z}}_{p}\rightarrow\overline{\mathbb{Q}}$
so that $\hat{\mu}\in B\left(\hat{\mathbb{Z}}_{p},K_{\mathfrak{z}}\right)$
for all $\mathfrak{z}\in U_{\textrm{non}}\left(\mathcal{F}\right)$;
that is: 
\begin{equation}
\left\Vert \hat{\mu}\right\Vert _{p,K_{\mathfrak{z}}}<\infty,\textrm{ }\forall\mathfrak{z}\in U_{\textrm{non}}\left(\mathcal{F}\right)\label{eq:Definition of an F-measure}
\end{equation}
where, recall $\left\Vert \hat{\mu}\right\Vert _{p,K_{\mathfrak{z}}}$
denotes $\sup_{t\in\hat{\mathbb{Z}}_{p}}\left|\hat{\mu}\left(t\right)\right|_{K_{\mathfrak{z}}}$.

Next, observe that for every $\mathfrak{z}\in U_{\textrm{non}}\left(\mathcal{F}\right)$,
the map: 
\begin{equation}
f\in C\left(\mathbb{Z}_{p},K_{\mathfrak{z}}\right)\mapsto\sum_{t\in\hat{\mathbb{Z}}_{p}}\hat{f}\left(t\right)\hat{\mu}\left(-t\right)\in K_{\mathfrak{z}}\label{eq:Definition of the action of an F-measure}
\end{equation}
then defines an element of $C\left(\mathbb{Z}_{p},K_{\mathfrak{z}}\right)^{\prime}$.
As such, we will identify $M\left(\mathcal{F}\right)$ with elements
of: 
\begin{equation}
\bigcap_{\mathfrak{z}\in U_{\textrm{non}}\left(\mathcal{F}\right)}C\left(\mathbb{Z}_{p},K_{\mathfrak{z}}\right)^{\prime}
\end{equation}
and refer to elements of $M\left(\mathcal{F}\right)$ as \textbf{$\mathcal{F}$-measures}\index{mathcal{F}-@$\mathcal{F}$-!measure},
writing them as $d\mu$. Thus, the statement ``$d\mu\in M\left(\mathcal{F}\right)$''
means that $d\mu$ is an algebraic measure on $C\left(\mathbb{Z}_{p},K_{\mathfrak{z}}\right)$
for all $\mathfrak{z}\in U_{\textrm{non}}\left(\mathcal{F}\right)$,
which necessarily forces $\hat{\mu}$ to be bounded on $\hat{\mathbb{Z}}_{p}$
in $K_{\mathfrak{z}}$-absolute-value for all $\mathfrak{z}\in U_{\textrm{non}}\left(\mathcal{F}\right)$. 
\end{defn}
\begin{rem}
If $\mathcal{F}$ is simple and there is a prime $q\neq p$ so that
$K_{\mathfrak{z}}=\mathbb{C}_{q}$ for all $\mathfrak{z}\in U_{\textrm{non}}\left(\mathcal{F}\right)$,
observe that $M\left(\mathcal{F}\right)$ is then precisely the set
of all algebraic $\left(p,q\right)$-adic measures. This is also the
case for the standard $\left(p,q\right)$-adic frame. 
\end{rem}
\begin{rem}
$M\left(\mathcal{F}\right)$ is a linear space over $\overline{\mathbb{Q}}$. 
\end{rem}
\vphantom{}

Next, we introduce definitions to link frames with the limits of Fourier
series, by way of the convolution of a measure with the $p$-adic
Dirichlet Kernels. 
\begin{defn}[\textbf{$\mathcal{F}$-rising measure}]
Given a $p$-adic frame $\mathcal{F}$, we say $d\mu\in M\left(\mathcal{F}\right)$
\textbf{rises over $\mathcal{F}$} / is \textbf{$\mathcal{F}$-rising}
whenever\index{measure!mathcal{F}-rising@$\mathcal{F}$-rising} \index{frame!rising measure}the
sequence $\left\{ \tilde{\mu}_{N}\right\} _{N\geq0}=\left\{ D_{p:N}*d\mu\right\} _{N\geq0}$
is $\mathcal{F}$-convergent: 
\begin{equation}
\lim_{N\rightarrow\infty}\tilde{\mu}_{N}\left(\mathfrak{z}\right)\textrm{ exists in }K_{\mathfrak{z}},\textrm{ }\forall\mathfrak{z}\in U_{\mathcal{F}}\label{eq:Definition of an F-rising measure}
\end{equation}
In the case of the standard $\left(p,q\right)$-adic frame, this means:

\vphantom{}

I. For every $\mathfrak{z}\in\mathbb{Z}_{p}^{\prime}$, $\lim_{N\rightarrow\infty}\tilde{\mu}_{N}\left(\mathfrak{z}\right)$
converges in $\mathbb{C}_{q}$.

\vphantom{}

II. For every $\mathfrak{z}\in\mathbb{N}_{0}$, $\lim_{N\rightarrow\infty}\tilde{\mu}_{N}\left(\mathfrak{z}\right)$
converges in $\mathbb{C}$.

\vphantom{}

In all cases, the convergence is point-wise. We then write $M_{\textrm{rise}}\left(\mathcal{F}\right)$\nomenclature{$M_{\textrm{rise}}\left(\mathcal{F}\right)$}{the set of $\mathcal{F}$-rising measures}
to denote the set of all $\mathcal{F}$-rising measures. Note that
$M_{\textrm{rise}}\left(\mathcal{F}\right)$ forms a linear space
over $\overline{\mathbb{Q}}$.

Given a $p$-adic frame $\mathcal{F}$ and an $\mathcal{F}$-rising
measure $d\mu\in M_{\textrm{rise}}\left(\mathcal{F}\right)$, the
\textbf{$\mathcal{F}$-derivative }of \index{measure!mathcal{F}-derivative@$\mathcal{F}$-derivative}
$d\mu$ (or simply \textbf{derivative }of $d\mu$), denoted $\tilde{\mu}$,
is the function on $U_{\mathcal{F}}$ defined by the $\mathcal{F}$-limit
of the $\tilde{\mu}_{N}$s: 
\begin{equation}
\tilde{\mu}\left(\mathfrak{z}\right)\overset{\mathcal{F}}{=}\lim_{N\rightarrow\infty}\tilde{\mu}_{N}\left(\mathfrak{z}\right)\label{eq:Definition of mu-twiddle / the derivative of mu}
\end{equation}
That is, for each $\mathfrak{z}\in U_{\mathcal{F}}$, $\tilde{\mu}\left(\mathfrak{z}\right)$
is defined by the limit of the $\tilde{\mu}_{N}\left(\mathfrak{z}\right)$s
in $K_{\mathfrak{z}}$.

When $\mathcal{F}$ is the standard $\left(p,q\right)$-adic frame,
we have that:

\begin{equation}
\tilde{\mu}\left(\mathfrak{z}\right)=\begin{cases}
\lim_{N\rightarrow\infty}\tilde{\mu}_{N}\left(\mathfrak{z}\right)\textrm{ in }\mathbb{C}_{q} & \textrm{if }\mathfrak{z}\in\mathbb{Z}_{p}^{\prime}\\
\lim_{N\rightarrow\infty}\tilde{\mu}_{N}\left(\mathfrak{z}\right)\textrm{ in }\mathbb{C} & \textrm{if }\mathfrak{z}\in\mathbb{N}_{0}
\end{cases}\label{eq:F-limit in the standard frame}
\end{equation}
We then write $\left(\tilde{\mu}\right)_{N}\left(\mathfrak{z}\right)$
to denote the $N$th truncation of the derivative of $d\mu$: 
\begin{equation}
\left(\tilde{\mu}\right)_{N}\left(\mathfrak{z}\right)\overset{\textrm{def}}{=}\sum_{n=0}^{p^{N}-1}\tilde{\mu}\left(n\right)\left[\mathfrak{z}\overset{p^{N}}{\equiv}n\right]\label{eq:Definition of the Nth truncation of the derivative of a rising measure}
\end{equation}
More generally, we have: 
\begin{equation}
\left(\tilde{\mu}_{M}\right)_{N}\left(\mathfrak{z}\right)\overset{\textrm{def}}{=}\sum_{n=0}^{p^{N}-1}\tilde{\mu}_{M}\left(n\right)\left[\mathfrak{z}\overset{p^{N}}{\equiv}n\right]\label{eq:Nth truncation of the Mth partial sum of mu's derivative}
\end{equation}
is the $N$th truncation of $\tilde{\mu}_{M}$. 
\end{defn}
\begin{defn}
We say $d\mu\in\mathcal{M}_{\textrm{rise}}\left(\mathcal{F}_{p,q}\right)$
is \textbf{rising-continuous }if\index{measure!rising-continuous}
it\index{rising-continuous!measure} whenever its $\mathcal{F}_{p,q}$-derivative
$\tilde{\mu}\left(\mathfrak{z}\right)$ is a rising-continuous function. 
\end{defn}
\vphantom{}

Next up: \emph{degenerate }measures.
\begin{defn}[\textbf{$\mathcal{F}$-degenerate measure}]
\label{def:degenerate measure}Given a $p$-adic frame $\mathcal{F}$,
an $\mathcal{F}$-rising measure $d\mu\in M_{\textrm{rise}}\left(\mathcal{F}\right)$
is said to be\textbf{ degenerate }whenever\index{measure!degenerate}
its\index{mathcal{F}-@$\mathcal{F}$-!degenerate measure} $\mathcal{F}$-derivative\index{frame!degenerate measure}
is identically zero. We drop the modifier $\mathcal{F}$ when there
is no confusion.

We write $M_{\textrm{dgen}}\left(\mathcal{F}\right)$\nomenclature{$M_{\textrm{dgen}}\left(\mathcal{F}\right)$}{the set of $\mathcal{F}$-degenerate measures}
to denote the set of all $\mathcal{F}$-degenerate measures. Note
that $M_{\textrm{dgen}}\left(\mathcal{F}\right)$ is a linear subspace
of $M_{\textrm{rise}}\left(\mathcal{F}\right)$. 
\end{defn}
\begin{defn}
Given a $p$-adic frame $\mathcal{F}$, an $\mathcal{F}$-rising measure
$d\mu\in M_{\textrm{rise}}\left(\mathcal{F}\right)$ is\index{measure!mathcal{F}-proper@$\mathcal{F}$-proper}
said\index{mathcal{F}-@$\mathcal{F}$-!proper measure} to be \textbf{proper
/ $\mathcal{F}$-proper }whenever: 
\begin{equation}
\lim_{N\rightarrow\infty}\left(\tilde{\mu}_{N}-\left(\tilde{\mu}\right)_{N}\right)\overset{\mathcal{F}}{=}0\label{eq:Definition of a proper measure}
\end{equation}
If this property fails to hold, we call $d\mu$ \textbf{improper /
$\mathcal{F}$-improper}.\textbf{ }We then write $M_{\textrm{prop}}\left(\mathcal{F}\right)$
to denote the set of all $\mathcal{F}$-proper measures. Note that
$M_{\textrm{prop}}\left(\mathcal{F}\right)$ is a linear subspace
of $M_{\textrm{rise}}\left(\mathcal{F}\right)$. 
\end{defn}
\begin{question}
\label{que:3.2}Is the notion of a non-proper measure a vacuous concept,
or does there exist a frame $\mathcal{F}$ and an $\mathcal{F}$-rising
measure $d\mu$ which is improper? 
\end{question}
\vphantom{}

Before we move on to the next subsection, we record one last result: 
\begin{prop}
Let $d\mu$ be a rising-continuous $\left(p,q\right)$-adic measure.
Then, the derivative $\tilde{\mu}$ is a rising-continuous function
if and only if:

\vphantom{}

I. $d\mu$ is proper.

\vphantom{}

II. $S_{p}\left\{ \tilde{\mu}\right\} \left(\mathfrak{z}\right)$
converges in $\mathbb{C}_{q}$ for every $\mathfrak{z}\in\mathbb{Z}_{p}$. 
\end{prop}
Proof: We begin by defining a function $\mu^{\prime}$ by: 
\begin{equation}
\mu^{\prime}\left(\mathfrak{z}\right)\overset{\textrm{def}}{=}S_{p}\left\{ \tilde{\mu}\right\} \left(\mathfrak{z}\right)\overset{\textrm{def}}{=}\sum_{n=0}^{\infty}c_{n}\left(\tilde{\mu}\right)\left[\mathfrak{z}\overset{p^{\lambda_{p}\left(n\right)}}{\equiv}n\right]\label{eq:Chi prime}
\end{equation}
That is to say, $\mu^{\prime}$ is obtained by computing the values
$\tilde{\mu}\left(m\right)$ via the limits of $\tilde{\mu}_{N}\left(m\right)$,
and then using those values to construct a van der Put series. On
the other hand, $\tilde{\mu}\left(\mathfrak{z}\right)$ is constructed
by computing the $\mathcal{F}_{p,q}$-limits of $\tilde{\mu}_{N}\left(\mathfrak{z}\right)$
for each $\mathfrak{z}\in\mathbb{Z}_{p}$. To complete the proof,
we just need to show that $\tilde{\mu}$ and $\mu^{\prime}$ are the
same function. In fact, we need only show that they are equal to one
another over $\mathbb{Z}_{p}^{\prime}$; that they are equal over
$\mathbb{N}_{0}$ is because $\mu^{\prime}$ is the function represented
by the van der Put series of $\tilde{\mu}$, and a van der Put series
of a function \emph{always} converges on $\mathbb{N}_{0}$ to said
function.

To do this, we write: 
\begin{equation}
\tilde{\mu}\left(\mathfrak{z}\right)-\mu^{\prime}\left(\mathfrak{z}\right)=\tilde{\mu}\left(\mathfrak{z}\right)-\tilde{\mu}_{N}\left(\mathfrak{z}\right)+\tilde{\mu}_{N}\left(\mathfrak{z}\right)-\left(\tilde{\mu}\right)_{N}\left(\mathfrak{z}\right)+\left(\tilde{\mu}\right)_{N}\left(\mathfrak{z}\right)-\mu^{\prime}\left(\mathfrak{z}\right)
\end{equation}
Since: 
\begin{equation}
S_{p:N}\left\{ \tilde{\mu}\right\} \left(\mathfrak{z}\right)=\tilde{\mu}\left(\left[\mathfrak{z}\right]_{p^{N}}\right)\overset{\textrm{def}}{=}\left(\tilde{\mu}\right)_{N}\left(\mathfrak{z}\right)
\end{equation}
we then have: 
\begin{equation}
\tilde{\mu}\left(\mathfrak{z}\right)-\mu^{\prime}\left(\mathfrak{z}\right)=\underbrace{\tilde{\mu}\left(\mathfrak{z}\right)-\tilde{\mu}_{N}\left(\mathfrak{z}\right)}_{\textrm{A}}+\underbrace{\tilde{\mu}_{N}\left(\mathfrak{z}\right)-\left(\tilde{\mu}\right)_{N}\left(\mathfrak{z}\right)}_{\textrm{B}}+\underbrace{S_{p:N}\left\{ \tilde{\mu}\right\} \left(\mathfrak{z}\right)-\mu^{\prime}\left(\mathfrak{z}\right)}_{\textrm{C}}\label{eq:A,B,C,}
\end{equation}

i. Suppose $d\mu$ is proper and that $S_{p}\left\{ \tilde{\mu}\right\} \left(\mathfrak{z}\right)\overset{\mathbb{C}_{q}}{=}\lim_{N\rightarrow\infty}S_{p:N}\left\{ \tilde{\mu}\right\} \left(\mathfrak{z}\right)$
converges for every $\mathfrak{z}\in\mathbb{Z}_{p}$. Since $d\mu$
is rising-continuous, $\lim_{N\rightarrow\infty}\left|\tilde{\mu}\left(\mathfrak{z}\right)-\tilde{\mu}_{N}\left(\mathfrak{z}\right)\right|_{q}=0$
for all $\mathfrak{z}\in\mathbb{Z}_{p}^{\prime}$, which proves that
(A) tends to zero $q$-adically. Additionally, because $d\mu$ is
proper, the fact that $d\mu$ is proper then guarantees that $\lim_{N\rightarrow\infty}\left|\tilde{\mu}_{N}\left(\mathfrak{z}\right)-\left(\tilde{\mu}\right)_{N}\left(\mathfrak{z}\right)\right|_{q}=0$
for all $\mathfrak{z}\in\mathbb{Z}_{p}^{\prime}$, seeing as $\mathbb{Z}_{p}^{\prime}$
is the domain of non-archimedean convergence for the standard $\left(p,q\right)$-adic
frame. This shows that (B) tends to zero $q$-adically.

Finally, (C) tends to zero $q$-adically as $N\rightarrow\infty$
for all $\mathfrak{z}\in\mathbb{Z}_{p}^{\prime}$ because firstly,
by definition, $\mu^{\prime}$ is the limit of $S_{p:N}\left\{ \tilde{\mu}\right\} \left(\mathfrak{z}\right)$
as $N\rightarrow\infty$, and, secondly, because $S_{p:N}\left\{ \tilde{\mu}\right\} \left(\mathfrak{z}\right)$
converges to a limit $q$-adically for every $\mathfrak{z}\in\mathbb{Z}_{p}$.
This shows that the right hand side of (\ref{eq:A,B,C,}) is identically
zero for all $\mathfrak{z}\in\mathbb{Z}_{p}^{\prime}$. So, $\tilde{\mu}=\mu^{\prime}$
on $\mathbb{Z}_{p}^{\prime}$, as desired.

\vphantom{}

ii. Suppose $\tilde{\mu}$ is rising-continuous. Then, $\mu^{\prime}$\textemdash the
van der Put series of $\tilde{\mu}$\textemdash is necessarily equal
to $\tilde{\mu}$, and converges everywhere, thereby forcing (II)
to hold true.

Next, by the van der Put identity (\textbf{Proposition \ref{prop:vdP identity}}),
we can write: 
\begin{equation}
\tilde{\mu}\left(\mathfrak{z}\right)=\mu^{\prime}\left(\mathfrak{z}\right)\overset{\mathbb{C}_{q}}{=}\lim_{N\rightarrow\infty}\left(\tilde{\mu}\right)_{N}\left(\mathfrak{z}\right),\textrm{ }\forall\mathfrak{z}\in\mathbb{Z}_{p}
\end{equation}
Moreover, $\left(\tilde{\mu}\right)_{N}\left(\mathfrak{z}\right)=\left(\tilde{\mu}\right)_{N+1}\left(\mathfrak{z}\right)$
for all sufficiently large $N$ whenever $\mathfrak{z}\in\mathbb{N}_{0}$,
so $\left(\tilde{\mu}\right)_{N}\left(\mathfrak{z}\right)$ converges
to $\tilde{\mu}\left(\mathfrak{z}\right)$ in $\mathbb{C}$ for $\mathfrak{z}\in\mathbb{N}_{0}$
and in $\mathbb{C}_{q}$ for $\mathfrak{z}\in\mathbb{Z}_{p}^{\prime}$.
This shows that $\left(\tilde{\mu}\right)_{N}$ is $\mathcal{F}_{p,q}$-convergent
to $\tilde{\mu}$. Since, $d\mu$ is rising-continuous, by definition,
the $\tilde{\mu}_{N}$s are $\mathcal{F}_{p,q}$-convergent to $\tilde{\mu}$.
Since $\left(\tilde{\mu}\right)_{N}$ and $\tilde{\mu}_{N}$ then
$\mathcal{F}_{p,q}$-converge to the same function, they then necessarily
$\mathcal{F}_{p,q}$-converge to one another. This proves $d\mu$
is proper, and hence, proves (I).

Q.E.D.

\subsection{\label{subsec:3.3.4 Toward-a-Taxonomy}Toward a Taxonomy of $\left(p,q\right)$-adic
Measures}

Although frames lubricate our discussion of the convergence behavior
of Fourier series generated by the Fourier-Stieltjes transforms of
$\left(p,q\right)$-adic measures, we still do not have many practical
examples of these measures. In the present subsection, I have gathered
the few results about $\left(p,q\right)$-adic measures which I \emph{have
}been able to prove, in the hopes of being able to give more substantial
descriptions of them. Somewhat vexingly, the methods in this subsection
depend nigh-entirely on the structural properties of whichever type
of measure we happen to be considering. Without assuming such properties
hold, it becomes difficult to say anything meaningful about the measures
and the Fourier series they generate\textemdash this is one of the
main challenges we have to deal with when doing theory with frames.
In order to keep this subsection text from drowning in adjectives,
we will use the following abbreviations: 
\begin{defn}[\textbf{Measure Type Codes}]
\ 

\vphantom{}

I.\textbf{ }(\textbf{ED1}) We say \index{measure!elementary type D1}$d\mu\in M_{\textrm{alg}}\left(\mathbb{Z}_{p},\mathbb{C}_{q}\right)$
is an \textbf{elementary D1 measure }(ED1) whenever it can be written
as $d\mu=d\nu*d\eta$ for $d\nu,d\eta\in M_{\textrm{alg}}\left(\mathbb{Z}_{p},\mathbb{C}_{q}\right)$,
where $d\nu$ is radially symmetric and where $d\eta$ is $\left(p,q\right)$-adically
regular magnitudinal measure. That is to say, $\hat{\nu}\left(t\right)$
is a $\overline{\mathbb{Q}}$-valued function which depends only on
$\left|t\right|_{p}$, and $\hat{\eta}\left(t\right)$ is of the form:
\begin{equation}
\hat{\eta}\left(t\right)=\begin{cases}
\kappa\left(0\right) & \textrm{if }t=0\\
\sum_{m=0}^{\left|t\right|_{p}-1}\kappa\left(m\right)e^{-2\pi imt} & \textrm{else}
\end{cases},\textrm{ }\forall t\in\hat{\mathbb{Z}}_{p}
\end{equation}
for some function $\kappa:\mathbb{N}_{0}\rightarrow\overline{\mathbb{Q}}$
such that: 
\begin{equation}
\lim_{n\rightarrow\infty}\left|\kappa\left(\left[\mathfrak{z}\right]_{p^{n}}\right)\right|_{q}=0,\textrm{ }\forall\mathfrak{z}\in\mathbb{Z}_{p}^{\prime}
\end{equation}
and, moreover, for which there are constants $a_{1},\ldots,a_{p-1}\in\overline{\mathbb{Q}}$
so that: 
\[
\kappa\left(\left[\mathfrak{z}\right]_{p^{n}}+jp^{n}\right)=a_{j}\kappa\left(\left[\mathfrak{z}\right]_{p^{n}}\right),\textrm{ }\forall\mathfrak{z}\in\mathbb{Z}_{p},\textrm{ }\forall n\geq1,\textrm{ }\forall j\in\left\{ 1,\ldots,p-1\right\} 
\]
Note that this then implies $\kappa\left(0\right)=1$.

\vphantom{}

II. (\textbf{D1}) We say $d\mu\in M_{\textrm{alg}}\left(\mathbb{Z}_{p},\mathbb{C}_{q}\right)$
a \textbf{type D1 measure }whenever\index{measure!type D1} it is
a $\overline{\mathbb{Q}}$-linear combination of finitely many ED1
measures.

\vphantom{}

III. (\textbf{ED2}) We say $d\mu\in M_{\textrm{alg}}\left(\mathbb{Z}_{p},\mathbb{C}_{q}\right)$
an \textbf{elementary D2 measure} whenever it is ED1, with\index{measure!elementary type D2}:
\begin{equation}
v_{n}=\left(1+A\right)v_{n+1},\textrm{ }\forall n\geq0\label{eq:ED2 measure}
\end{equation}
where, recall: 
\begin{equation}
A=\sum_{j=1}^{p-1}a_{j}=\sum_{j=1}^{p-1}\kappa\left(j\right)
\end{equation}
is the sum of the constants from the structural equation of $\kappa$
and $v_{n}$ denotes $\hat{\nu}\left(p^{-n}\right)$.

\vphantom{}

IV. (\textbf{D2}) We say $d\mu\in M_{\textrm{alg}}\left(\mathbb{Z}_{p},\mathbb{C}_{q}\right)$
a \textbf{type D2 measure }whenever it is a $\overline{\mathbb{Q}}$-linear
combination of finitely many ED2 measures. \index{measure!type D2} 
\end{defn}
\vphantom{}

Obviously, there are many other variants we could consider. For example,
w measures of the form: 
\begin{equation}
\hat{\mu}\left(t\right)=\begin{cases}
\kappa\left(0\right) & \textrm{if }t=0\\
\sum_{n=0}^{-v_{p}\left(t\right)-1}\kappa\left(n\right)e^{-2\pi int} & \textrm{else}
\end{cases},\textrm{ }\forall t\in\hat{\mathbb{Z}}_{p}
\end{equation}
for some $\left(p,q\right)$-adically regular $\kappa:\mathbb{N}_{0}\rightarrow\overline{\mathbb{Q}}$;
also, spaces generated by convolutions of measures with this type,
or any combination of convolutions of radially-symmetric and/or magnitudinal
measures, ED1 measures for which: 
\begin{equation}
\lim_{n\rightarrow\infty}\hat{\nu}\left(p^{-N}\right)p^{N}\overset{\mathbb{C}}{=}0
\end{equation}
and so on and so forth. 
\begin{prop}
\label{prop:ED1 measures are zero whenever kappa of 0 is 0}An\textbf{
}ED1\textbf{ }measure is the zero measure whenever $\kappa\left(0\right)=0$. 
\end{prop}
Proof: If $\kappa\left(0\right)=0$, then \textbf{Lemma \ref{lem:structural equations uniquely determine p-adic structured functions}}
applies, telling us that $\kappa$ is identically zero. This makes
$\hat{\eta}$ identically zero, which makes $\hat{\mu}=\hat{\nu}\cdot\hat{\eta}$
identically zero. This proves that $d\mu$ is the zero measure.

Q.E.D. 
\begin{lem}
\label{lem:degenerate rad sym measures are zero}The only degenerate
radially symmetric measure is the zero measure. 
\end{lem}
Proof: Let $d\mu$ be a degenerate radially-symmetric measure. Using
the Radial Fourier Resummation formula (\ref{eq:Radial Fourier Resummation Lemma}),
the degeneracy of $d\mu$ allows us to write: 
\begin{equation}
0\overset{\mathbb{C}_{q}}{=}\lim_{N\rightarrow\infty}\tilde{\mu}_{N}\left(\mathfrak{z}\right)\overset{\overline{\mathbb{Q}}}{=}\hat{\mu}\left(0\right)-\left|\mathfrak{z}\right|_{p}^{-1}\hat{\mu}\left(\frac{\left|\mathfrak{z}\right|_{p}}{p}\right)+\left(1-\frac{1}{p}\right)\sum_{n=1}^{v_{p}\left(\mathfrak{z}\right)}\hat{\mu}\left(\frac{1}{p^{n}}\right)p^{n},\textrm{ }\forall\mathfrak{z}\in\mathbb{Z}_{p}\backslash\left\{ 0\right\} \label{eq:Degeneracy Lemma - Algebraic Vanishing}
\end{equation}
Hence, it must be that the algebraic number: 
\begin{equation}
\hat{\mu}\left(0\right)-\left|\mathfrak{z}\right|_{p}^{-1}\hat{\mu}\left(\frac{\left|\mathfrak{z}\right|_{p}}{p}\right)+\left(1-\frac{1}{p}\right)\sum_{n=1}^{v_{p}\left(\mathfrak{z}\right)}\hat{\mu}\left(\frac{1}{p^{n}}\right)p^{n}
\end{equation}
is $0$ for all non-zero $p$-adic integers $\mathfrak{z}$.

Letting $\left|\mathfrak{z}\right|_{p}=p^{k}$ for $k\geq1$, and
letting $\hat{\mu}\left(p^{-n}\right)=c_{n}$, the above becomes:

\begin{equation}
0\overset{\overline{\mathbb{Q}}}{=}c_{0}-c_{k+1}p^{k}+\left(1-\frac{1}{p}\right)\sum_{n=1}^{k}c_{n}p^{n},\textrm{ }\forall k\geq1
\end{equation}
Next, we subtract the $\left(k-1\right)$th case from the $k$th case.
\begin{equation}
0\overset{\overline{\mathbb{Q}}}{=}c_{k}p^{k-1}-c_{k+1}p^{k}+\left(1-\frac{1}{p}\right)c_{k}p^{k},\textrm{ }\forall k\geq2
\end{equation}
and hence: 
\begin{equation}
c_{k+1}\overset{\overline{\mathbb{Q}}}{=}c_{k},\textrm{ }\forall k\geq2
\end{equation}
Meanwhile, when $k=1$, we get: 
\begin{equation}
c_{2}\overset{\overline{\mathbb{Q}}}{=}\frac{1}{p}c_{0}+\left(1-\frac{1}{p}\right)c_{1}
\end{equation}
As for the $\mathfrak{z}=0$ case, using the \textbf{Radial Fourier
Resummation Lemma }(\textbf{Lemma \ref{lem:Rad-Sym Fourier Resum Lemma}}),
the degeneracy of $d\mu$ tell us that:

\[
0\overset{\mathbb{C}}{=}\lim_{N\rightarrow\infty}\tilde{\mu}_{N}\left(0\right)=c_{0}+\left(p-1\right)\lim_{N\rightarrow\infty}\sum_{k=1}^{N}c_{k}p^{k-1}\overset{\mathbb{C}}{=}c_{0}+\left(1-\frac{1}{p}\right)\sum_{k=1}^{\infty}c_{k}p^{k}
\]
This then forces $c_{n}p^{n}$ to tend to $0$ in $\mathbb{C}$ as
$n\rightarrow\infty$. However, we just showed that $c_{k}=c_{k+1}$
for all $k\geq2$. So, the only way convergence occurs is if $c_{k}=0$
for all $k\geq2$. Consequently, $\hat{\mu}\left(t\right)=0$ for
all $\left|t\right|_{p}\geq p^{2}$.

Using the \textbf{Radial Fourier Resummation Lemma} once more, we
obtain: 
\begin{equation}
\tilde{\mu}\left(\mathfrak{z}\right)\overset{\overline{\mathbb{Q}}}{=}\begin{cases}
\hat{\mu}\left(0\right)-\hat{\mu}\left(\frac{1}{p}\right) & \textrm{if }\left|\mathfrak{z}\right|_{p}=1\\
\hat{\mu}\left(0\right)+\left(p-1\right)\hat{\mu}\left(\frac{1}{p}\right) & \textrm{if }\left|\mathfrak{z}\right|_{p}\leq\frac{1}{p}
\end{cases}
\end{equation}
Since $d\mu$ is degenerate, $\tilde{\mu}\left(\mathfrak{z}\right)$
is identically zero. This forces $\hat{\mu}\left(0\right)=\hat{\mu}\left(\frac{1}{p}\right)$.
So, we can write: 
\begin{equation}
0=\hat{\mu}\left(0\right)+\left(p-1\right)\hat{\mu}\left(\frac{1}{p}\right)=p\hat{\mu}\left(0\right)
\end{equation}
which forces $\hat{\mu}\left(0\right)=\hat{\mu}\left(\frac{1}{p}\right)=\hat{\mu}\left(\frac{1}{p^{2}}\right)=\cdots=0$.
This proves that $d\mu$ is the zero measure.

Q.E.D. 
\begin{prop}
Let $d\mu=d\nu*d\eta$ be an ED2 measure. Then, $d\mu$ is degenerate
if and only if: 
\begin{equation}
\lim_{N\rightarrow\infty}\hat{\nu}\left(p^{-N}\right)p^{N}\overset{\mathbb{C}}{=}0\label{eq:Degeneracy criterion for elementary D4 measures}
\end{equation}
\end{prop}
Proof: Since $d\mu$ is elementary and non-zero, the fact that every
ED2 measure is ED1 means that \textbf{Proposition \ref{prop:ED1 measures are zero whenever kappa of 0 is 0}
}applies, guaranteeing that $\kappa\left(0\right)\neq0$.

Now, supposing $d\mu$ satisfies (\ref{eq:Degeneracy criterion for elementary D4 measures}),
the \textbf{Radially-Magnitudinal Fourier Resummation Lemma }(\textbf{Lemma
\ref{lem:Radially-Mag Fourier Resummation Lemma}}) becomes:

\begin{equation}
\tilde{\mu}_{N}\left(\mathfrak{z}\right)=\sum_{\left|t\right|_{p}\leq p^{N}}\hat{\mu}\left(t\right)e^{2\pi i\left\{ t\mathfrak{z}\right\} _{p}}\overset{\overline{\mathbb{Q}}}{=}\hat{\nu}\left(p^{-N}\right)\kappa\left(\left[\mathfrak{z}\right]_{p^{N}}\right)p^{N}
\end{equation}
The $\left(p,q\right)$-adic regularity of $\kappa$ and along with
the assumptions then guarantee that $\tilde{\mu}_{N}\left(\mathfrak{z}\right)$
tends to $0$ in $\mathbb{C}$ (resp. $\mathbb{C}_{q}$) for $\mathfrak{z}\in\mathbb{N}_{0}$
(resp. $\mathbb{Z}_{p}^{\prime}$) as $N\rightarrow\infty$, which
shows that $d\mu$ is then degenerate.

Conversely, suppose that the non-zero ED2 measure $d\mu$ is degenerate.
Using \textbf{Lemma \ref{lem:Radially-Mag Fourier Resummation Lemma}}
once again, we write:

\begin{align}
\sum_{\left|t\right|_{p}\leq p^{N}}\hat{\mu}\left(t\right)e^{2\pi i\left\{ t\mathfrak{z}\right\} _{p}} & \overset{\mathbb{\overline{\mathbb{Q}}}}{=}v_{N}\kappa\left(\left[\mathfrak{z}\right]_{p^{N}}\right)p^{N}+\sum_{n=0}^{N-1}\left(v_{n}-\left(1+A\right)v_{n+1}\right)\kappa\left(\left[\mathfrak{z}\right]_{p^{n}}\right)p^{n}\label{eq:Radial Magnitudinal identity for the Degeneracy Lemma}
\end{align}
where, recall, $v_{n}=\hat{\nu}\left(p^{-n}\right)$. Since $d\mu$
is ED2, $v_{n}=Av_{n+1}$ for all $n\geq0$, and hence: 
\begin{equation}
\sum_{\left|t\right|_{p}\leq p^{N}}\hat{\mu}\left(t\right)e^{2\pi i\left\{ t\mathfrak{z}\right\} _{p}}\overset{\mathbb{\overline{\mathbb{Q}}}}{=}v_{N}\kappa\left(\left[\mathfrak{z}\right]_{p^{N}}\right)p^{N}
\end{equation}
Since $d\mu$ is ED2, it is also ED1, and as such, the assumption
that $d\mu$ is non-zero guarantees that $\kappa\left(0\right)\neq0$.
So, setting $\mathfrak{z}=0$, the degeneracy of $d\mu$ tells us
that: 
\begin{equation}
0\overset{\mathbb{C}}{=}\lim_{N\rightarrow\infty}\sum_{\left|t\right|_{p}\leq p^{N}}\hat{\mu}\left(t\right)\overset{\mathbb{\overline{\mathbb{Q}}}}{=}\kappa\left(0\right)\lim_{N\rightarrow\infty}v_{N}p^{N}
\end{equation}
Since $\kappa\left(0\right)\neq0$, this forces $\lim_{N\rightarrow\infty}v_{N}p^{N}$
to converge to $0$ in $\mathbb{C}$, which is exactly what we wished
to show (equation (\ref{eq:Degeneracy criterion for elementary D4 measures})).

Q.E.D. 
\begin{prop}
Let $d\mu=d\nu*d\eta$ be an ED1 measure. Then: 
\end{prop}
\begin{equation}
\tilde{\mu}_{N}\left(\mathfrak{z}\right)\overset{\overline{\mathbb{Q}}}{=}\hat{\nu}\left(0\right)\kappa\left(0\right)+\sum_{n=1}^{N}\hat{\nu}\left(p^{-n}\right)\left(\kappa\left(\left[\mathfrak{z}\right]_{p^{n}}\right)-\kappa\left(\left[\mathfrak{z}\right]_{p^{n-1}}\right)\right)p^{n}\label{eq:radial-magnitude Fourier resummation w/ SbP}
\end{equation}

Proof: Apply summation by parts to equation (\ref{eq:Radial-Magnitude p-adic distributed Fourier Resummation Identity, simplified})
from \textbf{Lemma \ref{lem:Radially-Mag Fourier Resummation Lemma}}.

Q.E.D. 
\begin{prop}
\label{prop:non-zero degenerate ED1}Let $d\mu=d\nu*d\eta$ be a non-zero
degenerate ED1 measure. Then $\left|\kappa\left(j\right)\right|_{q}<1$
for all $j\in\left\{ 1,\ldots,p-1\right\} $. In particular, there
is a $j\in\left\{ 1,\ldots,p-1\right\} $ so that $\left|\kappa\left(j\right)\right|_{q}>0$. 
\end{prop}
Proof: Suppose that the non-zero ED1\textbf{ }measure $d\mu$ is degenerate.
Then, by equation (\ref{eq:Radial Magnitudinal identity for the Degeneracy Lemma}):

\begin{align}
\sum_{\left|t\right|_{p}\leq p^{N}}\hat{\mu}\left(t\right)e^{2\pi i\left\{ t\mathfrak{z}\right\} _{p}} & \overset{\mathbb{\overline{\mathbb{Q}}}}{=}v_{N}\kappa\left(\left[\mathfrak{z}\right]_{p^{N}}\right)p^{N}+\sum_{n=0}^{N-1}\left(v_{n}-\left(1+A\right)v_{n+1}\right)\kappa\left(\left[\mathfrak{z}\right]_{p^{n}}\right)p^{n}\label{eq:Radial Magnitudinal identity for the Degeneracy Lemma-1}
\end{align}
where, again, $v_{n}=\hat{\nu}\left(p^{-n}\right)$.

Since $p$ and $q$ are distinct primes, and since $d\nu$ is a $\left(p,q\right)$-adic
measure, we have that: 
\begin{equation}
\sup_{N\geq1}\left|v_{N}p^{N}\right|_{q}<\infty
\end{equation}
Also, since $\kappa$ is $\left(p,q\right)$-adically regular, letting
$\mathfrak{z}\in\mathbb{Z}_{p}^{\prime}$ be arbitrary, taking the
limit of (\ref{eq:Radial Magnitudinal identity for the Degeneracy Lemma-1})
in $\mathbb{C}_{q}$ as $N\rightarrow\infty$ gives: 
\begin{equation}
\sum_{n=0}^{\infty}\left(v_{n}-\left(1+A\right)v_{n+1}\right)\kappa\left(\left[\mathfrak{z}\right]_{p^{n}}\right)p^{n}\overset{\mathbb{C}_{q}}{=}0,\textrm{ }\forall\mathfrak{z}\in\mathbb{Z}_{p}^{\prime}\label{eq:Degeneracy lemma, equation 1-1}
\end{equation}
which is equal to $0$, by the presumed degeneracy of $d\mu$. Now,
letting $j\in\left\{ 1,\ldots,p-1\right\} $, observe that $\left[-j\right]_{p^{n}}$
has exactly $n$ $p$-adic digits, all of which are $p-j$. That is:
\begin{equation}
\left[-j\right]_{p^{n+1}}=\left[-j\right]_{p^{n}}+\left(p-j\right)p^{n},\textrm{ }\forall n\geq0
\end{equation}
Because $\kappa$ is $\left(p,q\right)$-adically regular and non-zero,
we can then write: 
\[
\kappa\left(\left[-j\right]_{p^{n+1}}\right)=\kappa\left(\left[-j\right]_{p^{n}}+\left(p-j\right)p^{n}\right)=\kappa\left(p-j\kappa\right)\left(\left[-j\right]_{p^{n}}\right),\textrm{ }\forall n\geq0
\]
Thus, by induction: 
\begin{equation}
\kappa\left(\left[-j\right]_{p^{n}}\right)=\left(\kappa\left(p-j\right)\right)^{n}\kappa\left(0\right),\textrm{ }\forall j\in\left\{ 1,\ldots,p-1\right\} ,\textrm{ }\forall n\geq0
\end{equation}
Moreover, because $\left|\kappa\left(\left[\mathfrak{z}\right]_{p^{n}}\right)\right|_{q}\rightarrow0$
in $\mathbb{R}$ for all $\mathfrak{z}\in\mathbb{Z}_{p}^{\prime}$,
it must be that: 
\begin{equation}
0\overset{\mathbb{R}}{=}\lim_{n\rightarrow\infty}\left|\kappa\left(\left[-j\right]_{p^{n}}\right)\right|_{q}=\lim_{n\rightarrow\infty}\left|\kappa\left(p-j\right)\right|_{q}^{n}\left|\kappa\left(0\right)\right|_{q}
\end{equation}
Since $\kappa\left(0\right)\neq0$, this forces $\left|\kappa\left(p-j\right)\right|_{q}<1$
for all $j\in\left\{ 1,\ldots p-1\right\} $.

So, for $\mathfrak{z}=-j$, equation (\ref{eq:Degeneracy lemma, equation 1-1})
becomes: 
\begin{equation}
\kappa\left(0\right)\sum_{n=0}^{\infty}\left(v_{n}-\left(1+A\right)v_{n+1}\right)\left(\kappa\left(p-j\right)\right)^{n}p^{n}\overset{\mathbb{C}_{q}}{=}0,\textrm{ }\forall j\in\left\{ 1,\ldots,p-1\right\} 
\end{equation}
Since $\kappa\left(0\right)\neq0$, it must be that the series is
zero. Re-indexing the $j$s then gives us: 
\begin{equation}
\sum_{n=0}^{\infty}\left(v_{n}-\left(1+A\right)v_{n+1}\right)\left(\kappa\left(j\right)\right)^{n}p^{n}\overset{\mathbb{C}_{q}}{=}0,\textrm{ }\forall j\in\left\{ 1,\ldots,p-1\right\} 
\end{equation}

\begin{claim}
There exists a $j\in\left\{ 1,\ldots,p-1\right\} $ so that $\left|\kappa\left(j\right)\right|_{q}\in\left(0,1\right)$.

Proof of claim: We already know that all of the $\kappa\left(j\right)$s
have $q$-adic absolute value $<1$. So by way of contradiction, suppose
$\kappa\left(j\right)=0$ for all $j\in\left\{ 1,\ldots,p-1\right\} $.
This tells us that $\kappa\left(n\right)$ is $1$ when $n=0$ and
is $0$ for all $n\geq1$. Consequently, $\hat{\eta}\left(t\right)=\kappa\left(0\right)$
for all $t\in\hat{\mathbb{Z}}_{p}$. With this, we see that the degenerate
measure $d\mu$ must be radially-symmetric, and so, $\hat{\mu}\left(t\right)=\kappa\left(0\right)\hat{\nu}\left(t\right)=\hat{\nu}\left(t\right)$
for some radially-symmetric $\hat{\nu}$. Using \textbf{Lemma \ref{lem:degenerate rad sym measures are zero}},
$d\mu$ is then the zero measure, which contradicts the fact that
$d\mu$ was given to be non-zero. This proves the claim. 
\end{claim}
\vphantom{}

This shows that $\left|\kappa\left(j\right)\right|_{q}<1$ for all
$j\in\left\{ 1,\ldots,p-1\right\} $, and that there is one such $j$
for which $0<\left|\kappa\left(j\right)\right|_{q}<1$.

Q.E.D. 
\begin{prop}
Let $d\mu=d\nu*d\eta$ be a non-zero degenerate ED1 measure. Then
$A\neq p-1$, where, recall, $A=\sum_{j=1}^{p-1}\kappa\left(j\right)$. 
\end{prop}
Proof: Let $d\mu$ be as given. By way of contradiction, suppose $A=p-1$.
Note that since $d\mu$ is non-zero, we have that $\kappa\left(0\right)=1$.

Now, letting $\mathfrak{z}=0$, in (\ref{eq:radial-magnitude Fourier resummation w/ SbP})
we get: 
\begin{equation}
\tilde{\mu}_{N}\left(\mathfrak{z}\right)\overset{\overline{\mathbb{Q}}}{=}v_{0}\kappa\left(0\right)+\sum_{n=1}^{N}v_{n}\left(\kappa\left(0\right)-\kappa\left(0\right)\right)p^{n}=v_{0}\kappa\left(0\right)
\end{equation}
Letting $N\rightarrow\infty$, we have: 
\begin{equation}
0=v_{0}\kappa\left(0\right)
\end{equation}
Since $\kappa\left(0\right)=1$, this then implies $v_{0}=0$.

Next, by \textbf{Proposition \ref{prop:non-zero degenerate ED1}},
we can fix a $j\in\left\{ 1,\ldots,p-1\right\} $ so that $0<\left|\kappa\left(j\right)\right|_{q}<1$.
Since $\kappa$ is $\left(p,q\right)$-adically regular, we have:
\begin{equation}
\kappa\left(\left[j+jp+\cdots+jp^{m-1}\right]_{p^{n}}\right)=\kappa\left(j\right)\kappa\left(\left[j+jp+\cdots+jp^{m-2}\right]_{p^{n}}\right)
\end{equation}
As such: 
\begin{equation}
\kappa\left(\left[j\frac{p^{m}-1}{p-1}\right]_{p^{n}}\right)=\kappa\left(\left[j+jp+\cdots+jp^{m-1}\right]_{p^{n}}\right)=\left(\kappa\left(j\right)\right)^{\min\left\{ n,m\right\} }\label{eq:kappa of j times (p to the m -1 over p -1)}
\end{equation}
and so, letting $\mathfrak{z}=j\frac{p^{m}-1}{p-1}\in\mathbb{N}_{1}$
for $m\geq1$, (\ref{eq:radial-magnitude Fourier resummation w/ SbP})
becomes:

\begin{align*}
\tilde{\mu}_{N}\left(j\frac{p^{m}-1}{p-1}\right) & \overset{\overline{\mathbb{Q}}}{=}v_{0}+\sum_{n=1}^{N}v_{n}\left(\kappa\left(\left[j\frac{p^{m}-1}{p-1}\right]_{p^{n}}\right)-\kappa\left(\left[j\frac{p^{m}-1}{p-1}\right]_{p^{n-1}}\right)\right)p^{n}\\
 & =v_{0}+\sum_{n=1}^{N}v_{n}\left(\kappa\left(j\right)^{\min\left\{ n,m\right\} }-\left(\kappa\left(j\right)\right)^{\min\left\{ n-1,m\right\} }\right)p^{n}\\
 & =v_{0}+\sum_{n=1}^{m-1}v_{n}\left(\kappa\left(j\right)\right)^{n}p^{n}-\sum_{n=1}^{m}v_{n}\left(\kappa\left(j\right)\right)^{n-1}p^{n}\\
 & +\left(\kappa\left(j\right)\right)^{m}\left(\sum_{n=m}^{N}v_{n}p^{n}-\sum_{n=m+1}^{N}v_{n}p^{n}\right)\\
 & =v_{0}+\sum_{n=1}^{m}v_{n}\left(\left(\kappa\left(j\right)\right)^{n}-\left(\kappa\left(j\right)\right)^{n-1}\right)p^{n}
\end{align*}
Letting $N\rightarrow\infty$ gives: 
\begin{align}
v_{0}+\sum_{n=1}^{m}v_{n}\left(\left(\kappa\left(j\right)\right)^{n}-\left(\kappa\left(j\right)\right)^{n-1}\right)p^{n} & \overset{\overline{\mathbb{Q}}}{=}0,\textrm{ }\forall m\geq1\label{eq:m cases}
\end{align}
Subtracting the $\left(m-1\right)$th case from the $m$th case for
$m\geq2$ leaves us with: 
\begin{equation}
v_{m}\left(\left(\kappa\left(j\right)\right)^{m}-\left(\kappa\left(j\right)\right)^{m-1}\right)p^{m}\overset{\overline{\mathbb{Q}}}{=}0,\textrm{ }\forall m\geq2
\end{equation}
Multiplying through by the non-zero algebraic number $\left(\kappa\left(j\right)\right)^{1-m}$
yields: 
\begin{equation}
v_{m}\left(\kappa\left(j\right)-1\right)p^{m}\overset{\overline{\mathbb{Q}}}{=}0,\textrm{ }\forall m\geq2
\end{equation}
Here, $0<\left|\kappa\left(j\right)\right|_{q}<1$ guarantees that
$\left(\kappa\left(j\right)-1\right)p^{m}\neq0$, and hence, it must
be that $v_{m}=0$ for all $m\geq2$.

So, only the $m=1$ case of (\ref{eq:m cases}) is left: 
\begin{equation}
v_{0}+v_{1}\left(\kappa\left(j\right)-1\right)p\overset{\overline{\mathbb{Q}}}{=}0
\end{equation}
Since $v_{0}=0$, this gives: 
\begin{equation}
v_{1}\left(\kappa\left(j\right)-1\right)p\overset{\overline{\mathbb{Q}}}{=}0
\end{equation}
Since $0<\left|\kappa\left(j\right)\right|_{q}<1$, this then shows
that $\left(\kappa\left(j\right)-1\right)p\neq0$, and hence that
$v_{1}=0$. So, $v_{n}=0$ for all $n\geq0$. This means the radially
symmetric measure $d\nu$ is identically zero. Thus, so is $d\mu$\textemdash but
$d\mu$ was given to be non-zero. There's our contradiction.

Thus, if $d\mu$ is non-zero, it must be that $A\neq p-1$.

Q.E.D. 
\begin{prop}
Let $d\nu\in M_{\textrm{dgen}}\left(\mathbb{Z}_{p},K\right)$. Then
$e^{2\pi i\left\{ \tau\mathfrak{z}\right\} _{p}}d\nu\left(\mathfrak{z}\right)$
is a degenerate measure for all $\tau\in\hat{\mathbb{Z}}_{p}$. 
\end{prop}
Proof: The Fourier-Stieltjes transform of $d\mu\left(\mathfrak{z}\right)\overset{\textrm{def}}{=}e^{2\pi i\left\{ \tau\mathfrak{z}\right\} _{p}}d\nu\left(\mathfrak{z}\right)$
is $\hat{\nu}\left(t-\tau\right)$. Then, letting $N$ be large enough
so that $\left|\tau\right|_{p}\leq p^{N}$, we have that the map $t\mapsto t-\tau$
is then a bijection of the set $\left\{ t\in\hat{\mathbb{Z}}_{p}:\left|t\right|_{p}\leq p^{N}\right\} $.
For all such $N$: 
\begin{align*}
\tilde{\mu}_{N}\left(\mathfrak{z}\right) & =\sum_{\left|t\right|_{p}\leq p^{N}}\hat{\nu}\left(t-\tau\right)e^{2\pi i\left\{ t\mathfrak{z}\right\} _{p}}\\
 & =\sum_{\left|t\right|_{p}\leq p^{N}}\hat{\nu}\left(t\right)e^{2\pi i\left\{ \left(t+\tau\right)\mathfrak{z}\right\} _{p}}\\
 & =e^{2\pi i\left\{ \tau\mathfrak{z}\right\} _{p}}\tilde{\nu}_{N}\left(\mathfrak{z}\right)
\end{align*}
Since $e^{2\pi i\left\{ \tau\mathfrak{z}\right\} _{p}}$ is bounded
in $\mathbb{C}$ for all $\mathfrak{z}\in\mathbb{N}_{0}$, the degeneracy
of $d\nu$ guarantees that $\tilde{\mu}_{N}\left(\mathfrak{z}\right)\overset{\mathbb{C}}{\rightarrow}0$
as $N\rightarrow\infty$ for all $\mathfrak{z}\in\mathbb{N}_{0}$.
Likewise, since $e^{2\pi i\left\{ \tau\mathfrak{z}\right\} _{p}}$
is bounded in $\mathbb{C}_{q}$ for all $\mathfrak{z}\in\mathbb{Z}_{p}^{\prime}$,
the degeneracy of $d\nu$ guarantees that $\tilde{\mu}_{N}\left(\mathfrak{z}\right)\overset{\mathbb{C}_{q}}{\rightarrow}0$
as $N\rightarrow\infty$ for all $\mathfrak{z}\in\mathbb{Z}_{p}^{\prime}$.
Hence, $d\mu=e^{2\pi i\left\{ \tau\mathfrak{z}\right\} _{p}}d\nu\left(\mathfrak{z}\right)$
is degenerate.

Q.E.D.

\vphantom{}

For our purposes, the following proposition will be of great import. 
\begin{prop}
Let $d\mu$ be a degenerate ED2 measure. Then, either $d\mu$ is the
zero measure, or $\hat{\mu}\left(0\right)\neq0$. 
\end{prop}
Proof: Let $d\mu$ be as given. If $d\mu$ is the zero measure, then
$\hat{\mu}\left(0\right)=0$. Thus, to complete the proof, we need
only show that $\hat{\mu}\left(0\right)\neq0$ whenever $d\mu$ is
a non-zero degenerate ED2 measure.

First, let $d\mu$ be elementary, and suppose $\hat{\mu}\left(0\right)=0$.
Since $\hat{\mu}\left(0\right)=\hat{\nu}\left(0\right)\hat{\eta}\left(0\right)=\hat{\nu}\left(1\right)\kappa\left(0\right)$,
there are two possibilities: either $\kappa\left(0\right)=0$ or $\hat{\nu}\left(1\right)=0$.
If $\kappa\left(0\right)=0$, \textbf{Lemma \ref{lem:structural equations uniquely determine p-adic structured functions}}
shows that $d\mu$ is then the zero measure.

Alternatively, suppose $\hat{\nu}\left(1\right)=0$. Then, since $d\nu$
is an ED2 measure, we have that:

\begin{equation}
\hat{\nu}\left(p^{-n}\right)=\left(1+A\right)\hat{\nu}\left(p^{-\left(n+1\right)}\right),\textrm{ }\forall n\geq0
\end{equation}
If $A=-1$, then $\hat{\nu}\left(p^{-n}\right)$ is zero for all $n\geq1$.
Since $\hat{\nu}\left(1\right)=\hat{\nu}\left(p^{-0}\right)$ was
assumed to be $0$, the radial symmetry of $\hat{\nu}$ forces $\hat{\nu}\left(t\right)=0$
for all $t\in\hat{\mathbb{Z}}_{p}$. Since $\hat{\mu}\left(t\right)=\hat{\nu}\left(t\right)\hat{\eta}\left(t\right)$,
this shows that $d\mu$ is the zero measure.

On the other hand, if $A\neq-1$, then: 
\begin{equation}
\hat{\nu}\left(p^{-\left(n+1\right)}\right)=\frac{1}{A+1}\hat{\nu}\left(p^{-n}\right),\textrm{ }\forall n\geq0
\end{equation}
and thus, the assumption that $\hat{\nu}\left(1\right)=0$ then forces
$\hat{\nu}\left(p^{-n}\right)=0$ for all $n\geq0$. Like with the
$A=-1$ case, this forces $d\mu$ to be the zero measure.

Thus, either $\hat{\mu}\left(0\right)\neq0$ or $d\mu$ is the zero
measure.

Q.E.D. 
\begin{question}
\label{que:3.3}Although we shall not explore it here, there is also
the question of how spaces of $\left(p,q\right)$-adic measures behave
under pre-compositions by maps $\phi:\mathbb{Z}_{p}\rightarrow\mathbb{Z}_{p}$.
For example, given a rising-continuous $d\mu$, consider the measure
$d\nu$ defined by: 
\begin{equation}
\hat{\nu}\left(t\right)\overset{\textrm{def}}{=}\hat{\mu}\left(t\right)e^{-2\pi i\left\{ t\mathfrak{a}\right\} _{p}}
\end{equation}
where $\mathfrak{a}\in\mathbb{Z}_{p}$ is a constant (i.e., $d\nu\left(\mathfrak{z}\right)=d\mu\left(\mathfrak{z}-\mathfrak{a}\right)$).
Will $d\nu$ still be rising continuous? How will the frame of convergence
for $\lim_{N\rightarrow\infty}\tilde{\mu}_{N}\left(\mathfrak{z}\right)$
be affected by the map $d\mu\left(\mathfrak{z}\right)\mapsto d\mu\left(\mathfrak{z}-\mathfrak{a}\right)$?
The most general case of this kind of question would be to consider
the map $d\mu\left(\mathfrak{z}\right)\mapsto d\mu_{\phi}\left(\mathfrak{z}\right)$
where $\phi:\mathbb{Z}_{p}\rightarrow\mathbb{Z}_{p}$ is any sufficiently
nice map, as in \textbf{\emph{Definition \ref{def:pullback measure}}}\emph{
on page \pageref{def:pullback measure}}. These all appear to be worth
exploring. 
\end{question}

\subsection{\label{subsec:3.3.5 Quasi-Integrability}Quasi-Integrability}
\begin{rem}
Like in Subsection \ref{subsec:3.3.3 Frames}, this subsection exists
primarily to give a firm foundation for my concept of quasi-integrability.
With regard to the analysis of $\chi_{H}$ to be performed in Chapter
4, the only thing the reader needs to know is that I say a function
$\chi:\mathbb{Z}_{p}\rightarrow\mathbb{C}_{q}$ is \textbf{quasi-integrable
}with\index{quasi-integrability} respect to the standard $\left(p,q\right)$-adic
frame whenever there exists a function $\hat{\chi}:\hat{\mathbb{Z}}_{p}\rightarrow\overline{\mathbb{Q}}$
so that: 
\begin{equation}
\chi\left(\mathfrak{z}\right)\overset{\mathcal{F}_{p,q}}{=}\lim_{N\rightarrow\infty}\sum_{\left|t\right|_{p}\leq p^{N}}\hat{\chi}\left(t\right)e^{2\pi i\left\{ t\mathfrak{z}\right\} _{p}},\textrm{ }\forall\mathfrak{z}\in\mathbb{Z}_{p}
\end{equation}
where, recall, the $\mathcal{F}_{p,q}$ means that we interpret the
right-hand side in the topology of $\mathbb{C}$ when $\mathfrak{z}\in\mathbb{N}_{0}$
and in the topology of $\mathbb{C}_{q}$ whenever $\mathfrak{z}\in\mathbb{Z}_{p}^{\prime}$.
We call the function $\hat{\chi}$ \emph{a }\textbf{Fourier transform}\footnote{The reason we cannot say \emph{the }Fourier transform of $\chi$ is
because of the existence of degenerate measures. For those who did
not read all of Subsection \ref{subsec:3.3.3 Frames}, a $\left(p,q\right)$-adic
measure $d\mu$ is said to \index{measure!degenerate}be \textbf{degenerate
}(with respect to the standard frame) whenever its Fourier-Stieltjes
transform $\hat{\mu}$ takes values in $\overline{\mathbb{Q}}$ and
satisfies the limit condition: 
\begin{equation}
\lim_{N\rightarrow\infty}\underbrace{\sum_{\left|t\right|_{p}\leq p^{N}}\hat{\mu}\left(t\right)e^{2\pi i\left\{ t\mathfrak{z}\right\} _{p}}}_{\tilde{\mu}_{N}\left(\mathfrak{z}\right)}\overset{\mathcal{F}_{p,q}}{=}0,\textrm{ }\forall\mathfrak{z}\in\mathbb{Z}_{p}
\end{equation}
As a result, given any Fourier transform $\hat{\chi}$ of the quasi-integrable
function $\chi$, the function $\hat{\chi}+\hat{\mu}$ will \emph{also
}be a Fourier transform of $\chi$ for any degenerate measure $d\mu$.
This prevents $\hat{\chi}$ from being unique.}\textbf{ }of $\chi$. Finally, given a Fourier transform $\hat{\chi}$
of a quasi-integrable function $\chi$, we write $\tilde{\chi}_{N}$
to denote the $N$th partial Fourier series generated by $\hat{\chi}$.
Interested readers can continue with this subsection to learn more
(although, in that case, they should probably read Subsection \ref{subsec:3.3.3 Frames},
assuming they haven't already done so). Otherwise, the reader can
safely proceed to the next subsection.

\vphantom{}
\end{rem}
HERE BEGINS THE DISCUSSION OF QUASI-INTEGRABILITY.

\vphantom{}

Given the issues discussed in Subsection \ref{subsec:3.3.3 Frames},
it would be a bit unreasonable to expect a comprehensive definition
of the phenomenon of quasi-integrability given the current novelty
of $\left(p,q\right)$-adic analysis. As such, the definition I am
about to give should come with an asterisk attached; it is in desperate
need of deeper examination, especially regarding interactions with
the many types of measures described in \ref{subsec:3.3.3 Frames}. 
\begin{defn}
Consider a $p$-adic frame $\mathcal{F}$.

We say a function $\chi\in C\left(\mathcal{F}\right)$ is\textbf{
quasi-integrable} \textbf{with respect to $\mathcal{F}$} / $\mathcal{F}$\textbf{-quasi-integrable
}whenever $\chi:U_{\mathcal{F}}\rightarrow I\left(\mathcal{F}\right)$
is the derivative of some $\mathcal{F}$-rising measure $d\mu\in M_{\textrm{rise}}\left(\mathcal{F}\right)$.
That is, there is a measure $d\mu$ whose Fourier-Stieltjes transform
$\hat{\mu}$ satisfies:

\vphantom{}

I. $\hat{\mu}\left(t\right)\in\overline{\mathbb{Q}}$, $\forall t\in\hat{\mathbb{Z}}_{p}$;

\vphantom{}

II. $\left\Vert \hat{\mu}\right\Vert _{p,K_{\mathfrak{z}}}<\infty$
all $\mathfrak{z}\in U_{\textrm{non}}\left(\mathcal{F}\right)$;

\vphantom{}

II. $\tilde{\mu}_{N}$ $\mathcal{F}$-converges to $\chi$: 
\begin{equation}
\lim_{N\rightarrow\infty}\left|\chi\left(\mathfrak{z}\right)-\tilde{\mu}_{N}\left(\mathfrak{z}\right)\right|_{K_{\mathfrak{z}}}\overset{\mathbb{R}}{=}0,\textrm{ }\forall\mathfrak{z}\in U_{\mathcal{F}}\label{eq:definition of quasi-integrability}
\end{equation}

We call any $d\mu$ satisfying these properties an\textbf{ $\mathcal{F}$-quasi-integral}\index{mathcal{F}-@$\mathcal{F}$-!quasi-integral}\textbf{
}of $\chi$, dropping the $\mathcal{F}$ when there is no confusion.
More generally, we say a function $\chi$ is quasi-integrable\index{quasi-integrability}\textbf{
}if it is quasi-integrable with respect to some frame $\mathcal{F}$.
We then write \nomenclature{$\tilde{L}^{1}\left(\mathcal{F}\right)$}{set of $\mathcal{F}$-quasi-integrable functions}$\tilde{L}^{1}\left(\mathcal{F}\right)$
to denote the set of all $\mathcal{F}$-quasi-integrable functions.
(Note that this is a vector space over $\overline{\mathbb{Q}}$.)

When $\mathcal{F}$ is the standard $\left(p,q\right)$-adic frame,
this definition becomes:

\vphantom{}

i. $\chi$ is a function from $\mathbb{Z}_{p}\rightarrow\mathbb{C}_{q}$
so that $\chi\left(\mathfrak{z}\right)\in\mathbb{C}$ for all $\mathfrak{z}\in\mathbb{N}_{0}$
and $\chi\left(\mathfrak{z}\right)\in\mathbb{C}_{q}$ for all $\mathfrak{z}\in\mathbb{Z}_{p}^{\prime}$.

\vphantom{}

ii. For each $\mathfrak{z}\in\mathbb{N}_{0}$, as $N\rightarrow\infty$,
$\tilde{\mu}_{N}\left(\mathfrak{z}\right)$ converges to $\chi\left(\mathfrak{z}\right)$
in the topology of $\mathbb{C}$.

\vphantom{}

iii. For each $\mathfrak{z}\in\mathbb{Z}_{p}^{\prime}$, as $N\rightarrow\infty$,
$\tilde{\mu}_{N}\left(\mathfrak{z}\right)$ converges to $\chi\left(\mathfrak{z}\right)$
in the topology of $\mathbb{C}_{q}$.

\vphantom{}

We then write \nomenclature{$\tilde{L}^{1}\left(\mathcal{F}_{p,q}\right)$}{ }$\tilde{L}^{1}\left(\mathcal{F}_{p,q}\right)$
and $\tilde{L}^{1}\left(\mathbb{Z}_{p},\mathbb{C}_{q}\right)$ \nomenclature{$\tilde{L}^{1}\left(\mathbb{Z}_{p},\mathbb{C}_{q}\right)$}{ }
to denote the space of all functions $\mathbb{Z}_{p}\rightarrow\mathbb{C}_{q}$
which are quasi-integrable with respect to the standard frame $\left(p,q\right)$-adic
frame.
\end{defn}
\vphantom{}

Using quasi-integrability, we can define a generalization of the $\left(p,q\right)$-adic
Fourier transform. The proposition which guarantees this is as follows: 
\begin{prop}
Let $\chi\in\tilde{L}^{1}\left(\mathcal{F}\right)$. Then, the quasi-integral
of $\chi$ is unique modulo $M_{\textrm{dgen}}\left(\mathcal{F}\right)$.
That is, if two measures $d\mu$ and $d\nu$ are both quasi-integrals
of $\chi$, then the measure $d\mu-d\nu$ is $\mathcal{F}$-degenerate\index{measure!mathcal{F}-degenerate@$\mathcal{F}$-degenerate}\index{mathcal{F}-@$\mathcal{F}$-!degenerate}.
Equivalently, we have an isomorphism of $\overline{\mathbb{Q}}$-linear
spaces: 
\begin{equation}
\tilde{L}^{1}\left(\mathcal{F}\right)\cong M_{\textrm{rise}}\left(\mathcal{F}\right)/M_{\textrm{dgen}}\left(\mathcal{F}\right)\label{eq:L1 twiddle isomorphism to M alg / M dgen}
\end{equation}
where the isomorphism associates a given $\chi\in\tilde{L}^{1}\left(\mathcal{F}\right)$
to the equivalence class of $\mathcal{F}$-rising measures which are
quasi-integrals of $\chi$. 
\end{prop}
Proof: Let $\mathcal{F}$, $\chi$, $d\mu$, and $d\nu$ be as given.
Then, by the definition of what it means to be a quasi-integral of
$\chi$, for the measure $d\eta\overset{\textrm{def}}{=}d\mu-d\nu$,
we have that $\tilde{\eta}_{N}\left(\mathfrak{z}\right)$ tends to
$0$ in $K_{\mathfrak{z}}$ for all $\mathfrak{z}\in U_{\mathcal{F}}$.
Hence, $d\mu-d\nu$ is indeed $\mathcal{F}$-degenerate.

Q.E.D. 
\begin{defn}
Let $\chi$ be quasi-integrable with respect to $\mathcal{F}$. Then,
\textbf{a} \textbf{Fourier Transform} of $\chi$ is a function $\hat{\mu}:\hat{\mathbb{Z}}_{p}\rightarrow\overline{\mathbb{Q}}$
which is the Fourier-Stieltjes transform of some $\mathcal{F}$-quasi-integral
$d\mu$ of $\chi$. We write these functions as $\hat{\chi}$ and
write the associated measures as $\chi\left(\mathfrak{z}\right)d\mathfrak{z}$. 
\end{defn}
\vphantom{}

Because there exist degenerate measures which are not the zero measure,
our use of the notation $\chi\left(\mathfrak{z}\right)d\mathfrak{z}$
to denote a given quasi-integral $d\mu\left(\mathfrak{z}\right)$
of $\chi$ is \emph{not} well-defined, being dependent on our particular
choice of $\hat{\chi}$. So, provided we are working with a fixed
choice of $\hat{\chi}$, this then allows us to define the integral
$\int_{\mathbb{Z}_{p}}f\left(\mathfrak{z}\right)\chi\left(\mathfrak{z}\right)d\mathfrak{z}$
for any $f\in C\left(\mathbb{Z}_{p},\mathbb{C}_{q}\right)$ by interpreting
it as the image of $f$ under integration against the measure $\chi\left(\mathfrak{z}\right)d\mathfrak{z}$:
\begin{equation}
\int_{\mathbb{Z}_{p}}f\left(\mathfrak{z}\right)\chi\left(\mathfrak{z}\right)d\mathfrak{z}\overset{\mathbb{C}_{q}}{=}\sum_{t\in\hat{\mathbb{Z}}_{p}}\hat{f}\left(-t\right)\hat{\chi}\left(t\right),\textrm{ }\forall f\in C\left(\mathbb{Z}_{p},\mathbb{C}_{q}\right)\label{eq:Integral of a quasi-integrable function}
\end{equation}

More pedantically, we could say that ``the'' Fourier transform of
$\chi$ is an equivalence class of functions $\hat{\mathbb{Z}}_{p}\rightarrow\overline{\mathbb{Q}}$
so that the difference of any two functions in the equivalence class
is the Fourier-Stieltjes transform of an $\mathcal{F}$-degenerate
algebraic measure.

Additionally, because the Fourier transform of $\chi$ is unique only
modulo degenerate measures, there is (at least at present) no \emph{canonical
}choice for what $\hat{\chi}$ should be. In this dissertation, the
symbol $\hat{\chi}$ will be used in one of two senses: either to
designate an \emph{arbitrary }choice of a Fourier transform for a
quasi-integrable function $\chi$\textemdash as is the case in Chapters
3 \& 5\textemdash OR, to designate a \emph{specific }choice of a Fourier
transform for a quasi-integrable function $\chi$, as is the case
in Chapters 4 \& 6. 
\begin{defn}
Given a choice of Fourier transform $\hat{\chi}$ for \nomenclature{$\tilde{\chi}_{N}\left(\mathfrak{z}\right)$}{$N$th partial sum of Fourier series generated by $\hat{\chi}$.}
$\chi\in\tilde{L}^{1}\left(\mathcal{F}\right)$, we write: 
\begin{equation}
\tilde{\chi}_{N}\left(\mathfrak{z}\right)\overset{\textrm{def}}{=}\sum_{\left|t\right|_{p}\leq p^{N}}\hat{\chi}\left(t\right)e^{2\pi i\left\{ t\mathfrak{z}\right\} _{p}},\textrm{ }\forall\mathfrak{z}\in\mathbb{Z}_{p}\label{eq:Definition of Chi_N twiddle}
\end{equation}
to denote the \textbf{$N$th partial sum of the Fourier series generated
by $\hat{\chi}$}.
\end{defn}
\vphantom{}

Working with quasi-integrable functions in the abstract is particularly
difficult because it is not at all clear how to actually \emph{prove
}anything without being able to draw upon identities like the Fourier
Resummation Lemmata. Moreover, for \emph{specific }quasi-integrable
functions\textemdash those given by a particular formula\textemdash one
can often easily prove many results which seem nearly impossible to
address in the abstract. Conversely, those results which \emph{can
}be proven in the abstract are generally unimpressive; case in point: 
\begin{prop}
\label{prop:Chi_N and Chi_twiddle N F converge to one another}Let
$\chi\in\tilde{L}^{1}\left(\mathcal{F}\right)$. Then, for any Fourier
transform $\hat{\chi}$ of $\chi$, the difference $\chi_{N}-\tilde{\chi}_{N}$
(where $\chi_{N}$ is the $N$th truncation of $\chi$) $\mathcal{F}$-converges
to $0$. 
\end{prop}
Proof: Immediate from the definitions of $\hat{\chi}$, $\tilde{\chi}_{N}$,
and $\chi_{N}$.

Q.E.D.

\vphantom{}

Both because it is not my purpose to do so here, and as a simple matter
of time constraints, there are quite a few issues in the theory of
quasi-integrable functions that I have not yet been able to resolve
to my satisfaction. Below, I have outlined what I feel are the most
outstanding issues. 
\begin{question}
\label{que:3.4}Fix a $\left(p,q\right)$-adic frame $\mathcal{F}$.
Given a function $\chi:\mathbb{Z}_{p}\rightarrow\mathbb{C}_{q}$,
is there a practical criterion for determining whether or not $\chi$
is $\mathcal{F}$-quasi-integrable? Similarly, given a function $\hat{\chi}:\hat{\mathbb{Z}}_{p}\rightarrow\mathbb{C}_{q}$,
is there a practical criterion for determining whether or not $\hat{\chi}$
is the Fourier transform of an $\mathcal{F}$-quasi-integrable function
(other than directly summing the Fourier series generated by $\hat{\chi}$)? 
\end{question}
\vphantom{}

Question \ref{que:3.4} probably bothers me more than any other outstanding
issue\footnote{The verification of when the convolution of quasi-integrable functions
is, itself, quasi-integrable comes in at a close second.} in this fledgling theory. At present, the only means I have of verifying
quasi-integrability is by performing the rather intensive computational
regimen of Chapters 4 and 6 on a given $\chi$ in the hopes of being
able to find a function $\hat{\chi}$ which serves as the Fourier
transform of $\chi$ with respect to some frame. Arguably even worse
is the fact that there currently doesn't appear to be any reasonable
criterion for showing that a function is \emph{not }quasi-integrable.
In particular, despite being nearly certain that the following conjecture
is true, I have no way to prove it: 
\begin{conjecture}
Let $\mathfrak{a}\in\mathbb{Z}_{p}$. Then, the one-point function
$\mathbf{1}_{\mathfrak{a}}\left(\mathfrak{z}\right)$ is not quasi-integrable
with respect to any $\left(p,q\right)$-adic quasi-integrability frame. 
\end{conjecture}
\vphantom{}

The methods of Chapters 4 and 6 do not seem applicable here. These
methods are best described as a kind of ``asymptotic analysis''
in which we truncate the function $\chi$ being studied and then explicitly
compute the Fourier transform $\hat{\chi}_{N}\left(t\right)$ of the
locally constant (and hence, continuous) $\left(p,q\right)$-adic
function $\chi_{N}\left(\mathfrak{z}\right)$. Careful computational
manipulations allow us to ``untangle'' the dependence of $N$ and
$t$ in $\hat{\chi}_{N}\left(t\right)$ and to express $\hat{\chi}_{N}\left(t\right)$
in the form: 
\begin{equation}
\hat{\chi}_{N}\left(t\right)=\hat{F}_{N}\left(t\right)+\mathbf{1}_{0}\left(p^{N}t\right)\hat{G}\left(t\right)
\end{equation}
where $\hat{F}_{N},\hat{G}:\hat{\mathbb{Z}}_{p}\rightarrow\overline{\mathbb{Q}}$
are $q$-adically bounded, and with the Fourier series generated by
$\hat{F}_{N}$ converging to zero in the standard frame as $N\rightarrow\infty$.
$\hat{G}\left(t\right)$ can be thought of as the ``fine structure''
hidden beneath the tumult of the divergent $\hat{F}_{N}\left(t\right)$
term. This suggests that $\hat{G}$ should be a Fourier transform
of $\hat{\chi}$\textemdash and, with some work, this intuition will
be proved correct. Even though the Fourier series generated by $\hat{F}_{N}$
only $\mathcal{F}_{p,q}$-converges to $0$ when $\chi$ satisfies
certain arithmetical properties\footnote{For example, of the $\chi_{q}$s (the numina of the shortened $qx+1$
maps), the $\hat{F}_{N}$ term's Fourier series is $\mathcal{F}_{2,q}$-convergent
to $0$ if and only if $q=3$.}, we can bootstrap this special case and use it to tackle the more
general case where the Fourier series generated by $\hat{F}_{N}$
\emph{doesn't} $\mathcal{F}_{p,q}$-converge to $0$.

Applying this procedure to, say, $\mathbf{1}_{0}\left(\mathfrak{z}\right)$,
the $N$th truncation of $\mathbf{1}_{0}\left(\mathfrak{z}\right)$
is: 
\begin{equation}
\mathbf{1}_{0,N}\left(\mathfrak{z}\right)=\left[\mathfrak{z}\overset{p^{N}}{\equiv}0\right]\label{eq:Nth truncation of 1_0}
\end{equation}
which has the Fourier transform: 
\begin{equation}
\hat{\mathbf{1}}_{0,N}\left(t\right)=\frac{1}{p^{N}}\mathbf{1}_{0}\left(p^{N}t\right)=\begin{cases}
\frac{1}{p^{N}} & \textrm{if }\left|t\right|_{p}\leq p^{N}\\
0 & \textrm{if }\left|t\right|_{p}>p^{N}
\end{cases}\label{eq:Fourier transform of Nth truncation of 1_0}
\end{equation}
On the one hand, there appears to be no way to replicate Chapter 4's
argument to obtain a candidate for a Fourier transform of $\mathbf{1}_{0}$
and thereby establish its quasi-integrability. But how might we make
this \emph{rigorous}? At present, I have no idea. Only time will tell
if it will stay that way.

The next major issue regards the convolution of quasi-integrable functions. 
\begin{question}
Fix a frame $\mathcal{F}$, and let $\chi,\kappa\in\tilde{L}^{1}\left(\mathcal{F}\right)$.
The natural definition for the \textbf{convolution }of $\chi$ and
$\kappa$ would be: 
\begin{equation}
\lim_{N\rightarrow\infty}\sum_{\left|t\right|_{p}\leq p^{N}}\hat{\chi}\left(t\right)\hat{\kappa}\left(t\right)e^{2\pi i\left\{ t\mathfrak{z}\right\} _{p}}\label{eq:Limit formula for the convolution of two quasi-integrable functions}
\end{equation}
for a choice $\hat{\chi}$ and $\hat{\kappa}$ of Fourier transforms
of $\chi$ and $\kappa$, respectively. However, this immediately
leads to difficulties:

(1) As defined, $\hat{\chi}$ and $\hat{\kappa}$ are unique only
up to $\mathcal{F}$-degenerate measures. However, as is shown in
\emph{(\ref{eq:Limit of Fourier sum of v_p A_H hat}) }from \emph{Chapter
4's} \textbf{\emph{Lemma \ref{lem:v_p A_H hat summation formulae}}}
\emph{(page \pageref{lem:v_p A_H hat summation formulae})}, there
exists a degenerate measure $d\mu$ and a bounded $\hat{\nu}:\hat{\mathbb{Z}}_{p}\rightarrow\overline{\mathbb{Q}}$
so that: 
\[
\lim_{N\rightarrow\infty}\sum_{\left|t\right|_{p}\leq p^{N}}\hat{\nu}\left(t\right)\hat{\mu}\left(t\right)e^{2\pi i\left\{ t\mathfrak{z}\right\} _{p}}
\]
is \textbf{\emph{not}}\emph{ }$\mathcal{F}$-convergent to $0$ everywhere.
As a result, it appears that \emph{(\ref{eq:Limit formula for the convolution of two quasi-integrable functions})}
will not be well defined modulo $\mathcal{F}$-degenerate measures;
that is, replacing $\hat{\chi}$ in \emph{(\ref{eq:Limit formula for the convolution of two quasi-integrable functions})}
with $\hat{\chi}+\hat{\mu}$ where $d\mu$ is $\mathcal{F}$-degenerate
may cause the $\mathcal{F}$-limit of \emph{(\ref{eq:Limit formula for the convolution of two quasi-integrable functions})}
to change. This prompts the question: there a way to \uline{uniquely}
define the convolution of two quasi-integrable functions? Or will
we have to make everything contingent upon the particular choice of
Fourier transforms in order to work with quasi-integrable functions?

(2) As an extension of the present lack of a criterion for determining
whether or not a given $\chi$ is quasi-integrable (or whether a given
$\hat{\chi}$ is a Fourier transform of a quasi-integrable function)
are there any useful conditions we can impose on $\chi$, $\kappa$,
$\hat{\chi}$, and $\hat{\kappa}$ to ensure that \emph{(\ref{eq:Limit formula for the convolution of two quasi-integrable functions})}
exists and is quasi-integrable? The one obvious answer is when one
or both of $\hat{\chi}$ or $\hat{\kappa}$ is an element of $c_{0}\left(\hat{\mathbb{Z}}_{p},\mathbb{C}_{q}\right)$,
in which case $\chi*\kappa$ will then be a continuous function, due
to the $q$-adic decay of $\hat{\chi}\cdot\hat{\kappa}$. However,
that would be quite an unsatisfying state of affairs! Worse, the computations
done in \emph{Chapter 4} involve convolving quasi-integrable functions;
for example \emph{(\ref{eq:Limit of Fourier sum of v_p A_H hat} }from
\textbf{\emph{Lemma \ref{lem:v_p A_H hat summation formulae}}} on
page \emph{\pageref{lem:v_p A_H hat summation formulae})} arises
in the analysis of $\chi_{3}$, in which we convolve $dA_{3}$ with
the quasi-integrable function obtained by letting $N\rightarrow\infty$
in \emph{(\ref{eq:Fourier sum of v_p of t})}. So, there are clearly
cases where convolution \emph{\uline{is}}\emph{ }valid! 
\end{question}
\begin{question}
\label{que:3.6}How, if at all, can we work with functions $\chi$
and $\kappa$ which are quasi-integrable with respect to \uline{different
frames}? 
\end{question}

\subsection{\label{subsec:3.3.6 L^1 Convergence}$L^{1}$ Convergence}

When doing real or complex analysis, few settings are as fundamental
as $L^{p}$ spaces, with $L^{1}$ and $L^{2}$ possessing a mutual
claim to primacy. As\index{non-archimedean!$L^{1}$ theory} was mentioned
in Subsection \ref{subsec:3.1.2 Banach-Spaces-over}, however, neither
$L^{1}$ nor $L^{2}$\textemdash nor any $L^{p}$ space, for that
matter\textemdash have much significance when working with functions
taking values in a non-archimedean valued field. The magic of the
complex numbers is that they are a finite-dimensional field extension
of $\mathbb{R}$ which is nonetheless algebraically complete. Because
$\mathbb{C}_{q}$ is an infinite dimensional extension of $\mathbb{Q}_{q}$,
the Galois-theoretic identity $z\overline{z}=\left|z\right|^{2}$
from complex analysis does not hold, and the elegant duality of $L^{2}$
spaces gets lost as a result in the $q$-adic context.

From an analytical standpoint, the issue with non-archimedean $L^{1}$
spaces is arguably even more fundamental. For infinite series of complex
numbers, we need absolute convergence in order to guarantee the invariance
of the sum under groupings and rearrangements of its terms. In the
ultrametric setting, however, this is no longer the case: all that
matters is that the $n$th term of the series tend to $0$ in non-archimedean
absolute value. Because of this, as we saw in our discussion of Monna-Springer
Integration in Subsection \ref{subsec:3.1.6 Monna-Springer-Integration},
when working with non-archimedean absolute values, there is no longer
a meaningful connection between $\left|\int f\right|$ and $\int\left|f\right|$.

Despite this, it appears that a non-archimedean $L^{1}$ theory \emph{is}
implicated in $\left(p,q\right)$-adic analysis. Rather than begin
with theory or definitions, let us start with a concrete, Collatz-adjacent
example involving $\chi_{5}$\index{chi{5}@$\chi_{5}$}, the numen
of \index{$5x+1$ map} the shortened $5x+1$ map. 
\begin{example}
For a given $\mathcal{F}$-quasi-integrable $\chi$, \textbf{Proposition
\ref{prop:Chi_N and Chi_twiddle N F converge to one another}} tells
us that $\chi_{N}-\tilde{\chi}_{N}$ is $\mathcal{F}$-convergent
to zero. Nevertheless, for any given choice of $\hat{\chi}$, the
limit of the differences $\hat{\chi}_{N}\left(t\right)-\hat{\chi}\left(t\right)$
as $N\rightarrow\infty$ need not be well-behaved. As we will see
in Chapter 4, $\chi_{5}\in\tilde{L}^{1}\left(\mathbb{Z}_{2},\mathbb{C}_{5}\right)$
(that is, $\chi_{5}$ is the derivative of a rising-continuous $\left(2,5\right)$-adic
measure); moreover, one possible choice for a Fourier transform of
$\chi_{5}$ is: 
\begin{equation}
\hat{\chi}_{5}\left(t\right)=-\frac{1}{4}\mathbf{1}_{0}\left(t\right)-\frac{1}{4}\hat{A}_{5}\left(t\right),\textrm{ }\forall t\in\hat{\mathbb{Z}}_{2}\label{eq:Choice of Fourier Transform for Chi_5}
\end{equation}
where $\hat{A}_{5}$ is as defined in (\ref{eq:Definition of A_q hat}),
and where $\mathbf{1}_{0}$ is the indicator function of $0$ in $\hat{\mathbb{Z}}_{2}$.

Letting $\hat{\chi}_{5,N}$ denote the Fourier transform of the $N$th
truncation of $\chi_{5}$ ($\chi_{5,N}$) it can be shown that for
this choice of $\hat{\chi}_{5}$: 
\begin{equation}
\hat{\chi}_{5,N}\left(t\right)-\hat{\chi}_{5}\left(t\right)=\begin{cases}
\frac{1}{2}\left(\frac{3}{2}\right)^{N+\min\left\{ v_{2}\left(t\right),0\right\} }\hat{A}_{5}\left(t\right) & \textrm{if }\left|t\right|_{2}\leq2^{N}\\
\frac{1}{4}\hat{A}_{5}\left(t\right) & \textrm{if }\left|t\right|_{2}>2^{N}
\end{cases},\textrm{ }\forall t\in\hat{\mathbb{Z}}_{2}\label{eq:Difference between Chi_5,N hat and Chi_5 hat}
\end{equation}
For any given $t$, as $N\rightarrow\infty$, the limit of this difference
exists neither in $\mathbb{C}_{5}$ nor in $\mathbb{C}$, despite
the fact that the difference of $\chi_{5,N}$ and $\tilde{\chi}_{5,N}$
(the $N$th partial sum of the Fourier series generated by $\hat{\chi}_{5}$)
is $\mathcal{F}$-convergent to $0$. Nevertheless, going through
the motions and applying Fourier inversion to the above (multiplying
by $e^{2\pi i\left\{ t\mathfrak{z}\right\} _{2}}$ and summing over
$\left|t\right|_{2}\leq2^{N}$), we obtain: 
\begin{equation}
\chi_{5,N}\left(\mathfrak{z}\right)-\tilde{\chi}_{5,N}\left(\mathfrak{z}\right)=\frac{5^{\#_{1}\left(\left[\mathfrak{z}\right]_{2^{N}}\right)}}{2^{N+1}}\label{eq:Fourier inversion of Chi_5,N hat - Chi_5 hat}
\end{equation}
which, as expected, converges to zero $5$-adically as $N\rightarrow\infty$
for all $\mathfrak{z}\in\mathbb{Z}_{2}^{\prime}$ and converges to
zero in $\mathbb{C}$ for all $\mathfrak{z}\in\mathbb{N}_{0}$. Note
that this convergence is \emph{not }uniform over $\mathbb{Z}_{2}^{\prime}$:
we can make the rate of convergence with respect to $N$ arbitrarily
slow by choosing a $\mathfrak{z}\in\mathbb{Z}_{2}\backslash\mathbb{N}_{0}$
whose $1$s digits in its $2$-adic expansion are spaced sufficiently
far apart.

Despite the lack of the \emph{uniform} convergence of $\chi_{5,N}\left(\mathfrak{z}\right)-\tilde{\chi}_{5,N}\left(\mathfrak{z}\right)$
to zero, we can get $L^{1}$ convergence for certain cases, and in
the classical sense, no less. We have: 
\begin{align*}
\int_{\mathbb{Z}_{2}}\left|\chi_{5,N}\left(\mathfrak{z}\right)-\tilde{\chi}_{5,N}\left(\mathfrak{z}\right)\right|_{5}d\mathfrak{z} & \overset{\mathbb{R}}{=}\int_{\mathbb{Z}_{2}}\left(\frac{1}{5}\right)^{\#_{1}\left(\left[\mathfrak{z}\right]_{2^{N}}\right)}d\mathfrak{z}\\
 & =\frac{1}{2^{N}}\sum_{n=0}^{2^{N}-1}\left(\frac{1}{5}\right)^{\#_{1}\left(n\right)}\\
 & =\frac{1}{2^{N}}\underbrace{\prod_{n=0}^{N-1}\left(1+\frac{1}{5}\left(1^{2^{n}}\right)\right)}_{\left(6/5\right)^{N}}\\
 & =\left(\frac{3}{5}\right)^{N}
\end{align*}
which tends to $0$ in $\mathbb{R}$ as $N\rightarrow\infty$.

By using the non-trivial explicit formula for $\chi_{5}$ and $\chi_{5,N}$
from Subsection \ref{subsec:4.2.2}'s\textbf{ Corollary \ref{cor:F-series for Chi_H for p equals 2}}
and \textbf{Theorem \ref{thm:F-series for Nth truncation of Chi_H, alpha is 1}},\textbf{
}we obtain: 
\begin{equation}
\chi_{5}\left(\mathfrak{z}\right)-\chi_{5,N}\left(\mathfrak{z}\right)=-\frac{1}{4}\frac{5^{\#_{1}\left(\left[\mathfrak{z}\right]_{2^{N}}\right)}}{2^{N}}+\frac{1}{8}\sum_{n=N}^{\infty}\frac{5^{\#_{1}\left(\left[\mathfrak{z}\right]_{2^{n}}\right)}}{2^{n}}\label{eq:Difference between Chi_5 and Chi_5,N}
\end{equation}
where the convergence is in $\mathbb{C}$ for $\mathfrak{z}\in\mathbb{N}_{0}$
and is in $\mathbb{C}_{5}$ for $\mathfrak{z}\in\mathbb{Z}_{2}^{\prime}$.
Again, regardless of the topology, the convergence is point-wise.
Nevertheless, using the standard $5$-adic absolute value estimate
on the right-hand side of (\ref{eq:Difference between Chi_5 and Chi_5,N}),
we have: 
\begin{align*}
\int_{\mathbb{Z}_{2}}\left|\chi_{5}\left(\mathfrak{z}\right)-\chi_{5,N}\left(\mathfrak{z}\right)\right|_{5}d\mathfrak{z} & \overset{\mathbb{R}}{\leq}\int_{\mathbb{Z}_{2}}\max_{n\geq N}\left(5^{-\#_{1}\left(\left[\mathfrak{z}\right]_{2^{n}}\right)}\right)d\mathfrak{z}\\
 & =\int_{\mathbb{Z}_{2}}\left(\frac{1}{5}\right)^{\#_{1}\left(\left[\mathfrak{z}\right]_{2^{N}}\right)}d\mathfrak{z}\\
 & =\left(\frac{3}{5}\right)^{N}
\end{align*}
which, once again, tends to $0$ in $\mathbb{R}$ as $N\rightarrow\infty$.

Finally, when working with the $\left(2,5\right)$-adic integral of
the difference, we obtain: 
\begin{align*}
\int_{\mathbb{Z}_{2}}\left(\chi_{5,N}\left(\mathfrak{z}\right)-\tilde{\chi}_{5,N}\left(\mathfrak{z}\right)\right)d\mathfrak{z} & \overset{\mathbb{C}_{5}}{=}\int_{\mathbb{Z}_{2}}\frac{5^{\#_{1}\left(\left[\mathfrak{z}\right]_{2^{N}}\right)}}{2^{N+1}}d\mathfrak{z}\\
 & \overset{\mathbb{C}_{5}}{=}\frac{1}{2^{N}}\sum_{n=0}^{2^{N}-1}\frac{5^{\#_{1}\left(n\right)}}{2^{N+1}}\\
 & \overset{\mathbb{C}_{5}}{=}\frac{1/2}{4^{N}}\prod_{n=0}^{N-1}\left(1+5\left(1^{2^{n}}\right)\right)\\
 & \overset{\mathbb{C}_{5}}{=}\frac{1}{2}\left(\frac{3}{2}\right)^{N}
\end{align*}
which, just like $\hat{\chi}_{5,N}\left(t\right)-\hat{\chi}_{5}\left(t\right)$,
fails to converge $5$-adically as $N\rightarrow\infty$. Indeed,
we actually have that both of these quantities fail to converge $5$-adically
in precisely the same way: 
\[
\left(\frac{2}{3}\right)^{N}\int_{\mathbb{Z}_{2}}\left(\chi_{5,N}\left(\mathfrak{z}\right)-\tilde{\chi}_{5,N}\left(\mathfrak{z}\right)\right)d\mathfrak{z}\overset{\mathbb{C}_{5}}{=}\frac{1}{2},\textrm{ }\forall N\geq1
\]
\[
\lim_{N\rightarrow\infty}\left(\frac{2}{3}\right)^{N}\left(\hat{\chi}_{5,N}\left(t\right)-\hat{\chi}_{5}\left(t\right)\right)\overset{\mathbb{C}_{5}}{=}\frac{1}{2}\left(\frac{3}{2}\right)^{N+\min\left\{ v_{2}\left(t\right),0\right\} }\hat{A}_{5}\left(t\right),\textrm{ }\forall t\in\hat{\mathbb{Z}}_{2}
\]
and hence: 
\[
\lim_{N\rightarrow\infty}\frac{\hat{\chi}_{5,N}\left(t\right)-\hat{\chi}_{5}\left(t\right)}{\int_{\mathbb{Z}_{2}}\left(\chi_{5,N}\left(\mathfrak{z}\right)-\tilde{\chi}_{5,N}\left(\mathfrak{z}\right)\right)d\mathfrak{z}}\overset{\mathbb{C}_{5}}{=}\left(\frac{3}{2}\right)^{\min\left\{ v_{2}\left(t\right),0\right\} }\hat{A}_{5}\left(t\right),\textrm{ }\forall t\in\hat{\mathbb{Z}}_{2}
\]
\end{example}
\vphantom{}

At present, the stringent restrictions demanded by $\mathcal{F}$-convergence
seems to prevent us from equipping spaces such as $M\left(\mathcal{F}\right)$
or $\tilde{L}^{1}\left(\mathcal{F}\right)$ with a norm to induce
a topology. Non-archimedean $L^{1}$ theory\index{non-archimedean!$L^{1}$ theory}
might very well provided a viable work-around. 
\begin{defn}
We write$\mathcal{L}_{\mathbb{R}}^{1}\left(\mathbb{Z}_{p},\mathbb{C}_{q}\right)$
to denote the set of all functions $f:\mathbb{Z}_{p}\rightarrow\mathbb{C}_{q}$
so that the real-valued function $\mathfrak{z}\in\mathbb{Z}_{p}\mapsto\left|f\left(\mathfrak{z}\right)\right|_{q}\in\mathbb{R}$
is integrable with respect to the real-valued Haar probability measure
on $\mathbb{Z}_{p}$. As is customary, we then define an equivalence
relation $\sim$ on $\mathcal{L}_{\mathbb{R}}^{1}\left(\mathbb{Z}_{p},\mathbb{C}_{q}\right)$,
and say that $f\sim g$ whenever $\left|f\left(\mathfrak{z}\right)-g\left(\mathfrak{z}\right)\right|_{q}=0$
for almost every $\mathfrak{z}\in\mathbb{Z}_{p}$. We then write $L_{\mathbb{R}}^{1}\left(\mathbb{Z}_{p},\mathbb{C}_{q}\right)$\nomenclature{$L_{\mathbb{R}}^{1}\left(\mathbb{Z}_{p},\mathbb{C}_{q}\right)$}{set of $f:\mathbb{Z}_{p}\rightarrow\mathbb{C}_{q}$ so that $\left|f\right|_{q}\in L^{1}\left(\mathbb{Z}_{p},\mathbb{C}\right)$}
to denote the space of equivalence classes:
\begin{equation}
L_{\mathbb{R}}^{1}\left(\mathbb{Z}_{p},\mathbb{C}_{q}\right)\overset{\textrm{def}}{=}\mathcal{L}_{\mathbb{R}}^{1}\left(\mathbb{Z}_{p},\mathbb{C}_{q}\right)/\sim\label{eq:Definition of L^1_R}
\end{equation}
and we make this an\emph{ archimedean}\footnote{That is, $\left\Vert \mathfrak{a}f+\mathfrak{b}g\right\Vert _{1}\leq\left|\mathfrak{a}\right|_{q}\left\Vert f\right\Vert _{1}+\left|\mathfrak{b}\right|_{q}\left\Vert g\right\Vert _{1}$
for all $f,g\in L_{\mathbb{R}}^{1}\left(\mathbb{Z}_{p},\mathbb{C}_{q}\right)$
and all $\mathfrak{a},\mathfrak{b}\in\mathbb{C}_{q}$.}\emph{ }normed linear space over $\mathbb{C}_{q}$ by defining: 
\begin{equation}
\left\Vert f\right\Vert _{1}\overset{\textrm{def}}{=}\int_{\mathbb{Z}_{p}}\left|f\left(\mathfrak{z}\right)\right|_{q}d\mathfrak{z}\in\mathbb{R},\textrm{ }\forall f\in L_{\mathbb{R}}^{1}\left(\mathbb{Z}_{p},\mathbb{C}_{q}\right)\label{eq:Definition of L^1 norm}
\end{equation}
\end{defn}
\begin{lem}
\label{lem:Egorov lemma}Let $\left\{ f_{n}\right\} _{n\geq1}$ be
a sequence of functions $f_{n}:\mathbb{Z}_{p}\rightarrow\mathbb{C}_{q}$
so that:

\vphantom{}

i. The real-valued functions $\mathfrak{z}\mapsto\left|f_{n}\left(\mathfrak{z}\right)\right|_{q}$
are measurable with respect to the real-valued Haar probability measure
on $\mathbb{Z}_{p}$;

\vphantom{}

ii. There is a constant $C\in\mathbb{R}>0$ so that for each $n$,
$\left|f_{n}\left(\mathfrak{z}\right)\right|_{q}<C$ holds almost
everywhere.

\vphantom{}

Then:

\vphantom{}

I. If $\lim_{n\rightarrow\infty}\left|f_{n}\left(\mathfrak{z}\right)\right|_{q}\overset{\mathbb{R}}{=}0$
holds for almost every $\mathfrak{z}\in\mathbb{Z}_{p}$, then $\lim_{n\rightarrow\infty}\int_{\mathbb{Z}_{p}}\left|f_{n}\left(\mathfrak{z}\right)\right|_{q}d\mathfrak{z}\overset{\mathbb{R}}{=}0$.

\vphantom{}

II. If $\lim_{n\rightarrow\infty}\int_{\mathbb{Z}_{p}}\left|f_{n}\left(\mathfrak{z}\right)\right|_{q}d\mathfrak{z}\overset{\mathbb{R}}{=}0$
then there is a subsequence $\left\{ f_{n_{j}}\right\} _{j\geq1}$
so that $\lim_{j\rightarrow0}\left|f_{n_{j}}\left(\mathfrak{z}\right)\right|_{q}\overset{\mathbb{R}}{=}0$
holds for almost every $\mathfrak{z}\in\mathbb{Z}_{p}$. 
\end{lem}
Proof:

I. Suppose the $\left|f_{n}\right|_{q}$s converge point-wise to $0$
almost everywhere. Then, letting $\epsilon>0$, by \textbf{Egorov's
Theorem}\footnote{Let $\left(X,\Omega,\mu\right)$ be a measure space with $\mu\left(X\right)<\infty$,
and let $f$ and $\left\{ f_{N}\right\} _{N\geq0}$ be measurable,
complex-valued functions on $X$ so that $f_{N}$ converges to $f$
almost everywhere. Then, for every $\epsilon>0$, there is a measurable
set $E\subseteq X$ with $\mu\left(E\right)<\epsilon$ such that $f_{N}\rightarrow f$
uniformly on $X\backslash E$. \cite{Folland - real analysis}}, there is a measurable set $E\subseteq\mathbb{Z}_{p}$ of real-valued
Haar measure $<\epsilon$ so that the $\left|f_{n}\right|_{q}$s converge
uniformly to $0$ on $\mathbb{Z}_{p}\backslash E$. 
\begin{align*}
\int_{\mathbb{Z}_{p}}\left|f_{n}\left(\mathfrak{z}\right)\right|_{q}d\mathfrak{z} & =\int_{E}\left|f_{n}\left(\mathfrak{z}\right)\right|_{q}d\mathfrak{z}+\int_{\mathbb{Z}_{p}\backslash E}\left|f_{n}\left(\mathfrak{z}\right)\right|_{q}d\mathfrak{z}\\
\left(\left|f_{n}\right|_{q}<C\textrm{ a.e.}\right); & \leq C\epsilon+\int_{\mathbb{Z}_{p}\backslash E}\left|f_{n}\left(\mathfrak{z}\right)\right|_{q}d\mathfrak{z}
\end{align*}
Since the $\left|f_{n}\right|_{q}$s converge uniformly to $0$ on
$\mathbb{Z}_{p}\backslash E$, $\int_{\mathbb{Z}_{p}\backslash E}\left|f_{n}\left(\mathfrak{z}\right)\right|_{q}d\mathfrak{z}$
tends to $0$ in $\mathbb{R}$ as $n\rightarrow\infty$, and we are
left with: 
\begin{equation}
\limsup_{n\rightarrow\infty}\int_{\mathbb{Z}_{p}}\left|f_{n}\left(\mathfrak{z}\right)\right|_{q}d\mathfrak{z}\ll\epsilon
\end{equation}
Since $\epsilon$ was arbitrary, we conclude that: 
\begin{equation}
\lim_{n\rightarrow\infty}\int_{\mathbb{Z}_{p}}\left|f_{n}\left(\mathfrak{z}\right)\right|_{q}d\mathfrak{z}=0
\end{equation}

\vphantom{}

II. This is a standard result of measure theory. See Chapter 2.4 in
(\cite{Folland - real analysis}) for details.

Q.E.D. 
\begin{thm}
\label{thm:L^1_R is an archimedean Banach space}$L_{\mathbb{R}}^{1}\left(\mathbb{Z}_{p},\mathbb{C}_{q}\right)$
is an archimedean Banach space over $\mathbb{C}_{q}$. 
\end{thm}
Proof: To show that the archimedean normed vector space $L_{\mathbb{R}}^{1}\left(\mathbb{Z}_{p},\mathbb{C}_{q}\right)$
is a Banach space, we use \textbf{Theorem 5.1 }from \cite{Folland - real analysis}
(the archimedean analogue of \textbf{Proposition \ref{prop:series characterization of a Banach space}}):
it suffices to show that every absolutely convergent series in $L_{\mathbb{R}}^{1}\left(\mathbb{Z}_{p},\mathbb{C}_{q}\right)$
converges to a limit in $L_{\mathbb{R}}^{1}\left(\mathbb{Z}_{p},\mathbb{C}_{q}\right)$.

To do this, let $\left\{ f_{n}\right\} _{n\geq0}$ be a sequence in
$L_{\mathbb{R}}^{1}\left(\mathbb{Z}_{p},\mathbb{C}_{q}\right)$ so
that: 
\begin{equation}
\sum_{n=0}^{\infty}\left\Vert f_{n}\right\Vert _{1}<\infty
\end{equation}
That is: 
\begin{equation}
\lim_{N\rightarrow\infty}\sum_{n=0}^{N-1}\int_{\mathbb{Z}_{p}}\left|f_{n}\left(\mathfrak{z}\right)\right|_{q}d\mathfrak{z}
\end{equation}
exists in $\mathbb{R}$. Consequently: 
\begin{equation}
\lim_{N\rightarrow\infty}\sum_{n=0}^{N-1}\int_{\mathbb{Z}_{p}}\left|f_{n}\left(\mathfrak{z}\right)\right|_{q}d\mathfrak{z}\overset{\mathbb{R}}{=}\int_{\mathbb{Z}_{p}}\lim_{N\rightarrow\infty}\sum_{n=0}^{N-1}\left|f_{n}\left(\mathfrak{z}\right)\right|_{q}d\mathfrak{z}
\end{equation}
and hence, the function defined by $\lim_{N\rightarrow\infty}\sum_{n=0}^{N-1}\left|f_{n}\left(\mathfrak{z}\right)\right|_{q}$
is in $L^{1}\left(\mathbb{Z}_{p},\mathbb{R}\right)$, and the point-wise
limit exists $\sum_{n=0}^{\infty}\left|f_{n}\left(\mathfrak{z}\right)\right|_{q}$
for almost every $\mathfrak{z}\in\mathbb{Z}_{p}$. By the ordinary
$q$-adic triangle inequality (not the ultrametric inequality, just
the ordinary triangle inequality), this then shows that $\sum_{n=0}^{\infty}f_{n}\left(\mathfrak{z}\right)$
converges in $\mathbb{C}_{q}$ for almost every $\mathfrak{z}\in\mathbb{Z}_{p}$.

As such, let $F:\mathbb{Z}_{p}\rightarrow\mathbb{C}_{q}$ be defined
by: 
\begin{equation}
F\left(\mathfrak{z}\right)\overset{\textrm{def}}{=}\begin{cases}
\sum_{n=0}^{\infty}f_{n}\left(\mathfrak{z}\right) & \textrm{if \ensuremath{\sum_{n=0}^{\infty}f_{n}\left(\mathfrak{z}\right)} converges in }\mathbb{C}_{q}\\
0 & \textrm{if \ensuremath{\sum_{n=0}^{\infty}f_{n}\left(\mathfrak{z}\right)} does not converge in }\mathbb{C}_{q}
\end{cases},\textrm{ }\forall\mathfrak{z}\in\mathbb{Z}_{p}
\end{equation}
Then: 
\begin{equation}
\left|F\left(\mathfrak{z}\right)-\sum_{n=0}^{N-1}f_{n}\left(\mathfrak{z}\right)\right|_{q}=\left|\sum_{n=N}^{\infty}f_{n}\left(\mathfrak{z}\right)\right|_{q}
\end{equation}
holds for almost every $\mathfrak{z}\in\mathbb{Z}_{p}$, and so: 
\begin{align*}
\lim_{N\rightarrow\infty}\left\Vert F-\sum_{n=0}^{N-1}f_{n}\right\Vert _{1} & \overset{\mathbb{R}}{=}\lim_{N\rightarrow\infty}\int_{\mathbb{Z}_{p}}\left|F\left(\mathfrak{z}\right)-\sum_{n=0}^{N-1}f_{n}\left(\mathfrak{z}\right)\right|_{q}d\mathfrak{z}\\
 & =\lim_{N\rightarrow\infty}\int_{\mathbb{Z}_{p}}\left|\sum_{n=N}^{\infty}f_{n}\left(\mathfrak{z}\right)\right|_{q}d\mathfrak{z}\\
 & \leq\lim_{N\rightarrow\infty}\int_{\mathbb{Z}_{p}}\sum_{n=N}^{\infty}\left|f_{n}\left(\mathfrak{z}\right)\right|_{q}d\mathfrak{z}\\
 & =\lim_{N\rightarrow\infty}\sum_{n=N}^{\infty}\left\Vert f_{n}\right\Vert _{1}\\
\left(\sum_{n=0}^{\infty}\left\Vert f_{n}\right\Vert _{1}<\infty\right); & =0
\end{align*}
which shows that $F$ is the limit of the series $\sum_{n=0}^{N-1}f_{n}\left(\mathfrak{z}\right)$
in $L_{\mathbb{R}}^{1}\left(\mathbb{Z}_{p},\mathbb{C}_{q}\right)$.
Thus, every absolutely convergent series in $L_{\mathbb{R}}^{1}\left(\mathbb{Z}_{p},\mathbb{C}_{q}\right)$
converges, which proves that $L_{\mathbb{R}}^{1}\left(\mathbb{Z}_{p},\mathbb{C}_{q}\right)$
is a Banach space.

Q.E.D. 
\begin{prop}
Let $\chi:\mathbb{Z}_{p}\rightarrow\mathbb{C}_{q}$ be a function
so that: 
\begin{equation}
\chi\left(\mathfrak{z}\right)\overset{\mathbb{C}_{q}}{=}\lim_{N\rightarrow\infty}\chi\left(\left[\mathfrak{z}\right]_{p^{N}}\right)
\end{equation}
holds for almost every $\mathfrak{z}\in\mathbb{Z}_{p}$. Then, $\chi\in L_{\mathbb{R}}^{1}\left(\mathbb{Z}_{p},\mathbb{C}_{q}\right)$
whenever: 
\begin{equation}
\sum_{N=1}^{\infty}\frac{1}{p^{N}}\sum_{n=p^{N-1}}^{p^{N}-1}\left|c_{n}\left(\chi\right)\right|_{q}<\infty\label{eq:van der Put criterion for real L^1 integrability}
\end{equation}
\end{prop}
Proof: Using the truncated van der Put identity (\ref{eq:truncated van der Put identity}),
we have: 
\begin{align*}
\chi_{N}\left(\mathfrak{z}\right)-\chi_{N-1}\left(\mathfrak{z}\right) & \overset{\mathbb{C}_{q}}{=}\chi\left(\left[\mathfrak{z}\right]_{p^{N}}\right)-\chi\left(\left[\mathfrak{z}\right]_{p^{N-1}}\right)\\
 & =\sum_{n=p^{N-1}}^{p^{N}-1}c_{n}\left(\chi\right)\left[\mathfrak{z}\overset{p^{\lambda_{p}\left(n\right)}}{\equiv}n\right]\\
\left(\lambda_{p}\left(n\right)=N\Leftrightarrow p^{N-1}\leq n\leq p^{N}-1\right); & =\sum_{n=p^{N-1}}^{p^{N}-1}c_{n}\left(\chi\right)\left[\mathfrak{z}\overset{p^{N}}{\equiv}n\right]
\end{align*}
Now, fixing $N$, note that for any $\mathfrak{z}\in\mathbb{Z}_{p}$:
\begin{equation}
\sum_{n=p^{N-1}}^{p^{N}-1}c_{n}\left(\chi\right)\left[\mathfrak{z}\overset{p^{N}}{\equiv}n\right]\overset{\mathbb{R}}{=}\begin{cases}
c_{\left[\mathfrak{z}\right]_{p^{N}}}\left(\chi\right) & \textrm{if }p^{N-1}\leq\left[\mathfrak{z}\right]_{p^{N}}\leq p^{N}-1\\
0 & \textrm{else}
\end{cases}
\end{equation}
Hence: 
\begin{align*}
\left|\chi\left(\left[\mathfrak{z}\right]_{p^{N}}\right)-\chi\left(\left[\mathfrak{z}\right]_{p^{N-1}}\right)\right|_{q} & \overset{\mathbb{R}}{=}\left|c_{\left[\mathfrak{z}\right]_{p^{N}}}\left(\chi\right)\right|_{q}\left[p^{N-1}\leq\left[\mathfrak{z}\right]_{p^{N}}\leq p^{N}-1\right]\\
 & =\left|c_{\left[\mathfrak{z}\right]_{p^{N}}}\left(\chi\right)\right|_{q}\sum_{n=p^{N-1}}^{p^{N}-1}\left[\mathfrak{z}\overset{p^{N}}{\equiv}n\right]\\
 & =\sum_{n=p^{N-1}}^{p^{N}-1}\left|c_{n}\left(\chi\right)\right|_{q}\left[\mathfrak{z}\overset{p^{N}}{\equiv}n\right]
\end{align*}
Since the integral of $\left[\mathfrak{z}\overset{p^{N}}{\equiv}n\right]$
is $p^{-N}$, this yields: 
\begin{equation}
\left\Vert \chi_{N}-\chi_{N-1}\right\Vert _{1}\overset{\mathbb{R}}{=}\frac{1}{p^{N}}\sum_{n=p^{N-1}}^{p^{N}-1}\left|c_{n}\left(\chi\right)\right|_{q}\label{eq:Real Integral of difference of Nth truncations of Chi}
\end{equation}
If the right-hand side is summable in $\mathbb{R}$ with respect to
$N$, then: 
\begin{equation}
\sum_{N=1}^{\infty}\left\Vert \chi_{N}-\chi_{N-1}\right\Vert _{1}<\infty
\end{equation}
and hence, the series: 
\[
\sum_{N=1}^{\infty}\left(\chi_{N}\left(\mathfrak{z}\right)-\chi_{N-1}\left(\mathfrak{z}\right)\right)\overset{\mathbb{C}_{q}}{=}\lim_{N\rightarrow\infty}\chi_{N}\left(\mathfrak{z}\right)-\chi\left(0\right)\overset{\mathbb{C}_{q}}{=}\lim_{N\rightarrow\infty}\chi\left(\left[\mathfrak{z}\right]_{p^{N}}\right)-\chi\left(0\right)
\]
converges in $L_{\mathbb{R}}^{1}\left(\mathbb{Z}_{p},\mathbb{C}_{q}\right)$.
As such, since $\chi\left(\mathfrak{z}\right)\overset{\mathbb{C}_{q}}{=}\lim_{N\rightarrow\infty}\chi\left(\left[\mathfrak{z}\right]_{p^{N}}\right)$
almost everywhere, it then follows that $\chi\left(\mathfrak{z}\right)-\chi\left(0\right)$
is an element of $L_{\mathbb{R}}^{1}\left(\mathbb{Z}_{p},\mathbb{C}_{q}\right)$,
which shows that $\chi\in L_{\mathbb{R}}^{1}\left(\mathbb{Z}_{p},\mathbb{C}_{q}\right)$,
since the constant function $\chi\left(0\right)$ is in $L_{\mathbb{R}}^{1}\left(\mathbb{Z}_{p},\mathbb{C}_{q}\right)$.

Q.E.D. 
\begin{prop}
Let $d\mu\in C\left(\mathbb{Z}_{p},\mathbb{C}_{q}\right)^{\prime}$.
Then: 
\begin{equation}
\lim_{N\rightarrow\infty}\int_{\mathbb{Z}_{p}}\left|\tilde{\mu}_{N}\left(\mathfrak{z}\right)\right|_{q}d\mathfrak{z}\overset{\mathbb{R}}{=}0
\end{equation}
is a necessary\emph{ }condition for all of the following:

\vphantom{}

I. For $d\mu$ to be the zero measure;

\vphantom{}

II. For $d\mu$ to be a degenerate rising-continuous measure;

\vphantom{}

III. For $d\mu$ to be $\mathcal{F}$-rising, where the non-archimedean
domain of $\mathcal{F}$ is a set of full measure in $\mathbb{Z}_{p}$
(that is, its complement has measure zero in $\mathbb{Z}_{p}$). 
\end{prop}
Proof: If $d\mu$ is the zero measure, then, by \textbf{Proposition
\ref{prop:Criterion for zero measure in terms of partial sums of Fourier series}},
$\left\Vert \tilde{\mu}_{N}\right\Vert _{p,q}\rightarrow0$, and hence,
the $\tilde{\mu}_{N}$s converge $q$-adically to $0$ as $N\rightarrow\infty$
uniformly over $\mathbb{Z}_{p}$. Consequently, they are uniformly
bounded in $q$-adic absolute value. Moreover, since the $\tilde{\mu}_{N}$s
are $\left(p,q\right)$-adically continuous, the functions $\mathfrak{z}\mapsto\left|\tilde{\mu}_{N}\left(\mathfrak{z}\right)\right|_{q}$
are continuous, and hence measurable. Thus, by \textbf{Lemma \ref{lem:Egorov lemma}},
$\lim_{N\rightarrow\infty}\int_{\mathbb{Z}_{p}}\left|\tilde{\mu}_{N}\left(\mathfrak{z}\right)\right|_{q}d\mathfrak{z}\overset{\mathbb{R}}{=}0$.
So, this limit is necessary for $d\mu$ to be the zero measure.

(II) is a more specific case of (III). (III) holds true because if
the non-archimedean domain of $\mathcal{F}$ has full measure in $\mathbb{Z}_{p}$,
then the continuous (and hence measurable) real-valued functions $\mathfrak{z}\mapsto\left|\tilde{\mu}_{N}\left(\mathfrak{z}\right)\right|_{q}$
converge to $0$ almost everywhere on $\mathbb{Z}_{p}$. Moreover,
because $d\mu$ is a measure, we have that: 
\[
\left|\tilde{\mu}_{N}\left(\mathfrak{z}\right)\right|_{q}\leq\max_{\left|t\right|_{p}\leq p^{N}}\left|\hat{\mu}\left(t\right)\right|_{q}\leq\sup_{t\in\hat{\mathbb{Z}}_{p}}\left|\hat{\mu}\left(t\right)\right|_{q}<\infty,\textrm{ }\forall N\geq0
\]
As such, we can apply\textbf{ Lemma \ref{lem:Egorov lemma}} to obtain
$\lim_{N\rightarrow\infty}\int_{\mathbb{Z}_{p}}\left|\tilde{\mu}_{N}\left(\mathfrak{z}\right)\right|_{q}d\mathfrak{z}\overset{\mathbb{R}}{=}0$.

Q.E.D.

\vphantom{}

I find it interesting that everything I would like to have true for
quasi-integrable rising-continuous functions is easily verified for
\emph{continuous }$\left(p,q\right)$-adic functions. To show this,
we need the following formula:\index{van der Put!coefficients!Fourier coeffs.} 
\begin{prop}
\label{prop:3.78}\ 
\begin{equation}
\hat{f}_{N}\left(t\right)=\sum_{n=\frac{\left|t\right|_{p}}{p}}^{p^{N}-1}\frac{c_{n}\left(f\right)}{p^{\lambda_{p}\left(n\right)}}e^{-2\pi int},\textrm{ }\forall t\in\hat{\mathbb{Z}}_{p}\label{eq:Fourier transform of Nth truncation in terms of vdP coefficients}
\end{equation}
\end{prop}
Proof: Identical to the computation in \textbf{Theorem \ref{thm:formula for Fourier series}},
but with the $N$th partial van der Put series of $f$ ($S_{p:N}\left\{ f\right\} $)
used in place of the full van der Put series of $f$ ($S_{p}\left\{ f\right\} $).

Q.E.D.

\vphantom{}Before proceeding, we note that the above identity also
holds for rising-continuous functions. 
\begin{prop}
Let $\chi\in\tilde{C}\left(\mathbb{Z}_{p},\mathbb{C}_{q}\right)$.
Then, the Fourier transform of the $N$th truncation of $\chi$ is
given by: 
\end{prop}
\begin{equation}
\hat{\chi}_{N}\left(t\right)=\sum_{n=\frac{\left|t\right|_{p}}{p}}^{p^{N}-1}\frac{c_{n}\left(\chi\right)}{p^{\lambda_{p}\left(n\right)}}e^{-2\pi int},\textrm{ }\forall t\in\hat{\mathbb{Z}}_{p}\label{eq:Fourier transform of Nth truncation in terms of vdP coefficients for a rising continuous function}
\end{equation}
where the sum is defined to be $0$ whenever $\frac{\left|t\right|_{p}}{p}>p^{N}-1$
(i.e., whenever $\left|t\right|_{p}\geq p^{N+1}$). Note that this
equality holds in $\overline{\mathbb{Q}}$ whenever $\chi\mid_{\mathbb{N}_{0}}$
takes values in $\overline{\mathbb{Q}}$.

Proof: Since $\chi_{N}\in C\left(\mathbb{Z}_{p},\mathbb{C}_{q}\right)$
for all $N$, this immediately follows by \textbf{Proposition \ref{prop:3.78}}.

Q.E.D.

\vphantom{}

Here are the aforementioned properties provable for continuous $\left(p,q\right)$-adic
functions: 
\begin{lem}
Let $f\in C\left(\mathbb{Z}_{p},\mathbb{C}_{q}\right)$. Recall that
$f_{N}$ is the $N$th truncation of $f$, $\hat{f}_{N}$ is the Fourier
transform of $f_{N}$, $\hat{f}$ is the Fourier transform of $f$,
and write: 
\begin{equation}
\tilde{f}_{N}\left(\mathfrak{z}\right)\overset{\textrm{def}}{=}\sum_{\left|t\right|_{p}\leq p^{N}}\hat{f}\left(t\right)e^{2\pi i\left\{ t\mathfrak{z}\right\} _{p}}\label{eq:Definition of f_N twiddle}
\end{equation}
Then:

\vphantom{}

I. 
\begin{equation}
\hat{f}\left(t\right)-\hat{f}_{N}\left(t\right)=\sum_{n=p^{N}}^{\infty}\frac{c_{n}\left(f\right)}{p^{\lambda_{p}\left(n\right)}}e^{-2\pi int},\textrm{ }\forall t\in\hat{\mathbb{Z}}_{p}\label{eq:Difference between f-hat and f_N-hat}
\end{equation}
with: 
\begin{equation}
\left\Vert \hat{f}-\hat{f}_{N}\right\Vert _{p,q}\leq\sup_{n\geq p^{N}}\left|c_{n}\left(f\right)\right|_{q},\textrm{ }\forall N\geq0\label{eq:Infinit norm of f-hat - f_N-hat}
\end{equation}

\vphantom{}

II. 
\begin{equation}
\tilde{f}_{N}\left(\mathfrak{z}\right)-f_{N}\left(\mathfrak{z}\right)=p^{N}\sum_{n=p^{N}}^{\infty}\frac{c_{n}\left(f\right)}{p^{\lambda_{p}\left(n\right)}}\left[\mathfrak{z}\overset{p^{N}}{\equiv}n\right],\textrm{ }\forall\mathfrak{z}\in\mathbb{Z}_{p},\textrm{ }\forall N\geq0\label{eq:Difference between f_N-twiddle and f_N}
\end{equation}
with: 
\begin{equation}
\left\Vert \tilde{f}_{N}-f_{N}\right\Vert _{p,q}\leq\sup_{n\geq p^{N}}\left|c_{n}\left(f\right)\right|_{q}\label{eq:Infinity norm of difference between f_N twiddle and f_N}
\end{equation}
and: 
\begin{equation}
\int_{\mathbb{Z}_{p}}\left|\tilde{f}_{N}\left(\mathfrak{z}\right)-f_{N}\left(\mathfrak{z}\right)\right|_{q}d\mathfrak{z}\leq\sup_{n\geq p^{N}}\left|c_{n}\left(f\right)\right|_{q}\label{eq:L^1 norm of the difference between f_N-twiddle and f_N}
\end{equation}

\vphantom{}

III. 
\begin{equation}
\left\Vert f-\tilde{f}_{N}\right\Vert _{p,q}\leq\sup_{\left|t\right|_{p}>p^{N}}\left|\hat{f}\left(t\right)\right|_{q}\leq\sup_{n\geq p^{N}}\left|c_{n}\left(f\right)\right|_{q}\label{eq:Supremum norm of difference between f and f_N twiddle}
\end{equation}
and: 
\begin{equation}
\int_{\mathbb{Z}_{p}}\left|f\left(\mathfrak{z}\right)-\tilde{f}_{N}\left(\mathfrak{z}\right)\right|_{q}d\mathfrak{z}\leq\sup_{\left|t\right|_{p}>p^{N}}\left|\hat{f}\left(t\right)\right|_{q}\leq\sup_{n\geq p^{N}}\left|c_{n}\left(f\right)\right|_{q}\label{eq:L^1 norm of difference between f and f_N twiddle}
\end{equation}
\end{lem}
Proof: Since the Banach space $C\left(\mathbb{Z}_{p},\mathbb{C}_{q}\right)$
is isometrically isomorphic to $c_{0}\left(\hat{\mathbb{Z}}_{p},\mathbb{C}_{q}\right)$
by way of the $\left(p,q\right)$-adic Fourier transform (\ref{eq:Fourier transform of Nth truncation in terms of vdP coefficients}),
we have that: 
\begin{align*}
\hat{f}\left(t\right)-\hat{f}_{N}\left(t\right) & \overset{\mathbb{C}_{q}}{=}\sum_{n=\frac{\left|t\right|_{p}}{p}}^{\infty}\frac{c_{n}\left(f\right)}{p^{\lambda_{p}\left(n\right)}}e^{-2\pi int}-\sum_{n=\frac{\left|t\right|_{p}}{p}}^{p^{N}-1}\frac{c_{n}\left(f\right)}{p^{\lambda_{p}\left(n\right)}}e^{-2\pi int}\\
 & \overset{\mathbb{C}_{q}}{=}\sum_{n=p^{N}}^{\infty}\frac{c_{n}\left(f\right)}{p^{\lambda_{p}\left(n\right)}}e^{-2\pi int}
\end{align*}
This proves (\ref{eq:Difference between f-hat and f_N-hat}). Taking
$\left(p,q\right)$-adic supremum norm over $\hat{\mathbb{Z}}_{p}$
yields (\ref{eq:Infinit norm of f-hat - f_N-hat}).

Next, performing Fourier inversion, we have: 
\begin{align*}
\tilde{f}_{N}\left(\mathfrak{z}\right)-f_{N}\left(\mathfrak{z}\right) & =\sum_{\left|t\right|_{p}\leq p^{N}}\left(\hat{f}\left(t\right)-\hat{f}_{N}\left(t\right)\right)e^{2\pi i\left\{ t\mathfrak{z}\right\} _{p}}\\
 & =\sum_{\left|t\right|_{p}\leq p^{N}}\sum_{n=p^{N}}^{\infty}\frac{c_{n}\left(f\right)}{p^{\lambda_{p}\left(n\right)}}e^{-2\pi int}e^{2\pi i\left\{ t\mathfrak{z}\right\} _{p}}\\
 & =p^{N}\sum_{n=p^{N}}^{\infty}\frac{c_{n}\left(f\right)}{p^{\lambda_{p}\left(n\right)}}\left[\mathfrak{z}\overset{p^{N}}{\equiv}n\right]
\end{align*}
Hence: 
\[
\tilde{f}_{N}\left(\mathfrak{z}\right)-f_{N}\left(\mathfrak{z}\right)=p^{N}\sum_{n=p^{N}}^{\infty}\frac{c_{n}\left(f\right)}{p^{\lambda_{p}\left(n\right)}}\left[\mathfrak{z}\overset{p^{N}}{\equiv}n\right]
\]
which proves (\ref{eq:Difference between f_N-twiddle and f_N}). Taking
$\left(p,q\right)$-adic supremum norm over $\mathbb{Z}_{p}$ yields
(\ref{eq:Infinity norm of difference between f_N twiddle and f_N}).
On the other hand, taking $q$-adic absolute values, integrating gives
us: 
\begin{align*}
\int_{\mathbb{Z}_{p}}\left|\tilde{f}_{N}\left(\mathfrak{z}\right)-f_{N}\left(\mathfrak{z}\right)\right|_{q}d\mathfrak{z} & =\int_{\mathbb{Z}_{p}}\left|p^{N}\sum_{n=p^{N}}^{\infty}\frac{c_{n}\left(f\right)}{p^{\lambda_{p}\left(n\right)}}\left[\mathfrak{z}\overset{p^{N}}{\equiv}n\right]\right|_{q}d\mathfrak{z}\\
 & \leq\int_{\mathbb{Z}_{p}}\sup_{n\geq p^{N}}\left|\frac{c_{n}\left(f\right)}{p^{\lambda_{p}\left(n\right)}}\left[\mathfrak{z}\overset{p^{N}}{\equiv}n\right]\right|_{q}d\mathfrak{z}\\
 & \leq\int_{\mathbb{Z}_{p}}\sup_{n\geq p^{N}}\left(\left|c_{n}\left(f\right)\right|_{q}\left[\mathfrak{z}\overset{p^{N}}{\equiv}n\right]\right)d\mathfrak{z}\\
 & \leq\left(\sup_{n\geq p^{N}}\left|c_{n}\left(f\right)\right|_{q}\right)\int_{\mathbb{Z}_{p}}d\mathfrak{z}\\
 & =\sup_{n\geq p^{N}}\left|c_{n}\left(f\right)\right|_{q}
\end{align*}
which proves (\ref{eq:L^1 norm of the difference between f_N-twiddle and f_N}).

Lastly: 
\begin{align*}
f\left(\mathfrak{z}\right)-\tilde{f}_{N}\left(\mathfrak{z}\right) & =\sum_{t\in\hat{\mathbb{Z}}_{p}}\hat{f}\left(t\right)e^{2\pi i\left\{ t\mathfrak{z}\right\} _{p}}-\sum_{\left|t\right|_{p}\leq p^{N}}\hat{f}\left(t\right)e^{2\pi i\left\{ t\mathfrak{z}\right\} _{p}}\\
 & =\sum_{\left|t\right|_{p}>p^{N}}\hat{f}\left(t\right)e^{2\pi i\left\{ t\mathfrak{z}\right\} _{p}}
\end{align*}
and so: 
\begin{align*}
\left\Vert f-\tilde{f}_{N}\right\Vert _{p,q} & \leq\sup_{\left|t\right|_{p}>p^{N}}\left|\hat{f}\left(t\right)\right|_{q}\\
 & \leq\sup_{\left|t\right|_{p}>p^{N}}\left|\sum_{n=\frac{\left|t\right|_{p}}{p}}^{\infty}\frac{c_{n}\left(f\right)}{p^{\lambda_{p}\left(n\right)}}e^{-2\pi int}\right|_{q}\\
 & \leq\sup_{\left|t\right|_{p}>p^{N}}\sup_{n\geq\frac{\left|t\right|_{p}}{p}}\left|c_{n}\left(f\right)\right|_{q}\\
 & \leq\sup_{\left|t\right|_{p}>p^{N}}\sup_{n\geq\frac{p^{N+1}}{p}}\left|c_{n}\left(f\right)\right|_{q}\\
 & =\sup_{n\geq p^{N}}\left|c_{n}\left(f\right)\right|_{q}
\end{align*}
This proves (\ref{eq:Supremum norm of difference between f and f_N twiddle}).
If we instead integrate, we get: 
\begin{align*}
\int_{\mathbb{Z}_{p}}\left|f\left(\mathfrak{z}\right)-\tilde{f}_{N}\left(\mathfrak{z}\right)\right|_{q}d\mathfrak{z} & =\int_{\mathbb{Z}_{p}}\left|\sum_{\left|t\right|_{p}>p^{N}}\hat{f}\left(t\right)e^{2\pi i\left\{ t\mathfrak{z}\right\} _{p}}\right|_{q}d\mathfrak{z}\\
 & \leq\int_{\mathbb{Z}_{p}}\left(\sup_{\left|t\right|_{p}>p^{N}}\left|\hat{f}\left(t\right)\right|_{q}\right)d\mathfrak{z}\\
 & =\sup_{\left|t\right|_{p}>p^{N}}\left|\hat{f}\left(t\right)\right|_{q}\\
 & \leq\sup_{n\geq p^{N}}\left|c_{n}\left(f\right)\right|_{q}
\end{align*}

Q.E.D.

\vphantom{}

The last result of this subsection is a sufficient condition for a
rising-continuous function to be an element of $L_{\mathbb{R}}^{1}\left(\mathbb{Z}_{p},\mathbb{C}_{q}\right)$,
formulated as a decay condition on the $q$-adic absolute value of
the function's van der Put coefficient. We need two propositions first. 
\begin{prop}
Let $\chi:\mathbb{Z}_{p}\rightarrow\mathbb{C}_{q}$ be any function.
Then, for integers $N,k\geq1$ the difference between the $\left(N+k\right)$th
and $N$th truncations of $\chi$ is: 
\begin{align}
\chi_{N+k}\left(\mathfrak{z}\right)-\chi_{N}\left(\mathfrak{z}\right) & =\sum_{j=1}^{p^{k}-1}\sum_{n=0}^{p^{N}-1}\left(\chi\left(n+jp^{N}\right)-\chi\left(n\right)\right)\left[\mathfrak{z}\overset{p^{N+k}}{\equiv}n+jp^{N}\right]\label{eq:Difference between N+kth and Nth truncations of Chi}
\end{align}
\end{prop}
Proof: Note that: 
\begin{equation}
\left[\mathfrak{z}\overset{p^{N}}{\equiv}n\right]=\sum_{j=0}^{p^{k}-1}\left[\mathfrak{z}\overset{p^{N+k}}{\equiv}n+p^{N}j\right],\textrm{ }\forall N,k\geq0,\textrm{ }\forall n\in\left\{ 0,\ldots,p^{N}-1\right\} ,\textrm{ }\forall\mathfrak{z}\in\mathbb{Z}_{p}
\end{equation}
Hence: 
\begin{align*}
\chi_{N+k}\left(\mathfrak{z}\right)-\chi_{N}\left(\mathfrak{z}\right) & =\sum_{n=0}^{p^{N+k}-1}\chi\left(n\right)\left[\mathfrak{z}\overset{p^{N+k}}{\equiv}n\right]-\sum_{n=0}^{p^{N}-1}\chi\left(n\right)\left[\mathfrak{z}\overset{p^{N}}{\equiv}n\right]\\
 & =\sum_{n=0}^{p^{N+k}-1}\chi\left(n\right)\left[\mathfrak{z}\overset{p^{N+k}}{\equiv}n\right]-\sum_{n=0}^{p^{N}-1}\chi\left(n\right)\sum_{j=0}^{p^{k}-1}\left[\mathfrak{z}\overset{p^{N+k}}{\equiv}n+p^{N}j\right]\\
 & =\sum_{n=0}^{p^{N}-1}\chi\left(n\right)\underbrace{\left(\left[\mathfrak{z}\overset{p^{N+k}}{\equiv}n\right]-\sum_{j=0}^{p^{k}-1}\left[\mathfrak{z}\overset{p^{N+k}}{\equiv}n+p^{N}j\right]\right)}_{j=0\textrm{ term cancels}}\\
 & +\sum_{n=p^{N}}^{p^{N+k}-1}\chi\left(n\right)\left[\mathfrak{z}\overset{p^{N+k}}{\equiv}n\right]\\
 & =\sum_{n=p^{N}}^{p^{N+k}-1}\chi\left(n\right)\left[\mathfrak{z}\overset{p^{N+k}}{\equiv}n\right]-\sum_{n=0}^{p^{N}-1}\chi\left(n\right)\sum_{j=1}^{p^{k}-1}\left[\mathfrak{z}\overset{p^{N+k}}{\equiv}n+p^{N}j\right]
\end{align*}
Here: 
\begin{align*}
\sum_{n=p^{N}}^{p^{N+k}-1}\chi\left(n\right)\left[\mathfrak{z}\overset{p^{N+k}}{\equiv}n\right] & =\sum_{n=0}^{\left(p^{k}-1\right)p^{N}-1}\chi\left(n+p^{N}\right)\left[\mathfrak{z}\overset{p^{N+k}}{\equiv}n+p^{N}\right]\\
 & =\sum_{j=1}^{p^{k}-1}\sum_{n=\left(j-1\right)p^{N}}^{jp^{N}-1}\chi\left(n+p^{N}\right)\left[\mathfrak{z}\overset{p^{N+k}}{\equiv}n+p^{N}\right]\\
 & =\sum_{j=1}^{p^{k}-1}\sum_{n=0}^{p^{N}-1}\chi\left(n+jp^{N}\right)\left[\mathfrak{z}\overset{p^{N+k}}{\equiv}n+jp^{N}\right]
\end{align*}
and so: 
\begin{align*}
\chi_{N+k}\left(\mathfrak{z}\right)-\chi_{N}\left(\mathfrak{z}\right) & =\sum_{n=p^{N}}^{p^{N+k}-1}\chi\left(n\right)\left[\mathfrak{z}\overset{p^{N+k}}{\equiv}n\right]-\sum_{n=0}^{p^{N}-1}\chi\left(n\right)\sum_{j=1}^{p^{k}-1}\left[\mathfrak{z}\overset{p^{N+k}}{\equiv}n+jp^{N}\right]\\
 & =\sum_{j=1}^{p^{k}-1}\sum_{n=0}^{p^{N}-1}\left(\chi\left(n+jp^{N}\right)-\chi\left(n\right)\right)\left[\mathfrak{z}\overset{p^{N+k}}{\equiv}n+jp^{N}\right]
\end{align*}

Q.E.D. 
\begin{prop}
Let $\chi:\mathbb{Z}_{p}\rightarrow\mathbb{C}_{q}$ be any function.
Then, for all integers $N,k\geq1$:

\begin{equation}
\int_{\mathbb{Z}_{p}}\left|\chi_{N+k}\left(\mathfrak{z}\right)-\chi_{N}\left(\mathfrak{z}\right)\right|_{q}d\mathfrak{z}=\sum_{j=1}^{p^{k}-1}\sum_{n=0}^{p^{N}-1}\frac{\left|\chi\left(n+jp^{N}\right)-\chi\left(n\right)\right|_{q}}{p^{N+k}}\label{eq:L^1_R norm of Chi_N+k minus Chi_N}
\end{equation}
\end{prop}
Proof: Let $N$ and $k$ be arbitrary. Note that the brackets in (\ref{eq:Difference between N+kth and Nth truncations of Chi})
have supports which are pair-wise disjoint with respect to $n$ and
$j$, and that each of the brackets has Haar measure $1/p^{N+k}$.
As such: 
\[
\int_{\mathbb{Z}_{p}}\left|\chi_{N+k}\left(\mathfrak{z}\right)-\chi_{N}\left(\mathfrak{z}\right)\right|_{q}d\mathfrak{z}=\sum_{j=1}^{p^{k}-1}\sum_{n=0}^{p^{N}-1}\frac{\left|\chi\left(n+jp^{N}\right)-\chi\left(n\right)\right|_{q}}{p^{N+k}}
\]

Q.E.D.

\vphantom{}

Here, then, is the $L_{\mathbb{R}}^{1}$ criterion: 
\begin{lem}
Let\index{van der Put!coefficients!L_{mathbb{R}}^{1} criterion@$L_{\mathbb{R}}^{1}$ criterion}
$\chi\in\tilde{C}\left(\mathbb{Z}_{p},\mathbb{C}_{q}\right)$. If
there exists an $r\in\left(0,1\right)$ so that for all sufficiently
large integers $N\geq0$: 
\begin{equation}
\max_{0\leq n<p^{N}}\left|c_{n+jp^{N}}\left(\chi\right)\right|_{q}\ll r^{j},\textrm{ }\forall j\geq1\label{eq:van der Put coefficient asymptotic condition for L^1_R integrability}
\end{equation}
Then $\chi\in L_{\mathbb{R}}^{1}\left(\mathbb{Z}_{p},\mathbb{C}_{q}\right)$. 
\end{lem}
Proof: Letting $\chi$ and $r$ be as given, note that: 
\begin{equation}
c_{n+jp^{N}}\left(\chi\right)=\chi\left(n+jp^{N}\right)-\chi\left(n\right),\textrm{ }\forall n<p^{N},\textrm{ }\forall j\geq1
\end{equation}
By (\ref{eq:L^1_R norm of Chi_N+k minus Chi_N}), we have that: 
\begin{align*}
\int_{\mathbb{Z}_{p}}\left|\chi_{N+k}\left(\mathfrak{z}\right)-\chi_{N}\left(\mathfrak{z}\right)\right|_{q}d\mathfrak{z} & \ll\sum_{n=0}^{p^{N}-1}\sum_{j=1}^{p^{k}-1}\frac{r^{j}}{p^{N+k}}\\
 & =\sum_{n=0}^{p^{N}-1}\frac{1}{p^{N+k}}\left(-1+\frac{1-r^{p^{k}}}{1-r}\right)\\
 & =\frac{1}{p^{k}}\frac{r-r^{p^{k}}}{1-r}
\end{align*}
which tends to $0$ in $\mathbb{R}$ as $k\rightarrow\infty$. Since
$k$ was arbitrary, it follows that for any $\epsilon>0$, choosing
$N$ and $M=N+k$ sufficiently large, we can guarantee that: 
\begin{equation}
\int_{\mathbb{Z}_{p}}\left|\chi_{M}\left(\mathfrak{z}\right)-\chi_{N}\left(\mathfrak{z}\right)\right|_{q}d\mathfrak{z}<\epsilon
\end{equation}
which shows that the $\chi_{N}$s are Cauchy in $L_{\mathbb{R}}^{1}\left(\mathbb{Z}_{p},\mathbb{C}_{q}\right)$
norm. Since $L_{\mathbb{R}}^{1}\left(\mathbb{Z}_{p},\mathbb{C}_{q}\right)$
is complete, and since, as a rising-continuous function, the $\chi_{N}$s
converge point-wise to $\chi$, it follows that $\chi$ is then an
element of $L_{\mathbb{R}}^{1}\left(\mathbb{Z}_{p},\mathbb{C}_{q}\right)$.

Q.E.D.

\vphantom{}

We end this subsection with two questions: 
\begin{question}
\label{que:3.7}Is every function in $\tilde{C}\left(\mathbb{Z}_{p},\mathbb{C}_{q}\right)\cap L_{\mathbb{R}}^{1}\left(\mathbb{Z}_{p},\mathbb{C}_{q}\right)$
and element of $\tilde{L}^{1}\left(\mathcal{F}_{p,q}\right)$? 
\end{question}
\begin{question}
Let $\mathcal{F}$ be a $p$-aduc frame. Is $\tilde{C}\left(\mathbb{Z}_{p},\mathbb{C}_{q}\right)\cap\tilde{L}^{1}\left(\mathcal{F}\right)\subseteq L_{\mathbb{R}}^{1}\left(\mathbb{Z}_{p},\mathbb{C}_{q}\right)$?
Is $\tilde{L}^{1}\left(\mathcal{F}\right)\subseteq L_{\mathbb{R}}^{1}\left(\mathbb{Z}_{p},\mathbb{C}_{q}\right)$? 
\end{question}

\subsection{\label{subsec:3.3.7 -adic-Wiener-Tauberian}$\left(p,q\right)$-adic
Wiener Tauberian Theorems}

One of the great mathematical accomplishments of the twentieth century
was the transfiguration of the study of divergent series from the
realm of foolish fancy to one of serious, rigorous study. Asymptotic
analysis\textemdash be it Karamata's theory of functions of regular
variation\index{regular variation}\footnote{Even though it has little to do with the topic at hand, I cannot give
too strong of a recommendation of Bingham, Goldie, and Teugels' book
\emph{Regular Variation} \cite{Regular Variation}. It should be required
reading for anyone working in analysis. It is filled with a thousand
and one useful beauties that you never knew you wanted to know.}, or Abelian or Tauberian summability theorems\textemdash are powerful
and widely useful in number theory, probability theory, and nearly
everything in between.

The \index{Wiener!Tauberian Theorem}Wiener Tauberian Theorem (WTT)
is a capstone of these investigations. This theorem has a protean
diversity of forms and restatements, all of which are interesting
in their own particular way. Jacob Korevaar\index{Korevaar, Jacob}'s
book on Tauberian theory gives an excellent exposition of most of
them \cite{Korevaar}. The version of the WTT that matters to us sits
at the boarder between Fourier analysis and Functional analysis, where
it characterizes the relationship between a $1$-periodic function
$\phi:\mathbb{R}/\mathbb{Z}\rightarrow\mathbb{C}$ (where $\mathbb{R}/\mathbb{Z}$
is identified with $\left[0,1\right)$ with addition modulo $1$),
its Fourier coefficients $\hat{\phi}\left(n\right)$, and its reciprocal
$1/\phi$. 
\begin{defn}
The \index{Wiener!algebra}\textbf{Wiener algebra (on the circle)},
denote $A\left(\mathbb{T}\right)$, is the $\mathbb{C}$-linear space
of all absolutely convergent Fourier series: 
\begin{equation}
A\left(\mathbb{T}\right)\overset{\textrm{def}}{=}\left\{ \sum_{n\in\mathbb{Z}}\hat{\phi}\left(n\right)e^{2\pi int}:\hat{\phi}\in\ell^{1}\left(\mathbb{Z},\mathbb{C}\right)\right\} \label{eq:Definition of the Wiener algebra}
\end{equation}
where: 
\begin{equation}
\ell^{1}\left(\mathbb{Z},\mathbb{C}\right)\overset{\textrm{def}}{=}\left\{ \hat{\phi}:\mathbb{Z}\rightarrow\mathbb{C}:\sum_{n\in\mathbb{Z}}\left|\hat{\phi}\left(n\right)\right|<\infty\right\} \label{eq:Definition of ell 1 of Z, C}
\end{equation}
$A\left(\mathbb{T}\right)$ is then made into a Banach space with
the norm: 
\begin{equation}
\left\Vert \phi\right\Vert _{A}\overset{\textrm{def}}{=}\left\Vert \hat{\phi}\right\Vert _{1}\overset{\textrm{def}}{=}\sum_{n\in\mathbb{Z}}\left|\hat{\phi}\left(n\right)\right|,\textrm{ }\forall\phi\in A\left(\mathbb{T}\right)\label{eq:Definition of Wiener algebra norm}
\end{equation}
\end{defn}
\vphantom{}

As defined, observe that $A\left(\mathbb{T}\right)$ is then isometrically
isomorphic to the Banach space $\ell^{1}\left(\mathbb{Z},\mathbb{C}\right)$.
Because $\left(\mathbb{Z},+\right)$ is a locally compact abelian
group, $\ell^{1}\left(\mathbb{Z},\mathbb{C}\right)$ can be made into
a Banach algebra (the \textbf{group algebra}\index{group algebra}\textbf{
}of $\mathbb{Z}$) by equipping it with the convolution operation:
\begin{equation}
\left(\hat{\phi}*\hat{\psi}\right)\left(n\right)\overset{\textrm{def}}{=}\sum_{m\in\mathbb{Z}}\hat{\phi}\left(n-m\right)\hat{\psi}\left(m\right)\label{eq:ell 1 convolution}
\end{equation}
In this way, $A\left(\mathbb{T}\right)$ becomes a Banach algebra
under point-wise multiplication, thanks to \textbf{Young's Convolution
Inequality} and the fact that the Fourier transform turns multiplication
into convolution: 
\begin{equation}
\left\Vert \phi\cdot\psi\right\Vert _{A}=\left\Vert \widehat{\phi\cdot\psi}\right\Vert _{1}=\left\Vert \hat{\phi}*\hat{\psi}\right\Vert _{1}\overset{\textrm{Young}}{\leq}\left\Vert \hat{\phi}\right\Vert _{1}\left\Vert \hat{\psi}\right\Vert _{1}=\left\Vert \phi\right\Vert _{A}\left\Vert \psi\right\Vert _{A}
\end{equation}

With this terminology, we can state three versions of the WTT. 
\begin{thm}[\textbf{Wiener's Tauberian Theorem, Ver. 1}]
Let $\phi\in A\left(\mathbb{T}\right)$. Then, $1/\phi\in A\left(\mathbb{T}\right)$
if and only if\footnote{That is to say, the units of the Wiener algebra are precisely those
functions $\phi\in A\left(\mathbb{T}\right)$ which have no zeroes.} $\phi\left(t\right)\neq0$ for all $t\in\mathbb{R}/\mathbb{Z}$. 
\end{thm}
At first glance, this statement might seem almost tautological, until
you realize that we are requiring more than just the well-definedness
of $1/\phi$: we are also demanding it to have an absolutely convergent
Fourier series. Because the Fourier transform turns multiplication
into convolution, the existence of a multiplicative inverse for $\phi\in A\left(\mathbb{T}\right)$
is equivalent to the existence of a \index{convolution!inverse}\textbf{convolution
inverse}\emph{ }for $\hat{\phi}\in\ell^{1}\left(\mathbb{Z},\mathbb{C}\right)$,
by which I mean a function $\hat{\phi}^{-1}\in\ell^{1}\left(\mathbb{Z},\mathbb{C}\right)$
so that: 
\begin{equation}
\left(\hat{\phi}*\hat{\phi}^{-1}\right)\left(n\right)=\left(\hat{\phi}^{-1}*\hat{\phi}\right)\left(n\right)=\mathbf{1}_{0}\left(n\right)=\begin{cases}
1 & \textrm{if }n=0\\
0 & \textrm{else}
\end{cases}\label{eq:Definition of a convolution inverse}
\end{equation}
because the function $\mathbf{1}_{0}\left(n\right)$ which is $1$
at $n=0$ and $0$ for all other $n$ is identity element of the convolution
operation on $\ell^{1}\left(\mathbb{Z},\mathbb{C}\right)$. This gives
us a second version of the WTT. 
\begin{thm}[\textbf{Wiener's Tauberian Theorem, Ver. 2}]
Let $\hat{\phi}\in\ell^{1}\left(\mathbb{Z},\mathbb{C}\right)$. Then
$\hat{\phi}$ has a convolution inverse $\hat{\phi}^{-1}\in\ell^{1}\left(\mathbb{Z},\mathbb{C}\right)$
if and only if $\phi\left(t\right)\neq0$ for all $t\in\mathbb{R}/\mathbb{Z}$. 
\end{thm}
\vphantom{}

To get a more analytic statement of the theorem, observe that the
existence of a convolution inverse $\hat{\phi}^{-1}\in\ell^{1}\left(\mathbb{Z},\mathbb{C}\right)$
for $\hat{\phi}\in\ell^{1}\left(\mathbb{Z},\mathbb{C}\right)$ means
that for any $\hat{\psi}\in\ell^{1}\left(\mathbb{Z},\mathbb{C}\right)$,
we can convolve $\hat{\phi}$ with a function in $\ell^{1}\left(\mathbb{Z},\mathbb{C}\right)$
to obtain $\hat{\psi}$: 
\begin{equation}
\hat{\phi}*\left(\hat{\phi}^{-1}*\hat{\psi}\right)=\hat{\psi}
\end{equation}
Writing $\hat{\chi}$ to denote $\hat{\phi}^{-1}*\hat{\psi}$, observe
the limit:
\begin{equation}
\left(\hat{\phi}*\hat{\chi}\right)\left(n\right)\overset{\mathbb{C}}{=}\lim_{M\rightarrow\infty}\sum_{\left|m\right|\leq M}\hat{\chi}\left(m\right)\hat{\phi}\left(n-m\right)
\end{equation}
where, the right-hand side is a linear combination of translates of
$\hat{\phi}$. In fact, the convergence is in $\ell^{1}$-norm: 
\begin{equation}
\lim_{M\rightarrow\infty}\sum_{n\in\mathbb{Z}}\left|\hat{\psi}\left(n\right)-\sum_{\left|m\right|\leq M}\hat{\chi}\left(m\right)\hat{\phi}\left(n-m\right)\right|\overset{\mathbb{C}}{=}0
\end{equation}
This gives us a third version of the WTT: 
\begin{thm}[\textbf{Wiener's Tauberian Theorem, Ver. 3}]
Let $\phi\in A\left(\mathbb{T}\right)$. Then, the translates of
$\hat{\phi}$ are dense in $\ell^{1}\left(\mathbb{Z},\mathbb{C}\right)$
if and only if $\phi\left(t\right)\neq0$ for any $t\in\mathbb{R}/\mathbb{Z}$. 
\end{thm}
\vphantom{}

Because the \textbf{Correspondence Principle }tells us that periodic
points of a Hydra map $H$ are precisely those $x\in\mathbb{Z}$ which
lie in the image of $\chi_{H}$, a $\left(p,q\right)$-adic analogue
of WTT will allow us to study the multiplicative inverse of $\mathfrak{z}\mapsto\chi_{H}\left(\mathfrak{z}\right)-x$
(for any fixed $x$) by using Fourier transforms of $\chi_{H}$. And
while this hope is eventually borne out, it does not come to pass
without incident. In this subsection, I establish two realizations
of a $\left(p,q\right)$-adic WTT: one for continuous functions; another
for measures $d\mu$.

First, however, a useful notation for translates: 
\begin{defn}
For $s\in\hat{\mathbb{Z}}_{p}$, $\mathfrak{a}\in\mathbb{Z}_{p}$,
$f:\mathbb{Z}_{p}\rightarrow\mathbb{C}_{q}$, and $\hat{f}:\mathbb{Z}_{p}\rightarrow\mathbb{C}_{q}$,
we write: 
\begin{equation}
\tau_{s}\left\{ \hat{f}\right\} \left(t\right)\overset{\textrm{def}}{=}\hat{f}\left(t+s\right)\label{eq:definition of the translate of f hat}
\end{equation}
and: 
\begin{equation}
\tau_{\mathfrak{a}}\left\{ f\right\} \left(\mathfrak{z}\right)\overset{\textrm{def}}{=}f\left(\mathfrak{z}+\mathfrak{a}\right)\label{eq:Definition of the translate of f}
\end{equation}
\end{defn}
\begin{thm}[\textbf{Wiener Tauberian Theorem for Continuous $\left(p,q\right)$-adic
Functions}]
\index{Wiener!Tauberian Theorem!left(p,qright)-adic@$\left(p,q\right)$-adic}\label{thm:pq WTT for continuous functions}Let
$f\in C\left(\mathbb{Z}_{p},\mathbb{C}_{q}\right)$. Then, the following
are equivalent:

\vphantom{}

I. $\frac{1}{f}\in C\left(\mathbb{Z}_{p},\mathbb{C}_{q}\right)$;

\vphantom{}

II. $\hat{f}$ has a convolution inverse in $c_{0}\left(\hat{\mathbb{Z}}_{p},\mathbb{C}_{q}\right)$;

\vphantom{}

III. $\textrm{span}_{\mathbb{C}_{q}}\left\{ \tau_{s}\left\{ \hat{f}\right\} \left(t\right):s\in\hat{\mathbb{Z}}_{p}\right\} $
is dense in $c_{0}\left(\hat{\mathbb{Z}}_{p},\mathbb{C}_{q}\right)$;

\vphantom{}

IV. $f$ has no zeroes. 
\end{thm}
Proof:

\textbullet{} ($\textrm{I}\Rightarrow\textrm{II}$) Suppose $\frac{1}{f}$
is continuous. Then, since the $\left(p,q\right)$-adic Fourier transform
is an isometric isomorphism of the Banach algebra $C\left(\mathbb{Z}_{p},\mathbb{C}_{q}\right)$
onto the Banach algebra $c_{0}\left(\hat{\mathbb{Z}}_{p},\mathbb{C}_{q}\right)$,
it follows that: 
\begin{align*}
f\left(\mathfrak{z}\right)\cdot\frac{1}{f\left(\mathfrak{z}\right)} & =1,\textrm{ }\forall\mathfrak{z}\in\mathbb{Z}_{p}\\
\left(\textrm{Fourier transform}\right); & \Updownarrow\\
\left(\hat{f}*\widehat{\left(\frac{1}{f}\right)}\right)\left(t\right) & =\mathbf{1}_{0}\left(t\right),\textrm{ }\forall t\in\hat{\mathbb{Z}}_{p}
\end{align*}
where both $\hat{f}$ and $\widehat{\left(1/f\right)}$ are in $c_{0}$.
$\hat{f}$ has a convolution inverse in $c_{0}$, and this inverse
is $\widehat{\left(1/f\right)}$.

\vphantom{}

\textbullet{} ($\textrm{II}\Rightarrow\textrm{III}$) Suppose $\hat{f}$
has a convolution inverse $\hat{f}^{-1}\in c_{0}\left(\hat{\mathbb{Z}}_{p},\mathbb{C}_{q}\right)$.
Then, letting $\hat{g}\in c_{0}\left(\hat{\mathbb{Z}}_{p},\mathbb{C}_{q}\right)$
be arbitrary, we have that: 
\begin{equation}
\left(\hat{f}*\left(\hat{f}^{-1}*\hat{g}\right)\right)\left(t\right)=\left(\left(\hat{f}*\hat{f}^{-1}\right)*\hat{g}\right)\left(t\right)=\left(\mathbf{1}_{0}*\hat{g}\right)\left(t\right)=\hat{g}\left(t\right),\textrm{ }\forall t\in\hat{\mathbb{Z}}_{p}
\end{equation}
In particular, letting $\hat{h}$ denote $\hat{f}^{-1}*\hat{g}$,
we have that: 
\begin{equation}
\hat{g}\left(t\right)=\left(\hat{f}*\hat{h}\right)\left(t\right)\overset{\mathbb{C}_{q}}{=}\lim_{N\rightarrow\infty}\sum_{\left|s\right|_{p}\leq p^{N}}\hat{h}\left(s\right)\hat{f}\left(t-s\right)
\end{equation}
Since $\hat{f}$ and $\hat{h}$ are in $c_{0}$, note that: 
\begin{align*}
\sup_{t\in\hat{\mathbb{Z}}_{p}}\left|\sum_{s\in\hat{\mathbb{Z}}_{p}}\hat{h}\left(s\right)\hat{f}\left(t-s\right)-\sum_{\left|s\right|_{p}\leq p^{N}}\hat{h}\left(s\right)\hat{f}\left(t-s\right)\right|_{q} & \leq\sup_{t\in\hat{\mathbb{Z}}_{p}}\sup_{\left|s\right|_{p}>p^{N}}\left|\hat{h}\left(s\right)\hat{f}\left(t-s\right)\right|_{q}\\
\left(\left|\hat{f}\right|_{q}<\infty\right); & \leq\sup_{\left|s\right|_{p}>p^{N}}\left|\hat{h}\left(s\right)\right|_{q}
\end{align*}
and hence: 
\begin{equation}
\lim_{N\rightarrow\infty}\sup_{t\in\hat{\mathbb{Z}}_{p}}\left|\sum_{s\in\hat{\mathbb{Z}}_{p}}\hat{h}\left(s\right)\hat{f}\left(t-s\right)-\sum_{\left|s\right|_{p}\leq p^{N}}\hat{h}\left(s\right)\hat{f}\left(t-s\right)\right|_{q}\overset{\mathbb{R}}{=}\lim_{N\rightarrow\infty}\sup_{\left|s\right|_{p}>p^{N}}\left|\hat{h}\left(s\right)\right|_{q}\overset{\mathbb{R}}{=}0
\end{equation}
which shows that the $q$-adic convergence of $\sum_{\left|s\right|_{p}\leq p^{N}}\hat{h}\left(s\right)\hat{f}\left(t-s\right)$
to $\hat{g}\left(t\right)=\sum_{s\in\hat{\mathbb{Z}}_{p}}\hat{h}\left(s\right)\hat{f}\left(t-s\right)$
is uniform in $t$. Hence: 
\begin{equation}
\lim_{N\rightarrow\infty}\sup_{t\in\hat{\mathbb{Z}}_{p}}\left|\hat{g}\left(t\right)-\sum_{\left|s\right|_{p}\leq p^{N}}\hat{h}\left(s\right)\hat{f}\left(t-s\right)\right|_{q}=0
\end{equation}
which is precisely the definition of what it means for the sequence
$\left\{ \sum_{\left|s\right|_{p}\leq p^{N}}\hat{h}\left(s\right)\hat{f}\left(t-s\right)\right\} _{N\geq0}$
to converge in $c_{0}\left(\hat{\mathbb{Z}}_{p},\mathbb{C}_{q}\right)$
to $\hat{g}$. Since $\hat{g}$ was arbitrary, and since this sequence
is in the span of the translates of $\hat{f}$ over $\mathbb{C}_{q}$,
we see that said span is dense in $c_{0}\left(\hat{\mathbb{Z}}_{p},\mathbb{C}_{q}\right)$.

\vphantom{}

\textbullet{} ($\textrm{III}\Rightarrow\textrm{IV}$) Suppose $\textrm{span}_{\mathbb{C}_{q}}\left\{ \tau_{s}\left\{ \hat{f}\right\} \left(t\right):s\in\hat{\mathbb{Z}}_{p}\right\} $
is dense in $c_{0}\left(\hat{\mathbb{Z}}_{p},\mathbb{C}_{q}\right)$.
Since $\mathbf{1}_{0}\left(t\right)\in c_{0}\left(\hat{\mathbb{Z}}_{p},\mathbb{C}_{q}\right)$,
given any $\epsilon\in\left(0,1\right)$, we can then choose constants
$\mathfrak{c}_{1},\ldots,\mathfrak{c}_{M}\in\mathbb{C}_{q}$ and $t_{1},\ldots,t_{M}\in\hat{\mathbb{Z}}_{p}$
so that: 
\begin{equation}
\sup_{t\in\hat{\mathbb{Z}}_{p}}\left|\mathbf{1}_{0}\left(t\right)-\sum_{m=1}^{M}\mathfrak{c}_{m}\hat{f}\left(t-t_{m}\right)\right|_{q}<\epsilon
\end{equation}
Now, letting $N\geq\max\left\{ -v_{p}\left(t_{1}\right),\ldots,-v_{p}\left(t_{M}\right)\right\} $
be arbitrary (note that $-v_{p}\left(t\right)=-\infty$ when $t=0$),
the maps $t\mapsto t+t_{m}$ are then bijections of the set $\left\{ t\in\hat{\mathbb{Z}}_{p}:\left|t\right|_{p}\leq p^{N}\right\} $.
Consequently: 
\begin{align*}
\sum_{\left|t\right|_{p}\leq p^{N}}\left(\sum_{m=1}^{M}\mathfrak{c}_{m}\hat{f}\left(t-t_{m}\right)\right)e^{2\pi i\left\{ t\mathfrak{z}\right\} _{p}} & =\sum_{m=1}^{M}\mathfrak{c}_{m}\sum_{\left|t\right|_{p}\leq p^{N}}\hat{f}\left(t\right)e^{2\pi i\left\{ \left(t+t_{m}\right)\mathfrak{z}\right\} _{p}}\\
 & =\left(\sum_{m=1}^{M}\mathfrak{c}_{m}e^{2\pi i\left\{ t_{m}\mathfrak{z}\right\} _{p}}\right)\sum_{\left|t\right|_{p}\leq p^{N}}\hat{f}\left(t\right)e^{2\pi i\left\{ t\mathfrak{z}\right\} _{p}}
\end{align*}
Letting $N\rightarrow\infty$, we obtain: 
\begin{equation}
\lim_{N\rightarrow\infty}\sum_{\left|t\right|_{p}\leq p^{N}}\left(\sum_{m=1}^{M}\mathfrak{c}_{m}\hat{f}\left(t-t_{m}\right)\right)e^{2\pi i\left\{ t\mathfrak{z}\right\} _{p}}\overset{\mathbb{C}_{q}}{=}g_{m}\left(\mathfrak{z}\right)f\left(\mathfrak{z}\right)
\end{equation}
where $g_{m}:\mathbb{Z}_{p}\rightarrow\mathbb{C}_{q}$ is defined
by: 
\begin{equation}
g_{m}\left(\mathfrak{z}\right)\overset{\textrm{def}}{=}\sum_{m=1}^{M}\mathfrak{c}_{m}e^{2\pi i\left\{ t_{m}\mathfrak{z}\right\} _{p}},\textrm{ }\forall\mathfrak{z}\in\mathbb{Z}_{p}\label{eq:Definition of g_m}
\end{equation}
Moreover, this convergence is uniform with respect to $\mathfrak{z}$.

Consequently: 
\begin{align*}
\left|1-g_{m}\left(\mathfrak{z}\right)f\left(\mathfrak{z}\right)\right|_{q} & \overset{\mathbb{R}}{=}\lim_{N\rightarrow\infty}\left|\sum_{\left|t\right|_{p}\leq p^{N}}\left(\mathbf{1}_{0}\left(t\right)-\sum_{m=1}^{M}\mathfrak{c}_{m}\hat{f}\left(t-t_{m}\right)\right)e^{2\pi i\left\{ t\mathfrak{z}\right\} _{p}}\right|_{q}\\
 & \leq\sup_{t\in\hat{\mathbb{Z}}_{p}}\left|\mathbf{1}_{0}\left(t\right)-\sum_{m=1}^{M}\mathfrak{c}_{m}\hat{f}\left(t-t_{m}\right)\right|_{q}\\
 & <\epsilon
\end{align*}
for all $\mathfrak{z}\in\mathbb{Z}_{p}$.

If $f\left(\mathfrak{z}_{0}\right)=0$ for some $\mathfrak{z}_{0}\in\mathbb{Z}_{p}$,
we would then have: 
\begin{equation}
\epsilon>\left|1-g_{m}\left(\mathfrak{z}_{0}\right)\cdot0\right|_{q}=1
\end{equation}
which would contradict the fact that $\epsilon\in\left(0,1\right)$.
As such, $f$ cannot have any zeroes whenever the span of $\hat{f}$'s
translates are dense in $c_{0}\left(\hat{\mathbb{Z}}_{p},\mathbb{C}_{q}\right)$.

\vphantom{}

\textbullet{} ($\textrm{IV}\Rightarrow\textrm{I}$) Suppose $f$ has
no zeroes. Since $\mathfrak{y}\mapsto\frac{1}{\mathfrak{y}}$ is a
continuous map on $\mathbb{C}_{q}\backslash\left\{ 0\right\} $, and
since compositions of continuous maps are continuous, to show that
$1/f$ is continuous, it suffices to show that $\left|f\left(\mathfrak{z}\right)\right|_{q}$
is bounded away from $0$. To see this, suppose by way of contradiction
that there was a sequence $\left\{ \mathfrak{z}_{n}\right\} _{n\geq0}\subseteq\mathbb{Z}_{p}$
such that for all $\epsilon>0$, $\left|f\left(\mathfrak{z}_{n}\right)\right|_{q}<\epsilon$
holds for all sufficiently large $n$\textemdash say, for all $n\geq N_{\epsilon}$.
Since $\mathbb{Z}_{p}$ is compact, the $\mathfrak{z}_{n}$s have
a subsequence $\mathfrak{z}_{n_{k}}$ which converges in $\mathbb{Z}_{p}$
to some limit $\mathfrak{z}_{\infty}$. The continuity of $f$ forces
$f\left(\mathfrak{z}_{\infty}\right)=0$, which contradicts the hypothesis
that $f$ was given to have no zeroes.

Thus, if $f$ has no zeroes, $\left|f\left(\mathfrak{z}\right)\right|_{q}$
is bounded away from zero. This prove that $1/f$ is $\left(p,q\right)$-adically
continuous.

Q.E.D.

\vphantom{}

The proof of WTT for $\left(p,q\right)$-adic measures is significantly
more involved than the continuous case. While proving the non-vanishing
of the limit of $\tilde{\mu}_{N}\left(\mathfrak{z}\right)$ when $\hat{\mu}$'s
translates have dense span is as simple as in the continuous case,
the other direction requires a more intricate argument. The idea is
this: by way of contradiction, suppose there is a $\mathfrak{z}_{0}\in\mathbb{Z}_{p}$
so that $\lim_{N\rightarrow\infty}\tilde{\mu}_{N}\left(\mathfrak{z}_{0}\right)$
converges in $\mathbb{C}_{q}$ to zero, yet the span of the translates
of $\hat{\mu}$ \emph{are} dense in $c_{0}\left(\hat{\mathbb{Z}}_{p},\mathbb{C}_{q}\right)$.
We get a contradiction by showing that, with these assumptions, we
can construct two functions $c_{0}\left(\hat{\mathbb{Z}}_{p},\mathbb{C}_{q}\right)$,
both of which produce interesting results when convolved with $\hat{\mu}$.
The first function is constructed so as to give us something close
(in the sense of $c_{0}\left(\hat{\mathbb{Z}}_{p},\mathbb{C}_{q}\right)$'s
norm, the $\left(p,q\right)$-adic sup norm on $\hat{\mathbb{Z}}_{p}$)
to the constant function $0$ when we convolve it with $\hat{\mu}$,
whereas the second is constructed so as to give us something close
to $\mathbf{1}_{0}$ upon convolution with $\hat{\mu}$. Applications
of the ultrametric inequality (a.k.a., the strong triangle inequality)
show that these two estimates tear each other to pieces, and thus
cannot \emph{both }be true; this yields the desired contradiction.

While the idea is straightforward, the approach is somewhat technical,
because our argument will require restricting our attention to bounded
neighborhoods of zero in $\hat{\mathbb{Z}}_{p}$. The contradictory
part of estimates we need follow from delicate manipulations of convolutions
in tandem with these domain restrictions. As such, it will be very
helpful to have a notation for the supremum norm of a function $\hat{\mathbb{Z}}_{p}\rightarrow\mathbb{C}_{q}$
taken over a bounded neighborhood of zero, rather than all of $\hat{\mathbb{Z}}_{p}$.
\begin{defn}
\label{def:norm notation definition}For any integer $n\geq1$ and
any function $\hat{\chi}:\hat{\mathbb{Z}}_{p}\rightarrow\mathbb{C}_{q}$,
we write $\left\Vert \hat{\chi}\right\Vert _{p^{n},q}$ to denote
the non-archimedean norm:
\begin{equation}
\left\Vert \hat{\chi}\right\Vert _{p^{n},q}\overset{\textrm{def}}{=}\sup_{\left|t\right|_{p}\leq p^{n}}\left|\hat{\chi}\left(t\right)\right|_{q}\label{eq:Definition of truncated norm}
\end{equation}
\end{defn}
\begin{rem}
\textbf{WARNING }\textendash{} For the statement and proof of \textbf{Theorem
\ref{thm:pq WTT for measures}} and \textbf{Lemma \ref{lem:three convolution estimate}},
we will write $\left\Vert \cdot\right\Vert _{p^{\infty},q}$ to denote
the $\left(p,q\right)$-adic supremum norm for functions $\hat{\mathbb{Z}}_{p}\rightarrow\mathbb{C}_{q}$.
writing ``$\left\Vert \cdot\right\Vert _{p,q}$'' would potentially
cause confusion with $\left\Vert \cdot\right\Vert _{p^{1},q}$.
\end{rem}
\begin{thm}[\textbf{Wiener Tauberian Theorem for $\left(p,q\right)$-adic Measures}]
\index{Wiener!Tauberian Theorem!left(p,qright)-adic@$\left(p,q\right)$-adic}\label{thm:pq WTT for measures}Let
$d\mu\in C\left(\mathbb{Z}_{p},\mathbb{C}_{q}\right)^{\prime}$. Then,
$\textrm{span}_{\mathbb{C}_{q}}\left\{ \tau_{s}\left\{ \hat{\mu}\right\} \left(t\right):s\in\hat{\mathbb{Z}}_{p}\right\} $
is dense in $c_{0}\left(\hat{\mathbb{Z}}_{p},\mathbb{C}_{q}\right)$
if and only if, for any $\mathfrak{z}\in\mathbb{Z}_{p}$ for which
the limit $\lim_{N\rightarrow\infty}\tilde{\mu}_{N}\left(\mathfrak{z}\right)$
converges in $\mathbb{C}_{q}$, said limit is necessarily non-zero. 
\end{thm}
\vphantom{}

While much of the intricate business involving the intermingling of
restrictions and convolutions occurs in the various claims that structure
our proof of \textbf{Theorem \ref{thm:pq WTT for measures}}, one
particular result in this vein\textemdash the heart of \textbf{Theorem
\ref{thm:pq WTT for measures}}'s proof\textemdash is sufficiently
non-trivial as to merit separate consideration, so as to avoid cluttering
the flow of the argument. This is detailed in the following lemma:
\begin{lem}
\label{lem:three convolution estimate}Let $M,N\in\mathbb{N}_{0}$
be arbitrary, and let $\hat{\phi}\in B\left(\hat{\mathbb{Z}}_{p},\mathbb{C}_{q}\right)$
be supported on $\left|t\right|_{p}\leq p^{N}$. Then, for any $\hat{f},\hat{g}:\hat{\mathbb{Z}}_{p}\rightarrow\mathbb{C}_{q}$
where $\hat{g}$ is supported on $\left|t\right|_{p}\leq p^{M}$,
we have:
\begin{equation}
\left\Vert \hat{\phi}*\hat{f}\right\Vert _{p^{M},q}\leq\left\Vert \hat{\phi}\right\Vert _{p^{\infty},q}\left\Vert \hat{f}\right\Vert _{p^{\max\left\{ M,N\right\} },q}\label{eq:phi_N hat convolve f hat estimate-1}
\end{equation}
\begin{equation}
\left\Vert \hat{\phi}*\hat{f}*\hat{g}\right\Vert _{p^{M},q}\leq\left\Vert \hat{\phi}*\hat{f}\right\Vert _{p^{\max\left\{ M,N\right\} },q}\left\Vert \hat{g}\right\Vert _{p^{\infty},q}\label{eq:phi_N hat convolve f hat convolve g hat estimate-1}
\end{equation}
\end{lem}
Proof: (\ref{eq:phi_N hat convolve f hat estimate-1}) is the easier
of the two:
\begin{align*}
\left\Vert \hat{\phi}*\hat{f}\right\Vert _{p^{M},q} & =\sup_{\left|t\right|_{p}\leq p^{M}}\left|\sum_{s\in\hat{\mathbb{Z}}_{p}}\hat{\phi}\left(s\right)\hat{f}\left(t-s\right)\right|_{q}\\
\left(\hat{\phi}\left(s\right)=0,\textrm{ }\forall\left|s\right|_{p}>p^{N}\right); & \leq\sup_{\left|t\right|_{p}\leq p^{M}}\left|\sum_{\left|s\right|_{p}\leq p^{N}}\hat{\phi}\left(s\right)\hat{f}\left(t-s\right)\right|_{q}\\
\left(\textrm{ultrametric ineq.}\right); & \leq\sup_{\left|t\right|_{p}\leq p^{M}}\sup_{\left|s\right|_{p}\leq p^{N}}\left|\hat{\phi}\left(s\right)\hat{f}\left(t-s\right)\right|_{q}\\
 & \leq\sup_{\left|s\right|_{p}\leq p^{N}}\left|\hat{\phi}\left(s\right)\right|_{q}\sup_{\left|t\right|_{p}\leq p^{\max\left\{ M,N\right\} }}\left|\hat{f}\left(t\right)\right|_{q}\\
 & =\left\Vert \hat{\phi}\right\Vert _{p^{\infty},q}\cdot\left\Vert \hat{f}\right\Vert _{p^{\max\left\{ M,N\right\} },q}
\end{align*}
and we are done.

Proving (\ref{eq:phi_N hat convolve f hat convolve g hat estimate-1})
is similar, but more involved. We start by writing out the convolution
of $\hat{\phi}*\hat{f}$ and $\hat{g}$: 
\begin{align*}
\left\Vert \hat{\phi}*\hat{f}*\hat{g}\right\Vert _{p^{M},q} & =\sup_{\left|t\right|_{p}\leq p^{M}}\left|\sum_{s\in\hat{\mathbb{Z}}_{p}}\left(\hat{\phi}*\hat{f}\right)\left(t-s\right)\hat{g}\left(s\right)\right|_{q}\\
 & \leq\sup_{\left|t\right|_{p}\leq p^{M}}\sup_{s\in\hat{\mathbb{Z}}_{p}}\left|\left(\hat{\phi}*\hat{f}\right)\left(t-s\right)\hat{g}\left(s\right)\right|_{q}\\
\left(\hat{g}\left(s\right)=0,\textrm{ }\forall\left|s\right|_{p}>p^{M}\right); & \leq\sup_{\left|t\right|_{p}\leq p^{M}}\sup_{\left|s\right|_{p}\leq p^{M}}\left|\left(\hat{\phi}*\hat{f}\right)\left(t-s\right)\hat{g}\left(s\right)\right|_{q}\\
 & \leq\left\Vert \hat{g}\right\Vert _{p^{\infty},q}\sup_{\left|t\right|_{p},\left|s\right|_{p}\leq p^{M}}\left|\left(\hat{\phi}*\hat{f}\right)\left(t-s\right)\right|_{q}\\
\left(\textrm{write out }\hat{\phi}*\hat{f}\right); & =\left\Vert \hat{g}\right\Vert _{p^{\infty},q}\sup_{\left|t\right|_{p},\left|s\right|_{p}\leq p^{M}}\left|\sum_{\tau\in\hat{\mathbb{Z}}_{p}}\hat{\phi}\left(t-s-\tau\right)\hat{f}\left(\tau\right)\right|_{q}\\
\left(\textrm{let }u=s+\tau\right); & =\left\Vert \hat{g}\right\Vert _{p^{\infty},q}\sup_{\left|t\right|_{p},\left|s\right|_{p}\leq p^{M}}\left|\sum_{u-s\in\hat{\mathbb{Z}}_{p}}\hat{\phi}\left(t-u\right)\hat{f}\left(u-s\right)\right|_{q}\\
\left(s+\hat{\mathbb{Z}}_{p}=\hat{\mathbb{Z}}_{p}\right); & =\left\Vert \hat{g}\right\Vert _{p^{\infty},q}\sup_{\left|t\right|_{p},\left|s\right|_{p}\leq p^{M}}\left|\sum_{u\in\hat{\mathbb{Z}}_{p}}\hat{\phi}\left(t-u\right)\hat{f}\left(u-s\right)\right|_{q}
\end{align*}

Now we use the fact that $\hat{\phi}\left(t-u\right)$ vanishes for
all $\left|t-u\right|_{p}>p^{N}$. Because $t$ is restricted to $\left|t\right|_{p}\leq p^{M}$,
observe that for $\left|u\right|_{p}>p^{\max\left\{ M,N\right\} }$,
the ultrametric inequality allows us to write: 
\begin{equation}
\left|t-u\right|_{p}=\max\left\{ \left|t\right|_{p},\left|u\right|_{p}\right\} >p^{\max\left\{ M,N\right\} }>p^{N}
\end{equation}
So, for $\left|t\right|_{p},\left|s\right|_{p}\leq p^{M}$, the summand
$\hat{\phi}\left(t-u\right)\hat{f}\left(u-s\right)$ vanishes whenever
$\left|u\right|_{p}>p^{\max\left\{ M,N\right\} }$. This gives us:
\begin{equation}
\left\Vert \hat{\phi}*\hat{f}*\hat{g}\right\Vert _{p^{M},q}\leq\left\Vert \hat{g}\right\Vert _{p^{\infty},q}\sup_{\left|t\right|_{p},\left|s\right|_{p}\leq p^{M}}\left|\sum_{\left|u\right|_{p}\leq p^{\max\left\{ M,N\right\} }}\hat{\phi}\left(t-u\right)\hat{f}\left(u-s\right)\right|_{q}
\end{equation}

Next, we expand the range of $\left|s\right|_{p}$ and $\left|t\right|_{p}$
from $\leq p^{M}$ in $p$-adic absolute value to $\leq p^{\max\left\{ M,N\right\} }$
in $p$-adic absolute value: 
\begin{equation}
\left\Vert \hat{\phi}*\hat{f}*\hat{g}\right\Vert _{p^{M},q}\leq\left\Vert \hat{g}\right\Vert _{p^{\infty},q}\sup_{\left|t\right|_{p},\left|s\right|_{p}\leq p^{\max\left\{ M,N\right\} }}\left|\sum_{\left|u\right|_{p}\leq p^{\max\left\{ M,N\right\} }}\hat{\phi}\left(t-u\right)\hat{f}\left(u-s\right)\right|_{q}
\end{equation}
In doing so, we have put $s$, $t$, and $u$ all in the same $p$-adic
neighborhood of $0$: 
\begin{equation}
\left\{ x\in\hat{\mathbb{Z}}_{p}:\left|x\right|_{p}\leq p^{\max\left\{ M,N\right\} }\right\} 
\end{equation}
We do this because this set is closed under addition: for any $\left|s\right|_{p}\leq p^{\max\left\{ M,N\right\} }$,
our $u$-sum is invariant under the change of variables $u\mapsto u+s$,
and we obtain: 
\begin{align*}
\left\Vert \hat{\phi}*\hat{f}*\hat{g}\right\Vert _{p^{M},q} & \leq\left\Vert \hat{g}\right\Vert _{p^{\infty},q}\sup_{\left|t\right|_{p},\left|s\right|_{p}\leq p^{\max\left\{ M,N\right\} }}\left|\sum_{\left|u\right|_{p}\leq p^{\max\left\{ M,N\right\} }}\hat{\phi}\left(t-\left(u+s\right)\right)\hat{f}\left(u\right)\right|_{q}\\
 & =\left\Vert \hat{g}\right\Vert _{p^{\infty},q}\sup_{\left|t\right|_{p},\left|s\right|_{p}\leq p^{\max\left\{ M,N\right\} }}\left|\sum_{\left|u\right|_{p}\leq p^{\max\left\{ M,N\right\} }}\hat{\phi}\left(t-s-u\right)\hat{f}\left(u\right)\right|_{q}
\end{align*}

Finally, observing that: 
\begin{equation}
\left\{ t-s:\left|t\right|_{p},\left|s\right|_{p}\leq p^{\max\left\{ M,N\right\} }\right\} =\left\{ t:\left|t\right|_{p}\leq p^{\max\left\{ M,N\right\} }\right\} 
\end{equation}
we can write: 
\begin{align*}
\left\Vert \hat{\phi}*\hat{f}*\hat{g}\right\Vert _{p^{M},q} & \leq\left\Vert \hat{g}\right\Vert _{p^{\infty},q}\sup_{\left|t\right|_{p}\leq p^{\max\left\{ M,N\right\} }}\left|\sum_{\left|u\right|_{p}\leq p^{\max\left\{ M,N\right\} }}\hat{\phi}\left(t-u\right)\hat{f}\left(u\right)\right|_{q}\\
 & \leq\left\Vert \hat{g}\right\Vert _{p^{\infty},q}\sup_{\left|t\right|_{p}\leq p^{\max\left\{ M,N\right\} }}\left|\underbrace{\sum_{u\in\hat{\mathbb{Z}}_{p}}\hat{\phi}\left(t-u\right)\hat{f}\left(u\right)}_{\hat{\phi}*\hat{f}}\right|_{q}\\
 & =\left\Vert \hat{g}\right\Vert _{p^{\infty},q}\sup_{\left|t\right|_{p}\leq p^{\max\left\{ M,N\right\} }}\left|\left(\hat{\phi}*\hat{f}\right)\left(u\right)\right|_{q}\\
\left(\textrm{by definition}\right); & =\left\Vert \hat{g}\right\Vert _{p^{\infty},q}\left\Vert \hat{\phi}*\hat{f}\right\Vert _{p^{\max\left\{ M,N\right\} },q}
\end{align*}
This proves the desired estimate, and with it, the rest of the Lemma.

Q.E.D.

\vphantom{}

\textbf{Proof} \textbf{of Theorem \ref{thm:pq WTT for measures}}:
We start with the simpler of the two directions.

I. Suppose the span of the translates of $\hat{\mu}$ are dense. Just
as in the proof of the WTT for continuous $\left(p,q\right)$-adic
functions, we let $\epsilon\in\left(0,1\right)$ and then choose $\mathfrak{c}_{m}$s
and $t_{m}$s so that: 
\begin{equation}
\sup_{t\in\hat{\mathbb{Z}}_{p}}\left|\mathbf{1}_{0}\left(t\right)-\sum_{m=1}^{M}\mathfrak{c}_{m}\hat{\mu}\left(t-t_{m}\right)\right|_{q}<\epsilon
\end{equation}
Picking sufficiently large $N$, we obtain: 
\begin{align*}
\left|1-\left(\sum_{m=1}^{M}\mathfrak{c}_{m}e^{2\pi i\left\{ t_{m}\mathfrak{z}\right\} _{p}}\right)\tilde{\mu}_{N}\left(\mathfrak{z}\right)\right|_{q} & \leq\max_{\left|t\right|_{p}\leq p^{N}}\left|\mathbf{1}_{0}\left(t\right)-\sum_{m=1}^{M}\mathfrak{c}_{m}\hat{\mu}\left(t-t_{m}\right)\right|_{q}\\
 & \leq\sup_{t\in\hat{\mathbb{Z}}_{p}}\left|\mathbf{1}_{0}\left(t\right)-\sum_{m=1}^{M}\mathfrak{c}_{m}\hat{\mu}\left(t-t_{m}\right)\right|_{q}\\
 & <\epsilon
\end{align*}

Now, let $\mathfrak{z}_{0}\in\mathbb{Z}_{p}$ be a point so that $\mathfrak{L}\overset{\textrm{def}}{=}\lim_{N\rightarrow\infty}\tilde{\mu}_{N}\left(\mathfrak{z}_{0}\right)$
converges in $\mathbb{C}_{q}$. We need to show $\mathfrak{L}\neq0$.
To do this, plugging $\mathfrak{z}=\mathfrak{z}_{0}$ into the above
yields: 
\[
\epsilon>\lim_{N\rightarrow\infty}\left|1-\left(\sum_{m=1}^{M}\mathfrak{c}_{m}e^{2\pi i\left\{ t_{m}\mathfrak{z}_{0}\right\} _{p}}\right)\tilde{\mu}_{N}\left(\mathfrak{z}_{0}\right)\right|_{q}=\left|1-\left(\sum_{m=1}^{M}\mathfrak{c}_{m}e^{2\pi i\left\{ t_{m}\mathfrak{z}_{0}\right\} _{p}}\right)\cdot\mathfrak{L}\right|_{q}
\]
If $\mathfrak{L}=0$, the right-most expression will be $1$, and
hence, we get $\epsilon>1$, but this is impossible; $\epsilon$ was
given to be less than $1$. So, if $\lim_{N\rightarrow\infty}\tilde{\mu}_{N}\left(\mathfrak{z}_{0}\right)$
converges in $\mathbb{C}_{q}$ to $\mathfrak{L}$, $\mathfrak{L}$
must be non-zero.

\vphantom{}

II. Let $\mathfrak{z}_{0}\in\mathbb{Z}_{p}$ be a zero of $\tilde{\mu}\left(\mathfrak{z}\right)$
such that $\tilde{\mu}_{N}\left(\mathfrak{z}_{0}\right)\rightarrow0$
in $\mathbb{C}_{q}$ as $N\rightarrow\infty$. Then, by way of contradiction,
suppose the span of the translates of $\hat{\mu}$ is dense in $c_{0}\left(\hat{\mathbb{Z}}_{p},\mathbb{C}_{q}\right)$,
despite the zero of $\tilde{\mu}$ at $\mathfrak{z}_{0}$. Thus, we
can use linear combinations of translates of $\hat{\mu}$ approximate
any function in $c_{0}\left(\hat{\mathbb{Z}}_{p},\mathbb{C}_{q}\right)$'s
sup norm. In particular, we choose to approximate $\mathbf{1}_{0}$:
letting $\epsilon\in\left(0,1\right)$ be arbitrary, there is then
a choice of $\mathfrak{c}_{k}$s in $\mathbb{C}_{q}$ and $t_{k}$s
in $\hat{\mathbb{Z}}_{p}$ so that: 
\begin{equation}
\sup_{t\in\hat{\mathbb{Z}}_{p}}\left|\mathbf{1}_{0}\left(t\right)-\sum_{k=1}^{K}\mathfrak{c}_{k}\hat{\mu}\left(t-t_{k}\right)\right|_{q}<\epsilon
\end{equation}
Letting: 
\begin{equation}
\hat{\eta}_{\epsilon}\left(t\right)\overset{\textrm{def}}{=}\sum_{k=1}^{K}\mathfrak{c}_{k}\mathbf{1}_{t_{k}}\left(t\right)\label{eq:Definition of eta_epsilon hat}
\end{equation}
we can express the above linear combination as the convolution:

\begin{equation}
\left(\hat{\mu}*\hat{\eta}_{\epsilon}\right)\left(t\right)=\sum_{\tau\in\hat{\mathbb{Z}}_{p}}\hat{\mu}\left(t-\tau\right)\sum_{k=1}^{K}\mathfrak{c}_{k}\mathbf{1}_{t_{k}}\left(\tau\right)=\sum_{k=1}^{K}\mathfrak{c}_{k}\hat{\mu}\left(t-t_{k}\right)
\end{equation}
So, we get:
\begin{equation}
\left\Vert \mathbf{1}_{0}-\hat{\mu}*\hat{\eta}_{\epsilon}\right\Vert _{p^{\infty},q}\overset{\textrm{def}}{=}\sup_{t\in\hat{\mathbb{Z}}_{p}}\left|\mathbf{1}_{0}\left(t\right)-\left(\hat{\mu}*\hat{\eta}_{\epsilon}\right)\left(t\right)\right|_{q}<\epsilon\label{eq:Converse WTT - eq. 1}
\end{equation}
Before proceeding any further, it is vital to note that we can (and
must) assume that $\hat{\eta}_{\epsilon}$ is not identically\footnote{This must hold whenever $\epsilon\in\left(0,1\right)$, because\textemdash were
$\hat{\eta}_{\epsilon}$ identically zero\textemdash the supremum:
\[
\sup_{t\in\hat{\mathbb{Z}}_{p}}\left|\mathbf{1}_{0}\left(t\right)-\left(\hat{\mu}*\hat{\eta}_{\epsilon}\right)\left(t\right)\right|_{q}
\]
would then be equal to $1$, rather than $<\epsilon$.} $0$.

That being done, equation (\ref{eq:Converse WTT - eq. 1}) shows the
assumption we want to contradict: the existence of a function which
produces something close to $\mathbf{1}_{0}$ after convolution with
$\hat{\mu}$. We will arrive at our contradiction by showing that
the zero of $\tilde{\mu}$ at $\mathfrak{z}_{0}$ allows us to construct
a second function which yields something close to $0$ after convolution
with $\hat{\mu}$. By convolving \emph{both} of these functions with
$\hat{\mu}$, we will end up with something which is close to both
$\mathbf{1}_{0}$ \emph{and }close to $0$, which is, of course, impossible.

Our make-things-close-to-zero-by-convolution function is going to
be: 
\begin{equation}
\hat{\phi}_{N}\left(t\right)\overset{\textrm{def}}{=}\mathbf{1}_{0}\left(p^{N}t\right)e^{-2\pi i\left\{ t\mathfrak{z}_{0}\right\} _{p}},\textrm{ }\forall t\in\hat{\mathbb{Z}}_{p}\label{eq:Definition of Phi_N hat}
\end{equation}
Note that $\hat{\phi}_{N}\left(t\right)$ is only supported for $\left|t\right|_{p}\leq p^{N}$.
Now, as defined, we have that: 
\begin{align*}
\left(\hat{\mu}*\hat{\phi}_{N}\right)\left(\tau\right) & =\sum_{s\in\hat{\mathbb{Z}}_{p}}\hat{\mu}\left(\tau-s\right)\hat{\phi}_{N}\left(s\right)\\
\left(\hat{\phi}_{N}\left(s\right)=0,\textrm{ }\forall\left|s\right|_{p}>p^{N}\right); & =\sum_{\left|s\right|_{p}\leq p^{N}}\hat{\mu}\left(\tau-s\right)\hat{\phi}_{N}\left(s\right)\\
 & =\sum_{\left|s\right|_{p}\leq p^{N}}\hat{\mu}\left(\tau-s\right)e^{-2\pi i\left\{ s\mathfrak{z}_{0}\right\} _{p}}
\end{align*}
Fixing $\tau$, observe that the map $s\mapsto\tau-s$ is a bijection
of the set $\left\{ s\in\hat{\mathbb{Z}}_{p}:\left|s\right|_{p}\leq p^{N}\right\} $
whenever $\left|\tau\right|_{p}\leq p^{N}$. So, for any $N\geq-v_{p}\left(\tau\right)$,
we obtain: 
\begin{align*}
\left(\hat{\mu}*\hat{\phi}_{N}\right)\left(\tau\right) & =\sum_{\left|s\right|_{p}\leq p^{N}}\hat{\mu}\left(\tau-s\right)e^{-2\pi i\left\{ s\mathfrak{z}_{0}\right\} _{p}}\\
 & =\sum_{\left|s\right|_{p}\leq p^{N}}\hat{\mu}\left(s\right)e^{-2\pi i\left\{ \left(\tau-s\right)\mathfrak{z}_{0}\right\} _{p}}\\
 & =e^{-2\pi i\left\{ \tau\mathfrak{z}_{0}\right\} _{p}}\sum_{\left|s\right|_{p}\leq p^{N}}\hat{\mu}\left(s\right)e^{2\pi i\left\{ s\mathfrak{z}_{0}\right\} _{p}}\\
 & =e^{-2\pi i\left\{ \tau\mathfrak{z}_{0}\right\} _{p}}\tilde{\mu}_{N}\left(\mathfrak{z}_{0}\right)
\end{align*}
Since this holds for all $N\geq-v_{p}\left(\tau\right)$, upon letting
$N\rightarrow\infty$, $\tilde{\mu}_{N}\left(\mathfrak{z}_{0}\right)$
converges to $0$ in \emph{$\mathbb{C}_{q}$}, by our assumption.
So, for any $\epsilon^{\prime}>0$, there exists an $N_{\epsilon^{\prime}}$
so that $\left|\tilde{\mu}_{N}\left(\mathfrak{z}_{0}\right)\right|_{q}<\epsilon^{\prime}$
for all $N\geq N_{\epsilon^{\prime}}$. Combining this with the above
computation (after taking $q$-adic absolute values), we have established
the following:
\begin{claim}
\label{claim:phi_N hat claim}Let $\epsilon^{\prime}>0$ be arbitrary.
Then, there exists an $N_{\epsilon^{\prime}}\geq0$ (depending only
on $\hat{\mu}$ and $\epsilon^{\prime}$) so that, for all $\tau\in\hat{\mathbb{Z}}_{p}$:
\begin{equation}
\left|\left(\hat{\mu}*\hat{\phi}_{N}\right)\left(\tau\right)\right|_{q}=\left|e^{-2\pi i\left\{ \tau\mathfrak{z}_{0}\right\} _{p}}\tilde{\mu}_{N}\left(\mathfrak{z}_{0}\right)\right|_{q}<\epsilon^{\prime},\textrm{ }\forall N\geq\max\left\{ N_{\epsilon^{\prime}},-v_{p}\left(\tau\right)\right\} ,\tau\in\hat{\mathbb{Z}}_{p}\label{eq:WTT - First Claim}
\end{equation}
\end{claim}
\vphantom{}

As stated, the idea is to convolve $\hat{\mu}*\hat{\phi}_{N}$ with
$\hat{\eta}_{\epsilon}$ so as to obtain a function (via the associativity
of convolution) which is both close to $0$ and close to $\mathbf{1}_{0}$.
However, our present situation is less than ideal because the lower
bound on $N$ in \textbf{Claim \ref{claim:phi_N hat claim}} depends
on $\tau$, and so, the convergence of $\left|\left(\hat{\mu}*\hat{\phi}_{N}\right)\left(\tau\right)\right|_{q}$
to $0$ as $N\rightarrow\infty$ will \emph{not }be uniform in $\tau$.
This is where the difficulty of this direction of the proof lies.
To overcome this obstacle, instead of convolving $\hat{\mu}*\hat{\phi}_{N}$
with $\hat{\eta}_{\epsilon}$, we will convolve $\hat{\mu}*\hat{\phi}_{N}$
with a truncated version of $\hat{\eta}_{\epsilon}$, whose support
has been restricted to a finite subset of $\hat{\mathbb{Z}}_{p}$.
This is the function $\hat{\eta}_{\epsilon,M}:\hat{\mathbb{Z}}_{p}\rightarrow\mathbb{C}_{q}$
given by: 
\begin{equation}
\hat{\eta}_{\epsilon,M}\left(t\right)\overset{\textrm{def}}{=}\mathbf{1}_{0}\left(p^{M}t\right)\hat{\eta}_{\epsilon}\left(t\right)=\begin{cases}
\hat{\eta}_{\epsilon}\left(t\right) & \textrm{if }\left|t\right|_{p}\leq p^{M}\\
0 & \textrm{else}
\end{cases}\label{eq:Definition of eta epsilon M hat}
\end{equation}
where $M\geq0$ is arbitrary. With this, we can attain the desired
``close to both $0$ and $\mathbf{1}_{0}$'' contradiction for $\hat{\mu}*\hat{\phi}_{N}*\hat{\eta}_{\epsilon,M}$.
The proof will be completed upon demonstrating that this contradiction
will remain even as $M\rightarrow\infty$.

The analogue of the estimate (\ref{eq:Converse WTT - eq. 1}) for
this truncated case is:
\begin{align*}
\left(\hat{\mu}*\hat{\eta}_{\epsilon,M}\right)\left(t\right) & =\sum_{\tau\in\hat{\mathbb{Z}}_{p}}\hat{\mu}\left(t-\tau\right)\mathbf{1}_{0}\left(p^{M}\tau\right)\sum_{k=1}^{K}\mathfrak{c}_{k}\mathbf{1}_{t_{k}}\left(\tau\right)\\
 & =\sum_{\left|\tau\right|_{p}\leq p^{M}}\hat{\mu}\left(t-\tau\right)\sum_{k=1}^{K}\mathfrak{c}_{k}\mathbf{1}_{t_{k}}\left(\tau\right)\\
 & =\sum_{k:\left|t_{k}\right|_{p}\leq p^{M}}\mathfrak{c}_{k}\hat{\mu}\left(t-t_{k}\right)
\end{align*}
where, \emph{note}, instead of summing over all the $t_{k}$s that
came with $\hat{\eta}_{\epsilon}$, we only sum over those $t_{k}$s
in the set $\left\{ t\in\hat{\mathbb{Z}}_{p}:\left|t\right|_{p}\leq p^{M}\right\} $.
Because the $t_{k}$s came with the un-truncated $\hat{\eta}_{\epsilon}$,
observe that we can make $\hat{\mu}*\hat{\eta}_{\epsilon,M}$ \emph{equal}
to $\hat{\mu}*\hat{\eta}_{\epsilon}$ by simply choosing $M$ to be
large enough so that all the $t_{k}$s lie in $\left\{ t\in\hat{\mathbb{Z}}_{p}:\left|t\right|_{p}\leq p^{M}\right\} $.
The lower bound on such $M$s is given by: 
\begin{equation}
M_{0}\overset{\textrm{def}}{=}\max\left\{ -v_{p}\left(t_{1}\right),\ldots,-v_{p}\left(t_{K}\right)\right\} \label{eq:WTT - Choice for M_0}
\end{equation}
Then, we have: 
\begin{equation}
\left(\hat{\mu}*\hat{\eta}_{\epsilon,M}\right)\left(t\right)=\sum_{k=1}^{K}\mathfrak{c}_{k}\hat{\mu}\left(t-t_{k}\right)=\left(\hat{\mu}*\hat{\eta}_{\epsilon}\right)\left(t\right),\textrm{ }\forall M\geq M_{0},\textrm{ }\forall t\in\hat{\mathbb{Z}}_{p}\label{eq:Effect of M bigger than M0}
\end{equation}
So, applying $\left\Vert \cdot\right\Vert _{p^{M},q}$ norm: 
\begin{align*}
\left\Vert \mathbf{1}_{0}-\hat{\mu}*\hat{\eta}_{\epsilon,M}\right\Vert _{p^{M},q} & =\sup_{\left|t\right|_{p}\leq p^{M}}\left|\mathbf{1}_{0}\left(t\right)-\left(\hat{\mu}*\hat{\eta}_{\epsilon,M}\right)\left(t\right)\right|_{q}\\
\left(\textrm{if }M\geq M_{0}\right); & =\sup_{\left|t\right|_{p}\leq p^{M}}\left|\mathbf{1}_{0}\left(t\right)-\left(\hat{\mu}*\hat{\eta}_{\epsilon}\right)\left(t\right)\right|_{q}\\
 & \leq\left\Vert \mathbf{1}_{0}-\left(\hat{\mu}*\hat{\eta}_{\epsilon}\right)\right\Vert _{p^{\infty},q}\\
 & <\epsilon
\end{align*}
Finally, note that\emph{ $M_{0}$ }depends on\emph{ only $\hat{\eta}_{\epsilon}$
}and\emph{ $\epsilon$}. \textbf{Claim \ref{claim:truncated convolution estimates for 1d WTT}},
given below, summarizes these findings:
\begin{claim}
\label{claim:truncated convolution estimates for 1d WTT}Let $\epsilon>0$
be arbitrary. Then, there exists an integer $M_{0}$ depending only
on $\hat{\eta}_{\epsilon}$ and $\epsilon$, so that:

I. $\hat{\mu}*\hat{\eta}_{\epsilon,M}=\hat{\mu}*\hat{\eta}_{\epsilon},\textrm{ }\forall M\geq M_{0}$.

\vphantom{}

II. $\left\Vert \mathbf{1}_{0}-\hat{\mu}*\hat{\eta}_{\epsilon,M}\right\Vert _{p^{\infty},q}<\epsilon,\textrm{ }\forall M\geq M_{0}$.

\vphantom{}

III. $\left\Vert \mathbf{1}_{0}-\hat{\mu}*\hat{\eta}_{\epsilon,M}\right\Vert _{p^{M},q}<\epsilon,\textrm{ }\forall M\geq M_{0}$.
\end{claim}
\vphantom{}

The next step is to refine our ``make $\hat{\mu}$ close to zero''
estimate by taking into account $\left\Vert \cdot\right\Vert _{p^{m},q}$. 
\begin{claim}
\label{claim:p^m, q norm of convolution of mu-hat and phi_N hat }Let
$\epsilon^{\prime}>0$. Then, there exists $N_{\epsilon^{\prime}}$
depending only on $\epsilon^{\prime}$ and $\hat{\mu}$ so that: 
\begin{equation}
\left\Vert \hat{\phi}_{N}*\hat{\mu}\right\Vert _{p^{m},q}<\epsilon^{\prime},\textrm{ }\forall m\geq1,\textrm{ }\forall N\geq\max\left\{ N_{\epsilon^{\prime}},m\right\} \label{eq:WTT - eq. 4}
\end{equation}
Proof of claim: Let $\epsilon^{\prime}>0$. \textbf{Claim \ref{claim:phi_N hat claim}}
tells us that there is an $N_{\epsilon^{\prime}}$ (depending on $\epsilon^{\prime}$,
$\hat{\mu}$) so that: 
\begin{equation}
\left|\left(\hat{\phi}_{N}*\hat{\mu}\right)\left(\tau\right)\right|_{q}<\epsilon^{\prime},\textrm{ }\forall N\geq\max\left\{ N_{\epsilon^{\prime}},-v_{p}\left(\tau\right)\right\} ,\textrm{ }\forall\tau\in\hat{\mathbb{Z}}_{p}
\end{equation}
So, letting $m\geq1$ be arbitrary, note that $\left|\tau\right|_{p}\leq p^{m}$
implies $-v_{p}\left(\tau\right)\leq m$. As such, we can make the
result of \textbf{Claim \ref{claim:phi_N hat claim}} hold for all
$\left|\tau\right|_{p}\leq p^{m}$ by choosing $N\geq\max\left\{ N_{\epsilon^{\prime}},m\right\} $:

\begin{equation}
\underbrace{\sup_{\left|\tau\right|_{p}\leq p^{m}}\left|\left(\hat{\phi}_{N}*\hat{\mu}\right)\left(\tau\right)\right|_{q}}_{\left\Vert \hat{\mu}*\hat{\phi}_{N}\right\Vert _{p^{m},q}}<\epsilon^{\prime},\textrm{ }\forall N\geq\max\left\{ N_{\epsilon^{\prime}},m\right\} 
\end{equation}
This proves the claim. 
\end{claim}
\vphantom{}

Using \textbf{Lemma \ref{lem:three convolution estimate}}, we can
now set up the string of estimates we need to arrive at the desired
contradiction. First, let us choose an $\epsilon\in\left(0,1\right)$
and a function $\hat{\eta}_{\epsilon}:\hat{\mathbb{Z}}_{p}\rightarrow\mathbb{C}_{q}$
which is not identically zero, so that: 
\begin{equation}
\left\Vert \mathbf{1}_{0}-\hat{\mu}*\hat{\eta}_{\epsilon}\right\Vert _{p^{\infty},q}<\epsilon
\end{equation}
Then, by \textbf{Claim \ref{claim:truncated convolution estimates for 1d WTT}},
we can choose a $M_{\epsilon}$ depending only on $\epsilon$ and
$\hat{\eta}_{\epsilon}$ so that: 
\begin{equation}
M\geq M_{\epsilon}\Rightarrow\left\Vert \mathbf{1}_{0}-\hat{\mu}*\hat{\eta}_{\epsilon,M}\right\Vert _{p^{M},q}<\epsilon\label{eq:Thing to contradict}
\end{equation}
This shows that $\hat{\mu}$ can be convolved to become close to $\mathbf{1}_{0}$.
The contradiction will follow from convolving this with $\hat{\phi}_{N}$,
where $N$, at this point, is arbitrary: 
\begin{equation}
\left\Vert \hat{\phi}_{N}*\left(\mathbf{1}_{0}-\left(\hat{\mu}*\hat{\eta}_{\epsilon,M}\right)\right)\right\Vert _{p^{M},q}=\left\Vert \hat{\phi}_{N}-\left(\hat{\phi}_{N}*\hat{\mu}*\hat{\eta}_{\epsilon,M}\right)\right\Vert _{p^{M},q}\label{eq:WTT - Target of attack}
\end{equation}
Our goal here is to show that (\ref{eq:WTT - Target of attack}) is
close to both $0$ and $1$ simultaneously.

First, writing: 
\begin{equation}
\left\Vert \hat{\phi}_{N}*\left(\mathbf{1}_{0}-\left(\hat{\mu}*\hat{\eta}_{\epsilon,M}\right)\right)\right\Vert _{p^{M},q}
\end{equation}
the fact that $\hat{\phi}_{N}\left(t\right)$ is supported on $\left|t\right|_{p}\leq p^{N}$
allows us to apply equation (\ref{eq:phi_N hat convolve f hat estimate})
from \textbf{Lemma \ref{lem:three convolution estimate}}. This gives
the estimate:
\begin{equation}
\left\Vert \hat{\phi}_{N}*\left(\mathbf{1}_{0}-\left(\hat{\mu}*\hat{\eta}_{\epsilon,M}\right)\right)\right\Vert _{p^{M},q}\leq\underbrace{\left\Vert \hat{\phi}_{N}\right\Vert _{p^{\infty},q}}_{1}\left\Vert \mathbf{1}_{0}-\left(\hat{\mu}*\hat{\eta}_{\epsilon,M}\right)\right\Vert _{p^{\max\left\{ M,N\right\} },q}
\end{equation}
for all $M$ and $N$. Letting $M\geq M_{\epsilon}$, we can apply
\textbf{Claim \ref{claim:truncated convolution estimates for 1d WTT}}
and write $\hat{\mu}*\hat{\eta}_{\epsilon,M}=\hat{\mu}*\hat{\eta}_{\epsilon}$.
So, for $M\geq M_{\epsilon}$ and arbitrary $N$, we have: 
\begin{align*}
\left\Vert \hat{\phi}_{N}*\left(\mathbf{1}_{0}-\left(\hat{\mu}*\hat{\eta}_{\epsilon,M}\right)\right)\right\Vert _{p^{M},q} & \leq\left\Vert \mathbf{1}_{0}-\left(\hat{\mu}*\hat{\eta}_{\epsilon}\right)\right\Vert _{p^{\max\left\{ M,N\right\} },q}\\
 & \leq\left\Vert \mathbf{1}_{0}-\left(\hat{\mu}*\hat{\eta}_{\epsilon}\right)\right\Vert _{p^{\infty},q}\\
\left(\textrm{\textbf{Claim \ensuremath{\ref{claim:truncated convolution estimates for 1d WTT}}}}\right); & <\epsilon
\end{align*}
Thus, the \emph{left-hand side} of (\ref{eq:WTT - Target of attack})
is $<\epsilon$. This is the first estimate.

Keeping $M\geq M_{\epsilon}$ and $N$ arbitrary, we will obtain the
desired contradiction by showing that the \emph{right-hand side} of
(\ref{eq:WTT - Target of attack}) is $>\epsilon$. Since $\left\Vert \cdot\right\Vert _{p^{M},q}$
is a non-archimedean norm, it satisfies the ultrametric inequality.
Applying this to the right-hand side of (\ref{eq:WTT - Target of attack})
yields:
\begin{equation}
\left\Vert \hat{\phi}_{N}-\left(\hat{\phi}_{N}*\hat{\mu}*\hat{\eta}_{\epsilon,M}\right)\right\Vert _{p^{M},q}\leq\max\left\{ \left\Vert \hat{\phi}_{N}\right\Vert _{p^{M},q},\left\Vert \hat{\phi}_{N}*\hat{\mu}*\hat{\eta}_{\epsilon,M}\right\Vert _{p^{M},q}\right\} \label{eq:WTT - Ultrametric inequality}
\end{equation}
Because $\hat{\eta}_{\epsilon,M}\left(t\right)$ and $\hat{\phi}_{N}$
are supported on $\left|t\right|_{p}\leq p^{M}$ and $\left|t\right|_{p}\leq p^{N}$,
respectively, we can apply (\ref{eq:phi_N hat convolve f hat convolve g hat estimate})
from \textbf{Lemma \ref{lem:three convolution estimate}} and write:
\begin{equation}
\left\Vert \hat{\phi}_{N}*\hat{\mu}*\hat{\eta}_{\epsilon,M}\right\Vert _{p^{M},q}\leq\left\Vert \hat{\phi}_{N}*\hat{\mu}\right\Vert _{p^{\max\left\{ M,N\right\} },q}\cdot\left\Vert \hat{\eta}_{\epsilon}\right\Vert _{p^{\infty},q}\label{eq:Ready for epsilon prime}
\end{equation}

Since $\hat{\eta}_{\epsilon}$ was given to \emph{not }be identically
zero, the quantity $\left\Vert \hat{\eta}_{\epsilon}\right\Vert _{p^{\infty},q}$
must be positive. Consequently, for our given $\epsilon\in\left(0,1\right)$,
let us define $\epsilon^{\prime}$ by: 
\begin{equation}
\epsilon^{\prime}=\frac{\epsilon}{2\left\Vert \hat{\eta}_{\epsilon}\right\Vert _{p^{\infty},q}}
\end{equation}
So far, $N$ is still arbitrary. By \textbf{Claim \ref{claim:p^m, q norm of convolution of mu-hat and phi_N hat }},
for this $\epsilon^{\prime}$, we know there exists an $N_{\epsilon^{\prime}}$
(depending only on $\epsilon$, $\hat{\eta}_{\epsilon}$, and $\hat{\mu}$)
so that:\textbf{
\begin{equation}
\left\Vert \hat{\phi}_{N}*\hat{\mu}\right\Vert _{p^{m},q}<\epsilon^{\prime},\textrm{ }\forall m\geq1,\textrm{ }\forall N\geq\max\left\{ N_{\epsilon^{\prime}},m\right\} 
\end{equation}
}Choosing $m=N_{\epsilon^{\prime}}$ gives us: 
\begin{equation}
N\geq N_{\epsilon^{\prime}}\Rightarrow\left\Vert \hat{\phi}_{N}*\hat{\mu}\right\Vert _{p^{N},q}<\epsilon^{\prime}
\end{equation}
So, choose $N\geq\max\left\{ N_{\epsilon^{\prime}},M\right\} $. Then,
$\max\left\{ M,N\right\} =N$, and so (\ref{eq:Ready for epsilon prime})
becomes:
\begin{equation}
\left\Vert \hat{\phi}_{N}*\hat{\mu}*\hat{\eta}_{\epsilon,M}\right\Vert _{p^{M},q}\leq\left\Vert \hat{\phi}_{N}*\hat{\mu}\right\Vert _{p^{N},q}\cdot\left\Vert \hat{\eta}_{\epsilon}\right\Vert _{p^{\infty},q}<\frac{\epsilon}{2\left\Vert \hat{\eta}_{\epsilon}\right\Vert _{p^{\infty},q}}\cdot\left\Vert \hat{\eta}_{\epsilon}\right\Vert _{p^{\infty},q}=\frac{\epsilon}{2}
\end{equation}
Since $\left\Vert \hat{\phi}_{N}\right\Vert _{p^{M},q}=1$ for all
$M,N\geq0$, this shows that: 
\begin{equation}
\left\Vert \hat{\phi}_{N}*\hat{\mu}*\hat{\eta}_{\epsilon,M}\right\Vert _{p^{M},q}<\frac{\epsilon}{2}<1=\left\Vert \hat{\phi}_{N}\right\Vert _{p^{M},q}
\end{equation}

Now comes the hammer: by the ultrametric inequality, (\ref{eq:WTT - Ultrametric inequality})
is an \emph{equality }whenever one of $\left\Vert \hat{\phi}_{N}\right\Vert _{p^{M},q}$
or $\left\Vert \hat{\phi}_{N}*\hat{\mu}*\hat{\eta}_{\epsilon,M}\right\Vert _{p^{M},q}$
is strictly greater than the other. Having proved that to be the case,
(\ref{eq:WTT - Target of attack}) becomes: 
\begin{equation}
\epsilon>\left\Vert \hat{\phi}_{N}*\left(\mathbf{1}_{0}-\left(\hat{\mu}*\hat{\eta}_{\epsilon,M}\right)\right)\right\Vert _{p^{M},q}=\left\Vert \hat{\phi}_{N}-\left(\hat{\phi}_{N}*\hat{\mu}*\hat{\eta}_{\epsilon,M}\right)\right\Vert _{p^{M},q}=1>\epsilon
\end{equation}
for all $M\geq M_{\epsilon}$ and all $N\geq\max\left\{ N_{\epsilon}^{\prime},M\right\} $.
The left-hand side is our first estimate, and the right-hand side
is our second. Since $\epsilon<1$, this is clearly impossible.

This proves that the existence of the zero $\mathfrak{z}_{0}$ precludes
the translates of $\hat{\mu}$ from being dense in $c_{0}\left(\hat{\mathbb{Z}}_{p},\mathbb{C}_{q}\right)$.

Q.E.D.

\subsubsection{A Matter of \label{subsec:A-Matter-of}Matrices}

Along with the \textbf{Correspondence Principle }and the derivation
of a Formula for the Fourier transforms of the $\chi_{H}$s, \textbf{Theorem
\ref{thm:pq WTT for measures}} (and its multi-dimensional analogue
in Subsection \ref{subsec:5.4.4More-Fourier-Resummation}) is one
of this dissertation's central results. Combined, these three theorems
can (and will) be used to show that an integer $x$ is a periodic
point of a contracting, semi-basic $p$-Hydra map $H$ if and only
if the translates of $\hat{\chi}_{H}\left(t\right)-x\mathbf{1}_{0}\left(t\right)$
are dense in $c_{0}\left(\mathbb{Z}_{p},\mathbb{C}_{q_{H}}\right)$.
Aside from the mere aesthetic gratification of this equivalence, it
also appears to be fertile ground for future investigations into the
periodic points of Hydra maps, thanks to the algebraic structure of
$\hat{\mathbb{Z}}_{p}$ and a spoonful of linear algebra.

Identifying the group $\mathbb{Z}/p^{N}\mathbb{Z}$ with the set:
\begin{equation}
\left\{ 0,\frac{1}{p^{N}},\frac{2}{p^{N}}\ldots,\frac{p^{N}-1}{p^{N}}\right\} 
\end{equation}
equipped with addition modulo $1$, it is a standard fact of higher
algebra that $\hat{\mathbb{Z}}_{p}$ is the so-called \index{direct limit}\textbf{direct
limit }of the $\mathbb{Z}/p^{N}\mathbb{Z}$s. In terms of functions,
this is just the fact that any function $\hat{\mathbb{Z}}_{p}\rightarrow\mathbb{C}_{q}$
can be restricted to a function $\mathbb{Z}/p^{N}\mathbb{Z}\rightarrow\mathbb{C}_{q}$
by multiplication by $\mathbf{1}_{0}\left(p^{N}t\right)$. This makes
it feasible to study functions $\hat{\mathbb{Z}}_{p}\rightarrow\mathbb{C}_{q}$
by considering functions $\mathbb{Z}/p^{N}\mathbb{Z}\rightarrow\mathbb{C}_{q}$
as $N\rightarrow\infty$. This is particularly nice, seeing as the
question of dense translations over $\mathbb{Z}/p^{N}\mathbb{Z}$
reduces to a matter of matrices. 
\begin{defn}
\label{def:Conv inv}Let $\hat{\chi}:\hat{\mathbb{Z}}_{p}\rightarrow\mathbb{C}_{q}$.
Then, a \index{convolution!inverse}\textbf{convolution inverse }of
$\hat{\chi}$ is a function $\hat{\chi}^{-1}:\hat{\mathbb{Z}}_{p}\rightarrow\mathbb{C}_{q}$
so that: 
\begin{equation}
\underbrace{\lim_{N\rightarrow\infty}\sum_{\left|s\right|_{p}\leq p^{N}}\hat{\chi}\left(t-s\right)\hat{\chi}^{-1}\left(s\right)}_{\overset{\textrm{def}}{=}\left(\hat{\chi}*\hat{\chi}^{-1}\right)\left(t\right)}\overset{\mathbb{C}_{q}}{=}\mathbf{1}_{0}\left(t\right),\textrm{ }\forall t\in\hat{\mathbb{Z}}_{p}\label{eq:Definition of Convolution Inverse}
\end{equation}
Note that $\mathbf{1}_{0}$ is the identity element with respect to
convolution of functions $\hat{\mathbb{Z}}_{p}\rightarrow\mathbb{C}_{q}$.

Given $N\geq0$, we say that $\hat{\chi}^{-1}$ is an \textbf{$N$th
partial convolution inverse }of $\hat{\chi}$ if: 
\begin{equation}
\sum_{\left|s\right|_{p}\leq p^{N}}\hat{\chi}\left(t-s\right)\hat{\chi}^{-1}\left(s\right)\overset{\mathbb{C}_{q}}{=}\mathbf{1}_{0}\left(t\right),\textrm{ }\forall\left|t\right|_{p}\leq p^{N}\label{eq:Definition of Nth Partial Convolution Inverse}
\end{equation}
\end{defn}
\begin{prop}
Let $\hat{\chi}:\hat{\mathbb{Z}}_{p}\rightarrow\mathbb{C}_{q}$, and
suppose there is a function $\hat{\chi}^{-1}:\hat{\mathbb{Z}}_{p}\rightarrow\mathbb{C}_{q}$
which is an $N$th partial convolution inverse of $\hat{\chi}$ for
all sufficiently large $N$. Then, $\hat{\chi}^{-1}$ is a convolution
inverse of $\hat{\chi}$. 
\end{prop}
Proof: Take limits of (\ref{eq:Definition of Nth Partial Convolution Inverse})
as $N\rightarrow\infty$. (Note: the convergence is only guaranteed
to be point-wise with respect to $t$.)

Q.E.D.

\vphantom{}

Now, given any $N$, observe that the equation (\ref{eq:Definition of Nth Partial Convolution Inverse})
defining an $N$th partial convolution inverse of $\hat{\chi}$ is
actually a system of linear equations with the values of $\hat{\chi}^{-1}\left(s\right)$
on $\left|s\right|_{p}\leq p^{N}$ as its unknowns: 
\begin{align*}
1 & =\sum_{\left|s\right|_{p}\leq p^{N}}\hat{\chi}\left(0-s\right)\hat{\chi}^{-1}\left(s\right)\\
0 & =\sum_{\left|s\right|_{p}\leq p^{N}}\hat{\chi}\left(\frac{1}{p^{N}}-s\right)\hat{\chi}^{-1}\left(s\right)\\
0 & =\sum_{\left|s\right|_{p}\leq p^{N}}\hat{\chi}\left(\frac{2}{p^{N}}-s\right)\hat{\chi}^{-1}\left(s\right)\\
 & \vdots\\
0 & =\sum_{\left|s\right|_{p}\leq p^{N}}\hat{\chi}\left(\frac{p^{N}-1}{p^{N}}-s\right)\hat{\chi}^{-1}\left(s\right)
\end{align*}
We can express this linear system as a matrix. In fact, we can do
so in two noteworthy ways. 
\begin{defn}
Let $\hat{\chi}:\hat{\mathbb{Z}}_{p}\rightarrow\mathbb{C}_{q}$ and
let $N\geq1$.

\vphantom{}

I. We write $\mathbf{M}_{N}\left(\hat{\chi}\right)$ to denote the
$p^{N}\times p^{N}$matrix: 
\begin{equation}
\mathbf{M}_{N}\left(\hat{\chi}\right)\overset{\textrm{def}}{=}\left(\begin{array}{ccccc}
\hat{\chi}\left(0\right) & \hat{\chi}\left(\frac{p^{N}-1}{p^{N}}\right) & \cdots & \hat{\chi}\left(\frac{2}{p^{N}}\right) & \hat{\chi}\left(\frac{1}{p^{N}}\right)\\
\hat{\chi}\left(\frac{1}{p^{N}}\right) & \hat{\chi}\left(0\right) & \hat{\chi}\left(\frac{p^{N}-1}{p^{N}}\right) &  & \hat{\chi}\left(\frac{2}{p^{N}}\right)\\
\vdots & \hat{\chi}\left(\frac{1}{p^{N}}\right) & \hat{\chi}\left(0\right) & \ddots & \vdots\\
\hat{\chi}\left(\frac{p^{N}-2}{p^{N}}\right) &  & \ddots & \ddots & \hat{\chi}\left(\frac{p^{N}-1}{p^{N}}\right)\\
\hat{\chi}\left(\frac{p^{N}-1}{p^{N}}\right) & \hat{\chi}\left(\frac{p^{N}-2}{p^{N}}\right) & \cdots & \hat{\chi}\left(\frac{1}{p^{N}}\right) & \hat{\chi}\left(0\right)
\end{array}\right)\label{eq:Definition of bold M_N of Chi hat}
\end{equation}

\vphantom{}

II. By the \textbf{radial enumeration of $\mathbb{Z}/p^{N}\mathbb{Z}$},
we mean the sequence: 
\begin{equation}
0,\frac{1}{p},\ldots,\frac{p-1}{p},\frac{1}{p^{2}},\ldots,\frac{p-1}{p^{2}},\frac{p+1}{p^{2}},\ldots,\frac{p^{2}-1}{p^{2}},\frac{1}{p^{3}},\ldots,\frac{p^{N}-1}{p^{N}}\label{eq:Definition of the radial enumeration of Z mod p^N Z}
\end{equation}
that is to say, the radial enumeration of $\mathbb{Z}/p^{N}\mathbb{Z}$
starts with $0$, which is then followed by the list of all $t\in\hat{\mathbb{Z}}_{p}$
with $\left|t\right|_{p}=p$ in order of increasing value in the numerator,
which is then followed by the list of all $t\in\hat{\mathbb{Z}}_{p}$
with $\left|t\right|_{p}=p^{2}$ in order of increasing value in the
numerator, so on and so forth, until we have listed all $t$ with
$\left|t\right|_{p}\leq p^{N}$. By the \textbf{standard enumeration
of }$\mathbb{Z}/p^{N}\mathbb{Z}$, on the other hand, we mean the
enumeration: 
\begin{equation}
0,\frac{1}{p^{N}},\ldots,\frac{p^{N}-1}{p^{N}}\label{eq:Definition of the standard enumeration of Z mod p^N Z}
\end{equation}

\vphantom{}

III. We write $\mathbf{R}_{N}\left(\hat{\chi}\right)$ to denote the
$p^{N}\times p^{N}$ matrix: 
\[
\mathbf{R}_{N}\left(\hat{\chi}\right)\overset{\textrm{def}}{=}\left(\begin{array}{cccccc}
\hat{\chi}\left(0-0\right) & \hat{\chi}\left(0-\frac{1}{p}\right) & \cdots & \hat{\chi}\left(0-\frac{p-1}{p}\right) & \cdots & \hat{\chi}\left(0-\frac{p^{N}-1}{p^{N}}\right)\\
\hat{\chi}\left(\frac{1}{p}-0\right) & \hat{\chi}\left(0\right) & \cdots & \hat{\chi}\left(\frac{1}{p}-\frac{p-1}{p}\right) & \cdots & \hat{\chi}\left(\frac{1}{p}-\frac{p^{N}-1}{p^{N}}\right)\\
\vdots & \vdots & \ddots & \vdots &  & \vdots\\
\hat{\chi}\left(\frac{p-1}{p}-0\right) & \hat{\chi}\left(\frac{p-1}{p}-\frac{1}{p}\right) & \cdots & \hat{\chi}\left(0\right) & \cdots & \hat{\chi}\left(\frac{p-1}{p}-\frac{p^{N}-1}{p^{N}}\right)\\
\vdots & \vdots &  & \vdots & \ddots & \vdots\\
\hat{\chi}\left(\frac{p^{N}-1}{p^{N}}-0\right) & \hat{\chi}\left(\frac{p^{N}-1}{p^{N}}-\frac{1}{p}\right) & \cdots & \hat{\chi}\left(\frac{p^{N}-1}{p^{N}}-\frac{p-1}{p}\right) & \cdots & \hat{\chi}\left(0\right)
\end{array}\right)
\]
That is, the entries of the $n$th row of $\mathbf{R}_{N}\left(\hat{\chi}\right)$,
where $n\in\left\{ 0,\ldots,p^{N}-1\right\} $ are: 
\[
\hat{\chi}\left(t_{n}-s_{0}\right),\hat{\chi}\left(t_{n}-s_{1}\right),\hat{\chi}\left(t_{n}-s_{2}\right)\ldots
\]
where $t_{n}$ and $s_{k}$ are the $n$th and $k$th elements of
the radial enumeration of $\mathbb{Z}/p^{N}\mathbb{Z}$, respectively
($s_{0}=t_{0}=0$, $s_{1}=t_{1}=1/p$, $s_{2}=t_{2}=2/p$, etc.). 
\end{defn}
\begin{rem}
If $\hat{\chi}$ is a Fourier transform of some function $\chi:\mathbb{Z}_{p}\rightarrow\mathbb{C}_{q}$,
we will also write $\mathbf{M}_{N}\left(\chi\right)$ and $\mathbf{R}_{N}\left(\chi\right)$
to denote $\mathbf{M}_{N}\left(\hat{\chi}\right)$ and $\mathbf{R}_{N}\left(\hat{\chi}\right)$,
respectively. If $\chi$ and/or $\hat{\chi}$ are not in question,
we will just write $\mathbf{M}_{N}$ and $\mathbf{R}_{N}$. 
\end{rem}
\begin{rem}
As defined, we have that (\ref{eq:Definition of Nth Partial Convolution Inverse})
can be written in matrix form as: 
\begin{equation}
\mathbf{M}_{N}\left(\hat{\chi}\right)\left(\begin{array}{c}
\hat{\chi}^{-1}\left(0\right)\\
\hat{\chi}^{-1}\left(\frac{1}{p^{N}}\right)\\
\vdots\\
\hat{\chi}^{-1}\left(\frac{p^{N}-1}{p^{N}}\right)
\end{array}\right)=\underbrace{\left(\begin{array}{c}
1\\
0\\
\vdots\\
0
\end{array}\right)}_{\mathbf{e}_{1}}\label{eq:Matrix Equation - Circulant Form}
\end{equation}
and as: 
\begin{equation}
\mathbf{R}_{N}\left(\hat{\chi}\right)\left(\begin{array}{c}
\hat{\chi}^{-1}\left(0\right)\\
\hat{\chi}^{-1}\left(\frac{1}{p}\right)\\
\vdots\\
\hat{\chi}^{-1}\left(\frac{p-1}{p}\right)\\
\vdots\\
\hat{\chi}^{-1}\left(\frac{p^{N}-1}{p^{N}}\right)
\end{array}\right)=\mathbf{e}_{1}\label{eq:Matrix Equation - Radial Form}
\end{equation}
where the column of $\hat{\chi}^{-1}$ is listed according to the
radial enumeration of $\hat{\mathbb{Z}}_{p}$. 
\end{rem}
\begin{prop}
\label{prop:permutation similarity}For each $N$ for which at least
one of $\mathbf{R}_{N}\left(\hat{\chi}\right)$ and $\mathbf{M}_{N}\left(\hat{\chi}\right)$
is invertible, there exists a permutation matrix $\mathbf{E}_{N}$
(depending only on $N$ and $p$) so that: 
\begin{equation}
\mathbf{R}_{N}\left(\hat{\chi}\right)=\mathbf{E}_{N}\mathbf{M}_{N}\left(\hat{\chi}\right)\mathbf{E}_{N}^{-1}\label{eq:Permutation Similarity of M and R}
\end{equation}
That is to say, $\mathbf{M}_{N}\left(\hat{\chi}\right)$ and $\mathbf{R}_{N}\left(\hat{\chi}\right)$
are \textbf{permutation similar}. 
\end{prop}
Proof: Let $\mathbf{E}_{N}$ be the permutation matrix which sends
the standard enumeration of $\mathbb{Z}/p^{N}\mathbb{Z}$ to the radial
enumeration of $\mathbb{Z}/p^{N}\mathbb{Z}$: 
\begin{equation}
\mathbf{E}_{N}\left(\begin{array}{c}
0\\
\frac{1}{p^{N}}\\
\vdots\\
\frac{p^{N}-1}{p^{N}}
\end{array}\right)=\left(\begin{array}{c}
0\\
\frac{1}{p}\\
\vdots\\
\frac{p-1}{p}\\
\vdots\\
\frac{p^{N}-1}{p^{N}}
\end{array}\right)\label{eq:Definition of bold E_N}
\end{equation}
Consequently, $\mathbf{E}_{N}^{-1}$ sends the radial enumeration
to the standard enumeration.

Now, without loss of generality, suppose $\mathbf{M}_{N}\left(\hat{\chi}\right)$
is invertible. The proof for when $\mathbf{E}_{N}\left(\hat{\chi}\right)$
is nearly identical, just replace every instance of $\mathbf{E}_{N}^{-1}$
below with $\mathbf{E}_{N}$, and vice-versa, and replace every standard
enumeration with the radial enumeration, and vice-versa. Then, taking
(\ref{eq:Matrix Equation - Radial Form}), letting $\mathbf{x}$ denote
the column vector with $\hat{\chi}^{-1}$ in the standard enumeration,
and letting $\mathbf{y}$ denote the column vector with $\hat{\chi}^{-1}$
in the radial enumeration, observe that: 
\begin{equation}
\mathbf{R}_{N}\left(\hat{\chi}\right)\mathbf{y}=\mathbf{R}_{N}\left(\hat{\chi}\right)\mathbf{E}_{N}\mathbf{x}=\mathbf{e}_{1}
\end{equation}
Left-multiplying both sides by $\mathbf{E}_{N}^{-1}$ yields: 
\begin{equation}
\mathbf{E}_{N}^{-1}\mathbf{R}_{N}\left(\hat{\chi}\right)\mathbf{E}_{N}\mathbf{x}=\mathbf{E}_{N}^{-1}\mathbf{e}_{1}=\mathbf{e}_{1}
\end{equation}
where the right-most equality follows from the fact that the $1$
at the top of $\mathbf{e}_{1}$ is fixed by the permutation encoded
by $\mathbf{E}_{N}^{-1}$. Since $\mathbf{M}_{N}\left(\hat{\chi}\right)$
was assumed to be invertible, (\ref{eq:Matrix Equation - Circulant Form})
shows that $\mathbf{M}_{N}\left(\hat{\chi}\right)$ is the \emph{unique
}matrix so that: 
\[
\mathbf{M}_{N}\left(\hat{\chi}\right)\mathbf{x}=\mathbf{e}_{1}
\]
The above shows that $\left(\mathbf{E}_{N}^{-1}\mathbf{R}_{N}\left(\hat{\chi}\right)\mathbf{E}_{N}\right)\mathbf{x}=\mathbf{e}_{1}$,
which then forces $\mathbf{E}_{N}^{-1}\mathbf{R}_{N}\left(\hat{\chi}\right)\mathbf{E}_{N}=\mathbf{M}_{N}\left(\hat{\chi}\right)$,
and hence, (\ref{eq:Permutation Similarity of M and R}).

Q.E.D.

\vphantom{}

Observe that $\mathbf{M}_{N}\left(\hat{\chi}\right)$ is a matrix
of the form: 
\begin{equation}
\left(\begin{array}{ccccc}
c_{0} & c_{p^{N}-1} & \cdots & c_{2} & c_{1}\\
c_{1} & c_{0} & c_{p^{N}-1} &  & c_{2}\\
\vdots & c_{1} & c_{0} & \ddots & \vdots\\
c_{p^{N}-2} &  & \ddots & \ddots & c_{p^{N}-1}\\
c_{p^{N}-1} & c_{p^{N}-2} & \cdots & c_{1} & c_{0}
\end{array}\right)\label{eq:Bold M_N as a p^N by p^N Circulant Matrix}
\end{equation}
Our analysis is greatly facilitated by the elegant properties of properties
satisfied by matrices of the form (\ref{eq:Bold M_N as a p^N by p^N Circulant Matrix}),
which are called \textbf{circulant matrices}\index{circulant matrix}.
A standard reference for this subject is \cite{Circulant Matrices}. 
\begin{defn}[\textbf{Circulant matrices}]
Let $\mathbb{F}$ be an algebraically closed field of characteristic
$0$, and let $N$ be an integer $\geq1$. Then, an\textbf{ $N\times N$
circulant matrix with entries in $\mathbb{F}$} is a matrix $\mathbf{M}$
of the form: 
\begin{equation}
\left(\begin{array}{ccccc}
c_{0} & c_{N-1} & \cdots & c_{2} & c_{1}\\
c_{1} & c_{0} & c_{N-1} &  & c_{2}\\
\vdots & c_{1} & c_{0} & \ddots & \vdots\\
c_{N-2} &  & \ddots & \ddots & c_{N-1}\\
c_{N-1} & c_{N-2} & \cdots & c_{1} & c_{0}
\end{array}\right)\label{eq:Definition of a Circulant Matrix-1}
\end{equation}
where $c_{0},\ldots,c_{N-1}\in\mathbb{F}$. The polynomial $f_{\mathbf{M}}:\mathbb{F}\rightarrow\mathbb{F}$
defined by: 
\begin{equation}
f_{\mathbf{M}}\left(x\right)\overset{\textrm{def}}{=}\sum_{n=0}^{N-1}c_{n}x^{n}\label{eq:Definition of the associated polynomial of a circulant matrix}
\end{equation}
is called the \textbf{associated polynomial}\index{circulant matrix!associated polynomial}\textbf{
}of the matrix $\mathbf{M}$. 
\end{defn}
\begin{rem}
$\mathbf{M}_{N}\left(\hat{\chi}\right)$ is a circulant matrix for
all $\hat{\chi}:\hat{\mathbb{Z}}_{p}\rightarrow\mathbb{C}_{q}$. 
\end{rem}
\begin{defn}
Let $\mathbb{F}$ be a field of characteristic $0$. For a function
$\hat{\chi}:\hat{\mathbb{Z}}_{p}\rightarrow\mathbb{F}$, we then write
$f_{\hat{\chi},N}:\mathbb{F}\rightarrow\mathbb{F}$ to denote the
associated polynomial of $\mathbf{M}_{N}\left(\hat{\chi}\right)$:
\begin{equation}
f_{\hat{\chi},N}\left(z\right)=\sum_{n=0}^{p^{N}-1}\hat{\chi}\left(\frac{n}{p^{N}}\right)z^{n}\label{eq:Definition/Notation for the associated polynomial of bold M_N of Chi hat}
\end{equation}
\end{defn}
\begin{rem}
Like with $\mathbf{M}_{N}$ and $\mathbf{R}_{N}$, if $\hat{\chi}$
is a Fourier transform of some function $\chi$, we will write $f_{\chi,N}$
to denote $f_{\hat{\chi},N}$. 
\end{rem}
\vphantom{}

The chief properties of circulant matrices stem from their intimate
connection with convolution inverses and the finite Fourier transform. 
\begin{thm}
\label{thm:3.52}Let $\mathbf{M}$ be an $N\times N$ circulant matrix\index{circulant matrix!determinant}
with coefficients in an algebraically closed field $\mathbb{F}$,
and let $\omega\in\mathbb{F}$ be any primitive $N$th root of unity.
Then:

\vphantom{}

I. 
\begin{equation}
\det\mathbf{M}=\prod_{n=0}^{N-1}f_{\mathbf{M}}\left(\omega^{n}\right)\label{eq:Determinant of a Circulant Matrix}
\end{equation}

\vphantom{}

II. The\index{circulant matrix!eigenvalues} eigenvalues of $\mathbf{M}$
are $f_{\mathbf{M}}\left(\omega^{n}\right)$ for $n\in\left\{ 0,\ldots,N-1\right\} $. 
\end{thm}
\vphantom{}

We employ this result to compute the characteristic polynomial of
$\mathbf{M}_{N}\left(\hat{\chi}\right)$ for any $\hat{\chi}$. 
\begin{prop}
\label{prop:3.71}Let $\hat{\chi}:\hat{\mathbb{Z}}_{p}\rightarrow\mathbb{C}_{q}$
be any function. Then: 
\begin{equation}
\det\left(\mathbf{M}_{N}\left(\hat{\chi}\right)-\lambda\mathbf{I}_{p^{N}}\right)\overset{\mathbb{C}_{q}}{=}\prod_{n=0}^{p^{N}-1}\left(\tilde{\chi}_{N}\left(n\right)-\lambda\right),\textrm{ }\forall\lambda\in\mathbb{C}_{q}\label{eq:Characteristic polynomial of bold M_N of Chi hat}
\end{equation}
where, $\mathbf{I}_{p^{N}}$ is a $p^{N}\times p^{N}$ identity matrix,
and where, as usual, $\tilde{\chi}_{N}\left(\mathfrak{z}\right)=\sum_{\left|t\right|_{p}\leq p^{N}}\hat{\chi}\left(t\right)e^{2\pi i\left\{ t\mathfrak{z}\right\} _{p}}$. 
\end{prop}
Proof: First, observe that: 
\begin{equation}
\mathbf{M}_{N}\left(\hat{\chi}\right)-\lambda\mathbf{I}_{p^{N}}=\mathbf{M}_{N}\left(\hat{\chi}-\lambda\mathbf{1}_{0}\right)
\end{equation}
where: 
\begin{equation}
\hat{\chi}\left(t\right)-\lambda\mathbf{1}_{0}\left(t\right)=\begin{cases}
\hat{\chi}\left(0\right)-\lambda & \textrm{if }t=0\\
\hat{\chi}\left(t\right) & \textrm{else}
\end{cases},\textrm{ }\forall t\in\hat{\mathbb{Z}}_{p}
\end{equation}
Here, $\mathbf{M}_{N}\left(\hat{\chi}-\lambda\mathbf{1}_{0}\right)$
is a circulant matrix with the associated polynomial: 
\begin{align*}
f_{\hat{\chi}-\lambda\mathbf{1}_{0},N}\left(\mathfrak{z}\right) & =\sum_{n=0}^{p^{N}-1}\left(\hat{\chi}\left(\frac{n}{p^{N}}\right)-\lambda\mathbf{1}_{0}\left(\frac{n}{p^{N}}\right)\right)\mathfrak{z}^{n}\\
\left(\mathbf{1}_{0}\left(\frac{n}{p^{N}}\right)=\begin{cases}
1 & \textrm{if }n=0\\
0 & \textrm{else}
\end{cases}\right); & =-\lambda+\sum_{n=0}^{p^{N}-1}\hat{\chi}\left(\frac{n}{p^{N}}\right)\mathfrak{z}^{n}\\
 & =f_{\hat{\chi},N}\left(\mathfrak{z}\right)-\lambda
\end{align*}
By (\ref{eq:Determinant of a Circulant Matrix}), since $\mathbf{M}_{N}\left(\hat{\chi}-\lambda\mathbf{1}_{0}\right)=\mathbf{M}_{N}\left(\hat{\chi}\right)-\lambda\mathbf{I}_{p^{N}}$
is an $p^{N}\times p^{N}$ circulant matrix, we have that: 
\begin{align*}
\det\left(\mathbf{M}_{N}\left(\hat{\chi}\right)-\lambda\mathbf{I}_{p^{N}}\right) & \overset{\mathbb{C}_{q}}{=}\prod_{n=0}^{p^{N}-1}\left(f_{\hat{\chi},N}\left(e^{2\pi in/p^{N}}\right)-\lambda\right)\\
 & \overset{\mathbb{C}_{q}}{=}\prod_{n=0}^{p^{N}-1}\left(\sum_{k=0}^{p^{N}-1}\hat{\chi}\left(\frac{k}{p^{N}}\right)e^{\frac{2\pi ikn}{p^{N}}}-\lambda\right)\\
 & \overset{\mathbb{C}_{q}}{=}\prod_{n=0}^{p^{N}-1}\left(\sum_{\left|t\right|_{p}\leq p^{N}}\hat{\chi}\left(t\right)e^{2\pi itn}-\lambda\right)\\
\left(e^{2\pi itn}=e^{2\pi i\left\{ tn\right\} _{p}}\right); & =\prod_{n=0}^{p^{N}-1}\left(\tilde{\chi}_{N}\left(n\right)-\lambda\right)
\end{align*}

Q.E.D. 
\begin{thm}[\textbf{$\chi_{H}$ as an eigenvalue problem}]
\label{thm:Eigenvalue problem}Let $\hat{\chi}:\hat{\mathbb{Z}}_{p}\rightarrow\mathbb{C}_{q}$
be any function. Then:

\vphantom{}

I. For any $\lambda\in\mathbb{C}_{q}$, $\mathbf{M}_{N}-\lambda\mathbf{I}_{p^{N}}$
is invertible if and only if $\tilde{\chi}_{N}\left(n\right)\neq\lambda$
for any $n\in\left\{ 0,\ldots,p^{N}-1\right\} $.

\vphantom{}

II. The eigenvalues\index{matrix!eigenvalues} of $\mathbf{M}_{N}$
are $\left\{ \tilde{\chi}_{N}\left(n\right)\right\} _{0\leq n\leq p^{N}-1}$.

\vphantom{}

III. Let $N\geq1$ satisfy $\tilde{\chi}_{N}\left(n\right)\neq0$
for any $n\in\left\{ 0,\ldots,p^{N}-1\right\} $. Then $\mathbf{M}_{N}$
and $\mathbf{R}_{N}$ are permutation similar, and (\ref{eq:Permutation Similarity of M and R})
holds true. 
\end{thm}
Proof: Use \textbf{Proposition \ref{prop:permutation similarity}},
\textbf{Theorem \ref{thm:3.52}}, and \textbf{Proposition \ref{prop:3.71}}.

Q.E.D.

\chapter{\label{chap:4 A-Study-of}A Study of $\chi_{H}$ - The One-Dimensional
Case}

\includegraphics[scale=0.45]{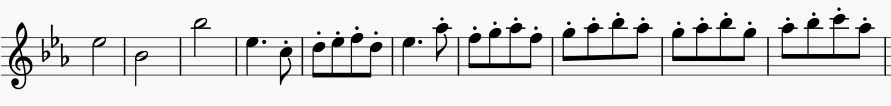}

\vphantom{}

THROUGHOUT THIS CHAPTER, UNLESS STATED OTHERWISE, WE ASSUME $p$ IS
A PRIME NUMBER, AND THAT $H$ IS A CONTRACTING, SEMI-BASIC $p$-HYDRA
MAP WHICH FIXES $0$.

\vphantom{}

Assuming the reader hasn't already figured it out, the musical epigrams
accompanying each of this dissertation's chapters are excerpts from
the fourth movement of Beethoven's Symphony No. 3\textemdash the \emph{Eroica}.
The excerpt for Chapter 4 is from one of the movement's fugal sections.
For the uninitiated, a fugue is a form of western classical composition
based around strict, rigorous, academic applications of imitative
counterpoint among multiple interacting musical voices, and it is
the perfect analogy for the work we are about to undertake.

It is this chapter that we will finally take the techniques recounted
or developed in Chapter 3 and deploy them to analyze $\chi_{H}$.
The main goal is to obtain a meaningful formula for a Fourier transform
($\hat{\chi}_{H}$) of $\chi_{H}$, and thereby show that $\chi_{H}$
is quasi-integrable. With a formula for $\hat{\chi}_{H}$, we can
utilize the Correspondence Principle in conjunction with our $\left(p,q\right)$-adic
generalization of Wiener's Tauberian Theorem (WTT) to reformulate
questions about periodic points and divergent points of $H$ into
Fourier-analytic questions about $\hat{\chi}_{H}$. Roughly speaking,
the equivalences are as follows:

\begin{eqnarray*}
 & x\in\mathbb{Z}\backslash\left\{ 0\right\} \textrm{ is a periodic point or divergent point of }H\\
 & \Updownarrow\\
 & \chi_{H}\left(\mathfrak{z}\right)-x\textrm{ vanishes at some }\mathfrak{z}\in\mathbb{Z}_{p}^{\prime}\\
 & \Updownarrow\\
 & \textrm{translates of }\hat{\chi}_{H}\left(t\right)-x\mathbf{1}_{0}\left(t\right)\textrm{ are not dense in }c_{0}\left(\hat{\mathbb{Z}}_{p},\mathbb{C}_{q}\right)
\end{eqnarray*}
where the upper $\Updownarrow$ is the Correspondence Principle and
the lower $\Updownarrow$ is the WTT.

I conjecture that \emph{both} of the $\Updownarrow$s are, in fact,
rigorous if-and-only-if equivalences. At present, however, each $\Updownarrow$
is only three-quarters true, and requires that $H$ satisfy the hypotheses
of \textbf{Theorem \ref{thm:Divergent trajectories come from irrational z}}
(page \pageref{thm:Divergent trajectories come from irrational z}).
For such an $H$, the actual implications of the Correspondence Principal
are:
\begin{eqnarray*}
 & x\in\mathbb{Z}\backslash\left\{ 0\right\} \textrm{ is a periodic point of }H\\
 & \Updownarrow\\
 & \chi_{H}\left(\mathfrak{z}\right)-x\textrm{ vanishes at some }\mathfrak{z}\in\mathbb{Q}\cap\mathbb{Z}_{p}^{\prime}
\end{eqnarray*}
and:

\begin{eqnarray*}
 & x\in\mathbb{Z}\backslash\left\{ 0\right\} \textrm{ is a divergent point of }H\\
 & \Uparrow\\
 & \chi_{H}\left(\mathfrak{z}\right)-x\textrm{ vanishes at some }\mathfrak{z}\in\mathbb{Z}_{p}\backslash\mathbb{Q}_{p}
\end{eqnarray*}
With the WTT and the help of a a condition on $H$ called \textbf{non-singularity
}we will prove the \textbf{Tauberian Spectral Theorem }for $\chi_{H}$
(\textbf{Theorem \ref{thm:Periodic Points using WTT}}, page \pageref{thm:Periodic Points using WTT})
then yields:
\begin{eqnarray*}
 & \chi_{H}\left(\mathfrak{z}\right)-x\textrm{ vanishes at some }\mathfrak{z}\in\mathbb{Z}_{p}^{\prime}\\
 & \Updownarrow\\
 & \textrm{translates of }\hat{\chi}_{H}\left(t\right)-x\mathbf{1}_{0}\left(t\right)\textrm{ are not dense in }c_{0}\left(\hat{\mathbb{Z}}_{p},\mathbb{C}_{q}\right)
\end{eqnarray*}

To keep the computations in this chapter as simple and manageable
as possible, Section \ref{sec:4.1 Preparatory-Work--} contains some
additional notational conventions (definitions of functions, constants,
etc.) which will streamline the presentation. Section \ref{sec:4.1 Preparatory-Work--}
also includes some identities and asymptotic estimates which will
be used throughout the rest of the chapter.

\section{\label{sec:4.1 Preparatory-Work--}Preparatory Work \textendash{}
Conventions, Identities, etc.}

We begin by introducing some new notations. 
\begin{defn}
\label{def:alpha, beta, gamma 1D}We define \nomenclature{$\alpha_{H}\left(t\right)$}{ }\nomenclature{$\beta_{H}\left(t\right)$}{ }$\alpha_{H},\beta_{H}:\hat{\mathbb{Z}}_{p}\rightarrow\overline{\mathbb{Q}}$
by: 
\begin{equation}
\alpha_{H}\left(t\right)\overset{\textrm{def}}{=}\frac{1}{p}\sum_{j=0}^{p-1}H_{j}^{\prime}\left(0\right)e^{-2\pi ijt}=\frac{1}{p}\sum_{j=0}^{p-1}\frac{\mu_{j}}{p}e^{-2\pi ijt}=\frac{1}{p}\sum_{j=0}^{p-1}\frac{a_{j}}{d_{j}}e^{-2\pi ijt}\label{eq:Definition of alpha_H}
\end{equation}
\begin{equation}
\beta_{H}\left(t\right)\overset{\textrm{def}}{=}\frac{1}{p}\sum_{j=0}^{p-1}H_{j}\left(0\right)e^{-2\pi ijt}=\frac{1}{p}\sum_{j=0}^{p-1}\frac{b_{j}}{d_{j}}e^{-2\pi ijt}\label{eq:Definition of beta_H}
\end{equation}
where, recall, $H_{j}$ denotes the $j$th branch of $H$. We also
adopt the notation: \nomenclature{$\gamma_{H}\left(t\right)$}{ }
\begin{equation}
\gamma_{H}\left(t\right)\overset{\textrm{def}}{=}\frac{\beta_{H}\left(t\right)}{\alpha_{H}\left(t\right)}\label{eq:Definition of gamma_H}
\end{equation}
Additionally, we write $\sigma_{H}$ to denote the constant: 
\begin{equation}
\sigma_{H}\overset{\textrm{def}}{=}\log_{p}\left(\sum_{j=0}^{p-1}H_{j}^{\prime}\left(0\right)\right)=\log_{p}\left(\frac{1}{p}\sum_{j=0}^{p-1}\mu_{j}\right)=1+\log_{p}\left(\alpha_{H}\left(0\right)\right)\label{eq:Definition of sigma_H}
\end{equation}
where $\log_{p}x=\frac{\ln x}{\ln p}$. Equivalently: 
\begin{equation}
\alpha_{H}\left(0\right)=p^{\sigma_{H}-1}\label{eq:alpha_H of 0 in terms of sigma_H}
\end{equation}
Finally, we define the function \nomenclature{$\kappa_{H}\left(n\right)$}{ }$\kappa_{H}:\mathbb{N}_{0}\rightarrow\mathbb{Q}$
by:\index{alpha{H}left(tright)@$\alpha_{H}\left(t\right)$}\index{beta{H}left(tright)@$\beta_{H}\left(t\right)$}\index{gamma{H}left(tright)@$\gamma_{H}\left(t\right)$}
\begin{equation}
\kappa_{H}\left(n\right)\overset{\textrm{def}}{=}M_{H}\left(n\right)\left(\frac{\mu_{0}}{p}\right)^{-\lambda_{p}\left(n\right)}\label{eq:Definition of Kappa_H}
\end{equation}
\end{defn}
\begin{defn}
We say $H$ is \textbf{non-singular }whenever\index{Hydra map!non-singular}
$\alpha_{H}\left(j/p\right)\neq0$ for any $j\in\mathbb{Z}/p\mathbb{Z}$. 
\end{defn}
\vphantom{}

The most important identities for this Chapter involve generating
functions for the $\#_{p:j}\left(n\right)$; recall that these functions
tell us the number of $j$s in the $p$-adic digits of $n$. 
\begin{prop}[\textbf{A Generating Function Identity}]
\label{prop:Generating function identities}Let $\mathbb{F}$ be
a field of characteristic zero, let $p$ be an integer $\geq2$, and
let $c_{0},\ldots,c_{p-1}$ be any non-zero elements of $\mathbb{F}$.
Then, for all $z\in\mathbb{F}$, we have the identities:

\vphantom{}

I. For all $n\geq1$: 
\begin{equation}
\prod_{m=0}^{n-1}\left(\sum_{j=0}^{p-1}c_{j}z^{jp^{m}}\right)=c_{0}^{n}\sum_{k=0}^{p^{n}-1}c_{0}^{-\lambda_{p}\left(k\right)}\left(\prod_{j=0}^{p-1}c_{j}^{\#_{p:j}\left(k\right)}\right)z^{k}\label{eq:M_H partial sum generating identity}
\end{equation}

\vphantom{}

II. If $c_{0}=1$, for all we have the identity: 
\begin{equation}
\sum_{k=p^{n-1}}^{p^{n}-1}\left(\prod_{j=0}^{p-1}c_{j}^{\#_{p:j}\left(k\right)}\right)z^{k}=\left(\sum_{k=1}^{p-1}c_{k}z^{kp^{n-1}}\right)\prod_{m=0}^{n-2}\left(\sum_{j=0}^{p-1}c_{j}z^{jp^{m}}\right)\label{eq:Lambda-restricted partial M_H sum}
\end{equation}
The $n$-product is defined to be $1$ whenever $n=1$. 
\end{prop}
Proof:

I. Each term on the right of (\ref{eq:M_H partial sum generating identity})
is obtained by taking a product of the form: 
\begin{equation}
\prod_{\ell=0}^{n-1}c_{j_{\ell}}z^{j_{\ell}p^{\ell}}=\left(\prod_{h=0}^{n-1}c_{j_{h}}\right)z^{\sum_{\ell=0}^{n-1}j_{\ell}p^{\ell}}
\end{equation}
for some choice of $j_{1},\ldots,j_{n-1}\in\mathbb{Z}/p\mathbb{Z}$.
Note, however, that if there is an $m\in\left\{ 0,\ldots,n-1\right\} $
so that $j_{\ell}=0$ for all $\ell\in\left\{ m,\ldots,n-1\right\} $,
the $c_{j_{\ell}}$s will still keep contributing to the product even
for those values of $\ell$. So, in order to use our $\#_{p:j}\left(k\right)$
notation\textemdash which stops counting the $p$-adic digits of $k$
after the last \emph{non-zero} digit in $k$'s $p$-adic expansion\textemdash we
need to modify the above product to take into account the number of
extra $0$s that occur past the last non-zero digit in $k$'s $p$-adic
expansion.

So, let $k\geq1$, and let: 
\begin{equation}
k=\sum_{\ell=0}^{n-1}j_{\ell}p^{\ell}
\end{equation}
Then, $j_{\ell}=0$ for all $\ell\geq\lambda_{p}\left(k\right)$.
As such: 
\begin{equation}
\prod_{h=0}^{n-1}c_{j_{h}}=c_{0}^{\left|\left\{ h\in\left\{ 0,\ldots,n-1\right\} :j_{h}=0\right\} \right|}\times\prod_{\ell=1}^{p-1}c_{\ell}^{\#_{p:\ell}\left(k\right)}
\end{equation}
Now, $\#_{p:0}\left(k\right)$ is the number of $\ell\leq\lambda_{p}\left(k\right)-1$
for which $j_{\ell}=0$. Since $j_{\ell}=0$ for all $\ell\geq\lambda_{p}\left(k\right)$,
and because there are $n-\lambda_{p}\left(k\right)$ values of $h\in\left\{ 0,\ldots,n-1\right\} $
which are in the interval $\left[\lambda_{p}\left(k\right),n-1\right]$,
we can write: 
\begin{equation}
\left|\left\{ h\in\left\{ 0,\ldots,n-1\right\} :j_{h}=0\right\} \right|=\#_{p:0}\left(k\right)+n-\lambda_{p}\left(k\right)
\end{equation}
So: 
\begin{equation}
\prod_{h=0}^{n-1}c_{j_{h}}=c_{0}^{\#_{p:0}\left(k\right)+n-\lambda_{p}\left(k\right)}\times\prod_{\ell=1}^{p-1}c_{\ell}^{\#_{p:\ell}\left(k\right)}=c_{0}^{n-\lambda_{p}\left(k\right)}\prod_{\ell=0}^{p-1}c_{\ell}^{\#_{p:\ell}\left(k\right)}
\end{equation}
Hence, every term on the right hand side of (\ref{eq:M_H partial sum generating identity})
is of the form: 
\begin{equation}
\prod_{\ell=0}^{n-1}c_{j_{\ell}}z^{j_{\ell}p^{\ell}}=\left(\prod_{h=0}^{n-1}c_{j_{h}}\right)z^{\sum_{\ell=0}^{n-1}j_{\ell}p^{\ell}}=c_{0}^{n-\lambda_{p}\left(k\right)}\left(\prod_{\ell=0}^{p-1}c_{\ell}^{\#_{p:\ell}\left(k\right)}\right)z^{k}
\end{equation}
where: 
\begin{equation}
k=\sum_{\ell=0}^{n-1}j_{\ell}p^{\ell}\leq\sum_{\ell=0}^{n-1}\left(p-1\right)p^{\ell}=p^{n}-1
\end{equation}
Summing over $k\in\left\{ 0,\ldots,p^{n}-1\right\} $ then proves
the identity in (I).

\vphantom{}

II. Suppose $c_{0}=1$. Then (I) can be written as: 
\[
\prod_{m=0}^{n-1}\left(\sum_{j=0}^{p-1}c_{j}z^{jp^{m}}\right)=\sum_{k=0}^{p^{n}-1}\left(\prod_{\ell=0}^{p-1}c_{\ell}^{\#_{p:\ell}\left(k\right)}\right)z^{k}
\]
Thus, subtracting the $\left(n-2\right)$nd case from the $\left(n-1\right)$st
case gives: 
\begin{align*}
\sum_{k=p^{n-1}}^{p^{n}-1}\left(\prod_{\ell=0}^{p-1}c_{\ell}^{\#_{p:\ell}\left(k\right)}\right)z^{k} & =\prod_{m=0}^{n-1}\left(\sum_{j=0}^{p-1}c_{j}z^{jp^{m}}\right)-\prod_{m=0}^{n-2}\left(\sum_{j=0}^{p-1}c_{j}z^{jp^{m}}\right)\\
 & =\left(\sum_{k=0}^{p-1}c_{k}z^{kp^{n-1}}-1\right)\prod_{m=0}^{n-2}\left(\sum_{j=0}^{p-1}c_{j}z^{jp^{m}}\right)
\end{align*}
Note that the $k=0$ term is $c_{0}z^{0\cdot p^{n-1}}=1\cdot1=1$;
this cancels out the $-1$, thereby proving (II).

Q.E.D.

\vphantom{}

Next, we will need some additional formulae, properties, and estimates
for $\kappa_{H}$. 
\begin{prop}[\textbf{Properties of} $\kappa_{H}$]
\label{prop:Properties of Kappa_H}\ 

\vphantom{}

I. \index{kappa{H}@$\kappa_{H}$!properties} 
\begin{equation}
\kappa_{H}\left(n\right)=\prod_{k=1}^{p-1}\left(\frac{\mu_{k}}{\mu_{0}}\right)^{\#_{p:k}\left(n\right)},\textrm{ }\forall n\in\mathbb{N}_{0}\label{eq:Kappa_H explicit formula}
\end{equation}

\vphantom{}

II. $\kappa_{H}$ is the unique function $\mathbb{N}_{0}\rightarrow\mathbb{Q}$
satisfying the functional equations \index{kappa{H}@$\kappa_{H}$!functional equation}:
\begin{equation}
\kappa_{H}\left(pn+j\right)=\frac{\mu_{j}}{\mu_{0}}\kappa_{H}\left(n\right),\textrm{ }\forall n\geq0,\textrm{ }\forall j\in\mathbb{Z}/p\mathbb{Z}\label{eq:Kappa_H functional equations}
\end{equation}
subject to the initial condition $\kappa_{H}\left(0\right)=1$.

\vphantom{}

III. $\kappa_{H}$ has $p$-adic structure, with: 
\begin{equation}
\kappa_{H}\left(m+jp^{n}\right)=\frac{\mu_{j}}{\mu_{0}}\kappa_{H}\left(m\right),\textrm{ }\forall n\in\mathbb{N}_{0},\textrm{ }\forall j\in\mathbb{Z}/p\mathbb{Z},\textrm{ }\forall m\in\mathbb{Z}/p^{n}\mathbb{Z}\label{eq:Kappa_H is rho-adically distributed}
\end{equation}

\vphantom{}

IV. If $H$ is semi-basic, then $\kappa_{H}$ is $\left(p,q_{H}\right)$-adically
regular, with: 
\begin{equation}
\lim_{n\rightarrow\infty}\left|\kappa_{H}\left(\left[\mathfrak{z}\right]_{p^{n}}\right)\right|_{q_{H}}\overset{\mathbb{R}}{=}0,\textrm{ }\forall\mathfrak{z}\in\mathbb{Z}_{p}^{\prime}\label{eq:Semi-basic q-adic decay for Kappa_H}
\end{equation}

\vphantom{}

V. If $H$ is contracting\index{Hydra map!contracting} ($\mu_{0}/p<1$),
then: 
\begin{equation}
\lim_{N\rightarrow\infty}\left(\frac{\mu_{0}}{p}\right)^{N}\kappa_{H}\left(\left[\mathfrak{z}\right]_{p^{N}}\right)\overset{\mathbb{R}}{=}0,\textrm{ }\forall\mathfrak{z}\in\mathbb{N}_{0}\label{eq:Kappa_H decay when H is contracting}
\end{equation}
If $H$ is expanding ($\mu_{0}/p>1$), then: 
\begin{equation}
\lim_{N\rightarrow\infty}\left(\frac{\mu_{0}}{p}\right)^{N}\kappa_{H}\left(\left[\mathfrak{z}\right]_{p^{N}}\right)\overset{\mathbb{R}}{=}+\infty,\textrm{ }\forall\mathfrak{z}\in\mathbb{N}_{0}\label{eq:Kappa_H growth when H is expanding}
\end{equation}
\end{prop}
Proof:

I. Using \textbf{Proposition \ref{prop:Explicit Formulas for M_H}
}(page \pageref{prop:Explicit Formulas for M_H}), we have: 
\begin{equation}
\prod_{j=0}^{p-1}\mu_{j}^{\#_{p:j}\left(\left[\mathfrak{z}\right]_{p^{N}}\right)}=p^{\lambda_{p}\left(\left[\mathfrak{z}\right]_{p^{N}}\right)}M_{H}\left(\left[\mathfrak{z}\right]_{p^{N}}\right)=\frac{M_{H}\left(\left[\mathfrak{z}\right]_{p^{N}}\right)}{p^{-\lambda_{p}\left(\left[\mathfrak{z}\right]_{p^{N}}\right)}}
\end{equation}
Multiplying through by: 
\begin{equation}
\left(\frac{\mu_{0}}{p}\right)^{N}\mu_{0}^{-\lambda_{p}\left(\left[\mathfrak{z}\right]_{p^{N}}\right)}
\end{equation}
makes the right-hand side into: 
\begin{equation}
\left(\frac{\mu_{0}}{p}\right)^{N}\kappa_{H}\left(\left[\mathfrak{z}\right]_{p^{N}}\right)
\end{equation}
So: 
\begin{align*}
\left(\frac{\mu_{0}}{p}\right)^{N}\kappa_{H}\left(\left[\mathfrak{z}\right]_{p^{N}}\right) & =\left(\frac{\mu_{0}}{p}\right)^{N}\mu_{0}^{-\lambda_{p}\left(\left[\mathfrak{z}\right]_{p^{N}}\right)}\prod_{j=0}^{p-1}\mu_{j}^{\#_{p:j}\left(\left[\mathfrak{z}\right]_{p^{N}}\right)}\\
 & =\left(\frac{\mu_{0}}{p}\right)^{N}\mu_{0}^{\#_{p:0}\left(\left[\mathfrak{z}\right]_{p^{N}}\right)-\lambda_{p}\left(\left[\mathfrak{z}\right]_{p^{N}}\right)}\prod_{j=1}^{p-1}\mu_{j}^{\#_{p:j}\left(\left[\mathfrak{z}\right]_{p^{N}}\right)}\\
\left(\lambda_{p}\left(k\right)=\sum_{j=0}^{p-1}\#_{p:j}\left(k\right)\right); & =\left(\frac{\mu_{0}}{p}\right)^{N}\mu_{0}^{-\sum_{j=1}^{p-1}\#_{p:j}\left(\left[\mathfrak{z}\right]_{p^{N}}\right)}\prod_{j=1}^{p-1}\mu_{j}^{\#_{p:j}\left(\left[\mathfrak{z}\right]_{p^{N}}\right)}\\
 & =\left(\frac{\mu_{0}}{p}\right)^{N}\prod_{j=1}^{p-1}\left(\frac{\mu_{j}}{\mu_{0}}\right)^{\#_{p:j}\left(\left[\mathfrak{z}\right]_{p^{N}}\right)}
\end{align*}
Now, divide by $\left(\mu_{0}/p\right)^{N}$, pick $\mathfrak{z}=n\in\mathbb{N}_{0}$,
and let $N\geq\lambda_{p}\left(n\right)$. This yields (I).

\vphantom{}

II \& III. Let $n\geq1$ and let $j\in\mathbb{Z}/p\mathbb{Z}$. Then:
\begin{align*}
\lambda_{p}\left(pn+j\right) & =\lambda_{p}\left(n\right)+1\\
M_{H}\left(pn+j\right) & =\frac{\mu_{j}}{p}M_{H}\left(n\right)
\end{align*}
Consequently, we have: 
\begin{align*}
\kappa_{H}\left(pn+j\right) & =\left(\frac{p}{\mu_{0}}\right)^{\lambda_{p}\left(pn+j\right)}M_{H}\left(pn+j\right)\\
 & =\frac{p}{\mu_{0}}\cdot\frac{\mu_{j}}{p}\cdot\left(\frac{p}{\mu_{0}}\right)^{\lambda_{p}\left(n\right)}M_{H}\left(n\right)\\
 & =\frac{\mu_{j}}{\mu_{0}}\kappa_{H}\left(n\right)
\end{align*}
Next: 
\begin{align*}
\kappa_{H}\left(j\right) & =\begin{cases}
\left(\frac{p}{\mu_{0}}\right)^{\lambda_{p}\left(0\right)}M_{H}\left(0\right) & \textrm{if }j=0\\
\left(\frac{p}{\mu_{0}}\right)^{\lambda_{p}\left(pn+j\right)}M_{H}\left(j\right) & \textrm{if }j\in\left\{ 1,\ldots,p-1\right\} 
\end{cases}\\
 & =\begin{cases}
1 & \textrm{if }j=0\\
\frac{\mu_{j}}{\mu_{0}} & \textrm{if }j\in\left\{ 1,\ldots,p-1\right\} 
\end{cases}
\end{align*}
So, $\kappa_{H}\left(0\right)=1$. Consequently 
\begin{equation}
\kappa_{H}\left(j\right)=\frac{\mu_{j}}{\mu_{0}}\kappa_{H}\left(0\right),\textrm{ }\forall j\in\mathbb{Z}/p\mathbb{Z}
\end{equation}
Combining this with the other cases computed above proves (\ref{eq:Kappa_H functional equations}).

As for (III), replace $n$ in (\ref{eq:Kappa_H explicit formula})
with $m+jp^{n}$, where $m\in\left\{ 0,\ldots,p^{n}-1\right\} $.
This gives us: 
\begin{align*}
\kappa_{H}\left(m+jp^{n}\right) & =\prod_{k=1}^{p-1}\left(\frac{\mu_{k}}{\mu_{0}}\right)^{\#_{p:k}\left(m+jp^{n}\right)}\\
\left(\#_{p:k}\left(m+jp^{n}\right)=\#_{p:k}\left(m\right)+\left[j=k\right]\right); & =\frac{\mu_{j}}{\mu_{0}}\prod_{k=1}^{p-1}\left(\frac{\mu_{k}}{\mu_{0}}\right)^{\#_{p:k}\left(m\right)}\\
 & =\frac{\mu_{j}}{\mu_{0}}\kappa_{H}\left(m\right)
\end{align*}
So:

\[
\kappa_{H}\left(pn+j\right)=\frac{\mu_{j}}{\mu_{0}}\kappa_{H}\left(n\right),\textrm{ }\forall n\geq0
\]
\[
\kappa_{H}\left(jp^{n}+m\right)=\frac{\mu_{m}}{\mu_{0}}\kappa_{H}\left(n\right),\textrm{ }\forall n\geq0
\]
Finally, by \textbf{Proposition \ref{prop:formula for functions with p-adic structure}},
we then have that $\kappa_{H}\left(n\right)$ is uniquely determined
by the $\mu_{j}$s and $\kappa_{H}\left(0\right)$. This proves that
any function $\mathbb{N}_{0}\rightarrow\mathbb{Q}$ satisfying (\ref{eq:Kappa_H functional equations})
and $\kappa_{H}\left(0\right)=1$ is necessarily equal to $\kappa_{H}$.
This shows the desired uniqueness, thereby proving the rest of (II).

\vphantom{}

IV. Let $\mathfrak{z}\in\mathbb{Z}_{p}^{\prime}$. If $H$ is semi-basic,
we have that $\mu_{0}$ and $p$ are both co-prime to $q_{H}$. Taking
$q$-adic absolute values of (I) yields: 
\begin{equation}
\left|\kappa_{H}\left(\left[\mathfrak{z}\right]_{p^{N}}\right)\right|_{q}=\left|\prod_{j=1}^{p-1}\left(\frac{\mu_{j}}{\mu_{0}}\right)^{\#_{p:j}\left(\left[\mathfrak{z}\right]_{p^{N}}\right)}\right|_{q}=\prod_{j=1}^{p-1}\left|\mu_{j}\right|_{q}^{\#_{p:j}\left(\left[\mathfrak{z}\right]_{p^{N}}\right)}
\end{equation}
With $H$ being semi-basic, $\mu_{j}$ is a multiple of $q_{H}$ for
all non-zero $j$. As such, by the pigeonhole principle, since $\mathfrak{z}\in\mathbb{Z}_{p}^{\prime}$,
there is at least one $j\in\left\{ 1,\ldots,p-1\right\} $ so that
infinitely many of $\mathfrak{z}$'s $p$-adic digits are $j$. This
makes $\#_{p:j}\left(\left[\mathfrak{z}\right]_{p^{N}}\right)\rightarrow\infty$
in $\mathbb{R}$ as $N\rightarrow\infty$, and demonstrates that (\ref{eq:Semi-basic q-adic decay for Kappa_H})
holds. This proves (IV).

\vphantom{}

V. Let $\mathfrak{z}\in\mathbb{N}_{0}$. Then, $\left[\mathfrak{z}\right]_{p^{N}}=\mathfrak{z}$
for all $N\geq\lambda_{p}\left(\mathfrak{z}\right)$, and so: 
\begin{equation}
\left(\frac{\mu_{0}}{p}\right)^{N}\kappa_{H}\left(\left[\mathfrak{z}\right]_{p^{N}}\right)=\left(\frac{\mu_{0}}{p}\right)^{N}\kappa_{H}\left(\mathfrak{z}\right),\textrm{ }\forall N\geq\lambda_{p}\left(\mathfrak{z}\right)
\end{equation}
where $\kappa_{H}\left(\mathfrak{z}\right)\neq0$ is then a non-zero
rational constant $c$. Thus, as $N\rightarrow\infty$ the limit in
$\mathbb{R}$ of the left-hand side will be $0$ if and only if $H$
is contracting ($\mu_{0}/p<1$) and will be $+\infty$ if and only
if $H$ is expanding ($\mu_{0}/p>1$). This proves (V).

Q.E.D.

\vphantom{}

Next, we compute the value of the summatory function\index{chi{H}@$\chi_{H}$!summatory function}
$\sum_{n=0}^{N}\chi_{H}\left(n\right)$ via a recursive method. 
\begin{lem}[\textbf{Summatory function of} $\chi_{H}$]
\label{lem:Summatory function of Chi_H}Let $H$ be a $p$-Hydra
map (we \emph{do not }require that $\mu_{0}=1$). Then: 
\begin{equation}
\sum_{n=0}^{p^{N}-1}\chi_{H}\left(n\right)=\begin{cases}
\beta_{H}\left(0\right)Np^{N} & \textrm{if }\sigma_{H}=1\\
\frac{\beta_{H}\left(0\right)}{\alpha_{H}\left(0\right)-1}\left(p^{\sigma_{H}N}-p^{N}\right) & \textrm{else}
\end{cases},\textrm{ }\forall N\geq1\label{eq:rho-N minus 1 sum of Chi_H}
\end{equation}
\end{lem}
Proof: For any $N\geq0$, define: 
\begin{equation}
S_{H}\left(N\right)\overset{\textrm{def}}{=}\frac{1}{p^{N}}\sum_{n=0}^{p^{N}-1}\chi_{H}\left(n\right)\label{eq:Definition of S_H of N}
\end{equation}
Letting $N\geq1$, note that: 
\begin{align*}
p^{N}S_{H}\left(N\right) & =\sum_{n=0}^{p^{N}-1}\chi_{H}\left(n\right)\\
 & =\sum_{j=0}^{p-1}\sum_{n=0}^{p^{N-1}-1}\chi_{H}\left(pn+j\right)\\
\left(\textrm{use }\chi_{H}\textrm{'s functional eqs.}\right); & =\sum_{j=0}^{p-1}\sum_{n=0}^{p^{N-1}-1}\frac{a_{j}\chi_{H}\left(n\right)+b_{j}}{d_{j}}\\
 & =\sum_{j=0}^{p-1}\left(\frac{a_{j}}{d_{j}}\sum_{n=0}^{p^{N-1}-1}\chi_{H}\left(n\right)+\frac{b_{j}}{d_{j}}p^{N-1}\right)\\
 & =\left(\sum_{j=0}^{p-1}\frac{a_{j}}{d_{j}}\right)p^{N-1}S_{H}\left(N-1\right)+\left(\sum_{j=0}^{p-1}\frac{b_{j}}{d_{j}}\right)p^{N-1}
\end{align*}
Dividing through by $p^{N}$, note that: 
\begin{align*}
\alpha_{H}\left(0\right) & =\frac{1}{p}\sum_{j=0}^{p-1}\frac{a_{j}}{d_{j}}=\frac{1}{p}\left(\frac{1}{p}\sum_{j=0}^{p-1}\mu_{j}\right)=p^{\sigma_{H}-1}\\
\beta_{H}\left(0\right) & =\frac{1}{p}\sum_{j=0}^{p-1}\frac{b_{j}}{d_{j}}
\end{align*}
As such, we can write: 
\begin{equation}
S_{H}\left(N\right)=p^{\sigma_{H}-1}S_{H}\left(N-1\right)+\beta_{H}\left(0\right)\label{eq:Recursive formula for S_H}
\end{equation}
With this, we see that $S_{H}\left(N\right)$ is the image of $S_{H}\left(0\right)$
under $N$ iterates of the affine linear map: 
\begin{equation}
x\mapsto p^{\sigma_{H}-1}x+\beta_{H}\left(0\right)\label{eq:Affine linear map generating S_H}
\end{equation}
This leaves us with two cases to investigate:

\vphantom{}

i. Suppose $\sigma_{H}=1$. Then, (\ref{eq:Affine linear map generating S_H})
reduces to a translation: 
\begin{equation}
x\mapsto x+\beta_{H}\left(0\right)
\end{equation}
in which case:

\begin{equation}
S_{H}\left(N\right)=S_{H}\left(0\right)+\beta_{H}\left(0\right)N\label{eq:Explicit formula for S_H of N when sigma_H =00003D00003D 1}
\end{equation}

\vphantom{}

ii. Suppose $\sigma_{H}\neq1$. Then, (\ref{eq:Affine linear map generating S_H})
is not a translation, and so, using the explicit formula for this
$N$th iterate, we obtain: 
\begin{equation}
S_{H}\left(N\right)=p^{\left(\sigma_{H}-1\right)N}S_{H}\left(0\right)+\frac{p^{\left(\sigma_{H}-1\right)N}-1}{p^{\sigma_{H}-1}-1}\beta_{H}\left(0\right)\label{eq:Explicit formula for S_H of N}
\end{equation}
Noting that: 
\begin{align}
S_{H}\left(0\right) & =\sum_{n=0}^{p^{0}-1}\chi_{H}\left(n\right)=\chi_{H}\left(0\right)=0
\end{align}
we then have: 
\begin{equation}
\frac{1}{p^{N}}\sum_{n=0}^{p^{N}-1}\chi_{H}\left(n\right)=S_{H}\left(N\right)=\begin{cases}
\beta_{H}\left(0\right)N & \textrm{if }\sigma_{H}=1\\
\frac{p^{\left(\sigma_{H}-1\right)N}-1}{p^{\sigma_{H}-1}-1}\beta_{H}\left(0\right) & \textrm{else}
\end{cases}
\end{equation}
Multiplying through on all sides by $p^{N}$ and re-writing $p^{\sigma_{H}-1}-1$
as $\alpha_{H}\left(0\right)-1$ then yields: 
\begin{equation}
\sum_{n=0}^{p^{N}-1}\chi_{H}\left(n\right)=\begin{cases}
\beta_{H}\left(0\right)Np^{N} & \textrm{if }\sigma_{H}=1\\
\frac{\beta_{H}\left(0\right)}{\alpha_{H}\left(0\right)-1}\left(p^{\sigma_{H}N}-p^{N}\right) & \textrm{else}
\end{cases}
\end{equation}
which proves (\ref{eq:rho-N minus 1 sum of Chi_H}).

Q.E.D.

\vphantom{}

Although not of any particular use (at least at the time of writing),
we can use the formula for the summatory function to obtain simple
upper and lower (archimedean) bounds. We do this with help of the
following trick: 
\begin{prop}[\textbf{Replacing $N$ with }$\log_{p}\left(N+1\right)$]
Let $\left\{ a_{n}\right\} _{n\geq0}$ be a sequence of non-negative
real numbers, and let $p$ be an integer $\geq2$. Suppose there is
a non-decreasing function $A:\left[0,\infty\right)\rightarrow\left[0,\infty\right)$
such that: 
\[
\sum_{n=0}^{p^{N}-1}a_{n}=A\left(N\right),\textrm{ }\forall N\in\mathbb{N}_{1}
\]
Then: 
\begin{equation}
A\left(\log_{p}N\right)\leq\sum_{n=0}^{N}a_{n}\leq2A\left(\log_{p}\left(N\right)+1\right)-A\left(\log_{p}\left(N\right)-1\right),\textrm{ }\forall N\geq1\label{eq:Replacement Lemma}
\end{equation}
\end{prop}
Proof: As just remarked, for each $k\in\mathbb{N}_{1}$, $\left\{ p^{k-1},p^{k-1}+1,\ldots,p^{k}-1\right\} $
is the set of all positive integers $n$ for which $\lambda_{p}\left(n\right)=k$.
As such, letting $N\in\mathbb{N}_{1}$ we have that $p^{\lambda_{p}\left(N\right)-1}$
is the smallest positive integer whose $p$-adic representation has
the same number of digits as the $p$-adic representation of $N$.
Consequently, we can write:

\begin{equation}
\sum_{n=0}^{N}a_{n}=\sum_{n=0}^{p^{\lambda_{p}\left(N\right)-1}-1}a_{n}+\sum_{n=p^{\lambda_{p}\left(N\right)-1}}^{N}a_{n}=A\left(\lambda_{p}\left(N\right)\right)+\sum_{n=p^{\lambda_{p}\left(N\right)-1}}^{N}a_{n}\label{eq:Lambda decomposition of Nth partial sum}
\end{equation}
Since the $a_{n}$s are non-negative, this then yields the lower bound:
\begin{equation}
\sum_{n=0}^{N}a_{n}\geq A\left(\lambda_{p}\left(N\right)\right)
\end{equation}
Since: 
\begin{equation}
\lambda_{p}\left(N\right)=1+\left\lfloor \log_{p}N\right\rfloor \geq\log_{p}N
\end{equation}
and since $A$ is a non-decreasing function, this then implies: 
\begin{equation}
\sum_{n=0}^{N}a_{n}\geq A\left(\log_{p}N\right)
\end{equation}

As for an upper bound, we note that the non-negativity of the $a_{n}$s
allows us to write: 
\begin{align*}
\sum_{n=p^{\lambda_{p}\left(N\right)-1}}^{N}a_{n} & \leq\sum_{n=p^{\lambda_{p}\left(N\right)-1}}^{p^{\lambda_{p}\left(N\right)}-1}a_{n}\\
 & =\sum_{n=0}^{p^{\lambda_{p}\left(N\right)}-1}a_{n}-\sum_{n=0}^{p^{\lambda_{p}\left(N\right)-1}-1}a_{n}\\
 & =A\left(\lambda_{p}\left(N\right)\right)-A\left(\lambda_{p}\left(N\right)-1\right)
\end{align*}
seeing as $p^{\lambda_{p}\left(N\right)}-1$ is the largest positive
integer with the same number of $p$-adic digits as $N$. With this,
(\ref{eq:Lambda decomposition of Nth partial sum}) becomes: 
\begin{align}
\sum_{n=0}^{N}a_{n} & \leq2A\left(\lambda_{p}\left(N\right)\right)-A\left(\lambda_{p}\left(N\right)-1\right)
\end{align}
Since $A$ is non-decreasing, and since:

\begin{equation}
\log_{p}N\leq1+\left\lfloor \log_{p}N\right\rfloor =\lambda_{p}\left(N\right)=1+\left\lfloor \log_{p}N\right\rfloor \leq1+\log_{p}N
\end{equation}
we then have that: 
\begin{align}
\sum_{n=0}^{N}a_{n} & \leq2A\left(\log_{p}\left(N\right)+1\right)-A\left(\log_{p}\left(N\right)-1\right)
\end{align}

Q.E.D. 
\begin{prop}[\textbf{Archimedean Bounds for $\chi_{H}$'s Summatory Function}]
\label{prop:archimedean bounds on Chi_H summatory function}Let $H$
be a $p$-Hydra map (we \emph{do not }require that $\mu_{0}=1$).
Then\index{chi{H}@$\chi_{H}$!summatory function!estimates}:

\vphantom{}

I. If $\sigma_{H}=1$: 
\begin{equation}
\beta_{H}\left(0\right)N\log_{p}N\leq\sum_{n=0}^{N}\chi_{H}\left(n\right)\leq C_{H,1}N+\frac{2p^{2}-1}{p}\beta_{H}\left(0\right)N\log_{p}N\label{eq:Chi_H partial sum asympotitcs for when sigma_H is 1}
\end{equation}
where: 
\begin{equation}
C_{H,1}\overset{\textrm{def}}{=}\frac{2p^{2}+1}{p}\beta_{H}\left(0\right)
\end{equation}

\vphantom{}

II. If $\sigma_{H}>1$: 
\begin{equation}
C_{H,2}N^{\sigma_{H}}-\gamma_{H}\left(0\right)N\leq\sum_{n=0}^{N}\chi_{H}\left(n\right)\leq\frac{2p^{2\sigma_{H}}-1}{p^{\sigma_{H}}}C_{H,2}N^{\sigma_{H}}-\frac{2p^{2}-1}{p}\gamma_{H}\left(0\right)N\label{eq:Chi_H partial sum asympotitcs for when sigma_H is greater than 1}
\end{equation}
where: 
\begin{equation}
C_{H,2}\overset{\textrm{def}}{=}\gamma_{H}\left(0\right)
\end{equation}
\begin{equation}
\gamma_{H}\left(0\right)\overset{\textrm{def}}{=}\frac{\beta_{H}\left(0\right)}{\alpha_{H}\left(0\right)-1}
\end{equation}
\end{prop}
Proof: When $\sigma_{H}\geq1$, the formulae in (\ref{eq:rho-N minus 1 sum of Chi_H})
are non-decreasing functions of $N$. As such, we may use (\ref{eq:Replacement Lemma})
which, upon simplification, yields (\ref{eq:Chi_H partial sum asympotitcs for when sigma_H is 1})
and (\ref{eq:Chi_H partial sum asympotitcs for when sigma_H is greater than 1}).

Q.E.D.

\vphantom{}

Finally, we have the truncations of $\chi_{H}$: 
\begin{notation}
We write $q$ to denote $q_{H}$. We write $\chi_{H,N}:\mathbb{Z}_{p}\rightarrow\mathbb{Q}\subset\mathbb{Q}_{q}$
to denote the $N$th truncation\index{chi{H}@$\chi_{H}$!$N$th truncation}
of $\chi_{H}$: 
\begin{equation}
\chi_{H,N}\left(\mathfrak{z}\right)\overset{\textrm{def}}{=}\sum_{n=0}^{p^{N}-1}\chi_{H}\left(n\right)\left[\mathfrak{z}\overset{p^{N}}{\equiv}n\right],\textrm{ }\forall\mathfrak{z}\in\mathbb{Z}_{p},\textrm{ }\forall N\in\mathbb{N}_{0}\label{eq:Definition of Nth truncation of Chi_H}
\end{equation}
Recall that: 
\begin{equation}
\chi_{H,N}\left(n\right)=\chi_{H}\left(n\right),\textrm{ }\forall n\in\left\{ 0,\ldots,p^{N}-1\right\} 
\end{equation}
In particular: 
\begin{equation}
\chi_{H}\left(n\right)=\chi_{H,\lambda_{p}\left(n\right)}\left(n\right),\textrm{ }\forall n\in\mathbb{N}_{0}
\end{equation}
Because $\chi_{H,N}$ is a rational-valued function which is a linear
combination of finitely many indicator functions for clopen subsets
of $\mathbb{Z}_{p}$, $\chi_{H,N}$ is continuous both as a function
$\mathbb{Z}_{p}\rightarrow\mathbb{C}_{q}$ and as a function $\mathbb{Z}_{p}\rightarrow\mathbb{C}$.
As a result of this, in writing $\hat{\chi}_{H,N}:\hat{\mathbb{Z}}_{p}\rightarrow\overline{\mathbb{Q}}$
to denote the Fourier Transform of $\chi_{H,N}$: 
\begin{equation}
\hat{\chi}_{H,N}\left(t\right)\overset{\textrm{def}}{=}\int_{\mathbb{Z}_{p}}\chi_{H,N}\left(\mathfrak{z}\right)e^{-2\pi i\left\{ t\mathfrak{z}\right\} _{p}}d\mathfrak{z},\textrm{ }\forall t\in\hat{\mathbb{Z}}_{p}\label{eq:Definition of the Fourier Coefficients of Chi_H,N}
\end{equation}
observe that this integral is convergent in both the topology of $\mathbb{C}$
and the topology of $\mathbb{C}_{q}$. In fact, because $\chi_{H,N}$
is locally constant, the integral will reduce to a finite sum (see
(\ref{eq:Chi_H,N hat transform formula - sum form})): 
\begin{equation}
\hat{\chi}_{H,N}\left(t\right)\overset{\overline{\mathbb{Q}}}{=}\frac{\mathbf{1}_{0}\left(p^{N}t\right)}{p^{N}}\sum_{n=0}^{p^{N}-1}\chi_{H}\left(n\right)e^{-2\pi int},\textrm{ }\forall t\in\hat{\mathbb{Z}}_{p}
\end{equation}
where: \nomenclature{$\mathbf{1}_{0}\left(t\right)$}{  } 
\begin{equation}
\mathbf{1}_{0}\left(p^{N}t\right)=\left[\left|t\right|_{p}\leq p^{N}\right],\textrm{ }\forall t\in\hat{\mathbb{Z}}_{p}
\end{equation}
As such, all of the computations and formal manipulations we will
perform hold simultaneously in $\mathbb{C}$ and in $\mathbb{C}_{q}$;
they both occur in $\overline{\mathbb{Q}}\subset\mathbb{C}\cap\mathbb{C}_{q}$.
The difference between the archimedean and non-archimedean topologies
will only emerge when we consider what happens as $N$ tends to $\infty$.

Additionally, note that because of $\chi_{H,N}$'s local constant-ness
and the finitude of its range, we have that, for each $N$, $\hat{\chi}_{H,N}$
has finite support ($\left|t\right|_{p}\geq p^{N+1}$): 
\begin{equation}
\hat{\chi}_{H,N}\left(t\right)=0,\textrm{ }\forall\left|t\right|_{p}\geq p^{N+1},\textrm{ }\forall N\in\mathbb{N}_{0}\label{eq:Vanishing of Chi_H,N hat for all t with sufficiently large denominators}
\end{equation}
and thus, that the Fourier series: 
\begin{equation}
\chi_{H,N}\left(\mathfrak{z}\right)=\sum_{t\in\hat{\mathbb{Z}}_{p}}\hat{\chi}_{H,N}\left(t\right)e^{2\pi i\left\{ t\mathfrak{z}\right\} _{p}}\label{eq:Fourier series for Chi_H,N}
\end{equation}
will be uniformly convergent in both $\mathbb{C}$ and $\mathbb{C}_{q}$
with respect to $\mathfrak{z}\in\mathbb{Z}_{p}$, reducing to a finite
sum in all cases.

Finally, with regard to algebraic-number-theoretic issues of picking
embeddings of $\overline{\mathbb{Q}}$ in $\mathbb{C}$ and $\mathbb{C}_{q}$,
there is no need to worry. This is because all of our work hinges
on the Fourier series identity:
\begin{equation}
\left[\mathfrak{z}\overset{p^{N}}{\equiv}n\right]\overset{\overline{\mathbb{Q}}}{=}\frac{1}{p^{N}}\sum_{\left|t\right|_{p}\leq p^{N}}e^{2\pi i\left\{ t\left(\mathfrak{z}-n\right)\right\} _{p}}
\end{equation}
which\textemdash \emph{crucially}\textemdash is invariant under the
action of elements of the Galois group $\textrm{Gal}\left(\overline{\mathbb{Q}}/\mathbb{Q}\right)$,
seeing as the sum evaluates to a rational number (the indicator function,
which is either $0$ or $1$). As such, all finite operations (sums,
multiplication, series rearrangements) involving this identity are
\emph{independent }of our choices of embeddings of $\overline{\mathbb{Q}}$
in $\mathbb{C}$ and $\mathbb{C}_{q}$.
\end{notation}
\vphantom{}

Lastly, we need to make is the following lemma regarding the interaction
between truncation and functional equations. This result will play
a critical role in actually \emph{proving }that the functions $\hat{\mathbb{Z}}_{p}\rightarrow\overline{\mathbb{Q}}$
which we will derive are, in fact, Fourier transforms of $\chi_{H}$. 
\begin{lem}
\label{lem:Functional equations and truncation}Let $\chi:\mathbb{N}_{0}\rightarrow\mathbb{Q}$,
and suppose that for $j\in\left\{ 0,\ldots,p-1\right\} $ there are
functions $\Phi_{j}:\mathbb{N}_{0}\times\mathbb{Q}\rightarrow\mathbb{Q}$
such that: 
\begin{equation}
\chi\left(pn+j\right)=\Phi_{j}\left(n,\chi\left(n\right)\right),\textrm{ }\forall n\in\mathbb{N}_{0},\textrm{ }\forall j\in\mathbb{Z}/p\mathbb{Z}\label{eq:Relation between truncations and functional equations - Hypothesis}
\end{equation}
Then, the $N$th truncations $\chi_{N}$ satisfy the functional equations\index{functional equation}:
\begin{equation}
\chi_{N}\left(pn+j\right)=\Phi_{j}\left(\left[n\right]_{p^{N-1}},\chi_{N-1}\left(n\right)\right),\textrm{ }\forall n\in\mathbb{N}_{0},\textrm{ }\forall j\in\mathbb{Z}/p\mathbb{Z}\label{eq:Relation between truncations and functional equations, version 1}
\end{equation}
Equivalently: 
\begin{equation}
\chi\left(\left[pn+j\right]_{p^{N}}\right)=\Phi_{j}\left(\left[n\right]_{p^{N-1}},\chi\left(\left[n\right]_{p^{N-1}}\right)\right),\textrm{ }\forall n\in\mathbb{N}_{0},\textrm{ }\forall j\in\mathbb{Z}/p\mathbb{Z}\label{eq:Relation between truncations and functional equations, version 2}
\end{equation}
\end{lem}
Proof: Fix $N\geq0$, $n\in\mathbb{N}_{0}$, and $j\in\mathbb{Z}/p\mathbb{Z}$.
Then: 
\begin{align*}
\chi\left(\left[pn+j\right]_{p^{N}}\right) & =\chi_{N}\left(pn+j\right)\\
 & =\sum_{m=0}^{p^{N}-1}\chi\left(m\right)\left[pn+j\overset{p^{N}}{\equiv}m\right]\\
\left(\textrm{split }m\textrm{ mod }p\right); & =\sum_{\ell=0}^{p^{N-1}-1}\sum_{k=0}^{p-1}\chi\left(p\ell+k\right)\left[pn+j\overset{p^{N}}{\equiv}p\ell+k\right]\\
 & =\sum_{\ell=0}^{p^{N-1}-1}\sum_{k=0}^{p-1}\Phi_{k}\left(\ell,\chi\left(\ell\right)\right)\underbrace{\left[n\overset{p^{N-1}}{\equiv}\ell+\frac{k-j}{p}\right]}_{0\textrm{ }\forall n\textrm{ if }k\neq j}\\
 & =\sum_{\ell=0}^{p^{N-1}-1}\Phi_{j}\left(\ell,\chi\left(\ell\right)\right)\left[n\overset{p^{N-1}}{\equiv}\ell\right]\\
 & =\Phi_{j}\left(\left[n\right]_{p^{N-1}},\chi\left(\left[n\right]_{p^{N-1}}\right)\right)
\end{align*}

Q.E.D. 
\begin{rem}
As an aside, it is worth mentioning that $\sigma_{H}$ is implicated
in the abscissa of convergence of the Dirichlet series\index{Dirichlet series}\index{chi{H}@$\chi_{H}$!Dirichlet series}:
\begin{equation}
\zeta_{M_{H}}\left(s\right)\overset{\textrm{def}}{=}\sum_{n=0}^{\infty}\frac{M_{H}\left(n\right)}{\left(n+1\right)^{s}}
\end{equation}
and: 
\begin{equation}
\zeta_{\chi_{H}}\left(s\right)\overset{\textrm{def}}{=}\sum_{n=0}^{\infty}\frac{\chi_{H}\left(n\right)}{\left(n+1\right)^{s}}
\end{equation}
For instance, using \textbf{Abel's Summation Formula }in the standard
way yields: 
\begin{equation}
\sum_{n=0}^{p^{N}-1}\frac{\chi_{H}\left(n\right)}{\left(1+n\right)^{s}}=\frac{1}{p^{Ns}}\sum_{n=0}^{p^{N}-1}\chi_{H}\left(n\right)+s\int_{0}^{p^{N}-1}\frac{\sum_{n=0}^{\left\lfloor x\right\rfloor }\chi_{H}\left(n\right)}{\left(1+x\right)^{s+1}}dx
\end{equation}
So, by \textbf{Proposition \ref{lem:Summatory function of Chi_H}},
when $\sigma_{H}=1$, we obtain: 
\[
\sum_{n=0}^{p^{N}-1}\frac{\chi_{H}\left(n\right)}{\left(1+n\right)^{s}}=\frac{N\beta_{H}\left(0\right)}{p^{N\left(s-1\right)}}+s\int_{0}^{p^{N}-1}\frac{\sum_{n=0}^{\left\lfloor x\right\rfloor }\chi_{H}\left(n\right)}{\left(1+x\right)^{s+1}}dx
\]
whereas for $\sigma_{H}\neq1$, we obtain:

\[
\sum_{n=0}^{p^{N}-1}\frac{\chi_{H}\left(n\right)}{\left(1+n\right)^{s}}=\frac{\beta_{H}\left(0\right)}{\alpha_{H}\left(0\right)-1}\left(\frac{1}{p^{N\left(s-\sigma_{H}\right)}}-\frac{1}{p^{N\left(s-1\right)}}\right)+s\int_{0}^{p^{N}-1}\frac{\sum_{n=0}^{\left\lfloor x\right\rfloor }\chi_{H}\left(n\right)}{\left(1+x\right)^{s+1}}dx
\]
which then gives: 
\begin{equation}
\zeta_{\chi_{H}}\left(s\right)=s\int_{0}^{\infty}\frac{\sum_{n=0}^{\left\lfloor x\right\rfloor }\chi_{H}\left(n\right)}{\left(1+x\right)^{s+1}}dx,\textrm{ }\forall\textrm{Re}\left(s\right)>\max\left\{ 1,\sigma_{H}\right\} 
\end{equation}
The (archimedean) bounds on $\chi_{H}$'s summatory function given
in \textbf{Proposition \ref{prop:archimedean bounds on Chi_H summatory function}
}then show that $\zeta_{\chi_{H}}\left(s\right)$ converges absolutely
for all $\textrm{Re}\left(s\right)>\sigma_{H}$ for all $H$ with
$\sigma_{H}\geq1$. As mentioned at the end of Chapter 2, this Dirichlet
series can be analytically continued to meromorphic function on $\mathbb{C}$
with a half-lattice of poles in the half-plane $\textrm{Re}\left(s\right)\leq\sigma_{H}$,
however, its growth rate as $\textrm{Re}\left(s\right)\rightarrow-\infty$
is hyper-exponential, making it seemingly ill-suited for analysis
via techniques of contour integration.
\end{rem}
\newpage{}

\section{\label{sec:4.2 Fourier-Transforms-=00003D000026}Fourier Transforms
and Quasi-Integrability}

THROUGHOUT THIS SECTION, WE ALSO ASSUME $H$ IS NON-SINGULAR.

\vphantom{}

This section is the core of our analysis of $\chi_{H}$. Understandably,
it is also the most computationally intensive portion of this dissertation\textemdash at
least until Chapter 6. Our work will be broken up into three parts.

In Subsection \ref{subsec:4.2.1}, we use the functional equations
(\ref{eq:Functional Equations for Chi_H over the rho-adics}) to establish
a recursion relation between the Fourier transforms of the $N$th
truncations of $\chi_{H}$. Solving this yields an explicit formula
for $\hat{\chi}_{H,N}\left(t\right)$ (\ref{eq:Explicit formula for Chi_H,N hat})
which is a sum of products of $\beta_{H}$ and a partial product involving
$\alpha_{H}$. The partial products involving $\alpha_{H}$ depend
on both $N$ and $t$, and by carefully manipulating the product,
we can extract from the product the part that depends on $N$. Doing
this yields a function $\hat{A}_{H}:\hat{\mathbb{Z}}_{p}\rightarrow\overline{\mathbb{Q}}$
(\vref{eq:Definition of A_H hat}) whose properties we then study
in detail, establishing various formulae along the way. Subsection
\ref{subsec:4.2.1} concludes with greatly simplified expressions
for $\hat{\chi}_{H,N}\left(t\right)$ given by\textbf{ Theorem} \textbf{\ref{thm:(N,t) asymptotic decomposition of Chi_H,N hat}}.
These can be considered explicit ``asymptotic formulae'' for $\hat{\chi}_{H,N}\left(t\right)$
in terms of $t$ and $N$ as $N\rightarrow\infty$.

The asymptotic formulae of \textbf{Theorem} \textbf{\ref{thm:(N,t) asymptotic decomposition of Chi_H,N hat}}
demonstrate the sensitive dependence of $\hat{\chi}_{H,N}\left(t\right)$
on $p$ and, most of all, on $\alpha_{H}\left(0\right)$. The dependence
on $p$ is a consequence of the fact that the multiplicative group
$\left(\mathbb{Z}/p\mathbb{Z}\right)^{\times}$ of multiplicatively
invertible integers modulo $p$ satisfies $\left|\left(\mathbb{Z}/p\mathbb{Z}\right)^{\times}\right|=1$
if and only if $p=2$. Much more significant, however, is the dependence
on $\alpha_{H}\left(0\right)$. The essence of this dependence was
already on display in the formula for $\chi_{H}$'s summatory function
as given in \textbf{Lemma \ref{lem:Summatory function of Chi_H}},
where we saw two distinct behaviors, one for $\sigma_{H}=1$ (equivalently,
for $\alpha_{H}\left(0\right)=1$) and one for all other values. This
dependency will appear again in our analysis of $\hat{A}_{H}$, where
it yields two cases: $\alpha_{H}\left(0\right)=1$ and $\alpha_{H}\left(0\right)\neq1$.

In Subsection \ref{subsec:4.2.2}, we first show that $\chi_{H}$
is quasi-integrable with respect to the standard $\left(p,q_{H}\right)$-adic
frame in the $\alpha_{H}\left(0\right)=1$ case. We will then exploit
the functional equations (\ref{eq:Functional Equations for Chi_H over the rho-adics})
to extend this result to arbitrary values of $\alpha_{H}\left(0\right)$.
As a consequence, we will be able to provide Fourier transforms for
$\chi_{H}$\textemdash thereby proving $\chi_{H}$'s quasi-integrability\textemdash along
with interesting non-trivial explicit series expressions for $\chi_{H}$
itself.

\subsection{\label{subsec:4.2.1}$\hat{\chi}_{H,N}$ and $\hat{A}_{H}$}

We begin by computing a recursive formula\index{hat{chi}{H,N}@$\hat{\chi}_{H,N}$!recursive formula}
for $\hat{\chi}_{H,N}$ in terms of $\hat{\chi}_{H,N-1}$. Solving
this produces an explicit formula for $\hat{\chi}_{H,N}$ in terms
of $t$ and $N$. 
\begin{prop}[\textbf{Formulae for }$\hat{\chi}_{H,N}$]
\label{prop:Computation of Chi_H,N hat}\ 

\vphantom{}

I. 
\begin{equation}
\hat{\chi}_{H,N}\left(t\right)\overset{\overline{\mathbb{Q}}}{=}\frac{\mathbf{1}_{0}\left(p^{N}t\right)}{p^{N}}\sum_{n=0}^{p^{N}-1}\chi_{H}\left(n\right)e^{-2\pi int}\label{eq:Chi_H,N hat transform formula - sum form}
\end{equation}

\vphantom{}

II. 
\begin{equation}
\hat{\chi}_{H,N}\left(t\right)=\mathbf{1}_{0}\left(p^{N}t\right)\alpha_{H}\left(t\right)\hat{\chi}_{H,N-1}\left(pt\right)+\mathbf{1}_{0}\left(pt\right)\beta_{H}\left(t\right),\textrm{ }\forall N\geq1,\textrm{ }\forall t\in\hat{\mathbb{Z}}_{p}\label{eq:Chi_H,N hat functional equation}
\end{equation}
the nesting of which yields: 
\begin{equation}
\hat{\chi}_{H,N}\left(t\right)=\sum_{n=0}^{N-1}\mathbf{1}_{0}\left(p^{n+1}t\right)\beta_{H}\left(p^{n}t\right)\prod_{m=0}^{n-1}\alpha_{H}\left(p^{m}t\right)\label{eq:Explicit formula for Chi_H,N hat}
\end{equation}
where the $m$-product is defined to be $1$ when $n=0$. In particular,
note that $\hat{\chi}_{H,N}\left(t\right)=0$ for all $t\in\hat{\mathbb{Z}}_{p}$
with $\left|t\right|_{p}>p^{N}$. 
\end{prop}
Proof: We use the Fourier series for $\left[\mathfrak{z}\overset{p^{N}}{\equiv}n\right]$
to re-write $\chi_{H,N}$ as a Fourier series:
\begin{align*}
\chi_{H,N}\left(\mathfrak{z}\right) & =\sum_{n=0}^{p^{N}-1}\chi_{H}\left(n\right)\left[\mathfrak{z}\overset{p^{N}}{\equiv}n\right]\\
 & =\sum_{n=0}^{p^{N}-1}\chi_{H}\left(n\right)\frac{1}{p^{N}}\sum_{\left|t\right|_{p}\leq p^{N}}e^{2\pi i\left\{ t\left(\mathfrak{z}-n\right)\right\} _{p}}\\
 & =\sum_{\left|t\right|_{p}\leq p^{N}}\left(\frac{1}{p^{N}}\sum_{n=0}^{p^{N}-1}\chi_{H}\left(n\right)e^{-2\pi int}\right)e^{2\pi i\left\{ t\mathfrak{z}\right\} _{p}}
\end{align*}
which is the Fourier series representation of $\chi_{H,N}$, with
$\hat{\chi}_{H,N}\left(t\right)$ being the coefficient of $e^{2\pi i\left\{ t\mathfrak{z}\right\} _{p}}$
in the series. As such: 
\[
\hat{\chi}_{H,N}\left(t\right)=\frac{\mathbf{1}_{0}\left(p^{N}t\right)}{p^{N}}\sum_{n=0}^{p^{N}-1}\chi_{H}\left(n\right)e^{-2\pi int}
\]
which proves (I).

Next, we split the $n$-sum modulo $p$ and utilize $\chi_{H}$'s
functional equation (here $\left|t\right|_{p}\leq p^{N}$): 
\begin{align*}
\hat{\chi}_{H,N}\left(t\right) & =\frac{\mathbf{1}_{0}\left(p^{N}t\right)}{p^{N}}\sum_{j=0}^{p-1}\sum_{n=0}^{p^{N-1}-1}\chi_{H}\left(pn+j\right)e^{-2\pi i\left(pn+j\right)t}\\
 & =\frac{\mathbf{1}_{0}\left(p^{N}t\right)}{p^{N}}\sum_{j=0}^{p-1}\sum_{n=0}^{p^{N-1}-1}\frac{a_{j}\chi_{H}\left(n\right)+b_{j}}{d_{j}}e^{-2\pi i\left(pn+j\right)t}\\
 & =\frac{\mathbf{1}_{0}\left(p^{N}t\right)}{p^{N}}\sum_{n=0}^{p^{N-1}-1}\left(\sum_{j=0}^{p-1}\frac{a_{j}\chi_{H}\left(n\right)+b_{j}}{d_{j}}e^{-2\pi ijt}\right)e^{-2\pi in\left(pt\right)}
\end{align*}
The right-hand side can be re-written like so:
\[
\frac{\mathbf{1}_{0}\left(p^{N}t\right)}{p^{N-1}}\sum_{n=0}^{p^{N-1}-1}\left(\underbrace{\frac{1}{p}\sum_{j=0}^{p-1}\frac{a_{j}}{d_{j}}e^{-2\pi ijt}}_{\alpha_{H}\left(t\right)}\chi_{H}\left(n\right)+\underbrace{\frac{1}{p}\sum_{j=0}^{p-1}\frac{b_{j}}{d_{j}}e^{-2\pi ijt}}_{\beta_{H}\left(t\right)}\right)e^{-2\pi in\left(pt\right)}
\]
which leaves us with:
\begin{align*}
\hat{\chi}_{H,N}\left(t\right) & =\frac{\mathbf{1}_{0}\left(p^{N}t\right)}{p^{N-1}}\sum_{n=0}^{p^{N-1}-1}\left(\alpha_{H}\left(t\right)\chi_{H}\left(n\right)+\beta_{H}\left(t\right)\right)e^{-2\pi in\left(pt\right)}\\
 & =\mathbf{1}_{0}\left(p^{N}t\right)\left(\underbrace{\frac{\alpha_{H}\left(t\right)}{p^{N-1}}\sum_{n=0}^{p^{N-1}-1}\chi_{H}\left(n\right)e^{-2\pi in\left(pt\right)}}_{\alpha_{H}\left(t\right)\hat{\chi}_{H,N-1}\left(pt\right)}+\frac{\beta_{H}\left(t\right)}{p^{N-1}}\sum_{n=0}^{p^{N-1}-1}e^{-2\pi in\left(pt\right)}\right)
\end{align*}
and so: 
\begin{equation}
\hat{\chi}_{H,N}\left(t\right)=\mathbf{1}_{0}\left(p^{N}t\right)\left(\alpha_{H}\left(t\right)\hat{\chi}_{H,N-1}\left(pt\right)+\frac{\beta_{H}\left(t\right)}{p^{N-1}}\sum_{n=0}^{p^{N-1}-1}e^{-2\pi in\left(pt\right)}\right)\label{eq:Chi_H N hat, almost ready to nest}
\end{equation}

Simplifying the remaining $n$-sum, we find that: 
\begin{align*}
\frac{1}{p^{N-1}}\sum_{n=0}^{p^{N-1}-1}e^{-2\pi in\left(pt\right)} & =\begin{cases}
\frac{1}{p^{N-1}}\sum_{n=0}^{p^{N-1}-1}1 & \textrm{if }\left|t\right|_{p}\leq p\\
\frac{1}{p^{N-1}}\frac{\left(e^{-2\pi i\left(pt\right)}\right)^{p^{N-1}}-1}{e^{-2\pi i\left(pt\right)}-1} & \textrm{if }\left|t\right|_{p}\geq p^{2}
\end{cases}\\
 & =\begin{cases}
1 & \textrm{if }\left|t\right|_{p}\leq p\\
0 & \textrm{if }\left|t\right|_{p}\geq p^{2}
\end{cases}\\
 & =\mathbf{1}_{0}\left(pt\right)
\end{align*}
Hence: 
\begin{align*}
\hat{\chi}_{H,N}\left(t\right) & =\mathbf{1}_{0}\left(p^{N}t\right)\left(\alpha_{H}\left(t\right)\hat{\chi}_{H,N-1}\left(pt\right)+\beta_{H}\left(t\right)\mathbf{1}_{0}\left(pt\right)\right)\\
 & =\mathbf{1}_{0}\left(p^{N}t\right)\alpha_{H}\left(t\right)\hat{\chi}_{H,N-1}\left(pt\right)+\mathbf{1}_{0}\left(pt\right)\beta_{H}\left(t\right)
\end{align*}
This proves (II).

With (II) proven, we can then nest (\ref{eq:Chi_H,N hat functional equation})
to derive an explicit formula for $\hat{\chi}_{H,N}\left(t\right)$:
\begin{align*}
\hat{\chi}_{H,N}\left(t\right) & =\mathbf{1}_{0}\left(p^{N}t\right)\alpha_{H}\left(t\right)\hat{\chi}_{H,N-1}\left(pt\right)+\mathbf{1}_{0}\left(pt\right)\beta_{H}\left(t\right)\\
 & =\mathbf{1}_{0}\left(p^{N}t\right)\alpha_{H}\left(t\right)\left(\mathbf{1}_{0}\left(p^{N-1}pt\right)\alpha_{H}\left(pt\right)\hat{\chi}_{H,N-2}\left(p^{2}t\right)+\mathbf{1}_{0}\left(p^{2}t\right)\beta_{H}\left(pt\right)\right)+\mathbf{1}_{0}\left(pt\right)\beta_{H}\left(t\right)\\
 & =\mathbf{1}_{0}\left(p^{N}t\right)\alpha_{H}\left(t\right)\alpha_{H}\left(pt\right)\hat{\chi}_{H,N-2}\left(p^{2}t\right)+\mathbf{1}_{0}\left(p^{2}t\right)\alpha_{H}\left(t\right)\beta_{H}\left(pt\right)+\mathbf{1}_{0}\left(pt\right)\beta_{H}\left(t\right)\\
 & \vdots\\
 & =\mathbf{1}_{0}\left(p^{N}t\right)\left(\prod_{n=0}^{N-2}\alpha_{H}\left(p^{n}t\right)\right)\hat{\chi}_{H,1}\left(p^{N-1}t\right)+\sum_{n=0}^{N-2}\beta_{H}\left(p^{n}t\right)\mathbf{1}_{0}\left(p^{n+1}t\right)\prod_{m=0}^{n-1}\alpha_{H}\left(p^{m}t\right)
\end{align*}
where the $m$-product is $1$ whenever $n=0$.

Finally: 
\begin{align*}
\hat{\chi}_{H,1}\left(t\right) & =\frac{\mathbf{1}_{0}\left(pt\right)}{p}\sum_{j=0}^{p-1}\chi_{H}\left(j\right)e^{-2\pi ijt}
\end{align*}
Since $\chi_{H}\left(0\right)=0$, $\chi_{H}$'s functional equation
gives: 
\[
\chi_{H}\left(j\right)=\chi_{H}\left(p\cdot0+j\right)=\frac{a_{j}\chi_{H}\left(0\right)+b_{j}}{d_{j}}=\frac{b_{j}}{d_{j}}
\]
So: 
\begin{align*}
\hat{\chi}_{H,1}\left(t\right) & =\frac{\mathbf{1}_{0}\left(pt\right)}{p}\sum_{j=0}^{p-1}\chi_{H}\left(j\right)e^{-2\pi ijt}=\frac{\mathbf{1}_{0}\left(pt\right)}{p}\sum_{j=0}^{p-1}\frac{b_{j}}{d_{j}}e^{-2\pi ijt}=\mathbf{1}_{0}\left(pt\right)\beta_{H}\left(t\right)
\end{align*}
Consequently: 
\begin{align*}
\hat{\chi}_{H,N}\left(t\right) & =\mathbf{1}_{0}\left(p^{N}t\right)\left(\prod_{n=0}^{N-2}\alpha_{H}\left(p^{n}t\right)\right)\hat{\chi}_{H,1}\left(p^{N-1}t\right)+\sum_{n=0}^{N-2}\mathbf{1}_{0}\left(p^{n+1}t\right)\beta_{H}\left(p^{n}t\right)\prod_{m=0}^{n-1}\alpha_{H}\left(p^{m}t\right)\\
 & =\mathbf{1}_{0}\left(p^{N}t\right)\beta_{H}\left(p^{N-1}t\right)\prod_{n=0}^{N-2}\alpha_{H}\left(p^{n}t\right)+\sum_{n=0}^{N-2}\mathbf{1}_{0}\left(p^{n+1}t\right)\beta_{H}\left(p^{n}t\right)\prod_{m=0}^{n-1}\alpha_{H}\left(p^{m}t\right)\\
 & =\sum_{n=0}^{N-1}\mathbf{1}_{0}\left(p^{n+1}t\right)\beta_{H}\left(p^{n}t\right)\prod_{m=0}^{n-1}\alpha_{H}\left(p^{m}t\right)
\end{align*}
This proves (\ref{eq:Explicit formula for Chi_H,N hat}).

Q.E.D.

\vphantom{}

Because $\chi_{H}$ is rising-continuous, rather than continuous,
for any fixed $t\in\hat{\mathbb{Z}}_{p}$, the limit of $\hat{\chi}_{H,N}\left(t\right)$
as $N\rightarrow\infty$ need not converge in $\mathbb{C}_{q}$. To
overcome this, we perform a kind of asymptotic analysis, excising
as much as possible the dependence of (\ref{eq:Explicit formula for Chi_H,N hat})
on $N$. Instead, we \emph{reverse} the dependence of $t$ and $n$
or $N$ in (\ref{eq:Explicit formula for Chi_H,N hat}). That this
can be done hinges on the observation that, for any fixed $t\in\hat{\mathbb{Z}}_{p}$,
the congruence $p^{m}t\overset{1}{\equiv}0$ will be satisfied for
all sufficiently large $m$\textemdash namely, all $m\geq-v_{p}\left(t\right)$.
So, as $n\rightarrow\infty$, the $m$th term of the product: 
\begin{equation}
\prod_{m=0}^{n-1}\alpha_{H}\left(p^{m}t\right)
\end{equation}
will be $\alpha_{H}\left(0\right)$ for all $m\geq-v_{p}\left(t\right)+1$.
This suggests that we can untangle the dependency of this product
on $n$ by breaking it up into the partial product from $m=0$ to
$m=-v_{p}\left(t\right)-1$, and a remainder product which will be
independent of $t$. This partial product is precisely $\hat{A}_{H}\left(t\right)$. 
\begin{defn}[$\hat{A}_{H}\left(t\right)$]
We\index{hat{A}{H}left(tright)@$\hat{A}_{H}\left(t\right)$} define
the function \nomenclature{$\hat{A}_{H}\left(t\right)$}{ }$\hat{A}_{H}:\hat{\mathbb{Z}}_{p}\rightarrow\overline{\mathbb{Q}}$
by: 
\begin{equation}
\hat{A}_{H}\left(t\right)\overset{\textrm{def}}{=}\prod_{m=0}^{-v_{p}\left(t\right)-1}\alpha_{H}\left(p^{m}t\right),\textrm{ }\forall t\in\hat{\mathbb{Z}}_{p}\label{eq:Definition of A_H hat}
\end{equation}
where the $m$-product is defined to be $1$ whenever $t=0$; that
is, $\hat{A}_{H}\left(0\right)\overset{\textrm{def}}{=}1$. 
\end{defn}
\vphantom{}

Our next proposition actually does the work of untangling $\prod_{m=0}^{n-1}\alpha_{H}\left(p^{m}t\right)$. 
\begin{prop}[\textbf{$\alpha_{H}$-product in terms of $\hat{A}_{H}$}]
\label{prop:alpha product in terms of A_H hat}Fix $t\in\hat{\mathbb{Z}}_{p}\backslash\left\{ 0\right\} $.
Then: 
\begin{align}
\mathbf{1}_{0}\left(p^{n+1}t\right)\prod_{m=0}^{n-1}\alpha_{H}\left(p^{m}t\right) & =\begin{cases}
0 & \textrm{if }n\leq-v_{p}\left(t\right)-2\\
\frac{\hat{A}_{H}\left(t\right)}{\alpha_{H}\left(\frac{t\left|t\right|_{p}}{p}\right)} & \textrm{if }n=-v_{p}\left(t\right)-1\\
\left(\alpha_{H}\left(0\right)\right)^{n+v_{p}\left(t\right)}\hat{A}_{H}\left(t\right) & \textrm{if }n\geq-v_{p}\left(t\right)
\end{cases}\label{eq:alpha product in terms of A_H hat}
\end{align}
\end{prop}
Proof: Fix $t\in\hat{\mathbb{Z}}_{p}\backslash\left\{ 0\right\} $.
Then, we can write: 
\begin{align*}
\mathbf{1}_{0}\left(p^{n+1}t\right)\prod_{m=0}^{n-1}\alpha_{H}\left(p^{m}t\right) & =\begin{cases}
0 & \textrm{if }\left|t\right|_{p}\geq p^{n+2}\\
\prod_{m=0}^{n-1}\alpha_{H}\left(p^{m}t\right) & \textrm{if }\left|t\right|_{p}\leq p^{n+1}
\end{cases}\\
 & =\begin{cases}
0 & \textrm{if }n\leq-v_{p}\left(t\right)-2\\
\prod_{m=0}^{n-1}\alpha_{H}\left(p^{m}t\right) & \textrm{if }n\geq-v_{p}\left(t\right)-1
\end{cases}
\end{align*}
As remarked above $m\geq-v_{p}\left(t\right)$ makes the congruence
$p^{m}t\overset{1}{\equiv}0$ hold true, guaranteeing that $\alpha_{H}\left(p^{m}t\right)=\alpha_{H}\left(0\right)$
for all $m\geq-v_{p}\left(t\right)$. So, for fixed $t$, we prepare
our decomposition like so: 
\[
\prod_{m=0}^{n-1}\alpha_{H}\left(p^{m}t\right)=\begin{cases}
\prod_{m=0}^{-v_{p}\left(t\right)-2}\alpha_{H}\left(p^{m}t\right) & \textrm{if }n=-v_{p}\left(t\right)-1\\
\prod_{m=0}^{-v_{p}\left(t\right)-1}\alpha_{H}\left(p^{m}t\right) & \textrm{if }n=-v_{p}\left(t\right)\\
\prod_{m=0}^{-v_{p}\left(t\right)-1}\alpha_{H}\left(p^{m}t\right)\times\prod_{k=-v_{p}\left(t\right)}^{n-1}\alpha_{H}\left(p^{k}t\right) & \textrm{if }n\geq-v_{p}\left(t\right)+1
\end{cases}
\]
Since $\alpha_{H}\left(p^{k}t\right)=\alpha_{H}\left(0\right)$ for
all terms of the $k$-product on the bottom line, we can write: 
\begin{align*}
\prod_{m=0}^{n-1}\alpha_{H}\left(p^{m}t\right) & =\begin{cases}
\prod_{m=0}^{-v_{p}\left(t\right)-2}\alpha_{H}\left(p^{m}t\right) & \textrm{if }n=-v_{p}\left(t\right)-1\\
\left(\alpha_{H}\left(0\right)\right)^{n+v_{p}\left(t\right)}\prod_{m=0}^{-v_{p}\left(t\right)-1}\alpha_{H}\left(p^{m}t\right) & \textrm{if }n\geq-v_{p}\left(t\right)
\end{cases}\\
\left(\times\&\div\textrm{ by }\alpha_{H}\left(p^{-v_{p}\left(t\right)-1}t\right)\right); & =\begin{cases}
\frac{\hat{A}_{H}\left(t\right)}{\alpha_{H}\left(p^{-v_{p}\left(t\right)-1}t\right)} & \textrm{if }n=-v_{p}\left(t\right)-1\\
\left(\alpha_{H}\left(0\right)\right)^{n+v_{p}\left(t\right)}\hat{A}_{H}\left(t\right) & \textrm{if }n\geq-v_{p}\left(t\right)
\end{cases}
\end{align*}
which gives us the desired formula.

Q.E.D.

\vphantom{}

Now, we learn how to express the $\alpha_{H}$ product as a series. 
\begin{prop}[$\alpha_{H}$\textbf{ series formula}]
\label{prop:alpha product series expansion} 
\end{prop}
\begin{equation}
\prod_{m=0}^{n-1}\alpha_{H}\left(p^{m}t\right)=\left(\frac{\mu_{0}}{p^{2}}\right)^{n}\sum_{m=0}^{p^{n}-1}\kappa_{H}\left(m\right)e^{-2\pi imt},\textrm{ }\forall n\geq0,\textrm{ }\forall t\in\hat{\mathbb{Z}}_{p}\label{eq:alpha_H product expansion}
\end{equation}

Proof: We start by writing: 
\[
\prod_{m=0}^{n-1}\alpha_{H}\left(p^{m}t\right)=\prod_{m=0}^{n-1}\left(\sum_{j=0}^{p-1}\frac{\mu_{j}}{p^{2}}e^{-2\pi ij\left(p^{m}t\right)}\right)
\]
and then apply (\ref{eq:M_H partial sum generating identity}) from
\textbf{Proposition \ref{prop:Generating function identities}} with
$c_{j}=\mu_{j}/p^{2}$ and $z=e^{-2\pi it}$: 
\begin{align*}
\prod_{m=0}^{n-1}\alpha_{H}\left(p^{m}t\right) & =\left(\frac{\mu_{0}}{p^{2}}\right)^{n}\sum_{m=0}^{p^{n}-1}\left(\frac{\mu_{0}}{p^{2}}\right)^{-\lambda_{p}\left(m\right)}\left(\prod_{j=0}^{p-1}\left(\frac{\mu_{j}}{p^{2}}\right)^{\#_{p:j}\left(m\right)}\right)e^{-2\pi imt}\\
\left(\sum_{j=0}^{p-1}\#_{p:j}\left(m\right)=\lambda_{p}\left(m\right)\right); & =\left(\frac{\mu_{0}}{p^{2}}\right)^{n}\sum_{m=0}^{p^{n}-1}\left(\frac{\mu_{0}}{p^{2}}\right)^{-\lambda_{p}\left(m\right)}\left(\frac{\prod_{j=0}^{p-1}\mu_{j}^{\#_{p:j}\left(m\right)}}{p^{2\lambda_{p}\left(m\right)}}\right)e^{-2\pi imt}\\
\left(M_{H}\left(m\right)=\frac{\prod_{j=0}^{p-1}\mu_{j}^{\#_{p:j}\left(m\right)}}{p^{\lambda_{p}\left(m\right)}}\right); & =\left(\frac{\mu_{0}}{p^{2}}\right)^{n}\sum_{m=0}^{p^{n}-1}\left(\frac{\mu_{0}}{p^{2}}\right)^{-\lambda_{p}\left(m\right)}\frac{M_{H}\left(m\right)}{p^{\lambda_{p}\left(m\right)}}e^{-2\pi imt}\\
 & =\left(\frac{\mu_{0}}{p^{2}}\right)^{n}\sum_{m=0}^{p^{n}-1}\left(\frac{p}{\mu_{0}}\right)^{\lambda_{p}\left(m\right)}M_{H}\left(m\right)e^{-2\pi imt}\\
 & =\left(\frac{\mu_{0}}{p^{2}}\right)^{n}\sum_{m=0}^{p^{n}-1}\kappa_{H}\left(m\right)e^{-2\pi imt}
\end{align*}

Q.E.D.

\vphantom{}

This formula reveals $\hat{A}_{H}$ as being a radially-magnitudinal
function. A simple estimate shows that that $\hat{A}_{H}$ is the
Fourier-Stieltjes transform of a $\left(p,q\right)$-adic measure. 
\begin{prop}
$\hat{A}_{H}\left(t\right)$ is the Fourier-Stieltjes transform of
a $\left(p,q\right)$-adic measure. 
\end{prop}
Proof: Let $H$ be semi-basic. Then, $p$ and the $d_{j}$s are co-prime
to both $q$ and the $a_{j}$s. So, we get: 
\[
\left|\alpha_{H}\left(t\right)\right|_{q}=\left|\sum_{j=0}^{p-1}\frac{\mu_{j}}{p^{2}}e^{-2\pi ijt}\right|_{q}=\left|\sum_{j=0}^{p-1}\frac{a_{j}}{pd_{j}}e^{-2\pi ijt}\right|_{q}\leq\max_{0\leq j\leq p-1}\left|\frac{a_{j}}{pd_{j}}\right|_{q}\leq1
\]
Hence: 
\begin{align*}
\sup_{t\in\hat{\mathbb{Z}}_{p}}\left|\hat{A}_{H}\left(t\right)\right|_{q} & \leq\max\left\{ 1,\sup_{t\in\hat{\mathbb{Z}}_{p}\backslash\left\{ 0\right\} }\prod_{n=0}^{-v_{p}\left(t\right)-1}\left|\alpha_{H}\left(p^{n}t\right)\right|_{q}\right\} \\
 & \leq\max\left\{ 1,\sup_{t\in\hat{\mathbb{Z}}_{p}\backslash\left\{ 0\right\} }\prod_{n=0}^{-v_{p}\left(t\right)-1}1\right\} \\
 & =1
\end{align*}
This shows $\hat{A}_{H}$ is then in $B\left(\hat{\mathbb{Z}}_{p},\mathbb{C}_{q}\right)$.
As such, the map: 
\begin{equation}
f\in C\left(\mathbb{Z}_{p},\mathbb{C}_{q}\right)\mapsto\sum_{t\in\hat{\mathbb{Z}}_{p}}\hat{f}\left(-t\right)\hat{A}_{H}\left(t\right)\in\mathbb{C}_{q}
\end{equation}
is a continuous linear functional on $C\left(\mathbb{Z}_{p},\mathbb{C}_{q}\right)$\textemdash that
is, a $\left(p,q\right)$-adic measure.

Q.E.D.

\vphantom{}

Here, we introduce the notation for the partial sums of the Fourier
series generated by $\hat{A}_{H}$. 
\begin{defn}
We\index{$dA_{H}$} write \nomenclature{$dA_{H}$}{ }$dA_{H}$ to
denote the $\left(p,q\right)$-adic measure whose Fourier-Stieltjes
transform is $\hat{A}_{H}\left(t\right)$. We then have \nomenclature{$\tilde{A}_{H,N}\left(\mathfrak{z}\right)$}{$N$th partial Fourier series generated by $\hat{A}_{H}\left(t\right)$}:
\begin{equation}
\tilde{A}_{H,N}\left(\mathfrak{z}\right)\overset{\textrm{def}}{=}\sum_{\left|t\right|_{p}\leq p^{N}}\hat{A}_{H}\left(t\right)e^{2\pi i\left\{ t\mathfrak{z}\right\} _{p}}\label{eq:Definition of A_H,N twiddle}
\end{equation}
\end{defn}
\vphantom{}

The theorem given below summarizes the main properties of $dA_{H}$: 
\begin{thm}[\textbf{Properties of $dA_{H}$}]
\label{thm:Properties of dA_H}\index{$dA_{H}$!properties of} \ 

\vphantom{}

I. ($dA_{H}$ is radially-magnitudinal and $\left(p,q\right)$-adically
regular) 
\begin{equation}
\hat{A}_{H}\left(t\right)=\begin{cases}
1 & \textrm{if }t=0\\
\frac{1}{\left|t\right|_{p}^{2-\log_{p}\mu_{0}}}\sum_{m=0}^{\left|t\right|_{p}-1}\kappa_{H}\left(m\right)e^{-2\pi imt} & \textrm{else}
\end{cases},\textrm{ }\forall t\in\hat{\mathbb{Z}}_{p}\label{eq:A_H hat as the product of radially symmetric and magnitude-dependent measures}
\end{equation}

\vphantom{}

II. (Formula for $\tilde{A}_{H,N}\left(\mathfrak{z}\right)$) For
any $N\in\mathbb{N}_{0}$ and $\mathfrak{z}\in\mathbb{Z}_{p}$: 
\begin{equation}
\tilde{A}_{H,N}\left(\mathfrak{z}\right)\overset{\overline{\mathbb{Q}}}{=}\left(\frac{\mu_{0}}{p}\right)^{N}\kappa_{H}\left(\left[\mathfrak{z}\right]_{p^{N}}\right)+\left(1-\alpha_{H}\left(0\right)\right)\sum_{n=0}^{N-1}\left(\frac{\mu_{0}}{p}\right)^{n}\kappa_{H}\left(\left[\mathfrak{z}\right]_{p^{n}}\right)\label{eq:Convolution of dA_H and D_N}
\end{equation}

\vphantom{}

III. ($q$-adic convergence of $\tilde{A}_{H,N}\left(\mathfrak{z}\right)$)
As $N\rightarrow\infty$, \emph{(\ref{eq:Convolution of dA_H and D_N})}
converges in $\mathbb{C}_{q}$ (in fact, $\mathbb{Z}_{q}$) for all
$\mathfrak{z}\in\mathbb{Z}_{p}^{\prime}$, with: 
\begin{equation}
\lim_{N\rightarrow\infty}\tilde{A}_{H,N}\left(\mathfrak{z}\right)\overset{\mathbb{C}_{q}}{=}\left(1-\alpha_{H}\left(0\right)\right)\sum_{n=0}^{\infty}\left(\frac{\mu_{0}}{p}\right)^{n}\kappa_{H}\left(\left[\mathfrak{z}\right]_{p^{n}}\right),\textrm{ }\forall\mathfrak{z}\in\mathbb{Z}_{p}^{\prime}\label{eq:Derivative of dA_H on Z_rho prime}
\end{equation}
This convergence is point-wise.

\vphantom{}

IV. (Convergence of $\tilde{A}_{H,N}\left(\mathfrak{z}\right)$ in
$\mathbb{C}$) As $N\rightarrow\infty$, \emph{(\ref{eq:Convolution of dA_H and D_N})}
either converges in $\mathbb{C}$ for all $\mathfrak{z}\in\mathbb{N}_{0}$,
which occurs if and only if $H$ is contracting ($\mu_{0}/p<1$);
or diverges in $\mathbb{C}$ for all $\mathfrak{z}\in\mathbb{N}_{0}$,
which occurs if and only if $H$ is expanding ($\mu_{0}/p>1$). The
convergence/divergence in both cases is point-wise. Moreover, if convergence
occurs, the limit in $\mathbb{C}$ is given by 
\begin{equation}
\lim_{N\rightarrow\infty}\tilde{A}_{H,N}\left(m\right)\overset{\mathbb{C}}{=}\left(1-\alpha_{H}\left(0\right)\right)\sum_{n=0}^{\infty}\left(\frac{\mu_{0}}{p}\right)^{n}\kappa_{H}\left(\left[m\right]_{p^{n}}\right),\textrm{ }\forall m\in\mathbb{N}_{0}\label{eq:derivative of dA_H on N_0}
\end{equation}

\vphantom{}

V. (Characterization of degeneracy)\index{$dA_{H}$!degeneracy} If
$H$ is contracting, then $dA_{H}$ is degenerate if and only if $\alpha_{H}\left(0\right)=1$. 
\end{thm}
Proof:

I. By \textbf{Proposition \ref{prop:alpha product series expansion}},
we have that: 
\begin{equation}
\hat{A}_{H}\left(t\right)=\prod_{m=0}^{-v_{p}\left(t\right)-1}\alpha_{H}\left(p^{m}t\right)=\left(\frac{\mu_{0}}{p^{2}}\right)^{-v_{p}\left(t\right)}\sum_{m=0}^{p^{-v_{p}\left(t\right)}-1}\kappa_{H}\left(m\right)e^{-2\pi imt}
\end{equation}
and hence: 
\begin{equation}
\hat{A}_{H}\left(t\right)=\prod_{m=0}^{-v_{p}\left(t\right)-1}\alpha_{H}\left(p^{m}t\right)=\frac{1}{\left|t\right|_{p}^{2-\log_{p}\mu_{0}}}\sum_{m=0}^{\left|t\right|_{p}-1}\kappa_{H}\left(m\right)e^{-2\pi imt}
\end{equation}

\vphantom{}

II. Since:

\[
\hat{A}_{H}\left(t\right)=\prod_{n=0}^{-v_{p}\left(t\right)-1}\alpha_{H}\left(p^{n}t\right)=\left(\frac{\mu_{0}}{p^{2}}\right)^{-v_{p}\left(t\right)}\sum_{m=0}^{\left|t\right|_{p}-1}\kappa_{H}\left(m\right)e^{-2\pi imt}
\]
we can use the Radial-Magnitudinal Fourier Resummation Lemma (\textbf{Lemma
\ref{lem:Radially-Mag Fourier Resummation Lemma}}) to write: 
\begin{align*}
\sum_{\left|t\right|_{p}\leq p^{N}}\hat{A}_{H}\left(t\right)e^{2\pi i\left\{ t\mathfrak{z}\right\} _{p}} & \overset{\overline{\mathbb{Q}}}{=}\hat{A}_{H}\left(0\right)+\sum_{n=1}^{N}\left(\frac{\mu_{0}}{p^{2}}\right)^{n}\left(p^{n}\kappa_{H}\left(\left[\mathfrak{z}\right]_{p^{n}}\right)-p^{n-1}\kappa_{H}\left(\left[\mathfrak{z}\right]_{p^{n-1}}\right)\right)\\
 & -\sum_{j=1}^{p-1}\sum_{n=1}^{N}p^{n-1}\left(\frac{\mu_{0}}{p^{2}}\right)^{n}\kappa_{H}\left(\left[\mathfrak{z}\right]_{p^{n-1}}+jp^{n-1}\right)
\end{align*}
By the $\left(p,q\right)$-adic regularity of $\kappa_{H}$, we then
have:

\[
\sum_{\left|t\right|_{p}\leq p^{N}}\hat{A}_{H}\left(t\right)e^{2\pi i\left\{ t\mathfrak{z}\right\} _{p}}\overset{\overline{\mathbb{Q}}}{=}1+\sum_{n=1}^{N}\left(\frac{\mu_{0}}{p}\right)^{n}\kappa_{H}\left(\left[\mathfrak{z}\right]_{p^{n}}\right)-\underbrace{\left(\frac{1}{p}\sum_{j=0}^{p-1}\frac{\mu_{j}}{p}\right)}_{\alpha_{H}\left(0\right)}\sum_{n=0}^{N-1}\left(\frac{\mu_{0}}{p}\right)^{n}\kappa_{H}\left(\left[\mathfrak{z}\right]_{p^{n}}\right)
\]
and hence (since $\kappa_{H}\left(0\right)=1$): 
\begin{equation}
\tilde{A}_{H,N}\left(\mathfrak{z}\right)\overset{\overline{\mathbb{Q}}}{=}\left(\frac{\mu_{0}}{p}\right)^{N}\kappa_{H}\left(\left[\mathfrak{z}\right]_{p^{N}}\right)+\left(1-\alpha_{H}\left(0\right)\right)\sum_{n=0}^{N-1}\left(\frac{\mu_{0}}{p}\right)^{n}\kappa_{H}\left(\left[\mathfrak{z}\right]_{p^{n}}\right)
\end{equation}
as desired.

\vphantom{}

III. Fix $\mathfrak{z}\in\mathbb{Z}_{p}^{\prime}$. Then, since $H$
is semi-basic, we have: 
\begin{equation}
\left|\left(\frac{\mu_{0}}{p}\right)^{n}\kappa_{H}\left(\left[\mathfrak{z}\right]_{p^{n}}\right)\right|_{q}=\left|\kappa_{H}\left(\left[\mathfrak{z}\right]_{p^{n}}\right)\right|_{q}
\end{equation}
So, by \textbf{Proposition \ref{prop:Properties of Kappa_H}}, $\kappa_{H}\left(\left[\mathfrak{z}\right]_{p^{n}}\right)$
converges $q$-adically to $0$ as $n\rightarrow\infty$.The ultrametric
structure of $\mathbb{C}_{q}$ then guarantees the convergence of
(\ref{eq:Derivative of dA_H on Z_rho prime}).

\vphantom{}

IV. Let $\mathfrak{z}\in\mathbb{N}_{0}$. Then, $\kappa_{H}\left(\mathfrak{z}\right)=\kappa_{H}\left(\left[\mathfrak{z}\right]_{p^{n}}\right)=\kappa_{H}\left(\left[\mathfrak{z}\right]_{p^{n+1}}\right)\in\mathbb{Q}\backslash\left\{ 0\right\} $
for all $n\geq\lambda_{p}\left(\mathfrak{z}\right)$. So, let $c$
denote the value of $M_{H}\left(\left[\mathfrak{z}\right]_{p^{n}}\right)$
for all such $n$. Examining the tail of the $n$-series in (\ref{eq:Convolution of dA_H and D_N}),
$N\geq\lambda_{p}\left(\mathfrak{z}\right)+1$ implies: 
\begin{align*}
\sum_{n=0}^{N-1}\left(\frac{\mu_{0}}{p}\right)^{n}\kappa_{H}\left(\left[\mathfrak{z}\right]_{p^{n}}\right) & \overset{\mathbb{R}}{=}\sum_{n=0}^{\lambda_{p}\left(\mathfrak{z}\right)-1}\left(\frac{\mu_{0}}{p}\right)^{n}\kappa_{H}\left(\left[\mathfrak{z}\right]_{p^{n}}\right)+\sum_{n=\lambda_{p}\left(\mathfrak{z}\right)}^{N-1}\left(\frac{\mu_{0}}{p}\right)^{n}\kappa_{H}\left(\left[\mathfrak{z}\right]_{p^{n}}\right)\\
 & \overset{\mathbb{R}}{=}O\left(1\right)+\kappa_{H}\left(\mathfrak{z}\right)\sum_{n=\lambda_{p}\left(\mathfrak{z}\right)}^{N-1}\left(\frac{\mu_{0}}{p}\right)^{n}
\end{align*}
The series on the bottom line is a \emph{geometric} series, and hence,
is convergent in $\mathbb{C}$ if and only if $\mu_{0}/p<1$. Likewise,
by \textbf{Proposition \ref{prop:Properties of Kappa_H}},\textbf{
}$\left(\frac{\mu_{0}}{p}\right)^{N}\kappa_{H}\left(\left[\mathfrak{z}\right]_{p^{N}}\right)$
tends to $0$ as $N\rightarrow\infty$ if and only if $\mu_{0}/p<1$,
and tends to $\infty$ if and only if $\mu_{0}/p>1$.

\vphantom{}

V. Letting $H$ be semi-basic and contracting, (\ref{eq:Derivative of dA_H on Z_rho prime})
shows that $\lim_{N\rightarrow\infty}\tilde{A}_{H,N}$ will be identically
zero if and only if $\alpha_{H}\left(0\right)=1$.

Q.E.D.

\vphantom{}

Using the above, we can give the formulae for integrating $\left(p,q\right)$-adic
functions against $dA_{H}$. 
\begin{cor}[\textbf{$dA_{H}$ Integration Formulae}]
\ 

\vphantom{}

I. \index{$dA_{H}$!integration}($dA_{H}$ is a probability measure\footnote{In the sense of \cite{Probabilities taking values in non-archimedean fields,Measure-theoretic approach to p-adic probability theory}.}):
\begin{equation}
\int_{\mathbb{Z}_{p}}dA_{H}\left(\mathfrak{z}\right)\overset{\mathbb{C}_{q}}{=}1\label{eq:dA_H is a probabiity measure}
\end{equation}
Also, for all $n\in\mathbb{N}_{1}$ and all $k\in\left\{ 0,\ldots,p^{n}-1\right\} $:
\begin{align}
\int_{\mathbb{Z}_{p}}\left[\mathfrak{z}\overset{p^{n}}{\equiv}k\right]dA_{H}\left(\mathfrak{z}\right) & \overset{\overline{\mathbb{Q}}}{=}\left(\frac{\mu_{0}}{p}\right)^{n}\kappa_{H}\left(k\right)+\left(1-\alpha_{H}\left(0\right)\right)\sum_{m=0}^{n-1}\left(\frac{\mu_{0}}{p}\right)^{m}\kappa_{H}\left(\left[k\right]_{p^{m}}\right)\label{eq:integration of indicator functions dA_H}
\end{align}

\vphantom{}

II. Let $f\in C\left(\mathbb{Z}_{p},\mathbb{C}_{q}\right)$. Then:
\begin{align}
\int_{\mathbb{Z}_{p}}f\left(\mathfrak{z}\right)dA_{H}\left(\mathfrak{z}\right) & \overset{\mathbb{C}_{q}}{=}\lim_{N\rightarrow\infty}\left(\left(\frac{\mu_{0}}{p^{2}}\right)^{N}\sum_{n=0}^{p^{N}-1}\kappa_{H}\left(n\right)f\left(n\right)\right.\label{eq:Truncation-based formula for integral of f dA_H}\\
 & +\left.\frac{1-\alpha_{H}\left(0\right)}{p^{N}}\sum_{n=0}^{p^{N}-1}\sum_{m=0}^{N-1}\left(\frac{\mu_{0}}{p}\right)^{m}\kappa_{H}\left(\left[n\right]_{p^{m}}\right)f\left(n\right)\right)\nonumber 
\end{align}
\end{cor}
\begin{rem}
One can also compute $\int_{\mathbb{Z}_{p}}f\left(\mathfrak{z}\right)dA_{H}\left(\mathfrak{z}\right)$
by expressing $f$ as a van der Put series and then integrating term-by-term;
this yields formulae nearly identical to the ones given above, albeit
in terms of $c_{n}\left(f\right)$ (the $n$th van der Put coefficient
of $f$), rather than $f\left(n\right)$. 
\end{rem}
Proof: (I) follows by using (\ref{eq:Convolution of dA_H and D_N})
along with the fact that, by definition: 
\begin{equation}
\int_{\mathbb{Z}_{p}}\left[\mathfrak{z}\overset{p^{n}}{\equiv}k\right]dA_{H}\left(\mathfrak{z}\right)=\frac{1}{p^{n}}\sum_{\left|t\right|_{p}\leq p^{n}}\hat{A}_{H}\left(t\right)e^{-2\pi i\left\{ tk\right\} _{p}}
\end{equation}
and that $\left[k\right]_{p^{n}}=k$ if and only if $k\in\left\{ 0,\ldots,p^{n}-1\right\} $.
As for (II), let $f\in C\left(\mathbb{Z}_{p},\mathbb{C}_{q}\right)$
be arbitrary. Then, the truncations $f_{N}$ converge $q$-adically
to $f$ uniformly over $\mathbb{Z}_{p}$. Consequently: 
\begin{align*}
\int_{\mathbb{Z}_{p}}f\left(\mathfrak{z}\right)dA_{H}\left(\mathfrak{z}\right) & \overset{\mathbb{C}_{q}}{=}\lim_{N\rightarrow\infty}\sum_{n=0}^{p^{N}-1}f\left(n\right)\int_{\mathbb{Z}_{p}}\left[\mathfrak{z}\overset{p^{N}}{\equiv}n\right]dA_{H}\left(\mathfrak{z}\right)\\
\left(\textrm{Use (I)}\right); & =\lim_{N\rightarrow\infty}\sum_{n=0}^{p^{N}-1}f\left(n\right)\left(\frac{1}{p^{N}}\left(\frac{\mu_{0}}{p}\right)^{N}\kappa_{H}\left(n\right)\right)\\
 & +\lim_{N\rightarrow\infty}\sum_{n=0}^{p^{N}-1}f\left(n\right)\frac{1-\alpha_{H}\left(0\right)}{p^{N}}\sum_{m=0}^{N-1}\left(\frac{\mu_{0}}{p}\right)^{m}\kappa_{H}\left(\left[n\right]_{p^{m}}\right)
\end{align*}

Q.E.D.

\vphantom{}

Given the import of the shortened $qx+1$ map\index{$qx+1$ map}s,
let me give the details for what happens with their associated $A_{H}$s.
\begin{cor}[\textbf{$dA_{H}$ for $qx+1$}]
\label{cor:Formula for Nth partial sum of the Fourier series generated by A_q hat}Let
$q$ be an odd integer $\geq3$, and let $A_{q}$ denote $A_{H}$
where $H=T_{q}$ is the shortened $qx+1$ map. Finally, write: 
\begin{equation}
\tilde{A}_{q,N}\left(\mathfrak{z}\right)\overset{\textrm{def}}{=}\sum_{\left|t\right|_{2}\leq2^{N}}\hat{A}_{q}\left(t\right)e^{2\pi i\left\{ t\mathfrak{z}\right\} _{2}}=\left(D_{2:N}*dA_{q}\right)\left(\mathfrak{z}\right)\label{eq:Definition of A_q,n twiddle}
\end{equation}
Then: 
\begin{equation}
\tilde{A}_{q,N}\left(\mathfrak{z}\right)=\frac{q^{\#_{1}\left(\left[\mathfrak{z}\right]_{2^{N}}\right)}}{2^{N}}+\left(1-\frac{q+1}{4}\right)\sum_{n=0}^{N-1}\frac{q^{\#_{1}\left(\left[\mathfrak{z}\right]_{2^{n}}\right)}}{2^{n}}\label{eq:Convolution of dA_q and D_N}
\end{equation}
where, for any integer $m\geq0$, $\#_{1}\left(m\right)$ is the number
of $1$s in the binary representation of $m$ (a shorthand for $\#_{2:1}\left(m\right)$).
Integration against $dA_{q}$ is then given by: 
\begin{align}
\int_{\mathbb{Z}_{2}}f\left(\mathfrak{z}\right)dA_{q}\left(\mathfrak{z}\right) & \overset{\mathbb{C}_{q}}{=}\lim_{N\rightarrow\infty}\frac{1}{4^{N}}\sum_{n=0}^{2^{N}-1}f\left(n\right)q^{\#_{1}\left(n\right)}\\
 & +\frac{3-q}{4}\lim_{N\rightarrow\infty}\frac{1}{2^{N}}\sum_{n=0}^{2^{N}-1}f\left(n\right)\sum_{m=0}^{N-1}\frac{q^{\#_{1}\left(\left[n\right]_{2^{m}}\right)}}{2^{m}}\nonumber 
\end{align}

In particular, for $q=3$ (the Shortened Collatz map\index{Collatz!map}),
we have: 
\begin{equation}
\tilde{A}_{3,N}\left(\mathfrak{z}\right)=\frac{3^{\#_{1}\left(\left[\mathfrak{z}\right]_{2^{N}}\right)}}{2^{N}}\label{eq:A_3 measure of a coset in Z_2}
\end{equation}
and hence: 
\begin{equation}
\int_{\mathbb{Z}_{2}}f\left(\mathfrak{z}\right)dA_{3}\left(\mathfrak{z}\right)\overset{\mathbb{C}_{3}}{=}\lim_{N\rightarrow\infty}\frac{1}{4^{N}}\sum_{n=0}^{2^{N}-1}3^{\#_{1}\left(n\right)}f\left(n\right),\textrm{ }\forall f\in C\left(\mathbb{Z}_{2},\mathbb{C}_{3}\right)\label{eq:Integral of f dA_3}
\end{equation}
\end{cor}
Proof: Use (\ref{eq:Truncation-based formula for integral of f dA_H}).

Q.E.D.

\vphantom{}

By using $\hat{A}_{H}$ and \textbf{Proposition \ref{prop:alpha product in terms of A_H hat}},
we can re-write our explicit formula for $\hat{\chi}_{H,N}$\textemdash equation
(\ref{eq:Explicit formula for Chi_H,N hat}) from \textbf{Proposition
\ref{prop:Computation of Chi_H,N hat}}\textemdash in a way that separates
the dependence of $t$ and $N$. Like in previous computations, the
formulae for \index{hat{chi}{H,N}@$\hat{\chi}_{H,N}$!left(N,tright) asymptotic decomposition@$\left(N,t\right)$
asymptotic decomposition}$\hat{\chi}_{H,N}$ come in two forms, based on whether or not $\alpha_{H}\left(0\right)=1$. 
\begin{thm}[\textbf{$\left(N,t\right)$-Asymptotic Decomposition of $\hat{\chi}_{H,N}$}]
\label{thm:(N,t) asymptotic decomposition of Chi_H,N hat}Recall
that we write $\gamma_{H}\left(t\right)$ to denote $\beta_{H}\left(t\right)/\alpha_{H}\left(t\right)$.

\vphantom{}

I. If $\alpha_{H}\left(0\right)=1$, then: 
\begin{equation}
\hat{\chi}_{H,N}\left(t\right)=\begin{cases}
\beta_{H}\left(0\right)N\hat{A}_{H}\left(t\right) & \textrm{if }t=0\\
\left(\gamma_{H}\left(\frac{t\left|t\right|_{p}}{p}\right)+\beta_{H}\left(0\right)\left(N+v_{p}\left(t\right)\right)\right)\hat{A}_{H}\left(t\right) & \textrm{if }0<\left|t\right|_{p}<p^{N}\\
\gamma_{H}\left(\frac{t\left|t\right|_{p}}{p}\right)\hat{A}_{H}\left(t\right) & \textrm{if }\left|t\right|_{p}=p^{N}\\
0 & \textrm{if }\left|t\right|_{p}>p^{N}
\end{cases},\textrm{ }\forall t\in\hat{\mathbb{Z}}_{p}\label{eq:Fine Structure Formula for Chi_H,N hat when alpha is 1}
\end{equation}

\vphantom{}

II. If $\alpha_{H}\left(0\right)\neq1$. Then: 
\begin{equation}
\hat{\chi}_{H,N}\left(t\right)=\begin{cases}
\beta_{H}\left(0\right)\frac{\left(\alpha_{H}\left(0\right)\right)^{N}-1}{\alpha_{H}\left(0\right)-1}\hat{A}_{H}\left(t\right) & \textrm{if }t=0\\
\left(\gamma_{H}\left(\frac{t\left|t\right|_{p}}{p}\right)+\beta_{H}\left(0\right)\frac{\left(\alpha_{H}\left(0\right)\right)^{N+v_{p}\left(t\right)}-1}{\alpha_{H}\left(0\right)-1}\right)\hat{A}_{H}\left(t\right) & \textrm{if }0<\left|t\right|_{p}<p^{N}\\
\gamma_{H}\left(\frac{t\left|t\right|_{p}}{p}\right)\hat{A}_{H}\left(t\right) & \textrm{if }\left|t\right|_{p}=p^{N}\\
0 & \textrm{if }\left|t\right|_{p}>p^{N}
\end{cases},\textrm{ }\forall t\in\hat{\mathbb{Z}}_{p}\label{eq:Fine Structure Formula for Chi_H,N hat when alpha is not 1}
\end{equation}
\end{thm}
Proof: For brevity, we write: 
\begin{equation}
\hat{A}_{H,n}\left(t\right)\overset{\textrm{def}}{=}\begin{cases}
1 & \textrm{if }n=0\\
\mathbf{1}_{0}\left(p^{n+1}t\right)\prod_{m=0}^{n-1}\alpha_{H}\left(p^{m}t\right) & \textrm{if }n\geq1
\end{cases}\label{eq:Definition of A_H,n+1 hat}
\end{equation}
so that (\ref{eq:Explicit formula for Chi_H,N hat}) becomes: 
\[
\hat{\chi}_{H,N}\left(t\right)=\sum_{n=0}^{N-1}\beta_{H}\left(p^{n}t\right)\hat{A}_{H,n}\left(t\right)
\]
Then, letting $t\in\hat{\mathbb{Z}}_{p}$ be non-zero and satisfy
$\left|t\right|_{p}\leq p^{N}$, we use \textbf{Proposition \ref{prop:alpha product in terms of A_H hat}}
to write: 
\begin{equation}
\hat{\chi}_{H,N}\left(t\right)\overset{\overline{\mathbb{Q}}}{=}\sum_{n=0}^{N-1}\beta_{H}\left(p^{n}t\right)\hat{A}_{H,n}\left(t\right)=\sum_{n=-v_{p}\left(t\right)-1}^{N-1}\beta_{H}\left(p^{n}t\right)\hat{A}_{H,n}\left(t\right)\label{eq:Chi_H,N hat as Beta_H plus A_H,n hat - ready for t,n analysis}
\end{equation}
Here, we have used the fact that: 
\begin{equation}
\hat{A}_{H,n}\left(t\right)=0,\textrm{ }\forall n\leq-v_{p}\left(t\right)-2
\end{equation}
because $\mathbf{1}_{0}\left(p^{n+1}t\right)$ vanishes whenever $n\leq-v_{p}\left(t\right)-2$.
This leaves us with two cases.

\vphantom{}

I. First, suppose $\left|t\right|_{p}=p^{N}$. Then $N=-v_{p}\left(t\right)$,
and so, (\ref{eq:Chi_H,N hat as Beta_H plus A_H,n hat - ready for t,n analysis})
becomes: 
\begin{align*}
\hat{\chi}_{H,N}\left(t\right) & =\sum_{n=-v_{p}\left(t\right)-1}^{N-1}\beta_{H}\left(p^{n}t\right)\hat{A}_{H,n}\left(t\right)\\
 & =\beta_{H}\left(p^{-v_{p}\left(t\right)-1}t\right)\hat{A}_{H,-v_{p}\left(t\right)-1}\left(t\right)\\
 & =\beta_{H}\left(p^{-v_{p}\left(t\right)-1}t\right)\cdot\mathbf{1}_{0}\left(p^{-v_{p}\left(t\right)-1+1}t\right)\prod_{m=0}^{-v_{p}\left(t\right)-1-1}\alpha_{H}\left(p^{m}t\right)\\
 & =\beta_{H}\left(p^{-v_{p}\left(t\right)-1}t\right)\cdot\frac{\mathbf{1}_{0}\left(p^{-v_{p}\left(t\right)}t\right)}{\alpha_{H}\left(p^{-v_{p}\left(t\right)-1}t\right)}\prod_{m=0}^{-v_{p}\left(t\right)-1}\alpha_{H}\left(p^{m}t\right)
\end{align*}
Since $p^{-v_{p}\left(t\right)}t=\left|t\right|_{p}t$ is the integer
in the numerator of $t$, $\mathbf{1}_{0}\left(p^{-v_{p}\left(t\right)}t\right)=1$,
and we are left with: 
\begin{align*}
\hat{\chi}_{H,N}\left(t\right) & =\frac{\beta_{H}\left(p^{-v_{p}\left(t\right)-1}t\right)}{\alpha_{H}\left(p^{-v_{p}\left(t\right)-1}t\right)}\underbrace{\prod_{m=0}^{-v_{p}\left(t\right)-1}\alpha_{H}\left(p^{m}t\right)}_{\hat{A}_{H}\left(t\right)}\\
 & =\frac{\beta_{H}\left(\frac{t\left|t\right|_{p}}{p}\right)}{\alpha_{H}\left(\frac{t\left|t\right|_{p}}{p}\right)}\hat{A}_{H}\left(t\right)\\
 & =\gamma_{H}\left(\frac{t\left|t\right|_{p}}{p}\right)\hat{A}_{H}\left(t\right)
\end{align*}

\vphantom{}

II. Second, suppose $0<\left|t\right|_{p}<p^{N}$. Then $-v_{p}\left(t\right)-1$
is strictly less than $N-1$, and so $p^{n}t\overset{1}{\equiv}0$
for all $n\geq-v_{p}\left(t\right)$. With this, (\ref{eq:Chi_H,N hat as Beta_H plus A_H,n hat - ready for t,n analysis})
becomes:

\begin{align*}
\hat{\chi}_{H,N}\left(t\right) & =\beta_{H}\left(p^{-v_{p}\left(t\right)-1}t\right)\hat{A}_{H,-v_{p}\left(t\right)-1}\left(t\right)+\sum_{n=-v_{p}\left(t\right)}^{N-1}\beta_{H}\left(p^{n}t\right)\hat{A}_{H,n}\left(t\right)\\
\left(p^{n}t\overset{1}{\equiv}0\textrm{ }\forall n\geq-v_{p}\left(t\right)\right); & =\gamma_{H}\left(\frac{t\left|t\right|_{p}}{p}\right)\hat{A}_{H}\left(t\right)+\beta_{H}\left(0\right)\sum_{n=-v_{p}\left(t\right)}^{N-1}\hat{A}_{H,n}\left(t\right)
\end{align*}
Using \textbf{Proposition \ref{prop:alpha product in terms of A_H hat}},
we get: 
\begin{equation}
\hat{A}_{H,n}\left(t\right)=\left(\alpha_{H}\left(0\right)\right)^{n+v_{p}\left(t\right)}\hat{A}_{H}\left(t\right),\textrm{ }\forall n\geq-v_{p}\left(t\right)
\end{equation}
and hence: 
\begin{align*}
\hat{\chi}_{H,N}\left(t\right) & =\gamma_{H}\left(\frac{t\left|t\right|_{p}}{p}\right)\hat{A}_{H}\left(t\right)+\beta_{H}\left(0\right)\sum_{n=-v_{p}\left(t\right)}^{N-1}\left(\alpha_{H}\left(0\right)\right)^{n+v_{p}\left(t\right)}\hat{A}_{H}\left(t\right)\\
 & =\gamma_{H}\left(\frac{t\left|t\right|_{p}}{p}\right)\hat{A}_{H}\left(t\right)+\beta_{H}\left(0\right)\left(\sum_{n=0}^{N+v_{p}\left(t\right)-1}\left(\alpha_{H}\left(0\right)\right)^{n}\right)\hat{A}_{H}\left(t\right)\\
 & =\begin{cases}
\gamma_{H}\left(\frac{t\left|t\right|_{p}}{p}\right)\hat{A}_{H}\left(t\right)+\beta_{H}\left(0\right)\left(N+v_{p}\left(t\right)\right)\hat{A}_{H}\left(t\right) & \textrm{if }\alpha_{H}\left(0\right)=1\\
\gamma_{H}\left(\frac{t\left|t\right|_{p}}{p}\right)\hat{A}_{H}\left(t\right)+\beta_{H}\left(0\right)\frac{\left(\alpha_{H}\left(0\right)\right)^{N+v_{p}\left(t\right)}-1}{\alpha_{H}\left(0\right)-1}\hat{A}_{H}\left(t\right) & \textrm{if }\alpha_{H}\left(0\right)\neq1
\end{cases}
\end{align*}

Finally, when $t=0$: 
\begin{align*}
\hat{\chi}_{H,N}\left(0\right) & =\sum_{n=0}^{N-1}\beta_{H}\left(0\right)\hat{A}_{H,n}\left(0\right)\\
 & =\beta_{H}\left(0\right)+\sum_{n=1}^{N-1}\beta_{H}\left(0\right)\mathbf{1}_{0}\left(0\right)\prod_{m=0}^{n-1}\alpha_{H}\left(0\right)\\
 & =\beta_{H}\left(0\right)+\sum_{n=1}^{N-1}\beta_{H}\left(0\right)\left(\alpha_{H}\left(0\right)\right)^{n}\\
 & =\beta_{H}\left(0\right)\sum_{n=0}^{N-1}\left(\alpha_{H}\left(0\right)\right)^{n}\\
 & =\begin{cases}
\beta_{H}\left(0\right)N & \textrm{if }\alpha_{H}\left(0\right)=1\\
\beta_{H}\left(0\right)\frac{\left(\alpha_{H}\left(0\right)\right)^{N}-1}{\alpha_{H}\left(0\right)-1} & \textrm{if }\alpha_{H}\left(0\right)\neq1
\end{cases}
\end{align*}

Since $\hat{A}_{H}\left(0\right)=1$, we have that $\beta_{H}\left(0\right)N=\beta_{H}\left(0\right)N\hat{A}_{H}\left(0\right)$
for the $\alpha_{H}\left(0\right)=1$ case and $\beta_{H}\left(0\right)\frac{\left(\alpha_{H}\left(0\right)\right)^{N}-1}{\alpha_{H}\left(0\right)-1}=\beta_{H}\left(0\right)\frac{\left(\alpha_{H}\left(0\right)\right)^{N}-1}{\alpha_{H}\left(0\right)-1}\hat{A}_{H}\left(0\right)$
for the $\alpha_{H}\left(0\right)\neq1$ case.

Q.E.D.

\vphantom{}

Subtracting $\beta_{H}\left(0\right)N\hat{A}_{H}\left(t\right)\mathbf{1}_{0}\left(p^{N-1}t\right)$
(which is $0$ for all $\left|t\right|_{p}\geq p^{N}$) from the $\alpha_{H}\left(0\right)=1$
case (equation (\ref{eq:Fine Structure Formula for Chi_H,N hat when alpha is 1}))\index{hat{chi}{H,N}@$\hat{\chi}_{H,N}$!fine structure}
we obtain: 
\begin{equation}
\hat{\chi}_{H,N}\left(t\right)-\beta_{H}\left(0\right)N\hat{A}_{H}\left(t\right)\mathbf{1}_{0}\left(p^{N-1}t\right)=\begin{cases}
0 & \textrm{if }t=0\\
\left(\gamma_{H}\left(\frac{t\left|t\right|_{p}}{p}\right)+\beta_{H}\left(0\right)v_{p}\left(t\right)\right)\hat{A}_{H}\left(t\right) & \textrm{if }0<\left|t\right|_{p}<p^{N}\\
\gamma_{H}\left(\frac{t\left|t\right|_{p}}{p}\right)\hat{A}_{H}\left(t\right) & \textrm{if }\left|t\right|_{p}=p^{N}\\
0 & \textrm{if }\left|t\right|_{p}>p^{N}
\end{cases}\label{eq:Fine structure of Chi_H,N hat when alpha is 1}
\end{equation}
Even though the ranges of $t$ assigned to the individual pieces of
this formula depends on $N$, note that the actual \emph{values} of
the pieces on the right-hand side are independent of $N$ for all
$\left|t\right|_{p}<p^{N}$. So, while the Fourier transform $\hat{\chi}_{H}\left(t\right)$
might not exist in a natural way, by having subtracted off the ``divergent''
$\beta_{H}\left(0\right)N\hat{A}_{H}\left(t\right)\mathbf{1}_{0}\left(p^{N-1}t\right)$
term from $\hat{\chi}_{H,N}\left(t\right)$, we have discovered fine-scale
behavior of $\hat{\chi}_{H,N}\left(t\right)$ which was hidden beneath
the tumult. Indeed, for fixed $t$, the $q$-adic limit of the left-hand
side of (\ref{eq:Fine structure of Chi_H,N hat when alpha is 1})
as $N\rightarrow\infty$ is: 
\begin{equation}
\begin{cases}
0 & \textrm{if }t=0\\
\left(\gamma_{H}\left(\frac{t\left|t\right|_{p}}{p}\right)+\beta_{H}\left(0\right)v_{p}\left(t\right)\right)\hat{A}_{H}\left(t\right) & \textrm{if }\left|t\right|_{p}>0
\end{cases}\label{eq:Point-wise limit of the left-hand side of the fine structure formula when alpha is 1}
\end{equation}
Moreover, for any given $t$, the left-hand side of (\ref{eq:Fine structure of Chi_H,N hat when alpha is 1})
is actually \emph{equal} to the above limit provided that $N>-v_{p}\left(t\right)$.
This strongly suggests that (\ref{eq:Point-wise limit of the left-hand side of the fine structure formula when alpha is 1})
is, or is very close to a ``valid'' formula for a Fourier transform
of $\chi_{H}$. The next subsection is dedicated to making this intuition
rigorous.

\subsection{\label{subsec:4.2.2}$\hat{\chi}_{H}$ and $\tilde{\chi}_{H,N}$}

We begin by computing the Fourier series generated by $v_{p}\left(t\right)\hat{A}_{H}\left(t\right)$.
\begin{lem}[\textbf{$v_{p}\hat{A}_{H}$ summation formulae}]
\label{lem:v_p A_H hat summation formulae}\ 

\vphantom{}

I. 
\begin{equation}
\sum_{0<\left|t\right|_{p}\leq p^{N}}v_{p}\left(t\right)\hat{A}_{H}\left(t\right)e^{2\pi i\left\{ t\mathfrak{z}\right\} _{p}}\overset{\overline{\mathbb{Q}}}{=}-N\left(\frac{\mu_{0}}{p}\right)^{N}\kappa_{H}\left(\left[\mathfrak{z}\right]_{p^{N}}\right)+\sum_{n=0}^{N-1}\left(\alpha_{H}\left(0\right)\left(n+1\right)-n\right)\left(\frac{\mu_{0}}{p}\right)^{n}\kappa_{H}\left(\left[\mathfrak{z}\right]_{p^{n}}\right)\label{eq:Fourier sum of A_H hat v_p}
\end{equation}

\vphantom{}

II. If $H$ is semi-basic and contracting: 
\begin{equation}
\sum_{t\in\hat{\mathbb{Z}}_{p}\backslash\left\{ 0\right\} }v_{p}\left(t\right)\hat{A}_{H}\left(t\right)e^{2\pi i\left\{ t\mathfrak{z}\right\} _{p}}\overset{\mathbb{C}_{q}}{=}\sum_{n=0}^{\infty}\left(\alpha_{H}\left(0\right)\left(n+1\right)-n\right)\left(\frac{\mu_{0}}{p}\right)^{n}\kappa_{H}\left(\left[\mathfrak{z}\right]_{p^{n}}\right),\textrm{ }\forall\mathfrak{z}\in\mathbb{Z}_{p}^{\prime}\label{eq:Limit of Fourier sum of v_p A_H hat}
\end{equation}
where the convergence is point-wise. 
\end{lem}
Proof: For (I), use (\ref{eq:Convolution of dA_H and D_N}) from \textbf{Theorem
\ref{thm:Properties of dA_H}}, along with \textbf{Proposition \ref{prop:v_p of t times mu hat sum}}.
If $H$ is semi-basic and contracting, the decay estimates on $\kappa_{H}$
and $\left(\frac{\mu_{0}}{p}\right)^{N}\kappa_{H}\left(\left[\mathfrak{z}\right]_{p^{N}}\right)$
given by \textbf{Proposition \ref{prop:Properties of Kappa_H}} then
guarantee the convergence of (I) to (II) in $\mathbb{C}_{q}$ for
$\mathfrak{z}\in\mathbb{Z}_{p}^{\prime}$ as $N\rightarrow\infty$.

Q.E.D.

\vphantom{}

Next, we need to deal with how $p$ affects the situation. The primary
complication in our computations comes from $\gamma_{H}\left(t\left|t\right|_{p}/p\right)$,
whose behavior bifurcates at $p=2$. The function: 
\begin{equation}
t\in\hat{\mathbb{Z}}_{p}\mapsto\frac{t\left|t\right|_{p}}{p}\in\hat{\mathbb{Z}}_{p}\label{eq:Projection onto 1 over rho}
\end{equation}
is a projection which sends every non-zero element $t=k/p^{n}$ in
$\hat{\mathbb{Z}}_{p}$ and outputs the fraction $\left[k\right]_{p}/p$
obtained by reducing $k$ mod $p$ and then sticking it over $p$.
When $p=2$, $\left[k\right]_{p}=1$ for all $t\in\hat{\mathbb{Z}}_{2}$,
making (\ref{eq:Projection onto 1 over rho}) constant on $\hat{\mathbb{Z}}_{2}\backslash\left\{ 0\right\} $,
where it takes the value $1/2$. When $p\geq3$, however, (\ref{eq:Projection onto 1 over rho})
is no longer constant on $\hat{\mathbb{Z}}_{p}\backslash\left\{ 0\right\} $,
which is going to significantly complicate our computations. This
can be made much more manageable, provided the reader will put up
with another bit of notation. 
\begin{defn}[$\varepsilon_{n}\left(\mathfrak{z}\right)$]
\nomenclature{$\varepsilon_{n}\left(\mathfrak{z}\right)$}{ }For
each $n\in\mathbb{N}_{0}$, we define $\varepsilon_{n}:\mathbb{Z}_{p}\rightarrow\overline{\mathbb{Q}}$
by: 
\begin{equation}
\varepsilon_{n}\left(\mathfrak{z}\right)\overset{\textrm{def}}{=}e^{\frac{2\pi i}{p^{n+1}}\left(\left[\mathfrak{z}\right]_{p^{n+1}}-\left[\mathfrak{z}\right]_{p^{n}}\right)}=e^{2\pi i\left\{ \frac{\mathfrak{z}}{p^{n+1}}\right\} _{p}}e^{-\frac{2\pi i}{p}\left\{ \frac{\mathfrak{z}}{p^{n}}\right\} _{p}}\label{eq:Definition of epsilon_n}
\end{equation}
\end{defn}
\vphantom{}

As with nearly every other noteworthy function in this dissertation,
the $\varepsilon_{n}$s satisfy functional equations, which we record
below.
\begin{prop}[\textbf{Properties of $\varepsilon_{n}$}]
\ 

\vphantom{}

I.

\begin{equation}
\varepsilon_{0}\left(\mathfrak{z}\right)=e^{2\pi i\left\{ \frac{\mathfrak{z}}{p}\right\} _{p}},\textrm{ }\forall\mathfrak{z}\in\mathbb{Z}_{p}\label{eq:Epsilon 0 of z}
\end{equation}

\begin{equation}
\varepsilon_{n}\left(j\right)=1,\textrm{ }\forall j\in\mathbb{Z}/p\mathbb{Z},\textrm{ }\forall n\geq1\label{eq:Epsilon_n of j}
\end{equation}

\vphantom{}

II. 
\begin{equation}
\varepsilon_{n}\left(pm+j\right)=\begin{cases}
\varepsilon_{0}\left(j\right) & \textrm{if }n=0\\
\varepsilon_{n-1}\left(m\right) & \textrm{if }n\geq1
\end{cases},\textrm{ }\forall m\in\mathbb{N}_{0},\textrm{ }\forall j\in\mathbb{Z}/p\mathbb{Z},\textrm{ }\forall n\geq1\label{eq:epsilon_n functional equations}
\end{equation}

\vphantom{}

III. Let $\mathfrak{z}\neq0$. Then $\varepsilon_{n}\left(\mathfrak{z}\right)=1\textrm{ }$
for all $n<v_{p}\left(\mathfrak{z}\right)$. 
\end{prop}
Proof:

I. The identity (\ref{eq:Epsilon 0 of z}) is immediate from the definition
of $\varepsilon_{n}$. As for the other identity, note that for $j\in\mathbb{Z}/p\mathbb{Z}$,
$\left[j\right]_{p^{n}}=j$ for all $n\geq1$. Hence, for $n\geq1$:
\begin{equation}
\varepsilon_{n}\left(j\right)=e^{\frac{2\pi i}{p^{n+1}}\left(\left[j\right]_{p^{n+1}}-\left[j\right]_{p^{n}}\right)}=e^{\frac{2\pi i}{p^{n+1}}\cdot0}=1
\end{equation}

\vphantom{}

II. 
\begin{align*}
\varepsilon_{n}\left(pm+j\right) & =e^{2\pi i\left\{ \frac{pm+j}{p^{n+1}}\right\} _{p}}e^{-\frac{2\pi i}{p}\left\{ \frac{pm+j}{p^{n}}\right\} _{p}}\\
 & =e^{2\pi i\left\{ \frac{j}{p^{n+1}}\right\} _{p}}e^{-\frac{2\pi i}{p}\left\{ \frac{j}{p^{n}}\right\} _{p}}\cdot e^{2\pi i\left\{ \frac{m}{p^{n}}\right\} _{p}}e^{-\frac{2\pi i}{p}\left\{ \frac{m}{p^{n-1}}\right\} _{p}}\\
 & =\varepsilon_{n}\left(j\right)\varepsilon_{n-1}\left(m\right)\\
\left(\textrm{by (I)}\right); & =\begin{cases}
\varepsilon_{0}\left(j\right) & \textrm{if }n=0\\
\varepsilon_{n-1}\left(m\right) & \textrm{if }n\geq1
\end{cases}
\end{align*}

\vphantom{}

III. Let $\mathfrak{z}$ be non-zero. When $n<v_{p}\left(\mathfrak{z}\right)$,
we have that $p^{-n}\mathfrak{z}$ and $p^{-\left(n+1\right)}\mathfrak{z}$
are then $p$-adic integers, and hence: 
\begin{equation}
\varepsilon_{n}\left(\mathfrak{z}\right)=e^{2\pi i\left\{ \frac{\mathfrak{z}}{p^{n+1}}\right\} _{p}}e^{-\frac{2\pi i}{p}\left\{ \frac{\mathfrak{z}}{p^{n}}\right\} _{p}}=e^{0}\cdot e^{-0}=1
\end{equation}

Q.E.D.

\vphantom{}

Now we compute the sum of the Fourier series generated by $\gamma_{H}\left(t\left|t\right|_{p}/p\right)\hat{A}_{H}\left(t\right)$. 
\begin{lem}[\textbf{$\gamma_{H}\hat{A}_{H}$ summation formulae}]
\label{lem:1D gamma formula}Let $p\geq2$. Then:

\vphantom{}

I.

\begin{align}
\sum_{0<\left|t\right|_{p}\leq p^{N}}\gamma_{H}\left(\frac{t\left|t\right|_{p}}{p}\right)\hat{A}_{H}\left(t\right)e^{2\pi i\left\{ t\mathfrak{z}\right\} _{p}} & \overset{\overline{\mathbb{Q}}}{=}\sum_{n=0}^{N-1}\left(\sum_{j=1}^{p-1}\beta_{H}\left(\frac{j}{p}\right)\varepsilon_{n}^{j}\left(\mathfrak{z}\right)\right)\left(\frac{\mu_{0}}{p}\right)^{n}\kappa_{H}\left(\left[\mathfrak{z}\right]_{p^{n}}\right)\label{eq:Gamma formula}
\end{align}

\vphantom{}

II. If $H$ is a contracting semi-basic $p$-Hydra map, then, as $N\rightarrow\infty$,
\emph{(\ref{eq:Gamma formula})} is $\mathcal{F}_{p,q_{H}}$ convergent
to: 
\begin{equation}
\sum_{t\in\hat{\mathbb{Z}}_{p}\backslash\left\{ 0\right\} }\gamma_{H}\left(\frac{t\left|t\right|_{p}}{p}\right)\hat{A}_{H}\left(t\right)e^{2\pi i\left\{ t\mathfrak{z}\right\} _{p}}\overset{\mathcal{F}_{p,q_{H}}}{=}\sum_{n=0}^{\infty}\left(\sum_{j=1}^{p-1}\beta_{H}\left(\frac{j}{p}\right)\varepsilon_{n}^{j}\left(\mathfrak{z}\right)\right)\left(\frac{\mu_{0}}{p}\right)^{n}\kappa_{H}\left(\left[\mathfrak{z}\right]_{p^{n}}\right)\label{eq:F limit of Gamma_H A_H hat Fourier series when alpha is 1}
\end{equation}
whenever $H$ is semi-basic and contracting.
\end{lem}
\begin{rem}
Note that the non-singularity of $H$ ($\alpha_{H}\left(j/p\right)\neq0$
for all $j\in\mathbb{Z}/p\mathbb{Z}$) is \emph{essential} for this
result. 
\end{rem}
Proof:

I. Note that the map $t\in\hat{\mathbb{Z}}_{p}\mapsto\frac{t\left|t\right|_{p}}{p}\in\hat{\mathbb{Z}}_{p}$
takes fractions $k/p^{n}$ and sends them to $\left[k\right]_{p}/p$.
Now, for brevity, let: 
\begin{align}
\gamma_{j} & \overset{\textrm{def}}{=}\gamma_{H}\left(\frac{j}{p}\right)\\
F_{N}\left(\mathfrak{z}\right) & \overset{\textrm{def}}{=}\sum_{0<\left|t\right|_{p}\leq p^{N}}\gamma_{H}\left(\frac{t\left|t\right|_{p}}{p}\right)\hat{A}_{H}\left(t\right)e^{2\pi i\left\{ t\mathfrak{z}\right\} _{p}}
\end{align}
Observing that: 
\begin{equation}
\left\{ t\in\hat{\mathbb{Z}}_{p}:\left|t\right|_{p}=p^{n}\right\} =\left\{ \frac{pk+j}{p^{n}}:k\in\left\{ 0,\ldots,p^{n-1}-1\right\} ,\textrm{ }j\in\left\{ 1,\ldots,p-1\right\} \right\} \label{eq:Decomposition of level sets}
\end{equation}
we can then express $F_{N}$ as a sum involving the $\gamma_{j}$s,
like so: 
\begin{align*}
F_{N}\left(\mathfrak{z}\right) & =\sum_{n=1}^{N}\sum_{\left|t\right|_{p}=p^{n}}\gamma_{H}\left(\frac{t\left|t\right|_{p}}{p}\right)\hat{A}_{H}\left(t\right)e^{2\pi i\left\{ t\mathfrak{z}\right\} _{p}}\\
\left(\textrm{use }(\ref{eq:Decomposition of level sets})\right); & =\sum_{n=1}^{N}\sum_{j=1}^{p-1}\sum_{k=0}^{p^{n-1}-1}\gamma_{H}\left(\frac{j}{p}\right)\hat{A}_{H}\left(\frac{pk+j}{p^{n}}\right)e^{2\pi i\left\{ \frac{pk+j}{p^{n}}\mathfrak{z}\right\} _{p}}
\end{align*}
Using the formal identity: 
\begin{equation}
\sum_{k=0}^{p^{n-1}-1}f\left(\frac{pk+j}{p^{n}}\right)=\sum_{\left|t\right|_{p}\leq p^{n-1}}f\left(t+\frac{j}{p^{n}}\right)
\end{equation}
we can then write:

\begin{align}
F_{N}\left(\mathfrak{z}\right) & =\sum_{n=1}^{N}\sum_{\left|t\right|_{p}\leq p^{n-1}}\sum_{j=1}^{p-1}\gamma_{j}e^{2\pi i\left\{ \frac{j\mathfrak{z}}{p^{n}}\right\} _{p}}\hat{A}_{H}\left(t+\frac{j}{p^{n}}\right)e^{2\pi i\left\{ t\mathfrak{z}\right\} _{p}}\label{eq:Halfway through the gamma computation}
\end{align}

To deal with the $j$-sum, we express $\hat{A}_{H}$ in product form,
changing: 
\begin{equation}
\sum_{j=1}^{p-1}\gamma_{j}e^{2\pi i\left\{ \frac{j\mathfrak{z}}{p^{n}}\right\} _{p}}\hat{A}_{H}\left(t+\frac{j}{p^{n}}\right)e^{2\pi i\left\{ t\mathfrak{z}\right\} _{p}}
\end{equation}
into:

\begin{equation}
\sum_{j=1}^{p-1}\gamma_{j}e^{2\pi i\left\{ \frac{j\mathfrak{z}}{p^{n}}\right\} _{p}}\left(\prod_{m=0}^{n-1}\alpha_{H}\left(p^{m}\left(t+\frac{j}{p^{n}}\right)\right)\right)e^{2\pi i\left\{ t\mathfrak{z}\right\} _{p}}
\end{equation}
Using\textbf{ Proposition \ref{prop:alpha product series expansion}}
to write the $\alpha_{H}$-product out as a series, the above becomes:
\begin{equation}
\sum_{j=1}^{p-1}\gamma_{j}e^{2\pi i\left\{ \frac{j\mathfrak{z}}{p^{n}}\right\} _{p}}\left(\left(\frac{\mu_{0}}{p^{2}}\right)^{n}\sum_{m=0}^{p^{n}-1}\kappa_{H}\left(m\right)e^{-2\pi im\left(t+\frac{j}{p^{n}}\right)}\right)e^{2\pi i\left\{ t\mathfrak{z}\right\} _{p}}
\end{equation}
Hence:

\begin{equation}
\sum_{j=1}^{p-1}\gamma_{j}\left(\frac{\mu_{0}}{p^{2}}\right)^{n}\sum_{m=0}^{p^{n}-1}\kappa_{H}\left(m\right)e^{2\pi i\left\{ \frac{j\left(\mathfrak{z}-m\right)}{p^{n}}\right\} _{p}}e^{2\pi i\left\{ t\left(\mathfrak{z}-m\right)\right\} _{p}}
\end{equation}
Summing over $\left|t\right|_{p}\leq p^{n-1}$, and using: 
\begin{equation}
\sum_{\left|t\right|_{p}\leq p^{n-1}}e^{2\pi i\left\{ t\left(\mathfrak{z}-m\right)\right\} _{p}}=p^{n-1}\left[\mathfrak{z}\overset{p^{n-1}}{\equiv}m\right]
\end{equation}
we obtain: 
\begin{equation}
\sum_{j=1}^{p-1}\gamma_{j}\left(\frac{\mu_{0}}{p^{2}}\right)^{n}\sum_{m=0}^{p^{n}-1}\kappa_{H}\left(m\right)e^{2\pi i\left\{ \frac{j\left(\mathfrak{z}-m\right)}{p^{n}}\right\} _{p}}p^{n-1}\left[\mathfrak{z}\overset{p^{n-1}}{\equiv}m\right]
\end{equation}
In summary, so far, we have; 
\begin{eqnarray}
 & \sum_{\left|t\right|_{p}\leq p^{n-1}}\sum_{j=1}^{p-1}\gamma_{j}e^{2\pi i\left\{ \frac{j\mathfrak{z}}{p^{n}}\right\} _{p}}\hat{A}_{H}\left(t+\frac{j}{p^{n}}\right)e^{2\pi i\left\{ t\mathfrak{z}\right\} _{p}}\nonumber \\
 & =\label{eq:2/3rds of the way through the gamma computation}\\
 & \sum_{j=1}^{p-1}\frac{\gamma_{j}}{p}\left(\frac{\mu_{0}}{p}\right)^{n}\sum_{m=0}^{p^{n}-1}\kappa_{H}\left(m\right)e^{2\pi i\left\{ \frac{j\left(\mathfrak{z}-m\right)}{p^{n}}\right\} _{p}}\left[\mathfrak{z}\overset{p^{n-1}}{\equiv}m\right]\nonumber 
\end{eqnarray}

Next, using the formal identity: 
\begin{equation}
\sum_{m=0}^{p^{n}-1}f\left(m\right)=\sum_{k=0}^{p-1}\sum_{m=0}^{p^{n-1}-1}f\left(m+kp^{n-1}\right)\label{eq:rho to the n formal identity}
\end{equation}
and the functional equation identity for $\kappa_{H}$ (equation (\ref{eq:Kappa_H functional equations})
from \textbf{Proposition \ref{prop:Properties of Kappa_H}}), we have:
\begin{eqnarray*}
 & \sum_{m=0}^{p^{n}-1}\kappa_{H}\left(m\right)e^{2\pi i\left\{ \frac{j\left(\mathfrak{z}-m\right)}{p^{n}}\right\} _{p}}\left[\mathfrak{z}\overset{p^{n-1}}{\equiv}m\right]\\
 & =\\
 & \sum_{k=0}^{p-1}\sum_{m=0}^{p^{n-1}-1}\frac{\mu_{k}}{\mu_{0}}\kappa_{H}\left(m\right)e^{2\pi i\left\{ \frac{j\left(\mathfrak{z}-m-kp^{n-1}\right)}{p^{n}}\right\} _{p}}\left[\mathfrak{z}\overset{p^{n-1}}{\equiv}m\right]
\end{eqnarray*}
Here, note that $\left[\mathfrak{z}\right]_{p^{n-1}}$ is the unique
integer $m\in\left\{ 0,\ldots,p^{n-1}-1\right\} $ satisfying $\mathfrak{z}\overset{p^{n-1}}{\equiv}m$.
This leaves us with: 
\begin{align*}
\sum_{m=0}^{p^{n}-1}\kappa_{H}\left(m\right)e^{2\pi i\left\{ \frac{j\left(\mathfrak{z}-m\right)}{p^{n}}\right\} _{p}}\left[\mathfrak{z}\overset{p^{n-1}}{\equiv}m\right] & =\kappa_{H}\left(\left[\mathfrak{z}\right]_{p^{n-1}}\right)\left(\varepsilon_{n-1}\left(\mathfrak{z}\right)\right)^{j}\overbrace{\sum_{k=0}^{p-1}\frac{\mu_{k}}{\mu_{0}}e^{-2\pi i\frac{jk}{p}}}^{\frac{p^{2}}{\mu_{0}}\alpha_{H}\left(\frac{j}{p}\right)}\\
 & =\frac{p^{2}}{\mu_{0}}\alpha_{H}\left(\frac{j}{p}\right)\kappa_{H}\left(\left[\mathfrak{z}\right]_{p^{n-1}}\right)\left(\varepsilon_{n-1}\left(\mathfrak{z}\right)\right)^{j}
\end{align*}

Returning with this to (\ref{eq:2/3rds of the way through the gamma computation}),
the equation 
\begin{equation}
\sum_{\left|t\right|_{p}\leq p^{n-1}}\sum_{j=1}^{p-1}\gamma_{j}e^{2\pi i\left\{ \frac{j\mathfrak{z}}{p^{n}}\right\} _{p}}\hat{A}_{H}\left(t+\frac{j}{p^{n}}\right)e^{2\pi i\left\{ t\mathfrak{z}\right\} _{p}}
\end{equation}
transforms into: 
\begin{equation}
\sum_{j=1}^{p-1}\frac{\gamma_{j}}{p}\left(\frac{\mu_{0}}{p}\right)^{n}\frac{p^{2}}{\mu_{0}}\alpha_{H}\left(\frac{j}{p}\right)\left(\varepsilon_{n-1}\left(\mathfrak{z}\right)\right)^{j}\kappa_{H}\left(\left[\mathfrak{z}\right]_{p^{n-1}}\right)
\end{equation}
Almost finished, recall that: 
\begin{equation}
\gamma_{j}=\gamma_{H}\left(\frac{j}{p}\right)=\frac{\beta_{H}\left(\frac{j}{p}\right)}{\alpha_{H}\left(\frac{j}{p}\right)}
\end{equation}
With this, we can write the previous line as: 
\begin{equation}
\left(\sum_{j=1}^{p-1}\beta_{H}\left(\frac{j}{p}\right)\left(\varepsilon_{n-1}\left(\mathfrak{z}\right)\right)^{j}\right)\left(\frac{\mu_{0}}{p}\right)^{n-1}\kappa_{H}\left(\left[\mathfrak{z}\right]_{p^{n-1}}\right)
\end{equation}
Applying this to the right-hand side of (\ref{eq:Halfway through the gamma computation})
gives us: 
\begin{align}
F_{N}\left(\mathfrak{z}\right) & =\sum_{n=1}^{N}\left(\sum_{j=1}^{p-1}\beta_{H}\left(\frac{j}{p}\right)\left(\varepsilon_{n-1}\left(\mathfrak{z}\right)\right)^{j}\right)\left(\frac{\mu_{0}}{p}\right)^{n-1}\kappa_{H}\left(\left[\mathfrak{z}\right]_{p^{n-1}}\right)
\end{align}
Re-indexing $n$ by a shift of $1$ produces equation (\ref{eq:Gamma formula}).

\vphantom{}

II. For each $\mathfrak{z}\in\mathbb{Z}_{p}$, the algebraic number
$\sum_{j=1}^{p-1}\beta_{H}\left(\frac{j}{p}\right)\left(\varepsilon_{n}\left(\mathfrak{z}\right)\right)^{j}$
is uniformly bounded with respect to $n\in\mathbb{N}_{0}$ in both
$\mathbb{C}$ and $\mathbb{C}_{q}$. Like in the $p=2$ case, since
$H$ is semi-basic, the sequence $\left(\frac{\mu_{0}}{p}\right)^{n}\kappa_{H}\left(\left[\mathfrak{z}\right]_{p^{n}}\right)$
is then $\mathcal{F}_{p,q_{H}}$-convergent to $0$. The $q$-adic
convergence to $0$ over $\mathbb{Z}_{p}^{\prime}$ guarantees the
$q$-adic convergence of (\ref{eq:Gamma formula}) as $N\rightarrow\infty$
for $\mathfrak{z}\in\mathbb{Z}_{p}^{\prime}$. On the other hand,
for $\mathfrak{z}\in\mathbb{N}_{0}$, as we have seen: 
\begin{equation}
\left(\frac{\mu_{0}}{p}\right)^{n}\kappa_{H}\left(\left[\mathfrak{z}\right]_{p^{n}}\right)=\left(\frac{\mu_{0}}{p}\right)^{n}\kappa_{H}\left(\mathfrak{z}\right),\textrm{ }\forall n\geq\lambda_{p}\left(\mathfrak{z}\right)
\end{equation}
and the fact that $H$ is contracting then guarantees (by the ratio
test, no less) that (\ref{eq:Gamma formula}) is convergent in $\mathbb{C}$
for $\mathfrak{z}\in\mathbb{N}_{0}$. This proves (II).

Q.E.D.

\vphantom{}

With these formulae, we can sum the Fourier series generated by (\ref{eq:Fine structure of Chi_H,N hat when alpha is 1})
to obtain a non-trivial formula for $\chi_{H,N}$. 
\begin{thm}
\label{thm:F-series for Nth truncation of Chi_H, alpha is 1}Suppose\index{chi{H}@$\chi_{H}$!$N$th truncation}
$\alpha_{H}\left(0\right)=1$. Then, for all $N\geq1$ and all $\mathfrak{z}\in\mathbb{Z}_{p}$:
\begin{align}
\chi_{H,N}\left(\mathfrak{z}\right) & \overset{\overline{\mathbb{Q}}}{=}\sum_{n=0}^{N-1}\left(\sum_{j=0}^{p-1}\beta_{H}\left(\frac{j}{p}\right)\left(\varepsilon_{n}\left(\mathfrak{z}\right)\right)^{j}\right)\left(\frac{\mu_{0}}{p}\right)^{n}\kappa_{H}\left(\left[\mathfrak{z}\right]_{p^{n}}\right)\label{eq:Chi_H,N when alpha is 1 and rho is arbitrary}
\end{align}
In particular, when $p=2$: 
\begin{equation}
\chi_{H,N}\left(\mathfrak{z}\right)\overset{\overline{\mathbb{Q}}}{=}-\gamma_{H}\left(\frac{1}{2}\right)+\gamma_{H}\left(\frac{1}{2}\right)\left(\frac{\mu_{0}}{2}\right)^{N}\kappa_{H}\left(\left[\mathfrak{z}\right]_{2^{N}}\right)+\beta_{H}\left(0\right)\sum_{n=0}^{N-1}\left(\frac{\mu_{0}}{2}\right)^{n}\kappa_{H}\left(\left[\mathfrak{z}\right]_{2^{n}}\right)\label{eq:Chi_H,N when alpha is 1 and rho is 2}
\end{equation}
\end{thm}
Proof: We start by multiplying (\ref{eq:Fine structure of Chi_H,N hat when alpha is 1})
by $e^{2\pi i\left\{ t\mathfrak{z}\right\} _{p}}$ and then summing
over all $\left|t\right|_{p}\leq p^{N}$. The left-hand side of (\ref{eq:Fine structure of Chi_H,N hat when alpha is 1})
becomes: 
\begin{equation}
\chi_{H,N}\left(\mathfrak{z}\right)-\beta_{H}\left(0\right)N\left(D_{p:N-1}*dA_{H}\right)\left(\mathfrak{z}\right)
\end{equation}
whereas the right-hand side is: 
\begin{align*}
\sum_{0<\left|t\right|_{p}\leq p^{N-1}}\left(\gamma_{H}\left(\frac{t\left|t\right|_{p}}{p}\right)+\beta_{H}\left(0\right)v_{p}\left(t\right)\right)\hat{A}_{H}\left(t\right)e^{2\pi i\left\{ t\mathfrak{z}\right\} _{p}}\\
+\sum_{\left|t\right|_{p}=p^{N}}\gamma_{H}\left(\frac{t\left|t\right|_{p}}{p}\right)\hat{A}_{H}\left(t\right)e^{2\pi i\left\{ t\mathfrak{z}\right\} _{p}}
\end{align*}
Simplifying produces: 
\begin{align}
\chi_{H,N}\left(\mathfrak{z}\right) & \overset{\overline{\mathbb{Q}}}{=}\beta_{H}\left(0\right)N\left(D_{p:N-1}*dA_{H}\right)\left(\mathfrak{z}\right)\label{eq:Chi_H,N rho not equal to 2, ready to simplify}\\
 & +\beta_{H}\left(0\right)\sum_{0<\left|t\right|_{p}\leq p^{N-1}}v_{p}\left(t\right)\hat{A}_{H}\left(t\right)e^{2\pi i\left\{ t\mathfrak{z}\right\} _{p}}\nonumber \\
 & +\sum_{0<\left|t\right|_{p}\leq p^{N}}\gamma_{H}\left(\frac{t\left|t\right|_{p}}{p}\right)\hat{A}_{H}\left(t\right)e^{2\pi i\left\{ t\mathfrak{z}\right\} _{p}}\nonumber 
\end{align}
Now we call upon our legion of formulae: (\ref{eq:Convolution of dA_H and D_N}),
\textbf{Lemma \ref{lem:v_p A_H hat summation formulae}}, and \textbf{Lemma
\ref{lem:1D gamma formula}}. Using them (whilst remembering that
$\alpha_{H}\left(0\right)=1$) makes (\ref{eq:Chi_H,N rho not equal to 2, ready to simplify})
into:

\begin{align*}
\chi_{H,N}\left(\mathfrak{z}\right) & \overset{\overline{\mathbb{Q}}}{=}\beta_{H}\left(0\right)N\left(\frac{\mu_{0}}{p}\right)^{N-1}\kappa_{H}\left(\left[\mathfrak{z}\right]_{p^{N-1}}\right)\\
 & -\beta_{H}\left(0\right)\left(N-1\right)\left(\frac{\mu_{0}}{p}\right)^{N-1}\kappa_{H}\left(\left[\mathfrak{z}\right]_{p^{N-1}}\right)\\
 & +\beta_{H}\left(0\right)\sum_{n=0}^{N-2}\left(\frac{\mu_{0}}{p}\right)^{n}\kappa_{H}\left(\left[\mathfrak{z}\right]_{p^{n}}\right)\\
 & \sum_{n=0}^{N-1}\left(\sum_{j=1}^{p-1}\beta_{H}\left(\frac{j}{p}\right)\left(\varepsilon_{n}\left(\mathfrak{z}\right)\right)^{j}\right)\left(\frac{\mu_{0}}{p}\right)^{n}\kappa_{H}\left(\left[\mathfrak{z}\right]_{p^{n}}\right)\\
 & \overset{\overline{\mathbb{Q}}}{=}\beta_{H}\left(0\right)\sum_{n=0}^{N-1}\left(\frac{\mu_{0}}{p}\right)^{n}\kappa_{H}\left(\left[\mathfrak{z}\right]_{p^{n}}\right)\\
 & +\sum_{n=0}^{N-1}\left(\sum_{j=1}^{p-1}\beta_{H}\left(\frac{j}{p}\right)\left(\varepsilon_{n}\left(\mathfrak{z}\right)\right)^{j}\right)\left(\frac{\mu_{0}}{p}\right)^{n}\kappa_{H}\left(\left[\mathfrak{z}\right]_{p^{n}}\right)
\end{align*}
The bottom two lines combine to form a single sum by noting that the
upper most of the two is the $j=0$ case of the bottom-most of the
two.

Finally, when $p=2$, rather than simplify (\ref{eq:Chi_H,N when alpha is 1 and rho is arbitrary}),
it will actually be easier to compute it from scratch all over again.
Multiplying (\ref{eq:Fine structure of Chi_H,N hat when alpha is 1})
by $e^{2\pi i\left\{ t\mathfrak{z}\right\} _{2}}$ and summing over
all $\left|t\right|_{2}\leq2^{N}$ gives: 
\begin{align*}
\chi_{H,N}\left(\mathfrak{z}\right)-\beta_{H}\left(0\right)N\left(D_{2:N-1}*dA_{H}\right)\left(\mathfrak{z}\right) & \overset{\overline{\mathbb{Q}}}{=}\gamma_{H}\left(\frac{1}{2}\right)\sum_{0<\left|t\right|_{2}\leq2^{N-1}}\hat{A}_{H}\left(t\right)e^{2\pi i\left\{ t\mathfrak{z}\right\} _{2}}\\
 & +\beta_{H}\left(0\right)\sum_{0<\left|t\right|_{2}\leq2^{N-1}}v_{2}\left(t\right)\hat{A}_{H}\left(t\right)e^{2\pi i\left\{ t\mathfrak{z}\right\} _{2}}\\
 & +\gamma_{H}\left(\frac{1}{2}\right)\sum_{\left|t\right|_{2}=2^{N}}\hat{A}_{H}\left(t\right)e^{2\pi i\left\{ t\mathfrak{z}\right\} _{2}}
\end{align*}
This simplifies to: 
\begin{align*}
\chi_{H,N}\left(\mathfrak{z}\right)-\beta_{H}\left(0\right)N\left(D_{2:N-1}*dA_{H}\right)\left(\mathfrak{z}\right) & \overset{\overline{\mathbb{Q}}}{=}\gamma_{H}\left(\frac{1}{2}\right)\sum_{0<\left|t\right|_{2}\leq2^{N}}\hat{A}_{H}\left(t\right)e^{2\pi i\left\{ t\mathfrak{z}\right\} _{2}}\\
 & +\beta_{H}\left(0\right)\sum_{0<\left|t\right|_{2}\leq2^{N-1}}v_{2}\left(t\right)\hat{A}_{H}\left(t\right)e^{2\pi i\left\{ t\mathfrak{z}\right\} _{2}}\\
 & \overset{\overline{\mathbb{Q}}}{=}-\gamma_{H}\left(\frac{1}{2}\right)\hat{A}_{H}\left(0\right)\\
 & +\gamma_{H}\left(\frac{1}{2}\right)\underbrace{\sum_{\left|t\right|_{2}\leq2^{N}}\hat{A}_{H}\left(t\right)e^{2\pi i\left\{ t\mathfrak{z}\right\} _{2}}}_{\left(D_{2:N}*dA_{H}\right)\left(\mathfrak{z}\right)}\\
 & +\beta_{H}\left(0\right)\sum_{0<\left|t\right|_{2}\leq2^{N-1}}v_{2}\left(t\right)\hat{A}_{H}\left(t\right)e^{2\pi i\left\{ t\mathfrak{z}\right\} _{2}}
\end{align*}
Applying (\ref{eq:Convolution of dA_H and D_N}) and \textbf{Lemma
\ref{lem:v_p A_H hat summation formulae}}, the above becomes: 
\begin{align*}
\chi_{H,N}\left(\mathfrak{z}\right)-\beta_{H}\left(0\right)N\left(\frac{\mu_{0}}{2}\right)^{N-1}\kappa_{H}\left(\left[\mathfrak{z}\right]_{2^{N-1}}\right) & \overset{\overline{\mathbb{Q}}}{=}-\gamma_{H}\left(\frac{1}{2}\right)\\
 & +\gamma_{H}\left(\frac{1}{2}\right)\left(\frac{\mu_{0}}{2}\right)^{N}\kappa_{H}\left(\left[\mathfrak{z}\right]_{2^{N}}\right)\\
 & -\beta_{H}\left(0\right)\left(N-1\right)\left(\frac{\mu_{0}}{2}\right)^{N-1}\kappa_{H}\left(\left[\mathfrak{z}\right]_{2^{N-1}}\right)\\
 & +\beta_{H}\left(0\right)\sum_{n=0}^{N-2}\left(\frac{\mu_{0}}{2}\right)^{n}\kappa_{H}\left(\left[\mathfrak{z}\right]_{2^{n}}\right)
\end{align*}
Simplifying yields: 
\begin{equation}
\chi_{H,N}\left(\mathfrak{z}\right)\overset{\overline{\mathbb{Q}}}{=}-\gamma_{H}\left(\frac{1}{2}\right)+\gamma_{H}\left(\frac{1}{2}\right)\left(\frac{\mu_{0}}{2}\right)^{N}\kappa_{H}\left(\left[\mathfrak{z}\right]_{2^{N}}\right)+\beta_{H}\left(0\right)\sum_{n=0}^{N-1}\left(\frac{\mu_{0}}{2}\right)^{n}\kappa_{H}\left(\left[\mathfrak{z}\right]_{2^{n}}\right)
\end{equation}

Q.E.D. 
\begin{cor}[\textbf{$\mathcal{F}$-series for $\chi_{H}$ when $\alpha_{H}\left(0\right)=1$}]
\label{cor:Chi_H, F-convergent formula}Suppose\index{chi{H}@$\chi_{H}$!mathcal{F}-series@$\mathcal{F}$-series}
For the given $p$-Hydra map $H$ (semi-basic, contracting, non-singular,
fixes $0$), suppose $\alpha_{H}\left(0\right)=1$. Then: 
\begin{equation}
\chi_{H}\left(\mathfrak{z}\right)\overset{\mathcal{F}_{p,q_{H}}}{=}\sum_{n=0}^{\infty}\left(\sum_{j=0}^{p-1}\beta_{H}\left(\frac{j}{p}\right)\left(\varepsilon_{n}\left(\mathfrak{z}\right)\right)^{j}\right)\left(\frac{\mu_{0}}{p}\right)^{n}\kappa_{H}\left(\left[\mathfrak{z}\right]_{p^{n}}\right),\textrm{ }\forall\mathfrak{z}\in\mathbb{Z}_{p}\label{eq:Explicit Formula for Chi_H when alpha is 1 and rho is arbitrary}
\end{equation}
with the special case: 
\begin{equation}
\chi_{H}\left(\mathfrak{z}\right)\overset{\mathcal{F}_{2,q_{H}}}{=}-\gamma_{H}\left(\frac{1}{2}\right)+\beta_{H}\left(0\right)\sum_{n=0}^{\infty}\left(\frac{\mu_{0}}{2}\right)^{n}\kappa_{H}\left(\left[\mathfrak{z}\right]_{2^{n}}\right),\textrm{ }\forall\mathfrak{z}\in\mathbb{Z}_{2}\label{eq:Explicit Formula for Chi_H when alpha is 1 and rho is 2}
\end{equation}
for when $p=2$.

In both cases, as indicated, the topology of convergence of the series
on the right-hand side is that of the standard $\left(p,q_{H}\right)$-adic
frame: the topology of $\mathbb{C}$ when $\mathfrak{z}\in\mathbb{N}_{0}$,
and the topology of $\mathbb{C}_{q_{H}}$ when $\mathfrak{z}\in\mathbb{Z}_{p}^{\prime}$. 
\end{cor}
Proof: The given properties of $H$ tell us that $\left(\frac{\mu_{0}}{p}\right)^{N}\kappa_{H}\left(\left[\mathfrak{z}\right]_{p^{N}}\right)$
has $0$ as a $\mathcal{F}_{p,q_{H}}$-limit. Appealing to the limits
formulae in \textbf{Lemmata \ref{lem:1D gamma formula}} and \textbf{\ref{lem:v_p A_H hat summation formulae}}
then allows us to take the $\mathcal{F}_{p,q_{H}}$-limits of (\ref{eq:Chi_H,N when alpha is 1 and rho is arbitrary})
and (\ref{eq:Chi_H,N when alpha is 1 and rho is 2}) as $N\rightarrow\infty$,
which produces (\ref{eq:Explicit Formula for Chi_H when alpha is 1 and rho is arbitrary})
and (\ref{eq:Explicit Formula for Chi_H when alpha is 1 and rho is 2}),
respectively.

Q.E.D.

\vphantom{}

Next, we sum in $\mathbb{C}$ for $\mathfrak{z}\in\mathbb{N}_{0}$. 
\begin{cor}
\label{cor:Chi_H explicit formula on N_0}If $\alpha_{H}\left(0\right)=1$,
then: 
\begin{equation}
\chi_{H}\left(n\right)\overset{\mathbb{C}}{=}\sum_{k=0}^{\lambda_{p}\left(n\right)-1}\left(\sum_{j=0}^{p-1}\beta_{H}\left(\frac{j}{p}\right)\left(\varepsilon_{k}\left(n\right)\right)^{j}\right)\left(\frac{\mu_{0}}{p}\right)^{k}\kappa_{H}\left(\left[n\right]_{p^{k}}\right),\textrm{ }\forall n\in\mathbb{N}_{0}\label{eq:archimedean Chi_H when rho is arbitrary and alpha_H of 0 is 1}
\end{equation}
with the special case: 
\begin{equation}
\chi_{H}\left(n\right)\overset{\mathbb{C}}{=}-\gamma_{H}\left(\frac{1}{2}\right)+\frac{2\beta_{H}\left(0\right)}{2-\mu_{0}}M_{H}\left(n\right)+\beta_{H}\left(0\right)\sum_{k=0}^{\lambda_{2}\left(n\right)-1}\left(\frac{\mu_{0}}{2}\right)^{k}\kappa_{H}\left(\left[n\right]_{2^{k}}\right),\textrm{ }\forall n\in\mathbb{N}_{0}\label{eq:archimedean Chi_H when rho is 2 and alpha_H of 0 is 1}
\end{equation}
when $p=2$. Regardless of the value of $p$, the $k$-sums are defined
to be $0$ when $n=0$. 
\end{cor}
Proof: Let $n\in\mathbb{N}_{0}$. Since $\kappa_{H}\left(\left[n\right]_{p^{k}}\right)=\kappa_{H}\left(n\right)$
and $\varepsilon_{k}\left(n\right)=1$ for all $k\geq\lambda_{p}\left(n\right)$,
the equation (\ref{eq:Explicit Formula for Chi_H when alpha is 1 and rho is arbitrary})
becomes: 
\begin{align*}
\chi_{H}\left(n\right) & \overset{\mathbb{C}}{=}\sum_{k=0}^{\lambda_{p}\left(n\right)-1}\left(\sum_{j=0}^{p-1}\beta_{H}\left(\frac{j}{p}\right)\left(\varepsilon_{k}\left(n\right)\right)^{j}\right)\left(\frac{\mu_{0}}{p}\right)^{k}\kappa_{H}\left(\left[n\right]_{p^{k}}\right)\\
 & +\left(\sum_{k=\lambda_{p}\left(n\right)}^{\infty}\left(\sum_{j=0}^{p-1}\beta_{H}\left(\frac{j}{p}\right)\cdot1\right)\left(\frac{\mu_{0}}{p}\right)^{k}\right)\kappa_{H}\left(n\right)
\end{align*}
Since $H$ is contracting, the geometric series: 
\begin{equation}
\sum_{k=\lambda_{p}\left(n\right)}^{\infty}\left(\frac{\mu_{0}}{p}\right)^{k}
\end{equation}
converges to a limit in $\mathbb{C}$. However, this does not matter
much, because: 
\[
\sum_{j=0}^{p-1}\beta_{H}\left(\frac{j}{p}\right)=\sum_{j=0}^{p-1}\frac{1}{p}\sum_{\ell=0}^{p-1}\frac{b_{\ell}}{d_{\ell}}e^{-2\pi i\ell\frac{j}{p}}=\sum_{\ell=0}^{p-1}\frac{b_{\ell}}{d_{\ell}}\sum_{j=0}^{p-1}\frac{1}{p}e^{-2\pi i\ell\frac{j}{p}}=\frac{b_{\ell}}{d_{\ell}}\left[\ell\overset{p}{\equiv}0\right]=\frac{b_{0}}{d_{0}}
\]
Since $b_{0}=H\left(0\right)=0$, we are then left with: 
\begin{align*}
\chi_{H}\left(n\right) & \overset{\mathbb{C}}{=}\sum_{k=0}^{\lambda_{p}\left(n\right)-1}\left(\sum_{j=0}^{p-1}\beta_{H}\left(\frac{j}{p}\right)\left(\varepsilon_{k}\left(n\right)\right)^{j}\right)\left(\frac{\mu_{0}}{p}\right)^{k}\kappa_{H}\left(\left[n\right]_{p^{k}}\right)
\end{align*}
which is the desired formula.

As for the $p=2$ case, applying the same argument given above to
(\ref{eq:Explicit Formula for Chi_H when alpha is 1 and rho is 2})
yields: 
\begin{align*}
\chi_{H}\left(n\right) & \overset{\mathbb{C}}{=}-\gamma_{H}\left(\frac{1}{2}\right)+\beta_{H}\left(0\right)\sum_{k=0}^{\lambda_{2}\left(n\right)-1}\left(\frac{\mu_{0}}{2}\right)^{n}\kappa_{H}\left(\left[n\right]_{2^{k}}\right)\\
 & +\beta_{H}\left(0\right)\kappa_{H}\left(n\right)\sum_{k=\lambda_{2}\left(n\right)}^{\infty}\left(\frac{\mu_{0}}{2}\right)^{k}\\
\left(H\textrm{ is contracting}\right); & \overset{\mathbb{C}}{=}-\gamma_{H}\left(\frac{1}{2}\right)+\beta_{H}\left(0\right)\sum_{k=0}^{\lambda_{2}\left(n\right)-1}\left(\frac{\mu_{0}}{2}\right)^{k}\kappa_{H}\left(\left[n\right]_{2^{k}}\right)\\
 & +\beta_{H}\left(0\right)\frac{\left(\frac{\mu_{0}}{2}\right)^{\lambda_{2}\left(n\right)}\kappa_{H}\left(n\right)}{1-\frac{\mu_{0}}{2}}
\end{align*}
Finally, since: 
\[
\kappa_{H}\left(n\right)=\left(\frac{2}{\mu_{0}}\right)^{\lambda_{2}\left(n\right)}M_{H}\left(n\right)
\]
we then obtain:

\begin{align*}
\chi_{H}\left(n\right) & \overset{\mathbb{C}}{=}-\gamma_{H}\left(\frac{1}{2}\right)+\beta_{H}\left(0\right)\frac{M_{H}\left(n\right)}{1-\frac{\mu_{0}}{2}}+\beta_{H}\left(0\right)\sum_{k=0}^{\lambda_{2}\left(n\right)-1}\left(\frac{\mu_{0}}{2}\right)^{k}\kappa_{H}\left(\left[n\right]_{2^{k}}\right)
\end{align*}

Q.E.D.

\vphantom{}

Taken together, these two corollaries then establish the quasi-integrability
of $\chi_{H}$\index{chi{H}@$\chi_{H}$!quasi-integrability} for $p\geq2$
and $\alpha_{H}\left(0\right)=1$. 
\begin{cor}[\textbf{Quasi-Integrability of $\chi_{H}$ when $\alpha_{H}\left(0\right)=1$}]
\label{cor:Quasi-integrability of Chi_H for alpha equals 1}If $\alpha_{H}\left(0\right)=1$,
then, $\chi_{H}$ is quasi-integrable with respect to the standard
$\left(p,q_{H}\right)$-adic frame.

In\index{chi{H}@$\chi_{H}$!Fourier transform} particular, when $p=2$,
the function $\hat{\chi}_{H}:\hat{\mathbb{Z}}_{2}\rightarrow\overline{\mathbb{Q}}$
defined by: 
\begin{equation}
\hat{\chi}_{H}\left(t\right)\overset{\textrm{def}}{=}\begin{cases}
-\gamma_{H}\left(\frac{1}{2}\right) & \textrm{if }t=0\\
\beta_{H}\left(0\right)v_{2}\left(t\right)\hat{A}_{H}\left(t\right) & \textrm{else }
\end{cases},\textrm{ }\forall t\in\hat{\mathbb{Z}}_{2}\label{eq:Formula for Chi_H hat when rho is 2 and alpha is 1}
\end{equation}
is then a Fourier transform of $\chi_{H}$. In this case, the function
defined by \emph{(\ref{eq:Point-wise limit of the left-hand side of the fine structure formula when alpha is 1})}
is also a Fourier transform of $\chi_{H}$, differing from the $\hat{\chi}_{H}$
given above by $\gamma_{H}\left(\frac{1}{2}\right)\hat{A}_{H}\left(t\right)$,
which, by \textbf{\emph{Theorem \ref{thm:Properties of dA_H}}}, is
a degenerate measure, seeing as $\alpha_{H}\left(0\right)=1$.

For $p\geq3$, we can obtain a Fourier transform for $\chi_{H}$ by
defining a function $\hat{\chi}_{H}:\hat{\mathbb{Z}}_{p}\rightarrow\overline{\mathbb{Q}}$
by: 
\begin{equation}
\hat{\chi}_{H}\left(t\right)\overset{\textrm{def}}{=}\begin{cases}
0 & \textrm{if }t=0\\
\left(\gamma_{H}\left(\frac{t\left|t\right|_{p}}{p}\right)+\beta_{H}\left(0\right)v_{p}\left(t\right)\right)\hat{A}_{H}\left(t\right) & \textrm{else}
\end{cases}\label{eq:Chi_H hat when rho is not 2 and when alpha is 1}
\end{equation}
\end{cor}
Proof: Corollaries\textbf{ \ref{cor:Chi_H explicit formula on N_0}}
and\textbf{ \ref{cor:Chi_H, F-convergent formula}} show that the
$N$th partial sums of the Fourier series generated by (\ref{eq:Point-wise limit of the left-hand side of the fine structure formula when alpha is 1})
are $\mathcal{F}_{p,q_{H}}$-convergent to (\ref{eq:Formula for Chi_H hat when rho is 2 and alpha is 1})
and (\ref{eq:Chi_H hat when rho is not 2 and when alpha is 1}) for
$p=2$ and $p\geq3$, respectively, thereby establishing the quasi-integrability
of $\chi_{H}$ with respect to the standard $\left(p,q_{H}\right)$-adic
frame.

Finally, letting $\hat{\chi}_{H}^{\prime}\left(t\right)$ denote (\ref{eq:Point-wise limit of the left-hand side of the fine structure formula when alpha is 1}),
observe that when $\alpha_{H}\left(0\right)=1$ and $p=2$: 
\begin{equation}
\hat{\chi}_{H}^{\prime}\left(t\right)\overset{\overline{\mathbb{Q}}}{=}\begin{cases}
0 & \textrm{if }t=0\\
\left(\gamma_{H}\left(\frac{1}{2}\right)+\beta_{H}\left(0\right)v_{2}\left(t\right)\right)\hat{A}_{H}\left(t\right) & \textrm{if }\left|t\right|_{2}>0
\end{cases}
\end{equation}
Since $\hat{A}_{H}\left(0\right)=1$, we have that: 
\begin{equation}
\hat{\chi}_{H}^{\prime}\left(t\right)-\gamma_{H}\left(\frac{1}{2}\right)\hat{A}_{H}\left(t\right)\overset{\overline{\mathbb{Q}}}{=}\begin{cases}
-\gamma_{H}\left(\frac{1}{2}\right) & \textrm{if }t=0\\
\beta_{H}\left(0\right)v_{2}\left(t\right)\hat{A}_{H}\left(t\right) & \textrm{if }\left|t\right|_{2}>0
\end{cases}=\hat{\chi}_{H}\left(t\right)
\end{equation}
which shows that $\hat{\chi}_{H}^{\prime}\left(t\right)$ and $\hat{\chi}_{H}\left(t\right)$
differ by a factor of $\gamma_{H}\left(\frac{1}{2}\right)\hat{A}_{H}\left(t\right)$,
which is a degenerate measure because $\alpha_{H}\left(0\right)=1$.

Q.E.D.

\vphantom{}

This completes the $\alpha_{H}\left(0\right)=1$ case. By using $\chi_{H}$'s
functional equations (\ref{eq:Functional Equations for Chi_H over the rho-adics}),
we can extend all of the above work to cover all $p$-Hydra maps,
regardless of the value of $\alpha_{H}\left(0\right)$. The key to
this are \textbf{Corollary \ref{cor:Chi_H explicit formula on N_0}
}and \textbf{Lemma \ref{lem:Functional equations and truncation}}.
First, however, a definition: 
\begin{defn}[\textbf{Little Psi-$H$ \& Big Psi-$H$}]
\ 

\vphantom{}

I. We define \nomenclature{$\psi_{H}\left(m\right)$}{ }$\psi_{H}:\mathbb{N}_{0}\rightarrow\mathbb{Q}$
(``Little Psi-$H$'') by: 
\begin{equation}
\psi_{H}\left(m\right)\overset{\textrm{def}}{=}\frac{M_{H}\left(m\right)}{1-\frac{\mu_{0}}{p}}+\sum_{n=0}^{\lambda_{p}\left(m\right)-1}\left(\frac{\mu_{0}}{p}\right)^{n}\kappa_{H}\left(\left[m\right]_{p^{n}}\right),\textrm{ }\forall m\in\mathbb{N}_{0}\label{eq:Definition of Little Psi_H}
\end{equation}
where the $n$-sum is defined to be $0$ when $m=0$.

\vphantom{}

II. We define \nomenclature{$\Psi_{H}\left(m\right)$}{ }$\Psi_{H}:\mathbb{N}_{0}\rightarrow\mathbb{Q}$
(``Big Psi-$H$'') by: 
\begin{equation}
\Psi_{H}\left(m\right)\overset{\textrm{def}}{=}-\beta_{H}\left(0\right)\frac{M_{H}\left(m\right)}{1-\frac{\mu_{0}}{p}}+\sum_{n=0}^{\lambda_{p}\left(m\right)-1}\left(\sum_{j=1}^{p-1}\beta_{H}\left(\frac{j}{p}\right)\varepsilon_{n}^{j}\left(m\right)\right)\left(\frac{\mu_{0}}{p}\right)^{n}\kappa_{H}\left(\left[m\right]_{p^{n}}\right)\label{eq:Definition of Big Psi_H}
\end{equation}
\end{defn}
\vphantom{}

As a simple computation will immediately verify, (\ref{eq:Definition of Little Psi_H})
and (\ref{eq:Definition of Big Psi_H}) are merely the sums in $\mathbb{C}$
of (\ref{eq:Limit of Fourier sum of v_p A_H hat}) and (\ref{eq:F limit of Gamma_H A_H hat Fourier series when alpha is 1})
from \textbf{Lemmata \ref{lem:v_p A_H hat summation formulae}} and
\textbf{\ref{lem:1D gamma formula}}, respectively, for $\mathfrak{z}\in\mathbb{N}_{0}$.
To derive Fourier transforms and quasi-integrability for $\chi_{H}$
when $\alpha_{H}\left(0\right)\neq1$, we will first show that $\psi_{H}$
and $\Psi_{H}$ are rising-continuable to standard-frame-quasi-integrable
$\left(p,q_{H}\right)$-adic functions, and that these continuations
are characterized by systems of functional equations very to near
the ones satisfied by $\chi_{H}$ (\ref{eq:Functional Equations for Chi_H over the rho-adics}).
With a bit of linear algebra, we can then \emph{solve }for the appropriate
linear combinations of $\psi_{H}$ and $\Psi_{H}$ to needed to obtain
$\chi_{H}$, regardless of the value of $\alpha_{H}\left(0\right)$. 
\begin{lem}[\textbf{Rising-Continuability and Functional Equations for $\psi_{H}$
\& $\Psi_{H}$}]
\label{lem:Rising continuation and uniquenss of the psi_Hs}For $H$
as given at the start of \emph{Section \pageref{sec:4.2 Fourier-Transforms-=00003D000026}}:

\vphantom{}

I. $\psi_{H}$ is rising-continuable to a $\left(p,q_{H}\right)$-adic
function $\psi_{H}:\mathbb{Z}_{p}\rightarrow\mathbb{Z}_{q_{H}}$ given
by: 
\begin{equation}
\psi_{H}\left(\mathfrak{z}\right)\overset{\mathcal{F}_{p,q_{H}}}{=}\sum_{n=0}^{\infty}\left(\frac{\mu_{0}}{p}\right)^{n}\kappa_{H}\left(\left[\mathfrak{z}\right]_{p^{n}}\right),\textrm{ }\forall\mathfrak{z}\in\mathbb{Z}_{p}\label{eq:Rising-continuation of Little Psi_H}
\end{equation}
Moreover, $\psi_{H}\left(\mathfrak{z}\right)$ is the unique rising-continuous
$\left(p,q_{H}\right)$-adic function satisfying the system of \index{functional equation!psi_{H}@$\psi_{H}$}functional
equations: 
\begin{equation}
\psi_{H}\left(p\mathfrak{z}+j\right)=\frac{\mu_{j}}{p}\psi_{H}\left(\mathfrak{z}\right)+1,\textrm{ }\forall\mathfrak{z}\in\mathbb{Z}_{p}\textrm{ \& }\forall j\in\mathbb{Z}/p\mathbb{Z}\label{eq:Little Psi_H functional equations}
\end{equation}

\vphantom{}

II. \index{functional equation!Psi_{H}@$\Psi_{H}$}$\Psi_{H}$ is
rising-continuable to a $\left(p,q_{H}\right)$-adic function $\Psi_{H}:\mathbb{Z}_{p}\rightarrow\mathbb{C}_{q_{H}}$
given by: 
\begin{equation}
\Psi_{H}\left(\mathfrak{z}\right)\overset{\mathcal{F}_{p,q_{H}}}{=}\sum_{n=0}^{\infty}\left(\sum_{j=1}^{p-1}\beta_{H}\left(\frac{j}{p}\right)\left(\varepsilon_{n}\left(\mathfrak{z}\right)\right)^{j}\right)\left(\frac{\mu_{0}}{p}\right)^{n}\kappa_{H}\left(\left[\mathfrak{z}\right]_{p^{n}}\right),\textrm{ }\forall\mathfrak{z}\in\mathbb{Z}_{p}\label{eq:Rising-continuation of Big Psi_H}
\end{equation}
Moreover, $\Psi_{H}\left(\mathfrak{z}\right)$ is the unique rising-continuous
$\left(p,q_{H}\right)$-adic function satisfying the system of functional
equations: 
\begin{equation}
\Psi_{H}\left(p\mathfrak{z}+j\right)=\frac{\mu_{j}}{p}\Psi_{H}\left(\mathfrak{z}\right)+H_{j}\left(0\right)-\beta_{H}\left(0\right),\textrm{ }\forall\mathfrak{z}\in\mathbb{Z}_{p}\textrm{ \& }\forall j\in\mathbb{Z}/p\mathbb{Z}\label{eq:Big Psi_H functional equations}
\end{equation}
\end{lem}
\begin{rem}
(\ref{eq:Little Psi_H functional equations}) is the most significant
part of this result. It shows that $\psi_{H}$ is a shifted form of
$\chi_{H}$, and\textemdash crucially\textemdash that this relation
between $\chi_{H}$ and $\psi_{H}$ is \emph{independent }of the value
of $\alpha_{H}\left(0\right)$. Indeed, it is solely a consequence
of the properties of the expression (\ref{eq:Definition of Little Psi_H}). 
\end{rem}
Proof: For both parts, we use (\ref{eq:Relation between truncations and functional equations, version 2})
from \textbf{Lemma \ref{lem:Functional equations and truncation}}.
With this, observe that the functional equations (\ref{eq:Kappa_H functional equations}):
\begin{equation}
\kappa_{H}\left(pm+j\right)=\frac{\mu_{j}}{\mu_{0}}\kappa_{H}\left(m\right)
\end{equation}
then imply that: 
\begin{equation}
\kappa_{H}\left(\left[pm+j\right]_{p^{n}}\right)=\frac{\mu_{j}}{\mu_{0}}\kappa_{H}\left(\left[m\right]_{p^{n-1}}\right),\textrm{ }\forall m\in\mathbb{N}_{0},\textrm{ }\forall n\in\mathbb{N}_{1},\textrm{ }\forall j\in\mathbb{Z}/p\mathbb{Z}\label{eq:functional equation truncation Lemma applied to kappa_H}
\end{equation}
where the function $\Phi_{j}$ from (\ref{eq:Relation between truncations and functional equations, version 1})
is, in this case: 
\begin{equation}
\Phi_{j}\left(m,n\right)=\frac{\mu_{j}}{\mu_{0}}n
\end{equation}

The rising-continuability of $\psi_{H}$ and $\psi_{H}$ to the given
series follow by the givens on $H$, which guarantee that, for each
$\mathfrak{z}\in\mathbb{Z}_{p}$, $M_{H}\left(\left[\mathfrak{z}\right]_{p^{n}}\right)=\left(\frac{\mu_{0}}{p}\right)^{n}\kappa_{H}\left(\left[\mathfrak{z}\right]_{p^{n}}\right)$
tends to $0$ in the standard $\left(p,q_{H}\right)$-adic frame as
$n\rightarrow\infty$, and convergence is then guaranteed by the same
arguments used for (\ref{eq:Limit of Fourier sum of v_p A_H hat})
and (\ref{eq:F limit of Gamma_H A_H hat Fourier series when alpha is 1})
from \textbf{Lemmata \ref{lem:v_p A_H hat summation formulae}} and
\textbf{\ref{lem:1D gamma formula}}, respectively.

All that remains is to verify the functional equations. \textbf{Theorem
\ref{thm:rising-continuability of Generic H-type functional equations}}
from Subsection \ref{subsec:3.2.1 -adic-Interpolation-of} then guarantees
the uniqueness of $\psi_{H}$ and $\Psi_{H}$ as rising-continuous
solutions of their respective systems of functional equations.

\vphantom{}

I. We pull out the $n=0$ term from $\psi_{H}\left(pm+j\right)$:

\begin{align*}
\psi_{H}\left(pm+j\right) & =\frac{\mu_{j}}{p}\frac{M_{H}\left(m\right)}{1-\frac{\mu_{0}}{p}}+\overbrace{\kappa_{H}\left(0\right)}^{1}+\sum_{n=1}^{\lambda_{p}\left(m\right)}\left(\frac{\mu_{0}}{p}\right)^{n}\kappa_{H}\left(\left[pm+j\right]_{p^{n}}\right)\\
 & =\frac{\mu_{j}}{p}\frac{M_{H}\left(m\right)}{1-\frac{\mu_{0}}{p}}+1+\sum_{n=1}^{\lambda_{p}\left(m\right)}\left(\frac{\mu_{0}}{p}\right)^{n}\frac{\mu_{j}}{p}\kappa_{H}\left(\left[m\right]_{p^{n-1}}\right)\\
 & =1+\frac{\mu_{j}}{p}\left(\frac{M_{H}\left(m\right)}{1-\frac{\mu_{0}}{p}}+\sum_{n=1}^{\lambda_{p}\left(m\right)}\left(\frac{\mu_{0}}{p}\right)^{n-1}\kappa_{H}\left(\left[m\right]_{p^{n-1}}\right)\right)\\
 & =1+\frac{\mu_{j}}{p}\underbrace{\left(\frac{M_{H}\left(m\right)}{1-\frac{\mu_{0}}{p}}+\sum_{n=0}^{\lambda_{p}\left(m\right)-1}\left(\frac{\mu_{0}}{p}\right)^{n}\kappa_{H}\left(\left[m\right]_{p^{n}}\right)\right)}_{\psi_{H}\left(m\right)}\\
 & =1+\frac{\mu_{j}}{p}\psi_{H}\left(m\right)
\end{align*}
Consequently: 
\begin{equation}
\psi_{H}\left(pm+j\right)=\frac{\mu_{j}}{p}\psi_{H}\left(m\right)+1,\textrm{ }\forall m\geq0\textrm{ \& }\forall j\in\mathbb{Z}/p\mathbb{Z}\label{eq:Little Psi_H functional equation on the integers}
\end{equation}
shows that (\ref{eq:Little Psi_H functional equation on the integers})
extends to hold for the rising-continuation of $\psi_{H}$, and that
this rising-continuation is the \emph{unique }$\left(p,q_{H}\right)$-adic
function satisfying (\ref{eq:Little Psi_H functional equations}).

Finally, letting $m\in\mathbb{N}_{0}$ and setting $\mathfrak{z}=m$,
the right-hand side of (\ref{eq:Rising-continuation of Little Psi_H})
becomes: 
\begin{align*}
\sum_{n=0}^{\infty}\left(\frac{\mu_{0}}{p}\right)^{n}\kappa_{H}\left(\left[m\right]_{p^{n}}\right) & =\sum_{n=0}^{\lambda_{p}\left(m\right)-1}\left(\frac{\mu_{0}}{p}\right)^{n}\kappa_{H}\left(\left[m\right]_{p^{n}}\right)+\sum_{n=\lambda_{p}\left(m\right)}^{\infty}\left(\frac{\mu_{0}}{p}\right)^{n}\kappa_{H}\left(m\right)\\
 & \overset{\mathbb{C}}{=}\sum_{n=0}^{\lambda_{p}\left(m\right)-1}\left(\frac{\mu_{0}}{p}\right)^{n}\kappa_{H}\left(\left[m\right]_{p^{n}}\right)+\frac{\kappa_{H}\left(m\right)\left(\frac{\mu_{0}}{p}\right)^{\lambda_{p}\left(m\right)}}{1-\frac{\mu_{0}}{p}}\\
 & =\frac{M_{H}\left(m\right)}{1-\frac{\mu_{0}}{p}}+\sum_{n=0}^{\lambda_{p}\left(m\right)-1}\left(\frac{\mu_{0}}{p}\right)^{n}\kappa_{H}\left(\left[m\right]_{p^{n}}\right)\\
 & =\psi_{H}\left(m\right)
\end{align*}
Hence, (\ref{eq:Rising-continuation of Little Psi_H}) converges to
$\psi_{H}$ in the standard frame.

\vphantom{}

II. Pulling out $n=0$ from (\ref{eq:Definition of Big Psi_H}) yields:
\begin{align*}
\Psi_{H}\left(m\right) & =-\beta_{H}\left(0\right)\frac{M_{H}\left(m\right)}{1-\frac{\mu_{0}}{p}}+\sum_{k=1}^{p-1}\beta_{H}\left(\frac{k}{p}\right)\varepsilon_{0}^{k}\left(m\right)\\
 & +\sum_{n=1}^{\lambda_{p}\left(m\right)-1}\left(\sum_{k=1}^{p-1}\beta_{H}\left(\frac{k}{p}\right)\varepsilon_{n}^{k}\left(m\right)\right)\left(\frac{\mu_{0}}{p}\right)^{n}\kappa_{H}\left(\left[m\right]_{p^{n}}\right)
\end{align*}
Then:

\begin{align*}
\sum_{k=0}^{p-1}\beta_{H}\left(\frac{k}{p}\right)\varepsilon_{0}^{k}\left(m\right) & =\sum_{k=0}^{p-1}\left(\frac{1}{p}\sum_{j=0}^{p-1}H_{j}\left(0\right)e^{-\frac{2\pi ijk}{p}}\right)\left(e^{\frac{2\pi i}{p}\left(\left[m\right]_{p}-\left[m\right]_{p^{0}}\right)}\right)^{k}\\
 & =\sum_{j=0}^{p-1}H_{j}\left(0\right)\frac{1}{p}\sum_{k=0}^{p-1}e^{-\frac{2\pi ijk}{p}}e^{\frac{2\pi ikm}{p}}\\
 & =\sum_{j=0}^{p-1}H_{j}\left(0\right)\frac{1}{p}\sum_{k=0}^{p-1}e^{\frac{2\pi ik\left(m-j\right)}{p}}\\
 & =\sum_{j=0}^{p-1}H_{j}\left(0\right)\left[j\overset{p}{\equiv}m\right]\\
 & =H_{\left[m\right]_{p}}\left(0\right)
\end{align*}
and so: 
\begin{equation}
\sum_{k=1}^{p-1}\beta_{H}\left(\frac{k}{p}\right)\varepsilon_{0}^{k}\left(m\right)=H_{\left[m\right]_{p}}\left(0\right)-\beta_{H}\left(0\right)\label{eq:Fourier sum (but not with 0) of beta_H}
\end{equation}
Consequently: 
\begin{align*}
\Psi_{H}\left(m\right) & =-\beta_{H}\left(0\right)\frac{M_{H}\left(m\right)}{1-\frac{\mu_{0}}{p}}+H_{\left[m\right]_{p}}\left(0\right)-\beta_{H}\left(0\right)\\
 & +\sum_{n=1}^{\lambda_{p}\left(m\right)-1}\left(\sum_{k=1}^{p-1}\beta_{H}\left(\frac{k}{p}\right)\varepsilon_{n}^{k}\left(m\right)\right)\left(\frac{\mu_{0}}{p}\right)^{n}\kappa_{H}\left(\left[m\right]_{p^{n}}\right)
\end{align*}

Replacing $m$ with $pm+j$ (where at least one of $m$ and $j$ is
non-zero), we use (\ref{eq:functional equation truncation Lemma applied to kappa_H})
and the functional equations for $M_{H}$, $\varepsilon_{n}$, and
$\lambda_{p}$ to obtain: 
\begin{align*}
\Psi_{H}\left(m\right) & =-\beta_{H}\left(0\right)\frac{\mu_{j}}{p}\frac{M_{H}\left(m\right)}{1-\frac{\mu_{0}}{p}}+H_{j}\left(0\right)-\beta_{H}\left(0\right)\\
 & +\sum_{n=1}^{\lambda_{p}\left(m\right)}\left(\sum_{k=1}^{p-1}\beta_{H}\left(\frac{k}{p}\right)\varepsilon_{n-1}^{k}\left(m\right)\varepsilon_{n}^{k}\left(j\right)\right)\left(\frac{\mu_{0}}{p}\right)^{n}\frac{\mu_{j}}{\mu_{0}}\kappa_{H}\left(\left[m\right]_{p^{n-1}}\right)
\end{align*}
Because $j\in\left\{ 0,\ldots,p-1\right\} $, note that $\left[j\right]_{p^{n}}=j$
for all $n\geq1$. As such: 
\[
\varepsilon_{n}\left(j\right)=e^{\frac{2\pi i}{p^{n+1}}\left(\left[j\right]_{p^{n+1}}-\left[j\right]_{p^{n}}\right)}=e^{\frac{2\pi i}{p^{n+1}}\left(j-j\right)}=1,\textrm{ }\forall n\geq1
\]
Re-indexing the $n$-sum then gives us: 
\begin{align*}
\Psi_{H}\left(m\right) & =-\beta_{H}\left(0\right)\frac{\mu_{j}}{p}\frac{M_{H}\left(m\right)}{1-\frac{\mu_{0}}{p}}+H_{j}\left(0\right)-\beta_{H}\left(0\right)\\
 & +\sum_{n=0}^{\lambda_{p}\left(m\right)-1}\left(\sum_{k=1}^{p-1}\beta_{H}\left(\frac{k}{p}\right)\varepsilon_{n}^{k}\left(m\right)\right)\left(\frac{\mu_{0}}{p}\right)^{n}\frac{\mu_{0}}{p}\frac{\mu_{j}}{\mu_{0}}\kappa_{H}\left(\left[m\right]_{p^{n}}\right)\\
 & =H_{j}\left(0\right)-\beta_{H}\left(0\right)\\
 & +\frac{\mu_{j}}{p}\underbrace{\left(-\beta_{H}\left(0\right)\frac{M_{H}\left(m\right)}{1-\frac{\mu_{0}}{p}}+\sum_{n=0}^{\lambda_{p}\left(m\right)-1}\left(\sum_{k=1}^{p-1}\beta_{H}\left(\frac{k}{p}\right)\varepsilon_{n}^{k}\left(m\right)\right)\left(\frac{\mu_{0}}{p}\right)^{n}\kappa_{H}\left(\left[m\right]_{p^{n}}\right)\right)}_{\Psi_{H}\left(m\right)}\\
 & =\frac{\mu_{j}}{p}\Psi_{H}\left(m\right)+H_{j}\left(0\right)-\beta_{H}\left(0\right)
\end{align*}

Finally, letting $\mathfrak{z}=m$ where $m\in\mathbb{N}_{0}$, the
right-hand side of (\ref{eq:Rising-continuation of Big Psi_H}) becomes:
\begin{align*}
\sum_{n=0}^{\lambda_{p}\left(m\right)}\left(\sum_{j=1}^{p-1}\beta_{H}\left(\frac{j}{p}\right)\left(\varepsilon_{n}\left(m\right)\right)^{j}\right)\left(\frac{\mu_{0}}{p}\right)^{n}\kappa_{H}\left(\left[m\right]_{p^{n}}\right)\\
+\sum_{n=\lambda_{p}\left(m\right)}^{\infty}\left(\sum_{j=1}^{p-1}\beta_{H}\left(\frac{j}{p}\right)\left(\varepsilon_{n}\left(m\right)\right)^{j}\right)\left(\frac{\mu_{0}}{p}\right)^{n}\kappa_{H}\left(m\right)
\end{align*}
Here: 
\begin{equation}
\varepsilon_{n}\left(m\right)=e^{\frac{2\pi i}{p^{n+1}}\left(\left[m\right]_{p^{n+1}}-\left[m\right]_{p^{n}}\right)}=e^{\frac{2\pi i}{p^{n+1}}\left(m-m\right)}=1,\textrm{ }\forall n\geq\lambda_{p}\left(m\right)
\end{equation}
So: 
\begin{align*}
\sum_{n=0}^{\lambda_{p}\left(m\right)}\left(\sum_{j=1}^{p-1}\beta_{H}\left(\frac{j}{p}\right)\left(\varepsilon_{n}\left(m\right)\right)^{j}\right)\left(\frac{\mu_{0}}{p}\right)^{n}\kappa_{H}\left(\left[m\right]_{p^{n}}\right)\\
+\left(\sum_{j=1}^{p-1}\beta_{H}\left(\frac{j}{p}\right)\right)\sum_{n=\lambda_{p}\left(m\right)}^{\infty}\left(\frac{\mu_{0}}{p}\right)^{n}\kappa_{H}\left(m\right)
\end{align*}
Summing this in the topology of $\mathbb{C}$ yields: 
\begin{align*}
\sum_{n=0}^{\lambda_{p}\left(m\right)}\left(\sum_{j=1}^{p-1}\beta_{H}\left(\frac{j}{p}\right)\left(\varepsilon_{n}\left(m\right)\right)^{j}\right)\left(\frac{\mu_{0}}{p}\right)^{n}\kappa_{H}\left(\left[m\right]_{p^{n}}\right)\\
+\left(\sum_{j=1}^{p-1}\beta_{H}\left(\frac{j}{p}\right)\right)\frac{M_{H}\left(m\right)}{1-\frac{\mu_{0}}{p}}
\end{align*}
Using (\ref{eq:Fourier sum (but not with 0) of beta_H}) with $m=0$
(which makes $\varepsilon_{0}^{k}\left(m\right)$ identically equal
to $1$) gives us: 
\begin{equation}
\sum_{j=1}^{p-1}\beta_{H}\left(\frac{j}{p}\right)=H_{0}\left(0\right)-\beta_{H}\left(0\right)=-\beta_{H}\left(0\right)
\end{equation}
and so, we are left with: 
\begin{equation}
-\beta_{H}\left(0\right)\frac{M_{H}\left(m\right)}{1-\frac{\mu_{0}}{p}}+\sum_{n=0}^{\lambda_{p}\left(m\right)}\left(\sum_{j=1}^{p-1}\beta_{H}\left(\frac{j}{p}\right)\left(\varepsilon_{n}\left(m\right)\right)^{j}\right)\left(\frac{\mu_{0}}{p}\right)^{n}\kappa_{H}\left(\left[m\right]_{p^{n}}\right)
\end{equation}
This, of course, is precisely $\Psi_{H}\left(m\right)$. This proves
that (\ref{eq:Rising-continuation of Big Psi_H}) converges to $\Psi_{H}$
in the standard frame.

Q.E.D.
\begin{prop}[\textbf{Quasi-Integrability of $\psi_{H}$ \& $\Psi_{H}$}]
Let $H$ be as given at the start of \emph{Section \pageref{sec:4.2 Fourier-Transforms-=00003D000026}}.
Then:

\vphantom{}

I. 
\begin{equation}
\psi_{H}\left(\mathfrak{z}\right)\overset{\mathcal{F}_{p,q_{H}}}{=}\begin{cases}
\lim_{N\rightarrow\infty}\sum_{0<\left|t\right|_{p}\leq p^{N}}v_{p}\left(t\right)\hat{A}_{H}\left(t\right)e^{2\pi i\left\{ t\mathfrak{z}\right\} _{p}} & \textrm{if }\alpha_{H}\left(0\right)=1\\
\lim_{N\rightarrow\infty}\sum_{\left|t\right|_{p}\leq p^{N}}\frac{\hat{A}_{H}\left(t\right)}{1-\alpha_{H}\left(0\right)}e^{2\pi i\left\{ t\mathfrak{z}\right\} _{p}} & \textrm{if }\alpha_{H}\left(0\right)\neq1
\end{cases},\textrm{ }\forall\mathfrak{z}\in\mathbb{Z}_{p}\label{eq:Little Psi_H as an F-limit}
\end{equation}
As such, $\psi_{H}$ is quasi-integrable with respect to the standard
$\left(p,q_{H}\right)$-adic quasi-integrability frame, and the function
$\hat{\psi}_{H}:\hat{\mathbb{Z}}_{p}\rightarrow\overline{\mathbb{Q}}$
defined by: 
\begin{equation}
\hat{\psi}_{H}\left(t\right)\overset{\textrm{def}}{=}\begin{cases}
\begin{cases}
0 & \textrm{if }t=0\\
v_{p}\left(t\right)\hat{A}_{H}\left(t\right) & \textrm{if }t\neq0
\end{cases} & \textrm{if }\alpha_{H}\left(0\right)=1\\
\frac{\hat{A}_{H}\left(t\right)}{1-\alpha_{H}\left(0\right)} & \textrm{if }\alpha_{H}\left(0\right)\neq1
\end{cases},\textrm{ }\forall t\in\hat{\mathbb{Z}}_{p}\label{eq:Fourier Transform of Little Psi_H}
\end{equation}
is a Fourier transform of $\psi_{H}$. Hence: 
\begin{equation}
\hat{\psi}_{H,N}\left(\mathfrak{z}\right)\overset{\overline{\mathbb{Q}}}{=}-N\left(\frac{\mu_{0}}{p}\right)^{N}\kappa_{H}\left(\left[\mathfrak{z}\right]_{p^{N}}\right)+\sum_{n=0}^{N-1}\left(\frac{\mu_{0}}{p}\right)^{n}\kappa_{H}\left(\left[\mathfrak{z}\right]_{p^{n}}\right)\label{eq:Little Psi H N twiddle when alpha is 1}
\end{equation}
when $\alpha_{H}\left(0\right)=1$ and:

\begin{equation}
\tilde{\psi}_{H,N}\left(\mathfrak{z}\right)\overset{\overline{\mathbb{Q}}}{=}\left(\frac{\mu_{0}}{p}\right)^{N}\frac{\kappa_{H}\left(\left[\mathfrak{z}\right]_{p^{N}}\right)}{1-\alpha_{H}\left(0\right)}+\sum_{n=0}^{N-1}\left(\frac{\mu_{0}}{p}\right)^{n}\kappa_{H}\left(\left[\mathfrak{z}\right]_{p^{n}}\right)\label{eq:Little Psi H N twiddle when alpha is not 1}
\end{equation}
when $\alpha_{H}\left(0\right)\neq1$.

\vphantom{}

II. 
\begin{equation}
\Psi_{H}\left(\mathfrak{z}\right)\overset{\mathcal{F}_{p,q_{H}}}{=}\lim_{N\rightarrow\infty}\sum_{0<\left|t\right|_{p}\leq p^{N}}\gamma_{H}\left(\frac{t\left|t\right|_{p}}{p}\right)\hat{A}_{H}\left(t\right)e^{2\pi i\left\{ t\mathfrak{z}\right\} _{p}},\textrm{ }\forall\mathfrak{z}\in\mathbb{Z}_{p}\label{eq:Big Psi_H as an F-limit}
\end{equation}
As such, $\Psi_{H}$ is quasi-integrable with respect to the standard
$\left(p,q_{H}\right)$-adic quasi-integrability frame, and the function
$\hat{\Psi}_{H}:\hat{\mathbb{Z}}_{p}\rightarrow\overline{\mathbb{Q}}$
defined by:
\begin{equation}
\hat{\Psi}_{H}\left(t\right)\overset{\textrm{def}}{=}\begin{cases}
0 & \textrm{if }t=0\\
\gamma_{H}\left(\frac{t\left|t\right|_{p}}{p}\right)\hat{A}_{H}\left(t\right) & \textrm{if }t\neq0
\end{cases},\textrm{ }\forall t\in\hat{\mathbb{Z}}_{p}\label{eq:Fourier Transform of Big Psi_H}
\end{equation}
is a Fourier transform of $\Psi_{H}$. Hence: 
\begin{equation}
\tilde{\Psi}_{H,N}\left(\mathfrak{z}\right)\overset{\overline{\mathbb{Q}}}{=}\sum_{n=0}^{N-1}\left(\sum_{j=1}^{p-1}\beta_{H}\left(\frac{j}{p}\right)\varepsilon_{n}^{j}\left(\mathfrak{z}\right)\right)\left(\frac{\mu_{0}}{p}\right)^{n}\kappa_{H}\left(\left[\mathfrak{z}\right]_{p^{n}}\right)\label{eq:Big Psi H N twiddle}
\end{equation}
\end{prop}
Proof:

I. When $\alpha_{H}\left(0\right)=1$, (\ref{eq:Little Psi_H as an F-limit})
follows from the limit (\ref{eq:Limit of Fourier sum of v_p A_H hat})
from \textbf{Lemma \ref{lem:v_p A_H hat summation formulae}}. So,
the $\alpha_{H}\left(0\right)=1$ case of (\ref{eq:Fourier Transform of Little Psi_H})
is then indeed a Fourier transform of $\psi_{H}$. This establishes
the $\mathcal{F}_{p,q_{H}}$-quasi-integrability of $\psi_{H}$ when
$\alpha_{H}\left(0\right)=1$.

When $\alpha_{H}\left(0\right)\neq1$, comparing (\ref{eq:Convolution of dA_H and D_N})
and (\ref{eq:Rising-continuation of Little Psi_H}) we observe that
$\left(1-\alpha_{H}\left(0\right)\right)\psi_{H}\left(\mathfrak{z}\right)$
then takes the values: 
\[
\lim_{N\rightarrow\infty}\sum_{\left|t\right|_{p}\leq p^{N}}\hat{A}_{H}\left(t\right)e^{2\pi i\left\{ t\mathfrak{z}\right\} _{p}}
\]
at every $\mathfrak{z}\in\mathbb{Z}_{p}^{\prime}$. Since $H$ is
semi-basic, applying our standard argument involving the $q$-adic
decay of: 
\[
\left(\frac{\mu_{0}}{p}\right)^{n}\kappa_{H}\left(\left[\mathfrak{z}\right]_{p^{n}}\right)
\]
to (\ref{eq:Little Psi_H functional equation on the integers}) shows
that $\psi_{H}\left(\left[\mathfrak{z}\right]_{p^{N}}\right)$ converges
in $\mathbb{Z}_{q_{H}}$ as $N\rightarrow\infty$ for all $\mathfrak{z}\in\mathbb{Z}_{p}$.
This establishes the rising-continuity of $\psi_{H}$. Consequently,
the function $A_{H}:\mathbb{Z}_{p}\rightarrow\mathbb{Z}_{q_{H}}$
defined by:

\begin{equation}
A_{H}\left(\mathfrak{z}\right)\overset{\textrm{def}}{=}\begin{cases}
\left(1-\alpha_{H}\left(0\right)\right)\psi_{H}\left(\mathfrak{z}\right) & \textrm{if }\mathfrak{z}\in\mathbb{N}_{0}\\
\left(1-\alpha_{H}\left(0\right)\right)\lim_{N\rightarrow\infty}\sum_{\left|t\right|_{p}\leq p^{N}}\hat{A}_{H}\left(t\right)e^{2\pi i\left\{ t\mathfrak{z}\right\} _{p}} & \textrm{if }\mathfrak{z}\in\mathbb{Z}_{p}^{\prime}
\end{cases}\label{eq:definition of A_H}
\end{equation}
is then rising-continuous, and its restriction to $\mathbb{N}_{0}$
is equal to $\left(1-\alpha_{H}\left(0\right)\right)\psi_{H}\left(\mathfrak{z}\right)$.
By (\ref{eq:Little Psi_H as an F-limit}), it then follows that $\hat{A}_{H}$
is a Fourier transform of $A_{H}$ with respect to the standard $\left(p,q_{H}\right)$-adic
frame, and thus, that $A_{H}$ is $\mathcal{F}_{p,q_{H}}$-quasi-integrable.
So, when $\alpha_{H}\left(0\right)\neq1$, $\psi_{H}\left(\mathfrak{z}\right)=\left(1-\alpha_{H}\left(0\right)\right)A_{H}\left(\mathfrak{z}\right)$
is also $\mathcal{F}_{p,q_{H}}$-quasi-integrable, and the $\alpha_{H}\left(0\right)\neq1$
case of (\ref{eq:Fourier Transform of Little Psi_H}) is then a Fourier
transform of $\psi_{H}$. As for formulae (\ref{eq:Little Psi H N twiddle when alpha is 1})
and (\ref{eq:Little Psi H N twiddle when alpha is not 1}), these
are nothing more than restatements of \textbf{Lemma \ref{lem:v_p A_H hat summation formulae}}
and (\ref{eq:Convolution of dA_H and D_N}), respectively.

\vphantom{}

II. (\ref{eq:Big Psi_H as an F-limit}) and (\ref{eq:Fourier Transform of Big Psi_H})
are exactly what we proved in \textbf{Lemma \ref{lem:1D gamma formula}};
(\ref{eq:Big Psi H N twiddle}) is just (\ref{eq:Gamma formula}).

Q.E.D.

\vphantom{}

Now that we know that $\psi_{H}$ and $\Psi_{H}$ are quasi-integrable
\emph{regardless} of the value of $\alpha_{H}\left(0\right)$ (or
$p$, for that matter), we can then reverse-engineer formulae for
the Fourier transform of $\chi_{H}$ \emph{regardless }of the value
of $\alpha_{H}\left(0\right)$; we need only solve a system of linear
equations to find the constants which relate $\psi_{H}$, $\Psi_{H}$
and $\chi_{H}$. We do the $p=2$ and $p\geq3$ cases separately.
Both of the following theorems can then be immediately applied to
obtain formulae for $\hat{\chi}_{H}$, $\tilde{\chi}_{H,N}\left(\mathfrak{z}\right)$
\emph{and} $\chi_{H}\left(\mathfrak{z}\right)$. 
\begin{thm}[\textbf{Quasi-Integrability of $\chi_{H}$ for a $2$-Hydra map}]
\label{thm:Quasi-integrability of an arbitrary 2-Hydra map}Let $H$
be any contracting non-singular semi-basic $2$-Hydra map which fixes
$0$. Then:\index{chi{H}@$\chi_{H}$!quasi-integrability} 
\begin{equation}
\chi_{H}\left(\mathfrak{z}\right)\overset{\mathcal{F}_{2,q_{H}}}{=}-\gamma_{H}\left(\frac{1}{2}\right)+\left(1-H^{\prime}\left(0\right)\right)\gamma_{H}\left(\frac{1}{2}\right)\psi_{H}\left(\mathfrak{z}\right),\textrm{ }\forall\mathfrak{z}\in\mathbb{Z}_{2}\label{eq:Chi_H in terms of Little Psi_H for a 2-Hydra map}
\end{equation}
In particular, by \emph{(\ref{eq:Fourier Transform of Little Psi_H})},
this shows that $\chi_{H}$ is then quasi-integrable with respect
to the standard $\left(2,q_{H}\right)$-adic quasi-integrability frame,
with the function $\hat{\chi}_{H}:\hat{\mathbb{Z}}_{2}\rightarrow\overline{\mathbb{Q}}$
defined by:

\index{chi{H}@$\chi_{H}$!Fourier transform}

\begin{equation}
\hat{\chi}_{H}\left(t\right)\overset{\textrm{def}}{=}\begin{cases}
\begin{cases}
-\gamma_{H}\left(\frac{1}{2}\right) & \textrm{if }t=0\\
\left(1-H^{\prime}\left(0\right)\right)\gamma_{H}\left(\frac{1}{2}\right)v_{2}\left(t\right)\hat{A}_{H}\left(t\right) & \textrm{if }t\neq0
\end{cases} & \textrm{if }\alpha_{H}\left(0\right)=1\\
-\gamma_{H}\left(\frac{1}{2}\right)\mathbf{1}_{0}\left(t\right)+\gamma_{H}\left(\frac{1}{2}\right)\frac{1-H^{\prime}\left(0\right)}{1-\alpha_{H}\left(0\right)}\hat{A}_{H}\left(t\right) & \textrm{if }\alpha_{H}\left(0\right)\neq1
\end{cases},\textrm{ }\forall t\in\hat{\mathbb{Z}}_{2}\label{eq:Fourier Transform of Chi_H for a contracting semi-basic 2-Hydra map}
\end{equation}
being a Fourier transform of $\chi_{H}$. 
\end{thm}
Proof: Let $A,B$ be undetermined constants, and let $\psi\left(\mathfrak{z}\right)=A\chi_{H}\left(\mathfrak{z}\right)+B$.
Then, multiplying $\chi_{H}$'s functional equations (\ref{eq:Functional Equations for Chi_H over the rho-adics})
by $A$ and adding $B$ gives: 
\begin{align*}
A\chi_{H}\left(2\mathfrak{z}+j\right)+B & =\frac{Aa_{j}\chi_{H}\left(\mathfrak{z}\right)+Ab_{j}}{d_{j}}+B\\
 & \Updownarrow\\
\psi\left(2\mathfrak{z}+j\right) & =\frac{a_{j}\left(\psi\left(\mathfrak{z}\right)-B\right)+Ab_{j}}{d_{j}}+B\\
 & =\frac{\mu_{j}}{2}\psi\left(\mathfrak{z}\right)+A\frac{b_{j}}{d_{j}}+B\left(1-\frac{a_{j}}{d_{j}}\right)
\end{align*}
for $j\in\left\{ 0,1\right\} $. Since: 
\[
\psi_{H}\left(2\mathfrak{z}+j\right)=\frac{\mu_{j}}{2}\psi\left(\mathfrak{z}\right)+1
\]
it then follows that to make $\psi=\psi_{H}$, we need for $A$ and
$B$ to satisfy: 
\begin{align*}
\frac{b_{0}}{d_{0}}A+\left(1-\frac{a_{0}}{d_{0}}\right)B & =1\\
\frac{b_{1}}{d_{1}}A+\left(1-\frac{a_{1}}{d_{1}}\right)B & =1
\end{align*}

Now, since: 
\begin{align*}
\psi_{H}\left(\mathfrak{z}\right) & =A\chi_{H}\left(\mathfrak{z}\right)+B\\
 & \Updownarrow\\
\chi_{H}\left(\mathfrak{z}\right) & =\frac{1}{A}\psi_{H}\left(\mathfrak{z}\right)-\frac{B}{A}
\end{align*}
upon letting $C=\frac{1}{A}$ and $D=-\frac{B}{A}$, we have that
$A=\frac{1}{C}$, $B=-\frac{D}{C}$, and hence, our linear system
becomes: 
\begin{align*}
d_{0}C+\left(d_{0}-a_{0}\right)D & =b_{0}\\
d_{1}C+\left(d_{1}-a_{1}\right)D & =b_{1}
\end{align*}
where: 
\begin{equation}
\chi_{H}\left(\mathfrak{z}\right)=C\psi_{H}\left(\mathfrak{z}\right)+D\label{eq:Chi_H =00003D00003D C Psi_H + D}
\end{equation}

The determinant of the coefficient matrix is: 
\begin{equation}
a_{0}d_{1}-a_{1}d_{0}
\end{equation}
So, as long as $a_{0}d_{1}\neq a_{1}d_{0}$, the system admits a unique
solution, which is easily computed to be: 
\begin{align*}
C & =\frac{a_{0}b_{1}-a_{1}b_{0}+b_{0}d_{1}-b_{1}d_{0}}{a_{0}d_{1}-a_{1}d_{0}}\\
D & =\frac{b_{1}d_{0}-d_{1}b_{0}}{a_{0}d_{1}-a_{1}d_{0}}
\end{align*}
Finally, note that: 
\begin{align*}
a_{0}d_{1} & \neq a_{1}d_{0}\\
 & \Updownarrow\\
0 & \neq\frac{a_{0}}{d_{0}}-\frac{a_{1}}{d_{1}}\\
 & =2\alpha_{H}\left(\frac{1}{2}\right)
\end{align*}
and so: 
\[
D=\frac{b_{1}d_{0}-d_{1}b_{0}}{a_{0}d_{1}-a_{1}d_{0}}=-\frac{\frac{1}{2}\left(\frac{b_{0}}{d_{0}}-\frac{b_{1}}{d_{1}}\right)}{\frac{1}{2}\left(\frac{a_{0}}{d_{0}}-\frac{a_{1}}{d_{1}}\right)}=-\frac{\beta_{H}\left(\frac{1}{2}\right)}{\alpha_{H}\left(\frac{1}{2}\right)}=-\gamma_{H}\left(\frac{1}{2}\right)
\]
As for $C$, note that our requirement that $H\left(0\right)=0$ then
forces $b_{0}=0$, and so: 
\begin{align*}
C & =\frac{a_{0}b_{1}-a_{1}b_{0}+b_{0}d_{1}-b_{1}d_{0}}{a_{0}d_{1}-a_{1}d_{0}}\\
 & =\frac{a_{0}-d_{0}}{d_{0}}\frac{\frac{b_{1}}{d_{1}}}{\frac{a_{0}}{d_{0}}-\frac{a_{1}}{d_{1}}}\\
\left(\begin{array}{c}
\alpha_{H}\left(\frac{1}{2}\right)=\frac{1}{2}\left(\frac{a_{0}}{a_{0}}-\frac{a_{1}}{a_{1}}\right)\\
\beta_{H}\left(\frac{1}{2}\right)=\frac{1}{2}\left(\frac{b_{0}}{d_{0}}-\frac{b_{1}}{d_{1}}\right)=-\frac{1}{2}\frac{b_{1}}{d_{1}}
\end{array}\right); & =\frac{a_{0}-d_{0}}{d_{0}}\frac{-2\beta_{H}\left(\frac{1}{2}\right)}{2\alpha_{H}\left(\frac{1}{2}\right)}\\
 & =\frac{d_{0}-a_{0}}{d_{0}}\gamma_{H}\left(\frac{1}{2}\right)\\
\left(H^{\prime}\left(0\right)=\frac{a_{0}}{d_{0}}\right); & =\left(1-H^{\prime}\left(0\right)\right)\gamma_{H}\left(\frac{1}{2}\right)
\end{align*}
Plugging these values for $C$ and $D$ into (\ref{eq:Chi_H =00003D00003D C Psi_H + D})
produces (\ref{eq:Chi_H in terms of Little Psi_H for a 2-Hydra map}).

Q.E.D.

\vphantom{}

Next up, our non-trivial formulae for $\tilde{\chi}_{H,N}$. 
\begin{cor}
\label{cor:Non-trivial formula for Chi_H,N twiddle for a 2-hydra map, alpha arbitrary}Let\index{chi{H}@$\chi_{H}$!$N$th partial Fourier series}
$H$ be a $2$-Hydra map like in \textbf{\emph{Theorem \ref{thm:Quasi-integrability of an arbitrary 2-Hydra map}}},
with $\hat{\chi}_{H}$ being as defined by \emph{(\ref{eq:Fourier Transform of Chi_H for a contracting semi-basic 2-Hydra map})}.
Then:

\begin{align}
\tilde{\chi}_{H,N}\left(\mathfrak{z}\right) & \overset{\overline{\mathbb{Q}}}{=}-\gamma_{H}\left(\frac{1}{2}\right)-\left(1-H^{\prime}\left(0\right)\right)\gamma_{H}\left(\frac{1}{2}\right)N\left(\frac{\mu_{0}}{2}\right)^{N}\kappa_{H}\left(\left[\mathfrak{z}\right]_{2^{N}}\right)\label{eq:Explicit Formula for Chi_H,N twiddle when rho is 2 and alpha is 1}\\
 & +\left(1-H^{\prime}\left(0\right)\right)\gamma_{H}\left(\frac{1}{2}\right)\sum_{n=0}^{N-1}\left(\frac{\mu_{0}}{2}\right)^{n}\kappa_{H}\left(\left[\mathfrak{z}\right]_{2^{n}}\right)\nonumber 
\end{align}
if $\alpha_{H}\left(0\right)=1$, and: 
\begin{align}
\tilde{\chi}_{H,N}\left(\mathfrak{z}\right) & \overset{\overline{\mathbb{Q}}}{=}-\gamma_{H}\left(\frac{1}{2}\right)+\gamma_{H}\left(\frac{1}{2}\right)\frac{1-H^{\prime}\left(0\right)}{1-\alpha_{H}\left(0\right)}\left(\frac{\mu_{0}}{2}\right)^{N}\kappa_{H}\left(\left[\mathfrak{z}\right]_{2^{N}}\right)\label{eq:Explicit Formula for Chi_H,N twiddle when rho is 2 and alpha is not 1}\\
 & +\gamma_{H}\left(\frac{1}{2}\right)\left(1-H^{\prime}\left(0\right)\right)\sum_{n=0}^{N-1}\left(\frac{\mu_{0}}{2}\right)^{n}\kappa_{H}\left(\left[\mathfrak{z}\right]_{2^{n}}\right)\nonumber 
\end{align}
if $\alpha_{H}\left(0\right)\neq1$. 
\end{cor}
Proof: Take (\ref{eq:Fourier Transform of Chi_H for a contracting semi-basic 2-Hydra map})
and use formulas from \textbf{Lemma \ref{lem:1D gamma formula}},
\textbf{Lemma (\ref{lem:v_p A_H hat summation formulae}}, and \textbf{Theorem
\ref{thm:Properties of dA_H}}.

Q.E.D.
\begin{cor}[\textbf{$\mathcal{F}$-series for $\chi_{H}$ when $H$ is a $2$-Hydra
map}]
Let \label{cor:F-series for Chi_H for p equals 2}\index{chi{H}@$\chi_{H}$!mathcal{F}-series@$\mathcal{F}$-series}$H$
be as given in \textbf{\emph{Theorem \ref{thm:Quasi-integrability of an arbitrary 2-Hydra map}}}.
Then: 
\begin{equation}
\chi_{H}\left(\mathfrak{z}\right)\overset{\mathcal{F}_{2,q_{H}}}{=}-\gamma_{H}\left(\frac{1}{2}\right)+\left(1-H^{\prime}\left(0\right)\right)\gamma_{H}\left(\frac{1}{2}\right)\sum_{n=0}^{\infty}\left(\frac{\mu_{0}}{2}\right)^{n}\left(\frac{a_{1}/d_{1}}{a_{0}/d_{0}}\right)^{\#_{1}\left(\left[\mathfrak{z}\right]_{2^{n}}\right)},\textrm{ }\forall\mathfrak{z}\in\mathbb{Z}_{2}\label{eq:Chi_H non-trivial series formula, rho =00003D00003D 2}
\end{equation}
\end{cor}
Proof: Write out (\ref{eq:Chi_H in terms of Little Psi_H for a 2-Hydra map}),
then use (\ref{eq:Definition of Kappa_H}) and $M_{H}$'s explicit
formula (\textbf{Proposition \ref{prop:Explicit Formulas for M_H}})
to write: 
\begin{align*}
\left(\frac{\mu_{0}}{2}\right)^{n}\kappa_{H}\left(m\right) & =\left(\frac{\mu_{0}}{2}\right)^{n}\left(\frac{2}{\mu_{0}}\right)^{\lambda_{2}\left(m\right)}M_{H}\left(m\right)\\
 & =\left(\frac{\mu_{0}}{2}\right)^{n}\left(\frac{2}{\mu_{0}}\right)^{\lambda_{2}\left(m\right)}\frac{\mu_{0}^{\#_{2:0}\left(m\right)}\mu_{1}^{\#_{1}\left(m\right)}}{2^{\lambda_{2}\left(m\right)}}\\
\left(\#_{2:0}\left(m\right)=\lambda_{2}\left(m\right)-\#_{1}\left(m\right)\right); & =\left(\frac{\mu_{0}}{2}\right)^{n}\left(\frac{\mu_{1}}{\mu_{0}}\right)^{\#_{1}\left(m\right)}\\
 & =\left(\frac{\mu_{0}}{2}\right)^{n}\left(\frac{a_{1}/d_{1}}{a_{0}/d_{0}}\right)^{\#_{1}\left(m\right)}
\end{align*}

Q.E.D.

\vphantom{}

Now we move to the general case. 
\begin{thm}[\textbf{Quasi-Integrability of $\chi_{H}$ for a $p$-Hydra map}]
\index{chi{H}@$\chi_{H}$!mathcal{F}-series@$\mathcal{F}$-series}\label{thm:F-series for an arbitrary 1D Chi_H}Let
$H$ be an arbitrary contracting, non-singular, semi-basic $p$-Hydra
map which fixes $0$, where $p\geq2$. Suppose additionally that $H_{j}\left(0\right)\neq0$
for any $j\in\left\{ 1,\ldots,p-1\right\} $. Then\footnote{I call the following identity for $\chi_{H}$ the \textbf{$\mathcal{F}$-series}
of $\chi_{H}$.}: 
\begin{equation}
\chi_{H}\left(\mathfrak{z}\right)\overset{\mathcal{F}_{p,q_{H}}}{=}\beta_{H}\left(0\right)\psi_{H}\left(\mathfrak{z}\right)+\Psi_{H}\left(\mathfrak{z}\right),\textrm{ }\forall\mathfrak{z}\in\mathbb{Z}_{p}\label{eq:Chi_H in terms of Little Psi_H and Big Psi_H for arbitrary rho}
\end{equation}
In particular, using \emph{(\ref{eq:Fourier Transform of Little Psi_H})},
this shows that $\chi_{H}$ is then quasi-integrable with respect
to the standard $\left(p,q_{H}\right)$-adic quasi-integrability frame,
with the function $\hat{\chi}_{H}:\hat{\mathbb{Z}}_{p}\rightarrow\overline{\mathbb{Q}}$
defined by:

\index{chi{H}@$\chi_{H}$!Fourier transform}

\begin{equation}
\hat{\chi}_{H}\left(t\right)\overset{\textrm{def}}{=}\begin{cases}
\begin{cases}
0 & \textrm{if }t=0\\
\left(\beta_{H}\left(0\right)v_{p}\left(t\right)+\gamma_{H}\left(\frac{t\left|t\right|_{p}}{p}\right)\right)\hat{A}_{H}\left(t\right) & \textrm{if }t\neq0
\end{cases} & \textrm{if }\alpha_{H}\left(0\right)=1\\
\frac{\beta_{H}\left(0\right)\hat{A}_{H}\left(t\right)}{1-\alpha_{H}\left(0\right)}+\begin{cases}
0 & \textrm{if }t=0\\
\gamma_{H}\left(\frac{t\left|t\right|_{p}}{p}\right)\hat{A}_{H}\left(t\right) & \textrm{if }t\neq0
\end{cases} & \textrm{if }\alpha_{H}\left(0\right)\neq1
\end{cases},\textrm{ }\forall t\in\hat{\mathbb{Z}}_{p}\label{eq:Fourier Transform of Chi_H for a contracting semi-basic rho-Hydra map}
\end{equation}
being a Fourier transform of $\chi_{H}$. 
\end{thm}
Proof: We make the ansatz:

\begin{equation}
\chi_{H}\left(\mathfrak{z}\right)=A\psi_{H}\left(\mathfrak{z}\right)+B\Psi_{H}\left(\mathfrak{z}\right)
\end{equation}
for constants $A$ and $B$. Using the functional equations (\ref{eq:Little Psi_H functional equations})
and (\ref{eq:Big Psi_H functional equations}), we then have that:

\begin{align*}
\chi_{H}\left(p\mathfrak{z}+j\right) & =A\psi_{H}\left(p\mathfrak{z}+j\right)+B\Psi_{H}\left(p\mathfrak{z}+j\right)\\
 & =A\frac{\mu_{j}}{p}\psi_{H}\left(\mathfrak{z}\right)+B\frac{\mu_{j}}{p}\Psi_{H}\left(\mathfrak{z}\right)+A+B\left(H_{j}\left(0\right)-\beta_{H}\left(0\right)\right)
\end{align*}
and: 
\[
\chi_{H}\left(p\mathfrak{z}+j\right)=\frac{\mu_{j}}{p}\chi_{H}\left(\mathfrak{z}\right)+H_{j}\left(0\right)=A\frac{\mu_{j}}{p}\psi_{H}\left(\mathfrak{z}\right)+B\frac{\mu_{j}}{p}\Psi_{H}\left(\mathfrak{z}\right)+H_{j}\left(0\right)
\]
Combining these two formulae for $\chi_{H}\left(p\mathfrak{z}+j\right)$
yields: 
\[
A\frac{\mu_{j}}{p}\psi_{H}\left(\mathfrak{z}\right)+B\frac{\mu_{j}}{p}\Psi_{H}\left(\mathfrak{z}\right)+A+BH_{j}\left(0\right)-B\beta_{H}\left(0\right)=A\frac{\mu_{j}}{p}\psi_{H}\left(\mathfrak{z}\right)+B\frac{\mu_{j}}{p}\Psi_{H}\left(\mathfrak{z}\right)+H_{j}\left(0\right)
\]
and hence: 
\begin{equation}
A+\left(H_{j}\left(0\right)-\beta_{H}\left(0\right)\right)B=H_{j}\left(0\right),\textrm{ }\forall j\in\left\{ 0,\ldots,p-1\right\} \label{eq:Seemingly over-determined system}
\end{equation}
Note that this is a system of $j$ linear equations in \emph{two }unknowns
($A$ and $B$). Most such systems are usually overdetermined and
admit no solution. Not here, however.

Since $H$ fixes $0$, $H_{0}\left(0\right)=0$, and the $j=0$ equation
becomes: 
\begin{equation}
A-\beta_{H}\left(0\right)B=0
\end{equation}
Letting $j\in\left\{ 1,\ldots,p-1\right\} $, we plug $A=\beta_{H}\left(0\right)B$
into (\ref{eq:Seemingly over-determined system}) to obtain: 
\begin{align*}
\beta_{H}\left(0\right)B+\left(H_{j}\left(0\right)-\beta_{H}\left(0\right)\right)B & =H_{j}\left(0\right)\\
 & \Updownarrow\\
H_{j}\left(0\right)B & =H_{j}\left(0\right)\\
\left(\textrm{if }H_{j}\left(0\right)\neq0\textrm{ for any }j\in\left\{ 1,\ldots,p-1\right\} \right); & \Updownarrow\\
B & =1
\end{align*}
Hence, our condition on the $H_{j}\left(0\right)$s guarantees that
(\ref{eq:Seemingly over-determined system}) has $A=\beta_{H}\left(0\right)$
and $B=1$ as its unique solution. Using (\ref{eq:Fourier Transform of Little Psi_H})
and (\ref{eq:Fourier Transform of Big Psi_H}) for the Fourier transforms
of $\psi_{H}$ and $\Psi_{H}$ yields (\ref{eq:Chi_H in terms of Little Psi_H and Big Psi_H for arbitrary rho}).

Q.E.D. 
\begin{cor}
\label{cor:Chi_H,N twiddle explicit formula, arbitrary p, arbitrary alpha}Let
$H$ be as given in \textbf{\emph{Theorem \ref{thm:F-series for an arbitrary 1D Chi_H}}},
with $p\geq2$. Then:

\begin{align}
\tilde{\chi}_{H,N}\left(\mathfrak{z}\right) & \overset{\overline{\mathbb{Q}}}{=}-\beta_{H}\left(0\right)N\left(\frac{\mu_{0}}{p}\right)^{N}\kappa_{H}\left(\left[\mathfrak{z}\right]_{p^{N}}\right)+\sum_{n=0}^{N-1}\left(\sum_{j=0}^{p-1}\beta_{H}\left(\frac{j}{p}\right)\varepsilon_{n}^{j}\left(\mathfrak{z}\right)\right)\left(\frac{\mu_{0}}{p}\right)^{n}\kappa_{H}\left(\left[\mathfrak{z}\right]_{p^{n}}\right)\label{eq:Explicit Formula for Chi_H,N twiddle for arbitrary rho and alpha equals 1}
\end{align}
if $\alpha_{H}\left(0\right)=1$, and: 
\begin{align}
\tilde{\chi}_{H,N}\left(\mathfrak{z}\right) & \overset{\overline{\mathbb{Q}}}{=}\frac{\beta_{H}\left(0\right)}{1-\alpha_{H}\left(0\right)}\left(\frac{\mu_{0}}{p}\right)^{N}\kappa_{H}\left(\left[\mathfrak{z}\right]_{p^{N}}\right)+\sum_{n=0}^{N-1}\left(\sum_{j=0}^{p-1}\beta_{H}\left(\frac{j}{p}\right)\varepsilon_{n}^{j}\left(\mathfrak{z}\right)\right)\left(\frac{\mu_{0}}{p}\right)^{n}\kappa_{H}\left(\left[\mathfrak{z}\right]_{p^{n}}\right)\label{eq:Explicit Formula for Chi_H,N twiddle for arbitrary rho and alpha not equal to 1}
\end{align}
if $\alpha_{H}\left(0\right)\neq1$. 
\end{cor}
Proof: Take (\ref{eq:Fourier Transform of Chi_H for a contracting semi-basic rho-Hydra map})
and use the formulas from \textbf{Lemmata \ref{lem:1D gamma formula}}
and \textbf{\ref{lem:v_p A_H hat summation formulae}}.

Q.E.D. 
\begin{cor}[\textbf{$\mathcal{F}$-series for $\chi_{H}$}]
\label{cor:F-series for Chi_H, arbitrary p and alpha}Let\index{chi{H}@$\chi_{H}$!mathcal{F}-series@$\mathcal{F}$-series}
$H$ be as given in \textbf{\emph{Theorem \ref{thm:F-series for an arbitrary 1D Chi_H}}}.
Then: 
\begin{equation}
\chi_{H}\left(\mathfrak{z}\right)\overset{\mathcal{F}_{p,q_{H}}}{=}\sum_{n=0}^{\infty}\left(\sum_{j=0}^{p-1}\beta_{H}\left(\frac{j}{p}\right)\left(\varepsilon_{n}\left(\mathfrak{z}\right)\right)^{j}\right)\left(\frac{\mu_{0}}{p}\right)^{n}\kappa_{H}\left(\left[\mathfrak{z}\right]_{p^{n}}\right)\label{eq:Chi_H non-trivial series formula, rho not equal to 2}
\end{equation}
\end{cor}
Proof: Take the $\mathcal{F}_{p,q_{H}}$-limit of the formulae in
\textbf{Corollary \ref{cor:Chi_H,N twiddle explicit formula, arbitrary p, arbitrary alpha}}.

Q.E.D.

\vphantom{}

Finally, using the Wiener Tauberian Theorem for $\left(p,q\right)$-adic
measures (\textbf{Theorem \ref{thm:pq WTT for measures}}), we can
prove the third chief result of this dissertation: the \textbf{Tauberian
Spectral Theorem }for $p$-Hydra maps. 
\begin{thm}[\textbf{Tauberian Spectral Theorem for $p$-Hydra Maps}]
\index{Hydra map!periodic points}\index{Hydra map!Tauberian Spectral Theorem}\label{thm:Periodic Points using WTT}\index{Wiener!Tauberian Theorem!left(p,qright)-adic@$\left(p,q\right)$-adic}\index{Hydra map!divergent trajectories}Let
$H$ be as given in \textbf{\emph{Theorem \ref{thm:F-series for an arbitrary 1D Chi_H}}}.
Let $x\in\mathbb{Z}\backslash\left\{ 0\right\} $, and let:

\begin{equation}
\hat{\chi}_{H}\left(t\right)=\begin{cases}
\begin{cases}
0 & \textrm{if }t=0\\
\left(\beta_{H}\left(0\right)v_{p}\left(t\right)+\gamma_{H}\left(\frac{t\left|t\right|_{p}}{p}\right)\right)\hat{A}_{H}\left(t\right) & \textrm{if }t\neq0
\end{cases} & \textrm{if }\alpha_{H}\left(0\right)=1\\
\frac{\beta_{H}\left(0\right)\hat{A}_{H}\left(t\right)}{1-\alpha_{H}\left(0\right)}+\begin{cases}
0 & \textrm{if }t=0\\
\gamma_{H}\left(\frac{t\left|t\right|_{p}}{p}\right)\hat{A}_{H}\left(t\right) & \textrm{if }t\neq0
\end{cases} & \textrm{if }\alpha_{H}\left(0\right)\neq1
\end{cases},\textrm{ }\forall t\in\hat{\mathbb{Z}}_{p}
\end{equation}

\vphantom{}

I. If $x$ is a periodic point of $H$, then the translates of the
function $\hat{\chi}_{H}\left(t\right)-x\mathbf{1}_{0}\left(t\right)$
are \emph{not }dense in $c_{0}\left(\hat{\mathbb{Z}}_{p},\mathbb{C}_{q_{H}}\right)$.

\vphantom{}

II. Suppose in addition that $H$ is integral\footnote{Propriety is defined alongside Hydra maps themselves on \pageref{def:p-Hydra map};
all the shortened $qx+1$ maps are integral.}, and that $\left|H_{j}\left(0\right)\right|_{q_{H}}=1$ for all $j\in\left\{ 1,\ldots,p-1\right\} $.
If the translates of the function $\hat{\chi}_{H}\left(t\right)-x\mathbf{1}_{0}\left(t\right)$
are \emph{not }dense in $c_{0}\left(\hat{\mathbb{Z}}_{p},\mathbb{C}_{q_{H}}\right)$,
then either $x$ is a periodic point of $H$ or $x$ belongs to a
divergent trajectory of $H$.

When $\alpha_{H}\left(0\right)=1$ and $p=2$, we can also work with:
\begin{equation}
\hat{\chi}_{H}\left(t\right)=\begin{cases}
-\gamma_{H}\left(\frac{1}{2}\right) & \textrm{if }t=0\\
\beta_{H}\left(0\right)v_{2}\left(t\right)\hat{A}_{H}\left(t\right) & \textrm{else }
\end{cases},\textrm{ }\forall t\in\hat{\mathbb{Z}}_{2}
\end{equation}
\end{thm}
Proof: Let $H$ and $\hat{\chi}_{H}$ be as given. In particular,
since $H$ is semi-basic, $\chi_{H}$ exists and the \index{Correspondence Principle}
\textbf{Correspondence Principle }applies. Moreover, since $H$ is
contracting and $H_{j}\left(0\right)\neq0$ for any $j\in\left\{ 1,\ldots,p-1\right\} $,
$\chi_{H}$ is quasi-integrable with respect to the standard $\left(p,q\right)$-adic
frame, and $\hat{\chi}_{H}$ is a Fourier transform of $\chi_{H}$.
Note that the alternative formula for $\hat{\chi}_{H}$ when $p=2$
and $\alpha_{H}\left(0\right)=1$ follows from the fact that it and
the first formula differ by the Fourier-Stieltjes transform of a degenerate
measure, as was shown in \textbf{Corollary \ref{cor:Quasi-integrability of Chi_H for alpha equals 1}}.

Now, let $x\in\mathbb{Z}\backslash\left\{ 0\right\} $. Since the
constant function $x$ has $x\mathbf{1}_{0}\left(t\right)$ as its
Fourier transform, the quasi-integrability of $\chi_{H}$ tells us
that the difference $\chi_{H}\left(\mathfrak{z}\right)-x$ is quasi-integrable,
and that $\hat{\chi}_{H}\left(t\right)-x\mathbf{1}_{0}\left(t\right)$
is then a Fourier transform of $\chi_{H}\left(\mathfrak{z}\right)-x$.
So, $\hat{\chi}_{H}\left(t\right)-x\mathbf{1}_{0}\left(t\right)$
is then the Fourier-Stieltjes transform of a $\left(p,q\right)$-adic
measure, and the Fourier series it generates converges in $\mathbb{C}_{q_{H}}$
(to $\chi_{H}\left(\mathfrak{z}\right)-x$) if and only if $\mathfrak{z}\in\mathbb{Z}_{p}^{\prime}$.
As such, the \textbf{Wiener Tauberian Theorem for $\left(p,q\right)$-adic
measures} (\textbf{Theorem \ref{thm:pq WTT for measures}}) guarantees
that the translates of $\hat{\chi}_{H}\left(t\right)-x\mathbf{1}_{0}\left(t\right)$
will be dense in $c_{0}\left(\hat{\mathbb{Z}}_{p},\mathbb{C}_{q_{H}}\right)$
if and only if $\chi_{H}\left(\mathfrak{z}\right)-x\neq0$ for every
$\mathfrak{z}\in\mathbb{Z}_{p}^{\prime}$. So:

\vphantom{}

I. Suppose $x$ is a periodic point of $H$. By the \textbf{Correspondence
Principle} (specifically, \textbf{Corollary \ref{cor:CP v4}}), there
exists a $\mathfrak{z}_{0}\in\mathbb{Q}\cap\mathbb{Z}_{p}^{\prime}$
so that $\chi_{H}\left(\mathfrak{z}_{0}\right)=x$. Then, the WTT
tells us that the translates of $\hat{\chi}_{H}\left(t\right)-x\mathbf{1}_{0}\left(t\right)$
are \emph{not }dense in $c_{0}\left(\hat{\mathbb{Z}}_{p},\mathbb{C}_{q}\right)$.

\vphantom{}

II. Let $H$ be integral and let $\left|H_{j}\left(0\right)\right|_{q_{H}}=1$
for all $j\in\left\{ 1,\ldots,p-1\right\} $. Since $H$ is integral,
this implies $H$ is proper (\textbf{Lemma \ref{lem:integrality lemma}}
from page \pageref{lem:integrality lemma}). Next, suppose the translates
of $\hat{\chi}_{H}\left(t\right)-x\mathbf{1}_{0}\left(t\right)$ are
not dense in $c_{0}\left(\hat{\mathbb{Z}}_{p},\mathbb{C}_{q}\right)$.
Then, there is a $\mathfrak{z}_{0}\in\mathbb{Z}_{p}^{\prime}$ so
that $\chi_{H}\left(\mathfrak{z}_{0}\right)=x$. If the $p$-adic
digits of $\mathfrak{z}_{0}$ are eventually periodic (i.e., if $\mathfrak{z}_{0}\in\mathbb{Q}\cap\mathbb{Z}_{p}^{\prime}$),
since $H$ is proper, \textbf{Corollary \ref{cor:CP v4}} guarantees
that $x$ is a periodic point of $H$. If the $p$-adic digits of
$\mathfrak{z}_{0}$ are not eventually periodic ($\mathfrak{z}_{0}\in\mathbb{Z}_{p}\backslash\mathbb{Q}$),
\textbf{Theorem \ref{thm:Divergent trajectories come from irrational z}}\textemdash applicable
because of the hypotheses placed on $H$\textemdash then guarantees
that $x$ belongs to a divergent trajectory of $H$.

Q.E.D.

\vphantom{}

Finally, for the reader's benefit, here are the formulae associated
to $\chi_{3}$\index{chi{3}@$\chi_{3}$!mathcal{F}-series@$\mathcal{F}$-series}.
\index{chi{3}@$\chi_{3}$}We have its $\mathcal{F}$-series: 
\begin{equation}
\chi_{3}\left(\mathfrak{z}\right)\overset{\mathcal{F}_{2,3}}{=}-\frac{1}{2}+\frac{1}{4}\sum_{k=0}^{\infty}\frac{3^{\#_{1}\left(\left[\mathfrak{z}\right]_{2^{k}}\right)}}{2^{k}},\textrm{ }\forall\mathfrak{z}\in\mathbb{Z}_{2}\label{eq:Explicit formula for Chi_3}
\end{equation}
Plugging $n\in\mathbb{N}_{0}$ in for $\mathfrak{z}$ and summing
the resultant geometric series in $\mathbb{R}$ yields a formula for
$\chi_{3}$ on the integers: 
\begin{equation}
\chi_{3}\left(n\right)=-\frac{1}{2}+\frac{1}{4}\frac{3^{\#_{1}\left(n\right)}}{2^{\lambda_{2}\left(n\right)}}+\frac{1}{4}\sum_{k=0}^{\lambda_{2}\left(n\right)}\frac{3^{\#_{1}\left(\left[n\right]_{2^{k}}\right)}}{2^{k}},\textrm{ }\forall n\in\mathbb{N}_{0}\label{eq:Explicit Formula for Chi_3 on the integers}
\end{equation}
A nice choice of Fourier transform for $\chi_{3}$ is\index{chi{3}@$\chi_{3}$!Fourier transform}:
\begin{equation}
\hat{\chi}_{3}\left(t\right)=\begin{cases}
-\frac{1}{2} & \textrm{if }t=0\\
\frac{1}{4}v_{2}\left(t\right)\hat{A}_{3}\left(t\right) & \textrm{if }t\neq0
\end{cases},\textrm{ }\forall t\in\hat{\mathbb{Z}}_{2}\label{eq:Formula for Chi_3 hat}
\end{equation}
Since it is the most natural case for comparison, here are corresponding
formulae for $\chi_{5}$:\index{chi{5}@$\chi_{5}$} 
\begin{equation}
\hat{\chi}_{5}\left(t\right)=\begin{cases}
-\frac{1}{2} & \textrm{if }t=0\\
-\frac{1}{4}\hat{A}_{5}\left(t\right) & \textrm{if }t\neq0
\end{cases}\label{eq:Chi_5 Fourier transform}
\end{equation}
\index{chi{5}@$\chi_{5}$!Fourier transform}\index{chi{5}@$\chi_{5}$!mathcal{F}-series@$\mathcal{F}$-series}
\begin{equation}
\chi_{5}\left(\mathfrak{z}\right)\overset{\mathcal{F}_{2,5}}{=}-\frac{1}{4}+\frac{1}{8}\sum_{k=0}^{\infty}\frac{5^{\#_{1}\left(\left[\mathfrak{z}\right]_{2^{k}}\right)}}{2^{k}}\label{eq:Chi_5 F-series}
\end{equation}

\begin{equation}
\chi_{5}\left(n\right)=-\frac{1}{4}+\frac{1}{4}\frac{5^{\#_{1}\left(n\right)}}{2^{\lambda_{2}\left(n\right)}}+\frac{1}{8}\sum_{k=0}^{\lambda_{2}\left(n\right)}\frac{5^{\#_{1}\left(\left[n\right]_{2^{k}}\right)}}{2^{k}},\textrm{ }\forall n\in\mathbb{N}_{0}\label{eq:Explicit Formula for Chi_5 on the integers}
\end{equation}

\index{chi{q}@$\chi_{q}$}\index{chi{q}@$\chi_{q}$!Fourier transform}\index{chi{q}@$\chi_{q}$!mathcal{F}-series@$\mathcal{F}$-series}
More generally, for any odd prime $q\geq5$, letting $\chi_{q}$ be
the numen of the Shortened $qx+1$ map, we have that:
\begin{equation}
\hat{\chi}_{q}\left(t\right)=\begin{cases}
-\frac{1}{q-3} & \textrm{if }t=0\\
-\frac{2\hat{A}_{q}\left(t\right)}{\left(q-1\right)\left(q-3\right)} & \textrm{if }t\neq0
\end{cases},\textrm{ }\forall t\in\hat{\mathbb{Z}}_{2}\label{eq:Fourier transform of the numen of the shortened qx+1 map}
\end{equation}
with the $\mathcal{F}$-series:
\begin{equation}
\chi_{q}\left(\mathfrak{z}\right)\overset{\mathcal{F}_{2,q}}{=}-\frac{1}{q-1}+\frac{1}{2\left(q-1\right)}\sum_{k=0}^{\infty}\frac{q^{\#_{1}\left(\left[\mathfrak{z}\right]_{2^{k}}\right)}}{2^{k}},\textrm{ }\forall\mathfrak{z}\in\mathbb{Z}_{2}\label{eq:F-series for Chi_q}
\end{equation}
and:
\begin{equation}
\chi_{q}\left(n\right)\overset{\mathbb{Q}}{=}-\frac{1}{q-1}+\frac{1}{2\left(q-1\right)}\frac{q^{\#_{1}\left(n\right)}}{2^{\lambda_{2}\left(n\right)}}+\frac{1}{2\left(q-1\right)}\sum_{k=0}^{\lambda_{2}\left(n\right)}\frac{q^{\#_{1}\left(\left[n\right]_{2^{k}}\right)}}{2^{k}},\textrm{ }\forall n\in\mathbb{N}_{0}\label{eq:Chi_q on N_0}
\end{equation}
Even though (\ref{eq:Fourier transform of the numen of the shortened qx+1 map})
is not the correct formula for $\hat{\chi}_{q}$ when $q=3$, it seems
impossible that it is merely coincidental that the right-hand side
of (\ref{eq:Fourier transform of the numen of the shortened qx+1 map})
is singular when $q\in\left\{ 1,3\right\} $. When $q=1$, it is easy
to show that the shortened $qx+1$ map iterates every positive integer
to $1$. Meanwhile, when $q=3$, it is conjectured that every positive
integer goes to $1$. This suggests that for a given family $\mathcal{H}$
of Hydra maps with quasi-integrable numina $\hat{\chi}_{H}$ for each
$H\in\mathcal{H}$, we have something along the lines of ``a given
map $T\in\mathcal{H}$ has finitely many orbit classes (possibly all
of which contain cycles) in $\mathbb{Z}$ if and only if the map $H\mapsto\hat{\chi}_{H}$
is 'singular' at $T$.''

\newpage{}

\section{\label{sec:4.3 Salmagundi}Salmagundi}

FOR THIS SUBSECTION, UNLESS STATED OTHERWISE, $H$ IS ALSO ASSUMED
TO BE INTEGRAL AND NON-SINGULAR.

\vphantom{}

The title of this section is a wonderful word which my dictionary
informs me, means either ``a dish of chopped meat, anchovies, eggs,
onions, and seasoning'', or ``a general mixture; a miscellaneous
collection''. It is the latter meaning that we are going to invoke
here. This section contains a general mixture of results on $\chi_{H}$.
Subsection \ref{subsec:4.3.1 p-adic-Estimates} demonstrates how the
$\mathcal{F}$-series formula for $\chi_{H}$ (\textbf{Corollaries
\ref{cor:F-series for Chi_H for p equals 2}} and \textbf{\ref{cor:F-series for Chi_H, arbitrary p and alpha}})
can be used to obtain a crude lower bound on the $q_{H}$-adic absolute
value of $\chi_{H}\left(\mathfrak{z}\right)-\lambda$. It would be
interesting to see if (and, if, then \emph{how}) this approach to
estimation could be refined. \ref{subsec:4.3.2 A-Wisp-of} takes us
back to the content of Subsection \ref{subsec:3.3.6 L^1 Convergence}\textemdash $L^{1}$
integrability of $\left(p,q\right)$-adic functions\textemdash and
shows that the real-valued function $\left|\chi_{H}\left(\mathfrak{z}\right)\right|_{q_{H}}$
will be integrable with respect to the real-valued Haar probability
measure on $\mathbb{Z}_{p}$ whenever $\gcd\left(\mu_{j},q_{H}\right)>1$
for at least one $j\in\left\{ 0,\ldots,p-1\right\} $ (\textbf{Theorem
\ref{thm:L^1 criterion for Chi_H}}). I also show that: 
\[
\lim_{N\rightarrow\infty}\int_{\mathbb{Z}_{p}}\left|\chi_{H,N}\left(\mathfrak{z}\right)-\tilde{\chi}_{H,N}\left(\mathfrak{z}\right)\right|_{q_{H}}d\mathfrak{z}\overset{\mathbb{R}}{=}0
\]
whenever $H$ is semi-basic and contracting (\textbf{Theorem \ref{thm:L^1 convergence of Chi_H,N minus Chi_H,N twiddle}}).
In Subsection \ref{subsec:4.3.3 Quick-Approach-of}, we revisit the
circle of ideas around the \textbf{Square Root Lemma }(page \pageref{lem:square root lemma}).
Finally, Subsection \ref{subsec:4.3.4 Archimedean-Estimates} establishes
what I have termed the ``$L^{1}$-method'', a means of using geometric
series universality to straddle the chasm between archimedean and
non-archimedean topologies so as to obtain upper bounds on the ordinary,
archimedean absolute value of $\chi_{H}\left(\mathfrak{z}\right)$
for appropriately chosen values of $\mathfrak{z}$ at which $\chi_{H}\left(\left[\mathfrak{z}\right]_{p^{N}}\right)$
converges to $\chi_{H}\left(\mathfrak{z}\right)$ in $\mathbb{Z}_{q_{H}}$
\emph{and }$\mathbb{R}$ as $N\rightarrow\infty$ (\textbf{Theorem
\ref{thm:The L^1 method}}). I also sketch a sample application of
this method to $\chi_{3}$.

\subsection{\label{subsec:4.3.1 p-adic-Estimates}$p$-adic Estimates}

While it is my belief that the most fruitful investigations into $\chi_{H}$
will come from studying $\hat{\chi}_{H}$ in the context of the $\left(p,q\right)$-adic
Wiener Tauberian Theorem\textemdash assuming, of course, that a proper
converse for the case of $\left(p,q\right)$-adic measures / quasi-integrable
functions can be formulated and proven\textemdash the explicit series
formulae from \textbf{Corollaries \ref{cor:F-series for Chi_H for p equals 2}}
and \textbf{\ref{cor:F-series for Chi_H, arbitrary p and alpha}}
may be of use in their own right. With sufficient knowledge of the
$p$-adic digits of a given $\mathfrak{z}$\textemdash equivalently,
with sufficient knowledge on the pattern of applications of branches
of $H$ needed to ensure that a trajectory cycles\textemdash one can
use the explicit series formulae to obtain $\left(p,q\right)$-adic
estimates for $\chi_{H}$ which can be used to rule out periodic points.
The more we know about $\mathfrak{z}$, the more refined the estimates
can be. In that respect, even though the estimates detailed below
are relatively simple, the method may be of use, so I might as well
give the details.

First, two notations to make our life easier: 
\begin{defn}[$u_{p}\left(\mathfrak{z}\right)$ \textbf{and} $K_{H}$]
\ 

\vphantom{}

I. We write \nomenclature{$u_{p}\left(\mathfrak{z}\right)$}{$\mathfrak{z}\left|\mathfrak{z}\right|_{p}$ }$u_{p}:\mathbb{Z}_{p}\rightarrow\mathbb{Z}_{p}^{\times}\cup\left\{ 0\right\} $
to denote the function: 
\begin{equation}
u_{p}\left(\mathfrak{z}\right)\overset{\textrm{def}}{=}\mathfrak{z}\left|\mathfrak{z}\right|_{p}\label{eq:Definition of u_p}
\end{equation}

\vphantom{}

II. Define:\nomenclature{$K_{H}$}{  } 
\begin{equation}
K_{H}\overset{\textrm{def}}{=}\max_{j\in\left\{ 1,\ldots,p-1\right\} }\left|\kappa_{H}\left(j\right)\right|_{q_{H}}=\max_{j\in\left\{ 1,\ldots,p-1\right\} }\left|\frac{\mu_{j}}{\mu_{0}}\right|_{q_{H}}\label{eq:Definition of K_H}
\end{equation}
Note that since $H$ is semi-basic, $\mu_{j}/\mu_{0}$ is a $q_{H}$-adic
integer with $q_{H}$-adic absolute value $\leq1/q_{H}$. 
\end{defn}
\begin{prop}
We have: 
\begin{equation}
\left|\kappa_{H}\left(n\right)\right|_{q_{H}}\leq K_{H}\leq\frac{1}{q_{H}}<1,\textrm{ }\forall n\in\mathbb{N}_{1}\label{eq:Kappa_H K_H bound}
\end{equation}
\end{prop}
Proof: Use (\ref{eq:Kappa_H is rho-adically distributed}).

Q.E.D.

\vphantom{}

Using the explicit formulae from \textbf{Proposition \ref{prop:formula for functions with p-adic structure}}
we can now establish the promised estimates. There isn't much thought
behind these estimates; it's mostly a matter of computation. So, to
keep things from getting too unwieldy, I think it appropriate to show
one concrete example\textemdash say, for $\chi_{3}$\textemdash before
diving into the details. 
\begin{example}
Using \textbf{Corollary \ref{cor:F-series for Chi_H for p equals 2}},
for the case of $\chi_{3}$, we have that: 
\begin{equation}
\tilde{\chi}_{3,N}\left(\mathfrak{z}\right)\overset{\overline{\mathbb{Q}}}{=}-\frac{1}{2}+\frac{N}{4}\frac{3^{\#_{1}\left(\left[\mathfrak{z}\right]_{2^{N}}\right)}}{2^{N}}+\frac{1}{4}\sum_{n=0}^{N-1}\frac{3^{\#_{1}\left(\left[\mathfrak{z}\right]_{2^{n}}\right)}}{2^{n}},\textrm{ }\forall N\in\mathbb{N}_{1},\textrm{ }\forall\mathfrak{z}\in\mathbb{Z}_{2}\label{eq:Explicit formula for Chi_3,N twiddle}
\end{equation}
Now, observe that $\left[\mathfrak{z}\right]_{2^{n}}=0$ for all $n\leq v_{2}\left(\mathfrak{z}\right)$;
ex. $\left[8\right]_{2}=\left[8\right]_{4}=\left[8\right]_{8}=0$.
As such, $n=v_{2}\left(\mathfrak{z}\right)+1$ is the smallest non-negative
integer for which $\left[\mathfrak{z}\right]_{2^{n}}>0$, and hence,
is the smallest integer for which $\kappa_{3}\left(\left[\mathfrak{z}\right]_{2^{n}}\right)=3^{\#_{1}\left(\left[\mathfrak{z}\right]_{2^{n}}\right)}$
is $>1$. Indeed, $\kappa_{3}\left(\left[\mathfrak{z}\right]_{2^{v_{2}\left(\mathfrak{z}\right)+1}}\right)=3$
for all $\mathfrak{z}\in\mathbb{Z}_{2}\backslash\left\{ 0\right\} $.
So, $3^{\#_{1}\left(\left[\mathfrak{z}\right]_{2^{n}}\right)}$ will
en be a multiple of $3$ for all $n\geq v_{2}\left(\mathfrak{z}\right)+1$.
This tells us we can pull out terms from the $n$-sum which have a
$3$-adic magnitude of $1$: 
\begin{align*}
\sum_{n=0}^{N-1}\frac{3^{\#_{1}\left(\left[\mathfrak{z}\right]_{2^{n}}\right)}}{2^{n}} & =\sum_{n=0}^{v_{2}\left(\mathfrak{z}\right)}\frac{1}{2^{n}}+\sum_{n=v_{2}\left(\mathfrak{z}\right)+1}^{N-1}\frac{3^{\#_{1}\left(\left[\mathfrak{z}\right]_{2^{n}}\right)}}{2^{n}}\\
\left(3^{\#_{1}\left(\left[\mathfrak{z}\right]_{2^{n}}\right)}=1\textrm{ }\forall n\leq v_{2}\left(\mathfrak{z}\right)\right); & =\sum_{n=0}^{v_{2}\left(\mathfrak{z}\right)}\frac{1}{2^{n}}+\sum_{n=v_{2}\left(\mathfrak{z}\right)+1}^{N-1}\frac{3^{\#_{1}\left(\left[\mathfrak{z}\right]_{2^{n}}\right)}}{2^{n}}\\
 & =\frac{1-\frac{1}{2^{v_{2}\left(\mathfrak{z}\right)+1}}}{1-\frac{1}{2}}+\sum_{n=v_{2}\left(\mathfrak{z}\right)+1}^{N-1}\frac{3^{\#_{1}\left(\left[\mathfrak{z}\right]_{2^{n}}\right)}}{2^{n}}\\
\left(\frac{1}{2^{v_{2}\left(\mathfrak{z}\right)}}=\left|\mathfrak{z}\right|_{2}\right); & =2-\left|\mathfrak{z}\right|_{2}+\sum_{n=v_{2}\left(\mathfrak{z}\right)+1}^{N-1}\frac{3^{\#_{1}\left(\left[\mathfrak{z}\right]_{2^{n}}\right)}}{2^{n}}
\end{align*}
which, for any $\mathfrak{z}\in\mathbb{Z}_{2}\backslash\left\{ 0\right\} $,
holds for all $N-1\geq v_{2}\left(\mathfrak{z}\right)+1$. So, (\ref{eq:Explicit formula for Chi_3,N twiddle})
becomes: 
\begin{equation}
\tilde{\chi}_{3,N}\left(\mathfrak{z}\right)\overset{\overline{\mathbb{Q}}}{=}\frac{N}{4}\frac{3^{\#_{1}\left(\left[\mathfrak{z}\right]_{2^{N}}\right)}}{2^{N}}-\frac{1}{4}\left|\mathfrak{z}\right|_{2}+\frac{1}{4}\sum_{n=v_{2}\left(\mathfrak{z}\right)+1}^{N-1}\frac{3^{\#_{1}\left(\left[\mathfrak{z}\right]_{2^{n}}\right)}}{2^{n}}\label{eq:Explicit formula for Chi_3,N twiddle with z pulled out}
\end{equation}

In terms of strings and composition sequence of trajectories, if $\mathfrak{z}\in\mathbb{Z}_{2}^{\prime}\cap\mathbb{Q}$
made $\chi_{3}\left(\mathfrak{z}\right)$ into a periodic point (belonging
to the cycle $\Omega$) of the shortened Collatz map $T_{3}$ , then
$v_{2}\left(\mathfrak{z}\right)$ would be the number of times we
had to divide by $2$ to get from the penultimate element of $\Omega$
and return to $\chi_{3}\left(\mathfrak{z}\right)$. In this light,
(\ref{eq:Explicit formula for Chi_3,N twiddle with z pulled out})
shows how we can obtain greater detail by using our knowledge of the
number of times we divide by $2$ on the final ``fall'' in a given
cycle. Channeling the Correspondence Principle, let $\mathfrak{z}\in\mathbb{Q}\cap\mathbb{Z}_{2}^{\prime}$.
Then, (\ref{eq:Explicit formula for Chi_3,N twiddle with z pulled out})
converges in $\mathbb{Z}_{3}$ to $\chi_{3}\left(\mathfrak{z}\right)$
as $N\rightarrow\infty$, and we obtain: 
\begin{equation}
\chi_{3}\left(\mathfrak{z}\right)\overset{\mathbb{Z}_{3}}{=}-\frac{1}{4}\left|\mathfrak{z}\right|_{2}+\frac{1}{4}\sum_{n=v_{2}\left(\mathfrak{z}\right)+1}^{\infty}\frac{3^{\#_{1}\left(\left[\mathfrak{z}\right]_{2^{n}}\right)}}{2^{n}},\textrm{ }\forall\mathfrak{z}\in\mathbb{Z}_{2}^{\prime}
\end{equation}
If we knew that the first $1$ in the sequence of $\mathfrak{z}$'s
$2$-adic digits was followed by $m$ consecutive zeroes, we would
then be able to pull out all of the terms of the $n$-sum of $3$-adic
magnitude $1/3$, and so on and so forth\textemdash but doing already
takes us beyond the simple estimate I had in mind. \emph{That }(admittedly
crude) estimate consists of taking the above, subtracting $\lambda\in\mathbb{Z}$,
and then applying the $3$-adic absolute value: 
\[
\left|\chi_{3}\left(\mathfrak{z}\right)-\lambda\right|_{3}\overset{\mathbb{Z}_{3}}{=}\left|\lambda+\frac{1}{4}\left|\mathfrak{z}\right|_{2}-\frac{1}{4}\underbrace{\sum_{n=v_{2}\left(\mathfrak{z}\right)+1}^{\infty}\frac{3^{\#_{1}\left(\left[\mathfrak{z}\right]_{2^{n}}\right)}}{2^{n}}}_{\textrm{has }\left|\cdot\right|_{3}\leq1/3}\right|_{3}
\]
By the $3$-adic ultrametric inequality, observe that $\left|\lambda+\frac{1}{4}\left|\mathfrak{z}\right|_{2}\right|_{3}=1$
then forces: 
\[
\left|\chi_{3}\left(\mathfrak{z}\right)-\lambda\right|_{3}=\left|\lambda+\frac{1}{4}\left|\mathfrak{z}\right|_{2}\right|_{3}=1
\]
So, we have that $\chi_{3}\left(\mathfrak{z}\right)\neq\lambda$ for
any $\lambda\in\mathbb{Z}$ and $\mathfrak{z}\in\mathbb{Z}_{2}^{\prime}$
such that $\left|\lambda+\frac{1}{4}\left|\mathfrak{z}\right|_{2}\right|_{3}=1$.

For example, if $\left|\mathfrak{z}\right|_{2}=2^{-n}$, then $\left|2^{n+2}\lambda+1\right|_{3}=1$
implies $\chi_{3}\left(\mathfrak{z}\right)\neq\lambda$. Hence, $2^{n+2}\lambda+1\overset{3}{\equiv}0$
is necessary in order for $\lambda$ to be a non-zero periodic point
of Collatz with $\chi_{3}\left(\mathfrak{z}\right)=\lambda$ where
$\left|\mathfrak{z}\right|_{2}=2^{-n}$. 
\end{example}
\vphantom{}

The following Corollary of the explicit formulae from Section \ref{sec:4.2 Fourier-Transforms-=00003D000026}
gives the above result for all of the kinds of Hydra maps we considered
in that section. 
\begin{cor}
\label{cor:Non-archimedean estimates of Chi_H}Let $\lambda\in\mathbb{Z}_{q_{H}}$
and $\mathfrak{z}\in\mathbb{Z}_{p}^{\prime}$.\index{chi{H}@$\chi_{H}$!non-archimedean estimates}

Then: 
\begin{equation}
\left|\chi_{H}\left(\mathfrak{z}\right)-\lambda\right|_{q_{H}}=\left|\chi_{H}\left(\left[u_{p}\left(\mathfrak{z}\right)\right]_{p}\right)\left|\mathfrak{z}\right|_{p}^{1-\log_{p}\mu_{0}}-\lambda\right|_{q_{H}}>0\label{eq:Crude Estimate}
\end{equation}
whenever: 
\begin{equation}
\left|\chi_{H}\left(\left[u_{p}\left(\mathfrak{z}\right)\right]_{p}\right)\left|\mathfrak{z}\right|_{p}^{1-\log_{p}\mu_{0}}-\lambda\right|_{q_{H}}>K_{H}\sup_{n>v_{p}\left(\mathfrak{z}\right)}\left|\sum_{j=0}^{p-1}\beta_{H}\left(\frac{j}{p}\right)\varepsilon_{n}^{j}\left(\mathfrak{z}\right)\right|_{q_{H}}\label{eq:Crude Estimate Prerequisite}
\end{equation}
In particular, $\chi_{H}\left(\mathfrak{z}\right)\neq\lambda$ whenever
\emph{(\ref{eq:Crude Estimate Prerequisite})} holds true. 
\end{cor}
Proof: Fix $\mathfrak{z}\in\mathbb{Z}_{p}^{\prime}$ and recall that:
\begin{equation}
\kappa_{H}\left(\left[\mathfrak{z}\right]_{p^{n}}\right)=\kappa_{H}\left(0\right)=1,\textrm{ }\forall n\leq v_{p}\left(\mathfrak{z}\right)
\end{equation}
\begin{equation}
\varepsilon_{n}\left(\mathfrak{z}\right)=e^{2\pi i\left\{ \frac{\mathfrak{z}}{p^{n+1}}\right\} _{p}}e^{-\frac{2\pi i}{p}\left\{ \frac{\mathfrak{z}}{p^{n}}\right\} _{p}}=1,\textrm{ }\forall n<v_{p}\left(\mathfrak{z}\right)
\end{equation}

Now, let $N-1\geq v_{p}\left(\mathfrak{z}\right)+1$. Then, there
are two cases.

\vphantom{}

I. Suppose $\alpha_{H}\left(0\right)=1$. Then, (\ref{eq:Explicit Formula for Chi_H,N twiddle for arbitrary rho and alpha equals 1})
becomes: 
\begin{align*}
\tilde{\chi}_{H,N}\left(\mathfrak{z}\right) & \overset{\overline{\mathbb{Q}}}{=}-\beta_{H}\left(0\right)N\left(\frac{\mu_{0}}{p}\right)^{N}\kappa_{H}\left(\left[\mathfrak{z}\right]_{p^{N}}\right)+\beta_{H}\left(0\right)\sum_{n=0}^{v_{p}\left(\mathfrak{z}\right)}\left(\frac{\mu_{0}}{p}\right)^{n}\kappa_{H}\left(0\right)\\
 & +\sum_{n=0}^{v_{p}\left(\mathfrak{z}\right)-1}\left(\sum_{j=1}^{p-1}\beta_{H}\left(\frac{j}{p}\right)\left(1\right)^{j}\right)\left(\frac{\mu_{0}}{p}\right)^{n}\kappa_{H}\left(0\right)\\
 & +\left(\sum_{j=1}^{p-1}\beta_{H}\left(\frac{j}{p}\right)\varepsilon_{v_{p}\left(\mathfrak{z}\right)}^{j}\left(\mathfrak{z}\right)\right)\left(\frac{\mu_{0}}{p}\right)^{v_{p}\left(\mathfrak{z}\right)}\kappa_{H}\left(0\right)\\
 & +\beta_{H}\left(0\right)\sum_{n=v_{p}\left(\mathfrak{z}\right)+1}^{N-1}\left(\frac{\mu_{0}}{p}\right)^{n}\kappa_{H}\left(\left[\mathfrak{z}\right]_{p^{n}}\right)\\
 & +\sum_{n=v_{p}\left(\mathfrak{z}\right)+1}^{N-1}\left(\sum_{j=1}^{p-1}\beta_{H}\left(\frac{j}{p}\right)\varepsilon_{n}^{j}\left(\mathfrak{z}\right)\right)\left(\frac{\mu_{0}}{p}\right)^{n}\kappa_{H}\left(\left[\mathfrak{z}\right]_{p^{n}}\right)
\end{align*}
Here: 
\begin{align*}
\sum_{j=0}^{p-1}\beta_{H}\left(\frac{j}{p}\right) & =\sum_{j=0}^{p-1}\left(\frac{1}{p}\sum_{k=0}^{p-1}H_{k}\left(0\right)e^{-\frac{2\pi ijk}{p}}\right)\\
 & =\sum_{k=0}^{p-1}H_{k}\left(0\right)\frac{1}{p}\sum_{j=0}^{p-1}e^{-\frac{2\pi ijk}{p}}\\
 & =\sum_{k=0}^{p-1}H_{k}\left(0\right)\left[k\overset{p}{\equiv}0\right]\\
 & =H_{0}\left(0\right)\\
 & =0
\end{align*}
So: 
\begin{equation}
\sum_{j=1}^{p-1}\beta_{H}\left(\frac{j}{p}\right)=-\beta_{H}\left(0\right)+\underbrace{\sum_{j=0}^{p-1}\beta_{H}\left(\frac{j}{p}\right)}_{0}=-\beta_{H}\left(0\right)\label{eq:Sum of beta H of j/p for j in 1,...,p-1}
\end{equation}
On the other hand: 
\begin{align*}
\sum_{j=0}^{p-1}\beta_{H}\left(\frac{j}{p}\right)\varepsilon_{v_{p}\left(\mathfrak{z}\right)}^{j}\left(\mathfrak{z}\right) & =\sum_{j=0}^{p-1}\left(\frac{1}{p}\sum_{k=0}^{p-1}H_{k}\left(0\right)e^{-\frac{2\pi ijk}{p}}\right)\left(e^{2\pi i\left\{ \frac{\mathfrak{z}}{p^{v_{p}\left(\mathfrak{z}\right)+1}}\right\} _{p}-\frac{2\pi i}{p}\left\{ \frac{\mathfrak{z}}{p^{v_{p}\left(\mathfrak{z}\right)}}\right\} _{p}}\right)^{j}\\
 & =\sum_{j=0}^{p-1}\left(\frac{1}{p}\sum_{k=0}^{p-1}H_{k}\left(0\right)e^{-\frac{2\pi ijk}{p}}\right)e^{2\pi i\frac{j\left[\left|\mathfrak{z}\right|_{p}\mathfrak{z}\right]_{p}}{p}}\\
 & =\frac{1}{p}\sum_{k=0}^{p-1}H_{k}\left(0\right)\sum_{j=0}^{p-1}e^{\frac{2\pi ij}{p}\left(\left[u_{p}\left(\mathfrak{z}\right)\right]_{p}-k\right)}\\
 & =\sum_{k=0}^{p-1}H_{k}\left(0\right)\left[\left|\mathfrak{z}\right|_{p}\mathfrak{z}\overset{p}{\equiv}k\right]\\
 & =H_{\left[u_{p}\left(\mathfrak{z}\right)\right]_{p}}\left(0\right)
\end{align*}
and so: 
\begin{equation}
\sum_{j=1}^{p-1}\beta_{H}\left(\frac{j}{p}\right)\varepsilon_{v_{p}\left(\mathfrak{z}\right)}^{j}\left(\mathfrak{z}\right)=-\beta_{H}\left(0\right)+\sum_{j=0}^{p-1}\beta_{H}\left(\frac{j}{p}\right)\varepsilon_{v_{p}\left(\mathfrak{z}\right)}^{j}\left(\mathfrak{z}\right)=H_{\left[u_{p}\left(\mathfrak{z}\right)\right]_{p}}\left(0\right)-\beta_{H}\left(0\right)\label{eq:Sum of beta H of j/p times the epsilon for j in 1,...,p-1}
\end{equation}
With this, $\tilde{\chi}_{H,N}\left(\mathfrak{z}\right)$ becomes:
\begin{align*}
\tilde{\chi}_{H,N}\left(\mathfrak{z}\right) & \overset{\overline{\mathbb{Q}}}{=}H_{\left[u_{p}\left(\mathfrak{z}\right)\right]_{p}}\left(0\right)\left(\frac{\mu_{0}}{p}\right)^{v_{p}\left(\mathfrak{z}\right)}-\beta_{H}\left(0\right)N\left(\frac{\mu_{0}}{p}\right)^{N}\kappa_{H}\left(\left[\mathfrak{z}\right]_{p^{N}}\right)\\
 & +\beta_{H}\left(0\right)\sum_{n=v_{p}\left(\mathfrak{z}\right)+1}^{N-1}\left(\frac{\mu_{0}}{p}\right)^{n}\kappa_{H}\left(\left[\mathfrak{z}\right]_{p^{n}}\right)\\
 & +\sum_{n=v_{p}\left(\mathfrak{z}\right)+1}^{N-1}\left(\sum_{j=1}^{p-1}\beta_{H}\left(\frac{j}{p}\right)\varepsilon_{n}^{j}\left(\mathfrak{z}\right)\right)\left(\frac{\mu_{0}}{p}\right)^{n}\kappa_{H}\left(\left[\mathfrak{z}\right]_{p^{n}}\right)
\end{align*}
Letting $N\rightarrow\infty$, the above converges in $\mathbb{C}_{q_{H}}$
to: 
\begin{align*}
\chi_{H}\left(\mathfrak{z}\right) & \overset{\mathbb{C}_{q_{H}}}{=}H_{\left[u_{p}\left(\mathfrak{z}\right)\right]_{p}}\left(0\right)\left(\frac{\mu_{0}}{p}\right)^{v_{p}\left(\mathfrak{z}\right)}+\beta_{H}\left(0\right)\sum_{n=v_{p}\left(\mathfrak{z}\right)+1}^{\infty}\left(\frac{\mu_{0}}{p}\right)^{n}\kappa_{H}\left(\left[\mathfrak{z}\right]_{p^{n}}\right)\\
 & +\sum_{n=v_{p}\left(\mathfrak{z}\right)+1}^{\infty}\left(\sum_{j=1}^{p-1}\beta_{H}\left(\frac{j}{p}\right)\varepsilon_{n}^{j}\left(\mathfrak{z}\right)\right)\left(\frac{\mu_{0}}{p}\right)^{n}\kappa_{H}\left(\left[\mathfrak{z}\right]_{p^{n}}\right)
\end{align*}
The $n$-sum on the top line is the $j=0$th term of the sum on the
bottom line. This leaves us with:
\[
\chi_{H}\left(\mathfrak{z}\right)\overset{\mathbb{C}_{q_{H}}}{=}H_{\left[u_{p}\left(\mathfrak{z}\right)\right]_{p}}\left(0\right)\left(\frac{\mu_{0}}{p}\right)^{v_{p}\left(\mathfrak{z}\right)}+\sum_{n=v_{p}\left(\mathfrak{z}\right)+1}^{\infty}\left(\sum_{j=0}^{p-1}\beta_{H}\left(\frac{j}{p}\right)\varepsilon_{n}^{j}\left(\mathfrak{z}\right)\right)\left(\frac{\mu_{0}}{p}\right)^{n}\kappa_{H}\left(\left[\mathfrak{z}\right]_{p^{n}}\right)
\]
Because $\mathfrak{z}\in\mathbb{Z}_{p}^{\prime}$, observe that $\left|\kappa_{H}\left(\left[\mathfrak{z}\right]_{p^{n}}\right)\right|_{q_{H}}\leq K_{H}$
for all $n>v_{p}\left(\mathfrak{z}\right)$. As such, subtracting
$\lambda\in\mathbb{Z}_{q_{H}}$ from both sides, we can apply the
$q_{H}$-adic ultrametric inequality. This gives us the equality:
\begin{equation}
\left|\chi_{H}\left(\mathfrak{z}\right)-\lambda\right|_{q_{H}}=\left|H_{\left[u_{p}\left(\mathfrak{z}\right)\right]_{p}}\left(0\right)\left(\frac{\mu_{0}}{p}\right)^{v_{p}\left(\mathfrak{z}\right)}-\lambda\right|_{q_{H}}
\end{equation}
provided that the condition: 
\begin{equation}
\left|H_{\left[u_{p}\left(\mathfrak{z}\right)\right]_{p}}\left(0\right)\left(\frac{\mu_{0}}{p}\right)^{v_{p}\left(\mathfrak{z}\right)}-\lambda\right|_{q_{H}}>K_{H}\sup_{n>v_{p}\left(\mathfrak{z}\right)}\left|\sum_{j=0}^{p-1}\beta_{H}\left(\frac{j}{p}\right)\varepsilon_{n}^{j}\left(\mathfrak{z}\right)\right|_{q_{H}}
\end{equation}
is satisfied. Finally, writing: 
\begin{equation}
\left(\frac{\mu_{0}}{p}\right)^{v_{p}\left(\mathfrak{z}\right)}=\left(p^{-v_{p}\left(\mathfrak{z}\right)}\right)^{1-\log_{p}\mu_{0}}=\left|\mathfrak{z}\right|_{p}^{1-\log_{p}\mu_{0}}
\end{equation}
we obtain (I) by noting that: 
\begin{equation}
H_{\left[u_{p}\left(\mathfrak{z}\right)\right]_{p}}\left(0\right)=\chi_{H}\left(\left[u_{p}\left(\mathfrak{z}\right)\right]_{p}\right)
\end{equation}

\vphantom{}

II. Suppose $\alpha_{H}\left(0\right)\neq1$. Then (\ref{eq:Explicit Formula for Chi_H,N twiddle for arbitrary rho and alpha not equal to 1})
becomes: 
\begin{align*}
\tilde{\chi}_{H,N}\left(\mathfrak{z}\right) & \overset{\overline{\mathbb{Q}}}{=}\frac{\beta_{H}\left(0\right)}{1-\alpha_{H}\left(0\right)}\left(\frac{\mu_{0}}{p}\right)^{N}\kappa_{H}\left(\left[\mathfrak{z}\right]_{p^{N}}\right)+\beta_{H}\left(0\right)\sum_{n=0}^{N-1}\left(\frac{\mu_{0}}{p}\right)^{n}\kappa_{H}\left(\left[\mathfrak{z}\right]_{p^{n}}\right)\\
 & +\sum_{n=0}^{N-1}\left(\sum_{j=1}^{p-1}\beta_{H}\left(\frac{j}{p}\right)\varepsilon_{n}^{j}\left(\mathfrak{z}\right)\right)\left(\frac{\mu_{0}}{p}\right)^{n}\kappa_{H}\left(\left[\mathfrak{z}\right]_{p^{n}}\right)
\end{align*}
The only term in this equation that differs from the $\alpha_{H}\left(0\right)=1$
case is the initial term: 
\begin{equation}
\frac{\beta_{H}\left(0\right)}{1-\alpha_{H}\left(0\right)}\left(\frac{\mu_{0}}{p}\right)^{N}\kappa_{H}\left(\left[\mathfrak{z}\right]_{p^{N}}\right)
\end{equation}
The corresponding term of the $\alpha_{H}\left(0\right)=1$ case was:
\begin{equation}
-\beta_{H}\left(0\right)N\left(\frac{\mu_{0}}{p}\right)^{N}\kappa_{H}\left(\left[\mathfrak{z}\right]_{p^{N}}\right)
\end{equation}
Since both of these terms converge $q_{H}$-adically to zero as $N\rightarrow\infty$
whenever $\mathfrak{z}\in\mathbb{Z}_{p}^{\prime}$, all of the computations
done for the $\alpha_{H}\left(0\right)=1$ case will apply here\textemdash the
only term that differs between the two cases ends up vanishing as
$N\rightarrow\infty$.

Q.E.D.

\subsection{\label{subsec:4.3.2 A-Wisp-of}A Wisp of $L^{1}$\index{non-archimedean!$L^{1}$ theory}}

One of my goals for future research is to pursue $\left(p,q\right)$-adic
analysis in the setting of $L_{\mathbb{R}}^{1}$. To that end, in
this subsection, a simple sufficient condition is establish for determining
when $\chi_{H}\in L_{\mathbb{R}}^{1}\left(\mathbb{Z}_{p},\mathbb{C}_{q_{H}}\right)$.

We begin by computing the van der Put series for $\chi_{H}$. When
$p=2$, this is easily done. For $p\geq3$, we will have to settle
on a mutated sort of van der Put series, seeing as the van der Put
coefficients of $\chi_{H}$ are not as well-behaved for $p\geq3$
as they are for $p=2$. 
\begin{prop}[\textbf{``van der Put'' series for $\chi_{H}$}]
\index{chi{H}@$\chi_{H}$!van der Put series}\label{prop:van der Put series for Chi_H}When
$p=2$: 
\begin{equation}
\chi_{H}\left(\mathfrak{z}\right)\overset{\mathbb{Z}_{q_{H}}}{=}\frac{b_{1}}{a_{1}}\sum_{n=1}^{\infty}M_{H}\left(n\right)\left[\mathfrak{z}\overset{2^{\lambda_{2}\left(n\right)}}{\equiv}n\right],\textrm{ }\forall\mathfrak{z}\in\mathbb{Z}_{2}\label{eq:Simplified vdP series for Chi_H when rho is 2}
\end{equation}
More generally, for any $p$: 
\begin{equation}
\chi_{H}\left(\mathfrak{z}\right)\overset{\mathbb{Z}_{q_{H}}}{=}\chi_{H}\left(\left[\mathfrak{z}\right]_{p}\right)+\sum_{n=1}^{\infty}\left(\sum_{j=1}^{p-1}H_{j}\left(0\right)\left[u_{p}\left(\mathfrak{z}-n\right)\overset{p}{\equiv}j\right]\right)M_{H}\left(n\right)\left[\mathfrak{z}\overset{p^{\lambda_{p}\left(n\right)}}{\equiv}n\right],\textrm{ }\forall\mathfrak{z}\in\mathbb{Z}_{p}\label{eq:Quasi vdP series for Chi__H}
\end{equation}
\end{prop}
Proof: Let $n\in\mathbb{N}_{1}$ be arbitrary, and let $\mathbf{j}=\left(j_{1},\ldots,j_{\lambda_{p}\left(n\right)}\right)$
be the shortest string in $\textrm{String}\left(p\right)$ representing
$n$. Since $n\neq0$, observe that $j_{\lambda_{p}\left(n\right)}\in\left\{ 1,\ldots,p-1\right\} $.
Using \textbf{Proposition \ref{prop:Explicit formula for Chi_H of bold j}
}(page \pageref{prop:Explicit formula for Chi_H of bold j}), we have
that: 
\begin{equation}
\chi_{H}\left(n\right)=\sum_{m=1}^{\lambda_{p}\left(n\right)}\frac{b_{j_{m}}}{a_{j_{m}}}\left(\prod_{k=1}^{m}\mu_{j_{k}}\right)p^{-m}
\end{equation}
Note that $n_{-}$ is the integer represented by the string $\mathbf{j}^{\prime}=\left(j_{1},\ldots,j_{\lambda_{p}\left(n\right)-1}\right)$.
Nevertheless, since $H\left(0\right)=0$, $\chi_{H}\left(\mathbf{j}^{\prime}\right)=\chi_{H}\left(n_{-}\right)$,
since removing zeroes at the right end of a string does not affect
the value of $\chi_{H}$. Consequently: 
\begin{equation}
\chi_{H}\left(n_{-}\right)=\sum_{m=1}^{\lambda_{p}\left(n\right)-1}\frac{b_{j_{m}}}{a_{j_{m}}}\left(\prod_{k=1}^{m}\mu_{j_{k}}\right)p^{-m}
\end{equation}
Thus: 
\begin{align*}
c_{n}\left(\chi_{H}\right) & =\chi_{H}\left(n\right)-\chi_{H}\left(n_{-}\right)=\frac{b_{j_{\lambda_{p}\left(n\right)}}}{a_{j_{\lambda_{p}\left(n\right)}}}\left(\prod_{k=1}^{\lambda_{p}\left(n\right)}\mu_{j_{k}}\right)p^{-\lambda_{p}\left(n\right)}=\frac{b_{j_{\lambda_{p}\left(n\right)}}}{a_{j_{\lambda_{p}\left(n\right)}}}M_{H}\left(n\right)
\end{align*}
which tends to $0$ $q_{H}$-adically as the number of non-zero $p$-adic
digits in $n$ tends to $\infty$, which gives another proof of the
rising-continuity of $\chi_{H}$; our first proof of this fact was
implicit in \textbf{Lemma \ref{lem:Unique rising continuation and p-adic functional equation of Chi_H}}.
Since $\chi_{H}$ is rising-continuous, it is represented everywhere
by its van der Put series: 
\begin{equation}
\chi_{H}\left(\mathfrak{z}\right)\overset{\mathbb{Z}_{q_{H}}}{=}\sum_{n=1}^{\infty}\frac{b_{j_{\lambda_{p}\left(n\right)}}}{a_{j_{\lambda_{p}\left(n\right)}}}M_{H}\left(n\right)\left[\mathfrak{z}\overset{p^{\lambda_{p}\left(n\right)}}{\equiv}n\right],\textrm{ }\forall\mathfrak{z}\in\mathbb{Z}_{p}\label{eq:Unsimplified van der Put series for Chi_H}
\end{equation}
where the $n=0$ term is $0$ because $c_{0}\left(\chi_{H}\right)=\chi_{H}\left(0\right)=0$.

When $p=2$, note that $j_{\lambda_{p}\left(n\right)}=1$ for all
$n\geq1$; that is to say, the right-most $2$-adic digit of every
non-zero integer is $1$. This allows us to write: 
\begin{equation}
\chi_{H}\left(\mathfrak{z}\right)\overset{\mathbb{Z}_{q_{H}}}{=}\frac{b_{1}}{a_{1}}\sum_{n=1}^{\infty}M_{H}\left(n\right)\left[\mathfrak{z}\overset{2^{\lambda_{2}\left(n\right)}}{\equiv}n\right],\textrm{ }\forall\mathfrak{z}\in\mathbb{Z}_{2}
\end{equation}
When $p\geq3$, we use the formal identity: 
\begin{equation}
\sum_{n=1}^{p^{N}-1}f\left(n\right)=\sum_{n=1}^{p^{N-1}-1}f\left(n\right)+\sum_{\ell=1}^{p-1}\sum_{n=1}^{p^{N-1}-1}f\left(n+\ell p^{N-1}\right)
\end{equation}
When applying this to $\chi_{H}$'s van der Put series, observe that
$j_{\lambda_{p}\left(n+\ell p^{N-1}\right)}$, the right-most $p$-adic
digit of $n+\ell p^{N-1}$, is then equal to $\ell$. Also: 
\begin{equation}
M_{H}\left(n+\ell p^{N-1}\right)=\frac{\mu_{\ell}}{p}M_{H}\left(n\right)
\end{equation}
and: 
\begin{equation}
\left[\mathfrak{z}\overset{p^{\lambda_{p}\left(n+\ell p^{N-1}\right)}}{\equiv}n+\ell p^{N-1}\right]=\left[\mathfrak{z}\overset{p^{N}}{\equiv}n+\ell p^{N-1}\right]
\end{equation}
These hold for all $n\in\left\{ 1,\ldots,p^{N-1}-1\right\} $ and
all $\ell\in\left\{ 1,\ldots,p-1\right\} $. As such: 
\begin{align*}
\sum_{n=1}^{p^{N}-1}\frac{b_{j_{\lambda_{p}\left(n\right)}}}{a_{j_{\lambda_{p}\left(n\right)}}}M_{H}\left(n\right)\left[\mathfrak{z}\overset{p^{\lambda_{p}\left(n\right)}}{\equiv}n\right] & =\sum_{n=1}^{p^{N-1}-1}\frac{b_{j_{\lambda_{p}\left(n\right)}}}{a_{j_{\lambda_{p}\left(n\right)}}}M_{H}\left(n\right)\left[\mathfrak{z}\overset{p^{\lambda_{p}\left(n\right)}}{\equiv}n\right]\\
 & +\sum_{\ell=1}^{p-1}\sum_{n=1}^{p^{N-1}-1}\frac{b_{\ell}}{a_{\ell}}\frac{\mu_{\ell}}{p}M_{H}\left(n\right)\left[\mathfrak{z}\overset{p^{N}}{\equiv}n+\ell p^{N-1}\right]\\
 & =\sum_{n=1}^{p^{N-1}-1}\frac{b_{j_{\lambda_{p}\left(n\right)}}}{a_{j_{\lambda_{p}\left(n\right)}}}M_{H}\left(n\right)\left[\mathfrak{z}\overset{p^{\lambda_{p}\left(n\right)}}{\equiv}n\right]\\
 & +\sum_{\ell=1}^{p-1}\frac{b_{\ell}}{d_{\ell}}\sum_{n=1}^{p^{N-1}-1}M_{H}\left(n\right)\left[\mathfrak{z}\overset{p^{N}}{\equiv}n+\ell p^{N-1}\right]
\end{align*}
Noting that: 
\begin{equation}
\sum_{n=1}^{p^{N}-1}\frac{b_{j_{\lambda_{p}\left(n\right)}}}{a_{j_{\lambda_{p}\left(n\right)}}}M_{H}\left(n\right)\left[\mathfrak{z}\overset{p^{\lambda_{p}\left(n\right)}}{\equiv}n\right]=S_{p:N}\left\{ \chi_{H}\right\} \left(\mathfrak{z}\right)
\end{equation}
we can simplify the previous equation to: 
\[
S_{p:N}\left\{ \chi_{H}\right\} \left(\mathfrak{z}\right)=S_{p:N-1}\left\{ \chi_{H}\right\} \left(\mathfrak{z}\right)+\sum_{\ell=1}^{p-1}\frac{b_{\ell}}{d_{\ell}}\sum_{n=1}^{p^{N-1}-1}M_{H}\left(n\right)\left[\mathfrak{z}\overset{p^{N}}{\equiv}n+\ell p^{N-1}\right]
\]
Next, we apply the truncated van der Put identity (\textbf{Proposition
\ref{prop:vdP identity}}) to note that: 
\begin{equation}
S_{p:N}\left\{ \chi_{H}\right\} \left(\mathfrak{z}\right)-S_{p:N-1}\left\{ \chi_{H}\right\} \left(\mathfrak{z}\right)=\chi_{H}\left(\left[\mathfrak{z}\right]_{p^{N}}\right)-\chi_{H}\left(\left[\mathfrak{z}\right]_{p^{N-1}}\right)
\end{equation}
So: 
\begin{align*}
\chi_{H}\left(\left[\mathfrak{z}\right]_{p^{K}}\right)-\chi_{H}\left(\left[\mathfrak{z}\right]_{p}\right) & =\sum_{N=2}^{K}\left(\chi_{H}\left(\left[\mathfrak{z}\right]_{p^{N}}\right)-\chi_{H}\left(\left[\mathfrak{z}\right]_{p^{N-1}}\right)\right)\\
 & =\sum_{N=2}^{K}\sum_{\ell=1}^{p-1}\frac{b_{\ell}}{d_{\ell}}\sum_{n=1}^{p^{N-1}-1}M_{H}\left(n\right)\left[\mathfrak{z}\overset{p^{N}}{\equiv}n+\ell p^{N-1}\right]\\
 & =\sum_{\ell=1}^{p-1}\frac{b_{\ell}}{d_{\ell}}\sum_{N=1}^{K-1}\sum_{n=1}^{p^{N}-1}M_{H}\left(n\right)\left[\mathfrak{z}\overset{p^{N+1}}{\equiv}n+\ell p^{N}\right]\\
 & =\sum_{\ell=1}^{p-1}\frac{b_{\ell}}{d_{\ell}}\sum_{n=1}^{p^{K-1}-1}\sum_{N=\lambda_{p}\left(n\right)}^{K-1}M_{H}\left(n\right)\left[\mathfrak{z}\overset{p^{N+1}}{\equiv}n+\ell p^{N}\right]\\
\left(H_{\ell}\left(0\right)=\frac{b_{\ell}}{d_{\ell}}\right); & =\sum_{\ell=1}^{p-1}H_{\ell}\left(0\right)\sum_{n=1}^{p^{K-1}-1}M_{H}\left(n\right)\sum_{N=\lambda_{p}\left(n\right)+1}^{K}\left[\mathfrak{z}\overset{p^{N}}{\equiv}n+\ell p^{N-1}\right]
\end{align*}
Letting $K\rightarrow\infty$ yields: 
\begin{equation}
\chi_{H}\left(\mathfrak{z}\right)-\chi_{H}\left(\left[\mathfrak{z}\right]_{p}\right)\overset{\mathbb{Z}_{q_{H}}}{=}\sum_{j=1}^{p-1}H_{j}\left(0\right)\sum_{n=1}^{\infty}M_{H}\left(n\right)\sum_{N=\lambda_{p}\left(n\right)+1}^{\infty}\left[\mathfrak{z}\overset{p^{N}}{\equiv}n+jp^{N-1}\right]\label{eq:Nearly finished with van-der-put computation for Chi_H}
\end{equation}

\[
\chi_{H}\left(\mathfrak{z}\right)-\chi_{H}\left(\left[\mathfrak{z}\right]_{p}\right)\overset{\mathbb{Z}_{q_{H}}}{=}\sum_{j=1}^{p-1}H_{j}\left(0\right)\sum_{n=1}^{\infty}M_{H}\left(n\right)\sum_{N=\lambda_{p}\left(n\right)+1}^{\infty}\left[\mathfrak{z}\overset{p^{N}}{\equiv}n+jp^{N-1}\right]
\]

Now, fix $n$ and $j$. Then: 
\[
\sum_{N=\lambda_{p}\left(n\right)+1}^{\infty}\left[\mathfrak{z}\overset{p^{N}}{\equiv}n+jp^{N-1}\right]=\sum_{N=\lambda_{p}\left(n\right)+1}^{\infty}\left[\mathfrak{z}\overset{p^{N}}{\equiv}n+jp^{N-1}\right]\left[\mathfrak{z}\overset{p^{N-1}}{\equiv}n\right]
\]
Since $j$ is co-prime to $p$, observe that the congruence $\mathfrak{z}\overset{p^{N}}{\equiv}n+jp^{N-1}$
is equivalent to the following pair of conditions:

i. $\mathfrak{z}\overset{p^{N-1}}{\equiv}n$;

ii. $v_{p}\left(\frac{\mathfrak{z}-n}{p^{N}}-\frac{j}{p}\right)\geq0$.

So: 
\begin{equation}
\left[\mathfrak{z}\overset{p^{N}}{\equiv}n+jp^{N-1}\right]=\left[\mathfrak{z}\overset{p^{N-1}}{\equiv}n\right]\left[v_{p}\left(\frac{\mathfrak{z}-n}{p^{N}}-\frac{j}{p}\right)\geq0\right]
\end{equation}
Here: 
\begin{align*}
v_{p}\left(\frac{\mathfrak{z}-n}{p^{N}}-\frac{j}{p}\right) & \geq0\\
\left(1+v_{p}\left(\mathfrak{y}\right)=v_{p}\left(p\mathfrak{y}\right)\right); & \Updownarrow\\
v_{p}\left(\frac{\mathfrak{z}-n}{p^{N-1}}-j\right) & \geq1\\
 & \Updownarrow\\
\left|\frac{\mathfrak{z}-n}{p^{N-1}}-j\right|_{p} & \leq\frac{1}{p}
\end{align*}
Because $\left|j\right|_{p}=1$, observe that if $\left|\frac{\mathfrak{z}-n}{p^{N-1}}\right|_{p}<1$,
we can then use $p$-adic ultrametric inequality to write: 
\begin{equation}
\frac{1}{p}\geq\left|\frac{\mathfrak{z}-n}{p^{N-1}}-j\right|_{p}=\max\left\{ \left|\frac{\mathfrak{z}-n}{p^{N-1}}\right|_{p},\left|j\right|_{p}\right\} =1
\end{equation}
Of course, this equality is impossible. So, if we have $\mathfrak{z}\overset{p^{N}}{\equiv}n+jp^{N-1}$,
it \emph{must} be that $\left|\frac{\mathfrak{z}-n}{p^{N-1}}\right|_{p}=1$.
This latter condition occurs if and only if $N=v_{p}\left(\mathfrak{z}-n\right)+1$.
Putting it all together, (\ref{eq:Nearly finished with van-der-put computation for Chi_H})
becomes: 
\begin{align*}
\chi_{H}\left(\mathfrak{z}\right)-\chi_{H}\left(\left[\mathfrak{z}\right]_{p}\right) & \overset{\mathbb{Z}_{q_{H}}}{=}\sum_{j=1}^{p-1}H_{j}\left(0\right)\sum_{n=1}^{\infty}M_{H}\left(n\right)\sum_{N=\lambda_{p}\left(n\right)+1}^{\infty}\left[\mathfrak{z}\overset{p^{N}}{\equiv}n+jp^{N-1}\right]\\
 & \overset{\mathbb{Z}_{q_{H}}}{=}\sum_{j=1}^{p-1}H_{j}\left(0\right)\sum_{n=1}^{\infty}M_{H}\left(n\right)\left[\lambda_{p}\left(n\right)\leq v_{p}\left(\mathfrak{z}-n\right)\right]\left[\mathfrak{z}\overset{p^{v_{p}\left(\mathfrak{z}-n\right)+1}}{\equiv}n+jp^{v_{p}\left(\mathfrak{z}-n\right)}\right]\\
 & \overset{\mathbb{Z}_{q_{H}}}{=}\sum_{j=1}^{p-1}H_{j}\left(0\right)\sum_{n=1}^{\infty}M_{H}\left(n\right)\left[\left|\mathfrak{z}-n\right|_{p}\leq p^{-\lambda_{p}\left(n\right)}\right]\left[\mathfrak{z}\overset{p^{v_{p}\left(\mathfrak{z}-n\right)+1}}{\equiv}n+jp^{v_{p}\left(\mathfrak{z}-n\right)}\right]\\
 & \overset{\mathbb{Z}_{q_{H}}}{=}\sum_{j=1}^{p-1}H_{j}\left(0\right)\sum_{n=1}^{\infty}M_{H}\left(n\right)\left[\mathfrak{z}\overset{p^{\lambda_{p}\left(n\right)}}{\equiv}n\right]\left[\mathfrak{z}\overset{p^{v_{p}\left(\mathfrak{z}-n\right)+1}}{\equiv}n+jp^{v_{p}\left(\mathfrak{z}-n\right)}\right]
\end{align*}
Here: 
\begin{align*}
\mathfrak{z} & \overset{p^{v_{p}\left(\mathfrak{z}-n\right)+1}}{\equiv}n+jp^{v_{p}\left(\mathfrak{z}-n\right)}\\
 & \Updownarrow\\
\mathfrak{z}-n & \in jp^{v_{p}\left(\mathfrak{z}-n\right)}+p^{v_{p}\left(\mathfrak{z}-n\right)+1}\mathbb{Z}_{p}\\
 & \Updownarrow\\
\left(\mathfrak{z}-n\right)\left|\mathfrak{z}-n\right|_{p} & \in j+p\mathbb{Z}_{p}\\
 & \Updownarrow\\
u_{p}\left(\mathfrak{z}-n\right) & \overset{p}{\equiv}j
\end{align*}
and so: 
\begin{align}
\chi_{H}\left(\mathfrak{z}\right)-\chi_{H}\left(\left[\mathfrak{z}\right]_{p}\right) & \overset{\mathbb{Z}_{q_{H}}}{=}\sum_{n=1}^{\infty}M_{H}\left(n\right)\left[\mathfrak{z}\overset{p^{\lambda_{p}\left(n\right)}}{\equiv}n\right]\sum_{j=1}^{p-1}H_{j}\left(0\right)\left[u_{p}\left(\mathfrak{z}-n\right)\overset{p}{\equiv}j\right]
\end{align}

Q.E.D.

\vphantom{}

With these formulae, we can directly estimate the $L_{\mathbb{R}}^{1}\left(\mathbb{Z}_{p},\mathbb{C}_{q_{H}}\right)$
norm of $\chi_{H}$. 
\begin{prop}
\label{prop:L_R^1 norm of Chi_H}\index{chi{H}@$\chi_{H}$!L{mathbb{R}}{1}-norm@$L_{\mathbb{R}}^{1}$-norm}
\begin{equation}
\int_{\mathbb{Z}_{p}}\left|\chi_{H}\left(\mathfrak{z}\right)\right|_{q_{H}}d\mathfrak{z}\leq\frac{1}{p}\left(\sum_{j=1}^{p-1}\left|\chi_{H}\left(j\right)\right|_{q_{H}}\right)\sum_{n=0}^{\infty}\frac{\left|M_{H}\left(n\right)\right|_{q_{H}}}{p^{\lambda_{p}\left(n\right)}}\label{eq:real L^1 estimate for Chi_H}
\end{equation}
\end{prop}
Proof: For brevity, let $q=q_{H}$. We use \textbf{Proposition \ref{prop:van der Put series for Chi_H}}.
This gives: 
\begin{align*}
\int_{\mathbb{Z}_{p}}\left|\chi_{H}\left(\mathfrak{z}\right)\right|_{q}d\mathfrak{z} & \leq\sum_{n=1}^{\infty}\sum_{j=1}^{p-1}\left|H_{j}\left(0\right)M_{H}\left(n\right)\right|_{q}\int_{\mathbb{Z}_{p}}\left[u_{p}\left(\mathfrak{z}-n\right)\overset{p}{\equiv}j\right]\left[\mathfrak{z}\overset{p^{\lambda_{p}\left(n\right)}}{\equiv}n\right]d\mathfrak{z}\\
 & +\int_{\mathbb{Z}_{p}}\left|\chi_{H}\left(\left[\mathfrak{z}\right]_{p}\right)\right|_{q}d\mathfrak{z}
\end{align*}
Using the translation invariance of $d\mathfrak{z}$, we can write:
\begin{align*}
\int_{\mathbb{Z}_{p}}\left[u_{p}\left(\mathfrak{z}-n\right)\overset{p}{\equiv}j\right]\left[\mathfrak{z}\overset{p^{\lambda_{p}\left(n\right)}}{\equiv}n\right]d\mathfrak{z} & \overset{\mathbb{R}}{=}\int_{\mathbb{Z}_{p}}\left[u_{p}\left(\mathfrak{z}\right)\overset{p}{\equiv}j\right]\left[\mathfrak{z}+n\overset{p^{\lambda_{p}\left(n\right)}}{\equiv}n\right]d\mathfrak{z}\\
 & =\int_{\mathbb{Z}_{p}}\left[u_{p}\left(\mathfrak{z}\right)\overset{p}{\equiv}j\right]\left[\mathfrak{z}\overset{p^{\lambda_{p}\left(n\right)}}{\equiv}0\right]d\mathfrak{z}\\
 & =\int_{p^{\lambda_{p}\left(n\right)}\mathbb{Z}_{p}}\left[u_{p}\left(\mathfrak{z}\right)\overset{p}{\equiv}j\right]d\mathfrak{z}\\
\left(\mathfrak{y}=\frac{\mathfrak{z}}{p^{\lambda_{p}\left(n\right)}}\right); & =\frac{1}{p^{\lambda_{p}\left(n\right)}}\int_{\mathbb{Z}_{p}}\left[u_{p}\left(p^{\lambda_{p}\left(n\right)}\mathfrak{y}\right)\overset{p}{\equiv}j\right]d\mathfrak{z}\\
\left(u_{p}\left(p\mathfrak{y}\right)=u_{p}\left(\mathfrak{y}\right)\right); & =\frac{1}{p^{\lambda_{p}\left(n\right)}}\int_{\mathbb{Z}_{p}}\left[u_{p}\left(\mathfrak{y}\right)\overset{p}{\equiv}j\right]d\mathfrak{z}\\
\left(j\in\left\{ 1,\ldots,p-1\right\} :u_{p}\left(\mathfrak{y}\right)\overset{p}{\equiv}j\Leftrightarrow\mathfrak{y}\overset{p}{\equiv}j\right); & =\frac{1}{p^{\lambda_{p}\left(n\right)}}\int_{\mathbb{Z}_{p}}\left[\mathfrak{y}\overset{p}{\equiv}j\right]d\mathfrak{z}\\
 & =\frac{1}{p^{\lambda_{p}\left(n\right)+1}}
\end{align*}
Thus: 
\begin{align*}
\int_{\mathbb{Z}_{p}}\left|\chi_{H}\left(\mathfrak{z}\right)\right|_{q}d\mathfrak{z} & \leq\sum_{n=1}^{\infty}\sum_{j=1}^{p-1}\frac{\left|H_{j}\left(0\right)M_{H}\left(n\right)\right|_{q}}{p^{\lambda_{p}\left(n\right)+1}}+\int_{\mathbb{Z}_{p}}\left|\chi_{H}\left(\left[\mathfrak{z}\right]_{p}\right)\right|_{q}d\mathfrak{z}\\
 & =\sum_{n=1}^{\infty}\sum_{j=1}^{p-1}\frac{\left|H_{j}\left(0\right)M_{H}\left(n\right)\right|_{q}}{p^{\lambda_{p}\left(n\right)+1}}+\frac{1}{p}\sum_{j=0}^{p-1}\left|\chi_{H}\left(j\right)\right|_{q}\\
\left(\chi_{H}\left(0\right)=0\right); & =\frac{1}{p}\sum_{j=1}^{p-1}\left(\left|\chi_{H}\left(j\right)\right|_{q}+\sum_{n=1}^{\infty}\left|H_{j}\left(0\right)\right|_{q}\frac{\left|M_{H}\left(n\right)\right|_{q}}{p^{\lambda_{p}\left(n\right)}}\right)
\end{align*}
Finally, since $\chi_{H}\left(j\right)=H_{j}\left(0\right)$ for all
$j\in\left\{ 0,\ldots,p-1\right\} $, we have: 
\begin{align*}
\int_{\mathbb{Z}_{p}}\left|\chi_{H}\left(\mathfrak{z}\right)\right|_{q}d\mathfrak{z} & \leq\frac{1}{p}\sum_{j=1}^{p-1}\left(\left|\chi_{H}\left(j\right)\right|_{q}+\sum_{n=1}^{\infty}\left|H_{j}\left(0\right)\right|_{q}\frac{\left|M_{H}\left(n\right)\right|_{q}}{p^{\lambda_{p}\left(n\right)}}\right)\\
 & =\frac{1}{p}\left(\sum_{j=1}^{p-1}\left|\chi_{H}\left(j\right)\right|_{q}\right)\left(1+\sum_{n=1}^{\infty}\frac{\left|M_{H}\left(n\right)\right|_{q}}{p^{\lambda_{p}\left(n\right)}}\right)\\
\left(M_{H}\left(0\right)=1\right); & =\frac{1}{p}\left(\sum_{j=1}^{p-1}\left|\chi_{H}\left(j\right)\right|_{q}\right)\sum_{n=0}^{\infty}\frac{\left|M_{H}\left(n\right)\right|_{q}}{p^{\lambda_{p}\left(n\right)}}
\end{align*}

Q.E.D.

\vphantom{}

With this, we can now state the sufficient condition for $\chi_{H}$
to be in $L_{\mathbb{R}}^{1}$. 
\begin{thm}[\textbf{Criterion for $\chi_{H}\in L_{\mathbb{R}}^{1}$}]
\label{thm:L^1 criterion for Chi_H}$\chi_{H}$ is an element of
$L_{\mathbb{R}}^{1}\left(\mathbb{Z}_{p},\mathbb{C}_{q_{H}}\right)$
whenever: 
\begin{equation}
\sum_{j=0}^{p-1}\left|\mu_{j}\right|_{q_{H}}<p\label{eq:Criterion for Chi_H to be real L^1 integrable}
\end{equation}
In particular, since the $\mu_{j}$s are integers, this shows that
$\chi_{H}\in L_{\mathbb{R}}^{1}\left(\mathbb{Z}_{p},\mathbb{C}_{q_{H}}\right)$
occurs whenever $\gcd\left(\mu_{j},q_{H}\right)>1$ for at least one
$j\in\left\{ 0,\ldots,p-1\right\} $. 
\end{thm}
Proof: For brevity, let $q=q_{H}$. By \textbf{Proposition \ref{prop:L_R^1 norm of Chi_H}},
to show that $\chi_{H}$ has finite $L_{\mathbb{R}}^{1}\left(\mathbb{Z}_{p},\mathbb{C}_{q}\right)$
norm, it suffices to prove that: 
\begin{equation}
\sum_{n=0}^{\infty}\frac{\left|M_{H}\left(n\right)\right|_{q}}{p^{\lambda_{p}\left(n\right)}}<\infty
\end{equation}
To do this, we use a lambda decomposition: 
\begin{equation}
\sum_{n=0}^{\infty}\frac{\left|M_{H}\left(n\right)\right|_{q}}{p^{\lambda_{p}\left(n\right)}}\overset{\mathbb{R}}{=}1+\sum_{m=1}^{\infty}\frac{1}{p^{m}}\sum_{n=p^{m-1}}^{p^{m}-1}\left|M_{H}\left(n\right)\right|_{q}\label{eq:M_H q-adic absolute value lambda decomposition}
\end{equation}
Once again, we proceed recursively: 
\begin{equation}
S_{m}\overset{\textrm{def}}{=}\sum_{n=0}^{p^{m}-1}\left|M_{H}\left(n\right)\right|_{q}
\end{equation}
Then, using $M_{H}$'s functional equations: 
\begin{align*}
\sum_{n=0}^{p^{m}-1}\left|M_{H}\left(n\right)\right|_{q} & =\sum_{j=0}^{p-1}\left|M_{H}\left(j\right)\right|_{q}+\sum_{j=0}^{p-1}\sum_{n=1}^{p^{m-1}-1}\left|M_{H}\left(pn+j\right)\right|_{q}\\
 & =\sum_{j=0}^{p-1}\left|M_{H}\left(j\right)\right|_{q}+\sum_{j=0}^{p-1}\sum_{n=1}^{p^{m-1}-1}\left|\frac{\mu_{j}}{p}M_{H}\left(n\right)\right|_{q}\\
\left(\gcd\left(p,q\right)=1\right); & =\sum_{j=0}^{p-1}\left|M_{H}\left(j\right)\right|_{q}+\sum_{n=1}^{p^{m-1}-1}\left(\sum_{j=0}^{p-1}\left|\mu_{j}\right|_{q}\right)\left|M_{H}\left(n\right)\right|_{q}\\
\left(M_{H}\left(j\right)=\mu_{j}\textrm{ }\forall j\in\left\{ 1,\ldots,p-1\right\} \right); & =1+\sum_{j=1}^{p-1}\left|\mu_{j}\right|_{q}+\left(\sum_{j=0}^{p-1}\left|\mu_{j}\right|_{q}\right)\sum_{n=1}^{p^{m-1}-1}\left|M_{H}\left(n\right)\right|_{q}\\
 & =1-\left|\mu_{0}\right|_{q}+\left(\sum_{j=0}^{p-1}\left|\mu_{j}\right|_{q}\right)\left(1+\sum_{n=1}^{p^{m-1}-1}\left|M_{H}\left(n\right)\right|_{q}\right)\\
 & =1-\left|\mu_{0}\right|_{q}+\left(\sum_{j=0}^{p-1}\left|\mu_{j}\right|_{q}\right)\sum_{n=0}^{p^{m-1}-1}\left|M_{H}\left(n\right)\right|_{q}
\end{align*}
So: 
\begin{equation}
S_{m}=1-\left|\mu_{0}\right|_{q}+\left(\sum_{j=0}^{p-1}\left|\mu_{j}\right|_{q}\right)S_{m-1},\textrm{ }\forall m\geq1
\end{equation}
where: 
\begin{equation}
S_{0}=\left|M_{H}\left(0\right)\right|_{q}=1
\end{equation}
Finally, defining $A$ and $B$ by:
\begin{align}
A & \overset{\textrm{def}}{=}1-\left|\mu_{0}\right|_{q}\\
B & \overset{\textrm{def}}{=}\sum_{j=0}^{p-1}\left|\mu_{j}\right|_{q}
\end{align}
we can write:
\begin{align*}
S_{m} & =A+BS_{m-1}\\
 & =A+AB+B^{2}S_{m-2}\\
 & =A+AB+AB^{2}+B^{3}S_{m-3}\\
 & \vdots\\
 & =A\sum_{k=0}^{m-1}B^{k}+B^{m}S_{0}\\
 & =B^{m}+A\sum_{k=0}^{m-1}B^{k}
\end{align*}
Using this, our lambda decomposition (\ref{eq:M_H q-adic absolute value lambda decomposition})
becomes: 
\begin{align*}
\sum_{n=0}^{\infty}\frac{\left|M_{H}\left(n\right)\right|_{q}}{p^{\lambda_{p}\left(n\right)}} & \overset{\mathbb{R}}{=}1+\sum_{m=1}^{\infty}\frac{1}{p^{m}}\left(\sum_{n=0}^{p^{m}-1}\left|M_{H}\left(n\right)\right|_{q}-\sum_{n=0}^{p^{m-1}-1}\left|M_{H}\left(n\right)\right|_{q}\right)\\
 & =1+\sum_{m=1}^{\infty}\frac{1}{p^{m}}\left(S_{m}-S_{m-1}\right)\\
 & =1+\sum_{m=1}^{\infty}\frac{S_{m}}{p^{m}}-p\sum_{m=1}^{\infty}\frac{S_{m-1}}{p^{m-1}}\\
 & =1-p\frac{S_{0}}{p^{0}}+\sum_{m=1}^{\infty}\frac{S_{m}}{p^{m}}-p\sum_{m=1}^{\infty}\frac{S_{m}}{p^{m}}\\
 & =1-p+\left(1-p\right)\sum_{m=1}^{\infty}\frac{S_{m}}{p^{m}}\\
 & =1-p+\left(1-p\right)\sum_{m=1}^{\infty}\frac{\left(B^{m}+A\sum_{k=0}^{m-1}B^{k}\right)}{p^{m}}\\
 & \overset{\mathbb{R}}{=}\begin{cases}
1-p+\left(1-p\right)\sum_{m=1}^{\infty}\frac{Am+1}{p^{m}} & \textrm{if }B=1\\
1-p+\left(1-p\right)\sum_{m=1}^{\infty}\frac{B^{m}+A\frac{B^{m}-1}{B-1}}{p^{m}} & \textrm{else}
\end{cases}
\end{align*}
The $B=1$ case is finite since $p\geq2$. As for $B\neq1$, summing
the geometric series yields: 
\begin{equation}
\sum_{m=1}^{\infty}\frac{B^{m}+A\frac{B^{m}-1}{B-1}}{p^{m}}\overset{\mathbb{R}}{=}-\frac{A}{B-1}\frac{1}{p-1}+\frac{A+B-1}{B-1}\sum_{m=1}^{\infty}\left(\frac{B}{p}\right)^{m}
\end{equation}
Hence, $0\leq B<p$ is sufficient to guarantee convergence in the
topology of $\mathbb{R}$.

Q.E.D.

\vphantom{}

We end our wisp of $L_{\mathbb{R}}^{1}$ by proving in general what
was shown for $\chi_{5}$ in Subsection \ref{subsec:3.3.6 L^1 Convergence}:
that $\tilde{\chi}_{H,N}-\chi_{H,N}$ converges to $0$ in $L_{\mathbb{R}}^{1}$
as $N\rightarrow\infty$.

The first step is another computation\textemdash this time outsourced
to previous work. 
\begin{prop}
\label{prop:Convolution of alpha and A_H hat}Suppose $\alpha_{H}\left(0\right)\notin\left\{ 0,1\right\} $.
Then:

\begin{equation}
\sum_{0<\left|t\right|_{p}\leq p^{N}}\left(\alpha_{H}\left(0\right)\right)^{v_{p}\left(t\right)}\hat{A}_{H}\left(t\right)e^{2\pi i\left\{ t\mathfrak{z}\right\} _{p}}\overset{\overline{\mathbb{Q}}}{=}\left(\frac{\mu_{0}}{p\alpha_{H}\left(0\right)}\right)^{N}\kappa_{H}\left(\left[\mathfrak{z}\right]_{p^{N}}\right)-1\label{eq:Partial Fourier Sum of A_H hat times alpha to the v_p(t)}
\end{equation}
\end{prop}
Proof: Use (\ref{eq:Convolution of dA_H and D_N}) from \textbf{Theorem
\ref{thm:Properties of dA_H}}.

Q.E.D.

\vphantom{}

With this, we can compute simple closed-form expressions for $\tilde{\chi}_{H,N}\left(\mathfrak{z}\right)-\chi_{H,N}\left(\mathfrak{z}\right)$.
It should be noted that the formulae given below only hold for the
indicated choice of $\hat{\chi}_{H}$. Choosing a different Fourier
transform of $\hat{\chi}_{H}$ will alter the formula for $\tilde{\chi}_{H,N}\left(\mathfrak{z}\right)-\chi_{H,N}\left(\mathfrak{z}\right)$. 
\begin{thm}[\textbf{$\chi_{H,N}\left(\mathfrak{z}\right)-\tilde{\chi}_{H,N}\left(\mathfrak{z}\right)$}]
\label{thm:Chi_H,N / Chi_H,N twiddle error}Choose $\hat{\chi}_{H}\left(t\right)$
to be the Fourier Transform of $\chi_{H}$ defined by \emph{(\ref{eq:Fourier Transform of Chi_H for a contracting semi-basic rho-Hydra map})}:
\[
\hat{\chi}_{H}\left(t\right)\overset{\textrm{def}}{=}\begin{cases}
\begin{cases}
0 & \textrm{if }t=0\\
\left(\beta_{H}\left(0\right)v_{p}\left(t\right)+\gamma_{H}\left(\frac{t\left|t\right|_{p}}{p}\right)\right)\hat{A}_{H}\left(t\right) & \textrm{if }t\neq0
\end{cases} & \textrm{if }\alpha_{H}\left(0\right)=1\\
\frac{\beta_{H}\left(0\right)\hat{A}_{H}\left(t\right)}{1-\alpha_{H}\left(0\right)}+\begin{cases}
0 & \textrm{if }t=0\\
\gamma_{H}\left(\frac{t\left|t\right|_{p}}{p}\right)\hat{A}_{H}\left(t\right) & \textrm{if }t\neq0
\end{cases} & \textrm{if }\alpha_{H}\left(0\right)\neq1
\end{cases},\textrm{ }\forall t\in\hat{\mathbb{Z}}_{p}
\]
Then:

\vphantom{}

I. If $\alpha_{H}\left(0\right)=1$: 
\begin{equation}
\chi_{H,N}\left(\mathfrak{z}\right)-\tilde{\chi}_{H,N}\left(\mathfrak{z}\right)\overset{\overline{\mathbb{Q}}}{=}\beta_{H}\left(0\right)\left(\frac{\mu_{0}}{p}\right)^{N-1}\left(\frac{2\mu_{0}N}{p}\kappa_{H}\left(\left[\mathfrak{z}\right]_{p^{N}}\right)-\left(N-1\right)\kappa_{H}\left(\left[\mathfrak{z}\right]_{p^{N-1}}\right)\right)\label{eq:Chi_H,N / Chi_H,N twiddle identity for alpha equals 1}
\end{equation}

\vphantom{}

II. If $\alpha_{H}\left(0\right)\neq1$: 
\begin{equation}
\chi_{H,N}\left(\mathfrak{z}\right)-\tilde{\chi}_{H,N}\left(\mathfrak{z}\right)\overset{\overline{\mathbb{Q}}}{=}\frac{\beta_{H}\left(0\right)}{\alpha_{H}\left(0\right)-1}\left(\frac{\mu_{0}}{p}\right)^{N-1}\left(\frac{2\mu_{0}}{p}\kappa_{H}\left(\left[\mathfrak{z}\right]_{p^{N}}\right)-\alpha_{H}\left(0\right)\kappa_{H}\left(\left[\mathfrak{z}\right]_{p^{N-1}}\right)\right)\label{eq:Chi_H,N / Chi_H,N twiddle identity for alpha is not 1}
\end{equation}
\end{thm}
Proof:

I. Let $\alpha_{H}\left(0\right)=1$. Then, we have that: 
\begin{equation}
\hat{\chi}_{H}\left(t\right)\overset{\overline{\mathbb{Q}}}{=}\begin{cases}
0 & \textrm{if }t=0\\
\left(\beta_{H}\left(0\right)v_{p}\left(t\right)+\gamma_{H}\left(\frac{t\left|t\right|_{p}}{p}\right)\right)\hat{A}_{H}\left(t\right) & \textrm{if }t\neq0
\end{cases}
\end{equation}
Turning to \textbf{Theorem \ref{thm:(N,t) asymptotic decomposition of Chi_H,N hat}}),
the associated ``fine-structure'' formula for $\hat{\chi}_{H,N}$
in this case is: 
\[
\hat{\chi}_{H,N}\left(t\right)-\beta_{H}\left(0\right)N\hat{A}_{H}\left(t\right)\mathbf{1}_{0}\left(p^{N-1}t\right)=\begin{cases}
0 & \textrm{if }t=0\\
\left(\gamma_{H}\left(\frac{t\left|t\right|_{p}}{p}\right)+\beta_{H}\left(0\right)v_{p}\left(t\right)\right)\hat{A}_{H}\left(t\right) & \textrm{if }0<\left|t\right|_{p}<p^{N}\\
\gamma_{H}\left(\frac{t\left|t\right|_{p}}{p}\right)\hat{A}_{H}\left(t\right) & \textrm{if }\left|t\right|_{p}=p^{N}\\
0 & \textrm{if }\left|t\right|_{p}>p^{N}
\end{cases}
\]
Using use our chosen $\hat{\chi}_{H}$ (and remembering that $\hat{A}_{H}\left(0\right)=1$),
this can be rewritten as: 
\begin{align}
\hat{\chi}_{H,N}\left(t\right)-\beta_{H}\left(0\right)N\hat{A}_{H}\left(t\right)\mathbf{1}_{0}\left(p^{N-1}t\right) & \overset{\overline{\mathbb{Q}}}{=}\mathbf{1}_{0}\left(p^{N-1}t\right)\hat{\chi}_{H}\left(t\right)\label{eq:Chi_H,N / Chi_H,N twiddle proof: eq 1}\\
 & +\left[\left|t\right|_{p}=p^{N}\right]\gamma_{H}\left(\frac{t\left|t\right|_{p}}{p}\right)\hat{A}_{H}\left(t\right),\textrm{ }\forall\left|t\right|_{p}\leq p^{N}\nonumber 
\end{align}
which is to say: 
\begin{align}
\hat{\chi}_{H,N}\left(t\right)-\beta_{H}\left(0\right)N\hat{A}_{H}\left(t\right)\mathbf{1}_{0}\left(p^{N-1}t\right) & \overset{\overline{\mathbb{Q}}}{=}\mathbf{1}_{0}\left(p^{N}t\right)\hat{\chi}_{H}\left(t\right)\label{eq:Chi_H,N / Chi_H,N twiddle proof: eq 2}\\
 & +\left[\left|t\right|_{p}=p^{N}\right]\left(\gamma_{H}\left(\frac{t\left|t\right|_{p}}{p}\right)\hat{A}_{H}\left(t\right)-\hat{\chi}_{H}\left(t\right)\right),\textrm{ }\forall\left|t\right|_{p}\leq p^{N}\nonumber 
\end{align}
Multiplying through by $e^{2\pi i\left\{ t\mathfrak{z}\right\} _{p}}$
and summing over $\left|t\right|_{p}\leq p^{N}$ yields: 
\begin{align}
\chi_{H,N}\left(\mathfrak{z}\right)-\beta_{H}\left(0\right)N\tilde{A}_{H,N-1}\left(\mathfrak{z}\right) & \overset{\overline{\mathbb{Q}}}{=}\tilde{\chi}_{H,N}\left(\mathfrak{z}\right)\label{eq:Chi_H,N / Chi_H,N twiddle proof: eq 3}\\
 & +\sum_{\left|t\right|_{p}=p^{N}}\left(\gamma_{H}\left(\frac{t\left|t\right|_{p}}{p}\right)\hat{A}_{H}\left(t\right)-\hat{\chi}_{H}\left(t\right)\right)e^{2\pi i\left\{ t\mathfrak{z}\right\} _{p}}\nonumber 
\end{align}
Using \textbf{Lemma \ref{lem:1D gamma formula}}, we have: 
\begin{align*}
\sum_{\left|t\right|_{p}=p^{N}}\gamma_{H}\left(\frac{t\left|t\right|_{p}}{p}\right)\hat{A}_{H}\left(t\right)e^{2\pi i\left\{ t\mathfrak{z}\right\} _{p}} & =\sum_{0<\left|t\right|_{p}\leq p^{N}}\gamma_{H}\left(\frac{t\left|t\right|_{p}}{p}\right)\hat{A}_{H}\left(t\right)e^{2\pi i\left\{ t\mathfrak{z}\right\} _{p}}\\
 & -\sum_{0<\left|t\right|_{p}\leq p^{N-1}}\gamma_{H}\left(\frac{t\left|t\right|_{p}}{p}\right)\hat{A}_{H}\left(t\right)e^{2\pi i\left\{ t\mathfrak{z}\right\} _{p}}\\
 & =\sum_{n=0}^{N-1}\left(\sum_{j=1}^{p-1}\beta_{H}\left(\frac{j}{p}\right)\varepsilon_{n}^{j}\left(\mathfrak{z}\right)\right)\left(\frac{\mu_{0}}{p}\right)^{n}\kappa_{H}\left(\left[\mathfrak{z}\right]_{p^{n}}\right)\\
 & -\sum_{n=0}^{N-2}\left(\sum_{j=1}^{p-1}\beta_{H}\left(\frac{j}{p}\right)\varepsilon_{n}^{j}\left(\mathfrak{z}\right)\right)\left(\frac{\mu_{0}}{p}\right)^{n}\kappa_{H}\left(\left[\mathfrak{z}\right]_{p^{n}}\right)\\
 & =\left(\sum_{j=1}^{p-1}\beta_{H}\left(\frac{j}{p}\right)\varepsilon_{N-1}^{j}\left(\mathfrak{z}\right)\right)\left(\frac{\mu_{0}}{p}\right)^{N-1}\kappa_{H}\left(\left[\mathfrak{z}\right]_{p^{N-1}}\right)
\end{align*}
Using (\ref{eq:Explicit Formula for Chi_H,N twiddle for arbitrary rho and alpha equals 1})
from \textbf{Corollary \ref{cor:Chi_H,N twiddle explicit formula, arbitrary p, arbitrary alpha}}
gives us:

\begin{align*}
\sum_{\left|t\right|_{p}=p^{N}}\hat{\chi}_{H}\left(t\right)e^{2\pi i\left\{ t\mathfrak{z}\right\} _{p}} & =\tilde{\chi}_{H,N}\left(\mathfrak{z}\right)-\tilde{\chi}_{H,N-1}\left(\mathfrak{z}\right)\\
 & \overset{\overline{\mathbb{Q}}}{=}-\beta_{H}\left(0\right)N\left(\frac{\mu_{0}}{p}\right)^{N}\kappa_{H}\left(\left[\mathfrak{z}\right]_{p^{N}}\right)+\beta_{H}\left(0\right)\sum_{n=0}^{N-1}\left(\frac{\mu_{0}}{p}\right)^{n}\kappa_{H}\left(\left[\mathfrak{z}\right]_{p^{n}}\right)\\
 & +\sum_{n=0}^{N-1}\left(\sum_{j=1}^{p-1}\beta_{H}\left(\frac{j}{p}\right)\varepsilon_{n}^{j}\left(\mathfrak{z}\right)\right)\left(\frac{\mu_{0}}{p}\right)^{n}\kappa_{H}\left(\left[\mathfrak{z}\right]_{p^{n}}\right)\\
 & +\beta_{H}\left(0\right)\left(N-1\right)\left(\frac{\mu_{0}}{p}\right)^{N-1}\kappa_{H}\left(\left[\mathfrak{z}\right]_{p^{N-1}}\right)\\
 & -\beta_{H}\left(0\right)\sum_{n=0}^{N-2}\left(\frac{\mu_{0}}{p}\right)^{n}\kappa_{H}\left(\left[\mathfrak{z}\right]_{p^{n}}\right)\\
 & -\sum_{n=0}^{N-2}\left(\sum_{j=1}^{p-1}\beta_{H}\left(\frac{j}{p}\right)\varepsilon_{n}^{j}\left(\mathfrak{z}\right)\right)\left(\frac{\mu_{0}}{p}\right)^{n}\kappa_{H}\left(\left[\mathfrak{z}\right]_{p^{n}}\right)\\
 & \overset{\overline{\mathbb{Q}}}{=}\beta_{H}\left(0\right)\left(N-1\right)\left(\frac{\mu_{0}}{p}\right)^{N-1}\kappa_{H}\left(\left[\mathfrak{z}\right]_{p^{N-1}}\right)\\
 & -\beta_{H}\left(0\right)N\left(\frac{\mu_{0}}{p}\right)^{N}\kappa_{H}\left(\left[\mathfrak{z}\right]_{p^{N}}\right)\\
 & +\left(\sum_{j=0}^{p-1}\beta_{H}\left(\frac{j}{p}\right)\varepsilon_{N-1}^{j}\left(\mathfrak{z}\right)\right)\left(\frac{\mu_{0}}{p}\right)^{N-1}\kappa_{H}\left(\left[\mathfrak{z}\right]_{p^{N-1}}\right)
\end{align*}
Finally, employing our familiar formula for for $\tilde{A}_{H,N}\left(\mathfrak{z}\right)$
(equation (\ref{eq:Convolution of dA_H and D_N}) from \textbf{Theorem
\ref{thm:Properties of dA_H}}) and remembering that $\alpha_{H}\left(0\right)=1$,
(\ref{eq:Chi_H,N / Chi_H,N twiddle proof: eq 3}) becomes: 
\begin{align*}
\chi_{H,N}\left(\mathfrak{z}\right)-\beta_{H}\left(0\right)N\left(\frac{\mu_{0}}{p}\right)^{N}\kappa_{H}\left(\left[\mathfrak{z}\right]_{p^{N}}\right) & \overset{\overline{\mathbb{Q}}}{=}\tilde{\chi}_{H,N}\left(\mathfrak{z}\right)\\
 & -\beta_{H}\left(0\right)\left(N-1\right)\left(\frac{\mu_{0}}{p}\right)^{N-1}\kappa_{H}\left(\left[\mathfrak{z}\right]_{p^{N-1}}\right)\\
 & +\beta_{H}\left(0\right)N\left(\frac{\mu_{0}}{p}\right)^{N}\kappa_{H}\left(\left[\mathfrak{z}\right]_{p^{N}}\right)
\end{align*}
and hence: 
\[
\chi_{H,N}\left(\mathfrak{z}\right)-\tilde{\chi}_{H,N}\left(\mathfrak{z}\right)\overset{\overline{\mathbb{Q}}}{=}\beta_{H}\left(0\right)\left(\frac{\mu_{0}}{p}\right)^{N-1}\left(\frac{2\mu_{0}N}{p}\kappa_{H}\left(\left[\mathfrak{z}\right]_{p^{N}}\right)-\left(N-1\right)\kappa_{H}\left(\left[\mathfrak{z}\right]_{p^{N-1}}\right)\right)
\]
which proves (I).

\vphantom{}

II. Suppose $\alpha_{H}\left(0\right)\neq1$. The ``fine-structure''
equation from our $\left(N,t\right)$-asymptotic analysis of $\hat{\chi}_{H,N}$
((\ref{eq:Fine Structure Formula for Chi_H,N hat when alpha is not 1})
from \textbf{Theorem \ref{thm:(N,t) asymptotic decomposition of Chi_H,N hat}})
is: 
\[
\hat{\chi}_{H,N}\left(t\right)=\begin{cases}
\beta_{H}\left(0\right)\frac{\left(\alpha_{H}\left(0\right)\right)^{N}-1}{\alpha_{H}\left(0\right)-1} & \textrm{if }t=0\\
\left(\gamma_{H}\left(\frac{t\left|t\right|_{p}}{p}\right)+\beta_{H}\left(0\right)\frac{\left(\alpha_{H}\left(0\right)\right)^{N+v_{p}\left(t\right)}-1}{\alpha_{H}\left(0\right)-1}\right)\hat{A}_{H}\left(t\right) & \textrm{if }0<\left|t\right|_{p}<p^{N}\\
\gamma_{H}\left(\frac{t\left|t\right|_{p}}{p}\right)\hat{A}_{H}\left(t\right) & \textrm{if }\left|t\right|_{p}=p^{N}\\
0 & \textrm{if }\left|t\right|_{p}>p^{N}
\end{cases},\textrm{ }\forall t\in\hat{\mathbb{Z}}_{p}
\]
Meanwhile, $\hat{\chi}_{H}\left(t\right)$ is: 
\[
\hat{\chi}_{H}\left(t\right)=\begin{cases}
\frac{\beta_{H}\left(0\right)}{1-\alpha_{H}\left(0\right)} & \textrm{if }t=0\\
\left(\frac{\beta_{H}\left(0\right)}{1-\alpha_{H}\left(0\right)}+\gamma_{H}\left(\frac{t\left|t\right|_{p}}{p}\right)\right)\hat{A}_{H}\left(t\right) & \textrm{if }t\neq0
\end{cases}
\]
Subtracting yields: 
\begin{equation}
\hat{\chi}_{H,N}\left(t\right)-\hat{\chi}_{H}\left(t\right)=\begin{cases}
\frac{\beta_{H}\left(0\right)\left(\alpha_{H}\left(0\right)\right)^{N}}{\alpha_{H}\left(0\right)-1} & \textrm{if }t=0\\
\frac{\beta_{H}\left(0\right)\left(\alpha_{H}\left(0\right)\right)^{N+v_{p}\left(t\right)}}{\alpha_{H}\left(0\right)-1}\hat{A}_{H}\left(t\right) & \textrm{if }0<\left|t\right|_{p}<p^{N}\\
\frac{\beta_{H}\left(0\right)}{\alpha_{H}\left(0\right)-1}\hat{A}_{H}\left(t\right) & \textrm{if }\left|t\right|_{p}=p^{N}\\
\left(\frac{\beta_{H}\left(0\right)}{\alpha_{H}\left(0\right)-1}-\gamma_{H}\left(\frac{t\left|t\right|_{p}}{p}\right)\right)\hat{A}_{H}\left(t\right) & \textrm{if }\left|t\right|_{p}>p^{N}
\end{cases},\textrm{ }\forall t\in\hat{\mathbb{Z}}_{p}\label{eq:Chi_H,N / Chi_H,N twiddle proof: eq 4}
\end{equation}
Since $\hat{\chi}_{H,N}\left(t\right)$ is supported on $\left|t\right|_{p}\leq p^{N}$,
when summing our Fourier series, we need only sum over $\left|t\right|_{p}\leq p^{N}$.
Doing so yields: 
\begin{align*}
\chi_{H,N}\left(\mathfrak{z}\right)-\tilde{\chi}_{H,N}\left(\mathfrak{z}\right) & =\frac{\beta_{H}\left(0\right)\left(\alpha_{H}\left(0\right)\right)^{N}}{\alpha_{H}\left(0\right)-1}+\frac{\beta_{H}\left(0\right)}{\alpha_{H}\left(0\right)-1}\sum_{\left|t\right|_{p}=p^{N}}\hat{A}_{H}\left(t\right)e^{2\pi i\left\{ t\mathfrak{z}\right\} _{p}}\\
 & +\frac{\beta_{H}\left(0\right)\left(\alpha_{H}\left(0\right)\right)^{N}}{\alpha_{H}\left(0\right)-1}\sum_{0<\left|t\right|_{p}<p^{N}}\left(\alpha_{H}\left(0\right)\right)^{v_{p}\left(t\right)}\hat{A}_{H}\left(t\right)e^{2\pi i\left\{ t\mathfrak{z}\right\} _{p}}
\end{align*}
Using (\ref{eq:Convolution of dA_H and D_N}) and \textbf{Proposition
\ref{prop:Convolution of alpha and A_H hat}}, we obtain: 
\begin{align*}
\chi_{H,N}\left(\mathfrak{z}\right)-\tilde{\chi}_{H,N}\left(\mathfrak{z}\right) & =\frac{\beta_{H}\left(0\right)\left(\alpha_{H}\left(0\right)\right)^{N}}{\alpha_{H}\left(0\right)-1}\\
 & +\frac{\beta_{H}\left(0\right)}{\alpha_{H}\left(0\right)-1}\left(\left(\frac{\mu_{0}}{p}\right)^{N}\kappa_{H}\left(\left[\mathfrak{z}\right]_{p^{N}}\right)-\alpha_{H}\left(0\right)\left(\frac{\mu_{0}}{p}\right)^{N-1}\kappa_{H}\left(\left[\mathfrak{z}\right]_{p^{N-1}}\right)\right)\\
 & +\frac{\beta_{H}\left(0\right)\left(\alpha_{H}\left(0\right)\right)^{N}}{\alpha_{H}\left(0\right)-1}\left(\left(\frac{\mu_{0}}{p\alpha_{H}\left(0\right)}\right)^{N}\kappa_{H}\left(\left[\mathfrak{z}\right]_{p^{N}}\right)-1\right)\\
 & =\frac{\beta_{H}\left(0\right)}{\alpha_{H}\left(0\right)-1}\left(\frac{\mu_{0}}{p}\right)^{N-1}\left(\frac{2\mu_{0}}{p}\kappa_{H}\left(\left[\mathfrak{z}\right]_{p^{N}}\right)-\alpha_{H}\left(0\right)\kappa_{H}\left(\left[\mathfrak{z}\right]_{p^{N-1}}\right)\right)
\end{align*}

Q.E.D.

\vphantom{}

The result just proved is of independent interest, seeing as it provides
us with an explicit form for the error between $\chi_{H,N}\left(\mathfrak{z}\right)$
and $\tilde{\chi}_{H,N}\left(\mathfrak{z}\right)$. A significant
issue can be seen by examining (\ref{eq:Chi_H,N / Chi_H,N twiddle identity for alpha equals 1})
and (\ref{eq:Chi_H,N / Chi_H,N twiddle identity for alpha is not 1})
for $\mathfrak{z}=0$. Since $\chi_{H,N}\left(0\right)=\chi_{H}\left(\left[0\right]_{p^{N}}\right)=0$
for all $N\geq0$, we obtain: 
\begin{equation}
\tilde{\chi}_{H,N}\left(0\right)\overset{\overline{\mathbb{Q}}}{=}\begin{cases}
\beta_{H}\left(0\right)\left(\left(\frac{2\mu_{0}}{p}-1\right)N+1\right)\left(\frac{\mu_{0}}{p}\right)^{N-1} & \textrm{if }\alpha_{H}\left(0\right)=1\\
\frac{\beta_{H}\left(0\right)}{\alpha_{H}\left(0\right)-1}\left(\frac{2\mu_{0}}{p}-\alpha_{H}\left(0\right)\right)\left(\frac{\mu_{0}}{p}\right)^{N-1} & \textrm{if }\alpha_{H}\left(0\right)\neq1
\end{cases},\textrm{ }\forall N\geq1\label{eq:Formula for Chi_H,N twiddle at 0}
\end{equation}
When $\alpha_{H}\left(0\right)=1$, observe that $\tilde{\chi}_{H,N}\left(0\right)$
is non-zero for all $N\geq1$, with it tending to $0$ in $\mathbb{R}$
as $N\rightarrow\infty$. The same is true when $\alpha_{H}\left(0\right)\neq1$,
unless $2\mu_{0}=p\alpha_{H}\left(0\right)$. Because the $N$th truncation
$\chi_{H,N}\left(\mathfrak{z}\right)$ is continuous and vanishes
for $\mathfrak{z}\overset{p^{N}}{\equiv}0$, the WTT for continuous
$\left(p,q\right)$-adic functions tells us that $\hat{\chi}_{H,N}\left(t\right)$
does not have a convolution inverse. However, this statement need
not be true for $\tilde{\chi}_{H,N}\left(\mathfrak{z}\right)$.

That being said, using the previous proposition, we can establish
the $L_{\mathbb{R}}^{1}$ convergence of $\tilde{\chi}_{H,N}$ to
$\chi_{H,N}$ and vice-versa. We just need one last computation under
our belts: the $L_{\mathbb{R}}^{1}$-norm of $\kappa_{H,N}\left(\mathfrak{z}\right)$ 
\begin{prop}[\textbf{$L_{\mathbb{R}}^{1}$-norm of $\kappa_{H,N}$}]
\label{prop:L^1_R norm of kapp_H,N}\index{kappa{H}@$\kappa_{H}$!truncation and L{mathbb{R}}{1}-norm@truncation
and $L_{\mathbb{R}}^{1}$-norm}Let $H$ be a contracting semi-basic $p$-Hydra map, and let $q=q_{H}$.
Then: 
\begin{equation}
\int_{\mathbb{Z}_{p}}\left|\kappa_{H,N}\left(\mathfrak{z}\right)\right|_{q}d\mathfrak{z}=\int_{\mathbb{Z}_{p}}\left|\kappa_{H}\left(\left[\mathfrak{z}\right]_{p^{N}}\right)\right|_{q}d\mathfrak{z}\leq\left(\frac{p+q-1}{pq}\right)^{N},\textrm{ }\forall N\geq1\label{eq:L^1_R norm estimate of Nth truncation of Kappa_H}
\end{equation}
In particular, since $p+q-1<pq$ for all $p,q\geq2$, we have that
the $N$th truncations of $\kappa_{H}$ converge to $0$ in $L_{\mathbb{R}}^{1}\left(\mathbb{Z}_{p},\mathbb{C}_{q}\right)$
as $N\rightarrow\infty$. 
\end{prop}
Proof: Since $\left|\kappa_{H}\left(\left[\mathfrak{z}\right]_{p^{N}}\right)\right|_{q}$
is a locally constant real valued function of $\mathfrak{z}$ whose
output is determined by the value of $\mathfrak{z}$ mod $p^{N}$,
we have that: 
\begin{equation}
\int_{\mathbb{Z}_{p}}\left|\kappa_{H}\left(\left[\mathfrak{z}\right]_{p^{N}}\right)\right|_{q}d\mathfrak{z}=\frac{1}{p^{N}}\sum_{n=0}^{p^{N}-1}\left|\kappa_{H}\left(n\right)\right|_{q}
\end{equation}
Once again, we proceed recursively. Let: 
\begin{equation}
S_{N}=\frac{1}{p^{N}}\sum_{n=0}^{p^{N}-1}\left|\kappa_{H}\left(n\right)\right|_{q}
\end{equation}
where, note, $S_{0}=\left|\kappa_{H}\left(0\right)\right|_{q}=\left|1\right|_{q}=1$.
Then, splitting the $n$-sum modulo $p$: 
\begin{align*}
S_{N} & =\frac{1}{p^{N}}\sum_{n=0}^{p^{N-1}-1}\sum_{j=0}^{p-1}\left|\kappa_{H}\left(pn+j\right)\right|_{q}\\
 & =\frac{1}{p^{N}}\sum_{n=0}^{p^{N-1}-1}\sum_{j=0}^{p-1}\left|\frac{\mu_{j}}{\mu_{0}}\kappa_{H}\left(n\right)\right|_{q}\\
 & =\left(\sum_{j=0}^{p-1}\left|\frac{\mu_{j}}{\mu_{0}}\right|_{q}\right)\frac{1}{p^{N}}\sum_{n=0}^{p^{N-1}-1}\left|\kappa_{H}\left(n\right)\right|_{q}\\
 & =\left(\frac{1}{p}\sum_{j=0}^{p-1}\left|\frac{\mu_{j}}{\mu_{0}}\right|_{q}\right)S_{N-1}
\end{align*}
Hence: 
\begin{equation}
S_{N}=\left(\frac{1}{p}\sum_{j=0}^{p-1}\left|\frac{\mu_{j}}{\mu_{0}}\right|_{q}\right)^{N}S_{0}=\left(\frac{1}{p}\sum_{j=0}^{p-1}\left|\frac{\mu_{j}}{\mu_{0}}\right|_{q}\right)^{N}
\end{equation}
Since $H$ is given to be semi-basic, we have that: 
\begin{equation}
\left|\frac{\mu_{j}}{\mu_{0}}\right|_{q}\leq\frac{1}{q},\textrm{ }\forall j\in\left\{ 1,\ldots,q-1\right\} 
\end{equation}
and so: 
\[
S_{N}\leq\left(\frac{1}{p}\left(1+\sum_{j=1}^{p-1}\frac{1}{q}\right)\right)^{N}=\left(\frac{p+q-1}{pq}\right)^{N}
\]

Q.E.D.

\vphantom{}

Here, then, is the desired theorem:
\begin{thm}
\label{thm:L^1 convergence of Chi_H,N minus Chi_H,N twiddle}Let $H$
be a contracting semi-basic $p$-Hydra map. Then: 
\begin{equation}
\lim_{N\rightarrow\infty}\int_{\mathbb{Z}_{p}}\left|\chi_{H,N}\left(\mathfrak{z}\right)-\tilde{\chi}_{H,N}\left(\mathfrak{z}\right)\right|_{q}d\mathfrak{z}\overset{\mathbb{R}}{=}0\label{eq:L^1 converges of the difference between Chi_H,N and Chi_H,N twiddle}
\end{equation}
\end{thm}
Proof: Regardless of whether or not $\alpha_{H}\left(0\right)=1$,
the co-primality of $p$ and $q$ along with an application of the
\emph{ordinary} $q$-adic triangle inequality, (\ref{eq:Chi_H,N / Chi_H,N twiddle identity for alpha equals 1})
and (\ref{eq:Chi_H,N / Chi_H,N twiddle identity for alpha is not 1})
from \textbf{Theorem \ref{thm:Chi_H,N / Chi_H,N twiddle error}} both
become: 
\begin{align*}
\int_{\mathbb{Z}_{p}}\left|\chi_{H,N}\left(\mathfrak{z}\right)-\tilde{\chi}_{H,N}\left(\mathfrak{z}\right)\right|_{q}d\mathfrak{z} & \ll\int_{\mathbb{Z}_{p}}\left|\kappa_{H}\left(\left[\mathfrak{z}\right]_{p^{N}}\right)\right|_{q}d\mathfrak{z}+\int_{\mathbb{Z}_{p}}\left|\kappa_{H}\left(\left[\mathfrak{z}\right]_{p^{N-1}}\right)\right|_{q}d\mathfrak{z}
\end{align*}
By \textbf{Proposition \ref{prop:L^1_R norm of kapp_H,N}}, both of
the terms in the upper bound go to $0$ as $N\rightarrow\infty$.

Q.E.D.

\subsection{\label{subsec:4.3.3 Quick-Approach-of}\index{quick approach}Quick
Approach of $\chi_{H}$}

In Subsection \ref{subsec:3.2.2 Truncations-=00003D000026-The}, we
proved the \textbf{Square Root Lemma }(page \pageref{eq:Square Root Lemma}),
which\textemdash recall\textemdash asserted that for any $\chi\in\tilde{C}\left(\mathbb{Z}_{p},K\right)$,
any $\mathfrak{z}\in\mathbb{Z}_{p}^{\prime}$, and any $\mathfrak{c}\in K\backslash\chi\left(\mathbb{N}_{0}\right)$,
the equality $\chi\left(\mathfrak{z}\right)=\mathfrak{c}$ held if
and only if: 
\begin{equation}
\liminf_{n\rightarrow\infty}\frac{\left|\chi\left(\left[\mathfrak{z}\right]_{p^{n}}\right)-\mathfrak{c}\right|_{q}}{\left|\nabla_{p^{n}}\left\{ \chi\right\} \left(\mathfrak{z}\right)\right|_{q}^{1/2}}<\infty
\end{equation}
With that in mind, we then said the pair $\left(\mathfrak{z},\mathfrak{c}\right)$
was quickly (resp. slowly) approached by $\chi$ whenever the above
$\liminf$ was $0$ (resp. $\in\left(0,\infty\right)$). Below, we
prove that the pair $\left(\mathfrak{z},\chi_{H}\left(\mathfrak{z}\right)\right)$
is approached quickly by $\chi_{H}$ for all $\mathfrak{z}\in\mathbb{Z}_{p}$
whenever $H$ is a contracting, semi-basic $p$-Hydra map which fixes
$0$. 
\begin{thm}[\textbf{\textit{Quick Approach of $\chi_{H}$ on $\mathbb{Z}_{p}^{\prime}$}}]
Let $H$ be a contracting, semi-basic $p$-Hydra map which fixes
$0$, and write $q=q_{H}$. Then: 
\begin{equation}
\lim_{n\rightarrow\infty}\frac{\left|\chi_{H}\left(\left[\mathfrak{z}\right]_{p^{n}}\right)-\chi_{H}\left(\mathfrak{z}\right)\right|_{q}}{\left|M_{H}\left(\left[\mathfrak{z}\right]_{p^{n}}\right)\right|_{q}^{1/2}}=0,\textrm{ }\forall\mathfrak{z}\in\mathbb{Z}_{p}^{\prime}\label{eq:Quickness of a p-Hydra map on Z_p prime}
\end{equation}
\end{thm}
Proof: Suppose that $\alpha_{H}\left(0\right)=1$, and fix $\mathfrak{z}\in\mathbb{Z}_{p}^{\prime}$.
We wish to investigate $\chi_{H}\left(\mathfrak{z}\right)-\chi_{H}\left(\left[\mathfrak{z}\right]_{p^{N}}\right)$,
or\textemdash which is the same\textemdash $\chi_{H}\left(\mathfrak{z}\right)-\chi_{H,N}\left(\mathfrak{z}\right)$.
To do this, we write: 
\begin{equation}
\chi_{H}\left(\mathfrak{z}\right)-\chi_{H,N}\left(\mathfrak{z}\right)=\left(\chi_{H}\left(\mathfrak{z}\right)-\tilde{\chi}_{H,N}\left(\mathfrak{z}\right)\right)-\left(\chi_{H,N}\left(\mathfrak{z}\right)-\tilde{\chi}_{H,N}\left(\mathfrak{z}\right)\right)
\end{equation}
Subtracting \textbf{Corollary \ref{cor:Chi_H,N twiddle explicit formula, arbitrary p, arbitrary alpha}}'s
formula for $\tilde{\chi}_{H,N}\left(\mathfrak{z}\right)$) from \textbf{Corollary
\ref{cor:F-series for Chi_H, arbitrary p and alpha}}'s $\mathcal{F}$-series
formula for $\chi_{H}\left(\mathfrak{z}\right)$) gives: 
\begin{align*}
\chi_{H}\left(\mathfrak{z}\right)-\tilde{\chi}_{H,N}\left(\mathfrak{z}\right) & \overset{\mathcal{F}_{p,q_{H}}}{=}\beta_{H}\left(0\right)N\left(\frac{\mu_{0}}{p}\right)^{N}\kappa_{H}\left(\left[\mathfrak{z}\right]_{p^{N}}\right)\\
 & +\beta_{H}\left(0\right)\sum_{n=N}^{\infty}\left(\frac{\mu_{0}}{p}\right)^{n}\kappa_{H}\left(\left[\mathfrak{z}\right]_{p^{n}}\right)\\
 & +\sum_{n=N}^{\infty}\left(\sum_{j=1}^{p-1}\beta_{H}\left(\frac{j}{p}\right)\left(\varepsilon_{n}\left(\mathfrak{z}\right)\right)^{j}\right)\left(\frac{\mu_{0}}{p}\right)^{n}\kappa_{H}\left(\left[\mathfrak{z}\right]_{p^{n}}\right)
\end{align*}
Since $\mathfrak{z}\in\mathbb{Z}_{p}^{\prime}$, the topology of convergence
as $N\rightarrow\infty$ is the $q$-adic topology, and as such, estimating
with the $q$-adic absolute value gives: 
\begin{equation}
\left|\chi_{H}\left(\mathfrak{z}\right)-\tilde{\chi}_{H,N}\left(\mathfrak{z}\right)\right|_{q}\ll\left|\kappa_{H}\left(\left[\mathfrak{z}\right]_{p^{N}}\right)\right|_{q}
\end{equation}
On the other hand, because $\alpha_{H}\left(0\right)=1$,\textbf{
Theorem \ref{thm:Chi_H,N / Chi_H,N twiddle error}} tells us that:
\begin{equation}
\chi_{H,N}\left(\mathfrak{z}\right)-\tilde{\chi}_{H,N}\left(\mathfrak{z}\right)=\beta_{H}\left(0\right)\left(\frac{\mu_{0}}{p}\right)^{N-1}\left(\frac{2\mu_{0}N}{p}\kappa_{H}\left(\left[\mathfrak{z}\right]_{p^{N}}\right)-\left(N-1\right)\kappa_{H}\left(\left[\mathfrak{z}\right]_{p^{N-1}}\right)\right)
\end{equation}
and so: 
\begin{equation}
\left|\chi_{H,N}\left(\mathfrak{z}\right)-\tilde{\chi}_{H,N}\left(\mathfrak{z}\right)\right|_{q}\ll\left|\kappa_{H}\left(\left[\mathfrak{z}\right]_{p^{N-1}}\right)\right|_{q}
\end{equation}
Since $\mathfrak{z}\in\mathbb{Z}_{p}^{\prime}$, $\left|\kappa_{H}\left(\left[\mathfrak{z}\right]_{p^{N-1}}\right)\right|_{q}\geq\left|\kappa_{H}\left(\left[\mathfrak{z}\right]_{p^{N}}\right)\right|_{q}$,
which leaves us with: 
\[
\max\left\{ \left|\chi_{H}\left(\mathfrak{z}\right)-\tilde{\chi}_{H,N}\left(\mathfrak{z}\right)\right|_{q},\left|\chi_{H,N}\left(\mathfrak{z}\right)-\tilde{\chi}_{H,N}\left(\mathfrak{z}\right)\right|_{q}\right\} \ll\left|\kappa_{H}\left(\left[\mathfrak{z}\right]_{p^{N-1}}\right)\right|_{q}
\]
Hence, by the ultrametric inequality: 
\begin{equation}
\left|\chi_{H}\left(\mathfrak{z}\right)-\chi_{H,N}\left(\mathfrak{z}\right)\right|_{q}\ll\left|\kappa_{H}\left(\left[\mathfrak{z}\right]_{p^{N-1}}\right)\right|_{q}
\end{equation}
Since: 
\begin{equation}
\kappa_{H}\left(m\right)=M_{H}\left(m\right)\left(\frac{\mu_{0}}{p}\right)^{-\lambda_{p}\left(m\right)}
\end{equation}
we have that: 
\begin{equation}
\left|\kappa_{H}\left(m\right)\right|_{q}=\left|M_{H}\left(m\right)\left(\frac{\mu_{0}}{p}\right)^{-\lambda_{p}\left(m\right)}\right|_{q}=\left|M_{H}\left(m\right)\right|_{q}
\end{equation}
and so: 
\begin{equation}
\left|\chi_{H}\left(\mathfrak{z}\right)-\chi_{H,N}\left(\mathfrak{z}\right)\right|_{q}\ll\left|M_{H}\left(\left[\mathfrak{z}\right]_{p^{N-1}}\right)\right|_{q}
\end{equation}
Dividing by $\left|M_{H}\left(\left[\mathfrak{z}\right]_{p^{N}}\right)\right|_{q}^{1/2}$
gives: 
\begin{equation}
\frac{\left|\chi_{H}\left(\mathfrak{z}\right)-\chi_{H,N}\left(\mathfrak{z}\right)\right|_{q}}{\left|M_{H}\left(\left[\mathfrak{z}\right]_{p^{N}}\right)\right|_{q}^{1/2}}\ll\frac{\left|M_{H}\left(\left[\mathfrak{z}\right]_{p^{N-1}}\right)\right|_{q}}{\left|M_{H}\left(\left[\mathfrak{z}\right]_{p^{N}}\right)\right|_{q}^{1/2}}
\end{equation}
Since the string of $p$-adic digits of $\left[\mathfrak{z}\right]_{p^{N}}$
is either equal to that of $\left[\mathfrak{z}\right]_{p^{N-1}}$
or differs by a single digit on the right, $M_{H}$'s functional equations
(\textbf{Proposition \ref{prop:M_H concatenation identity}}) tell
us that: 
\begin{equation}
M_{H}\left(\left[\mathfrak{z}\right]_{p^{N}}\right)=\begin{cases}
M_{H}\left(\left[\mathfrak{z}\right]_{p^{N-1}}\right) & \textrm{if }\left[\mathfrak{z}\right]_{p^{N-1}}=\left[\mathfrak{z}\right]_{p^{N}}\\
M_{H}\left(\left[\mathfrak{z}\right]_{p^{N-1}}\right)\cdot M_{H}\left(\textrm{coefficient of }p^{N}\textrm{ in }\mathfrak{z}\right) & \textrm{if }\left[\mathfrak{z}\right]_{p^{N-1}}\neq\left[\mathfrak{z}\right]_{p^{N}}
\end{cases}
\end{equation}
Seeing $\left[\mathfrak{z}\right]_{p^{N-1}}\neq\left[\mathfrak{z}\right]_{p^{N}}$
occurs if and only if the coefficient of $p^{N}$ in $\mathfrak{z}$
is non-zero, the semi-basicness of $H$ then guarantees that: 
\begin{equation}
\left|M_{H}\left(\textrm{coefficient of }p^{N}\textrm{ in }\mathfrak{z}\right)\right|_{p}\leq\max_{1\leq j\leq p-1}\left|M_{H}\left(j\right)\right|_{p}=p^{-\nu}
\end{equation}
occurs for some integer constant $\nu\geq1$. Hence: 
\begin{equation}
\left|M_{H}\left(\left[\mathfrak{z}\right]_{p^{N}}\right)\right|_{q}\geq p^{-\nu}\left|M_{H}\left(\left[\mathfrak{z}\right]_{p^{N-1}}\right)\right|_{q}
\end{equation}
and so:
\begin{align*}
\frac{\left|\chi_{H}\left(\mathfrak{z}\right)-\chi_{H,N}\left(\mathfrak{z}\right)\right|_{q}}{\left|M_{H}\left(\left[\mathfrak{z}\right]_{p^{N}}\right)\right|_{q}^{1/2}} & \ll\frac{\left|M_{H}\left(\left[\mathfrak{z}\right]_{p^{N-1}}\right)\right|_{q}}{\left|M_{H}\left(\left[\mathfrak{z}\right]_{p^{N}}\right)\right|_{q}^{1/2}}\\
 & \leq\frac{\left|M_{H}\left(\left[\mathfrak{z}\right]_{p^{N-1}}\right)\right|_{q}}{p^{-\frac{\nu}{2}}\left|M_{H}\left(\left[\mathfrak{z}\right]_{p^{N-1}}\right)\right|_{q}^{1/2}}\\
 & =p^{\nu/2}\left|M_{H}\left(\left[\mathfrak{z}\right]_{p^{N-1}}\right)\right|_{q}^{1/2}
\end{align*}
which tends to $0$ as $N\rightarrow\infty$, since $\mathfrak{z}\in\mathbb{Z}_{p}^{\prime}$.

The case where $\alpha_{H}\left(0\right)\neq1$ works in almost exactly
the same way.

Q.E.D.

\subsection{\label{subsec:4.3.4 Archimedean-Estimates}Archimedean Estimates}

As we discussed in Chapter 1, mixing up convergence in $\mathbb{Q}_{p}$
or $\mathbb{C}_{p}$ with convergence in $\mathbb{R}$ or $\mathbb{C}$
is, quite literally, the oldest mistake in the big book of $p$-adic
analysis. As a rule, mixing topologies in this way is fraught with
danger. So far, however, the formalism of frames has been a helpful
guard rail, keeping us on track. However, in this subsection, we are
going to simply dive off the edge of the cliff head-first.

Let $\hat{\eta}:\hat{\mathbb{Z}}_{p}\rightarrow\overline{\mathbb{Q}}$
be an element of $c_{0}\left(\hat{\mathbb{Z}}_{p},\mathbb{C}_{q}\right)$,
so that the $\left(p,q\right)$-adic Fourier series: 
\begin{equation}
\eta\left(\mathfrak{z}\right)=\sum_{t\in\hat{\mathbb{Z}}_{p}}\hat{\eta}\left(t\right)e^{2\pi i\left\{ t\mathfrak{z}\right\} _{p}}\label{eq:Getting Ready to Jump off the Cliff}
\end{equation}
converges in $\mathbb{C}_{q}$ uniformly with respect to $\mathfrak{z}\in\mathbb{Z}_{p}$.
In particular, for any $\mathfrak{z}$, the ultrametric structure
of $\mathbb{C}_{q}$ allows us to freely group and rearrange the terms
of (\ref{eq:Getting Ready to Jump off the Cliff}) however we please
without affecting the convergence. Because every term of the series
is, technically, an element of $\overline{\mathbb{Q}}$, note that
the $N$th partial sums of this series are perfectly well-defined\emph{
}complex-valued functions (complex, not $q$-adic complex) on $\mathbb{Z}_{p}$.
Moreover, it may just so happen that there are $\mathfrak{z}_{0}\in\mathbb{Z}_{p}$
so that the partial sums: 
\begin{equation}
\tilde{\eta}_{N}\left(\mathfrak{z}_{0}\right)=\sum_{\left|t\right|_{p}\leq p^{N}}\hat{\eta}\left(t\right)e^{2\pi i\left\{ t\mathfrak{z}_{0}\right\} _{p}}\label{eq:Jumping off the Cliff}
\end{equation}
converge to a limit in $\mathbb{C}$ as $N\rightarrow\infty$, in
addition to the known limit to which they converge in $\mathbb{C}_{q}$.
As we saw when discussing Hensel's error, \textbf{\emph{we cannot,
in general, assume that the limit of }}(\ref{eq:Jumping off the Cliff})\textbf{\emph{
in $\mathbb{C}$ has, in any way, shape, or form, a meaningful relation
to its limit in $\mathbb{C}_{q}$}}\textemdash \textbf{except}, of
course, when (\ref{eq:Jumping off the Cliff}) happens to be a \textbf{\emph{geometric
series}}, thanks to the \emph{universality }of the geometric series\index{geometric series universality}
(page \pageref{fact:Geometric series universality}). Specifically,
for any $\mathfrak{z}_{0}$, we can characterize the limiting behavior
in $\mathbb{C}$ (not $\mathbb{C}_{q}$!) of (\ref{eq:Jumping off the Cliff})
as $N\rightarrow\infty$ like so: 
\begin{fact}
Let $\mathfrak{z}_{0}\in\mathbb{Z}_{p}$. Then, exactly one of the
following occurs for \emph{(\ref{eq:Jumping off the Cliff})} as $N\rightarrow\infty$:

\vphantom{}

I. \emph{(\ref{eq:Jumping off the Cliff})} fails to converge to a
limit in $\mathbb{C}$.

\vphantom{}

II. (\emph{\ref{eq:Jumping off the Cliff})} converges to a limit
in $\mathbb{C}$, and \emph{(\ref{eq:Getting Ready to Jump off the Cliff})}\textbf{
}\textbf{\emph{can}} be rearranged into a geometric series of the
form $\sum_{n=0}^{\infty}c_{n}r^{n}$ for $c_{n},r\in\overline{\mathbb{Q}}$,
where $r$ has archimedean absolute value $<1$.

\vphantom{}

III. \emph{(\ref{eq:Jumping off the Cliff})} converges to a limit
in $\mathbb{C}$, but \emph{(\ref{eq:Getting Ready to Jump off the Cliff})}\textbf{
}\textbf{\emph{cannot}} be rearranged into a geometric series of the
form $\sum_{n=0}^{\infty}c_{n}r^{n}$ for $c_{n},r\in\overline{\mathbb{Q}}$,
where $r$ has archimedean absolute value $<1$. 
\end{fact}
\vphantom{}

By the universality of the geometric series, (II) is \emph{precisely
}the case where we \emph{can }conclude that the limit of (\ref{eq:Jumping off the Cliff})
in $\mathbb{C}$ as $N\rightarrow\infty$ is the same number as the
limit of (\ref{eq:Jumping off the Cliff}) in $\mathbb{C}_{q}$ as
$N\rightarrow\infty$. For any $\mathfrak{z}$ satisfying (II), the
value of (\ref{eq:Getting Ready to Jump off the Cliff}) is the same
regardless of the topology (archimedean or $q$-adic) we use to make
sense of it. This then leads to the following observation: let $U$
be the set of all $\mathfrak{z}\in\mathbb{Z}_{p}$ for which either
(II) or (III) holds true. Then, we can define a function $\eta^{\prime}:U\rightarrow\mathbb{C}$
by the limit of (\ref{eq:Jumping off the Cliff}) in $\mathbb{C}$
(not $\mathbb{C}_{q}$!). For any $\mathfrak{z}_{0}$ in $U$ satisfying
(III), there is no reason for the complex number $\eta^{\prime}\left(\mathfrak{z}_{0}\right)$
to have anything to do with the $q$-adic complex number $\eta\left(\mathfrak{z}_{0}\right)$
defined by the limit of (\ref{eq:Jumping off the Cliff}) in $\mathbb{C}_{q}$.
However$\ldots$ for $\mathfrak{z}_{0}\in U$ satisfying (II), the
universality of the geometric series \emph{guarantees} that $\eta^{\prime}\left(\mathfrak{z}_{0}\right)=\eta\left(\mathfrak{z}_{0}\right)$,
where the equality holds \emph{simultaneously }in $\mathbb{C}$ and
in $\mathbb{C}_{q}$.

The general implementation of this idea is as following. First\textemdash of
course\textemdash some terminology: 
\begin{notation}
We write $L^{1}\left(\hat{\mathbb{Z}}_{p},\mathbb{C}\right)$ to denote
the Banach space of all complex-valued functions (not $q$-adic complex
valued!) $\hat{\eta}:\hat{\mathbb{Z}}_{p}\rightarrow\mathbb{C}$ with
finite $L^{1}$ norm: 
\begin{equation}
\left\Vert \hat{\eta}\right\Vert _{L^{1}}\overset{\textrm{def}}{=}\sum_{t\in\hat{\mathbb{Z}}_{p}}\left|\hat{\eta}\left(t\right)\right|\label{eq:Definition of complex L^1 norm on Z_p hat}
\end{equation}
For $r\in\left[1,2\right]$, we write $L^{r}\left(\mathbb{Z}_{p},\mathbb{C}\right)$
to denote the Banach space of all complex-valued functions (not $q$-adic
complex valued!) $\eta:\mathbb{Z}_{p}\rightarrow\mathbb{C}$ with
finite $L^{r}$ norm: 
\begin{equation}
\left\Vert \eta\right\Vert _{L^{r}}\overset{\textrm{def}}{=}\left(\int_{\mathbb{Z}_{p}}\left|\eta\left(\mathfrak{z}\right)\right|^{r}d\mathfrak{z}\right)^{1/r}\label{eq:Definition of complex L^r norm on Z_p}
\end{equation}
\end{notation}
\begin{defn}
Let $\chi\in\tilde{C}$$\left(\mathbb{Z}_{p},\mathbb{C}_{q}\right)$,
and suppose that $\chi\left(\mathbb{N}_{0}\right)\subseteq\overline{\mathbb{Q}}$
(so that every van der Put coefficient of $\chi$ is in $\overline{\mathbb{Q}}$).
Let $\mathfrak{z}_{0}\in\mathbb{Z}_{p}$.

\vphantom{}

I. We\index{double convergence} say $\chi$ is \textbf{doubly convergent
}at $\mathfrak{z}_{0}$ whenever the van der Put series $S_{p}\left\{ \chi\right\} \left(\mathfrak{z}_{0}\right)$
can be rearranged into a geometric series of the form $\sum_{n=0}^{\infty}c_{n}r^{n}$
where $r$ and the $c_{n}$s are elements of $\overline{\mathbb{Q}}$
so that both the archimedean and $q$-adic absolute values of $r$
are $<1$. We call $\mathfrak{z}_{0}$ a \textbf{double convergence
point }of $\chi$ whenever $\chi$ is doubly convergent at $\mathfrak{z}_{0}$.

\vphantom{}

II. We say $\chi$ is \textbf{spurious }at $\mathfrak{z}_{0}$ whenever
$S_{p}\left\{ \chi\right\} \left(\mathfrak{z}_{0}\right)$ converges
in $\mathbb{C}$, yet cannot be rearranged into a geometric series
which is $q$-adically convergent. Likewise, we call $\mathfrak{z}_{0}$
a \textbf{spurious point }of $\chi$ whenever $\chi$ is spurious
at $\mathfrak{z}_{0}$.

\vphantom{}

III. We say $\chi$ is \textbf{mixed }at $\mathfrak{z}_{0}$ whenever
$\chi$ is either doubly convergent or spurious at $\mathfrak{z}_{0}$.
Likewise, we call $\mathfrak{z}_{0}$ a \textbf{mixed point }of $\chi$
whenever $\chi$ is either doubly convergent or spurious at $\mathfrak{z}_{0}$. 
\end{defn}
\begin{rem}
Observe that the complex number value attained by $S_{p}\left\{ \chi\right\} \left(\mathfrak{z}_{0}\right)$
at a mixed point $\mathfrak{z}_{0}$ is equal to the $q$-adic value
attained by $\chi$ at $\mathfrak{z}_{0}$ \emph{if and only if} $\mathfrak{z}_{0}$
is a point of double convergence. 
\end{rem}
\vphantom{}

The next proposition shows that, however audacious this approach might
be, it is nevertheless sound in mind and body.
\begin{prop}
\label{prop:measurability proposition}$\chi\in\tilde{C}\left(\mathbb{Z}_{p},\mathbb{C}_{q}\right)$,
suppose that $\chi\left(\mathbb{N}_{0}\right)\subseteq\overline{\mathbb{Q}}$,
and let $\phi:\mathbb{Z}_{p}\rightarrow\mathbb{Z}_{p}$ be a function
which is measurable with respect to the \textbf{\emph{real-valued}}\emph{
}Haar probability measure on $\mathbb{Z}_{p}$ so that every $\mathfrak{z}\in\phi\left(\mathbb{Z}_{p}\right)$
is a mixed point of $\chi$. Then, the function $\eta:\mathbb{Z}_{p}\rightarrow\mathbb{C}$
given by: 
\begin{equation}
\eta\left(\mathfrak{z}\right)\overset{\textrm{def}}{=}\lim_{N\rightarrow\infty}S_{p,N}\left\{ \chi\right\} \left(\phi\left(\mathfrak{z}\right)\right)
\end{equation}
is measurable with respect to the \textbf{\emph{real-valued}}\emph{
}Haar probability measure on $\mathbb{Z}_{p}$. Moreover, we have
that $\eta\left(\mathfrak{z}\right)=\chi\left(\phi\left(\mathfrak{z}\right)\right)$
for all $\mathfrak{z}$ for which $\chi$ is doubly convergent at
$\phi\left(\mathfrak{z}\right)$. 
\end{prop}
Proof: Let $\chi$, $\phi$, and $\eta$ be as given. Then, for each
$N$, the $N$th partial van der Put series $S_{p,N}\left\{ \chi\right\} \left(\mathfrak{z}\right)$
is a locally-constant $\overline{\mathbb{Q}}$-valued function on
$\mathbb{Z}_{p}$, and, as such, is a \emph{continuous }$\overline{\mathbb{Q}}$-valued
function on $\mathbb{Z}_{p}$. Since $\phi$ is measurable, this means
that $\mathfrak{z}\mapsto S_{p,N}\left\{ \chi\right\} \left(\phi\left(\mathfrak{z}\right)\right)$
is a measurable $\overline{\mathbb{Q}}$-valued function on $\mathbb{Z}_{p}$.
Since $\phi\left(\mathfrak{z}\right)$ is a mixed point of $\chi$
for every $\mathfrak{z}$, the van der Put series $S_{p,N}\left\{ \chi\right\} \left(\phi\left(\mathfrak{z}\right)\right)$
converge point-wise in $\mathbb{C}$ to limit, and this limit is,
by definition, $\eta\left(\mathfrak{z}\right)$. Since $\eta$ is
the point-wise limit of a sequence of measurable complex-valued functions
on $\mathbb{Z}_{p}$, $\eta$ itself is then measurable.

Finally, note that for any $\mathfrak{z}$ for which $\phi\left(\mathfrak{z}\right)$
is a double convergent point of $\chi$, the limit of $S_{p,N}\left\{ \chi\right\} \left(\phi\left(\mathfrak{z}\right)\right)$
is a geometric series which converges simultaneously in $\mathbb{C}$
and $\mathbb{C}_{q}$. The universality of the geometric series then
guarantees that the limit of the series in $\mathbb{C}$ is equal
to the limit of the series in $\mathbb{C}_{q}$. Since $\chi$ is
rising-continuous, the $N$th partial van der Put series $S_{p,N}\left\{ \chi\right\} \left(\phi\left(\mathfrak{z}\right)\right)$
converges in $\mathbb{C}_{q}$ to $\chi\left(\phi\left(\mathfrak{z}\right)\right)$
for all $\mathfrak{z}$. So, $\eta\left(\mathfrak{z}\right)=\chi\left(\phi\left(\mathfrak{z}\right)\right)$
for all $\mathfrak{z}$ for which $\phi\left(\mathfrak{z}\right)$
is a double convergent point of $\chi$.

Q.E.D. 
\begin{thm}[\textbf{$L^{1}$ Method for Archimedean Estimates}]
\label{thm:The L^1 method} Let\index{$L^{1}$ method} $\chi\in\tilde{C}$$\left(\mathbb{Z}_{p},\mathbb{C}_{q}\right)$,
suppose that $\chi\left(\mathbb{N}_{0}\right)\subseteq\overline{\mathbb{Q}}$,
and let $\phi:\mathbb{Z}_{p}\rightarrow\mathbb{Z}_{p}$ be a function
which is measurable with respect to the \textbf{\emph{real-valued}}\emph{
}Haar probability measure on $\mathbb{Z}_{p}$ so that:

\vphantom{}

I. Every $\mathfrak{z}\in\phi\left(\mathbb{Z}_{p}\right)$ is a mixed
point of $\chi$.

\vphantom{}

II. The measurable function $\eta:\mathbb{Z}_{p}\rightarrow\mathbb{C}$
defined by: 
\begin{equation}
\eta\left(\mathfrak{z}\right)\overset{\textrm{def}}{=}\lim_{N\rightarrow\infty}S_{p,N}\left\{ \chi\right\} \left(\phi\left(\mathfrak{z}\right)\right)
\end{equation}
is in $L^{r}\left(\mathbb{Z}_{p},\mathbb{C}_{q}\right)$ for some
$r\in\left[1,2\right]$, so that the Fourier transform $\hat{\eta}:\hat{\mathbb{Z}}_{p}\rightarrow\mathbb{C}$:
\begin{equation}
\hat{\eta}\left(t\right)\overset{\textrm{def}}{=}\int_{\mathbb{Z}_{p}}\eta\left(\mathfrak{z}\right)e^{2\pi i\left\{ t\mathfrak{z}\right\} _{p}}d\mathfrak{z}
\end{equation}
is an element of $L^{1}\left(\mathbb{\hat{Z}}_{p},\mathbb{C}_{q}\right)$.

Then: 
\begin{equation}
\left|\chi\left(\phi\left(\mathfrak{z}\right)\right)\right|\leq\left\Vert \hat{\eta}\right\Vert _{L^{1}}\label{eq:The L^1 Method}
\end{equation}
for all $\mathfrak{z}\in\mathbb{Z}_{p}$ so that $\chi$ is doubly
convergent at $\phi\left(\mathfrak{z}\right)$. 
\end{thm}
\begin{rem}
For any $\mathfrak{z}\in\mathbb{Z}_{p}$ so that $\chi$ is doubly
convergent at $\phi\left(\mathfrak{z}\right)$, it necessarily follows
that the value of $\chi\left(\mathfrak{y}\right)$ at $\mathfrak{y}=\phi\left(\mathfrak{z}\right)$
is then an element of $\mathbb{C}_{q}\cap\overline{\mathbb{Q}}$.
As such, the $L^{1}$ method allows us to extract \emph{archimedean
upper bounds }on outputs of $\chi$ that happen to be both $q$-adic
complex numbers and ordinary algebraic numbers. 
\end{rem}
Proof: Let everything be as given. By \textbf{Proposition \ref{prop:measurability proposition}},
$\eta$ is measurable. It is a well-known fact that the Fourier transform
of $\eta$ (a complex-valued function on the compact abelian group
$\mathbb{Z}_{p}$) is defined whenever $\eta$ is in $L^{r}\left(\mathbb{Z}_{p},\mathbb{C}\right)$
for some $r\in\left[1,2\right]$ (see \cite{Folland - harmonic analysis},
for example). Consequently, the hypothesis that $\hat{\eta}$ be in
$L^{1}$ tells us that $\eta\left(\mathfrak{z}\right)$'s Fourier
series is absolutely convergent, and hence, that: 
\begin{equation}
\left|\eta\left(\mathfrak{z}\right)\right|=\left|\sum_{t\in\hat{\mathbb{Z}}_{p}}\hat{\eta}\left(t\right)e^{2\pi i\left\{ t\mathfrak{z}\right\} _{p}}\right|\leq\sum_{t\in\hat{\mathbb{Z}}_{p}}\left|\hat{\eta}\left(t\right)\right|=\left\Vert \hat{\eta}\right\Vert _{L^{1}}
\end{equation}
Since $\eta\left(\mathfrak{z}\right)=\chi\left(\phi\left(\mathfrak{z}\right)\right)$
for all $\mathfrak{z}\in\mathbb{Z}_{p}$ for which $\chi$ is doubly
convergent at $\phi\left(\mathfrak{z}\right)$, the result follows.

Q.E.D.

\vphantom{}

Originally, I employed this method to obtain archimedean bounds on
$\chi_{q}\left(\mathfrak{z}\right)$ for odd $q\geq3$. However, the
bounds were not very good, in that my choice of $\phi$ was relatively
poor\textemdash the exact meaning of this will be clear in a moment.
Like with most of my original work in this dissertation, the bounds
followed only after a lengthy, involved Fourier analytic computation.
As such, rather than continue to draw out this three-hundred-plus-page-long
behemoth, I hope the reader will forgive me for merely sketching the
argument.

The idea for $\phi$ comes from how the $2$-adic digits of a given
$\mathfrak{z}\in\mathbb{Z}_{2}$ affect the value of $\chi_{q}\left(\mathfrak{z}\right)$.
Consider, for instance the explicit formula for $\chi_{3}\left(\mathfrak{z}\right)$
(equation (\ref{eq:Explicit formula for Chi_3})): 
\begin{equation}
\chi_{3}\left(\mathfrak{z}\right)\overset{\mathcal{F}_{2,3}}{=}-\frac{1}{2}+\frac{1}{4}\sum_{n=0}^{\infty}\frac{3^{\#_{1}\left(\left[\mathfrak{z}\right]_{2^{n}}\right)}}{2^{n}},\textrm{ }\forall\mathfrak{z}\in\mathbb{Z}_{2}
\end{equation}
For fixed $\mathfrak{z}\in\mathbb{Z}_{2}^{\prime}$, as $n$ increases
to $n+1$, $3^{\#_{1}\left(\left[\mathfrak{z}\right]_{2^{n}}\right)}$
will increase only if the coefficient of $2^{n}$ in the $2$-adic
expansion of $\mathfrak{z}$ is $1$. If that coefficient is $0$,
the denominator $2^{n}$ in the geometric series will increase, but
the numerator will not increase. A simple analysis shows that, for
$\chi_{3}$, we can guarantee that the series will become a convergent
geometric series for $\mathfrak{z}\in\mathbb{Q}\cap\mathbb{Z}_{2}^{\prime}$
by requiring that any two $1$s in the $2$-adic digits of $\mathfrak{z}$
are separated by at least one $0$.

There is, in fact, a function $\phi:\mathbb{Z}_{2}\rightarrow\mathbb{Z}_{2}$
which can guarantee this for us, and\textemdash amusingly\textemdash we
have actually already encountered this particular $\phi$ before:
it was the function from \textbf{Example \ref{exa:p-adic differentiation is crazy}}
(page \pageref{exa:p-adic differentiation is crazy}) which altered
the $2$-adic expansion of a given $2$-adic integer like so: 
\begin{equation}
\phi_{3}\left(\sum_{n=0}^{\infty}c_{n}2^{n}\right)\overset{\textrm{def}}{=}\sum_{n=0}^{\infty}c_{n}2^{2n}\label{eq:Definition of Phi_2}
\end{equation}
As we saw in Subsection \ref{subsec:3.1.1 Some-Historical-and}, $\phi_{2}$
is injective, continuous (and hence, measurable), differentiable,
and in fact \emph{continuously} differentiable, and its derivative
is identically $0$. Observe that for a $\mathfrak{z}$ with a $2$-adic
digit sequence of $\mathfrak{z}=\centerdot c_{0}c_{1}c_{2}\ldots$,
we have: 
\begin{align*}
\phi_{2}\left(\mathfrak{z}\right) & =c_{0}+c_{1}2^{2}+c_{2}2^{4}+\cdots\\
 & =c_{0}+0\cdot2^{1}+c_{1}2^{2}+0\cdot2^{3}+c_{2}2^{4}+\cdots\\
 & =\centerdot c_{0}0c_{1}0c_{2}0\ldots
\end{align*}
Thus, $\phi_{2}\left(\mathfrak{z}\right)$ is a mixed point of $\chi_{3}$
for all $\mathfrak{z}\in\mathbb{Z}_{2}$.

To use the $L^{1}$ method, like with nearly everything else in this
dissertation, we once again appeal to functional equations\index{functional equation}.
Specifically: 
\begin{prop}
Let $d$ be an integer $\geq2$, and let $\phi_{d}:\mathbb{Z}_{p}\rightarrow\mathbb{Z}_{p}$
be defined by: 
\begin{equation}
\phi_{d}\left(\sum_{n=0}^{\infty}c_{n}p^{n}\right)\overset{\textrm{def}}{=}\sum_{n=0}^{\infty}c_{n}p^{dn}\label{eq:Definition of phi_d}
\end{equation}
for all $p$-adic integers $\sum_{n=0}^{\infty}c_{n}p^{n}$. Then:
\begin{equation}
\phi_{d}\left(p\mathfrak{z}+j\right)=p^{d}\phi_{d}\left(\mathfrak{z}\right)+j,\textrm{ }\forall\mathfrak{z}\in\mathbb{Z}_{p},\textrm{ }\forall j\in\mathbb{Z}/p\mathbb{Z}\label{eq:phi_d functional equations}
\end{equation}
\end{prop}
Proof: Let $j\in\left\{ 0,\ldots,p-1\right\} $, and let $\mathfrak{z}=\sum_{n=0}^{\infty}c_{n}p^{n}$.
Then: 
\begin{align*}
\phi_{d}\left(p\mathfrak{z}+j\right) & =\phi_{d}\left(j+\sum_{n=0}^{\infty}c_{n}p^{n+1}\right)\\
 & =j+\sum_{n=0}^{\infty}c_{n}p^{d\left(n+1\right)}\\
 & =j+p^{d}\sum_{n=0}^{\infty}c_{n}p^{dn}\\
 & =j+p^{d}\phi_{d}\left(\mathfrak{z}\right)
\end{align*}

Q.E.D.

\vphantom{}

The $\phi_{d}$s are just one example of a possible choice of $\phi$
for which we can apply the $L^{1}$ method. The goal is to choose
a $\phi$ satisfying desirable functional equations like these which
will then allow us to exploit $\chi_{H}$'s functional equation to
show that $\eta=\chi\circ\phi$ satisfies a functional equation as
well. This will allow us to recursively solve for an explicit formula
of $\hat{\text{\ensuremath{\eta}}}$ whose $L^{1}$ convergence can
then be proven by direct estimation. 
\begin{example}[\textbf{A sample application of the $L^{1}$ method to $\chi_{3}$}]
\label{exa:L^1 method example}Let $\eta=\chi_{3}\circ\phi_{2}$.
Then, combining the functional equations of $\chi_{3}$ and $\phi_{2}$
yields: 
\begin{equation}
\eta\left(2\mathfrak{z}\right)=\frac{1}{4}\eta\left(\mathfrak{z}\right)
\end{equation}
\begin{equation}
\eta\left(2\mathfrak{z}+1\right)=\frac{3\eta\left(\mathfrak{z}\right)+2}{4}
\end{equation}
Since $\eta$ will be integrable over $\mathbb{Z}_{2}$ whenever $\hat{\eta}$
is in $L^{1}$, we can proceed by formally computing $\hat{\eta}$
and then showing that it is in $L^{1}$. Formally: 
\begin{align*}
\hat{\eta}\left(t\right) & \overset{\mathbb{C}}{=}\int_{\mathbb{Z}_{2}}\eta\left(\mathfrak{z}\right)e^{-2\pi i\left\{ t\mathfrak{z}\right\} _{2}}d\mathfrak{z}\\
 & =\int_{2\mathbb{Z}_{2}}\eta\left(\mathfrak{z}\right)e^{-2\pi i\left\{ t\mathfrak{z}\right\} _{2}}d\mathfrak{z}+\int_{2\mathbb{Z}_{2}+1}\eta\left(\mathfrak{z}\right)e^{-2\pi i\left\{ t\mathfrak{z}\right\} _{2}}d\mathfrak{z}\\
\left(\mathfrak{y}=\frac{\mathfrak{z}}{2},\frac{\mathfrak{z}-1}{2}\right); & =\frac{1}{2}\int_{\mathbb{Z}_{2}}\eta\left(2\mathfrak{y}\right)e^{-2\pi i\left\{ 2t\mathfrak{y}\right\} _{2}}d\mathfrak{y}+\frac{1}{2}\int_{\mathbb{Z}_{2}}\eta\left(2\mathfrak{y}+1\right)e^{-2\pi i\left\{ t\left(2\mathfrak{y}+1\right)\right\} _{2}}d\mathfrak{y}\\
 & =\frac{1}{8}\int_{\mathbb{Z}_{2}}\eta\left(\mathfrak{y}\right)e^{-2\pi i\left\{ 2t\mathfrak{y}\right\} _{2}}d\mathfrak{y}+\frac{1}{8}\int_{\mathbb{Z}_{2}}\left(3\eta\left(\mathfrak{z}\right)+2\right)e^{-2\pi i\left\{ t\left(2\mathfrak{y}+1\right)\right\} _{2}}d\mathfrak{y}\\
 & =\frac{1}{8}\hat{\eta}\left(2t\right)+\frac{3e^{-2\pi i\left\{ t\right\} _{2}}}{8}\hat{\eta}\left(2t\right)+\frac{e^{-2\pi i\left\{ t\right\} _{2}}}{4}\mathbf{1}_{0}\left(2t\right)
\end{align*}
and hence: 
\begin{equation}
\hat{\eta}\left(t\right)\overset{\mathbb{C}}{=}\frac{1+3e^{-2\pi i\left\{ t\right\} _{2}}}{8}\hat{\eta}\left(2t\right)+\frac{e^{-2\pi i\left\{ t\right\} _{2}}}{4}\mathbf{1}_{0}\left(2t\right)
\end{equation}
Just like we did with $\hat{\chi}_{H,N}\left(t\right)$, nesting this
formula will lead to an explicit formula for $\hat{\eta}\left(t\right)$,
which can then be used to prove $\left\Vert \hat{\eta}\right\Vert _{L^{1}}<\infty$
and thereby obtain a bound on $\left|\chi_{3}\left(\phi_{2}\left(\mathfrak{z}\right)\right)\right|$
via \textbf{Theorem \ref{thm:The L^1 method}}.

The reason why this particular choice of $\phi$ was relatively poor
is because by inserting a $0$ between \emph{every }two consecutive
$2$-adic digits of $\mathfrak{z}$, it forces us to consider $\mathfrak{z}$s
associated to composition sequences of the branches of the Shortened
Collatz Map where every application of $\frac{3x+1}{2}$ is followed
by at least one application of $\frac{x}{2}$. By experimenting with
more flexible choices of $\phi$ (say, those satisfying approximate
functional equations / inequalities, rather than \emph{exact} ones),
it may be possible to obtain non-trivial archimedean bounds on $\chi_{3}$\textemdash and,
on $\chi_{H}$s in general. 
\end{example}

\chapter{\label{chap:5 The-Multi-Dimensional-Case}The Multi-Dimensional Case}

\includegraphics[scale=0.45]{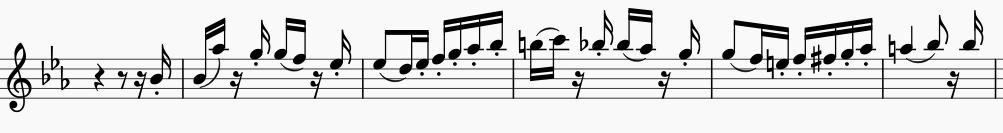}

\vphantom{}

As stated in the Introduction (Chapter 1), Chapter 5 is essentially
a compactification of the contents of Chapters 2 and 3, presenting
most of their contents as their occur in the context of multi-dimensional
Hydra maps and the techniques of multi-variable $\left(p,q\right)$-adic
analysis employed to study them. Section \ref{sec:5.1 Hydra-Maps-on}
introduces multi-dimensional Hydra maps and the two equivalent contexts
we use to study them: $\mathbb{Z}^{d}$ and $\mathcal{O}_{\mathbb{F}}$.

\section{\label{sec:5.1 Hydra-Maps-on}Hydra Maps on Lattices and Number Rings}

When attempting to generalize Hydra maps to ``higher-dimensional''
spaces, there are two possible approaches. The first is to view the
higher-dimensional analogue of $\mathbb{Z}$ as being the ring of
integers (a.k.a. \textbf{number ring}) of a given \textbf{number field
}$\mathbb{F}$\textemdash a degree $d<\infty$ extension of $\mathbb{Q}$.
As is traditional in algebraic number theory, we write $\mathcal{O}_{\mathbb{F}}$
to denote the number ring of integers of $\mathbb{F}$ (``$\mathbb{F}$-integers'').
The second approach is to work with the lattice $\mathbb{Z}^{d}$
instead of $\mathcal{O}_{\mathbb{F}}$. In theory, these two approaches
are ultimately equivalent, with $\mathcal{O}_{\mathbb{F}}$ being
$\mathbb{Z}$-module-isomorphic to $\mathbb{Z}^{d}$ after making
a choice of a basis. In practice, however, the way of number rings
introduces an irksome complication by way of the \textbf{Structure
Theorem for Finitely Generated Modules over a Principal Ideal Domain}.

Before we begin, it is informative to consider how we might translate
the one-dimensional case into the number ring setting. In the one-dimensional
case, a $p$-Hydra map $H$ was obtained by stitching together a finite
number of affine-linear maps $H_{0},\ldots,H_{p-1}$ on $\mathbb{Q}$
by defining $H\left(n\right)$ to be the image of the integer $n$
under $H_{j}$, where $j$ was the value of $n$ modulo $p$. The
essential features of this construction were: 
\begin{itemize}
\item $\mathbb{F}$, a number field; previously, this was $\mathbb{Q}$. 
\item $\mathcal{O}_{\mathbb{F}}$, the ring of integers of $\mathbb{F}$;
previously, this was $\mathbb{Z}$. 
\item \nomenclature{$\mathfrak{I}$}{a non-zero proper ideal of $\mathcal{O}_{\mathbb{F}}$}A
proper ideal\index{ideal!of mathcal{O}{mathbb{F}}@of $\mathcal{O}_{\mathbb{F}}$}
of $\mathcal{O}_{\mathbb{F}}$; previously, this was $p\mathbb{Z}$. 
\item The different co-sets of $\mathfrak{I}$ in $\mathcal{O}_{\mathbb{F}}$;
previously, these were the different equivalence classes modulo $p$. 
\item For each co-set of $\mathfrak{I}$, an affine linear map which sends
elements of that particular co-set to some (possibly different) co-set
of $\mathfrak{I}$; previously, these were the $H_{j}$s. 
\item Congruences which tell us the co-set of $\mathfrak{I}$ to which a
given element of $\mathcal{O}_{\mathbb{F}}$ belongs; previously,
these were $\overset{p}{\equiv}j$ for $j\in\mathbb{Z}/p\mathbb{Z}$. 
\end{itemize}
\begin{example}
\label{exa:root 2 example}Let $\mathbb{F}=\mathbb{Q}\left(\sqrt{2}\right)$,
$\mathcal{O}_{\mathbb{F}}=\mathbb{Z}\left[\sqrt{2}\right]$, and consider
the ideal $\mathfrak{I}=\left\langle \sqrt{2}\right\rangle $ generated
by $\sqrt{2}$ along with the quotient ring $\mathbb{Z}\left[\sqrt{2}\right]/\left\langle \sqrt{2}\right\rangle $.
We $\left\{ 1,\sqrt{2}\right\} $ as our basis for $\mathbb{Z}\left[\sqrt{2}\right]$.

Since $\left\langle 2\right\rangle $ is a proper ideal of $\left\langle \sqrt{2}\right\rangle $,
it follows that everything in $\mathbb{Z}\left[\sqrt{2}\right]$ which
is congruent to $0$ mod $2$ will also be congruent to $0$ mod $\sqrt{2}$.
Additionally, since $2\in\mathfrak{I}$, it follows that $a+b\sqrt{2}\overset{\sqrt{2}}{\equiv}a$
for all $a,b\in\mathbb{Z}$. If $a$ is even, then $a\overset{\sqrt{2}}{\equiv}0$.
If $a$ is odd, then $a\overset{2}{\equiv}1$. Moreover, since $\sqrt{2}$
divides $2$, note that $a\overset{2}{\equiv}1$ forces $a\overset{\sqrt{2}}{\equiv}1$.

Consequently, given any number $a+b\sqrt{2}\in\mathbb{Z}\left[\sqrt{2}\right]$,
the equivalence class in $\mathbb{Z}\left[\sqrt{2}\right]/\left\langle \sqrt{2}\right\rangle $
to which belongs\textemdash a.k.a, the value of $a+b\sqrt{2}$ ``mod
$\sqrt{2}$''\textemdash is then uniquely determined by the value
of $a$ mod $2$. In particular, every $a+b\sqrt{2}\in\mathbb{Z}\left[\sqrt{2}\right]$
is congruent to either $0$ or $1$ mod $\sqrt{2}$, with: 
\[
a+b\sqrt{2}\overset{\sqrt{2}}{\equiv}0\Leftrightarrow a\overset{2}{\equiv}0
\]
\[
a+b\sqrt{2}\overset{\sqrt{2}}{\equiv}1\Leftrightarrow a\overset{2}{\equiv}1
\]
As such, if we write $\mathbb{Z}\left[\sqrt{2}\right]$ in $\left\{ 1,\sqrt{2}\right\} $-coordinates,
the $2$-tuple $\mathbf{m}=\left(a,b\right)\in\mathbb{Z}^{2}$ corresponding
to $a+b\sqrt{2}\in\mathbb{Z}\left[\sqrt{2}\right]$ represents a number
congruent to $0$ (resp. $1$) mod $\sqrt{2}$ if and only if the
first component of $\mathbf{m}$ ($a$) is congruent to $0$ (resp.
$1$) mod $2$.
\end{example}
\vphantom{}

The reason \textbf{Example \ref{exa:root 2 example}} works out so
nicely is because quotienting the number ring $\mathbb{Z}\left[\sqrt{2}\right]$
by the ideal $\left\langle \sqrt{2}\right\rangle $ yields a ring
isomorphic to $\mathbb{Z}/2\mathbb{Z}$, which is cyclic as an additive
group. However, as any algebraic number theorist will readily tell
you, for a general number field $\mathbb{F}$ and an arbitrary non-zero
proper ideal $\mathfrak{I}$ of $\mathcal{O}_{\mathbb{F}}$, there
is no guarantee that the ring $\mathcal{O}_{\mathbb{F}}/\mathfrak{I}$
will be cyclic as an additive group. Rather, the Structure Theorem
for Finitely-Generated Modules over a Principal Ideal Domain tells
us that the most we can expect is for there to be an isomorphism of
additive groups:
\begin{equation}
\mathcal{O}_{\mathbb{F}}/\mathfrak{I}\cong\mathfrak{C}_{p_{1}}\times\cdots\times\mathfrak{C}_{p_{r}}\label{eq:Direct Product Representation of O_F / I}
\end{equation}
where $\mathfrak{C}_{p_{n}}$ is the cyclic group of order $p_{n}$
\nomenclature{$\mathfrak{C}_{n}$}{cyclic group of order $n$ \nopageref},
where $r\in\left\{ 1,\ldots,\left|\mathcal{O}_{\mathbb{F}}/\mathfrak{I}\right|\right\} $,
and where the $p_{n}$s are integers $\geq2$ so that $p_{n}\mid p_{n+1}$
for all $n\in\left\{ 1,\ldots,r-1\right\} $. In analogy to the one-dimensional
case, $\mathcal{O}_{\mathbb{F}}/\mathfrak{I}$ was $\mathbb{Z}/p\mathbb{Z}$,
the set to which the branch-determining parameter $j$ belonged. As
(\ref{eq:Direct Product Representation of O_F / I}) shows, however,
for the multi-dimensional case, we cannot assume that the branch-determining
parameter will take integer values like in \textbf{Example \ref{exa:root 2 example}},
where the parameter was $j\in\left\{ 0,1\right\} $. Instead, the
branch-determining parameter we be an \textbf{$r$-tuple of integers
}$\mathbf{j}\in\left(\mathbb{Z}/p_{1}\mathbb{Z}\right)\times\cdots\times\left(\mathbb{Z}/p_{r}\mathbb{Z}\right)$.
In the one-dimensional case, we obtained an extension of $\chi_{H}$
by considering infinitely long lists (strings) of parameters. In the
multi-dimensional case, an infinitely long list of $r$-tuples $\mathbf{j}\in\left(\mathbb{Z}/p_{1}\mathbb{Z}\right)\times\cdots\times\left(\mathbb{Z}/p_{r}\mathbb{Z}\right)$
would then be identified with an $r$-tuple $\mathbf{z}\in\mathbb{Z}_{p_{1}}\times\cdots\times\mathbb{Z}_{p_{r}}$
whose $n$th entry was a $p_{n}$-adic integer corresponding to the
string whose elements are the $n$th entries of the $\mathbf{j}$s.

Although this situation creates no problems when performing Fourier
analysis on multi-dimensional $\chi_{H}$, it causes trouble for our
proof of the multi-dimensional analogue of the Correspondence Principle,
because it potentially makes it impossible for us to define a $p$-adic
extension of the multi-dimensional Hydra map $H$. As such, for Subsection
\ref{sec:5.2 The-Numen-of}, we will need to restrict to the case
where every element of $P$ is a single prime, $p$.

\subsection{\label{subsec:5.1.1. Algebraic-Conventions}Algebraic Conventions}

To minimize any cumbersome aspects of notation, we introduce the following
conventions\footnote{The one exception to this will be for multi-dimensional analogues
of functions from the one-dimensional case, such as $\chi_{H}$, $\hat{A}_{H}$,
$\alpha_{H}$, and the rest. For them, I retain as much of the one-dimensional
notation as possible to emphasize the parallels between the two cases.
Indeed, most of the computations will be repeated nearly verbatim.}: 
\begin{itemize}
\item \textbf{Bold, lower case letters }(ex: $\mathbf{j}$, $\mathbf{a}$,
$\mathbf{n}$, $\mathbf{x}$, $\mathbf{z}$, etc.) are reserved to
denote tuples of finite length. In computations involving matrices,
such tuples will always be treated as column vectors. Row vectors
are written with a superscript $T$ to denote the transpose of column
vector (ex: \textbf{$\mathbf{j}^{T}$}). 
\item \textbf{BOLD, UPPER CASE LETTERS }(ex: $\mathbf{A},\mathbf{D},\mathbf{P}$,
etc.) are reserved to denote $d\times d$ matrices. In particular,
$\mathbf{D}$ will be used for \textbf{diagonal matrices}, while $\mathbf{P}$
will be used for \textbf{permutation matrices}\textemdash matrices
(whose entries are $0$s and $1$s) that give the so-called ``defining
representation'' of the symmetric group $\mathfrak{S}_{d}$ acting
on $\mathbb{R}^{d}$. 
\end{itemize}
It is not an exaggeration to say that the primary challenge of the
multi-dimensional case is its notation. As a rule of thumb, whenever
vectors are acting on vectors, reader should assume that the operations
are being done entry-wise, unless state otherwise. That being said,
there is quite a lot of notation we will have to cover. In order to
avoid symbol overload, the notation will only be introduced as needed.
The primary chunks of notation occur in the segment given below, and
then again at the start of Subsections \ref{subsec:5.3.1 Tensor-Products}
and \ref{subsec:5.4.1 Multi-Dimensional--adic-Fourier}. 
\begin{notation}[\textbf{Multi-Dimensional Notational Conventions} \index{multi-dimensional!notation}]
\label{nota:First MD notation batch}Let $d$ be an integer, let
$\mathbb{F}$ be a field, let $s\in\mathbb{F}$, and let $\mathbf{a}=\left(a_{1},\ldots,a_{d}\right)$,
$\mathbf{b}=\left(b_{1},\ldots,b_{d}\right)$, and $\mathbf{c}=\left(c_{1},\ldots,c_{d}\right)$
be elements of $\mathbb{F}^{d}$. Then:\index{multi-dimensional!notation}

\vphantom{}

I. $\mathbf{a}+\mathbf{b}\overset{\textrm{def}}{=}\left(a_{1}+b_{1},\ldots,a_{d}+b_{d}\right)$

\vphantom{}

II. $\mathbf{a}\mathbf{b}\overset{\textrm{def}}{=}\left(a_{1}b_{1},\ldots,a_{d}b_{d}\right)$

\vphantom{}

III. $\frac{\mathbf{a}}{\mathbf{b}}\overset{\textrm{def}}{=}\left(\frac{a_{1}}{b_{1}},\ldots,\frac{a_{d}}{b_{d}}\right)$

\vphantom{}

IV. $s\mathbf{a}\overset{\textrm{def}}{=}\left(sa_{1},\ldots,sa_{d}\right)$

\vphantom{}

V. $\mathbf{a}+s\overset{\textrm{def}}{=}\left(a_{1}+s,\ldots,a_{d}+s\right)$

\vphantom{}

VI. 
\[
\sum_{\mathbf{k}=\mathbf{a}}^{\mathbf{b}}\overset{\textrm{def}}{=}\sum_{k_{1}=a_{1}}^{b_{1}}\cdots\sum_{k_{d}=a_{d}}^{b_{d}}
\]

\vphantom{}

VII. $\mathbf{I}_{d}$ \nomenclature{$\mathbf{I}_{d}$}{$d\times d$ identity matrix \nopageref}
denotes the $d\times d$ identity matrix. $\mathbf{O}_{d}$\nomenclature{$\mathbf{O}_{d}$}{$d\times d$ zero matrix \nopageref}
denotes the $d\times d$ zero matrix. $\mathbf{0}$\nomenclature{$\mathbf{0}$}{zero vector \nopageref}
denotes a column vector of $0$s (a.k.a., the zero vector). Unfortunately,
the length of this $\mathbf{0}$ will usually depend on context.

\vphantom{}

VIII. Given $d\times d$ matrices $\mathbf{A}$ and $\mathbf{B}$
with entries in a field $\mathbb{F}$, if $\mathbf{B}$ is invertible,
we write: \nomenclature{$\frac{\mathbf{A}}{\mathbf{B}}$}{$\mathbf{B}^{-1}\mathbf{A}$ \nopageref}
\begin{equation}
\frac{\mathbf{A}}{\mathbf{B}}\overset{\textrm{def}}{=}\mathbf{B}^{-1}\mathbf{A}\label{eq:Matrix fraction convention}
\end{equation}
\end{notation}
\begin{defn}
We write\nomenclature{$P$}{$\left(p_{1},\ldots,p_{r}\right)$}: 
\begin{equation}
P\overset{\textrm{def}}{=}\left(p_{1},\ldots,p_{r}\right)\label{eq:Definition of Big P}
\end{equation}
be to denote $r$-tuple of integers (all of which are $\geq2$), where
$r$ is a positive integer, and where $p_{n}\mid p_{n+1}$ for all
$n\in\left\{ 1,\ldots,r-1\right\} $. We then write\nomenclature{$\mathbb{Z}^{r}/P\mathbb{Z}^{r}$}{$\overset{\textrm{def}}{=}\prod_{m=1}^{r}\left(\mathbb{Z}/p_{m}\mathbb{Z}\right)$ }
to denote the direct product of rings: 
\begin{equation}
\mathbb{Z}^{r}/P\mathbb{Z}^{r}\overset{\textrm{def}}{=}\prod_{m=1}^{r}\left(\mathbb{Z}/p_{m}\mathbb{Z}\right)
\end{equation}
That is, $\mathbb{Z}^{r}/P\mathbb{Z}^{r}$ is set of all $r$-tuples
$\mathbf{j}=\left(j_{1},\ldots,j_{r}\right)$ so that $j_{m}\in\mathbb{Z}/p_{m}\mathbb{Z}$
for all $m$, equipped with component-wise addition and multiplication
operations on the $\mathbb{Z}/p_{m}\mathbb{Z}$s. We also write $\mathbb{Z}^{r}/p\mathbb{Z}^{r}\overset{\textrm{def}}{=}\left(\mathbb{Z}/p\mathbb{Z}\right)^{r}$
to denote the case where $p_{n}=p$ (for some integer $p\geq2$) for
all $n\in\left\{ 1,\ldots,r\right\} $.

Given a $\mathbf{j}=\left(j_{1},\ldots,j_{r}\right)\in\mathbb{Z}^{r}/P\mathbb{Z}^{r}$
and any tuple $\mathbf{x}=\left(x_{1},\ldots,x_{\left|\mathbf{x}\right|}\right)$
of length $\geq1$, we write:
\begin{align*}
\mathbf{x} & \overset{P}{\equiv}\mathbf{j}\\
 & \Updownarrow\\
x_{n} & \overset{p_{n}}{\equiv}j_{n},\textrm{ }\forall n\in\left\{ 1,\ldots,\min\left\{ r,\left|\mathbf{x}\right|\right\} \right\} 
\end{align*}
When there is an integer $p\geq2$ so that $p_{n}=p$ for all $n$,
we write:
\begin{align*}
\mathbf{x} & \overset{p}{\equiv}\mathbf{j}\\
 & \Updownarrow\\
x_{n} & \overset{p}{\equiv}j_{n},\textrm{ }\forall n\in\left\{ 1,\ldots,\min\left\{ r,\left|\mathbf{x}\right|\right\} \right\} 
\end{align*}
\end{defn}

\subsection{Multi-Dimensional Hydra Maps \label{subsec:5.1.2 Co=00003D0000F6rdinates,-Half-Lattices,-and}}
\begin{defn}
Fix a number field $\mathbb{F}$ of degree $d$ over $\mathbb{Q}$,
and a non-zero proper ideal $\mathfrak{I}\subset\mathcal{O}_{\mathbb{F}}$.
\nomenclature{$\iota$}{$\left|\mathcal{O}_{\mathbb{F}}/\mathfrak{I}\right|$}We
write: 
\begin{equation}
\iota\overset{\textrm{def}}{=}\left|\mathcal{O}_{\mathbb{F}}/\mathfrak{I}\right|\label{eq:definition of iota_I}
\end{equation}
In algebraic terminology, $\iota$ is the \textbf{index }of\index{ideal!index of a}
$\mathfrak{I}$ in\emph{ $\mathcal{O}_{\mathbb{F}}$. }As discussed
above, there is an isomorphism of additive groups:
\begin{equation}
\mathcal{O}_{\mathbb{F}}/\mathfrak{I}\cong\mathfrak{C}_{p_{1}}\times\cdots\times\mathfrak{C}_{p_{r}}
\end{equation}
where $\mathfrak{C}_{p_{i}}$ is the cyclic group of order $p_{i}$,
where $r\in\left\{ 1,\ldots,\iota\right\} $, where $p_{n}\mid p_{n+1}$
for all $n\in\left\{ 1,\ldots,r-1\right\} $ and: 
\begin{equation}
\prod_{n=1}^{r}p_{n}=\iota\label{eq:Relation between the rho_is and iota_I}
\end{equation}
I call $r$ the \textbf{depth}\footnote{From an algebraic perspective, the decomposition (\ref{eq:Direct Product Representation of O_F / I-1})
is relatively unremarkable, being a specific case of the structure
theorem for finitely generated modules over a principal ideal domain.
To an algebraist, $r_{\mathfrak{I}}$ would be ``the number of invariant
factors of the $\mathbb{Z}$-module $\mathcal{O}_{\mathbb{F}}/\mathfrak{I}$'';
``depth'', though not standard terminology, is certainly pithier.}\textbf{ }of\index{ideal!depth of a} $\mathfrak{I}$.

In addition to the above, note that there is a set $\mathcal{B}\subset\mathcal{O}_{\mathbb{F}}$
which is a $\mathbb{Z}$-basis of the group $\left(\mathcal{O}_{\mathbb{F}},+\right)$,
so that, for any $d$-tuple $\mathbf{x}=\left(x_{1},\ldots,x_{d}\right)\in\mathbb{Z}^{d}$
representing an element $z\in\mathcal{O}_{\mathbb{F}}$, the equivalence
class of $\mathcal{O}_{\mathbb{F}}/\mathfrak{I}$ to which $z$ belongs
is then completely determined by the values $x_{1}$ mod $p_{1}$,
$x_{2}$ mod $p_{2}$, $\ldots$, $x_{r_{\mathfrak{I}}}$ mod $p_{r}$.
\emph{From here on out, we will work with a fixed choice of such a
basis $\mathcal{B}$}.

We write $\mathbf{R}$ to denote the $d\times d$ diagonal matrix:
\begin{equation}
\mathbf{R}\overset{\textrm{def}}{=}\left[\begin{array}{cccccc}
p_{1}\\
 & \ddots\\
 &  & p_{r}\\
 &  &  & 1\\
 &  &  &  & \ddots\\
 &  &  &  &  & 1
\end{array}\right]\label{eq:Definition of bold R}
\end{equation}
where there are $d-r$ $1$s after $p_{r}$. Lastly, when adding vectors
of unequal length, we do everything from left to right: 
\begin{equation}
\mathbf{v}+\mathbf{j}\overset{\textrm{def}}{=}\left(v_{1}+j_{1},\ldots,v_{r}+j_{r},v_{r+1}+0,\ldots,v_{d}+0\right),\textrm{ }\forall\mathbf{v}\in\mathbb{Z}^{d},\textrm{ }\forall\mathbf{j}\in\mathbb{Z}^{r}/P\mathbb{Z}^{r}\label{eq:Definition of the sum of bold_m and bold_j}
\end{equation}
\end{defn}
\begin{defn}
Letting everything be as described above, enumerate the elements of
the basis $\mathcal{B}$ as $\left\{ \gamma_{1},\ldots,\gamma_{d}\right\} $.
We then write $\mathcal{O}_{\mathbb{F},\mathcal{B}}^{+}$ to denote
the \textbf{half-lattice}\index{half-lattice}\textbf{ of }$\mathcal{B}$\textbf{-positive
$\mathbb{F}$-integers} (or ``positive half-lattice''), defined
by: 
\begin{equation}
\mathcal{O}_{\mathbb{F},\mathcal{B}}^{+}\overset{\textrm{def}}{=}\left\{ \sum_{\ell=1}^{d}x_{\ell}\gamma_{\ell}:x_{1},\ldots,x_{d}\in\mathbb{N}_{0}\right\} \label{eq:Def of O-plus_F,B}
\end{equation}
The \textbf{half-lattice of $\mathcal{B}$-negative $\mathbb{F}$-integers
}(or ``negative half-lattice''), denoted $\mathcal{O}_{\mathbb{F},\mathcal{B}}^{-}$,
is defined by: 
\begin{equation}
\mathcal{O}_{\mathbb{F},\mathcal{B}}^{-}\overset{\textrm{def}}{=}\left\{ -\sum_{\ell=1}^{d}x_{\ell}\gamma_{\ell}:x_{1},\ldots,x_{d}\in\mathbb{N}_{0}\right\} \label{eq:Def of O_F,B minus}
\end{equation}
Finally, we write \nomenclature{$\varphi_{\mathcal{B}}$}{ }$\varphi_{\mathcal{B}}:\mathbb{Q}^{d}\rightarrow\mathbb{F}$
to denote the vector space isomorphism that sends each $\mathbf{x}\in\mathbb{Q}^{d}$
to the unique element of $\mathbb{F}$ whose representation in $\mathcal{B}$-coördinates
is $\mathbf{x}$, and which satisfies the property that the restriction
$\varphi_{\mathcal{B}}\mid_{\mathbb{Z}^{d}}$ of $\varphi_{\mathcal{B}}$
to $\mathbb{Z}^{d}$ is an isomorphism of the groups $\left(\mathbb{Z}^{d},+\right)$
and $\left(\mathcal{O}_{\mathbb{F}},+\right)$. We write $\varphi_{\mathcal{B}}^{-1}$
to denote the inverse of $\varphi_{\mathcal{B}}$; $\varphi_{\mathcal{B}}^{-1}$
outputs the $d$-tuple $\mathcal{B}$-coördinate representation of
the inputted $z\in\mathbb{F}$.

We then write $\mathbb{Z}_{\mathbb{F},\mathcal{B}}^{+}$ \nomenclature{$\mathbb{Z}_{\mathbb{F},\mathcal{B}}^{+}$}{$\varphi_{\mathcal{B}}\left(\mathcal{O}_{\mathbb{F},\mathcal{B}}^{+}\right)$}
and $\mathbb{Z}_{\mathbb{F},\mathcal{B}}^{-}$ \nomenclature{$\mathbb{Z}_{\mathbb{F},\mathcal{B}}^{-}$}{$\varphi_{\mathcal{B}}\left(\mathcal{O}_{\mathbb{F},\mathcal{B}}^{-}\right)$}
to denote $\varphi_{\mathcal{B}}\left(\mathcal{O}_{\mathbb{F},\mathcal{B}}^{+}\right)\subset\mathbb{Z}^{d}$
and $\varphi_{\mathcal{B}}\left(\mathcal{O}_{\mathbb{F},\mathcal{B}}^{-}\right)\subset\mathbb{Z}^{d}$,
respectively. Note that, as subsets of $\mathbb{Z}^{d}$, we have:
\begin{align*}
\mathbb{Z}_{\mathbb{F},\mathcal{B}}^{+} & =\mathbb{N}_{0}^{d}\\
\mathbb{Z}_{\mathbb{F},\mathcal{B}}^{-} & =-\mathbb{N}_{0}^{d}
\end{align*}
\end{defn}
\begin{defn}[\textbf{$\left(\mathbb{F},\mathfrak{I},\mathcal{B}\right)$-Hydra
maps}]
\label{def:(F,I,B)-Hydra map}Let $\mathbb{F}$ be a number field
of dimension $d$, let $\mathfrak{I}$ be a non-zero proper ideal
of $\mathcal{O}_{\mathbb{F}}$ of index $\iota$, and let $\mathcal{B}$
be a basis as discussed above. We write $\mathfrak{I}_{0},\ldots,\mathfrak{I}_{\iota}$
to denote the co-sets of $\mathfrak{I}$ in $\mathcal{O}_{\mathbb{F}}$,
with $\mathfrak{I}_{0}$ denoting $\mathfrak{I}$ itself. Then, a
\textbf{$\left(\mathbb{F},\mathfrak{I},\mathcal{B}\right)$-Hydra
map} \textbf{on $\mathcal{O}_{\mathbb{F}}$} is a surjective map $\tilde{H}:\mathcal{O}_{\mathbb{F}}\rightarrow\mathcal{O}_{\mathbb{F}}$
of the form: 
\begin{equation}
\tilde{H}\left(z\right)=\begin{cases}
\frac{a_{0}z+b_{0}}{d_{0}} & \textrm{if }z\in\mathfrak{I}_{0}\\
\vdots & \vdots\\
\frac{a_{\iota-1}z+b_{\iota-1}}{d_{\iota-1}} & \textrm{if }z\in\mathfrak{I}_{\iota}
\end{cases}\label{eq:Definition of I-hydra map}
\end{equation}
We call $r$ the \textbf{depth }of $\tilde{H}$. Here, the $a_{j}$s,
$b_{j}$s, and $d_{j}$s are elements of $\mathcal{O}_{\mathbb{F}}$
so that:

\vphantom{}

I. $a_{j},d_{j}\neq0$ for all $j\in\left\{ 0,\ldots,\iota-1\right\} $.

\vphantom{}

II. $\gcd\left(a_{j},d_{j}\right)=1$ for all $j\in\left\{ 0,\ldots,\iota-1\right\} $.

\vphantom{}

III. $H\left(\mathcal{O}_{\mathbb{F},\mathcal{B}}^{+}\right)\subseteq H\left(\mathcal{O}_{\mathbb{F},\mathcal{B}}^{+}\right)$
and $H\left(\mathcal{O}_{\mathbb{F},\mathcal{B}}^{-}\backslash\mathcal{O}_{\mathbb{F},\mathcal{B}}^{+}\right)\subseteq H\left(\mathcal{O}_{\mathbb{F},\mathcal{B}}^{-}\backslash\mathcal{O}_{\mathbb{F},\mathcal{B}}^{+}\right)$,
where $\mathcal{O}_{\mathbb{F},\mathcal{B}}^{-}\backslash\mathcal{O}_{\mathbb{F},\mathcal{B}}^{+}$
is the set of all elements of $\mathcal{O}_{\mathbb{F},\mathcal{B}}^{-}$
which are not elements of $\mathcal{O}_{\mathbb{F},\mathcal{B}}^{+}$.

\vphantom{}

IV. For all $j\in\left\{ 0,\ldots,\iota-1\right\} $, the ideal $\left\langle d_{j}\right\rangle _{\mathcal{O}_{\mathbb{F}}}$
in $\mathcal{O}_{\mathbb{F}}$ generated by $d_{j}$ is contained
in $\mathfrak{I}$.

\vphantom{}

V. For all $j\in\left\{ 0,\ldots,\iota-1\right\} $, the the matrix
representation in $\mathcal{B}$-coordinates on $\mathbb{Q}^{d}$
of the ``multiplication by $a_{j}/d_{j}$'' map on $\mathbb{F}$
is of the form: 
\begin{equation}
\frac{\mathbf{A}}{\mathbf{D}}\overset{\textrm{def}}{=}\mathbf{D}^{-1}\mathbf{A}\label{eq:A / D notation}
\end{equation}
where $\mathbf{A},\mathbf{D}$ are invertible $d\times d$ matrices
so that:

\vphantom{}

V-i. $\mathbf{A}=\tilde{\mathbf{A}}\mathbf{P}$, where $\mathbf{P}$
is a permutation matrix\footnote{That is, a matrix of $0$s and $1$s which is a representation of
the action on $\mathbb{Z}^{d}$ of an element of the symmetric group
on $d$ objects by way of a permutation of the coordinate entries
of the $d$-tuples in $\mathbb{Z}^{d}$.} and where $\mathbf{\tilde{\mathbf{A}}}$ is a diagonal matrix whose
non-zero entries are positive integers.

\vphantom{}

V-ii. $\mathbf{D}$ is a diagonal matrix whose non-zero entries are
positive integers such that every entry on the diagonal of $\mathbf{R}\mathbf{D}^{-1}$
is a positive integer. Note that this forces the $\left(r+1\right)$th
through $d$th diagonal entries of $\mathbf{D}$ to be equal to $1$.
Moreover, for each $\ell\in\left\{ 1,\ldots,r\right\} $, this forces
the $\ell$th entry of the diagonal of $\mathbf{D}$ to be a divisor
of $p_{\ell}$.

\vphantom{}

V-iii. Every non-zero element of $\mathbf{A}$ is co-prime to every
non-zero element of $\mathbf{D}$.

\vphantom{}

Additionally, we say that $\tilde{H}$ is \textbf{integral }if it
satisfies\index{Hydra map!integral}:

\vphantom{}

VI. For all $j\in\left\{ 0,\ldots,\iota-1\right\} $, the number $\frac{a_{j}z+b_{j}}{d_{j}}$
is an element of $\mathcal{O}_{\mathbb{F}}$ if and only if $z\in\mathcal{O}_{\mathbb{F}}\cap\mathfrak{I}_{j}$.

\vphantom{}

If $\tilde{H}$ does not satisfy (VI), we say $\tilde{H}$ is \textbf{non-integral}.
\index{Hydra map!non-integral}
\end{defn}
\begin{defn}[\textbf{$P$-Hydra maps}]
\label{def:P-Hydra map}Let $\mathbb{F}$ be a number field of dimension
$d$, let $\mathfrak{I}$ be a non-zero proper ideal of $\mathcal{O}_{\mathbb{F}}$
of index $\iota$, and let $\mathcal{B}$ be a basis as discussed
above. A $P$\textbf{-Hydra map} \textbf{on $\mathbb{Z}^{d}$}\index{$P$-Hydra map}\textbf{
}is a map $H:\mathbb{Z}^{d}\rightarrow\mathbb{Z}^{d}$ so that: 
\begin{equation}
H=\varphi_{\mathcal{B}}\circ\tilde{H}\circ\varphi_{\mathcal{B}}^{-1}\label{eq:Definition of a Hydra map on Zd}
\end{equation}
for some $\left(\mathbb{F},\mathfrak{I},\mathcal{B}\right)$-Hydra
map $\tilde{H}$ on $\mathcal{O}_{\mathbb{F}}$ for which $P=\left(p_{1},\ldots,p_{r}\right)$.
We\index{Hydra map!field analogue} \index{Hydra map!lattice analogue}call
$\tilde{H}$ the \textbf{field analogue }of $H$, and call $H$ the
\textbf{lattice analogue }of $\tilde{H}$. We say $H$ has \textbf{depth}
$r$ whenever $\tilde{H}$ has depth $r$.\index{Hydra map!depth}
\end{defn}
\begin{prop}[Formula for $P$-Hydra maps]
Let $\tilde{H}:\mathcal{O}_{\mathbb{F}}\rightarrow\mathcal{O}_{\mathbb{F}}$
be an $\left(\mathbb{F},\mathfrak{I},\mathcal{B}\right)$-Hydra map,
and let $H$ be its lattice analogue (a $P$-Hydra map). Then, there
is a unique collection of integer-entry matrices: 
\[
\left\{ \mathbf{A}_{\mathbf{j}}\right\} _{\mathbf{j}\in\mathbb{Z}^{r}/P\mathbb{Z}^{r}},\left\{ \mathbf{D}_{\mathbf{j}}\right\} _{\mathbf{j}\in\mathbb{Z}^{r}/P\mathbb{Z}^{r}}\subseteq\textrm{GL}_{d}\left(\mathbb{Q}\right)
\]
satisfying conditions \emph{(V-i)}, \emph{(V-ii), and (V-iii) }from\textbf{
}\emph{(\ref{def:(F,I,B)-Hydra map})} for each $\mathbf{j}$, respectively,
and a unique collection of $d\times1$ column vectors $\left\{ \mathbf{b}_{\mathbf{j}}\right\} _{\mathbf{j}\in\mathbb{Z}^{r}/P\mathbb{Z}^{r}}\subseteq\mathbb{Z}^{d}$
so that: \nomenclature{$\mathbf{A}_{\mathbf{j}}$}{ } \nomenclature{$\mathbf{D}_{\mathbf{j}}$}{ }
\nomenclature{$\mathbf{b}_{\mathbf{j}}$}{ } \index{Hydra map}\index{Hydra map!on mathbb{Z}{d}@on $\mathbb{Z}^{d}$}
\index{$P$-Hydra map} \index{multi-dimensional!Hydra map} 
\begin{equation}
H\left(\mathbf{x}\right)=\sum_{\mathbf{j}\in\mathbb{Z}^{r}/P\mathbb{Z}^{r}}\left[\mathbf{x}\overset{P}{\equiv}\mathbf{j}\right]\frac{\mathbf{A}_{\mathbf{j}}\mathbf{x}+\mathbf{b}_{\mathbf{j}}}{\mathbf{D}_{\mathbf{j}}},\textrm{ }\forall\mathbf{x}\in\mathbb{Z}^{d}\label{eq:MD Hydra Map Formula}
\end{equation}
In particular, for all $\mathbf{j}\in\mathbb{Z}^{r}/P\mathbb{Z}^{r}$,
the $d$-tuple $\mathbf{D}_{\mathbf{j}}^{-1}\left(\mathbf{A}_{\mathbf{j}}\mathbf{x}+\mathbf{b}_{\mathbf{j}}\right)$
will be an element of $\mathbb{Z}_{\mathbb{F},\mathcal{B}}^{+}$ for
all $\mathbf{x}\in\mathbb{Z}_{\mathbb{F},\mathcal{B}}^{+}$ for which
$\mathbf{x}\overset{P}{\equiv}\mathbf{j}$. 
\end{prop}
\begin{rem}
Recall that the congruence $\mathbf{x}\overset{P}{\equiv}\mathbf{j}$
means that for each $n\in\left\{ 1,\ldots,r\right\} $, the $n$th
entry of $\mathbf{x}$ is congruent mod $p_{n}$ to the $n$th entry
of $\mathbf{j}$. Moreover, $\mathbf{x}\overset{P}{\equiv}\mathbf{j}$
is completely independent of the $\left(r+1\right)$th through $d$th
entries of $\mathbf{x}$.
\end{rem}
Proof: Let everything be as described above. Fix $z=\sum_{n=1}^{d}c_{n}\gamma_{n}\in\mathcal{O}_{\mathbb{F},\mathcal{B}}^{+}$,
where the $c_{n}$s are non-negative integers. Letting $\mathfrak{I}_{j}$
denote the unique equivalence class of $\mathfrak{I}$ in $\mathcal{O}_{\mathbb{F}}$
to which $z$ belongs, it follows by definition of $\tilde{H}$ that:
\begin{equation}
\tilde{H}\left(z\right)=\frac{a_{j}z+b_{j}}{d_{j}}
\end{equation}
Now, consider $\mathbb{F}$ as a $d$-dimensional linear space over
$\mathbb{Q}$, equipped with the coordinate system given by the basis
$\mathcal{B}=\left\{ \gamma_{1},\ldots,\gamma_{d}\right\} $. In these
coordinates, the ``multiplication by $a_{j}/d_{j}$'' map on $\mathbb{F}$
can be uniquely represented as left-multiplication by some $\mathbf{C}\in\textrm{GL}_{d}\left(\mathbb{Q}\right)$
with rational entries. So, letting $\mathbf{x}$ denote the coordinate
$d$-tuple representing $z$ (that is, $\mathbf{x}=\varphi_{\mathcal{B}}^{-1}\left(z\right)$),
it follows that: 
\begin{equation}
\varphi_{\mathcal{B}}^{-1}\left(\frac{a_{j}z}{d_{j}}\right)=\mathbf{C}\mathbf{x}
\end{equation}
Letting $\mathbf{k}=\varphi_{\mathcal{B}}^{-1}\left(b_{j}/d_{j}\right)$
be unique the coordinate $d$-tuple representing $b_{j}/d_{j}$, we
then have that: 
\begin{equation}
\varphi_{\mathcal{B}}^{-1}\left(\tilde{H}\left(z\right)\right)=\varphi_{\mathcal{B}}^{-1}\left(\frac{a_{j}z+b_{j}}{d_{j}}\right)=\varphi_{\mathcal{B}}^{-1}\left(\frac{a_{j}z}{d_{j}}\right)+\varphi_{\mathcal{B}}^{-1}\left(\frac{b_{j}}{d_{j}}\right)=\mathbf{C}\mathbf{x}+\mathbf{k}
\end{equation}
where we used the fact that, as defined, $\varphi_{\mathcal{B}}$
is an isomorphism of the linear spaces $\mathbb{F}$ and $\mathbb{Q}^{d}$.
By (\ref{eq:A / D notation}), we can write $\mathbf{C}$ as: 
\begin{equation}
\mathbf{C}=\mathbf{D}^{-1}\mathbf{A}
\end{equation}
for $\mathbf{D},\mathbf{A}$ as described in (V) of (\ref{def:(F,I,B)-Hydra map}).
We can make $\mathbf{D}$ and $\mathbf{A}$ unique by choosing them
so as to make their non-zero entires (all of which are positive integers)
as small as possible. Consequently, the vector: 
\begin{equation}
\mathbf{b}\overset{\textrm{def}}{=}\mathbf{D}\mathbf{k}
\end{equation}
will have integer entries, and we can express the action of the $j$th
branch of $\tilde{H}$ on an arbitrary $\mathbf{x}$ by: 
\begin{equation}
\varphi_{\mathcal{B}}^{-1}\left(\frac{a_{j}z+b_{j}}{d_{j}}\right)=\frac{\mathbf{A}\mathbf{x}+\mathbf{b}}{\mathbf{D}}
\end{equation}

Since $\mathcal{O}_{\mathbb{F}}/\mathfrak{I}\cong\mathbb{Z}^{r}/P\mathbb{Z}^{r}$
for each $\mathbf{j}\in\mathbb{Z}^{r}/P\mathbb{Z}^{r}$, by the argument
given in the previous paragraph, there are unique invertible $\mathbf{A}_{\mathbf{j}},\mathbf{D}_{\mathbf{j}}$
satisfying (V-i), (V-ii), and (V-iii) respectively and a unique $d\times1$
column vector $\mathbf{b}_{\mathbf{j}}$ with integer entries so that:
\begin{equation}
\varphi_{\mathcal{B}}^{-1}\left(\tilde{H}\left(z\right)\right)=\frac{\mathbf{A}_{\mathbf{j}}\varphi_{\mathcal{B}}^{-1}\left(z\right)+\mathbf{b}_{\mathbf{j}}}{\mathbf{D}_{\mathbf{j}}}
\end{equation}
holds for all $z\in\mathcal{O}_{\mathbb{F}}$ that are congruent to
$\varphi_{\mathcal{B}}\left(\mathbf{j}\right)$ mod $\mathfrak{I}$.
Multiplying by the Iverson bracket $\left[\varphi_{\mathcal{B}}^{-1}\left(z\right)\overset{P}{\equiv}\mathbf{j}\right]$
and summing over all $\mathbf{j}\in\mathbb{Z}^{r}/P\mathbb{Z}^{r}$
gives: 
\begin{eqnarray*}
\sum_{\mathbf{j}\in\mathbb{Z}^{r}/P\mathbb{Z}^{r}}\left[\varphi_{\mathcal{B}}^{-1}\left(z\right)\overset{P}{\equiv}\mathbf{j}\right]\left(\frac{\mathbf{A}_{\mathbf{j}}\varphi_{\mathcal{B}}^{-1}\left(z\right)+\mathbf{b}_{\mathbf{j}}}{\mathbf{D}_{\mathbf{j}}}\right) & = & \sum_{\mathbf{j}\in\mathbb{Z}^{r}/P\mathbb{Z}^{r}}\left[\varphi_{\mathcal{B}}^{-1}\left(z\right)\overset{P}{\equiv}\mathbf{j}\right]\varphi_{\mathcal{B}}^{-1}\left(\tilde{H}\left(z\right)\right)\\
 & = & \sum_{\mathbf{j}\in\mathbb{Z}^{r}/P\mathbb{Z}^{r}}\underbrace{\left[z\overset{\mathfrak{I}}{\equiv}\varphi_{\mathcal{B}}\left(\mathbf{j}\right)\right]}_{\in\left\{ 0,1\right\} }\varphi_{\mathcal{B}}^{-1}\left(\tilde{H}\left(z\right)\right)\\
 & = & \varphi_{\mathcal{B}}^{-1}\left(\sum_{\mathbf{j}\in\mathcal{O}_{\mathbb{F}}/\mathfrak{I}}\left[z\overset{\mathfrak{I}}{\equiv}\varphi_{\mathcal{B}}\left(\mathbf{j}\right)\right]\tilde{H}\left(z\right)\right)\\
\left(\varphi_{\mathcal{B}}\left(\mathbf{j}\right)\textrm{s partition }\mathcal{O}_{\mathbb{F},\mathcal{B}}^{+}\right); & = & \varphi_{\mathcal{B}}^{-1}\left(\tilde{H}\left(z\right)\right)
\end{eqnarray*}
Hence: 
\begin{equation}
\varphi_{\mathcal{B}}^{-1}\left(\tilde{H}\left(z\right)\right)=\sum_{\mathbf{j}\in\mathbb{Z}^{r}/P\mathbb{Z}^{r}}\left[\varphi_{\mathcal{B}}^{-1}\left(z\right)\overset{P}{\equiv}\mathbf{j}\right]\frac{\mathbf{A}_{\mathbf{j}}\varphi_{\mathcal{B}}^{-1}\left(z\right)+\mathbf{b}_{\mathbf{j}}}{\mathbf{D}_{\mathbf{j}}},\forall z\in\mathcal{O}_{\mathbb{F},\mathcal{B}}^{+}
\end{equation}
Replacing $z$ with $\varphi_{\mathcal{B}}\left(\mathbf{x}\right)$
(where $\mathbf{x}\in\mathbb{Z}_{\mathbb{F},\mathcal{B}}^{+}$) then
gives the desired formula for $H\left(\mathbf{x}\right)$: 
\begin{equation}
\underbrace{\left(\varphi_{\mathcal{B}}^{-1}\circ\tilde{H}\circ\varphi_{\mathcal{B}}\right)\left(\mathbf{x}\right)}_{H\left(\mathbf{x}\right)}=\sum_{\mathbf{j}\in\mathbb{Z}^{r}/P\mathbb{Z}^{r}}\left[\mathbf{x}\overset{P}{\equiv}\mathbf{j}\right]\frac{\mathbf{A}_{\mathbf{j}}\mathbf{x}+\mathbf{b}_{\mathbf{j}}}{\mathbf{D}_{\mathbf{j}}},\forall\mathbf{x}\in\mathbb{Z}_{\mathbb{F},\mathcal{B}}^{+}
\end{equation}

Q.E.D.

\vphantom{}

Finally, we introduce notation to take on the roles played by $\mu_{j}$
and $p$ in the one-dimensional case. 
\begin{defn}
Let $H$ be a $P$-Hydra map on $\mathbb{Z}^{d}$.

\vphantom{}

I. For each $\mathbf{j}\in\mathbb{Z}^{r}/P\mathbb{Z}^{r}$, we write:
\begin{equation}
\mathbf{M}_{\mathbf{j}}\overset{\textrm{def}}{=}\mathbf{R}\frac{\mathbf{A}_{\mathbf{j}}}{\mathbf{D}_{\mathbf{j}}}\overset{\textrm{def}}{=}\mathbf{R}\mathbf{D}_{\mathbf{j}}^{-1}\mathbf{A}_{\mathbf{j}}\label{eq:Definition of bold M bold j}
\end{equation}

\vphantom{}

II. \nomenclature{$H_{\mathbf{j}}\left(\mathbf{x}\right)$}{$\mathbf{j}$th branch of a multi-dimensional Hydra map}For
each $\mathbf{j}\in\mathbb{Z}^{r}/P\mathbb{Z}^{r}$, we call the affine
linear map $H_{\mathbf{j}}:\mathbb{Q}^{d}\rightarrow\mathbb{Q}^{d}$
defined by: 
\begin{equation}
H_{\mathbf{j}}\left(\mathbf{x}\right)\overset{\textrm{def}}{=}\frac{\mathbf{A}_{\mathbf{j}}\mathbf{x}+\mathbf{b}_{\mathbf{j}}}{\mathbf{D}_{\mathbf{j}}}\label{eq:Definition of the bold jth branch of H}
\end{equation}
the \textbf{$\mathbf{j}$th branch }of $H$. Note that this overrides
the definition used for $H_{\mathbf{j}}$ in Chapter 2.

Using $\mathbf{M}_{\mathbf{j}}$, we can write: 
\begin{equation}
H_{\mathbf{j}}\left(\mathbf{x}\right)=\frac{\mathbf{M}_{\mathbf{j}}}{\mathbf{R}}\mathbf{x}+H_{\mathbf{j}}\left(\mathbf{0}\right)\label{eq:bold jth branch of H in terms of bold M bold j}
\end{equation}
We then define $H_{\mathbf{j}}^{\prime}\left(\mathbf{0}\right)$ as:
\begin{equation}
H_{\mathbf{j}}^{\prime}\left(\mathbf{0}\right)\overset{\textrm{def}}{=}\frac{\mathbf{A}_{\mathbf{j}}}{\mathbf{D}_{\mathbf{j}}}=\frac{\mathbf{M}_{\mathbf{j}}}{\mathbf{R}}\label{eq:Definition of H_bold j prime of bold 0}
\end{equation}
Also, note then that $H^{\prime}\left(\mathbf{0}\right)=H_{\mathbf{0}}^{\prime}\left(\mathbf{0}\right)$.

\vphantom{}

III. We say $H$ is integral (resp. non-integral) if its field analogue
is integral (resp. non-integral).
\end{defn}
\begin{example}
Consider the map $\tilde{H}:\mathbb{Z}\left[\sqrt{3}\right]\rightarrow\mathbb{Z}\left[\sqrt{3}\right]$
defined by: 
\begin{equation}
\tilde{H}\left(z\right)\overset{\textrm{def}}{=}\begin{cases}
\frac{z}{\sqrt{3}} & \textrm{if }z=0\mod\sqrt{3}\\
\frac{z-1}{\sqrt{3}} & \textrm{if }z=1\mod\sqrt{3}\\
\frac{4z+1}{\sqrt{3}} & \textrm{if }z=2\mod\sqrt{3}
\end{cases}
\end{equation}
Here $\gamma_{1}=1$ and $\gamma_{2}=\sqrt{3}$, $\mathfrak{I}=\left\langle \sqrt{3}\right\rangle $.
Since: 
\begin{equation}
\tilde{H}\left(a+b\sqrt{3}\right)=\begin{cases}
b+\frac{a}{3}\sqrt{3} & \textrm{if }a=0\mod3\\
b+\frac{a-1}{3}\sqrt{3} & \textrm{if }a=1\mod3\\
4b+\frac{4a+1}{3}\sqrt{3} & \textrm{if }a=2\mod3
\end{cases}
\end{equation}
we have: 
\begin{align*}
H\left(\left[\begin{array}{c}
v_{1}\\
v_{2}
\end{array}\right]\right) & =\begin{cases}
\left[\begin{array}{c}
v_{2}\\
\frac{v_{1}}{3}
\end{array}\right] & \textrm{if }v_{1}=0\mod3\\
\left[\begin{array}{c}
v_{2}\\
\frac{v_{1}-1}{3}
\end{array}\right] & \textrm{if }v_{1}=1\mod3\\
\left[\begin{array}{c}
4v_{2}\\
\frac{4v_{1}+1}{3}
\end{array}\right] & \textrm{if }v_{1}=2\mod3
\end{cases}\\
 & =\begin{cases}
\left[\begin{array}{cc}
0 & 1\\
\frac{1}{3} & 0
\end{array}\right]\left[\begin{array}{c}
v_{1}\\
v_{2}
\end{array}\right] & \textrm{if }v_{1}=0\mod3\\
\left[\begin{array}{cc}
0 & 1\\
\frac{1}{3} & 0
\end{array}\right]\left[\begin{array}{c}
v_{1}\\
v_{2}
\end{array}\right]+\left[\begin{array}{c}
0\\
-\frac{1}{3}
\end{array}\right] & \textrm{if }v_{1}=1\mod3\\
\left[\begin{array}{cc}
0 & 4\\
\frac{4}{3} & 0
\end{array}\right]\left[\begin{array}{c}
v_{1}\\
v_{2}
\end{array}\right]+\left[\begin{array}{c}
0\\
\frac{1}{3}
\end{array}\right] & \textrm{if }v_{1}=2\mod3
\end{cases}
\end{align*}
and so: 
\begin{align*}
H\left(\mathbf{v}\right) & =\left[\mathbf{v}\overset{3}{\equiv}\left[\begin{array}{c}
0\\
0
\end{array}\right]\right]\left[\begin{array}{cc}
0 & 1\\
\frac{1}{3} & 0
\end{array}\right]\mathbf{v}+\left[\mathbf{v}\overset{3}{\equiv}\left[\begin{array}{c}
1\\
0
\end{array}\right]\right]\left(\left[\begin{array}{cc}
0 & 1\\
\frac{1}{3} & 0
\end{array}\right]\mathbf{v}+\left[\begin{array}{c}
0\\
-\frac{1}{3}
\end{array}\right]\right)\\
 & +\left[\mathbf{v}\overset{3}{\equiv}\left[\begin{array}{c}
2\\
0
\end{array}\right]\right]\left(\left[\begin{array}{cc}
0 & 4\\
\frac{4}{3} & 0
\end{array}\right]\mathbf{v}+\left[\begin{array}{c}
0\\
\frac{1}{3}
\end{array}\right]\right)
\end{align*}
where, in all three branches, the permutation is: 
\begin{equation}
\left[\begin{array}{cc}
0 & 1\\
1 & 0
\end{array}\right]
\end{equation}
\end{example}
\vphantom{}

Lastly, we will need to distinguish those $P$-Hydra maps for which
$p_{n}=p$ for all $n\in\left\{ 1,\ldots,r\right\} $ for some prime
$p$.
\begin{defn}
\index{ideal!smooth}I say a non-zero proper ideal $\mathfrak{I}$
of $\mathcal{O}_{\mathbb{F}}$ is \textbf{smooth }if the $p_{n}$s
from the factorization:
\[
\mathcal{O}_{\mathbb{F}}/\mathfrak{I}\cong\mathfrak{C}_{p_{1}}\times\cdots\times\mathfrak{C}_{p_{r}}
\]
are all equal, with there being some integer $p\geq2$ so that $p_{n}=p$
for all $n\in\left\{ 1,\ldots,r\right\} $. I call a Multi-Dimensional
Hydra map \textbf{smooth}\index{Hydra map!smooth}\index{Hydra map!$p$-smooth}\textbf{
}(or, more specifically, \textbf{$p$-smooth})\textbf{ }whenever either
of the equivalent conditions is satisfied:

\vphantom{}

i $H$ is an $\left(\mathbb{F},\mathcal{B},\mathfrak{I}\right)$-Hydra
map and $\mathfrak{I}$ is smooth.

\vphantom{}

ii. $H$ is a $d$-dimensional $P$-Hydra map of depth $r$, and there
is an integer $p\geq2$ so that $p_{n}=p$ for all $n\in\left\{ 1,\ldots,r\right\} $.
In this case, I refer to $H$ is a \textbf{$d$-dimensional $p$-Hydra
map of depth $r$}.

\vphantom{}

Finally, I say $H$ (equivalently, $\tilde{H}$) is \textbf{prime
}whenever it is $p$-smooth for some prime $p$.
\end{defn}
\newpage{}

\section{\label{sec:5.2 The-Numen-of}The Numen of a Multi-Dimensional Hydra
Map}

THROUGHOUT THIS SECTION, UNLESS STATED OTHERWISE, $H$ DENOTES A PRIME
($p$-SMOOTH) $d$-DIMENSIONAL HYDRA MAP OF DEPTH $r$ WHICH FIXES
$\mathbf{0}$.

\vphantom{}

In this section, we will do for multi-dimensional Hydra maps what
we did for their one-dimensional cousins in Chapter 2. The exposition
here will be more briskly paced than in that earlier chapter. Throughout,
we will keep the one-dimensional case as a guiding analogy. Of course,
we will pause as needed to discuss deviations particular to the multi-dimensional
setting.

\subsection{\label{subsec:5.2.1 Notation-and-Preliminary}Notation and Preliminary
Definitions}

Much like as in the one-dimensional case, we will restrict our attention
to $\mathbf{m}\in\mathbb{N}_{0}^{d}$, the multi-dimensional analogue
of $\mathbb{N}_{0}$. Whereas our composition sequences in the one-dimensional
case were of the form: 
\begin{equation}
H_{j_{1}}\circ H_{j_{2}}\circ\cdots
\end{equation}
for integers $j_{1},j_{2},\ldots\in\mathbb{Z}/p\mathbb{Z}$, in the
multi-dimensional case, our composition sequences will be of the form:
\begin{equation}
H_{\mathbf{j}_{1}}\circ H_{\mathbf{j}_{2}}\circ\cdots
\end{equation}
where $\mathbf{j}_{1},\mathbf{j}_{2},\ldots\in\mathbb{Z}^{r}/p\mathbb{Z}^{r}$
are \emph{tuples}, rather than integers.\emph{ }What we denoted by
$H_{\mathbf{j}}$ in the one-dimensional case will now be written
as $H_{\mathbf{J}}$. As indicated by the capital letter, $\mathbf{J}$
denotes a \emph{matrix} whose columns are $\mathbf{j}_{1},\mathbf{j}_{2},\ldots$.
The particular definitions are as follows:
\begin{defn}[\textbf{Block strings}]
\ 

\vphantom{}

I. A \textbf{$p$-block string }of \textbf{depth }$r$ is a\index{block string}\index{block string!depth}
finite tuple $\mathbf{J}=\left(\mathbf{j}_{1},\ldots,\mathbf{j}_{\left|\mathbf{J}\right|}\right)^{T}$,
where:

\vphantom{}

i. \nomenclature{$\left|\mathbf{J}\right|$}{Length of the block string $\mathbf{J}$}$\left|\mathbf{J}\right|$
is an integer $\geq1$ which we call the \textbf{length }of $\mathbf{J}$.

\vphantom{}

ii. For each $m\in\left\{ 1,\ldots,\left|\mathbf{J}\right|\right\} $,
$\mathbf{j}_{m}$ is an $r$-tuple: 
\begin{equation}
\mathbf{j}_{m}=\left(j_{m,1},j_{m,2},\ldots,j_{m,r}\right)\in\mathbb{Z}^{r}/p\mathbb{Z}^{r}
\end{equation}
That is, for each $m$, the $n$th entry of $\mathbf{j}_{m}$ is an
element of $\mathbb{Z}/p\mathbb{Z}$.

We can also view $\mathbf{J}$ as a $\left|\mathbf{J}\right|\times r$
matrix, where the $m$th row of $\mathbf{J}$ consists of the entries
of $\mathbf{j}_{m}$: 
\begin{equation}
\mathbf{J}=\left(\mathbf{j}_{1},\ldots,\mathbf{j}_{\left|\mathbf{J}\right|}\right)^{T}=\left(\begin{array}{cccc}
j_{1,1} & j_{1,2} & \cdots & j_{1,r}\\
j_{2,1} & j_{2,2} & \cdots & j_{2,r}\\
\vdots & \vdots & \ddots & \vdots\\
j_{\left|\mathbf{J}\right|,1} & j_{\left|\mathbf{J}\right|,2} & \cdots & j_{\left|\mathbf{J}\right|,r}
\end{array}\right)\label{eq:J as a matrix}
\end{equation}

\vphantom{}

II. We write $\textrm{String}^{r}\left(p\right)$\nomenclature{$\textrm{String}^{r}\left(p\right)$}{ }
to denote the set of all finite length $p$-block strings of depth
$r$.

\vphantom{}

III. Given any $\mathbf{J}\in\textrm{String}^{r}\left(p\right)$,
the \textbf{composition}\index{composition sequence!multi-dimensional}\textbf{
sequence }$H_{\mathbf{J}}:\mathbb{Q}^{d}\rightarrow\mathbb{Q}^{d}$
\nomenclature{$H_{\mathbf{J}}\left(\mathbf{x}\right)$}{$\overset{\textrm{def}}{=}\left(H_{\mathbf{j}_{1}}\circ\cdots\circ H_{\mathbf{j}_{\left|\mathbf{J}\right|}}\right)\left(\mathbf{x}\right) $}is
the map: 
\begin{equation}
H_{\mathbf{J}}\left(\mathbf{x}\right)\overset{\textrm{def}}{=}\left(H_{\mathbf{j}_{1}}\circ\cdots\circ H_{\mathbf{j}_{\left|\mathbf{J}\right|}}\right)\left(\mathbf{x}\right),\textrm{ }\forall\mathbf{J}\in\textrm{String}^{r}\left(p\right),\textrm{ }\forall\mathbf{x}\in\mathbb{R}^{d}\label{eq:Def of composition sequence-1}
\end{equation}

\vphantom{}

IV. We write $\textrm{String}_{\infty}^{r}\left(p\right)$ \nomenclature{$\textrm{String}_{\infty}^{r}\left(p\right)$}{ }to
denote the set of all $p$-block strings of finite or infinite length.
We also include the empty set as an element of $\textrm{String}_{\infty}^{r}\left(p\right)$,
and call it the \textbf{empty $p$-block string}. 
\end{defn}
\begin{rem}
As $\left|\mathbf{J}\right|\rightarrow\infty$, we will have that
the sequence consisting of the $m$th entry of $\mathbf{j}_{1}$ followed
by the $m$th entry of $\mathbf{j}_{2}$, followed by $m$th entry
of $\mathbf{j}_{3}$, and so on will be the digits of a $p_{m}$-adic
integer. 
\end{rem}
\begin{defn}[\textbf{Concatenation}]
The\index{block string!concatenation} \textbf{concatenation operator}
$\wedge:\textrm{String}^{r}\left(p\right)\times\textrm{String}_{\infty}^{r}\left(p\right)\rightarrow\textrm{String}_{\infty}^{r}\left(p\right)$
is defined by: \index{concatenation!operation!multi-dimensional}
\begin{equation}
\mathbf{J}\wedge\mathbf{K}=\left(\mathbf{j}_{1},\ldots,\mathbf{j}_{\left|\mathbf{J}\right|}\right)\wedge\left(\mathbf{k}_{1},\ldots,\mathbf{k}_{\left|\mathbf{K}\right|}\right)\overset{\textrm{def}}{=}\left(\mathbf{j}_{1},\ldots,\mathbf{j}_{\left|\mathbf{J}\right|},\mathbf{k}_{1},\ldots,\mathbf{k}_{\left|\mathbf{K}\right|}\right)\label{eq:Definition of Concatenation-1}
\end{equation}
for all $\mathbf{J}\in\textrm{String}^{r}\left(p\right)$ and $\mathbf{K}\in\textrm{String}_{\infty}^{r}\left(p\right)$.
Also, for any integer $m\geq1$, if $\mathbf{J}\in\textrm{String}^{r}\left(p\right)$,
we write $\mathbf{J}^{\wedge m}$ to denote the concatenation of $m$
copies of $\mathbf{J}$: 
\begin{equation}
\mathbf{J}^{\wedge m}\overset{\textrm{def}}{=}\left(\underbrace{\mathbf{j}_{1},\ldots,\mathbf{j}_{\left|\mathbf{J}\right|},\mathbf{j}_{1},\ldots,\mathbf{j}_{\left|\mathbf{J}\right|},\ldots,\mathbf{j}_{1},\ldots,\mathbf{j}_{\left|\mathbf{J}\right|}}_{m\textrm{ times}}\right)^{T}\label{eq:Definition of block string concatenation}
\end{equation}
We shall only use this notation when $\mathbf{J}$ has finite length.
Note that the empty block string is the identity element with respect
to $\wedge$: 
\begin{equation}
\mathbf{J}\wedge\varnothing=\varnothing\wedge\mathbf{J}=\mathbf{J}
\end{equation}
\end{defn}
\begin{defn}
We write $\left\Vert \cdot\right\Vert _{p}$ to denote the non-archimedean
norm:\nomenclature{$\left\Vert \mathbf{z}\right\Vert _{p}$}{$\max\left\{ \left|\mathfrak{z}_{1}\right|_{p},\ldots,\left|\mathfrak{z}_{r}\right|_{p}\right\}$ }
\begin{equation}
\left\Vert \mathbf{z}\right\Vert _{p}\overset{\textrm{def}}{=}\max\left\{ \left|\mathfrak{z}_{1}\right|_{p},\ldots,\left|\mathfrak{z}_{r}\right|_{p}\right\} ,\textrm{ }\forall\mathbf{z}\in\mathbb{Z}_{p}^{r}\label{eq:Definition of p norm}
\end{equation}
which outputs the maximum of the $p$-adic absolute values of the
entries of $\mathbf{z}$. Here
\begin{equation}
\left\Vert \mathbf{z}-\mathbf{w}\right\Vert _{p}\leq\max\left\{ \left\Vert \mathbf{z}\right\Vert _{p},\left\Vert \mathbf{w}\right\Vert _{p}\right\} \label{eq:p norm ultrametric inequality}
\end{equation}
with equality whenever $\left\Vert \mathbf{z}\right\Vert _{p}\neq\left\Vert \mathbf{w}\right\Vert _{p}$,
and: 
\begin{equation}
\left\Vert \mathbf{z}\mathbf{w}\right\Vert _{p}\leq\left\Vert \mathbf{z}\right\Vert _{p}\left\Vert \mathbf{w}\right\Vert _{p}\label{eq:Product property of p norm}
\end{equation}
\end{defn}
\vphantom{}

Next, we have notations for projections mod $p^{k}$ and for the multi-dimensional
analogue of $\mathbb{Z}_{p}^{\prime}\backslash\mathbb{N}_{0}$.
\begin{defn}
For any $\mathbf{z}=\left(\mathfrak{z}_{1},\ldots,\mathfrak{z}_{r}\right)\in\mathbb{Z}_{p}^{r}$,
we write: \nomenclature{$\left[\mathbf{z}\right]_{p^{k}}$}{$\overset{\textrm{def}}{=}\left(\left[\mathfrak{z}_{1}\right]_{p^{k}},\ldots,\left[\mathfrak{z}_{r}\right]_{p^{k}}\right)$ \nopageref}
\begin{equation}
\left[\mathbf{z}\right]_{p^{k}}\overset{\textrm{def}}{=}\left(\left[\mathfrak{z}_{1}\right]_{p^{k}},\ldots,\left[\mathfrak{z}_{r}\right]_{p^{k}}\right)\label{eq:definition of z-bar mod P^k}
\end{equation}
We also use the notation:\nomenclature{$\left(\mathbb{Z}_{p}^{r}\right)^{\prime}$}{$\overset{\textrm{def}}{=}\mathbb{Z}_{p}^{r}\backslash\mathbb{N}_{0}^{r}$ \nopageref}
\begin{equation}
\left(\mathbb{Z}_{p}^{r}\right)^{\prime}\overset{\textrm{def}}{=}\mathbb{Z}_{p}^{r}\backslash\mathbb{N}_{0}^{r}\label{eq:Definition of Z_P prime}
\end{equation}
\end{defn}
\begin{defn}
Given $\mathbf{J}\in\textrm{String}_{\infty}^{r}\left(p\right)$ and
$\mathbf{z}=\left(\mathfrak{z}_{1},\ldots,\mathfrak{z}_{r}\right)\in\mathbb{Z}_{p}^{r}$
we say $\mathbf{J}$ \textbf{represents }(or \textbf{is} \textbf{associated
to})\textbf{ }$\mathbf{z}$ and vice-versa\textemdash written $\mathbf{J}\sim\mathbf{z}$
or $\mathbf{z}\sim\mathbf{J}$\textemdash whenever the entries of
the $k$th column of $\mathbf{J}$ represent the $p$-adic integer
$\mathfrak{z}_{k}$ for all $k\in\left\{ 1,\ldots,r\right\} $ in
the sense of Chapter 2: that is, for each $k$, the entries of the
$k$th column of $\mathbf{J}$ are the sequence of the $p$-adic digits
of $\mathfrak{z}_{k}$: 
\begin{equation}
\mathbf{J}\sim\mathbf{z}\Leftrightarrow\mathbf{z}\sim\mathbf{J}\Leftrightarrow\mathfrak{z}_{k}=j_{1,k}+j_{2,k}p+j_{3,k}p^{2}+\cdots\textrm{ }\forall k\in\left\{ 1,\ldots,r\right\} \label{eq:Definition of n-bold-j correspondence.-1}
\end{equation}
We extend the relation $\sim$ to a binary relation on $\textrm{String}_{\infty}^{r}\left(p\right)$
by writing $\mathbf{J}\sim\mathbf{K}$ if and only if both $\mathbf{J}$
and $\mathbf{K}$ represent the same $d$-tuple $\mathbf{z}\in\mathbb{Z}_{p}^{r}$,
and write \nomenclature{$\textrm{String}_{\infty}^{r}\left(p\right)/\sim$}{ }
$\textrm{String}_{\infty}^{r}\left(p\right)/\sim$ to denote the set
of all equivalence classes of $\textrm{String}_{\infty}^{r}\left(p\right)$
with respect to $\sim$. In particular, we identify any $p$-block
string $\mathbf{J}$ of depth $r$ whose entries are all $0$s (such
a string is called a \textbf{zero $p$-block string}) with the empty
set, and view this particular equivalence class of $\textrm{String}_{\infty}^{r}\left(p\right)/\sim$
as representing the zero vector in $\mathbb{Z}_{p}^{r}$.
\end{defn}
\vphantom{}

Similar to what we did in Chapter 2, every function  on $\mathbb{Z}_{p}^{r}$
can be defined as a function of block strings. 
\begin{defn}
For each $\mathbf{k}\in\mathbb{Z}^{r}/p\mathbb{Z}^{r}$, we define
$\#_{\mathfrak{I}:\mathbf{k}}:\textrm{String}^{r}\left(p\right)\rightarrow\mathbb{N}_{0}$
by: 
\begin{equation}
\#_{\mathfrak{I}:\mathbf{k}}\left(\mathbf{J}\right)\overset{\textrm{def}}{=}\left|\left\{ \mathbf{j}\in\mathbf{J}:\mathbf{j}=\mathbf{k}\right\} \right|,\textrm{ }\forall\mathbf{J}\in\textrm{String}^{r}\left(p\right)\label{eq:Definition of pound of big bold j}
\end{equation}
We also define $\#_{\mathfrak{I}:\mathbf{j}}$ as a function on $\mathbb{N}_{0}^{r}$
by writing: 
\begin{equation}
\#_{\mathfrak{I}:\mathbf{k}}\left(\mathbf{n}\right)\overset{\textrm{def}}{=}\#_{\mathfrak{I}:\mathbf{k}}\left(\mathbf{J}\right)\label{eq:Definition of pound bold n}
\end{equation}
where $\mathbf{J}$ is the shortest element of $\textrm{String}^{r}\left(p\right)$
representing $\mathbf{n}$. 
\end{defn}
\begin{defn}
We define the function \nomenclature{$\lambda_{p}\left(\mathbf{n}\right)$}{ }
$\lambda_{p}:\mathbb{N}_{0}^{r}\rightarrow\mathbb{N}_{0}$ by:
\begin{equation}
\lambda_{p}\left(\mathbf{n}\right)\overset{\textrm{def}}{=}\max_{1\leq\ell\leq r}\lambda_{p}\left(n_{\ell}\right)\label{eq:Definition of lambda P of bold n}
\end{equation}
Note that $\lambda_{p}\left(\mathbf{n}\right)$ is then the length
of the shortest $p$-block string representing $\mathbf{n}$.
\end{defn}
\begin{prop}
For any $\mathbf{z}\in\mathbb{Z}_{p}^{r}$ there exists a unique $\mathbf{J}\in\textrm{String}_{\infty}^{r}\left(p\right)$
of minimal length which represents $\mathbf{z}$. This $\mathbf{J}$
has finite length if and only if $\mathbf{z}=\mathbf{n}\in\mathbb{N}_{0}^{r}$,
in which case, $\left|\mathbf{J}\right|=\min\left\{ 1,\lambda_{p}\left(\mathbf{n}\right)\right\} $
(the minimum is necessary to account for the case when $\mathbf{n}=\mathbf{0}$). 
\end{prop}
Proof: We begin with the following special case: 
\begin{claim}
Let $\mathbf{n}\in\mathbb{N}_{0}^{r}$. Then, there is a unique $\mathbf{J}\in\textrm{String}_{\infty}^{r}\left(p\right)$
of minimal length which represents $\mathbf{n}$.

Proof of claim:\textbf{ }If $\mathbf{n}$ is the zero vector, then
we can choose $\mathbf{J}$ to have length $1$ and contain only $0$s.
So, suppose $\mathbf{n}$ is not the zero vector. Then $\lambda_{p}\left(\mathbf{n}\right)=\max_{1\leq k\leq r}\lambda_{p}\left(n_{k}\right)$
is positive. Using this, we construct the desired $\mathbf{J}$ like
so: for any $k$ for which $\lambda_{p}\left(n_{k}\right)=L$, we
write out in the $k$th column of $\mathbf{J}$ (starting from the
top, proceeding downward) all of the $p$-adic digits of $n_{k}$.
Next, for any $k$ for which $\lambda_{p}\left(n_{k}\right)<L$, we
have that $j_{\ell,k}=0$ for all $\ell\in\left\{ \lambda_{p}\left(n_{k}\right)+1,\ldots,L\right\} $;
that is, after we have written out all of the $p$-adic digits of
$n_{k}$ in the $k$th column of $\mathbf{J}$, if there are still
spaces further down in that column without entries, we fill those
spaces with $0$s.

As constructed, this $\mathbf{J}$ represents $\mathbf{n}$. Moreover,
its length is minimal: if $\mathbf{J}^{\prime}$ is a block string
of length $\left|\mathbf{J}\right|-1=L-1$, then, for any $k\in\left\{ 1,\ldots,r\right\} $
for which $\lambda_{p}\left(n_{k}\right)=L$, the $k$th column of
$\mathbf{J}^{\prime}$ will fail to represent $n_{k}$. This also
establishes $\mathbf{J}$'s uniqueness. This proves the claim.

\vphantom{}
\end{claim}
Having proved the claim, using the density of $\mathbb{N}_{0}^{r}$
in $\mathbb{Z}_{p}^{r}$, given any $\mathbf{z}\in\mathbb{Z}_{p}^{r}$,
let $\mathbf{J}$ be the limit (in the projective sense) of the sequence
$\mathbf{J}_{1},\mathbf{J}_{2},\ldots\in\textrm{String}_{\infty}^{r}\left(p\right)$,
where $\mathbf{J}_{m}$ is the unique minimal-length block string
representing $\left[\mathbf{z}\right]_{p^{m}}$, the $r$-tuple whose
entries are $\left[\mathfrak{z}_{k}\right]_{p^{m}}$. $\mathbf{J}$
will then be the unique block string representing $\mathbf{z}$.

Q.E.D.

\vphantom{}

We can restate the above as: 
\begin{prop}
The map which sends each $\mathbf{z}\in\mathbb{Z}_{p}^{r}$ to the
unique string $\mathbf{J}\in\textrm{String}_{\infty}^{r}\left(p\right)$
of minimal length which represents $\mathbf{z}$ is a bijection from
$\mathbb{Z}_{p}^{r}$ onto $\textrm{String}_{\infty}^{r}\left(p\right)/\sim$. 
\end{prop}
\begin{rem}
We identify the zero vector in $\mathbb{Z}_{p}^{r}$ with the empty
string. On the other hand, for any non-zero vector $\mathbf{z}=\left(\mathfrak{z}_{1},\ldots,\mathfrak{z}_{r}\right)$
so that $\mathfrak{z}_{k}=0$ for some $k$, we have that the $k$th
column of the unique $\mathbf{J}$ representing $\mathbf{z}$ will
contain only $0$s as entries.
\end{rem}
\begin{prop}
\label{prop:MD lambda and digit-number functional equations}$\lambda_{p}$
and the $\#_{\mathfrak{I}:\mathbf{k}}$s satisfy the \index{functional equation!lambda{p}@$\lambda_{p}$}functional
equations: 
\begin{equation}
\lambda_{p}\left(p\mathbf{n}+\mathbf{j}\right)=\lambda_{p}\left(\mathbf{n}\right)+1,\textrm{ }\forall\mathbf{n}\in\mathbb{N}_{0}^{r},\textrm{ }\forall\mathbf{j}\in\mathbb{Z}^{r}/p\mathbb{Z}^{r}\label{eq:lambda_P functional equation}
\end{equation}

\begin{equation}
\#_{\mathfrak{I}:\mathbf{k}}\left(p\mathbf{n}+\mathbf{j}\right)=\begin{cases}
\#_{\mathfrak{I}:\mathbf{k}}\left(\mathbf{n}\right)+1 & \textrm{if }\mathbf{j}=\mathbf{k}\\
\#_{\mathfrak{I}:\mathbf{k}}\left(\mathbf{n}\right) & \textrm{else}
\end{cases},\textrm{ }\forall\mathbf{n}\in\mathbb{N}_{0}^{r},\textrm{ }\forall\mathbf{j}\in\mathbb{Z}^{r}/p\mathbb{Z}^{r}\label{eq:Pound of bold n functional equation}
\end{equation}
More generally, for any $k\geq1$, any $\mathbf{m},\mathbf{n}\in\mathbb{N}_{0}^{r}$
with $\lambda_{p}\left(\mathbf{n}\right)\leq k$: 
\begin{equation}
\lambda_{p}\left(\mathbf{m}+p^{k}\mathbf{n}\right)=\lambda_{p}\left(\mathbf{m}\right)+\lambda_{p}\left(\mathbf{n}\right)\label{eq:lambda_P pseudoconcatenation identity}
\end{equation}
\begin{equation}
\#_{\mathfrak{I}:\mathbf{k}}\left(\mathbf{m}+p^{k}\mathbf{n}\right)=\#_{\mathfrak{I}:\mathbf{k}}\left(\mathbf{m}\right)+\#_{\mathfrak{I}:\mathbf{k}}\left(\mathbf{n}\right)\label{eq:Pound of bold n pseudoconcatenation identity}
\end{equation}
\end{prop}

\subsection{\label{subsec:5.2.2 Construction-of-the}Construction of the Numen}

For any affine linear map $L:\mathbb{R}^{d}\rightarrow\mathbb{R}^{d}$
defined by: 
\begin{equation}
L\left(\mathbf{x}\right)=\mathbf{A}\mathbf{x}+\mathbf{b}
\end{equation}
for a $d\times d$ matrix $\mathbf{A}$ and a $d\times1$ column vector
$\mathbf{b}$ (both with real entries), we will define the derivative
of $L$ to be the matrix $\mathbf{A}$. Thus, given any such $L$,
we will write $L^{\prime}\left(\mathbf{0}\right)$ to denote the derivative
of $L$. $L^{\prime}\left(\mathbf{0}\right)$ is then the unique $d\times d$
matrix so that: 
\begin{equation}
L\left(\mathbf{x}\right)=L^{\prime}\left(\mathbf{0}\right)\mathbf{x}+\mathbf{b}
\end{equation}
In fact, we can write: 
\begin{equation}
L\left(\mathbf{x}\right)=L^{\prime}\left(\mathbf{0}\right)\mathbf{x}+L\left(\mathbf{0}\right)
\end{equation}

Applying this formalism to a composition sequence $H_{\mathbf{J}}$\textemdash which
is valid, since, as a composition of a sequence of the affine linear
maps $H_{\mathbf{j}}:\mathbb{R}^{d}\rightarrow\mathbb{R}^{d}$\textemdash we
see that $H_{\mathbf{J}}$ is an affine linear map on $\mathbb{R}^{d}$,
with:

\begin{equation}
H_{\mathbf{J}}\left(\mathbf{x}\right)=H_{\mathbf{J}}^{\prime}\left(\mathbf{0}\right)\mathbf{x}+H_{\mathbf{J}}\left(\mathbf{0}\right),\textrm{ }\forall\mathbf{J}\in\textrm{String}^{r}\left(p\right)\label{eq:ax+b formula for H_bold_J}
\end{equation}
If $\mathbf{J}$ is a $p$-block string for which $H_{\mathbf{J}}\left(\mathbf{x}\right)=\mathbf{x}$,
this becomes: 
\begin{equation}
\mathbf{x}=H_{\mathbf{J}}^{\prime}\left(\mathbf{0}\right)\mathbf{x}+H_{\mathbf{J}}\left(\mathbf{0}\right)\label{eq:affine formula for bold x}
\end{equation}
So: 
\begin{equation}
\mathbf{x}=\left(\mathbf{I}_{d}-H_{\mathbf{J}}^{\prime}\left(\mathbf{0}\right)\right)^{-1}H_{\mathbf{J}}\left(\mathbf{0}\right)\label{eq:Formula for bold x in terms of bold J}
\end{equation}
where $\mathbf{I}_{d}$ is a $d\times d$ identity matrix.

Now for the definitions: 
\begin{defn}[\textbf{Multi-Dimensional $\chi_{H}$}]
\label{def:MD Chi_H} \index{chi{H}@$\chi_{H}$!multi-dimensional}\nomenclature{$\chi_{H}\left(\mathbf{J}\right)$}{ }\nomenclature{$\chi_{H}\left(\mathbf{n}\right)$}{ }\ 

\vphantom{}

I. We define $\chi_{H}:\textrm{String}^{r}\left(p\right)\rightarrow\mathbb{Q}^{d}$
by: 
\begin{equation}
\chi_{H}\left(\mathbf{J}\right)\overset{\textrm{def}}{=}H_{\mathbf{J}}\left(\mathbf{0}\right),\textrm{ }\forall\mathbf{J}\in\textrm{String}^{r}\left(p\right)\label{eq:MD Definition of Chi_H of bold J}
\end{equation}

\vphantom{}

II. For any $\mathbf{n}\in\mathbb{N}_{0}^{r}$, we write: 
\begin{equation}
\chi_{H}\left(\mathbf{n}\right)\overset{\textrm{def}}{=}\chi_{H}\left(\mathbf{J}\right)=\begin{cases}
H_{\mathbf{J}}\left(\mathbf{0}\right) & \textrm{if }\mathbf{n}\neq\mathbf{0}\\
\mathbf{0} & \textrm{else}
\end{cases}\label{eq:MD Definition of Chi_H of n}
\end{equation}
where $\mathbf{J}$ is any element of $\textrm{String}^{r}\left(p\right)$
which represents $\mathbf{n}$. 
\end{defn}
\begin{rem}
Much like in the one-dimensional case, since $H_{\mathbf{0}}\left(\mathbf{0}\right)=\mathbf{0}$,
it then follows that $\chi_{H}$ is then a well-defined function from
$\textrm{String}^{r}\left(p\right)/\sim$ (equivalently, $\mathbb{N}_{0}^{r}$)
to $\mathbb{Q}^{d}$.
\end{rem}
\begin{defn}[\textbf{Multi-Dimensional $M_{H}$}]
We define $M_{H}:\mathbb{N}_{0}^{r}\rightarrow\textrm{GL}_{d}\left(\mathbb{Q}\right)$
by: \nomenclature{$M_{H}\left(\mathbf{J}\right)$}{$\overset{\textrm{def}}{=}H_{\mathbf{J}}^{\prime}\left(\mathbf{0}\right)$  }\nomenclature{$M_{H}\left(\mathbf{n}\right)$}{ }\index{$M_{H}$!multi-dimensional}
\begin{equation}
M_{H}\left(\mathbf{n}\right)\overset{\textrm{def}}{=}M_{H}\left(\mathbf{J}\right)\overset{\textrm{def}}{=}H_{\mathbf{J}}^{\prime}\left(\mathbf{0}\right),\textrm{ }\forall\mathbf{n}\in\mathbb{N}_{0}^{r}\backslash\left\{ \mathbf{0}\right\} \label{eq:MD Definition of M_H}
\end{equation}
where $\mathbf{J}\in\textrm{String}^{r}\left(p\right)$ is the unique
minimal length block string representing $\mathbf{n}$. We define
$M_{H}\left(\mathbf{0}\right)$ and $M_{H}\left(\varnothing\right)$
to be the $d\times d$ identity matrix $\mathbf{I}_{d}$, where $\mathbf{0}$
is the $r\times1$ zero vector and where $\varnothing$ is the empty
block string. 
\end{defn}
\vphantom{}

In the one-dimensional case, the $q$-adic convergence of $\chi_{H}$
over $\mathbb{Z}_{p}$ was due to the $q$-adic decay of $M_{H}\left(\mathbf{j}\right)$
as $\left|\mathbf{j}\right|\rightarrow\infty$ with infinitely many
non-zero entries. We will do the same in the multi-dimensional case,
except with matrices. Unlike the one-dimensional case, however, the
multi-dimensional case has to deal with the various ways in which
the $\mathbf{A}_{\mathbf{j}}$s might interact with respect to multiplication.

Like in the one-dimensional case, $M_{H}$ and $\chi_{H}$ will be
characterized by functional equations and concatenation identities. 
\begin{prop}[\textbf{Functional Equations for Multi-Dimensional $M_{H}$}]
\label{prop:MD M_H functional equations}\index{functional equation!$M_{H}$!multi-dimensional}\ 
\begin{equation}
M_{H}\left(p\mathbf{n}+\mathbf{j}\right)=\mathbf{D}_{\mathbf{j}}^{-1}\mathbf{A}_{\mathbf{j}}M_{H}\left(\mathbf{n}\right),\textrm{ }\forall\mathbf{n}\in\mathbb{N}_{0}^{r},\textrm{ }\forall\mathbf{j}\in\mathbb{Z}^{r}/p\mathbb{Z}^{r}:p\mathbf{n}+\mathbf{j}\neq\mathbf{0}\label{eq:MD M_H functional equation}
\end{equation}
In terms of block strings and concatenations, we have: 
\begin{equation}
M_{H}\left(\mathbf{J}\wedge\mathbf{K}\right)=M_{H}\left(\mathbf{J}\right)M_{H}\left(\mathbf{K}\right),\textrm{ }\forall\mathbf{J},\mathbf{K}\in\textrm{String}^{r}\left(p\right)\label{eq:MD M_H concatenation identity}
\end{equation}
\end{prop}
\begin{rem}
Because $M_{H}$ is matrix-valued, $M_{H}\left(\mathbf{J}\right)M_{H}\left(\mathbf{K}\right)$
is not generally equal to $M_{H}\left(\mathbf{K}\right)M_{H}\left(\mathbf{J}\right)$.
However, if $H$ is $\mathbf{A}$-scalar, then: 
\begin{equation}
M_{H}\left(\mathbf{J}\right)M_{H}\left(\mathbf{K}\right)=M_{H}\left(\mathbf{K}\right)M_{H}\left(\mathbf{J}\right)
\end{equation}
does hold for all $\mathbf{J},\mathbf{K}$. Consequently, \emph{when
$H$ is $\mathbf{A}$-scalar, permuting the rows of $\mathbf{J}$
(a $\left|\mathbf{J}\right|\times r$ matrix) does not change the
value of $M_{H}\left(\mathbf{J}\right)$}.\index{concatenation!identities!multi-dimensional}
\end{rem}
\begin{lem}[\textbf{Functional Equations for Multi-Dimensional $\chi_{H}$}]
\label{lem:MD Chi_H functional equations and characterization}$\chi_{H}$
is the unique function $\mathbb{N}_{0}^{r}\rightarrow\mathbb{Q}^{d}$
satisfying the system of functional equations\index{functional equation!chi_{H}@$\chi_{H}$!multi-dimensional}:
\begin{equation}
\chi_{H}\left(p\mathbf{n}+\mathbf{j}\right)=\frac{\mathbf{A_{j}}\chi_{H}\left(\mathbf{n}\right)+\mathbf{b}_{\mathbf{j}}}{\mathbf{D}_{\mathbf{j}}},\textrm{ }\forall\mathbf{n}\in\mathbb{N}_{0}^{r},\textrm{ }\forall\mathbf{j}\in\mathbb{Z}^{r}/p\mathbb{Z}^{r}\label{eq:MD Chi_H functional equation}
\end{equation}
Equivalently: 
\begin{equation}
\chi_{H}\left(p\mathbf{n}+\mathbf{j}\right)=H_{\mathbf{j}}\left(\chi_{H}\left(\mathbf{n}\right)\right),\textrm{ }\forall\mathbf{n}\in\mathbb{N}_{0}^{r},\textrm{ }\forall\mathbf{j}\in\mathbb{Z}^{r}/p\mathbb{Z}^{r}\label{eq:MD Chi_H functional equation, alternate statement}
\end{equation}
In terms of concatenations of block strings, the above are written
as: 
\begin{equation}
\chi_{H}\left(\mathbf{J}\wedge\mathbf{K}\right)=H_{\mathbf{J}}\left(\chi_{H}\left(\mathbf{K}\right)\right),\textrm{ }\forall\mathbf{J},\mathbf{K}\in\textrm{String}^{r}\left(p\right)\label{eq:MD Chi_H concatenation identity}
\end{equation}
\end{lem}
Proof: For once, I'll leave this to the reader as an exercise.

Q.E.D.

\vphantom{}

Like in the one-dimensional case, we will also need explicit formulae
for $M_{H}$ and $\chi_{H}$. 
\begin{prop}
\label{prop:H_bold J in terms of M_H and Chi_H, explicitly}For all
$\mathbf{x}\in\mathbb{N}_{0}^{d}$ and all $\mathbf{J}\in\textrm{String}^{r}\left(p\right)$:
\begin{equation}
H_{\mathbf{J}}\left(\mathbf{x}\right)=\left(\prod_{k=1}^{\left|\mathbf{J}\right|}\frac{\mathbf{M}_{\mathbf{j}_{k}}}{\mathbf{R}}\right)\mathbf{x}+\sum_{n=1}^{\left|\mathbf{J}\right|}\left(\prod_{k=1}^{n-1}\frac{\mathbf{M}_{\mathbf{j}_{k}}}{\mathbf{R}}\right)H_{\mathbf{j}_{k}}\left(\mathbf{0}\right)\label{eq:MD H composition sequence formula}
\end{equation}
where the products are defined to be $\mathbf{I}_{d}$ whenever their
upper limits of multiplication are $0$. 
\end{prop}
Proof: Induction.

Q.E.D.

\vphantom{}

Next up, the explicit, string-based formulae for $M_{H}$ and $\chi_{H}$. 
\begin{prop}
\label{prop:Bold J formulae for M_H and Chi_H}\ 

\vphantom{}

I. For all non-zero $\mathbf{n}\in\mathbb{N}_{0}^{r}$ and any $\mathbf{J}\in\textrm{String}^{r}\left(p\right)$
representing $\mathbf{n}$: 
\begin{equation}
M_{H}\left(\mathbf{n}\right)=M_{H}\left(\mathbf{J}\right)=\prod_{k=1}^{\left|\mathbf{J}\right|}\frac{\mathbf{M}_{\mathbf{j}_{k}}}{\mathbf{R}}=\prod_{k=1}^{\left|\mathbf{J}\right|}\frac{\mathbf{A}_{\mathbf{j}_{k}}}{\mathbf{D}_{\mathbf{j}_{k}}}\label{eq:MD M_H formula}
\end{equation}

\vphantom{}

II. For all $\mathbf{J}\in\textrm{String}^{r}\left(p\right)$: 
\begin{equation}
\chi_{H}\left(\mathbf{J}\right)=\sum_{n=1}^{\left|\mathbf{J}\right|}\left(\prod_{k=1}^{n-1}\frac{\mathbf{M}_{\mathbf{j}_{k}}}{\mathbf{R}}\right)H_{\mathbf{j}_{n}}\left(\mathbf{0}\right)\label{eq:MD Chi_H formula}
\end{equation}
\end{prop}
\begin{rem}
All the above equalities hold in $\mathbb{Q}^{d,d}$, the space of
$d\times d$ matrices with rational entries.
\end{rem}
Proof: Use (\ref{eq:MD H composition sequence formula}) along with
the fact that: 
\begin{equation}
H_{\mathbf{J}}\left(\mathbf{x}\right)=H_{\mathbf{J}}^{\prime}\left(\mathbf{0}\right)\mathbf{x}+H_{\mathbf{J}}\left(\mathbf{0}\right)=M_{H}\left(\mathbf{J}\right)\mathbf{x}+\chi_{H}\left(\mathbf{J}\right)
\end{equation}

Q.E.D.

\vphantom{}

Now we deal with the rising-continuation of $\chi_{H}$. Like in the
one-dimensional case, $\chi_{H}$ will end up being a rising-continuous
function, albeit multi-dimensionally so. The specifics of multi-dimensional
rising-continuous functions are dealt with in Section \ref{sec:5.3 Rising-Continuity-in-Multiple}.
Our current goal is to state the multi-dimensional generalizations
of conditions like semi-basicness and basicness which were needed
to guaranteed the rising-continuability of $\chi_{H}$ in the one-dimensional
case.

To motivate these, let us examine (\ref{eq:MD Chi_H formula}). Recall
that the \emph{columns }of $\mathbf{J}$ are going to represent $p$-adic
integers as $\left|\mathbf{J}\right|\rightarrow\infty$; that is:
\begin{equation}
\mathbf{J}\rightarrow\left.\left[\begin{array}{cccc}
\mid & \mid &  & \mid\\
\mathfrak{z}_{1} & \mathfrak{z}_{2} & \cdots & \mathfrak{z}_{r}\\
\mid & \mid &  & \mid
\end{array}\right]\right\} \textrm{ infinitely many rows}
\end{equation}
where the $m$th column on the right contains the digits of the $p$-adic
integer $\mathfrak{z}_{m}$; the $k$th row on the right, recall,
is $\mathbf{j}_{k}$. In particular, the length of $\mathbf{J}$ is
infinite if and only if at least one of the $\mathfrak{z}_{m}$s is
an element of $\mathbb{Z}_{p}^{\prime}$. This, in turn, will occur
if and only if said $\mathfrak{z}_{m}$ has \emph{infinitely many
non-zero $p$-adic digits}, which occurs if and only if $\lambda_{p}\left(\mathbf{J}\right)\rightarrow\infty$.
Applying $\chi_{H}$ to $\mathbf{J}$, we see that:

\begin{equation}
\chi_{H}\left(\left[\begin{array}{cccc}
\mid & \mid &  & \mid\\
\mathfrak{z}_{1} & \mathfrak{z}_{2} & \cdots & \mathfrak{z}_{r}\\
\mid & \mid &  & \mid
\end{array}\right]\right)=\sum_{n=1}^{\left|\mathbf{J}\right|}\left(\prod_{k=1}^{n-1}\frac{\mathbf{M}_{\mathbf{j}_{k}}}{\mathbf{R}}\right)H_{\mathbf{j}_{k}}\left(\mathbf{0}\right)
\end{equation}
So, $\chi_{H}$ will extend to a function on $\mathbb{Z}_{p}^{r}$,
with the convergence of the series on the right occurring because
at least one of the $\mathfrak{z}_{m}$s has infinitely many non-zero
$p$-adic digits. To deal with the convergence of sums of products
of matrices, let briefly recall how convergence issues like this work
in spaces of matrices. 
\begin{prop}[\textbf{Convergence for Matrix Sums and Matrix Products}]
Let $K$ be a metrically complete non-archimedean valued field, let
$d\geq2$, and let $\left\{ \mathbf{M}_{n}\right\} _{n\geq1}$ be
a sequence in $K^{d,d}$ ($d\times d$ matrices with entries in $K$).
Then, the partial products: 
\begin{equation}
\prod_{n=1}^{N}\mathbf{M}_{n}
\end{equation}
converge in $K^{d,d}$ to the $d\times d$ zero matrix $\mathbf{O}_{d}$
whenever $\prod_{n=1}^{N}\left\Vert \mathbf{M}_{n}\right\Vert _{K}\rightarrow0$
in $\mathbb{R}$ as $N\rightarrow\infty$, where $\left\Vert \mathbf{M}_{n}\right\Vert _{K}$
is the maximum of the $K$-absolute-values of the entries of $\mathbf{M}_{n}$.\index{matrix!infinite product}\index{series!of matrices}\index{matrix!infinite series}

Likewise, since $K$ is non-archimedean, we have that a matrix series:
\begin{equation}
\sum_{n=1}^{\infty}\mathbf{M}_{n}
\end{equation}
converges entry-wise in $K$ to an element of $K^{d,d}$ if and only
if $\left\Vert \mathbf{M}_{n}\right\Vert _{K}\rightarrow0$. 
\end{prop}
Proof: $\left(K^{d,d},\left\Vert \cdot\right\Vert _{K}\right)$ is
a finite-dimensional Banach space over $K$.

Q.E.D.

\vphantom{}

To cleanly state the multi-dimensional analogue of what it means for
$H$ to be contracting, we will need to use the proposition given
below. Recall that we write $\left\Vert \cdot\right\Vert _{\infty}$
to denote the maximum of the ordinary, complex absolute value of the
entries of a matrix or vector. 
\begin{prop}
\label{prop:MD contracting proposition}Let $\mathbf{D}$ and $\mathbf{P}$
be invertible $d\times d$ matrices with rational entries, with $\mathbf{D}$
being diagonal and $\mathbf{P}$ being a permutation matrix, representing
a permutation of order $\omega$. Then, $\lim_{n\rightarrow\infty}\left\Vert \left(\mathbf{DP}\right)^{n}\right\Vert _{\infty}=0$
whenever $\left\Vert \left(\mathbf{D}\mathbf{P}\right)^{\omega}\right\Vert _{\infty}<1$. 
\end{prop}
Proof: As given, $\omega$ is the smallest positive integer so that
$\mathbf{P}^{\omega}=\mathbf{I}_{d}$. Now, letting $\mathbf{v}\in\mathbb{Q}^{d}$
be arbitrary observe that there will be an invertible diagonal matrix
$\mathbf{D}_{\omega}$ so that: 
\begin{equation}
\left(\mathbf{D}\mathbf{P}\right)^{\omega}\mathbf{v}=\mathbf{D}_{\omega}\mathbf{v}
\end{equation}
This is because each time we apply $\mathbf{D}\mathbf{P}$ to $\mathbf{v}$,
we permute the entries of $\mathbf{v}$ by $\mathbf{P}$ and then
multiply by $\mathbf{D}$. Since $\mathbf{P}$ has order $\omega$,
if we label the entries of $\mathbf{v}$ as $v_{1},\ldots,v_{d}$,
the positions of the $v_{j}$s in $\mathbf{v}$ will return to $v_{1},\ldots,v_{d}$
only after we have multiplied $\mathbf{v}$ by $\mathbf{P}$ $\omega$
times. In each application of $\mathbf{D}\mathbf{P}$, the $v_{j}$s
accrue a non-zero multiplicative factor\textemdash{} courtesy of $\mathbf{D}$.
So, after $\omega$ applications of $\mathbf{D}\mathbf{P}$, all of
the entries of $\mathbf{v}$ will be returned to the original places,
albeit with their original values having been multiplied by constants
chosen from the diagonal entries of $\mathbf{D}$. These constants
are then precisely the diagonal entries of $\mathbf{D}_{\omega}$.
Since $\mathbf{v}$ was arbitrary, we conclude that $\left(\mathbf{D}\mathbf{P}\right)^{\omega}=\mathbf{D}_{\omega}$.

Now, letting the integer $n\geq0$ be arbitrary, we can write $n$
as $\omega\left\lfloor \frac{n}{\omega}\right\rfloor +\left[n\right]_{\omega}$,
where $\left\lfloor \frac{n}{\omega}\right\rfloor $ is the largest
integer $\leq n/\omega$. As such: 
\begin{equation}
\left\Vert \left(\mathbf{D}\mathbf{P}\right)^{n}\right\Vert _{\infty}\leq\left\Vert \underbrace{\left(\mathbf{D}\mathbf{P}\right)^{\omega}}_{\mathbf{D}_{\omega}}\right\Vert _{\infty}^{\left\lfloor \frac{n}{\omega}\right\rfloor }\left\Vert \mathbf{D}\right\Vert _{\infty}^{\left[n\right]_{\omega}}\underbrace{\left\Vert \mathbf{P}\right\Vert _{\infty}^{\left[n\right]_{\omega}}}_{1}
\end{equation}
Since $\left[n\right]_{\omega}$ only takes values in $\left\{ 0,\ldots,\omega-1\right\} $,
the quantity $\left\Vert \mathbf{D}\right\Vert _{\infty}^{\left[n\right]_{\omega}}$
is uniformly bounded with respect to $n$, we can write: 
\begin{equation}
\left\Vert \left(\mathbf{D}\mathbf{P}\right)^{n}\right\Vert _{\infty}\ll\left\Vert \mathbf{D}_{\omega}\right\Vert _{\infty}^{\left\lfloor \frac{n}{\omega}\right\rfloor }
\end{equation}
The upper bounds tends to zero as $n\rightarrow\infty$ whenever $\left\Vert \mathbf{D}_{\omega}\right\Vert _{\infty}<1$.

Q.E.D.

\vphantom{}

Below, we formulate the analogues of properties such as monogenicity,
degeneracy, simplicity, and the like for multi-dimensional Hydra maps. 
\begin{defn}[\textbf{Qualitative Terminology for Multi-Dimensional Hydra Maps}]
\label{def:qualitative definitions for MD hydra maps}Let $H$ be
a $P$-Hydra map, where $P=\left(p_{1},\ldots,p_{r}\right)$ where
the $p_{\ell}$s are possibly distinct, and possibly non-prime.

\vphantom{}

I. The \textbf{ground }of $H$ is\index{Hydra map!ground} the set:
\begin{equation}
\textrm{Gr}\left(H\right)\overset{\textrm{def}}{=}\left\{ p\textrm{ prime}:p\mid p_{m}\textrm{ for some }m\in\left\{ 1,\ldots,r\right\} \right\} \label{eq:Definition of the Ground of H}
\end{equation}

\vphantom{}

II. I say $H$ is \textbf{rough }if $H$ is not smooth. In particular,
$H$ is called:

i. \textbf{Weakly rough }if $H$ is rough and $\left|\textrm{Gr}\left(H\right)\right|=1$;

\vphantom{}

ii. \textbf{Strongly rough }if $H$ is rough and $\left|\textrm{Gr}\left(H\right)\right|\geq2$.

\vphantom{}

III.\index{Hydra map!degenerate} We say $H$ is \textbf{degenerate
}if $\mathbf{A}_{\mathbf{j}}$ is a permutation matrix (including
possibly the identity matrix) for some $\mathbf{j}\in\left(\mathbb{Z}^{r}/P\mathbb{Z}^{r}\right)\backslash\left\{ \mathbf{0}\right\} $.
We say $H$ is \textbf{non-degenerate }whenever it is not degenerate.

\vphantom{}

IV. We say $H$ is \textbf{semi-simple }\index{Hydra map!semi-simple}whenever,
for any $\mathbf{j},\mathbf{k}\in\mathbb{Z}^{r}/P\mathbb{Z}^{r}$,
every non-zero entry in $\mathbf{A}_{\mathbf{j}}$ is co-prime to
every non-zero entry in $\mathbf{D}_{\mathbf{k}}$. If, in addition:
\begin{equation}
\mathbf{D}_{\mathbf{j}}=\left(\begin{array}{cccccc}
p_{1}\\
 & \ddots\\
 &  & p_{r}\\
 &  &  & 1\\
 &  &  &  & \ddots\\
 &  &  &  &  & 1
\end{array}\right),\textrm{ }\forall\mathbf{j}\in\mathbb{Z}^{r}/P\mathbb{Z}^{r}\label{eq:MD Definition of simplicity}
\end{equation}
we then say $H$ is \textbf{simple}.

\vphantom{}

V. We say $H$ is \textbf{monogenic }\index{Hydra map!monogenic}if
there exists a prime $q\notin\textrm{Gr}\left(H\right)$ so that $\left\Vert \mathbf{A}_{\mathbf{j}}\right\Vert _{q}\leq1/q$
for all $\mathbf{j}\in\left(\mathbb{Z}^{r}/P\mathbb{Z}^{r}\right)\backslash\left\{ \mathbf{0}\right\} $.
We write $q_{H}$ to denote the smallest such prime. We say $H$ \index{Hydra map!polygenic}
is \textbf{polygenic }if $H$ is not monogenic, but there exists a
set $Q$ primes, with $Q\cap\textrm{Gr}\left(H\right)=\varnothing$
so that for each $\mathbf{j}\in\left(\mathbb{Z}^{r}/P\mathbb{Z}^{r}\right)\backslash\left\{ \mathbf{0}\right\} $,
there exists a $q\in Q$ for which $\left\Vert \mathbf{A}_{\mathbf{j}}\right\Vert _{q}\leq1/q$.

\vphantom{}

VI.\footnote{This definition will not be relevant until Chapter 6.}
For each $\mathbf{j}\in\mathbb{Z}^{r}/P\mathbb{Z}^{r}$, we write
$\omega_{\mathbf{j}}$ to denote the order of the permutation represented
by the permutation matrix $\mathbf{P}_{\mathbf{j}}$, where $\mathbf{A}_{\mathbf{j}}=\tilde{\mathbf{A}}_{\mathbf{j}}\mathbf{P}_{\mathbf{j}}$,
where $\tilde{\mathbf{A}}_{\mathbf{j}}$ is diagonal. We say $H$
\index{Hydra map!contracting} is \textbf{contracting }whenever $\lim_{n\rightarrow\infty}\left\Vert \left(H^{\prime}\left(\mathbf{0}\right)\right)^{n}\right\Vert _{\infty}=0$;
recall that $H^{\prime}\left(\mathbf{0}\right)=\mathbf{D}_{\mathbf{0}}^{-1}\mathbf{A}_{\mathbf{0}}$.
By \textbf{Proposition \ref{prop:MD contracting proposition}}, it
follows that $H$ is contracting whenever $\left\Vert \left(H^{\prime}\left(\mathbf{0}\right)\right)^{\omega_{\mathbf{0}}}\right\Vert _{\infty}<1$.

\vphantom{}

VII. We say $H$ \index{Hydra map!basic}is \textbf{basic }if $H$
is simple, non-degenerate, and monogenic.

\vphantom{}

VIII. We say $H$ is \textbf{semi-basic }if \index{Hydra map!semi-basic}$H$
is semi-simple, non-degenerate, and monogenic.
\end{defn}
\begin{rem}
For brevity, we will sometimes write $q$ to denote $q_{H}$. More
generally, in any context where a Hydra map is being worked with\textemdash unless
explicitly stated otherwise\textemdash $q$ will \emph{always }denote
$q_{H}$.
\end{rem}
\begin{rem}
Like in the one-dimensional case, everything given in this subsection
holds when $q_{H}$ is replaced by any prime $q\neq p$ such that
$\left\Vert \mathbf{A}_{\mathbf{j}}\right\Vert _{q}\leq1/q$ for all
$\mathbf{j}\in\left(\mathbb{Z}^{r}/p\mathbb{Z}^{r}\right)\backslash\left\{ \mathbf{0}\right\} $.
\end{rem}
\begin{rem}
Like with polygenicity in the one-dimensional case, I would like to
believe that frames can be used to provide a setting where we can
speak meaningfully about the $p$-adic interpolation of $\chi_{H}$
in the case of a $p$-smooth polygenic multi-dimensional Hydra map.
\end{rem}
\begin{prop}
\label{prop:MD Contracting H proposition}If $H$ is contracting,
the matrix $\mathbf{I}_{d}-H^{\prime}\left(\mathbf{0}\right)$ is
invertible.
\end{prop}
Proof: If $H$ is contracting, there is an integer $\omega_{\mathbf{0}}\geq1$
so that the entries of the matrix $\left(H^{\prime}\left(\mathbf{0}\right)\right)^{\omega_{\mathbf{0}}}$
are all rational numbers with archimedean absolute value $<1$. Consequently,
letting $\mathbf{X}$ denote $H^{\prime}\left(\mathbf{0}\right)$:
\begin{align*}
\sum_{n=0}^{\infty}\mathbf{X}^{n} & =\sum_{k=0}^{\omega_{\mathbf{0}}-1}\sum_{n=0}^{\infty}\mathbf{X}^{\omega_{\mathbf{0}}n+k}\\
 & =\sum_{k=0}^{\omega_{\mathbf{0}}-1}\mathbf{X}^{k}\sum_{n=0}^{\infty}\mathbf{X}^{\omega_{\mathbf{0}}n}\\
\left(\left\Vert \mathbf{X}^{\omega_{\mathbf{0}}}\right\Vert _{\infty}<1\right); & \overset{\mathbb{R}^{d,d}}{=}\sum_{k=0}^{\omega_{\mathbf{0}}-1}\mathbf{X}^{k}\left(\mathbf{I}_{d}-\mathbf{X}^{\omega_{\mathbf{0}}}\right)^{-1}
\end{align*}
As such, the series $\sum_{n=0}^{\infty}\mathbf{X}^{n}$ converges
in $\mathbb{R}^{d,d}$, and therefore defines a matrix $\mathbf{M}\in\mathbb{R}^{d,d}$.
Since:
\[
\left(\mathbf{I}_{d}-\mathbf{X}\right)\mathbf{M}=\mathbf{M}\left(\mathbf{I}_{d}-\mathbf{X}\right)=\mathbf{I}_{d}
\]
this proves that $\mathbf{M}=\left(\mathbf{I}_{d}-\mathbf{X}\right)^{-1}$,
showing that $\mathbf{I}_{d}-H^{\prime}\left(\mathbf{0}\right)$ is
invertible.

Q.E.D.

\vphantom{}

Like in the one-dimensional case, monogenicity will be needed to ensure
the existence of a rising-continuation of $\chi_{H}$. 
\begin{prop}
\label{prop:MD M_H decay estimate}Let $H$ be monogenic. Then:

\textup{ 
\begin{equation}
\left\Vert M_{H}\left(\mathbf{J}\right)\right\Vert _{q_{H}}=\left\Vert \prod_{k=1}^{\left|\mathbf{J}\right|}\frac{\mathbf{M}_{\mathbf{j}_{k}}}{\mathbf{R}}\right\Vert _{q_{H}}\leq q_{H}^{-\#\left(\mathbf{J}\right)},\textrm{ }\forall\mathbf{J}\in\textrm{String}^{r}\left(p\right)\label{eq:Decay estimate on M_H of big bold j}
\end{equation}
}This result also holds when we replace $\mathbf{J}$ by the unique
$\mathbf{n}\in\mathbb{N}_{0}^{r}$ for which $\mathbf{J}\sim\mathbf{n}$. 
\end{prop}
Proof: Since $\mathbf{M}_{\mathbf{j}}=\mathbf{R}\mathbf{D}_{\mathbf{j}}^{-1}\mathbf{A}_{\mathbf{j}}$,
and since the entries of $\mathbf{R}$ and $\mathbf{D}_{\mathbf{j}}$
are co-prime to $q_{H}$, we have that: 
\begin{equation}
\left\Vert \frac{\mathbf{M}_{\mathbf{j}}}{\mathbf{R}}\right\Vert _{q_{H}}\leq\left\Vert \mathbf{A}_{\mathbf{j}}\right\Vert _{q_{H}}
\end{equation}
Since $H$ is monogenic, $\left\Vert \mathbf{A}_{\mathbf{j}}\right\Vert _{q_{H}}\leq1/q_{H}$
whenever $\mathbf{j}\neq\mathbf{0}$. Since $M_{H}\left(\mathbf{J}\right)=\prod_{k=1}^{\left|\mathbf{J}\right|}\frac{\mathbf{M}_{\mathbf{j}_{k}}}{\mathbf{R}}$,
(\ref{eq:Decay estimate on M_H of big bold j}) then holds.

Q.E.D.

\vphantom{}

Now, the all-important characterization of $\chi_{H}$'s interpolation
as the unique rising-continuous solution of a certain system of functional
equations. 
\begin{lem}[\textbf{$\left(p,q\right)$-adic Characterization of $\chi_{H}$}]
\label{lem:MD rising-continuation and functional equations of Chi_H}Let
$H$ be monogenic. Then, the limit:\nomenclature{$\chi_{H}\left(\mathbf{z}\right)$}{ }
\begin{equation}
\chi_{H}\left(\mathbf{z}\right)\overset{\mathbb{Z}_{q_{H}}^{d}}{=}\lim_{n\rightarrow\infty}\chi_{H}\left(\left[\mathbf{z}\right]_{p^{n}}\right)\label{eq:MD Rising Continuity Formula for Chi_H}
\end{equation}
exists for all $\mathbf{z}\in\mathbb{Z}_{p}^{r}$, and thereby defines
a continuation of $\chi_{H}$ to a function $\mathbb{Z}_{p}^{r}\rightarrow\mathbb{Z}_{q_{H}}^{d}$.
Moreover:\index{functional equation!chi_{H}@$\chi_{H}$!multi-dimensional}

\vphantom{}

I. 
\begin{equation}
\chi_{H}\left(p\mathbf{z}+\mathbf{j}\right)\overset{\mathbb{Z}_{q_{H}}^{d}}{=}\frac{\mathbf{A_{j}}\chi_{H}\left(\mathbf{z}\right)+\mathbf{b}_{\mathbf{j}}}{\mathbf{D}_{\mathbf{j}}},\textrm{ }\forall\mathbf{z}\in\mathbb{Z}_{p}^{r},\textrm{ }\forall\mathbf{j}\in\mathbb{Z}^{r}/p\mathbb{Z}^{r}\label{eq:MD Functional Equations for Chi_H over the rho-bar-adics}
\end{equation}

\vphantom{}

II. \emph{(\ref{eq:MD Rising Continuity Formula for Chi_H})} is the
unique function $\mathbf{f}:\mathbb{Z}_{p}^{r}\rightarrow\mathbb{Z}_{q_{H}}^{d}$
satisfying both\emph{: 
\begin{equation}
\mathbf{f}\left(p\mathbf{z}+\mathbf{j}\right)\overset{\mathbb{Z}_{q_{H}}^{d}}{=}\frac{\mathbf{A_{j}}\mathbf{f}\left(\mathbf{z}\right)+\mathbf{b}_{\mathbf{j}}}{\mathbf{D}_{\mathbf{j}}},\textrm{ }\forall\mathbf{z}\in\mathbb{Z}_{p}^{r},\textrm{ }\forall\mathbf{j}\in\mathbb{Z}^{r}/p\mathbb{Z}^{r}\label{eq:MD extension Lemma for Chi_H - functional equation}
\end{equation}
and: 
\begin{equation}
\mathbf{f}\left(\mathbf{z}\right)\overset{\mathbb{Z}_{q_{H}}^{d}}{=}\lim_{n\rightarrow\infty}\mathbf{f}\left(\left[\mathbf{z}\right]_{p^{n}}\right)\label{eq:MD extension Lemma for Chi_H - rising continuity}
\end{equation}
} 
\end{lem}
Proof: Let $\mathbf{z}\in\mathbb{Z}_{p}^{r}$ and $N\geq0$ be arbitrary.
Then, choose: 
\begin{equation}
\mathbf{J}=\left[\begin{array}{ccc}
\mathfrak{z}_{1,0} & \cdots & \mathfrak{z}_{r,0}\\
\mathfrak{z}_{1,1} & \cdots & \mathfrak{z}_{r,1}\\
\vdots &  & \vdots\\
\mathfrak{z}_{1,N-1} & \cdots & \mathfrak{z}_{r,N-1}
\end{array}\right]=\left[\begin{array}{ccc}
\mid &  & \mid\\
\left[\mathfrak{z}_{1}\right]_{p^{N}} & \cdots & \left[\mathfrak{z}_{r}\right]_{p^{N}}\\
\mid &  & \mid
\end{array}\right]\sim\left[\mathbf{z}\right]_{^{N}}
\end{equation}
where the $m$th column's entries are the first $N$ $p$-adic digits
of $\mathfrak{z}_{m}$. Then, taking (\ref{eq:MD Chi_H formula})
from \textbf{Proposition \ref{prop:Bold J formulae for M_H and Chi_H}}:
\begin{equation}
\chi_{H}\left(\mathbf{J}\right)=\sum_{n=1}^{\left|\mathbf{J}\right|}\left(\prod_{k=1}^{n-1}\frac{\mathbf{M}_{\mathbf{j}_{k}}}{\mathbf{R}}\right)H_{\mathbf{j}_{k}}\left(\mathbf{0}\right)
\end{equation}
we obtain: 
\begin{equation}
\chi_{H}\left(\left[\mathbf{z}\right]_{p^{N}}\right)=\sum_{n=1}^{N}\left(\prod_{k=1}^{n-1}\frac{\mathbf{M}_{\left(\mathfrak{z}_{1,k},\ldots,\mathfrak{z}_{r,k}\right)}}{\mathbf{R}}\right)H_{\left(\mathfrak{z}_{1,n},\ldots,\mathfrak{z}_{r,n}\right)}\left(\mathbf{0}\right)\label{eq:MD Chi_H formula for truncation of z-arrow}
\end{equation}

If $\mathbf{z}\in\mathbb{N}_{0}^{r}$, then $\left(\mathfrak{z}_{1,n},\ldots,\mathfrak{z}_{r,n}\right)$
must be a tuple of $0$s for all sufficiently large $n$. Since $H$
fixes $\mathbf{0}$, we know that $\mathbf{b}_{\mathbf{0}}=\mathbf{0}$.
As such, the $n$th terms of (\ref{eq:MD Chi_H formula for truncation of z-arrow})
are then $\mathbf{0}$ for all sufficiently large $n$, which shows
that (\ref{eq:MD Rising Continuity Formula for Chi_H}) holds for
$\mathbf{z}\in\mathbb{N}_{0}^{r}$.

On the other hand, suppose at least one entry of $\mathbf{z}$ is
not in $\mathbb{N}_{0}$. Then, as $n\rightarrow\infty$, in the matrix
product: 
\begin{equation}
\prod_{k=1}^{n-1}\frac{\mathbf{M}_{\mathbf{j}_{k}}}{\mathbf{R}}
\end{equation}
there will be infinitely many values of $k$ for which the $r$-tuple
$\mathbf{j}=\left(\mathfrak{z}_{1,k},\ldots,\mathfrak{z}_{r,k}\right)$
has a non-zero entry. Since $H$ is monogenic, the decay estimate
given by \textbf{Proposition \ref{prop:MD M_H decay estimate}} guarantees
that: 
\begin{equation}
\left\Vert \prod_{k=1}^{n-1}\frac{\mathbf{M}_{\mathbf{j}_{k}}}{\mathbf{R}}\right\Vert _{q_{H}}\leq q_{H}^{-\left|\left\{ k\in\left\{ 1,\ldots,n-1\right\} :\mathbf{j}_{k}\neq0\right\} \right|}
\end{equation}
which tends to $0$ as $n\rightarrow\infty$. This shows that (\ref{eq:MD Chi_H formula for truncation of z-arrow})
converges point-wise in $\mathbb{Z}_{q_{H}}^{d}$ to a limit.

Since (\ref{eq:MD Functional Equations for Chi_H over the rho-bar-adics})
holds for all $\mathbf{z}\in\mathbb{N}_{0}^{r}$, the limit (\ref{eq:MD Functional Equations for Chi_H over the rho-bar-adics})
then shows that (\ref{eq:MD Functional Equations for Chi_H over the rho-bar-adics})
in fact holds for all $\mathbf{z}\in\mathbb{Z}_{p}^{r}$, thereby
proving (I).

Finally, \textbf{Lemma \ref{lem:MD Chi_H functional equations and characterization}
}shows that any solution $\mathbf{f}$ of (\ref{eq:MD extension Lemma for Chi_H - functional equation})
must be equal to $\chi_{H}$ on $\mathbb{N}_{0}^{r}$. If $\mathbf{f}$
also satisfies (\ref{eq:MD extension Lemma for Chi_H - rising continuity}),
this then shows that, just like $\chi_{H}$'s extension to $\mathbb{Z}_{p}^{r}$,
the extension of $\mathbf{f}$ to $\mathbb{Z}_{p}^{r}$ is uniquely
determined by $\mathbf{f}$'s values on $\mathbb{N}_{0}^{r}$. Since
$\mathbf{f}=\chi_{H}$ on $\mathbb{N}_{0}^{r}$, this then forces
$\mathbf{f}=\chi_{H}$ on $\mathbb{Z}_{p}^{r}$, proving the uniqueness
of $\chi_{H}$.

Q.E.D.

\subsection{\label{subsec:5.2.3 The-Correspondence-Principle}The Correspondence
Principle Revisited}

THROUGHOUT THIS SUBSECTION, WE ALSO ASSUME $H$ IS INTEGRAL, UNLESS
STATED OTHERWISE.

\vphantom{}

Like in the one-dimensional case, we start with a self-concatenation
identity.
\begin{prop}
\label{prop:MD self-concatenation identity}Let $\mathbf{n}\in\mathbb{N}_{0}^{r}$
and let $\mathbf{J}\in\textrm{String}^{r}\left(p\right)$ be the shortest
block string representing. Then, for all $k\in\mathbb{N}_{1}$: 
\begin{equation}
\left(\mathbf{j}_{1}^{\wedge k},\ldots,\mathbf{j}_{\left|\mathbf{J}\right|}^{\wedge k}\right)^{T}=\mathbf{J}^{\wedge k}\sim\frac{1-p^{k\lambda_{p}\left(\mathbf{n}\right)}}{1-p^{\lambda_{p}\left(\mathbf{n}\right)}}\mathbf{n}\overset{\textrm{def}}{=}\left(n_{1}\frac{1-p^{k\lambda_{p}\left(n_{1}\right)}}{1-p^{\lambda_{p_{}}\left(n_{1}\right)}},\ldots,n_{r}\frac{1-p^{k\lambda_{p}\left(n_{r}\right)}}{1-p^{\lambda_{p}\left(n_{r}\right)}}\right)\label{eq:MD Self-Concatentation Identity}
\end{equation}
\end{prop}
Proof: Apply \textbf{Proposition \ref{prop:Concatenation exponentiation}}
to each row of $\mathbf{J}^{\wedge k}$.

Q.E.D.

\vphantom{}Next come the multi-dimensional analogues of $B_{p}$
and the functional equation satisfied by $\chi_{H}\circ B_{p}$. 
\begin{defn}[$\mathbf{B}_{p}$]
We define $\mathbf{B}_{p}:\mathbb{Z}_{p}^{r}\rightarrow\mathbb{Z}_{p}^{r}$
by:\nomenclature{$\mathbf{B}_{p}\left(\mathbf{z}\right)$}{$\overset{\textrm{def}}{=}\left(B_{p}\left(\mathfrak{z}_{1}\right),\ldots,B_{p}\left(\mathfrak{z}_{m}\right)\right)$}
\begin{equation}
\mathbf{B}_{p}\left(\mathbf{z}\right)\overset{\textrm{def}}{=}\left(B_{p}\left(\mathfrak{z}_{1}\right),\ldots,B_{p}\left(\mathfrak{z}_{m}\right)\right)\label{eq:MD Definition of B projection function}
\end{equation}
where, recall: 
\begin{equation}
B_{p}\left(\mathfrak{z}\right)=\begin{cases}
\mathfrak{z} & \textrm{if }\mathfrak{z}\in\mathbb{Z}_{p}^{\prime}\\
\frac{\mathfrak{z}}{1-p^{\lambda_{p}\left(\mathfrak{z}\right)}} & \textrm{if }\mathfrak{z}\in\mathbb{N}_{0}
\end{cases}
\end{equation}
\end{defn}
\begin{lem}[\textbf{Functional Equation for $\chi_{H}\circ\mathbf{B}_{p}$}]
\label{lem:MD Chi_H o B functional equation}Let\index{functional equation!chi{H}circmathbf{B}{p}@$\chi_{H}\circ\mathbf{B}_{p}$}
$H$ be monogenic. Then: 
\begin{equation}
\chi_{H}\left(\mathbf{B}_{p}\left(\mathbf{n}\right)\right)\overset{\mathbb{Z}_{q_{H}}^{d}}{=}\left(\mathbf{I}_{d}-M_{H}\left(\mathbf{n}\right)\right)^{-1}\chi_{H}\left(\mathbf{n}\right),\textrm{ }\forall\mathbf{n}\in\mathbb{N}_{0}^{r}\backslash\left\{ \mathbf{0}\right\} \label{eq:MD Chi_H B functional equation}
\end{equation}
Moreover, both sides of the above identity are $d\times1$ column
vectors with entries in $\mathbb{Q}$. 
\end{lem}
Proof: Let $\mathbf{n}\in\mathbb{N}_{0}^{r}\backslash\left\{ \mathbf{0}\right\} $,
and let $\mathbf{J}\in\textrm{String}^{r}\left(p\right)$ be the shortest
block string representing $\mathbf{n}$. Then: 
\begin{equation}
H_{\mathbf{J}^{\wedge k}}\left(\mathbf{m}\right)=H_{\mathbf{J}}^{\circ k}\left(\mathbf{m}\right),\textrm{ }\forall\mathbf{m}\in\mathbb{Z}^{d},\textrm{ }\forall k\geq1
\end{equation}
Since: 
\begin{equation}
H_{\mathbf{J}}\left(\mathbf{m}\right)=M_{H}\left(\mathbf{J}\right)\mathbf{m}+\chi_{H}\left(\mathbf{J}\right)
\end{equation}
observe that: 
\begin{align*}
H_{\mathbf{J}^{\wedge k}}\left(\mathbf{m}\right) & \overset{\mathbb{Q}^{d}}{=}H_{\mathbf{J}}^{\circ k}\left(\mathbf{m}\right)\\
 & =\left(M_{H}\left(\mathbf{J}\right)\right)^{k}\mathbf{m}+\left(\sum_{\ell=0}^{k-1}\left(M_{H}\left(\mathbf{J}\right)\right)^{\ell}\right)\chi_{H}\left(\mathbf{J}\right)\\
 & =\left(M_{H}\left(\mathbf{n}\right)\right)^{k}\mathbf{m}+\left(\sum_{\ell=0}^{k-1}\left(M_{H}\left(\mathbf{n}\right)\right)^{\ell}\right)\chi_{H}\left(\mathbf{n}\right)
\end{align*}
Since $\mathbf{n}\neq\mathbf{0}$, there is an $m\in\left\{ 1,\ldots,r\right\} $
so that the $m$th entry of $\mathbf{n}$ has at least one non-zero
$p$-adic digit. Because $H$ is monogenic, we then have that at least
\emph{one} of the matrices in the product $M_{H}\left(\mathbf{n}\right)$
has $q_{H}$-adic norm $<1$. Because all the matrices in that product
have $q$-adic norm $\leq1$, this shows that $\left\Vert M_{H}\left(\mathbf{n}\right)\right\Vert _{q_{H}}<1$.
Consequently, the series: 
\begin{equation}
\sum_{\ell=0}^{\infty}\left(M_{H}\left(\mathbf{n}\right)\right)^{\ell}
\end{equation}
converges in $\mathbb{Q}_{q_{H}}^{d,d}$, thereby defining the inverse
of $\mathbf{I}_{d}-M_{H}\left(\mathbf{n}\right)$. In fact, because
$\mathbf{I}_{d}-M_{H}\left(\mathbf{n}\right)$ has entries in $\mathbb{Q}$,
note that that $\left(\mathbf{I}_{d}-M_{H}\left(\mathbf{n}\right)\right)^{-1}\in\mathbb{Q}_{q_{H}}^{d,d}$.
Moreover, we can use the geometric series formula to write: 
\begin{equation}
\left(\mathbf{I}_{d}-M_{H}\left(\mathbf{n}\right)\right)\sum_{\ell=0}^{k-1}\left(M_{H}\left(\mathbf{n}\right)\right)^{\ell}=\sum_{\ell=0}^{k-1}\left(M_{H}\left(\mathbf{n}\right)\right)^{\ell}-\sum_{\ell=0}^{k-1}\left(M_{H}\left(\mathbf{n}\right)\right)^{\ell+1}=\mathbf{I}_{d}-\left(M_{H}\left(\mathbf{n}\right)\right)^{k}
\end{equation}
So: 
\begin{equation}
\sum_{\ell=0}^{k-1}\left(M_{H}\left(\mathbf{n}\right)\right)^{\ell}=\frac{\mathbf{I}_{d}-\left(M_{H}\left(\mathbf{n}\right)\right)^{k}}{\mathbf{I}_{d}-M_{H}\left(\mathbf{n}\right)}\overset{\textrm{def}}{=}\left(\mathbf{I}_{d}-\left(M_{H}\left(\mathbf{n}\right)\right)^{k}\right)^{-1}\left(\mathbf{I}_{d}-\left(M_{H}\left(\mathbf{n}\right)\right)^{k}\right)
\end{equation}

Consequently: 
\begin{equation}
H_{\mathbf{J}^{\wedge k}}\left(\mathbf{m}\right)\overset{\mathbb{Q}^{d}}{=}\left(M_{H}\left(\mathbf{n}\right)\right)^{k}\mathbf{m}+\frac{\mathbf{I}_{d}-\left(M_{H}\left(\mathbf{n}\right)\right)^{k}}{\mathbf{I}_{d}-M_{H}\left(\mathbf{n}\right)}\chi_{H}\left(\mathbf{n}\right)
\end{equation}
and so: 
\begin{equation}
M_{H}\left(\mathbf{J}^{\wedge k}\right)\mathbf{m}+\chi_{H}\left(\mathbf{J}^{\wedge k}\right)\overset{\mathbb{Q}^{d}}{=}\left(M_{H}\left(\mathbf{n}\right)\right)^{k}\mathbf{m}+\frac{\mathbf{I}_{d}-\left(M_{H}\left(\mathbf{n}\right)\right)^{k}}{\mathbf{I}_{d}-M_{H}\left(\mathbf{n}\right)}\chi_{H}\left(\mathbf{n}\right)
\end{equation}
Now, using $M_{H}$'s concatenation identity (\textbf{Proposition
\ref{prop:MD M_H functional equations}}): 
\begin{equation}
M_{H}\left(\mathbf{J}^{\wedge k}\right)\mathbf{m}=\left(M_{H}\left(\mathbf{J}\right)\right)^{k}\mathbf{m}=\left(M_{H}\left(\mathbf{n}\right)\right)^{k}\mathbf{m}
\end{equation}
we have: 
\begin{equation}
\left(M_{H}\left(\mathbf{n}\right)\right)^{k}\mathbf{m}+\chi_{H}\left(\mathbf{J}^{\wedge k}\right)\overset{\mathbb{Q}^{d}}{=}\left(M_{H}\left(\mathbf{n}\right)\right)^{k}\mathbf{m}+\frac{\mathbf{I}_{d}-\left(M_{H}\left(\mathbf{n}\right)\right)^{k}}{\mathbf{I}_{d}-M_{H}\left(\mathbf{n}\right)}\chi_{H}\left(\mathbf{n}\right)
\end{equation}
and so: 
\begin{equation}
\chi_{H}\left(\mathbf{J}^{\wedge k}\right)=\frac{\mathbf{I}_{d}-\left(M_{H}\left(\mathbf{n}\right)\right)^{k}}{\mathbf{I}_{d}-M_{H}\left(\mathbf{n}\right)}\chi_{H}\left(\mathbf{n}\right)\label{eq:MD Chi_H B functional equation, ready to take limits}
\end{equation}
Using \textbf{Proposition \ref{prop:MD self-concatenation identity}}
gives us: 
\begin{align*}
\chi_{H}\left(\mathbf{J}^{\wedge k}\right) & =\chi_{H}\left(\mathbf{n}\frac{1-p^{k\lambda_{p}\left(\mathbf{n}\right)}}{1-p^{\lambda_{p}\left(\mathbf{n}\right)}}\right)\\
 & =\chi_{H}\left(n_{1}\frac{1-p^{k\lambda_{p}\left(n_{1}\right)}}{1-p^{\lambda_{p}\left(n_{1}\right)}},\ldots,n_{r}\frac{1-p^{k\lambda_{p}\left(n_{r}\right)}}{1-p^{\lambda_{p}\left(n_{r}\right)}}\right)
\end{align*}
$\chi_{H}$'s rising-continuity property (\ref{eq:MD Rising Continuity Formula for Chi_H})
allows us to take the limit as $k\rightarrow\infty$: 
\begin{equation}
\lim_{k\rightarrow\infty}\chi_{H}\left(\mathbf{J}^{\wedge k}\right)\overset{\mathbb{Z}_{q_{H}}^{d}}{=}\chi_{H}\left(\frac{n_{1}}{1-p^{\lambda_{p}\left(n_{1}\right)}},\ldots,\frac{n_{r}}{1-p^{\lambda_{p}\left(n_{r}\right)}}\right)=\chi_{H}\left(\mathbf{B}_{p}\left(\mathbf{n}\right)\right)
\end{equation}
On the other hand, (\ref{eq:MD Chi_H B functional equation, ready to take limits})
yields; 
\begin{align*}
\lim_{k\rightarrow\infty}\chi_{H}\left(\mathbf{J}^{\wedge k}\right) & \overset{\mathbb{Z}_{q_{H}}^{d}}{=}\lim_{k\rightarrow\infty}\frac{\mathbf{I}_{d}-\left(M_{H}\left(\mathbf{n}\right)\right)^{k}}{\mathbf{I}_{d}-M_{H}\left(\mathbf{n}\right)}\chi_{H}\left(\mathbf{n}\right)\\
\left(\left\Vert M_{H}\left(\mathbf{n}\right)\right\Vert _{q}<1\right); & \overset{\mathbb{Z}_{q_{H}}^{d}}{=}\frac{\mathbf{I}_{d}-\mathbf{O}_{d}}{\mathbf{I}_{d}-M_{H}\left(\mathbf{n}\right)}\chi_{H}\left(\mathbf{n}\right)\\
 & =\frac{\chi_{H}\left(\mathbf{n}\right)}{\mathbf{I}_{d}-M_{H}\left(\mathbf{n}\right)}\\
 & =\left(\mathbf{I}_{d}-M_{H}\left(\mathbf{n}\right)\right)^{-1}\chi_{H}\left(\mathbf{n}\right)
\end{align*}
Putting the two together gives: 
\begin{equation}
\chi_{H}\left(\mathbf{B}_{p}\left(\mathbf{n}\right)\right)\overset{\mathbb{Z}_{q_{H}}^{d}}{=}\left(\mathbf{I}_{d}-M_{H}\left(\mathbf{n}\right)\right)^{-1}\chi_{H}\left(\mathbf{n}\right)
\end{equation}

Lastly, because the right-hand side is the product of $\left(\mathbf{I}_{d}-M_{H}\left(\mathbf{n}\right)\right)^{-1}$
(a $d\times d$ matrix with entries in $\mathbb{Q}$) and $\chi_{H}\left(\mathbf{n}\right)$
(a $d\times1$ column vector with entries in $\mathbb{Q}$), we see
that both sides of (\ref{eq:MD Chi_H B functional equation}) are
$d\times1$ column vectors with entries in $\mathbb{Q}$.

Q.E.D.

\vphantom{}Next up: the multi-dimensional notion of a ``wrong value''.
But first, let us establish the $p$-adic extendibility of $H$.
\begin{prop}
\label{prop:p-adic extension of H-1}Let $H:\mathbb{Z}^{d}\rightarrow\mathbb{Z}^{d}$
be any $p$-Hydra map of depth $r$, not necessarily integral. Then,
$H$ admits an extension to a continuous map $\mathbb{Z}_{p}^{d}\rightarrow\mathbb{Z}_{p}^{d}$
defined by:
\begin{equation}
H\left(\mathbf{z}\right)\overset{\textrm{def}}{=}\sum_{\mathbf{j}\in\mathbb{Z}^{r}/p\mathbb{Z}^{r}}\left[\mathbf{z}\overset{p}{\equiv}\mathbf{j}\right]\mathbf{D}_{\mathbf{j}}^{-1}\left(\mathbf{A}_{\mathbf{j}}\mathbf{z}+\mathbf{b}_{\mathbf{j}}\right),\textrm{ }\forall\mathbf{z}\in\mathbb{Z}_{p}^{r}\label{eq:MD p-adic extension of H}
\end{equation}
Moreover, the function $\mathbf{f}:\mathbb{Z}_{p}^{d}\rightarrow\mathbb{Z}_{p}^{d}$
defined by:
\begin{equation}
\mathbf{f}\left(\mathbf{z}\right)=\sum_{\mathbf{j}\in\mathbb{Z}^{r}/p\mathbb{Z}^{r}}\left[\mathbf{z}\overset{p}{\equiv}\mathbf{j}\right]\mathbf{D}_{\mathbf{j}}^{-1}\left(\mathbf{A}_{\mathbf{j}}\mathbf{z}+\mathbf{b}_{\mathbf{j}}\right),\textrm{ }\forall\mathbf{z}\in\mathbb{Z}_{p}^{r}
\end{equation}
is the unique continuous function $\mathbb{Z}_{p}^{d}\rightarrow\mathbb{Z}_{p}^{d}$
whose restriction to $\mathbb{Z}^{d}$ is equal to $H$.
\end{prop}
Proof: Immediate, just like the one-dimensional case.

Q.E.D.
\begin{defn}[\textbf{Wrong Values and Propriety}]
\label{def:MD Wrong values and propriety}\index{Hydra map!wrong value}\index{Hydra map!proper}Let
$H$ be any $p$-Hydra map, not necessarily integral.

\vphantom{}

I. We say $\mathbf{y}=\left(\mathfrak{y}_{1},\ldots,\mathfrak{y}_{d}\right)\in\mathbb{Q}_{p}^{d}$
is a \textbf{wrong value for $H$ }whenever there are $\mathbf{z}\in\mathbb{Z}_{p}^{d}$
and $\mathbf{J}\in\textrm{String}^{r}\left(p\right)$ so that $\mathbf{y}=H_{\mathbf{J}}\left(\mathbf{z}\right)$
and $H_{\mathbf{J}}\left(\mathbf{z}\right)\neq H^{\circ\left|\mathbf{J}\right|}\left(\mathbf{z}\right)$.
We call $\mathbf{z}$ a \textbf{seed }of $\mathbf{y}$.

\vphantom{}

II.\index{Hydra map!proper} We say $H$ is \textbf{proper }whenever
every wrong value of $H$ is an element of $\mathbb{Q}_{p}^{d}\backslash\mathbb{Z}_{p}^{d}$;
that is, every wrong value of $H$ is a $d$-tuple of $p$-adic rational
numbers at least one of which is \emph{not} a $p$-adic integer.
\end{defn}
\vphantom{}

We now repeat the lead-up to the Correspondence Principle from the
one-dimensional case.
\begin{prop}
\label{prop:MD Q_p / Z_p prop}Let $H$ be any $p$-Hydra map of dimension
$d$ and depth $r$. Then, $H_{\mathbf{j}}\left(\mathbb{Q}_{p}^{d}\backslash\mathbb{Z}_{p}^{d}\right)\subseteq\mathbb{Q}_{p}^{d}\backslash\mathbb{Z}_{p}^{d}$
for all $\mathbf{j}\in\mathbb{Z}^{r}/p\mathbb{Z}^{r}$.
\end{prop}
\begin{rem}
Here, we are viewing the $H_{\mathbf{j}}$s as functions $\mathbb{Q}_{p}^{d}\rightarrow\mathbb{Q}_{p}^{d}$.
\end{rem}
Proof: Let $\mathbf{y}\in\mathbb{Q}_{p}^{d}\backslash\mathbb{Z}_{p}^{d}$
and $\mathbf{j}\in\mathbb{Z}^{r}/p\mathbb{Z}^{r}$ be arbitrary. Note
that as a $p$-Hydra map, $p$ divides every element of $\mathbf{D}_{\mathbf{j}}$
which is not $1$. Moreover, as a multi-dimensional $p$-Hydra map,
every non-zero element of $\mathbf{A}_{\mathbf{j}}$ is co-prime to
every non-zero element of $\mathbf{D}_{\mathbf{j}}$. Consequently,
$p$ does not divide any non-zero element of $\mathbf{A}_{\mathbf{j}}$,
and so, $\left\Vert \mathbf{A}_{\mathbf{j}}\right\Vert _{p}=1$. Thus:
\begin{equation}
\left\Vert \mathbf{A}_{\mathbf{j}}\mathbf{y}\right\Vert _{p}=\left\Vert \mathbf{y}\right\Vert _{p}>1
\end{equation}
So, $\mathbf{A}_{\mathbf{j}}\mathbf{y}\in\mathbb{Q}_{p}^{d}\backslash\mathbb{Z}_{p}^{d}$.
Since every entry of the $d$-tuple $\mathbf{b}_{\mathbf{j}}$ is
a rational integer, the $p$-adic ultrametric inequality guarantees
that $\left\Vert \mathbf{A}_{\mathbf{j}}\mathbf{y}+\mathbf{b}_{\mathbf{j}}\right\Vert >1$,
and thus, $\mathbf{A}_{\mathbf{j}}\mathbf{y}+\mathbf{b}_{\mathbf{j}}$
is in $\mathbb{Q}_{p}^{d}\backslash\mathbb{Z}_{p}^{d}$. Finally,
$p$ dividing the diagonal elements of $\mathbf{D}_{\mathbf{j}}$
which are \emph{not }$1$ implies $\left\Vert D_{\mathbf{j}}\right\Vert _{p}<1$,
and so:
\begin{equation}
\left\Vert H_{\mathbf{j}}\left(\mathbf{y}\right)\right\Vert _{p}=\left\Vert \frac{\mathbf{A}_{\mathbf{j}}\mathbf{y}+\mathbf{b}_{\mathbf{j}}}{\mathbf{D}_{\mathbf{j}}}\right\Vert _{p}\geq\left\Vert \mathbf{A}_{\mathbf{j}}\mathbf{y}+\mathbf{b}_{\mathbf{j}}\right\Vert _{p}>1
\end{equation}
which shows that $H_{\mathbf{j}}\left(\mathbf{y}\right)\in\mathbb{Q}_{p}^{d}\backslash\mathbb{Z}_{p}^{d}$.

Q.E.D.
\begin{lem}
\label{lem:MD integrality lemma}Let $H$ be a semi-basic $p$-Hydra
map of dimension $d$ and depth $r$, where $p$ is prime; we do not
require $H$ to be integral. Then, $H$ is proper if and only if $H$
is integral.
\end{lem}
Proof: Let $H$ be semi-basic.

I. (Proper implies integral) Suppose $H$ is proper, and let $\mathbf{n}\in\mathbb{Z}^{d}$
be arbitrary. Then, clearly, when $\mathbf{j}\in\mathbb{Z}^{r}/p\mathbb{Z}^{r}$
satisfies $\mathbf{j}=\left[\mathbf{n}\right]_{p}$ (recall, this
means the first $r$ entries of $\mathbf{n}$ are congruent mod $p$
to the first $r$ entries of $\mathbf{j}$), we have that $H_{\mathbf{j}}\left(\mathbf{n}\right)\in\mathbb{Z}^{d}$.
So, suppose instead that $\mathbf{j}\neq\left[\mathbf{n}\right]_{p}$.
Since $H$ is proper, the fact that $\mathbf{n}\in\mathbb{Z}^{d}\subseteq\mathbb{Z}_{p}^{d}$
tells us that $\left\Vert H_{\mathbf{j}}\left(\mathbf{n}\right)\right\Vert _{p}>1$.
So, at least one entry of $H_{\mathbf{j}}\left(\mathbf{n}\right)$
is not a $p$-adic integer, and thus, is not a rational integer, either.
This proves $H$ is integral.

\vphantom{}

II. (Integral implies proper) Suppose $H$ is integral, and\textemdash by
way of contradiction\textemdash suppose $H$ is \emph{not}\textbf{
}proper. Then, there is a $\mathbf{z}\in\mathbb{Z}_{p}^{d}$ and a
$\mathbf{j}\in\mathbb{Z}^{r}/p\mathbb{Z}^{r}$ with $\mathbf{j}\neq\left[\mathbf{z}\right]_{p}$
(again, only the first $r$ entries are affected) so that $\left\Vert H_{\mathbf{j}}\left(\mathbf{z}\right)\right\Vert _{p}\leq1$.

Now, writing 
\[
\mathbf{z}=\sum_{m=1}^{d}\sum_{n=0}^{\infty}c_{m,n}p^{n}\mathbf{e}_{m}
\]
where $\mathbf{e}_{m}$ is the $m$th standard basis vector in $\mathbb{Z}_{p}^{d}$:
\begin{align*}
H_{\mathbf{j}}\left(\mathbf{z}\right) & =\frac{\mathbf{A}_{\mathbf{j}}\left(\sum_{m=1}^{d}\sum_{n=0}^{\infty}c_{m,n}p^{n}\mathbf{e}_{m}\right)+\mathbf{b}_{\mathbf{j}}}{\mathbf{D}_{\mathbf{j}}}\\
 & =\frac{\mathbf{A}_{\mathbf{j}}\left(\sum_{m=1}^{d}c_{m,0}\mathbf{e}_{m}\right)+\mathbf{b}_{\mathbf{j}}}{\mathbf{D}_{\mathbf{j}}}+\sum_{m=1}^{d}\sum_{n=1}^{\infty}c_{m,n}p^{n-1}\left(p\mathbf{D}_{\mathbf{j}}^{-1}\right)\mathbf{e}_{m}\\
\left(\sum_{m=1}^{d}c_{m,0}\mathbf{e}_{m}=\left[\mathbf{z}\right]_{p}\right); & =H_{\mathbf{j}}\left(\left[\mathbf{z}\right]_{p}\right)+\sum_{m=1}^{d}\sum_{n=1}^{\infty}c_{m,n}p^{n-1}\left(p\mathbf{D}_{\mathbf{j}}^{-1}\right)\mathbf{e}_{m}
\end{align*}
Because $H$ is a $p$-Hydra map, every diagonal entry of $\mathbf{D}_{\mathbf{j}}$
which is not $1$ must be a divisor of $p$$d_{j}$ must be a divisor
of $p$. The primality of $p$ then forces $d_{j}$ to be either $1$
or $p$. In either case, we have that $p\mathfrak{y}/d_{j}$ is an
element of $\mathbb{Z}_{p}$. Our contradictory assumption $\left|H_{j}\left(\mathfrak{z}\right)\right|_{p}\leq1$
tells us that $H_{j}\left(\mathfrak{z}\right)$ is also in $\mathbb{Z}_{p}$,
and so:
\begin{equation}
H_{j}\left(\left[\mathfrak{z}\right]_{p}\right)=H_{j}\left(\mathfrak{z}\right)-\frac{p\mathfrak{y}}{d_{j}}\in\mathbb{Z}_{p}
\end{equation}
Since $H$ was given to be integral, $j\neq\left[\mathfrak{z}\right]_{p}$
implies $H_{j}\left(\left[\mathfrak{z}\right]_{p}\right)\notin\mathbb{Z}$.
As such, $H_{j}\left(\left[\mathfrak{z}\right]_{p}\right)$ is a $p$-adic
integer which is \emph{not }a rational integer. Since $H_{j}\left(\left[\mathfrak{z}\right]_{p}\right)$
is a non-integer rational number which is a $p$-adic integer, the
denominator of $H_{j}\left(\left[\mathfrak{z}\right]_{p}\right)=\frac{a_{j}\left[\mathfrak{z}\right]_{p}+b_{j}}{d_{j}}$
must be divisible by some prime $q\neq p$, and hence, $q\mid d_{j}$.
However, we saw that $p$ being prime forced $d_{j}\in\left\{ 1,p\right\} $\textemdash this
is impossible!

Thus, it must be that $H$ is proper.

Q.E.D.

\vphantom{}

Next, we use \textbf{Proposition \ref{prop:MD Q_p / Z_p prop} }to
establish a result about the wrong values of $H$ when $H$ is proper.
\begin{lem}
\label{lem:MD wrong values lemma}Let $H$ be any proper $p$-Hydra
map, not necessarily integral or semi-basic. Then, all wrong values
of $H$ are elements of $\mathbb{Q}_{p}^{d}\backslash\mathbb{Z}_{p}^{d}$.
\end{lem}
Proof: Let $H$ be as given, let $\mathbf{z}\in\mathbb{Z}_{p}^{d}$,
and let $\mathbf{i}\in\mathbb{Z}^{r}/p\mathbb{Z}^{r}$ be such that
$\left[\mathbf{z}\right]_{p}\neq\mathbf{i}$. Then, by definition
of properness, the quantity: 
\begin{equation}
H_{\mathbf{i}}\left(\mathbf{z}\right)=\frac{\mathbf{A}_{\mathbf{i}}\mathbf{z}+\mathbf{b}_{\mathbf{i}}}{\mathbf{D}_{\mathbf{i}}}
\end{equation}
has $p$-adic norm $>1$. By \textbf{Proposition \ref{prop:MD Q_p / Z_p prop}},
this then forces $H_{\mathbf{J}}\left(H_{\mathbf{i}}\left(\mathbf{z}\right)\right)$
to be an element of $\mathbb{Q}_{p}^{d}\backslash\mathbb{Z}_{p}^{d}$
for all $\mathbf{J}\in\textrm{String}^{r}\left(p\right)$. Since every
wrong value with seed $\mathbf{z}$ is of the form $H_{\mathbf{J}}\left(H_{\mathbf{i}}\left(\mathbf{z}\right)\right)$
for some $\mathbf{J}\in\textrm{String}^{r}\left(p\right)$, some $\mathbf{z}\in\mathbb{Z}_{p}^{d}$,
and some $\mathbf{i}\in\mathbb{Z}^{r}/p\mathbb{Z}^{r}$ for which
$\left[\mathbf{z}\right]_{p}\neq\mathbf{i}$, this shows that every
wrong value of $H$ is in $\mathbb{Q}_{p}^{d}\backslash\mathbb{Z}_{p}^{d}$.
So, $H$ is proper.

Q.E.D.
\begin{lem}
\label{lem:MD properness lemma}Let $H$ be any proper $p$-Hydra
map, not necessarily semi-basic or integral. Let $\mathbf{z}\in\mathbb{Z}_{p}^{d}$,
and let $\mathbf{J}\in\mathrm{String}^{r}\left(p\right)$. If $H_{\mathbf{J}}\left(\mathbf{z}\right)=\mathbf{z}$,
then $\ensuremath{H^{\circ\left|\mathbf{J}\right|}\left(\mathbf{z}\right)=\mathbf{z}}$.
\end{lem}
Proof: Let $H$, $\mathbf{z}$, and $\mathbf{J}$ be as given. By
way of contradiction, suppose $H^{\circ\left|\mathbf{J}\right|}\left(\mathbf{z}\right)\neq\mathbf{z}$.
But then $\mathbf{z}=H_{\mathbf{J}}\left(\mathbf{z}\right)$ implies
$H^{\circ\left|\mathbf{J}\right|}\left(\mathbf{z}\right)\neq H_{\mathbf{J}}\left(\mathbf{z}\right)$.
Hence, $H_{\mathbf{J}}\left(\mathbf{z}\right)$ is a wrong value of
$H$ with seed $\mathbf{z}$. \textbf{Lemma} \ref{lem:MD wrong values lemma}
then forces $H_{\mathbf{J}}\left(\mathbf{z}\right)\in\mathbb{Q}_{p}^{d}\backslash\mathbb{Z}_{p}^{d}$.
However, $H_{\mathbf{J}}\left(\mathbf{z}\right)=\mathbf{z}$, and
$\mathbf{z}$ was given to be in $\mathbb{Z}_{p}^{d}$. This is impossible!

Consequently, $H_{\mathbf{J}}\left(\mathbf{z}\right)=\mathbf{z}$
implies $H^{\circ\left|\mathbf{J}\right|}\left(\mathbf{z}\right)=\mathbf{z}$.

Q.E.D.

\vphantom{}

Now comes the multi-dimensional Correspondence Principle. Like in
the one-dimensional case, we start by proving a correspondence between
cycles of $H$ and rational-integer values attained by $\chi_{H}$.
The task of refining this correspondence to one involving the individual
periodic points of $H$ will be dealt with in a Corollary. 
\begin{thm}[\textbf{Multi-Dimensional Correspondence Principle, Ver. 1}]
\label{thm:MD CP v1}Let $H$ be a semi-basic $p$-Hydra map of dimension
$d$ and depth $r$ that fixes $\mathbf{0}$. \index{Correspondence Principle}
Then:

\vphantom{}

I. Let $\Omega$ be any cycle of $H$ in $\mathbb{Z}^{d}$, with $\left|\Omega\right|\geq2$.
Then, there exist $\mathbf{x}\in\Omega$ and an $\mathbf{n}\in\mathbb{N}_{0}^{r}\backslash\left\{ \mathbf{0}\right\} $
so that: 
\begin{equation}
\chi_{H}\left(\mathbf{B}_{p}\left(\mathbf{n}\right)\right)\overset{\mathbb{Z}_{q_{H}}^{d}}{=}\mathbf{x}
\end{equation}

\vphantom{}

II. Suppose also that $H$ is integral, and let $\mathbf{n}\in\mathbb{N}_{0}^{r}$.
If: 
\begin{equation}
\chi_{H}\left(\mathbf{B}_{p}\left(\mathbf{n}\right)\right)\in\mathbb{Z}_{q_{H}}^{d}\cap\mathbb{Z}^{d}
\end{equation}
then $\chi_{H}\left(\mathbf{B}_{p}\left(\mathbf{n}\right)\right)$
is a periodic point of $H$ in $\mathbb{Z}^{d}$. 
\end{thm}
Proof:

I. Let $\Omega$ be a cycle of $H$ in $\mathbb{Z}^{d}$ with $\left|\Omega\right|\geq2$.
Observe that for any $\mathbf{x}\in\Omega$ and any $\mathbf{J}\in\textrm{String}^{r}\left(p\right)$
for which $\left|\mathbf{J}\right|=\left|\Omega\right|\geq2$ and
$H_{\mathbf{J}}\left(\mathbf{x}\right)=\mathbf{x}$, it must be that
$\mathbf{J}$ contains at least one non-zero $r$-tuple; that is,
in writing $\mathbf{J}=\left(\mathbf{j}_{1},\ldots,\mathbf{j}_{\left|\mathbf{J}\right|}\right)$,
there is at least one $\mathbf{j}_{k}$ which is not $\mathbf{0}$.
As in the one-dimensional case, we can assume without loss of generality
that we have chosen $\mathbf{x}\in\Omega$ for which the known non-zero
$r$-tuple of $\mathbf{J}$ occurs in the bottom-most row of $\mathbf{J}$
(that is, $\mathbf{j}_{\left|\mathbf{J}\right|}\neq\mathbf{0}$),
viewing $\mathbf{J}$ as a matrix.

Writing the affine map $H_{\mathbf{J}}$ in affine form gives: 
\begin{equation}
\mathbf{x}=H_{\mathbf{J}}\left(\mathbf{x}\right)=H_{\mathbf{J}}^{\prime}\left(\mathbf{0}\right)\mathbf{x}+H_{\mathbf{J}}\left(\mathbf{0}\right)=M_{H}\left(\mathbf{J}\right)\mathbf{x}+\chi_{H}\left(\mathbf{J}\right)
\end{equation}
where all the equalities are in $\mathbb{Q}^{d}$. Then, noting the
implication: 
\begin{equation}
H_{\mathbf{J}}\left(\mathbf{x}\right)=\mathbf{x}\Rightarrow H_{\mathbf{J}^{\wedge m}}\left(\mathbf{x}\right)=\mathbf{x}
\end{equation}
we have: 
\begin{equation}
\mathbf{x}=H_{\mathbf{J}^{\wedge m}}\left(\mathbf{x}\right)=M_{H}\left(\mathbf{J}^{\wedge m}\right)\mathbf{x}+\chi_{H}\left(\mathbf{J}^{\wedge m}\right)=\left(M_{H}\left(\mathbf{J}\right)\right)^{m}\mathbf{x}+\chi_{H}\left(\mathbf{J}^{\wedge m}\right)
\end{equation}
Since $\mathbf{J}$ contains a tuple which is \emph{not} the zero
tuple, the decay estimate from \textbf{Proposition \ref{prop:MD M_H decay estimate}}
tells us that that $\left\Vert \left(M_{H}\left(\mathbf{J}\right)\right)^{m}\right\Vert _{q_{H}}$
tends to $0$ in $\mathbb{R}$ as $m\rightarrow\infty$. Moreover,
the $\mathbf{n}\in\mathbb{N}_{0}^{r}$ represented by $\mathbf{J}$
is necessarily non-$\mathbf{0}$. So, we can apply \textbf{Proposition
\ref{prop:MD self-concatenation identity}}: 
\begin{equation}
\mathbf{J}^{\wedge m}\sim\frac{1-p^{m\lambda_{p}\left(\mathbf{n}\right)}}{1-p^{\lambda_{p}\left(\mathbf{n}\right)}}\mathbf{n}
\end{equation}
and, like in the proof of \textbf{Lemma \ref{lem:MD Chi_H o B functional equation}},
we can write: 
\begin{align*}
\lim_{m\rightarrow\infty}\chi_{H}\left(\mathbf{J}^{\wedge m}\right) & \overset{\mathbb{Z}_{q_{H}}^{d}}{=}\lim_{m\rightarrow\infty}\chi_{H}\left(\frac{1-p^{m\lambda_{p}\left(\mathbf{n}\right)}}{1-p^{\lambda_{p}\left(\mathbf{n}\right)}}\mathbf{n}\right)\\
 & \overset{\mathbb{Z}_{q_{H}}^{d}}{=}\chi_{H}\left(\frac{\mathbf{n}}{1-p^{\lambda_{p}\left(\mathbf{n}\right)}}\right)\\
 & =\chi_{H}\left(\mathbf{B}_{p}\left(\mathbf{n}\right)\right)
\end{align*}
Finally, letting $m\rightarrow\infty$ in the equation:

\begin{equation}
\mathbf{x}=\left(M_{H}\left(\mathbf{J}\right)\right)^{m}\mathbf{x}+\chi_{H}\left(\mathbf{J}^{\wedge m}\right)=\left(M_{H}\left(\mathbf{J}\right)\right)^{m}\mathbf{x}+\chi_{H}\left(\frac{1-p^{m\lambda_{p}\left(\mathbf{n}\right)}}{1-p^{\lambda_{p}\left(\mathbf{n}\right)}}\mathbf{n}\right)
\end{equation}
we obtain the equality: 
\begin{equation}
\mathbf{x}\overset{\mathbb{Z}_{q_{H}}^{d}}{=}\chi_{H}\left(\frac{\mathbf{n}}{1-p^{\lambda_{p}\left(\mathbf{n}\right)}}\right)=\chi_{H}\left(\mathbf{B}_{p}\left(\mathbf{n}\right)\right)
\end{equation}
This proves the existence of the desired $\mathbf{x}$ and $\mathbf{n}$.

\vphantom{}

II. Let $H$ be integral; then, by \textbf{Lemma \ref{lem:MD integrality lemma}},
$H$ is proper. Next, let $\mathbf{n}\in\mathbb{N}_{0}^{r}\backslash\left\{ \mathbf{0}\right\} $
be such that $\chi_{H}\left(\mathbf{B}_{p}\left(\mathbf{n}\right)\right)$
is in $\mathbb{Z}_{q_{H}}^{d}\cap\mathbb{Z}^{d}$, and let $\mathbf{J}\in\textrm{String}^{r}\left(p\right)$
be the shortest block string representing $\mathbf{n}$. Then, by
$\chi_{H}$'s concatenation identity (equation \ref{eq:MD Chi_H concatenation identity}
from \textbf{Lemma \ref{lem:MD Chi_H functional equations and characterization}}),
we can write: 
\begin{equation}
\chi_{H}\left(\frac{1-p^{k\lambda_{p}\left(\mathbf{n}\right)}}{1-p^{\lambda_{p}\left(\mathbf{n}\right)}}\mathbf{n}\right)=\chi_{H}\left(\mathbf{J}^{\wedge k}\right)=H\left(\chi_{H}\left(\mathbf{J}^{\wedge\left(k-1\right)}\right)\right)=\cdots=H^{\circ\left(k-1\right)}\left(\chi_{H}\left(\mathbf{J}\right)\right)\label{eq:MD self-concatenation identity for Chi_H}
\end{equation}
Letting $k\rightarrow\infty$, \textbf{Lemma \ref{lem:MD Chi_H functional equations and characterization}}
along with the convergence of $\frac{1-p^{k\lambda_{p}\left(\mathbf{n}\right)}}{1-p^{\lambda_{p}\left(\mathbf{n}\right)}}\mathbf{n}$
to $\frac{\mathbf{n}}{1-p^{\lambda_{p}\left(\mathbf{n}\right)}}$
in $\mathbb{Z}_{p}^{r}$ then yields the $q$-adic equality: 
\begin{equation}
\lim_{k\rightarrow\infty}H_{\mathbf{J}}^{\circ\left(k-1\right)}\left(\chi_{H}\left(\mathbf{n}\right)\right)\overset{\mathbb{Z}_{q_{H}}^{d}}{=}\chi_{H}\left(\frac{\mathbf{n}}{1-p^{\lambda_{p}\left(\mathbf{n}\right)}}\right)=\chi_{H}\left(\mathbf{B}_{p}\left(\mathbf{n}\right)\right)\label{eq:Iterating H_bold-J on Chi_H}
\end{equation}

Because $H$ is semi-basic, the non-zero entries of any $\mathbf{A}_{\mathbf{j}}$
are co-prime to the non-zero entries of any $\mathbf{D}_{\mathbf{k}}$.
Moreover, since $\left\Vert \mathbf{A}_{\mathbf{j}}\right\Vert _{q_{H}}\leq1/q_{H}$
for all $\mathbf{j}\in\left(\mathbb{Z}^{r}/p\mathbb{Z}^{r}\right)\backslash\left\{ \mathbf{0}\right\} $,
this guarantees that for any finite length block string such as our
$\mathbf{J}$, the matrix $H_{\mathbf{J}}^{\prime}\left(\mathbf{0}\right)$
and the $d$-tuple $H_{\mathbf{J}}\left(\mathbf{0}\right)$ have entries
in $\mathbb{Q}$ whose $q_{H}$-adic absolute values are $\leq1$;
that is, their entries are in $\mathbb{Q}\cap\mathbb{Z}_{q_{H}}$.
Consequently, the affine linear map: 
\begin{equation}
\mathbf{w}\in\mathbb{Z}_{q_{H}}^{d}\mapsto H_{\mathbf{J}}\left(\mathbf{w}\right)=H_{\mathbf{J}}^{\prime}\left(\mathbf{0}\right)\mathbf{w}+H_{\mathbf{J}}\left(\mathbf{0}\right)\in\mathbb{Z}_{q_{H}}^{d}
\end{equation}
defines a continuous self-map of $\mathbb{Z}_{q_{H}}^{d}$. With this,
we can then write: 
\begin{align*}
H_{\mathbf{J}}\left(\chi_{H}\left(\mathbf{B}_{p}\left(\mathbf{n}\right)\right)\right) & \overset{\mathbb{Z}_{q_{H}}^{d}}{=}H_{\mathbf{j}}\left(\lim_{k\rightarrow\infty}H_{\mathbf{J}}^{\circ\left(k-1\right)}\left(\chi_{H}\left(\mathbf{n}\right)\right)\right)\\
 & \overset{\mathbb{Z}_{q_{H}}^{d}}{=}\lim_{k\rightarrow\infty}H_{\mathbf{J}}^{\circ k}\left(\chi_{H}\left(\mathbf{n}\right)\right)\\
 & \overset{\mathbb{Z}_{q_{H}}^{d}}{=}\lim_{k\rightarrow\infty}H_{\mathbf{J}}^{\circ\left(k-1\right)}\left(\chi_{H}\left(\mathbf{n}\right)\right)\\
 & \overset{\mathbb{Z}_{q_{H}}^{d}}{=}\chi_{H}\left(\mathbf{B}_{p}\left(\mathbf{n}\right)\right)
\end{align*}
This proves that the rational $d$-tuple $\chi_{H}\left(\mathbf{B}_{p}\left(\mathbf{n}\right)\right)$
is a fixed point of the composition sequence $H_{\mathbf{J}}$ as
an element of $\mathbb{Z}_{q_{H}}^{d}$.

Now, since we assumed that $\chi_{H}\left(\mathbf{B}_{p}\left(\mathbf{n}\right)\right)$
was in $\mathbb{Z}^{d}$, the fact that $H_{\mathbf{J}}$ is a map
on $\mathbb{Q}^{d}$, tells us that the equality $H_{\mathbf{J}}\left(\chi_{H}\left(\mathbf{B}_{p}\left(\mathbf{n}\right)\right)\right)=\chi_{H}\left(\mathbf{B}_{p}\left(\mathbf{n}\right)\right)$\textemdash nominally
occurring in $\mathbb{Z}_{q_{H}}^{d}$\textemdash is actually valid
in $\mathbb{Q}^{d}$. Consequently, this equality is also valid in
$\mathbb{Q}_{p}^{d}$, and so, $\chi_{H}\left(\mathbf{B}_{p}\left(\mathbf{n}\right)\right)$
is a fixed point of the map $H_{\mathbf{J}}:\mathbb{Q}_{p}^{d}\rightarrow\mathbb{Q}_{p}^{d}$.

Next, because the assumed integrality of $H$ makes $H$ proper, we
can apply \textbf{Lemma \ref{lem:MD properness lemma}} to conclude
that $\chi_{H}\left(\mathbf{B}_{p}\left(\mathbf{n}\right)\right)\in\mathbb{Z}^{d}\subseteq\mathbb{Z}_{p}^{d}$,
as a fixed point of $H_{\mathbf{J}}$, is in fact a fixed point of
$H^{\circ\left|\mathbf{J}\right|}$: 
\begin{equation}
H^{\circ\left|\mathbf{J}\right|}\left(\chi_{H}\left(\mathbf{B}_{p}\left(\mathbf{n}\right)\right)\right)=\chi_{H}\left(\mathbf{B}_{p}\left(\mathbf{n}\right)\right)
\end{equation}
Lastly, because $\chi_{H}\left(\mathbf{B}_{p}\left(\mathbf{n}\right)\right)$
is in $\mathbb{Z}^{d}$, the fact that $H:\mathbb{Z}^{d}\rightarrow\mathbb{Z}^{d}$
is equal to the restriction to $\mathbb{Z}^{d}$ of $H:\mathbb{Z}_{p}^{d}\rightarrow\mathbb{Z}_{p}^{d}$
necessarily forces $\chi_{H}\left(\mathbf{B}_{p}\left(\mathbf{n}\right)\right)$
to be a periodic point of $H$ in $\mathbb{Z}^{d}$, as desired.

Q.E.D.

\vphantom{}

Like with the one-dimensional case, we now have several alternative
versions.
\begin{cor}[\textbf{Multi-Dimensional Correspondence Principle, Ver. 2}]
\label{cor:MD CP v2}Suppose that $H$ is semi-basic\index{Correspondence Principle}
$p$-Hydra map of dimension $d$ and depth $r$ that fixes $\mathbf{0}$.
Then:

\vphantom{}

I. For every cycle $\Omega\subseteq\mathbb{Z}^{d}$ of $H$, viewing
$\Omega\subseteq\mathbb{Z}^{d}$ as a subset of $\mathbb{Z}_{q_{H}}^{d}$,
the intersection $\chi_{H}\left(\mathbb{Z}_{p}^{r}\right)\cap\Omega$
is non-empty. Moreover, for every $\mathbf{x}\in\chi_{H}\left(\mathbb{Z}_{p}^{r}\right)\cap\Omega$,
there is an $\mathbf{n}\in\mathbb{N}_{0}^{r}\backslash\left\{ \mathbf{0}\right\} $
so that: 
\[
\mathbf{x}=\left(\mathbf{I}_{d}-M_{H}\left(\mathbf{n}\right)\right)^{-1}\chi_{H}\left(\mathbf{n}\right)=\chi_{H}\left(\mathbf{B}_{p}\left(\mathbf{n}\right)\right)
\]
where the equality occurs in $\mathbb{Q}^{d}$.

\vphantom{}

II. Suppose in addition that $H$ is integral. Let $\mathbf{n}\in\mathbb{N}_{0}^{r}\backslash\left\{ \mathbf{0}\right\} $.
If the tuple $\mathbf{x}\in\mathbb{Z}_{q_{H}}^{d}$ given by: 
\[
\mathbf{x}=\left(\mathbf{I}_{d}-M_{H}\left(\mathbf{n}\right)\right)^{-1}\chi_{H}\left(\mathbf{n}\right)=\chi_{H}\left(\mathbf{B}_{p}\left(\mathbf{n}\right)\right)
\]
is an element of $\mathbb{Z}^{d}$, then $\mathbf{x}$ is necessarily
a periodic point of $H$. 
\end{cor}
Proof: Re-write the results of \textbf{Theorem \ref{thm:MD CP v1}}
using \textbf{Lemma \ref{lem:MD Chi_H o B functional equation}}.

Q.E.D. 
\begin{cor}[\textbf{Multi-Dimensional Correspondence Principle, Ver. 3}]
\label{cor:MD CP v3}Suppose that $H$ is an integral semi-basic
\index{Correspondence Principle}$p$-Hydra map of dimension $d$
and depth $r$ that fixes $\mathbf{0}$. Then, the set of non-zero
periodic points of $H$ is equal to $\mathbb{Z}^{d}\cap\chi_{H}\left(\mathbb{Q}^{r}\cap\left(\mathbb{Z}_{p}^{r}\right)^{\prime}\right)$
(where, recall, $\left(\mathbb{Z}_{p}^{r}\right)^{\prime}=\mathbb{Z}_{p}^{\prime}\backslash\mathbb{N}_{0}^{r}$).
\end{cor}
\begin{rem}
The implication ``If $\mathbf{x}\in\mathbb{Z}^{d}\backslash\left\{ \mathbf{0}\right\} $
is a periodic point, then $\mathbf{x}\in\chi_{H}\left(\mathbb{Q}^{r}\cap\left(\mathbb{Z}_{p}^{r}\right)^{\prime}\right)$''
\emph{does not }require $H$ to be integral.
\end{rem}
Proof: Before we do anything else, note that because $H$ is integral
and semi-basic, \textbf{Lemma \ref{lem:MD integrality lemma} }guarantees
$H$ will be proper.

I. Let $\mathbf{x}$ be a non-zero periodic point of $H$, and let
$\Omega$ be the unique cycle of $H$ in $\mathbb{Z}$ which contains
$\mathbf{x}$. By Version 1 of the Correspondence Principle (\textbf{Theorem
\ref{thm:MD CP v1}}), there exists a $\mathbf{y}\in\Omega$ and a
$\mathbf{z}=B_{p}\left(\mathbf{n}\right)\subset\mathbb{Z}_{p}^{r}$
(for some $\mathbf{n}\in\mathbb{N}_{0}^{r}$ containing at least one
non-zero entry) so that $\chi_{H}\left(\mathbf{z}\right)=\mathbf{y}$.
Because $\mathbf{y}$ is in $\Omega$, there is an $k\geq1$ so that
$\mathbf{x}=H^{\circ k}\left(\mathbf{y}\right)$. As such, there is
a \emph{unique} length $k$ block-string $\mathbf{I}\in\textrm{String}^{r}\left(p\right)$
so that $H_{\mathbf{I}}\left(\mathbf{y}\right)=H^{\circ k}\left(\mathbf{y}\right)=\mathbf{x}$.

Now, let $\mathbf{J}\in\textrm{String}_{\infty}^{r}\left(p\right)$
represent $\mathbf{z}$; note that $\mathbf{J}$ is infinite and that
its columns are periodic. Using $\chi_{H}$'s concatenation identity
(\textbf{Lemma \ref{eq:MD Chi_H concatenation identity}}), we have:
\begin{equation}
\mathbf{x}=H_{\mathbf{I}}\left(\mathbf{y}\right)=H_{\mathbf{I}}\left(\chi_{H}\left(\mathbf{z}\right)\right)=\chi_{H}\left(\mathbf{I}\wedge\mathbf{J}\right)
\end{equation}
Next, let $\mathbf{w}\in\mathbb{Z}_{p}^{r}$ be the $r$-tuple of
$p$-adic integers represented by the infinite block-string $\mathbf{I}\wedge\mathbf{J}$;
note that $\mathbf{w}$ is \emph{not} an element of $\mathbb{N}_{0}^{r}$.
By the above, $\chi_{H}\left(\mathbf{w}\right)=\mathbf{x}$, and hence
that $\mathbf{x}\in\mathbb{Z}^{d}\cap\chi_{H}\left(\mathbb{Z}_{p}^{r}\right)$.

Finally, since $\mathbf{z}=B_{p}\left(\mathbf{n}\right)$ and $\mathbf{n}\neq\mathbf{0}$,
the $p$-adic digits of $\mathbf{z}$'s entries are periodic, which
forces $\mathbf{z}$ to be an element of $\mathbb{Q}^{r}\cap\left(\mathbb{Z}_{p}^{r}\right)^{\prime}$.
So, letting $\mathbf{m}$ be the rational integer $r$-tuple represented
by length-$k$ block-string $\mathbf{I}$, we have: 
\begin{equation}
\mathbf{w}\sim\mathbf{I}\wedge\mathbf{J}\sim\mathbf{m}+p^{\lambda_{p}\left(\mathbf{m}\right)}\mathbf{z}
\end{equation}
This shows that $\mathbf{w}\in\mathbb{Q}^{r}\cap\left(\mathbb{Z}_{p}^{r}\right)^{\prime}$,
and hence, that: 
\[
\mathbf{x}=\chi_{H}\left(\mathbf{w}\right)\in\mathbb{Z}^{d}\cap\chi_{H}\left(\mathbb{Q}^{r}\cap\left(\mathbb{Z}_{p}^{r}\right)^{\prime}\right)
\]

\vphantom{}

II. Suppose $\mathbf{x}\in\mathbb{Z}^{d}\cap\chi_{H}\left(\mathbb{Q}^{r}\cap\left(\mathbb{Z}_{p}^{r}\right)^{\prime}\right)$,
with $\mathbf{x}=\chi_{H}\left(\mathbf{z}\right)$ for some $\mathbf{z}\in\mathbb{Q}^{r}\cap\left(\mathbb{Z}_{p}^{r}\right)^{\prime}$.
As a $r$-tuple of rational numbers are all $p$-adic integers which
are \emph{not} elements of $\mathbb{N}_{0}$, the $p$-adic digits
of the entries of $\mathbf{z}$ are all eventually periodic. As such,
there are $\mathbf{m},\mathbf{n}\in\mathbb{N}_{0}^{r}$ (with $\mathbf{n}\neq\mathbf{0}$)
so that: 
\begin{equation}
\mathbf{z}=\mathbf{m}+p^{\lambda_{p}\left(\mathbf{m}\right)}B_{p}\left(\mathbf{n}\right)
\end{equation}
Here, the $p$-adic digits of the entries of $\mathbf{n}$ generate
the periodic parts of the digits of $\mathbf{z}$'s entries. On the
other hand, the $p$-adic digits of $\mathbf{m}$'s entries are the
finite-length sequences in the digits of the entries of $\mathbf{z}$
that occur before the periodicity sets in.

Now, let $\mathbf{I}$ be the finite block string representing $\mathbf{m}$,
and let $\mathbf{J}$ be the infinite block-string representing $B_{p}\left(\mathbf{n}\right)$.
Then, $\mathbf{z}=\mathbf{I}\wedge\mathbf{J}$. Using \textbf{Lemmata
\ref{eq:MD Chi_H concatenation identity}} and \textbf{\ref{lem:MD Chi_H o B functional equation}}
(the concatenation identity and $\chi_{H}\circ\mathbf{B}_{p}$ functional
equation, respectively), we have: 
\begin{equation}
\mathbf{x}=\chi_{H}\left(\mathbf{I}\wedge\mathbf{J}\right)=H_{\mathbf{I}}\left(\chi_{H}\left(\mathbf{J}\right)\right)=H_{\mathbf{I}}\left(\chi_{H}\left(B_{p}\left(\mathbf{n}\right)\right)\right)=H_{\mathbf{I}}\left(\left(\mathbf{I}_{d}-M_{H}\left(\mathbf{n}\right)\right)^{-1}\chi_{H}\left(\mathbf{n}\right)\right)
\end{equation}
Here, let $\mathbf{y}\overset{\textrm{def}}{=}\left(\mathbf{I}_{d}-M_{H}\left(\mathbf{n}\right)\right)^{-1}\chi_{H}\left(\mathbf{n}\right)$
is a $d$-tuple of rational numbers.
\begin{claim}
$\left\Vert \mathbf{y}\right\Vert _{p}\leq1$.

Proof of claim: First, since $\mathbf{y}\in\mathbb{Q}^{d}$, it also
lies in $\mathbb{Q}_{p}^{d}$. So, by way of contradiction, suppose
$\left\Vert \mathbf{y}\right\Vert _{p}>1$. Since $H$ is proper,
\textbf{Lemma \ref{lem:MD wrong values lemma}} tells us that no matter
which branches of $H$ are specified by $\mathbf{I}$, the composition
sequence $H_{\mathbf{I}}$ maps $\mathbf{y}\in\mathbb{Q}_{p}^{d}\backslash\mathbb{Z}_{p}^{d}$
to $H_{\mathbf{I}}\left(\mathbf{y}\right)\in\mathbb{Q}_{p}^{d}\backslash\mathbb{Z}_{p}^{d}$,
and so $\left\Vert H_{\mathbf{i}}\left(\mathbf{y}\right)\right\Vert _{p}>1$.
However, $H_{\mathbf{i}}\left(\mathbf{y}\right)=\mathbf{x}$, and
$\mathbf{x}\in\mathbb{N}_{0}^{d}$; hence, $1<\left\Vert H_{\mathbf{i}}\left(\mathbf{y}\right)\right\Vert _{p}=\left\Vert \mathbf{x}\right\Vert _{p}\leq1$.
This is impossible! So, it must be that $\left\Vert \mathbf{y}\right\Vert _{p}\leq1$.
This proves the claim.
\end{claim}
\begin{claim}
\label{claim:5}$\mathbf{x}=H^{\circ\left|\mathbf{I}\right|}\left(\mathbf{y}\right)$

Proof of claim: Suppose the equality failed. Then $H_{\mathbf{I}}\left(\mathbf{y}\right)=\mathbf{x}\neq H^{\circ\left|\mathbf{I}\right|}\left(\mathbf{y}\right)$,
and so $\mathbf{x}=H_{\mathbf{I}}\left(\mathbf{y}\right)$ is a wrong
value of $H$ with seed $\mathbf{y}$. Because $H$ is proper,\textbf{
Lemma \ref{lem:MD wrong values lemma}} forces $\left\Vert \mathbf{x}\right\Vert _{p}=\left\Vert H_{\mathbf{I}}\left(\mathbf{y}\right)\right\Vert _{p}>1$.
However, like in the one-dimensional case, this is just as impossible
as it was in the previous paragraph, seeing as $\left\Vert \mathbf{x}\right\Vert _{p}\leq1$.
This proves the claim.
\end{claim}
\vphantom{}

Finally, let $\mathbf{V}$ be the shortest block-string representing
$\mathbf{n}$, so that $\mathbf{J}$ (the block-string representing
$B_{p}\left(\mathbf{n}\right)$) is obtained by concatenating infinitely
many copies of $\mathbf{v}$. Because $q_{H}$ is co-prime to all
the non-zero entries of $\mathbf{D}_{\mathbf{j}}$ for all $\mathbf{j}\in\mathbb{Z}^{r}/p\mathbb{Z}^{r}$,
\emph{note that $H_{\mathbf{V}}$ is continuous on $\mathbb{Z}_{q_{H}}^{d}$}.
As such:
\begin{equation}
\chi_{H}\left(B_{p}\left(\mathbf{n}\right)\right)\overset{\mathbb{Z}_{q_{H}}^{d}}{=}\lim_{k\rightarrow\infty}\chi_{H}\left(\mathbf{V}^{\wedge k}\right)
\end{equation}
implies:
\begin{align*}
H_{\mathbf{V}}\left(\chi_{H}\left(B_{p}\left(\mathbf{n}\right)\right)\right) & \overset{\mathbb{Z}_{q_{H}}^{d}}{=}\lim_{k\rightarrow\infty}H_{\mathbf{V}}\left(\chi_{H}\left(\mathbf{V}^{\wedge k}\right)\right)\\
 & \overset{\mathbb{Z}_{q_{H}}^{d}}{=}\lim_{k\rightarrow\infty}\chi_{H}\left(\mathbf{V}^{\wedge\left(k+1\right)}\right)\\
 & \overset{\mathbb{Z}_{q_{H}}^{d}}{=}\chi_{H}\left(B_{p}\left(\mathbf{n}\right)\right)
\end{align*}
Hence, $H_{\mathbf{V}}\left(\mathbf{y}\right)=\mathbf{y}$. Since
we showed that $\left\Vert \mathbf{y}\right\Vert _{p}\leq1$, the
propriety of $H$ allows us to apply \textbf{Lemma \ref{lem:properness lemma}}
and conclude that $H_{\mathbf{V}}\left(\mathbf{y}\right)=H^{\circ\left|\mathbf{V}\right|}\left(\mathbf{y}\right)$.
Thus, $\mathbf{y}$ is a periodic point of $H:\mathbb{Z}_{p}^{d}\rightarrow\mathbb{Z}_{p}^{d}$.

By \textbf{Claim \ref{claim:5}}, $H$ iterates $\mathbf{y}$ to $\mathbf{x}$.
Since $\mathbf{y}$ is a periodic point of $H$ in $\mathbb{Z}_{p}^{d}$,
this forces $\mathbf{x}$ and $\mathbf{y}$ to belong to the same
cycle of $H$ in $\mathbb{Z}_{p}^{d}$, with $\mathbf{x}$ being a
periodic point of $H$. As such, just as $H$ iterates $\mathbf{y}$
to $\mathbf{x}$, so too does $H$ iterate $\mathbf{x}$ to $\mathbf{y}$.
Finally, since $\mathbf{x}\in\mathbb{N}_{0}^{d}$, so too is $\mathbf{y}$.
integer as well. Thus, $\mathbf{x}$ belongs to a cycle of $H$ in
$\mathbb{Z}^{d}$, as desired.

Q.E.D.

\vphantom{}

Like in the one-dimensional case, the fact that the orbit classes
of $H$'s attracting cycles, isolated cycles, and divergent trajectories
constitute a partition of $\mathbb{Z}^{d}$ allows us to use the Correspondence
Principle to make conclusions about the divergent trajectories. 
\begin{prop}
\label{prop:MD Preparing for application to divergent trajectories}Let
$H$ be as given in \textbf{\emph{Corollary \ref{cor:MD CP v3}}}.
Additionally, suppose that $\left\Vert H_{\mathbf{j}}\left(\mathbf{0}\right)\right\Vert _{q_{H}}=1$
for all $\mathbf{j}\in\left(\mathbb{Z}^{r}/p\mathbb{Z}^{r}\right)\backslash\left\{ \mathbf{0}\right\} $.
Then $\chi_{H}\left(\mathbf{z}\right)\neq\mathbf{0}$ for any $\mathbf{z}\in\mathbb{Z}_{p}^{r}\backslash\mathbb{Q}^{r}$. 
\end{prop}
Proof: Let $H$ as given. By way of contradiction, suppose that $\chi_{H}\left(\mathbf{z}\right)=\mathbf{0}$
for some $\mathbf{z}=\left(\mathfrak{z}_{1},\ldots,\mathfrak{z}_{r}\right)\in\mathbb{Z}_{p}^{r}\backslash\mathbb{Q}^{r}$.
Let $\mathbf{j}\in\mathbb{Z}^{r}/p\mathbb{Z}^{r}$ denote $\left[\mathbf{z}\right]_{p}$,
and let $\mathbf{z}^{\prime}$ denote: 
\begin{equation}
\mathbf{z}^{\prime}\overset{\textrm{def}}{=}\frac{\mathbf{z}-\mathbf{j}}{p}=\left(\frac{\mathfrak{z}_{1}-j_{1}}{p},\ldots,\frac{\mathfrak{z}_{r}-j_{r}}{p}\right)
\end{equation}
Then, we can write: 
\[
\mathbf{0}=\chi_{H}\left(\mathbf{z}\right)=\chi_{H}\left(p\mathbf{z}^{\prime}+\mathbf{j}\right)=H_{\mathbf{j}}\left(\chi_{H}\left(\mathbf{z}^{\prime}\right)\right)
\]
Next, suppose $\mathbf{j}=\mathbf{0}$. Because $H\left(\mathbf{0}\right)=\mathbf{0}$,
we have: 
\[
\mathbf{0}=H_{\mathbf{0}}\left(\chi_{H}\left(\mathbf{z}^{\prime}\right)\right)=\mathbf{D}_{\mathbf{0}}^{-1}\mathbf{A}_{\mathbf{0}}\chi_{H}\left(\mathbf{z}^{\prime}\right)
\]
Since $\mathbf{A}_{\mathbf{0}}$ and $\mathbf{D}_{\mathbf{0}}$ are
invertible, this forces $\chi_{H}\left(\mathbf{z}^{\prime}\right)=\mathbf{0}$.
In this way, if the the first $n$ $p$-adic digits of $\mathfrak{z}_{m}$
are $0$ for all $m\in\left\{ 1,\ldots,r\right\} $, we can pull those
digits out from $\mathbf{z}$ to obtain a $p$-adic integer tuple
of unit $p$-adic norm ($\left\Vert \cdot\right\Vert _{p}$). If \emph{all
}of the digits of all the $\mathfrak{z}_{m}$s are $0$, this then
makes $\mathbf{z}=\mathbf{0}$, contradicting that $\mathbf{z}$ was
given to be an element of $\mathbb{Z}_{p}^{r}\backslash\mathbb{Q}^{r}$.
So, without loss of generality, we can assume that $\mathbf{j}\neq\mathbf{0}$.

Hence: 
\begin{align*}
\mathbf{0} & =H_{\mathbf{j}}\left(\chi_{H}\left(\mathbf{z}^{\prime}\right)\right)\\
 & =H_{\mathbf{j}}^{\prime}\left(\mathbf{0}\right)\chi_{H}\left(\mathbf{z}^{\prime}\right)+H_{\mathbf{j}}\left(\mathbf{0}\right)\\
\left(H_{\mathbf{j}}^{\prime}\left(\mathbf{0}\right)=\mathbf{D}_{\mathbf{j}}^{-1}\mathbf{A}_{\mathbf{j}}\right); & \Updownarrow\\
\mathbf{D}_{\mathbf{j}}^{-1}\mathbf{A}_{\mathbf{j}}\chi_{H}\left(\mathbf{z}^{\prime}\right) & =-H_{\mathbf{j}}\left(\mathbf{0}\right)
\end{align*}
Because $H$ is semi-basic, $\mathbf{j}\neq\mathbf{0}$ tells us that
$\left\Vert \mathbf{A}_{\mathbf{j}}\right\Vert _{q_{H}}<1$ and $\left\Vert \mathbf{D}_{\mathbf{j}}\right\Vert _{q_{H}}=1$.
So: 
\begin{align*}
\left\Vert H_{\mathbf{j}}\left(\mathbf{0}\right)\right\Vert _{q_{H}} & =\left\Vert \mathbf{D}_{\mathbf{j}}^{-1}\mathbf{A}_{\mathbf{j}}\chi_{H}\left(\mathbf{z}^{\prime}\right)\right\Vert _{q_{H}}\\
 & \leq\underbrace{\left\Vert \mathbf{D}_{\mathbf{j}}^{-1}\right\Vert _{q_{H}}\left\Vert \mathbf{A}_{\mathbf{j}}\right\Vert _{q_{H}}}_{<1}\left\Vert \chi_{H}\left(\mathbf{z}^{\prime}\right)\right\Vert _{q_{H}}\\
 & <\left\Vert \chi_{H}\left(\mathbf{z}^{\prime}\right)\right\Vert _{q_{H}}
\end{align*}
Since $\chi_{H}\left(\mathbb{Z}_{p}^{r}\right)\subseteq\mathbb{Z}_{q_{H}}^{d}$,
we see that $\left\Vert \chi_{H}\left(\mathbf{z}^{\prime}\right)\right\Vert _{q_{H}}\leq1$,
which forces $\left\Vert H_{\mathbf{j}}\left(\mathbf{0}\right)\right\Vert _{q_{H}}<1$
for our non-zero $\mathbf{j}$. However, we were given that $\left\Vert H_{\mathbf{j}}\left(\mathbf{0}\right)\right\Vert _{q_{H}}=1$
for all $\mathbf{j}\in\left(\mathbb{Z}^{r}/p\mathbb{Z}^{r}\right)\backslash\left\{ \mathbf{0}\right\} $.
This is a contradiction!

So, for each $m$, $\mathfrak{z}_{m}$ has no non-zero $p$-adic digits.
This forces all of the $\mathfrak{z}_{m}$s to be $0$, which forces
$\mathbf{z}=\mathbf{0}$. But $\mathbf{z}$ was given to be an element
of $\mathbb{Z}_{p}^{r}\backslash\mathbb{Q}^{r}$, yet $\mathbf{0}\in\mathbb{Q}^{r}$\textemdash and
that's our contradiction.

So, it must be that $\chi_{H}\left(\mathbf{z}\right)\neq\mathbf{0}$.

Q.E.D. 
\begin{thm}
\label{thm:MD Divergent trajectories come from irrational z}Let\index{Hydra map@\emph{Hydra map!divergent trajectories}}
$H$ be as given in \textbf{\emph{Corollary \ref{cor:MD CP v3}}}.
Additionally, suppose that:

\vphantom{}

I. $H$ is contracting;

\vphantom{}

II. $\left\Vert H_{\mathbf{j}}\left(\mathbf{0}\right)\right\Vert _{q_{H}}=1$
for all $\mathbf{j}\in\left(\mathbb{Z}^{r}/p\mathbb{Z}^{r}\right)\backslash\left\{ \mathbf{0}\right\} $.

\vphantom{}

Under these hypotheses, let $\mathbf{z}\in\mathbb{Z}_{p}^{r}\backslash\mathbb{Q}^{r}$.
If $\chi_{H}\left(\mathbf{z}\right)\in\mathbb{Z}^{d}$, then $\chi_{H}\left(\mathbf{z}\right)$
belongs to a divergent orbit class of $H$. 
\end{thm}
Proof: Let $H$ and $\mathbf{z}$ be as given. By \textbf{Proposition
\ref{prop:MD Preparing for application to divergent trajectories}},
$\chi_{H}\left(\mathbf{z}\right)\neq\mathbf{0}$.

Now, by way of contradiction, suppose that $\chi_{H}\left(\mathbf{z}\right)$
did not belong to a divergent orbit class of $H$. By \textbf{Theorem
\ref{thm:orbit classes partition domain}}, every element of $\mathbb{Z}^{d}$
belongs to either a divergent orbit class of $H$, or to an orbit
class of $H$ which contains a cycle of $H$. This forces $\chi_{H}\left(\mathbf{z}\right)$
to be a pre-periodic point of $H$. As such, there is an $n\geq0$
so that $H^{\circ n}\left(\chi_{H}\left(\mathbf{z}\right)\right)$
is a periodic point of $H$. Using the functional equations of $\chi_{H}$
(\textbf{Lemma \ref{lem:MD Chi_H functional equations and characterization}}),
we can write $H^{\circ n}\left(\chi_{H}\left(\mathbf{z}\right)\right)=\chi_{H}\left(\mathbf{y}\right)$
where $\mathbf{y}=\mathbf{m}+p^{\lambda_{p}\left(\mathbf{m}\right)}\mathbf{z}$,
and where $\mathbf{m}\in\mathbb{N}_{0}^{r}$ is the unique $r$-tuple
of non-negative integers which represents the composition sequence
of branches used in iterating $H$ $n$ times at $\chi_{H}\left(\mathbf{z}\right)$
to produce $\chi_{H}\left(\mathbf{y}\right)$.

Because $\chi_{H}\left(\mathbf{z}\right)\neq\mathbf{0}$, we can use
the integrality of $H$ along with the fact that $\mathbf{0}=H\left(\mathbf{0}\right)=H_{\mathbf{0}}\left(\mathbf{0}\right)$
to conclude that $\left\{ \mathbf{0}\right\} $ is an isolated cycle
of $H$. This then guarantees that the periodic point $\chi_{H}\left(\mathbf{y}\right)$
is non-zero. Consequently, the multi-dimensional \textbf{Correspondence
Principle }tells us that $\mathbf{y}$ must be an element of $\mathbb{Q}^{r}\cap\left(\mathbb{Z}_{p}^{r}\right)^{\prime}$.
So, \emph{for each} $\ell\in\left\{ 1,\ldots,r\right\} $, the $\ell$th
entry of $\mathbf{y}$ is a rational number; therefore, its sequence
of $p$-adic digits is eventually periodic. However, since $\mathbf{y}=\mathbf{m}+p^{\lambda_{p}\left(\mathbf{m}\right)}\mathbf{z}$
where $\mathbf{z}\in\mathbb{Z}_{p}^{r}\backslash\mathbb{Q}^{r}$,
there is at least one $\ell\in\left\{ 1,\ldots,r\right\} $ so that
$\mathfrak{z}_{\ell}$'s $p$-adic digits are \emph{aperiodic}. This
is a contradiction!

So, it must be that $\chi_{H}\left(\mathbf{z}\right)$ is a divergent
point of $H$.

Q.E.D.

\vphantom{}

Like with the one-dimensional case, it would be highly desirable to
prove the following conjecture true: 
\begin{conjecture}[\textbf{Correspondence Principle for Multi-Dimensional Divergent Points?}]
\label{conj:MD correspondence theorem for divergent trajectories}Provided
that $H$ satisfies certain prerequisites such as the hypotheses of
\textbf{\emph{Theorem \ref{thm:MD Divergent trajectories come from irrational z}}},
$\mathbf{x}\in\mathbb{Z}^{d}$ belongs to a divergent trajectory under
$H$ if and only if there is a $\mathbf{z}\in\mathbb{Z}_{p}^{r}\backslash\mathbb{Q}^{r}$
so that $\chi_{H}\left(\mathbf{z}\right)\in\mathbb{Z}^{d}$.
\end{conjecture}
\newpage{}

\section{\label{sec:5.3 Rising-Continuity-in-Multiple}Rising-Continuity in
Multiple Dimensions}

IN THIS SECTION, WE FIX INTEGERS $r,d\geq1$ (WITH $r\leq d$) ALONG
WITH TWO DISTINCT PRIME NUMBERS $p$ AND $q$. UNLESS SAID OTHERWISE,
$\mathbb{K}$ DENOTES A METRICALLY COMPLETE VALUED FIELD, POSSIBLY
ARCHIMEDEAN; $K$, MEANWHILE, DENOTES A METRICALLY COMPLETE VALUED
\emph{NON-ARCHIMEDEAN} FIELD, UNLESS STATED OTHERWISE.

\subsection{\label{subsec:5.3.1 Tensor-Products}Multi-Dimensional Notation and
Tensor Products}

\begin{notation}
\index{multi-dimensional!notation}\ 

\vphantom{}

I. Like with $\left\Vert \mathbf{z}\right\Vert _{p}$ for $\mathbf{z}\in\mathbb{Z}_{p}^{r}$,
for $\hat{\mathbb{Z}}_{p}^{r}$, we write: 
\begin{equation}
\left\Vert \mathbf{t}\right\Vert _{p}\overset{\textrm{def}}{=}\max\left\{ p^{-v_{p}\left(t_{1}\right)},\ldots,p^{-v_{p}\left(t_{m}\right)}\right\} ,\textrm{ }\forall\mathbf{t}\in\hat{\mathbb{Z}}_{p}^{r}\label{eq:Definition of bold t p norm}
\end{equation}
\begin{equation}
v_{p}\left(\mathbf{t}\right)\overset{\textrm{def}}{=}\min\left\{ v_{p}\left(t_{1}\right),\ldots,v_{p}\left(t_{m}\right)\right\} ,\textrm{ }\forall\mathbf{t}\in\hat{\mathbb{Z}}_{p}^{r}\label{eq:Definition of v_P of bold t}
\end{equation}
so that $\left\Vert \mathbf{t}\right\Vert _{p}=p^{-v_{p}\left(\mathbf{t}\right)}$.
In this notation, $\left\{ \mathbf{t}:\left\Vert \mathbf{t}\right\Vert _{p}=p^{n}\right\} $
(equivalently $\left\{ \mathbf{t}:v_{p}\left(\mathbf{t}\right)=-n\right\} $)
tells us that for each $\mathbf{t}$, the inequality $\left|t_{m}\right|_{p}\leq p^{n}$
occurs for all $m$, with equality actually being achieved for \emph{at
least} one $m\in\left\{ 1,\ldots,r\right\} $.

$\left\{ \mathbf{t}:\left\Vert \mathbf{t}\right\Vert _{p}\leq p^{N}\right\} $,
meanwhile, tells us that, for each $\mathbf{t}$, $\left|t_{m}\right|_{p}\leq p^{N}$
occurs for all $m$. In particular: 
\begin{equation}
\sum_{\left\Vert \mathbf{t}\right\Vert _{p}\leq p^{N}}=\sum_{\left|t_{1}\right|_{p}\leq p^{N}}\cdots\sum_{\left|t_{r}\right|_{p}\leq p^{N}}\label{eq:Definition of sum of bold t norm _p lessthanorequalto p^N}
\end{equation}
That being said, the reader should take care to notice that:

\begin{equation}
\sum_{\left\Vert \mathbf{t}\right\Vert _{p}=p^{n}}\neq\sum_{\left|t_{1}\right|_{p}=p^{n}}\cdots\sum_{\left|t_{r}\right|_{p}=p^{n}}\label{eq:Warning about summing bold t norm _p equal to p^N}
\end{equation}

\vphantom{}

II. We write: 
\begin{equation}
\sum_{\mathbf{n}=\mathbf{0}}^{p^{N}-1}\overset{\textrm{def}}{=}\sum_{n_{1}=0}^{p^{N}-1}\cdots\sum_{n_{r}=0}^{p^{N}-1}\label{eq:Definition of sum from bold n eq zero to P^N minus 1}
\end{equation}
and define: 
\begin{equation}
\mathbf{n}\leq p^{m}-1\label{eq:Definition of bold n lessthanoreqto P^m minus 1}
\end{equation}
to indicate that $0\leq n_{\ell}\leq p^{m}-1$ for all $\ell\in\left\{ 1,\ldots,r\right\} $.
Note then that: 
\begin{equation}
\sum_{\mathbf{j}=\mathbf{0}}^{p-1}=\sum_{\mathbf{j}\in\mathbb{Z}^{r}/p\mathbb{Z}^{r}}\label{eq:Definition of sum from bold j is 0 to P minus 1}
\end{equation}
and also: 
\begin{equation}
\sum_{\mathbf{j}>\mathbf{0}}^{p-1}\overset{\textrm{def}}{=}\sum_{\mathbf{j}\in\left(\mathbb{Z}^{r}/p\mathbb{Z}^{r}\right)\backslash\left\{ \mathbf{0}\right\} }\label{eq:Definition of sum from bold j greaterthan 0 to P minus 1}
\end{equation}
Consequently, $\forall\mathbf{j}\in\mathbb{Z}^{r}/\mathbb{Z}^{r}$
and $\forall\mathbf{j}\leq p-1$ denote the same sets of $\mathbf{j}$s.

\vphantom{}

III. $\mathbf{1}_{\mathbf{0}}\left(\mathbf{t}\right)$\nomenclature{$\mathbf{1}_{\mathbf{0}}\left(\mathbf{t}\right)$}{  }
denotes the function on $\hat{\mathbb{Z}}_{p}^{r}$ which is $1$
if $\mathbf{t}\overset{1}{\equiv}\mathbf{0}$ and is $0$ otherwise.
For $\mathbf{s}\in\hat{\mathbb{Z}}_{p}^{r}$, I write $\mathbf{1}_{\mathbf{s}}\left(\mathbf{t}\right)\overset{\textrm{def}}{=}\mathbf{1}_{\mathbf{0}}\left(\mathbf{t}-\mathbf{s}\right)$.

\vphantom{}

IV. For a matrix $\mathbf{A}$ with entries in $\mathbb{Z}_{q}$,
recall we write $\left\Vert \mathbf{A}\right\Vert _{q}$ to denote
the maximum of the $q$-adic absolute values of $\mathbf{A}$'s entries.
When $\mathbf{A}$ has entries in $\overline{\mathbb{Q}}$, we write
$\left\Vert \mathbf{A}\right\Vert _{\infty}$ to denote the maximum
of the complex absolute values of $\mathbf{A}$'s entries.

\vphantom{}

V. For any $\mathbf{t}\in\hat{\mathbb{Z}}_{p}^{r}$, we write \nomenclature{$\left|\mathbf{t}\right|_{p}$}{$\left(p^{-v_{p}\left(\mathbf{t}\right)},\ldots,p^{-v_{p}\left(\mathbf{t}\right)}\right)$}$\left|\mathbf{t}\right|_{p}$
to denote the $r$\emph{-tuple}: 
\begin{equation}
\left|\mathbf{t}\right|_{p}\overset{\textrm{def}}{=}\left(p^{-v_{p}\left(\mathbf{t}\right)},\ldots,p^{-v_{p}\left(\mathbf{t}\right)}\right)\label{eq:Definition of single bars of bold t _P}
\end{equation}
Consequently: 
\begin{equation}
\frac{\mathbf{t}\left|\mathbf{t}\right|_{p}}{p}=p^{-v_{p}\left(\mathbf{t}\right)-1}\mathbf{t}=\left(p^{-v_{p}\left(\mathbf{t}\right)-1}t_{1},\ldots,p^{-v_{p}\left(\mathbf{t}\right)-1}t_{r}\right),\textrm{ }\forall\mathbf{t}\in\hat{\mathbb{Z}}_{p}^{r}\label{eq:Definition of bold t projection onto first Z_P hat level set}
\end{equation}
\end{notation}
\vphantom{}

Recall our earlier notation for writing $\mathbb{K}^{\rho,c}$ to
denote the $\mathbb{K}$-linear space of all $\rho\times c$ matrices
with entries in $\mathbb{K}$ for integers $\rho,c\geq1$. When $c=1$,
we identify $\mathbb{K}^{\rho,1}$ with $\mathbb{K}^{\rho}$.
\begin{defn}
\label{nota:Second batch}For a function $\mathbf{F}:\mathbb{Z}_{p}^{r}\rightarrow\mathbb{K}^{\rho,c}$,
we write \nomenclature{$\left\Vert \mathbf{F}\right\Vert _{p,\mathbb{K}}$}{ }
$\left\Vert \mathbf{F}\right\Vert _{p,\mathbb{K}}$ to denote: 
\begin{equation}
\left\Vert \mathbf{F}\right\Vert _{p,\mathbb{K}}\overset{\textrm{def}}{=}\sup_{\mathbf{z}\in\mathbb{Z}_{p}^{r}}\left\Vert \mathbf{F}\left(\mathbf{z}\right)\right\Vert _{\mathbb{K}}\label{eq:Definition of P,K norm}
\end{equation}
Here, $\left\Vert \mathbf{F}\left(\mathbf{z}\right)\right\Vert _{\mathbb{K}}$
is the maximum of the $\mathbb{K}$ absolute values of the entries
of the matrix $\mathbf{F}\left(\mathbf{z}\right)$. We write $\left\Vert \cdot\right\Vert _{p,q}$
in the case where $\mathbb{K}$ is a $q$-adic field, and write $\left\Vert \cdot\right\Vert _{p,\infty}$
when $\mathbb{K}$ is $\mathbb{C}$. We \emph{also} use $\left\Vert \cdot\right\Vert _{p,\mathbb{K}}$
to denote the corresponding norm of a matrix-valued function on $\hat{\mathbb{Z}}_{p}^{r}$. 
\end{defn}
\begin{defn}
Let\index{matrix!Hadamard product} $\mathbf{A}=\left\{ a_{j,k}\right\} _{1\leq j\leq\rho,1\leq k\leq c}$
and $\mathbf{B}=\left\{ b_{j,k}\right\} _{1\leq j\leq\rho,1\leq k\leq c}$
be elements of $\mathbb{K}^{\rho,c}$. Then, we write \nomenclature{$\mathbf{A}\odot\mathbf{B}$}{the Hadamard product of matrices}$\mathbf{A}\odot\mathbf{B}$
to denote the matrix: 
\begin{equation}
\mathbf{A}\odot\mathbf{B}\overset{\textrm{def}}{=}\left\{ a_{j,k}\cdot b_{j,k}\right\} _{1\leq j\leq\rho,1\leq k\leq c}\label{eq:Definition of the Hadamard product of two matrices}
\end{equation}
This is called the \textbf{Hadamard product }of $\mathbf{A}$ and
$\mathbf{B}$. Given $\mathbf{A}_{1},\ldots,\mathbf{A}_{N}$, where:
\begin{equation}
\mathbf{A}_{n}=\left\{ a_{j,k}\left(n\right)\right\} _{1\leq j\leq\rho,1\leq k\leq c}
\end{equation}
for each $n$, we then define write: 
\begin{equation}
\bigodot_{n=1}^{N}\mathbf{A}_{n}\overset{\textrm{def}}{=}\left\{ \prod_{n=1}^{N}a_{j,k}\left(n\right)\right\} _{1\leq j\leq\rho,1\leq k\leq c}\label{eq:Definition of the tensor product of n matrices}
\end{equation}
Note that $\odot$ reduces to standard multiplication when $\rho=c=1$. 
\end{defn}
\begin{defn}
We say a function $\mathbf{F}$ (resp. $\hat{\mathbf{F}}$) defined
on $\mathbb{Z}_{p}^{r}$ (resp. $\hat{\mathbb{Z}}_{p}^{r}$) taking
values in $\mathbb{K}^{\rho,c}$ is \textbf{elementary }if\index{elementary!function}
there are functions $\mathbf{F}_{1},\ldots,\mathbf{F}_{r}$ (resp.
$\hat{\mathbf{F}}_{1},\ldots,\hat{\mathbf{F}}_{r}$) defined on $\mathbb{Z}_{p}$
(resp. $\hat{\mathbb{Z}}_{p}$) so that: 
\begin{equation}
\mathbf{F}\left(\mathbf{z}\right)=\bigodot_{m=1}^{r}\mathbf{F}\left(\mathfrak{z}_{m}\right),\textrm{ }\forall\mathbf{z}\in\mathbb{Z}_{p}^{r}
\end{equation}
and: 
\begin{equation}
\hat{\mathbf{F}}\left(\mathbf{t}\right)=\bigodot_{m=1}^{r}\hat{\mathbf{F}}\left(t_{m}\right),\textrm{ }\forall\mathbf{t}\in\hat{\mathbb{Z}}_{p}^{r}
\end{equation}
respectively. We write \nomenclature{$E\left(\mathbb{Z}_{p}^{r},\mathbb{K}^{\left(\rho,c\right)}\right)$}{set of elementary tensors $\mathbb{Z}_{p}^{r}\rightarrow\mathbb{K}^{\rho,c}$}$E\left(\mathbb{Z}_{p}^{r},\mathbb{K}^{\rho,c}\right)$
(resp. \nomenclature{$E\left(\hat{\mathbb{Z}}_{p}^{r},\mathbb{K}^{\left(\rho,c\right)}\right)$}{set of elementary tensors $\hat{\mathbb{Z}}_{p}^{r}\rightarrow\mathbb{K}^{\rho,c}$}$E\left(\hat{\mathbb{Z}}_{p}^{r},\mathbb{K}^{\rho,c}\right)$)
to denote the $\mathbb{K}$-span of all elementary functions (a.k.a.,
\textbf{elementary tensors}\index{elementary!tensors}) on $\mathbb{Z}_{p}^{r}$
(resp. $\hat{\mathbb{Z}}_{p}^{r}$).

Note that the Hadamard products reduce to entry-wise products of vectors
when $c=1$.
\end{defn}
\begin{defn}[\textbf{Multi-Dimensional Function Spaces}]
\ 

\vphantom{}

I. We write \nomenclature{$B\left(\mathbb{Z}_{p}^{r},K\right)$}{set of bounded functions $\mathbb{Z}_{p}^{r}\rightarrow K$}$B\left(\mathbb{Z}_{p}^{r},K\right)$
to denote the $K$-linear space of all bounded functions $f:\mathbb{Z}_{p}^{r}\rightarrow K$.
This is a non-archimedean Banach space under the norm: 
\begin{equation}
\left\Vert f\right\Vert _{p,K}\overset{\textrm{def}}{=}\sup_{\mathbf{z}\in\mathbb{Z}_{p}^{r}}\left|f\left(\mathbf{z}\right)\right|_{K}=\sup_{\mathfrak{z}_{1},\ldots,\mathfrak{z}_{r}\in\mathbb{Z}_{p}}\left|f\left(\mathfrak{z}_{1},\ldots,\mathfrak{z}_{r}\right)\right|_{K}\label{eq:Definition of scalar valued C Z_P norm}
\end{equation}
We write \nomenclature{$C\left(\mathbb{Z}_{p}^{r},K\right)$}{set of continuous functions $\mathbb{Z}_{p}^{r}\rightarrow K$ }$C\left(\mathbb{Z}_{p}^{r},K\right)$
to denote the space of all continuous $f:\mathbb{Z}_{p}^{r}\rightarrow K$.

\vphantom{}

II. We write \nomenclature{$B\left(\mathbb{Z}_{p}^{r},K^{\left(\rho,c\right)}\right)$}{set of bounded functions $\mathbb{Z}_{p}^{r}\rightarrow K^{\rho,c}$  }$B\left(\mathbb{Z}_{p}^{r},K^{\rho,c}\right)$
to denote the $K$-linear space of all bounded functions $\mathbf{F}:\mathbb{Z}_{p}^{r}\rightarrow K^{\rho,c}$,
where: 
\begin{equation}
\mathbf{F}\left(\mathbf{z}\right)=\left\{ F_{j,k}\left(\mathbf{z}\right)\right\} _{1\leq j\leq\rho,1\leq k\leq c}\in K^{\rho,c},\textrm{ }\forall\mathbf{z}\in\mathbb{Z}_{p}^{r}
\end{equation}
where, for each $\left(j,k\right)$, $F_{j,k}\in B\left(\mathbb{Z}_{p}^{r},K\right)$.
$B\left(\mathbb{Z}_{p}^{r},K^{\rho,c}\right)$ becomes a non-archimedean
Banach space under the norm $\left\Vert \mathbf{F}\right\Vert _{p,K}$
(the maximum of the supremum-over-$\mathbb{Z}_{p}^{r}$ of the $K$-adic
absolute values of each entry of $\mathbf{F}$). We write \nomenclature{$C\left(\mathbb{Z}_{p}^{r},K^{\left(\rho,c\right)}\right)$}{set of continuous functions $\mathbb{Z}_{p}^{r}\rightarrow K^{\rho,c}$   }$C\left(\mathbb{Z}_{p}^{r},K^{\rho,c}\right)$
to denote the space of all continuous functions $\mathbf{F}:\mathbb{Z}_{p}^{r}\rightarrow K^{\rho,c}$.

\vphantom{}

III. We write \nomenclature{$B\left(\hat{\mathbb{Z}}_{p}^{r},K\right)$}{set of bounded functions $\hat{\mathbb{Z}}_{p}^{r}\rightarrow K$}$B\left(\hat{\mathbb{Z}}_{p}^{r},K\right)$
to denote the $K$-linear space of all bounded functions $\hat{f}:\hat{\mathbb{Z}}_{p}^{r}\rightarrow K$.
This is a non-archimedean Banach space under the norm: 
\begin{equation}
\left\Vert \hat{f}\right\Vert _{p,K}\overset{\textrm{def}}{=}\sup_{\mathbf{t}\in\hat{\mathbb{Z}}_{p}^{r}}\left|\hat{f}\left(\mathbf{t}\right)\right|_{K}=\sup_{t_{1},\ldots,t_{r}\in\hat{\mathbb{Z}}_{p}}\left|\hat{f}\left(\hat{t}_{1},\ldots,\hat{t}_{r}\right)\right|_{K}\label{eq:Definition of scalar valued C Z_P hat norm}
\end{equation}
We write $c_{0}\left(\hat{\mathbb{Z}}_{p}^{r},K\right)$\nomenclature{$c_{0}\left(\hat{\mathbb{Z}}_{p}^{r},K\right)$}{set of $f\in B\left(\hat{\mathbb{Z}}_{p}^{r},K\right)$ so that $\lim_{\left\Vert \mathbf{t}\right\Vert _{p}\rightarrow\infty}\left|\hat{f}\left(\mathbf{t}\right)\right|_{K}=0$}
to denote the subspace of $B\left(\hat{\mathbb{Z}}_{p}^{r},K\right)$
consisting of all $\hat{f}$ for which: 
\begin{equation}
\lim_{\left\Vert \mathbf{t}\right\Vert _{p}\rightarrow\infty}\left|\hat{f}\left(\mathbf{t}\right)\right|_{K}=0
\end{equation}

\vphantom{}

IV. We write \nomenclature{$B\left(\hat{\mathbb{Z}}_{p}^{r},K^{\left(\rho,c\right)}\right)$}{set of bounded functions $\hat{\mathbb{Z}}_{p}^{r}\rightarrow K^{\rho,c}$}$B\left(\hat{\mathbb{Z}}_{p}^{r},K^{\rho,c}\right)$
to denote the $K$-linear space of all bounded functions $\hat{\mathbf{F}}:\hat{\mathbb{Z}}_{p}^{r}\rightarrow K^{\rho,c}$,
where: 
\begin{equation}
\hat{\mathbf{F}}\left(\mathbf{t}\right)=\left\{ \hat{F}_{j,k}\left(\mathbf{t}\right)\right\} _{1\leq j\leq\rho,1\leq k\leq c}\in K^{\rho,c},\textrm{ }\forall\mathbf{t}\in\hat{\mathbb{Z}}_{p}^{r}
\end{equation}
where, for each $\left(j,k\right)$, $\hat{F}_{m}\in B\left(\hat{\mathbb{Z}}_{p}^{r},K^{\rho,c}\right)$.
$B\left(\hat{\mathbb{Z}}_{p}^{r},K^{\rho,c}\right)$ becomes a non-archimedean
Banach space under the norm $\left\Vert \hat{\mathbf{F}}\right\Vert _{p,K}^{r}$.
We then write \nomenclature{$c_{0}\left(\hat{\mathbb{Z}}_{p}^{r},K^{\left(\rho,c\right)}\right)$}{set of bounded functions $\hat{\mathbb{Z}}_{p}^{r}\rightarrow K^{\rho,c}$ so that $\lim_{\left\Vert \mathbf{t}\right\Vert _{p}^{r}\rightarrow\infty}\left\Vert \hat{\mathbf{F}}\left(\mathbf{t}\right)\right\Vert _{K}=0$}$c_{0}\left(\hat{\mathbb{Z}}_{p}^{r},K^{\rho,c}\right)$
to denote the Banach subspace of $B\left(\hat{\mathbb{Z}}_{p}^{r},K^{\rho,c}\right)$
consisting of all those $\hat{\mathbf{F}}$ so that: 
\begin{equation}
\lim_{\left\Vert \mathbf{t}\right\Vert _{p}^{r}\rightarrow\infty}\left\Vert \hat{\mathbf{F}}\left(\mathbf{t}\right)\right\Vert _{K}=0
\end{equation}
\end{defn}
\begin{rem}
The Fourier transform will send functions $\mathbb{Z}_{p}^{r}\rightarrow K$
to functions $\hat{\mathbb{Z}}_{p}^{r}\rightarrow K$. Fourier transforms
on $C\left(\mathbb{Z}_{p}^{r},K^{\rho,c}\right)$ are done component-wise.
\end{rem}
\vphantom{}

Now, we need to introduce the tensor product. Fortunately, in all
of the cases we shall be working with, the tensor product is nothing
more than point-wise multiplication of functions. A reader desirous
for greater abstract detail can refer to \cite{van Rooij - Non-Archmedean Functional Analysis}
to slake their thirst. The fourth chapter of van Rooij's book has
an \emph{entire section} dedicated to the tensor product\index{tensor!product}
($\otimes$\nomenclature{$\otimes$}{tensor product}). 
\begin{thm}
\cite{van Rooij - Non-Archmedean Functional Analysis} Let $E$, $F$,
and $G$ be Banach spaces over metrically complete non-archimedean
fields. Then:
\begin{equation}
\left\Vert f\otimes g\right\Vert =\left\Vert f\right\Vert _{E}\left\Vert g\right\Vert _{F},\textrm{ }\forall f\in E,\textrm{ }\forall g\in F\label{eq:Norm of a tensor product}
\end{equation}
Additionally, for any linear operator $T:E\otimes F\rightarrow G$,
the operator norm\index{operator norm} of $T$ is given by: 
\begin{equation}
\left\Vert T\right\Vert \overset{\textrm{def}}{=}\sup_{\left\Vert f\right\Vert _{E}\leq1,\left\Vert g\right\Vert _{F}\leq1}\frac{\left\Vert T\left\{ f\otimes g\right\} \right\Vert }{\left\Vert f\right\Vert _{E}\left\Vert g\right\Vert _{F}}\label{eq:Definition of the norm of a linear operator on a tensor product}
\end{equation}
\end{thm}
\begin{thm}
\cite{van Rooij - Non-Archmedean Functional Analysis} Let $E_{1},F_{1},E_{2},F_{2}$
be Banach spaces over metrically complete non-archimedean fields.
Given continuous linear operators $T_{j}:E_{j}\rightarrow F_{j}$
for $j\in\left\{ 1,2\right\} $, the operator $T_{1}\otimes T_{2}:E_{1}\otimes F_{1}\rightarrow E_{2}\otimes F_{2}$
is given by: 
\begin{equation}
\left(T_{1}\otimes T_{2}\right)\left\{ f\otimes g\right\} =T_{1}\left\{ f\right\} \otimes T_{2}\left\{ g\right\} \label{eq:Action of tensor products of linear operators}
\end{equation}
\end{thm}
\begin{cor}
\cite{van Rooij - Non-Archmedean Functional Analysis} Let $E$ and
$F$ be Banach spaces over a metrically complete non-archimedean field
$K$. Then, the dual of $E\otimes F$ is $E^{\prime}\otimes F^{\prime}$,
with elementary tensors $\mu\otimes\nu\in E^{\prime}\otimes F^{\prime}$
acting on elementary tensors $f\otimes g\in E\otimes F$ by: 
\begin{equation}
\left(\mu\otimes\nu\right)\left(f\otimes g\right)\overset{K}{=}\mu\left(f\right)\cdot\nu\left(g\right)\label{eq:action of tensor product of linear functionals}
\end{equation}
\end{cor}
\begin{thm}
\cite{van Rooij - Non-Archmedean Functional Analysis}\ 

\vphantom{}

I. 
\begin{equation}
B\left(\mathbb{Z}_{p}^{r},K\right)=\bigotimes_{m=1}^{r}B\left(\mathbb{Z}_{p},K\right)
\end{equation}
\index{tensor!elementary}Elementary tensors in $B\left(\mathbb{Z}_{p}^{r},K\right)$
are of the form: 
\begin{equation}
\left(\bigotimes_{m=1}^{r}f_{m}\right)\left(\mathbf{z}\right)=\left(\bigotimes_{m=1}^{r}f_{m}\right)\left(\mathfrak{z}_{1},\ldots,\mathfrak{z}_{r}\right)\overset{\textrm{def}}{=}\prod_{m=1}^{r}f_{m}\left(\mathfrak{z}_{m}\right)
\end{equation}
This also holds for $C\left(\mathbb{Z}_{p}^{r},K\right)$.

\vphantom{}

II. 
\begin{equation}
B\left(\hat{\mathbb{Z}}_{p}^{r},K\right)=\bigotimes_{m=1}^{r}B\left(\hat{\mathbb{Z}}_{p},K\right)
\end{equation}
Elementary tensors in $B\left(\hat{\mathbb{Z}}_{p}^{r},K\right)$
are of the form: 
\begin{equation}
\left(\bigotimes_{m=1}^{r}\hat{f}_{m}\right)\left(\mathbf{t}\right)=\left(\bigotimes_{m=1}^{r}\hat{f}_{m}\right)\left(t_{1},\ldots,t_{r}\right)\overset{\textrm{def}}{=}\prod_{m=1}^{r}\hat{f}_{m}\left(t_{m}\right)
\end{equation}
This also holds for $c_{0}\left(\hat{\mathbb{Z}}_{p}^{r},K\right)$.

\vphantom{}

III. 
\begin{equation}
B\left(\mathbb{Z}_{p}^{r},K^{\rho,c}\right)=\bigotimes_{m=1}^{r}B\left(\mathbb{Z}_{p},K^{\rho,c}\right)
\end{equation}
where an elementary tensor in $B\left(\mathbb{Z}_{p}^{r},K^{\rho,c}\right)$
is of the form: 
\begin{equation}
\mathbf{F}\left(\mathbf{z}\right)=\bigodot_{m=1}^{r}\mathbf{F}_{m}\left(\mathfrak{z}_{m}\right)
\end{equation}
where $\mathbf{F}_{m}$ is in $B\left(\mathbb{Z}_{p},K^{\rho,c}\right)$
for each $m$. This also holds for $C\left(\mathbb{Z}_{p}^{r},K^{\rho,c}\right)$.

\vphantom{}

IV.

\begin{equation}
B\left(\hat{\mathbb{Z}}_{p}^{r},K^{\rho,c}\right)=\bigotimes_{m=1}^{r}B\left(\hat{\mathbb{Z}}_{p},K^{\rho,c}\right)
\end{equation}
where an elementary tensor in $B\left(\hat{\mathbb{Z}}_{p}^{r},K^{\rho,c}\right)$
is of the form: 
\begin{equation}
\hat{\mathbf{F}}\left(\mathbf{t}\right)=\bigodot_{m=1}^{r}\hat{\mathbf{F}}_{m}\left(t_{m}\right)
\end{equation}
where $\hat{\mathbf{F}}_{m}$ is in $B\left(\hat{\mathbb{Z}}_{p},K^{\rho,c}\right)$
for each $m$. This also holds for $c_{0}\left(\hat{\mathbb{Z}}_{p}^{r},K^{\rho,c}\right)$. 
\end{thm}
\begin{rem}
Every function in $C\left(\mathbb{Z}_{p}^{r},K\right)$ can be written
as a linear combination of (possibly infinitely many) elementary functions.
If the linear combination is infinite, then elements of the linear
combination can be enumerated in a sequence which converges to $0$
in norm, so as to guarantee the convergence of the resulting infinite
sum. 
\end{rem}
\begin{defn}[\textbf{Multi-Dimensional van der Put Basis}]
Given a prime $p$, we write $\mathcal{B}_{p}$ to denote the set
of van der Put basis functions for $\mathbb{Z}_{p}$: 
\begin{equation}
\mathcal{B}_{p}\overset{\textrm{def}}{=}\left\{ \left[\mathfrak{z}\overset{p^{\lambda_{p}\left(n\right)}}{\equiv}n\right]:n\in\mathbb{N}_{0}\right\} \label{eq:Definition of the van der Put basis}
\end{equation}
Next, for any $\mathbf{n}\in\mathbb{N}_{0}^{r}$, we write: 
\begin{equation}
\left[\mathbf{z}\overset{p^{\lambda_{p}\left(\mathbf{n}\right)}}{\equiv}\mathbf{n}\right]\overset{\textrm{def}}{=}\prod_{\ell=1}^{r}\left[\mathfrak{z}_{\ell}\overset{p^{\lambda_{p}\left(n_{\ell}\right)}}{\equiv}n_{\ell}\right],\textrm{ }\forall\mathbf{z}\in\mathbb{Z}_{p}^{r}\label{eq:Notation for elementary P-adic van der Put tensors}
\end{equation}
We call (\ref{eq:Notation for elementary P-adic van der Put tensors})
an \textbf{elementary $p$-adic van der Put tensor}. \nomenclature{$\mathcal{B}_{p}^{\otimes r}$}{$\left\{ \left[\mathbf{z}\overset{p^{\lambda_{p}\left(\mathbf{n}\right)}}{\equiv}\mathbf{n}\right]:\mathbf{n}\in\mathbb{N}_{0}^{r}\right\}$}We
then define the \textbf{$p$-adic van der Put basis of depth $r$
}as the set: 
\begin{equation}
\mathcal{B}_{p}^{\otimes r}\overset{\textrm{def}}{=}\left\{ \left[\mathbf{z}\overset{p^{\lambda_{p}\left(\mathbf{n}\right)}}{\equiv}\mathbf{n}\right]:\mathbf{n}\in\mathbb{N}_{0}^{r}\right\} \label{eq:definition of elementary tensor van der Put basis}
\end{equation}
\end{defn}
\begin{thm}
$\mathcal{B}_{p}^{\otimes r}$ is a basis for $C\left(\mathbb{Z}_{p}^{r},K\right)$. 
\end{thm}
Proof: $\mathcal{B}_{p}$ is a basis for $C\left(\mathbb{Z}_{p},K\right)$,
so, tensoring the $\mathcal{B}_{p}$s to get $\mathcal{B}_{p}^{\otimes r}$
yields a basis for $C\left(\mathbb{Z}_{p}^{r},K\right)$.

Q.E.D. 
\begin{defn}[\textbf{Multi-Dimensional van der Put Series}]
We write \nomenclature{$\textrm{vdP}\left(\mathbb{Z}_{p}^{r},K\right)$}{the set of formal $\left(p,K\right)$-adic van der Put series of depth $r$ }
to denote the $K$-linear space of all\index{van der Put!series!formal}
\textbf{formal $\left(p,K\right)$-adic van der Put series of depth
$r$}. The elements of $\textrm{vdP}\left(\mathbb{Z}_{p}^{r},K\right)$
are formal sums: 
\begin{equation}
\sum_{\mathbf{n}\in\mathbb{N}_{0}^{r}}\mathfrak{a}_{\mathbf{n}}\left[\mathbf{z}\overset{p^{\lambda_{p}\left(\mathbf{n}\right)}}{\equiv}\mathbf{n}\right]\label{eq:Definition of a formal (P,K)-adic van der Put series}
\end{equation}
where the $\mathfrak{a}_{\mathbf{n}}$s are elements of $K$. The
\textbf{domain of convergence }of a formal van der Put series is the
set of all $\mathbf{z}\in\mathbb{Z}_{p}^{r}$ for which the series
converges\footnote{Note: since $\mathbf{z}\overset{p^{\lambda_{p}\left(\mathbf{n}\right)}}{\equiv}\mathbf{n}$
can hold true for at most finitely many $\mathbf{n}$ whenever $\mathbf{z}\in\mathbb{N}_{0}^{r}$,
observe that every formal van der Put series necessarily converges
in $K$ at every $\mathbf{z}\in\mathbb{N}_{0}^{r}$.} in $K$. We call $\textrm{vdP}\left(\mathbb{Z}_{p}^{r},\mathbb{C}_{q}\right)$
the space of \textbf{formal $\left(p,q\right)$-adic van der Put series
of depth $r$}.
\end{defn}
\begin{defn}
Recall that for an integer $n\geq0$ and a prime $p$, we write: 
\begin{equation}
n_{-}=\begin{cases}
0 & \textrm{if }n=0\\
n-n_{\lambda_{p}\left(n\right)-1}p^{\lambda_{p}\left(n\right)-1} & \textrm{if }n\geq1
\end{cases}
\end{equation}
where $n_{\lambda_{p}\left(n\right)-1}$ is the coefficient of the
largest power of $p$ present in the $p$-adic expansion of $n$.

\vphantom{}

I. Given $\mathbf{n}\in\mathbb{N}_{0}^{d}$, we write \nomenclature{$\mathbf{n}_{-}$}{multi-dimensional analogue of $n_{-}$}$\mathbf{n}_{-}$
to denote the $d$-tuple obtained by applying the subscript $-$ operation
to each entry of $\mathbf{n}$: $\mathbf{n}_{-}\overset{\textrm{def}}{=}\left(\left(n_{1}\right)_{-},\ldots,\left(n_{d}\right)_{-}\right)$.

\vphantom{}

II. Let $f:\mathbb{Z}_{p}^{r}\rightarrow K$ be an elementary function,
with $f=\prod_{m=1}^{r}f_{m}$. Then, we define the \textbf{van der
Put coefficients}\index{van der Put!coefficients}\textbf{ }\nomenclature{$c_{\mathbf{n}}\left(f\right)$}{$\mathbf{n}$th van der Put coefficient of $f$}$c_{\mathbf{n}}\left(f\right)$
of $f$ by: 
\begin{equation}
c_{\mathbf{n}}\left(f\right)\overset{\textrm{def}}{=}\prod_{m=1}^{r}c_{n_{m}}\left(f_{m}\right),\textrm{ }\forall\mathbf{n}\in\mathbb{N}_{0}^{r}\label{eq:Def of the van der put coefficients of an elementary function}
\end{equation}
\end{defn}
\begin{prop}
\label{prop:vdP series of an elementary function}Every elementary\index{elementary!function}
$f\in C\left(\mathbb{Z}_{p}^{r},K\right)$ is uniquely representable
as a $\left(p,K\right)$-adic van der Put series of depth $r$: 
\begin{equation}
f\left(\mathbf{z}\right)\overset{K}{=}\sum_{\mathbf{n}\in\mathbb{N}_{0}^{r}}c_{\mathbf{n}}\left(f\right)\left[\mathbf{z}\overset{p^{\lambda_{p}\left(\mathbf{n}\right)}}{\equiv}\mathbf{n}\right],\textrm{ }\mathbf{z}\in\mathbb{Z}_{p}^{r}
\end{equation}
Moreover, the $K$-convergence of this series is uniform over $\mathbb{Z}_{p}^{r}$. 
\end{prop}
Proof: Let $f$ be elementary. Then, $f$ is a product of the form:
\begin{equation}
f\left(\mathbf{z}\right)=f\left(\mathfrak{z}_{1},\ldots,\mathfrak{z}_{r}\right)=\prod_{m=1}^{r}f_{m}\left(\mathfrak{z}_{m}\right)
\end{equation}
for $f_{m}\in C\left(\mathbb{Z}_{p},K\right)$. As such, each $f_{m}$
is uniquely representable by a uniformly convergent van der Put series:
\begin{equation}
f_{m}\left(\mathfrak{z}_{m}\right)=\sum_{n_{m}=0}^{\infty}c_{n_{m}}\left(f_{m}\right)\left[\mathfrak{z}_{m}\overset{p^{\lambda_{p}\left(n_{m}\right)}}{\equiv}n_{m}\right]
\end{equation}
with $\left|c_{n_{m}}\left(f_{m}\right)\right|_{K}\rightarrow0$ as
$n_{m}\rightarrow\infty$. So: 
\begin{align*}
f\left(\mathfrak{z}_{1},\ldots,\mathfrak{z}_{r}\right) & =\prod_{m=1}^{r}\left(\sum_{n_{m}=0}^{\infty}c_{n_{m}}\left(f_{m}\right)\left[\mathfrak{z}_{m}\overset{p^{\lambda_{p}\left(n_{m}\right)}}{\equiv}n_{m}\right]\right)\\
 & =\sum_{n_{1}=0}^{\infty}\cdots\sum_{n_{r}=0}^{\infty}\left(\prod_{m=1}^{r}c_{n_{m}}\left(f_{m}\right)\left[\mathfrak{z}_{m}\overset{p^{\lambda_{p}\left(n_{m}\right)}}{\equiv}n_{m}\right]\right)\\
 & =\sum_{\mathbf{n}\in\mathbb{N}_{0}^{r}}c_{\mathbf{n}}\left(f\right)\left[\mathbf{z}\overset{p^{\lambda_{p}\left(\mathbf{n}\right)}}{\equiv}\mathbf{n}\right]
\end{align*}

Q.E.D. 
\begin{prop}
Let $g_{1},\ldots g_{N}$ be elementary functions in $C\left(\mathbb{Z}_{p}^{r},K\right)$.
Then, for any $\alpha_{1},\ldots,\alpha_{N}\in K$, the linear combination
$\sum_{k=1}^{N}\alpha_{k}g_{k}$ is represented by a uniformly convergent
$\left(p,K\right)$-adic van der Put series of depth $r$, the coefficients
of which are given by: 
\begin{equation}
c_{\mathbf{n}}\left(\sum_{k=1}^{N}\alpha_{k}g_{k}\right)=\sum_{k=1}^{N}\alpha_{k}c_{\mathbf{n}}\left(g_{k}\right)\label{eq:Linear extension of MD van der Put coefficients}
\end{equation}
where $c_{\mathbf{n}}\left(g_{k}\right)$ is computed by \emph{(\ref{eq:Def of the van der put coefficients of an elementary function})}.
This result also extends to infinite linear combinations, provided
that they converge in norm. 
\end{prop}
Proof: Each $g_{k}$ is representable by a uniformly convergent $\left(p,K\right)$-adic
van der Put series of depth $r$. The proposition follows from the
fact that the space of van der Put series is closed under linear combinations.

Q.E.D. 
\begin{rem}
Consequently, we can use (\ref{eq:Linear extension of MD van der Put coefficients})
to compute the van der Put coefficients of \emph{any }function in
$C\left(\mathbb{Z}_{p}^{r},K\right)$, seeing as the $K$-span of
$C\left(\mathbb{Z}_{p}^{r},K\right)$'s elementary functions is dense
in $C\left(\mathbb{Z}_{p}^{r},K\right)$. 
\end{rem}
\begin{lem}
\label{lem:van der Put series for an arbitrary K-valued function on Z_P}Let
$f\in C\left(\mathbb{Z}_{p}^{r},K\right)$. Then, there are unique
constants $\left\{ c_{\mathbf{n}}\left(f\right)\right\} _{\mathbf{n}\in\mathbb{N}_{0}^{r}}$
so that: 
\begin{equation}
f\left(\mathbf{z}\right)\overset{K}{=}\sum_{\mathbf{n}\in\mathbb{N}_{0}^{r}}c_{\mathbf{n}}\left(f\right)\left[\mathbf{z}\overset{p^{\lambda_{p}\left(\mathbf{n}\right)}}{\equiv}\mathbf{n}\right]\label{eq:van der Put series for an arbitrary K-valued function on Z_P}
\end{equation}
where the convergence is uniform in $\mathbf{z}$; specifically: 
\begin{equation}
\lim_{\left\Vert \mathbf{n}\right\Vert _{\infty}\rightarrow\infty}\left|c_{\mathbf{n}}\left(f\right)\right|_{K}=0
\end{equation}
where, recall (see \emph{(\ref{eq:Definition of infinity norm})}
on page \emph{\pageref{eq:Definition of infinity norm}}):
\begin{equation}
\left\Vert \mathbf{n}\right\Vert _{\infty}=\max\left\{ n_{1},\ldots,n_{r}\right\} 
\end{equation}

In particular, letting $\left\{ g_{k}\right\} _{k\geq1}$ and $\left\{ \alpha_{k}\right\} _{k\geq1}$
be any sequences of \index{elementary!function}elementary functions
in $C\left(\mathbb{Z}_{p}^{r},K\right)$ and scalars in $K$, respectively,
so that: 
\begin{equation}
\lim_{N\rightarrow\infty}\left\Vert f-\sum_{k=1}^{N}\alpha_{k}g_{k}\right\Vert _{p,q}=0
\end{equation}
then: 
\begin{equation}
c_{\mathbf{n}}\left(f\right)\overset{K}{=}\lim_{N\rightarrow\infty}c_{\mathbf{n}}\left(\sum_{k=1}^{N}\alpha_{k}g_{k}\right)=\lim_{N\rightarrow\infty}\sum_{k=1}^{N}\alpha_{k}c_{\mathbf{n}}\left(g_{k}\right)\label{eq:Limit formula for the van der Put coefficients of an abitrary continuous (P,K)-adic function}
\end{equation}
We call the $c_{\mathbf{n}}\left(f\right)$s the \index{van der Put!coefficients}\textbf{van
der Put coefficients }of $f$. 
\end{lem}
Proof:

I. (Existence) Since the $C\left(\mathbb{Z}_{p}^{r},K\right)$-elementary
functions are dense in $C\left(\mathbb{Z}_{p}^{r},K\right)$, given
any $f\in C\left(\mathbb{Z}_{p}^{r},K\right)$, we can choose $\left\{ g_{k}\right\} _{k\geq1}$
and $\left\{ \alpha_{k}\right\} _{k\geq1}$ as described above so
that the norm convergence occurs. Writing: 
\begin{equation}
\sum_{k=1}^{N}\alpha_{k}g_{k}\left(\mathbf{z}\right)\overset{K}{=}\sum_{\mathbf{n}\in\mathbb{N}_{0}^{r}}\left(\sum_{k=1}^{N}\alpha_{k}c_{\mathbf{n}}\left(g_{k}\right)\right)\left[\mathbf{z}\overset{p^{\lambda_{p}\left(\mathbf{n}\right)}}{\equiv}\mathbf{n}\right]
\end{equation}
\textemdash note, the convergence is uniform in $\mathbf{z}$\textemdash the
assumed norm convergence allows us to interchange limits and sums:
\begin{equation}
f\left(\mathbf{z}\right)\overset{K}{=}\lim_{N\rightarrow\infty}\sum_{k=1}^{N}\alpha_{k}g_{k}\left(\mathbf{z}\right)\overset{K}{=}\sum_{\mathbf{n}\in\mathbb{N}_{0}^{r}}\underbrace{\left(\lim_{N\rightarrow\infty}\sum_{k=1}^{N}\alpha_{k}c_{\mathbf{n}}\left(g_{k}\right)\right)}_{c_{\mathbf{n}}\left(f\right)}\left[\mathbf{z}\overset{p^{\lambda_{p}\left(\mathbf{n}\right)}}{\equiv}\mathbf{n}\right]
\end{equation}

\vphantom{}

II. (Uniqueness) Let $\left\{ g_{k}\right\} _{k\geq1}$ and $\left\{ \alpha_{k}\right\} _{k\geq1}$
be one choice of $g$s and $\alpha$s whose linear combinations converge
in norm to $f$; let $\left\{ g_{k}^{\prime}\right\} _{k\geq1}$ and
$\left\{ \alpha_{k}^{\prime}\right\} _{k\geq1}$ be another choice
of $g$s and $\alpha$s satisfying the same. Then: 
\begin{align*}
0 & \overset{K}{=}f\left(\mathbf{z}\right)-f\left(\mathbf{z}\right)\\
 & \overset{K}{=}\lim_{N\rightarrow\infty}\left(\sum_{k=1}^{N}\left(\alpha_{k}g_{k}\left(\mathbf{z}\right)-\alpha_{k}^{\prime}g_{k}^{\prime}\left(\mathbf{z}\right)\right)\right)\\
 & \overset{K}{=}\sum_{\mathbf{n}\in\mathbb{N}_{0}^{r}}\lim_{N\rightarrow\infty}\left(\sum_{k=1}^{N}\left(\alpha_{k}c_{\mathbf{n}}\left(g_{k}\right)-\alpha_{k}^{\prime}c_{\mathbf{n}}\left(g_{k}^{\prime}\right)\right)\right)\left[\mathbf{z}\overset{p^{\lambda_{p}\left(\mathbf{n}\right)}}{\equiv}\mathbf{n}\right]
\end{align*}
Since $0$ is an elementary function, by \textbf{Proposition \ref{prop:vdP series of an elementary function},}
it is uniquely represented by the van der Put series whose coefficients
are all $0$. This forces: 
\begin{equation}
\lim_{N\rightarrow\infty}\sum_{k=1}^{N}\alpha_{k}c_{\mathbf{n}}\left(g_{k}\right)\overset{K}{=}\lim_{N\rightarrow\infty}\sum_{k=1}^{N}\alpha_{k}^{\prime}c_{\mathbf{n}}\left(g_{k}^{\prime}\right)
\end{equation}
which shows that $c_{\mathbf{n}}\left(f\right)$ is, indeed, unique.

Q.E.D. 
\begin{defn}
\label{def:MD S_p}We write $F\left(\mathbb{Z}_{p}^{r},K\right)$\nomenclature{$F\left(\mathbb{Z}_{p}^{r},K\right)$}{$K$-valued functions on $\mathbb{Z}_{p}^{r}$},
to denote the space of all $K$-valued functions on $\mathbb{Z}_{p}^{r}$.
We then define the linear operator $S_{p}:F\left(\mathbb{Z}_{p}^{r},K\right)\rightarrow\textrm{vdP}\left(\mathbb{Z}_{p}^{r},K\right)$
by: 
\begin{equation}
S_{p}\left\{ f\right\} \left(\mathbf{z}\right)\overset{\textrm{def}}{=}\sum_{\mathbf{n}\in\mathbb{N}_{0}^{r}}c_{\mathbf{n}}\left(f\right)\left[\mathbf{z}\overset{p^{\lambda_{p}\left(\mathbf{n}\right)}}{\equiv}\mathbf{n}\right],\textrm{ }\forall f\in F\left(\mathbb{Z}_{p}^{r},K\right)\label{eq:MD Definition of S_p of f}
\end{equation}
We also define partial sum operators: $S_{p:N}:F\left(\mathbb{Z}_{p}^{r},K\right)\rightarrow C\left(\mathbb{Z}_{p}^{r},K\right)$
by: 
\begin{equation}
S_{p:N}\left\{ f\right\} \left(\mathbf{z}\right)\overset{\textrm{def}}{=}\sum_{\mathbf{n}=\mathbf{0}}^{p^{N}-1}c_{\mathbf{n}}\left(f\right)\left[\mathbf{z}\overset{p^{\lambda_{p}\left(\mathbf{n}\right)}}{\equiv}\mathbf{n}\right],\textrm{ }\forall f\in F\left(\mathbb{Z}_{p}^{r},K\right)\label{eq:MD Definition of S_p N of f}
\end{equation}
Recall, the sum here is taken over all $\mathbf{n}\in\mathbb{N}_{0}^{r}$
so that, for each $m\in\left\{ 1,\ldots,r\right\} $, $0\leq n_{m}\leq p^{N}-1$. 
\end{defn}
\begin{prop}
We can extend $S_{p}$ from $F\left(\mathbb{Z}_{p}^{r},K\right)$
to $B\left(\mathbb{Z}_{p}^{r},K\right)$ by taking limits in $F\left(\mathbb{Z}_{p}^{r},K\right)$
with respect to $B\left(\mathbb{Z}_{p}^{r},K\right)$'s norm. We can
then compute $c_{\mathbf{n}}\left(f\right)$ for any $f\in B\left(\mathbb{Z}_{p}^{r},K\right)$
by using the construction given in \textbf{\emph{Lemma \ref{lem:van der Put series for an arbitrary K-valued function on Z_P}}},
albeit with elementary functions in $B\left(\mathbb{Z}_{p}^{r},K\right)$,
rather than in $C\left(\mathbb{Z}_{p}^{r},K\right)$. 
\end{prop}
Proof: Essentially the same as the proof of \textbf{Lemma \ref{lem:van der Put series for an arbitrary K-valued function on Z_P}}.

Q.E.D. 
\begin{lem}[\textbf{Multi-Dimensional van der Put Identities}]
\label{lem:MD vdP identities}\ 

\vphantom{}

I. For any $f:\mathbb{Z}_{p}^{r}\rightarrow K$ and any $N\in\mathbb{N}_{0}$:

\begin{equation}
S_{p:N}\left\{ f\right\} \left(\mathbf{z}\right)=\sum_{\mathbf{n}=\mathbf{0}}^{p^{N}-1}c_{\mathbf{n}}\left(f\right)\left[\mathbf{z}\overset{p^{\lambda_{p}\left(\mathbf{n}\right)}}{\equiv}\mathbf{n}\right]\overset{K}{=}f\left(\left[\mathbf{z}\right]_{p^{N}}\right),\textrm{ }\forall\mathbf{z}\in\mathbb{Z}_{p}^{r}\label{eq:MD truncated vdP identity}
\end{equation}

\vphantom{}

II. For any $f\in B\left(\mathbb{Z}_{p}^{r},K\right)$: 
\begin{equation}
S_{p}\left\{ f\right\} \left(\mathbf{z}\right)=\sum_{\mathbf{n}\in\mathbb{N}_{0}^{r}}c_{\mathbf{n}}\left(f\right)\left[\mathbf{z}\overset{p^{\lambda_{p}\left(\mathbf{n}\right)}}{\equiv}\mathbf{n}\right]\overset{K}{=}\lim_{k\rightarrow\infty}f\left(\left[\mathbf{z}\right]_{p^{k}}\right),\textrm{ }\forall\mathbf{z}\in\mathbb{Z}_{p}^{r}\label{eq:MD vdP identity}
\end{equation}
\end{lem}
Proof:

I. Let $f:\mathbb{Z}_{p}^{r}\rightarrow K$ be elementary, with $f=\prod_{m=1}^{r}f_{m}$.
Then, by the one-dimensional truncated van der Put identity (\ref{eq:truncated van der Put identity}):
\begin{align*}
f\left(\left[\mathbf{z}\right]_{p^{N}}\right) & =\prod_{m=1}^{r}f_{m}\left(\left[\mathfrak{z}_{m}\right]_{p^{N}}\right)\\
 & =\prod_{m=1}^{r}S_{p:N}\left\{ f_{m}\right\} \left(\mathfrak{z}_{m}\right)\\
\left(\textrm{Computation from \textbf{Prop. \ref{prop:vdP series of an elementary function}})}\right); & =S_{p:N}\left\{ \prod_{m=1}^{r}f_{m}\right\} \left(\mathfrak{z}_{1},\ldots,\mathfrak{z}_{r}\right)\\
 & =S_{p:N}\left\{ f\right\} \left(\mathbf{z}\right)
\end{align*}
Using the linearity of $S_{p:N}$ then gives the desired result.

\vphantom{}

II. Take limits with (I).

Q.E.D. 
\begin{thm}[\textbf{Multi-Dimensional van der Put Basis Theorem}]
\label{thm:MD vdP basis theorem}$\mathcal{B}_{p}^{\otimes r}$ is
a basis for $C\left(\mathbb{Z}_{p}^{r},K\right)$. In particular,
for any $f:\mathbb{Z}_{p}^{r}\rightarrow K$, the following are equivalent:

\vphantom{}

I. $f$ is continuous.

\vphantom{}

II. $\lim_{\left\Vert \mathbf{n}\right\Vert _{\infty}\rightarrow\infty}\left|c_{\mathbf{n}}\left(f\right)\right|_{K}=0$.

\vphantom{}

III. $\lim_{N\rightarrow\infty}\left\Vert f-S_{p:N}\left\{ f\right\} \right\Vert _{p,K}=0$ 
\end{thm}
Proof: Formally equivalent to its one-dimensional analogue\textemdash \textbf{Theorem
\ref{thm:vdP basis theorem}}.

Q.E.D. 
\begin{cor}
If $K$ is non-archimedean, then the Banach space $C\left(\mathbb{Z}_{p}^{r},K\right)$
is isometrically isomorphic to $c_{0}\left(\mathbb{N}_{0}^{r},K\right)$
(the space of sequences $\left\{ c_{\mathbf{n}}\right\} _{\mathbf{n}\in\mathbb{N}_{0}^{r}}$
in $K$ that converge to $0$ in $K$ as $\left\Vert \mathbf{n}\right\Vert _{\infty}\rightarrow\infty$). 
\end{cor}
Proof: Use the tensor product along with the isometric isomorphism
$C\left(\mathbb{Z}_{p},K\right)\cong c_{0}\left(K\right)$ proven
in \textbf{Theorem \ref{thm:C(Z_p,K) is iso to c_0 K}} (page \pageref{thm:C(Z_p,K) is iso to c_0 K}).

Q.E.D. 
\begin{rem}
Before we move on to the next subsection, note that everything we
have done so far extends in the obvious way to vector\emph{-} or matrix-valued\emph{
}functions on $\mathbb{Z}_{p}^{r}$ by working component-wise. To
that end, for $\mathbf{F}:\mathbb{Z}_{p}^{r}\rightarrow K^{\rho,c}$
given by $\mathbf{F}\left(\mathbf{z}\right)=\left\{ F_{j,k}\right\} _{j,k}$,
\nomenclature{$c_{\mathbf{n}}\left(\mathbf{F}\right)$}{$\mathbf{n}$th van der Put coefficient of $\mathbf{F}$}
we write: 
\begin{equation}
c_{\mathbf{n}}\left(\mathbf{F}\right)\overset{\textrm{def}}{=}\left\{ c_{\mathbf{n}}\left(F_{j,k}\right)\right\} _{j,k}
\end{equation}
\begin{equation}
S_{p:N}\left\{ \mathbf{F}\right\} \left(\mathbf{z}\right)\overset{\textrm{def}}{=}\left\{ S_{p:N}\left\{ F_{j,k}\right\} \left(\mathbf{z}\right)\right\} _{j,k}
\end{equation}
\end{rem}

\subsection{\label{subsec:5.3.2. Interpolation-Revisited}Interpolation Revisited}

We begin with rising sequences and rising-continuous functions. 
\begin{defn}
A sequence $\left\{ \mathbf{z}_{n}\right\} _{n\geq0}$ in $\mathbb{Z}_{p}^{r}$
(where $\mathbf{z}_{n}=\left(\mathfrak{z}_{n,1},\ldots,\mathfrak{z}_{n,r}\right)$)
is said to be \textbf{($p$-adically) rising }if there exists an $\ell\in\left\{ 1,\ldots,r\right\} $
so that the number of non-zero $p$-adic digits in $\mathfrak{z}_{n,\ell}$
tends to $\infty$ as $n\rightarrow\infty$. 
\end{defn}
\begin{defn}
\label{def:MD rising-continuity}A function $\chi:\mathbb{Z}_{p}^{r}\rightarrow K^{d}$
is said to be \textbf{($\left(p,K\right)$-adically)} \textbf{rising-continuous}\index{rising-continuous!left(p,Kright)-adically@$\left(p,K\right)$-adically}\textbf{
}\index{rising-continuous!function}whenever: 
\begin{equation}
\lim_{n\rightarrow\infty}\chi\left(\left[\mathbf{z}\right]_{p^{n}}\right)\overset{K^{d}}{=}\chi\left(\mathbf{z}\right),\textrm{ }\forall\mathbf{z}\in\mathbb{Z}_{p}^{r}\label{eq:MD Definition of a rising-continuous function}
\end{equation}
where the convergence is point-wise. We write \nomenclature{$\tilde{C}\left(\mathbb{Z}_{p}^{r},K^{d}\right)$}{set of $K^{d}$-valued rising-continuous functions on $\mathbb{Z}_{p}^{r}$ }$\tilde{C}\left(\mathbb{Z}_{p}^{r},K^{d}\right)$
to denote the $K$-linear space of all rising-continuous functions
$\mathbb{Z}_{p}^{r}\rightarrow K^{d}$.
\end{defn}
\begin{prop}
\label{prop:MD vdP criterion for rising continuity}Let $\chi\in B\left(\mathbb{Z}_{p}^{r},\mathbb{C}_{q}^{d}\right)$
be any function. Then, $S_{p}\left\{ \chi\right\} $ (the van der
Put series\index{van der Put!series} of $\chi$) converges at $\mathbf{z}\in\mathbb{Z}_{p}^{r}$
if and only if: 
\begin{equation}
\lim_{k\rightarrow\infty}c_{\left[\mathbf{z}\right]_{p^{k}}}\left(\chi\right)\left[\lambda_{p}\left(\left[\mathbf{z}\right]_{p^{k}}\right)=k\right]\overset{\mathbb{C}_{q}^{d}}{=}\mathbf{0}\label{eq:MD vdP criterion for rising-continuity}
\end{equation}
where $c_{\left[\mathbf{z}\right]_{p^{k}}}\left(\chi\right)$ is the
$\left[\mathbf{z}\right]_{p^{k}}=\left(\left[\mathfrak{z}_{1}\right]_{p^{k}},\ldots,\left[\mathfrak{z}_{r}\right]_{p^{k}}\right)$
th van der Put coefficient of $\chi$. 
\end{prop}
Proof: We start by writing: 
\begin{align*}
S_{p}\left\{ \chi\right\} \left(\mathbf{z}\right) & \overset{\mathbb{C}_{q}^{d}}{=}\sum_{\mathbf{n}\in\mathbb{N}_{0}^{r}}c_{\mathbf{n}}\left(\chi\right)\left[\mathbf{z}\overset{p^{\lambda_{p}\left(\mathbf{n}\right)}}{\equiv}\mathbf{n}\right]\\
 & =c_{\mathbf{0}}\left(\chi\right)+\sum_{k=1}^{\infty}\sum_{\mathbf{n}:\lambda_{p}\left(\mathbf{n}\right)=k}c_{\mathbf{n}}\left(\chi\right)\left[\mathbf{z}\overset{p^{k}}{\equiv}\mathbf{n}\right]
\end{align*}
For each $\mathbf{n}$, observe that if $\lambda_{p}\left(\left[\mathbf{z}\right]_{p^{k}}\right)=k$,
then: 
\[
\mathbf{n}\overset{\textrm{def}}{=}\left[\mathbf{z}\right]_{p^{k}}=\left(\left[\mathfrak{z}_{1}\right]_{p^{k}},\ldots,\left[\mathfrak{z}_{r}\right]_{p^{k}}\right)
\]
 is the unique element of $\mathbb{N}_{0}^{r}$ satisfying both $\lambda_{p}\left(\mathbf{n}\right)=k$
and $\mathbf{z}\overset{p^{k}}{\equiv}\mathbf{n}$. Hence: 
\begin{align}
S_{p}\left\{ \chi\right\} \left(\mathbf{z}\right) & \overset{\mathbb{C}_{q}^{d}}{=}c_{\mathbf{0}}\left(\chi\right)+\sum_{k=1}^{\infty}c_{\left[\mathbf{z}\right]_{p^{k}}}\left(\chi\right)\left[\lambda_{p}\left(\left[\mathbf{z}\right]_{p^{k}}\right)=k\right]
\end{align}
The ultrametric properties of $\mathbb{C}_{q}^{d}$ tell us that the
$q$-adic convergence of this series at any given $\mathbf{z}\in\mathbb{Z}_{p}^{r}$
is equivalent to:

\begin{equation}
\lim_{k\rightarrow\infty}c_{\left[\mathbf{z}\right]_{p^{k}}}\left(\chi\right)\left[\lambda_{p}\left(\left[\mathbf{z}\right]_{p^{k}}\right)=k\right]\overset{\mathbb{C}_{q}^{d}}{=}\mathbf{0}
\end{equation}

Q.E.D. 
\begin{thm}
The operator $S_{p}$ which sends a function to its formal van der
Put series is a isomorphism of $K$-linear spaces. This isomorphism
map $\tilde{C}\left(\mathbb{Z}_{p}^{r},K^{d}\right)$ onto the subspace
of $\textrm{vdP}\left(\mathbb{Z}_{p}^{r},K^{d}\right)$ consisting
of all van der Put series which converge at every $\mathbf{z}\in\mathbb{Z}_{p}^{r}$.

Additionally, for every $\chi\in\tilde{C}\left(\mathbb{Z}_{p}^{r},K^{d}\right)$:

\vphantom{}

I. $\chi=S_{p}\left\{ \chi\right\} $;

\vphantom{}

II. $\chi$ is uniquely represented by its van der Put series: 
\begin{equation}
\chi\left(\mathbf{z}\right)\overset{K^{d}}{=}\sum_{\mathbf{n}\in\mathbb{N}_{0}^{r}}c_{\mathbf{n}}\left(\chi\right)\left[\mathbf{z}\overset{p^{\lambda_{p}\left(\mathbf{n}\right)}}{\equiv}\mathbf{n}\right],\textrm{ }\forall\mathbf{z}\in\mathbb{Z}_{p}^{r}\label{eq:MD Chi vdP series}
\end{equation}
where the convergence is point-wise. 
\end{thm}
Proof: Let $\chi\in\tilde{C}\left(\mathbb{Z}_{p}^{r},K^{d}\right)$
be arbitrary. By the truncated van der Put identity (\ref{eq:MD truncated vdP identity}),
$\chi\left(\left[\mathbf{z}\right]_{p^{N}}\right)=S_{p:N}\left\{ \chi\right\} \left(\mathbf{z}\right)$.
Here, the rising-continuity of $\chi$ guarantees the point-wise convergence
of $S_{p:N}\left\{ \chi\right\} \left(\mathbf{z}\right)$ in $K^{d}$
to $\chi\left(\mathbf{z}\right)$ as $N\rightarrow\infty$. By \textbf{Proposition
\ref{prop:MD vdP criterion for rising continuity}}, this implies
the van der Put coefficients of $\chi$ satisfy (\ref{eq:MD vdP criterion for rising-continuity})
for all $\mathbf{z}\in\mathbb{Z}_{p}^{r}$. Consequently, the van
der Put series $S_{p}\left\{ \chi\right\} \left(\mathbf{z}\right)$
converges at every $\mathbf{z}\in\mathbb{Z}_{p}^{r}$, where it is
equal to $\chi\left(\mathbf{z}\right)$. This proves (I).

As for (II), the uniqueness specified therein is equivalent to demonstrating
that $S_{p}$ is an isomorphism in the manner described above. We
do this below:
\begin{itemize}
\item (Surjectivity) Let $V\left(\mathbf{z}\right)\in\textrm{vdP}\left(\mathbb{Z}_{p}^{r},K^{d}\right)$
be any formal van der Put series which converges $q$-adically at
every $\mathbf{z}\in\mathbb{Z}_{p}^{r}$. Letting: 
\begin{equation}
\chi\left(\mathbf{z}\right)\overset{\textrm{def}}{=}\lim_{N\rightarrow\infty}S_{p:N}\left\{ V\right\} \left(\mathbf{z}\right)
\end{equation}
we have $V\left(\mathbf{m}\right)=\chi\left(\mathbf{m}\right)$ for
all $\mathbf{m}\in\mathbb{N}_{0}^{r}$, and hence, $V\left(\mathbf{z}\right)=S_{p}\left\{ \chi\right\} $.
Thus, $S_{p:N}\left\{ \chi\right\} \left(\mathbf{z}\right)=\chi\left(\left[\mathbf{z}\right]_{p^{N}}\right)$.
Since $\chi\left(\mathbf{z}\right)$ is defined by $\lim_{N\rightarrow\infty}S_{p:N}\left\{ V\right\} \left(\mathbf{z}\right)$,
this gives: 
\begin{equation}
\chi\left(\mathbf{z}\right)=\lim_{N\rightarrow\infty}S_{p:N}\left\{ V\right\} \left(\mathbf{z}\right)=\chi\left(\left[\mathbf{z}\right]_{p^{N}}\right)
\end{equation}
which establishes the rising-continuity of $\chi$. This proves $V=S_{p}\left\{ \chi\right\} $,
and thus, that $S_{p}$ is surjective. 
\item (Injectivity) Let $\chi_{1},\chi_{2}\in\tilde{C}\left(\mathbb{Z}_{p}^{r},K^{d}\right)$
and suppose $S_{p}\left\{ \chi_{1}\right\} =S_{p}\left\{ \chi_{2}\right\} $.
Then, by (I): 
\[
\chi_{1}\left(\mathbf{z}\right)\overset{\textrm{(I)}}{=}S_{p}\left\{ \chi_{1}\right\} \left(\mathbf{z}\right)=S_{p}\left\{ \chi_{2}\right\} \left(\mathbf{z}\right)\overset{\textrm{(I)}}{=}\chi_{2}\left(\mathbf{z}\right),\textrm{ }\forall\mathbf{z}\in\mathbb{Z}_{p}^{r}
\]
This proves $\chi_{1}=\chi_{2}$, which establishes $S_{p}$'s injectivity.
\end{itemize}
Thus, $S_{p}$ is an isomorphism.

Q.E.D.
\begin{defn}
Let $\mathbb{F}$ be $\mathbb{Q}$ or a field extension thereof, and
let $\chi:\mathbb{N}_{0}^{r}\rightarrow\mathbb{F}^{d}$ be a function.
We say $\chi$ has (or ``admits'') a \textbf{$\left(p,q\right)$-adic
rising-continuation} \textbf{(to $K^{d}$)} whenever\index{rising-continuation}
there is a metrically complete $q$-adic field extension $K$ of $\mathbb{F}$
and a rising-continuous function $\chi^{\prime}:\mathbb{Z}_{p}^{r}\rightarrow K^{d}$
so that $\chi^{\prime}\left(\mathbf{n}\right)=\chi\left(\mathbf{n}\right)$
for all $\mathbf{n}\in\mathbb{N}_{0}^{r}$. We call any $\chi^{\prime}$
satisfying this property a \textbf{($\left(p,q\right)$-adic)}\emph{
}\textbf{rising-continuation }of $\chi$ (to $K^{d}$). 
\end{defn}
\begin{prop}
\label{prop:Uniqueness of rising-continuation, MD}Let $\chi:\mathbb{N}_{0}^{r}\rightarrow\mathbb{F}^{d}$
be a function admitting a $\left(p,q\right)$-adic rising-continuation
to $K^{d}$. Then:

\vphantom{}

I. The rising-continuation of $\chi$ is unique. Consequently, we
write $\chi^{\prime}$ to denote the rising continuation of $\chi$.

\vphantom{}

II. 
\begin{equation}
\chi^{\prime}\left(\mathbf{z}\right)\overset{K^{d}}{=}\lim_{k\rightarrow\infty}\chi\left(\left[\mathbf{z}\right]_{p^{k}}\right)=\sum_{\mathbf{n}\in\mathbb{N}_{0}^{r}}c_{\mathbf{n}}\left(\chi\right)\left[\mathbf{z}\overset{p^{\lambda_{p}\left(\mathbf{n}\right)}}{\equiv}\mathbf{n}\right],\textrm{ }\forall\mathbf{z}\in\mathbb{Z}_{p}^{r}\label{eq:MD Rising continuation limit formula}
\end{equation}
\end{prop}
Proof:

I. Suppose $\chi$ admits two (possibly distinct) rising-continuations,
$\chi^{\prime}$ and $\chi^{\prime\prime}$. To see that $\chi^{\prime}$
and $\chi^{\prime\prime}$ must be the same, we note that since the
restrictions of both $\chi^{\prime}$ and $\chi^{\prime\prime}$ to
$\mathbb{N}_{0}^{r}$ are, by definition, equal to $\chi$, these
restrictions must be equal to one another.

So, let $\mathbf{z}$ be an arbitrary element of $\left(\mathbb{Z}_{p}^{r}\right)^{\prime}$.
Then, we note that at least one entry of $\mathbf{z}$ (say, the $\ell$th
entry) necessarily has infinitely many non-zero $p$-adic digits.
Consequently, $\left\{ \left[\mathbf{z}\right]_{p^{k}}\right\} _{k\geq1}$
is a rising sequence of non-negative integer $r$-tuples converging
to $\mathbf{z}$. As such, by the rising-continuity of $\chi^{\prime}$
and $\chi^{\prime\prime}$: 
\begin{equation}
\lim_{k\rightarrow\infty}\chi^{\prime}\left(\left[\mathbf{z}\right]_{p^{k}}\right)\overset{K^{d}}{=}\chi^{\prime}\left(\mathbf{z}\right)
\end{equation}
\begin{equation}
\lim_{k\rightarrow\infty}\chi^{\prime\prime}\left(\left[\mathbf{z}\right]_{p^{k}}\right)\overset{K^{d}}{=}\chi^{\prime\prime}\left(\mathbf{z}\right)
\end{equation}
Because the $\left[\mathbf{z}\right]_{p^{k}}$s are $r$-tuples of
integers, we can then write: 
\begin{equation}
\chi^{\prime}\left(\mathbf{z}\right)=\lim_{k\rightarrow\infty}\chi^{\prime}\left(\left[\mathbf{z}\right]_{p^{k}}\right)=\lim_{k\rightarrow\infty}\chi\left(\left[\mathbf{z}\right]_{p^{k}}\right)=\lim_{k\rightarrow\infty}\chi^{\prime\prime}\left(\left[\mathbf{z}\right]_{p^{k}}\right)=\chi^{\prime\prime}\left(\mathbf{z}\right)
\end{equation}
Since $\mathbf{z}$ was arbitrary, we conclude $\chi^{\prime}\left(\mathbf{z}\right)=\chi^{\prime\prime}\left(\mathbf{z}\right)$
for all $\mathbf{z}\in\left(\mathbb{Z}_{p}^{r}\right)^{\prime}$.
Thus, $\chi^{\prime}$ and $\chi^{\prime\prime}$ are equal to one
another on both $\left(\mathbb{Z}_{p}^{r}\right)^{\prime}$ and $\mathbb{N}_{0}^{r}$,
which is all of $\mathbb{Z}_{p}^{r}$. So, $\chi^{\prime}$ and $\chi^{\prime\prime}$
are, in fact, the same function. This proves the uniqueness of $\chi$'s
rising-continuation.

\vphantom{}

II. As a rising-continuous function, $\chi^{\prime}\left(\mathbf{z}\right)$
is uniquely determined by its values on $\mathbb{N}_{0}^{r}$. Since
$\chi^{\prime}\mid_{\mathbb{N}_{0}^{r}}=\chi$, $\chi^{\prime}$ and
$\chi$ then have the same van der Put coefficients. Applying the
van der Put identity (\ref{eq:van der Put identity}) then yields
(\ref{eq:MD Rising continuation limit formula}).

Q.E.D. 
\begin{thm}
Let $\mathbb{F}$ be $\mathbb{Q}$ or a field extension thereof, let
$\chi:\mathbb{N}_{0}^{r}\rightarrow\mathbb{F}^{d}$ be a function,
and let $K$ be a metrically complete $q$-adic field extension of
$\mathbb{F}$, where $q$ is prime. Then, the following are equivalent:

\vphantom{}

I. $\chi$ admits a $\left(p,q\right)$-adic rising-continuation to
$K^{d}$.

\vphantom{}

II. For each $\mathbf{z}\in\mathbb{Z}_{p}^{r}$, $\chi\left(\left[\mathbf{z}\right]_{p^{n}}\right)$
converges to a limit in $K^{d}$ as $n\rightarrow\infty$.
\end{thm}
Proof:

i. Suppose (I) holds. Then, by \textbf{Proposition \ref{prop:Uniqueness of rising-continuation, MD}}
we have that (\ref{eq:MD Rising continuation limit formula}) holds,
which shows that (II) is true.

\vphantom{}

ii. Conversely, suppose (II) holds. Then, by the van der Put identity
(\ref{eq:van der Put identity}), $S_{p}\left\{ \chi\right\} $ is
a rising-continuous function whose restriction to $\mathbb{N}_{0}^{r}$
is equal to $\chi$. So, $S_{p}\left\{ \chi\right\} $ is the rising-continuation
of $\chi$, and hence, $\chi$ is rising-continuable.

This shows the equivalence of (II) and (I).

Q.E.D.
\begin{rem}
Like in the one-dimensional case, we will now identify a function
$\chi:\mathbb{N}_{0}^{r}\rightarrow K^{d}$ with its rising-continuation
$\chi^{\prime}$.
\end{rem}
\vphantom{}

Finally, we have the analogues of the functional equation results
from the end of Subsection \ref{subsec:3.2.1 -adic-Interpolation-of}. 
\begin{thm}
\label{thm:generic Chi_type functional equations, MD}Let $H$ be
a semi-basic $p$-smooth $d$-dimensional depth-$r$ Hydra map which
fixes $\mathbf{0}$, and consider the system of functional equations\index{functional equation!rising-continuability}:
\begin{equation}
\mathbf{f}\left(p\mathbf{n}+\mathbf{j}\right)=H_{\mathbf{j}}^{\prime}\left(\mathbf{0}\right)\mathbf{f}\left(\mathbf{n}\right)+\mathbf{c}_{\mathbf{j}},\textrm{ }\forall\mathbf{j}\in\mathbb{Z}^{r}/p\mathbb{Z}^{r},\textrm{ }\forall\mathbf{n}\in\mathbb{N}_{0}^{r}\label{eq:MD Generic H-type functional equations}
\end{equation}
where $\left\{ \mathbf{c}_{\mathbf{j}}\right\} _{\mathbf{j}\in\mathbb{Z}^{r}/p\mathbb{Z}^{r}}$
are vector constants in $\overline{\mathbb{Q}}^{d}$. Additionally,
suppose $\mathbf{I}_{d}-H^{\prime}\left(\mathbf{0}\right)$ is invertible.
Then:

\vphantom{}

I. There is a unique function $\chi:\mathbb{N}_{0}^{r}\rightarrow\overline{\mathbb{Q}}^{d}$
such that $\mathbf{f}=\chi$ is a solution of \emph{(\ref{eq:MD Generic H-type functional equations})}.

\vphantom{}

II. The solution $\chi$ \emph{(\ref{eq:MD Generic H-type functional equations})}
is rising-continuable to a function $\chi:\mathbb{Z}_{p}^{r}\rightarrow\mathbb{C}_{q_{H}}^{d}$
which satisfies: 
\begin{equation}
\chi\left(p\mathbf{z}+\mathbf{j}\right)=H_{\mathbf{j}}^{\prime}\left(\mathbf{0}\right)\chi\left(\mathbf{n}\right)+\mathbf{c}_{\mathbf{j}},\textrm{ }\forall\mathbf{j}\in\mathbb{Z}^{r}/p\mathbb{Z}^{r},\textrm{ }\forall\mathbf{z}\in\mathbb{Z}_{p}^{r}\label{eq:MD Rising-continuation Generic H-type functional equations}
\end{equation}

\vphantom{}

III. The function $\chi:\mathbb{Z}_{p}^{r}\rightarrow\mathbb{C}_{q_{H}}^{d}$
described in \emph{(III)} is the unique rising-continuous function
$\mathbf{f}:\mathbb{Z}_{p}^{r}\rightarrow\mathbb{C}_{q_{H}}^{d}$
satisfying:
\begin{equation}
\mathbf{f}\left(p\mathbf{z}+\mathbf{j}\right)=H_{\mathbf{j}}^{\prime}\left(\mathbf{0}\right)\mathbf{f}\left(\mathbf{n}\right)+\mathbf{c}_{\mathbf{j}},\textrm{ }\forall\mathbf{j}\in\mathbb{Z}^{r}/p\mathbb{Z}^{r},\textrm{ }\forall\mathbf{z}\in\mathbb{Z}_{p}^{r}
\end{equation}
\end{thm}
Proof:

I. Let $\mathbf{f}:\mathbb{N}_{0}^{r}\rightarrow\overline{\mathbb{Q}}^{d}$
be any solution of (\ref{eq:MD Generic H-type functional equations}).
Setting $\mathbf{n}=\mathbf{j}=\mathbf{0}$ yields: 
\begin{align*}
\mathbf{f}\left(\mathbf{0}\right) & =H^{\prime}\left(\mathbf{0}\right)\mathbf{f}\left(\mathbf{0}\right)+\mathbf{c}_{\mathbf{0}}\\
 & \Updownarrow\\
\mathbf{f}\left(\mathbf{0}\right) & =\left(\mathbf{I}_{d}-H^{\prime}\left(\mathbf{0}\right)\right)^{-1}\mathbf{c}_{\mathbf{0}}
\end{align*}
which is well-defined since $\mathbf{I}_{d}-H^{\prime}\left(\mathbf{0}\right)$
was given to be invertible. Then, we have that: 
\begin{equation}
\mathbf{f}\left(\mathbf{j}\right)=H_{\mathbf{j}}^{\prime}\left(\mathbf{0}\right)\mathbf{f}\left(\mathbf{0}\right)+\mathbf{c}_{\mathbf{j}},\textrm{ }\forall\mathbf{j}\in\mathbb{Z}^{r}/p\mathbb{Z}^{r}
\end{equation}
and, more generally: 
\begin{equation}
\mathbf{f}\left(\mathbf{m}\right)=H_{\left[\mathbf{m}\right]_{p}}^{\prime}\left(\mathbf{0}\right)\mathbf{f}\left(\frac{\mathbf{m}-\left[\mathbf{m}\right]_{p}}{p}\right)+\mathbf{c}_{\left[\mathbf{m}\right]_{p}},\textrm{ }\forall\mathbf{m}\in\mathbb{N}_{0}^{r}\label{eq:MD f,m, digit shifting}
\end{equation}
Next, because the map $\mathbf{m}\mapsto\frac{\mathbf{m}-\left[\mathbf{m}\right]_{p}}{p}$
sends the non-negative integer tuple: 
\begin{equation}
\mathbf{m}=\left(\sum_{k=0}^{\lambda_{p}\left(m_{1}\right)-1}m_{1,k}p^{k},\ldots,\sum_{k=0}^{\lambda_{p}\left(m_{r}\right)-1}m_{r,k}p^{k}\right)
\end{equation}
to the integer tuple: 
\begin{equation}
\mathbf{m}=\left(\sum_{k=0}^{\lambda_{p}\left(m_{1}\right)-2}m_{1,k}p^{k},\ldots,\sum_{k=0}^{\lambda_{p}\left(m_{r}\right)-2}m_{r,k}p^{k}\right)
\end{equation}
it follows that $\mathbf{m}\mapsto\frac{\mathbf{m}-\left[\mathbf{m}\right]_{p}}{p}$
eventually iterates every $\mathbf{m}\in\mathbb{N}_{0}^{r}$ to $\mathbf{0}$.
So, (\ref{eq:MD f,m, digit shifting}) implies that, for every $\mathbf{m}\in\mathbb{N}_{0}^{r}$,
$\mathbf{f}\left(\mathbf{m}\right)$ is entirely determined by $\mathbf{f}\left(\mathbf{0}\right)$
and the $\mathbf{c}_{\mathbf{j}}$s. Since $\mathbf{f}\left(\mathbf{0}\right)$
is uniquely determined by $\mathbf{c}_{\mathbf{0}}$ and $H$, (\ref{eq:MD Generic H-type functional equations})
then possesses exactly one solution, which we shall denote by $\chi$.

\vphantom{}

II. Rehashing the argument for existence of the numen $\chi_{H}$
when $H$ is semi-basic and fixes $\mathbf{0}$, because $H$ is here
semi-basic, any $\mathbf{z}\in\left(\mathbb{Z}_{p}^{r}\right)^{\prime}$
will have at least one entry infinitely many non-zero $p$-adic digits,
and hence, the product of $H_{\mathbf{j}}^{\prime}\left(\mathbf{0}\right)$
taken over all the digits of such a $\mathbf{z}$ will converge $q_{H}$-adically
to zero. Then, using (\ref{eq:MD Generic H-type functional equations}),
we see that $\chi\left(\mathbf{z}\right)$ will be a sum of the form:
\begin{equation}
\beta_{0}+\alpha_{1}\beta_{1}+\alpha_{1}\alpha_{2}\beta_{2}+\alpha_{1}\alpha_{2}\alpha_{3}\beta_{3}+\cdots
\end{equation}
where, for each $n$, $\beta_{n}$ is one of the $\mathbf{c}_{\mathbf{j}}$s
and $\alpha_{n}$ is $H_{\mathbf{j}_{n}}^{\prime}\left(\mathbf{0}\right)$,
where $\mathbf{j}_{n}$ is the $r$-tuple of the $n$th $p$-adic
digit of $\mathfrak{z}_{1}$ through the $n$th $p$-adic digit of
$\mathfrak{z}_{r}$, where $\mathbf{z}=\left(\mathfrak{z}_{1},\ldots,\mathfrak{z}_{r}\right)$.
Consequently, this series converges in $\mathbb{C}_{q_{H}}^{d}$ for
all $\mathbf{z}\in\left(\mathbb{Z}_{p}^{r}\right)^{\prime}$, seeing
as $\left\Vert H_{\mathbf{j}}^{\prime}\left(\mathbf{0}\right)\right\Vert _{q_{H}}<1$
for all $\mathbf{j}\in\left(\mathbb{Z}^{r}/p\mathbb{Z}^{r}\right)\backslash\left\{ \mathbf{0}\right\} $
because $H$ is semi-basic. This guarantees the rising-continuability
of $\chi$.

\vphantom{}

III. Because $\chi$ admits a rising continuation, this continuation
is given at every $\mathbf{z}\in\mathbb{Z}_{p}^{r}$ by the van der
Put series $S_{p}\left\{ \chi\right\} \left(\mathbf{z}\right)$ .
As such: 
\begin{equation}
\chi\left(\mathbf{z}\right)\overset{\mathbb{C}_{q_{H}}^{d}}{=}S_{p}\left\{ \chi\right\} \left(\mathbf{z}\right)\overset{\mathbb{C}_{q_{H}}^{d}}{=}\lim_{N\rightarrow\infty}\chi\left(\left[\mathbf{z}\right]_{p^{N}}\right)
\end{equation}
Because $\chi$ satisfies (\ref{eq:MD Generic H-type functional equations}),
we can write: 
\begin{equation}
\chi\left(p\left[\mathbf{z}\right]_{p^{N}}+\mathbf{j}\right)\overset{\overline{\mathbb{Q}}^{d}}{=}H_{\mathbf{j}}^{\prime}\left(\mathbf{0}\right)\chi\left(\left[\mathbf{z}\right]_{p^{N}}\right)+\mathbf{c}_{\mathbf{j}}
\end{equation}
for all $\mathbf{z}\in\mathbb{Z}_{p}^{r}$, all $\mathbf{j}\in\mathbb{Z}^{r}/p\mathbb{Z}^{r}$,
and all $N\geq0$. So, for $\mathbf{z}\in\left(\mathbb{Z}_{p}^{r}\right)^{\prime}$,
letting $N\rightarrow\infty$ yields: 
\begin{equation}
\chi\left(p\mathbf{z}+\mathbf{j}\right)\overset{\mathbb{C}_{q_{H}}^{d}}{=}H_{\mathbf{j}}^{\prime}\left(\mathbf{0}\right)\chi\left(\mathbf{z}\right)+\mathbf{c}_{\mathbf{j}},\textrm{ }\forall\mathbf{j}\in\mathbb{Z}^{r}/p\mathbb{Z}^{r},\textrm{ }\forall\mathbf{z}\in\left(\mathbb{Z}_{p}^{r}\right)^{\prime}
\end{equation}
Note that these identities hold automatically for $\mathbf{z}\in\mathbb{N}_{0}^{r}$,
because those are the cases governed by (\ref{eq:MD Generic H-type functional equations}).
The uniqueness of $\chi$ as a solution of (\ref{eq:MD Rising-continuation Generic H-type functional equations})
occurs because any $\mathbf{f}:\mathbb{Z}_{p}^{r}\rightarrow\mathbb{C}_{q_{H}}^{d}$
satisfying (\ref{eq:MD Rising-continuation Generic H-type functional equations})
has a restriction to $\mathbb{N}_{0}^{r}$ which satisfies (\ref{eq:MD Generic H-type functional equations}),
thereby forcing $\mathbf{f}=\chi$.

Q.E.D.

\vphantom{}A slightly more general version of this type of argument
is as follows: 
\begin{lem}
Fix an integer $q\geq2$, let $K$ be a metrically complete $q$-adic
field, and let $\Phi_{\mathbf{j}}:\mathbb{Z}_{p}^{r}\times K^{d}\rightarrow K^{d}$
be continuous for $\mathbf{j}\in\mathbb{Z}^{r}/p\mathbb{Z}^{r}$.
Suppose $\chi:\mathbb{N}_{0}^{r}\rightarrow K^{d}$ is $\left(p,K\right)$-adically
rising-continuable. If $\chi$ satisfies the functional equations:
\begin{equation}
\chi\left(p\mathbf{n}+\mathbf{j}\right)=\Phi_{\mathbf{j}}\left(\mathbf{n},\chi\left(\mathbf{n}\right)\right),\textrm{ }\forall\mathbf{n}\in\mathbb{N}_{0}^{r},\textrm{ }\forall\mathbf{j}\in\mathbb{Z}^{r}/p\mathbb{Z}^{r}
\end{equation}
then the rising-continuation of $\chi$ satisfies: 
\begin{equation}
\chi\left(p\mathbf{z}+\mathbf{j}\right)=\Phi_{\mathbf{j}}\left(\mathbf{z},\chi\left(\mathbf{z}\right)\right),\textrm{ }\forall\mathbf{z}\in\mathbb{Z}_{p},\textrm{ }\forall\mathbf{j}\in\mathbb{Z}^{r}/p\mathbb{Z}^{r}
\end{equation}
\end{lem}
Proof: Let everything be as given. Since: 
\begin{equation}
\chi\left(p\left[\mathbf{z}\right]_{p^{N}}+\mathbf{j}\right)=\Phi_{\mathbf{j}}\left(\left[\mathbf{z}\right]_{p^{N}},\chi\left(\left[\mathbf{z}\right]_{p^{N}}\right)\right)
\end{equation}
holds true for all $\mathbf{z}\in\mathbb{Z}_{p}$ and all $N\geq0$,
the rising-continuability of $\chi$ and the continuity of $\Phi_{\mathbf{j}}$
guarantee that: 
\begin{align*}
\chi\left(p\mathbf{z}+\mathbf{j}\right) & \overset{K^{d}}{=}\lim_{N\rightarrow\infty}\chi\left(p\left[\mathbf{z}\right]_{p^{N}}+\mathbf{j}\right)\\
 & \overset{K^{d}}{=}\lim_{N\rightarrow\infty}\Phi_{\mathbf{j}}\left(\left[\mathbf{z}\right]_{p^{N}},\chi\left(\left[\mathbf{z}\right]_{p^{N}}\right)\right)\\
 & \overset{K^{d}}{=}\Phi_{\mathbf{j}}\left(\mathbf{z},\chi\left(\mathbf{z}\right)\right)
\end{align*}
as desired.

Q.E.D.

\vphantom{}
\begin{rem}
Like in Subsection \ref{subsec:5.3.2. Interpolation-Revisited}, all
of the above holds when functions taking values in $\mathbb{C}_{q}^{d}$
(vectors) are replaced with functions taking values in $\mathbb{C}_{q}^{\rho,c}$
(matrices).
\end{rem}
\newpage{}

\section{\label{sec:5.4. Quasi-Integrability-in-Multiple}Quasi-Integrability
in Multiple Dimensions}

I see little reason to cover the multi-dimensional analogue of the
results on the Banach algebra of rising-continuous functions, seeing
as we will not use them in our analysis of the multi-dimensional $\chi_{H}$.
As such, we will instead proceed directly to the multi-dimensional
$\left(p,q\right)$-adic Fourier transform, and from there to frames
and quasi-integrability.

\subsection{Multi-Dimensional $\left(p,q\right)$-adic Fourier Analysis and Thick
Measures\label{subsec:5.4.1 Multi-Dimensional--adic-Fourier}}

With regard to multi-dimensional Fourier analysis, I will be less
pedantic than the reader has probably come to expect me to be. So
long as the reader understands how to do $\left(p,q\right)$-adic
Fourier analysis in one dimension, the formalism for Fourier analysis
of scalar or vector valued functions over $\mathbb{R}^{n}$ will provide
enough information for the reader to do multi-dimensional $\left(p,q\right)$-adic
Fourier analysis\index{multi-dimensional!Fourier transform}, with
the help of the notation provided below, of course.
\begin{notation}[\textbf{Formalism for Multi-Dimensional $\left(p,q\right)$-adic Fourier
Analysis}]
\label{nota:third batch}\ 

\vphantom{}

I. Given $\mathbf{x}\in\mathbb{Q}_{p}^{r}$, where $\mathbf{x}=\left(\mathfrak{x}_{1},\ldots,\mathfrak{x}_{r}\right)$,
we write: 
\begin{equation}
\left\{ \mathbf{x}\right\} _{p}\overset{\textrm{def}}{=}\sum_{m=1}^{r}\left\{ \mathfrak{x}_{m}\right\} _{p}\label{eq:Definition of P-adic fractional part (MD)}
\end{equation}
Consequently, for $\mathbf{t}\in\hat{\mathbb{Z}}_{p}^{r}$ and $\mathbf{z}\in\mathbb{Z}_{p}^{r}$:
\begin{equation}
\left\{ \mathbf{t}\cdot\mathbf{z}\right\} _{p}=\left\{ \mathbf{t}\mathbf{z}\right\} _{p}=\sum_{m=1}^{r}\left\{ t_{m}\mathfrak{z}_{m}\right\} _{p}\label{eq:P-adic fractional part of bold t times bold z (MD)}
\end{equation}

\vphantom{}

II. Let $f\in C\left(\mathbb{Z}_{p}^{r},\mathbb{C}_{q}\right)$ be
elementary, with $f\left(\mathbf{z}\right)=\prod_{m=1}^{r}f_{m}\left(\mathfrak{z}_{m}\right)$.
Then, the Fourier transform of $f$ is defined by: 
\begin{equation}
\hat{f}\left(\mathbf{t}\right)\overset{\textrm{def}}{=}\int_{\mathbb{Z}_{p}^{r}}f\left(\mathbf{z}\right)e^{-2\pi i\left\{ \mathbf{t}\mathbf{z}\right\} _{p}}d\mathbf{z}=\prod_{m=1}^{r}\int_{\mathbb{Z}_{p}^{r}}f\left(\mathfrak{z}_{m}\right)e^{-2\pi i\left\{ t_{m}\mathfrak{z}_{m}\right\} _{p}}d\mathfrak{z}_{m}\label{eq:Definition of Fourier transform of an elementary scalar function on Z_P}
\end{equation}
We then define $\hat{f}\left(\mathbf{t}\right)$ for any $C\left(\mathbb{Z}_{p}^{r},\mathbb{C}_{q}\right)$
by extending via linearity.

\vphantom{}

III. Let $\mathbf{f}\in C\left(\mathbb{Z}_{p},\mathbb{C}_{q}^{d}\right)$,
with $\mathbf{f}\left(\mathfrak{z}\right)=\left(f_{1}\left(\mathfrak{z}\right),\ldots,f_{d}\left(\mathfrak{z}\right)\right)$
for $f_{1},\ldots,f_{d}\in C\left(\mathbb{Z}_{p},\mathbb{C}_{q}\right)$.
Then, we define $\hat{\mathbf{f}}\left(t\right)$ as the $d$-tuple
of the Fourier transforms of the $f_{m}$s: 
\begin{equation}
\hat{\mathbf{f}}\left(t\right)\overset{\textrm{def}}{=}\left(\int_{\mathbb{Z}_{p}}f_{1}\left(\mathfrak{z}\right)e^{-2\pi i\left\{ t\mathfrak{z}\right\} _{p}}d\mathfrak{z},\ldots,\int_{\mathbb{Z}_{p}}f_{d}\left(\mathfrak{z}\right)e^{-2\pi i\left\{ t\mathfrak{z}\right\} _{p}}d\mathfrak{z}\right)\label{eq:Definition of the Fourier transform of a vector-valued function on Z_p}
\end{equation}

\vphantom{}

IV. Let $\mathbf{f}\in C\left(\mathbb{Z}_{p}^{r},\mathbb{C}_{q}^{d}\right)$,
with $\mathbf{f}\left(\mathbf{z}\right)=\left(f_{1}\left(\mathbf{z}\right),\ldots,f_{d}\left(\mathbf{z}\right)\right)$
with $f_{1},\ldots,f_{d}\in C\left(\mathbb{Z}_{p}^{r},\mathbb{C}_{q}\right)$.
Then, we define: 
\begin{equation}
\hat{\mathbf{f}}\left(\mathbf{t}\right)\overset{\textrm{def}}{=}\left(\int_{\mathbb{Z}_{p}^{r}}f_{1}\left(\mathbf{z}\right)e^{-2\pi i\left\{ \mathbf{t}\mathbf{z}\right\} _{p}}d\mathbf{z},\ldots,\int_{\mathbb{Z}_{p}^{r}}f_{d}\left(\mathbf{z}\right)e^{-2\pi i\left\{ \mathbf{t}\mathbf{z}\right\} _{p}}d\mathbf{z}\right)\label{eq:Definition of the Fourier transform of a vector-valued function on Z_P}
\end{equation}
This will also be written as: 
\begin{equation}
\hat{\mathbf{f}}\left(\mathbf{t}\right)=\int_{\mathbb{Z}_{p}^{r}}\mathbf{f}\left(\mathbf{z}\right)e^{-2\pi i\left\{ \mathbf{t}\mathbf{z}\right\} _{p}}d\mathbf{z}\label{eq:Non-tuple formula for the Fourier transform of a vector-valued function on Z_P}
\end{equation}

\vphantom{}

V. Most generally, let $\mathbf{F}\in C\left(\mathbb{Z}_{p}^{r},\mathbb{C}_{q}^{\rho,c}\right)$,
with $\mathbf{F}\left(\mathbf{z}\right)=\left\{ F_{j,k}\left(\mathbf{z}\right)\right\} _{1\leq j\leq\rho,1\leq k\leq c}$
where each $F_{j,k}$ is in $C\left(\mathbb{Z}_{p}^{r},\mathbb{C}_{q}\right)$.
Then: 
\begin{equation}
\hat{\mathbf{F}}\left(\mathbf{t}\right)\overset{\textrm{def}}{=}\left\{ \int_{\mathbb{Z}_{p}^{r}}F_{j,k}\left(\mathbf{z}\right)e^{-2\pi i\left\{ \mathbf{t}\mathbf{z}\right\} _{p}}d\mathbf{z}\right\} _{1\leq j\leq\rho,1\leq k\leq c}\label{eq:Definition of the Fourier transform of a matrix-valued function on Z_P}
\end{equation}
This will also be written as: 
\begin{equation}
\hat{\mathbf{F}}\left(\mathbf{t}\right)=\int_{\mathbb{Z}_{p}^{r}}\mathbf{F}\left(\mathbf{z}\right)e^{-2\pi i\left\{ \mathbf{t}\mathbf{z}\right\} _{p}}d\mathbf{z}\label{eq:Non-tuple formula for the Fourier transform of a matrix-valued function on Z_P}
\end{equation}
\end{notation}
Like in the one-dimensional case, the most important identity for
us is the Fourier series for an indicator function for the clopen
set $\mathbf{n}+p^{N}\mathbb{Z}_{p}^{r}$:
\begin{prop}[\textbf{Fourier Series for the Multi-Dimensional Indicator Function}]
\label{prop:Multi-Dimensional indicator function Fourier Series}For
$r$-tuples $\mathbf{z}$ and $\mathbf{n}$:
\begin{equation}
\left[\mathbf{z}\overset{p^{N}}{\equiv}\mathbf{n}\right]=\frac{1}{p^{Nr}}\sum_{\left\Vert \mathbf{t}\right\Vert _{p}\leq p^{N}}e^{2\pi i\left\{ \mathbf{t}\left(\mathbf{z}-\mathbf{n}\right)\right\} _{p}}\label{eq:MD Fourier series for indicator function}
\end{equation}
\end{prop}
Proof: 
\begin{align*}
\left[\mathbf{z}\overset{p^{N}}{\equiv}\mathbf{n}\right] & =\prod_{\ell=1}^{r}\left[\mathfrak{z}_{\ell}\overset{p^{N}}{\equiv}n_{\ell}\right]\\
 & =\prod_{\ell=1}^{r}\frac{1}{p^{N}}\sum_{k_{\ell}=0}^{p^{N}-1}e^{2\pi i\left\{ \frac{k_{\ell}}{p^{N}}\left(\mathfrak{z}_{\ell}-n_{\ell}\right)\right\} _{p}}\\
 & =\frac{1}{p^{Nr}}\prod_{\ell=1}^{r}\sum_{\left|t_{\ell}\right|_{p}\leq p^{N}}e^{2\pi i\left\{ t_{\ell}\left(\mathfrak{z}_{\ell}-n_{\ell}\right)\right\} _{p}}\\
 & =\frac{1}{p^{Nr}}\sum_{\left\Vert \mathbf{t}\right\Vert _{p}\leq p^{N}}e^{2\pi i\left\{ \mathbf{t}\left(\mathbf{z}-\mathbf{n}\right)\right\} _{p}}
\end{align*}

Q.E.D.

\vphantom{}

Just to be extra cautious, let us prove that the following formal
summation identities \emph{do }in fact hold true: 
\begin{equation}
\sum_{\left\Vert \mathbf{t}\right\Vert _{p}\leq p^{N}}f\left(\mathbf{t}\right)=f\left(\mathbf{0}\right)+\sum_{n=1}^{N}\sum_{\left\Vert \mathbf{t}\right\Vert _{p}=p^{n}}f\left(\mathbf{t}\right)\label{eq:First MD level set summation identity}
\end{equation}
\begin{equation}
\sum_{\left\Vert \mathbf{t}\right\Vert _{p}=p^{n}}f\left(\mathbf{t}\right)=\sum_{\left\Vert \mathbf{t}\right\Vert _{p}\leq p^{n}}f\left(\mathbf{t}\right)-\sum_{\left\Vert \mathbf{t}\right\Vert _{p}\leq p^{n-1}}f\left(\mathbf{t}\right)\label{eq:Second MD level set summation identity}
\end{equation}
so as to justify all this hullabaloo over notation.
\begin{prop}
For all $n\geq1$: 
\begin{equation}
\left\{ \left\Vert \mathbf{t}\right\Vert _{p}=p^{n}\right\} =\left\{ \left\Vert \mathbf{t}\right\Vert _{p}\leq p^{n}\right\} \backslash\left\{ \left\Vert \mathbf{t}\right\Vert _{p}\leq p^{n-1}\right\} \label{eq:MD Level set decomposition}
\end{equation}
\end{prop}
Proof:

I. If $\left\Vert \mathbf{t}\right\Vert _{p}=p^{n}$, then $\left\Vert \mathbf{t}\right\Vert _{p}\leq p^{n}$.
However, there is at least one $m\in\left\{ 1,\ldots,r\right\} $
so that $-v_{p}\left(t_{m}\right)=n>n-1$. As such, $\left\Vert \mathbf{t}\right\Vert _{p}>p^{n-1}$,
and so:
\begin{equation}
\left\{ \left\Vert \mathbf{t}\right\Vert _{p}=p^{n}\right\} \subseteq\left\{ \left\Vert \mathbf{t}\right\Vert _{p}\leq p^{n}\right\} \backslash\left\{ \left\Vert \mathbf{t}\right\Vert _{p}\leq p^{n-1}\right\} 
\end{equation}

\vphantom{}

II. If $\left\Vert \mathbf{t}\right\Vert _{p}\leq p^{n}$ and $\left\Vert \mathbf{t}\right\Vert _{p}\nleq p^{n-1}$,
then $-v_{p}\left(t_{m}\right)\leq n$ occurs for all $m\in\left\{ 1,\ldots,r\right\} $.
Moreover, there is at least one such $m$ for which $-v_{p}\left(t_{m}\right)=n$.
This means $\left\Vert \mathbf{t}\right\Vert _{p}=p^{n}$, and so:
\begin{equation}
\left\{ \left\Vert \mathbf{t}\right\Vert _{p}\leq p^{n}\right\} \backslash\left\{ \left\Vert \mathbf{t}\right\Vert _{p}\leq p^{n-1}\right\} \subseteq\left\{ \left\Vert \mathbf{t}\right\Vert _{p}=p^{n}\right\} 
\end{equation}

Hence, the two sets are equal.

Q.E.D.

\vphantom{}Consequently, we have: 
\begin{prop}
\label{prop:MD ramanujan sum} 
\begin{equation}
\sum_{\left\Vert \mathbf{t}\right\Vert _{p}=p^{n}}e^{2\pi i\left\{ \mathbf{t}\mathbf{z}\right\} _{p}}=p^{rn}\left[\mathbf{z}\overset{p^{n}}{\equiv}\mathbf{0}\right]-p^{r\left(n-1\right)}\left[\mathbf{z}\overset{p^{n-1}}{\equiv}\mathbf{0}\right],\textrm{ }\forall n\geq1,\textrm{ }\forall\mathbf{z}\in\mathbb{Z}_{p}^{r}\label{eq:MD Ramanujan sum}
\end{equation}
\end{prop}
Proof: Use (\ref{eq:Second MD level set summation identity}) and
\textbf{Proposition \ref{prop:Multi-Dimensional indicator function Fourier Series}}.

Q.E.D. 
\begin{thm}
Let $\mathbf{F}\left(\mathbf{t}\right)=\bigodot_{m=1}^{r}\mathbf{F}_{m}\left(t_{m}\right)$
be an elementary function in $C\left(\mathbb{Z}_{p}^{r},\mathbb{C}_{q}^{\rho,c}\right)$.
Then: 
\begin{equation}
\hat{\mathbf{F}}\left(\mathbf{t}\right)=\hat{\mathbf{F}}\left(t_{1},\ldots,t_{r}\right)=\bigodot_{m=1}^{r}\hat{\mathbf{F}}_{m}\left(t_{m}\right)\label{eq:Fourier transform mixes with tensoring}
\end{equation}
\end{thm}
Proof: The complex exponentials turn multiplication into addition.

Q.E.D.
\begin{thm}
The $\mathbb{C}_{q}$-valued functions $\left\{ \mathbf{z}\mapsto e^{2\pi i\left\{ \mathbf{t}\mathbf{z}\right\} _{p}}\right\} _{\mathbf{t}\in\hat{\mathbb{Z}}_{p}^{r}}$
form a basis of $C\left(\mathbb{Z}_{p}^{r},\mathbb{C}_{q}\right)$.
More generally, the vectors: 
\[
\left[\begin{array}{c}
e^{2\pi i\left\{ \mathbf{t}\mathbf{z}\right\} _{p}}\\
0\\
\vdots\\
0
\end{array}\right],\left[\begin{array}{c}
0\\
e^{2\pi i\left\{ \mathbf{t}\mathbf{z}\right\} _{p}}\\
\vdots\\
0
\end{array}\right],\ldots,\left[\begin{array}{c}
0\\
0\\
\vdots\\
e^{2\pi i\left\{ \mathbf{t}\mathbf{z}\right\} _{p}}
\end{array}\right]
\]
form a basis of $C\left(\mathbb{Z}_{p}^{r},\mathbb{C}_{q}^{d}\right)$.
\end{thm}
Proof: $C\left(\mathbb{Z}_{p}^{r},\mathbb{C}_{q}^{d}\right)$ is the
closure in $\left\Vert \cdot\right\Vert _{p,q}$-norm of the set of
linear combinations of elementary continuous functions $\mathbf{f}:\mathbb{Z}_{p}^{r}\rightarrow\mathbb{C}_{q}^{d}$.
By the properties of the tensor product and the Fourier transform,
every such elementary function $\mathbf{f}\left(\mathbf{z}\right)=\prod_{m=1}^{r}\mathbf{f}_{m}\left(\mathfrak{z}_{m}\right)$
is then uniquely expressible as: 
\begin{align*}
\mathbf{f}\left(\mathbf{z}\right) & =\prod_{m=1}^{r}\sum_{t_{m}\in\hat{\mathbb{Z}}_{p}}\hat{\mathbf{f}}\left(t_{m}\right)e^{2\pi i\left\{ t_{m}\mathfrak{z}_{m}\right\} _{p}}\\
 & =\sum_{t_{1}\in\hat{\mathbb{Z}}_{p}}\cdots\sum_{t_{r}\in\hat{\mathbb{Z}}_{p}}\left(\prod_{m=1}^{r}\hat{\mathbf{f}}\left(t_{m}\right)\right)e^{2\pi i\sum_{\ell=1}^{r}\left\{ t_{\ell}\mathfrak{z}_{\ell}\right\} _{p}}\\
 & =\sum_{\mathbf{t}\in\hat{\mathbb{Z}}_{p}}\hat{\mathbf{f}}\left(\mathbf{t}\right)e^{2\pi i\left\{ \mathbf{t}\mathbf{z}\right\} _{p}}
\end{align*}
The theorem follows by the density of these elementary functions in
$C\left(\mathbb{Z}_{p}^{r},\mathbb{C}_{q}^{d}\right)$.

Q.E.D.
\begin{thm}
The Fourier transform $C\left(\mathbb{Z}_{p}^{r},\mathbb{C}_{q}^{\rho,c}\right)\rightarrow c_{0}\left(\hat{\mathbb{Z}}_{p}^{r},\mathbb{C}_{q}^{\rho,c}\right)$
is an isometric isomorphism of non-archimedean Banach spaces.
\end{thm}
Proof: Use linearity and tensors to extend the one-dimensional case
(\textbf{Corollary \ref{cor:pq adic Fourier transform is an isometric isomorphism}}).

Q.E.D.
\begin{defn}
Given a function $\hat{\mathbf{F}}:\hat{\mathbb{Z}}_{p}^{r}\rightarrow\mathbb{C}_{q}^{\rho,c}$
and an $N\geq0$, we write $\tilde{\mathbf{F}}_{N}$ to denote the
function \nomenclature{$\tilde{\mathbf{F}}_{N}$}{$N$th partial Fourier series generated by $\hat{\mathbf{F}}\left(\mathbf{t}\right)$}$\tilde{\mathbf{F}}_{N}:\mathbb{Z}_{p}^{r}\rightarrow\mathbb{C}_{q}^{\rho,c}$
defined by: 
\begin{equation}
\tilde{\mathbf{F}}_{N}\left(\mathbf{z}\right)\overset{\textrm{def}}{=}\sum_{\left\Vert \mathbf{t}\right\Vert _{p}\leq p^{N}}\hat{\mathbf{F}}\left(\mathbf{t}\right)e^{2\pi i\left\{ \mathbf{t}\mathbf{z}\right\} _{p}}\label{eq:Definition of bold F_N twiddle}
\end{equation}
We call $\tilde{\mathbf{F}}_{N}$ the \textbf{$N$th partial Fourier
series }generated by $\hat{\mathbf{F}}$.
\end{defn}
\vphantom{}

At this point, were we in the one-dimensional case, we would introduce
measures. However, the situation is not as simple in the multi-dimensional
case. Let us review our set-up so far. $\chi_{H}$\textemdash the
function we want to study\textemdash goes from $\mathbb{Z}_{p}^{r}$
to $\mathbb{Z}_{q}^{d}$. As such, our setting will be functions $\mathbb{Z}_{p}^{r}\rightarrow\mathbb{C}_{q}^{d}$.
The goal is to interpret integration of $\chi_{H}$ as the action
of \emph{some sort }of linear map acting on $C\left(\mathbb{Z}_{p}^{r},\mathbb{C}_{q}^{d}\right)$.
In the one-dimensional case, this linear map was a continuous, \emph{scalar-valued}
linear operator $C\left(\mathbb{Z}_{p},K\right)\rightarrow K$ (a.k.a.
a \emph{measure}), which acted upon continuous functions by way of
the formula:

\begin{equation}
f\in C\left(\mathbb{Z}_{p},K\right)\mapsto\sum_{t\in\hat{\mathbb{Z}}_{p}}\hat{f}\left(-t\right)\hat{\mu}\left(t\right)\in K
\end{equation}
Because $\chi_{H}$ is vector-valued, in order to replicate this Parseval-Plancherel
identity, we will need to find a way for the vector $\hat{\chi}_{H}\left(\mathbf{t}\right)$
to act on $\hat{\mathbf{f}}\left(\mathbf{t}\right)\in C\left(\mathbb{Z}_{p}^{r},\mathbb{C}_{q}^{d}\right)$.
However, because we want to use this action to integrate $\chi_{H}\left(\mathbf{z}\right)$
entry-by-entry, the \emph{result }of the action cannot be scalar:
it needs to be a vector.

Thus, our linear operator will not send $\hat{\mathbf{f}}\left(\mathbf{t}\right)$
to a scalar, but to a \emph{vector}. So, instead of linear functionals,
we will work with linear operators on finite-dimensional vector spaces.
This is where a subtlety emerges: \emph{we will need to distinguish
between different realizations of linear operators}.\emph{ }For example,
we could have a linear operator which sends $\mathbf{f}\left(\mathbf{z}\right)\in C\left(\mathbb{Z}_{p}^{r},\mathbb{C}_{q}^{d}\right)$
to the sum: 
\begin{equation}
\sum_{\mathbf{t}\in\hat{\mathbb{Z}}_{p}^{r}}\hat{\mathbf{f}}\left(-\mathbf{t}\right)\hat{\mathbf{m}}\left(\mathbf{t}\right)
\end{equation}
where $\hat{\mathbf{m}}:\hat{\mathbb{Z}}_{p}^{r}\rightarrow\mathbb{C}_{q}^{d}$
is a column-vector valued function on $\hat{\mathbb{Z}}_{p}^{r}$
and the juxtaposition $\hat{\mathbf{f}}\left(-\mathbf{t}\right)\hat{\mathbf{m}}\left(\mathbf{t}\right)$
denotes the entry-wise product of $\hat{\mathbf{f}}\left(-\mathbf{t}\right)$
and $\mathbf{\hat{m}}\left(\mathbf{t}\right)$. Alternatively, we
could have a linear operator which sends $\mathbf{f}\left(\mathbf{z}\right)$
to: 
\begin{equation}
\sum_{\mathbf{t}\in\hat{\mathbb{Z}}_{p}^{r}}\hat{\mathbf{M}}\left(\mathbf{t}\right)\hat{\mathbf{f}}\left(-\mathbf{t}\right)
\end{equation}
where $\hat{\mathbf{M}}:\hat{\mathbb{Z}}_{p}^{r}\rightarrow\mathbb{C}_{q}^{d,d}$
is a $d\times d$-matrix valued function and the juxtaposition $\hat{\mathbf{M}}\left(\mathbf{t}\right)\hat{\mathbf{f}}\left(-\mathbf{t}\right)$
denotes ordinary left-multiplication of the vector $\hat{\mathbf{f}}\left(-\mathbf{t}\right)$
by the matrix $\hat{\mathbf{M}}\left(\mathbf{t}\right)$.

Note that $\hat{\mathbf{f}}\left(-\mathbf{t}\right)\hat{\mathbf{m}}\left(\mathbf{t}\right)$
could also be written as left-multiplication of $\hat{\mathbf{f}}\left(-\mathbf{t}\right)$
by the $d\times d$ diagonal matrix whose diagonal entries are precisely
the entries of the vector $\hat{\mathbf{m}}\left(-\mathbf{t}\right)$.
However, the object we would obtain by summing the Fourier series
generated by said matrix would be a $d\times d$ matrix rather than
a $d\times1$ vector. With all of this in mind, instead of forcing
the reader to keep different context-based identifications in mind,
I feel it will be simpler to distinguish between those linear operators
which are represented by entry-wise multiplication against a vector
and those which are represented by left-multiplication by a matrix.
These considerations motivate our next definition, along with those
elaborated in the examples that follow.
\begin{defn}[\textbf{Thick measures}]
 \index{measure!thick}\index{thick measure}A\textbf{ $\left(p,q\right)$-adic}
\textbf{thick measure of depth $r$ }is a continuous linear operator
$\mathcal{M}:C\left(\mathbb{Z}_{p}^{r},\mathbb{C}_{q}^{d}\right)\rightarrow\mathbb{C}_{q}^{d}$.
I write\footnote{In this way, we can reserve the notation $C\left(\mathbb{Z}_{p}^{r},\mathbb{C}_{q}^{d}\right)^{\prime}$
to denote the dual space of $C\left(\mathbb{Z}_{p}^{r},\mathbb{C}_{q}^{d}\right)$\textemdash continuous
\emph{scalar}-valued linear operators on $C\left(\mathbb{Z}_{p}^{r},\mathbb{C}_{q}^{d}\right)$.} \nomenclature{$C\left(\mathbb{Z}_{p}^{r},\mathbb{C}_{q}^{d}\right)^{*}$}{Thick measures on $C\left(\mathbb{Z}_{p}^{r},\mathbb{C}_{q}^{d}\right)$}$C\left(\mathbb{Z}_{p}^{r},\mathbb{C}_{q}^{d}\right)^{*}$
to denote the space of vector measures on $C\left(\mathbb{Z}_{p}^{r},\mathbb{C}_{q}^{d}\right)$.
Note that this is a $\mathbb{C}_{q}$-Banach space under operator
norm.

Finally, we will write $\mathcal{M}\left(\mathbf{f}\right)$ to denote
the $d\times1$ vector obtained by applying $\mathcal{M}$ to $\mathbf{f}\in C\left(\mathbb{Z}_{p}^{r},\mathbb{C}_{q}^{d}\right)$.
\end{defn}
\begin{rem}
I use the asterisk $*$ for the space of thick measures rather than
a $\prime$, because, for the sake of consistency, the space denoted
$C\left(\mathbb{Z}_{p}^{r},\mathbb{C}_{q}^{d}\right)^{\prime}$ should
be the dual of $C\left(\mathbb{Z}_{p}^{r},\mathbb{C}_{q}^{d}\right)$,
which is the space of continuous \emph{scalar-valued }linear maps
on $C\left(\mathbb{Z}_{p}^{r},\mathbb{C}_{q}^{d}\right)$, as opposed
to $C\left(\mathbb{Z}_{p}^{r},\mathbb{C}_{q}^{d}\right)^{*}$, which
is the space of continuous \emph{vector-valued }linear maps on $C\left(\mathbb{Z}_{p}^{r},\mathbb{C}_{q}^{d}\right)$.
\end{rem}
\begin{example}[\textbf{Thick measures of matrix type}]
\label{exa:thick matrix measures}Consider a bounded $d\times d$-matrix-valued
function $\hat{\mathbf{M}}\in B\left(\hat{\mathbb{Z}}_{p}^{r},\mathbb{C}_{q}^{d,d}\right)$.
Then, the map: 
\begin{equation}
\mathbf{f}\in C\left(\mathbb{Z}_{p}^{r},\mathbb{C}_{q}^{d}\right)\mapsto\sum_{\mathbf{t}\in\hat{\mathbb{Z}}_{p}^{r}}\hat{\mathbf{M}}\left(\mathbf{t}\right)\hat{\mathbf{f}}\left(-\mathbf{t}\right)\in\mathbb{C}_{q}^{d}\label{eq:Action of a matrix-valued multiplier}
\end{equation}
defines a continuous linear operator $\mathcal{M}:C\left(\mathbb{Z}_{p}^{r},\mathbb{C}_{q}^{d}\right)\rightarrow\mathbb{C}_{q}^{d}$,
and hence, a $\mathbb{C}_{q}^{d}$-valued thick measure. I call such
an $\mathcal{M}$ a \textbf{thick measure of} \textbf{matrix type}\index{thick measure!matrix type}.\textbf{
}Consequently, we define $\hat{\mathbf{M}}:\hat{\mathbb{Z}}_{p}^{r}\rightarrow\mathbb{C}_{q}^{d,d}$
as the \textbf{Fourier-Stieltjes transform} of $\mathcal{M}$, and
write: 
\begin{equation}
\hat{\mathcal{M}}\left(\mathbf{t}\right)\overset{\textrm{def}}{=}\hat{\mathbf{M}}\left(\mathbf{t}\right)\label{eq:Fourier-Stieltjes transform of a matrix type thick measure}
\end{equation}
\index{Fourier-Stieltjes transform} 
\end{example}
\begin{example}[\textbf{Thick measures of vector type}]
\label{exa:thick vector measures}Consider a bounded $d\times1$-vector-valued
function $\hat{\mathbf{m}}\in B\left(\hat{\mathbb{Z}}_{p}^{r},\mathbb{C}_{q}^{d}\right)$,
as well as the map: 
\begin{equation}
\mathbf{f}\in C\left(\mathbb{Z}_{p}^{r},\mathbb{C}_{q}^{d}\right)\mapsto\sum_{\mathbf{t}\in\hat{\mathbb{Z}}_{p}^{r}}\hat{\mathbf{m}}\left(\mathbf{t}\right)\hat{\mathbf{f}}\left(-\mathbf{t}\right)\in\mathbb{C}_{q}^{d}\label{eq:Action of a vector-valued multiplier}
\end{equation}
Here $\hat{\mathbf{m}}\left(\mathbf{t}\right)=\left(\hat{m}_{1}\left(\mathbf{t}\right),\ldots,\hat{m}_{d}\left(\mathbf{t}\right)\right)$
and we write: 
\begin{equation}
\hat{\mathbf{m}}\left(\mathbf{t}\right)\hat{\mathbf{f}}\left(-\mathbf{t}\right)=\left(\begin{array}{c}
\hat{m}_{1}\left(\mathbf{t}\right)\hat{f}_{1}\left(-\mathbf{t}\right)\\
\vdots\\
\hat{m}_{d}\left(\mathbf{t}\right)\hat{f}_{d}\left(-\mathbf{t}\right)
\end{array}\right)\in\mathbb{C}_{q}^{d}
\end{equation}
to denote the entry-wise product of the column vectors $\hat{\mathbf{m}}\left(\mathbf{t}\right)$
and $\hat{\mathbf{f}}\left(-\mathbf{t}\right)$. As written, (\ref{eq:Action of a vector-valued multiplier})
defines a a continuous linear operator $\mathcal{M}:C\left(\mathbb{Z}_{p}^{r},\mathbb{C}_{q}^{d}\right)\rightarrow\mathbb{C}_{q}^{d}$,
and hence, a $\mathbb{C}_{q}^{d}$-valued Banach measure. I call such
an $\mathcal{M}$ a \textbf{thick measure of vector type}\index{thick measure!vector type}.\textbf{
}The \textbf{Fourier-Stieltjes transform} of $\mathcal{M}$ is given
by: 
\begin{equation}
\hat{\mathcal{M}}\left(\mathbf{t}\right)\overset{\textrm{def}}{=}\hat{\mathbf{m}}\left(\mathbf{t}\right)\label{eq:Fourier transform of a vector-type multiplier}
\end{equation}
\end{example}
\begin{defn}
Given a thick measure $\mathcal{M}:C\left(\mathbb{Z}_{p}^{r},\mathbb{C}_{q}^{d}\right)\rightarrow\mathbb{C}_{q}^{d}$
of matrix or vector type, the \textbf{Fourier-Stieltjes transform
}of $\mathcal{M}$ is the function $\hat{\mathcal{M}}$ defined in
the above two examples, respectively. We say a thick measure $\mathcal{M}$
is \textbf{algebraic }whenever the entries of $\hat{\mathcal{M}}\left(\mathbf{t}\right)$
are elements of $\overline{\mathbb{Q}}$ for all $\mathbf{t}\in\hat{\mathbb{Z}}_{p}^{r}$.
\end{defn}
\begin{defn}
We write $M\left(\mathbb{Z}_{p}^{r},\mathbb{C}_{q}^{d}\right)$ and
$M\left(\mathbb{Z}_{p}^{r},\mathbb{C}_{q}^{d,d}\right)$ \nomenclature{$M\left(\mathbb{Z}_{p}^{r},\mathbb{C}_{q}^{d}\right)$}{thick vector-type measures}
\nomenclature{$M\left(\mathbb{Z}_{p}^{r},\mathbb{C}_{q}^{\left(d,d\right)}\right)$}{thick matrix-type measures}
to denote the sets (in fact, $\mathbb{C}_{q}$-linear spaces) of vector-
and matrix-type thick measures on $C\left(\mathbb{Z}_{p}^{r},\mathbb{C}_{q}^{d}\right)$,
respectively. By definition, $M\left(\mathbb{Z}_{p}^{r},\mathbb{C}_{q}^{d}\right)$
and $M\left(\mathbb{Z}_{p}^{r},\mathbb{C}_{q}^{d,d}\right)$ are the
operators induced by Fourier-Stieltjes transforms in $B\left(\hat{\mathbb{Z}}_{p}^{r},\mathbb{C}_{q}^{d}\right)$
and $B\left(\hat{\mathbb{Z}}_{p}^{r},\mathbb{C}_{q}^{d,d}\right)$,
respectively, in the manner described in \textbf{Examples \ref{exa:thick vector measures}}
and \textbf{\ref{exa:thick matrix measures}}.
\end{defn}
\begin{defn}[\textbf{$N$th partial kernel of a thick measure}]
Consider a thick $\left(p,q\right)$-adic measure $\mathcal{M}$
with Fourier-Stieltjes transform $\hat{\mathcal{M}}\left(\mathbf{t}\right)$.
For each $N\geq0$, we write \nomenclature{$\tilde{\mathcal{M}}_{N}\left(\mathbf{z}\right)$}{$N$th partial kernel of $\mathcal{M}$}$\tilde{\mathcal{M}}_{N}$
to denote the \textbf{$N$th partial kernel }of\index{thick measure!partial kernel}
$\mathcal{M}$, which is defined as the continuous (in fact, locally
constant)\emph{ }function:

\begin{equation}
\tilde{\mathcal{M}}_{N}\left(\mathbf{z}\right)\overset{\textrm{def}}{=}\sum_{\left\Vert \mathbf{t}\right\Vert _{p}\leq p^{N}}\hat{\mathcal{M}}\left(\mathbf{t}\right)e^{2\pi i\left\{ \mathbf{t}\mathbf{z}\right\} _{p}},\textrm{ }\forall\mathbf{z}\in\mathbb{Z}_{p}^{r}\label{eq:Concrete Definition of the Nth Partial Kernel of a Fourier Multiplier}
\end{equation}
This function is vector-valued when $\mathcal{M}$ is vector-type
and is matrix-valued when $\mathcal{M}$ is matrix-type. 
\end{defn}

\subsection{\label{subsec:5.4.2 Multi-Dimensional-Frames}Multi-Dimensional Frames}
\begin{rem}
Like with \ref{subsec:3.3.3 Frames}, the present subsection is going
to be heavy on abstract definitions. Practical-minded readers can
safely skip this section so long as the keep the following concepts
in mind:

Consider a function $\mathbf{F}:\mathbb{Z}_{p}^{r}\rightarrow\mathbb{C}_{q}^{\rho,c}$.
A sequence of functions $\left\{ \mathbf{F}_{n}\right\} _{n\geq1}$
on $\mathbb{Z}_{p}^{r}$ is said to \textbf{converge} \textbf{to $f$
with respect to the standard $\left(p,q\right)$-adic frame }if\index{frame!standard left(p,qright)-adic@standard $\left(p,q\right)$-adic}:

\vphantom{}

I. For all $n$, $\mathbf{F}_{n}\left(\mathbf{z}\right)\in\overline{\mathbb{Q}}^{\rho,c}$
for all $\mathbf{z}\in\mathbb{Z}_{p}^{r}$.

\vphantom{}

II. For all $\mathbf{z}\in\mathbb{N}_{0}^{r}$, $\mathbf{F}\left(\mathbf{z}\right)\in\overline{\mathbb{Q}}^{\rho,c}$
and $\lim_{n\rightarrow\infty}\mathbf{F}_{n}\left(\mathbf{z}\right)\overset{\mathbb{C}^{\rho,c}}{=}\mathbf{F}\left(\mathbf{z}\right)$
(meaning the convergence is in the topology of $\mathbb{C}$).

\vphantom{}

III. For all $\mathbf{z}\in\left(\mathbb{Z}_{p}^{r}\right)^{\prime}$,
$\mathbf{F}\left(\mathbf{z}\right)\in\mathbb{C}_{q}^{\rho,c}$ and
$\lim_{n\rightarrow\infty}\mathbf{F}_{n}\left(\mathbf{z}\right)\overset{\mathbb{C}_{q}^{\rho,c}}{=}\mathbf{F}\left(\mathbf{z}\right)$
(meaning the convergence is in the topology of $\mathbb{C}_{q}$).

\vphantom{}

We write $\mathcal{F}_{p,q}$ to denote the standard $\left(p,q\right)$-adic
frame. In particular, we write:

\vphantom{}

i. $\mathcal{F}_{p,q}^{d}$, to denote the case where $\mathbf{F}\left(\mathbf{z}\right)$
is a $d\times1$-column-vector-valued function (the $d$-dimensional
standard frame\index{frame!standard left(p,qright)-adic@standard $\left(p,q\right)$-adic!$d$-dimensional});

\vphantom{}

ii. $\mathcal{F}_{p,q}^{d,d}$ to denote the case where $\mathbf{F}\left(\mathbf{z}\right)$
is a $d\times d$-matrix-valued function (the $d,d$-dimensional standard
frame \index{frame!standard left(p,qright)-adic@standard $\left(p,q\right)$-adic!left(d,dright)
-dimensional@$d,d$-dimensional}).

\vphantom{}

We then write $\lim_{n\rightarrow\infty}\mathbf{f}_{n}\left(\mathbf{z}\right)\overset{\mathcal{F}_{p,q}^{d}}{=}\mathbf{f}\left(\mathbf{z}\right)$
and $\lim_{n\rightarrow\infty}\mathbf{F}_{n}\left(\mathbf{z}\right)\overset{\mathcal{F}_{p,q}^{d,d}}{=}\mathbf{F}\left(\mathbf{z}\right)$
to denote convergence with respect to the $d$-dimensional and $d,d$-dimensional
standard frames, respectively.
\end{rem}
\vphantom{}

HERE BEGINS THE DISCUSSION OF MULTI-DIMENSIONAL FRAMES

\vphantom{}

The set-up for the multi-dimensional case will be a tad bit more involved
than the one-dimensional case. The most significant difference is
that because our functions can now take vector or matrix values, instead
of assigning topological fields to each point of the domain, we will
assign topological \emph{vector spaces}, all of which contain a finite-dimensional
$\overline{\mathbb{Q}}$-vector-space of an appropriate dimension.
These will come in two flavors: those with the discrete topology,
and those which are Banach spaces with respect to a norm induced by
an absolute value.
\begin{defn}[\textbf{Multi-Dimensional Frames}]
\label{def:MD frame}Fix an integer $d\geq1$ . A $d$-dimensional
$p$-adic \index{frame!$p$-adic}\index{frame}(or just ``frame'',
for short) $\mathcal{F}$ of depth $r$ is the following collection
of data:

\vphantom{}

I. A set $U_{\mathcal{F}}\subseteq\mathbb{Z}_{p}^{r}$, called the
\textbf{domain }of $\mathcal{F}$.

\vphantom{}

II. For each $\mathbf{z}\in U_{\mathcal{F}}$, a topological field
$K_{\mathbf{z}}$ containing $\overline{\mathbb{Q}}$ and a $d$-dimensional
topological vector space $V_{\mathbf{z}}$ over $K_{\mathbf{z}}$.
\nomenclature{$K_{\mathbf{z}}^{d}$}{ }\nomenclature{$V_{\mathbf{z}}$}{ }
We allow for $K_{\mathbf{z}}=\overline{\mathbb{Q}}$. In particular,
for any $\mathbf{z}\in U_{\mathcal{F}}$, we require the topology
of $K_{\mathbf{z}}$ to be either the discrete topology\emph{ or}
the topology induced by an absolute value,\nomenclature{$K_{\mathbf{z}}$}{ },
denoted \nomenclature{$\left|\cdot\right|_{K_{\mathbf{z}}}$}{ }$\left|\cdot\right|_{K_{\mathbf{z}}}$.
For any $\mathbf{z}$ for which $K_{\mathbf{z}}$ is given the discrete
topology, we endow $V_{\mathbf{z}}$ with the discrete topology as
well. If $K_{\mathbf{z}}$ is given the topology induced by an absolute
value, we require $\left(K_{\mathbf{z}},\left|\cdot\right|_{K_{\mathbf{z}}}\right)$
to be a complete metric space, and then endow $V_{\mathbf{z}}$ with
the topology induced by the norm $\left\Vert \cdot\right\Vert _{K_{\mathbf{z}}}$\nomenclature{$\left\Vert \cdot\right\Vert _{K_{\mathbf{z}}}$}{ },
defined here as outputting the maximum of the $K_{\mathbf{z}}$-absolute
values of the entries of a given element of $V_{\mathbf{z}}$. This
makes $V_{\mathbf{z}}$ into a finite-dimensional Banach space.
\end{defn}
\begin{defn}
I adopt the convention of writing $\mathcal{F}^{d}$ when I mean that
the $V_{\mathbf{z}}$s are spaces of $d\times1$ column vectors. In
that case, I write $V_{\mathbf{z}}$ as $K_{\mathbf{z}}^{d}$. On
the other hand, when working with $\rho\times c$ matrices (where
$d=\rho c$) instead of $d\times1$ column vectors, I write $\mathcal{F}^{\rho,c}$
instead of $\mathcal{F}$, and say that this frame has dimension $\rho,c$.
For this case, I write $K_{\mathbf{z}}^{\rho,c}$ instead of $V_{\mathbf{z}}$.
\end{defn}
\vphantom{}

The two most important objects associated with a given frame are its
image and the space of compatible functions, defined below. 
\begin{defn}[\textbf{Image and Compatible Functions}]
\label{def:MD Frame terminology}Let $\mathcal{F}$ be a $p$-adic
frame.

\vphantom{}

I. The \textbf{image }of $\mathcal{F}$, denoted $I\left(\mathcal{F}\right)$,
is the \emph{set }defined by: 
\begin{equation}
I\left(\mathcal{F}\right)\overset{\textrm{def}}{=}\bigcup_{\mathbf{z}\in U_{\mathcal{F}}}V_{\mathbf{z}}\label{eq:MD The image of a frame}
\end{equation}
\index{frame!image}

\vphantom{}

II. A function $\chi:U_{\mathcal{F}}\rightarrow I\left(\mathcal{F}\right)$
is said to be \textbf{$\mathcal{F}$-compatible} / \textbf{compatible
}(\textbf{with $\mathcal{F}$}) whenever \index{mathcal{F}-@$\mathcal{F}$-!compatible}\index{frame!compatible functions}
$\chi\left(\mathfrak{z}\right)\in V_{\mathbf{z}}$ for all $\mathbf{z}\in U_{\mathcal{F}}$.
I write $C\left(\mathcal{F}\right)$ to denote the set of all $\mathcal{F}$-compatible
functions.
\end{defn}
\begin{rem}
Note that $C\left(\mathcal{F}\right)$ is a vector space over $\overline{\mathbb{Q}}$
with respect to point-wise addition of functions and scalar multiplication.
\end{rem}
\vphantom{}

In addition to this, we also have the following bits of terminology:
\begin{defn}[\textbf{Frame Terminology}]
Let $\mathcal{F}$ be a $p$-adic frame of dimension $d$ and depth
$r$.

\vphantom{}

I. I call $\mathbb{Z}_{p}^{r}\backslash U_{\mathcal{F}}$ the \textbf{set
of singularities }of the frame. I say that $\mathcal{F}$ is \textbf{non-singular
}whenever $U_{\mathcal{F}}=\mathbb{Z}_{p}^{r}$.\index{frame!non-singular}

\vphantom{}

II. Given a topology $\tau$ (so, $\tau$ could be $\textrm{dis}$,
$\infty$, $\textrm{non}$, or $p$), I writ $U_{\tau}\left(\mathcal{F}\right)$
to denote the set of $\mathbf{z}\in U_{\mathcal{F}}$ so that $K_{\mathbf{z}}$
has been equipped with $\tau$. I call the $U_{\tau}\left(\mathcal{F}\right)$\textbf{
$\tau$-convergence domain }or \textbf{domain of $\tau$ convergence
}of $\mathcal{F}$. In this way, we can speak of the \textbf{discrete
convergence domain}, the \textbf{archimedean convergence domain},
the \textbf{non-archimedean convergence domain}, and the \textbf{$p$-adic
convergence domain} of a given frame $\mathcal{F}$.

\vphantom{}

IV. I say $\mathcal{F}$ is \textbf{simple} \index{frame!simple}if
either $U_{\textrm{non}}\left(\mathcal{F}\right)$ is empty or there
exists a single metrically complete non-archimedean field extension
$K$ of $\overline{\mathbb{Q}}$ so that $K=K_{\mathbf{z}}$ for all
$\mathbf{z}\in U_{\textrm{non}}\left(\mathcal{F}\right)$. (That is
to say, for a simple frame, we use at most one non-archimedean topology.)

\vphantom{}

V. I say $\mathcal{F}$ is \textbf{proper }\index{frame!proper}proper
whenever $K_{\mathbf{z}}$ is not a $p$-adic field for any $\mathbf{z}\in U_{\textrm{non}}\left(\mathcal{F}\right)$. 
\end{defn}
\vphantom{}

UNLESS STATED OTHERWISE, ALL FRAMES ARE ASSUMED TO BE PROPER.

\vphantom{}

Like in the one-dimensional case, only one $p$-adic frame will ever
be used in this dissertation: the standard one.
\begin{defn}
The \textbf{standard $d$-dimensional ($\left(p,q\right)$-adic) frame}\index{frame!standard left(p,qright)-adic@standard $\left(p,q\right)$-adic},
denoted $\mathcal{F}_{p,q}$, is the frame for which the topology
of $\mathbb{C}$ is associated to $\mathbb{N}_{0}^{d}$ and the topology
of $\mathbb{C}_{q}$ is associated to $\left(\mathbb{Z}_{p}^{r}\right)^{\prime}$.
In particular, I write $\mathcal{F}_{p,q}^{d}$ to denote the case
where the topological vector spaces are spaces of $d\times1$ column
vectors; I write $\mathcal{F}_{p,q}^{\rho,c}$ to denote the case
where the topological vector spaces are spaces of $\rho\times c$
column vectors.
\end{defn}
\vphantom{}

Next we have the all-important notion of $\mathcal{F}$-convergence.
\begin{defn}[\textbf{$\mathcal{F}$-convergence}]
Given a frame a $\mathcal{F}$ and a function $\chi\in C\left(\mathcal{F}\right)$,
we say a sequence $\left\{ \chi_{n}\right\} _{n\geq1}$ in $C\left(\mathcal{F}\right)$
\textbf{converges to $\chi$ over $\mathcal{F}$ }(or \textbf{is $\mathcal{F}$-convergent
to $\chi$}) whenever, for each $\mathbf{z}\in U_{\mathcal{F}}$,
we have: 
\begin{equation}
\lim_{n\rightarrow\infty}\chi_{n}\left(\mathbf{z}\right)\overset{V_{\mathbf{z}}}{=}\chi\left(\mathbf{z}\right)
\end{equation}
Note that if $K_{\mathbf{z}}$ has the discrete topology, the limit
implies that $\chi_{n}\left(\mathbf{z}\right)=\chi\left(\mathbf{z}\right)$
for all sufficiently large $n$.\index{mathcal{F}-@$\mathcal{F}$-!convergence}\index{frame!convergence}
On the other hand, if $K_{\mathbf{z}}$ has the topology of a valued
field, we then require:
\begin{equation}
\lim_{n\rightarrow\infty}\left\Vert \chi\left(\mathbf{z}\right)-\chi_{n}\left(\mathbf{z}\right)\right\Vert _{K_{\mathbf{z}}}\overset{\mathbb{R}}{=}0
\end{equation}
I call $\chi$ the \index{mathcal{F}-@$\mathcal{F}$-!limit}\textbf{$\mathcal{F}$-limit
of the $\chi_{n}$s}. \index{frame!limit}More generally, I say $\left\{ \chi_{n}\right\} _{n\geq1}$
is \textbf{$\mathcal{F}$-convergent / converges over $\mathcal{F}$}
whenever there is a function $\chi\in C\left(\mathcal{F}\right)$
so that the $\chi_{n}$s are $\mathcal{F}$-convergent to $\chi$.

In symbols, I denote $\mathcal{F}$-convergence by: 
\begin{equation}
\lim_{n\rightarrow\infty}\chi_{n}\left(\mathbf{z}\right)\overset{\mathcal{F}}{=}\chi\left(\mathbf{z}\right),\textrm{ }\forall\mathbf{z}\in U_{\mathcal{F}}\label{eq:MD Definition-in-symbols of F-convergence}
\end{equation}
or simply:

\begin{equation}
\lim_{n\rightarrow\infty}\chi_{n}\overset{\mathcal{F}}{=}\chi\label{eq:MD Simplified version of expressing f as the F-limit of f_ns}
\end{equation}
\end{defn}
\begin{rem}
In an abuse of notation, we will sometimes write $\mathcal{F}$ convergence
as: 
\begin{equation}
\lim_{n\rightarrow\infty}\left\Vert \chi_{n}\left(\mathbf{z}\right)-\chi\left(\mathbf{z}\right)\right\Vert _{K_{\mathbf{z}}},\textrm{ }\forall\mathfrak{z}\in U_{\mathcal{F}}
\end{equation}
The abuse here is that the norm $\left\Vert \cdot\right\Vert _{K_{\mathbf{z}}}$
will not exist if $\mathbf{z}\in U_{\textrm{dis}}\left(\mathcal{F}\right)$.
As such, for any $\mathbf{z}\in U_{\textrm{dis}}\left(\mathcal{F}\right)$,
I define the expression:
\begin{equation}
\lim_{n\rightarrow\infty}\left\Vert \chi_{n}\left(\mathbf{z}\right)-\chi\left(\mathbf{z}\right)\right\Vert _{K_{\mathbf{z}}}=0
\end{equation}
as meaning that for all sufficiently large $n$, the equality $\chi_{n}\left(\mathbf{z}\right)=\chi_{n+1}\left(\mathbf{z}\right)$
holds in $K_{\mathbf{z}}^{d}$.
\end{rem}
\vphantom{}

Now that we have the language of frames at our disposal, we can begin
to define classes of $\left(p,q\right)$-adic measures which will
be of interest to us.
\begin{defn}[\textbf{Thick} \textbf{$\mathcal{F}$-measures}]
Given a $d$-dimensional depth $r$ $p$-adic frame $\mathcal{F}$,
let $V=\overline{\mathbb{Q}}^{d}$ if $\mathcal{F}=\mathcal{F}^{d}$
and let $V=\overline{\mathbb{Q}}^{\rho,c}$ if $\mathcal{F}=\mathcal{F}^{\rho,c}$.
We write \nomenclature{$M\left(\mathcal{F}\right)$}{$\mathcal{F}$-measures}$M\left(\mathcal{F}\right)$
to denote the set of all functions $\hat{\mathcal{M}}:\hat{\mathbb{Z}}_{p}^{r}\rightarrow V$
so that $\hat{\mathcal{M}}\in B\left(\hat{\mathbb{Z}}_{p}^{r},V_{\mathbf{z}}\right)$
for all $\mathbf{z}\in U_{\textrm{non}}\left(\mathcal{F}\right)$;
that is: 
\begin{equation}
\left\Vert \hat{\mathcal{M}}\right\Vert _{p,K_{\mathbf{z}}}<\infty,\textrm{ }\forall\mathbf{z}\in U_{\textrm{non}}\left(\mathcal{F}\right)\label{eq:Definition of a thick F-measure}
\end{equation}
where:
\[
\left\Vert \hat{\mathcal{M}}\right\Vert _{p,K_{\mathbf{z}}}=\sup_{\mathbf{t}\in\hat{\mathbb{Z}}_{p}^{r}}\left\Vert \hat{\mathcal{M}}\left(\mathbf{t}\right)\right\Vert _{K_{\mathbf{z}}}
\]
where $\left\Vert \cdot\right\Vert _{K_{\mathbf{z}}}$, recall, is
the norm on $V_{\mathbf{z}}$, being the maximum of the $K_{\mathbf{z}}$-absolute-values
of the entries of $\hat{\mathcal{M}}\left(\mathbf{t}\right)$.

\vphantom{}

Next observe that for every $\mathbf{z}\in U_{\textrm{non}}\left(\mathcal{F}\right)$,
the map: 
\begin{equation}
\mathbf{f}\in C\left(\mathbb{Z}_{p}^{r},K_{\mathbf{z}}\right)\mapsto\sum_{\mathbf{t}\in\hat{\mathbb{Z}}_{p}^{r}}\hat{\mathcal{M}}\left(\mathbf{t}\right)\hat{\mathbf{f}}\left(-\mathbf{t}\right)\in K_{\mathbf{z}}\label{eq:Definition of the action of a thick F-measure}
\end{equation}
then defines an element of $C\left(\mathbb{Z}_{p}^{r},K_{\mathbf{z}}^{d}\right)^{*}$,
the space of all continuous $K_{\mathbf{z}}^{d}$-valued linear maps
on the space of continuous functions $\mathbf{f}:\mathbb{Z}_{p}^{r}\rightarrow K_{\mathbf{z}}^{d}$.
As such, we will identify $M\left(\mathcal{F}\right)$ with elements
of $\bigcap_{\mathbf{z}\in U_{\textrm{non}}\left(\mathcal{F}\right)}C\left(\mathbb{Z}_{p}^{r},K_{\mathbf{z}}\right)^{*}$,
and refer to elements of $M\left(\mathcal{F}\right)$ as \textbf{thick}
\textbf{$\mathcal{F}$-measures}\index{mathcal{F}-@$\mathcal{F}$-!thick measure}\index{measure!thick}.

Note that the specific type of multiplication signified by the juxtaposition
of $\hat{\mathcal{M}}\left(\mathbf{t}\right)$ and $\hat{\mathbf{f}}\left(-\mathbf{t}\right)$
in (\ref{eq:Definition of the action of a thick F-measure}) depends
on whether $\hat{\mathcal{M}}$ is of vector type or matrix type.
We write $M\left(\mathcal{F}^{d}\right)$\nomenclature{$M\left(\mathcal{F}^{d}\right)$}{ }
and $M\left(\mathcal{F}^{\rho,c}\right)$\nomenclature{$M\left(\mathcal{F}^{\rho,c}\right)$}{ }
to denote the vector type and matrix type cases, respectively. Regardless,
we then identify $\hat{\mathcal{M}}$ with the Fourier-Stieltjes transform
of a thick measure $\mathcal{M}$. Thus, the statement ``$\mathcal{M}\in M\left(\mathcal{F}^{d}\right)$''
means that, for all $\mathbf{z}\in U_{\textrm{non}}\left(\mathcal{F}\right)$,
$\mathcal{M}$ is an element of $C\left(\mathbb{Z}_{p}^{r},K_{\mathbf{z}}\right)^{*}$
whose action on a function is given by $\sum_{\mathbf{t}\in\hat{\mathbb{Z}}_{p}^{r}}\hat{\mathcal{M}}\left(\mathbf{t}\right)\hat{\mathbf{f}}\left(-\mathbf{t}\right)$,
where $\hat{\mathcal{M}}\in B\left(\hat{\mathbb{Z}}_{p}^{r},K_{\mathbf{z}}^{d}\right)$.
``$\mathcal{M}\in M\left(\mathcal{F}^{\rho,c}\right)$'', meanwhile,
means that, for all $\mathbf{z}\in U_{\textrm{non}}\left(\mathcal{F}\right)$,
$\mathcal{M}$ is an element of $C\left(\mathbb{Z}_{p}^{r},K_{\mathbf{z}}\right)^{*}$
whose action on a function is given by $\sum_{\mathbf{t}\in\hat{\mathbb{Z}}_{p}^{r}}\hat{\mathcal{M}}\left(\mathbf{t}\right)\hat{\mathbf{f}}\left(-\mathbf{t}\right)$,
where $\hat{\mathcal{M}}\in B\left(\hat{\mathbb{Z}}_{p}^{r},K_{\mathbf{z}}^{d,d}\right)$.
Both cases force $\left\Vert \hat{\mathcal{M}}\right\Vert _{p,K_{\mathbf{z}}}$
to be finite for all $\mathbf{z}\in U_{\textrm{non}}\left(\mathcal{F}\right)$.
\end{defn}
\begin{prop}
Let $\hat{\mathcal{M}}\in M\left(\mathcal{F}\right)$. Then, the $N$th
partial sum of the Fourier series generated by $\hat{\mathcal{M}}$:
\begin{equation}
\tilde{\mathcal{M}}_{N}\left(\mathbf{z}\right)\overset{\textrm{def}}{=}\sum_{\left\Vert \mathbf{t}\right\Vert _{p}\leq p^{N}}\hat{\mathcal{M}}\left(\mathbf{t}\right)e^{2\pi i\left\{ \mathbf{t}\mathbf{z}\right\} _{p}},\textrm{ }\forall\mathbf{z}\in\mathbb{Z}_{p}^{r}
\end{equation}
is an element of $C\left(\mathcal{F}\right)$. 
\end{prop}
Proof: A straight-forward multi-dimensional analogue of the one-dimensional
case.

Q.E.D.
\begin{defn}[\textbf{$\mathcal{F}$-rising thick measure}]
Given a $p$-adic frame $\mathcal{F}$ of dimension $d$ and depth
$r$, we say $\mathcal{M}\in M\left(\mathcal{F}\right)$ \index{thick measure!mathcal{F}-rising@$\mathcal{F}$-rising}
\textbf{rises over $\mathcal{F}$} / is \textbf{$\mathcal{F}$-rising}
whenever the sequence $\left\{ \tilde{\mathcal{M}}_{N}\right\} _{N\geq0}$
is $\mathcal{F}$-convergent. That is:
\begin{equation}
\tilde{\mathcal{M}}\left(\mathbf{z}\right)\overset{\mathcal{F}}{=}\lim_{N\rightarrow\infty}\tilde{\mathcal{M}}_{N}\left(\mathbf{z}\right),\textrm{ }\forall\mathbf{z}\in U_{\mathcal{F}}
\end{equation}
which, recall, means:

\vphantom{}

I. For every $\mathbf{z}\in U_{\textrm{dis}}\left(\mathcal{F}\right)$,
$\tilde{\mathcal{M}}_{N}\left(\mathbf{z}\right)=\tilde{\mathcal{M}}_{N+1}\left(\mathbf{z}\right)$
for all sufficiently large $N$.

\vphantom{}

II. For every $\mathbf{z}\in U_{\textrm{arch}}\left(\mathcal{F}\right)\cup U_{\textrm{non}}\left(\mathcal{F}\right)$,
$\lim_{N\rightarrow\infty}\tilde{\mathcal{M}}_{N}\left(\mathbf{z}\right)$
converges to a limit in the topology of $V_{\mathbf{z}}$.

\vphantom{}

We write $\mathcal{M}_{\textrm{rise}}\left(\mathcal{F}\right)$ to
denote the space of thick $\mathcal{F}$-rising measures, with $\mathcal{M}_{\textrm{rise}}\left(\mathcal{F}^{d}\right)$
and $\mathcal{M}_{\textrm{rise}}\left(\mathcal{F}^{\rho,c}\right)$
being used to remind the reader that the thick measures in question
are of vector type and matrix type, respectively.

Additionally, we call a thick measure \textbf{rising-continuous }whenever
it rises over the standard frame. \index{thick measure!rising-continuous} 
\end{defn}
\begin{defn}[\textbf{Kernel of an $\mathcal{F}$-rising measure}]
Given a $p$-adic frame $\mathcal{F}$ of dimension $d$ and depth
$r$ along with an $\mathcal{F}$-rising thick measure $\mathcal{M}$,
the\index{mathcal{F}-@$\mathcal{F}$-!kernel} \textbf{kernel }of $\mathcal{M}$\textemdash denoted\index{thick measure!mathcal{F}-kernel@$\mathcal{F}$-kernel}
\nomenclature{$\tilde{\mathcal{M}}\left(\mathbf{z}\right)$}{kernel of $\mathcal{M}$}\textemdash is
the vector- or matrix-valued function on $U_{\mathcal{F}}$ defined
by the $\mathcal{F}$-limit of the $\tilde{\mathcal{M}}_{N}$s: 
\begin{equation}
\tilde{\mathcal{M}}\left(\mathbf{z}\right)\overset{\mathcal{F}}{=}\lim_{N\rightarrow\infty}\tilde{\mathcal{M}}_{N}\left(\mathbf{z}\right),\textrm{ }\forall\mathbf{z}\in U_{\mathcal{F}}\label{eq:Definition of M-twiddle / the kernel of M}
\end{equation}
We write $\left(\tilde{\mathcal{M}}\right)_{N}\left(\mathbf{z}\right)$\nomenclature{$\left(\tilde{\mathcal{M}}\right)_{N}\left(\mathbf{z}\right)$}{$N$th truncation of the kernel of $\mathcal{M}$}
to denote the $N$th truncation of the kernel of $\mathcal{M}$: 
\begin{equation}
\left(\tilde{\mathcal{M}}\right)_{N}\left(\mathbf{z}\right)\overset{\textrm{def}}{=}\sum_{\mathbf{n}=\mathbf{0}}^{p^{N}-1}\tilde{\mathcal{M}}\left(\mathbf{n}\right)\left[\mathbf{z}\overset{p^{N}}{\equiv}\mathbf{n}\right]\label{eq:Definition of the Nth truncation of the full kernel of a rising multiplier}
\end{equation}
More generally, we write:
\begin{equation}
\left(\tilde{\mathcal{M}}_{M}\right)_{N}\left(\mathbf{z}\right)\overset{\textrm{def}}{=}\sum_{\mathbf{n}=\mathbf{0}}^{p^{N}-1}\tilde{\mathcal{M}}_{M}\left(\mathbf{n}\right)\left[\mathbf{z}\overset{p^{N}}{\equiv}\mathbf{n}\right]\label{eq:Nth truncation of the Mth partial Kernel of M}
\end{equation}
to denote the $N$th truncation of $\tilde{\mathcal{M}}_{M}$. 
\end{defn}
\begin{rem}
In this terminology, observe that for $\mathbf{f}\in C\left(\mathbb{Z}_{p}^{r},\mathbb{C}_{q}^{d}\right)$,
we can express the image of $\mathbf{f}$ under $\mathcal{M}$ as:
\begin{equation}
\mathcal{M}\left(\mathbf{f}\right)\overset{\mathbb{C}_{q}^{d}}{=}\lim_{N\rightarrow\infty}\int_{\mathbb{Z}_{p}^{r}}\tilde{\mathcal{M}}_{N}\left(\mathbf{z}\right)\mathbf{f}\left(\mathbf{z}\right)d\mathbf{z}=\lim_{N\rightarrow\infty}\sum_{\left\Vert \mathbf{t}\right\Vert _{p}\leq p^{N}}\hat{\mathcal{M}}\left(\mathbf{t}\right)\hat{\mathbf{f}}\left(-\mathbf{t}\right)
\end{equation}
\end{rem}
\begin{defn}[\textbf{$\mathcal{F}$-degenerate thick measure}]
For a $p$-adic frame $\mathcal{F}$, a thick $\mathcal{F}$-rising
measure \index{thick measure!mathcal{F}-degenerate@$\mathcal{F}$-degenerate}$\mathcal{M}$
is said to be \textbf{($\mathcal{F}$-)degenerate }whenever its kernel
is identically zero. We write $M_{\textrm{dgen}}\left(\mathcal{F}\right)$
to denote the set of $\mathcal{F}$-degenerate thick measures. We
write this as $M_{\textrm{dgen}}\left(\mathcal{F}^{d}\right)$ and
$M_{\textrm{dgen}}\left(\mathcal{F}^{\rho,c}\right)$ when we wish
to single out vector-type and matrix-type thick measures, respectively.
\end{defn}

\subsection{\label{subsec:5.4.2 Quasi-Integrability-Revisited}Quasi-Integrability
and Thick Measures}
\begin{rem}
Like in \ref{subsec:3.3.5 Quasi-Integrability}, the purpose of Subsection
\ref{subsec:5.4.3 -adic-Wiener-Tauberian} is to provide a solid foundation
for understanding and discussing quasi-integrability in its multi-dimensional
incarnation. As it regards Chapter \ref{chap:6 A-Study-of}'s analysis
of multi-dimensional $\chi_{H}$, the main thing the reader needs
to know is that I say a function $\chi:\mathbb{Z}_{p}^{r}\rightarrow\mathbb{C}_{q}^{d}$
is \textbf{quasi-integrable }with\index{quasi-integrability} respect
to the standard $\left(p,q\right)$-adic frame whenever there exists
a function $\hat{\chi}:\hat{\mathbb{Z}}_{p}^{r}\rightarrow\overline{\mathbb{Q}}^{d}$
so that: 
\begin{equation}
\chi\left(\mathbf{z}\right)\overset{\mathcal{F}_{p,q}^{d}}{=}\lim_{N\rightarrow\infty}\sum_{\left\Vert \mathbf{t}\right\Vert _{p}\leq p^{N}}\hat{\chi}\left(\mathbf{t}\right)e^{2\pi i\left\{ \mathbf{t}\mathbf{z}\right\} _{p}},\textrm{ }\forall\mathbf{z}\in\mathbb{Z}_{p}^{r}
\end{equation}
We call the function $\hat{\chi}$ a\emph{ }\textbf{Fourier transform
}of $\chi$. Once again, this Fourier transform is only unique modulo
the Fourier-Stieltjes transform of a thick measure of vector type
which is degenerate with respect to the standard $\left(p,q\right)$-adic
frame. Like in the one-dimensional case, given a Fourier transform
$\hat{\chi}$ of a quasi-integrable function $\chi$, we write $\tilde{\chi}_{N}$
to denote the $N$th partial Fourier series generated by $\hat{\chi}$.
\end{rem}
\vphantom{}

HERE BEGINS THE DISCUSSION OF MULTI-DIMENSIONAL QUASI-INTEGRABILITY

\vphantom{}
\begin{defn}
Consider a $d$-dimensional $p$-adic frame $\mathcal{F}$ of depth
$r$.

\vphantom{}

We say a function $\chi\in C\left(\mathcal{F}\right)$ is\textbf{
quasi-integrable} \textbf{with respect to $\mathcal{F}$} / $\mathcal{F}$\textbf{-quasi-integrable
}whenever $\chi:U_{\mathcal{F}}\rightarrow I\left(\mathcal{F}\right)$
is the kernel of some $\mathcal{F}$-rising measure $\mathcal{M}\in M\left(\mathcal{F}^{d}\right)$
(that is, of \emph{vector }type). In other words:

\vphantom{}

I. There is a vector-type thick measure $\mathcal{M}$ with a Fourier-Stieltjes
transform $\hat{\mathcal{M}}:\hat{\mathbb{Z}}_{p}^{r}\rightarrow\overline{\mathbb{Q}}^{d}$
which is bounded in $\left\Vert \cdot\right\Vert _{p,K_{\mathbf{z}}}$-norm
for all $\mathbf{z}\in U_{\textrm{non}}\left(\mathcal{F}\right)$;

\vphantom{}

II. $\tilde{\mathcal{M}}_{N}$ $\mathcal{F}$-converges to $\chi$
as $N\rightarrow\infty$: 
\begin{equation}
\lim_{N\rightarrow\infty}\left\Vert \chi\left(\mathbf{z}\right)-\tilde{\mathcal{M}}_{N}\left(\mathbf{z}\right)\right\Vert _{K_{\mathbf{z}}}\overset{\mathbb{R}}{=}0,\textrm{ }\forall\mathbf{z}\in U_{\mathcal{F}}\label{eq:MD definition of quasi-integrability}
\end{equation}

\vphantom{}

I call any $\mathcal{M}$ satisfying these properties an\textbf{ $\mathcal{F}$-quasi-integral}\index{mathcal{F}-@$\mathcal{F}$-!quasi-integral}.
I write $\tilde{L}^{1}\left(\mathcal{F}\right)$ to denote the set
of all $\mathcal{F}$-quasi-integrable functions. (Note that $\tilde{L}^{1}\left(\mathcal{F}\right)$
is then a vector space over $\overline{\mathbb{Q}}$.) We write $\tilde{L}^{1}\left(\mathcal{F}^{d}\right)$
when we wish to remind our audience that the elements of $\tilde{L}^{1}\left(\mathcal{F}^{d}\right)$
are $d\times1$-vector-valued functions, rather than scalar valued
functions.

When $\mathcal{F}$ is the standard $d$-dimensional depth $r$ $\left(p,q\right)$-adic
frame, this definition becomes:

\vphantom{}

i. $\chi$ is a function from $\mathbb{Z}_{p}^{r}\rightarrow\mathbb{C}_{q}^{d}$
so that $\chi\left(\mathbf{z}\right)\in\mathbb{C}^{d}$ for all $\mathbf{z}\in\mathbb{N}_{0}^{r}$
and $\chi\left(\mathbf{z}\right)\in\mathbb{C}_{q}^{d}$ for all $\mathbf{z}\in\left(\mathbb{Z}_{p}^{r}\right)^{\prime}$;

\vphantom{}

ii. For each $\mathbf{z}\in\mathbb{N}_{0}^{r}$, as $N\rightarrow\infty$,
$\tilde{\mathcal{M}}_{N}\left(\mathbf{z}\right)$ converges to $\chi\left(\mathbf{z}\right)$
in $\mathbb{C}^{d}$.

\vphantom{}

iii. For each $\mathbf{z}\in\left(\mathbb{Z}_{p}^{r}\right)^{\prime}$,
as $N\rightarrow\infty$, $\tilde{\mathcal{M}}_{N}\left(\mathbf{z}\right)$
converges to $\chi\left(\mathbf{z}\right)$ in $\mathbb{C}_{q}^{d}$.

\vphantom{}

We then write \nomenclature{$\tilde{L}^{1}\left(\mathbb{Z}_{p}^{r},\mathbb{C}_{q}^{d}\right)$}{ }$\tilde{L}^{1}\left(\mathbb{Z}_{p}^{r},\mathbb{C}_{q}^{d}\right)$
and \nomenclature{$\tilde{L}^{1}\left(\mathcal{F}_{p,q}^{d}\right)$}{ }$\tilde{L}^{1}\left(\mathcal{F}_{p,q}^{d}\right)$
to denote the space of all functions $\mathbb{Z}_{p}^{r}\rightarrow\mathbb{C}_{q}^{d}$
which are quasi-integrable with respect to the standard frame. 
\end{defn}
\vphantom{}

Once again, we will end up with a generalization of the Fourier transform
of a function $\chi:\mathbb{Z}_{p}^{r}\rightarrow\mathbb{C}_{q}^{d}$
by treating $\chi$ as the kernel of some thick vector-type measure.
\begin{defn}
Let $\chi$ be $\mathcal{F}$-quasi-integrable. Then, \textbf{a} \textbf{Fourier
Transform} of $\chi$ is a function $\hat{\mathcal{M}}:\hat{\mathbb{Z}}_{p}^{r}\rightarrow\overline{\mathbb{Q}}^{d}$
which is the Fourier-Stieltjes transform of some $\mathcal{F}$-quasi-integral
$\mathcal{M}$ of $\chi$. We write these functions as $\hat{\chi}$
and write the associated thick measure as $\chi\left(\mathbf{z}\right)d\mathbf{z}$.
Given $\mathbf{f}\in C\left(\mathbb{Z}_{p}^{r},\mathbb{C}_{q}^{d}\right)$,
we then define: 
\begin{equation}
\int_{\mathbb{Z}_{p}^{r}}\mathbf{f}\left(\mathbf{z}\right)\chi\left(\mathbf{z}\right)d\mathbf{z}\overset{\textrm{def}}{=}\sum_{\mathbf{t}\in\hat{\mathbb{Z}}_{p}^{r}}\hat{\mathbf{f}}\left(-\mathbf{t}\right)\hat{\chi}\left(\mathbf{t}\right)\label{eq:MD Definition of the integral of a quasi-integrable function against a continuous function}
\end{equation}
Note that because $\mathcal{M}$ is a thick measure of vector type,
its Fourier-Stieltjes transform $\hat{\mathcal{M}}\left(\mathbf{t}\right)=\hat{\chi}\left(\mathbf{t}\right)$
is a $d\times1$ column vector with entries in $\overline{\mathbb{Q}}$.
The above integral is then the same thing as the image of $\mathbf{f}$
under $\mathcal{M}$: 
\begin{equation}
\mathcal{M}\left(\mathbf{f}\right)\overset{\mathbb{C}_{q}^{d}}{=}\sum_{\mathbf{t}\in\hat{\mathbb{Z}}_{p}^{r}}\hat{\mathbf{f}}\left(-\mathbf{t}\right)\hat{\mathcal{M}}\left(\mathbf{t}\right)\overset{\mathbb{C}_{q}^{d}}{=}\sum_{\mathbf{t}\in\hat{\mathbb{Z}}_{p}^{r}}\hat{\mathcal{M}}\left(\mathbf{t}\right)\hat{\mathbf{f}}\left(-\mathbf{t}\right)
\end{equation}
where the order of juxtaposition of the $d$-tuples $\hat{\mathbf{f}}$
and $\hat{\mathcal{M}}=\hat{\chi}$ is irrelevant because the juxtaposition
signifies \emph{entry-wise} multiplication\textemdash a commutative
operation.
\end{defn}
\begin{defn}
Given a choice of Fourier transform $\hat{\chi}$ for $\chi\in\tilde{L}^{1}\left(\mathcal{F}^{d}\right)$,
we write: 
\begin{equation}
\tilde{\chi}_{N}\left(\mathbf{z}\right)\overset{\textrm{def}}{=}\sum_{\left\Vert \mathbf{t}\right\Vert _{p}\leq p^{N}}\hat{\chi}\left(\mathbf{t}\right)e^{2\pi i\left\{ \mathbf{t}\mathbf{z}\right\} _{p}},\textrm{ }\forall\mathbf{z}\in\mathbb{Z}_{p}^{r}\label{eq:MD Definition of Chi_N twiddle}
\end{equation}
to denote the \textbf{$N$th partial sum of the Fourier series generated
by $\hat{\chi}$}.
\end{defn}
\vphantom{}

Just like with the one-dimensional case, Fourier transforms of quasi-integrable
functions are only unique modulo the Fourier-Stieltjes transform of
a degenerate thick measure.
\begin{prop}
Let $\chi\in\tilde{L}^{1}\left(\mathcal{F}^{d}\right)$. The quasi-integral
of $\chi$ is unique modulo $M_{\textrm{dgen}}\left(\mathcal{F}^{d}\right)$.
That is, if two thick measures $\mathcal{M}$ and $\mathcal{N}$ are
both quasi-integrals of $\chi$, then the thick measure $\mathcal{M}-\mathcal{N}$
is $\mathcal{F}$-degenerate\index{thick measure!mathcal{F}-degenerate@$\mathcal{F}$-degenerate}\index{mathcal{F}-@$\mathcal{F}$-!degenerate}.
Equivalently, we have an isomorphism of $\overline{\mathbb{Q}}$-linear
spaces: 
\begin{equation}
\tilde{L}^{1}\left(\mathcal{F}^{d}\right)\cong M_{\textrm{rise}}\left(\mathcal{F}^{d}\right)/M_{\textrm{dgen}}\left(\mathcal{F}^{d}\right)\label{eq:MD L1 twiddle isomorphism to M alg / M dgen}
\end{equation}
where the isomorphism associates a given $\chi\in\tilde{L}^{1}\left(\mathcal{F}^{d}\right)$
to the equivalence class of $\mathcal{F}$-rising thick measures which
are quasi-integrals of $\chi$. 
\end{prop}
Proof: Let $\mathcal{F}$, $\chi$, $\mathcal{M}$, and $\mathcal{N}$
be as given. Then, by the definition of what it means to be a quasi-integral
of $\chi$, for the thick vector-type measure $\mathcal{P}\overset{\textrm{def}}{=}\mathcal{M}-\mathcal{N}$,
we have that $\tilde{\mathcal{P}}_{N}\left(\mathbf{z}\right)$ tends
to $\mathbf{0}$ in $K_{\mathbf{z}}^{d}$ for all $\mathbf{z}\in U_{\mathcal{F}}$.
Hence, $\mathcal{P}=\mathcal{M}-\mathcal{N}$ is indeed $\mathcal{F}$-degenerate.

Q.E.D.
\begin{prop}
Let $\chi\in\tilde{L}^{1}\left(\mathcal{F}^{d}\right)$. Then, for
any Fourier transform $\hat{\chi}$ of $\chi$, the difference $\chi_{N}-\tilde{\chi}_{N}$
(where $\chi_{N}$ is the $N$th truncation of $\chi$) $\mathcal{F}$-converges
to $\mathbf{0}$. 
\end{prop}
Proof: Immediate from the definitions of $\hat{\chi}$, $\tilde{\chi}_{N}$,
and $\chi_{N}$.

Q.E.D.
\begin{rem}
Note that all of the main unresolved issues regarding one-dimensional
quasi-integrable functions (convolutions, a criterion for determining
non-quasi-integrability, etc.) also apply to the multi-dimensional
case.
\end{rem}

\subsection{\label{subsec:5.4.3 -adic-Wiener-Tauberian}$\left(p,q\right)$-adic
Wiener Tauberian Theorems Revisited}

Here, we quickly establish the multi-dimensional analogues of the
results of Subsection \ref{subsec:3.3.7 -adic-Wiener-Tauberian}. 
\begin{defn}
For \nomenclature{$\tau_{\mathbf{s}}$}{Multi-dimensional translation operator}$\mathbf{s}\in\hat{\mathbb{Z}}_{p}^{r}$,
and for a thick $\left(p,q\right)$-adic measure $\mathcal{M}$ with
Fourier-Stieltjes transform $\hat{\mathcal{M}}$, we write: 
\begin{equation}
\tau_{\mathbf{s}}\left\{ \hat{\mathcal{M}}\right\} \left(\mathbf{t}\right)\overset{\textrm{def}}{=}\hat{\mathcal{M}}_{\mathbf{t}+\mathbf{s}},\textrm{ }\forall\mathbf{t}\in\hat{\mathbb{Z}}_{p}^{r}\label{eq:definition of the translate of M hat}
\end{equation}
We also write: 
\begin{equation}
\tau_{\mathbf{s}}\left\{ \hat{\mathbf{f}}\right\} \left(\mathbf{t}\right)\overset{\textrm{def}}{=}\hat{\mathbf{f}}\left(\mathbf{t}+\mathbf{s}\right),\textrm{ }\forall\mathbf{t}\in\hat{\mathbb{Z}}_{p}^{r}\label{eq:MD Definition of the translate of f hat}
\end{equation}
\end{defn}
\begin{defn}
We write $\left(\mathbb{C}_{q}^{d}\right)^{\times}$\nomenclature{$\left(\mathbb{C}_{q}^{d}\right)^{\times}$}{$d\times1$ vectors whose entries are non-zero $q$-adic complex numbers }
to denote the set of all $d\times1$ vectors whose entries are non-zero
$q$-adic complex numbers. We make $\left(\mathbb{C}_{q}^{d}\right)^{\times}$
an abelian group with respect to the operation of entry-wise multiplication;
i.e.: 
\begin{equation}
\left(\begin{array}{c}
\mathfrak{a}_{1}\\
\mathfrak{a}_{2}
\end{array}\right)\left(\begin{array}{c}
\mathfrak{b}_{1}\\
\mathfrak{b}_{2}
\end{array}\right)=\left(\begin{array}{c}
\mathfrak{a}_{1}\mathfrak{b}_{1}\\
\mathfrak{a}_{2}\mathfrak{b}_{2}
\end{array}\right)
\end{equation}
We write \nomenclature{$\mathbf{1}$}{$d\times 1$ vector whose every entry is $1$}$\mathbf{1}$
to denote the identity element of $\left(\mathbb{C}_{q}^{d}\right)^{\times}$,
the $d\times1$ vector whose entries are all $1$s. We also write
$\vec{\mathbf{1}}_{\mathbf{0}}\left(\mathbf{t}\right)$\nomenclature{$\vec{\mathbf{1}}_{\mathbf{0}}\left(\mathbf{t}\right)$}{function which is $\mathbf{1}$ when $\mathbf{t}=\mathbf{0}$ and which is $\mathbf{0}$ otherwise.}
to denote the function $\hat{\mathbb{Z}}_{p}^{r}\rightarrow\mathbb{C}_{q}^{d}$
which is equal to $\mathbf{1}$ when $\mathbf{t}=\mathbf{0}$ and
which is equal to $\mathbf{0}$ for all other $\mathbf{t}$. Note
that for any $\hat{\mathbf{f}}:\hat{\mathbb{Z}}_{p}^{r}\rightarrow\mathbb{C}_{q}^{d}$:
\begin{equation}
\left(\hat{\mathbf{f}}*\vec{\mathbf{1}}_{\mathbf{0}}\right)\left(\mathbf{t}\right)=\sum_{\mathbf{s}\in\hat{\mathbb{Z}}_{p}^{r}}\hat{\mathbf{f}}\left(\mathbf{t}-\mathbf{s}\right)\vec{\mathbf{1}}_{\mathbf{0}}\left(\mathbf{s}\right)=\hat{\mathbf{f}}\left(\mathbf{t}\right),\textrm{ }\forall\mathbf{t}\in\hat{\mathbb{Z}}_{p}^{r}
\end{equation}
so $\vec{\mathbf{1}}_{\mathbf{0}}\left(\mathbf{t}\right)$ is the
identity element with respect to the convolution operation on $B\left(\hat{\mathbb{Z}}_{p}^{r},\mathbb{C}_{q}^{d}\right)$.
Incidentally, this also proves that the constant function $\mathbf{z}\mapsto\mathbf{1}$
and the function $\mathbf{t}\mapsto\vec{\mathbf{1}}_{\mathbf{0}}\left(\mathbf{t}\right)$
are a Fourier transform pair. More generally, we write $\vec{\mathbf{1}}_{\mathbf{s}}\left(\mathbf{t}\right)$
to denote $\vec{\mathbf{1}}_{\mathbf{0}}\left(\mathbf{t}-\mathbf{s}\right)$,
which is $\mathbf{1}$ for $\mathbf{t}\overset{\mathbf{1}}{\equiv}\mathbf{s}$
and $\mathbf{0}$ otherwise. The reason for introducing the vector
notation is to avoid confusing the vector-valued indicator function
for the point $\left\{ \mathbf{s}\right\} $ ($\vec{\mathbf{1}}_{\mathbf{s}}\left(\mathbf{t}\right)$)
and the \emph{scalar}-valued indicator function for the point $\left\{ \mathbf{s}\right\} $
($\mathbf{1}_{\mathbf{s}}\left(\mathbf{t}\right)$).
\end{defn}
\vphantom{}

As in the one-dimensional case, we establish two WTTs: one for continuous
functions, and another for thick measures.
\begin{thm}[\textbf{Multi-Dimensional Wiener Tauberian Theorem for Continuous
$\left(p,q\right)$-adic Functions}]
\label{thm:MD pq-adic WTT for continuous functions}Let\index{Wiener!Tauberian Theorem!left(p,qright)-adic@$\left(p,q\right)$-adic}
$\mathbf{f}=\left(f_{1}\left(\mathbf{z}\right),\ldots,f_{d}\left(\mathbf{z}\right)\right)\in C\left(\mathbb{Z}_{p}^{r},\mathbb{C}_{q}^{d}\right)$,
and let: 
\begin{equation}
\mathbf{f}^{-1}\left(\mathbf{z}\right)=\left(\frac{1}{f_{1}\left(\mathbf{z}\right)},\ldots,\frac{1}{f_{d}\left(\mathbf{z}\right)}\right)
\end{equation}
denote the multiplicative inverse\footnote{If $\mathbf{F}$ is matrix-valued, then $\mathbf{F}^{-1}$ is the
matrix inverse; if $\mathbf{F}$ is scalar-valued, then $\mathbf{F}^{-1}$
is the reciprocal of $\mathbf{F}$; if $\mathbf{F}$ is vector-valued,
then $\mathbf{F}^{-1}$ is the vector whose entries are the reciprocals
of the entries of $\mathbf{F}$.} of $\mathbf{f}$, if it exists. Then, the following are equivalent:

\vphantom{}

I. $\mathbf{f}^{-1}$ exists and is an element of $C\left(\mathbb{Z}_{p}^{r},\mathbb{C}_{q}^{d}\right)$;

\vphantom{}

II. $\hat{\mathbf{f}}$ has a convolution inverse in $c_{0}\left(\hat{\mathbb{Z}}_{p}^{r},\mathbb{C}_{q}^{d}\right)$;

\vphantom{}

III. $\textrm{span}_{\mathbb{C}_{q}}\left\{ \tau_{\mathbf{s}}\left\{ \hat{\mathbf{f}}\right\} \left(\mathbf{t}\right):\mathbf{s}\in\hat{\mathbb{Z}}_{p}^{r}\right\} $
is dense in $c_{0}\left(\hat{\mathbb{Z}}_{p}^{r},\mathbb{C}_{q}^{d}\right)$;

\vphantom{}

IV. $\mathbf{f}\left(\mathbf{z}\right)\in\left(\mathbb{C}_{q}^{d}\right)^{\times}$
for all $\mathbf{z}\in\mathbb{Z}_{p}^{r}$, . 
\end{thm}
Proof:

\textbullet{} ($\textrm{I}\Rightarrow\textrm{II}$) Suppose $\mathbf{f}^{-1}$
exists and is continuous. Because the $\left(p,q_{H}\right)$-adic
Fourier transform isometrically and isomorphically maps the Banach
algebra $C\left(\mathbb{Z}_{p}^{r},\mathbb{C}_{q}^{d}\right)$ onto
the Banach algebra $c_{0}\left(\hat{\mathbb{Z}}_{p}^{r},\mathbb{C}_{q}^{d}\right)$,
we can write: 
\begin{align*}
\mathbf{f}\left(\mathbf{z}\right)\cdot\mathbf{f}^{-1}\left(\mathbf{z}\right) & =\mathbf{1},\textrm{ }\forall\mathbf{z}\in\mathbb{Z}_{p}^{r}\\
\left(\textrm{Fourier transform}\right); & \Updownarrow\\
\left(\hat{\mathbf{f}}*\widehat{\mathbf{f}^{-1}}\right)\left(\mathbf{t}\right) & =\vec{\mathbf{1}}_{\mathbf{0}}\left(\mathbf{t}\right),\textrm{ }\forall\mathbf{t}\in\hat{\mathbb{Z}}_{p}^{r}
\end{align*}
where both $\hat{\mathbf{f}}$ and $\widehat{\mathbf{f}^{-1}}$ are
in $c_{0}$. This shows that $\widehat{\mathbf{f}^{-1}}$ is $\hat{\mathbf{f}}^{-1}$\textemdash the
convolution inverse of $\hat{\mathbf{f}}$.

\vphantom{}

\textbullet{} ($\textrm{II}\Rightarrow\textrm{III}$) Suppose $\hat{\mathbf{f}}$
has a convolution inverse $\hat{\mathbf{f}}^{-1}\in c_{0}\left(\hat{\mathbb{Z}}_{p}^{r},\mathbb{C}_{q}^{d}\right)$.
Then, letting $\hat{\mathbf{g}}\in c_{0}\left(\hat{\mathbb{Z}}_{p}^{r},\mathbb{C}_{q}^{d}\right)$
be arbitrary, we can write: 
\begin{equation}
\left(\hat{\mathbf{f}}*\left(\hat{\mathbf{f}}^{-1}*\hat{\mathbf{g}}\right)\right)\left(\mathbf{t}\right)=\left(\left(\hat{\mathbf{f}}*\hat{\mathbf{f}}^{-1}\right)*\hat{\mathbf{g}}\right)\left(\mathbf{t}\right)=\left(\mathbf{1}_{\mathbf{0}}*\hat{\mathbf{g}}\right)\left(\mathbf{t}\right)=\hat{\mathbf{g}}\left(\mathbf{t}\right),\textrm{ }\forall\mathbf{t}\in\hat{\mathbb{Z}}_{p}^{r}
\end{equation}
Letting $\hat{\mathbf{h}}$ denote $\hat{\mathbf{f}}^{-1}*\hat{\mathbf{g}}$,
we then have: 
\begin{equation}
\hat{\mathbf{g}}\left(\mathbf{t}\right)=\left(\hat{\mathbf{f}}*\hat{\mathbf{h}}\right)\left(\mathbf{t}\right)\overset{\mathbb{C}_{q}^{d}}{=}\lim_{N\rightarrow\infty}\sum_{\left\Vert \mathbf{s}\right\Vert _{p}\leq p^{N}}\hat{\mathbf{h}}\left(\mathbf{s}\right)\hat{\mathbf{f}}\left(\mathbf{t}-\mathbf{s}\right)
\end{equation}
Since $\hat{\mathbf{f}}$ and $\hat{\mathbf{h}}$ are in $c_{0}$,
note that: 
\begin{align*}
\sup_{\mathbf{t}\in\hat{\mathbb{Z}}_{p}^{r}}\left\Vert \sum_{\mathbf{s}\in\hat{\mathbb{Z}}_{p}^{r}}\hat{\mathbf{h}}\left(\mathbf{s}\right)\hat{\mathbf{f}}\left(\mathbf{t}-\mathbf{s}\right)-\sum_{\left\Vert \mathbf{s}\right\Vert _{p}\leq p^{N}}\hat{\mathbf{h}}\left(\mathbf{s}\right)\hat{\mathbf{f}}\left(\mathbf{t}-\mathbf{s}\right)\right\Vert _{q} & \leq\sup_{\mathbf{t}\in\hat{\mathbb{Z}}_{p}^{r}}\sup_{\left\Vert \mathbf{s}\right\Vert _{p}>p^{N}}\left\Vert \hat{\mathbf{h}}\left(\mathbf{s}\right)\hat{\mathbf{f}}\left(\mathbf{t}-\mathbf{s}\right)\right\Vert _{q}\\
\left(\left\Vert \hat{\mathbf{f}}\right\Vert _{p,q}<\infty\right); & \leq\sup_{\left\Vert \mathbf{s}\right\Vert _{p}>p^{N}}\left\Vert \hat{\mathbf{h}}\left(\mathbf{s}\right)\right\Vert _{q}
\end{align*}
and hence: 
\[
\lim_{N\rightarrow\infty}\sup_{\mathbf{t}\in\hat{\mathbb{Z}}_{p}^{r}}\left\Vert \sum_{\mathbf{s}\in\hat{\mathbb{Z}}_{p}^{r}}\hat{\mathbf{h}}\left(\mathbf{s}\right)\hat{\mathbf{f}}\left(\mathbf{t}-\mathbf{s}\right)-\sum_{\left\Vert \mathbf{s}\right\Vert _{p}\leq p^{N}}\hat{\mathbf{h}}\left(\mathbf{s}\right)\hat{\mathbf{f}}\left(\mathbf{t}-\mathbf{s}\right)\right\Vert _{q}
\]
is equal to $\lim_{N\rightarrow\infty}\sup_{\left\Vert \mathbf{s}\right\Vert _{p}>p^{N}}\left\Vert \hat{\mathbf{h}}\left(\mathbf{s}\right)\right\Vert _{q}$,
which tends to $0$. This shows that the $q$-adic convergence of
$\sum_{\left\Vert \mathbf{s}\right\Vert _{p}>p^{N}}\hat{\mathbf{h}}\left(\mathbf{s}\right)\hat{\mathbf{f}}\left(\mathbf{t}-\mathbf{s}\right)$
to $\hat{\mathbf{g}}\left(\mathbf{t}\right)=\sum_{\mathbf{s}\in\hat{\mathbb{Z}}_{p}^{r}}\hat{\mathbf{h}}\left(\mathbf{s}\right)\hat{\mathbf{f}}\left(\mathbf{t}-\mathbf{s}\right)$
is uniform in $\mathbf{t}$. Hence: 
\begin{equation}
\lim_{N\rightarrow\infty}\sup_{\mathbf{t}\in\hat{\mathbb{Z}}_{p}^{r}}\left\Vert \hat{\mathbf{g}}\left(\mathbf{t}\right)-\sum_{\left\Vert \mathbf{s}\right\Vert _{p}\leq p^{N}}\hat{\mathbf{h}}\left(\mathbf{s}\right)\hat{\mathbf{f}}\left(\mathbf{t}-\mathbf{s}\right)\right\Vert _{q}=0
\end{equation}
which is precisely the definition of what it means for the sequence:
\begin{equation}
\left\{ \sum_{\left\Vert \mathbf{s}\right\Vert _{p}\leq p^{N}}\hat{\mathbf{h}}\left(\mathbf{s}\right)\hat{\mathbf{f}}\left(\mathbf{t}-\mathbf{s}\right)\right\} _{N\geq0}
\end{equation}
to converge in $c_{0}\left(\hat{\mathbb{Z}}_{p}^{r},\mathbb{C}_{q}^{d}\right)$
to $\hat{\mathbf{g}}$. Because our sequence converging to the arbitrary
$\hat{\mathbf{g}}\in c_{0}\left(\hat{\mathbb{Z}}_{p}^{r},\mathbb{C}_{q}^{d}\right)$
is an element of: 
\begin{equation}
\textrm{span}_{\mathbb{C}_{q}}\left\{ \tau_{\mathbf{s}}\left\{ \hat{\mathbf{f}}\right\} \left(\mathbf{t}\right):\mathbf{s}\in\hat{\mathbb{Z}}_{p}^{r}\right\} 
\end{equation}
we have proven the density of $\textrm{span}_{\mathbb{C}_{q}}\left\{ \tau_{\mathbf{s}}\left\{ \hat{\mathbf{f}}\right\} \left(\mathbf{t}\right):\mathbf{s}\in\hat{\mathbb{Z}}_{p}^{r}\right\} $
in $c_{0}\left(\hat{\mathbb{Z}}_{p}^{r},\mathbb{C}_{q}^{d}\right)$.

\vphantom{}

\textbullet{} ($\textrm{III}\Rightarrow\textrm{IV}$) Suppose $\textrm{span}_{\mathbb{C}_{q}}\left\{ \tau_{\mathbf{s}}\left\{ \hat{\mathbf{f}}\right\} \left(\mathbf{t}\right):\mathbf{s}\in\hat{\mathbb{Z}}_{p}^{r}\right\} $
is dense in $c_{0}\left(\hat{\mathbb{Z}}_{p}^{r},\mathbb{C}_{q}^{d}\right)$.
Now, let $\epsilon>0$. Because $\vec{\mathbf{1}}_{\mathbf{0}}\left(\mathbf{t}\right)\in c_{0}\left(\hat{\mathbb{Z}}_{p}^{r},\mathbb{C}_{q}^{d}\right)$,
the assumed density of the span of translates guarantees the existence
of constant vectors $\mathbf{c}_{1},\ldots,\mathbf{c}_{M}\in\mathbb{C}_{q}^{d}$
and points $\mathbf{t}_{1},\ldots,\mathbf{t}_{M}\in\hat{\mathbb{Z}}_{p}^{r}$
so that: 
\begin{equation}
\sup_{\mathbf{t}\in\hat{\mathbb{Z}}_{p}^{r}}\left\Vert \vec{\mathbf{1}}_{\mathbf{0}}\left(\mathbf{t}\right)-\sum_{m=1}^{M}\mathbf{c}_{m}\hat{\mathbf{f}}\left(\mathbf{t}-\mathbf{t}_{m}\right)\right\Vert _{q}<\epsilon
\end{equation}
Here, note that $\mathbf{c}_{m}\hat{\mathbf{f}}\left(\mathbf{t}-\mathbf{t}_{m}\right)$
is the entry-wise product of the $d\times1$ vectors $\mathbf{c}_{m}$
and $\hat{\mathbf{f}}\left(\mathbf{t}-\mathbf{t}_{m}\right)$.

Next, letting $N\geq\max\left\{ -v_{p}\left(\mathbf{t}_{1}\right),\ldots,-v_{p}\left(\mathbf{t}_{M}\right)\right\} $
be arbitrary\footnote{Note that this choice of $N$ is larger than the negative $p$-adic
valuations of all of the entries of all of the $\mathbf{t}_{m}$s.}, we have that for all $m$, the map $\mathbf{t}\mapsto\mathbf{t}+\mathbf{t}_{m}$
is a bijection of $\left\{ \mathbf{t}\in\hat{\mathbb{Z}}_{p}^{r}:\left\Vert \mathbf{t}\right\Vert _{p}\leq p^{N}\right\} $.
Consequently: 
\begin{align*}
\sum_{\left\Vert \mathbf{t}\right\Vert _{p}\leq p^{N}}\left(\sum_{m=1}^{M}\mathbf{c}_{m}\hat{\mathbf{f}}\left(\mathbf{t}-\mathbf{t}_{m}\right)\right)e^{2\pi i\left\{ \mathbf{t}\mathbf{z}\right\} _{p}} & =\sum_{m=1}^{M}\mathbf{c}_{m}\sum_{\left\Vert \mathbf{t}\right\Vert _{p}\leq p^{N}}\hat{\mathbf{f}}\left(\mathbf{t}\right)e^{2\pi i\left\{ \left(\mathbf{t}+\mathbf{t}_{m}\right)\mathbf{z}\right\} _{p}}\\
 & =\left(\sum_{m=1}^{M}\mathbf{c}_{m}e^{2\pi i\left\{ \mathbf{t}_{m}\mathbf{z}\right\} _{p}}\right)\sum_{\left\Vert \mathbf{t}\right\Vert _{p}\leq p^{N}}\hat{\mathbf{f}}\left(\mathbf{t}\right)e^{2\pi i\left\{ \mathbf{t}\mathbf{z}\right\} _{p}}
\end{align*}
Since $N$ was arbitrary, we may let it tend to $\infty$. This gives:
\begin{equation}
\lim_{N\rightarrow\infty}\sum_{\left\Vert \mathbf{t}\right\Vert _{p}\leq p^{N}}\left(\sum_{m=1}^{M}\mathbf{c}_{m}\hat{\mathbf{f}}\left(\mathbf{t}-\mathbf{t}_{m}\right)\right)e^{2\pi i\left\{ \mathbf{t}\mathbf{z}\right\} _{p}}\overset{\mathbb{C}_{q}^{d}}{=}\mathbf{g}_{m}\left(\mathbf{z}\right)\mathbf{f}\left(\mathbf{z}\right)
\end{equation}
where $\mathbf{g}_{m}=\left(g_{m,1},\ldots,g_{m,d}\right):\mathbb{Z}_{p}^{r}\rightarrow\mathbb{C}_{q}^{d}$
is defined by: 
\begin{equation}
\mathbf{g}_{m}\left(\mathbf{z}\right)\overset{\textrm{def}}{=}\sum_{m=1}^{M}\mathbf{c}_{m}e^{2\pi i\left\{ \mathbf{t}_{m}\mathbf{z}\right\} _{p}},\textrm{ }\forall\mathbf{z}\in\mathbb{Z}_{p}^{r}\label{eq:MD Definition of G_m}
\end{equation}
Moreover, the convergence of the $N$-limit is uniform in $\mathbf{z}$.

Consequently: 
\begin{align*}
\left\Vert \mathbf{1}-\mathbf{g}_{m}\left(\mathbf{z}\right)\mathbf{f}\left(\mathbf{z}\right)\right\Vert _{q} & \overset{\mathbb{R}}{=}\lim_{N\rightarrow\infty}\left\Vert \sum_{\left\Vert \mathbf{t}\right\Vert _{p}\leq p^{N}}\left(\vec{\mathbf{1}}_{\mathbf{0}}\left(\mathbf{t}\right)-\sum_{m=1}^{M}\mathbf{c}_{m}\hat{\mathbf{f}}\left(\mathbf{t}-\mathbf{t}_{m}\right)\right)e^{2\pi i\left\{ \mathbf{t}\mathbf{z}\right\} _{p}}\right\Vert _{q}\\
 & \leq\sup_{\mathbf{t}\in\hat{\mathbb{Z}}_{p}^{r}}\left\Vert \vec{\mathbf{1}}_{\mathbf{0}}\left(\mathbf{t}\right)-\sum_{m=1}^{M}\mathbf{c}_{m}\hat{\mathbf{f}}\left(\mathbf{t}-\mathbf{t}_{m}\right)\right\Vert _{q}\\
 & <\epsilon
\end{align*}
for all $\mathbf{z}\in\mathbb{Z}_{p}^{r}$.

Finally, by way of contradiction, suppose there is an $\ell\in\left\{ 1,\ldots,d\right\} $
so that $f_{\ell}\left(\mathbf{z}_{0}\right)=0$. Then, the $\ell$th
entry of the $d\times1$ vector $\mathbf{1}-\mathbf{g}_{m}\left(\mathbf{z}_{0}\right)\mathbf{f}\left(\mathbf{z}_{0}\right)$
is: 
\[
1-g_{m,\ell}\left(\mathbf{z}_{0}\right)f_{\ell}\left(\mathbf{z}_{0}\right)=1-g_{m,\ell}\left(\mathbf{z}_{0}\right)\cdot0=1
\]
Because $\left\Vert \mathbf{1}-\mathbf{g}_{m}\left(\mathbf{z}_{0}\right)\mathbf{f}\left(\mathbf{z}_{0}\right)\right\Vert _{q}\geq\left\Vert 1-g_{m,k}\left(\mathbf{z}_{0}\right)f_{k}\left(\mathbf{z}_{0}\right)\right\Vert _{q}$
for all $k\in\left\{ 1,\ldots,d\right\} $, we then have:
\begin{equation}
\epsilon>\left\Vert \mathbf{1}-\mathbf{g}_{m}\left(\mathbf{z}_{0}\right)\mathbf{f}\left(\mathbf{z}_{0}\right)\right\Vert _{q}\geq\left\Vert 1-g_{m,\ell}\left(\mathbf{z}_{0}\right)\cdot0\right\Vert _{q}=1
\end{equation}
However, $\epsilon$ was given to be in $\left(0,1\right)$. This
is a contradiction!

As such, the entries of $\mathbf{f}$ cannot have any zeroes whenever
the span of $\hat{\mathbf{f}}$'s translates are dense in $c_{0}\left(\hat{\mathbb{Z}}_{p}^{r},\mathbb{C}_{q}^{d}\right)$.

\vphantom{}

\textbullet{} ($\textrm{IV}\Rightarrow\textrm{I}$) Suppose $\mathbf{f}\left(\mathbf{z}\right)\in\left(\mathbb{C}_{q}^{d}\right)^{\times}$
for all $\mathbf{z}\in\mathbb{Z}_{p}^{r}$; that is, suppose no entry
of $\mathbf{f}\left(\mathbf{z}\right)$ is zero for any $\mathbf{z}\in\mathbb{Z}_{p}^{r}$.
Because the map which sends $\mathbf{y}=\left(\mathfrak{y}_{1},\ldots,\mathfrak{y}_{d}\right)\in\left(\mathbb{C}_{q}^{d}\right)^{\times}$
to: 
\[
\frac{1}{\mathbf{y}}=\left(\frac{1}{\mathfrak{y}_{1}},\ldots,\frac{1}{\mathfrak{y}_{d}}\right)\in\left(\mathbb{C}_{q}^{d}\right)^{\times}
\]
is continuous on $\left(\mathbb{C}_{q}^{d}\right)^{\times}$, the
continuity of compositions of continuous maps tells us that we can
show $1/\mathbf{f}$ is continuous by establishing that show that
$\inf_{\mathbf{z}\in\mathbb{Z}_{p}^{r}}\left|f_{k}\left(\mathbf{z}\right)\right|_{q}>0$
for each $k\in\left\{ 1,\ldots,d\right\} $.

So, by way of contradiction, suppose there is a $k\in\left\{ 1,\ldots,d\right\} $
so that $\inf_{\mathbf{z}\in\mathbb{Z}_{p}^{r}}\left|f_{k}\left(\mathbf{z}\right)\right|_{q}=0$.
Then, there exists a sequence $\left\{ \mathbf{z}_{n}\right\} _{n\geq0}\subseteq\mathbb{Z}_{p}^{r}$
so that $\lim_{n\rightarrow\infty}\left|f_{k}\left(\mathbf{z}_{n}\right)\right|_{q}=0$.
By the compactness of $\mathbb{Z}_{p}^{r}$, the $\mathbf{z}_{n}$s
have a subsequence $\mathbf{z}_{n_{\ell}}$ which converges in $\mathbb{Z}_{p}^{r}$
to some limit $\mathbf{z}_{\infty}\in\mathbb{Z}_{p}^{r}$. Because
$\mathbf{f}$ is continuous, so is $f_{k}$; this forces $f_{k}\left(\mathbf{z}_{\infty}\right)=0$.
Hence, $\mathbf{z}_{\infty}$ is an element of $\mathbb{Z}_{p}^{r}$
for which $\mathbf{f}\left(\mathbf{z}_{\infty}\right)$ has a non-zero
entry\textemdash which contradicts our given hypothesis on $\mathbf{f}$.

Thus, if none of the entries of $\mathbf{f}$ are ever zero on $\mathbb{Z}_{p}^{r}$,
$\left|f_{k}\left(\mathbf{z}\right)\right|_{q}$ is bounded away from
zero for all $k$, which then shows that $1/\mathbf{f}$ is $\left(p,q\right)$-adically
continuous.

Q.E.D.
\begin{thm}[\textbf{Wiener Tauberian Theorem for Vector-Type Thick $\left(p,q\right)$-adic
Measures}]
\label{thm:MD WTT for pq adic vector type thick measures}Let $\mathcal{M}\in M\left(\mathbb{Z}_{p}^{r},\mathbb{C}_{q}^{d}\right)$
be a $\left(p,q\right)$-adic thick measure of vector type, with Fourier-Stieltjes
transform $\hat{\mathcal{M}}:\hat{\mathbb{Z}}_{p}^{r}\rightarrow\mathbb{C}_{q}^{d}$.
Then, $\textrm{span}_{\mathbb{C}_{q}}\left\{ \tau_{\mathbf{s}}\left\{ \hat{\mathcal{M}}\right\} \left(\mathbf{t}\right):\mathbf{s}\in\hat{\mathbb{Z}}_{p}^{r}\right\} $
is dense in $c_{0}\left(\hat{\mathbb{Z}}_{p}^{r},\mathbb{C}_{q}^{d}\right)$
if and only if for any $\mathbf{z}\in\mathbb{Z}_{p}^{r}$ for which
the limit $\lim_{N\rightarrow\infty}\tilde{\mathcal{M}}_{N}\left(\mathbf{z}\right)$
converges in $\mathbb{C}_{q}^{d}$, said limit is necessarily an element
of $\left(\mathbb{C}_{q}^{d}\right)^{\times}$. \index{Wiener!Tauberian Theorem!left(p,qright)-adic@$\left(p,q\right)$-adic}
\end{thm}
\begin{rem}
\textbf{WARNING \textendash{} }Like in the one-dimensional case (\textbf{Theorem
\ref{thm:pq WTT for measures}}), for the proof of this theorem\textemdash and
only this theorem\textemdash we will need to modify our notation for
the $\left(p,q\right)$-adic norm. Instead of writing $\left\Vert \cdot\right\Vert _{p,q}$
to denote the supremum over $\hat{\mathbb{Z}}_{p}^{r}$ of a vector-valued
function with entries in $\mathbb{C}_{q}$, we write that norm as
$\left\Vert \cdot\right\Vert _{p^{\infty},q}$. 
\end{rem}
Proof:

I. Suppose the span of the translates of $\hat{\mathcal{M}}$ are
dense. Like in the one-dimensional case, we let $\epsilon\in\left(0,1\right)$
and then, using the assumed density, choose constant vectors $\mathbf{c}_{m}\in\mathbb{C}_{q}^{d}$
and $\mathbf{t}_{m}\in\hat{\mathbb{Z}}_{p}^{r}$ so that: 
\begin{equation}
\sup_{\mathbf{t}\in\hat{\mathbb{Z}}_{p}^{r}}\left\Vert \vec{\mathbf{1}}_{\mathbf{0}}\left(\mathbf{t}\right)-\sum_{m=1}^{M}\mathbf{c}_{m}\hat{\mathcal{M}}\left(\mathbf{t}-\mathbf{t}_{m}\right)\right\Vert _{q}<\epsilon
\end{equation}
When $N$ is sufficiently large, we obtain: 
\begin{align*}
\left\Vert \mathbf{1}-\left(\sum_{m=1}^{M}\mathbf{c}_{m}e^{2\pi i\left\{ \mathbf{t}_{m}\mathbf{z}\right\} _{p}}\right)\tilde{\mathcal{M}}_{N}\left(\mathbf{z}\right)\right\Vert _{q} & \leq\max_{\left\Vert \mathbf{t}\right\Vert _{p}\leq p^{N}}\left\Vert \vec{\mathbf{1}}_{\mathbf{0}}\left(\mathbf{t}\right)-\sum_{m=1}^{M}\mathbf{c}_{m}\hat{\mathcal{M}}\left(\mathbf{t}-\mathbf{t}_{m}\right)\right\Vert _{q}\\
 & \leq\sup_{\mathbf{t}\in\hat{\mathbb{Z}}_{p}^{r}}\left\Vert \vec{\mathbf{1}}_{\mathbf{0}}\left(\mathbf{t}\right)-\sum_{m=1}^{M}\mathbf{c}_{m}\hat{\mathcal{M}}\left(\mathbf{t}-\mathbf{t}_{m}\right)\right\Vert _{q}\\
 & <\epsilon
\end{align*}

So, let $\mathbf{z}_{0}\in\mathbb{Z}_{p}^{r}$ be a point so that
$\mathfrak{L}=\lim_{N\rightarrow\infty}\tilde{\mathcal{M}}_{N}\left(\mathbf{z}_{0}\right)$
converges in $\mathbb{C}_{q}^{d}$. Then, plugging $\mathbf{z}=\mathbf{z}_{0}$
in the above ends yields: 
\[
\epsilon>\lim_{N\rightarrow\infty}\left\Vert \mathbf{1}-\left(\sum_{m=1}^{M}\mathbf{c}_{m}e^{2\pi i\left\{ \mathbf{t}_{m}\mathbf{z}_{0}\right\} _{p}}\right)\tilde{\mathcal{M}}_{N}\left(\mathbf{z}_{0}\right)\right\Vert _{q}=\left\Vert \mathbf{1}-\left(\sum_{m=1}^{M}\mathbf{c}_{m}e^{2\pi i\left\{ \mathbf{t}_{m}\mathbf{z}_{0}\right\} _{p}}\right)\cdot\mathfrak{L}\right\Vert _{q}
\]
If $\mathfrak{L}$ had an entry which is $0$ (i.e., if $\mathfrak{L}\in\mathbb{C}_{q}^{d}\backslash\left(\mathbb{C}_{q}^{d}\right)^{\times}$),
then we end up with $\epsilon>1$, which contradicts the fact that
$\epsilon\in\left(0,1\right)$. Thus, if $\lim_{N\rightarrow\infty}\tilde{\mathcal{M}}_{N}\left(\mathbf{z}_{0}\right)$
converges in $\mathbb{C}_{q}^{d}$, it must converge to an element
of $\left(\mathbb{C}_{q}^{d}\right)^{\times}$.

\vphantom{}

II. Let $\mathbf{z}_{0}\in\mathbb{Z}_{p}^{r}$ be such that $\tilde{\mathcal{M}}\left(\mathbf{z}_{0}\right)\in\mathbb{C}_{q}^{d}\backslash\left(\mathbb{C}_{q}^{d}\right)^{\times}$
and that $\tilde{\mathcal{M}}_{N}\left(\mathbf{z}_{0}\right)\rightarrow\tilde{\mathcal{M}}\left(\mathbf{z}_{0}\right)$
in $\mathbb{C}_{q}^{d}$ as $N\rightarrow\infty$; in an abuse of
terminology, let us call this $\mathbf{z}_{0}$ a ``zero'' of $\tilde{\mathcal{M}}$.
Next, by way of contradiction, suppose the span of the translates
of $\hat{\mathcal{M}}$ is dense in $c_{0}\left(\hat{\mathbb{Z}}_{p}^{r},\mathbb{C}_{q}^{d}\right)$
in spite of the zero of $\tilde{\mathcal{M}}$ at $\mathbf{z}_{0}$.
So, letting $\epsilon\in\left(0,1\right)$ be arbitrary, by the assumed
density,, there exists a linear combination of translates of $\hat{\mathcal{M}}$
which approximates $\vec{\mathbf{1}}_{\mathbf{0}}$ in sup norm: 
\[
\sup_{\mathbf{t}\in\hat{\mathbb{Z}}_{p}^{r}}\left\Vert \vec{\mathbf{1}}_{\mathbf{0}}\left(\mathbf{t}\right)-\sum_{k=1}^{K}\mathbf{c}_{k}\hat{\mathcal{M}}\left(\mathbf{t}-\mathbf{t}_{k}\right)\right\Vert _{q}<\epsilon
\]
where the $\mathbf{c}_{k}$s and $\mathbf{t}_{k}$s are constants
in $\mathbb{C}_{q}^{d}$ and $\hat{\mathbb{Z}}_{p}^{r}$, respectively.
Writing: 
\begin{equation}
\hat{\mathbf{h}}_{\epsilon}\left(\mathbf{t}\right)\overset{\textrm{def}}{=}\sum_{k=1}^{K}\mathbf{c}_{k}\vec{\mathbf{1}}_{\mathbf{t}_{k}}\left(\mathbf{t}\right)\label{eq:Definition of bold h_epsilon hat}
\end{equation}
we can represent our approximating linear combination as a convolution:

\[
\left(\hat{\mathcal{M}}*\hat{\mathbf{h}}_{\epsilon}\right)\left(\mathbf{t}\right)=\sum_{\mathbf{u}\in\hat{\mathbb{Z}}_{p}^{r}}\hat{\mathcal{M}}\left(\mathbf{t}-\mathbf{u}\right)\sum_{k=1}^{K}\mathbf{c}_{k}\vec{\mathbf{1}}_{\mathbf{t}_{k}}\left(\mathbf{u}\right)=\sum_{k=1}^{K}\mathbf{c}_{k}\hat{\mathcal{M}}\left(\mathbf{t}-\mathbf{t}_{k}\right)
\]
and so: 
\begin{equation}
\left\Vert \vec{\mathbf{1}}_{\mathbf{0}}-\hat{\mathcal{M}}*\hat{\mathbf{h}}_{\epsilon}\right\Vert _{p^{\infty},q}\overset{\textrm{def}}{=}\sup_{\mathbf{t}\in\hat{\mathbb{Z}}_{p}^{r}}\left\Vert \vec{\mathbf{1}}_{\mathbf{0}}\left(\mathbf{t}\right)-\left(\hat{\mathcal{M}}*\hat{\mathbf{h}}_{\epsilon}\right)\left(\mathbf{t}\right)\right\Vert _{q}<\epsilon\label{eq:MD Converse WTT - eq. 1}
\end{equation}
Before we proceed, exactly like in the one-dimensional case, it is
\emph{vital }to note that we can and \emph{must} assume that $\hat{\mathbf{h}}_{\epsilon}$
is \emph{not identically zero}. In fact, we can assume this must hold
for any $\epsilon\in\left(0,1\right)$. Indeed, if there was an $\epsilon\in\left(0,1\right)$
for which $\hat{\mathbf{h}}_{\epsilon}$ was identically zero, the
supremum $\sup_{\mathbf{t}\in\hat{\mathbb{Z}}_{p}^{r}}\left\Vert \vec{\mathbf{1}}_{\mathbf{0}}\left(\mathbf{t}\right)-\left(\hat{\mathcal{M}}*\hat{\mathbf{h}}_{\epsilon}\right)\left(\mathbf{t}\right)\right\Vert _{q}$
would be equal to $1$, which is quite impossible, given that the
supremum is, by construction, less than $\epsilon$.

Exactly like before, we will obtain a contradiction by showing that
we can convolve $\hat{\mathcal{M}}$ by two different functions; this
first makes the resultant vector-valued function close to $\vec{\mathbf{1}}_{\mathbf{0}}$;
the second makes one of the entries of the resultant vector-valued
function close to $0$.

Our close-to-zero convolving function is going to be: 
\begin{equation}
\hat{\phi}_{N}\left(\mathbf{t}\right)\overset{\textrm{def}}{=}\vec{\mathbf{1}}_{\mathbf{0}}\left(p^{N}\mathbf{t}\right)e^{-2\pi i\left\{ \mathbf{t}\mathbf{z}_{0}\right\} _{p}}\label{eq:MD Definition of Phi_N hat}
\end{equation}
Note that we can then write: 
\begin{align*}
\left(\hat{\mathcal{M}}*\hat{\phi}_{N}\right)\left(\mathbf{u}\right) & =\sum_{\mathbf{s}\in\hat{\mathbb{Z}}_{p}^{r}}\hat{\mathcal{M}}\left(\mathbf{u}-\mathbf{s}\right)\hat{\phi}_{N}\left(\mathbf{s}\right)\\
\left(\hat{\phi}_{N}\left(\mathbf{s}\right)=\mathbf{0},\textrm{ }\forall\left\Vert \mathbf{s}\right\Vert _{p}>p^{N}\right); & =\sum_{\left\Vert \mathbf{s}\right\Vert _{p}\leq p^{N}}\hat{\mathcal{M}}\left(\mathbf{u}-\mathbf{s}\right)\hat{\phi}_{N}\left(\mathbf{s}\right)\\
 & =\sum_{\left\Vert \mathbf{s}\right\Vert _{p}\leq p^{N}}\hat{\mathcal{M}}\left(\mathbf{u}-\mathbf{s}\right)e^{-2\pi i\left\{ \mathbf{s}\mathbf{z}_{0}\right\} _{p}}
\end{align*}

Now, fixing $\mathbf{u}$ with $\left\Vert \mathbf{u}\right\Vert _{p}\leq p^{N}$,
observe that the map $\mathbf{s}\mapsto\mathbf{u}-\mathbf{s}$ is
a bijection of the set $\left\{ \mathbf{s}\in\hat{\mathbb{Z}}_{p}^{r}:\left\Vert \mathbf{s}\right\Vert _{p}\leq p^{N}\right\} $;
equivalently, this map is a bijection whenever $N\geq-v_{p}\left(\mathbf{u}\right)$.
So, picking $N\geq-v_{p}\left(\mathbf{u}\right)$, we can write:
\begin{align*}
\sum_{\left\Vert \mathbf{s}\right\Vert _{p}\leq p^{N}}\hat{\mathcal{M}}\left(\mathbf{u}-\mathbf{s}\right)e^{-2\pi i\left\{ \mathbf{s}\mathbf{z}_{0}\right\} _{p}} & =\sum_{\left\Vert \mathbf{s}\right\Vert _{p}\leq p^{N}}\hat{\mathcal{M}}\left(\mathbf{s}\right)e^{-2\pi i\left\{ \left(\mathbf{u}-\mathbf{s}\right)\mathbf{z}_{0}\right\} _{p}}\\
 & =e^{-2\pi i\left\{ \mathbf{u}\mathbf{z}_{0}\right\} _{p}}\sum_{\left\Vert \mathbf{s}\right\Vert _{p}\leq p^{N}}\hat{\mathcal{M}}\left(\mathbf{s}\right)e^{2\pi i\left\{ \mathbf{s}\mathbf{z}_{0}\right\} _{p}}\\
 & =e^{-2\pi i\left\{ \mathbf{u}\mathbf{z}_{0}\right\} _{p}}\tilde{\mathcal{M}}_{N}\left(\mathbf{z}_{0}\right)
\end{align*}
Denoting the entries of $\tilde{\mathcal{M}}_{N}\left(\mathbf{z}_{0}\right)$
as $\tilde{\mathcal{M}}_{N}\left(\mathbf{z}_{0}\right)=\left(\tilde{\mathcal{M}}_{N,1}\left(\mathbf{z}_{0}\right),\ldots,\tilde{\mathcal{M}}_{N,d}\left(\mathbf{z}_{0}\right)\right)$,
the given assumption about the zero at $\mathbf{z}_{0}$ tells us
there is an $\ell\in\left\{ 1,\ldots,d\right\} $ so that $\tilde{\mathcal{M}}_{N}\left(\mathbf{z}_{0}\right)$'s
$\ell$th entry ($\tilde{\mathcal{M}}_{N,\ell}\left(\mathbf{z}_{0}\right)$)
converges to $0$ in\emph{ $\mathbb{C}_{q}$} as $N\rightarrow\infty$.
Thus, for any $\epsilon^{\prime}>0$, there exists an $N_{\epsilon^{\prime}}$
so that $\left|\tilde{\mathcal{M}}_{N,\ell}\left(\mathbf{z}_{0}\right)\right|_{q}<\epsilon^{\prime}$
for all $N\geq N_{\epsilon^{\prime}}$.

Next, let $\mathbf{e}_{\ell}$ denote the $\ell$th standard basis
vector of $\mathbb{C}_{q}^{d}$. Then, for any $\hat{\mathbf{f}}=\left(\hat{f}_{1},\ldots,\hat{f}_{d}\right):\hat{\mathbb{Z}}_{p}^{r}\rightarrow\mathbb{C}_{q}^{d}$,
the entry-wise product $\mathbf{e}_{\ell}\hat{\mathbf{f}}$ is: 
\begin{equation}
\left(\mathbf{e}_{\ell}\hat{\mathbf{f}}\right)\left(\mathbf{t}\right)=\left(0,\ldots,0,\hat{f}_{\ell}\left(\mathbf{t}\right),0,\ldots,0\right)\label{eq:e_ell f hat}
\end{equation}
Moreover, since the convolution of two vector-valued functions on
$\hat{\mathbb{Z}}_{p}^{r}$ is done entry-wise: 
\begin{align*}
\left(\hat{\mathbf{f}}*\hat{\mathbf{g}}\right)\left(\mathbf{t}\right) & =\sum_{\mathbf{s}\in\hat{\mathbb{Z}}_{p}^{r}}\hat{\mathbf{f}}\left(\mathbf{t}-\mathbf{s}\right)\hat{\mathbf{g}}\left(\mathbf{s}\right)\\
 & =\left(\sum_{\mathbf{s}\in\hat{\mathbb{Z}}_{p}^{r}}\hat{f}_{1}\left(\mathbf{t}-\mathbf{s}\right)\hat{g}_{1}\left(\mathbf{s}\right),\ldots,\sum_{\mathbf{s}\in\hat{\mathbb{Z}}_{p}^{r}}\hat{f}_{d}\left(\mathbf{t}-\mathbf{s}\right)\hat{g}_{d}\left(\mathbf{s}\right)\right)
\end{align*}
we obtain: 
\begin{equation}
\left(\left(\mathbf{e}_{\ell}\hat{\mathbf{f}}\right)*\hat{\mathbf{g}}\right)\left(\mathbf{t}\right)=\left(\mathbf{e}_{\ell}\left(\hat{\mathbf{f}}*\hat{\mathbf{g}}\right)\right)\left(\mathbf{t}\right)=\left(\left(\hat{f}_{1}*\hat{g}_{1}\right)\left(\mathbf{t}\right),0,\ldots,0\right)\label{eq:e_ell f hat convolve with g hat}
\end{equation}
With this notation, we can write:

\begin{align*}
\left(\mathbf{e}_{\ell}\hat{\mathcal{M}}*\hat{\phi}_{N}\right)\left(\mathbf{u}\right) & =\mathbf{e}_{\ell}\left(\hat{\mathcal{M}}*\hat{\phi}_{N}\right)\left(\mathbf{u}\right)\\
 & =\mathbf{e}_{\ell}\left(e^{-2\pi i\left\{ \mathbf{u}\mathbf{z}_{0}\right\} _{p}}\tilde{\mathcal{M}}_{N}\left(\mathbf{z}_{0}\right)\right)\\
 & =e^{-2\pi i\left\{ \mathbf{u}\mathbf{z}_{0}\right\} _{p}}\left(\mathbf{e}_{\ell}\tilde{\mathcal{M}}_{N}\right)\left(\mathbf{z}_{0}\right)\\
 & =\left(0,\ldots,0,e^{-2\pi i\left\{ \mathbf{u}\mathbf{z}_{0}\right\} _{p}}\tilde{\mathcal{M}}_{N,\ell}\left(\mathbf{z}_{0}\right),0,\ldots,0\right)
\end{align*}
Taking $q$-adic absolute values, we have thereby proven the following
claim:
\begin{claim}
\label{claim:5.11}For any $\epsilon^{\prime}>0$, there is an $N_{\epsilon^{\prime}}\geq0$
(depending only on $\hat{\mathcal{M}}$, $\ell$ and $\epsilon^{\prime}$)
so that, for any $\mathbf{u}\in\hat{\mathbb{Z}}_{p}$: 
\end{claim}
\begin{equation}
\left\Vert \left(\mathbf{e}_{\ell}\hat{\mathcal{M}}*\hat{\phi}_{N}\right)\left(\mathbf{u}\right)\right\Vert _{p^{\infty},q}<\epsilon^{\prime},\textrm{ }\forall N\geq\max\left\{ N_{\epsilon^{\prime}},-v_{p}\left(\mathbf{u}\right)\right\} \label{eq:MD WTT - First Claim}
\end{equation}

\vphantom{}

For our next step, we restrict the support of $\hat{\mathbf{h}}_{\epsilon}$
by defining the function $\hat{\mathbf{h}}_{\epsilon,M}:\hat{\mathbb{Z}}_{p}^{r}\rightarrow\mathbb{C}_{q}^{d}$
as:
\begin{equation}
\hat{\mathbf{h}}_{\epsilon,M}\left(\mathbf{t}\right)\overset{\textrm{def}}{=}\vec{\mathbf{1}}_{\mathbf{0}}\left(p^{M}\mathbf{t}\right)\hat{\mathbf{h}}_{\epsilon}\left(\mathbf{t}\right)=\begin{cases}
\hat{\mathbf{h}}_{\epsilon}\left(\mathbf{t}\right) & \textrm{if }\left\Vert \mathbf{t}\right\Vert _{p}\leq p^{M}\\
\mathbf{0} & \textrm{else}
\end{cases}\label{eq:MD Definition of h epsilon M hat}
\end{equation}
Here $M$ is an arbitrary integer $\geq0$. Then: 
\begin{align*}
\left(\hat{\mathcal{M}}*\hat{\mathbf{h}}_{\epsilon,M}\right)\left(\mathbf{t}\right) & =\sum_{\mathbf{u}\in\mathbb{Z}_{p}}\hat{\mathcal{M}}\left(\mathbf{t}-\mathbf{u}\right)\vec{\mathbf{1}}_{\mathbf{0}}\left(p^{M}\mathbf{u}\right)\sum_{k=1}^{K}\mathbf{c}_{k}\vec{\mathbf{1}}_{\mathbf{t}_{k}}\left(\mathbf{u}\right)\\
 & =\sum_{\left\Vert \mathbf{u}\right\Vert _{p}\leq p^{M}}\hat{\mathcal{M}}\left(\mathbf{t}-\mathbf{u}\right)\sum_{k=1}^{K}\mathbf{c}_{k}\vec{\mathbf{1}}_{\mathbf{t}_{k}}\left(\mathbf{u}\right)\\
 & =\sum_{k:\left\Vert \mathbf{t}_{k}\right\Vert _{p}\leq p^{M}}\mathbf{c}_{k}\hat{\mathcal{M}}\left(\mathbf{t}-\mathbf{t}_{k}\right)
\end{align*}
Since there are finitely many $\mathbf{t}_{k}$s, we can choose the
non-negative integer: 
\begin{equation}
M_{\epsilon}\overset{\textrm{choose}}{=}\max\left\{ -v_{p}\left(\mathbf{t}_{1}\right),\ldots,-v_{p}\left(\mathbf{t}_{K}\right)\right\} \label{eq:MD WTT - Choice for M_e}
\end{equation}
In particular, note that\emph{ $M_{\epsilon}$ depends only on $\mathbf{h}_{\epsilon}$
and $\epsilon$. }Then: 
\begin{equation}
\left(\hat{\mathcal{M}}*\hat{\mathbf{h}}_{\epsilon,M}\right)\left(\mathbf{t}\right)=\sum_{k=1}^{K}\mathbf{c}_{k}\hat{\mathcal{M}}\left(\mathbf{t}-\mathbf{t}_{k}\right)=\left(\hat{\mathcal{M}}*\hat{\mathbf{h}}_{\epsilon}\right)\left(\mathbf{t}\right),\textrm{ }\forall M\geq M_{\epsilon}\label{eq:MD Effect of M bigger than M0}
\end{equation}
Next, for any $\hat{\mathbf{f}}:\hat{\mathbb{Z}}_{p}^{r}\rightarrow\mathbb{C}_{q}^{d}$
and any $n\in\mathbb{N}_{0}$, let us write: \nomenclature{$\left\Vert \hat{\mathbf{f}}\right\Vert _{p^{n},q}$}{ }

\begin{equation}
\left\Vert \hat{\mathbf{f}}\right\Vert _{p^{n},q}\overset{\textrm{def}}{=}\sup_{\left\Vert \mathbf{t}\right\Vert _{p}\leq p^{n}}\left\Vert \hat{\mathbf{f}}\left(\mathbf{t}\right)\right\Vert _{q}\label{eq:MD Definition of infinity,p^n norm}
\end{equation}
Then:

\begin{align*}
\left\Vert \vec{\mathbf{1}}_{\mathbf{0}}-\hat{\mathcal{M}}*\hat{\mathbf{h}}_{\epsilon,M}\right\Vert _{p^{M},q} & =\sup_{\left\Vert \mathbf{t}\right\Vert _{p}\leq p^{M}}\left\Vert \vec{\mathbf{1}}_{\mathbf{0}}\left(\mathbf{t}\right)-\left(\hat{\mathcal{M}}*\hat{\mathbf{h}}_{\epsilon,M}\right)\left(\mathbf{t}\right)\right\Vert _{q}\\
\left(\textrm{if }M\geq M_{\epsilon}\right); & =\sup_{\left\Vert \mathbf{t}\right\Vert _{p}\leq p^{M}}\left\Vert \vec{\mathbf{1}}_{\mathbf{0}}\left(\mathbf{t}\right)-\left(\hat{\mathcal{M}}*\hat{\mathbf{h}}_{\epsilon}\right)\left(\mathbf{t}\right)\right\Vert _{q}\\
 & \leq\left\Vert \vec{\mathbf{1}}_{\mathbf{0}}-\left(\hat{\mathcal{M}}*\hat{\mathbf{h}}_{\epsilon}\right)\right\Vert _{p^{\infty},q}\\
 & <\epsilon
\end{align*}
In summary, we have shown: 
\begin{claim}
\label{claim:5.12}Let $\epsilon>0$ be arbitrary. Then, there exists
an $M_{\epsilon}$ (depending only on $\hat{\mathcal{M}}$ and $\epsilon$
so that:

\vphantom{}

I. $\hat{\mathcal{M}}*\hat{\mathbf{h}}_{\epsilon,M}=\hat{\mathcal{M}}*\hat{\mathbf{h}}_{\epsilon},\textrm{ }\forall M\geq M_{\epsilon}$

\vphantom{}

II. $\left\Vert \vec{\mathbf{1}}_{\mathbf{0}}-\hat{\mathcal{M}}*\hat{\mathbf{h}}_{\epsilon,M}\right\Vert _{p^{\infty},q}<\epsilon,\textrm{ }\forall M\geq M_{\epsilon}$

\vphantom{}

III. $\left\Vert \vec{\mathbf{1}}_{\mathbf{0}}-\hat{\mathcal{M}}*\hat{\mathbf{h}}_{\epsilon,M}\right\Vert _{p^{M},q}<\epsilon,\textrm{ }\forall M\geq M_{\epsilon}$ 
\end{claim}
\vphantom{}

We now modify these estimates to take into account the norm $\left\Vert \cdot\right\Vert _{p^{m},q}$. 
\begin{claim}
\label{claim:5.13}Let $\epsilon^{\prime}>0$ and $m\in\mathbb{N}_{1}$
be arbitrary. Then, there exists $N_{\epsilon^{\prime}}$ (depending
only on $\epsilon^{\prime}$ and $\hat{\mathcal{M}}$) so that: 
\begin{equation}
\left\Vert \left(\mathbf{e}_{\ell}\hat{\mathcal{M}}\right)*\hat{\phi}_{N}\right\Vert _{p^{m},q}<\epsilon^{\prime},\textrm{ }\forall N\geq\max\left\{ N_{\epsilon^{\prime}},m\right\} \label{eq:MD WTT - eq. 4}
\end{equation}

Proof of claim: For arbitrary $\epsilon^{\prime}>0$, \textbf{Claim
\ref{claim:5.11}} tells us there is an $N_{\epsilon^{\prime}}$ with
the stated dependencies so that, for any $\mathbf{u}\in\hat{\mathbb{Z}}_{p}^{r}$:
\begin{equation}
\left\Vert \left(\left(\mathbf{e}_{\ell}\hat{\mathcal{M}}\right)*\hat{\phi}_{N}\right)\left(\mathbf{u}\right)\right\Vert _{p^{\infty},q}<\epsilon^{\prime},\textrm{ }\forall N\geq\max\left\{ N_{\epsilon},-v\left(\mathbf{u}\right)\right\} 
\end{equation}
So, letting $m\geq1$ be arbitrary, note that $\left\Vert \mathbf{u}\right\Vert _{p}\leq p^{m}$
implies $-v_{p}\left(\mathbf{u}\right)\leq m$. As such, by choosing
$N\geq\max\left\{ N_{\epsilon},m\right\} $, we can make the result
of \textbf{Claim \ref{claim:5.11}} hold for all $\left\Vert \mathbf{u}\right\Vert _{p}\leq p^{m}$:

\begin{equation}
\underbrace{\sup_{\left\Vert \mathbf{u}\right\Vert _{p}\leq p^{m}}\left\Vert \left(\left(\mathbf{e}_{\ell}\hat{\mathcal{M}}\right)*\hat{\phi}_{N}\right)\left(\mathbf{u}\right)\right\Vert _{q}}_{\left\Vert \left(\mathbf{e}_{\ell}\hat{\mathcal{M}}\right)*\hat{\phi}_{N}\right\Vert _{p^{m},q}}<\epsilon^{\prime},\textrm{ }\forall N\geq\max\left\{ N_{\epsilon},m\right\} 
\end{equation}
This proves the claim. 
\end{claim}
\vphantom{}

We need two more estimates before we can wrap up the argument. 
\begin{claim}
\label{claim:5.14}Let $M,N\in\mathbb{N}_{0}$ be arbitrary. Then,
for any $\hat{\varphi},\hat{\mathbf{f}},\hat{\mathbf{g}}:\hat{\mathbb{Z}}_{p}^{r}\rightarrow\mathbb{C}_{q}^{d}$
satisfying:

\vphantom{}

i. $\hat{\varphi}\left(\mathbf{t}\right)=\mathbf{0}$ for all $\left\Vert \mathbf{t}\right\Vert _{p}>p^{N}$;

\vphantom{}

ii. $\left\Vert \hat{\varphi}\right\Vert _{p^{\infty},q}\leq1$;

\vphantom{}

iii. $\hat{\mathbf{g}}\left(\mathbf{t}\right)=\mathbf{0}$ for all
$\left\Vert \mathbf{t}\right\Vert _{p}>p^{M}$;

\vphantom{}

the following estimates hold:

\begin{equation}
\left\Vert \hat{\varphi}*\hat{\mathbf{f}}\right\Vert _{p^{M},q}\leq\left\Vert \hat{\mathbf{f}}\right\Vert _{p^{\max\left\{ M,N\right\} },q}\label{eq:MD phi_N hat convolve f hat estimate}
\end{equation}
\begin{equation}
\left\Vert \hat{\varphi}*\hat{\mathbf{f}}*\hat{\mathbf{g}}\right\Vert _{p^{M},q}\leq\left\Vert \hat{\mathbf{g}}\right\Vert _{p^{\infty},q}\left\Vert \hat{\varphi}*\hat{\mathbf{f}}\right\Vert _{p^{\max\left\{ M,N\right\} },q}\label{eq:MD phi_N hat convolve f hat convolve g hat estimate}
\end{equation}

Proof of claim: As in the one-dimensional case, we start with the
top-most inequality:

\begin{align*}
\left\Vert \hat{\varphi}*\hat{\mathbf{f}}\right\Vert _{p^{M},q} & =\sup_{\left\Vert \mathbf{t}\right\Vert _{p}\leq p^{M}}\left|\sum_{\mathbf{s}\in\hat{\mathbb{Z}}_{p}^{r}}\hat{\varphi}\left(\mathbf{s}\right)\hat{\mathbf{f}}\left(\mathbf{t}-\mathbf{s}\right)\right|_{q}\\
\left(\hat{\varphi}\left(\mathbf{s}\right)=\mathbf{0},\textrm{ }\forall\left\Vert \mathbf{s}\right\Vert _{p}>p^{N}\right); & \leq\sup_{\left\Vert \mathbf{t}\right\Vert _{p}\leq p^{M}}\left\Vert \sum_{\left\Vert \mathbf{s}\right\Vert _{p}\leq p^{N}}\hat{\varphi}\left(\mathbf{s}\right)\hat{\mathbf{f}}\left(\mathbf{t}-\mathbf{s}\right)\right\Vert _{q}\\
\left(\textrm{ultrametric ineq.}\right); & \leq\sup_{\left\Vert \mathbf{t}\right\Vert _{p}\leq p^{M}}\sup_{\left\Vert \mathbf{s}\right\Vert _{p}\leq p^{N}}\left\Vert \hat{\varphi}\left(\mathbf{s}\right)\hat{\mathbf{f}}\left(\mathbf{t}-\mathbf{s}\right)\right\Vert _{q}\\
 & \leq\sup_{\left\Vert \mathbf{s}\right\Vert _{p}\leq p^{N}}\left\Vert \hat{\varphi}\left(\mathbf{s}\right)\right\Vert _{q}\sup_{\left\Vert \mathbf{t}\right\Vert _{p}\leq p^{\max\left\{ M,N\right\} }}\left\Vert \hat{\mathbf{f}}\left(\mathbf{t}\right)\right\Vert _{q}\\
 & =\left\Vert \hat{\varphi}\right\Vert _{p^{\infty},q}\cdot\left\Vert \hat{\mathbf{f}}\right\Vert _{p^{\max\left\{ M,N\right\} },q}\\
\left(\left\Vert \hat{\varphi}\right\Vert _{p^{\infty},q}\leq1\right); & \leq\left\Vert \hat{\mathbf{f}}\right\Vert _{p^{\max\left\{ M,N\right\} },q}
\end{align*}
This proves (\ref{eq:MD phi_N hat convolve f hat estimate}).

Next, we deal with (\ref{eq:MD phi_N hat convolve f hat convolve g hat estimate}).
For that, we start by writing out the convolution of $\hat{\varphi}*\hat{\mathbf{f}}$
with $\hat{\mathbf{g}}$: 
\begin{align*}
\left\Vert \hat{\varphi}*\hat{\mathbf{f}}*\hat{\mathbf{g}}\right\Vert _{p^{M},q} & =\sup_{\left\Vert \mathbf{t}\right\Vert _{p}\leq p^{M}}\left\Vert \sum_{\mathbf{s}\in\hat{\mathbb{Z}}_{p}^{r}}\left(\hat{\varphi}*\hat{\mathbf{f}}\right)\left(\mathbf{t}-\mathbf{s}\right)\hat{\mathbf{g}}\left(\mathbf{s}\right)\right\Vert _{q}\\
\left(\textrm{ultrametric ineq.}\right); & \leq\sup_{\left\Vert \mathbf{t}\right\Vert _{p}\leq p^{M}}\sup_{\mathbf{s}\in\hat{\mathbb{Z}}_{p}^{r}}\left\Vert \left(\hat{\varphi}*\hat{\mathbf{f}}\right)\left(\mathbf{t}-\mathbf{s}\right)\hat{\mathbf{g}}\left(\mathbf{s}\right)\right\Vert _{q}\\
\left(\hat{\mathbf{g}}\left(\mathbf{s}\right)=\mathbf{0},\textrm{ }\forall\left\Vert \mathbf{s}\right\Vert _{p}>p^{M}\right); & \leq\sup_{\left\Vert \mathbf{t}\right\Vert _{p}\leq p^{M}}\sup_{\left\Vert \mathbf{s}\right\Vert _{p}\leq p^{M}}\left\Vert \left(\hat{\varphi}*\hat{\mathbf{f}}\right)\left(\mathbf{t}-\mathbf{s}\right)\hat{\mathbf{g}}\left(\mathbf{s}\right)\right\Vert _{q}\\
 & \leq\left\Vert \hat{\mathbf{g}}\right\Vert _{p^{\infty},q}\sup_{\left\Vert \mathbf{t}\right\Vert _{p},\left\Vert \mathbf{s}\right\Vert _{p}\leq p^{M}}\left\Vert \left(\hat{\varphi}*\hat{\mathbf{f}}\right)\left(\mathbf{t}-\mathbf{s}\right)\right\Vert _{q}
\end{align*}
Next, we write out the convolution $\hat{\varphi}*\hat{\mathbf{f}}$.
This gives us:
\begin{align*}
\left\Vert \hat{\varphi}*\hat{\mathbf{f}}*\hat{\mathbf{g}}\right\Vert _{p^{M},q} & \leq\left\Vert \hat{\mathbf{g}}\right\Vert _{p^{\infty},q}\sup_{\left\Vert \mathbf{t}\right\Vert _{p},\left\Vert \mathbf{s}\right\Vert _{p}\leq p^{M}}\left\Vert \sum_{\mathbf{v}\in\hat{\mathbb{Z}}_{p}}\hat{\varphi}\left(\mathbf{t}-\mathbf{s}-\mathbf{v}\right)\hat{\mathbf{f}}\left(\mathbf{v}\right)\right\Vert _{q}\\
\left(\textrm{let }\mathbf{u}=\mathbf{s}+\mathbf{v}\right); & =\left\Vert \hat{\mathbf{g}}\right\Vert _{p^{\infty},q}\sup_{\left\Vert \mathbf{t}\right\Vert _{p},\left\Vert \mathbf{s}\right\Vert _{p}\leq p^{M}}\left\Vert \sum_{\mathbf{u}-\mathbf{s}\in\hat{\mathbb{Z}}_{p}^{r}}\hat{\varphi}\left(\mathbf{t}-\mathbf{u}\right)\hat{\mathbf{f}}\left(\mathbf{u}-\mathbf{s}\right)\right\Vert _{q}\\
\left(\mathbf{s}+\hat{\mathbb{Z}}_{p}^{r}=\hat{\mathbb{Z}}_{p}^{r}\right); & =\left\Vert \hat{\mathbf{g}}\right\Vert _{p^{\infty},q}\sup_{\left\Vert \mathbf{t}\right\Vert _{p},\left\Vert \mathbf{s}\right\Vert _{p}\leq p^{M}}\left\Vert \sum_{\mathbf{u}\in\hat{\mathbb{Z}}_{p}^{r}}\hat{\varphi}\left(\mathbf{t}-\mathbf{u}\right)\hat{\mathbf{f}}\left(\mathbf{u}-\mathbf{s}\right)\right\Vert _{q}
\end{align*}

We now invoke the vanishing of $\hat{\varphi}\left(\mathbf{t}-\mathbf{u}\right)$
for all $\left\Vert \mathbf{t}-\mathbf{u}\right\Vert _{p}>p^{N}$.
Indeed, because $\mathbf{t}$ is restricted to $\left\Vert \mathbf{t}\right\Vert _{p}\leq p^{M}$,
observe that when $\left\Vert \mathbf{u}\right\Vert _{p}>p^{\max\left\{ M,N\right\} }$,
the ultrametric inequality yields: 
\begin{equation}
\left\Vert \mathbf{t}-\mathbf{u}\right\Vert _{p}=\max\left\{ \left\Vert \mathbf{t}\right\Vert _{p},\left\Vert \mathbf{u}\right\Vert _{p}\right\} >p^{\max\left\{ M,N\right\} }>p^{N}
\end{equation}
Hence, for $\left\Vert \mathbf{t}\right\Vert _{p},\left\Vert \mathbf{s}\right\Vert _{p}\leq p^{M}$,
the summand $\hat{\varphi}\left(\mathbf{t}-\mathbf{u}\right)\hat{\mathbf{f}}\left(\mathbf{u}-\mathbf{s}\right)$
vanishes whenever $\left\Vert \mathbf{u}\right\Vert >p^{\max\left\{ M,N\right\} }$.
This gives us: 
\begin{equation}
\left\Vert \hat{\varphi}*\hat{\mathbf{f}}*\hat{\mathbf{g}}\right\Vert _{p^{M},q}\leq\left\Vert \hat{\mathbf{g}}\right\Vert _{p^{\infty},q}\sup_{\left\Vert \mathbf{t}\right\Vert _{p},\left\Vert \mathbf{s}\right\Vert _{p}\leq p^{M}}\left\Vert \sum_{\left\Vert \mathbf{u}\right\Vert _{p}\leq p^{\max\left\{ M,N\right\} }}\hat{\varphi}\left(\mathbf{t}-\mathbf{u}\right)\hat{\mathbf{f}}\left(\mathbf{u}-\mathbf{s}\right)\right\Vert _{q}
\end{equation}
Next, we enlarge things slightly by expanding the range of $\left\Vert \mathbf{s}\right\Vert _{p}$
and $\left\Vert \mathbf{t}\right\Vert _{p}$ from $\leq p^{M}$ to
$\leq p^{\max\left\{ M,N\right\} }$:
\begin{equation}
\left\Vert \hat{\varphi}*\hat{\mathbf{f}}*\hat{\mathbf{g}}\right\Vert _{p^{M},q}\leq\left\Vert \hat{\mathbf{g}}\right\Vert _{p^{\infty},q}\sup_{\left\Vert \mathbf{t}\right\Vert _{p},\left\Vert \mathbf{s}\right\Vert _{p}\leq p^{\max\left\{ M,N\right\} }}\left\Vert \sum_{\left\Vert \mathbf{u}\right\Vert _{p}\leq p^{\max\left\{ M,N\right\} }}\hat{\varphi}\left(\mathbf{t}-\mathbf{u}\right)\hat{\mathbf{f}}\left(\mathbf{u}-\mathbf{s}\right)\right\Vert _{q}
\end{equation}
Like in the one-dimensional case, everything is now on the same level:
$\mathbf{t}$, $\mathbf{s}$, and $\mathbf{u}$ are all bounded in
$p$-adic norm by $p^{\max\left\{ M,N\right\} }$. Exploiting the
closure of the set: 
\[
\left\{ \mathbf{x}\in\hat{\mathbb{Z}}_{p}^{r}:\left\Vert \mathbf{x}\right\Vert _{p}\leq p^{\max\left\{ M,N\right\} }\right\} 
\]
under addition allows us to note that our $\mathbf{u}$-sum is invariant
under the change of variables $\mathbf{u}\mapsto\mathbf{u}+\mathbf{s}$
for any $\left\Vert \mathbf{s}\right\Vert _{p}\leq p^{\max\left\{ M,N\right\} }$.
Making this change of variables then gives us:
\begin{align*}
\left\Vert \hat{\varphi}*\hat{\mathbf{f}}*\hat{\mathbf{g}}\right\Vert _{p^{M},q} & \leq\left\Vert \hat{\mathbf{g}}\right\Vert _{p^{\infty},q}\sup_{\left\Vert \mathbf{t}\right\Vert _{p},\left\Vert \mathbf{s}\right\Vert _{p}\leq p^{\max\left\{ M,N\right\} }}\left\Vert \sum_{\left\Vert \mathbf{u}\right\Vert _{p}\leq p^{\max\left\{ M,N\right\} }}\hat{\varphi}\left(\mathbf{t}-\left(\mathbf{u}+\mathbf{s}\right)\right)\hat{\mathbf{f}}\left(\mathbf{u}\right)\right\Vert _{q}\\
 & =\left\Vert \hat{\mathbf{g}}\right\Vert _{p^{\infty},q}\sup_{\left\Vert \mathbf{t}\right\Vert _{p},\left\Vert \mathbf{s}\right\Vert _{p}\leq p^{\max\left\{ M,N\right\} }}\left\Vert \sum_{\left\Vert \mathbf{u}\right\Vert _{p}\leq p^{\max\left\{ M,N\right\} }}\hat{\varphi}\left(\mathbf{t}-\mathbf{s}-\mathbf{u}\right)\hat{\mathbf{f}}\left(\mathbf{u}\right)\right\Vert _{q}
\end{align*}
Because: 
\begin{equation}
\left\{ \mathbf{t}-\mathbf{s}:\left\Vert \mathbf{t}\right\Vert _{p},\left\Vert \mathbf{s}\right\Vert _{p}\leq p^{\max\left\{ M,N\right\} }\right\} =\left\{ \mathbf{t}:\left\Vert \mathbf{t}\right\Vert _{p}\leq p^{\max\left\{ M,N\right\} }\right\} 
\end{equation}
we then have: 
\begin{align*}
\left\Vert \hat{\varphi}*\hat{\mathbf{f}}*\hat{\mathbf{g}}\right\Vert _{p^{M},q} & \leq\left\Vert \hat{\mathbf{g}}\right\Vert _{p^{\infty},q}\sup_{\left\Vert \mathbf{t}\right\Vert _{p}\leq p^{\max\left\{ M,N\right\} }}\left\Vert \sum_{\left\Vert \mathbf{u}\right\Vert _{p}\leq p^{\max\left\{ M,N\right\} }}\hat{\varphi}\left(\mathbf{t}-\mathbf{u}\right)\hat{\mathbf{f}}\left(\mathbf{u}\right)\right\Vert _{q}\\
 & \leq\left\Vert \hat{\mathbf{g}}\right\Vert _{p^{\infty},q}\sup_{\left\Vert \mathbf{t}\right\Vert _{p}\leq p^{\max\left\{ M,N\right\} }}\left\Vert \sum_{\mathbf{u}\in\hat{\mathbb{Z}}_{p}^{r}}\hat{\varphi}\left(\mathbf{t}-\mathbf{u}\right)\hat{\mathbf{f}}\left(\mathbf{u}\right)\right\Vert _{q}\\
 & =\left\Vert \hat{\mathbf{g}}\right\Vert _{p^{\infty},q}\sup_{\left\Vert \mathbf{t}\right\Vert _{p}\leq p^{\max\left\{ M,N\right\} }}\left\Vert \left(\hat{\varphi}*\hat{\mathbf{f}}\right)\left(\mathbf{u}\right)\right\Vert _{q}\\
\left(\textrm{by definition}\right); & =\left\Vert \hat{\mathbf{g}}\right\Vert _{p^{\infty},q}\left\Vert \hat{\varphi}*\hat{\mathbf{f}}\right\Vert _{p^{\max\left\{ M,N\right\} },q}
\end{align*}
which proves the desired estimate. 
\end{claim}
\vphantom{}

Now, let's line everything up: first, we choose $\epsilon\in\left(0,1\right)$
and a non-identically-zero $\hat{\mathbf{h}}_{\epsilon}$ so that:
\begin{equation}
\left\Vert \vec{\mathbf{1}}_{\mathbf{0}}-\hat{\mathcal{M}}*\hat{\mathbf{h}}_{\epsilon}\right\Vert _{p^{\infty},q}<\epsilon
\end{equation}
Then, by \textbf{Claim \ref{claim:5.12}}, we can choose an integer
$M_{\epsilon}$ (depending only on $\epsilon$ and $\hat{\mathbf{h}}_{\epsilon}$)
so that: 
\begin{equation}
\left\Vert \vec{\mathbf{1}}_{\mathbf{0}}-\hat{\mathcal{M}}*\hat{\mathbf{h}}_{\epsilon,M}\right\Vert _{p^{M},q}<\epsilon,\textrm{ }\forall M\geq M_{\epsilon}
\end{equation}
Since the norm is the maximum of the components, multiplying the left-hand
side entry-wise by the standard basis vector $\mathbf{e}_{\ell}$
can only decrease the norm: 
\begin{equation}
\left\Vert \mathbf{e}_{\ell}\left(\vec{\mathbf{1}}_{\mathbf{0}}-\hat{\mathcal{M}}*\hat{\mathbf{h}}_{\epsilon,M}\right)\right\Vert _{p^{M},q}\leq\left\Vert \vec{\mathbf{1}}_{\mathbf{0}}-\hat{\mathcal{M}}*\hat{\mathbf{h}}_{\epsilon,M}\right\Vert _{p^{M},q}<\epsilon,\textrm{ }\forall M\geq M_{\epsilon}
\end{equation}

Next, we observe that: 
\begin{equation}
\left\Vert \mathbf{e}_{\ell}\left(\hat{\phi}_{N}*\left(\vec{\mathbf{1}}_{\mathbf{0}}-\left(\hat{\mathcal{M}}*\hat{\mathbf{h}}_{\epsilon,M}\right)\right)\right)\right\Vert _{p^{M},q}=\left\Vert \mathbf{e}_{\ell}\hat{\phi}_{N}*\mathbf{e}_{\ell}\left(\vec{\mathbf{1}}_{\mathbf{0}}-\left(\hat{\mathcal{M}}*\hat{\mathbf{h}}_{\epsilon,M}\right)\right)\right\Vert _{p^{M},q}\label{eq:MD WTT - Target of attack}
\end{equation}
Like in the one-dimensional case, we will estimate this in two ways.
First, we apply (\ref{eq:MD phi_N hat convolve f hat estimate}) from
\textbf{Claim \ref{claim:5.14}} with $\hat{\varphi}=\mathbf{e}_{\ell}\hat{\phi}_{N}$
and $\hat{\mathbf{f}}=\mathbf{e}_{\ell}\left(\vec{\mathbf{1}}_{\mathbf{0}}-\left(\hat{\mathcal{M}}*\hat{\mathbf{h}}_{\epsilon,M}\right)\right)$.
This gives us: 
\begin{align*}
\left\Vert \mathbf{e}_{\ell}\hat{\phi}_{N}*\mathbf{e}_{\ell}\left(\vec{\mathbf{1}}_{\mathbf{0}}-\left(\hat{\mathcal{M}}*\hat{\mathbf{h}}_{\epsilon,M}\right)\right)\right\Vert _{p^{M},q} & \leq\left\Vert \mathbf{e}_{\ell}\left(\vec{\mathbf{1}}_{\mathbf{0}}-\left(\hat{\mathcal{M}}*\hat{\mathbf{h}}_{\epsilon,M}\right)\right)\right\Vert _{p^{\max\left\{ M,N\right\} },q}\\
 & \leq\left\Vert \vec{\mathbf{1}}_{\mathbf{0}}-\left(\hat{\mathcal{M}}*\hat{\mathbf{h}}_{\epsilon,M}\right)\right\Vert _{p^{\max\left\{ M,N\right\} },q}\\
\left(M\geq M_{\epsilon}\Rightarrow\textrm{use \textbf{Claim \ref{claim:5.12}}}\right); & <\epsilon
\end{align*}
So, the \emph{left-hand side} of (\ref{eq:MD WTT - Target of attack})
is $<\epsilon$.

To contradict this, we use the ultrametric inequality for the non-archimedean
norm $\left\Vert \cdot\right\Vert _{p^{M},q}$ on the \emph{right-hand
side} of (\ref{eq:MD WTT - Target of attack}). This yields:
\begin{equation}
\max\left\{ \left\Vert \mathbf{e}_{\ell}\hat{\phi}_{N}\right\Vert _{p^{M},q},\left\Vert \hat{\phi}_{N}*\mathbf{e}_{\ell}\hat{\mathcal{M}}*\hat{\mathbf{h}}_{\epsilon,M}\right\Vert _{p^{M},q}\right\} \label{eq:WTT - Ultrametric inequality-1}
\end{equation}
as an upper bound for: 
\begin{equation}
\left\Vert \mathbf{e}_{\ell}\left(\hat{\phi}_{N}-\left(\hat{\phi}_{N}*\hat{\mathcal{M}}*\hat{\mathbf{h}}_{\epsilon,M}\right)\right)\right\Vert _{p^{M},q}
\end{equation}
(Note that we could have also written $\mathbf{e}_{\ell}$ next to
$\hat{\phi}_{N}$ and $\hat{\mathbf{h}}_{\epsilon,M}$, though it
would have been redundant.) Using (\ref{eq:MD phi_N hat convolve f hat convolve g hat estimate})
from \textbf{Claim \ref{claim:5.14}}, we get: 
\begin{align*}
\left\Vert \hat{\phi}_{N}*\mathbf{e}_{\ell}\hat{\mathcal{M}}*\hat{\mathbf{h}}_{\epsilon,M}\right\Vert _{p^{M},q} & \leq\left\Vert \hat{\phi}_{N}*\mathbf{e}_{\ell}\hat{\mathcal{M}}\right\Vert _{p^{M},q}\left\Vert \hat{\mathbf{h}}_{\epsilon,M}\right\Vert _{p^{M},q}\\
\left(M\geq M_{\epsilon}\right); & \leq\left\Vert \hat{\phi}_{N}*\mathbf{e}_{\ell}\hat{\mathcal{M}}\right\Vert _{p^{M},q}\cdot\left\Vert \hat{\mathbf{h}}_{\epsilon}\right\Vert _{p^{\infty},q}
\end{align*}

Now, let $\epsilon^{\prime}=\frac{\epsilon}{2\left\Vert \hat{\mathbf{h}}_{\epsilon}\right\Vert _{p^{\infty},q}}$;
we are allowed to do this, thanks to our assumption that $\hat{\mathbf{h}}_{\epsilon}$
was not identically zero. Seeing as $N$ is still arbitrary, let $N$
be larger than both $N_{\epsilon^{\prime}}$ and our choice for $M$;
\textbf{Claim \ref{claim:5.13}} tells us that $\left\Vert \mathbf{e}_{\ell}\hat{\mathcal{M}}*\hat{\phi}_{N}\right\Vert _{p^{M},q}<\epsilon^{\prime}$.
Since the size of the $\ell$th component of $\mathbf{e}_{\ell}\hat{\mathcal{M}}*\hat{\phi}_{N}$
is: 
\begin{equation}
\left\Vert \mathbf{e}_{\ell}\hat{\mathcal{M}}*\hat{\phi}_{N}\right\Vert _{p^{M},q}<\epsilon^{\prime}
\end{equation}
we can write:
\begin{align*}
\left\Vert \hat{\phi}_{N}*\mathbf{e}_{\ell}\hat{\mathcal{M}}*\hat{\mathbf{h}}_{\epsilon,M}\right\Vert _{p^{M},q} & \leq\left\Vert \hat{\phi}_{N}*\mathbf{e}_{\ell}\hat{\mathcal{M}}\right\Vert _{p^{M},q}\cdot\left\Vert \hat{\mathbf{h}}_{\epsilon}\right\Vert _{p^{\infty},q}\\
 & <\epsilon^{\prime}\cdot\left\Vert \hat{\mathbf{h}}_{\epsilon}\right\Vert _{p^{\infty},q}\\
 & =\frac{\epsilon}{2\left\Vert \hat{\mathbf{h}}_{\epsilon}\right\Vert _{p^{\infty},q}}\cdot\left\Vert \hat{\mathbf{h}}_{\epsilon}\right\Vert _{p^{\infty},q}\\
 & =\frac{\epsilon}{2}
\end{align*}
for our large choice of $N$. Because $\left\Vert \mathbf{e}_{\ell}\hat{\phi}_{N}\right\Vert _{p^{M},q}=1$
for all $M,N\geq0$, this shows that: 
\begin{equation}
\left\Vert \mathbf{e}_{\ell}\hat{\phi}_{N}\right\Vert _{p^{M},q}=1>\frac{\epsilon}{2}>\left\Vert \hat{\phi}_{N}*\mathbf{e}_{\ell}\hat{\mathcal{M}}*\hat{\mathbf{h}}_{\epsilon,M}\right\Vert _{p^{M},q}
\end{equation}
Hence, by the ultrametric inequality, the upper bound (\ref{eq:WTT - Ultrametric inequality-1})
is in fact an equality: 
\begin{align*}
\left\Vert \mathbf{e}_{\ell}\left(\hat{\phi}_{N}-\left(\hat{\phi}_{N}*\hat{\mathcal{M}}*\hat{\mathbf{h}}_{\epsilon,M}\right)\right)\right\Vert _{p^{M},q} & =\max\left\{ \left\Vert \mathbf{e}_{\ell}\hat{\phi}_{N}\right\Vert _{p^{M},q},\left\Vert \hat{\phi}_{N}*\mathbf{e}_{\ell}\hat{\mathcal{M}}*\hat{\mathbf{h}}_{\epsilon,M}\right\Vert _{p^{M},q}\right\} \\
 & =1
\end{align*}
With this, (\ref{eq:MD WTT - Target of attack}) becomes: 
\begin{align*}
\epsilon & >\left\Vert \mathbf{e}_{\ell}\left(\hat{\phi}_{N}*\left(\vec{\mathbf{1}}_{\mathbf{0}}-\left(\hat{\mathcal{M}}*\hat{\mathbf{h}}_{\epsilon,M}\right)\right)\right)\right\Vert _{p^{M},q}\\
 & =\left\Vert \mathbf{e}_{\ell}\hat{\phi}_{N}-\mathbf{e}_{\ell}\left(\hat{\phi}_{N}*\hat{\mathcal{M}}*\hat{\mathbf{h}}_{\epsilon,M}\right)\right\Vert _{p^{M},q}\\
 & =1
\end{align*}
But $\epsilon<1$!\textemdash we have arrived at our contradiction.

Consequently, the existence of the ``zero'' $\mathbf{z}_{0}$ precludes
the translates of $\hat{\mathcal{M}}$ from being dense in $c_{0}\left(\hat{\mathbb{Z}}_{p}^{r},\mathbb{C}_{q}^{d}\right)$.

Q.E.D.
\begin{rem}
Versions of \textbf{Theorem \ref{thm:MD pq-adic WTT for continuous functions}}
likely hold for arbitrary thick $\left(p,q\right)$-adic measures,
but I will not pursue the matter here.
\end{rem}

\subsection{\emph{\label{subsec:5.4.4More-Fourier-Resummation}More} Fourier
Resummation Lemmata}

As in the one-dimensional case, we will need to avail ourselves of
various Fourier Resummation Lemmata. First, of course, we need the
multi-dimensional definitions: 
\begin{defn}
A $\left(p,q\right)$-adic function $\hat{\mu}:\hat{\mathbb{Z}}_{p}^{r}\rightarrow\mathbb{C}_{q}^{\rho,c}$
is said to be:

\vphantom{}

I. \textbf{Radial} whenever:\index{thick measure!radial} 
\[
\hat{\mu}\left(\mathbf{t}\right)=\hat{\mu}\left(\left|t_{1}\right|_{p_{1}}^{-1},\ldots,\left|t_{r}\right|_{p_{r}}^{-1}\right),\textrm{ }\forall\mathbf{t}\in\hat{\mathbb{Z}}_{p}^{r}\backslash\left\{ \mathbf{0}\right\} 
\]

\vphantom{}

II. \textbf{Magnitudinal} \index{thick measure!magnitudinal}whenever
there is a function $\kappa:\mathbb{N}_{0}^{r}\rightarrow\mathbb{C}_{q}^{\rho,c}$
so that: 
\[
\hat{\mu}\left(\mathbf{t}\right)=\begin{cases}
\kappa\left(\mathbf{0}\right) & \textrm{if }\mathbf{t}=\mathbf{0}\\
\sum_{\mathbf{m}=\mathbf{0}}^{p^{-v_{p}\left(\mathbf{t}\right)}-1}\kappa\left(\mathbf{m}\right)e^{-2\pi i\left(\mathbf{m}\cdot\mathbf{t}\right)} & \textrm{else}
\end{cases},\textrm{ }\forall\mathbf{t}\in\hat{\mathbb{Z}}_{p}^{r}
\]

\vphantom{}

III. \textbf{Radially-Magnitudinal} \index{thick measure!radial-magnitudinal}
whenever there is a radial function $\hat{\nu}:\hat{\mathbb{Z}}_{p}^{r}\rightarrow\mathbb{C}_{q}^{\left(r,c\right)}$
and a magnitudinal function $\hat{\eta}:\hat{\mathbb{Z}}_{p}^{r}\rightarrow\mathbb{C}_{q}^{\left(\rho,r\right)}$
so that: 
\begin{equation}
\hat{\mu}\left(\mathbf{t}\right)=\hat{\eta}\left(\mathbf{t}\right)\hat{\nu}\left(\mathbf{t}\right),\textrm{ }\forall\mathbf{t}\in\hat{\mathbb{Z}}_{p}^{r}
\end{equation}
where the right-hand side, as a product of a $\rho\times r$ matrix-valued
function and an $r\times c$ matrix-valued function is a $\rho\times c$
matrix-valued function.
\end{defn}
\begin{defn}
We say a thick measure $\mathcal{M}$ is radial (resp. magnitudinal,
radially-magnitudinal) whenever $\hat{\mathcal{M}}\left(\mathbf{t}\right)$
is radial (resp. magnitudinal, radially-magnitudinal). 
\end{defn}
\begin{defn}
Consider a function $\kappa:\mathbb{N}_{0}^{r}\rightarrow\mathbb{C}_{q}^{d,d}$.

\vphantom{}

I. We say $\kappa$ is \textbf{$p$-adically structured} / \textbf{has
$p$-adic structure }whenever:\index{$p$-adic!structure} 
\begin{equation}
\kappa\left(p\mathbf{n}+\mathbf{j}\right)=\mathbf{M}_{\mathbf{j}}\kappa_{H}\left(\mathbf{n}\right)\mathbf{M}_{\mathbf{0}}^{-1},\textrm{ }\forall\mathbf{j}\in\mathbb{Z}^{r}/p\mathbb{Z}^{r},\textrm{ }\forall\mathbf{n}\in\mathbb{N}_{0}^{r}\label{eq:Definition of P-adic structure}
\end{equation}
where $\left\{ \mathbf{M}_{\mathbf{j}}\right\} _{\mathbf{j}\in\mathbb{Z}^{r}/p\mathbb{Z}^{r}}$
are invertible matrices in $\mathbb{C}_{q}^{d,d}$.

\vphantom{}

II. We say $\kappa$\textbf{ }is \index{tame}\textbf{$\left(p,K\right)$-adically
tame} on a set $X\subseteq\mathbb{Z}_{p}^{r}$ whenever $\lim_{n\rightarrow\infty}\left\Vert \kappa\left(\left[\mathbf{z}\right]_{p^{n}}\right)\right\Vert _{K}=0$
for all $\mathbf{z}\in X$. We do not mention $X$ when $X=\mathbb{Z}_{p}^{r}$.
If $K$ is a $q$-adic field, we say $\kappa$ is $\left(p,q\right)$-adically
tame on $X$; if $K$ is archimedean, we say $\kappa$ is $\left(p,\infty\right)$-adically
tame on $X$. 
\end{defn}
\begin{prop}
\label{prop:MD p-adic structure prop}A function $\kappa:\mathbb{N}_{0}^{r}\rightarrow\mathbb{C}_{q}^{d,d}$
is $p$-adically structured if and only if, for every $m\geq0$: 
\begin{equation}
\kappa\left(\mathbf{n}+\mathbf{j}p^{m}\right)=\mathbf{M}_{\mathbf{n}}\kappa_{H}\left(\mathbf{j}\right)\mathbf{M}_{\mathbf{0}}^{-1},\textrm{ }\forall\mathbf{j}\in\mathbb{Z}^{r}/p\mathbb{Z}^{r},\textrm{ }\forall\mathbf{n}\leq p^{m}-1
\end{equation}
where $\mathbf{M}_{\mathbf{n}}=\mathbf{M}_{\mathbf{j}_{1}}\times\cdots\times\mathbf{M}_{\mathbf{j}_{\lambda_{p}\left(\mathbf{n}\right)}}$,
where $\mathbf{J}=\left(\mathbf{j}_{1},\ldots,\mathbf{j}_{\lambda_{p}\left(\mathbf{n}\right)}\right)\in\textrm{String}^{r}\left(p\right)$
is the shortest block string representing $\mathbf{n}$. 
\end{prop}
Proof:

Write: 
\[
\mathbf{n}=\left(n_{1},\ldots,n_{r}\right)
\]
\[
\mathbf{j}=\left(j_{1},\ldots,j_{r}\right)
\]
where: 
\[
n_{\ell}=n_{\ell,0}+n_{\ell,1}p+\cdots+n_{\ell,\lambda_{p}\left(n_{\ell}\right)-1}p^{\lambda_{p}\left(n_{\ell}\right)-1},\textrm{ }\forall\ell\in\left\{ 1,\ldots,r\right\} 
\]
and: 
\[
\mathbf{n}+\mathbf{j}p^{m}=\left(n_{1}+j_{1}p^{m},\ldots,n_{r}+j_{r}p^{m}\right)
\]
where, for each $\ell$, $\lambda_{p}\left(n_{\ell}\right)-1<m$.
The rest is just like the one-dimensional case.

Q.E.D.

\vphantom{}Next up, the multi-dimensional resummation lemmata.
\begin{lem}
\label{lem:MD magnitude F Resum Lemma}Let $\hat{\mathcal{M}}:\hat{\mathbb{Z}}_{p}^{r}\rightarrow\mathbb{C}_{q}^{d,d}$
be the Fourier-Stieltjes transform of a magnitudinal thick measure,
and let $\kappa$ have $p$-adic structure. Then, for all $N\geq0$
and all $\mathbf{z}\in\mathbb{Z}_{p}^{r}$: 
\begin{align}
\tilde{\mathcal{M}}_{N}\left(\mathbf{z}\right) & =p^{rN}\kappa\left(\left[\mathbf{z}\right]_{p^{N}}\right)-\sum_{n=0}^{N-1}\sum_{\mathbf{j}>\mathbf{0}}^{p-1}p^{rn}\kappa\left(\left[\mathbf{z}\right]_{p^{n}}+\mathbf{j}p^{n+1}\right)\label{eq:MD magnitudinal resummation formula}
\end{align}
\end{lem}
Proof: We begin with: 
\[
\hat{\mathcal{M}}\left(\mathbf{t}\right)=\sum_{\mathbf{m}=\mathbf{0}}^{p^{-v_{p}\left(\mathbf{t}\right)}-1}\kappa\left(\mathbf{m}\right)e^{-2\pi i\left(\mathbf{m}\cdot\mathbf{t}\right)}
\]
Consequently: 
\begin{align*}
\tilde{\mathcal{M}}_{N}\left(\mathbf{z}\right) & =\sum_{\left\Vert \mathbf{t}\right\Vert _{p}\leq p^{N}}\hat{\mathcal{M}}\left(\mathbf{t}\right)e^{2\pi i\left\{ \mathbf{t}\mathbf{z}\right\} _{p}}\\
 & =\kappa\left(\mathbf{0}\right)+\sum_{n=1}^{N}\sum_{\left\Vert \mathbf{t}\right\Vert _{p}=p^{n}}\left(\sum_{\mathbf{m}=\mathbf{0}}^{p^{n}-1}\kappa\left(\mathbf{m}\right)e^{-2\pi i\left(\mathbf{m}\cdot\mathbf{t}\right)}\right)e^{2\pi i\left\{ \mathbf{t}\mathbf{z}\right\} _{p}}\\
 & =\kappa\left(\mathbf{0}\right)+\sum_{n=1}^{N}\sum_{\mathbf{m}=\mathbf{0}}^{p^{n}-1}\kappa\left(\mathbf{m}\right)\sum_{\left\Vert \mathbf{t}\right\Vert _{p}=p^{n}}e^{2\pi i\left\{ \mathbf{t}\left(\mathbf{z}-\mathbf{m}\right)\right\} _{p}}\\
 & =\kappa\left(\mathbf{0}\right)+\sum_{n=1}^{N}\sum_{\mathbf{m}=\mathbf{0}}^{p^{n}-1}\kappa\left(\mathbf{m}\right)\left(p^{rn}\left[\mathbf{z}\overset{p^{n}}{\equiv}\mathbf{m}\right]-p^{r\left(n-1\right)}\left[\mathbf{z}\overset{p^{n-1}}{\equiv}\mathbf{m}\right]\right)\\
 & =\kappa\left(\mathbf{0}\right)+\sum_{n=1}^{N}\left(p^{rn}\kappa\left(\left[\mathbf{z}\right]_{p^{n}}\right)-p^{r\left(n-1\right)}\sum_{\mathbf{m}=\mathbf{0}}^{p^{n}-1}\kappa\left(\mathbf{m}\right)\left[\mathbf{z}\overset{p^{n-1}}{\equiv}\mathbf{m}\right]\right)
\end{align*}
Because $\kappa$ has $p$-adic structure, we can then write:
\begin{align*}
\sum_{\mathbf{m}=\mathbf{0}}^{p^{n}-1}\kappa\left(\mathbf{m}\right)\left[\mathbf{z}\overset{p^{n-1}}{\equiv}\mathbf{m}\right] & =\sum_{\mathbf{j}=\mathbf{0}}^{p-1}\sum_{\mathbf{m}=\mathbf{0}}^{p^{n-1}-1}\kappa\left(\mathbf{m}+\mathbf{j}p^{n}\right)\left[\mathbf{z}\overset{p^{n-1}}{\equiv}\mathbf{m}+\mathbf{j}p^{n}\right]\\
 & =\sum_{\mathbf{j}=\mathbf{0}}^{p-1}\kappa\left(\left[\mathbf{z}\right]_{p^{n-1}}+\mathbf{j}p^{n}\right)
\end{align*}
and so: 
\begin{equation}
\sum_{\left\Vert \mathbf{t}\right\Vert _{p}=p^{n}}\sum_{\mathbf{m}=\mathbf{0}}^{p^{n}-1}\kappa\left(\mathbf{m}\right)e^{-2\pi i\mathbf{t}\cdot\mathbf{m}}e^{2\pi i\left\{ \mathbf{t}\mathbf{z}\right\} _{p}}=p^{rn}\kappa\left(\left[\mathbf{z}\right]_{p^{n}}\right)-p^{r\left(n-1\right)}\sum_{\mathbf{j}=\mathbf{0}}^{p-1}\kappa\left(\left[\mathbf{z}\right]_{p^{n-1}}+\mathbf{j}p^{n}\right)\label{eq:MD Level set Fourier sum of a magnitudinal multiplier}
\end{equation}
Consequently: 
\begin{align*}
\tilde{\mathcal{M}}_{N}\left(\mathbf{z}\right) & =\kappa\left(\mathbf{0}\right)+\sum_{n=1}^{N}\left(p^{rn}\kappa\left(\left[\mathbf{z}\right]_{^{n}}\right)-p^{r\left(n-1\right)}\sum_{\mathbf{m}=\mathbf{0}}^{p^{n}-1}\kappa\left(\mathbf{m}\right)\left[\mathbf{z}\overset{p^{n-1}}{\equiv}\mathbf{m}\right]\right)\\
 & =\kappa\left(\mathbf{0}\right)+\sum_{n=1}^{N}\left(p^{rn}\kappa\left(\left[\mathbf{z}\right]_{p^{n}}\right)-p^{r\left(n-1\right)}\sum_{\mathbf{j}=\mathbf{0}}^{p-1}\kappa\left(\left[\mathbf{z}\right]_{p^{n-1}}+\mathbf{j}p^{n}\right)\right)\\
 & =\sum_{n=0}^{N}p^{rn}\kappa\left(\left[\mathbf{z}\right]_{p^{n}}\right)-\sum_{n=0}^{N-1}\sum_{\mathbf{j}=\mathbf{0}}^{p-1}p^{rn}\kappa\left(\left[\mathbf{z}\right]_{p^{n}}+\mathbf{j}p^{n+1}\right)\\
 & =\sum_{n=0}^{N}p^{rn}\kappa\left(\left[\mathbf{z}\right]_{p^{n}}\right)-\sum_{n=0}^{N-1}p^{rn}\kappa\left(\left[\mathbf{z}\right]_{p^{n}}\right)-\sum_{n=0}^{N-1}\sum_{\mathbf{j}>\mathbf{0}}^{p-1}p^{rn}\kappa\left(\left[\mathbf{z}\right]_{p^{n}}+\mathbf{j}p^{n+1}\right)\\
 & =p^{rN}\kappa\left(\left[\mathbf{z}\right]_{p^{N}}\right)-\sum_{n=0}^{N-1}\sum_{\mathbf{j}>\mathbf{0}}^{p-1}p^{rn}\kappa\left(\left[\mathbf{z}\right]_{p^{n}}+\mathbf{j}p^{n+1}\right)
\end{align*}

Q.E.D. 
\begin{lem}
\label{lem:MD radial-magnitude F Resum Lemma}Let $\hat{\mathcal{M}}:\hat{\mathbb{Z}}_{p}^{r}\rightarrow\mathbb{C}_{q}^{d,d}$
be the Fourier-Stieltjes transform of a radially-magnitudinal thick
measure, and let $\kappa$ have $p$-adic structure. Then, for all
$N\geq0$ and all $\mathbf{z}\in\mathbb{Z}_{p}^{r}$: 
\begin{align}
\tilde{\mathcal{M}}_{N}\left(\mathbf{z}\right) & =p^{rN}\kappa\left(\left[\mathbf{z}\right]_{p^{N}}\right)\hat{\nu}\left(\frac{1}{p^{N}}\right)+\sum_{n=0}^{N-1}p^{rn}\kappa\left(\left[\mathbf{z}\right]_{p^{n}}\right)\left(\hat{\nu}\left(\frac{1}{p^{n}}\right)-\hat{\nu}\left(\frac{1}{p^{n+1}}\right)\right)\label{eq:MD radially-magnitudinal resummation formula}\\
 & -\sum_{n=0}^{N-1}p^{rn}\sum_{\mathbf{j}>\mathbf{0}}^{p-1}\kappa\left(\left[\mathbf{z}\right]_{p^{n}}+\mathbf{j}p^{n}\right)\hat{\nu}\left(\frac{1}{p^{n+1}}\right)\nonumber 
\end{align}
\end{lem}
Proof: Write: 
\[
\hat{\mathcal{M}}\left(\mathbf{t}\right)=\sum_{\mathbf{m}=\mathbf{0}}^{p^{-v_{p}\left(\mathbf{t}\right)}-1}\kappa\left(\mathbf{m}\right)\hat{\nu}\left(\mathbf{t}\right)e^{-2\pi i\left(\mathbf{m}\cdot\mathbf{t}\right)}
\]
Then: 
\begin{align*}
\tilde{\mathcal{M}}_{N}\left(\mathbf{z}\right) & =\sum_{\left\Vert \mathbf{t}\right\Vert _{p}\leq p^{N}}\left(\sum_{\mathbf{m}=\mathbf{0}}^{p^{-v_{p}\left(\mathbf{t}\right)}-1}\kappa\left(\mathbf{m}\right)\hat{\nu}\left(\mathbf{t}\right)e^{-2\pi i\left(\mathbf{m}\cdot\mathbf{t}\right)}\right)e^{2\pi i\left\{ \mathbf{t}\mathbf{z}\right\} _{p}}\\
 & =\kappa\left(\mathbf{0}\right)\hat{\nu}\left(\mathbf{0}\right)+\sum_{n=1}^{N}\sum_{\left\Vert \mathbf{t}\right\Vert _{p}=p^{n}}\left(\sum_{\mathbf{m}=\mathbf{0}}^{p^{-v_{p}\left(\mathbf{t}\right)}-1}\kappa\left(\mathbf{m}\right)\hat{\nu}\left(\mathbf{t}\right)e^{-2\pi i\left(\mathbf{m}\cdot\mathbf{t}\right)}\right)e^{2\pi i\left\{ \mathbf{t}\mathbf{z}\right\} _{p}}\\
\left(\textrm{use }(\ref{eq:First MD level set summation identity})\right); & =\kappa\left(\mathbf{0}\right)\hat{\nu}\left(\mathbf{0}\right)+\sum_{n=1}^{N}\sum_{\left\Vert \mathbf{t}\right\Vert _{p}=p^{n}}\left(\sum_{\mathbf{m}=\mathbf{0}}^{p^{n}-1}\kappa\left(\mathbf{m}\right)e^{-2\pi i\left(\mathbf{m}\cdot\mathbf{t}\right)}\right)e^{2\pi i\left\{ \mathbf{t}\mathbf{z}\right\} _{p}}\hat{\nu}\left(\frac{1}{p^{n}}\right)
\end{align*}
Using (\ref{eq:MD Level set Fourier sum of a magnitudinal multiplier})
from the proof of \textbf{Lemma \ref{lem:MD magnitude F Resum Lemma}}
gives us:
\begin{align*}
\tilde{\mathcal{M}}_{N}\left(\mathbf{z}\right) & =\kappa\left(\mathbf{0}\right)\hat{\nu}\left(\mathbf{0}\right)+\sum_{n=1}^{N}\sum_{\left\Vert \mathbf{t}\right\Vert _{p}=p^{n}}\left(\sum_{\mathbf{m}=\mathbf{0}}^{p^{n}-1}\kappa\left(\mathbf{m}\right)e^{-2\pi i\left(\mathbf{m}\cdot\mathbf{t}\right)}\right)e^{2\pi i\left\{ \mathbf{t}\mathbf{z}\right\} _{p}}\hat{\nu}\left(\frac{1}{p^{n}}\right)\\
 & =\kappa\left(\mathbf{0}\right)\hat{\nu}\left(\mathbf{0}\right)+\sum_{n=1}^{N}\left(p^{rn}\kappa\left(\left[\mathbf{z}\right]_{p^{n}}\right)-p^{r\left(n-1\right)}\sum_{\mathbf{j}=\mathbf{0}}^{p-1}\kappa\left(\left[\mathbf{z}\right]_{p^{n-1}}+\mathbf{j}p^{n-1}\right)\right)\hat{\nu}\left(\frac{1}{p^{n}}\right)\\
 & =\sum_{n=0}^{N}p^{rn}\kappa\left(\left[\mathbf{z}\right]_{p^{n}}\right)\hat{\nu}\left(\frac{1}{p^{n}}\right)-\sum_{n=0}^{N-1}p^{rn}\sum_{\mathbf{j}=\mathbf{0}}^{p-1}\kappa\left(\left[\mathbf{z}\right]_{p^{n}}+\mathbf{j}p^{n}\right)\hat{\nu}\left(\frac{1}{p^{n+1}}\right)\\
 & =\sum_{n=0}^{N}p^{rn}\kappa\left(\left[\mathbf{z}\right]_{p^{n}}\right)\hat{\nu}\left(\frac{1}{p^{n}}\right)-\sum_{n=0}^{N-1}p^{rn}\kappa\left(\left[\mathbf{z}\right]_{p^{n}}\right)\hat{\nu}\left(\frac{1}{p^{n+1}}\right)\\
 & -\sum_{n=0}^{N-1}p^{rn}\sum_{\mathbf{j}>\mathbf{0}}^{p-1}\kappa\left(\left[\mathbf{z}\right]_{p^{n}}+\mathbf{j}p^{n}\right)\hat{\nu}\left(\frac{1}{p^{n+1}}\right)
\end{align*}
and so: 
\begin{align*}
\tilde{\mathcal{M}}_{N}\left(\mathbf{z}\right) & =p^{rN}\kappa\left(\left[\mathbf{z}\right]_{p^{N}}\right)\hat{\nu}\left(\frac{1}{p^{N}}\right)+\sum_{n=0}^{N-1}p^{rn}\kappa\left(\left[\mathbf{z}\right]_{p^{n}}\right)\left(\hat{\nu}\left(\frac{1}{p^{n}}\right)-\hat{\nu}\left(\frac{1}{p^{n+1}}\right)\right)\\
 & -\sum_{n=0}^{N-1}p^{rn}\sum_{\mathbf{j}>\mathbf{0}}^{p-1}\kappa\left(\left[\mathbf{z}\right]_{p^{n}}+\mathbf{j}p^{n}\right)\hat{\nu}\left(\frac{1}{p^{n+1}}\right)
\end{align*}

Q.E.D. 
\begin{prop}
\label{prop:MD convolution with v_p identity}Let $\hat{\mathcal{M}}:\hat{\mathbb{Z}}_{p}^{r}\rightarrow\mathbb{C}_{q}^{d,d}$
be any function. Then: 
\begin{equation}
\sum_{0<\left\Vert \mathbf{t}\right\Vert _{p}\leq p^{N}}v_{p}\left(\mathbf{t}\right)\hat{\mathcal{M}}\left(\mathbf{t}\right)e^{2\pi i\left\{ \mathbf{t}\mathbf{z}\right\} _{p}}\overset{\mathbb{C}_{q}^{d,d}}{=}-N\tilde{\mathcal{M}}_{N}\left(\mathbf{z}\right)+\sum_{n=0}^{N-1}\tilde{\mathcal{M}}_{n}\left(\mathbf{z}\right)\label{eq:MD Fourier sum of v_p times mu-hat}
\end{equation}
\end{prop}
Proof: 
\begin{align*}
\sum_{0<\left\Vert \mathbf{t}\right\Vert _{p}\leq p^{N}}v_{p}\left(\mathbf{t}\right)\hat{\mathcal{M}}\left(\mathbf{t}\right)e^{2\pi i\left\{ \mathbf{t}\mathbf{z}\right\} } & =\sum_{n=1}^{N}\sum_{\left\Vert \mathbf{t}\right\Vert _{p}=p^{n}}v_{p}\left(\mathbf{t}\right)\hat{\mathcal{M}}\left(\mathbf{t}\right)e^{2\pi i\left\{ \mathbf{t}\mathbf{z}\right\} _{p}}\\
 & =-\sum_{n=1}^{N}n\sum_{\left\Vert \mathbf{t}\right\Vert _{p}=p^{n}}\hat{\mathcal{M}}\left(\mathbf{t}\right)e^{2\pi i\left\{ \mathbf{t}\mathbf{z}\right\} _{p}}\\
 & =-\sum_{n=1}^{N}n\left(\tilde{\mathcal{M}}_{n}\left(\mathbf{z}\right)-\tilde{\mathcal{M}}_{n-1}\left(\mathbf{z}\right)\right)
\end{align*}
The rest of the proof is formally identical to its one-dimensional
counterpart (\textbf{Proposition \ref{prop:v_p of t times mu hat sum}}).

Q.E.D.

\chapter{\label{chap:6 A-Study-of}A Study of $\chi_{H}$ - The Multi-Dimensional
Case}

\includegraphics[scale=0.45]{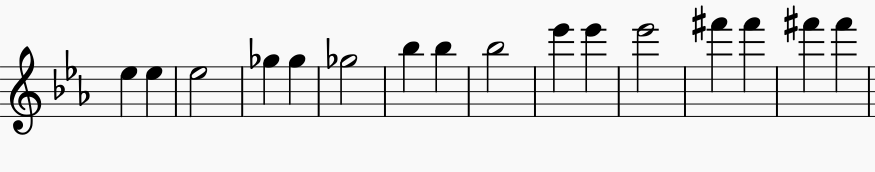}

\vphantom{}

IN THIS CHAPTER, UNLESS STATED OTHERWISE, WE ASSUME $p$ IS A PRIME
AND THAT $H$ IS A CONTRACTING, SEMI-BASIC $p$-SMOOTH $d$-DIMENSIONAL
DEPTH $r$ HYDRA MAP WHICH FIXES $\mathbf{0}$.

\vphantom{}

The layout of Chapter 6 is much like that of its one-dimensional predecessor\textemdash Chapter
4. After an initial string of notational definitions and preparatory
work, we investigate $\chi_{H}$ by squeezing out explicit ``asymptotic
formulae'' for the Fourier transform of the $N$th truncation of
$\chi_{H}$; this is the content of Subsection \ref{subsec:6.2.1 -and-},
with the formulae themselves occurring in \textbf{Theorem \ref{thm:MD N,t asympotics for Chi_H,N hat}}.
Just as $1-\alpha_{H}\left(0\right)$ was the key determinative quantity
for the one-dimensional case, so too will $\mathbf{I}_{d}-\alpha_{H}\left(\mathbf{0}\right)$
be a key quantity for the multi-dimensional case, where $\alpha_{H}$
(defined at the start of Section \ref{sec:6.1. Preparatory-Work-(Again)})
is the multi-dimensional generalization of the one-dimensional $\alpha_{H}$.
However, the non-commutativity of matrix multiplication will, unfortunately,
complicate matters somewhat. As a result, at the beginning of \ref{sec:6.1. Preparatory-Work-(Again)},
I introduce a qualitative notion of the \textbf{commutativity }of
a $p$-Hydra map $H$. This condition turns out to be exactly what
is needed in order to minimize the computational differences between
the one-dimensional case and the multi-dimensional case currently
under consideration.

The primary consequence of this issue is that I have only been able
to establish quasi-integrability results and formulae for Fourier
transforms of the multi-dimensional $\chi_{H}$ in the case where
$H$ is commutative. A treatment of this case is given in Subsection
\ref{subsec:6.2.3 Multi-Dimensional--=00003D000026}, following the
$\alpha_{H}\left(\mathbf{0}\right)=\mathbf{I}_{d}$ case dealt with
in \ref{subsec:6.2.2 Multi-Dimensional--=00003D000026}. That being
said, as a headache-preventing prophylaxis for future explorers of
this subject, \emph{I have taken the liberty of doing all the major
computations without relying on the assumption that $H$ is commutative}.
For ease of readability, alongside these non-commutative cases, I
also state the simpler commutative cases for all significant formulae.
Subsection \ref{subsec:6.2.4 Multi-Dimensional--=00003D000026} contains
additional computations leading all the way up to the doorstep of
a proof of the quasi-integrability of $\chi_{H}$ for non-commutative
$H$, leaving the final step as a conjecture to be tackled in future
work.

\section{\label{sec:6.1. Preparatory-Work-(Again)}Preparatory Work (Again)}

This section is just a multi-dimensional copy of its predecessor in
Section \pageref{sec:4.1 Preparatory-Work--}.
\begin{defn}
We define the functions $\alpha_{H}:\hat{\mathbb{Z}}_{p}^{r}\rightarrow\textrm{GL}_{d}\left(\overline{\mathbb{Q}}\right)$
and $\beta_{H}:\hat{\mathbb{Z}}_{p}^{r}\rightarrow\overline{\mathbb{Q}}^{d}$
by:\nomenclature{$\alpha_{H}\left(\mathbf{t}\right)$}{ }\nomenclature{$\beta_{H}\left(\mathbf{t}\right)$}{ }

\begin{equation}
\alpha_{H}\left(\mathbf{t}\right)\overset{\textrm{def}}{=}\frac{1}{p^{r}}\sum_{\mathbf{j}=\mathbf{0}}^{p-1}\mathbf{D}_{\mathbf{j}}^{-1}\mathbf{A}_{\mathbf{j}}e^{-2\pi i\mathbf{j}\cdot\mathbf{t}}\label{eq:MD definition of alpha_H}
\end{equation}
\begin{equation}
\beta_{H}\left(\mathbf{t}\right)\overset{\textrm{def}}{=}\frac{1}{p^{r}}\sum_{\mathbf{j}=\mathbf{0}}^{p-1}\mathbf{D}_{\mathbf{j}}^{-1}\mathbf{b}_{\mathbf{j}}e^{-2\pi i\mathbf{j}\cdot\mathbf{t}}\label{eq:MD definition of beta_H}
\end{equation}
We then write $\gamma_{H}:\hat{\mathbb{Z}}_{p}^{r}\rightarrow\overline{\mathbb{Q}}^{d}$
to denote:\nomenclature{$\gamma_{H}\left(\mathbf{t}\right)$}{ } 
\begin{equation}
\gamma_{H}\left(\mathbf{t}\right)\overset{\textrm{def}}{=}\left(\alpha_{H}\left(\mathbf{t}\right)\right)^{-1}\beta_{H}\left(\mathbf{t}\right)\label{eq:MD definition of gamma_H}
\end{equation}
\end{defn}
\begin{defn}
We say $H$ is \textbf{non-singular }whenever the matrix $\alpha_{H}\left(\mathbf{j}/p\right)$
is invertible for all $\mathbf{j}\in\mathbb{Z}^{r}/p\mathbb{Z}^{r}$.
\index{Hydra map!non-singular} 
\end{defn}
\vphantom{}

The primary distinction between the work we will do here and the one-dimensional
case comes from the non-commutativity of matrix multiplication. Despite
this, there is a simple\textemdash and not unreasonable\textemdash qualitative
condition we can place on $H$ to make the computational distinctions
between the one-dimensional and multi-dimensional cases all-but-trivial: 
\begin{defn}
We say $H$ is\index{Hydra map!commutative} \textbf{commutative }whenever
$\alpha_{H}\left(\mathbf{0}\right)$ commutes with $H^{\prime}\left(\mathbf{0}\right)$:
\begin{equation}
\alpha_{H}\left(\mathbf{0}\right)H^{\prime}\left(\mathbf{0}\right)=H^{\prime}\left(\mathbf{0}\right)\alpha_{H}\left(\mathbf{0}\right)\label{eq:Definition of H commutativity}
\end{equation}
\end{defn}
\vphantom{}

Next, we introduce the multi-dimensional analogue of $\kappa_{H}$,
here a matrix-valued function.
\begin{defn}
\nomenclature{$\kappa_{H}\left(\mathbf{n}\right)$}{ }We define the
function $\kappa_{H}:\mathbb{N}_{0}^{r}\rightarrow\textrm{GL}_{d}\left(\mathbb{Q}\right)$
by: 
\begin{equation}
\kappa_{H}\left(\mathbf{n}\right)\overset{\textrm{def}}{=}M_{H}\left(\mathbf{n}\right)\left(H^{\prime}\left(\mathbf{0}\right)\right)^{-\lambda_{p}\left(\mathbf{n}\right)},\textrm{ }\forall\mathbf{n}\in\mathbb{N}_{0}^{r}\label{eq:MD definition of Kappa_H}
\end{equation}
\end{defn}
\vphantom{}

\textbf{Proposition \ref{prop:MD generating function}}, given below,
contains the multi-dimensional analogue of the generating function
identities from Chapter \ref{chap:4 A-Study-of}'s \textbf{Proposition
\ref{prop:Generating function identities}}.
\begin{prop}[\textbf{A Generating Function Identity}]
\label{prop:MD generating function}For all $n\geq1$, all scalars
$z$, and all $\mathbf{t}\in\hat{\mathbb{Z}}_{p}^{r}$: 
\begin{equation}
\prod_{m=0}^{n-1}\left(\sum_{\mathbf{j}=\mathbf{0}}^{p-1}\frac{\mathbf{A}_{\mathbf{j}}}{\mathbf{D}_{\mathbf{j}}}z^{\mathbf{j}\cdot p^{m}\mathbf{t}}\right)=\sum_{\mathbf{m}=\mathbf{0}}^{p^{n}-1}M_{H}\left(\mathbf{m}\right)\left(H^{\prime}\left(\mathbf{0}\right)\right)^{n-\lambda_{p}\left(\mathbf{m}\right)}z^{\mathbf{m}\cdot\mathbf{t}}\label{eq:MD M_H partial sum generating identity}
\end{equation}
\end{prop}
Proof: Each term on the right of (\ref{eq:MD M_H partial sum generating identity})
is obtained by taking a product of the form: 
\begin{equation}
\prod_{\ell=0}^{n-1}\frac{\mathbf{A}_{\mathbf{j}_{\ell}}}{\mathbf{D}_{\mathbf{j}_{\ell}}}z^{\mathbf{j}_{\ell}\cdot p^{\ell}\mathbf{t}}=\left(\prod_{\ell=0}^{n-1}\frac{\mathbf{A}_{\mathbf{j}_{\ell}}}{\mathbf{D}_{\mathbf{j}_{\ell}}}\right)z^{\mathbf{t}\cdot\sum_{\ell=0}^{n-1}p^{\ell}\mathbf{j}_{\ell}}
\end{equation}
for some choice of $\mathbf{j}_{1},\ldots,\mathbf{j}_{n-1}\in\mathbb{Z}^{r}/p\mathbb{Z}^{r}$.

Here, writing $\mathbf{j}_{\ell}=\left(j_{\ell,1},\ldots,j_{\ell,r}\right)$
we have: 
\begin{equation}
\sum_{\ell=0}^{n-1}\mathbf{j}_{\ell}p^{\ell}=\left(\sum_{\ell=0}^{n-1}j_{\ell,1}p^{\ell},\ldots,\sum_{\ell=0}^{n-1}j_{\ell,r}p^{\ell}\right)
\end{equation}
So, writing $\mathbf{J}$ to denote the block string $\left(\mathbf{j}_{0},\ldots,\mathbf{j}_{n-1}\right)$,
we define $\mathbf{m}\in\mathbb{N}_{0}^{r}$ by:
\begin{equation}
\mathbf{m}\overset{\textrm{def}}{=}\sum_{\ell=0}^{n-1}\mathbf{j}_{\ell}p^{\ell}\label{eq:definition of bold m}
\end{equation}
so that $\mathbf{m}$ is then represented by $\mathbf{J}$. As such:
\begin{equation}
\prod_{m=0}^{n-1}\left(\sum_{\mathbf{j}=\mathbf{0}}^{-1}\frac{\mathbf{A}_{\mathbf{j}}}{\mathbf{D}_{\mathbf{j}}}z^{\mathbf{j}\cdot\left(p^{m}\mathbf{t}\right)}\right)=\sum_{\begin{array}{c}
\mathbf{J}\in\textrm{String}^{r}\left(p\right)\\
\left|\mathbf{J}\right|=n
\end{array}}\left(\frac{\mathbf{A}_{\mathbf{j}_{0}}}{\mathbf{D}_{\mathbf{j}_{0}}}\times\cdots\times\frac{\mathbf{A}_{\mathbf{j}_{n-1}}}{\mathbf{D}_{\mathbf{j}_{n-1}}}\right)z^{\mathbf{J}\cdot\mathbf{t}}
\end{equation}

As an illustrative example, let $n=4$, consider: 
\begin{equation}
\mathbf{J}=\left(\mathbf{j}_{0},\mathbf{j}_{1},\mathbf{0},\mathbf{0}\right)
\end{equation}
and let $\mathbf{m}$ be the tuple represented by this $\mathbf{J}$.
Note that this $\mathbf{J}$ is the unique \emph{length $4$ }block
string representing $\mathbf{m}$; any other length $4$ block string
will have a non-zero tuple in either the $\mathbf{j}_{2}$ slot or
the $\mathbf{j}_{3}$ slot. More generally, there will be a bijection
between the set of all length $n$ block strings and the set of all
$r$-tuples of integers \emph{representable }by said block strings.
This set of $r$-tuples of integers is precisely: 
\begin{equation}
\left\{ \left(m_{1},\ldots,m_{r}\right):0\leq m_{\ell}\leq p^{n}-1\textrm{ }\forall\ell\in\left\{ 1,\ldots,r\right\} \right\} 
\end{equation}
More compactly, this is exactly the range of tuples indicated by the
summation notation: 
\begin{equation}
\sum_{\mathbf{m}=\mathbf{0}}^{p^{N}-1}
\end{equation}
The problem here is that any terminal $\mathbf{0}$-tuples in $\mathbf{J}$
they affect the matrix product associated to $\mathbf{J}$: 
\begin{equation}
\frac{\mathbf{A}_{\mathbf{j}_{0}}}{\mathbf{D}_{\mathbf{j}_{0}}}\times\cdots\times\frac{\mathbf{A}_{\mathbf{j}_{n-1}}}{\mathbf{D}_{\mathbf{j}_{n-1}}}
\end{equation}
even though those $\mathbf{0}$-tuples do \emph{not} affect $z^{\mathbf{J}\cdot\mathbf{t}}$.
For example, for $\mathbf{J}=\left(\mathbf{j}_{0},\mathbf{j}_{1},\mathbf{0},\mathbf{0}\right)$,
the product would be: 
\begin{align*}
\frac{\mathbf{A}_{\mathbf{j}_{0}}}{\mathbf{D}_{\mathbf{j}_{0}}}\times\frac{\mathbf{A}_{\mathbf{j}_{1}}}{\mathbf{D}_{\mathbf{j}_{1}}}\times\frac{\mathbf{A}_{\mathbf{0}}}{\mathbf{D}_{\mathbf{0}}}\times\frac{\mathbf{A}_{\mathbf{0}}}{\mathbf{D}_{\mathbf{0}}} & =M_{H}\left(\left(\mathbf{j}_{0},\mathbf{j}_{1}\right)\right)\times\left(\frac{\mathbf{A}_{\mathbf{0}}}{\mathbf{D}_{\mathbf{0}}}\right)^{2}\\
\left(\mathbf{m}\textrm{ is represented by }\left(\mathbf{j}_{0},\mathbf{j}_{1}\right)\right); & =M_{H}\left(\mathbf{m}\right)\times\left(H^{\prime}\left(\mathbf{0}\right)\right)^{2}
\end{align*}

So, we need to modify the product to take into account terminal $\mathbf{0}$s
in $\mathbf{J}$. To do this, consider the worst-case scenario where
$\mathbf{J}$ is a length-$n$ block string whose every $\mathbf{j}$
is $\mathbf{0}$. Letting $\mathbf{m}$ be the unique $r$-tuple of
integers represented by an arbitrary length-$n$ $\mathbf{J}$, observe
that $\lambda_{p}\left(\mathbf{m}\right)$ is the number of tuples
in $\mathbf{J}$ which are \emph{not }part of a possible run of consecutive
terminal $\mathbf{0}$s. This is because $\lambda_{p}\left(\mathbf{m}\right)$
is the length of the shortest block string representing $\mathbf{m}$.
Consequently, for each $\mathbf{J}$, we have:
\begin{equation}
\frac{\mathbf{A}_{\mathbf{j}_{0}}}{\mathbf{D}_{\mathbf{j}_{0}}}\times\cdots\times\frac{\mathbf{A}_{\mathbf{j}_{n-1}}}{\mathbf{D}_{\mathbf{j}_{n-1}}}=M_{H}\left(\mathbf{m}\right)\left(H^{\prime}\left(\mathbf{0}\right)\right)^{n-\lambda\left(\mathbf{m}\right)}
\end{equation}
which leaves us with: 
\begin{align*}
\prod_{m=0}^{n-1}\left(\sum_{\mathbf{j}=\mathbf{0}}^{-1}\frac{\mathbf{A}_{\mathbf{j}}}{\mathbf{D}_{\mathbf{j}}}z^{\mathbf{j}\cdot\left(^{m}\mathbf{t}\right)}\right) & =\sum_{\begin{array}{c}
\mathbf{J}\in\textrm{String}^{r}\left(p\right)\\
\left|\mathbf{J}\right|=n
\end{array}}\left(\frac{\mathbf{A}_{\mathbf{j}_{0}}}{\mathbf{D}_{\mathbf{j}_{0}}}\times\cdots\times\frac{\mathbf{A}_{\mathbf{j}_{n-1}}}{\mathbf{D}_{\mathbf{j}_{n-1}}}\right)z^{\mathbf{J}\cdot\mathbf{t}}\\
 & =\sum_{\mathbf{m}=\mathbf{0}}^{^{n}-1}M_{H}\left(\mathbf{m}\right)\left(H^{\prime}\left(\mathbf{0}\right)\right)^{n-\lambda\left(\mathbf{m}\right)}z^{\mathbf{m}\cdot\mathbf{t}}
\end{align*}
which is, of course, our sought-after identity.

Q.E.D.

\vphantom{}

In the process of writing this chapter, I took pains to devise a notation
which would make multi-dimensional case's computations as near as
possible to verbatim repeats of their one-dimensional predecessors.
The principal difficulty our notation must overcome is the non-commutativity
of matrix multiplication. This manifests most strongly in the functional
equations for the multi-dimensional $\kappa_{H}$. In order to keep
our computations from spilling out into the margins of the page\textemdash or
beyond\textemdash we will need to introduce a notation for the linear
operator which conjugates $d\times d$ matrices by the $n$th power
of $H^{\prime}\left(\mathbf{0}\right)$.
\begin{defn}
For any $n\in\mathbb{N}_{0}$ and any $\mathbf{A}\in\textrm{GL}_{d}\left(\overline{\mathbb{Q}}\right)$,
$\textrm{GL}_{d}\left(\mathbb{C}\right)$, or $\textrm{GL}_{d}\left(\mathbb{C}_{q}\right)$,
we define: \nomenclature{$\mathcal{C}_{H}\left(\mathbf{A}:n\right)$}{$\overset{\textrm{def}}{=}\left(H^{\prime}\left(\mathbf{0}\right)\right)^{n}\mathbf{A}\left(H^{\prime}\left(\mathbf{0}\right)\right)^{-n}$}
\begin{equation}
\mathcal{C}_{H}\left(\mathbf{A}:n\right)\overset{\textrm{def}}{=}\left(H^{\prime}\left(\mathbf{0}\right)\right)^{n}\mathbf{A}\left(H^{\prime}\left(\mathbf{0}\right)\right)^{-n}\label{eq:Definition of script C_H}
\end{equation}
\end{defn}
\begin{rem}
The most important thing to note is that: 
\begin{equation}
\mathcal{C}_{H}\left(\mathbf{I}_{d}:n\right)=\mathbf{I}_{d}\label{eq:Script C_H of I_d}
\end{equation}
This property is responsible for the impact had on the dynamics of
$H$ when $\alpha_{H}\left(\mathbf{0}\right)=\mathbf{I}_{d}$ (the
multi-dimensional analogue of the ``$\alpha_{H}\left(0\right)=1$''
case.). Additionally, note that $\mathcal{C}_{H}\left(\alpha_{H}\left(\mathbf{0}\right):n\right)=\alpha_{H}\left(\mathbf{0}\right)$
whenever $H$ is commutative.
\end{rem}
\begin{lem}[\textbf{Properties of Multi-Dimensional $\kappa_{H}$}]
\label{lem:properties of MD kappa_H}\ 

\vphantom{}

I. 
\begin{equation}
\sum_{\mathbf{j}>\mathbf{0}}^{p-1}\kappa_{H}\left(\mathbf{j}\right)=\alpha_{H}\left(\mathbf{0}\right)\left(\frac{H^{\prime}\left(\mathbf{0}\right)}{p^{r}}\right)^{-1}-\mathbf{I}_{d}\label{eq:MD kappa H sum in terms of MD alpha}
\end{equation}

\vphantom{}

II. $\kappa_{H}$ is the unique function $\mathbb{N}_{0}^{r}\rightarrow\textrm{GL}_{d}\left(\mathbb{Q}\right)$
satisfying the functional equations: 
\begin{equation}
\kappa_{H}\left(p\mathbf{n}+\mathbf{j}\right)=H_{\mathbf{j}}^{\prime}\left(\mathbf{0}\right)\kappa_{H}\left(\mathbf{n}\right)\left(H^{\prime}\left(\mathbf{0}\right)\right)^{-1},\textrm{ }\forall\mathbf{n}\in\mathbb{N}_{0}^{d},\textrm{ }\forall\mathbf{j}\leq p-1\label{eq:MD Kappa_H functional equations}
\end{equation}
subject to the initial condition $\kappa_{H}\left(\mathbf{0}\right)=\mathbf{I}_{d}$.

\vphantom{}

III. $\kappa_{H}$ satisfies: 
\begin{equation}
\kappa_{H}\left(\mathbf{m}+p^{k}\mathbf{j}\right)=\kappa_{H}\left(\mathbf{m}\right)\left(H^{\prime}\left(\mathbf{0}\right)\right)^{\lambda_{p}\left(\mathbf{m}\right)}\kappa_{H}\left(\mathbf{j}\right)\left(H^{\prime}\left(\mathbf{0}\right)\right)^{-\lambda_{p}\left(\mathbf{m}\right)}\label{eq:MD Kappa_H has P-adic structure}
\end{equation}
for all $k\in\mathbb{N}_{1}$, all $\mathbf{m}\leq p^{k}-1$, and
all $\mathbf{j}\leq p-1$. Equivalently, for these parameters: 
\begin{equation}
\kappa_{H}\left(\mathbf{m}+p^{k}\mathbf{j}\right)=\kappa_{H}\left(\mathbf{m}\right)\mathcal{C}_{H}\left(\kappa_{H}\left(\mathbf{j}\right):\lambda_{p}\left(\mathbf{m}\right)\right)\label{eq:MD Kappa_H has P-adic structure, with script C_H}
\end{equation}

\vphantom{}

IV. If $H$ is semi-basic, then $\kappa_{H}$ is $\left(p,q_{H}\right)$-adically
regular, with: 
\begin{equation}
\lim_{n\rightarrow\infty}\left\Vert \kappa_{H}\left(\left[\mathbf{z}\right]_{p^{n}}\right)\right\Vert _{q_{H}}\overset{\mathbb{R}}{=}0,\textrm{ }\forall\mathbf{z}\in\left(\mathbb{Z}_{p}^{r}\right)^{\prime}\label{eq:MD Semi-basic q-adic decay for Kappa_H}
\end{equation}

\vphantom{}

V. If $H$ is contracting, then: 
\begin{equation}
\lim_{N\rightarrow\infty}\kappa_{H}\left(\left[\mathbf{z}\right]_{p^{N}}\right)\left(H^{\prime}\left(\mathbf{0}\right)\right)^{N}\overset{\mathbb{R}^{d,d}}{=}\mathbf{O}_{d},\textrm{ }\forall\mathbf{z}\in\mathbb{N}_{0}^{r}\label{eq:MD Kappa_H decay when H is contracting}
\end{equation}
\end{lem}
Proof:

I. 
\begin{equation}
\sum_{\mathbf{j}>\mathbf{0}}^{p-1}\kappa_{H}\left(\mathbf{j}\right)=\sum_{\mathbf{j}>\mathbf{0}}^{p-1}M_{H}\left(\mathbf{j}\right)\left(H^{\prime}\left(\mathbf{0}\right)\right)^{-\lambda_{p}\left(\mathbf{j}\right)}
\end{equation}
Since $\lambda_{p}\left(\mathbf{j}\right)=1$ for all $\mathbf{j}\in\left(\mathbb{Z}^{r}/p\mathbb{Z}^{r}\right)\backslash\left\{ \mathbf{0}\right\} $,
and since: 
\begin{equation}
M_{H}\left(\mathbf{j}\right)=\mathbf{D}_{\mathbf{j}}^{-1}\mathbf{A}_{\mathbf{j}}
\end{equation}
we have: 
\begin{align*}
\sum_{\mathbf{j}>\mathbf{0}}^{p-1}\kappa_{H}\left(\mathbf{j}\right) & =\sum_{\mathbf{j}>\mathbf{0}}^{p-1}\mathbf{D}_{\mathbf{j}}^{-1}\mathbf{A}_{\mathbf{j}}\left(H^{\prime}\left(\mathbf{0}\right)\right)^{-1}\\
 & =\frac{1}{p^{r}}\sum_{\mathbf{j}>\mathbf{0}}^{p-1}\mathbf{D}_{\mathbf{j}}^{-1}\mathbf{A}_{\mathbf{j}}\left(\frac{H^{\prime}\left(\mathbf{0}\right)}{p^{r}}\right)^{-1}\\
\left(\mathbf{D}_{\mathbf{0}}^{-1}\mathbf{A}_{\mathbf{0}}=H^{\prime}\left(\mathbf{0}\right)\right); & =\frac{1}{p^{r}}\left(-H^{\prime}\left(\mathbf{0}\right)+\sum_{\mathbf{j}=\mathbf{0}}^{-1}\mathbf{D}_{\mathbf{j}}^{-1}\mathbf{A}_{\mathbf{j}}\right)\left(\frac{H^{\prime}\left(\mathbf{0}\right)}{p^{r}}\right)^{-1}\\
 & =-\frac{H^{\prime}\left(\mathbf{0}\right)}{p^{r}}\left(\frac{H^{\prime}\left(\mathbf{0}\right)}{p^{r}}\right)^{-1}+\underbrace{\left(\frac{1}{p^{r}}\sum_{\mathbf{j}=\mathbf{0}}^{p-1}\mathbf{D}_{\mathbf{j}}^{-1}\mathbf{A}_{\mathbf{j}}\right)}_{\alpha_{H}\left(\mathbf{0}\right)}\left(\frac{H^{\prime}\left(\mathbf{0}\right)}{p^{r}}\right)^{-1}\\
 & =-\mathbf{I}_{d}+\alpha_{H}\left(\mathbf{0}\right)\left(\frac{H^{\prime}\left(\mathbf{0}\right)}{p^{r}}\right)^{-1}
\end{align*}

\vphantom{}

II. Let $\mathbf{n}\in\mathbb{N}_{0}^{r}\backslash\left\{ \mathbf{0}\right\} $
and let $\mathbf{j}\in\mathbb{Z}^{r}/p\mathbb{Z}^{r}$. Then: 
\begin{align*}
\lambda_{p}\left(p\mathbf{n}+\mathbf{j}\right) & =\lambda_{p}\left(\mathbf{n}\right)+1\\
M_{H}\left(p\mathbf{n}+\mathbf{j}\right) & =M_{H}\left(\mathbf{n}\right)\frac{\mathbf{A}_{\mathbf{j}}}{\mathbf{D}_{\mathbf{j}}}=M_{H}\left(\mathbf{n}\right)\mathbf{D}_{\mathbf{j}}^{-1}\mathbf{A}_{\mathbf{j}}
\end{align*}
As such: 
\begin{align*}
\kappa_{H}\left(p\mathbf{n}+\mathbf{j}\right) & =M_{H}\left(p\mathbf{n}+\mathbf{j}\right)\left(H^{\prime}\left(\mathbf{0}\right)\right)^{-\lambda_{p}\left(p\mathbf{n}+\mathbf{j}\right)}\\
 & =\left(\frac{\mathbf{A}_{\mathbf{j}}}{\mathbf{D}_{\mathbf{j}}}M_{H}\left(\mathbf{n}\right)\right)\left(H^{\prime}\left(\mathbf{0}\right)\right)^{-\lambda_{p}\left(\mathbf{n}\right)-1}\\
 & =\mathbf{D}_{\mathbf{j}}^{-1}\mathbf{A}_{\mathbf{j}}\overbrace{\left(M_{H}\left(\mathbf{n}\right)\left(H^{\prime}\left(\mathbf{0}\right)\right)^{-\lambda_{p}\left(\mathbf{n}\right)}\right)}^{\kappa_{H}\left(\mathbf{n}\right)}\left(H^{\prime}\left(\mathbf{0}\right)\right)^{-1}\\
\left(H^{\prime}\left(\mathbf{0}\right)=\mathbf{D}_{\mathbf{0}}^{-1}\mathbf{A}_{\mathbf{0}}\right); & =\mathbf{D}_{\mathbf{j}}^{-1}\mathbf{A}_{\mathbf{j}}\kappa_{H}\left(\mathbf{n}\right)\mathbf{A}_{\mathbf{0}}^{-1}\mathbf{D}_{\mathbf{0}}
\end{align*}
Next: 
\begin{align*}
\kappa_{H}\left(\mathbf{j}\right) & =\begin{cases}
M_{H}\left(\mathbf{0}\right)\left(H^{\prime}\left(\mathbf{0}\right)\right)^{-0} & \textrm{if }\mathbf{j}=\mathbf{0}\\
M_{H}\left(\mathbf{j}\right)\left(H^{\prime}\left(\mathbf{0}\right)\right)^{-1} & \textrm{if }\mathbf{j}\in\left(\mathbb{Z}^{r}/p\mathbb{Z}^{r}\right)\backslash\left\{ \mathbf{0}\right\} 
\end{cases}\\
 & =\begin{cases}
\mathbf{I}_{d} & \textrm{if }\mathbf{j}=\mathbf{0}\\
\frac{\mathbf{A}_{\mathbf{j}}}{\mathbf{D}_{\mathbf{j}}}\frac{\mathbf{A}_{\mathbf{0}}}{\mathbf{D}_{\mathbf{0}}} & \textrm{if }\mathbf{j}\in\left(\mathbb{Z}^{r}/p\mathbb{Z}^{r}\right)\backslash\left\{ \mathbf{0}\right\} 
\end{cases}
\end{align*}
where: 
\begin{equation}
\frac{\mathbf{A}_{\mathbf{j}}}{\mathbf{D}_{\mathbf{j}}}\frac{\mathbf{A}_{\mathbf{0}}}{\mathbf{D}_{\mathbf{0}}}=\mathbf{D}_{\mathbf{j}}^{-1}\mathbf{A}_{\mathbf{j}}\mathbf{D}_{\mathbf{0}}^{-1}\mathbf{A}_{\mathbf{0}}
\end{equation}
Since this shows $\kappa_{H}\left(\mathbf{0}\right)=\mathbf{I}_{d}$,
we can then write: 
\begin{equation}
\kappa_{H}\left(\mathbf{j}\right)=\frac{\mathbf{A}_{\mathbf{j}}}{\mathbf{D}_{\mathbf{j}}}\kappa_{H}\left(\mathbf{0}\right),\textrm{ }\forall\mathbf{j}\leq p-1
\end{equation}
Combining this last equation with the other cases computed above yields
(\ref{eq:MD Kappa_H functional equations}). The uniqueness follows
by an inductive argument, showing that the initial condition $\kappa_{H}\left(\mathbf{0}\right)=\mathbf{I}_{d}$
uniquely determines solutions of the functional equation (\ref{eq:MD Kappa_H functional equations})
on $\mathbb{N}_{0}^{r}$. Finally, we put things in the desired form
by recalling that $H_{\mathbf{j}}^{\prime}\left(\mathbf{0}\right)=\mathbf{D}_{\mathbf{j}}^{-1}\mathbf{A}_{\mathbf{j}}$
for all $\mathbf{j}$, and that $H^{\prime}\left(\mathbf{0}\right)=H_{\mathbf{0}}^{\prime}\left(\mathbf{0}\right)$.

\vphantom{}

III. Let $\mathbf{j}\in\mathbb{Z}^{r}/p\mathbb{Z}^{r}$ and $k\geq1$
be arbitrary, and let $\mathbf{m}\in\mathbb{N}_{0}^{r}$ satisfy $\lambda\left(\mathbf{m}\right)\leq k$
By definition of multi-dimensional $\kappa_{H}$ (equation (\ref{eq:MD definition of Kappa_H})
and using the functional equations for $M_{H}$ (\textbf{Proposition
\ref{prop:MD M_H functional equations}}) and $\lambda_{p}$ (\textbf{Proposition
\ref{prop:MD lambda and digit-number functional equations}})
\begin{align*}
\kappa_{H}\left(\mathbf{m}+p^{k}\mathbf{j}\right) & =M_{H}\left(\mathbf{m}+p^{k}\mathbf{j}\right)\left(H^{\prime}\left(\mathbf{0}\right)\right)^{-\lambda_{p}\left(\mathbf{m}+p^{k}\mathbf{j}\right)}\\
 & =M_{H}\left(\mathbf{m}\right)M_{H}\left(p^{k}\mathbf{j}\right)\left(H^{\prime}\left(\mathbf{0}\right)\right)^{-\lambda_{p}\left(\mathbf{m}\right)-\lambda_{p}\left(\mathbf{j}\right)}\\
\left(M_{H}\left(p^{k}\mathbf{j}\right)=M_{H}\left(\mathbf{j}\right)\right); & =M_{H}\left(\mathbf{m}\right)M_{H}\left(\mathbf{j}\right)\left(H^{\prime}\left(\mathbf{0}\right)\right)^{-\lambda_{p}\left(\mathbf{m}\right)-\lambda_{p}\left(\mathbf{j}\right)}\\
 & =M_{H}\left(\mathbf{m}\right)\underbrace{M_{H}\left(\mathbf{j}\right)\left(H^{\prime}\left(\mathbf{0}\right)\right)^{-\lambda_{p}\left(\mathbf{j}\right)}}_{\kappa_{H}\left(\mathbf{j}\right)}\left(H^{\prime}\left(\mathbf{0}\right)\right)^{-\lambda_{p}\left(\mathbf{m}\right)}\\
 & =M_{H}\left(\mathbf{m}\right)\kappa_{H}\left(\mathbf{j}\right)\left(H^{\prime}\left(\mathbf{0}\right)\right)^{-\lambda_{p}\left(\mathbf{m}\right)}
\end{align*}
Inserting $\mathbf{I}_{d}=\left(H^{\prime}\left(\mathbf{0}\right)\right)^{-\lambda_{p}\left(\mathbf{m}\right)}\left(H^{\prime}\left(\mathbf{0}\right)\right)^{\lambda_{p}\left(\mathbf{m}\right)}$
in between $M_{H}\left(\mathbf{m}\right)$ and $\kappa_{H}\left(\mathbf{j}\right)$
yields:
\begin{align*}
\kappa_{H}\left(\mathbf{m}+p^{k}\mathbf{j}\right) & =\underbrace{M_{H}\left(\mathbf{m}\right)\left(H^{\prime}\left(\mathbf{0}\right)\right)^{-\lambda_{p}\left(\mathbf{m}\right)}}_{\kappa_{H}\left(\mathbf{m}\right)}\left(H^{\prime}\left(\mathbf{0}\right)\right)^{\lambda_{p}\left(\mathbf{m}\right)}\kappa_{H}\left(\mathbf{j}\right)\left(H^{\prime}\left(\mathbf{0}\right)\right)^{-\lambda_{p}\left(\mathbf{m}\right)}\\
 & =\kappa_{H}\left(\mathbf{m}\right)\left(H^{\prime}\left(\mathbf{0}\right)\right)^{\lambda_{p}\left(\mathbf{m}\right)}\kappa_{H}\left(\mathbf{j}\right)\left(H^{\prime}\left(\mathbf{0}\right)\right)^{-\lambda_{p}\left(\mathbf{m}\right)}
\end{align*}
as desired.

\vphantom{}

IV. Applying $q$-adic norm to (\ref{eq:MD Kappa_H has P-adic structure})
gives us:
\begin{align*}
\left\Vert \kappa_{H}\left(\mathbf{m}+p^{k}\mathbf{j}\right)\right\Vert _{q} & \leq\left\Vert \kappa_{H}\left(\mathbf{m}\right)\right\Vert _{q}\left\Vert H^{\prime}\left(\mathbf{0}\right)\right\Vert _{q}^{\lambda_{p}\left(\mathbf{m}\right)}\left\Vert \kappa_{H}\left(\mathbf{j}\right)\right\Vert _{q}\left\Vert H^{\prime}\left(\mathbf{0}\right)\right\Vert _{q}^{-\lambda_{p}\left(\mathbf{m}\right)}\\
 & =\left\Vert \kappa_{H}\left(\mathbf{m}\right)\right\Vert _{q}\left\Vert \kappa_{H}\left(\mathbf{j}\right)\right\Vert _{q}\\
 & =\left\Vert M_{H}\left(\mathbf{m}\right)\left(H^{\prime}\left(\mathbf{0}\right)\right)^{-\lambda_{p}\left(\mathbf{m}\right)}\right\Vert _{q}\left\Vert M_{H}\left(\mathbf{j}\right)\left(H^{\prime}\left(\mathbf{0}\right)\right)^{-\lambda_{p}\left(\mathbf{j}\right)}\right\Vert _{q}\\
 & \leq\left\Vert M_{H}\left(\mathbf{m}\right)\right\Vert _{q}\left\Vert M_{H}\left(\mathbf{j}\right)\right\Vert _{q}\left\Vert H^{\prime}\left(\mathbf{0}\right)\right\Vert _{q}^{-\lambda_{p}\left(\mathbf{m}\right)-\lambda_{p}\left(\mathbf{j}\right)}
\end{align*}
Since $H$ is semi-basic, $\left\Vert H^{\prime}\left(\mathbf{0}\right)\right\Vert _{q}=1$,
and so: 
\begin{equation}
\left\Vert \kappa_{H}\left(\mathbf{m}+p^{k}\mathbf{j}\right)\right\Vert _{q}\leq\left\Vert M_{H}\left(\mathbf{m}\right)\right\Vert _{q}\left\Vert M_{H}\left(\mathbf{j}\right)\right\Vert _{q}
\end{equation}
From this, it follows that for a block string $\mathbf{J}$, the $q$-adic
norm $\left\Vert \kappa_{H}\left(\mathbf{J}\right)\right\Vert _{q}$
tends to $0$ as $\mathbf{J}$ converges ($\left|\mathbf{J}\right|\rightarrow\infty$)
to any infinite block string representing an element of $\left(\mathbb{Z}_{p}^{r}\right)^{\prime}$.
This proves (IV).

\vphantom{}

V. Let $\mathbf{z}\in\mathbb{N}_{0}^{r}$. Then, $\left[\mathbf{z}\right]_{p^{N}}=\mathbf{z}$
for all $N\geq\lambda_{p}\left(\mathbf{z}\right)$, and so: 
\begin{equation}
\kappa_{H}\left(\left[\mathbf{z}\right]_{p^{N}}\right)\left(H^{\prime}\left(\mathbf{0}\right)\right)^{N}=\kappa_{H}\left(\mathbf{z}\right)\left(H^{\prime}\left(\mathbf{0}\right)\right)^{N},\textrm{ }\forall N\geq\lambda_{p}\left(\mathbf{z}\right)
\end{equation}
the right-hand side of which tends to the zero matrix $\mathbf{O}_{d}$
in the topology of $\mathbb{R}^{d,d}$ as $N\rightarrow\infty$, seeing
as $H$ is contracting. This proves (V).

Q.E.D. 
\begin{prop}[\textbf{Summatory Function of $\chi_{H}$}]
\label{prop:MD summatory function formula for Chi_H}
\begin{equation}
\sum_{\mathbf{n}=\mathbf{0}}^{p^{N}-1}\chi_{H}\left(\mathbf{n}\right)=\begin{cases}
\beta_{H}\left(\mathbf{0}\right)Np^{rN} & \textrm{if }\alpha_{H}\left(\mathbf{0}\right)=\mathbf{I}_{d}\\
\frac{\left(\alpha_{H}\left(\mathbf{0}\right)\right)^{N}-1}{\alpha_{H}\left(\mathbf{0}\right)-1}\beta_{H}\left(\mathbf{0}\right)p^{rN} & \textrm{else}
\end{cases}\label{eq:MD Summatory formula for Chi_H}
\end{equation}
\end{prop}
Proof: Let: 
\[
S\left(N\right)=\frac{1}{p^{rN}}\sum_{\mathbf{n}=\mathbf{0}}^{p^{N}-1}\chi_{H}\left(\mathbf{n}\right)
\]
Then, using $\chi_{H}$'s functional equation (\textbf{Lemma \ref{lem:MD Chi_H functional equations and characterization}}):

\begin{align*}
p^{rN}S\left(N\right) & =\sum_{\mathbf{n}=\mathbf{0}}^{p^{N}-1}\chi_{H}\left(\mathbf{n}\right)\\
 & =\sum_{\mathbf{j}=\mathbf{0}}^{p-1}\sum_{\mathbf{n}=\mathbf{0}}^{p^{N-1}-1}\chi_{H}\left(p\mathbf{n}+\mathbf{j}\right)\\
 & =\sum_{\mathbf{j}=\mathbf{0}}^{p-1}\sum_{\mathbf{n}=\mathbf{0}}^{p^{N-1}-1}\frac{\mathbf{A_{j}}\chi_{H}\left(\mathbf{n}\right)+\mathbf{b}_{\mathbf{j}}}{\mathbf{D}_{\mathbf{j}}}\\
 & =\sum_{\mathbf{j}=\mathbf{0}}^{p-1}\left(p^{r\left(N-1\right)}\mathbf{D}_{\mathbf{j}}^{-1}\mathbf{A}_{\mathbf{j}}\frac{1}{p^{r\left(N-1\right)}}\sum_{\mathbf{n}=\mathbf{0}}^{p^{N-1}-1}\chi_{H}\left(\mathbf{n}\right)+p^{r\left(N-1\right)}\mathbf{D}_{\mathbf{j}}^{-1}\mathbf{b}_{\mathbf{j}}\right)\\
 & =p^{r\left(N-1\right)}\left(\sum_{\mathbf{j}=\mathbf{0}}^{p-1}\mathbf{D}_{\mathbf{j}}^{-1}\mathbf{A}_{\mathbf{j}}\right)S\left(N-1\right)+p^{rN}\sum_{\mathbf{j}=\mathbf{0}}^{p-1}\mathbf{D}_{\mathbf{j}}^{-1}\mathbf{b}_{\mathbf{j}}
\end{align*}
Dividing through by $p^{rN}$ gives us: 
\begin{equation}
S_{H}\left(N\right)=\alpha_{H}\mathbf{\left(0\right)}S_{H}\left(N-1\right)+\beta_{H}\left(\mathbf{0}\right)\label{eq:MD Recursive formula for S_H}
\end{equation}
Much like the one-dimensional case, this shows that $S_{H}\left(N\right)$
is the image of $S_{H}\left(0\right)$ under $N$ iterates of an affine
linear map (this time, on $\mathbb{Q}^{d}$); specifically: 
\begin{equation}
\mathbf{x}\mapsto\alpha_{H}\left(\mathbf{0}\right)\mathbf{x}+\beta_{H}\left(\mathbf{0}\right)\label{eq:Affine linear map generating S_H-1}
\end{equation}
This leaves us with two cases:

\vphantom{}

i. Suppose $\alpha_{H}\left(\mathbf{0}\right)=\mathbf{I}_{d}$. Then,
(\ref{eq:Affine linear map generating S_H-1}) reduces to a translation:
\begin{equation}
\mathbf{x}\mapsto\mathbf{x}+\beta_{H}\left(\mathbf{0}\right)
\end{equation}
in which case:

\begin{equation}
S_{H}\left(N\right)=S_{H}\left(0\right)+\beta_{H}\left(\mathbf{0}\right)N\label{eq:MD Explicit formula for S_H of N when alpha =00003D00003D 1}
\end{equation}

\vphantom{}

ii. Suppose $\alpha_{H}\left(\mathbf{0}\right)\neq\mathbf{I}_{d}$,
so that (\ref{eq:Affine linear map generating S_H-1}) is not a translation.
Using the explicit formula for this $N$th iterate, we obtain: 
\begin{equation}
S_{H}\left(N\right)=\left(\alpha_{H}\left(\mathbf{0}\right)\right)^{N}S_{H}\left(0\right)+\left(\alpha_{H}\left(\mathbf{0}\right)-1\right)^{-1}\left(\left(\alpha_{H}\left(\mathbf{0}\right)\right)^{N}-1\right)\beta_{H}\left(\mathbf{0}\right)\label{eq:MD Explicit formula for S_H of N}
\end{equation}
Noting that: 
\begin{align}
S_{H}\left(0\right) & =\sum_{\mathbf{n}=\mathbf{0}}^{p^{0}-1}\chi_{H}\left(\mathbf{n}\right)=\chi_{H}\left(\mathbf{0}\right)=\mathbf{0}
\end{align}
we then have: 
\begin{equation}
\frac{1}{p^{rN}}\sum_{\mathbf{n}=\mathbf{0}}^{p^{N}-1}\chi_{H}\left(\mathbf{n}\right)=S_{H}\left(N\right)=\begin{cases}
\beta_{H}\left(\mathbf{0}\right)N & \textrm{if }\alpha_{H}\left(\mathbf{0}\right)=\mathbf{I}_{d}\\
\frac{\left(\alpha_{H}\left(\mathbf{0}\right)\right)^{N}-1}{\alpha_{H}\left(\mathbf{0}\right)-1}\beta_{H}\left(\mathbf{0}\right) & \textrm{else}
\end{cases}
\end{equation}

Q.E.D. 
\begin{notation}
Just as a reminder, we write $q$ to denote $q_{H}$. Like in the
one-dimensional case, we write $\chi_{H,N}:\mathbb{Z}_{p}^{r}\rightarrow\mathbb{Q}^{d}\subset\mathbb{Q}_{q}^{d}$
to denote the $N$th truncation of $\chi_{H}$: 
\begin{equation}
\chi_{H,N}\left(\mathbf{z}\right)\overset{\textrm{def}}{=}\sum_{\mathbf{n}=\mathbf{0}}^{p^{N}-1}\chi_{H}\left(\mathbf{n}\right)\left[\mathbf{z}\overset{p^{N}}{\equiv}\mathbf{n}\right],\textrm{ }\forall\mathbf{z}\in\mathbb{Z}_{p}^{r},\textrm{ }\forall N\in\mathbb{N}_{0}\label{eq:MD Definition of Nth truncation of Chi_H}
\end{equation}
Recall that: 
\begin{equation}
\chi_{H,N}\left(\mathbf{n}\right)\overset{\mathbb{Q}^{d}}{=}\chi_{H}\left(\mathbf{n}\right),\textrm{ }\forall\mathbf{n}\leq p^{N}-1
\end{equation}
In particular: 
\begin{equation}
\chi_{H}\left(\mathbf{n}\right)=\chi_{H,\lambda_{p}\left(\mathbf{n}\right)}\left(\mathbf{n}\right),\textrm{ }\forall\mathbf{n}\in\mathbb{N}_{0}^{r}
\end{equation}
Here, $\chi_{H,N}$ is a $\mathbb{Q}^{d}$-valued function which is
a linear combination of finitely many indicator functions for clopen
subsets of $\mathbb{Z}_{p}^{r}$. As such, $\chi_{H,N}$ is continuous
both as a function $\mathbb{Z}_{p}^{r}\rightarrow\mathbb{C}_{q}^{d}$
and as a function $\mathbb{Z}_{p}^{r}\rightarrow\mathbb{C}^{d}$;
in fact, $\chi_{H,N}$ is continuous on $K^{d}$ for any topological
field $K$ containing $\mathbb{Q}$. As a result, of this, in writing
$\hat{\chi}_{H,N}:\hat{\mathbb{Z}}_{p}^{r}\rightarrow\overline{\mathbb{Q}}^{d}$
to denote the Fourier Transform of $\chi_{H,N}$: 
\begin{equation}
\hat{\chi}_{H,N}\left(\mathbf{t}\right)\overset{\textrm{def}}{=}\int_{\mathbb{Z}_{p}^{r}}\chi_{H,N}\left(\mathbf{z}\right)e^{-2\pi i\left\{ \mathbf{t}\mathbf{z}\right\} _{p}}d\mathbf{z},\textrm{ }\forall\mathbf{t}\in\hat{\mathbb{Z}}_{p}^{r}\label{eq:MD Definition of the Fourier Coefficients of Chi_H,N}
\end{equation}
observe that this integral is convergent in both the topology of $\mathbb{C}$
\emph{and} the topology of $\mathbb{C}_{q}$. As before, because $\chi_{H,N}$
is locally constant, it will reduce to a finite sum (see \textbf{Proposition
\ref{prop:MD Recursive formula for Chi_H,N hat}} in Subsection \ref{subsec:6.2.1 -and-},
below)): 
\begin{equation}
\hat{\chi}_{H,N}\left(\mathbf{t}\right)\overset{\overline{\mathbb{Q}}^{d}}{=}\frac{\mathbf{1}_{\mathbf{0}}\left(p^{N}\mathbf{t}\right)}{p^{rN}}\sum_{\mathbf{n}=\mathbf{0}}^{p^{N}-1}\chi_{H}\left(\mathbf{n}\right)e^{-2\pi i\mathbf{n}\cdot\mathbf{t}},\textrm{ }\forall\mathbf{t}\in\hat{\mathbb{Z}}_{p}^{r}
\end{equation}
where: 
\begin{equation}
\mathbf{1}_{\mathbf{0}}\left(p^{N}\mathbf{t}\right)=\left[\left\Vert \mathbf{t}\right\Vert _{p}\leq p^{N}\right]=\prod_{\ell=1}^{r}\left[\left|t_{\ell}\right|_{p}\leq p^{N}\right],\textrm{ }\forall\mathbf{t}\in\hat{\mathbb{Z}}_{p}^{r}
\end{equation}
is the \emph{scalar-valued }indicator function for the set $\left\{ \left\Vert \mathbf{t}\right\Vert _{p}\leq p^{N}\right\} $.
Once more, all of the computations and formal manipulations we will
perform hold simultaneously in the topologies of $\mathbb{C}$ and
in $\mathbb{C}_{q}$; they both occur in $\overline{\mathbb{Q}}\subset\mathbb{C}\cap\mathbb{C}_{q}$.
The difference between the archimedean and non-archimedean topologies
will only emerge when we consider what happens as $N$ tends to $\infty$.
Much like the one-dimensional case, the reason this all works because
of the invariance of the Fourier series formula for $\left[\mathbf{z}\overset{p^{N}}{\equiv}\mathbf{n}\right]$
under the action of elements of the Galois group $\textrm{Gal}\left(\overline{\mathbb{Q}}/\mathbb{Q}\right)$.

Additionally, note that because of $\chi_{H,N}$'s local constant-ness
and the finitude of its range, we have that, for each $N$, $\hat{\chi}_{H,N}$
has finite support: 
\begin{equation}
\hat{\chi}_{H,N}\left(\mathbf{t}\right)=\mathbf{0},\textrm{ }\forall\left\Vert \mathbf{t}\right\Vert _{p}\geq p^{N+1},\textrm{ }\forall N\in\mathbb{N}_{0}\label{eq:MD Vanishing of Chi_H,N hat for all t with sufficiently large denominators}
\end{equation}
Consequently, the Fourier series: 
\begin{equation}
\chi_{H,N}\left(\mathbf{z}\right)\overset{\overline{\mathbb{Q}}^{d}}{=}\sum_{\mathbf{t}\in\hat{\mathbb{Z}}_{p}^{r}}\hat{\chi}_{H,N}\left(\mathbf{t}\right)e^{2\pi i\left\{ \mathbf{t}\mathbf{z}\right\} _{p}}\label{eq:MD Fourier series for Chi_H,N}
\end{equation}
will converge in both $\mathbb{C}^{d}$ and $\mathbb{C}_{q}^{d}$
uniformly with respect to $\mathbf{z}\in\mathbb{Z}_{p}^{r}$, reducing
to a finite sum in all cases. 
\end{notation}
Lastly, like in the one-dimensional case, we will need to express
the interaction between $p$-adic structure functional equations and
$N$th truncations. For full generality (matrix-type and vector-type),
this result is given in terms of functions taking values in an arbitrary
linear space $V$ over $\mathbb{Q}$. 
\begin{lem}
\label{lem:MD functional equations and truncation lemma}Let $V$
be a $\mathbb{Q}$-linear space, and consider a function $\chi:\mathbb{N}_{0}^{r}\rightarrow V$.
Suppose that for all $\mathbf{j}\in\mathbb{Z}^{r}/p\mathbb{Z}^{r}$
there are functions $\Phi_{\mathbf{j}}:\mathbb{N}_{0}^{r}\times V\rightarrow V$
so that: 
\begin{equation}
\chi\left(p\mathbf{n}+\mathbf{j}\right)\overset{V}{=}\Phi_{\mathbf{j}}\left(\mathbf{n},\chi\left(\mathbf{n}\right)\right),\textrm{ }\forall\mathbf{n}\in\mathbb{N}_{0}^{r},\textrm{ }\forall\mathbf{j}\in\mathbb{Z}^{r}/p\mathbb{Z}^{r}\label{eq:MD Relation between truncations and functional equations - Hypothesis}
\end{equation}
Then, the $N$th truncations $\chi_{N}$ satisfy the functional equations:
\begin{equation}
\chi_{N}\left(p\mathbf{n}+\mathbf{j}\right)\overset{V}{=}\Phi_{\mathbf{j}}\left(\left[\mathbf{n}\right]_{p^{N-1}},\chi_{N-1}\left(\mathbf{n}\right)\right),\textrm{ }\forall\mathbf{n}\in\mathbb{N}_{0}^{r},\textrm{ }\forall\mathbf{j}\in\mathbb{Z}^{r}/p\mathbb{Z}^{r}\label{eq:MD Relation between truncations and functional equations, version 1}
\end{equation}
Equivalently: 
\begin{equation}
\chi\left(\left[p\mathbf{n}+\mathbf{j}\right]_{p^{N}}\right)\overset{V}{=}\Phi_{\mathbf{j}}\left(\left[\mathbf{n}\right]_{p^{N-1}},\chi\left(\left[\mathbf{n}\right]_{p^{N-1}}\right)\right),\textrm{ }\forall\mathbf{n}\in\mathbb{N}_{0}^{r},\textrm{ }\forall\mathbf{j}\in\mathbb{Z}^{r}/p\mathbb{Z}^{r}\label{eq:MD Relation between truncations and functional equations, version 2}
\end{equation}
\end{lem}
Proof: Fix $N\geq0$, $\mathbf{n}\in\mathbb{N}_{0}^{r}$, and $\mathbf{j}\in\mathbb{Z}^{r}/p\mathbb{Z}^{r}$.
Then: 
\begin{align*}
\chi\left(\left[p\mathbf{n}+\mathbf{j}\right]_{p^{N}}\right) & \overset{V}{=}\chi_{N}\left(p\mathbf{n}+\mathbf{j}\right)\\
 & =\sum_{\mathbf{m}=\mathbf{0}}^{p^{N}-1}\chi\left(\mathbf{m}\right)\left[p\mathbf{n}+\mathbf{j}\overset{p^{N}}{\equiv}\mathbf{m}\right]\\
\left(\textrm{split }\mathbf{m}\textrm{ mod }p\right); & =\sum_{\mathbf{m}=\mathbf{0}}^{p^{N-1}-1}\sum_{\mathbf{k}=\mathbf{0}}^{p-1}\chi\left(p\mathbf{m}+\mathbf{k}\right)\left[p\mathbf{n}+\mathbf{j}\overset{p^{N}}{\equiv}p\mathbf{m}+\mathbf{k}\right]\\
 & =\sum_{\mathbf{m}=\mathbf{0}}^{p^{N-1}-1}\sum_{\mathbf{k}=\mathbf{0}}^{p-1}\Phi_{\mathbf{k}}\left(\mathbf{m},\chi\left(\mathbf{m}\right)\right)\underbrace{\left[\mathbf{n}\overset{p^{N-1}}{\equiv}\mathbf{m}+\frac{\mathbf{k}-\mathbf{j}}{p}\right]}_{0\textrm{ }\forall\mathbf{n}\textrm{ if }\mathbf{k}\neq\mathbf{j}}\\
 & =\sum_{\mathbf{m}=\mathbf{0}}^{p^{N-1}-1}\Phi_{\mathbf{j}}\left(\mathbf{m},\chi\left(\mathbf{m}\right)\right)\left[\mathbf{n}\overset{p^{N-1}}{\equiv}\mathbf{m}\right]\\
 & =\Phi_{\mathbf{j}}\left(\left[\mathbf{n}\right]_{p^{N-1}},\chi\left(\left[\mathbf{n}\right]_{p^{N-1}}\right)\right)
\end{align*}

Q.E.D.

\newpage{}

\section{\label{sec:6.2 Fourier-Transforms-and}Fourier Transforms and Quasi-Integrability
(Again)}

IN THIS SECTION, WE ASSUME THAT $H$ IS NON-SINGULAR.

\vphantom{}

Section \pageref{sec:6.2 Fourier-Transforms-and} is mostly the same
as its one-dimensional predecessor in Section \pageref{sec:4.2 Fourier-Transforms-=00003D000026}.
As mentioned earlier, the primary distinction of the multi-dimensional
case is the non-commutativity of matrix multiplication. This eventually
leads to three distinct cases: those where $\alpha_{H}\left(\mathbf{0}\right)=\mathbf{I}_{d}$,
those where $\mathbf{I}_{d}-\alpha_{H}\left(\mathbf{0}\right)$ is
invertible and $H$ is commutative, and those where $\mathbf{I}_{d}-\alpha_{H}\left(\mathbf{0}\right)$
is invertible and $H$ is \emph{not }commutative. Although an $\mathcal{F}$-series
representation of $\chi_{H}$ is obtained for all three cases (\textbf{Theorem
\ref{thm:MD F-series for Chi_H}}), it is only for the first two cases
that this series can be used to deduce a closed-form expression for
a Fourier transform of $\hat{\chi}_{H}$. Subsection \pageref{subsec:6.2.4 Multi-Dimensional--=00003D000026}
takes the reader right up to the computational obstacle that appears
in the case of a non-commutative $H$.

That being said, for the two cases where we \emph{can }deduce formulae
for $\hat{\chi}_{H}$, our multi-dimensional $\left(p,q\right)$-adic
Wiener Tauberian Theorem can then be applied to yield a \textbf{Tauberian
Spectral Theorem }for suitable multi-dimensional Hydra maps, exactly
like in the one-dimensional case.

\subsection{\label{subsec:6.2.1 -and-}$\hat{\chi}_{H,N}$ and $\hat{A}_{H}$
(Again)}

We begin by computing a recursive formula for $\hat{\chi}_{H,N}$
in terms of $\hat{\chi}_{H,N-1}$, which we then solve to obtain an
explicit formula for $\hat{\chi}_{H,N}$. 
\begin{prop}[\textbf{Formulae for Multi-Dimensional} $\hat{\chi}_{H,N}$]
\label{prop:MD Recursive formula for Chi_H,N hat}\ 

\vphantom{}

I. 
\begin{equation}
\hat{\chi}_{H,N}\left(\mathbf{t}\right)\overset{\overline{\mathbb{Q}}^{d}}{=}\frac{\mathbf{1}_{\mathbf{0}}\left(p^{N}\mathbf{t}\right)}{p^{rN}}\sum_{\mathbf{n}=\mathbf{0}}^{p^{N}-1}\chi_{H}\left(\mathbf{n}\right)e^{-2\pi i\mathbf{n}\cdot\mathbf{t}}\label{eq:MD Chi_H,N hat transform formula - sum form}
\end{equation}

\vphantom{}

II. 
\begin{equation}
\mathbf{1}_{\mathbf{0}}\left(p^{N}\mathbf{t}\right)\alpha_{H}\left(\mathbf{t}\right)\hat{\chi}_{H,N-1}\left(p\mathbf{t}\right)+\mathbf{1}_{\mathbf{0}}\left(p\mathbf{t}\right)\beta_{H}\left(\mathbf{t}\right),\textrm{ }\forall N\geq1,\textrm{ }\forall\mathbf{t}\in\hat{\mathbb{Z}}_{p}^{r}\label{eq:MD Chi_H,N hat functional equation}
\end{equation}
the nesting of which yields: 
\begin{equation}
\hat{\chi}_{H,N}\left(\mathbf{t}\right)=\sum_{n=0}^{N-1}\mathbf{1}_{\mathbf{0}}\left(p^{n+1}\mathbf{t}\right)\left(\prod_{m=0}^{n-1}\alpha_{H}\left(p^{m}\mathbf{t}\right)\right)\beta_{H}\left(p^{n}\mathbf{t}\right)\label{eq:MD Explicit formula for Chi_H,N hat}
\end{equation}
where the $m$-product is defined to be $1$ when $n=0$. In particular,
note that $\hat{\chi}_{H,N}\left(\mathbf{t}\right)=\mathbf{0}$ for
all $\mathbf{t}\in\hat{\mathbb{Z}}_{p}^{r}$ with $\left\Vert \mathbf{t}\right\Vert _{p}>p^{N}$. 
\end{prop}
Proof: Using \textbf{Proposition \ref{prop:Multi-Dimensional indicator function Fourier Series}},
we have: 
\begin{align*}
\chi_{H,N}\left(\mathbf{z}\right) & =\sum_{\mathbf{n}=\mathbf{0}}^{p^{N}-1}\chi_{H}\left(\mathbf{n}\right)\left[\mathbf{z}\overset{p^{N}}{\equiv}\mathbf{n}\right]\\
 & =\sum_{\mathbf{n}=\mathbf{0}}^{p^{N}-1}\chi_{H}\left(\mathbf{n}\right)\frac{1}{p^{rN}}\sum_{\left\Vert \mathbf{t}\right\Vert _{p}\leq p^{N}}e^{2\pi i\left\{ \mathbf{t}\left(\mathbf{z}-\mathbf{n}\right)\right\} _{p}}\\
 & =\sum_{\left\Vert \mathbf{t}\right\Vert _{p}\leq p^{N}}\left(\frac{1}{p^{rN}}\sum_{\mathbf{n}=\mathbf{0}}^{p^{N}-1}\chi_{H}\left(\mathbf{n}\right)e^{-2\pi i\mathbf{n}\cdot\mathbf{t}}\right)e^{2\pi i\left\{ \mathbf{t}\mathbf{z}\right\} _{p}}
\end{align*}
This is the Fourier series representation of $\chi_{H,N}$; as such,
the coefficient of $e^{2\pi i\left\{ \mathbf{t}\mathbf{z}\right\} _{p}}$
in the series is precisely the formula for $\hat{\chi}_{H,N}\left(\mathbf{t}\right)$:
\[
\hat{\chi}_{H,N}\left(\mathbf{t}\right)=\frac{\mathbf{1}_{\mathbf{0}}\left(p^{N}\mathbf{t}\right)}{p^{rN}}\sum_{\mathbf{n}=\mathbf{0}}^{p^{N}-1}\chi_{H}\left(\mathbf{n}\right)e^{-2\pi i\mathbf{n}\cdot\mathbf{t}}
\]
This proves (I).

To prove (II), just like in the one-dimensional case, we split the
$\mathbf{n}$-sum modulo $p$ and utilize $\chi_{H}$'s functional
equation. Here\textbf{ }$\mathbf{t}$ satisfies $\left\Vert \mathbf{t}\right\Vert _{p}\leq p^{N}$.
\begin{align*}
\hat{\chi}_{H,N}\left(\mathbf{t}\right) & =\frac{\mathbf{1}_{\mathbf{0}}\left(p^{N}\mathbf{t}\right)}{p^{rN}}\sum_{\mathbf{j}=\mathbf{0}}^{p-1}\sum_{\mathbf{n}=\mathbf{0}}^{p^{N-1}-1}\chi_{H}\left(p\mathbf{n}+\mathbf{j}\right)e^{-2\pi i\left(p\mathbf{n}+\mathbf{j}\right)\cdot\mathbf{t}}\\
 & =\frac{\mathbf{1}_{\mathbf{0}}\left(p^{N}\mathbf{t}\right)}{p^{rN}}\sum_{\mathbf{j}=\mathbf{0}}^{p-1}\sum_{\mathbf{n}=\mathbf{0}}^{p^{N-1}-1}\frac{\mathbf{A}_{\mathbf{j}}\chi_{H}\left(\mathbf{n}\right)+\mathbf{b}_{\mathbf{j}}}{\mathbf{D}_{\mathbf{j}}}e^{-2\pi i\left(p\mathbf{n}+\mathbf{j}\cdot\mathbf{t}\right)}
\end{align*}
Interchanging the $\mathbf{j}$ and $\mathbf{n}$ sums and breaking
things up yields
\begin{align*}
\frac{\mathbf{1}_{\mathbf{0}}\left(p^{N}\mathbf{t}\right)}{p^{r\left(N-1\right)}}\sum_{\mathbf{n}=\mathbf{0}}^{p^{N-1}-1}\underbrace{\left(\frac{1}{p^{r}}\sum_{\mathbf{j}=\mathbf{0}}^{p-1}\mathbf{D}_{\mathbf{j}}^{-1}\mathbf{A}_{\mathbf{j}}e^{-2\pi i\mathbf{j}\cdot\mathbf{t}}\right)}_{\alpha_{H}\left(\mathbf{t}\right)}\chi_{H}\left(\mathbf{n}\right)e^{-2\pi i\left(p\mathbf{n}\right)\cdot\mathbf{t}}\\
+\frac{\mathbf{1}_{\mathbf{0}}\left(p^{N}\mathbf{t}\right)}{p^{r\left(N-1\right)}}\sum_{\mathbf{n}=\mathbf{0}}^{p^{N-1}-1}\underbrace{\left(\frac{1}{p^{r}}\sum_{\mathbf{j}=\mathbf{0}}^{p-1}\mathbf{D}_{\mathbf{j}}^{-1}\mathbf{b}_{\mathbf{j}}e^{-2\pi i\mathbf{j}\cdot\mathbf{t}}\right)}_{\beta_{H}\left(\mathbf{t}\right)}e^{-2\pi i\left(p\mathbf{n}\right)\cdot\mathbf{t}}
\end{align*}
Hence:
\begin{equation}
\hat{\chi}_{H,N}\left(\mathbf{t}\right)=\frac{\mathbf{1}_{\mathbf{0}}\left(p^{N}\mathbf{t}\right)}{p^{r\left(N-1\right)}}\sum_{\mathbf{n}=\mathbf{0}}^{p^{N-1}-1}\left(\alpha_{H}\left(\mathbf{t}\right)\chi_{H}\left(\mathbf{n}\right)+\beta_{H}\left(\mathbf{t}\right)\right)e^{-2\pi i\mathbf{n}\cdot\left(p\mathbf{t}\right)}
\end{equation}
Re-writing the right-hand side in terms of $\chi_{H,N-1}$ yields:
\begin{equation}
\hat{\chi}_{H,N}\left(\mathbf{t}\right)=\mathbf{1}_{\mathbf{0}}\left(p^{N}\mathbf{t}\right)\alpha_{H}\left(\mathbf{t}\right)\hat{\chi}_{H,N-1}\left(p\mathbf{t}\right)+\beta_{H}\left(\mathbf{t}\right)\frac{\mathbf{1}_{\mathbf{0}}\left(p^{N}\mathbf{t}\right)}{p^{r\left(N-1\right)}}\sum_{\mathbf{n}=\mathbf{0}}^{p^{N-1}-1}e^{-2\pi i\mathbf{n}\cdot\left(p\mathbf{t}\right)}\label{eq:Chi_H N hat, almost ready to nest-1}
\end{equation}
Simplifying the remaining $n$-sum, we obtain:
\begin{align*}
\frac{1}{p^{r\left(N-1\right)}}\sum_{\mathbf{n}=\mathbf{0}}^{p^{N-1}-1}e^{-2\pi i\mathbf{n}\cdot\left(p\mathbf{t}\right)} & =\begin{cases}
\frac{1}{p^{r\left(N-1\right)}}\sum_{\mathbf{n}=\mathbf{0}}^{p^{N-1}-1}1 & \textrm{if }\left\Vert \mathbf{t}\right\Vert _{p}\leq p\\
\frac{1}{p^{r\left(N-1\right)}}\prod_{m=1}^{r}\frac{\left(e^{-2\pi ipt_{m}}\right)^{p^{N-1}}-1}{e^{-2\pi ipt_{m}}-1} & \textrm{if }\left\Vert \mathbf{t}\right\Vert _{p}\geq p^{2}
\end{cases}\\
 & =\begin{cases}
1 & \textrm{if }\left\Vert \mathbf{t}\right\Vert _{p}\leq p\\
0 & \textrm{if }\left\Vert \mathbf{t}\right\Vert _{p}\geq p^{2}
\end{cases}\\
 & =\mathbf{1}_{\mathbf{0}}\left(p\mathbf{t}\right)
\end{align*}
Hence: 
\begin{align*}
\hat{\chi}_{H,N}\left(\mathbf{t}\right) & =\mathbf{1}_{\mathbf{0}}\left(p^{N}\mathbf{t}\right)\left(\alpha_{H}\left(\mathbf{t}\right)\hat{\chi}_{H,N-1}\left(p\mathbf{t}\right)+\beta_{H}\left(\mathbf{t}\right)\mathbf{1}_{\mathbf{0}}\left(p\mathbf{t}\right)\right)\\
 & =\mathbf{1}_{\mathbf{0}}\left(p^{N}\mathbf{t}\right)\alpha_{H}\left(\mathbf{t}\right)\hat{\chi}_{H,N-1}\left(p\mathbf{t}\right)+\mathbf{1}_{\mathbf{0}}\left(p\mathbf{t}\right)\beta_{H}\left(\mathbf{t}\right)
\end{align*}
which proves (II).

With (II) proven, we can then nest (\ref{eq:MD Chi_H,N hat functional equation})
to derive an explicit formula for $\hat{\chi}_{H,N}\left(\mathbf{t}\right)$:
\begin{align*}
\hat{\chi}_{H,N}\left(\mathbf{t}\right) & =\mathbf{1}_{\mathbf{0}}\left(p^{N}\mathbf{t}\right)\alpha_{H}\left(\mathbf{t}\right)\hat{\chi}_{H,N-1}\left(p\mathbf{t}\right)+\mathbf{1}_{\mathbf{0}}\left(p\mathbf{t}\right)\beta_{H}\left(\mathbf{t}\right)\\
 & =\mathbf{1}_{\mathbf{0}}\left(p^{N}\mathbf{t}\right)\alpha_{H}\left(\mathbf{t}\right)\left(\mathbf{1}_{\mathbf{0}}\left(p^{N-1}p\mathbf{t}\right)\alpha_{H}\left(p\mathbf{t}\right)\hat{\chi}_{H,N-2}\left(p^{2}\mathbf{t}\right)+\mathbf{1}_{0}\left(p^{2}\mathbf{t}\right)\beta_{H}\left(p\mathbf{t}\right)\right)\\
 & +\mathbf{1}_{0}\left(p\mathbf{t}\right)\beta_{H}\left(\mathbf{t}\right)\\
 & =\mathbf{1}_{\mathbf{0}}\left(p^{N}\mathbf{t}\right)\alpha_{H}\left(\mathbf{t}\right)\alpha_{H}\left(p\mathbf{t}\right)\hat{\chi}_{H,N-2}\left(p^{2}\mathbf{t}\right)\\
 & +\mathbf{1}_{\mathbf{0}}\left(p^{2}\mathbf{t}\right)\alpha_{H}\left(\mathbf{t}\right)\beta_{H}\left(p\mathbf{t}\right)+\mathbf{1}_{\mathbf{0}}\left(p\mathbf{t}\right)\beta_{H}\left(\mathbf{t}\right)\\
 & \vdots\\
 & =\mathbf{1}_{\mathbf{0}}\left(p^{N}\mathbf{t}\right)\left(\prod_{n=0}^{N-2}\alpha_{H}\left(p^{n}\mathbf{t}\right)\right)\hat{\chi}_{H,1}\left(p^{N-1}\mathbf{t}\right)\\
 & +\sum_{n=0}^{N-2}\left(\prod_{m=0}^{n-1}\alpha_{H}\left(p^{m}\mathbf{t}\right)\right)\beta_{H}\left(p^{n}\mathbf{t}\right)\mathbf{1}_{\mathbf{0}}\left(p^{n+1}\mathbf{t}\right)
\end{align*}
where the $m$-product is $1$ whenever $n=0$.

Finally, using (\ref{eq:MD Chi_H,N hat transform formula - sum form})
to compute $\hat{\chi}_{H,1}\left(\mathbf{t}\right)$, we have:
\begin{align}
\hat{\chi}_{H,1}\left(\mathbf{t}\right) & =\frac{\mathbf{1}_{\mathbf{0}}\left(p\mathbf{t}\right)}{p^{r}}\sum_{\mathbf{j}=\mathbf{0}}^{p-1}\chi_{H}\left(\mathbf{j}\right)e^{-2\pi i\mathbf{j}\cdot\mathbf{t}}
\end{align}
Since $\chi_{H}\left(\mathbf{0}\right)=\mathbf{0}$, $\chi_{H}$'s
functional equation gives: 
\begin{equation}
\chi_{H}\left(\mathbf{j}\right)=\chi_{H}\left(p\mathbf{0}+\mathbf{j}\right)=\frac{\mathbf{A}_{\mathbf{j}}\chi_{H}\left(\mathbf{0}\right)+\mathbf{b}_{\mathbf{j}}}{\mathbf{D}_{\mathbf{j}}}=\mathbf{D}_{\mathbf{j}}^{-1}\mathbf{b}_{\mathbf{j}}
\end{equation}
and so: 
\begin{align*}
\hat{\chi}_{H,1}\left(\mathbf{t}\right) & =\frac{\mathbf{1}_{\mathbf{0}}\left(p\mathbf{t}\right)}{p^{r}}\sum_{\mathbf{j}=\mathbf{0}}^{p-1}\chi_{H}\left(\mathbf{j}\right)e^{-2\pi i\mathbf{j}\cdot\mathbf{t}}\\
 & =\frac{\mathbf{1}_{\mathbf{0}}\left(p\mathbf{t}\right)}{p^{r}}\sum_{\mathbf{j}=\mathbf{0}}^{p-1}\mathbf{D}_{\mathbf{j}}^{-1}\mathbf{b}_{\mathbf{j}}e^{-2\pi i\mathbf{j}\cdot\mathbf{t}}\\
 & =\mathbf{1}_{\mathbf{0}}\left(p\mathbf{t}\right)\beta_{H}\left(\mathbf{t}\right)
\end{align*}
Consequently: 
\begin{align*}
\hat{\chi}_{H,N}\left(\mathbf{t}\right) & =\mathbf{1}_{\mathbf{0}}\left(p^{N}\mathbf{t}\right)\left(\prod_{n=0}^{N-2}\alpha_{H}\left(p^{n}\mathbf{t}\right)\right)\hat{\chi}_{H,1}\left(p^{N-1}\mathbf{t}\right)\\
 & +\sum_{n=0}^{N-2}\left(\prod_{m=0}^{n-1}\alpha_{H}\left(p^{m}\mathbf{t}\right)\right)\beta_{H}\left(p^{n}\mathbf{t}\right)\mathbf{1}_{\mathbf{0}}\left(p^{n+1}\mathbf{t}\right)\\
 & =\mathbf{1}_{\mathbf{0}}\left(p^{N}t\right)\beta_{H}\left(p^{N-1}\mathbf{t}\right)\prod_{n=0}^{N-2}\alpha_{H}\left(p^{n}\mathbf{t}\right)\\
 & +\sum_{n=0}^{N-2}\mathbf{1}_{\mathbf{0}}\left(p^{n+1}\mathbf{t}\right)\left(\prod_{m=0}^{n-1}\alpha_{H}\left(p^{m}\mathbf{t}\right)\right)\beta_{H}\left(p^{n}\mathbf{t}\right)\\
 & =\sum_{n=0}^{N-1}\mathbf{1}_{\mathbf{0}}\left(p^{n+1}\mathbf{t}\right)\left(\prod_{m=0}^{n-1}\alpha_{H}\left(p^{m}\mathbf{t}\right)\right)\beta_{H}\left(p^{n}\mathbf{t}\right)
\end{align*}
which proves (\ref{eq:MD Explicit formula for Chi_H,N hat}).

Q.E.D.

\vphantom{}

Next, we introduce our multi-dimensional analogue $\hat{A}_{H}$. 
\begin{defn}[\textbf{Multi-Dimensional $\hat{A}_{H}$}]
\nomenclature{$\hat{A}_{H}\left(\mathbf{t}\right)$}{ }\textbf{ }We
define the function $\hat{A}_{H}:\hat{\mathbb{Z}}_{p}^{r}\rightarrow\textrm{GL}_{d}\left(\overline{\mathbb{Q}}\right)$
by: 
\begin{equation}
\hat{A}_{H}\left(\mathbf{t}\right)\overset{\textrm{def}}{=}\prod_{m=0}^{-v_{p}\left(\mathbf{t}\right)-1}\alpha_{H}\left(p^{m}\mathbf{t}\right),\textrm{ }\forall\mathbf{t}\in\hat{\mathbb{Z}}_{p}^{r}\label{eq:MD Definition of A_H hat}
\end{equation}
where the $m$-product is defined to be $\mathbf{I}_{d}$ whenever
$\mathbf{t}=\mathbf{0}$\textemdash so, $\hat{A}_{H}\left(\mathbf{0}\right)\overset{\textrm{def}}{=}\mathbf{I}_{d}$. 
\end{defn}
\begin{prop}[\textbf{Matrix $\alpha_{H}$ product in terms of $\hat{A}_{H}$}]
\label{prop:MD alpha_H A_H hat product} Fix $\mathbf{t}\in\hat{\mathbb{Z}}_{p}^{r}\backslash\left\{ \mathbf{0}\right\} $.
Then: 
\begin{align}
\mathbf{1}_{\mathbf{0}}\left(p^{n+1}\mathbf{t}\right)\prod_{m=0}^{n-1}\alpha_{H}\left(p^{m}\mathbf{t}\right) & =\begin{cases}
\mathbf{O}_{d} & \textrm{if }n\leq-v_{p}\left(\mathbf{t}\right)-2\\
\hat{A}_{H}\left(\mathbf{t}\right)\left(\alpha_{H}\left(\frac{\mathbf{t}\left|\mathbf{t}\right|_{p}}{p}\right)\right)^{-1} & \textrm{if }n=-v_{p}\left(\mathbf{t}\right)-1\\
\hat{A}_{H}\left(\mathbf{t}\right)\left(\alpha_{H}\left(\mathbf{0}\right)\right)^{n+v_{p}\left(\mathbf{t}\right)} & \textrm{if }n\geq-v_{p}\left(\mathbf{t}\right)
\end{cases}\label{eq:MD alpha product in terms of A_H hat}
\end{align}
\end{prop}
Proof: Fix $\mathbf{t}\in\hat{\mathbb{Z}}_{p}^{r}\backslash\left\{ \mathbf{0}\right\} $.
Then, we can write: 
\begin{align*}
\mathbf{1}_{\mathbf{0}}\left(p^{n+1}\mathbf{t}\right)\prod_{m=0}^{n-1}\alpha_{H}\left(p^{m}\mathbf{t}\right) & =\begin{cases}
\mathbf{O}_{d} & \textrm{if }\left\Vert \mathbf{t}\right\Vert _{p}\geq p^{n+2}\\
\prod_{m=0}^{n-1}\alpha_{H}\left(p^{m}\mathbf{t}\right) & \textrm{if }\left\Vert \mathbf{t}\right\Vert _{p}\leq p^{n+1}
\end{cases}\\
 & =\begin{cases}
\mathbf{O}_{d} & \textrm{if }\exists\ell:n\leq-v_{p}\left(t_{\ell}\right)-2\\
\prod_{m=0}^{n-1}\alpha_{H}\left(p^{m}\mathbf{t}\right) & \textrm{if }n\geq-v_{p}\left(t_{\ell}\right)-1,\textrm{ }\forall\ell
\end{cases}
\end{align*}
So, we get the zero matrix whenever $n+1$ is not greater than the
negatives of the valuations of each of $\mathbf{t}$s components.

Fixing $\mathbf{t}$, we have that $p^{m}\mathbf{t}\overset{1}{\equiv}\mathbf{0}$
only when $m$ is larger than the negative valuations of each of the
components of $\mathbf{t}$: 
\begin{equation}
m\geq\max_{1\leq\ell\leq r}\left(-v_{p}\left(t_{\ell}\right)\right)=-\min_{1\leq\ell\leq r}v_{p}\left(t_{\ell}\right)=-v_{p}\left(\mathbf{t}\right)
\end{equation}
Hence, $\alpha_{H}\left(p^{m}\mathbf{t}\right)=\alpha_{H}\left(\mathbf{0}\right)$
for all $m\geq-v_{p}\left(\mathbf{t}\right)$. Thus, for fixed $\mathbf{t}$:
\begin{equation}
\prod_{m=0}^{n-1}\alpha_{H}\left(p^{m}\mathbf{t}\right)=\begin{cases}
\prod_{m=0}^{-v_{p}\left(\mathbf{t}\right)-2}\alpha_{H}\left(p^{m}\mathbf{t}\right) & \textrm{if }n=-v_{p}\left(\mathbf{t}\right)-1\\
\prod_{m=0}^{-v_{p}\left(\mathbf{t}\right)-1}\alpha_{H}\left(p^{m}\mathbf{t}\right) & \textrm{if }n=-v_{p}\left(\mathbf{t}\right)\\
\prod_{m=0}^{-v_{p}\left(\mathbf{t}\right)-1}\alpha_{H}\left(p^{m}\mathbf{t}\right)\times\prod_{k=-v_{p}\left(\mathbf{t}\right)}^{n-1}\alpha_{H}\left(p^{k}\mathbf{t}\right) & \textrm{if }n\geq-v_{p}\left(\mathbf{t}\right)+1
\end{cases}
\end{equation}
Because $\alpha_{H}\left(p^{k}\mathbf{t}\right)=\alpha_{H}\left(\mathbf{0}\right)$
for all $k\geq-v_{p}\left(\mathbf{t}\right)$, this can be simplified
to:
\begin{align*}
\prod_{m=0}^{n-1}\alpha_{H}\left(p^{m}\mathbf{t}\right) & =\begin{cases}
\prod_{m=0}^{-v_{p}\left(\mathbf{t}\right)-2}\alpha_{H}\left(p^{m}\mathbf{t}\right) & \textrm{if }n=-v_{p}\left(\mathbf{t}\right)-1\\
\left(\prod_{m=0}^{-v_{p}\left(\mathbf{t}\right)-1}\alpha_{H}\left(p^{m}\mathbf{t}\right)\right)\left(\alpha_{H}\left(\mathbf{0}\right)\right)^{n+v_{p}\left(\mathbf{t}\right)} & \textrm{if }n\geq-v_{p}\left(\mathbf{t}\right)
\end{cases}\\
\left(\times\&\div\textrm{ by }\alpha_{H}\left(p^{-v_{p}\left(\mathbf{t}\right)-1}\mathbf{t}\right)\right); & =\begin{cases}
\hat{A}_{H}\left(\mathbf{t}\right)\left(\alpha_{H}\left(p^{-v_{p}\left(\mathbf{t}\right)-1}\mathbf{t}\right)\right)^{-1} & \textrm{if }n=-v_{p}\left(\mathbf{t}\right)-1\\
\hat{A}_{H}\left(\mathbf{t}\right)\left(\alpha_{H}\left(\mathbf{0}\right)\right)^{n+v_{p}\left(\mathbf{t}\right)} & \textrm{if }n\geq-v_{p}\left(\mathbf{t}\right)
\end{cases}
\end{align*}
which gives us the desired formula.

Q.E.D.

\vphantom{}

Next, we express the $\alpha_{H}$ product as a series. 
\begin{prop}[\textbf{Matrix $\alpha_{H}$ series formula}]
\label{prop:MD alpha H series} 
\end{prop}
\begin{equation}
\prod_{m=0}^{n-1}\alpha_{H}\left(p^{m}\mathbf{t}\right)\overset{\overline{\mathbb{Q}}^{d,d}}{=}\sum_{\mathbf{m}=\mathbf{0}}^{p^{n}-1}\kappa_{H}\left(\mathbf{m}\right)\left(\frac{H^{\prime}\left(\mathbf{0}\right)}{p^{r}}\right)^{n}e^{-2\pi i\left(\mathbf{m}\cdot\mathbf{t}\right)},\textrm{ }\forall n\geq0,\textrm{ }\forall\mathbf{t}\in\hat{\mathbb{Z}}_{p}^{r}\label{eq:MD alpha_H product expansion}
\end{equation}

Proof: We start by writing: 
\begin{equation}
\prod_{m=0}^{n-1}\alpha_{H}\left(p^{m}\mathbf{t}\right)=\prod_{m=0}^{n-1}\left(\sum_{\mathbf{j}=\mathbf{0}}^{p-1}\frac{1}{p^{r}}\frac{\mathbf{A}_{\mathbf{j}}}{\mathbf{D}_{\mathbf{j}}}e^{-2\pi i\mathbf{j}\cdot\left(p^{m}\mathbf{t}\right)}\right)
\end{equation}
and then apply (\ref{eq:MD M_H partial sum generating identity})
from \textbf{Proposition \ref{prop:MD generating function} }with
$z^{\mathbf{j}\cdot\mathbf{t}}$ replaced by $e^{-2\pi i\left(\mathbf{j}\cdot\mathbf{t}\right)}$:
\begin{align*}
\prod_{m=0}^{n-1}\alpha_{H}\left(p^{m}\mathbf{t}\right) & =\prod_{m=0}^{n-1}\left(\sum_{\mathbf{j}=\mathbf{0}}^{p-1}\frac{1}{p^{r}}\frac{\mathbf{A}_{\mathbf{j}}}{\mathbf{D}_{\mathbf{j}}}e^{-2\pi i\mathbf{j}\cdot\left(p^{m}\mathbf{t}\right)}\right)\\
 & =\frac{1}{p^{rn}}\sum_{\mathbf{m}=\mathbf{0}}^{p^{n}-1}M_{H}\left(\mathbf{m}\right)\left(H^{\prime}\left(\mathbf{0}\right)\right)^{n-\lambda_{p}\left(\mathbf{m}\right)}e^{-2\pi i\left(\mathbf{m}\cdot\mathbf{t}\right)}\\
\left(\kappa_{H}\left(\mathbf{m}\right)=M_{H}\left(\mathbf{m}\right)\left(H^{\prime}\left(\mathbf{0}\right)\right)^{-\lambda_{p}\left(\mathbf{m}\right)}\right); & =\frac{1}{p^{rn}}\sum_{\mathbf{m}=\mathbf{0}}^{p^{n}-1}\kappa_{H}\left(\mathbf{m}\right)\left(H^{\prime}\left(\mathbf{0}\right)\right)^{n}e^{-2\pi i\left(\mathbf{m}\cdot\mathbf{t}\right)}\\
 & =\sum_{\mathbf{m}=\mathbf{0}}^{p^{n}-1}\kappa_{H}\left(\mathbf{m}\right)\left(\frac{H^{\prime}\left(\mathbf{0}\right)}{p^{r}}\right)^{n}e^{-2\pi i\left(\mathbf{m}\cdot\mathbf{t}\right)}
\end{align*}

Q.E.D.

\vphantom{}

Just like the one-dimensional case, we now use this formula to exhibit
$\hat{A}_{H}$ as a radially-magnitudinal function, and then\textemdash with
the help of our multi-dimensional Fourier Resummation Lemmata\textemdash show
that $\hat{A}_{H}$ is the Fourier-Stieltjes transform of a rising-continuous
$\left(p,q\right)$-adic thick measure of matrix-type.
\begin{prop}
Let $H$ be semi-basic. Then, $\hat{A}_{H}\left(\mathbf{t}\right)$
is a depth-$r$ $\left(p,q\right)$-adic thick measure of matrix type. 
\end{prop}
Proof: Let $H$ be semi-basic. Then, $p$ divides the diagonal entires
of the $\mathbf{D}_{\mathbf{j}}$s which are not $1$, and the entries
of any $\mathbf{D}_{\mathbf{j}}$ are co-prime to the entries of every
$\mathbf{A}_{\mathbf{j}}$. As such, taking $q_{H}$-adic matrix norms
yields: 
\begin{equation}
\left\Vert \alpha_{H}\left(\mathbf{t}\right)\right\Vert _{q_{H}}=\left\Vert \frac{1}{p^{r}}\sum_{\mathbf{j}=\mathbf{0}}^{p-1}\frac{\mathbf{A}_{\mathbf{j}}}{\mathbf{D}_{\mathbf{j}}}e^{-2\pi i\mathbf{j}\cdot\mathbf{t}}\right\Vert _{q_{H}}\leq\max_{\mathbf{j}\leq p-1}\left|\frac{1}{p^{r}}\frac{\mathbf{A}_{\mathbf{j}}}{\mathbf{D}_{\mathbf{j}}}\right|_{q_{H}}\leq1
\end{equation}
Hence: 
\begin{align*}
\sup_{\mathbf{t}\in\hat{\mathbb{Z}}_{p}^{r}}\left\Vert \hat{A}_{H}\left(\mathbf{t}\right)\right\Vert _{q_{H}} & \leq\max\left\{ 1,\sup_{\mathbf{t}\in\hat{\mathbb{Z}}_{p}^{r}\backslash\left\{ \mathbf{0}\right\} }\left\Vert \prod_{m=0}^{-v_{p}\left(\mathbf{t}\right)-1}\alpha_{H}\left(p^{m}\mathbf{t}\right)\right\Vert _{q_{H}}\right\} \\
 & \leq\max\left\{ 1,\sup_{\mathbf{t}\in\hat{\mathbb{Z}}_{p}^{r}\backslash\left\{ \mathbf{0}\right\} }\prod_{m=0}^{-v_{p}\left(\mathbf{t}\right)-1}1\right\} \\
 & =1
\end{align*}
and so, $\hat{A}_{H}$ is then in $B\left(\hat{\mathbb{Z}}_{p}^{r},\mathbb{C}_{q}^{d,d}\right)$.
As such, the map: 
\begin{equation}
\mathbf{f}\in C\left(\mathbb{Z}_{p}^{r},\mathbb{C}_{q}^{d}\right)\mapsto\sum_{\mathbf{t}\in\hat{\mathbb{Z}}_{p}^{r}}\hat{A}_{H}\left(\mathbf{t}\right)\hat{\mathbf{f}}\left(-\mathbf{t}\right)\in\mathbb{C}_{q}^{d}
\end{equation}
is then a continuous, defining a $\left(p,q\right)$-adic thick measure
of matrix type.

Q.E.D.

\vphantom{}

Next, some obligatory abbreviations:
\begin{defn}
\nomenclature{$\tilde{A}_{H,N}\left(\mathbf{z}\right)$}{$N$th partial sum of the Fourier series generated by $\hat{A}_{H}\left(\mathbf{t}\right)$.}For
each $N\in\mathbb{N}_{0}$, we write $\tilde{A}_{H,N}:\mathbb{Z}_{p}^{r}\rightarrow\overline{\mathbb{Q}}^{d,d}$
to denote the matrix-valued function: 
\begin{equation}
\tilde{A}_{H,N}\left(\mathbf{z}\right)\overset{\textrm{def}}{=}\sum_{\left\Vert \mathbf{t}\right\Vert _{p}\leq p^{N}}\hat{A}_{H}\left(\mathbf{t}\right)e^{2\pi i\left\{ \mathbf{t}\mathbf{z}\right\} _{p}}\label{eq:MD Definition of A_H,N twiddle}
\end{equation}
\end{defn}
\begin{defn}
We write $\mathbf{I}_{H,n}\left(n\right)$ \nomenclature{$\mathbf{I}_{H,n}\left(n\right)$}{$\overset{\textrm{def}}{=}\mathbf{I}_{d}-\left(H^{\prime}\left(\mathbf{0}\right)\right)^{n}\alpha_{H}\left(\mathbf{0}\right)\left(H^{\prime}\left(\mathbf{0}\right)\right)^{-n}$ \nopageref}to
denote the $d\times d$ matrix: 
\begin{align}
\mathbf{I}_{H}\left(n\right) & \overset{\textrm{def}}{=}\mathbf{I}_{d}-\mathcal{C}_{H}\left(\alpha_{H}\left(\mathbf{0}\right):n\right)\label{eq:Definition of I_H}\\
 & =\mathbf{I}_{d}-\left(H^{\prime}\left(\mathbf{0}\right)\right)^{n}\alpha_{H}\left(\mathbf{0}\right)\left(H^{\prime}\left(\mathbf{0}\right)\right)^{-n}\nonumber 
\end{align}
\end{defn}
\begin{rem}
Note that $\mathbf{I}_{H}\left(n\right)=\mathbf{I}_{d}-\alpha_{H}\left(\mathbf{0}\right)$
for all $n\in\mathbb{N}_{0}$ whenever $H$ is commutative. Also,
$\mathbf{I}_{H}\left(n\right)=\mathbf{O}_{d}$ for all $n\in\mathbb{N}_{0}$
whenever $\alpha_{H}\left(\mathbf{0}\right)=\mathbf{I}_{d}$. 
\end{rem}
Our next theorem details the properties of $A_{H}$: 
\begin{thm}[\textbf{Properties of Multi-Dimensional $\hat{A}_{H}$}]
\label{thm:MD properties of A_H hat}\ 

\vphantom{}

I. ($\hat{A}_{H}$ is radially-magnitudinal and $\left(p,q\right)$-adically
regular) 
\begin{equation}
\hat{A}_{H}\left(\mathbf{t}\right)\overset{\overline{\mathbb{Q}}^{d,d}}{=}\begin{cases}
\mathbf{I}_{d} & \textrm{if }\mathbf{t}=0\\
\left(\sum_{\mathbf{m}=\mathbf{0}}^{p^{-v_{p}\left(\mathbf{t}\right)}-1}\kappa_{H}\left(\mathbf{m}\right)e^{-2\pi i\left(\mathbf{m}\cdot\mathbf{t}\right)}\right)\left(\frac{H^{\prime}\left(\mathbf{0}\right)}{p^{r}}\right)^{-v_{p}\left(\mathbf{t}\right)} & \textrm{else}
\end{cases},\textrm{ }\forall\mathbf{t}\in\hat{\mathbb{Z}}_{p}^{r}\label{eq:A_H hat as the product of radially symmetric and magnitude-dependent measures-1}
\end{equation}
Here, in an abuse of notation, we write: 
\[
\left(\frac{H^{\prime}\left(\mathbf{0}\right)}{p^{r}}\right)^{-v_{p}\left(\mathbf{0}\right)}\overset{\textrm{def}}{=}\mathbf{I}_{d}
\]

\vphantom{}

II. (Formula for $\tilde{A}_{H,N}\left(\mathbf{z}\right)$) For any
$N\in\mathbb{N}_{0}$ and $\mathbf{z}\in\mathbb{Z}_{p}^{r}$: 
\begin{align}
\tilde{A}_{H,N}\left(\mathbf{z}\right)\overset{\overline{\mathbb{Q}}^{d,d}}{=} & \kappa_{H}\left(\left[\mathbf{z}\right]_{p^{N}}\right)\left(H^{\prime}\left(\mathbf{0}\right)\right)^{N}+\sum_{n=0}^{N-1}\kappa_{H}\left(\left[\mathbf{z}\right]_{p^{n}}\right)\mathbf{I}_{H}\left(\lambda_{p}\left(\left[\mathbf{z}\right]_{p^{n}}\right)\right)\left(H^{\prime}\left(\mathbf{0}\right)\right)^{n}\label{eq:MD Convolution of dA_H and D_N}
\end{align}
In particular, if $H$ is commutative: 
\begin{equation}
\tilde{A}_{H,N}\left(\mathbf{z}\right)\overset{\overline{\mathbb{Q}}^{d,d}}{=}\kappa_{H}\left(\left[\mathbf{z}\right]_{p^{N}}\right)\left(H^{\prime}\left(\mathbf{0}\right)\right)^{N}+\sum_{n=0}^{N-1}\kappa_{H}\left(\left[\mathbf{z}\right]_{p^{n}}\right)\left(H^{\prime}\left(\mathbf{0}\right)\right)^{n}\left(\mathbf{I}_{d}-\alpha_{H}\left(\mathbf{0}\right)\right)\label{eq:MD Convolution of dA_H and D_N when H is commutative}
\end{equation}

\vphantom{}

III. As $N\rightarrow\infty$, \emph{(\ref{eq:MD Convolution of dA_H and D_N})}
converges in the standard $\left(p,q_{H}\right)$-adic frame to:
\begin{equation}
\lim_{N\rightarrow\infty}\tilde{A}_{H,N}\left(\mathbf{z}\right)\overset{\mathcal{F}_{p,q_{H}}^{d,d}}{=}\sum_{n=0}^{\infty}\kappa_{H}\left(\left[\mathbf{z}\right]_{p^{n}}\right)\mathbf{I}_{H}\left(\lambda_{p}\left(\left[\mathbf{z}\right]_{p^{n}}\right)\right)\left(H^{\prime}\left(\mathbf{0}\right)\right)^{n},\textrm{ }\forall\mathbf{z}\in\mathbb{Z}_{p}^{r}\label{eq:MD Fourier Limit of A_H,N twiddle in standard frame}
\end{equation}
There are two particular cases:

i If $\alpha_{H}\left(\mathbf{0}\right)=\mathbf{I}_{d}$, then: 
\begin{equation}
\lim_{N\rightarrow\infty}\tilde{A}_{H,N}\left(\mathbf{z}\right)\overset{\mathcal{F}_{p,q_{H}}^{d,d}}{=}\mathbf{O}_{d},\textrm{ }\forall\mathbf{z}\in\mathbb{Z}_{p}^{r}\label{eq:Full kernel of dA_H on when alpha is 1}
\end{equation}

ii. If $H$ is commutative, then: 
\begin{equation}
\lim_{N\rightarrow\infty}\tilde{A}_{H,N}\left(\mathbf{z}\right)\overset{\mathcal{F}_{p,q_{H}}^{d,d}}{=}\sum_{n=0}^{\infty}\kappa_{H}\left(\left[\mathbf{z}\right]_{p^{n}}\right)\left(H^{\prime}\left(\mathbf{0}\right)\right)^{n}\left(\mathbf{I}_{d}-\alpha_{H}\left(\mathbf{0}\right)\right),\textrm{ }\forall\mathbf{z}\in\mathbb{Z}_{p}^{r}\label{eq:MD Fourier Limit of A_H,N twiddle in standard frame, commutative}
\end{equation}
\end{thm}
Proof:

I. By (\ref{eq:MD alpha_H product expansion}), we have that: 
\[
\hat{A}_{H}\left(\mathbf{t}\right)\overset{\overline{\mathbb{Q}}^{d,d}}{=}\prod_{m=0}^{-v_{p}\left(\mathbf{t}\right)-1}\alpha_{H}\left(p^{m}\mathbf{t}\right)=\left(\sum_{\mathbf{m}=\mathbf{0}}^{p^{-v_{p}\left(\mathbf{t}\right)}-1}\kappa_{H}\left(\mathbf{m}\right)e^{-2\pi i\left(\mathbf{m}\cdot\mathbf{t}\right)}\right)\left(\frac{H^{\prime}\left(\mathbf{0}\right)}{p^{r}}\right)^{-v_{p}\left(\mathbf{t}\right)}
\]

\vphantom{}

II. By (\ref{eq:MD radially-magnitudinal resummation formula}), we
have: 
\begin{align*}
\tilde{A}_{H,N}\left(\mathbf{z}\right) & =p^{rN}\kappa_{H}\left(\left[\mathbf{z}\right]_{p^{N}}\right)\left(\frac{H^{\prime}\left(\mathbf{0}\right)}{p^{r}}\right)^{-v_{p}\left(p^{-N}\right)}\\
 & +\sum_{n=0}^{N-1}p^{rn}\kappa_{H}\left(\left[\mathbf{z}\right]_{p^{n}}\right)\left(\left(\frac{H^{\prime}\left(\mathbf{0}\right)}{p^{r}}\right)^{-v_{p}\left(p^{-n}\right)}-\left(\frac{H^{\prime}\left(\mathbf{0}\right)}{p^{r}}\right)^{-v_{p}\left(p^{-\left(n+1\right)}\right)}\right)\\
 & -\sum_{n=0}^{N-1}p^{rn}\sum_{\mathbf{j}>\mathbf{0}}^{p-1}\kappa_{H}\left(\left[\mathbf{z}\right]_{p^{n}}+\mathbf{j}p^{n}\right)\left(\frac{H^{\prime}\left(\mathbf{0}\right)}{p^{r}}\right)^{-v_{p}\left(p^{-\left(n+1\right)}\right)}
\end{align*}
Here, $-v_{p}\left(p^{-n}\right)=n$. Since we are using the convention:
\begin{equation}
\left(\frac{H^{\prime}\left(\mathbf{0}\right)}{p^{r}}\right)^{-v_{p}\left(\frac{1}{p^{0}}\right)}=\left(\frac{H^{\prime}\left(\mathbf{0}\right)}{p^{r}}\right)^{-v_{p}\left(\mathbf{0}\right)}\overset{\textrm{def}}{=}\mathbf{I}_{d}
\end{equation}
and since $H^{\prime}\left(\mathbf{0}\right)/p^{r}$ is an invertible
$d\times d$ matrix, we have that: 
\begin{equation}
\left(\frac{H^{\prime}\left(\mathbf{0}\right)}{p^{r}}\right)^{-v_{p}\left(1/p^{n}\right)}=\left(\frac{H^{\prime}\left(\mathbf{0}\right)}{p^{r}}\right)^{n},\textrm{ }\forall n\in\mathbb{N}_{0}
\end{equation}
Consequently: 
\begin{align*}
\tilde{A}_{H,N}\left(\mathbf{z}\right) & =p^{rN}\kappa_{H}\left(\left[\mathbf{z}\right]_{p^{n}}\right)\left(\frac{H^{\prime}\left(\mathbf{0}\right)}{p^{r}}\right)^{N}\\
 & +\sum_{n=0}^{N-1}p^{rn}\kappa_{H}\left(\left[\mathbf{z}\right]_{p^{n}}\right)\left(\left(\frac{H^{\prime}\left(\mathbf{0}\right)}{p^{r}}\right)^{n}-\left(\frac{H^{\prime}\left(\mathbf{0}\right)}{p^{r}}\right)^{n+1}\right)\\
 & -\sum_{n=0}^{N-1}p^{rn}\sum_{\mathbf{j}>\mathbf{0}}^{p-1}\kappa_{H}\left(\left[\mathbf{z}\right]_{p^{n}}+\mathbf{j}p^{n}\right)\left(\frac{H^{\prime}\left(\mathbf{0}\right)}{p^{r}}\right)^{n+1}
\end{align*}
Here, (\ref{eq:MD Kappa_H has P-adic structure}) from \textbf{Lemma
\ref{lem:properties of MD kappa_H}} yields: 
\begin{align*}
\kappa_{H}\left(\left[\mathbf{z}\right]_{p^{n}}+p^{n}\mathbf{j}\right) & =\kappa_{H}\left(\left[\mathbf{z}\right]_{p^{n}}\right)\left(H^{\prime}\left(\mathbf{0}\right)\right)^{\lambda_{p}\left(\left[\mathbf{z}\right]_{p^{n}}\right)}\kappa_{H}\left(\mathbf{j}\right)\left(H^{\prime}\left(\mathbf{0}\right)\right)^{-\lambda_{p}\left(\left[\mathbf{z}\right]_{p^{n}}\right)}\\
 & =M_{H}\left(\left[\mathbf{z}\right]_{p^{n}}\right)\kappa_{H}\left(\mathbf{j}\right)\left(H^{\prime}\left(\mathbf{0}\right)\right)^{-\lambda_{p}\left(\left[\mathbf{z}\right]_{p^{n}}\right)}
\end{align*}
As such: 
\begin{align*}
\tilde{A}_{H,N}\left(\mathbf{z}\right) & =\kappa_{H}\left(\left[\mathbf{z}\right]_{p^{N}}\right)\left(H^{\prime}\left(\mathbf{0}\right)\right)^{N}\\
 & +\sum_{n=0}^{N-1}p^{rn}\kappa_{H}\left(\left[\mathbf{z}\right]_{p^{n}}\right)\left(\left(\frac{H^{\prime}\left(\mathbf{0}\right)}{p^{r}}\right)^{n}-\left(\frac{H^{\prime}\left(\mathbf{0}\right)}{p^{r}}\right)^{n+1}\right)\\
 & -\sum_{n=0}^{N-1}p^{rn}M_{H}\left(\left[\mathbf{z}\right]_{p^{n}}\right)\left(\sum_{\mathbf{j}>\mathbf{0}}^{p-1}\kappa_{H}\left(\mathbf{j}\right)\right)\left(H^{\prime}\left(\mathbf{0}\right)\right)^{-\lambda_{p}\left(\left[\mathbf{z}\right]_{p^{n}}\right)}\left(\frac{H^{\prime}\left(\mathbf{0}\right)}{p^{r}}\right)^{n+1}
\end{align*}
Using (\ref{eq:MD kappa H sum in terms of MD alpha}) gives: 
\begin{align*}
\tilde{A}_{H,N}\left(\mathbf{z}\right) & =\kappa_{H}\left(\left[\mathbf{z}\right]_{p^{N}}\right)\left(H^{\prime}\left(\mathbf{0}\right)\right)^{N}\\
 & +\sum_{n=0}^{N-1}p^{rn}\kappa_{H}\left(\left[\mathbf{z}\right]_{p^{n}}\right)\left(\left(\frac{H^{\prime}\left(\mathbf{0}\right)}{p^{r}}\right)^{n}-\left(\frac{H^{\prime}\left(\mathbf{0}\right)}{p^{r}}\right)^{n+1}\right)\\
 & -\sum_{n=0}^{N-1}p^{rn}M_{H}\left(\left[\mathbf{z}\right]_{p^{n}}\right)\left(\alpha_{H}\left(\mathbf{0}\right)\left(\frac{H^{\prime}\left(\mathbf{0}\right)}{p^{r}}\right)^{-1}-\mathbf{I}_{d}\right)\left(H^{\prime}\left(\mathbf{0}\right)\right)^{-\lambda_{p}\left(\left[\mathbf{z}\right]_{p^{n}}\right)}\left(\frac{H^{\prime}\left(\mathbf{0}\right)}{p^{r}}\right)^{n+1}
\end{align*}
Re-writing things to make some cancellations evident:
\begin{align*}
\tilde{A}_{H,N}\left(\mathbf{z}\right) & =\kappa_{H}\left(\left[\mathbf{z}\right]_{p^{N}}\right)\left(H^{\prime}\left(\mathbf{0}\right)\right)^{N}\\
 & +\sum_{n=0}^{N-1}p^{rn}\kappa_{H}\left(\left[\mathbf{z}\right]_{p^{n}}\right)\left(\frac{H^{\prime}\left(\mathbf{0}\right)}{p^{r}}\right)^{n}-\overbrace{\sum_{n=0}^{N-1}p^{rn}\kappa_{H}\left(\left[\mathbf{z}\right]_{p^{n}}\right)\left(\frac{H^{\prime}\left(\mathbf{0}\right)}{p^{r}}\right)^{n+1}}^{\textrm{these two}}\\
 & -\sum_{n=0}^{N-1}p^{rn}M_{H}\left(\left[\mathbf{z}\right]_{p^{n}}\right)\alpha_{H}\left(\mathbf{0}\right)\left(H^{\prime}\left(\mathbf{0}\right)\right)^{-\lambda_{p}\left(\left[\mathbf{z}\right]_{p^{n}}\right)}\left(\frac{H^{\prime}\left(\mathbf{0}\right)}{p^{r}}\right)^{n}\\
 & +\underbrace{\sum_{n=0}^{N-1}p^{rn}\overbrace{M_{H}\left(\left[\mathbf{z}\right]_{p^{n}}\right)\left(H^{\prime}\left(\mathbf{0}\right)\right)^{-\lambda_{p}\left(\left[\mathbf{z}\right]_{p^{n}}\right)}}^{\kappa_{H}\left(\left[\mathbf{z}\right]_{p^{n}}\right)}\left(\frac{H^{\prime}\left(\mathbf{0}\right)}{p^{r}}\right)^{n+1}}_{\textrm{cancel out}}
\end{align*}
we end up with:
\begin{align*}
\tilde{A}_{H,N}\left(\mathbf{z}\right) & =\kappa_{H}\left(\left[\mathbf{z}\right]_{p^{N}}\right)\left(H^{\prime}\left(\mathbf{0}\right)\right)^{N}+\sum_{n=0}^{N-1}p^{rn}\kappa_{H}\left(\left[\mathbf{z}\right]_{p^{n}}\right)\left(\frac{H^{\prime}\left(\mathbf{0}\right)}{p^{r}}\right)^{n}\\
 & -\sum_{n=0}^{N-1}p^{rn}M_{H}\left(\left[\mathbf{z}\right]_{p^{n}}\right)\alpha_{H}\left(\mathbf{0}\right)\left(H^{\prime}\left(\mathbf{0}\right)\right)^{-\lambda_{p}\left(\left[\mathbf{z}\right]_{p^{n}}\right)}\left(\frac{H^{\prime}\left(\mathbf{0}\right)}{p^{r}}\right)^{n}
\end{align*}

Now we utilize our $\mathcal{C}_{H}$ notation to make the lower line
manageable. First, recall that:
\begin{equation}
M_{H}\left(\left[\mathbf{z}\right]_{p^{n}}\right)=\kappa_{H}\left(\left[\mathbf{z}\right]_{p^{n}}\right)\left(H^{\prime}\left(\mathbf{0}\right)\right)^{\lambda_{p}\left(\left[\mathbf{z}\right]_{p^{n}}\right)}
\end{equation}
With this, the expression:
\begin{equation}
M_{H}\left(\left[\mathbf{z}\right]_{p^{n}}\right)\alpha_{H}\left(\mathbf{0}\right)\left(H^{\prime}\left(\mathbf{0}\right)\right)^{-\lambda_{p}\left(\left[\mathbf{z}\right]_{p^{n}}\right)}
\end{equation}
becomes: 
\begin{equation}
\kappa_{H}\left(\left[\mathbf{z}\right]_{p^{n}}\right)\left(H^{\prime}\left(\mathbf{0}\right)\right)^{\lambda_{p}\left(\left[\mathbf{z}\right]_{p^{n}}\right)}\alpha_{H}\left(\mathbf{0}\right)\left(H^{\prime}\left(\mathbf{0}\right)\right)^{-\lambda_{p}\left(\left[\mathbf{z}\right]_{p^{n}}\right)}
\end{equation}
which, using $\mathcal{C}_{H}$, is just:
\begin{equation}
\kappa_{H}\left(\left[\mathbf{z}\right]_{p^{n}}\right)\mathcal{C}_{H}\left(\alpha_{H}\left(\mathbf{0}\right):\lambda_{p}\left(\left[\mathbf{z}\right]_{p^{n}}\right)\right)
\end{equation}
Finally, writing out the formula for $\tilde{A}_{H,N}\left(\mathbf{z}\right)$,
we have:
\begin{align*}
\tilde{A}_{H,N}\left(\mathbf{z}\right) & =\kappa_{H}\left(\left[\mathbf{z}\right]_{p^{N}}\right)\left(H^{\prime}\left(\mathbf{0}\right)\right)^{N}+\sum_{n=0}^{N-1}p^{rn}\kappa_{H}\left(\left[\mathbf{z}\right]_{p^{n}}\right)\left(\frac{H^{\prime}\left(\mathbf{0}\right)}{p^{r}}\right)^{n}\\
 & -\sum_{n=0}^{N-1}p^{rn}\kappa_{H}\left(\left[\mathbf{z}\right]_{p^{n}}\right)\mathcal{C}_{H}\left(\alpha_{H}\left(\mathbf{0}\right):\lambda_{p}\left(\left[\mathbf{z}\right]_{p^{n}}\right)\right)\left(\frac{H^{\prime}\left(\mathbf{0}\right)}{p^{r}}\right)^{n}\\
 & =\kappa_{H}\left(\left[\mathbf{z}\right]_{p^{N}}\right)\left(H^{\prime}\left(\mathbf{0}\right)\right)^{N}\\
 & +\sum_{n=0}^{N-1}\kappa_{H}\left(\left[\mathbf{z}\right]_{p^{n}}\right)\underbrace{\left(\mathbf{I}_{d}-\mathcal{C}_{H}\left(\alpha_{H}\left(\mathbf{0}\right):\lambda_{p}\left(\left[\mathbf{z}\right]_{p^{n}}\right)\right)\right)}_{\mathbf{I}_{H}\left(\lambda_{p}\left(\left[\mathbf{z}\right]_{p^{n}}\right)\right)}\left(H^{\prime}\left(\mathbf{0}\right)\right)^{n}
\end{align*}
which is the desired identity.

\vphantom{}

III. Let $H$ be semi-basic, and fix $\mathbf{z}\in\left(\mathbb{Z}_{p}^{r}\right)^{\prime}$.
Then: 
\begin{equation}
\left\Vert \kappa_{H}\left(\left[\mathbf{z}\right]_{p^{N}}\right)\left(H^{\prime}\left(\mathbf{0}\right)\right)^{N}\right\Vert _{q}\leq\left\Vert \kappa_{H}\left(\left[\mathbf{z}\right]_{p^{N}}\right)\right\Vert _{q}
\end{equation}
Applying \textbf{Lemma \ref{lem:properties of MD kappa_H}}, $\kappa_{H}\left(\left[\mathbf{z}\right]_{p^{N}}\right)$
converges $q$-adically to $\mathbf{O}_{d}$ as $N\rightarrow\infty$.
By the ultrametric structure of $\mathbb{C}_{q}$, this then guarantees
the convergence of (\ref{eq:MD Fourier Limit of A_H,N twiddle in standard frame}).

Next, let $\mathfrak{\mathbf{z}}=\mathbf{m}\in\mathbb{N}_{0}^{r}$.
Then, for all $n\geq\lambda_{p}\left(\mathbf{m}\right)$:
\begin{equation}
\kappa_{H}\left(\mathbf{m}\right)=\kappa_{H}\left(\left[\mathbf{m}\right]_{p^{n}}\right)=\kappa_{H}\left(\left[\mathbf{m}\right]_{p^{n+1}}\right)\in\textrm{Gl}_{d}\left(\mathbb{Q}\right)\backslash\left\{ \mathbf{O}_{d}\right\} 
\end{equation}
Likewise, keeping $n\geq\lambda_{p}\left(\mathbf{m}\right)$, we have:
\begin{align*}
\mathbf{I}_{H}\left(\lambda_{p}\left(\left[\mathbf{m}\right]_{p^{n}}\right)\right) & =\mathbf{I}_{d}-\left(H^{\prime}\left(\mathbf{0}\right)\right)^{\lambda_{p}\left(\left[\mathbf{m}\right]_{p^{n}}\right)}\alpha_{H}\left(\mathbf{0}\right)\left(H^{\prime}\left(\mathbf{0}\right)\right)^{-\lambda_{p}\left(\left[\mathbf{m}\right]_{p^{n}}\right)}\\
 & =\mathbf{I}_{d}-\left(H^{\prime}\left(\mathbf{0}\right)\right)^{\lambda_{p}\left(\mathbf{m}\right)}\alpha_{H}\left(\mathbf{0}\right)\left(H^{\prime}\left(\mathbf{0}\right)\right)^{-\lambda_{p}\left(\mathbf{m}\right)}\\
 & =\mathbf{I}_{d}-\mathcal{C}_{H}\left(\alpha_{H}\left(\mathbf{0}\right):\lambda_{p}\left(\mathbf{m}\right)\right)\\
 & =\mathbf{I}_{H}\left(\lambda_{p}\left(\mathbf{m}\right)\right)
\end{align*}

Examining the tail of the $n$-series in (\ref{eq:MD Convolution of dA_H and D_N}),
note that: 
\[
\sum_{n=0}^{N-1}\kappa_{H}\left(\left[\mathbf{m}\right]_{p^{n}}\right)\mathbf{I}_{H}\left(\lambda_{p}\left(\mathbf{m}\right)\right)\left(H^{\prime}\left(\mathbf{0}\right)\right)^{n}\overset{\overline{\mathbb{Q}}^{d,d}}{=}O\left(\mathbf{I}_{d}\right)+\kappa_{H}\left(\mathbf{m}\right)\mathbf{I}_{H}\left(\lambda_{p}\left(\mathbf{m}\right)\right)\sum_{n=\lambda_{p}\left(\mathbf{m}\right)}^{N-1}\left(H^{\prime}\left(\mathbf{0}\right)\right)^{n}
\]
holds whenever $N\geq\lambda_{p}\left(\mathbf{m}\right)+1$. Since
the $n$-series on the right is geometric in the matrix $H^{\prime}\left(\mathbf{0}\right)$,
it will converge in the topology of in $\mathbb{C}$ as $N\rightarrow\infty$
if and only if $H$ is contracting. Likewise, by \textbf{Lemma \ref{lem:properties of MD kappa_H}},\textbf{
}$\kappa_{H}\left(\left[\mathbf{z}\right]_{p^{N}}\right)\left(H^{\prime}\left(\mathbf{0}\right)\right)^{N}$
will tend to $\mathbf{O}_{d}$ in the topology of $\mathbb{C}$ as
$N\rightarrow\infty$ if and only if $H$ is contracting, and will
diverge in the topology of $\mathbb{C}$ if and only if $H$ is expanding.

This shows that (\ref{eq:MD Fourier Limit of A_H,N twiddle in standard frame})
holds true.

Q.E.D.

\vphantom{}

Now we can show that $\hat{A}_{H}$ is the Fourier-Stieltjes transform
of a thick measure of matrix type.
\begin{cor}[\textbf{$\hat{A}_{H}$ Induces a Matrix-Type Thick Measure}]
\ 

\vphantom{}

I. The matrix-valued function $\hat{A}_{H}\left(\mathbf{t}\right)$
is the Fourier-Stieltjes transform of a $\left(p,q_{H}\right)$-adic
matrix-type thick measure of depth $r$. This measure acts on functions
by way of the formula:
\begin{equation}
\mathbf{f}\in C\left(\mathbb{Z}_{p}^{r},\mathbb{C}_{q}^{d}\right)\mapsto\sum_{\mathbf{t}\in\hat{\mathbb{Z}}_{p}^{r}}\hat{A}_{H}\left(-\mathbf{t}\right)\hat{\mathbf{f}}\left(\mathbf{t}\right)\in\mathbb{C}_{q}^{d}
\end{equation}
Moreover, this thick measure is degenerate with respect to the standard
$\left(p,q_{H}\right)$-adic frame if and only if $\alpha_{H}\left(\mathbf{0}\right)=\mathbf{I}_{d}$.

\vphantom{}

II. For any $\mathbf{v}\in\mathbb{C}_{q}^{d}$, the vector-valued
function $\hat{A}_{H}\left(\mathbf{t}\right)\mathbf{v}$ is the Fourier-Stieltjes
transform of a $\left(p,q_{H}\right)$-adic thick measure of vector
type. This measure acts on functions by way of the formula: 
\begin{equation}
\mathbf{f}\in C\left(\mathbb{Z}_{p}^{r},\mathbb{C}_{q}^{d}\right)\mapsto\sum_{\mathbf{t}\in\hat{\mathbb{Z}}_{p}^{r}}\left(\hat{A}_{H}\left(-\mathbf{t}\right)\mathbf{v}\right)\hat{\mathbf{f}}\left(\mathbf{t}\right)\in\mathbb{C}_{q}^{d}
\end{equation}
Moreover, this thick measure is degenerate with respect to the standard
$\left(p,q_{H}\right)$-adic frame if and only if $\alpha_{H}\left(\mathbf{0}\right)=\mathbf{I}_{d}$. 
\end{cor}
Proof:

I. Follows by (II) and (III) of \textbf{Theorem \ref{thm:MD properties of A_H hat}}.

\vphantom{}

II. (I) shows that the matrix-type thick measure associated to $\hat{A}_{H}\left(\mathbf{t}\right)$
is degenerate if and only if $\alpha_{H}\left(\mathbf{0}\right)=\mathbf{I}_{d}$.
Noting that the vector-valued Fourier series generated by $\hat{A}_{H}\left(\mathbf{t}\right)\mathbf{v}$
has an $\mathcal{F}_{p,q}$ limit of: 
\begin{equation}
\lim_{N\rightarrow\infty}\sum_{\left\Vert \mathbf{t}\right\Vert _{p}\leq p^{N}}\left(\hat{A}_{H}\left(\mathbf{t}\right)\mathbf{v}\right)e^{2\pi i\left\{ \mathbf{t}\mathbf{z}\right\} _{p}}
\end{equation}
we can write this as: 
\begin{equation}
\left(\lim_{N\rightarrow\infty}\sum_{\left\Vert \mathbf{t}\right\Vert _{p}\leq p^{N}}\hat{A}_{H}\left(\mathbf{t}\right)e^{2\pi i\left\{ \mathbf{t}\mathbf{z}\right\} _{p}}\right)\mathbf{v}
\end{equation}
Note that (I) is exactly the statement that the limit in parenthesis
(taken with respect to $\mathcal{F}_{p,q}$) is $\mathbf{O}_{d}$
for every $\mathbf{z}\in\mathbb{Z}_{p}^{r}$ if and only if $\alpha_{H}\left(\mathbf{0}\right)=\mathbf{I}_{d}$.
Hence: 
\begin{align*}
\lim_{N\rightarrow\infty}\sum_{\left\Vert \mathbf{t}\right\Vert _{p}\leq p^{N}}\left(\hat{A}_{H}\left(\mathbf{t}\right)\mathbf{v}\right)e^{2\pi i\left\{ \mathbf{t}\mathbf{z}\right\} _{p}} & =\left(\lim_{N\rightarrow\infty}\sum_{\left\Vert \mathbf{t}\right\Vert _{p}\leq p^{N}}\hat{A}_{H}\left(\mathbf{t}\right)e^{2\pi i\left\{ \mathbf{t}\mathbf{z}\right\} _{p}}\right)\mathbf{v}\\
\left(\textrm{iff }\alpha_{H}\left(\mathbf{0}\right)=\mathbf{I}_{d}\right); & =\mathbf{O}_{d}\mathbf{v}\\
 & =\mathbf{0}
\end{align*}
This shows that the vector-type thick measure induced by $\hat{A}_{H}\left(\mathbf{t}\right)\mathbf{v}$
is degenerate if and only if $\alpha_{H}\left(\mathbf{0}\right)=\mathbf{I}_{d}$.

Q.E.D.
\begin{defn}[$\tilde{A}_{H}\left(\mathbf{z}\right)$]
We write $\tilde{A}_{H}\left(\mathbf{z}\right)$ to denote the kernel
of the matrix-type thick measure induced by $\hat{A}_{H}$; that is,
$\tilde{A}_{H}\left(\mathbf{z}\right)$ is the $\mathcal{F}_{p,q_{H}}^{d,d}$-limit
of $\tilde{A}_{H,N}\left(\mathbf{z}\right)$ as $N\rightarrow\infty$.
\end{defn}
\vphantom{}

Like in the one-dimensional case, we now can write out our asymptotic
decomposition of $\hat{\chi}_{H,N}$.
\begin{thm}[\textbf{$\left(N,\mathbf{t}\right)$ Asymptotic Decomposition of $\hat{\chi}_{H,N}$}]
\label{thm:MD N,t asympotics for Chi_H,N hat}\ 

\vphantom{}

I. If $\alpha_{H}\left(\mathbf{0}\right)=\mathbf{I}_{d}$, then: 
\begin{equation}
\hat{\chi}_{H,N}\left(\mathbf{t}\right)\overset{\overline{\mathbb{Q}}^{d}}{=}\begin{cases}
N\hat{A}_{H}\left(\mathbf{t}\right)\beta_{H}\left(\mathbf{0}\right) & \textrm{if }\mathbf{t}=\mathbf{0}\\
\hat{A}_{H}\left(\mathbf{t}\right)\left(\left(N+v_{p}\left(\mathbf{t}\right)\right)\beta_{H}\left(\mathbf{0}\right)+\gamma_{H}\left(\frac{\mathbf{t}\left|\mathbf{t}\right|_{p}}{p}\right)\right) & \textrm{if }0<\left\Vert \mathbf{t}\right\Vert _{p}<p^{N}\\
\hat{A}_{H}\left(\mathbf{t}\right)\gamma_{H}\left(\frac{\mathbf{t}\left|\mathbf{t}\right|_{p}}{p}\right) & \textrm{if }\left\Vert \mathbf{t}\right\Vert _{p}=p^{N}\\
\mathbf{0} & \textrm{if }\left\Vert \mathbf{t}\right\Vert _{p}>p^{N}
\end{cases},\textrm{ }\forall\mathbf{t}\in\hat{\mathbb{Z}}_{p}^{r}\label{eq:MD Formula for Chi_H,N hat when alpha is 1}
\end{equation}

\vphantom{}

II. If $\mathbf{I}_{d}-\alpha_{H}\left(\mathbf{0}\right)$ is invertible,
then: 
\begin{equation}
\hat{\chi}_{H,N}\left(\mathbf{t}\right)\overset{\overline{\mathbb{Q}}^{d}}{=}\begin{cases}
\hat{A}_{H}\left(\mathbf{t}\right)\frac{\mathbf{I}_{d}-\left(\alpha_{H}\left(\mathbf{0}\right)\right)^{N}}{\mathbf{I}_{d}-\alpha_{H}\left(\mathbf{0}\right)}\beta_{H}\left(\mathbf{0}\right) & \textrm{if }\mathbf{t}=\mathbf{0}\\
\hat{A}_{H}\left(\mathbf{t}\right)\left(\gamma_{H}\left(\frac{\mathbf{t}\left|\mathbf{t}\right|_{p}}{p}\right)+\frac{\mathbf{I}_{d}-\left(\alpha_{H}\left(\mathbf{0}\right)\right)^{N+v_{p}\left(\mathbf{t}\right)}}{\mathbf{I}_{d}-\alpha_{H}\left(\mathbf{0}\right)}\beta_{H}\left(\mathbf{0}\right)\right) & \textrm{if }0<\left\Vert \mathbf{t}\right\Vert _{p}<p^{N}\\
\hat{A}_{H}\left(\mathbf{t}\right)\gamma_{H}\left(\frac{\mathbf{t}\left|\mathbf{t}\right|_{p}}{p}\right) & \textrm{if }\left\Vert \mathbf{t}\right\Vert _{p}=p^{N}\\
\mathbf{0} & \textrm{if }\left\Vert \mathbf{t}\right\Vert _{p}>p^{N}
\end{cases},\textrm{ }\forall\mathbf{t}\in\hat{\mathbb{Z}}_{p}^{r}\label{eq:MD Formula for Chi_H,N hat when alpha is not 1}
\end{equation}
\end{thm}
Proof: For brevity, we write: 
\begin{equation}
\hat{A}_{H,n}\left(\mathbf{t}\right)\overset{\textrm{def}}{=}\begin{cases}
\mathbf{I}_{d} & \textrm{if }n=0\\
\mathbf{1}_{\mathbf{0}}\left(p^{n+1}\mathbf{t}\right)\prod_{m=0}^{n-1}\alpha_{H}\left(p^{m}\mathbf{t}\right) & \textrm{if }n\geq1
\end{cases}\label{eq:MD Definition of A_H,n+1 hat}
\end{equation}
so that (\ref{eq:MD Chi_H,N hat functional equation}) becomes: 
\begin{equation}
\hat{\chi}_{H,N}\left(\mathbf{t}\right)=\sum_{n=0}^{N-1}\hat{A}_{H,n}\left(\mathbf{t}\right)\beta_{H}\left(p^{n}\mathbf{t}\right)
\end{equation}
Now, letting $\mathbf{t}\in\hat{\mathbb{Z}}_{p}^{r}$ be non-zero
and satisfy $\left\Vert \mathbf{t}\right\Vert _{p}\leq p^{N}$, we
use \textbf{Proposition \ref{prop:MD alpha_H A_H hat product}} to
write: 
\begin{equation}
\hat{A}_{H,n}\left(\mathbf{t}\right)=\begin{cases}
\mathbf{O}_{d} & \textrm{if }n\leq-v_{p}\left(\mathbf{t}\right)-2\\
\frac{\hat{A}_{H}\left(\mathbf{t}\right)}{\alpha_{H}\left(p^{-v_{p}\left(\mathbf{t}\right)-1}\mathbf{t}\right)} & \textrm{if }n=-v_{p}\left(\mathbf{t}\right)-1\\
\left(\alpha_{H}\left(\mathbf{0}\right)\right)^{n+v_{p}\left(\mathbf{t}\right)}\hat{A}_{H}\left(\mathbf{t}\right) & \textrm{if }n\geq-v_{p}\left(\mathbf{t}\right)
\end{cases}
\end{equation}
So: 
\begin{equation}
\hat{\chi}_{H,N}\left(\mathbf{t}\right)\overset{\overline{\mathbb{Q}}^{d}}{=}\sum_{n=0}^{N-1}\hat{A}_{H,n}\left(\mathbf{t}\right)\beta_{H}\left(p^{n}\mathbf{t}\right)=\sum_{n=-v_{p}\left(\mathbf{t}\right)-1}^{N-1}\hat{A}_{H,n}\left(\mathbf{t}\right)\beta_{H}\left(p^{n}\mathbf{t}\right)\label{eq:MD Chi_H,N hat as Beta_H plus A_H,n hat - ready for t,n analysis}
\end{equation}
Note that in obtaining (\ref{eq:MD Chi_H,N hat as Beta_H plus A_H,n hat - ready for t,n analysis}),
we have used the fact: 
\begin{equation}
\hat{A}_{H,n}\left(\mathbf{t}\right)=\mathbf{O}_{d},\textrm{ }\forall n\leq-v_{p}\left(\mathbf{t}\right)-2
\end{equation}
which is, itself, a consequence of the vanishing of $\mathbf{1}_{\mathbf{0}}\left(p^{n+1}\mathbf{t}\right)$
whenever $n\leq-v_{p}\left(\mathbf{t}\right)-2$. With that done,
we are left with two cases.

\vphantom{}

\textbullet{} First, suppose $\left\Vert \mathbf{t}\right\Vert _{p}=p^{N}$.
Then $N=-v_{p}\left(\mathbf{t}\right)$, and so, (\ref{eq:MD Chi_H,N hat as Beta_H plus A_H,n hat - ready for t,n analysis})
becomes: 
\begin{align*}
\hat{\chi}_{H,N}\left(\mathbf{t}\right) & =\sum_{n=-v_{p}\left(\mathbf{t}\right)-1}^{N-1}\hat{A}_{H,n}\left(\mathbf{t}\right)\beta_{H}\left(p^{n}\mathbf{t}\right)\\
 & =\hat{A}_{H,-v_{p}\left(\mathbf{t}\right)-1}\left(\mathbf{t}\right)\beta_{H}\left(p^{-v_{p}\left(\mathbf{t}\right)-1}\mathbf{t}\right)\\
 & =\left(\mathbf{1}_{\mathbf{0}}\left(p^{-v_{p}\left(\mathbf{t}\right)-1+1}\mathbf{t}\right)\prod_{m=0}^{-v_{p}\left(\mathbf{t}\right)-1-1}\alpha_{H}\left(p^{m}\mathbf{t}\right)\right)\beta_{H}\left(p^{-v_{p}\left(\mathbf{t}\right)-1}\mathbf{t}\right)\\
 & =\left(\mathbf{1}_{\mathbf{0}}\left(p^{-v_{p}\left(\mathbf{t}\right)}\mathbf{t}\right)\prod_{m=0}^{-v_{p}\left(\mathbf{t}\right)-1}\alpha_{H}\left(p^{m}\mathbf{t}\right)\right)\left(\alpha_{H}\left(p^{-v_{p}\left(\mathbf{t}\right)-1}\mathbf{t}\right)\right)^{-1}\beta_{H}\left(p^{-v_{p}\left(\mathbf{t}\right)-1}\mathbf{t}\right)
\end{align*}
Because $p^{-v_{p}\left(\mathbf{t}\right)}\mathbf{t}$ contains only
integer entries, note that $\mathbf{1}_{\mathbf{0}}\left(p^{-v_{p}\left(\mathbf{t}\right)}\mathbf{t}\right)=1$.
This leaves us with: 
\begin{align*}
\hat{\chi}_{H,N}\left(\mathbf{t}\right) & =\left(\prod_{m=0}^{-v_{p}\left(\mathbf{t}\right)-1}\alpha_{H}\left(p^{m}\mathbf{t}\right)\right)\left(\alpha_{H}\left(p^{-v_{p}\left(\mathbf{t}\right)-1}\mathbf{t}\right)\right)^{-1}\beta_{H}\left(p^{-v_{p}\left(\mathbf{t}\right)-1}\mathbf{t}\right)\\
 & =\hat{A}_{H}\left(\mathbf{t}\right)\gamma_{H}\left(\frac{\mathbf{t}\left|\mathbf{t}\right|_{p}}{p}\right)
\end{align*}

\vphantom{}

\textbullet{} Second, suppose $0<\left\Vert \mathbf{t}\right\Vert _{p}<p^{N}$.
Then $-v_{p}\left(\mathbf{t}\right)-1<N-1$, and so $p^{n}\mathbf{t}\overset{1}{\equiv}\mathbf{0}$
for all $n\geq-v_{p}\left(\mathbf{t}\right)$. Consequently, (\ref{eq:MD Chi_H,N hat as Beta_H plus A_H,n hat - ready for t,n analysis})
becomes:

\begin{align*}
\hat{\chi}_{H,N}\left(\mathbf{t}\right) & =\hat{A}_{H,-v_{p}\left(\mathbf{t}\right)-1}\left(\mathbf{t}\right)\beta_{H}\left(p^{-v_{p}\left(\mathbf{t}\right)-1}\mathbf{t}\right)+\sum_{n=-v_{p}\left(\mathbf{t}\right)}^{N-1}\hat{A}_{H,n}\left(\mathbf{t}\right)\beta_{H}\left(p^{n}\mathbf{t}\right)\\
\left(p^{n}\mathbf{t}\overset{1}{\equiv}\mathbf{0}\textrm{ }\forall n\geq-v_{p}\left(\mathbf{t}\right)\right); & =\hat{A}_{H}\left(\mathbf{t}\right)\gamma_{H}\left(\frac{\mathbf{t}\left|\mathbf{t}\right|_{p}}{p}\right)+\sum_{n=-v_{p}\left(\mathbf{t}\right)}^{N-1}\hat{A}_{H,n}\left(\mathbf{t}\right)\beta_{H}\left(\mathbf{0}\right)
\end{align*}
Now, using \textbf{Proposition \ref{prop:MD alpha_H A_H hat product}},
we write:
\begin{equation}
\hat{A}_{H,n}\left(\mathbf{t}\right)=\hat{A}_{H}\left(\mathbf{t}\right)\left(\alpha_{H}\left(\mathbf{0}\right)\right)^{n+v_{p}\left(\mathbf{t}\right)},\textrm{ }\forall n\geq-v_{p}\left(\mathbf{t}\right)
\end{equation}
and hence: 
\begin{align*}
\hat{\chi}_{H,N}\left(\mathbf{t}\right) & =\hat{A}_{H}\left(\mathbf{t}\right)\gamma_{H}\left(\frac{\mathbf{t}\left|\mathbf{t}\right|_{p}}{p}\right)+\sum_{n=-v_{p}\left(\mathbf{t}\right)}^{N-1}\hat{A}_{H}\left(\mathbf{t}\right)\left(\alpha_{H}\left(\mathbf{0}\right)\right)^{n+v_{p}\left(\mathbf{t}\right)}\beta_{H}\left(\mathbf{0}\right)\\
 & =\hat{A}_{H}\left(\mathbf{t}\right)\gamma_{H}\left(\frac{\mathbf{t}\left|\mathbf{t}\right|_{p}}{p}\right)+\hat{A}_{H}\left(\mathbf{t}\right)\sum_{n=0}^{N+v_{p}\left(\mathbf{t}\right)-1}\left(\alpha_{H}\left(\mathbf{0}\right)\right)^{n}\beta_{H}\left(\mathbf{0}\right)
\end{align*}
When $\alpha_{H}\left(\mathbf{0}\right)=\mathbf{I}_{d}$, this gives:
\begin{equation}
\hat{\chi}_{H,N}\left(\mathbf{t}\right)=\hat{A}_{H}\left(\mathbf{t}\right)\gamma_{H}\left(\frac{\mathbf{t}\left|\mathbf{t}\right|_{p}}{p}\right)+\left(N+v_{p}\left(\mathbf{t}\right)\right)\hat{A}_{H}\left(\mathbf{t}\right)\beta_{H}\left(\mathbf{0}\right)\label{eq:MD asymptotic analysis, nearly done}
\end{equation}
On the other hand, suppose that $\mathbf{I}_{d}-\alpha_{H}\left(\mathbf{0}\right)$
is invertible. Then, since: 
\begin{align}
\left(\mathbf{I}_{d}-\alpha_{H}\left(\mathbf{0}\right)\right)\sum_{m=0}^{M-1}\left(\alpha_{H}\left(\mathbf{0}\right)\right)^{m} & =\mathbf{I}_{d}-\left(\alpha_{H}\left(\mathbf{0}\right)\right)^{M}=\sum_{m=0}^{M-1}\left(\alpha_{H}\left(\mathbf{0}\right)\right)^{m}\left(\mathbf{I}_{d}-\alpha_{H}\left(\mathbf{0}\right)\right)
\end{align}
it follows that: 
\begin{align*}
\sum_{m=0}^{M-1}\left(\alpha_{H}\left(\mathbf{0}\right)\right)^{m} & =\left(\mathbf{I}_{d}-\alpha_{H}\left(\mathbf{0}\right)\right)^{-1}\left(\mathbf{I}_{d}-\left(\alpha_{H}\left(\mathbf{0}\right)\right)^{M}\right)\\
 & =\left(\mathbf{I}_{d}-\left(\alpha_{H}\left(\mathbf{0}\right)\right)^{M}\right)\left(\mathbf{I}_{d}-\alpha_{H}\left(\mathbf{0}\right)\right)^{-1}
\end{align*}
With this, the right-hand side of (\ref{eq:MD asymptotic analysis, nearly done})
becomes: 
\[
\begin{cases}
\hat{A}_{H}\left(\mathbf{t}\right)\gamma_{H}\left(\frac{\mathbf{t}\left|\mathbf{t}\right|_{p}}{p}\right)+\left(N+v_{p}\left(\mathbf{t}\right)\right)\hat{A}_{H}\left(\mathbf{t}\right)\beta_{H}\left(\mathbf{0}\right) & \textrm{if }\alpha_{H}\left(\mathbf{0}\right)=\mathbf{I}_{d}\\
\hat{A}_{H}\left(\mathbf{t}\right)\gamma_{H}\left(\frac{\mathbf{t}\left|\mathbf{t}\right|_{p}}{p}\right)+\hat{A}_{H}\left(\mathbf{t}\right)\frac{\mathbf{I}_{d}-\left(\alpha_{H}\left(\mathbf{0}\right)\right)^{N+v_{p}\left(\mathbf{t}\right)}}{\mathbf{I}_{d}-\alpha_{H}\left(\mathbf{0}\right)}\beta_{H}\left(\mathbf{0}\right) & \textrm{if }\det\left(\mathbf{I}_{d}-\alpha_{H}\left(\mathbf{0}\right)\right)\neq0
\end{cases}
\]

Finally, when $\mathbf{t}=\mathbf{0}$: 
\begin{align*}
\hat{\chi}_{H,N}\left(\mathbf{0}\right) & =\sum_{n=0}^{N-1}\hat{A}_{H,n}\left(\mathbf{0}\right)\beta_{H}\left(\mathbf{0}\right)\\
 & =\beta_{H}\left(\mathbf{0}\right)+\sum_{n=1}^{N-1}\mathbf{1}_{\mathbf{0}}\left(\mathbf{0}\right)\left(\prod_{m=0}^{n-1}\alpha_{H}\left(\mathbf{0}\right)\right)\beta_{H}\left(\mathbf{0}\right)\\
 & =\beta_{H}\left(\mathbf{0}\right)+\sum_{n=1}^{N-1}\left(\alpha_{H}\left(\mathbf{0}\right)\right)^{n}\beta_{H}\left(\mathbf{0}\right)\\
 & =\begin{cases}
\beta_{H}\left(\mathbf{0}\right)N & \textrm{if }\alpha_{H}\left(\mathbf{0}\right)=\mathbf{I}_{d}\\
\frac{\mathbf{I}_{d}-\left(\alpha_{H}\left(\mathbf{0}\right)\right)^{N}}{\mathbf{I}_{d}-\alpha_{H}\left(\mathbf{0}\right)}\beta_{H}\left(\mathbf{0}\right) & \textrm{if }\det\left(\mathbf{I}_{d}-\alpha_{H}\left(\mathbf{0}\right)\right)\neq0
\end{cases}
\end{align*}
Since $\hat{A}_{H}\left(\mathbf{0}\right)=\mathbf{I}_{d}$, we have
that $\beta_{H}\left(\mathbf{0}\right)N=\hat{A}_{H}\left(\mathbf{t}\right)\beta_{H}\left(\mathbf{0}\right)N$
when $\mathbf{t}=\mathbf{0}$ (in the $\alpha_{H}\left(\mathbf{0}\right)=\mathbf{I}_{d}$
case) and: 
\[
\frac{\mathbf{I}_{d}-\left(\alpha_{H}\left(\mathbf{0}\right)\right)^{N}}{\mathbf{I}_{d}-\alpha_{H}\left(\mathbf{0}\right)}\beta_{H}\left(\mathbf{0}\right)=\hat{A}_{H}\left(\mathbf{t}\right)\frac{\mathbf{I}_{d}-\left(\alpha_{H}\left(\mathbf{0}\right)\right)^{N}}{\mathbf{I}_{d}-\alpha_{H}\left(\mathbf{0}\right)}\beta_{H}\left(\mathbf{0}\right)
\]
when $\mathbf{t=0}$ in the $\det\left(\mathbf{I}_{d}-\alpha_{H}\left(\mathbf{0}\right)\right)\neq0$
case.

Q.E.D.

\vphantom{}

When $\alpha_{H}\left(\mathbf{0}\right)=\mathbf{I}_{d}$, we have
that 
\[
\hat{\chi}_{H,N}\left(\mathbf{t}\right)-N\mathbf{1}_{\mathbf{0}}\left(p^{N-1}\mathbf{t}\right)\hat{A}_{H}\left(\mathbf{t}\right)\beta_{H}\left(\mathbf{0}\right)
\]
is given by:
\begin{equation}
\begin{cases}
\mathbf{0} & \textrm{if }\mathbf{t}=\mathbf{0}\\
\hat{A}_{H}\left(\mathbf{t}\right)\left(v_{p}\left(\mathbf{t}\right)\beta_{H}\left(\mathbf{0}\right)+\gamma_{H}\left(\frac{\mathbf{t}\left|\mathbf{t}\right|_{p}}{p}\right)\right) & \textrm{if }0<\left\Vert \mathbf{t}\right\Vert _{p}<p^{N}\\
\hat{A}_{H}\left(\mathbf{t}\right)\gamma_{H}\left(\frac{\mathbf{t}\left|\mathbf{t}\right|_{p}}{p}\right) & \textrm{if }\left\Vert \mathbf{t}\right\Vert _{p}=p^{N}\\
\mathbf{0} & \textrm{if }\left\Vert \mathbf{t}\right\Vert _{p}>p^{N}
\end{cases}\label{eq:MD Fine structure of Chi_H,N hat when alpha is 1}
\end{equation}

\subsection{\label{subsec:6.2.2 Multi-Dimensional--=00003D000026}Multi-Dimensional
$\hat{\chi}_{H}$ and $\tilde{\chi}_{H,N}$ for $\alpha_{H}\left(\mathbf{0}\right)=\mathbf{I}_{d}$}

We begin by computing the Fourier series generated by $v_{p}\left(\mathbf{t}\right)\hat{A}_{H}\left(\mathbf{t}\right)$. 
\begin{lem}[\textbf{$v_{p}\hat{A}_{H}$ Summation Formulae}]
\label{lem:MD v_p A_H hat summation formulae}\ 

\vphantom{}

I.

\begin{align}
\sum_{0<\left\Vert \mathbf{t}\right\Vert _{p}\leq p^{N}}v_{p}\left(\mathbf{t}\right)\hat{A}_{H}\left(\mathbf{t}\right)e^{2\pi i\left\{ \mathbf{t}\mathbf{z}\right\} _{p}}\overset{\overline{\mathbb{Q}}^{d,d}}{=} & -N\kappa_{H}\left(\left[\mathbf{z}\right]_{p^{N}}\right)\left(H^{\prime}\left(\mathbf{0}\right)\right)^{N}\label{eq:MD Fourier sum of A_H hat v_rho}\\
 & +\sum_{n=0}^{N-1}\kappa_{H}\left(\left[\mathbf{z}\right]_{p^{n}}\right)\left(H^{\prime}\left(\mathbf{0}\right)\right)^{n}\nonumber \\
 & -\sum_{n=0}^{N-1}\left(n+1\right)\kappa_{H}\left(\left[\mathbf{z}\right]_{p^{n}}\right)\mathbf{I}_{H}\left(\lambda_{p}\left(\left[\mathbf{z}\right]_{p^{n}}\right)\right)\left(H^{\prime}\left(\mathbf{0}\right)\right)^{n}\nonumber 
\end{align}
Additionally, if $H$ is commutative: 
\begin{align}
\sum_{0<\left\Vert \mathbf{t}\right\Vert _{p}\leq p^{N}}v_{p}\left(\mathbf{t}\right)\hat{A}_{H}\left(\mathbf{t}\right)e^{2\pi i\left\{ \mathbf{t}\mathbf{z}\right\} _{p}}\overset{\overline{\mathbb{Q}}^{d,d}}{=} & -N\kappa_{H}\left(\left[\mathbf{z}\right]_{p^{N}}\right)\left(H^{\prime}\left(\mathbf{0}\right)\right)^{N}\label{eq:MD Fourier sum of A_H hat v_rho, commutative}\\
 & +\sum_{n=0}^{N-1}\kappa_{H}\left(\left[\mathbf{z}\right]_{p^{n}}\right)\left(H^{\prime}\left(\mathbf{0}\right)\right)^{n}\nonumber \\
 & -\sum_{n=0}^{N-1}\left(n+1\right)\kappa_{H}\left(\left[\mathbf{z}\right]_{p^{n}}\right)\left(H^{\prime}\left(\mathbf{0}\right)\right)^{n}\left(\mathbf{I}_{d}-\alpha_{H}\left(\mathbf{0}\right)\right)\nonumber 
\end{align}

\vphantom{}

II. 
\begin{align}
\sum_{\mathbf{t}\in\hat{\mathbb{Z}}_{p}^{r}\backslash\left\{ \mathbf{0}\right\} }v_{p}\left(\mathbf{t}\right)\hat{A}_{H}\left(\mathbf{t}\right)e^{2\pi i\left\{ \mathbf{t}\mathbf{z}\right\} _{p}}\overset{\mathcal{F}_{p,q}^{d,d}}{=} & \sum_{n=0}^{\infty}\kappa_{H}\left(\left[\mathbf{z}\right]_{p^{n}}\right)\left(H^{\prime}\left(\mathbf{0}\right)\right)^{n}\label{eq:MD Limit of Fourier sum of v_rho A_H hat}\\
 & -\sum_{n=0}^{\infty}\left(n+1\right)\kappa_{H}\left(\left[\mathbf{z}\right]_{p^{n}}\right)\mathbf{I}_{H}\left(\left[\mathbf{z}\right]_{p^{n}}\right)\left(H^{\prime}\left(\mathbf{0}\right)\right)^{n}\nonumber 
\end{align}
for all $\mathbf{z}\in\mathbb{Z}_{p}^{r}$, where the convergence
is point-wise. The right-hand side can also be written as: 
\begin{equation}
\sum_{n=0}^{\infty}\kappa_{H}\left(\left[\mathbf{z}\right]_{p^{n}}\right)\left(\left(n+1\right)\mathcal{C}_{H}\left(\alpha_{H}\left(\mathbf{0}\right):\lambda_{p}\left(\left[\mathbf{z}\right]_{p^{n}}\right)\right)-n\mathbf{I}_{d}\right)\left(H^{\prime}\left(\mathbf{0}\right)\right)^{n}\label{eq:MD Limit of Fourier sum of v_rho A_H hat, Alt}
\end{equation}
When $H$ is commutative, we get: 
\begin{equation}
\sum_{\mathbf{t}\in\hat{\mathbb{Z}}_{p}^{r}\backslash\left\{ \mathbf{0}\right\} }v_{p}\left(\mathbf{t}\right)\hat{A}_{H}\left(\mathbf{t}\right)e^{2\pi i\left\{ \mathbf{t}\mathbf{z}\right\} _{p}}\overset{\mathcal{F}_{p,q}^{d,d}}{=}\sum_{n=0}^{\infty}\kappa_{H}\left(\left[\mathbf{z}\right]_{p^{n}}\right)\left(\left(n+1\right)\alpha_{H}\left(\mathbf{0}\right)-n\mathbf{I}_{d}\right)\left(H^{\prime}\left(\mathbf{0}\right)\right)^{n}\label{eq:MD Limit of Fourier sum of v_rho A_H hat, Alt, commutative}
\end{equation}
\end{lem}
Proof: For (I), using (\ref{eq:MD Convolution of dA_H and D_N}) gives:
\begin{align*}
\sum_{0<p\left\Vert \mathbf{t}\right\Vert \leq p^{N}}v_{p}\left(\mathbf{t}\right)\hat{A}_{H}\left(\mathbf{t}\right)e^{2\pi i\left\{ \mathbf{t}\mathbf{z}\right\} _{p}} & =\sum_{n=0}^{N-1}\kappa_{H}\left(\left[\mathbf{z}\right]_{p^{n}}\right)\left(H^{\prime}\left(\mathbf{0}\right)\right)^{n}\\
 & -N\kappa_{H}\left(\left[\mathbf{z}\right]_{p^{N}}\right)\left(H^{\prime}\left(\mathbf{0}\right)\right)^{N}\\
 & -N\sum_{m=0}^{N-1}\kappa_{H}\left(\left[\mathbf{z}\right]_{p^{m}}\right)\mathbf{I}_{H}\left(\lambda_{p}\left(\left[\mathbf{z}\right]_{p^{m}}\right)\right)\left(H^{\prime}\left(\mathbf{0}\right)\right)^{m}\\
 & +\sum_{n=1}^{N-1}\sum_{m=0}^{n-1}\kappa_{H}\left(\left[\mathbf{z}\right]_{p^{m}}\right)\mathbf{I}_{H}\left(\lambda_{p}\left(\left[\mathbf{z}\right]_{p^{m}}\right)\right)\left(H^{\prime}\left(\mathbf{0}\right)\right)^{m}
\end{align*}
where we used the fact that $\tilde{A}_{H,0}\left(\mathbf{z}\right)=\hat{A}_{H}\left(\mathbf{0}\right)=\mathbf{I}_{d}$.
Next, we note the formal identity (provable by summation by parts):
\begin{equation}
\sum_{n=1}^{N-1}\sum_{m=0}^{n-1}f\left(m\right)=\sum_{m=0}^{N-2}\sum_{n=m+1}^{N-1}f\left(m\right)=\sum_{m=0}^{N-2}\left(N-1-m\right)f\left(m\right)
\end{equation}
Applying this to our sum yields:

\begin{align*}
\sum_{0<\left\Vert \mathbf{t}\right\Vert _{p}\leq p^{N}}v_{p}\left(\mathbf{t}\right)\hat{A}_{H}\left(\mathbf{t}\right)e^{2\pi i\left\{ \mathbf{t}\mathbf{z}\right\} _{p}} & =\sum_{n=0}^{N-1}\kappa_{H}\left(\left[\mathbf{z}\right]_{p^{n}}\right)\left(H^{\prime}\left(\mathbf{0}\right)\right)^{n}\\
 & -N\kappa_{H}\left(\left[\mathbf{z}\right]_{p^{N}}\right)\left(H^{\prime}\left(\mathbf{0}\right)\right)^{N}\\
 & -N\sum_{m=0}^{N-1}\kappa_{H}\left(\left[\mathbf{z}\right]_{p^{m}}\right)\mathbf{I}_{H}\left(\lambda_{p}\left(\left[\mathbf{z}\right]_{p^{m}}\right)\right)\left(H^{\prime}\left(\mathbf{0}\right)\right)^{m}\\
 & +\sum_{m=0}^{N-2}\left(N-1-m\right)\kappa_{H}\left(\left[\mathbf{z}\right]_{p^{m}}\right)\mathbf{I}_{H}\left(\lambda_{p}\left(\left[\mathbf{z}\right]_{p^{m}}\right)\right)\left(H^{\prime}\left(\mathbf{0}\right)\right)^{m}
\end{align*}
and hence: 
\begin{align*}
\sum_{0<\left\Vert \mathbf{t}\right\Vert _{p}\leq p^{N}}v_{p}\left(\mathbf{t}\right)\hat{A}_{H}\left(\mathbf{t}\right)e^{2\pi i\left\{ \mathbf{t}\mathbf{z}\right\} _{p}} & =-N\kappa_{H}\left(\left[\mathbf{z}\right]_{p^{N}}\right)\left(H^{\prime}\left(\mathbf{0}\right)\right)^{N}\\
 & +\sum_{n=0}^{N-1}\kappa_{H}\left(\left[\mathbf{z}\right]_{p^{n}}\right)\left(H^{\prime}\left(\mathbf{0}\right)\right)^{n}\\
 & -\sum_{n=0}^{N-1}\left(n+1\right)\kappa_{H}\left(\left[\mathbf{z}\right]_{p^{n}}\right)\mathbf{I}_{H}\left(\lambda_{p}\left(\left[\mathbf{z}\right]_{p^{n}}\right)\right)\left(H^{\prime}\left(\mathbf{0}\right)\right)^{n}
\end{align*}
The $\mathcal{F}_{p,q}^{d,d}$-convergence of this sum as $N\rightarrow\infty$
in the case where $H$ is semi-basic and contracting follows from
the decay properties ($q$-adic and archimedean, respectively) of
$\kappa_{H}$ and $\left(H^{\prime}\left(\mathbf{0}\right)\right)^{n}$,
as well as the boundedness of $\mathbf{I}_{H}\left(\lambda_{p}\left(\left[\mathbf{z}\right]_{p^{n}}\right)\right)$
in both topologies with respect to $n$. Finally, writing: 
\begin{equation}
\mathbf{I}_{H}\left(\lambda_{p}\left(\left[\mathbf{z}\right]_{p^{n}}\right)\right)=\mathbf{I}_{d}-\mathcal{C}_{H}\left(\alpha_{H}\left(\mathbf{0}\right):\lambda_{p}\left(\left[\mathbf{z}\right]_{p^{n}}\right)\right)
\end{equation}
yields (\ref{eq:MD Limit of Fourier sum of v_rho A_H hat, Alt}).

Q.E.D.

\vphantom{}

Like in the one-dimensional case, we now introduce $\varepsilon_{n}$. 
\begin{defn}[\textbf{Multi-Dimensional $\varepsilon_{n}$}]
For each $n\in\mathbb{N}_{0}$, \nomenclature{$\varepsilon_{n}\left(\mathbf{z}\right)$}{ }we
define $\varepsilon_{n}:\mathbb{Z}_{p}^{r}\rightarrow\overline{\mathbb{Q}}$
by: 
\begin{equation}
\varepsilon_{n}\left(\mathbf{z}\right)\overset{\textrm{def}}{=}e^{\frac{2\pi i}{p^{n+1}}\left(\left[\mathbf{z}\right]_{p^{n+1}}-\left[\mathbf{z}\right]_{p^{n}}\right)}=e^{2\pi i\left\{ \frac{\mathbf{z}}{p^{n+1}}\right\} _{p}}\cdot e^{-\frac{2\pi i}{p}\left\{ \frac{\mathbf{z}}{p^{n}}\right\} _{p}}\label{eq:MD Definition of epsilon_n}
\end{equation}
\end{defn}
Once again, we have functional equations, among other things. 
\begin{prop}[\textbf{Properties of Multi-Dimensional $\varepsilon_{n}$}]
\label{prop:properties of MD epsilon n}\ 

\vphantom{}

I.

\begin{equation}
\varepsilon_{0}\left(\mathbf{z}\right)=e^{2\pi i\left\{ \frac{\mathbf{z}}{p}\right\} _{p}},\textrm{ }\forall\mathbf{z}\in\mathbb{Z}_{p}^{r}\label{eq:MD Epsilon 0 of z}
\end{equation}

\begin{equation}
\varepsilon_{n}\left(\mathbf{j}\right)=1,\textrm{ }\forall\mathbf{j}\in\mathbb{Z}^{r}/p\mathbb{Z}^{r},\textrm{ }\forall n\geq1\label{eq:MD Epsilon_n of j}
\end{equation}

\vphantom{}

II. 
\begin{equation}
\varepsilon_{n}\left(p\mathbf{m}+\mathbf{j}\right)=\begin{cases}
\varepsilon_{0}\left(\mathbf{j}\right) & \textrm{if }n=0\\
\varepsilon_{n-1}\left(\mathbf{m}\right) & \textrm{if }n\geq1
\end{cases},\textrm{ }\forall\mathbf{m}\in\mathbb{N}_{0}^{r},\textrm{ }\forall\mathbf{j}\in\mathbb{Z}^{r}/p\mathbb{Z}^{r},\textrm{ }\forall n\geq1\label{eq:MD epsilon_n functional equations}
\end{equation}

\vphantom{}

III. Let $\mathbf{z}\neq\mathbf{0}$. Then $\varepsilon_{n}\left(\mathbf{z}\right)=1\textrm{ }$
for all $n<v_{p}\left(\mathbf{z}\right)$.
\end{prop}
Proof:

I. The identity (\ref{eq:MD Epsilon 0 of z}) is immediate from the
definition of $\varepsilon_{n}$. As for the other one, note that
$\left[\mathbf{j}\right]_{p^{n}}=\mathbf{j}$ for all $n\geq1$ and
all $\mathbf{j}\in\mathbb{Z}^{r}/p\mathbb{Z}^{r}$. So, letting $n\geq1$,
we get:
\begin{equation}
\varepsilon_{n}\left(\mathbf{j}\right)=e^{\frac{2\pi i}{p^{n+1}}\left(\left[\mathbf{j}\right]_{p^{n+1}}-\left[\mathbf{j}\right]_{p^{n}}\right)}=e^{\frac{2\pi i}{p^{n+1}}\cdot\mathbf{0}}=1
\end{equation}
as desired.

\vphantom{}

II. 
\begin{align*}
\varepsilon_{n}\left(p\mathbf{m}+\mathbf{j}\right) & =e^{2\pi i\left\{ \frac{p\mathbf{m}+\mathbf{j}}{p^{n+1}}\right\} _{p}}e^{-\frac{2\pi i}{p}\left\{ \frac{p\mathbf{m}+\mathbf{j}}{p^{n}}\right\} _{p}}\\
 & =e^{2\pi i\left\{ \frac{\mathbf{j}}{p^{n+1}}\right\} _{p}}e^{-\frac{2\pi i}{p}\left\{ \frac{\mathbf{j}}{p^{n}}\right\} _{p}}\cdot e^{2\pi i\left\{ \frac{\mathbf{m}}{p^{n}}\right\} _{p}}e^{-\frac{2\pi i}{p}\left\{ \frac{\mathbf{m}}{p^{n-1}}\right\} _{p}}\\
 & =\varepsilon_{n}\left(\mathbf{j}\right)\varepsilon_{n-1}\left(\mathbf{m}\right)\\
\left(\textrm{by (I)}\right); & =\begin{cases}
\varepsilon_{0}\left(\mathbf{j}\right) & \textrm{if }n=0\\
\varepsilon_{n-1}\left(\mathbf{m}\right) & \textrm{if }n\geq1
\end{cases}
\end{align*}

\vphantom{}

III. Let $\mathbf{z}$ be non-zero. When $n<v_{p}\left(\mathbf{z}\right)$,
we have that $p^{-n}\mathbf{z}$ and $p^{-\left(n+1\right)}\mathbf{z}$
are then $p$-adic integer tuples. As such:
\begin{equation}
\varepsilon_{n}\left(\mathbf{z}\right)=e^{2\pi i\left\{ \frac{\mathbf{z}}{p^{n+1}}\right\} _{p}}e^{-\frac{2\pi i}{p}\left\{ \frac{\mathbf{z}}{p^{n}}\right\} _{p}}=e^{0}\cdot e^{-0}=1
\end{equation}

Q.E.D.

\vphantom{}Now we compute the sum of the Fourier series generated
by $\hat{A}_{H}\left(\mathbf{t}\right)\gamma_{H}\left(\frac{\mathbf{t}\left|\mathbf{t}\right|_{p}}{p}\right)$.
Once again, the non-singularity of $H$ will be necessary for this
to work out, as will $H$'s semi-basicness and contracting-ness.
\begin{lem}[\textbf{Multi-Dimensional $\gamma_{H}\hat{A}_{H}$ Summation Formulae}]
\label{lem:MD gamma formulae}\ 

\vphantom{}

I. The column vector: 
\begin{equation}
\sum_{0<\left\Vert \mathbf{t}\right\Vert _{p}\leq p^{N}}\hat{A}_{H}\left(\mathbf{t}\right)\gamma_{H}\left(\frac{\mathbf{t}\left|\mathbf{t}\right|_{p}}{p}\right)e^{2\pi i\left\{ \mathbf{t}\mathbf{z}\right\} _{p}}\in\overline{\mathbb{Q}}^{d}\label{eq:Left Hand Side of MD Gamma Formulae (partial Fourier sum)}
\end{equation}
is given by:

\begin{equation}
\sum_{n=0}^{N-1}\kappa_{H}\left(\left[\mathbf{z}\right]_{p^{n}}\right)\sum_{\mathbf{j}>\mathbf{0}}^{p-1}\mathcal{C}_{H}\left(\alpha_{H}\left(\frac{\mathbf{j}}{p}\right)\varepsilon_{n}\left(\mathbf{j}\mathbf{z}\right):\lambda_{p}\left(\left[\mathbf{z}\right]_{p^{n}}\right)\right)\left(H^{\prime}\left(\mathbf{0}\right)\right)^{n}\gamma_{H}\left(\frac{\mathbf{j}}{p}\right)\label{eq:MD Gamma formula}
\end{equation}

\vphantom{}

II. As $N\rightarrow\infty$, \emph{(\ref{eq:Left Hand Side of MD Gamma Formulae (partial Fourier sum)})}
is $\mathcal{F}_{p,q}^{d}$ convergent to: 
\begin{equation}
\sum_{n=0}^{\infty}\kappa_{H}\left(\left[\mathbf{z}\right]_{p^{n}}\right)\sum_{\mathbf{j}>\mathbf{0}}^{p-1}\mathcal{C}_{H}\left(\alpha_{H}\left(\frac{\mathbf{j}}{p}\right)\varepsilon_{n}\left(\mathbf{j}\mathbf{z}\right):\lambda_{p}\left(\left[\mathbf{z}\right]_{p^{n}}\right)\right)\left(H^{\prime}\left(\mathbf{0}\right)\right)^{n}\gamma_{H}\left(\frac{\mathbf{j}}{p}\right)\label{eq:F limit of MD Gamma_H A_H hat Fourier series}
\end{equation}
since $H$ is semi-basic and contracting. 
\end{lem}
Proof: Throughout this proof, we write $\mathbf{X}$ to denote $H^{\prime}\left(\mathbf{0}\right)$.

\vphantom{}

I. Like in the one-dimensional case, we first note that the map: 
\begin{equation}
\mathbf{t}\in\hat{\mathbb{Z}}_{p}^{r}\mapsto\frac{\mathbf{t}\left|\mathbf{t}\right|_{p}}{p}=\left(\frac{t_{1}\left|t_{1}\right|_{p}}{p},\ldots,\frac{t_{r}\left|t_{r}\right|_{p}}{p}\right)\in\hat{\mathbb{Z}}_{p}^{r}
\end{equation}
takes tuples: 
\begin{equation}
\mathbf{t}=\left(\frac{k_{1}}{p^{n_{1}}},\ldots,\frac{k_{r}}{p^{n_{r}}}\right)
\end{equation}
and outputs: 
\begin{equation}
\frac{\mathbf{t}\left|\mathbf{t}\right|_{p}}{p}=\left(\frac{\left[k_{1}\right]_{p}}{p},\ldots,\frac{\left[k_{r}\right]_{p}}{p}\right)
\end{equation}
Now, for brevity, let: 
\begin{align}
\gamma_{\mathbf{j}} & \overset{\textrm{def}}{=}\gamma_{H}\left(\frac{\mathbf{j}}{p}\right)=\gamma_{H}\left(\frac{j_{1}}{p},\ldots,\frac{j_{r}}{p}\right),\textrm{ }\forall\mathbf{j}\in\mathbb{Z}^{r}/p\mathbb{Z}^{r}\\
F_{N}\left(\mathbf{z}\right) & \overset{\textrm{def}}{=}\sum_{0<\left\Vert \mathbf{t}\right\Vert _{p}\leq p^{N}}\hat{A}_{H}\left(\mathbf{t}\right)\gamma_{H}\left(\frac{\mathbf{t}\left|\mathbf{t}\right|_{p}}{p}\right)e^{2\pi i\left\{ \mathbf{t}\mathbf{z}\right\} _{p}}
\end{align}
Finally, for each $\mathbf{k}\leq p^{n-1}-1$, each $\mathbf{j}\in\mathbb{Z}^{r}/p\mathbb{Z}^{r}$,
and each $n\geq1$, we write: 
\begin{equation}
\frac{p\mathbf{k}+\mathbf{j}}{p^{n}}=\left(\frac{pk_{1}+j_{1}}{p^{n}},\ldots,\frac{pk_{r}+j_{r}}{p^{n}}\right)
\end{equation}

\begin{claim}
In the above notation, we have: 
\begin{equation}
\left\{ \mathbf{t}\in\hat{\mathbb{Z}}_{p}^{r}:\left\Vert \mathbf{t}\right\Vert _{p}=p^{n}\right\} =\left\{ \frac{p\mathbf{k}+\mathbf{j}}{p^{n}}:\mathbf{k}\leq p^{n-1}-1,\textrm{ }\mathbf{j}\in\left(\mathbb{Z}^{r}/p\mathbb{Z}^{r}\right)\backslash\left\{ \mathbf{0}\right\} \right\} \label{eq:MD Decomposition of level sets}
\end{equation}

Proof of claim:

\vphantom{}

I. Let $\mathbf{t}=\frac{p\mathbf{k}+\mathbf{j}}{p^{n}}$. Then, for
each $m\in\left\{ 1,\ldots,r\right\} $, $t_{m}=\frac{pk_{m}+j_{m}}{p^{n}}$
where $0\leq k_{m}\leq p^{n-1}-1$ and $0\leq j_{m}\leq p-1$; moreover,
there exists an $\ell\in\left\{ 1,\ldots,r\right\} $ so that $j_{\ell}>0$.
Since $pk_{\ell}+j_{\ell}$ is then co-prime to $p_{\ell}$, it follows
that $\left|t_{\ell}\right|_{p}=p^{n}$. Since $\left|t_{m}\right|_{p}\leq p^{n}$
for all $m$, this shows that $\left\Vert \mathbf{t}\right\Vert _{p}=p^{n}$.

\vphantom{}

II. Let $\left\Vert \mathbf{t}\right\Vert _{p}=p^{n}$. Then, for
every $m\in\left\{ 1,\ldots,r\right\} $, we can write $t_{m}=\frac{\nu_{m}}{p^{n}}$,
where $\nu_{m}\in\left\{ 0,\ldots,p^{n}-1\right\} $. Moreover, there
is an $\ell\in\left\{ 1,\ldots,r\right\} $ so that $\nu_{\ell}$
is co-prime to $p$, so as to guarantee that $\left|t_{\ell}\right|_{p}=p^{n}$
Each $\nu_{m}$ can be written as $pk_{m}-j_{m}$, where we have chosen
$j_{m}=\left[\nu_{m}\right]_{p}\leq p-1$ and $k_{m}=\left(\nu_{m}-j_{m}\right)/p$.
Since $\nu_{m}\leq p^{n}-1$, this forces $k_{m}\leq p^{n-1}-1$.
Consequently, for these $j_{m}$s, the tuple $\mathbf{j}=\left(j_{1},\ldots,j_{r}\right)$
is then an element of $\mathbb{Z}^{r}/p\mathbb{Z}^{r}$, and $\mathbf{k}$
is a tuple $\leq p^{n-1}-1$ so that $\mathbf{t}=\left(p\mathbf{k}+\mathbf{j}\right)/p^{n}$.
Finally, since $\nu_{\ell}$ is co-prime to $p$, $j_{\ell}\neq0$,
which shows that $\mathbf{j}\in\left(\mathbb{Z}^{r}/p\mathbb{Z}^{r}\right)\backslash\left\{ \mathbf{0}\right\} $.

\vphantom{}

This proves the claim.

\vphantom{} 
\end{claim}
Consequently, we can express $F_{N}$ as a sum involving the $\gamma_{\mathbf{j}}$s
like so: 
\begin{align*}
F_{N}\left(\mathbf{z}\right) & =\sum_{n=1}^{N}\sum_{\left\Vert \mathbf{t}\right\Vert _{p}=p^{n}}\hat{A}_{H}\left(\mathbf{t}\right)\gamma_{H}\left(\frac{\mathbf{t}\left|\mathbf{t}\right|_{p}}{p}\right)e^{2\pi i\left\{ \mathbf{t}\mathbf{z}\right\} _{p}}\\
\left(\textrm{use }(\ref{eq:MD Decomposition of level sets})\right); & =\sum_{n=1}^{N}\sum_{\mathbf{j}>\mathbf{0}}^{p-1}\sum_{\mathbf{k}=\mathbf{0}}^{p^{n-1}-1}\hat{A}_{H}\left(\frac{p\mathbf{k}+\mathbf{j}}{p^{n}}\right)\gamma_{H}\left(\frac{\mathbf{j}}{p}\right)e^{2\pi i\left\{ \frac{p\mathbf{k}+\mathbf{j}}{p^{n}}\mathbf{z}\right\} _{p}}
\end{align*}
Using the formal identity: 
\begin{equation}
\sum_{\mathbf{k}=\mathbf{0}}^{p^{n-1}-1}f\left(\frac{p\mathbf{k}+\mathbf{j}}{p^{n}}\right)=\sum_{\left\Vert \mathbf{t}\right\Vert _{p}\leq p^{n-1}}f\left(\mathbf{t}+\frac{\mathbf{j}}{p^{n}}\right)
\end{equation}
we can then write:

\begin{align}
F_{N}\left(\mathbf{z}\right) & =\sum_{n=1}^{N}\sum_{\left\Vert \mathbf{t}\right\Vert _{p}\leq p^{n-1}}\sum_{\mathbf{j}>\mathbf{0}}^{p-1}\hat{A}_{H}\left(\mathbf{t}+\frac{\mathbf{j}}{p^{n}}\right)\gamma_{\mathbf{j}}e^{2\pi i\left\{ \frac{\mathbf{j}\mathbf{z}}{p^{n}}\right\} _{p}}e^{2\pi i\left\{ \mathbf{t}\mathbf{z}\right\} _{p}}\label{eq:MD 1/3rd of the way through the gamma computation}
\end{align}

To deal with the $\mathbf{j}$-sum, we express $\hat{A}_{H}$ in product
form, changing: 
\[
\sum_{\mathbf{j}>\mathbf{0}}^{p-1}\hat{A}_{H}\left(\mathbf{t}+\frac{\mathbf{j}}{p^{n}}\right)\gamma_{\mathbf{j}}e^{2\pi i\left\{ \frac{\mathbf{j}\mathbf{z}}{p^{n}}\right\} _{p}}e^{2\pi i\left\{ \mathbf{t}\mathbf{z}\right\} _{p}}
\]
into:

\[
\sum_{\mathbf{j}>\mathbf{0}}^{p-1}\left(\prod_{\mathbf{m}=\mathbf{0}}^{n-1}\alpha_{H}\left(p^{m}\left(\mathbf{t}+\frac{\mathbf{j}}{p^{n}}\right)\right)\right)\gamma_{\mathbf{j}}e^{2\pi i\left\{ \frac{\mathbf{j}\mathbf{z}}{p^{n}}\right\} _{p}}e^{2\pi i\left\{ \mathbf{t}\mathbf{z}\right\} _{p}}
\]
Using \textbf{Proposition \ref{prop:MD alpha H series}} to write
the $\alpha_{H}$-product out as a series gives: 
\begin{equation}
\sum_{\mathbf{j}>\mathbf{0}}^{p-1}\left(\sum_{\mathbf{m}=\mathbf{0}}^{p^{n}-1}\kappa_{H}\left(\mathbf{m}\right)\left(\frac{\mathbf{X}}{p^{r}}\right)^{n}e^{-2\pi i\mathbf{m}\cdot\left(\mathbf{t}+\frac{\mathbf{j}}{p^{n}}\right)}\right)\gamma_{\mathbf{j}}e^{2\pi i\left\{ \frac{\mathbf{j}\mathbf{z}}{p^{n}}\right\} _{p}}e^{2\pi i\left\{ \mathbf{t}\mathbf{z}\right\} _{p}}
\end{equation}
and hence:

\begin{equation}
\sum_{\mathbf{j}>\mathbf{0}}^{p-1}\sum_{\mathbf{m}=\mathbf{0}}^{p^{n}-1}\kappa_{H}\left(\mathbf{m}\right)\left(\frac{\mathbf{X}}{p^{r}}\right)^{n}\gamma_{\mathbf{j}}e^{2\pi i\left\{ \frac{\mathbf{j}\left(\mathbf{z}-\mathbf{m}\right)}{p^{n}}\right\} _{p}}e^{2\pi i\left\{ \mathbf{t}\left(\mathbf{z}-\mathbf{m}\right)\right\} _{p}}
\end{equation}
Summing over $\left\Vert \mathbf{t}\right\Vert _{p}\leq p^{n-1}$
and using the Fourier series for $\left[\mathbf{z}\overset{p^{n-1}}{\equiv}\mathbf{m}\right]$:
\begin{equation}
\sum_{\left\Vert \mathbf{t}\right\Vert _{p}\leq p^{n-1}}e^{2\pi i\left\{ \mathbf{t}\left(\mathbf{z}-\mathbf{m}\right)\right\} _{p}}=p^{r\left(n-1\right)}\left[\mathbf{z}\overset{p^{n-1}}{\equiv}\mathbf{m}\right]
\end{equation}
we obtain: 
\begin{equation}
\sum_{\mathbf{j}>\mathbf{0}}^{p-1}\sum_{\mathbf{m}=\mathbf{0}}^{p^{n}-1}\kappa_{H}\left(\mathbf{m}\right)\left(\frac{\mathbf{X}}{p^{r}}\right)^{n}\gamma_{\mathbf{j}}e^{2\pi i\left\{ \frac{\mathbf{j}\left(\mathbf{z}-\mathbf{m}\right)}{p^{n}}\right\} _{p}}p^{r\left(n-1\right)}\left[\mathbf{z}\overset{p^{n-1}}{\equiv}\mathbf{m}\right]
\end{equation}

In summary: 
\begin{eqnarray}
 & \sum_{\left\Vert \mathbf{t}\right\Vert _{p}\leq p^{n-1}}\sum_{\mathbf{j}>\mathbf{0}}^{p-1}\hat{A}_{H}\left(\mathbf{t}+\frac{\mathbf{j}}{p^{n}}\right)\gamma_{\mathbf{j}}e^{2\pi i\left\{ \frac{\mathbf{j}\mathbf{z}}{p^{n}}\right\} _{p}}e^{2\pi i\left\{ \mathbf{t}\mathbf{z}\right\} _{p}}\nonumber \\
 & =\label{eq:MD 2/3rds of the way through the gamma computation}\\
 & \frac{1}{p^{r}}\sum_{\mathbf{j}>\mathbf{0}}^{p-1}\left(\sum_{\mathbf{m}=\mathbf{0}}^{p^{n}-1}\kappa_{H}\left(\mathbf{m}\right)e^{2\pi i\left\{ \frac{\mathbf{j}\left(\mathbf{z}-\mathbf{m}\right)}{p^{n}}\right\} _{p}}\left[\mathbf{z}\overset{p^{n-1}}{\equiv}\mathbf{m}\right]\right)\mathbf{X}^{n}\gamma_{\mathbf{j}}\nonumber 
\end{eqnarray}

Next, using the formal summation identity: 
\begin{equation}
\sum_{\mathbf{m}=\mathbf{0}}^{p^{n}-1}f\left(\mathbf{m}\right)=\sum_{\mathbf{k}=\mathbf{0}}^{p-1}\sum_{\mathbf{m}=\mathbf{0}}^{p^{n-1}-1}f\left(\mathbf{m}+p^{n-1}\mathbf{k}\right)\label{eq:MD rho to the n formal identity}
\end{equation}
the expression:
\begin{equation}
\sum_{\mathbf{m}=\mathbf{0}}^{p^{n}-1}\kappa_{H}\left(\mathbf{m}\right)e^{2\pi i\left\{ \frac{\mathbf{j}\left(\mathbf{z}-\mathbf{m}\right)}{p^{n}}\right\} _{p}}\left[\mathbf{z}\overset{p^{n-1}}{\equiv}\mathbf{m}\right]
\end{equation}
becomes: 
\begin{equation}
\sum_{\mathbf{k}=\mathbf{0}}^{p-1}\sum_{\mathbf{m}=\mathbf{0}}^{p^{n-1}-1}\kappa_{H}\left(\mathbf{m}+p^{n-1}\mathbf{k}\right)e^{-2\pi i\frac{\mathbf{j}\cdot\mathbf{k}}{p}}e^{2\pi i\left\{ \frac{\mathbf{j}\left(\mathbf{z}-\mathbf{m}\right)}{p^{n}}\right\} _{p}}\left[\mathbf{z}\overset{p^{n-1}}{\equiv}\mathbf{m}\right]
\end{equation}
Applying the functional equation identity for $\kappa_{H}$ (equation
(\ref{eq:MD Kappa_H has P-adic structure}) from \textbf{Lemma \ref{lem:properties of MD kappa_H}})
yields: 
\begin{equation}
\sum_{\mathbf{k}=\mathbf{0}}^{p-1}\sum_{\mathbf{m}=\mathbf{0}}^{p^{n-1}-1}\kappa_{H}\left(\mathbf{m}\right)\mathcal{C}_{H}\left(\kappa_{H}\left(\mathbf{k}\right):\lambda_{p}\left(\mathbf{m}\right)\right)e^{-2\pi i\frac{\mathbf{j}\cdot\mathbf{k}}{p}}e^{2\pi i\left\{ \frac{\mathbf{j}\left(\mathbf{z}-\mathbf{m}\right)}{p^{n}}\right\} _{p}}\left[\mathbf{z}\overset{p^{n-1}}{\equiv}\mathbf{m}\right]
\end{equation}
Here, note that $\left[\mathbf{z}\right]_{p^{n-1}}$ is the unique
integer $r$-tuple $\mathbf{m}$ so that $\mathbf{m}\leq p^{n-1}-1$
and $\mathbf{z}\overset{p^{n-1}}{\equiv}\mathbf{m}$. This leaves
us with: 
\begin{equation}
\sum_{\mathbf{k}=\mathbf{0}}^{p-1}\kappa_{H}\left(\left[\mathbf{z}\right]_{p^{n-1}}\right)\mathcal{C}_{H}\left(\kappa_{H}\left(\mathbf{k}\right):\lambda_{p}\left(\left[\mathbf{z}\right]_{p^{n-1}}\right)\right)e^{-2\pi i\frac{\mathbf{j}\cdot\mathbf{k}}{p}}e^{2\pi i\left\{ \frac{\mathbf{j}\left(\mathbf{z}-\left[\mathbf{z}\right]_{p^{n-1}}\right)}{p^{n}}\right\} _{p}}
\end{equation}
which is: 
\begin{equation}
\sum_{\mathbf{k}=\mathbf{0}}^{p-1}\kappa_{H}\left(\left[\mathbf{z}\right]_{p^{n-1}}\right)\mathcal{C}_{H}\left(\kappa_{H}\left(\mathbf{k}\right):\lambda_{p}\left(\left[\mathbf{z}\right]_{p^{n-1}}\right)\right)e^{-2\pi i\frac{\mathbf{j}\cdot\mathbf{k}}{p}}\varepsilon_{n-1}\left(\mathbf{j}\mathbf{z}\right)\label{eq:gamma proof - the above}
\end{equation}
Next, observe that for fixed $\left[\mathbf{z}\right]_{p^{n-1}}$,
the map $\mathbf{A}\mapsto\mathcal{C}_{H}\left(\mathbf{A}:\lambda_{p}\left(\left[\mathbf{z}\right]_{p^{n-1}}\right)\right)$
is linear. So, (\ref{eq:gamma proof - the above}) can be written
as: 
\begin{equation}
\kappa_{H}\left(\left[\mathbf{z}\right]_{p^{n-1}}\right)\mathcal{C}_{H}\left(\sum_{\mathbf{k}=\mathbf{0}}^{p-1}\kappa_{H}\left(\mathbf{k}\right)e^{-2\pi i\frac{\mathbf{j}\cdot\mathbf{k}}{p}}:\lambda_{p}\left(\left[\mathbf{z}\right]_{p^{n-1}}\right)\right)\varepsilon_{n-1}\left(\mathbf{j}\mathbf{z}\right)
\end{equation}
With this, (\ref{eq:MD 2/3rds of the way through the gamma computation})
becomes:
\begin{eqnarray}
 & \sum_{\left\Vert \mathbf{t}\right\Vert _{p}\leq p^{n-1}}\sum_{\mathbf{j}=\mathbf{0}}^{p-1}\hat{A}_{H}\left(\mathbf{t}+\frac{\mathbf{j}}{p^{n}}\right)\gamma_{\mathbf{j}}e^{2\pi i\left\{ \frac{\mathbf{j}\mathbf{z}}{p^{n}}\right\} _{p}}e^{2\pi i\left\{ \mathbf{t}\mathbf{z}\right\} _{p}}\nonumber \\
 & =\label{eq:MD 2/3rds and 1/4th of the way through the gamma computation}\\
 & \frac{1}{p^{r}}\sum_{\mathbf{j}>\mathbf{0}}^{p-1}\kappa_{H}\left(\left[\mathbf{z}\right]_{p^{n-1}}\right)\mathcal{C}_{H}\left(\sum_{\mathbf{k}=\mathbf{0}}^{p-1}\kappa_{H}\left(\mathbf{k}\right)e^{-2\pi i\frac{\mathbf{j}\cdot\mathbf{k}}{p}}:\lambda_{p}\left(\left[\mathbf{z}\right]_{p^{n-1}}\right)\right)\varepsilon_{n-1}\left(\mathbf{j}\mathbf{z}\right)\mathbf{X}^{n}\gamma_{\mathbf{j}}\nonumber 
\end{eqnarray}

Our next step is to simplify the expression with $\mathcal{C}_{H}$.
Here, we use $\kappa_{H}$'s functional equation ((\ref{eq:MD Kappa_H functional equations})
from \textbf{Lemma \ref{lem:properties of MD kappa_H}}), setting
$\mathbf{j}=\mathbf{k}$ and $\mathbf{n}=\mathbf{0}$ in it. This
gives us:
\begin{align*}
\sum_{\mathbf{k}=\mathbf{0}}^{p-1}\kappa_{H}\left(\mathbf{k}\right)e^{-2\pi i\frac{\mathbf{j}\cdot\mathbf{k}}{p}} & =\sum_{\mathbf{k}=\mathbf{0}}^{p-1}\mathbf{D}_{\mathbf{k}}^{-1}\mathbf{A}_{\mathbf{k}}\kappa_{H}\left(\mathbf{0}\right)\mathbf{A}_{\mathbf{0}}^{-1}\mathbf{D}_{\mathbf{0}}e^{-2\pi i\frac{\mathbf{j}\cdot\mathbf{k}}{p}}\\
\left(\kappa_{H}\left(\mathbf{0}\right)=\mathbf{I}_{d}\right); & =\left(\sum_{\mathbf{k}=\mathbf{0}}^{p-1}\mathbf{D}_{\mathbf{k}}^{-1}\mathbf{A}_{\mathbf{k}}e^{-2\pi i\frac{\mathbf{j}\cdot\mathbf{k}}{p}}\right)\mathbf{A}_{\mathbf{0}}^{-1}\mathbf{D}_{\mathbf{0}}\\
 & =p^{r}\alpha_{H}\left(\frac{\mathbf{j}}{p}\right)\mathbf{A}_{\mathbf{0}}^{-1}\mathbf{D}_{\mathbf{0}}\\
 & =p^{r}\alpha_{H}\left(\frac{\mathbf{j}}{p}\right)\mathbf{X}^{-1}
\end{align*}
Hence, (\ref{eq:MD 2/3rds and 1/4th of the way through the gamma computation})
becomes: 
\begin{eqnarray*}
 & \sum_{\left\Vert \mathbf{t}\right\Vert _{p}\leq p^{n-1}}\sum_{\mathbf{j}>\mathbf{0}}^{p-1}\hat{A}_{H}\left(\mathbf{t}+\frac{\mathbf{j}}{p^{n}}\right)\gamma_{\mathbf{j}}e^{2\pi i\left\{ \frac{\mathbf{j}\mathbf{z}}{p^{n}}\right\} _{p}}e^{2\pi i\left\{ \mathbf{t}\mathbf{z}\right\} _{p}}\\
 & =\\
 & \frac{1}{p^{r}}\sum_{\mathbf{j}>\mathbf{0}}^{p-1}\kappa_{H}\left(\left[\mathbf{z}\right]_{p^{n-1}}\right)\mathcal{C}_{H}\left(p\alpha_{H}\left(\frac{\mathbf{j}}{p}\right)\mathbf{X}^{-1}:\lambda_{p}\left(\left[\mathbf{z}\right]_{p^{n-1}}\right)\right)\underbrace{\varepsilon_{n-1}\left(\mathbf{j}\mathbf{z}\right)}_{\textrm{scalar}}\mathbf{X}^{n}\gamma_{\mathbf{j}}\\
 & =\\
 & \kappa_{H}\left(\left[\mathbf{z}\right]_{p^{n-1}}\right)\sum_{\mathbf{j}>\mathbf{0}}^{p-1}\varepsilon_{n-1}\left(\mathbf{j}\mathbf{z}\right)\mathcal{C}_{H}\left(\alpha_{H}\left(\frac{\mathbf{j}}{p}\right)\mathbf{X}^{-1}:\lambda_{p}\left(\left[\mathbf{z}\right]_{p^{n-1}}\right)\right)\mathbf{X}^{n}\gamma_{\mathbf{j}}
\end{eqnarray*}

Next, we note the following formal identity for matrices $\mathbf{A}$
and $\mathbf{X}$ and integers $m,n\geq0$:
\begin{align*}
\mathcal{C}_{H}\left(\mathbf{A}\mathbf{X}^{-1}:m\right)\mathbf{X}^{n} & =\left(\mathbf{X}^{m}\mathbf{A}\mathbf{X}^{-1}\mathbf{X}^{-m}\right)\mathbf{X}^{n}\\
 & =\mathbf{X}^{m}\mathbf{A}\mathbf{X}^{-m+n-1}\\
 & =\mathbf{X}^{-n+1}\mathbf{X}^{m-n+1}\mathbf{A}\mathbf{X}^{-\left(m-n+1\right)}\\
 & =\mathbf{X}^{-n+1}\mathcal{C}_{H}\left(\mathbf{A}:m-n+1\right)
\end{align*}
As such:

\begin{eqnarray*}
 & \sum_{\left\Vert \mathbf{t}\right\Vert _{p}\leq p^{n-1}}\sum_{\mathbf{j}>\mathbf{0}}^{p-1}\hat{A}_{H}\left(\mathbf{t}+\frac{\mathbf{j}}{p^{n}}\right)\gamma_{\mathbf{j}}e^{2\pi i\left\{ \frac{\mathbf{j}\cdot\mathbf{z}}{p^{n}}\right\} _{p}}e^{2\pi i\left\{ \mathbf{t}\cdot\mathbf{z}\right\} _{p}}\\
 & =\\
 & \kappa_{H}\left(\left[\mathbf{z}\right]_{p^{n-1}}\right)\sum_{\mathbf{j}>\mathbf{0}}^{p-1}\underbrace{\varepsilon_{n-1}\left(\mathbf{j}\mathbf{z}\right)}_{\textrm{a scalar}}\mathbf{X}^{-n+1}\mathcal{C}_{H}\left(\alpha_{H}\left(\frac{\mathbf{j}}{p}\right):\lambda_{p}\left(\left[\mathbf{z}\right]_{p^{n-1}}\right)-n+1\right)\gamma_{\mathbf{j}}\\
 & =\\
 & \kappa_{H}\left(\left[\mathbf{z}\right]_{p^{n-1}}\right)\mathbf{X}^{-n+1}\sum_{\mathbf{j}>\mathbf{0}}^{p-1}\mathcal{C}_{H}\left(\alpha_{H}\left(\frac{\mathbf{j}}{p}\right)\varepsilon_{n-1}\left(\mathbf{j}\mathbf{z}\right):\lambda_{p}\left(\left[\mathbf{z}\right]_{p^{n-1}}\right)-n+1\right)\gamma_{\mathbf{j}}
\end{eqnarray*}
Since: 
\begin{align*}
\mathbf{X}^{-n+1}\mathcal{C}_{H}\left(\mathbf{A}:m-\left(n-1\right)\right) & =\mathbf{X}^{-n+1}\left(\mathbf{X}^{m-\left(n-1\right)}\mathbf{A}\mathbf{X}^{-m+n-1}\right)\\
 & =\mathbf{X}^{-n+1}\mathbf{X}^{-\left(n-1\right)}\left(\mathbf{x}^{m}\mathbf{A}\mathbf{x}^{-m}\right)\mathbf{X}^{n-1}\\
 & =\mathcal{C}_{H}\left(\mathbf{A}:m\right)\mathbf{X}^{n-1}
\end{align*}
we then have:

\begin{eqnarray}
 & \sum_{\left\Vert \mathbf{t}\right\Vert _{p}\leq p^{n-1}}\sum_{\mathbf{j}>\mathbf{0}}^{p-1}\hat{A}_{H}\left(\mathbf{t}+\frac{\mathbf{j}}{p^{n}}\right)\gamma_{\mathbf{j}}e^{2\pi i\left\{ \frac{\mathbf{j}\mathbf{z}}{p^{n}}\right\} _{p}}e^{2\pi i\left\{ \mathbf{t}\mathbf{z}\right\} _{p}}\nonumber \\
 & =\label{eq:MD 2/3rds and 1/4th and a bit of the way through the gamma computation}\\
 & \kappa_{H}\left(\left[\mathbf{z}\right]_{p^{n-1}}\right)\sum_{\mathbf{j}>\mathbf{0}}^{p-1}\mathcal{C}_{H}\left(\alpha_{H}\left(\frac{\mathbf{j}}{p}\right)\varepsilon_{n-1}\left(\mathbf{j}\mathbf{z}\right):\lambda_{p}\left(\left[\mathbf{z}\right]_{p^{n-1}}\right)\right)\mathbf{X}^{n-1}\gamma_{\mathbf{j}}\nonumber 
\end{eqnarray}

Finally, returning with this to (\ref{eq:MD 1/3rd of the way through the gamma computation}),
we get:
\begin{align}
F_{N}\left(\mathbf{z}\right) & =\sum_{n=1}^{N}\kappa_{H}\left(\left[\mathbf{z}\right]_{p^{n-1}}\right)\sum_{\mathbf{j}>\mathbf{0}}^{p-1}\mathcal{C}_{H}\left(\alpha_{H}\left(\frac{\mathbf{j}}{p}\right)\varepsilon_{n-1}\left(\mathbf{j}\mathbf{z}\right):\lambda_{p}\left(\left[\mathbf{z}\right]_{p^{n-1}}\right)\right)\left(H^{\prime}\left(\mathbf{0}\right)\right)^{n-1}\gamma_{\mathbf{j}}
\end{align}
which gives (\ref{eq:MD Gamma formula}) after re-indexing $n$ by
a shift of $1$.

\vphantom{}

II. Taking the $n$th term from (\ref{eq:MD Gamma formula}), we first
apply $q$-adic norm to get and upper bound of: 
\begin{equation}
\left\Vert \kappa_{H}\left(\left[\mathbf{z}\right]_{p^{n}}\right)\right\Vert _{q}\left\Vert H^{\prime}\left(\mathbf{0}\right)\right\Vert _{q}^{n}\max_{\mathbf{0}<\mathbf{j}\leq p-1}\left\Vert \beta_{H}\left(\frac{\mathbf{j}}{p}\right)\right\Vert _{q}
\end{equation}
In this, we used: 
\begin{equation}
\gamma_{\mathbf{j}}=\gamma_{H}\left(\frac{\mathbf{j}}{p}\right)=\left(\alpha_{H}\left(\frac{\mathbf{j}}{p}\right)\right)^{-1}\beta_{H}\left(\frac{\mathbf{j}}{p}\right)
\end{equation}
Since the $q$-adic bound on the $\beta_{H}$s is independent of $n$,
and since $\left\Vert H^{\prime}\left(\mathbf{0}\right)\right\Vert _{q}=1$,
the fact that $\left\Vert \kappa_{H}\left(\left[\mathbf{z}\right]_{p^{n}}\right)\right\Vert _{q}$
tends to $0$ as $n\rightarrow\infty$ for all $\mathbf{z}\in\left(\mathbb{Z}_{p}^{r}\right)^{\prime}$
then guarantees the convergence of (\ref{eq:F limit of MD Gamma_H A_H hat Fourier series})
for $\mathbf{z}\in\left(\mathbb{Z}_{p}^{r}\right)^{\prime}$.

For $\mathbf{z}\in\mathbb{N}_{0}^{r}$, applying the archimedean complex
norm to the $n$th term of (\ref{eq:MD Gamma formula}) gives an upper
bound of: 
\begin{equation}
\sum_{n=0}^{N-1}\left\Vert \kappa_{H}\left(\left[\mathbf{z}\right]_{p^{n}}\right)\right\Vert _{\infty}\left\Vert H^{\prime}\left(\mathbf{0}\right)\right\Vert _{\infty}^{n}\sum_{\mathbf{j}>\mathbf{0}}^{p-1}\left\Vert \beta_{H}\left(\frac{\mathbf{j}}{p}\right)\right\Vert _{\infty}
\end{equation}
Here, $\sum_{\mathbf{j}>\mathbf{0}}^{p-1}\left\Vert \beta_{H}\left(\frac{\mathbf{j}}{p}\right)\right\Vert _{\infty}$
is uniformly bounded with respect to $n$. Moreover, since $\mathbf{z}\in\mathbb{N}_{0}^{r}$,
we have that $\kappa_{H}\left(\left[\mathbf{z}\right]_{p^{n}}\right)=\kappa_{H}\left(\mathbf{z}\right)$
for all sufficiently large $n$. So, if we let $N\rightarrow\infty$,
we obtain:
\begin{equation}
\sum_{n=0}^{N-1}\left\Vert \kappa_{H}\left(\left[\mathbf{z}\right]_{p^{n}}\right)\right\Vert _{\infty}\left\Vert H^{\prime}\left(\mathbf{0}\right)\right\Vert _{\infty}^{n}\sum_{\mathbf{j}>\mathbf{0}}^{p-1}\left\Vert \beta_{H}\left(\frac{\mathbf{j}}{p}\right)\right\Vert _{\infty}\ll O\left(1\right)+\sum_{n}\left\Vert H^{\prime}\left(\mathbf{0}\right)\right\Vert _{\infty}^{n}
\end{equation}
Here, the upper bound is a convergent geometric series because $H$
is contracting. This then guarantees the convergence of (\ref{eq:F limit of MD Gamma_H A_H hat Fourier series})
for $\mathbf{z}\in\mathbb{N}_{0}^{r}$.

Q.E.D.

\vphantom{}

With these formulae, we can now sum the Fourier series generated by
(\ref{eq:MD Fine structure of Chi_H,N hat when alpha is 1}) to obtain
a non-trivial formula for $\chi_{H,N}$. 
\begin{thm}
Suppose $\alpha_{H}\left(\mathbf{0}\right)=\mathbf{I}_{d}$. Then,
for all $N\geq1$ and all $\mathbf{z}\in\mathbb{Z}_{p}^{r}$: 
\begin{align}
\chi_{H,N}\left(\mathbf{z}\right) & \overset{\overline{\mathbb{Q}}^{d}}{=}\sum_{n=0}^{N-1}\kappa_{H}\left(\left[\mathbf{z}\right]_{p^{n}}\right)\sum_{\mathbf{j}=\mathbf{0}}^{p-1}\mathcal{C}_{H}\left(\alpha_{H}\left(\frac{\mathbf{j}}{p}\right)\varepsilon_{n}\left(\mathbf{j}\mathbf{z}\right):\lambda_{p}\left(\left[\mathbf{z}\right]_{p^{n}}\right)\right)\left(H^{\prime}\left(\mathbf{0}\right)\right)^{n}\gamma_{H}\left(\frac{\mathbf{j}}{p}\right)\label{eq:MD Chi_H,N when alpha is 1 and rho is arbitrary}
\end{align}
In particular, when $p=2$: 
\begin{equation}
\chi_{H,N}\left(\mathbf{z}\right)=-\gamma_{H}\left(\mathbf{\frac{1}{2}}\right)+\kappa_{H}\left(\left[\mathbf{z}\right]_{2^{N}}\right)\left(H^{\prime}\left(\mathbf{0}\right)\right)^{N}\gamma_{H}\left(\mathbf{\frac{1}{2}}\right)+\sum_{n=0}^{N-1}\kappa_{H}\left(\left[\mathbf{z}\right]_{2^{n}}\right)\left(H^{\prime}\left(\mathbf{0}\right)\right)^{n}\beta_{H}\left(\mathbf{0}\right)\label{eq:MD Chi_H,N when alpha is 1 and every prime in P is 2}
\end{equation}
where: 
\begin{equation}
\gamma_{H}\left(\mathbf{\frac{1}{2}}\right)\overset{\textrm{def}}{=}\gamma_{H}\left(\left(\frac{1}{2},\ldots,\frac{1}{2}\right)\right)\label{eq:Definition of MD gamma_H of 1/2}
\end{equation}
\end{thm}
Proof: We start by multiplying (\ref{eq:MD Fine structure of Chi_H,N hat when alpha is 1})
by $e^{2\pi i\left\{ \mathbf{t}\mathbf{z}\right\} _{p}}$ and then
summing over all $\left\Vert \mathbf{t}\right\Vert _{p}\leq p^{N}$.
The left-hand side of (\ref{eq:MD Fine structure of Chi_H,N hat when alpha is 1})
becomes: 
\begin{equation}
\chi_{H,N}\left(\mathbf{z}\right)-N\tilde{A}_{H,N-1}\left(\mathbf{z}\right)\beta_{H}\left(\mathbf{0}\right)
\end{equation}
while the right-hand side becomes: 
\begin{align*}
 & \sum_{0<\left\Vert \mathbf{t}\right\Vert _{p}\leq p^{N-1}}\hat{A}_{H}\left(\mathbf{t}\right)\left(v_{p}\left(\mathbf{t}\right)\beta_{H}\left(\mathbf{0}\right)+\gamma_{H}\left(\frac{\mathbf{t}\left|\mathbf{t}\right|_{p}}{p}\right)\right)e^{2\pi i\left\{ \mathbf{t}\mathbf{z}\right\} _{p}}\\
 & +\sum_{\left\Vert \mathbf{t}\right\Vert _{p}=p^{N}}\hat{A}_{H}\left(\mathbf{t}\right)\gamma_{H}\left(\frac{\mathbf{t}\left|\mathbf{t}\right|_{p}}{p}\right)e^{2\pi i\left\{ \mathbf{t}\mathbf{z}\right\} _{p}}
\end{align*}
Simplifying produces: 
\begin{align}
\chi_{H,N}\left(\mathbf{z}\right) & \overset{\overline{\mathbb{Q}}^{d}}{=}N\tilde{A}_{H,N-1}\left(\mathbf{z}\right)\beta_{H}\left(\mathbf{0}\right)\label{eq:MD Chi_H,N rho not equal to 2, ready to simplify}\\
 & +\sum_{0<\left\Vert \mathbf{t}\right\Vert _{p}\leq p^{N-1}}\hat{A}_{H}\left(\mathbf{t}\right)v_{p}\left(\mathbf{t}\right)\beta_{H}\left(\mathbf{0}\right)e^{2\pi i\left\{ \mathbf{t}\mathbf{z}\right\} _{p}}\nonumber \\
 & +\sum_{0<\left\Vert \mathbf{t}\right\Vert _{p}\leq p^{N}}\hat{A}_{H}\left(\mathbf{t}\right)\gamma_{H}\left(\frac{\mathbf{t}\left|\mathbf{t}\right|_{p}}{p}\right)e^{2\pi i\left\{ \mathbf{t}\mathbf{z}\right\} _{p}}\nonumber 
\end{align}
Once again, we call upon our legion of formulae: (\ref{eq:MD Convolution of dA_H and D_N}),
\textbf{Lemmata \ref{lem:MD v_p A_H hat summation formulae}}, and
\textbf{\ref{lem:MD gamma formulae}}. Using them (while remembering
that $\alpha_{H}\left(\mathbf{0}\right)=\mathbf{I}_{d}$ makes $\mathbf{I}_{H}$
identically $\mathbf{O}_{d}$) turns (\ref{eq:MD Chi_H,N rho not equal to 2, ready to simplify})
into:

\begin{align*}
\chi_{H,N}\left(\mathbf{z}\right) & \overset{\overline{\mathbb{Q}}^{d}}{=}N\kappa_{H}\left(\left[\mathbf{z}\right]_{p^{N-1}}\right)\left(H^{\prime}\left(\mathbf{0}\right)\right)^{N-1}\beta_{H}\left(\mathbf{0}\right)\\
 & -\left(N-1\right)\kappa_{H}\left(\left[\mathbf{z}\right]_{p^{N-1}}\right)\left(H^{\prime}\left(\mathbf{0}\right)\right)^{N-1}\beta_{H}\left(\mathbf{0}\right)\\
 & +\sum_{n=0}^{N-2}\kappa_{H}\left(\left[\mathbf{z}\right]_{p^{n}}\right)\left(H^{\prime}\left(\mathbf{0}\right)\right)^{n}\beta_{H}\left(\mathbf{0}\right)\\
 & +\sum_{n=0}^{N-1}\kappa_{H}\left(\left[\mathbf{z}\right]_{p^{n}}\right)\sum_{\mathbf{j}>\mathbf{0}}^{p-1}\mathcal{C}_{H}\left(\alpha_{H}\left(\frac{\mathbf{j}}{p}\right)\varepsilon_{n}\left(\mathbf{j}\mathbf{z}\right):\lambda_{p}\left(\left[\mathbf{z}\right]_{p^{n}}\right)\right)\left(H^{\prime}\left(\mathbf{0}\right)\right)^{n}\gamma_{H}\left(\frac{\mathbf{j}}{p}\right)
\end{align*}
Simplifying yields:

\begin{align*}
\chi_{H,N}\left(\mathbf{z}\right) & \overset{\overline{\mathbb{Q}}^{d}}{=}\sum_{n=0}^{N-1}\kappa_{H}\left(\left[\mathbf{z}\right]_{p^{n}}\right)\left(H^{\prime}\left(\mathbf{0}\right)\right)^{n}\beta_{H}\left(\mathbf{0}\right)\\
 & +\sum_{n=0}^{N-1}\kappa_{H}\left(\left[\mathbf{z}\right]_{p^{n}}\right)\sum_{\mathbf{j}>\mathbf{0}}^{p-1}\mathcal{C}_{H}\left(\alpha_{H}\left(\frac{\mathbf{j}}{p}\right)\varepsilon_{n}\left(\mathbf{j}\mathbf{z}\right):\lambda_{p}\left(\left[\mathbf{z}\right]_{p^{n}}\right)\right)\left(H^{\prime}\left(\mathbf{0}\right)\right)^{n}\gamma_{H}\left(\frac{\mathbf{j}}{p}\right)
\end{align*}
Now, because $\alpha_{H}\left(\mathbf{0}\right)=\mathbf{I}_{d}$,
the vector:
\begin{equation}
\mathcal{C}_{H}\left(\alpha_{H}\left(\frac{\mathbf{0}}{p}\right)\varepsilon_{n}\left(\mathbf{0}\mathbf{z}\right):\lambda_{p}\left(\left[\mathbf{z}\right]_{p^{n}}\right)\right)\left(H^{\prime}\left(\mathbf{0}\right)\right)^{n}\gamma_{H}\left(\frac{\mathbf{0}}{p}\right)
\end{equation}
becomes: 
\begin{equation}
\underbrace{\mathcal{C}_{H}\left(\mathbf{I}_{d}:\lambda_{p}\left(\left[\mathbf{z}\right]_{p^{n}}\right)\right)}_{\mathbf{I}_{d}}\left(H^{\prime}\left(\mathbf{0}\right)\right)^{n}\beta_{H}\left(\mathbf{0}\right)=\left(H^{\prime}\left(\mathbf{0}\right)\right)^{n}\beta_{H}\left(\mathbf{0}\right)
\end{equation}
So, the two lines on the right-hand side of the above formula for
$\chi_{H,N}\left(\mathbf{z}\right)$ combine to form a single sum
because the upper line is the $\mathbf{j}=\mathbf{0}$ case of the
bottom line.

Finally, when $p=2$, we can compute everything directly from (\ref{eq:MD Fine structure of Chi_H,N hat when alpha is 1})
by multiplying by $e^{2\pi i\left\{ \mathbf{t}\mathbf{z}\right\} _{2}}$
and summing over all $\left\Vert \mathbf{t}\right\Vert _{2}\leq2^{N}$.
This gives: 
\begin{align*}
\chi_{H,N}\left(\mathbf{z}\right)-N\tilde{A}_{H,N-1}\left(\mathbf{z}\right)\beta_{H}\left(\mathbf{0}\right)= & \sum_{0<\left\Vert \mathbf{t}\right\Vert _{2}\leq2^{N-1}}\hat{A}_{H}\left(\mathbf{t}\right)v_{2}\left(\mathbf{t}\right)e^{2\pi i\left\{ \mathbf{t}\mathbf{z}\right\} _{2}}\beta_{H}\left(\mathbf{0}\right)\\
 & +\sum_{0<\left\Vert \mathbf{t}\right\Vert _{2}\leq2^{N}}\hat{A}_{H}\left(\mathbf{t}\right)e^{2\pi i\left\{ \mathbf{t}\mathbf{z}\right\} _{2}}\gamma_{H}\left(\mathbf{\frac{1}{2}}\right)
\end{align*}
which simplifies to: 
\begin{align*}
\chi_{H,N}\left(\mathbf{z}\right)= & N\tilde{A}_{H,N-1}\left(\mathbf{z}\right)\beta_{H}\left(\mathbf{0}\right)-\gamma_{H}\left(\mathbf{\frac{1}{2}}\right)+\sum_{\left\Vert \mathbf{t}\right\Vert _{2}\leq2^{N}}\hat{A}_{H}\left(\mathbf{t}\right)e^{2\pi i\left\{ \mathbf{t}\mathbf{z}\right\} _{2}}\gamma_{H}\left(\mathbf{\frac{1}{2}}\right)\\
 & +\sum_{0<\left\Vert \mathbf{t}\right\Vert _{2}\leq2^{N-1}}\hat{A}_{H}\left(\mathbf{t}\right)v_{2}\left(\mathbf{t}\right)e^{2\pi i\left\{ \mathbf{t}\mathbf{z}\right\} _{2}}\beta_{H}\left(\mathbf{0}\right)
\end{align*}
Applying (\ref{eq:MD Convolution of dA_H and D_N}) and \textbf{Lemma
\ref{lem:MD v_p A_H hat summation formulae}}\textemdash and, again
remembering that $\alpha_{H}\left(\mathbf{0}\right)=\mathbf{I}_{d}$
makes $\mathbf{I}_{H}$ identically $\mathbf{O}_{d}$\textemdash the
above becomes: 
\begin{align*}
\chi_{H,N}\left(\mathbf{z}\right)= & N\kappa_{H}\left(\left[\mathbf{z}\right]_{2^{N-1}}\right)\left(H^{\prime}\left(\mathbf{0}\right)\right)^{N-1}\beta_{H}\left(\mathbf{0}\right)-\gamma_{H}\left(\mathbf{\frac{1}{2}}\right)\\
 & +\kappa_{H}\left(\left[\mathbf{z}\right]_{2^{N}}\right)\left(H^{\prime}\left(\mathbf{0}\right)\right)^{N}\gamma_{H}\left(\mathbf{\frac{1}{2}}\right)\\
 & -\left(N-1\right)\kappa_{H}\left(\left[\mathbf{z}\right]_{2^{N-1}}\right)\left(H^{\prime}\left(\mathbf{0}\right)\right)^{N-1}\beta_{H}\left(\mathbf{0}\right)\\
 & +\sum_{n=0}^{N-2}\kappa_{H}\left(\left[\mathbf{z}\right]_{2^{n}}\right)\left(H^{\prime}\left(\mathbf{0}\right)\right)^{n}\beta_{H}\left(\mathbf{0}\right)
\end{align*}
Simplifying then gives us what we wanted: 
\[
\chi_{H,N}\left(\mathbf{z}\right)=-\gamma_{H}\left(\mathbf{\frac{1}{2}}\right)+\kappa_{H}\left(\left[\mathbf{z}\right]_{2^{N}}\right)\left(H^{\prime}\left(\mathbf{0}\right)\right)^{N}\gamma_{H}\left(\mathbf{\frac{1}{2}}\right)+\sum_{n=0}^{N-1}\kappa_{H}\left(\left[\mathbf{z}\right]_{2^{n}}\right)\left(H^{\prime}\left(\mathbf{0}\right)\right)^{n}\beta_{H}\left(\mathbf{0}\right)
\]

Q.E.D. 
\begin{cor}[\textbf{$\mathcal{F}$-Series for $\chi_{H}$ when $\alpha_{H}\left(\mathbf{0}\right)=\mathbf{I}_{d}$}]
\label{cor:MD alpha is 1 case, F-series}Suppose $\alpha_{H}\left(\mathbf{0}\right)=\mathbf{I}_{d}$.
Then: 
\begin{equation}
\chi_{H}\left(\mathbf{z}\right)\overset{\mathcal{F}_{p,q_{H}}^{d}}{=}\sum_{n=0}^{\infty}\kappa_{H}\left(\left[\mathbf{z}\right]_{p^{n}}\right)\sum_{\mathbf{j}=\mathbf{0}}^{p-1}\mathcal{C}_{H}\left(\alpha_{H}\left(\frac{\mathbf{j}}{p}\right)\varepsilon_{n}\left(\mathbf{j}\mathbf{z}\right):\lambda_{p}\left(\left[\mathbf{z}\right]_{p^{n}}\right)\right)\left(H^{\prime}\left(\mathbf{0}\right)\right)^{n}\gamma_{H}\left(\frac{\mathbf{j}}{p}\right)\label{eq:MD Explicit Formula for Chi_H when alpha is 1 and rho is arbitrary}
\end{equation}
for all $\mathbf{z}\in\mathbb{Z}_{p}^{r}$, with the special case:
\begin{equation}
\chi_{H}\left(\mathbf{z}\right)\overset{\mathcal{F}_{2,q_{H}}^{d}}{=}-\gamma_{H}\left(\mathbf{\frac{1}{2}}\right)+\sum_{n=0}^{\infty}\kappa_{H}\left(\left[\mathbf{z}\right]_{2^{n}}\right)\left(H^{\prime}\left(\mathbf{0}\right)\right)^{n}\beta_{H}\left(\mathbf{0}\right)\label{eq:MD Explicit Formula for Chi_H when alpha is 1 and every prime in P is 2}
\end{equation}
for all $\mathbf{z}\in\mathbb{Z}_{2}$ when $p=2$. 
\end{cor}
Proof: Same as in the one-dimensional case, but with vector norms.

Q.E.D.

\vphantom{}

Next, we sum in $\mathbb{C}$ for $\mathbf{z}\in\mathbb{N}_{0}^{r}$. 
\begin{cor}
\label{cor:MD alpha is 1, F-series on N_0^r}Suppose that $\alpha_{H}\left(\mathbf{0}\right)=\mathbf{I}_{d}$.
Then: 
\begin{align}
\chi_{H}\left(\mathbf{n}\right) & \overset{\mathbb{C}^{d}}{=}\sum_{k=0}^{\lambda_{p}\left(\mathbf{n}\right)-1}\kappa_{H}\left(\left[\mathbf{n}\right]_{p^{k}}\right)\sum_{\mathbf{j}=\mathbf{0}}^{p-1}\mathcal{C}_{H}\left(\alpha_{H}\left(\frac{\mathbf{j}}{p}\right)\varepsilon_{k}\left(\mathbf{j}\mathbf{n}\right):\lambda_{p}\left(\left[\mathbf{n}\right]_{p^{k}}\right)\right)\left(H^{\prime}\left(\mathbf{0}\right)\right)^{k}\gamma_{H}\left(\frac{\mathbf{j}}{p}\right)\label{eq:MD archimedean Chi_H when rho is arbitrary and alpha_H of 0 is 1}\\
 & +M_{H}\left(\mathbf{n}\right)\sum_{\mathbf{j}=\mathbf{0}}^{p-1}\alpha_{H}\left(\frac{\mathbf{j}}{p}\right)\left(\mathbf{I}_{d}-H^{\prime}\left(\mathbf{0}\right)\right)^{-1}\gamma_{H}\left(\frac{\mathbf{j}}{p}\right)
\end{align}
for all $\mathbf{n}\in\mathbb{N}_{0}^{r}$, with the special case:
\begin{align}
\chi_{H}\left(\mathbf{n}\right) & \overset{\mathbb{C}^{d}}{=}-\gamma_{H}\left(\mathbf{\frac{1}{2}}\right)+M_{H}\left(\mathbf{n}\right)\left(\mathbf{I}_{d}-H^{\prime}\left(\mathbf{0}\right)\right)^{-1}\beta_{H}\left(\mathbf{0}\right)\label{eq:MD archimedean Chi_H when every prime in P is 2 and alpha_H of 0 is 1}\\
 & +\sum_{k=0}^{\lambda_{2}\left(\mathbf{n}\right)-1}\kappa_{H}\left(\left[\mathbf{n}\right]_{2^{k}}\right)\left(H^{\prime}\left(\mathbf{0}\right)\right)^{k}\beta_{H}\left(\mathbf{0}\right)\chi_{H}\left(\mathbf{n}\right)
\end{align}
when $p=2$.

Regardless of the value of $p$, the $k$-sums are defined to be $\mathbf{0}$
when $\mathbf{n}=\mathbf{0}$.
\end{cor}
Proof: Let $\mathbf{n}\in\mathbb{N}_{0}^{r}$. Since $\left[\mathbf{n}\right]_{p^{k}}=\mathbf{n}$
and $\varepsilon_{k}\left(\mathbf{n}\right)=1$ for all $k\geq\lambda_{p}\left(\mathbf{n}\right)$,
(\ref{eq:MD Explicit Formula for Chi_H when alpha is 1 and rho is arbitrary})
becomes: 
\begin{align*}
\chi_{H}\left(\mathbf{n}\right) & \overset{\mathbb{C}^{d}}{=}\sum_{k=0}^{\lambda_{p}\left(\mathbf{n}\right)-1}\kappa_{H}\left(\left[\mathbf{n}\right]_{p^{k}}\right)\sum_{\mathbf{j}=\mathbf{0}}^{p-1}\mathcal{C}_{H}\left(\alpha_{H}\left(\frac{\mathbf{j}}{p}\right)\varepsilon_{k}\left(\mathbf{j}\mathbf{n}\right):\lambda_{p}\left(\left[\mathbf{n}\right]_{p^{k}}\right)\right)\left(H^{\prime}\left(\mathbf{0}\right)\right)^{k}\gamma_{H}\left(\frac{\mathbf{j}}{p}\right)\\
 & +\kappa_{H}\left(\mathbf{n}\right)\sum_{k=\lambda_{p}\left(\mathbf{n}\right)}^{\infty}\sum_{\mathbf{j}=\mathbf{0}}^{p-1}\mathcal{C}_{H}\left(\alpha_{H}\left(\frac{\mathbf{j}}{p}\right):\lambda_{p}\left(\mathbf{n}\right)\right)\left(H^{\prime}\left(\mathbf{0}\right)\right)^{k}\gamma_{H}\left(\frac{\mathbf{j}}{p}\right)
\end{align*}
Here: 
\[
\sum_{k=\lambda_{p}\left(\mathbf{n}\right)}^{\infty}\sum_{\mathbf{j}=\mathbf{0}}^{p-1}\mathcal{C}_{H}\left(\alpha_{H}\left(\frac{\mathbf{j}}{p}\right):\lambda_{p}\left(\mathbf{n}\right)\right)\left(H^{\prime}\left(\mathbf{0}\right)\right)^{k}\gamma_{H}\left(\frac{\mathbf{j}}{p}\right)
\]
becomes: 
\[
\sum_{k=\lambda_{p}\left(\mathbf{n}\right)}^{\infty}\sum_{\mathbf{j}=\mathbf{0}}^{p-1}\left(H^{\prime}\left(\mathbf{0}\right)\right)^{\lambda_{p}\left(\mathbf{n}\right)}\alpha_{H}\left(\frac{\mathbf{j}}{p}\right)\left(H^{\prime}\left(\mathbf{0}\right)\right)^{k-\lambda_{p}\left(\mathbf{n}\right)}\gamma_{H}\left(\frac{\mathbf{j}}{p}\right)
\]
which simplifies to :

\[
\sum_{\mathbf{j}=\mathbf{0}}^{p-1}\left(H^{\prime}\left(\mathbf{0}\right)\right)^{\lambda_{p}\left(\mathbf{n}\right)}\alpha_{H}\left(\frac{\mathbf{j}}{p}\right)\sum_{k=0}^{\infty}\left(H^{\prime}\left(\mathbf{0}\right)\right)^{k}\gamma_{H}\left(\frac{\mathbf{j}}{p}\right)
\]
Since $H$ is contracting, we get a $\mathbb{C}$-convergent geometric
series: 
\[
\sum_{k=0}^{\infty}\left(H^{\prime}\left(\mathbf{0}\right)\right)^{k}\overset{\mathbb{C}}{=}\left(\mathbf{I}_{d}-H^{\prime}\left(\mathbf{0}\right)\right)^{-1}
\]
and so, we are left with: 
\[
\sum_{\mathbf{j}=\mathbf{0}}^{p-1}\left(H^{\prime}\left(\mathbf{0}\right)\right)^{\lambda_{p}\left(\mathbf{n}\right)}\alpha_{H}\left(\frac{\mathbf{j}}{p}\right)\left(\mathbf{I}_{d}-H^{\prime}\left(\mathbf{0}\right)\right)^{-1}\gamma_{H}\left(\frac{\mathbf{j}}{p}\right)
\]
Consequently: 
\begin{align*}
\chi_{H}\left(\mathbf{n}\right) & \overset{\mathbb{C}^{d}}{=}\sum_{k=0}^{\lambda_{p}\left(\mathbf{n}\right)-1}\kappa_{H}\left(\left[\mathbf{n}\right]_{p^{k}}\right)\sum_{\mathbf{j}=\mathbf{0}}^{p-1}\mathcal{C}_{H}\left(\alpha_{H}\left(\frac{\mathbf{j}}{p}\right)\varepsilon_{k}\left(\mathbf{j}\mathbf{n}\right):\lambda_{p}\left(\left[\mathbf{n}\right]_{p^{k}}\right)\right)\left(H^{\prime}\left(\mathbf{0}\right)\right)^{k}\gamma_{H}\left(\frac{\mathbf{j}}{p}\right)\\
 & +\kappa_{H}\left(\mathbf{n}\right)\left(H^{\prime}\left(\mathbf{0}\right)\right)^{\lambda_{p}\left(\mathbf{n}\right)}\sum_{\mathbf{j}=\mathbf{0}}^{p-1}\alpha_{H}\left(\frac{\mathbf{j}}{p}\right)\left(\mathbf{I}_{d}-H^{\prime}\left(\mathbf{0}\right)\right)^{-1}\gamma_{H}\left(\frac{\mathbf{j}}{p}\right)
\end{align*}
Since: 
\[
\kappa_{H}\left(\mathbf{n}\right)=M_{H}\left(\mathbf{n}\right)\left(H^{\prime}\left(\mathbf{0}\right)\right)^{-\lambda_{p}\left(\mathbf{n}\right)}
\]
the above can be simplified to: 
\begin{align*}
\chi_{H}\left(\mathbf{n}\right) & \overset{\mathbb{C}^{d}}{=}\sum_{k=0}^{\lambda_{p}\left(\mathbf{n}\right)-1}\kappa_{H}\left(\left[\mathbf{n}\right]_{p^{k}}\right)\sum_{\mathbf{j}=\mathbf{0}}^{p-1}\mathcal{C}_{H}\left(\alpha_{H}\left(\frac{\mathbf{j}}{p}\right)\varepsilon_{k}\left(\mathbf{j}\mathbf{n}\right):\lambda_{p}\left(\left[\mathbf{n}\right]_{p^{k}}\right)\right)\left(H^{\prime}\left(\mathbf{0}\right)\right)^{k}\gamma_{H}\left(\frac{\mathbf{j}}{p}\right)\\
 & +M_{H}\left(\mathbf{n}\right)\sum_{\mathbf{j}=\mathbf{0}}^{p-1}\alpha_{H}\left(\frac{\mathbf{j}}{p}\right)\left(\mathbf{I}_{d}-H^{\prime}\left(\mathbf{0}\right)\right)^{-1}\gamma_{H}\left(\frac{\mathbf{j}}{p}\right)
\end{align*}
which is the desired formula.

As for the case where $p=2$, applying the above argument to (\ref{eq:MD Explicit Formula for Chi_H when alpha is 1 and every prime in P is 2})
given produces: 
\begin{align*}
\chi_{H}\left(\mathbf{n}\right) & \overset{\mathbb{C}^{d}}{=}-\gamma_{H}\left(\mathbf{\frac{1}{2}}\right)+\sum_{k=0}^{\lambda_{2}\left(\mathbf{n}\right)-1}\kappa_{H}\left(\left[\mathbf{n}\right]_{2^{k}}\right)\left(H^{\prime}\left(\mathbf{0}\right)\right)^{k}\beta_{H}\left(\mathbf{0}\right)\\
 & +\sum_{k=\lambda_{2}\left(\mathbf{n}\right)}^{\infty}\kappa_{H}\left(\mathbf{n}\right)\left(H^{\prime}\left(\mathbf{0}\right)\right)^{k}\beta_{H}\left(\mathbf{0}\right)\\
\left(H\textrm{ is contracting}\right); & \overset{\mathbb{C}^{d}}{=}-\gamma_{H}\left(\mathbf{\frac{1}{2}}\right)+\sum_{k=0}^{\lambda_{2}\left(\mathbf{n}\right)-1}\kappa_{H}\left(\left[\mathbf{n}\right]_{2^{k}}\right)\left(H^{\prime}\left(\mathbf{0}\right)\right)^{k}\beta_{H}\left(\mathbf{0}\right)\\
 & +\underbrace{\kappa_{H}\left(\mathbf{n}\right)\left(H^{\prime}\left(\mathbf{0}\right)\right)^{\lambda_{2}\left(\mathbf{n}\right)}}_{M_{H}\left(\mathbf{n}\right)}\left(\mathbf{I}_{d}-H^{\prime}\left(\mathbf{0}\right)\right)^{-1}\beta_{H}\left(\mathbf{0}\right)
\end{align*}

Q.E.D.

\vphantom{}

Taken together, these two corollaries establish the quasi-integrability
of $\chi_{H}$ for arbitrary $p$, provided that $\alpha_{H}\left(\mathbf{0}\right)=\mathbf{I}_{d}$:
\begin{cor}[\textbf{Quasi-Integrability of $\chi_{H}$ When $\alpha_{H}\left(\mathbf{0}\right)=\mathbf{I}_{d}$}]
If $\alpha_{H}\left(\mathbf{0}\right)=\mathbf{I}_{d}$, then, $\chi_{H}$
is quasi-integrable with respect to the standard $\left(p,q_{H}\right)$-adic
frame.

In particular, when $p=2$, the function $\hat{\chi}_{H}:\hat{\mathbb{Z}}_{2}^{r}\rightarrow\overline{\mathbb{Q}}^{d}$
defined by: 
\begin{equation}
\hat{\chi}_{H}\left(\mathbf{t}\right)\overset{\textrm{def}}{=}\begin{cases}
-\gamma_{H}\left(\mathbf{\frac{1}{2}}\right) & \textrm{if }\mathbf{t}=\mathbf{0}\\
\hat{A}_{H}\left(\mathbf{t}\right)v_{2}\left(\mathbf{t}\right)\beta_{H}\left(\mathbf{0}\right) & \textrm{else }
\end{cases},\textrm{ }\forall\mathbf{t}\in\hat{\mathbb{Z}}_{2}\label{eq:MD Formula for Chi_H hat when every prime in P is 2 and alpha is 1}
\end{equation}
is then a Fourier transform of $\chi_{H}$. In this case, the function
defined by: 
\begin{equation}
\hat{\chi}_{H}\left(\mathbf{t}\right)=\begin{cases}
\mathbf{0} & \textrm{if }\mathbf{t}=\mathbf{0}\\
\hat{A}_{H}\left(\mathbf{t}\right)\left(v_{2}\left(\mathbf{t}\right)\beta_{H}\left(\mathbf{0}\right)+\gamma_{H}\left(\mathbf{\frac{1}{2}}\right)\right) & \textrm{if }\textrm{else}
\end{cases},\textrm{ }\forall\mathbf{t}\in\hat{\mathbb{Z}}_{2}^{r}\label{eq:MD Formula for Chi_H hat when every prime in P is 2 and alpha is 1, alt}
\end{equation}
is also a Fourier transform of $\chi_{H}$, differing from the $\hat{\chi}_{H}$
given above by $\hat{A}_{H}\left(\mathbf{t}\right)\gamma_{H}\left(\mathbf{\frac{1}{2}}\right)$.
\emph{By}\textbf{\emph{ Theorem \ref{thm:MD properties of A_H hat}}},\textbf{
}since $\alpha_{H}\left(\mathbf{0}\right)=\mathbf{I}_{d}$, the function
$\hat{A}_{H}\left(\mathbf{t}\right)\gamma_{H}\left(\mathbf{\frac{1}{2}}\right)$
is then the Fourier-Stieltjes transform of a degenerate thick measure
of vector type.

For odd primes $p$, we can obtain a Fourier transform for $\chi_{H}$
by defining a function $\hat{\chi}_{H}:\hat{\mathbb{Z}}_{p}^{r}\rightarrow\overline{\mathbb{Q}}^{d}$
by: 
\begin{equation}
\hat{\chi}_{H}\left(\mathbf{t}\right)\overset{\textrm{def}}{=}\begin{cases}
\mathbf{0} & \textrm{if }\mathbf{t}=\mathbf{0}\\
\hat{A}_{H}\left(\mathbf{t}\right)\left(v_{p}\left(\mathbf{t}\right)\beta_{H}\left(\mathbf{0}\right)+\gamma_{H}\left(\frac{\mathbf{t}\left|\mathbf{t}\right|_{p}}{p}\right)\right) & \textrm{else}
\end{cases},\textrm{ }\forall\mathbf{t}\in\hat{\mathbb{Z}}_{p}^{r}\label{eq:MD Chi_H hat when rho is not 2 and when alpha is 1}
\end{equation}
\end{cor}
Proof: \textbf{Corollaries \ref{cor:MD alpha is 1 case, F-series}}
and \textbf{\ref{cor:MD alpha is 1, F-series on N_0^r}} show that
the $N$th partial sums of the Fourier series generated by (\ref{eq:MD Formula for Chi_H hat when every prime in P is 2 and alpha is 1})
and (\ref{eq:MD Chi_H hat when rho is not 2 and when alpha is 1})
are $\mathcal{F}_{p,q_{H}}$-convergent to (\ref{eq:MD Explicit Formula for Chi_H when alpha is 1 and every prime in P is 2})
and (\ref{eq:MD Explicit Formula for Chi_H when alpha is 1 and rho is arbitrary})
for the case where $p=2$ and the case $p\geq3$, respectively, thereby
establishing the quasi-integrability of $\chi_{H}$ with respect to
the standard $\left(p,q_{H}\right)$-adic frame.

Finally, letting $\hat{\chi}_{H}^{\prime}\left(\mathbf{t}\right)$
denote (\ref{eq:MD Formula for Chi_H hat when every prime in P is 2 and alpha is 1, alt}),
observe that when $\alpha_{H}\left(\mathbf{0}\right)=\mathbf{I}_{d}$
and $p=2$: 
\begin{equation}
\hat{\chi}_{H}^{\prime}\left(\mathbf{t}\right)-\hat{A}_{H}\left(\mathbf{t}\right)\gamma_{H}\left(\mathbf{\frac{1}{2}}\right)\overset{\overline{\mathbb{Q}}^{d}}{=}\begin{cases}
-\gamma_{H}\left(\mathbf{\frac{1}{2}}\right) & \textrm{if }\mathbf{t}=\mathbf{0}\\
\hat{A}_{H}\left(\mathbf{t}\right)v_{2}\left(\mathbf{t}\right)\beta_{H}\left(\mathbf{0}\right) & \textrm{if }\left\Vert \mathbf{t}\right\Vert _{2}>0
\end{cases}=\hat{\chi}_{H}\left(\mathbf{t}\right)
\end{equation}
which shows that $\hat{\chi}_{H}^{\prime}\left(\mathbf{t}\right)$
and $\hat{\chi}_{H}\left(\mathbf{t}\right)$ differ by a factor of
$\hat{A}_{H}\left(\mathbf{t}\right)\gamma_{H}\left(\mathbf{\frac{1}{2}}\right)$,
which is the Fourier-Stieltjes transform of a degenerate vector-type
thick measure, by \textbf{Theorem \ref{thm:MD properties of A_H hat}},
since $\alpha_{H}\left(\mathbf{0}\right)=\mathbf{I}_{d}$.

Q.E.D.

\subsection{\label{subsec:6.2.3 Multi-Dimensional--=00003D000026}Multi-Dimensional
$\hat{\chi}_{H}$ and $\tilde{\chi}_{H,N}$ \textendash{} The Commutative
Case}

Like in the one-dimensional case, we now introduce $\psi_{H}$ and
$\Psi_{H}$, which we will use in conjunction with $\chi_{H}$'s functional
equations (\ref{eq:Functional Equations for Chi_H over the rho-adics})
to extend what we have done to cover all \textbf{\emph{commutative}}
$p$-Hydra maps, regardless of the value of $\alpha_{H}\left(\mathbf{0}\right)$.
\begin{defn}[\textbf{Multi-Dimensional Little Psi-$H$ \& Big Psi-$H$}]
\ 

\vphantom{}

I. We define $\psi_{H}:\mathbb{N}_{0}^{r}\rightarrow\overline{\mathbb{Q}}^{d,d}$
(``Little Psi-$H$'') by: 
\begin{equation}
\psi_{H}\left(\mathbf{n}\right)\overset{\textrm{def}}{=}M_{H}\left(\mathbf{n}\right)\left(\mathbf{I}_{d}-H^{\prime}\left(\mathbf{0}\right)\right)^{-1}+\sum_{k=0}^{\lambda_{p}\left(\mathbf{n}\right)-1}\kappa_{H}\left(\left[\mathbf{n}\right]_{p^{k}}\right)\left(H^{\prime}\left(\mathbf{0}\right)\right)^{k}\label{eq:MD Definition of Little Psi_H}
\end{equation}

\vphantom{}

II. We define $\Psi_{H}:\mathbb{N}_{0}^{r}\rightarrow\overline{\mathbb{Q}}^{d}$
(``Big Psi-$H$'') by: 
\begin{align}
\Psi_{H}\left(\mathbf{n}\right) & \overset{\overline{\mathbb{Q}}^{d}}{=}M_{H}\left(\mathbf{n}\right)\sum_{\mathbf{j}>\mathbf{0}}^{p-1}\alpha_{H}\left(\frac{\mathbf{j}}{p}\right)\left(\mathbf{I}_{d}-H^{\prime}\left(\mathbf{0}\right)\right)^{-1}\gamma_{H}\left(\frac{\mathbf{j}}{p}\right)\label{eq:MD Definition of Big Psi_H}\\
 & \sum_{k=0}^{\lambda_{p}\left(\mathbf{n}\right)-1}\kappa_{H}\left(\left[\mathbf{n}\right]_{p^{k}}\right)\sum_{\mathbf{j}>\mathbf{0}}^{p-1}\mathcal{C}_{H}\left(\alpha_{H}\left(\frac{\mathbf{j}}{p}\right)\varepsilon_{k}\left(\mathbf{j}\mathbf{n}\right):\lambda_{p}\left(\left[\mathbf{n}\right]_{p^{k}}\right)\right)\left(H^{\prime}\left(\mathbf{0}\right)\right)^{k}\gamma_{H}\left(\frac{\mathbf{j}}{p}\right)
\end{align}

In both cases, the $k$-sums are defined to be $\mathbf{0}$ when
$\mathbf{n}=\mathbf{0}$. 
\end{defn}
\begin{lem}[\textbf{Rising-Continuability and Functional Equations for Multi-Dimensional
$\psi_{H}$ \& $\Psi_{H}$}]
\label{lem:MD Rising-continuability of Psi_Hs}\ 

\vphantom{}

I. $\psi_{H}$ is rising-continuable to a $d\times d$-matrix-valued
$\left(p,q_{H}\right)$-adic function $\psi_{H}:\mathbb{Z}_{p}^{r}\rightarrow\mathbb{Z}_{q_{H}}^{d,d}$
given by: 
\begin{equation}
\psi_{H}\left(\mathbf{z}\right)\overset{\mathcal{F}_{p,q_{H}}^{d,d}}{=}\sum_{n=0}^{\infty}\kappa_{H}\left(\left[\mathbf{z}\right]_{p^{n}}\right)\left(H^{\prime}\left(\mathbf{0}\right)\right)^{n},\textrm{ }\forall\mathbf{z}\in\left(\mathbb{Z}_{p}^{r}\right)^{\prime}\label{eq:MD Rising-continuation of Little Psi_H}
\end{equation}
Moreover, $\psi_{H}\left(\mathbf{z}\right)$ is the unique rising-continuous
$d\times d$-matrix-valued $\left(p,q_{H}\right)$-adic function satisfying
the system of functional equations: 
\begin{equation}
\psi_{H}\left(p\mathbf{z}+\mathbf{j}\right)\overset{\mathbb{C}_{q_{H}}^{d,d}}{=}H_{\mathbf{j}}^{\prime}\left(\mathbf{0}\right)\psi_{H}\left(\mathbf{z}\right)+\mathbf{I}_{d},\textrm{ }\forall\mathbf{z}\in\mathbb{Z}_{p}^{r}\textrm{ \& }\forall\mathbf{j}\in\mathbb{Z}^{r}/p\mathbb{Z}^{r}\label{eq:MD Little Psi_H functional equations}
\end{equation}
\index{functional equation!psi_{H}@$\psi_{H}$!multi-dimensional}

\vphantom{}

II. $\Psi_{H}$ is rising-continuable to a $d\times1$-vector-valued
$\left(p,q_{H}\right)$-adic function $\Psi_{H}:\mathbb{Z}_{p}^{r}\rightarrow\mathbb{C}_{q_{H}}^{d}$
given by:\index{functional equation!Psi_{H}@$\Psi_{H}$!multi-dimensional}
\begin{equation}
\Psi_{H}\left(\mathbf{z}\right)\overset{\mathcal{F}_{p,q_{H}}^{d}}{=}\sum_{n=0}^{\infty}\kappa_{H}\left(\left[\mathbf{z}\right]_{p^{n}}\right)\sum_{\mathbf{j}>\mathbf{0}}^{p-1}\mathcal{C}_{H}\left(\alpha_{H}\left(\frac{\mathbf{j}}{p}\right)\varepsilon_{n}\left(\mathbf{j}\mathbf{z}\right):\lambda_{p}\left(\left[\mathbf{z}\right]_{p^{n}}\right)\right)\left(H^{\prime}\left(\mathbf{0}\right)\right)^{n}\gamma_{H}\left(\frac{\mathbf{j}}{p}\right)\label{eq:MD Rising-continuation of Big Psi_H}
\end{equation}
for all $\mathbf{z}\in\left(\mathbb{Z}_{p}^{r}\right)^{\prime}$.
Moreover, $\Psi_{H}\left(\mathbf{z}\right)$ is the unique rising-continuous
$d\times1$-vector-valued $\left(p,q_{H}\right)$-adic function satisfying
the system of functional equations: 
\begin{equation}
\Psi_{H}\left(p\mathbf{z}+\mathbf{j}\right)\overset{\mathbb{C}_{q_{H}}^{d}}{=}H_{\mathbf{j}}^{\prime}\left(\mathbf{0}\right)\Psi_{H}\left(\mathbf{z}\right)+H_{\mathbf{j}}\left(\mathbf{0}\right)-\beta_{H}\left(\mathbf{0}\right),\textrm{ }\forall\mathbf{z}\in\mathbb{Z}_{p}^{r},\textrm{ }\forall\mathbf{j}\in\mathbb{Z}^{r}/p\mathbb{Z}^{r}\label{eq:MD Big Psi_H functional equations}
\end{equation}
\end{lem}
Proof: For both parts, we use (\ref{eq:MD Relation between truncations and functional equations, version 2})
and (\ref{eq:MD Kappa_H functional equations}) from \textbf{Lemma
\ref{lem:MD functional equations and truncation lemma}}. With these,
observe that $\kappa_{H}$'s functional equations: 
\begin{equation}
\kappa_{H}\left(p\mathbf{m}+\mathbf{j}\right)=\underbrace{\mathbf{D}_{\mathbf{j}}^{-1}\mathbf{A}_{\mathbf{j}}}_{H_{\mathbf{j}}^{\prime}\left(\mathbf{0}\right)}\kappa_{H}\left(\mathbf{m}\right)\underbrace{\mathbf{A}_{\mathbf{0}}^{-1}\mathbf{D}_{\mathbf{0}}}_{\left(H_{\mathbf{j}}\left(\mathbf{0}\right)\right)^{-1}}
\end{equation}
imply that: 
\begin{equation}
\kappa_{H}\left(\left[p\mathbf{m}+\mathbf{j}\right]_{p^{n}}\right)=\mathbf{D}_{\mathbf{j}}^{-1}\mathbf{A}_{\mathbf{j}}\kappa_{H}\left(\left[\mathbf{m}\right]_{p^{n-1}}\right)\mathbf{A}_{\mathbf{0}}^{-1}\mathbf{D}_{\mathbf{0}},\textrm{ }\forall n\in\mathbb{N}_{1},\textrm{ }\forall\mathbf{j}\in\mathbb{Z}^{r}/p\mathbb{Z}^{r},\textrm{ }\forall\mathbf{m}\in\mathbb{N}_{0}^{r}\label{eq:MD functional equation truncation Lemma applied to kappa_H}
\end{equation}
In this case, the function $\Phi_{\mathbf{j}}$ from (\ref{eq:MD Relation between truncations and functional equations, version 1})
is, in this case: 
\begin{equation}
\Phi_{\mathbf{j}}\left(\mathbf{m},\mathbf{X}\right)=\mathbf{D}_{\mathbf{j}}^{-1}\mathbf{A}_{\mathbf{j}}\mathbf{X}\mathbf{A}_{\mathbf{0}}^{-1}\mathbf{D}_{\mathbf{0}}
\end{equation}
where $\mathbf{X}$ is a $d\times d$ matrix.

The rising-continuability of $\psi_{H}$ and $\psi_{H}$ to the given
series follow by the givens on $H$, which guarantee that, for each
$\mathbf{z}\in\mathbb{Z}_{p}^{r}$, $M_{H}\left(\left[\mathbf{z}\right]_{p^{n}}\right)=\kappa_{H}\left(\left[\mathbf{z}\right]_{p^{n}}\right)\left(H^{\prime}\left(\mathbf{0}\right)\right)^{\lambda_{p}\left(n\right)}$
tends to $0$ in the standard $\left(p,q_{H}\right)$-adic frame as
$n\rightarrow\infty$, and convergence is then guaranteed by the same
arguments used for \textbf{Lemma \ref{lem:MD v_p A_H hat summation formulae}}
and equation (\ref{eq:F limit of MD Gamma_H A_H hat Fourier series}).
All that remains is to verify the functional equations; \textbf{Theorem
\ref{thm:rising-continuability of Generic H-type functional equations}}
from Subsection \ref{subsec:5.3.2. Interpolation-Revisited} then
guarantees the uniqueness of $\psi_{H}$ and $\Psi_{H}$ as rising-continuous
solutions of their respective systems of functional equations. In
what follows, recall that since we have assumed at the start of this
chapter that $H$ is contracting, \textbf{Proposition \ref{prop:MD Contracting H proposition}}
(page \pageref{prop:MD Contracting H proposition}) shows that the
matrix $\mathbf{I}_{d}-H^{\prime}\left(\mathbf{0}\right)$ is invertible.

\vphantom{}

I. We pull out the $k=0$ term from $\psi_{H}\left(p\mathbf{n}+\mathbf{j}\right)$:

\begin{align*}
\psi_{H}\left(p\mathbf{n}+\mathbf{j}\right) & \overset{\overline{\mathbb{Q}}^{d,d}}{=}M_{H}\left(p\mathbf{n}+\mathbf{j}\right)\left(\mathbf{I}_{d}-H^{\prime}\left(\mathbf{0}\right)\right)^{-1}\mathbf{I}_{d}\\
 & +\sum_{k=1}^{\lambda_{p}\left(p\mathbf{n}+\mathbf{j}\right)-1}\kappa_{H}\left(\left[p\mathbf{n}+\mathbf{j}\right]_{p^{k}}\right)\left(H^{\prime}\left(\mathbf{0}\right)\right)^{k}\\
 & =\mathbf{D}_{\mathbf{j}}^{-1}\mathbf{A}_{\mathbf{j}}M_{H}\left(\mathbf{n}\right)\left(\mathbf{I}_{d}-H^{\prime}\left(\mathbf{0}\right)\right)^{-1}+\mathbf{I}_{d}\\
 & +\sum_{k=1}^{\lambda_{p}\left(\mathbf{n}\right)}\mathbf{D}_{\mathbf{j}}^{-1}\mathbf{A}_{\mathbf{j}}\kappa_{H}\left(\left[\mathbf{n}\right]_{p^{k-1}}\right)\underbrace{\mathbf{A}_{\mathbf{0}}^{-1}\mathbf{D}_{\mathbf{0}}}_{\left(H^{\prime}\left(\mathbf{0}\right)\right)^{-1}}\left(H^{\prime}\left(\mathbf{0}\right)\right)^{k}\\
 & =\mathbf{D}_{\mathbf{j}}^{-1}\mathbf{A}_{\mathbf{j}}\underbrace{\left(M_{H}\left(\mathbf{n}\right)\left(\mathbf{I}_{d}-H^{\prime}\left(\mathbf{0}\right)\right)^{-1}+\sum_{k=0}^{\lambda_{p}\left(\mathbf{n}\right)-1}\kappa_{H}\left(\left[\mathbf{n}\right]_{p^{k}}\right)\left(H^{\prime}\left(\mathbf{0}\right)\right)^{k}\right)}_{\psi_{H}\left(\mathbf{n}\right)}+\mathbf{I}_{d}
\end{align*}
Consequently: 
\begin{equation}
\psi_{H}\left(p\mathbf{n}+\mathbf{j}\right)\overset{\overline{\mathbb{Q}}^{d,d}}{=}\underbrace{\mathbf{D}_{\mathbf{j}}^{-1}\mathbf{A}_{\mathbf{j}}}_{H_{\mathbf{j}}^{\prime}\left(\mathbf{0}\right)}\psi_{H}\left(\mathbf{n}\right)+\mathbf{I}_{d},\textrm{ }\forall\mathbf{n}\in\mathbb{N}_{0}^{r}\textrm{ \& }\forall j\in\mathbb{Z}^{r}/p\mathbb{Z}^{r}\label{eq:MD Little Psi_H functional equation on the integer tuples}
\end{equation}
then shows that (\ref{eq:MD Little Psi_H functional equation on the integer tuples})
extends to hold for the rising-continuation of $\psi_{H}$, and that
this rising-continuation is then the \emph{unique }$\left(p,q_{H}\right)$-adic
function satisfying (\ref{eq:MD Little Psi_H functional equations}).

Finally, letting $\mathbf{m}\in\mathbb{N}_{0}^{r}$, setting $\mathbf{z}=\mathbf{m}$,
the right-hand side of (\ref{eq:MD Rising-continuation of Little Psi_H})
becomes: 
\begin{align*}
\sum_{n=0}^{\infty}\kappa_{H}\left(\left[\mathbf{m}\right]_{p^{n}}\right)\left(H^{\prime}\left(\mathbf{0}\right)\right)^{n} & \overset{\mathbb{C}^{d}}{=}\sum_{n=0}^{\lambda_{p}\left(\mathbf{m}\right)-1}\kappa_{H}\left(\left[\mathbf{m}\right]_{p^{n}}\right)\left(H^{\prime}\left(\mathbf{0}\right)\right)^{n}\\
 & +\kappa_{H}\left(\mathbf{m}\right)\left(H^{\prime}\left(\mathbf{0}\right)\right)^{\lambda_{p}\left(\mathbf{m}\right)}\left(\mathbf{I}_{d}-H^{\prime}\left(\mathbf{0}\right)\right)\\
 & =M_{H}\left(\mathbf{n}\right)\left(\mathbf{I}_{d}-H^{\prime}\left(\mathbf{0}\right)\right)^{-1}+\sum_{n=0}^{\lambda_{p}\left(\mathbf{m}\right)-1}\kappa_{H}\left(\left[\mathbf{m}\right]_{p^{n}}\right)\left(H^{\prime}\left(\mathbf{0}\right)\right)^{n}\\
 & =\psi_{H}\left(\mathbf{m}\right)
\end{align*}
Hence, (\ref{eq:MD Rising-continuation of Little Psi_H}) converges
to $\psi_{H}$ in the standard frame.

\vphantom{}

II. Pulling out $k=0$ from (\ref{eq:MD Definition of Big Psi_H})
yields: 
\begin{align*}
\Psi_{H}\left(\mathbf{n}\right) & =M_{H}\left(\mathbf{n}\right)\sum_{\mathbf{j}>\mathbf{0}}^{p-1}\alpha_{H}\left(\frac{\mathbf{j}}{p}\right)\left(\mathbf{I}_{d}-H^{\prime}\left(\mathbf{0}\right)\right)^{-1}\gamma_{H}\left(\frac{\mathbf{j}}{p}\right)\\
 & +\kappa_{H}\left(\mathbf{0}\right)\sum_{\mathbf{j}>\mathbf{0}}^{p-1}\mathcal{C}_{H}\left(\alpha_{H}\left(\frac{\mathbf{j}}{p}\right)\varepsilon_{0}\left(\mathbf{j}\mathbf{n}\right):0\right)\gamma_{H}\left(\frac{\mathbf{j}}{p}\right)\\
 & +\sum_{k=1}^{\lambda_{p}\left(\mathbf{n}\right)-1}\kappa_{H}\left(\left[\mathbf{n}\right]_{p^{k}}\right)\sum_{\mathbf{j}>\mathbf{0}}^{p-1}\mathcal{C}_{H}\left(\alpha_{H}\left(\frac{\mathbf{j}}{p}\right)\varepsilon_{k}\left(\mathbf{j}\mathbf{n}\right):\lambda_{p}\left(\left[\mathbf{n}\right]_{p^{k}}\right)\right)\left(H^{\prime}\left(\mathbf{0}\right)\right)^{k}\gamma_{H}\left(\frac{\mathbf{j}}{p}\right)
\end{align*}
Here: 
\begin{equation}
\kappa_{H}\left(\mathbf{0}\right)\sum_{\mathbf{j}>\mathbf{0}}^{p-1}\mathcal{C}_{H}\left(\alpha_{H}\left(\frac{\mathbf{j}}{p}\right)\varepsilon_{0}\left(\mathbf{j}\mathbf{n}\right):\lambda_{p}\left(0\right)\right)\gamma_{H}\left(\frac{\mathbf{j}}{p}\right)
\end{equation}
is: 
\begin{equation}
\sum_{\mathbf{j}>\mathbf{0}}^{p-1}\alpha_{H}\left(\frac{\mathbf{j}}{p}\right)\gamma_{H}\left(\frac{\mathbf{j}}{p}\right)\varepsilon_{0}\left(\mathbf{j}\mathbf{n}\right)
\end{equation}
Now:

\begin{align*}
\sum_{\mathbf{j}=\mathbf{0}}^{p-1}\alpha_{H}\left(\frac{\mathbf{j}}{p}\right)\gamma_{H}\left(\frac{\mathbf{j}}{p}\right)\varepsilon_{0}\left(\mathbf{j}\mathbf{n}\right) & =\sum_{\mathbf{j}=\mathbf{0}}^{p-1}\beta_{H}\left(\frac{\mathbf{j}}{p}\right)\varepsilon_{0}\left(\mathbf{j}\mathbf{n}\right)\\
 & =\sum_{\mathbf{k}=\mathbf{0}}^{p-1}\mathbf{D}_{\mathbf{k}}^{-1}\mathbf{b}_{\mathbf{k}}\frac{1}{p^{r}}\sum_{\mathbf{j}=\mathbf{0}}^{p-1}e^{2\pi i\frac{\mathbf{j}\cdot\left(\mathbf{n}-\mathbf{k}\right)}{p}}\\
 & =\sum_{\mathbf{k}=\mathbf{0}}^{p-1}\mathbf{D}_{\mathbf{k}}^{-1}\mathbf{b}_{\mathbf{k}}\left[\mathbf{n}\overset{p}{\equiv}\mathbf{k}\right]\\
 & =\mathbf{D}_{\left[\mathbf{n}\right]_{p}}^{-1}\mathbf{b}_{\left[\mathbf{n}\right]_{p}}\\
 & =H_{\left[\mathbf{n}\right]_{p}}\left(\mathbf{0}\right)
\end{align*}
and so: 
\begin{equation}
\sum_{\mathbf{j}>\mathbf{0}}^{p-1}\alpha_{H}\left(\frac{\mathbf{j}}{p}\right)\gamma_{H}\left(\frac{\mathbf{j}}{p}\right)\varepsilon_{0}\left(\mathbf{j}\mathbf{n}\right)=H_{\left[\mathbf{n}\right]_{p}}\left(\mathbf{0}\right)-\beta_{H}\left(\mathbf{0}\right)
\end{equation}
Consequently: 
\begin{align*}
\Psi_{H}\left(\mathbf{n}\right) & =M_{H}\left(\mathbf{n}\right)\sum_{\mathbf{j}>\mathbf{0}}^{p-1}\alpha_{H}\left(\frac{\mathbf{j}}{p}\right)\left(\mathbf{I}_{d}-H^{\prime}\left(\mathbf{0}\right)\right)^{-1}\gamma_{H}\left(\frac{\mathbf{j}}{p}\right)+H_{\left[\mathbf{n}\right]_{p}}\left(\mathbf{0}\right)-\beta_{H}\left(\mathbf{0}\right)\\
 & +\sum_{k=1}^{\lambda_{p}\left(\mathbf{n}\right)-1}\kappa_{H}\left(\left[\mathbf{n}\right]_{p^{k}}\right)\sum_{\mathbf{j}>\mathbf{0}}^{p-1}\mathcal{C}_{H}\left(\alpha_{H}\left(\frac{\mathbf{j}}{p}\right)\varepsilon_{k}\left(\mathbf{j}\mathbf{n}\right):\lambda_{p}\left(\left[\mathbf{n}\right]_{p^{k}}\right)\right)\left(H^{\prime}\left(\mathbf{0}\right)\right)^{k}\gamma_{H}\left(\frac{\mathbf{j}}{p}\right)
\end{align*}

Now, replacing $\mathbf{n}$ with $p\mathbf{n}+\mathbf{i}$ (where
at least one of $\mathbf{n}$ and $\mathbf{i}$ is not the zero tuple),
we use (\ref{eq:MD functional equation truncation Lemma applied to kappa_H})
and the functional equations for $M_{H}$, $\varepsilon_{n}$, and
$\lambda_{p}$ (along with using (\ref{eq:MD Relation between truncations and functional equations, version 2})
for $\lambda_{p}\left(\left[\mathbf{n}\right]_{p^{k}}\right)$) and
the definition of $\mathcal{C}_{H}$ to obtain: 
\begin{align*}
 & \mathbf{D}_{\mathbf{i}}^{-1}\mathbf{A}_{\mathbf{i}}M_{H}\left(\mathbf{n}\right)\sum_{\mathbf{j}>\mathbf{0}}^{p-1}\alpha_{H}\left(\frac{\mathbf{j}}{p}\right)\left(\mathbf{I}_{d}-H^{\prime}\left(\mathbf{0}\right)\right)^{-1}\gamma_{H}\left(\frac{\mathbf{j}}{p}\right)+H_{\mathbf{i}}\left(\mathbf{0}\right)-\beta_{H}\left(\mathbf{0}\right)\\
 & +\mathbf{D}_{\mathbf{i}}^{-1}\mathbf{A}_{\mathbf{i}}\sum_{k=0}^{\lambda_{p}\left(\mathbf{n}\right)-1}\kappa_{H}\left(\left[\mathbf{n}\right]_{p^{k}}\right)\sum_{\mathbf{j}>\mathbf{0}}^{p-1}\mathcal{C}_{H}\left(\alpha_{H}\left(\frac{\mathbf{j}}{p}\right)\varepsilon_{k}\left(\mathbf{j}\mathbf{n}\right):\lambda_{p}\left(\left[\mathbf{n}\right]_{p^{k}}\right)\right)\left(H^{\prime}\left(\mathbf{0}\right)\right)^{k}\gamma_{H}\left(\frac{\mathbf{j}}{p}\right)
\end{align*}
as the formula for $\Psi_{H}\left(p\mathbf{n}+\mathbf{i}\right)$.
This is: 
\[
\Psi_{H}\left(p\mathbf{n}+\mathbf{i}\right)=\underbrace{\mathbf{D}_{\mathbf{i}}^{-1}\mathbf{A}_{\mathbf{i}}}_{H_{\mathbf{i}}^{\prime}\left(\mathbf{0}\right)}\Psi_{H}\left(p\mathbf{n}+\mathbf{i}\right)+H_{\mathbf{i}}\left(\mathbf{0}\right)-\beta_{H}\left(\mathbf{0}\right)
\]

Finally, letting $\mathbf{X}$ denote $H^{\prime}\left(\mathbf{0}\right)$
and letting $\mathbf{m}\in\mathbb{N}_{0}^{r}$, we have that for $\mathbf{z}=\mathbf{m}$,
the right-hand side of (\ref{eq:MD Rising-continuation of Big Psi_H})
becomes: 
\begin{align*}
\sum_{n=0}^{\lambda_{p}\left(\mathbf{m}\right)-1}\kappa_{H}\left(\left[\mathbf{m}\right]_{p^{n}}\right)\sum_{\mathbf{j}>\mathbf{0}}^{p-1}\mathcal{C}_{H}\left(\alpha_{H}\left(\frac{\mathbf{j}}{p}\right)\varepsilon_{n}\left(\mathbf{j}\mathbf{m}\right):\lambda_{p}\left(\left[\mathbf{m}\right]_{p^{n}}\right)\right)\mathbf{X}^{n}\gamma_{H}\left(\frac{\mathbf{j}}{p}\right)\\
+\sum_{n=\lambda_{p}\left(\mathbf{m}\right)}^{\infty}\kappa_{H}\left(\mathbf{m}\right)\sum_{\mathbf{j}>\mathbf{0}}^{p-1}\mathcal{C}_{H}\left(\alpha_{H}\left(\frac{\mathbf{j}}{p}\right)\varepsilon_{n}\left(\mathbf{j}\mathbf{m}\right):\lambda_{p}\left(\mathbf{m}\right)\right)\mathbf{X}^{n}\gamma_{H}\left(\frac{\mathbf{j}}{p}\right)
\end{align*}
Here: 
\[
\varepsilon_{n}\left(\mathbf{j}\mathbf{m}\right)=e^{\frac{2\pi i}{p^{n+1}}\left(\left[\mathbf{j}\mathbf{m}\right]_{p^{n+1}}-\left[\mathbf{j}\mathbf{m}\right]_{p^{n}}\right)}=e^{\frac{2\pi i}{p^{n+1}}\left(\mathbf{j}\mathbf{m}-\mathbf{j}\mathbf{m}\right)}=1,\textrm{ }\forall n\geq\lambda_{p}\left(\mathbf{m}\right)
\]
and so: 
\begin{align*}
\sum_{n=0}^{\lambda_{p}\left(\mathbf{m}\right)-1}\kappa_{H}\left(\left[\mathbf{m}\right]_{p^{n}}\right)\sum_{\mathbf{j}>\mathbf{0}}^{p-1}\mathcal{C}_{H}\left(\alpha_{H}\left(\frac{\mathbf{j}}{p}\right)\varepsilon_{n}\left(\mathbf{j}\mathbf{m}\right):\lambda_{p}\left(\left[\mathbf{m}\right]_{p^{n}}\right)\right)\mathbf{X}^{n}\gamma_{H}\left(\frac{\mathbf{j}}{p}\right)\\
+\sum_{n=\lambda_{p}\left(\mathbf{m}\right)}^{\infty}\kappa_{H}\left(\mathbf{m}\right)\sum_{\mathbf{j}>\mathbf{0}}^{p-1}\mathbf{X}^{\lambda_{p}\left(\mathbf{m}\right)}\alpha_{H}\left(\frac{\mathbf{j}}{p}\right)\mathbf{X}^{-\lambda_{p}\left(\mathbf{m}\right)}\mathbf{X}^{n}\gamma_{H}\left(\frac{\mathbf{j}}{p}\right)
\end{align*}
Summed in the topology of $\mathbb{C}$, this becomes: 
\begin{align*}
\sum_{n=0}^{\lambda_{p}\left(\mathbf{m}\right)-1}\kappa_{H}\left(\left[\mathbf{m}\right]_{p^{n}}\right)\sum_{\mathbf{j}>\mathbf{0}}^{p-1}\mathcal{C}_{H}\left(\alpha_{H}\left(\frac{\mathbf{j}}{p}\right)\varepsilon_{n}\left(\mathbf{j}\mathbf{m}\right):\lambda_{p}\left(\left[\mathbf{m}\right]_{p^{n}}\right)\right)\mathbf{X}^{n}\gamma_{H}\left(\frac{\mathbf{j}}{p}\right)\\
+\kappa_{H}\left(\mathbf{m}\right)\sum_{\mathbf{j}>\mathbf{0}}^{p-1}\mathbf{X}^{\lambda_{p}\left(\mathbf{m}\right)}\alpha_{H}\left(\frac{\mathbf{j}}{p}\right)\mathbf{X}^{-\lambda_{p}\left(\mathbf{m}\right)}\mathbf{X}^{\lambda_{p}\left(\mathbf{m}\right)}\left(\mathbf{I}_{d}-\mathbf{X}\right)^{-1}\gamma_{H}\left(\frac{\mathbf{j}}{p}\right)
\end{align*}
which simplifies to: 
\begin{align*}
\sum_{n=0}^{\lambda_{p}\left(\mathbf{m}\right)-1}\kappa_{H}\left(\left[\mathbf{m}\right]_{p^{n}}\right)\sum_{\mathbf{j}>\mathbf{0}}^{p-1}\mathcal{C}_{H}\left(\alpha_{H}\left(\frac{\mathbf{j}}{p}\right)\varepsilon_{n}\left(\mathbf{j}\mathbf{m}\right):\lambda_{p}\left(\left[\mathbf{m}\right]_{p^{n}}\right)\right)\mathbf{X}^{n}\gamma_{H}\left(\frac{\mathbf{j}}{p}\right)\\
+M_{H}\left(\mathbf{m}\right)\sum_{\mathbf{j}>\mathbf{0}}^{p-1}\alpha_{H}\left(\frac{\mathbf{j}}{p}\right)\left(\mathbf{I}_{d}-\mathbf{X}\right)^{-1}\gamma_{H}\left(\frac{\mathbf{j}}{p}\right)
\end{align*}
which is precisely $\Psi_{H}\left(\mathbf{m}\right)$ as given in
(\ref{eq:MD Definition of Big Psi_H}). Hence, (\ref{eq:MD Rising-continuation of Big Psi_H})
converges to $\Psi_{H}$ in the standard frame.

Q.E.D. 
\begin{prop}[\textbf{Quasi-Integrability of $\psi_{H}$ and $\Psi_{H}$}]
\label{prop:quasi-integrability of MD Psis}Let $H$ be commutative.
Then:

\vphantom{}

I. $\psi_{H}$ is quasi-integrable with respect to the standard $\left(p,q_{H}\right)$-adic
frame whenever either $\alpha_{H}\left(\mathbf{0}\right)=\mathbf{I}_{d}$
or the matrix $\mathbf{I}_{d}-\alpha_{H}\left(\mathbf{0}\right)$
is invertible. For these conditions, the function $\hat{\psi}_{H}:\hat{\mathbb{Z}}_{p}^{r}\rightarrow\overline{\mathbb{Q}}^{d,d}$
defined by: 
\begin{equation}
\hat{\psi}_{H}\left(\mathbf{t}\right)\overset{\textrm{def}}{=}\begin{cases}
\begin{cases}
\mathbf{O}_{d} & \textrm{if }\mathbf{t}=\mathbf{0}\\
v_{p}\left(\mathbf{t}\right)\hat{A}_{H}\left(\mathbf{t}\right) & \textrm{if }\mathbf{t}\neq\mathbf{0}
\end{cases} & \textrm{if }\alpha_{H}\left(\mathbf{0}\right)=\mathbf{I}_{d}\\
\hat{A}_{H}\left(\mathbf{t}\right)\left(\mathbf{I}_{d}-\alpha_{H}\left(\mathbf{0}\right)\right)^{-1} & \textrm{if }\mathbf{I}_{d}-\alpha_{H}\left(\mathbf{0}\right)\textrm{ is invertible}
\end{cases},\textrm{ }\forall\mathbf{t}\in\hat{\mathbb{Z}}_{p}^{r}\label{eq:MD Fourier Transform of Little Psi_H, commutative}
\end{equation}
is then a Fourier transform of $\psi_{H}$.

Hence: 
\begin{equation}
\hat{\psi}_{H,N}\left(\mathbf{z}\right)\overset{\overline{\mathbb{Q}}^{d,d}}{=}-N\kappa_{H}\left(\left[\mathbf{z}\right]_{p^{N}}\right)\left(H^{\prime}\left(\mathbf{0}\right)\right)^{N}+\sum_{n=0}^{N-1}\kappa_{H}\left(\left[\mathbf{z}\right]_{p^{n}}\right)\left(H^{\prime}\left(\mathbf{0}\right)\right)^{n}\label{eq:MD Little Psi H N twiddle when alpha is 1, commutative}
\end{equation}
when $\alpha_{H}\left(\mathbf{0}\right)=\mathbf{I}_{d}$ and:

\begin{equation}
\tilde{\psi}_{H,N}\left(\mathbf{z}\right)\overset{\overline{\mathbb{Q}}^{d,d}}{=}\kappa_{H}\left(\left[\mathbf{z}\right]_{p^{N}}\right)\left(H^{\prime}\left(\mathbf{0}\right)\right)^{N}\left(\mathbf{I}_{d}-\alpha_{H}\left(\mathbf{0}\right)\right)^{-1}+\sum_{n=0}^{N-1}\kappa_{H}\left(\left[\mathbf{z}\right]_{p^{n}}\right)\left(H^{\prime}\left(\mathbf{0}\right)\right)^{n}\label{eq:MD Little Psi H N twiddle when 1 minus alpha is invertible, commutative}
\end{equation}
when $\mathbf{I}_{d}-\alpha_{H}\left(\mathbf{0}\right)$ is invertible.

\vphantom{}

II. $\Psi_{H}$ is quasi-integrable with respect to the standard $\left(p,q_{H}\right)$-adic
frame for all $\alpha_{H}\left(\mathbf{0}\right)$, and the function
$\hat{\Psi}_{H}:\hat{\mathbb{Z}}_{p}^{r}\rightarrow\overline{\mathbb{Q}}^{d}$
defined by: 
\begin{equation}
\hat{\Psi}_{H}\left(t\right)\overset{\textrm{def}}{=}\begin{cases}
0 & \textrm{if }\mathbf{t}=\mathbf{0}\\
\hat{A}_{H}\left(\mathbf{t}\right)\gamma_{H}\left(\frac{\mathbf{t}\left|\mathbf{t}\right|_{p}}{p}\right) & \textrm{if }\mathbf{t}\neq\mathbf{0}
\end{cases},\textrm{ }\forall\mathbf{t}\in\hat{\mathbb{Z}}_{p}^{r}\label{eq:MD Fourier Transform of Big Psi_H}
\end{equation}
is a Fourier transform of $\Psi_{H}$. Hence: 
\begin{equation}
\tilde{\Psi}_{H,N}\left(\mathbf{z}\right)\overset{\overline{\mathbb{Q}}^{d}}{=}\sum_{n=0}^{N-1}\kappa_{H}\left(\left[\mathbf{z}\right]_{p^{n}}\right)\sum_{\mathbf{j}>\mathbf{0}}^{p-1}\mathcal{C}_{H}\left(\alpha_{H}\left(\frac{\mathbf{j}}{p}\right)\varepsilon_{n}\left(\mathbf{j}\mathbf{z}\right):\lambda_{p}\left(\left[\mathbf{z}\right]_{p^{n}}\right)\right)\left(H^{\prime}\left(\mathbf{0}\right)\right)^{n}\gamma_{H}\left(\frac{\mathbf{j}}{p}\right)\label{eq:MD Big Psi H N twiddle}
\end{equation}
\end{prop}
Proof:

I. When $\alpha_{H}\left(\mathbf{0}\right)=\mathbf{I}_{d}$, (\ref{eq:MD Fourier Transform of Little Psi_H, commutative})
follows from \textbf{Lemma \ref{lem:MD v_p A_H hat summation formulae}},
and hence, the $\alpha_{H}\left(\mathbf{0}\right)=\mathbf{I}_{d}$
case of (\ref{eq:MD Fourier Transform of Little Psi_H, commutative})
is then indeed a Fourier transform of $\psi_{H}$, thus proving the
$\mathcal{F}_{p,q_{H}}$-quasi-integrability of $\psi_{H}$ when $\alpha_{H}\left(\mathbf{0}\right)=\mathbf{I}_{d}$.

When $\mathbf{I}_{d}-\alpha_{H}\left(\mathbf{0}\right)$ is invertible,
comparing (\ref{eq:MD Fourier Limit of A_H,N twiddle in standard frame, commutative}):
\begin{align*}
\tilde{A}_{H}\left(\mathbf{z}\right) & \overset{\mathcal{F}_{p,q_{H}}^{d,d}}{=}\sum_{n=0}^{\infty}\kappa_{H}\left(\left[\mathbf{z}\right]_{p^{n}}\right)\left(H^{\prime}\left(\mathbf{0}\right)\right)^{n}\left(\mathbf{I}_{d}-\alpha_{H}\left(\mathbf{0}\right)\right)
\end{align*}
and (\ref{eq:MD Rising-continuation of Little Psi_H}): 
\[
\psi_{H}\left(\mathbf{z}\right)\overset{\mathcal{F}_{p,q_{H}}^{d,d}}{=}\sum_{n=0}^{\infty}\kappa_{H}\left(\left[\mathbf{z}\right]_{p^{n}}\right)\left(H^{\prime}\left(\mathbf{0}\right)\right)^{n}
\]
we observe that: 
\[
\tilde{A}_{H}\left(\mathbf{z}\right)=\psi_{H}\left(\mathbf{z}\right)\left(\mathbf{I}_{d}-\alpha_{H}\left(\mathbf{0}\right)\right)
\]
Since \textbf{Theorem \ref{thm:MD properties of A_H hat}} shows that
$\tilde{A}_{H}\left(\mathbf{z}\right)$ is $\mathcal{F}_{p,q_{H}}$-quasi-integrable
when $\mathbf{I}_{d}-\alpha_{H}\left(\mathbf{0}\right)$ is invertible
(which necessarily forces $\alpha_{H}\left(\mathbf{0}\right)\neq\mathbf{I}_{d}$),
it then follows that: 
\[
\psi_{H}\left(\mathbf{z}\right)=\tilde{A}_{H}\left(\mathbf{z}\right)\left(\mathbf{I}_{d}-\alpha_{H}\left(\mathbf{0}\right)\right)^{-1}
\]
and hence, that: 
\[
\hat{\psi}_{H}\left(\mathbf{t}\right)\overset{\textrm{def}}{=}\hat{A}_{H}\left(\mathbf{t}\right)\left(\mathbf{I}_{d}-\alpha_{H}\left(\mathbf{0}\right)\right)^{-1}
\]
is a Fourier transform of $\psi_{H}$. This proves that $\psi_{H}$
is quasi-integrable with respect to the standard frame when $\mathbf{I}_{d}-\alpha_{H}\left(\mathbf{0}\right)$
is invertible, and that (\ref{eq:MD Fourier Transform of Little Psi_H, commutative})
is then a Fourier transform of $\psi_{H}$ in this case.

Finally, since $H$ is commutative, for the case $\alpha_{H}\left(\mathbf{0}\right)=\mathbf{I}_{d}$,
(\ref{eq:MD Fourier sum of A_H hat v_rho, commutative}) becomes:
\begin{align*}
\sum_{0<\left\Vert \mathbf{t}\right\Vert _{p}\leq p^{N}}v_{p}\left(\mathbf{t}\right)\hat{A}_{H}\left(\mathbf{t}\right)e^{2\pi i\left\{ \mathbf{t}\mathbf{z}\right\} _{p}} & \overset{\overline{\mathbb{Q}}^{d,d}}{=}-N\kappa_{H}\left(\left[\mathbf{z}\right]_{p^{N}}\right)\left(H^{\prime}\left(\mathbf{0}\right)\right)^{N}\\
 & +\sum_{n=0}^{N-1}\kappa_{H}\left(\left[\mathbf{z}\right]_{p^{n}}\right)\left(H^{\prime}\left(\mathbf{0}\right)\right)^{n}
\end{align*}
the left-hand side of which is exactly (\ref{eq:MD Little Psi H N twiddle when alpha is 1, commutative}).
On the other hand, for the case where $\mathbf{I}_{d}-\alpha_{H}\left(\mathbf{0}\right)$
is invertible, the commutativity of $H$ tells us that (\ref{eq:MD Convolution of dA_H and D_N when H is commutative})
holds: 
\[
\tilde{A}_{H,N}\left(\mathbf{z}\right)=\kappa_{H}\left(\left[\mathbf{z}\right]_{p^{N}}\right)\left(H^{\prime}\left(\mathbf{0}\right)\right)^{N}+\sum_{n=0}^{N-1}\kappa_{H}\left(\left[\mathbf{z}\right]_{p^{n}}\right)\left(H^{\prime}\left(\mathbf{0}\right)\right)^{n}\left(\mathbf{I}_{d}-\alpha_{H}\left(\mathbf{0}\right)\right)
\]
Right-multiplying by $\mathbf{I}_{d}-\alpha_{H}\left(\mathbf{0}\right)$
then yields (\ref{eq:MD Little Psi H N twiddle when 1 minus alpha is invertible, commutative}).

\vphantom{}

II. (\ref{eq:MD Fourier Transform of Big Psi_H}) is exactly what
we proved in \textbf{Lemma \ref{lem:MD Rising-continuability of Psi_Hs}}.
(\ref{eq:MD Big Psi H N twiddle}) is precisely what we proved in
(\ref{eq:MD Gamma formula}).

Q.E.D.

\vphantom{}

Like in the one-dimensional case, we can bootstrap a Fourier transform
for $\chi_{H}$ for the case where $H$ is commutative. However, this
follows from a more general result which tells us how to express $\chi_{H}$
in terms of $\psi_{H}$ and $\Psi_{H}$:
\begin{thm}[\textbf{$\mathcal{F}$-Series for Multi-Dimensional $\chi_{H}$}]
\label{thm:MD F-series for Chi_H}Let $H$ be a $p$-smooth $d$-dimensional
depth $r$ Hydra map dimension $d$ which is contracting, semi-basic,
non-singular, and fixes $\mathbf{0}$. Regardless of:

\vphantom{}

I. Whether or not $\alpha_{H}\mathbf{\left(0\right)}=\mathbf{I}_{d}$;

\vphantom{}

II. Whether or not $\mathbf{I}_{d}-\alpha_{H}\left(\mathbf{0}\right)$
is invertible;

\vphantom{}

III. Whether or not $H$ is commutative;

\vphantom{}

the numen $\chi_{H}$ admits the $\mathcal{F}$-series representation:
\begin{equation}
\chi_{H}\left(\mathbf{z}\right)\overset{\mathcal{F}_{p,q_{H}}^{d}}{=}\psi_{H}\left(\mathbf{z}\right)\beta_{H}\left(\mathbf{0}\right)+\Psi_{H}\left(\mathbf{z}\right),\textrm{ }\forall\mathbf{z}\in\mathbb{Z}_{p}^{r}\label{eq:MD Chi_H in terms of Little Psi_H and Big Psi_H}
\end{equation}
\end{thm}
Proof: Let $f\left(\mathbf{z}\right)$ denote the function $\psi_{H}\left(\mathbf{z}\right)\beta_{H}\left(\mathbf{0}\right)+\Psi_{H}\left(\mathbf{z}\right)$.
Then, using the functional equations for $\psi_{H}$ and $\Psi_{H}$:
\begin{equation}
\psi_{H}\left(p\mathbf{z}+\mathbf{j}\right)=H_{\mathbf{j}}^{\prime}\left(\mathbf{0}\right)\psi_{H}\left(\mathbf{z}\right)+\mathbf{I}_{d}
\end{equation}
\begin{equation}
\Psi_{H}\left(p\mathbf{z}+\mathbf{j}\right)=H_{\mathbf{j}}^{\prime}\left(\mathbf{0}\right)\Psi_{H}\left(\mathbf{z}\right)+H_{\mathbf{j}}\left(\mathbf{0}\right)-\beta_{H}\left(\mathbf{0}\right)
\end{equation}
we have that: 
\begin{align*}
f\left(p\mathbf{z}+\mathbf{j}\right) & =\psi_{H}\left(p\mathbf{z}+\mathbf{j}\right)\beta_{H}\left(\mathbf{0}\right)+\Psi_{H}\left(p\mathbf{z}+\mathbf{j}\right)\\
 & =\left(H_{\mathbf{j}}^{\prime}\left(\mathbf{0}\right)\psi_{H}\left(\mathbf{z}\right)+\mathbf{I}_{d}\right)\beta_{H}\left(\mathbf{0}\right)+\left(H_{\mathbf{j}}^{\prime}\left(\mathbf{0}\right)\Psi_{H}\left(\mathbf{z}\right)+H_{\mathbf{j}}\left(\mathbf{0}\right)-\beta_{H}\left(\mathbf{0}\right)\right)\\
 & =H_{\mathbf{j}}^{\prime}\left(\mathbf{0}\right)\psi_{H}\left(\mathbf{z}\right)\beta_{H}\left(\mathbf{0}\right)+H_{\mathbf{j}}^{\prime}\left(\mathbf{0}\right)\Psi_{H}\left(\mathbf{z}\right)+H_{\mathbf{j}}\left(\mathbf{0}\right)\\
 & =H_{\mathbf{j}}^{\prime}\left(\mathbf{0}\right)\left(\psi_{H}\left(\mathbf{z}\right)\beta_{H}\left(\mathbf{0}\right)+\Psi_{H}\left(\mathbf{z}\right)\right)+H_{\mathbf{j}}\left(\mathbf{0}\right)\\
 & =H_{\mathbf{j}}^{\prime}\left(\mathbf{0}\right)f\left(\mathbf{z}\right)+H_{\mathbf{j}}\left(\mathbf{0}\right)\\
 & =H_{\mathbf{j}}\left(f\left(\mathbf{z}\right)\right)
\end{align*}
Thus, $f$ satisfies $\chi_{H}$'s functional equation. Since $\chi_{H}$
is the unique $\left(p,q_{H}\right)$-adic function satisfying its
functional equation asa described in \textbf{Lemma \ref{lem:MD rising-continuation and functional equations of Chi_H}},
this forces $f=\chi_{H}$.

Q.E.D. 
\begin{cor}[\textbf{Quasi-Integrability of $\chi_{H}$ When $\mathbf{I}_{d}-\alpha_{H}\left(\mathbf{0}\right)$
is Invertible}]
\label{cor:MD Quasi-integrability of Chi_H for I minus alpha invertible}Let
$H$ be as given in \textbf{\emph{Theorem \ref{thm:MD F-series for Chi_H}}}.
In addition, suppose $H$ is commutative, and that $\mathbf{I}_{d}-\alpha_{H}\left(\mathbf{0}\right)$
is invertible. Then, $\chi_{H}$ is quasi-integrable with respect
to the standard $\left(p,q_{H}\right)$-adic frame, and the function
$\hat{\chi}_{H}:\hat{\mathbb{Z}}_{p}^{r}\rightarrow\overline{\mathbb{Q}}^{d}$
defined by: 
\begin{equation}
\hat{\chi}_{H}\left(\mathbf{t}\right)\overset{\textrm{def}}{=}\hat{A}_{H}\left(\mathbf{t}\right)\left(\mathbf{I}_{d}-\alpha_{H}\left(\mathbf{0}\right)\right)^{-1}\beta_{H}\left(\mathbf{0}\right)+\begin{cases}
\mathbf{0} & \textrm{if }\mathbf{t}=\mathbf{0}\\
\hat{A}_{H}\left(\mathbf{t}\right)\gamma_{H}\left(\frac{\mathbf{t}\left|\mathbf{t}\right|_{p}}{p}\right) & \textrm{else}
\end{cases},\textrm{ }\forall\mathbf{t}\in\hat{\mathbb{Z}}_{p}^{r}\label{eq:MD Chi_H hat for a contracting semi-basic commutative P Hydra map where alpha minus I is invertible}
\end{equation}
is a Fourier transform of $\chi_{H}$. 
\end{cor}
Proof: The given hypotheses tell us that both $\psi_{H}$ and $\Psi_{H}$
are quasi-integrable. As such, by the linearity of the Fourier transform,
(\ref{eq:MD Chi_H in terms of Little Psi_H and Big Psi_H}) becomes:
\begin{equation}
\hat{\chi}_{H}\left(\mathbf{t}\right)=\hat{\psi}_{H}\left(\mathbf{t}\right)\beta_{H}\left(\mathbf{0}\right)+\hat{\Psi}_{H}\left(\mathbf{t}\right)
\end{equation}
Using (\ref{eq:MD Fourier Transform of Little Psi_H, commutative})
and (\ref{eq:MD Fourier Transform of Big Psi_H}) from \textbf{Proposition
\ref{prop:quasi-integrability of MD Psis}} to express the right-hand
side of this yields (\ref{eq:MD Chi_H hat for a contracting semi-basic commutative P Hydra map where alpha minus I is invertible}).

Q.E.D. 
\begin{cor}[\textbf{Formulae for $\tilde{\chi}_{H,N}$}]
Let $H$ and $\hat{\chi}_{H}$ be as given in \textbf{\emph{Corollary
\ref{cor:MD Quasi-integrability of Chi_H for I minus alpha invertible}}}.

\vphantom{}

I. If $p=2$, then: 
\begin{align}
\tilde{\chi}_{H,N}\left(\mathbf{z}\right) & \overset{\overline{\mathbb{Q}}^{d}}{=}-\gamma_{H}\left(\mathbf{\frac{1}{2}}\right)+\kappa_{H}\left(\left[\mathbf{z}\right]_{2^{N}}\right)\left(H^{\prime}\left(\mathbf{0}\right)\right)^{N}\left(\left(\mathbf{I}_{d}-\alpha_{H}\left(\mathbf{0}\right)\right)^{-1}\beta_{H}\left(\mathbf{0}\right)+\gamma_{H}\left(\mathbf{\frac{1}{2}}\right)\right)\label{eq:MD Chi_H,N twiddle formula when P is 2 and alpha minus 1 is invertible, commutative}\\
 & +\sum_{n=0}^{N-1}\kappa_{H}\left(\left[\mathbf{z}\right]_{2^{n}}\right)\left(H^{\prime}\left(\mathbf{0}\right)\right)^{n}\left(\beta_{H}\left(\mathbf{0}\right)+\left(\mathbf{I}_{d}-\alpha_{H}\left(\mathbf{0}\right)\right)\gamma_{H}\left(\mathbf{\frac{1}{2}}\right)\right)\nonumber 
\end{align}
 Hence, passing to the $\mathcal{F}_{2,q_{H}}$-limit as $N\rightarrow\infty$:
\begin{equation}
\chi_{H}\left(\mathbf{z}\right)\overset{\mathcal{F}_{2,q_{H}}^{d}}{=}-\gamma_{H}\left(\mathbf{\frac{1}{2}}\right)+\sum_{n=0}^{\infty}\kappa_{H}\left(\left[\mathbf{z}\right]_{2^{n}}\right)\left(H^{\prime}\left(\mathbf{0}\right)\right)^{n}\left(\beta_{H}\left(\mathbf{0}\right)+\left(\mathbf{I}_{d}-\alpha_{H}\left(\mathbf{0}\right)\right)\gamma_{H}\left(\mathbf{\frac{1}{2}}\right)\right)\label{eq:MD Chi_H explicit formula when P is 2 and alpha minus 1 is invertible, commutative}
\end{equation}

\vphantom{}

II. If $p$ is odd: 
\begin{align}
\tilde{\chi}_{H,N}\left(\mathbf{z}\right) & \overset{\overline{\mathbb{Q}}^{d}}{=}\kappa_{H}\left(\left[\mathbf{z}\right]_{p^{N}}\right)\left(H^{\prime}\left(\mathbf{0}\right)\right)^{N}+\sum_{n=0}^{N-1}\kappa_{H}\left(\left[\mathbf{z}\right]_{p^{n}}\right)\left(H^{\prime}\left(\mathbf{0}\right)\right)^{n}\left(\mathbf{I}_{d}-\alpha_{H}\left(\mathbf{0}\right)\right)\label{eq:MD Chi_H,N twiddle formula for arbitrary P and alpha minus 1 is invertible, commutative}\\
 & +\sum_{n=0}^{N-1}\kappa_{H}\left(\left[\mathbf{z}\right]_{p^{n}}\right)\sum_{\mathbf{j}>\mathbf{0}}^{p-1}\mathcal{C}_{H}\left(\alpha_{H}\left(\frac{\mathbf{j}}{p}\right)\varepsilon_{n}\left(\mathbf{j}\mathbf{z}\right):\lambda_{p}\left(\left[\mathbf{z}\right]_{p^{n}}\right)\right)\left(H^{\prime}\left(\mathbf{0}\right)\right)^{n}\gamma_{H}\left(\frac{\mathbf{j}}{p}\right)\nonumber 
\end{align}
Hence, passing to the $\mathcal{F}_{p,q_{H}}$-limit as $N\rightarrow\infty$:
\begin{align}
\chi_{H}\left(\mathbf{z}\right) & \overset{\mathcal{F}_{p,q_{H}}^{d}}{=}\sum_{n=0}^{\infty}\kappa_{H}\left(\left[\mathbf{z}\right]_{p^{n}}\right)\left(H^{\prime}\left(\mathbf{0}\right)\right)^{n}\left(\mathbf{I}_{d}-\alpha_{H}\left(\mathbf{0}\right)\right)\label{eq:MD Chi_H explicit formula for arbitrary P and alpha minus 1 is invertible, commutative}\\
 & +\sum_{n=0}^{\infty}\kappa_{H}\left(\left[\mathbf{z}\right]_{p^{n}}\right)\sum_{\mathbf{j}>\mathbf{0}}^{p-1}\mathcal{C}_{H}\left(\alpha_{H}\left(\frac{\mathbf{j}}{p}\right)\varepsilon_{n}\left(\mathbf{j}\mathbf{z}\right):\lambda_{p}\left(\left[\mathbf{z}\right]_{p^{n}}\right)\right)\left(H^{\prime}\left(\mathbf{0}\right)\right)^{n}\gamma_{H}\left(\frac{\mathbf{j}}{p}\right)\nonumber 
\end{align}
\end{cor}
Proof:

I. Suppose $p=2$. Then, by \textbf{Corollary \ref{cor:MD Quasi-integrability of Chi_H for I minus alpha invertible}},
we have that: 
\begin{equation}
\hat{\chi}_{H}\left(\mathbf{t}\right)=\hat{A}_{H}\left(\mathbf{t}\right)\left(\mathbf{I}_{d}-\alpha_{H}\left(\mathbf{0}\right)\right)^{-1}\beta_{H}\left(\mathbf{0}\right)+\begin{cases}
\mathbf{0} & \textrm{if }\mathbf{t}=\mathbf{0}\\
\hat{A}_{H}\left(\mathbf{t}\right)\gamma_{H}\left(\mathbf{\frac{1}{2}}\right) & \textrm{else}
\end{cases}
\end{equation}
Hence, by (\ref{eq:MD Convolution of dA_H and D_N when H is commutative})
from\textbf{ Theorem \ref{thm:MD properties of A_H hat}}: 
\begin{align*}
\tilde{\chi}_{H,N}\left(\mathbf{z}\right) & =\sum_{\left\Vert \mathbf{t}\right\Vert _{2}\leq2^{N}}\hat{A}_{H}\left(\mathbf{t}\right)e^{2\pi i\left\{ \mathbf{t}\mathbf{z}\right\} _{2}}\left(\mathbf{I}_{d}-\alpha_{H}\left(\mathbf{0}\right)\right)^{-1}\beta_{H}\left(\mathbf{0}\right)\\
 & +\sum_{0<\left\Vert \mathbf{t}\right\Vert _{2}\leq2^{N}}\hat{A}_{H}\left(\mathbf{t}\right)e^{2\pi i\left\{ \mathbf{t}\mathbf{z}\right\} _{2}}\gamma_{H}\left(\mathbf{\frac{1}{2}}\right)\\
 & =-\gamma_{H}\left(\mathbf{\frac{1}{2}}\right)+\kappa_{H}\left(\left[\mathbf{z}\right]_{2^{N}}\right)\left(H^{\prime}\left(\mathbf{0}\right)\right)^{N}\left(\left(\mathbf{I}_{d}-\alpha_{H}\left(\mathbf{0}\right)\right)^{-1}\beta_{H}\left(\mathbf{0}\right)+\gamma_{H}\left(\frac{1}{2}\right)\right)\\
 & +\sum_{n=0}^{N-1}\kappa_{H}\left(\left[\mathbf{z}\right]_{2^{n}}\right)\left(H^{\prime}\left(\mathbf{0}\right)\right)^{n}\left(\beta_{H}\left(\mathbf{0}\right)+\left(\mathbf{I}_{d}-\alpha_{H}\left(\mathbf{0}\right)\right)\gamma_{H}\left(\mathbf{\frac{1}{2}}\right)\right)
\end{align*}
Letting $N\rightarrow\infty$, the above $\mathcal{F}_{2,q_{H}}^{d}$-converges
to: 
\begin{align*}
\chi_{H}\left(\mathbf{z}\right) & \overset{\mathcal{F}_{2,q_{H}}^{d}}{=}-\gamma_{H}\left(\mathbf{\frac{1}{2}}\right)+\sum_{n=0}^{\infty}\kappa_{H}\left(\left[\mathbf{z}\right]_{2^{n}}\right)\left(H^{\prime}\left(\mathbf{0}\right)\right)^{n}\left(\beta_{H}\left(\mathbf{0}\right)+\left(\mathbf{I}_{d}-\alpha_{H}\left(\mathbf{0}\right)\right)\gamma_{H}\left(\mathbf{\frac{1}{2}}\right)\right)
\end{align*}

\vphantom{}

II. Letting $p\geq3$, \textbf{Corollary \ref{cor:MD Quasi-integrability of Chi_H for I minus alpha invertible}}
gives us: 
\[
\hat{\chi}_{H}\left(\mathbf{t}\right)=\hat{A}_{H}\left(\mathbf{t}\right)\left(\mathbf{I}_{d}-\alpha_{H}\left(\mathbf{0}\right)\right)^{-1}\beta_{H}\left(\mathbf{0}\right)+\begin{cases}
\mathbf{0} & \textrm{if }\mathbf{t}=\mathbf{0}\\
\hat{A}_{H}\left(\mathbf{t}\right)\gamma_{H}\left(\frac{\mathbf{t}\left|\mathbf{t}\right|_{p}}{p}\right) & \textrm{else}
\end{cases}
\]
Using (\ref{eq:MD Convolution of dA_H and D_N when H is commutative})
and \textbf{Lemma \ref{lem:MD gamma formulae}}, we then have: 
\begin{align*}
\tilde{\chi}_{H,N}\left(\mathbf{z}\right) & =\kappa_{H}\left(\left[\mathbf{z}\right]_{p^{N}}\right)\left(H^{\prime}\left(\mathbf{0}\right)\right)^{N}+\sum_{n=0}^{N-1}\kappa_{H}\left(\left[\mathbf{z}\right]_{p^{n}}\right)\left(H^{\prime}\left(\mathbf{0}\right)\right)^{n}\left(\mathbf{I}_{d}-\alpha_{H}\left(\mathbf{0}\right)\right)\\
 & +\sum_{n=0}^{N-1}\kappa_{H}\left(\left[\mathbf{z}\right]_{p^{n}}\right)\sum_{\mathbf{j}>\mathbf{0}}^{p-1}\mathcal{C}_{H}\left(\alpha_{H}\left(\frac{\mathbf{j}}{p}\right)\varepsilon_{n}\left(\mathbf{j}\mathbf{z}\right):\lambda_{p}\left(\left[\mathbf{z}\right]_{p^{n}}\right)\right)\left(H^{\prime}\left(\mathbf{0}\right)\right)^{n}\gamma_{H}\left(\frac{\mathbf{j}}{p}\right)
\end{align*}
Taking limits as $N\rightarrow\infty$ gives the desired formula.

Q.E.D. 
\begin{thm}[\textbf{Tauberian Spectral Theorem for Multi-Dimensional $p$-Hydra
Maps}]
\index{Hydra map!Tauberian Spectral Theorem}\label{thm:MD Periodic Points using WTT}
\index{$p,q$-adic!Wiener Tauberian Theorem}\index{Hydra map!divergent trajectories}Let
$H$ be as given in \textbf{\emph{Theorem \ref{thm:MD F-series for Chi_H}}}.

If $\alpha_{H}\left(\mathbf{0}\right)=\mathbf{I}_{d}$, let: 
\begin{equation}
\hat{\chi}_{H}\left(\mathbf{t}\right)=\begin{cases}
-\gamma_{H}\left(\mathbf{\frac{1}{2}}\right) & \textrm{if }\mathbf{t}=\mathbf{0}\\
\hat{A}_{H}\left(\mathbf{t}\right)v_{2}\left(\mathbf{t}\right)\beta_{H}\left(\mathbf{0}\right) & \textrm{else }
\end{cases},\textrm{ }\forall\mathbf{t}\in\hat{\mathbb{Z}}_{2}^{r}
\end{equation}
if every $p=2$. If $p\geq3$, let: 
\begin{equation}
\hat{\chi}_{H}\left(\mathbf{t}\right)=\begin{cases}
\mathbf{0} & \textrm{if }\mathbf{t}=\mathbf{0}\\
\hat{A}_{H}\left(\mathbf{t}\right)\left(v_{p}\left(\mathbf{t}\right)\beta_{H}\left(\mathbf{0}\right)+\gamma_{H}\left(\frac{\mathbf{t}\left|\mathbf{t}\right|_{p}}{p}\right)\right) & \textrm{else}
\end{cases},\textrm{ }\forall\mathbf{t}\in\hat{\mathbb{Z}}_{p}^{r}
\end{equation}

Alternatively, if $\alpha_{H}\left(\mathbf{0}\right)\neq\mathbf{I}_{d}$,
suppose $\mathbf{I}_{d}-\alpha_{H}\left(\mathbf{0}\right)$ is invertible
and that $H$ is commutative. Then, let: 
\begin{equation}
\hat{\chi}_{H}\left(\mathbf{t}\right)=\hat{A}_{H}\left(\mathbf{t}\right)\left(\mathbf{I}_{d}-\alpha_{H}\left(\mathbf{0}\right)\right)^{-1}\beta_{H}\left(\mathbf{0}\right)+\begin{cases}
\mathbf{0} & \textrm{if }\mathbf{t}=\mathbf{0}\\
\hat{A}_{H}\left(\mathbf{t}\right)\gamma_{H}\left(\frac{\mathbf{t}\left|\mathbf{t}\right|_{p}}{p}\right) & \textrm{else}
\end{cases},\textrm{ }\forall\mathbf{t}\in\hat{\mathbb{Z}}_{p}^{r}
\end{equation}

Finally, let $\mathbf{x}\in\mathbb{Z}^{d}\backslash\left\{ \mathbf{0}\right\} $.
Then, for the formula for $\hat{\chi}_{H}$ chosen according to the
situations described above:

\vphantom{}

I. If $\mathbf{x}$ is a periodic point of $H$, the translates of
the function $\hat{\chi}_{H}\left(\mathbf{t}\right)-\mathbf{x}\mathbf{1}_{\mathbf{0}}\left(\mathbf{t}\right)$
are \emph{not }dense in $c_{0}\left(\hat{\mathbb{Z}}_{p}^{r},\mathbb{C}_{q_{H}}^{d}\right)$.

\vphantom{}

II. Suppose further that $H$ is integral, and that $\left\Vert H_{\mathbf{j}}\left(\mathbf{0}\right)\right\Vert _{q_{H}}=1$
for all $\mathbf{j}\in\left(\mathbb{Z}^{r}/p\mathbb{Z}^{r}\right)\backslash\left\{ \mathbf{0}\right\} $.
If the translates of the function $\hat{\chi}_{H}\left(\mathbf{t}\right)-\mathbf{x}\mathbf{1}_{\mathbf{0}}\left(\mathbf{t}\right)$
are \emph{not }dense in $c_{0}\left(\hat{\mathbb{Z}}_{p}^{r},\mathbb{C}_{q_{H}}^{d}\right)$,
then either $\mathbf{x}$ is a periodic point of $H$ or $\mathbf{x}$
belongs to a divergent trajectory of $H$. 
\end{thm}
Proof: Essentially identical to the one-dimensional case (\textbf{Theorem
\ref{thm:MD pq-adic WTT for continuous functions}}).

Q.E.D.

\subsection{\label{subsec:6.2.4 Multi-Dimensional--=00003D000026}Multi-Dimensional
$\hat{\chi}_{H}$ and $\tilde{\chi}_{H,N}$ \textendash{} The Non-Commutative
Case}

Suppose that $H$ is non-commutative. Then, the multi-dimensional
analogue of the factor $1-\alpha_{H}\left(0\right)$ is now: 
\begin{equation}
\mathbf{I}_{H}\left(\lambda_{p}\left(\left[\mathbf{z}\right]_{p^{n}}\right)\right)=\mathbf{I}_{d}-\left(H^{\prime}\left(\mathbf{0}\right)\right)^{\lambda_{p}\left(\left[\mathbf{z}\right]_{p^{n}}\right)}\alpha_{H}\left(\mathbf{0}\right)\left(H^{\prime}\left(\mathbf{0}\right)\right)^{-\lambda_{p}\left(\left[\mathbf{z}\right]_{p^{n}}\right)}
\end{equation}
In the multi-dimensional case, the $\mathcal{F}_{p,q_{H}}$-limit
of the Fourier series generated by $\hat{A}_{H}\left(t\right)$ is:
\begin{equation}
\tilde{A}_{H}\left(\mathbf{z}\right)\overset{\textrm{def}}{=}\lim_{N\rightarrow\infty}\tilde{A}_{H,N}\left(\mathbf{z}\right)\overset{\mathcal{F}_{p,q_{H}}^{d,d}}{=}\sum_{n=0}^{\infty}\kappa_{H}\left(\left[\mathbf{z}\right]_{p^{n}}\right)\mathbf{I}_{H}\left(\lambda_{p}\left(\left[\mathbf{z}\right]_{p^{n}}\right)\right)\left(H^{\prime}\left(\mathbf{0}\right)\right)^{n}\label{eq:MD Definition of A_H twiddle}
\end{equation}
(equation (\ref{eq:MD Fourier Limit of A_H,N twiddle in standard frame})
from \textbf{Theorem \ref{thm:MD properties of A_H hat}}). To speak
of the multi-dimensional case for a moment using the terminology of
the one-dimensional case, because $\psi_{H}\left(\mathfrak{z}\right)$
satisfied: 
\begin{equation}
\psi_{H}\left(\mathfrak{z}\right)\overset{\mathcal{F}_{p,q_{H}}}{=}\sum_{n=0}^{\infty}\kappa_{H}\left(\left[\mathfrak{z}\right]_{p^{n}}\right)\left(H^{\prime}\left(0\right)\right)^{n}
\end{equation}
(\textbf{Lemma \ref{eq:Rising-continuation of Little Psi_H}}) the
method of the one-dimensional case was to exploit the fact that the
one-dimensional counterpart of $\mathbf{I}_{H}\left(\lambda_{p}\left(\left[\mathbf{z}\right]_{p^{n}}\right)\right)$
reduced to a constant which could then be factored out of the $n$-sum.
To overcome this obstacle in the $\alpha_{H}\left(0\right)\neq1$
case, we ended up writing: 
\begin{align}
\hat{\psi}_{H}\left(t\right) & =\frac{\hat{A}_{H}\left(t\right)}{1-\alpha_{H}\left(0\right)}\\
\psi_{H}\left(\mathfrak{z}\right) & =\frac{A_{H}\left(\mathfrak{z}\right)}{1-\alpha_{H}\left(0\right)}
\end{align}
the right-hand sides of which were defined solely because $\alpha_{H}\left(0\right)\neq1$.
In particular, the non-vanishing of $\alpha_{H}\left(0\right)$ makes
multiplication by $1/\left(1-\alpha_{H}\left(0\right)\right)$ into
an \emph{invertible linear transformation}, with the linear transformation:
\begin{equation}
\mathcal{L}\left\{ f\right\} \left(\mathfrak{z}\right)\overset{\textrm{def}}{=}\left(1-\alpha_{H}\left(0\right)\right)f\left(\mathfrak{z}\right)
\end{equation}
then satisfying $\mathcal{L}\left\{ \psi_{H}\right\} \left(\mathfrak{z}\right)=A_{H}\left(\mathfrak{z}\right)$.

Turning now to the multi-dimensional case, this suggests that to replicate
this argument, we will need to find an invertible linear transformation
which sends $\psi_{H}\left(\mathbf{z}\right)$ to $\tilde{A}_{H}\left(\mathbf{z}\right)$.
In theory, by computing the effects of this linear transformation
on the Fourier transforms of $\psi_{H}$ and $\tilde{A}_{H}$, we
can then reverse-engineer a formula for $\hat{\psi}_{H}$ from our
formula for $\hat{A}_{H}$.
\begin{defn}
We define the following operators on the space of functions $\mathbb{Z}_{p}^{r}\rightarrow\mathbb{C}_{q}^{d,d}$:

\vphantom{}

I. For all $n\geq0$: 
\begin{equation}
\mathcal{T}_{n}\left\{ \mathbf{F}\right\} \left(\mathbf{z}\right)\overset{\textrm{def}}{=}\mathbf{F}\left(\left[\mathbf{z}\right]_{p^{n}}\right)\label{eq:Definition of the nth MD truncation operator}
\end{equation}

\vphantom{}

II. For all $n\geq0$: 
\begin{equation}
\mathcal{L}_{H,1,n}\left\{ \mathbf{F}\right\} \left(\mathbf{z}\right)\overset{\textrm{def}}{=}\mathcal{T}_{n}\left\{ \mathbf{F}\right\} \left(\mathbf{z}\right)\left(\mathbf{I}_{d}-H^{\prime}\left(\mathbf{0}\right)\right)\left(H^{\prime}\left(\mathbf{0}\right)\right)^{-n}\label{eq:Definition of L_H,1,n}
\end{equation}

\vphantom{}

III. For all $n\geq1$: 
\begin{equation}
\mathcal{L}_{H,2,n}\left\{ \mathbf{F}\right\} \left(\mathbf{z}\right)\overset{\textrm{def}}{=}\mathcal{L}_{H,1,n}\left\{ \mathbf{F}\right\} \left(\mathbf{z}\right)-\mathcal{L}_{H,1,n-1}\left\{ \mathbf{F}\right\} \left(\mathbf{z}\right)\left(H^{\prime}\left(\mathbf{0}\right)\right)^{-1}\label{eq:Definition of L_H,2,n}
\end{equation}

\vphantom{}

IV. 
\begin{equation}
\mathcal{E}_{H}\left\{ \mathbf{F}\right\} \left(\mathbf{z}\right)\overset{\textrm{def}}{=}\mathbf{F}\left(\mathbf{0}\right)\left(\mathbf{I}_{d}-H^{\prime}\left(\mathbf{0}\right)\right)\label{eq:Definition of E_H}
\end{equation}

\vphantom{}

V. For all $n\geq0$: 
\begin{equation}
\mathcal{L}_{H,3,n}\left\{ \mathbf{F}\right\} \left(\mathbf{z}\right)\overset{\textrm{def}}{=}\begin{cases}
\mathcal{E}_{0}\left\{ \mathbf{F}\right\} \left(\mathbf{z}\right) & \textrm{if }n=0\\
\mathcal{E}_{0}\left\{ \mathbf{F}\right\} \left(\mathbf{z}\right)+\sum_{m=1}^{n}\mathcal{L}_{H,2,m}\left\{ \mathbf{F}\right\} \left(\mathbf{z}\right) & \textrm{if }n\geq1
\end{cases}\label{eq:Definition of L_H,3,n}
\end{equation}

\vphantom{}

VI. For all $n\geq0$: 
\begin{equation}
\mathcal{L}_{H,4,n}\left\{ \mathbf{F}\right\} \left(\mathbf{z}\right)\overset{\textrm{def}}{=}\mathcal{L}_{H,3,n}\left\{ \mathbf{F}\right\} \left(\mathbf{z}\right)\left(\mathbf{I}_{H}\left(\lambda_{p}\left(\left[\mathbf{z}\right]_{p^{n}}\right)\right)\right)^{-1}\left(H^{\prime}\left(\mathbf{0}\right)\right)^{n}\label{eq:Definition of L_H,4,n}
\end{equation}

\vphantom{}

VII. 
\begin{equation}
\mathcal{L}_{H}\left\{ \mathbf{F}\right\} \left(\mathbf{z}\right)\overset{\textrm{def}}{=}\sum_{n=0}^{\infty}\mathcal{L}_{H,4,n}\left\{ \mathbf{F}\right\} \left(\mathbf{z}\right)\label{eq:Definition of L_H}
\end{equation}
\end{defn}
\begin{rem}
Note that all of these operators are \emph{linear}. 
\end{rem}
\begin{prop}
If $\mathbf{I}_{d}-\alpha_{H}\left(\mathbf{0}\right)$ is invertible,
then $\mathbf{I}_{H}\left(n\right)$ is invertible for all $n\geq0$. 
\end{prop}
Proof: 
\begin{align*}
\mathbf{I}_{H}\left(n\right) & =\mathbf{I}_{d}-\left(H^{\prime}\left(\mathbf{0}\right)\right)^{n}\alpha_{H}\left(\mathbf{0}\right)\left(H^{\prime}\left(\mathbf{0}\right)\right)^{-n}\\
 & =\left(H^{\prime}\left(\mathbf{0}\right)\right)^{n}\left(\mathbf{I}_{d}-\alpha_{H}\left(\mathbf{0}\right)\right)\left(H^{\prime}\left(\mathbf{0}\right)\right)^{-n}
\end{align*}
Since $H$ is a $p$-Hydra map, $H^{\prime}\left(\mathbf{0}\right)$
is invertible, and hence: 
\[
\det\mathbf{I}_{H}\left(n\right)=\det\left(\mathbf{I}_{d}-\alpha_{H}\left(\mathbf{0}\right)\right)
\]
which is non-zero if and only if $\mathbf{I}_{d}-\alpha_{H}\left(\mathbf{0}\right)$
is invertible.

Q.E.D. 
\begin{prop}
Suppose $\mathbf{I}_{d}-\alpha_{H}\left(\mathbf{0}\right)$ is invertible.
Then:

\begin{equation}
\mathcal{L}_{H}\left\{ \tilde{A}_{H}\right\} \left(\mathbf{z}\right)\overset{\mathcal{F}_{p,q_{H}}^{d,d}}{=}\psi_{H}\left(\mathbf{z}\right)\label{eq:L_H sends A_H twiddle to Little Psi_H}
\end{equation}
where, as indicated, we compute the value of either side by summing
the infinite series on either side in accordance with the standard
$\left(p,q_{H}\right)$-adic frame. 
\end{prop}
\begin{rem}
The linear operator which sends $\psi_{H}$ to $\tilde{A}_{H}$ is
nearly identical to $\mathcal{L}_{H}$; all you do is change the multiplication
by $\left(I_{H}\left(\lambda_{p}\left(\left[\mathbf{z}\right]_{p^{n}}\right)\right)\right)^{-1}$
in $\mathcal{L}_{H,4,n}$ to multiplication by $I_{H}\left(\lambda_{p}\left(\left[\mathbf{z}\right]_{p^{n}}\right)\right)$. 
\end{rem}
Proof: First, we observe that: 
\begin{align*}
\mathcal{T}_{m}\left\{ \tilde{A}_{H}\right\} \left(\mathbf{z}\right) & =\tilde{A}_{H}\left(\left[\mathbf{z}\right]_{p^{m}}\right)\\
 & \overset{\mathbb{C}^{d,d}}{=}\sum_{n=0}^{\infty}\kappa_{H}\left(\left[\left[\mathbf{z}\right]_{p^{m}}\right]_{p^{n}}\right)\mathbf{I}_{H}\left(\lambda_{p}\left(\left[\left[\mathbf{z}\right]_{p^{m}}\right]_{p^{n}}\right)\right)\left(H^{\prime}\left(\mathbf{0}\right)\right)^{n}\\
 & \overset{\mathbb{C}^{d,d}}{=}\kappa_{H}\left(\left[\mathbf{z}\right]_{p^{m}}\right)\mathbf{I}_{H}\left(\lambda_{p}\left(\left[\mathbf{z}\right]_{p^{m}}\right)\right)\sum_{n=m}^{\infty}\left(H^{\prime}\left(\mathbf{0}\right)\right)^{n}\\
 & +\sum_{n=0}^{m-1}\kappa_{H}\left(\left[\mathbf{z}\right]_{p^{n}}\right)\mathbf{I}_{H}\left(\lambda_{p}\left(\left[\mathbf{z}\right]_{p^{n}}\right)\right)\left(H^{\prime}\left(\mathbf{0}\right)\right)^{n}
\end{align*}
Summing the geometric series yields:
\begin{align*}
\mathcal{T}_{m}\left\{ \tilde{A}_{H}\right\} \left(\mathbf{z}\right) & \overset{\mathbb{C}^{d,d}}{=}\kappa_{H}\left(\left[\mathbf{z}\right]_{p^{m}}\right)\mathbf{I}_{H}\left(\lambda_{p}\left(\left[\mathbf{z}\right]_{p^{m}}\right)\right)\left(H^{\prime}\left(\mathbf{0}\right)\right)^{m}\left(\mathbf{I}_{d}-H^{\prime}\left(\mathbf{0}\right)\right)^{-1}\\
 & +\sum_{n=0}^{m-1}\kappa_{H}\left(\left[\mathbf{z}\right]_{p^{n}}\right)\mathbf{I}_{H}\left(\lambda_{p}\left(\left[\mathbf{z}\right]_{p^{n}}\right)\right)\left(H^{\prime}\left(\mathbf{0}\right)\right)^{n}
\end{align*}
Thus: 
\begin{align*}
\mathcal{L}_{H,1,m}\left\{ \tilde{A}_{H}\right\} \left(\mathbf{z}\right) & \overset{\overline{\mathbb{Q}}^{d,d}}{=}\kappa_{H}\left(\left[\mathbf{z}\right]_{p^{m}}\right)\mathbf{I}_{H}\left(\lambda_{p}\left(\left[\mathbf{z}\right]_{p^{m}}\right)\right)\\
 & +\sum_{n=0}^{m-1}\kappa_{H}\left(\left[\mathbf{z}\right]_{p^{n}}\right)\mathbf{I}_{H}\left(\lambda_{p}\left(\left[\mathbf{z}\right]_{p^{n}}\right)\right)\left(H^{\prime}\left(\mathbf{0}\right)\right)^{n}\left(\mathbf{I}_{d}-H^{\prime}\left(\mathbf{0}\right)\right)\left(H^{\prime}\left(\mathbf{0}\right)\right)^{-m}
\end{align*}
and so: 
\begin{equation}
\mathcal{L}_{H,2,m}\left\{ \tilde{A}_{H}\right\} \left(\mathbf{z}\right)\overset{\overline{\mathbb{Q}}^{d,d}}{=}\kappa_{H}\left(\left[\mathbf{z}\right]_{p^{m}}\right)\mathbf{I}_{H}\left(\lambda_{p}\left(\left[\mathbf{z}\right]_{p^{m}}\right)\right)-\kappa_{H}\left(\left[\mathbf{z}\right]_{p^{m-1}}\right)\mathbf{I}_{H}\left(\lambda_{p}\left(\left[\mathbf{z}\right]_{p^{m-1}}\right)\right)
\end{equation}
Next, because: 
\begin{equation}
\kappa_{H}\left(\mathbf{0}\right)=\mathbf{I}_{d}
\end{equation}
and: 
\begin{equation}
\mathbf{I}_{H}\left(0\right)=\mathbf{I}_{d}-\alpha_{H}\left(\mathbf{0}\right)
\end{equation}
we have: 
\[
\sum_{m=1}^{n}\mathcal{L}_{H,2,m}\left\{ \tilde{A}_{H}\right\} \left(\mathbf{z}\right)=\kappa_{H}\left(\left[\mathbf{z}\right]_{p^{n}}\right)\mathbf{I}_{H}\left(\lambda_{p}\left(\left[\mathbf{z}\right]_{p^{n}}\right)\right)-\left(\mathbf{I}_{d}-\alpha_{H}\left(\mathbf{0}\right)\right)
\]

Now: 
\begin{equation}
\mathcal{E}_{H}\left\{ \tilde{A}_{H}\right\} \left(\mathbf{z}\right)=\tilde{A}_{H}\left(\mathbf{0}\right)\left(\mathbf{I}_{d}-H^{\prime}\left(\mathbf{0}\right)\right)
\end{equation}
Since $\alpha_{H}\left(\mathbf{0}\right)\neq\mathbf{I}_{d}$: 
\begin{align*}
\tilde{A}_{H}\left(\mathbf{0}\right) & =\sum_{n=0}^{\infty}\kappa_{H}\left(\mathbf{0}\right)\mathbf{I}_{H}\left(0\right)\left(H^{\prime}\left(\mathbf{0}\right)\right)^{n}\\
 & =\left(\mathbf{I}_{d}-\alpha_{H}\left(\mathbf{0}\right)\right)\sum_{n=0}^{\infty}\left(H^{\prime}\left(\mathbf{0}\right)\right)^{n}\\
 & =\left(\mathbf{I}_{d}-\alpha_{H}\left(\mathbf{0}\right)\right)\left(\mathbf{I}_{d}-H^{\prime}\left(\mathbf{0}\right)\right)^{-1}
\end{align*}
and so: 
\begin{align*}
\mathcal{E}_{H}\left\{ \tilde{A}_{H}\right\} \left(\mathbf{z}\right) & =\tilde{A}_{H}\left(\mathbf{0}\right)\left(\mathbf{I}_{d}-H^{\prime}\left(\mathbf{0}\right)\right)\\
 & =\left(\mathbf{I}_{d}-\alpha_{H}\left(\mathbf{0}\right)\right)\left(\mathbf{I}_{d}-H^{\prime}\left(\mathbf{0}\right)\right)^{-1}\left(\mathbf{I}_{d}-H^{\prime}\left(\mathbf{0}\right)\right)\\
 & =\mathbf{I}_{d}-\alpha_{H}\left(\mathbf{0}\right)
\end{align*}
So, for $n\geq1$: 
\begin{align*}
\mathcal{L}_{H,3,n}\left\{ \tilde{A}_{H}\right\} \left(\mathbf{z}\right) & =\mathcal{E}_{H}\left\{ \tilde{A}_{H}\right\} \left(\mathbf{z}\right)+\sum_{m=1}^{n}\mathcal{L}_{H,2,m}\left\{ \tilde{A}_{H}\right\} \left(\mathbf{z}\right)\\
 & =\left(\mathbf{I}_{d}-\alpha_{H}\left(\mathbf{0}\right)\right)+\kappa_{H}\left(\left[\mathbf{z}\right]_{p^{n}}\right)\mathbf{I}_{H}\left(\lambda_{p}\left(\left[\mathbf{z}\right]_{p^{n}}\right)\right)-\left(\mathbf{I}_{d}-\alpha_{H}\left(\mathbf{0}\right)\right)\\
 & =\kappa_{H}\left(\left[\mathbf{z}\right]_{p^{n}}\right)\mathbf{I}_{H}\left(\lambda_{p}\left(\left[\mathbf{z}\right]_{p^{n}}\right)\right)
\end{align*}
Finally:

\begin{align*}
\mathcal{L}_{H}\left\{ \tilde{A}_{H}\right\} \left(\mathbf{z}\right) & \overset{\mathcal{F}_{p,q_{H}}^{d,d}}{=}\sum_{n=0}^{\infty}\mathcal{L}_{H,4,n}\left\{ \tilde{A}_{H}\right\} \left(\mathbf{z}\right)\\
 & \overset{\mathcal{F}_{p,q_{H}}^{d,d}}{=}\sum_{n=0}^{\infty}\mathcal{L}_{H,3,n}\left\{ \tilde{A}_{H}\right\} \left(\mathbf{z}\right)\left(\mathbf{I}_{H}\left(\lambda_{p}\left(\left[\mathbf{z}\right]_{p^{n}}\right)\right)\right)^{-1}\left(H^{\prime}\left(\mathbf{0}\right)\right)^{n}\\
 & \overset{\mathcal{F}_{p,q_{H}}^{d,d}}{=}\sum_{n=0}^{\infty}\kappa_{H}\left(\left[\mathbf{z}\right]_{p^{n}}\right)\mathbf{I}_{H}\left(\lambda_{p}\left(\left[\mathbf{z}\right]_{p^{n}}\right)\right)\left(\mathbf{I}_{H}\left(\lambda_{p}\left(\left[\mathbf{z}\right]_{p^{n}}\right)\right)\right)^{-1}\left(H^{\prime}\left(\mathbf{0}\right)\right)^{n}\\
 & \overset{\mathcal{F}_{p,q_{H}}^{d,d}}{=}\sum_{n=0}^{\infty}\kappa_{H}\left(\left[\mathbf{z}\right]_{p^{n}}\right)\left(H^{\prime}\left(\mathbf{0}\right)\right)^{n}\\
 & \overset{\mathcal{F}_{p,q_{H}}^{d,d}}{=}\psi_{H}\left(\mathbf{z}\right)
\end{align*}
where $\left(\mathbf{I}_{H}\left(\lambda_{p}\left(\left[\mathbf{z}\right]_{p^{n}}\right)\right)\right)^{-1}$
is defined for all $\mathbf{z}$ and all $n$ because the invertibility
of $\mathbf{I}_{d}-\alpha_{H}\left(\mathbf{0}\right)$ guarantees
the invertibility of $\mathbf{I}_{H}\left(m\right)$ for all $m\geq0$.

Q.E.D.

\vphantom{}

So, to obtain a Fourier transform for $\psi_{H}$, we just need only
compute the effect $\mathcal{L}_{H}$ has on $\tilde{A}_{H}$'s Fourier
transform ($\hat{A}_{H}$). This, however, is easier said than done.
I have broken down the computation into several steps; to begin, we
need to following definition: 
\begin{defn}
For each $m\in\mathbb{N}_{0}$, let $\hat{K}_{m}:\hat{\mathbb{Z}}_{p}^{r}\rightarrow\mathbb{C}_{q}$
be defined by: 
\begin{equation}
\hat{K}_{m}\left(\mathbf{t}\right)\overset{\textrm{def}}{=}\frac{1}{p^{rm}}\sum_{\mathbf{n}=\mathbf{0}}^{p^{m}-1}e^{-2\pi i\mathbf{n}\cdot\mathbf{t}}\label{eq:Definition of K_m hat}
\end{equation}
\end{defn}
\begin{prop}
\label{prop:A_H hat convolve K_H hat MD}Let $m\geq0$. Then, for
all $\mathbf{t}\in\hat{\mathbb{Z}}_{p}^{r}$: 
\begin{equation}
\left(\hat{A}_{H}*\hat{K}_{m}\right)\left(\mathbf{t}\right)=\sum_{n=0}^{m-1}\mathbf{1}_{\mathbf{0}}\left(p^{n}\mathbf{t}\right)\prod_{k=0}^{n-1}\alpha_{H}\left(p^{k}\mathbf{t}\right)+\left(\prod_{k=0}^{m-1}\alpha_{H}\left(p^{k}\mathbf{t}\right)\right)\left(\mathbf{I}_{d}-H^{\prime}\left(\mathbf{0}\right)\right)^{-1}\label{eq:Convolution of A_H hat and K_m hat}
\end{equation}
where the $k$ product is defined to be $\mathbf{I}_{d}$ whenever
its upper limit is $<0$. 
\end{prop}
Proof: We begin with the definition of $\hat{A}_{H}*\hat{K}_{m}$:
\begin{equation}
\left(\hat{A}_{H}*\hat{K}_{m}\right)\left(\mathbf{t}\right)=\sum_{\mathbf{s}\in\hat{\mathbb{Z}}_{p}^{r}}\hat{A}_{H}\left(\mathbf{s}\right)\frac{1}{p^{rm}}\sum_{\mathbf{n}=\mathbf{0}}^{p^{m}-1}e^{-2\pi i\left\{ \mathbf{n}\left(\mathbf{t}-\mathbf{s}\right)\right\} _{p}}
\end{equation}
We then pull out the $\mathbf{s}=\mathbf{0}$ term: 
\begin{align*}
\hat{A}_{H}\left(\mathbf{0}\right)\frac{1}{p^{rm}}\sum_{\mathbf{n}=\mathbf{0}}^{p^{m}-1}e^{-2\pi i\left\{ \mathbf{n}\left(\mathbf{t}-\mathbf{0}\right)\right\} _{p}} & =\frac{1}{p^{rm}}\sum_{\mathbf{n}=\mathbf{0}}^{p^{m}-1}e^{-2\pi i\left\{ \mathbf{n}\mathbf{t}\right\} _{p}}\\
 & =\frac{1}{p^{rm}}\sum_{\mathbf{n}=\mathbf{0}}^{p^{m}-1}e^{-2\pi i\left\{ \frac{\mathbf{n}}{p^{m}}p^{m}\mathbf{t}\right\} _{p}}\\
\left(\left\{ \mathbf{s}\in\hat{\mathbb{Z}}_{p}^{r}:\left\Vert \mathbf{s}\right\Vert _{p}\leq p^{m}\right\} =\left\{ \frac{\mathbf{n}}{p^{m}}:\mathbf{n}\leq p^{m}-1\right\} \right); & =\frac{1}{p^{rm}}\sum_{\left\Vert \mathbf{s}\right\Vert _{p}\leq p^{m}}e^{-2\pi i\left\{ \mathbf{s}p^{m}\mathbf{t}\right\} _{p}}\\
\left(\textrm{\textbf{proposition \ref{prop:Multi-Dimensional indicator function Fourier Series}}}\right); & =\left[p^{m}\mathbf{t}\overset{p^{m}}{\equiv}\mathbf{0}\right]
\end{align*}
Since congruence $p^{m}\mathbf{t}\overset{p^{m}}{\equiv}\mathbf{0}$
means that $p^{m}t_{\ell}\overset{p^{m}}{\equiv}0$ for all $\ell\in\left\{ 1,\ldots,r\right\} $,
note that this occurs if and only if $\left|t_{\ell}\right|_{p}\leq p^{m}$.
As such, our $\mathbf{s}=\mathbf{0}$ term becomes: 
\[
\hat{A}_{H}\left(\mathbf{0}\right)\frac{1}{p^{rm}}\sum_{\mathbf{n}=\mathbf{0}}^{p^{m}-1}e^{-2\pi i\left\{ \mathbf{n}\left(\mathbf{t}-\mathbf{0}\right)\right\} _{p}}=\left[p^{m}\mathbf{t}\overset{p^{m}}{\equiv}\mathbf{0}\right]=\mathbf{1}_{\mathbf{0}}\left(p^{m}\mathbf{t}\right)
\]

Next, using \textbf{Proposition \ref{prop:MD alpha H series}}, we
express $\hat{A}_{H}\left(\mathbf{s}\right)$ as a series: 
\begin{equation}
\hat{A}_{H}\left(\mathbf{s}\right)=\prod_{m=0}^{-v_{p}\left(\mathbf{s}\right)-1}\alpha_{H}\left(p^{m}\mathbf{s}\right)=\sum_{\mathbf{k}=\mathbf{0}}^{p^{-v_{p}\left(\mathbf{s}\right)}-1}\kappa_{H}\left(\mathbf{k}\right)\left(\frac{H^{\prime}\left(\mathbf{0}\right)}{p^{r}}\right)^{-v_{p}\left(\mathbf{s}\right)}e^{-2\pi i\left(\mathbf{k}\cdot\mathbf{s}\right)}
\end{equation}
and so: 
\begin{equation}
\left(\hat{A}_{H}*\hat{K}_{m}\right)\left(\mathbf{t}\right)=\sum_{\mathbf{s}\in\hat{\mathbb{Z}}_{p}^{r}}\sum_{\mathbf{k}=\mathbf{0}}^{p^{-v_{p}\left(\mathbf{s}\right)}-1}\kappa_{H}\left(\mathbf{k}\right)\left(\frac{H^{\prime}\left(\mathbf{0}\right)}{p^{r}}\right)^{-v_{p}\left(\mathbf{s}\right)}\frac{e^{-2\pi i\left(\mathbf{k}\cdot\mathbf{s}\right)}}{p^{rm}}\sum_{\mathbf{n}=\mathbf{0}}^{p^{m}-1}e^{-2\pi i\left\{ \mathbf{n}\left(\mathbf{t}-\mathbf{s}\right)\right\} _{p}}
\end{equation}
Like many times before, we now split the $\mathbf{s}$-sum into a
sum over level sets $\left\Vert \mathbf{s}\right\Vert _{p}=p^{h}$.
This gives: 
\begin{equation}
\mathbf{1}_{\mathbf{0}}\left(p^{m}\mathbf{t}\right)+\sum_{h=1}^{\infty}\sum_{\left\Vert \mathbf{s}\right\Vert _{p}=p^{h}}\sum_{\mathbf{k}=\mathbf{0}}^{p^{h}-1}\kappa_{H}\left(\mathbf{k}\right)\left(\frac{H^{\prime}\left(\mathbf{0}\right)}{p^{r}}\right)^{h}\frac{e^{-2\pi i\left(\mathbf{k}\cdot\mathbf{s}\right)}}{p^{rm}}\sum_{\mathbf{n}=\mathbf{0}}^{p^{m}-1}e^{-2\pi i\left\{ \mathbf{n}\left(\mathbf{t}-\mathbf{s}\right)\right\} _{p}}
\end{equation}
as our expression for $\left(\hat{A}_{H}*\hat{K}_{m}\right)\left(\mathbf{t}\right)$.

Pulling the $\left\Vert \mathbf{s}\right\Vert _{p}$-sum to the far
right, we have: 
\begin{align*}
e^{-2\pi i\left(\mathbf{k}\cdot\mathbf{s}\right)}\sum_{\left\Vert \mathbf{s}\right\Vert _{p}=p^{h}}e^{-2\pi i\left\{ \mathbf{n}\left(\mathbf{t}-\mathbf{s}\right)\right\} _{p}} & =e^{-2\pi i\mathbf{n}\cdot\mathbf{t}}\sum_{\left\Vert \mathbf{s}\right\Vert _{p}=p^{h}}e^{2\pi i\left(\mathbf{n}-\mathbf{k}\right)\cdot\mathbf{s}}\\
 & =p^{rh}e^{-2\pi i\mathbf{n}\cdot\mathbf{t}}\left[\mathbf{k}\overset{p^{h}}{\equiv}\mathbf{n}\right]
\end{align*}
and so: 
\begin{align*}
\left(\hat{A}_{H}*\hat{K}_{m}\right)\left(\mathbf{t}\right) & =\mathbf{1}_{\mathbf{0}}\left(p^{m}\mathbf{t}\right)+\sum_{h=1}^{\infty}\sum_{\mathbf{k}=\mathbf{0}}^{p^{h}-1}\frac{1}{p^{rm}}\sum_{\mathbf{n}=\mathbf{0}}^{p^{m}-1}\kappa_{H}\left(\mathbf{k}\right)\left(H^{\prime}\left(\mathbf{0}\right)\right)^{h}e^{-2\pi i\mathbf{n}\cdot\mathbf{t}}\left[\mathbf{k}\overset{p^{h}}{\equiv}\mathbf{n}\right]\\
 & =\mathbf{1}_{\mathbf{0}}\left(p^{m}\mathbf{t}\right)+\sum_{h=1}^{\infty}\frac{1}{p^{rm}}\sum_{\mathbf{n}=\mathbf{0}}^{p^{m}-1}\kappa_{H}\left(\left[\mathbf{n}\right]_{p^{h}}\right)\left(H^{\prime}\left(\mathbf{0}\right)\right)^{h}e^{-2\pi i\mathbf{n}\cdot\mathbf{t}}\\
 & =\sum_{h=0}^{\infty}\frac{1}{p^{rm}}\sum_{\mathbf{n}=\mathbf{0}}^{p^{m}-1}\kappa_{H}\left(\left[\mathbf{n}\right]_{p^{h}}\right)\left(H^{\prime}\left(\mathbf{0}\right)\right)^{h}e^{-2\pi i\mathbf{n}\cdot\mathbf{t}}
\end{align*}
Because $m$ is finite, $\mathbf{n}\leq p^{m}-1$ tells us that $\left[\mathbf{n}\right]_{p^{h}}=\mathbf{n}$
for all $h\geq m$. Consequently: 
\begin{align*}
\left(\hat{A}_{H}*K_{m}\right)\left(\mathbf{t}\right) & =\sum_{h=0}^{m-1}\frac{1}{p^{rm}}\sum_{\mathbf{n}=\mathbf{0}}^{p^{m}-1}\kappa_{H}\left(\left[\mathbf{n}\right]_{p^{h}}\right)\left(H^{\prime}\left(\mathbf{0}\right)\right)^{h}e^{-2\pi i\mathbf{n}\cdot\mathbf{t}}\\
 & +\frac{1}{p^{rm}}\sum_{\mathbf{n}=\mathbf{0}}^{p^{m}-1}\kappa_{H}\left(\mathbf{n}\right)e^{-2\pi i\mathbf{n}\cdot\mathbf{t}}\sum_{h=m}^{\infty}\left(H^{\prime}\left(\mathbf{0}\right)\right)^{h}
\end{align*}
Since $H$ is contracting, by \textbf{Proposition \ref{prop:MD Contracting H proposition}}
(page \pageref{prop:MD Contracting H proposition}) the series $\sum_{h=m}^{\infty}\left(H^{\prime}\left(\mathbf{0}\right)\right)^{h}$
converges in $\mathbb{R}$ to the matrix: 
\[
\sum_{h=m}^{\infty}\left(H^{\prime}\left(\mathbf{0}\right)\right)^{h}\overset{\mathbb{R}^{d,d}}{=}\left(H^{\prime}\left(\mathbf{0}\right)\right)^{m}\left(\mathbf{I}_{d}-H^{\prime}\left(\mathbf{0}\right)\right)^{-1}
\]
which has entries in $\mathbb{Q}$. Thus: 
\begin{align*}
\left(\hat{A}_{H}*K_{m}\right)\left(\mathbf{t}\right) & =\sum_{h=0}^{m-1}\frac{1}{p^{rm}}\sum_{\mathbf{n}=\mathbf{0}}^{p^{m}-1}\kappa_{H}\left(\left[\mathbf{n}\right]_{p^{h}}\right)e^{-2\pi i\mathbf{n}\cdot\mathbf{t}}\left(H^{\prime}\left(\mathbf{0}\right)\right)^{h}\\
 & +\frac{1}{p^{rm}}\sum_{\mathbf{n}=\mathbf{0}}^{p^{m}-1}\kappa_{H}\left(\mathbf{n}\right)e^{-2\pi i\mathbf{n}\cdot\mathbf{t}}\left(H^{\prime}\left(\mathbf{0}\right)\right)^{m}\left(\mathbf{I}_{d}-H^{\prime}\left(\mathbf{0}\right)\right)^{-1}
\end{align*}
Next, we split the $\mathbf{n}$-sum on the upper line modulo $p^{h}$:
\begin{align*}
\sum_{\mathbf{n}=\mathbf{0}}^{p^{m}-1}\kappa_{H}\left(\left[\mathbf{n}\right]_{p^{h}}\right)e^{-2\pi i\mathbf{n}\cdot\mathbf{t}} & =\sum_{\mathbf{k}=\mathbf{0}}^{p^{h}-1}\sum_{\mathbf{n}=\mathbf{0}}^{p^{m-h}-1}\kappa_{H}\left(\left[p^{h}\mathbf{n}+\mathbf{k}\right]_{p^{h}}\right)e^{-2\pi i\left(p^{h}\mathbf{n}+\mathbf{k}\right)\cdot\mathbf{t}}\\
 & =\sum_{\mathbf{n}=\mathbf{0}}^{p^{m-h}-1}\left(\sum_{\mathbf{k}=\mathbf{0}}^{p^{h}-1}\kappa_{H}\left(\mathbf{k}\right)e^{-2\pi i\mathbf{k}\cdot\mathbf{t}}\right)e^{-2\pi i\mathbf{n}\cdot p^{h}\mathbf{t}}
\end{align*}
and so: 
\begin{align}
\left(\hat{A}_{H}*K_{m}\right)\left(\mathbf{t}\right) & =\sum_{h=0}^{m-1}\frac{1}{p^{rm}}\sum_{\mathbf{n}=\mathbf{0}}^{p^{m-h}-1}\left(\sum_{\mathbf{k}=\mathbf{0}}^{p^{h}-1}\kappa_{H}\left(\mathbf{k}\right)e^{-2\pi i\mathbf{k}\cdot\mathbf{t}}\right)e^{-2\pi i\mathbf{n}\cdot p^{h}\mathbf{t}}\left(H^{\prime}\left(\mathbf{0}\right)\right)^{h}\label{eq:MD Halfway done with K_m convolution computation}\\
 & +\frac{1}{p^{rm}}\sum_{\mathbf{n}=\mathbf{0}}^{p^{m}-1}\kappa_{H}\left(\mathbf{n}\right)e^{-2\pi i\mathbf{n}\cdot\mathbf{t}}\left(H^{\prime}\left(\mathbf{0}\right)\right)^{m}\left(\mathbf{I}_{d}-H^{\prime}\left(\mathbf{0}\right)\right)^{-1}\nonumber 
\end{align}

To finish, we yet again apply the familiar recursive evaluation technique.
First, define: 
\begin{equation}
S_{h}\left(\mathbf{t}\right)\overset{\textrm{def}}{=}\sum_{\mathbf{k}=\mathbf{0}}^{p^{h}-1}\kappa_{H}\left(\mathbf{k}\right)e^{-2\pi i\left(\mathbf{k}\cdot\mathbf{t}\right)}\label{eq:S_h for MD K_m convolution computation}
\end{equation}
Then: 
\begin{align*}
S_{h}\left(\mathbf{t}\right) & =\sum_{\mathbf{k}=\mathbf{0}}^{p^{h}-1}\kappa_{H}\left(\mathbf{k}\right)e^{-2\pi i\left(\mathbf{k}\cdot\mathbf{t}\right)}\\
 & =\sum_{\mathbf{j}=\mathbf{0}}^{p-1}\sum_{\mathbf{k}=\mathbf{0}}^{p^{h-1}-1}\kappa_{H}\left(p\mathbf{k}+\mathbf{j}\right)e^{-2\pi i\left(\left(p\mathbf{k}+\mathbf{j}\right)\cdot\mathbf{t}\right)}\\
 & =\sum_{\mathbf{j}=\mathbf{0}}^{p-1}\sum_{\mathbf{k}=\mathbf{0}}^{p^{h-1}-1}H_{\mathbf{j}}^{\prime}\left(\mathbf{0}\right)\kappa_{H}\left(\mathbf{k}\right)\left(H^{\prime}\left(\mathbf{0}\right)\right)^{-1}e^{-2\pi i\left(\left(p\mathbf{k}+\mathbf{j}\right)\cdot\mathbf{t}\right)}\\
 & =\left(\sum_{\mathbf{j}=\mathbf{0}}^{p-1}H_{\mathbf{j}}^{\prime}\left(\mathbf{0}\right)e^{-2\pi i\mathbf{j}\cdot\mathbf{t}}\right)\sum_{\mathbf{k}=\mathbf{0}}^{p^{h-1}-1}\kappa_{H}\left(\mathbf{k}\right)e^{-2\pi i\left(\mathbf{k}\cdot p\mathbf{t}\right)}\left(H^{\prime}\left(\mathbf{0}\right)\right)^{-1}\\
 & =p^{r}\alpha_{H}\left(\mathbf{t}\right)S_{h-1}\left(p\mathbf{t}\right)\left(H^{\prime}\left(\mathbf{0}\right)\right)^{-1}
\end{align*}
This yields the recursion relation: 
\begin{equation}
S_{h}\left(\mathbf{t}\right)=p^{r}\alpha_{H}\left(\mathbf{t}\right)S_{h-1}\left(p\mathbf{t}\right)\left(H^{\prime}\left(\mathbf{0}\right)\right)^{-1}\label{eq:S_h recursion relation for K_m convolution computation}
\end{equation}
and so: 
\begin{equation}
S_{h}\left(\mathbf{t}\right)=p^{rh}\left(\prod_{k=0}^{h-1}\alpha_{H}\left(p^{k}\mathbf{t}\right)\right)S_{0}\left(p^{h}\mathbf{t}\right)\left(H^{\prime}\left(\mathbf{0}\right)\right)^{-h}
\end{equation}
Here: 
\begin{equation}
S_{0}\left(\mathbf{t}\right)=\kappa_{H}\left(\mathbf{0}\right)=\mathbf{I}_{d}
\end{equation}
which leaves us with: 
\begin{equation}
S_{h}\left(\mathbf{t}\right)=p^{hr}\left(\prod_{k=0}^{h-1}\alpha_{H}\left(p^{k}\mathbf{t}\right)\right)\left(H^{\prime}\left(\mathbf{0}\right)\right)^{-h}\label{eq:Explicit formula for S_h in K_m convolution computation}
\end{equation}

Returning to (\ref{eq:MD Halfway done with K_m convolution computation}),
we have: 
\begin{align*}
\left(\hat{A}_{H}*K_{m}\right)\left(\mathbf{t}\right) & =\sum_{h=0}^{m-1}\frac{1}{p^{rm}}\sum_{\mathbf{n}=\mathbf{0}}^{p^{m-h}-1}\left(p^{rh}\left(\prod_{k=0}^{h-1}\alpha_{H}\left(p^{k}\mathbf{t}\right)\right)\left(H^{\prime}\left(\mathbf{0}\right)\right)^{-h}\right)e^{-2\pi i\mathbf{n}\cdot p^{h}\mathbf{t}}\left(H^{\prime}\left(\mathbf{0}\right)\right)^{h}\\
 & +\frac{1}{p^{rm}}p^{rm}\left(\prod_{k=0}^{m-1}\alpha_{H}\left(p^{k}\mathbf{t}\right)\right)\left(H^{\prime}\left(\mathbf{0}\right)\right)^{-m}\left(H^{\prime}\left(\mathbf{0}\right)\right)^{m}\left(\mathbf{I}_{d}-H^{\prime}\left(\mathbf{0}\right)\right)^{-1}\\
 & =\sum_{h=0}^{m-1}\left(\prod_{k=0}^{h-1}\alpha_{H}\left(p^{k}\mathbf{t}\right)\right)\frac{1}{p^{r\left(m-h\right)}}\sum_{\mathbf{n}=\mathbf{0}}^{p^{m-h}-1}e^{-2\pi i\mathbf{n}\cdot p^{h}\mathbf{t}}\\
 & +\left(\prod_{k=0}^{m-1}\alpha_{H}\left(p^{k}\mathbf{t}\right)\right)\left(\mathbf{I}_{d}-H^{\prime}\left(\mathbf{0}\right)\right)^{-1}
\end{align*}
Here: 
\begin{align*}
\sum_{\mathbf{n}=\mathbf{0}}^{p^{m-h}-1}e^{-2\pi i\mathbf{n}\cdot p^{h}\mathbf{t}} & =\sum_{\mathbf{n}=\mathbf{0}}^{p^{m-h}-1}e^{-2\pi i\frac{\mathbf{n}}{p^{m-h}}\cdot p^{m}\mathbf{t}}\\
 & =\sum_{\left\Vert \mathbf{s}\right\Vert _{p}\leq p^{m-h}}e^{-2\pi i\left(\mathbf{s}\cdot p^{m}\mathbf{t}\right)}\\
 & =p^{r\left(m-h\right)}\left[p^{m}\mathbf{t}\overset{p^{m-h}}{\equiv}\mathbf{0}\right]
\end{align*}
The congruence says that $p^{m}\mathbf{t}\in p^{m-h}\mathbb{Z}^{r}$.
This implies $p^{h}\mathbf{t}\in\mathbb{Z}^{r}$, which means: 
\begin{equation}
\left[p^{m}\mathbf{t}\overset{p^{m-h}}{\equiv}\mathbf{0}\right]=\left[p^{h}\mathbf{t}\overset{1}{\equiv}\mathbf{0}\right]=\mathbf{1}_{\mathbf{0}}\left(p^{h}\mathbf{t}\right)
\end{equation}
Hence:
\begin{align*}
\left(\hat{A}_{H}*K_{m}\right)\left(\mathbf{t}\right) & =\sum_{h=0}^{m-1}\left(\prod_{k=0}^{h-1}\alpha_{H}\left(p^{k}\mathbf{t}\right)\right)\frac{1}{p^{r\left(m-h\right)}}p^{m-h}\mathbf{1}_{\mathbf{0}}\left(p^{h}\mathbf{t}\right)\\
 & +\left(\prod_{k=0}^{m-1}\alpha_{H}\left(p^{k}\mathbf{t}\right)\right)\left(\mathbf{I}_{d}-H^{\prime}\left(\mathbf{0}\right)\right)^{-1}\\
 & =\sum_{h=0}^{m-1}\mathbf{1}_{\mathbf{0}}\left(p^{h}\mathbf{t}\right)\left(\prod_{k=0}^{h-1}\alpha_{H}\left(p^{k}\mathbf{t}\right)\right)+\left(\prod_{k=0}^{m-1}\alpha_{H}\left(p^{k}\mathbf{t}\right)\right)\left(\mathbf{I}_{d}-H^{\prime}\left(\mathbf{0}\right)\right)^{-1}
\end{align*}

Q.E.D. 
\begin{lem}
\label{lem:T_m hat of A_H hat}Let $m\in\mathbb{N}_{0}$. Then: 
\begin{equation}
\hat{\mathcal{T}}_{m}\left\{ \hat{A}_{H}\right\} \left(\mathbf{t}\right)=\sum_{n=0}^{m-1}\mathbf{1}_{\mathbf{0}}\left(p^{n}\mathbf{t}\right)\prod_{k=0}^{n-1}\alpha_{H}\left(p^{k}\mathbf{t}\right)+\left(\prod_{k=0}^{m-1}\alpha_{H}\left(p^{k}\mathbf{t}\right)\right)\left(\mathbf{I}_{d}-H^{\prime}\left(\mathbf{0}\right)\right)^{-1}\label{eq:T_m hat of A_H hat}
\end{equation}
where the $k$ product is defined to be $\mathbf{I}_{d}$ whenever
its upper limit is $<0$, and where the $m$-sum is defined to be
$\mathbf{O}_{d}$ when $m=0$. 
\end{lem}
Proof: Note that: 
\begin{align*}
\mathcal{T}_{m}\left\{ \tilde{A}_{H}\right\} \left(\mathbf{z}\right) & =\sum_{\mathbf{n}=\mathbf{0}}^{p^{m}-1}\tilde{A}_{H}\left(\mathbf{n}\right)\left[\mathbf{z}\overset{p^{m}}{\equiv}\mathbf{n}\right]\\
 & =\frac{1}{p^{rm}}\sum_{\mathbf{n}=\mathbf{0}}^{p^{m}-1}\tilde{A}_{H}\left(\mathbf{n}\right)\sum_{\left\Vert \mathbf{s}\right\Vert _{p}\leq p^{m}}e^{2\pi i\left\{ \mathbf{s}\left(\mathbf{z}-\mathbf{n}\right)\right\} _{p}}\\
 & =\frac{1}{p^{rm}}\sum_{\mathbf{n}=\mathbf{0}}^{p^{m}-1}\sum_{\mathbf{t}\in\hat{\mathbb{Z}}_{p}^{r}}\hat{A}_{H}\left(\mathbf{t}\right)e^{2\pi i\left\{ \mathbf{t}\mathbf{n}\right\} _{p}}\sum_{\left\Vert \mathbf{s}\right\Vert _{p}\leq p^{m}}e^{2\pi i\left\{ \mathbf{s}\left(\mathbf{z}-\mathbf{n}\right)\right\} _{p}}\\
 & =\sum_{\left\Vert \mathbf{s}\right\Vert _{p}\leq p^{m}}\sum_{\mathbf{t}\in\hat{\mathbb{Z}}_{p}^{r}}\hat{A}_{H}\left(\mathbf{t}\right)\frac{1}{p^{rm}}\sum_{\mathbf{n}=\mathbf{0}}^{p^{m}-1}e^{-2\pi i\left\{ \mathbf{n}\left(\mathbf{s}-\mathbf{t}\right)\right\} _{p}}e^{2\pi i\left\{ \mathbf{s}\mathbf{z}\right\} _{p}}\\
 & =\sum_{\left\Vert \mathbf{s}\right\Vert _{p}\leq p^{m}}\left(\hat{A}_{H}*\hat{K}_{m}\right)\left(\mathbf{s}\right)e^{2\pi i\left\{ \mathbf{s}\mathbf{z}\right\} _{p}}
\end{align*}
Then use \textbf{Proposition \ref{prop:A_H hat convolve K_H hat MD}}.

Q.E.D.

\vphantom{}

The last ingredient is the proposition given below. However, we will
need another bit of notation to help us along the way: 
\begin{notation}
Let $\mathbf{J}=\left(\mathbf{j}_{1},\ldots,\mathbf{j}_{\left|\mathbf{J}\right|}\right)\in\textrm{String}^{r}\left(p\right)$
be a $p$-block string. Then, we define: 
\begin{equation}
\mathbf{t}\cdot\mathbf{J}\overset{\textrm{def}}{=}\mathbf{t}\cdot\sum_{k=1}^{\left|\mathbf{J}\right|}\mathbf{j}_{k}p^{k-1}\label{eq:Definition of t dot big bold J}
\end{equation}
Writing each $\mathbf{j}_{k}$ as $\left(j_{k,1},\ldots,j_{k,r}\right)$,
we have that: 
\begin{equation}
\sum_{k=1}^{\left|\mathbf{J}\right|}\mathbf{j}_{k}p^{k-1}=\left(j_{1,1}+\cdots+j_{\left|\mathbf{J}\right|,1}p^{\left|\mathbf{J}\right|-1},\ldots,j_{1,r}+\cdots+j_{\left|\mathbf{J}\right|,r}p^{\left|\mathbf{J}\right|-1}\right)
\end{equation}
In particular, $\sum_{k=1}^{\left|\mathbf{J}\right|}\mathbf{j}_{k}p^{k-1}=\mathbf{n}$,
where $\mathbf{n}\in\mathbb{N}_{0}^{r}$ is the $r$-tuple represented
by $\mathbf{J}$. Additionally, we write: 
\begin{equation}
\sum_{\left|\mathbf{J}\right|=n}=\sum_{\mathbf{j}_{1},\ldots,\mathbf{j}_{n}\in\mathbb{Z}^{r}/p\mathbb{Z}^{r}}
\end{equation}
to denote a sum taken over all $\mathbf{J}\in\textrm{String}^{r}\left(p\right)$
of length exactly $n$. 
\end{notation}
\begin{prop}
Let $n\geq0$. Then, the Fourier transform of the locally constant
$\overline{\mathbb{Q}}^{d,d}$-valued function: 
\begin{equation}
\mathbf{z}\mapsto\mathbf{I}_{H}\left(\lambda_{p}\left(\left[\mathbf{z}\right]_{p^{n}}\right)\right)
\end{equation}
is given by: 
\begin{align}
\int_{\mathbb{Z}_{p}^{r}}\mathbf{I}_{H}\left(\lambda_{p}\left(\left[\mathbf{z}\right]_{p^{n}}\right)\right)e^{-2\pi i\left\{ \mathbf{t}\mathbf{z}\right\} _{p}}d\mathbf{z} & =\mathbf{1}_{\mathbf{0}}\left(\mathbf{t}\right)\mathbf{I}_{d}-\mathcal{C}_{H}\left(\alpha_{H}\left(\mathbf{0}\right):n\right)\hat{K}_{n}\left(\mathbf{t}\right)\label{eq:Fourier transform of I_H of lambda_P of z mod P^n}
\end{align}
\end{prop}
Proof: To keep the computations from running off the right side of
the page, let us write $\mathbf{A}=\alpha_{H}\left(\mathbf{0}\right)$
and $\mathbf{X}=H^{\prime}\left(\mathbf{0}\right)$. Then, the integral:
\begin{equation}
\int_{\mathbb{Z}_{p}^{r}}\mathbf{I}_{H}\left(\lambda_{p}\left(\left[\mathbf{z}\right]_{p^{n}}\right)\right)e^{-2\pi i\left\{ \mathbf{t}\mathbf{z}\right\} _{p}}d\mathbf{z}
\end{equation}
is: 
\begin{equation}
\int_{\mathbb{Z}_{p}^{r}}\left(\mathbf{I}_{d}-\mathbf{X}^{\lambda_{p}\left(\left[\mathbf{z}\right]_{p^{n}}\right)}\mathbf{A}\mathbf{X}^{-\lambda_{p}\left(\left[\mathbf{z}\right]_{p^{n}}\right)}\right)e^{-2\pi i\left\{ \mathbf{t}\mathbf{z}\right\} _{p}}d\mathbf{z}
\end{equation}
Because the main part of the integrand is a constant with respect
to values of $\mathbf{z}$ modulo $p^{n}$, the integral reduces to
a finite sum:
\begin{equation}
\mathbf{1}_{\mathbf{0}}\left(\mathbf{t}\right)\mathbf{I}_{d}-\frac{1}{p^{rn}}\sum_{\mathbf{m}=\mathbf{0}}^{p^{n}-1}\mathbf{X}^{\lambda_{p}\left(\mathbf{m}\right)}\mathbf{A}\mathbf{X}^{-\lambda_{p}\left(\mathbf{m}\right)}e^{-2\pi i\mathbf{t}\cdot\mathbf{\mathbf{m}}}
\end{equation}
Hopefully, the reader should not be surprised that we evaluate the
$\mathbf{m}$-sum using the recursive method. This time around, we
define: 
\begin{equation}
S_{n}\left(\mathbf{t}\right)\overset{\textrm{def}}{=}\frac{1}{p^{rn}}\sum_{\mathbf{m}=\mathbf{0}}^{p^{n}-1}\mathbf{X}^{\lambda_{p}\left(\mathbf{m}\right)}\mathbf{A}\mathbf{X}^{-\lambda_{p}\left(\mathbf{m}\right)}e^{-2\pi i\mathbf{t}\cdot\mathbf{m}}\label{eq:Definition of S_n for I_H computation}
\end{equation}
Then: 
\begin{align*}
S_{n}\left(\mathbf{t}\right) & =\frac{1}{p^{rn}}\sum_{\mathbf{j}=\mathbf{0}}^{p-1}\sum_{\mathbf{m}=\mathbf{0}}^{p^{n-1}-1}\mathbf{X}^{\lambda_{p}\left(p\mathbf{m}+\mathbf{j}\right)}\mathbf{A}\mathbf{X}^{-\lambda_{p}\left(p\mathbf{m}+\mathbf{j}\right)}e^{-2\pi i\mathbf{t}\cdot\left(p\mathbf{m}+\mathbf{j}\right)}\\
 & =\mathbf{X}\underbrace{\left(\frac{1}{p^{r\left(n-1\right)}}\sum_{\mathbf{m}=\mathbf{0}}^{p^{n-1}-1}\mathbf{X}^{\lambda_{p}\left(\mathbf{m}\right)}\mathbf{A}\mathbf{X}^{-\lambda_{p}\left(\mathbf{m}\right)}e^{-2\pi i\left(\mathbf{m}\cdot p\mathbf{t}\right)}\right)}_{S_{n-1}\left(p\mathbf{t}\right)}\mathbf{X}^{-1}\frac{1}{p^{r}}\sum_{\mathbf{j}=\mathbf{0}}^{p-1}e^{-2\pi i\mathbf{t}\cdot\mathbf{j}}
\end{align*}
So: 
\begin{equation}
S_{n}\left(\mathbf{t}\right)=\mathbf{X}\left(S_{n-1}\left(p\mathbf{t}\right)\right)\mathbf{X}^{-1}\frac{1}{p^{r}}\sum_{\mathbf{j}=\mathbf{0}}^{p-1}e^{-2\pi i\mathbf{t}\cdot\mathbf{j}}\label{eq:S_n recursion identity for I_H computation}
\end{equation}
Nesting this yields: 
\begin{equation}
S_{n}\left(\mathbf{t}\right)=\mathbf{X}^{n}\left(S_{0}\left(p^{n}\mathbf{t}\right)\right)\mathbf{X}^{-n}\frac{1}{p^{rn}}\sum_{\mathbf{j}_{1},\ldots,\mathbf{j}_{n}\leq p-1}e^{-2\pi i\mathbf{t}\cdot\sum_{k=1}^{n}\mathbf{j}_{k}p^{k}}
\end{equation}
Because:
\begin{equation}
S_{0}\left(\mathbf{s}\right)=\mathbf{X}^{\lambda_{p}\left(\mathbf{0}\right)}\mathbf{A}\mathbf{X}^{-\lambda_{p}\left(\mathbf{0}\right)}e^{-2\pi i\mathbf{s}\cdot\mathbf{0}}=\mathbf{A}
\end{equation}
we are then left with: 
\begin{equation}
S_{n}\left(\mathbf{t}\right)=\frac{1}{p^{rn}}\mathbf{X}^{n}\mathbf{A}\mathbf{X}^{-n}\sum_{\left|\mathbf{J}\right|=n}e^{-2\pi i\mathbf{t}\cdot\mathbf{J}}
\end{equation}

Finally, recall that the set of all $p$-block strings of length $n$
is in a bijective correspondence with the set of all $r$-tuples $\mathbf{k}\in\mathbb{N}_{0}^{r}$
satisfying $\mathbf{k}\leq p^{n}-1$. As such, we can write:
\begin{equation}
\sum_{\left|\mathbf{J}\right|=n}e^{-2\pi i\mathbf{t}\cdot\mathbf{J}}=\sum_{\mathbf{k}=\mathbf{0}}^{p^{n}-1}e^{-2\pi i\mathbf{t}\cdot\mathbf{k}}
\end{equation}
and hence:

\begin{equation}
S_{n}\left(\mathbf{t}\right)=\frac{1}{p^{rn}}\mathbf{X}^{n}\mathbf{A}\mathbf{X}^{-n}\sum_{\mathbf{k}=\mathbf{0}}^{p^{n}-1}e^{-2\pi i\mathbf{t}\cdot\mathbf{k}}\label{eq:Explicit formula for S_n for I_H computation}
\end{equation}

Putting everything together yields: 
\begin{align*}
\int_{\mathbb{Z}_{p}^{r}}\mathbf{I}_{H}\left(\lambda_{p}\left(\left[\mathbf{z}\right]_{p^{n}}\right)\right)e^{-2\pi i\left\{ \mathbf{t}\mathbf{z}\right\} _{p}}d\mathbf{z} & =\mathbf{1}_{\mathbf{0}}\left(\mathbf{t}\right)\mathbf{I}_{d}-\mathbf{X}^{n}\mathbf{A}\mathbf{X}^{-n}\frac{1}{p^{rn}}\sum_{\mathbf{k}=\mathbf{0}}^{p^{n}-1}e^{-2\pi i\mathbf{t}\cdot\mathbf{k}}\\
 & =\mathbf{1}_{\mathbf{0}}\left(\mathbf{t}\right)\mathbf{I}_{d}-\mathbf{X}^{n}\mathbf{A}\mathbf{X}^{-n}\hat{K}_{n}\left(\mathbf{t}\right)\\
 & =\mathbf{1}_{\mathbf{0}}\left(\mathbf{t}\right)\mathbf{I}_{d}-\mathcal{C}_{H}\left(\mathbf{A}:n\right)\hat{K}_{n}\left(\mathbf{t}\right)
\end{align*}

Q.E.D.

\vphantom{}

The way forward is then like so: we start our computation with:

\begin{equation}
\hat{\mathcal{L}}_{H,1,n}\left\{ \hat{A}_{H}\right\} \left(\mathbf{t}\right)=\hat{\mathcal{T}}_{n}\left\{ \hat{A}_{H}\right\} \left(\mathbf{t}\right)\left(\mathbf{I}_{d}-H^{\prime}\left(\mathbf{0}\right)\right)\left(H^{\prime}\left(\mathbf{0}\right)\right)^{-n}
\end{equation}
Using \textbf{Lemma \ref{lem:T_m hat of A_H hat}}, we get: 
\begin{align*}
\hat{\mathcal{L}}_{H,1,n}\left\{ \hat{A}_{H}\right\} \left(\mathbf{t}\right) & =\sum_{m=0}^{n-1}\mathbf{1}_{\mathbf{0}}\left(p^{m}\mathbf{t}\right)\left(\prod_{k=0}^{m-1}\alpha_{H}\left(p^{k}\mathbf{t}\right)\right)\left(H^{\prime}\left(\mathbf{0}\right)\right)^{-n}\left(\mathbf{I}_{d}-H^{\prime}\left(\mathbf{0}\right)\right)\\
 & +\left(\prod_{k=0}^{n-1}\alpha_{H}\left(p^{k}\mathbf{t}\right)\right)\left(H^{\prime}\left(\mathbf{0}\right)\right)^{-n}
\end{align*}
Next: 
\begin{equation}
\hat{\mathcal{L}}_{H,2,n}\left\{ \hat{A}_{H}\right\} \left(\mathbf{t}\right)=\mathcal{L}_{H,1,n}\left\{ \hat{A}_{H}\right\} \left(\mathbf{t}\right)-\mathcal{L}_{H,1,n-1}\left\{ \hat{A}_{H}\right\} \left(\mathbf{t}\right)\left(H^{\prime}\left(\mathbf{0}\right)\right)^{-1}
\end{equation}
and so: 
\begin{align*}
\hat{\mathcal{L}}_{H,2,n}\left\{ \hat{A}_{H}\right\} \left(\mathbf{t}\right) & =\mathbf{1}_{\mathbf{0}}\left(p^{n-1}\mathbf{t}\right)\left(\prod_{k=0}^{n-2}\alpha_{H}\left(p^{k}\mathbf{t}\right)\right)\left(H^{\prime}\left(\mathbf{0}\right)\right)^{-n}\left(\mathbf{I}_{d}-H^{\prime}\left(\mathbf{0}\right)\right)\\
 & +\left(\prod_{k=0}^{n-1}\alpha_{H}\left(p^{k}\mathbf{t}\right)\right)\left(H^{\prime}\left(\mathbf{0}\right)\right)^{-n}-\left(\prod_{k=0}^{n-2}\alpha_{H}\left(p^{k}\mathbf{t}\right)\right)\left(H^{\prime}\left(\mathbf{0}\right)\right)^{-n}\\
 & =\mathbf{1}_{\mathbf{0}}\left(p^{n-1}\mathbf{t}\right)\left(\prod_{k=0}^{n-2}\alpha_{H}\left(p^{k}\mathbf{t}\right)\right)\left(H^{\prime}\left(\mathbf{0}\right)\right)^{-n}\left(\mathbf{I}_{d}-H^{\prime}\left(\mathbf{0}\right)\right)\\
 & +\left(\prod_{k=0}^{n-2}\alpha_{H}\left(p^{k}\mathbf{t}\right)\right)\left(\alpha_{H}\left(p^{n-1}\mathbf{t}\right)-\mathbf{I}_{d}\right)\left(H^{\prime}\left(\mathbf{0}\right)\right)^{-n}
\end{align*}
Then: 
\begin{equation}
\hat{\mathcal{L}}_{H,3,n}\left\{ \hat{A}_{H}\right\} \left(\mathbf{t}\right)=\hat{\mathcal{E}}_{0}\left\{ \hat{A}_{H}\right\} \left(\mathbf{t}\right)+\sum_{m=1}^{n}\hat{\mathcal{L}}_{H,2,m}\left\{ \hat{A}_{H}\right\} \left(\mathbf{t}\right)
\end{equation}
where the $m$-sum is $\mathbf{O}_{d}$ when $n=0$. Recalling that:
\begin{equation}
\mathcal{E}_{0}\left\{ \tilde{A}_{H}\right\} \left(\mathbf{z}\right)=\mathbf{I}_{d}-\alpha_{H}\left(\mathbf{0}\right)
\end{equation}
we then have:
\begin{equation}
\hat{\mathcal{E}}_{0}\left\{ \hat{A}_{H}\right\} \left(\mathbf{t}\right)=\left(\mathbf{I}_{d}-\alpha_{H}\left(\mathbf{0}\right)\right)\mathbf{1}_{\mathbf{0}}\left(\mathbf{t}\right)
\end{equation}
which gives: 
\begin{align*}
\hat{\mathcal{L}}_{H,3,n}\left\{ \hat{A}_{H}\right\} \left(\mathbf{t}\right) & =\left(\mathbf{I}_{d}-\alpha_{H}\left(\mathbf{0}\right)\right)\mathbf{1}_{\mathbf{0}}\left(\mathbf{t}\right)\\
 & +\sum_{m=1}^{n}\mathbf{1}_{\mathbf{0}}\left(p^{m-1}\mathbf{t}\right)\left(\prod_{k=0}^{m-2}\alpha_{H}\left(p^{k}\mathbf{t}\right)\right)\left(H^{\prime}\left(\mathbf{0}\right)\right)^{-m}\left(\mathbf{I}_{d}-H^{\prime}\left(\mathbf{0}\right)\right)\\
 & +\sum_{m=1}^{n}\left(\prod_{k=0}^{m-2}\alpha_{H}\left(p^{k}\mathbf{t}\right)\right)\left(\alpha_{H}\left(p^{m-1}\mathbf{t}\right)-\mathbf{I}_{d}\right)\left(H^{\prime}\left(\mathbf{0}\right)\right)^{-m}
\end{align*}
where the $k$-products are $\mathbf{I}_{d}$ whenever $m=1$. 
\begin{equation}
\prod_{k=0}^{n-1}\alpha_{H}\left(p^{k}\mathbf{t}\right)=\begin{cases}
\hat{A}_{H}\left(\mathbf{t}\right)\left(\alpha_{H}\left(\frac{\mathbf{t}\left|\mathbf{t}\right|_{p}}{p}\right)\right)^{-1} & \textrm{if }n=-v_{p}\left(\mathbf{t}\right)-1\\
\hat{A}_{H}\left(\mathbf{t}\right)\left(\alpha_{H}\left(\mathbf{0}\right)\right)^{n+v_{p}\left(\mathbf{t}\right)} & \textrm{if }n\geq-v_{p}\left(\mathbf{t}\right)
\end{cases}
\end{equation}

Finally, note that: 
\begin{equation}
\mathcal{L}_{H,4,n}\left\{ \tilde{A}_{H}\right\} \left(\mathbf{z}\right)=\mathcal{L}_{H,3,n}\left\{ \tilde{A}_{H}\right\} \left(\mathbf{z}\right)\left(\mathbf{I}_{H}\left(\lambda_{p}\left(\left[\mathbf{z}\right]_{p^{n}}\right)\right)\right)^{-1}\left(H^{\prime}\left(\mathbf{0}\right)\right)^{n}
\end{equation}
Since the Fourier transform turns multiplication into convolution,
we get: 
\begin{equation}
\hat{\mathcal{L}}_{H,4,n}\left\{ \hat{A}_{H}\right\} \left(\mathbf{t}\right)=\left(\hat{\mathcal{L}}_{H,3,n}\left\{ \hat{A}_{H}\right\} *\hat{\mathbf{G}}_{n}^{-1}\right)\left(\mathbf{t}\right)\left(H^{\prime}\left(\mathbf{0}\right)\right)^{n}
\end{equation}
where $\hat{\mathbf{G}}_{n}^{-1}$ denotes the convolution inverse
of the Fourier transform of: 
\begin{equation}
\mathbf{G}_{n}\left(\mathbf{z}\right)\overset{\textrm{def}}{=}\mathbf{I}_{H}\left(\lambda_{p}\left(\left[\mathbf{z}\right]_{p^{n}}\right)\right)\label{eq:Definition of bold G_n}
\end{equation}
It is here that we hit a rocky obstacle. Because $\hat{\mathbf{G}}_{n}^{-1}\left(\mathbf{t}\right)$
does not appear to have a manageably simple explicit formula, determining
whether or not the sum: 
\begin{equation}
\hat{\mathcal{L}}_{H}\left\{ \hat{A}_{H}\right\} \left(\mathbf{t}\right)=\sum_{n=0}^{\infty}\hat{\mathcal{L}}_{H,4,n}\left\{ \hat{A}_{H}\right\} \left(\mathbf{t}\right)\label{eq:L_H hat of A_H hat}
\end{equation}
converges with respect to $\mathcal{F}_{p,q_{H}}^{d,d}$ becomes problematic.
Without an explicit formula with which convergence of the right-hand
side of (\ref{eq:L_H hat of A_H hat}) could be directly ascertained,
it would seem we are at the mercy of foremost weakness of fledgling
frame theory: the lack of a method of \emph{indirectly }verifying
convergence with respect to a given frame by way of estimates, approximation
arguments, and the like.

So, while it is intuitively clear that (\ref{eq:L_H hat of A_H hat})
should then furnish a formula for $\hat{\psi}_{H}\left(\mathbf{t}\right)$
in the case of non-commutative $H$ for which $\mathbf{I}_{d}-\alpha_{H}\left(\mathbf{0}\right)$
is invertible, at the time of this writing\textemdash the third day
of Putin's invasion of Ukraine (26 February, 2022)\textemdash I do
not have an argument to prove this rigorously. If and when such an
argument is obtained, it would immediately follow by \textbf{Theorem
\ref{thm:MD F-series for Chi_H}}'s formula $\chi_{H}\left(\mathbf{z}\right)=\psi_{H}\left(\mathbf{z}\right)\beta_{H}\left(\mathbf{0}\right)+\Psi_{H}\left(\mathbf{z}\right)$
that $\chi_{H}$ would be quasi-integrable. As such, I end with a
conjecture: 
\begin{conjecture}[\textbf{Quasi-Integrability of $\chi_{H}$ in the Non-Commutative
Case}]
\label{conj:MD Non-Commutative case conjecture}Let $H$ be as given
in \textbf{\emph{Theorem \ref{thm:MD F-series for Chi_H}}}, and suppose
that $\mathbf{I}_{d}-\alpha_{H}\left(\mathbf{0}\right)$ is invertible.
Then,\emph{ (\ref{eq:L_H hat of A_H hat})} gives a formula $\hat{\mathbb{Z}}_{p}^{r}\rightarrow\overline{\mathbb{Q}}^{d,d}$
for a Fourier transform of $\psi_{H}$ with respect to the standard
$\left(p,q_{H}\right)$-adic frame.
\end{conjecture}

\chapter*{Coda}

\pagestyle{fancy}\fancyfoot{}\fancyhead[L]{\sl CODA}\fancyhead[R]{\thepage}\addcontentsline{toc}{chapter}{Coda}

\includegraphics[scale=0.4]{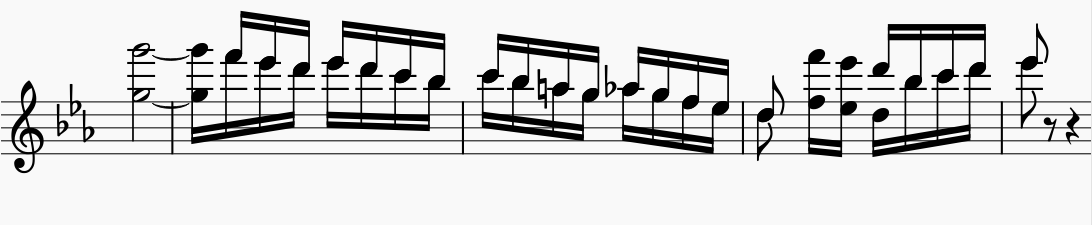}

\section*{Conclusion - A Call to Arms}

\addcontentsline{toc}{section}{Conclusion -   A Call to Arms}

Even though insight drives mathematical progress, the utility and
efficacy of mathematics is only knowable in hindsight. Will the methods
I have presented in my dissertation be of any use in solving the Collatz
Conjecture once and for all? The answer\textemdash if it is to ever
be found\textemdash lies in posterity. I would very much like to meet
it, if I can.

In writing this dissertation, the two biggest surprises were the connection
to eigenvalues of matrices in Sub-subsection \ref{subsec:A-Matter-of}\textemdash although,
in hindsight this feels obvious (of course)\textemdash and, most of
all, the result that divergent trajectories of $H$ were associated
with values in $\mathbb{Z}$ attained by $\chi_{H}$ over $\mathbb{Z}_{p}\backslash\mathbb{Q}$
(\textbf{Theorem \ref{thm:Divergent trajectories come from irrational z}}
on page \pageref{thm:Divergent trajectories come from irrational z},
and \textbf{Theorem \ref{thm:MD Divergent trajectories come from irrational z}}
on page \pageref{thm:MD Divergent trajectories come from irrational z}).
I only discovered $\chi_{H}$'s involvement in divergent trajectories
in February 2022 while I was beginning my polishing and editing runs
on this monograph. Prior to that, I was firmly convinced that my methods
had nothing to say about divergent trajectories. For once, I was \emph{thrilled}
to be proven wrong.

Regarding the theory of Hydra maps as I have presented it, the most
significant outstanding issue are the \textbf{Conjectures} \textbf{\ref{conj:correspondence theorem for divergent trajectories}}
(page \pageref{conj:correspondence theorem for divergent trajectories})
and \textbf{\ref{conj:MD correspondence theorem for divergent trajectories}}
(page \pageref{conj:MD correspondence theorem for divergent trajectories})
regarding a Correspondence Principle for divergent trajectories. If
these Conjectures could be proven true, the \textbf{Tauberian Spectral
Theorems} (pages \pageref{thm:Periodic Points using WTT} \& \pageref{thm:MD Periodic Points using WTT})
for $\chi_{H}$ would then \emph{completely} characterize the dynamics
of the Hydra maps under consideration\textemdash that is, both periodic
points \emph{and }divergent trajectories. The entire matter of the
Collatz Conjecture (not to mention most any comparable conjectures
for other Hydra maps) would then be reduced to the study of the density
of translates of $\hat{\chi}_{H}\left(t\right)-x\mathbf{1}_{0}\left(t\right)$
in $c_{0}\left(\hat{\mathbb{Z}}_{p},\mathbb{C}_{q}\right)$ and its
multi-dimensional analogue. At a lesser level, there are the matters
of one-dimensional $p^{n}$-Hydra maps (where $p$ is a prime and
$n\in\mathbb{N}_{1}$) and multi-dimensional $P$-Hydra maps where
$P=\left(p^{n_{1}},\ldots,p^{n_{r}}\right)$ where $p$ is prime and
$n_{1}\leq\ldots\leq n_{r}$ is a non-decreasing sequence of positive
integers.

While I certainly \emph{like} to think that my Tauberian Spectral
reformulation of Hydra maps' dynamics is aesthetically pleasing, I
have yet to explore whether or not it is actually \emph{useful}. The
question of the density of the translates of $\hat{\chi}_{H}\left(t\right)-x\mathbf{1}_{0}\left(t\right)$
might very well turn out to be just as intractable as the original
question of whether or not $x$ is a periodic point or divergent trajectory
of $H$. Although, obviously, I am biased in favor of believing that
my work \emph{will} turn out to be useful, I cannot help but feel
that the distinctions between, say, Fourier transforms of the Shortened
$3x+1$ map (for which $\alpha_{H}\left(0\right)=1$) and those of
the Shortened $5x+1$ map (for which $\alpha_{H}\left(0\right)\neq1$)
have the look of a smoking gun. For reference, valid choices of Fourier
transforms for these numina are:
\begin{equation}
\hat{\chi}_{3}\left(t\right)=\begin{cases}
-\frac{1}{2} & \textrm{if }t=0\\
\frac{1}{4}v_{2}\left(t\right)\hat{A}_{3}\left(t\right) & \textrm{if }t\neq0
\end{cases},\textrm{ }\forall t\in\hat{\mathbb{Z}}_{2}
\end{equation}
\begin{equation}
\hat{\chi}_{5}\left(t\right)=\begin{cases}
-\frac{1}{2} & \textrm{if }t=0\\
-\frac{1}{4}\hat{A}_{5}\left(t\right) & \textrm{if }t\neq0
\end{cases}=-\frac{1}{4}\mathbf{1}_{0}\left(t\right)-\frac{1}{4}\hat{A}_{5}\left(t\right),\textrm{ }\forall t\in\hat{\mathbb{Z}}_{2}
\end{equation}
Because $v_{2}\left(t\right)$ takes values in $\left\{ -1,-2,-3,\ldots\right\} $
for $t\in\hat{\mathbb{Z}}_{2}\backslash\left\{ 0\right\} $, the $3$-adic
absolute value of $\hat{\chi}_{3}\left(t\right)$ will become arbitrarily
small infinitely often:
\begin{equation}
\liminf_{\left|t\right|_{2}\rightarrow\infty}\left|\hat{\chi}_{3}\left(t\right)\right|_{3}=0
\end{equation}
whereas: 
\begin{equation}
\left|\hat{\chi}_{5}\left(t\right)\right|_{5}=1,\textrm{ }\forall t\in\hat{\mathbb{Z}}_{2}
\end{equation}
Unlike the probabilistic heuristics usually given to argue that $3x+1$
sends all positive integers to $1$ or that $5x+1$ sends almost all
positive integers to $\infty$, this observation about the behaviors
of $\hat{\chi}_{3}$ and $\hat{\chi}_{5}$ holds with \emph{absolute}
certainty. This is but one reason for my my conviction that $\hat{\chi}_{3}\left(t\right)$'s
failure to be $3$-adically bounded away from $0$ \emph{must }be
a key feature in any proof of the Collatz Conjecture, should one arise.
To that end, it would be wonderful if we could establish something
along the lines of:
\begin{conjecture}
\label{conj:Implication of Tauberian Spectral Theorem}Let $H:\mathbb{Z}\rightarrow\mathbb{Z}$
be a contracting, semi-basic $p$-Hydra map, and let $\hat{\chi}_{H}$
be a Fourier transform of $\chi_{H}$.

\vphantom{}

I. $H$ has finitely many periodic points if and only if $\liminf_{\left|t\right|_{p}\rightarrow\infty}\left|\hat{\chi}_{H}\left(t\right)\right|_{q_{H}}=0$.

\vphantom{}

II. If $\liminf_{\left|t\right|_{p}\rightarrow\infty}\left|\hat{\chi}_{H}\left(t\right)\right|_{q_{H}}>0$,
then $H$ has a divergent trajectory.
\end{conjecture}
\vphantom{}

While the particulars of the hypotheses and conclusions of this conjecture
might very well need to be tinkered with, my hope in making it will
be that questions about the dynamics of $H$ can be reduced to number-theoretic
($q$-adic) properties of the formulae for $\hat{\chi}_{H}$.
\begin{example}[\textbf{Matthews' Map, Revisited}]
\index{Matthews' Conjecture}Of special interest is the $3$-Hydra
map obtained by conjugating the map featured in \textbf{Matthews'
Conjecture} (\textbf{Conjecture \ref{conj:Matthews conjecture}} on
page \pageref{conj:Matthews conjecture}); specifically, the conjugated
version $\tilde{M}$ defined in equation (\ref{eq:Matthews' Conjecture Map, conjugated})
on page \pageref{exa:Matthews' map}:
\begin{equation}
\tilde{M}\left(n\right)=\begin{cases}
\frac{n}{3} & \textrm{if }n\overset{3}{\equiv}0\\
7n-3 & \textrm{if }n\overset{3}{\equiv}1\\
\frac{7n-2}{3} & \textrm{if }n\overset{3}{\equiv}2
\end{cases}
\end{equation}
The reason this map is of interest, recall, is because it is easily
proved that every integer congruent to $1$ mod $3$ belongs to a
divergent trajectory of $\tilde{M}$ (a simple modification of page
\pageref{prop:Matthews' map}'s \textbf{Proposition \ref{prop:Matthews' map}}).
Though \emph{non-integral}, this $3$-Hydra is nevertheless semi-basic
(with $\mu_{0}=1$, $\mu_{1}=21$, $\mu_{2}=7$, and $q=7$) and so
it possesses a numen $\chi_{\tilde{M}}$, and every periodic point
of $\tilde{M}$ is of the form $\chi_{\tilde{M}}\left(B_{3}\left(n\right)\right)$
for some $n\in\mathbb{N}_{1}$. However, because $\tilde{M}$ is non-integrable,
the versions of the Correspondence Principal established in this dissertation
do not apply to conclude that every value in $\mathbb{Z}$ attained
by $\chi_{\tilde{M}}$ over $\mathbb{Q}\cap\mathbb{Z}_{3}^{\prime}$
is a periodic point of $\tilde{M}$, or that $\mathfrak{z}\in\mathbb{Z}_{3}\backslash\mathbb{Q}$
for which $\chi_{\tilde{M}}\left(\mathfrak{z}\right)\in\mathbb{Z}$
necessarily make $\chi_{\tilde{M}}\left(\mathfrak{z}\right)$ into
a divergent point of $\tilde{M}$. Nevertheless, I have not explored
this issue in detail\textemdash my goal in my dissertation was to
obtain as broad a theory as possible, rather than whittle away at
any particular Hydra map\textemdash so there may still be a way to
make these results applicable to $\tilde{M}$.

Despite these obstacles, nothing prevents us from going through with
Chapter 4's Fourier analysis of $\chi_{\tilde{M}}$. Doing so, we
find that:
\begin{equation}
\alpha_{\tilde{M}}\left(t\right)\overset{\textrm{def}}{=}\frac{1+21e^{-2\pi it}+7e^{-4\pi it}}{9},\textrm{ }\forall t\in\hat{\mathbb{Z}}_{3}\label{eq:alpha for conjugated matthews map}
\end{equation}
\begin{equation}
\beta_{\tilde{M}}\left(t\right)\overset{\textrm{def}}{=}e^{-2\pi it}-\frac{2}{9}e^{-4\pi it},\textrm{ }\forall t\in\hat{\mathbb{Z}}_{3}\label{eq:Beta for conjugated Matthews map}
\end{equation}
This gives:
\begin{align}
\alpha_{\tilde{M}}\left(0\right) & =\frac{29}{9}\\
\beta_{\tilde{M}}\left(0\right) & =\frac{7}{9}
\end{align}
with:
\begin{equation}
\hat{A}_{\tilde{M}}\left(t\right)\overset{\textrm{def}}{=}\begin{cases}
1 & \textrm{if }t=0\\
\prod_{n=0}^{-v_{3}\left(t\right)-1}\frac{1+21e^{-2\pi i3^{m}t}+7e^{-4\pi i3^{m}t}}{9} & \textrm{else}
\end{cases},\textrm{ }\forall t\in\hat{\mathbb{Z}}_{3}\label{eq:A_H-hat for Matthews map}
\end{equation}
and:
\begin{equation}
\kappa_{\tilde{M}}\left(n\right)=\left(21\right)^{\#_{3:1}\left(n\right)}7^{\#_{3:2}\left(n\right)}\label{eq:Kappa_H for matthews map}
\end{equation}
Using \textbf{Theorem \ref{thm:F-series for an arbitrary 1D Chi_H}}
(page \pageref{thm:F-series for an arbitrary 1D Chi_H}), a Fourier
transform for $\chi_{\tilde{M}}:\mathbb{Z}_{3}\rightarrow\mathbb{Z}_{7}$
is given by:
\begin{equation}
\hat{\chi}_{\tilde{M}}\left(t\right)=\begin{cases}
-\frac{7}{20} & \textrm{if }t=0\\
\left(\frac{9e^{-2\pi i\frac{t\left|t\right|_{3}}{3}}-2e^{-4\pi i\frac{t\left|t\right|_{3}}{3}}}{1+21e^{-2\pi i\frac{t\left|t\right|_{3}}{3}}+7e^{-4\pi i\frac{t\left|t\right|_{3}}{3}}}-\frac{7}{20}\right)\hat{A}_{\tilde{M}}\left(t\right) & \textrm{if }t\neq0
\end{cases}\label{eq:Fourier Transform for Chi_H for Matthews map}
\end{equation}

Here, number theory comes into play. Noting that:
\begin{equation}
\alpha_{\tilde{M}}\left(t\right)=\frac{1+\overbrace{21e^{-2\pi it}+7e^{-4\pi it}}^{\textrm{a multiple of }7}}{9}\in\frac{1}{9}+7\mathbb{C}_{7},\textrm{ }\forall t\in\hat{\mathbb{Z}}_{3}
\end{equation}
the $7$-adic ultrametric inequality yields:
\begin{equation}
\left|\alpha_{\tilde{M}}\left(t\right)\right|_{7}=1,\textrm{ }\forall t\in\hat{\mathbb{Z}}_{3}
\end{equation}
As such:
\begin{equation}
\left|\hat{A}_{\tilde{M}}\left(t\right)\right|_{7}=\prod_{n=0}^{-v_{3}\left(t\right)-1}\left|\alpha_{\tilde{M}}\left(3^{n}t\right)\right|_{7}=1,\textrm{ }\forall t\in\hat{\mathbb{Z}}_{3}
\end{equation}
So, the $7$-adic absolute value of $\hat{\chi}_{\tilde{M}}\left(t\right)$
for non-zero $t$ is entirely determined by the function:
\begin{equation}
\frac{9e^{-2\pi i\frac{t\left|t\right|_{3}}{3}}-2e^{-4\pi i\frac{t\left|t\right|_{3}}{3}}}{1+21e^{-2\pi i\frac{t\left|t\right|_{3}}{3}}+7e^{-4\pi i\frac{t\left|t\right|_{3}}{3}}}-\frac{7}{20}\label{eq:7-adic magnitude determinator for Chi_H for Matthews map}
\end{equation}

Because $t\left|t\right|_{3}/3$ is congruent mod $1$ to either $1/3$
(if the numerator of $t$ is $1$ mod $3$) or $2/3$ (if the numerator
of $t$ is $2$ mod $3$) for any $t\in\hat{\mathbb{Z}}_{3}\backslash\left\{ 0\right\} $,
(\ref{eq:7-adic magnitude determinator for Chi_H for Matthews map})
takes precisely two values on $\hat{\mathbb{Z}}_{3}\backslash\left\{ 0\right\} $.
Letting $\xi$ denote $e^{2\pi i/6}$, and letting $\zeta$ denote
$\xi^{2}$, these are:
\begin{equation}
\frac{9e^{-2\pi i\frac{t\left|t\right|_{3}}{3}}-2e^{-4\pi i\frac{t\left|t\right|_{3}}{3}}}{1+21e^{-2\pi i\frac{t\left|t\right|_{3}}{3}}+7e^{-4\pi i\frac{t\left|t\right|_{3}}{3}}}-\frac{7}{20}=\begin{cases}
\frac{17-47\zeta^{2}}{1580} & \textrm{if }t\left|t\right|_{3}\overset{3}{\equiv}1\\
\frac{17-47\zeta}{1580} & \textrm{if }t\left|t\right|_{3}\overset{3}{\equiv}2
\end{cases}
\end{equation}
Recalling that the torsion subgroup of $\mathbb{Z}_{7}^{\times}$
is isomorphic to $\mathbb{Z}/6\mathbb{Z}$, all $6$th roots of unity
are occur naturally as $7$-adic integers. As such, to compute the
$7$-adic absolute values of these expressions, we will need to fall
back on our convention that, for an odd prime $p$, $\xi\in\mathbb{Z}_{p}$
denotes the unique primitive $\left(p-1\right)$th root of unity whose
first $p$-adic digit is the smallest primitive root of unity in $\mathbb{Z}/p\mathbb{Z}$.
Since $p=7$, that digit is $3$. Consequently, $\zeta=\xi^{2}$ has
$3^{2}\overset{7}{\equiv}2$ as its first $7$-adic digit and $\zeta^{2}=\xi^{4}$
has $3^{4}\overset{7}{\equiv}4$ as its first $7$-adic digit. Observing
the identity:
\begin{equation}
\left(17-47\zeta^{2}\right)\left(17-47\zeta\right)=\left(17\right)^{2}+\left(47\right)^{2}+17\cdot47=3297=3\cdot7\cdot157
\end{equation}
we have that:
\begin{equation}
17-47\zeta^{2}\overset{7}{\equiv}17-\left(47\right)\left(4\right)\overset{7}{\equiv}3-\left(5\right)\left(2\right)=-17\overset{7}{\equiv}4
\end{equation}
So, $\left|17-47\zeta^{2}\right|_{7}=1$. As such:
\begin{equation}
\left|\left(17-47\zeta^{2}\right)\left(17-47\zeta\right)\right|_{7}=\frac{1}{7}
\end{equation}
then forces: 
\begin{equation}
\left|17-47\zeta\right|_{7}=\frac{1}{7}
\end{equation}
Finally, noting that:
\begin{equation}
1580=2^{2}\cdot5\cdot79
\end{equation}
we have:
\begin{align}
\left|\frac{17-47\zeta}{1580}\right|_{7} & =\frac{1}{7}\\
\left|\frac{17-47\zeta^{2}}{1580}\right|_{7} & =1
\end{align}
Since $\left|\hat{\chi}_{\tilde{M}}\left(0\right)\right|_{7}=\left|-7/20\right|_{7}=1/7$,
we then conclude that:
\begin{equation}
\left|\hat{\chi}_{\tilde{M}}\left(t\right)\right|_{7}=\begin{cases}
\frac{1}{7} & \textrm{if }t=0\\
1 & \textrm{if }t\left|t\right|_{3}\overset{3}{\equiv}1\\
\frac{1}{7} & \textrm{if }t\left|t\right|_{3}\overset{3}{\equiv}2
\end{cases},\textrm{ }\forall t\in\hat{\mathbb{Z}}_{7}\label{eq:7-adic absolute value of the fourier transform of the numen of Matthews map}
\end{equation}
With this, we see that $\chi_{\tilde{M}}$ is more like $\chi_{q}$
for $q\geq5$, in that its $7$-adic absolute value is bounded away
from $0$. This is consistent with \textbf{Conjecture \ref{conj:Implication of Tauberian Spectral Theorem}}\textemdash or,
at least, with this particular formulation of the Conjecture. As mentioned
above, the precise statement of that Conjecture is still subject to
revision; more specific cases should be investigated to help see what
the ``correct'' conjecture might be.

Subject to the verification of \textbf{Conjecture \ref{conj:Implication of Tauberian Spectral Theorem}},
the method described above\textemdash checking $\hat{\chi}_{H}$ and
the values of $\gamma_{H}\left(j/p\right)$ for $j\in\left\{ 1,\ldots,p-1\right\} $
over the $q$-adics\textemdash would then be an extremely elegant
means for determining the dynamical properties of $H$.
\end{example}
\vphantom{}

The two other most pressing issues I can think of are the \textbf{exploration
of Hydra maps on fields of positive characteristic}, such as those
presented near the end of Matthews' slides \cite{Matthews' slides},
and the \textbf{exploration of the polygenic case} as discussed in
\textbf{Example \ref{exa:Polygenic example, part 1}} on page \pageref{exa:Polygenic example, part 1}.
As mentioned in the remark given after \textbf{Example \ref{exa:Polygenic example, part 1}},
I believe it will be possible to use frames to put $\chi_{H}$ on
rigorous footing in the polygenic case. Because the $N$th truncations
of $\chi_{H}$ take values in $\mathbb{Q}$, the asymptotic analysis
of $\chi_{H,N}$\textemdash untangling the link between $N$ and $t$\textemdash will
yield Fourier transforms for polygenic $\chi_{H}$. The only possible
issue I can foresee is in verifying that the resultant formulae for
$\hat{\chi}_{H}$ form continuous linear functionals\textemdash although,
there, the question becomes \emph{on what space?} Provided this can
be answered in a straight-forward matter, the $\left(p,q\right)$-adic
Wiener Tauberian Theorem for measures will almost certainly apply
as in the monogenic case, thereby furnishing a Tauberian Spectral
Theorem for polygenic $\chi_{H}$, seeing as how they already depend
on the points at which the measure's Fourier series converge $q$-adically.

Aside from these obvious next steps, it also seems worthwhile to explore
if $L_{\mathbb{R}}^{1}$ can, in fact, be used as a ``base of operations''
for analytic investigations of quasi-integrability and frames.

\subsection*{A Soapbox Moment}

Before I begin the obligatory bibliographic essay, I would like to
make a ``call to arms'' for fellow travelers in mathematical analysis
and its many subspecialties. Although I am pleased to have broken
new ground in non-archimedean analysis\textemdash without it, I would
not have attained my hard-earned PhD\textemdash it troubles me that
some of the phenomena I have discovered appear to haven't already
been discovered, especially considering their astonishing simplicity.
The identity (\ref{eq:Fourier sum of v_p of t}) (page \pageref{eq:Fourier sum of v_p of t}):
\[
\sum_{0<\left|t\right|_{p}\leq p^{N}}v_{p}\left(t\right)e^{2\pi i\left\{ t\mathfrak{z}\right\} _{p}}\overset{\mathbb{Q}}{=}\frac{p\left|\mathfrak{z}\right|_{p}^{-1}-1}{p-1},\textrm{ }\forall N>v_{p}\left(\mathfrak{z}\right)
\]
requires no extraordinary machinery to prove. The computation is as
near as mindless as can be, yet\textemdash to my bemusement\textemdash no
one seems to have noticed it until now. The mathematical community's
failure to notice this identity, as well as the implication it has
for $\left(p,q\right)$-adic integration ($\left(p,q\right)$-adic
Mellin transforms, distributional derivatives, etc.) cannot be due
to a lack of genius or ingenuity; it can only be attributed to a lack
of curiosity and an aversion to simplicity. It frustrates me to no
end that this state of affairs is allowed to stand.

I have abiding empathy for Vladimir Arnold's protests against contemporary
mathematical pedagogy \cite{Arnold}. As the venerable Russian analyst
put it: 
\begin{quote}
Attempts to create ``pure'' deductive-axiomatic mathematics have
led to the rejection of the scheme used in physics (observation -
model - investigation of the model - conclusions - testing by observations)
and its substitution by the scheme: definition - theorem - proof. 
\end{quote}
I could not agree with this sentiment more strongly. The deductive-axiomatic
scheme is a powerful organizational tool, but I feel it horribly misrepresents
the reality of mathematical inquiry. In the research that led to my
dissertation, I proceeded very much like a physicist of old, testing
different ideas and formulas in the hopes of finding something that
managed to hold water and say something useful about the objects of
my investigations. Make no mistake: I am not a mindless iconoclast.
The trend toward abstraction has and always will be a vital part of
mathematical inquiry. In clearing the road of debris and obstacles,
they make it easier to discern deep patterns and unexpected symmetries.
But, it is my conviction that flights of fancy such as these need
to be matched and tempered by good old-fashioned concrete play. We
need to be able to get our hands dirty and loamy with tedious, self-limiting
specificity. Euler spent a decade trying to prove his famous (and
beautiful) generating function identity for pentagonal numbers. If
only us moderns had the tenacity to give simple, specific questions
that same level of consideration.

Progress comes about when people look in places they shouldn't, and
when ideas take root far from home. Theories stagnate when all their
voices sound the same. Practical-minded analysts should not let topics
like non-archimedean analysis fall by the wayside and become sole
purview of the algebraists. The existence of the topological, open-set-based
definition of continuity does not diminish the $\epsilon$-$\delta$
definition, nor render it obsolete. The same ought to be true of non-archimedean
analysis and classical analysis. Just because a thing can be tied
to the cross of commutative algebra does not mean there is nothing
to be gained from an independent examination of the thing from a purely
concrete, analytical perspective.

Even if this dissertation of mine ends up being a mere tangent on
the long road to Collatz, if it can succeed at enticing other scholars
to consider topics or settings they might not have bothered to ever
explore, I can say I will be happy with it, with my efforts, and with
myself.

\vphantom{}
\begin{quote}
\emph{Oh me! Oh life! of the questions of these recurring,}

\emph{Of the endless trains of the faithless, of cities fill\textquoteright d
with the foolish,}

\emph{Of myself forever reproaching myself, (for who more foolish
than I, and who more faithless?)}

\emph{Of eyes that vainly crave the light, of the objects mean, of
the struggle ever renew\textquoteright d,}

\emph{Of the poor results of all, of the plodding and sordid crowds
I see around me,}

\emph{Of the empty and useless years of the rest, with the rest me
intertwined,}

\emph{The question, O me! so sad, recurring\textemdash What good amid
these, O me, O life?}
\begin{center}
Answer.
\par\end{center}
\emph{That you are here\textemdash that life exists and identity,}

\emph{That the powerful play goes on, and you may contribute a verse.}

\textemdash Walt Whitman (1892)
\end{quote}
\newpage{}

\section*{Bibliographic Essay}

\addcontentsline{toc}{section}{Bibliographic Essay}As explained in
the Preface, the $\left(p,q\right)$-adic analytic methods used in
this dissertation were not the first path of attack I tried on the
Collatz Conjecture, and I doubt they will be the last. To that end,
I have not exercised much restraint with regard to the references
catalogued in the Bibliography below; it is my hope that others will
find the listed sources as interesting as I did. Like most bibliographic
essays, the next few paragraphs are intended to assist the reader
in assessing the resources I have left for them. For ease of access,
I have organized my discussion by topic.

\subsubsection*{Collatz Studies}

At the time of writing, ``Collatz studies'' has yet to become a
well-established mathematical discipline, and the scattered, amorphous
state of the literature reflects this. Lagarias' \emph{Ultimate Challenge
}book \cite{Ultimate Challenge} is an wonderful start, as are his
annotated bibliographies of Collatz research \cite{Collatz Biography}.
As of the first quarter of the twenty-first century, Lagarias is generally
considered the top authority in Collatz studies, and to that end,
all of his publications on the topic are worth reading (\cite{Applegate and Lagarias - Trees,Applegate and Lagarias - Difference inequalities,3x+1 semigroup,natural boundaries,Lagarias-Kontorovich Paper,Lagarias' Survey,Lagarias - rational cycles for 3x+1}).
There are doubtless others to be found, should the reader spend the
time searching for them.

Matthews' \href{http://www.numbertheory.org/3x\%2B1/}{Collatz webpage}\footnote{http://www.numbertheory.org/3x+1/}
website contains links to his slides and related articles (\cite{Matthews' Conjecture,Matthews' Leigh Article,Matthews' slides,Matthews and Watts});
these focus primarily on the Markov-chain approach. Also of note,
the late Meinardus and Berg \cite{Meinardus and Berg,Meinardus} attempted
analytic approaches using functional equations, as did their colleague
G. Opfer (\cite{Berg =00003D000026 Opfer,Opfer}), with Opfer's 2011
paper \cite{Opfer} being notable as having been thought to have actually\emph{
proved }Collatz until a gap was discovered (see, for instance, \cite{de Weger on Opfer}).
There is also my own (albeit somewhat messy) work in this vein \cite{Dreancatchers for Hydra Maps}.
My bibliography only scratches the surface of the extant literature.

Eric Roosendaal maintains a $3x+1$ website dedicated to chronicling
various computational phenomena and statistics thereof found in the
iteration of the Collatz map \cite{Roosendaal's Website}. Then, of
course, there is Tao's recent work \cite{Tao Probability paper}.

The reader should be aware that most Collatz studies tend to focus
solely on the $3x+1$ map itself. One notable exception to this trend
is \cite{RCWAG}, which is part of a computational package for studying
more general Collatz-type maps\textemdash the ungainly (though accurately)
termed ``residue-class-wise affine maps'' (RCWA)\textemdash all
due to Stephan Kohl. \cite{RCWAG} is a freely available computational
software for studying Collatz type maps, written using the GAP programming
language. It strikes me as the sort of thing a person might fool around
with\footnote{Particularly if it can be used to create attractive, eye-catching
visualizations. As they say, a picture is worth a thousand words.
It also makes for a very effective tool for marketing and recruitment.} while lounging about at home, either on a rainy day or a golden Sunday
afternoon.

\subsubsection*{Non-Archimedean Analysis (Including $\left(p,q\right)$)}

W. M. Schikhof PhD dissertation \cite{Schikhof's Thesis} establishes
the fundamentals of harmonic analysis in the general non-archimedean
setting, including the $\left(p,q\right)$-adic case. A particularly
valuable\textemdash and accessible!\textemdash introduction to non-archimedean
analysis is Schikhof's \emph{Ultrametric Calculus }\cite{Ultrametric Calculus},
a classic of the subject. The bulk of it is focused on $\left(p,p\right)$-adic
analysis, though. A treatment of the Monna-Springer integral (though
not by that name) is given in one of \emph{Ultrametric Calculus'}
appendices. A \emph{much }more exhaustive collection of non-archimedean
material can be found in van Rooij's \emph{Non-Archimedean Functional
Analysis }\cite{van Rooij - Non-Archmedean Functional Analysis},
though\textemdash as I have mentioned elsewhere\textemdash this book
is out of print\footnote{It \emph{can}, however, be pirated.}. van
Rooij's vantage is quite general, to the point that his book is a
bit bewildering to work with. Those without photographic memories
are advised to keep a glossary of his definitions and notations.

Khrennikov's writings (especially the wild ride of a book in \cite{Quantum Paradoxes})
are very much worth the read, both for sheer entertainment value,
and because of their tendency to be more concrete than van Rooij's
presentation or that of other scholars: Aguayo (\cite{Aguayo 1,Aguayo 2,Aguayo and Moraga radon-nikodym theorem paper}),
and\textemdash most egregiously\textemdash Ludkovsky terribly dense
writing (see \cite{Ludkovsky on non-archimedean measures}, for example).
Khrennikov's co-authored paper \cite{Measure-theoretic approach to p-adic probability theory}
on non-archimedean probability theory is also an excellent, self-contained
introduction to the Monna-Springer integration technique presented
here, supplemented with additional considerations regarding the formulation
of measure-theoretic approaches to probability theory for probabilities
taking values in non-archimedean fields.

As a rule of thumb, much of the literature in non-archimedean analysis
(as opposed to $p$-adic analysis proper) is couched in the language
of abstract functional analysis. Some may find that register of presentation
elegant or high-brow; I find it frustratingly un-welcoming, especially
for newcomers. This is particularly important for the reader to keep
in mind, seeing as the nature of the work I have done in my dissertation
is foundationally \emph{concrete}.

As mentioned in Subsection \ref{subsec:3.1.1 Some-Historical-and},
a would-be non-archimedean analyst should be aware of 'false friends'
like \cite{Bosch lying title,Schneider awful book,More Schneider lies}
which are works of algebraic geometry prancing about with analysis-sounding
titles.

\subsubsection*{$\left(p,\infty\right)$-Adic Analysis}

An excellent reference for the methods of Fourier analysis for complex-valued
functions on local fields is Taibleson's book \cite{Taibleson - Fourier Analysis on Local Fields}.
I also recommend perusing the matter of harmonic analysis on locally
compact groups. Jordan Bell's lecture notes \cite{Bell - Harmonic Analysis on the p-adics,Bell - Pontryagin Duals of Q/Z and Q}
are a great resource for anyone who share my ``just shut up and tell
me how to compute stuff!'' attitude and the lack of patience implicit
therein. Folland's text on abstract harmonic analysis \cite{Folland - harmonic analysis}
makes for an excellent chaser to introduce the theory as a whole\textemdash Pontryagin
duality, connections to representation theory, and all. And although
\cite{Automorphic Representations} is a text on the representation
theory of the general linear group, it nevertheless presents a thorough
account of the tools of $\left(p,\infty\right)$-adic analysis\textemdash Fourier
transform \emph{and }Mellin transform\textemdash in its opening chapters.
Tate's thesis \cite{Tate's thesis}, is of course, a classic, and
a beautiful display of these analytical methods at work in a number-theoretic
application.

The connections with $p$-adic mathematical physics are also worth
exploring, particularly for the reader interested in learning more
about the nitty-gritty of integration and the theory of distributions.
Vladimirov's article \cite{Vladimirov - the big paper about complex-valued distributions over the p-adics}
contains a comprehensive account of $\left(p,\infty\right)$-adic
distributions and the accompanying notion of distributional derivatives,
going so far as to even solve some differential equations by this
method. \cite{Real and p-Adic Oscillatory integrals,p-adic van der Corput lemma}
deal with $\left(p,\infty\right)$-adic oscillatory integrals.

\subsubsection*{$\left(p,p\right)$-Adic Analysis}

After reading the first chapter of Gouvea's undergraduate-level text
\cite{Gouvea's introudction to p-adic numbers book} in $\left(p,p\right)$-adic
analysis, I used Robert's \emph{A course in $p$-adic analysis }(\cite{Robert's Book})
to acquaint myself with the ins and outs of the titular subject, along
with a little help from the seminar in algebraic number theory taught
at the University of Southern California by my eventual co-advisor\footnote{Along with Nicolai Haydn},
Prof. Sheldon Kamienny.

Robert's book gives an excellent, balanced treatment of the subject;
I never had the feeling that he had to tie and gag his inner number
theorist or algebraic geometer to keep it from trying to monopolize
the discussion. Curious undergraduates should probably peruse \cite{Gouvea's introudction to p-adic numbers book}
before attempting Robert, if only for the pleasure of reading Gouvea's
beautiful, lucid explanations. Koblitz's books \cite{Koblitz's book,Koblitz's other book}
are more biased toward number theoretical concerns, but are still
very accessible. More advanced\footnote{In other words, more \emph{algebraic}.}
references include Kedlaya's book on $p$-adic differential equations
\cite{Kedlaya}, and\textemdash diving straight in to the algebra
and number theory\textemdash Iwasawa's classic treatise on $p$-adic
$L$-functions \cite{Iwasawa} or Colmez's exposition of the same
\cite{p-adic L-functions paper}, the chapter on $p$-adic distributions
in Washington's \emph{Introduction to Cyclotomic Fields} \cite{Cyclotomic fields}.
Readers with a desire to force-feed themselves the theory of Berkovitch
spaces need only turn to \cite{Berkovich spaces and applications}
and \cite{Bosch lying title} for algebraic geometry in all its haunting,
malefic glory. On the analytical side, Anashi has a free e-text on
$p$-adic ergodic theory \cite{Anashi - p-adic ergodic theory} and
Adams' work on $p$-adic transcendental number theory contains an
exposition of the $p$-adic analogue of contour integration (the \index{integral!Shnirelmann}\textbf{Shnirelmann
integral}) \cite{Adams' paper on p-adic transcendental numbers with an appendix on Schnirelman integration}\textemdash though,
to be clear, this is \emph{not }the same as the $p$-adic line integral
used in rigid analytic spaces. Cherry's notes \cite{Cherry non-archimedean function theory notes}
give an excellent introduction to the theory of analytic and meromorphic
$\left(p,p\right)$-adic functions, especially for readers who don't
have the patience for Robert's more completionist presentation, or
Escassut more advanced, specialized treatment of the topic in \cite{value distribution in p-adic analysis}.
George Brom's thesis \cite{pp adic Fourier theory} covers Fourier
analysis in the $\left(p,p\right)$-adic setting, a fascinating oddity,
seeing as\textemdash recall\textemdash there is no non-trivial translation-invariant
$p$-adic-valued linear functional on the space of continuous functions
$\mathbb{Z}_{p}\rightarrow\mathbb{C}_{p}$.

\subsubsection*{Miscellany}

For completeness' sake, I have included in the bibliography much of
the material that I had on mind from 2017 to 2019, back when I was
still trying to use complex analytic methods and their $\left(p,\infty\right)$-adic
extensions to study Hydra maps. To reiterate, the guiding light of
behind my investigations of those two years was the functional equation
approach which I discovered independently of Meinardus and Berg. This
line of research had me trying to learn as much as I could about holomorphic
functions on the open unit disk and the boundary values of power series.
Garnett's\footnote{My undergraduate analysis professor!} book \cite{Bounded analytic functions}
gives a comprehensive account of harmonic analysis on the disk (and
also half-plane), including\textemdash but not limited to\textemdash the
Poisson Kernel and Hardy spaces. \cite{Ross et al} is an expository
text on the Cauchy transform; \cite{Fractional Cauchy transforms}
is a textbook on the Fractional Cauchy transform. Pavlovi\'{c}'s book\footnote{Not to be confused with the less-advanced but similarly-titled \cite{Baby Pavlovic},
also by Pavlovi\'{c}.} \cite{function classes on the unit disc} is an \emph{encyclopedic}
account of spaces of complex-valued functions on the open unit disk
in $\mathbb{C}$.

One of the central themes of my complex-analytic investigations was
the study of generating functions I called \textbf{set-series}\index{set-series};
these are of the form: 
\begin{equation}
\varsigma_{V}\left(z\right)\overset{\textrm{def}}{=}\sum_{v\in V}z^{v}\label{eq:Definition of a set-series}
\end{equation}
where $V\subseteq\mathbb{N}_{0}$ is a set of interest. For my work,
$V$ was an orbit class (or union thereof) of a given Hydra map on
$\mathbb{Z}$. H. S. Wilf's delightfully titled \emph{generatingfunctionology
}\cite{generatingfunctionology} is a thrilling read, as is Flajolet
and Sedgewick's \emph{Analytic Combinatorics} \cite{analytic combinatorics}.
For a thorough accounting of Tauberian methods used to study the coefficients
of functions represented by power series or integrals, I cannot give
too strong of a recommendation of Korevaar's book \cite{Korevaar}
or Bingham, Goldie, and Teugels' \emph{magisterial }tome, the \emph{Regular
Variation} \cite{Regular Variation}. For a more general background
in (divergent) infinite series, Hardy's \emph{Divergent Series} \cite{Hardy - Divergent Series}
is a classic. Flajolet's fascinating series of expository articles
on the Mellin transform (\cite{Flajolet - Mellin Transforms,Flajolet - Digital sums},
etc.) provide powerful analytical tools for anyone wishing to study
generating functions and power series asymptotics.

It is also noteworthy that generating functions have a fascinating
connection to transcendental number theory all their own. The study
of this connection often goes by the name of \textbf{Mahler Theory}\footnote{Kurt, not Gustav.}\textbf{
}\index{Mahler Theory}\index{Mahler, Kurt}, after the German-Australian
mathematician Kurt Mahler (of $p$-adic Mahler basis fame). Letting
$d$ be an integer $\geq2$, the function\footnote{The notation $\phi_{d}$ is my own.}:
\begin{equation}
\phi_{d}\left(z\right)\overset{\textrm{def}}{=}\sum_{n=0}^{\infty}z^{d^{n}}\label{eq:Definition of Mahler's phi_d}
\end{equation}
is a \textbf{lacunary}\footnote{Meaning ``gap''.}\textbf{ series
}\index{series!lacunary} with a natural boundary on the unit circle
in $\mathbb{C}$. This function satisfies the functional equation:
\begin{equation}
\phi_{d}\left(z^{d}\right)=\phi_{d}\left(z\right)-z\label{eq:Mahler functional equation}
\end{equation}
By using this equation, Mahler was able to prove that the complex
number $\phi_{d}\left(\alpha\right)$ is transcendental for any algebraic
number $\alpha$ with complex absolute value $0<\left|\alpha\right|<1$
\cite{Mahler theory}. Mahler Theory uses functions characterized
by polynomial functional equations like (\ref{eq:Mahler functional equation})
to show that said functions' take transcendental values at algebraic
numbers \cite{mahler functions,Loxton and van der Poorten - arithemtic properties of functional equations}.
Mahler theory is also connected with the theory of automata and automatic
sequences\textemdash see \cite{mahler functions} and the sources
cited therein.\newpage{}

\pagestyle{fancy}\fancyfoot{}\fancyhead[L]{\sl LIST OF SYMBOLS}\fancyhead[R]{\thepage}\printnomenclature{}

\addcontentsline{toc}{chapter}{Index}\pagestyle{fancy}\fancyfoot{}\fancyhead[L]{\sl INDEX}\fancyhead[R]{\thepage}

\printindex

\begin{thebibliography}{100}
\bibitem{mahler functions}\addcontentsline{toc}{chapter}{Bibliography}

\pagestyle{fancy}\fancyfoot{}\fancyhead[L]{\sl BIBLIOGRAPHY}\fancyhead[R]{\thepage}Adamczewski, Boris, and Jason Bell. ``A problem around Mahler functions.''
(2012). <https://hal.archives-ouvertes.fr/hal-00798303/document>

\bibitem{Adams' paper on p-adic transcendental numbers with an appendix on Schnirelman integration}Adams,
William Wells. ``Transcendental numbers in the $p$-adic domain''.
Columbia University, 1964.

\bibitem{Aguayo 1}Aguayo, J. ``Vector measures and integral operators''.
Ultrametric Functional Analysis, Contemp. Math., vol. 384, Amer. Math.
Soc., Providence, RI (2005), pp. 1-13.

\bibitem{Aguayo 2}Aguayo, J. and T. E. Gilsdorf. ``Non-Archimedean
vector measures and integral operators'' $p$-Adic Functional Analysis,
Lect. Notes Pure Appl. Math., vol. 222, Marcel Dekker, New York (2001),
pp. 1-11.

\bibitem{Aguayo and Moraga radon-nikodym theorem paper}José Aguayo
and Mirta Moraga. ``A radon nikodym theorem in the non-archimedean
setting.'' Proyecciones (Antofagasta) 20.3 (2001): 263-279.

\bibitem{key-10}Akin, Ethan (2004), ``Why is the $3x+1$ Problem
Hard?'', In: Chapel Hill Ergodic Theory Workshops (I. Assani, Ed.),
Contemp. Math. vol 356, Amer. Math. Soc. 2004, pp. 1\textendash 20.
(MR 2005f:37031).

\bibitem{Amice}Amice, Yvette. \emph{Les nombres $p$-adiques}. Vol.
14. Presses universitaires de France, 1975.

\bibitem{Anashi - p-adic ergodic theory}Anashi, Vladimir Sergeevich
(2014). \emph{The p-adic ergodic theory and applications}. <https://www.researchgate.net/publication/269571423\_The\_p-adic\_ergodic\_theory\_and\_applications>/

\bibitem{Andaloro - first paper on sufficient sets}Andaloro, Paul.
``On total stopping times under $3x+1$ iteration'', Fibonacci Quart.
38 (2000), no. 1, 73\textendash 78. MR 1738650.

\bibitem{Applegate and Lagarias - Trees}Applegate, David, and Jeffrey
C. Lagarias. ``Density bounds for the $3x+1$ problem. I. Tree-search
method.'' mathematics of computation 64.209 (1995): 411-426.

\bibitem{Applegate and Lagarias - Difference inequalities}Applegate,
David, and Jeffrey C. Lagarias. ``Density bounds for the $3x+1$
problem. II. Krasikov inequalities.'' Mathematics of computation
64.209 (1995): 427-438.

\bibitem{3x+1 semigroup}David Applegate, Jeffrey C. Lagarias, ``The
$3x+1$ semigroup'', Journal of Number Theory, Volume 117, Issue
1, 2006, Pages 146-159, ISSN 0022-314X. <https://www.sciencedirect.com/science/article/pii/S0022314X05001459>

\bibitem{Arens =00003D000026 Singer - Generalized Analytic Functions}Arens,
Richard and I. M. Singer (1956). ``Generalized Analytic Functions''.
Transactions of the American Mathematical Society, Vol. 81, No. 2
(Mar., 1956), pp. 379-393. Published by: American Mathematical Society.
<http://www.jstor.org/stable/1992923>.

\bibitem{Arnold}Arnold, V. I., ``On teaching mathematics''. Translated
from the original Russian by A.V. Goryunov. Published in: Uspekhi
Mat. Nauk 53 (1998), no. 1, 229-234; English translation: Russian
Math. Surveys 53 (1998), no. 1, 229-236. <http://www.ceremade.dauphine.fr/\textasciitilde msfr/articles/arnold/PRE\_anglais.ps>

\bibitem{key-2-1}Axelsson, Ekaterina Yurova. ``On the sub-coordinate
representation of $p$-adic functions''. \emph{Advances in Ultrametric
Analysis}, edited by Alain Escassut, et al., American Mathematical
Society, 2018. ProQuest Ebook Central, <https://ebookcentral.proquest.com/lib/socal/detail.action?docID=5347081>.

\bibitem{Baker's Transcendental Number Theory}Baker, Alan (1990).
\emph{Transcendental number theory}. Cambridge Mathematical Library
(2nd ed.), Cambridge University Press, ISBN 978-0-521-39791-9, MR
0422171.

\bibitem{p-adic proof for =00003D0003C0}Baker, Matt. ``A p-adic
proof that pi is transcendental''. \emph{Matt Baker's Math Blog}.
20 March. 2015. <https://mattbaker.blog/2015/03/20/a-p-adic-proof-that-pi-is-transcendental/>

\bibitem{Simon =00003D000026 Breuer - Natural Boundaries =00003D000026 Spectral Theory}Barry
Simon and Jonathan Breuer (2011).\emph{ }``Natural Boundaries and
Spectral Theory''.

\bibitem{natural boundaries}Bell, Jason P., and Jeffrey C. Lagarias.
``$3x+1$ inverse orbit generating functions almost always have natural
boundaries.'' arXiv preprint arXiv:1408.6884 (2014).

\bibitem{Bell - Harmonic Analysis on the p-adics}Bell, Jordan. (2016)
``Harmonic analysis on the $p$-adic numbers''.\emph{ }<https://pdfs.semanticscholar.org/9d4f/f3efad867fe867331bb00402807e0856170c.pdf>

\bibitem{Bell - Pontryagin Duals of Q/Z and Q}Bell, Jordan. (2015)
``The Pontryagin duals of $\mathbb{Q}/\mathbb{Z}$ and $\mathbb{Q}$''.
<https://pdfs.semanticscholar.org/bffd/643d972eab01c0b063c34ad1560791899d3a.pdf>

\bibitem{Color Picture Book}Alek Bellos and Edmund Harriss. \emph{Visions
of the Universe: A Coloring Journey Through Math's Great Mysteries}.\emph{
}Experiment; Clr Csm edition. 29 Nov. 2016. ISBN-10: 1615193677.

\bibitem{Berg =00003D000026 Opfer}Berg, Lothar, and Gerhard Opfer.
``An Analytic Approach to the Collatz $3n+1$ Problem for Negative
Start Values.'' Computational Methods and Function Theory 13.2 (2013):
225-236.

\bibitem{Berg =00003D000026 Meinardus}Berg, Lothar, and Günter Meinardus.
``The $3n+1$ Collatz Problem and Functional Equations.'' (1995).

\bibitem{Meinardus and Berg}Lothar Berg and Günter Meinardus. ``Functional
equations connected with the Collatz problem.'' Results in Mathematics
25, no. 1-2 (1994): 1-12.

\bibitem{Besikovitch - Almost Periodic Functions}Besikovitch, A.
S. (1954) \emph{Almost Periodic Functions}. Dover Publications, Inc.,
through special permission of Cambridge University Press.

\bibitem{Regular Variation}Bingham, Nicholas H., Charles M. Goldie,
and Jozef L. Teugels. \emph{Regular variation}. No. 27. Cambridge
university press, 1989.

\bibitem{pp adic Fourier theory}Borm, George Florimond. \emph{$p$-adic
Fourier theory}. Diss. {[}Sl: sn{]}, 1988.

\bibitem{Bosch lying title}S. Bosch, U. Güntzer, and R. Remmert,
Non-Archimedean Analysis, Grundlehren der Math. Wiss. 261, Springer-Verlag,
Berlin, 1984. <https://www.math.arizona.edu/\textasciitilde cais/scans/BGR-Non\_Archimedean\_Analysis.pdf>

\bibitem{Introduction to Langlands}Bump, Daniel, et al. \emph{An
introduction to the Langlands program}. Vol. 140. Springer Science
\& Business Media, 2003.

\bibitem{key-20}Chabrera, Víctor Martín. ``An algebraic fractal
approach to Collatz Conjecture''. BS thesis. Universitat Politècnica
de Catalunya, 2019.

\bibitem{Cherry non-archimedean function theory notes}Cherry, William.
``Lectures on non-archimedean function theory.'' arXiv preprint
arXiv:0909.4509 (2009). <https://arxiv.org/pdf/0909.4509.pdf>

\bibitem{Cohen Number Theory}Cohen, H. (2008). \emph{Number theory:
Volume II: Analytic and modern tools} (Vol. 240). Springer Science
and Business Media.

\bibitem{Collatz Letter}Collatz, Lothar. ``\href{http://www.cecm.sfu.ca/organics/papers/lagarias/paper/goodies/ubersetzung/html/ubersetzung.html}{Letter to Professor Mays}''.
Translated from the German by Marion Meudt and John Read\footnote{<http://www.cecm.sfu.ca/organics/papers/lagarias/paper/goodies/ubersetzung/html/ubersetzung.html>}.

\bibitem{p-adic L-functions paper}Colmez, Pierre. ``Fontaine's
rings and $p$-adic $L$-functions.'' Lecture notes 32 (2004): 33.
<http://people.math.jussieu.fr/\textasciitilde colmez/tsinghua.pdf>

\bibitem{Conrad on p-adic series}Conrad, K. ``Infinite Series in
$p$-Adic Fields''.

\bibitem{Conrad - Local Global Principle}Conrad, K. ``The Local-Global
Principle''.

\bibitem{Conway - Unpredictable Iterations}Conway, John H. (1972).
``Unpredictable Iterations''. Proc. 1972 Number Theory Conf., Univ.
Colorado, Boulder, pp. 49\textendash 52.

\bibitem{Bohm and Sontacchi}Corrado Böhm and Giovanna Sontacchi (1978).
``On the existence of cycles of given length in integer sequences
like $x_{n}+1=x_{n}/2$ if $x_{n}$ even, and $x_{n+1}=3x_{n}+1$
otherwise'', Atti Ac- cad. Naz. Lincei Rend. Cl. Sci. Fis. Mat. Natur.
64 (1978), 260\textendash 264. (MR 83h:10030)

\bibitem{Circulant Matrices}Davis, Philip J. \emph{Circulant matrices}.
American Mathematical Soc., 2013.

\bibitem{Dynamical Systems}Devaney, Robert L. \emph{An introduction
to chaotic dynamical systems}. CRC press, 2018.

\bibitem{Non-archimedean operator theory}Diagana, Toka, and François
Ramaroson. \emph{Non-Archimedean operator theory}. Springer, 2016.

\bibitem{First 30 years of p-adic mathematical physics}Dragovich,
B., et al. ``$p$-Adic mathematical physics: the first 30 years.''
P-Adic numbers, ultrametric analysis and applications 9.2 (2017):
87-121.

\bibitem{Adelic Harmonic Oscillator}Dragovich, Branko. ``Adelic
harmonic oscillator.'' International Journal of Modern Physics A
10.16 (1995): 2349-2365.

\bibitem{Berkovich spaces and applications}Ducros, Antoine, Charles
Favre, and Johannes Nicaise, eds. \emph{Berkovich spaces and applications}.
Springer, 2015.

\bibitem{Dwork Zeta Rationality}Dwork, Bernard. ``On the rationality
of the zeta function of an algebraic variety''. Amer. J. Math., 82:631\textendash 648,
1960.

\bibitem{value distribution in p-adic analysis}Escassut, Alain. \emph{Value
Distribution in p-adic Analysis}. World Scientific, 2015.

\bibitem{n-furcations}M. R. Feix, M. Muriel and J. L. Rouet. ``Statistical
Properties of an Iterated Arithmetic Mapping'', J. Stat. Phys., 76
(1994), 725-741

\bibitem{Geometric Measure Theory}Federer, Herbert. \emph{Geometric
measure theory}. Springer, 2014.

\bibitem{analytic combinatorics}Flajolet, Philippe, and Robert Sedgewick.
\emph{Analytic combinatorics}. Cambridge University press, 2009.

\bibitem{Flajolet - Mellin Transforms}Philippe Flajolet, Xavier Gourdon,
and Philippe Dumas. ``Mellin transforms and asymptotics: Harmonic
sums\textquotedblright . In: Theoretical Computer Science 144 (1995),
pp. 3\textendash 58.

\bibitem{Flajolet - Digital sums}Philippe Flajolet, et al. ``Mellin
transforms and asymptotics: digital sums\textquotedblright . In: Theoretical
Computer Science 123 (1994), pp. 291\textendash 314.

\bibitem{Folland - real analysis}Folland, Gerald B. \emph{Real analysis:
modern techniques and their applications}. Vol. 40. John Wiley and
Sons, 1999.

\bibitem{Folland - harmonic analysis}Folland, Gerald B. \emph{A course
in abstract harmonic analysis}. Vol. 29. CRC press, 2016.

\bibitem{ergodic theory of p-adic transformations}Furno, Joanna.
2013. \emph{Ergodic Theory of P-Adic Transformations}. PhD Thesis.
University of North Carolina at Chapel Hill. <https://doi.org/10.17615/8p5w-9v68>

\bibitem{Bounded analytic functions}Garnett, John. \emph{Bounded
analytic functions}. Vol. 236. Springer Science \& Business Media,
2007.

\bibitem{George collatz visualization}Georgiou, George M. ``Encoding
and Visualization in the Collatz Conjecture.'' arXiv preprint arXiv:1811.00384
(2018).

\bibitem{Automorphic Representations}Goldfeld, D., and Hundley, J.
(2011). \emph{Automorphic Representations and L-Functions for the
General Linear Group: Volume 1} (Vol. 129). Cambridge University Press.

\bibitem{Gouvea's introudction to p-adic numbers book}Gouvêa, Fernando
Quadros. \emph{$p$-adic Numbers: An Introduction}. Springer Science
and Business Media, 2012.

\bibitem{Gouvea's p-adic number history slides}Gouvêa, Fernando Quadros.
``Hensel's $p$-adic Numbers: early history'' <https://www-fourier.ujf-grenoble.fr/\textasciitilde panchish/Mag2009L3/GouveaHensel2.pdf>

\bibitem{Walsh}B. Golubov, A. Efimov, and V. Skvortsov. ``Walsh
Functions and Their Generalizations.'' \emph{Walsh Series and Transforms}.
Springer, Dordrecht, 1991. 1-33.

\bibitem{key-18}Hashimoto, Yukihiro. ``A fractal set associated
with the Collatz problem.'' Bull. of Aichi Univ. of Education, Natural
Science 56 (2007): 1-6.

\bibitem{key-21}Hassan, Sk. ``Discrete dynamics of one dimensional
Collatz like integral value transformations.'' Journal of Applied
Mathematics and Computing 49.1 (2015): 91-105.

\bibitem{Hardy - Divergent Series}Hardy, G. H. (1949). \emph{Divergent
Series}. Oxford University Press. Oxford: 1949.

\bibitem{Dirichlet Series (Hardy)}G. H. Hardy and Marcel Reisz. (1915).
\emph{The general theory of Dirichlet's series}. Cambridge {[}Eng.{]}
University Press, 1915. <http://hdl.handle.net/2027/mdp.39015064521282>

\bibitem{key-27}Helson, Henry (1955). ``On a theorem of Szeg\H{o}''.
Proceedings of the American Mathematical Society.

\bibitem{Hensel-s original article}Hensel, Kurt. \textquotedbl Ueber
die Zurückführung der Divisorensysteme auf eine reducirte Form.\textquotedbl{}
(1897): 234-250.

\bibitem{Fractional Cauchy transforms}Hibschweiler, Rita A., and
Thomas H. MacGregor. \emph{Fractional Cauchy Transforms}. Chapman
and Hall/CRC, 2020.

\bibitem{Iwasawa}Iwasawa, Kinkichi. \emph{Lectures on $P$-Adic $L$-Functions}.(AM-74),
Volume 74. Princeton University Press, 2016.

\bibitem{Jackson}Jackson, Glenn Micah Jr., ``Explorations of the
Collatz Conjecture (mod $m$)'' (2014). Honors College Theses. 70.
<https://digitalcommons.georgiasouthern.edu/honors-theses/70>

\bibitem{key-19}Joseph, L. Pe. ``The $3x+1$ fractal.'' Computers
\& Graphics 28.3 (2004): 431-435.

\bibitem{Non-archimedean function spaces and the dominated convergence theorem}J.
K\k{a}kol. C. Perez-Garcia. W. \'{S}liwa. ``Non-archimedean function
spaces and the Lebesgue dominated convergence theorem.'' Bull. Belg.
Math. Soc. Simon Stevin 19 (1) 173 - 184, march 2012. https://doi.org/10.36045/bbms/1331153417

\bibitem{Katsaras}Katsaras. A. K. ``On $p$-adic vector measure
spaces'', J. Math. Anal. Appl. ,365(2010), 342-357.

\bibitem{Kedlaya}Kedlaya, Kiran S. \emph{p-adic Differential Equations}.
Vol. 125. Cambridge University Press, 2010.

\bibitem{Quantum Paradoxes}Khrennikov, Andrei. \emph{Non-Archimedean
Analysis: Quantum Paradoxes, Dynamical Systems and Biological Models}.
Springer, Dordrecht, 1997. 101-129.

\bibitem{Probabilities taking values in non-archimedean fields}Khrennikov,
A. Yu. ``Generalized probabilities taking values in non-Archimedean
fields and in topological groups.''\textquotedbl{} Russian Journal
of Mathematical Physics 14.2 (2007): 142-159.

\bibitem{Measure-theoretic approach to p-adic probability theory}Khrennikov,
Andrei; Yamada, Shinichi; van Rooij, Arnoud. ``The measure-theoretical
approach to $p$-adic probability theory''. Annales Mathématiques
Blaise Pascal, Tome 6 (1999) no. 1, pp. 21-32. <http://www.numdam.org/item/AMBP\_1999\_\_6\_1\_21\_0/>

\bibitem{p-adic Tauberian}Khrennikov, A. Yu, and V. M. Shelkovich.
``Distributional asymptotics and $p$-adic Tauberian and Shannon\textendash Kotelnikov
theorems.'' Asymptotic Analysis 46.2 (2006): 163-187.

\bibitem{Koblitz's book}Koblitz, Neal. \emph{$p$-adic analysis:
A short course on recent work}. Vol. 46. Cambridge University Press,
1980.

\bibitem{Koblitz's other book}Koblitz, Neal. \emph{$p$-adic Numbers,
$p$-adic Analysis, and Zeta-Functions}. Vol. 58. Springer Science
\& Business Media, 2012.

\bibitem{Kochubei}Kochubei, Anatoly N. ``Non-Archimedean Duality:
Algebras, Groups, and Multipliers''. 2016 <https://arxiv.org/pdf/1510.06876.pdf>

\bibitem{RCWAG}Kohl, Stefan (2016).\emph{ RCWA: Residue-Class-Wise
Affine Groups}. Version 4.0.0 <https://www.gap-system.org/Manuals/pkg/rcwa/doc/manual.pdf>

\bibitem{dgh paper}A. V. Kontorovich and Y. G. Sinai. ``Structure
Theorem for $\left(d,g,h\right)$-Maps'', 2006 <https://arxiv.org/pdf/math/0601622.pdf>

\bibitem{Lagarias-Kontorovich Paper}Alex V. Kontorovich and Jeffery
C. Lagarias. (2009) ``Stochastic Models for the 3x+1 and 5x+1 Problems''.
<https://arxiv.org/pdf/0910.1944>

\bibitem{Korevaar}Korevaar, Jacob. \emph{Tauberian theory: a century
of developments}. Vol. 329. Springer Science and Business Media, 2013.

\bibitem{Krasikov}Krasikov, Ilia. ``How many numbers satisfy the
$3x+1$ conjecture?'' International Journal of Mathematics and Mathematical
Sciences 12.4 (1989): 791-796.

\bibitem{key-4}Kudla, Stephen S. ``Tate's thesis.'' An introduction
to the Langlands program. Birkhäuser, Boston, MA, 2004. 109-131.

\bibitem{key-4-1}Kyung Soo Rim. ``Gibbs phenomenon for Fourier partial
sums on $\mathbb{Z}_{p}$'', \emph{Journal of Mathematical Analysis
and Applications}, Volume 433, Issue 1, 2016, Pages 392-404, ISSN
0022-247X, <https://doi.org/10.1016/j.jmaa.2015.08.005>

\bibitem{Ultimate Challenge}Lagarias, Jeffrey C., ed. \emph{The ultimate
challenge: The $3x+1$ problem}. American Mathematical Soc., 2010.

\bibitem{Lagarias' Survey}Lagarias, J.C. ``The $3x+1$ problem and
its generalizations'', Amer. Math. Monthly 92 (1985), 3\textendash 23.

\bibitem{Lagarias - rational cycles for 3x+1}Lagarias, Jeffrey C.
``The set of rational cycles for the $3x+1$ problem'', Acta Arith.
56 (1990), no. 1, 33\textendash 53. MR 1067980, DOI <https://doi.org/10.4064/aa-56-1-33-53>

\bibitem{Lang - Algebraic Number Theory}Lang, Serge. \emph{Algebraic
number theory}. Vol. 110. Springer Science and Business Media, 2013.

\bibitem{key-9}Laarhoven, Thijs, and Benne de Weger. ``The Collatz
conjecture and De Bruijn graphs.'' Indagationes Mathematicae 24.4
(2013): 971-983.

\bibitem{Leigh's article}Leigh, G.M. ``A Markov process underlying
the generalized Syracuse algorithm'', Acta Arith. 46 (1986) 125-143.

\bibitem{Liepins}Liepins, Gunar E. (1976). ``A Paley-Wiener Theorem
for Locally Compact Abelian Groups''. Transactions of the American
Mathematical Society, Vol. 222 (Sep., 1976), pp. 193-210 Published
by: American Mathematical Society. <http://www.jstor.org/stable/1997665>

\bibitem{Loxton and van der Poorten - arithemtic properties of functional equations}J.
H. Loxton and A. J. van der Poorten (1981). ``Arithmetic properties
of the solutions of a class of functional equations''.

\bibitem{Ludkovsky on non-archimedean measures}Ludkovsky, S. V. ``Quasi-Invariant
and Pseudo-Differentiable Measures with Values in Non-Archimedean
Fields on a Non-Archimedean Banach Space''. J Math Sci 128, 3428\textendash 3460
(2005). <https://doi.org/10.1007/s10958-005-0280-2>

\bibitem{Mackey - Laplace Transform for Locally Compact Abelian groups}Mackey,
George W. (1948). ``The Laplace Transform for Locally Compact Abelian
Groups''. Proceedings of the National Academy of Sciences, Volume
34, Number 4, pp. 156-162. <http://www.pnas.org/content/34/4/156>

\bibitem{Mahler}Mahler, Kurt. \emph{Introduction to p-adic numbers
and their functions}. No. 64. CUP Archive, 1973.

\bibitem{Mahler Series}Mahler, Kurt. ``An Interpolation Series for
Continuous Functions of a $p$-adic Variable.'' (1958): 23-34.

\bibitem{Mandelbrojt's lacunary series paper}Mandelbrojt, Szolem.
``Lacunary series''. Rice Institute Pamphlet - Rice University Studies,
14, no. 4 (1927) Rice University: <https://hdl.handle.net/1911/8511>.

\bibitem{Mandelbrojt's Dirichlet Series book}Mandelbrojt, S. (1972)
\emph{Dirichlet Series: Principals and Methods}. D. Reidel Publishing
Company / Dordrecht-Holland.

\bibitem{The Berkovitch Space Paper}Maïnetti, Nicolas. ``Gefland
transform and spectral radius formulae for ultrametric Banach algebras''
Ultrametric Functional Analysis: Seventh International Conference
on P-Adic Functional Analysis, June 17-21, 2002, University of Nijmegen,
the Netherlands. Vol. 319. American Mathematical Soc., 2003.

\bibitem{Matthews' Conjecture}Matthews, K. R. ``A Generalized $3x+1$
Conjecture''. <https://www.numbertheory.org/gnubc/challenge>

\bibitem{Matthews' Leigh Article}Matthews, K. R. ``The generalized
$3x+1$ mapping: George Leigh's approach''. <http://www.numbertheory.org/keith/george.html>

\bibitem{Matthews' slides}Matthews, K. R. ``Generalizations of the
$3x+1$ problem and connections with Markov matrices and chains''
(2010) <http://www.numbertheory.org/PDFS/3x+1\_slides.pdf>

\bibitem{Matthews and Watts}Matthews, K. R., and Watts, A.M. (1985).
``A Markov approach to the generalized Syracuse algorithm''. Acta
Arithmetica, 45(1), 29-42.

\bibitem{Meinardus}Meinardus, Günter. ``Some Analytic Aspects Concerning
the Collatz Problem.'' (2001).

\bibitem{key-6-2-1}Memic, Nacima. ``Multiplicative Systems on Ultra-Metric
Spaces.'' Mathematica Balkanica (2010).

\bibitem{Moller's paper (german)}Möller, Herbert. ``Über Hasses
Verallgemeinerung des Syracuse-Algorithmus (Kakutanis Problem)'',
Acta Arith. 34 (1978) 219-226.

\bibitem{Monks' sufficiency paper for 3x+1}Monks, Kenneth. ``The
sufficiency of arithmetic progressions for the $3x+1$ Conjecture.''
Proceedings of the American mathematical society 134.10 (2006): 2861-2872.

\bibitem{The many Monks paper}Monks, Keenan, et al. ``On the distribution
of arithmetic sequences in the Collatz graph.'' <http://www.mathematicalgemstones.com/maria/slides/SlidesV4.pdf>

\bibitem{Monna and Springer}Monna, A. and Springer, T. ``Integration
non-Archimedienne'', 1,2. Indag. Math., 25, 634-653 (1963)

\bibitem{Mahler theory}Nishioka, Kumiko. \emph{Mahler functions and
transcendence}. Springer, 2006.

\bibitem{Schikhof memorial}Ochsenius, H., E. Olivos, and C. Perez-Garcia.
``Remembering WH Schikhof.'' P-Adic Numbers, Ultrametric Analysis,
and Applications 7.2 (2015): 81-95.

\bibitem{Collatz Biography}J. J. O'Connor and E. F. Robertson. ``Lothar
Collatz''. <https://mathshistory. st-andrews.ac.uk/Biographies/Collatz/>,
2006. St Andrews University School of Mathematics and Statistics,
Scotland.

\bibitem{Oliveira}T. Oliveira e Silva ``Maximum excursion and stopping
time record-holders for the problem: computational results''. Math.
Comp., 68 (1) (1999), pp. 371-384 (see website: <http://www.ieeta.pt/\textasciitilde tos/3x+1.html>
for additional computations)

\bibitem{Opfer}Opfer, Gerhard. ``An Analytic Approach to the Collatz
$3n+1$ Problem''. Hamburger Beiträge zur Angewandten Mathematik
Nr. 2011-09, May 2011.

\bibitem{Baby Pavlovic}Pavlovi\'{c}, Miroslav. \emph{Introduction
to function spaces on the disk}. Belgrade, 2004.

\bibitem{function classes on the unit disc}Pavlovi\'{c}, Miroslav.
\emph{Function classes on the unit disc}. de Gruyter, 2019.

\bibitem{collatz harmonic oscillator}Perelman, C. C., Carbó-Dorca,
R. ``The Collatz conjecture and the quantum mechanical harmonic oscillator''.
J Math Chem 60, 145\textendash 160 (2022). <https://doi.org/10.1007/s10910-021-01296-6>

\bibitem{key-4-2}Perez-Garcia, C. and W. H. Schikhof. ``The Orlicz-Pettis
property in $p$-adic analysis''. Collect. Math., 43 (1992), 225-233.

\bibitem{Journey throughout the history of p-adic numbers}Perrin,
Yvette. ``A journey throughout the history of $p$-adic numbers.''
\emph{Advances in Ultrametric Analysis}, edited by Alain Escassut,
et al., American Mathematical Society, 2018. ProQuest Ebook Central,
<https://ebookcentral.proquest.com/lib/socal/detail.action?docID=5347081>

\bibitem{Remmert}Remmert, Reinhold. \emph{Classical topics in complex
function theory}. Vol. 172. Springer Science and Business Media, 2013.

\bibitem{Robert's Book}Robert, A. M. (2013). \emph{A course in p-adic
analysis} (Vol. 198). Springer Science and Business Media. Chicago.

\bibitem{Real and p-Adic Oscillatory integrals}Rogers, Keith McKenzie.
\emph{Real and p-Adic Oscillatory Integrals}. Diss. The University
of New South Wales, 2004.

\bibitem{p-adic van der Corput lemma}Rogers, Keith M. ``A van der
Corput lemma for the $p$-adic numbers.'' Proceedings of the American
Mathematical Society (2005): 3525-3534.

\bibitem{van Rooij - Non-Archmedean Functional Analysis}van Rooij,
A.C.M. \emph{Non-Archimedean functional analysis.} Pure and Applied
Math., vol. 51, Marcel Dekker, New York, 1978.

\bibitem{van Rooij and Schikhof - "Non-archimedean analysis"}A. C.
M. van Rooij and W. H. Schikhof. ``Non-Archimedean analysis.\textquotedblright{}
Nieuw Arch. Wisk. (3)19, 120\textendash 160 (1971).

\bibitem{van Rooij and Schikhof "Non-archimedean integration theory"}Arnoud
C. M van Rooij, and Wim H. Schikhof. ``Non-Archimedean integration
theory.'' Indagationes Mathematicae (Proceedings). Vol. 72. No. 2.
North-Holland, 1969.

\bibitem{Roosendaal's Website}Roosendaal, E. ``On the $3x+1$ problem''
<www.ericr.nl/wondrous>

\bibitem{Ross et al}Ross, William T., Joseph A. Cima, and Alec L.
Matheson. ``The Cauchy Transform.'' (2006).

\bibitem{Parity Sequences}Rozier, O. (2018). ``Parity sequences
of the $3x+1$ map on the $2$-adic integers and Euclidean embedding''.
arXiv preprint arXiv:1805.00133.

\bibitem{Ultrametric Calculus}Schikhof, W. (1985). \emph{Ultrametric
Calculus: An Introduction to p-Adic Analysis} (Cambridge Studies in
Advanced Mathematics). Cambridge: Cambridge University Press. doi:10.1017/CBO9780511623844

\bibitem{Schikhof's Thesis}Schikhof, W. H. \emph{Non-Archimedean
Harmonic Analysis}, Ph.D. Thesis, pp. 1\textendash 80 (Catholic Univ.
of Nijmegen, The Netherlands, 1967).

\bibitem{Schikhof - Radon-Nikodym}Schikhof, W. H. ``A Radon-Nikodym
theorem for non-Archimedean integrals and absolutely continuous measures
on groups,\textquotedblright{} Indag. Math. 33, 78\textendash 85 (1971).

\bibitem{key-8-1}Schikhof, W. H. ``An Approach to the Ultrametric
Moment Problem'', Trudy Mat. Inst. Steklova, 2004, Volume 245, 251\textendash 256

\bibitem{Schikhof Banach Space Paper}Schikhof, W. H. ``Banach spaces
over non-Archimedean valued fields,\textquotedblright{} Topology Proc.
24, 547\textendash 581 (2001).

\bibitem{Schneider awful book}Schneider, Peter. \emph{Nonarchimedean
functional analysis}. Springer Science \& Business Media, 2001.

\bibitem{More Schneider lies}Schneider, Peter, and Jeremy Teitelbaum.
``$p$-adic Fourier theory.'' arXiv preprint math/0102012 (2001).

\bibitem{On Wiener's Lemma}Schulte, Eike. ``On Wiener's lemma for
locally compact abelian groups''. arXiv:2005.05212 {[}math.FA{]}
<https://arxiv.org/abs/2005.05212> (2020)

\bibitem{Local-Global Principle Failure}Selmer, E. ``The Diophantine
equation $ax^{3}+by^{3}+cz^{3}=0$'', Acta Arithmetica 85 (1951),
203\textendash 362.

\bibitem{Dreancatchers for Hydra Maps}Siegel, M. C. ``Conservation
of Singularities in Functional Equations Associated to Collatz-Type
Dynamical Systems; or, Dreamcatchers for Hydra Maps'' (2019) arXiv:1909.09733v3
{[}math.GM{]} <https://arxiv.org/abs/1909.09733>

\bibitem{Mellin transform paper}Siegel, M. C. ``Syracuse Random
Variables and the Periodic Points of Collatz-type maps'' (2020) <https://arxiv.org/abs/2007.15936>

\bibitem{Sultanow}Sultanow, Eldar, Christian Koch, and Sean Cox.
``Collatz sequences in the light of graph theory.'' (2020).

\bibitem{Taibleson - Fourier Analysis on Local Fields}Taibleson,
M. H. \emph{Fourier Analysis on Local Fields}. Notes, vol. 15, Princeton
Univ. Press (1975)

\bibitem{Tao Blog}Tao, Terence. (2011) ``The Collatz conjecture,
Littlewood-Offord theory, and powers of $2$ and $3$'' <https://terrytao.wordpress.com/2011/08/25/the-collatz-conjecture-littlewood-offord-theory-and-powers-of-2-and-3/>

\bibitem{Tao Fourier Transform Blog Post}Tao, Terence. (2009) ``245C,
Notes 2: The Fourier transform'' <https://terrytao.wordpress.com/2009/04/06/the-fourier-transform/>

\bibitem{Tao Probability paper}Tao, Terence. (2019) ``Almost all
orbits of the Collatz map attain almost bounded values.\textquotedblright{}
<https://arxiv.org/abs/1909.03562>

\bibitem{Tate's thesis}Tate, John (1950), \emph{Fourier analysis
in number fields and Hecke's zeta functions}, Princeton University
Ph.D. thesis under Emil Artin. Reprinted in Cassels, J. W. S.; Fröhlich,
Albrecht, eds. (1967), Algebraic number theory, London: Academic Press,
pp. 305\textendash 347, MR 0215665

\bibitem{Terras 76}Terras, Riho. (1976) ``A stopping time problem
on the positive integers''. Acta Arithmetica. XXX.

\bibitem{Terras 79}Terras, Riho. (1979) ``On the existence of a
density''. Acta Arithmetica. XXXV.

\bibitem{On the origins of p-adic analysis}Ullrich, Peter. ``On
the origins of p-adic analysis.'' Symposia Gaussiana, Proc. of the
2nd Gauss Symposium (München 1993), Conference A: Mathematics and
Theoretical Physics, ed. M. Behara, R. Fritsch und RG Lintz, (Berlin
1995). 1995.

\bibitem{Vladimirov - the big paper about complex-valued distributions over the p-adics}Vladimirov,
Vasilii S. ``Generalized functions over the field of $p$-adic numbers.''
Russian Mathematical Surveys 43.5 (1988): 19.

\bibitem{key-5-1}Volosivets, S. S. ``Fourier-Vilenkin series and
analogs of Besov and Sobolev classes.'' Annales Univ. Sci. Budapest,
Sect. Comp. Vol. 33. 2010.

\bibitem{Volo - p-adic Wiener}Volosivets, S. S. ``A Distributional
Proof of $p$-Adic Wiener Tauberian Theorem and Approximation by Translates
of a Function.'' p-Adic Numbers, Ultrametric Analysis and Applications
13.4 (2021): 308-315.

\bibitem{p-adic space-time}Volovich, Igor V.. ``$p$-adic space-time
and string theory.\textquotedblright{} Theoretical and Mathematical
Physics 71 (1987): 574-576.

\bibitem{Cyclotomic fields}Washington, Lawrence C. \emph{Introduction
to Cyclotomic Fields}. Second Edition. Vol. 83. New York: Springer
New York, 1997

\bibitem{de Weger on Opfer}de Weger, Benne. ``Comments on Opfer's
alleged proof of the $3n+1$ Conjecture''. Tech. rep. Eindhoven University
of Technology, 2011.

\bibitem{generatingfunctionology}Wilf, Herbert S. \emph{generatingfunctionology}.
CRC press, 2005.

\bibitem{Wirsching's book on 3n+1}Wirsching, Günther J. \emph{The
dynamical system generated by the $3n+1$ function}, Lecture Notes
in Mathematics, vol. 1681, Springer-Verlag, Berlin, 1998. MR 1612686\newpage{}
\end{thebibliography}
\end{document}